\documentclass[12pt]{tlpbook}
\usepackage{times,mathrsfs,amscd,amsmath,amsfonts,amssymb,graphicx,tikz,enumerate,xcolor,fancyhdr}
\usepackage[T1]{fontenc}
\usepackage{helvet}
\usepackage[most, listings, theorems]{tcolorbox}
\usepackage{babel}
\usepackage{thmtools}
\definecolor{grey}{rgb}{0.5,0.5,0.5}
\definecolor{orange}{rgb}{1,0.5,0}
\definecolor{purple}{rgb}{0.7,0,1}
\definecolor{darkgreen}{rgb}{0,0.5,0.2}
\definecolor{darkblue}{rgb}{0,0,1}
\definecolor{linkc}{rgb}{0.4,0.2,0.4}
\usepackage{hyperref}
\hypersetup{
    colorlinks=true,
    citecolor=darkblue,
    linkcolor=linkc,
    urlcolor=red,
    pdfpagemode=FullScreen,
    pdftitle={Mathematical Analysis\\Volume II},
    }

\newcounter{myquestion}
\newcommand{\atc}{\addtocounter{myquestion}{1}}

\newtcolorbox[auto counter, number within = chapter]
{theorem}[2][]{%
fonttitle= \bfseries\upshape, fontupper= \upshape,  arc=0mm,
colback=red!5,colframe=red!80!black,
title={Theorem \thetcbcounter\; #2},#1}

\newtcolorbox[use counter from = theorem, number within = chapter]
{proposition}[2][]{%
fonttitle= \bfseries\upshape, fontupper= \upshape, arc=0mm,
colback=red!5,colframe=red!80!black,
title={Proposition \thetcbcounter\; #2},#1}

\newtcolorbox[use counter from = theorem, number within = chapter]
{lemma}[2][]{%
fonttitle= \bfseries\upshape, fontupper= \upshape,
arc=0mm, colback=green!10,colframe=green!50!black,
title={Lemma \thetcbcounter\; #2},#1}

\newtcolorbox[use counter from = theorem, number within = chapter]
{corollary}[2][]{%
fonttitle= \bfseries\upshape, fontupper= \upshape, arc=0mm,
colback=green!10,colframe=green!50!black,
title={Corollary \thetcbcounter\; #2},#1}

 \newtcolorbox{myproof}[1]{%
fonttitle= \bfseries\upshape, fontupper= \upshape, arc=0mm,
colback=grey!10,colframe=grey!50!black, colbacktitle=grey!75!black,
enhanced, attach boxed title to top center={yshift=-2mm}, title={#1}}

\newtcolorbox[auto counter, number within = chapter]
{definition}[2][]{%
fonttitle= \bfseries\upshape, fontupper= \upshape, arc=0mm,
colback=blue!5,colframe=blue!75!black, 
title={Definition  \thetcbcounter\; #2},#1}

 \newtcolorbox[auto counter, number within = chapter]{remark}[2][]{%
fonttitle= \bfseries\upshape, fontupper= \upshape, arc=0mm,
colback=yellow!10,colframe=yellow!20!red,
title={Remark \thetcbcounter\; #2},#1}

\newtcolorbox[auto counter, number within = chapter]{example}[2][]{%
fonttitle= \bfseries\upshape, fontupper= \upshape, arc=0mm,
colback=orange!15,colframe=orange!55!black,
title={Example \thetcbcounter\; #2},#1}

 \newtcolorbox{example2}[1]{%
fonttitle= \bfseries\upshape, fontupper= \upshape, arc=0mm,
colback=orange!15,colframe=orange!55!black,  title=#1}

 \newtcolorbox{solution}[1]{%
fonttitle= \bfseries\upshape, fontupper= \upshape, arc=0mm,
colback=grey!10,colframe=grey!50!black, colbacktitle=grey!75!black,
enhanced, attach boxed title to top center={yshift=-2mm}, title={#1}}

 \newtcolorbox{highlight}[1]{%
fonttitle= \bfseries\upshape, fontupper= \upshape, arc=0mm,
colback=yellow!10,colframe=yellow!20!red, colbacktitle=yellow!20!red,
enhanced, attach boxed title to top center={yshift=-2mm}, title={#1}}

\newtcolorbox[auto counter, number within = section]{question}[2][]{%
fonttitle= \bfseries\upshape, fontupper= \upshape, arc=0mm,
colback=purple!5,colframe=purple!75!black, title={\linkt Question #2\linko},#1}

\definecolor{covercolor}{rgb}{0.7,0,0.2}
\definecolor{covertext}{rgb}{1,1,0.3}
\definecolor{darkblue}{rgb}{0,0,0.5}

\newcommand{\di}{\displaystyle}
\newcommand{\pa}{\partial}
\newcommand{\bookauthor}{Teo Lee Peng}
\newcommand{\booktitle}{Mathematical Analysis\\Volume II}
 
\newcommand{\bp}{\end{myproof}\begin{myproof}{}}
\newcommand{\bs}{\end{solution}\begin{solution}{}}
\newcommand{\be}{\end{example}\begin{example2}{}}
\newcommand{\vp}{\vfill\pagebreak}

\newcommand{\linkt}{\hypersetup{linkcolor=white}}
\newcommand{\linko}{\hypersetup{linkcolor=linkc}}
 \newcommand{\rvline}{\hspace*{-\arraycolsep}\vline\hspace*{-\arraycolsep}}
 \newcommand{\mf}{\mathbf}
 \newcommand{\mk}{\mathfrak}
  \newcommand{\mb}{\mathbb}
\begin{document}
  \begin{coverpage}
~\vspace{2cm}
\begin{center}
{\fontfamily{phv}\fontseries{mc}\fontsize{24}{28}\selectfont \textcolor{covertext}{  {\bfseries {\booktitle}} }

\vspace{6cm}

\textcolor{covertext}{{ \bfseries{\bookauthor}}}

}
\end{center}
 \end{coverpage}

\pagecolor{white}
\title{\uppercase{\booktitle}}

\author{\bookauthor}

\dedication{}

\date{\today}

\maketitle

 \frontmatter
\setcounter{page}{1}
\hypersetup{linkcolor=darkblue}

\tableofcontents
\linko

\chapter*{Preface}

Mathematical analysis is a standard course which introduces students to rigorous reasonings in mathematics, as well as the theories needed for advanced analysis courses. It is a compulsory course for all mathematics majors. It is also strongly recommended for students that major in computer science, physics, data science, financial analysis,  and other areas that require a lot of analytical skills. Some standard textbooks in mathematical analysis include the classical one by  Apostol \cite{Apostol} and Rudin \cite{Rudin}, and the modern one by Bartle \cite{Bartle}, Fitzpatrick \cite{Fitzpatrick}, Abbott \cite{Abbott}, Tao \cite{Tao_1, Tao_2} and Zorich \cite{Zorich_1, Zorich_2}.
 
 This book is the second volume of the textbooks intended for a one-year course in mathematical analysis.   We introduce the fundamental concepts in a pedagogical way. Lots of examples are given to illustrate the theories.  We assume that students are familiar with the material of calculus such as those in the book \cite{Stewart}.  Thus, we do not emphasize on the computation techniques. Emphasis is put on building up  analytical skills through rigorous reasonings. 

Besides calculus, it is also assumed that students have taken introductory courses in discrete mathematics and linear algebra, which covers topics such as logic, sets, functions, vector spaces, inner products, and quadratic forms. Whenever needed, these concepts would be briefly revised. 

In this book, we have defined all the mathematical terms we use carefully. While most of the terms have standard definitions, some of the terms may have definitions  defer from authors to authors. The readers are advised to check the definitions of the terms used in this book when they encounter them. This can be easily done by using the search function provided by any PDF viewer. The readers are also encouraged to fully utilize the hyper-referencing provided.

 \vspace{0.9cm}
~\hfill\bookauthor


\mainmatter

\chapter{Euclidean Spaces} \label{chapter1}

 In this second volume of mathematical analysis, we  study functions defined on   subsets of $\mathbb{R}^n$. For this, we need to study the structure and topology of $\mathbb{R}^n$ first. We start by a revision on $\mathbb{R}^n$ as a vector space.
 
 In the sequel, $n$ is  a fixed positive integer reserved to be used for  $\mathbb{R}^n$.
 
\section[The Euclidean Space $\mathbb{R}^n$ as a Vector Space]{The Euclidean Space $\pmb{\mathbb{R}^n}$ as a Vector Space} 
If $S_1$, $S_2$, $\ldots$, $S_n$ are sets, the cartesian product of these $n$ sets is defined as the set
\[S=S_1\times\cdots\times S_n=\prod_{i=1}^n S_i=\left\{(a_1, \ldots, a_n)\,|\, a_i\in S_i, 1\leq i\leq n\right\}\]that contains all $n$-tuples $(a_1, \ldots, a_n)$, where $a_i\in S_i$ for all $1\leq i\leq n$.

 The set $\mathbb{R}^n$ is the cartesian product of $n$ copies of $\mathbb{R}$. Namely,
\[\mathbb{R}^n=\left\{(x_1, x_2,\ldots, x_n)\,|\, x_1, x_2, \ldots, x_n\in\mathbb{R}\right\}.\]The point $(x_1, x_2, \ldots, x_n)$ is denoted as $\mathbf{x}$, whereas $x_1, x_2, \ldots, x_n$ are called the components of the point $\mathbf{x}$. We can define an addition and a scalar multiplication on $\mathbb{R}^n$. If $\mathbf{x}=(x_1, x_2, \ldots, x_n)$ and $\mathbf{y}=(y_1, y_2, \ldots, y_n)$ are in $\mathbb{R}^n$, the addition of $\mathbf{x}$ and $\mathbf{y}$ is defined as
\[\mathbf{x}+\mathbf{y}=(x_1+y_1, x_2+y_2, \ldots, x_n+y_n).\]
In other words, it is a componentwise addition. Given a real number $\alpha$, the scalar multiplication of $\alpha$ with $\mathbf{x}$ is given by the componentwise multiplication
\[\alpha\mathbf{x}=(\alpha x_1, \alpha x_2, \ldots, \alpha x_n).\]

The set $\mathbb{R}^n$ with the addition and scalar multiplication operations is a vector space. It satisfies the  10 axioms for a real vector space $V$.

\newcommand{\ve}{\vspace{0.4cm}}

\begin{highlight}{The 10 Axioms for a Real Vector Space $\pmb{V}$}

\vspace{0.2cm}
Let $V$ be a set that is equipped with two operations -- the addition   and the scalar multiplication. For any two vectors $\mathbf{u}$ and $\mathbf{v}$ in $V$, their addition is denoted by $\mathbf{u}+\mathbf{v}$. For a vector $\mathbf{u}$ in $V$ and a scalar $\alpha\in\mathbb{R}$, the scalar multiplication of $\mathbf{v}$ by $\alpha$ is denoted by $\alpha\mathbf{v}$. We say that $V$ with the addition and scalar multiplication is a real vector space provided that the following 10 axioms are satisfied for any $\mathbf{u}$, $\mathbf{v}$ and $\mathbf{w}$   in $V$, and any $\alpha$ and  $\beta$ in $\mathbb{R}$.

\ve
\textbf{\underline{Axiom 1}}
\quad
If $\mathbf{u}$ and $\mathbf{v}$ are in $V$, then $\mathbf{u}+\mathbf{v}$ is in $V$.

\ve
\textbf{\underline{Axiom 2}}
\quad $\mathbf{u}+\mathbf{v}   =\mathbf{v}+\mathbf{u}$.

\ve
\textbf{\underline{Axiom 3}}
\quad $\di (\mathbf{u}+ \mathbf{v})+\mathbf{w}= \mathbf{u}+(\mathbf{v}+\mathbf{w})$.

\ve
\textbf{\underline{Axiom 4}}
\quad There is a zero vector $\mathbf{0}$ in $V$ such that
\[\mathbf{0}+\mathbf{v}=\mathbf{v}=\mathbf{v}+\mathbf{0}\hspace{1cm}\text{for all}\;\mathbf{v}\in V.\]

\ve
\textbf{\underline{Axiom 5}}
\quad For any $\mathbf{v}$ in $V$, there is a vector $\mathbf{w}$ in $V$ such that
\[\mathbf{v}+\mathbf{w}=\mathbf{0}=\mathbf{w}+\mathbf{v}.\]The vector $\mathbf{w}$ satisfying this equation is called the {\it negative} of $\mathbf{v}$, and is denoted by $-\mathbf{v}$.

\ve
\textbf{\underline{Axiom 6}}
\quad For any $\mathbf{v}$ in $V$, and any $\alpha\in\mathbb{R}$, $\alpha\mathbf{v}$ is in $V$.

\ve
\textbf{\underline{Axiom 7}}
 \quad $\di\alpha(\mathbf{u}+\mathbf{v})=\alpha\mathbf{u}+\alpha \mathbf{v}$.
 
 \ve
\textbf{\underline{Axiom 8}}
\quad $\di(\alpha+\beta)\mathbf{v}=\alpha\mathbf{v}+\beta\mathbf{v}$.

\ve
\textbf{\underline{Axiom 9}}
\quad $\di\alpha(\beta\mathbf{v})=(\alpha\beta)\mathbf{v}$.

\ve
\textbf{\underline{Axiom 10}}
\quad $\di 1\mathbf{v}=\mathbf{v}$.
\end{highlight}
$\mathbb{R}^n$ is a real vector space. The zero vector is  the point $\mathbf{0}=(0, 0, \ldots, 0)$ with all components equal to 0. Sometimes we also call a point $\mathbf{x}=(x_1, \ldots, x_n)$ in $\mathbb{R}^n$ a vector, and identify it as the vector from the origin $\mathbf{0}$ to the point $\mathbf{x}$.

\begin{definition}{Standard Unit Vectors}
In $\mathbb{R}^n$, there are $n$ standard unit vectors $\mathbf{e}_1$, $\ldots$, $\mathbf{e}_n$  given by 
\[\mathbf{e}_1=(1, 0, \ldots, 0), \;\mathbf{e}_2=(0, 1, \ldots, 0), \;\cdots,\; \mathbf{e}_n=(0, \ldots, 0, 1).\]

\end{definition}

Let us  review some concepts from linear algebra which will be useful later.
Given that  $\mathbf{v}_1, \ldots, \mathbf{v}_k$ are vectors in a vector space $V$, a linear combination of  $\mathbf{v}_1, \ldots, \mathbf{v}_k$ is a vector $\mathbf{v}$ in $V$ of the form
\[\mathbf{v}=c_1\mathbf{v}_1+\cdots+c_k\mathbf{v}_k\] for some scalars $c_1, \ldots, c_k$, which are known as the coefficients of the linear combination.

A subspace of a vector space $V$  is a subset of $V$ that is itself a vector space. There is a simple way to construct subspaces.
\begin{proposition}{}
Let $V$ be a vector space, and let  $\mathbf{v}_1, \ldots, \mathbf{v}_k$ be vectors in $V$. The subset  \[W=\left\{c_1\mathbf{v}_1+\cdots+c_k\mathbf{v}_k\,|\, c_1, \ldots, c_k\in\mathbb{R}\right\}\] of $V$ that contains all linear combinations of  $\mathbf{v}_1, \ldots, \mathbf{v}_k$ is itself a vector space. It is called the subspace of $V$ {\it spanned} by  $\mathbf{v}_1, \ldots, \mathbf{v}_k$.
\end{proposition}

\begin{example}{}
In $\mathbb{R}^3$, the subspace spanned by the vectors $\mathbf{e}_1=(1,0,0)$ and $\mathbf{e}_3=(0,0,1)$ is the set $W$ that contains all points of the form
 \[x(1,0,0)+z(0,0,1)=(x,0,z),\] which is the $xz$-plane.
\end{example}

Next, we recall the concept of linear independence.
\begin{definition}{Linear Independence}
Let $V$ be a vector space, and let  $\mathbf{v}_1, \ldots, \mathbf{v}_k$ be vectors in $V$. 
We say that the set  $\{\mathbf{v}_1, \ldots, \mathbf{v}_k\}$ is a linearly independent set of vectors, or the vectors $\mathbf{v}_1, \ldots, \mathbf{v}_k$ are linearly independent, if the only $k$-tuple of real numbers $(c_1, \ldots, c_k)$ which satisfies
\[c_1\mathbf{v}_1+\cdots+c_k\mathbf{v}_k=\mathbf{0}\] is the {\it trivial} $k$-tuple $(c_1, \ldots, c_k)=(0,\ldots,0)$.
\end{definition}

 \begin{example}{}
In $\mathbb{R}^n$, the standard unit vectors $\mathbf{e}_1, \ldots, \mathbf{e}_n$ are linearly independent.
\end{example}

 \begin{example}{}
If $V$ is a vector space, a vector $\mathbf{v}$ in $V$ is linearly independent if and only if $\mathbf{v}\neq\mathbf{0}$.
\end{example}

 \begin{example}{}
Let $V$ be a vector space. Two vectors  $\mathbf{u}$ and $\mathbf{v}$ in $V$ are linearly independent if and only if $\mathbf{u}\neq \mathbf{0}$, $\mathbf{v}\neq \mathbf{0}$, and there does not exists a constant $\alpha$ such that $\mathbf{v}=\alpha\mathbf{u}$.
\end{example}

Let us recall the following definition for two vectors to be parallel.
\begin{definition}{Parallel Vectors}
 Let $V$ be a vector space. Two vectors  $\mathbf{u}$ and $\mathbf{v}$ in $V$ are parallel if either $\mathbf{u}= \mathbf{0}$ or there exists a constant $\alpha$ such that $\mathbf{v}=\alpha\mathbf{u}$.
\end{definition}In other words, two vectors $\mathbf{u}$ and $\mathbf{v}$ in $V$ are linearly independent if and only if they are not parallel.

\begin{example}{}
If $S=\{\mathbf{v}_1, \ldots, \mathbf{v}_k\}$ is a linearly independent set of vectors, then for any $S'\subset S$, $S'$ is also a linearly independent set of vectors.
\end{example}

Now we discuss the concept of dimension and basis.
\begin{definition}{Dimension and Basis}
Let $V$ be a vector space, and let $W$ be a subspace of $V$. If $W$ can be spanned by $k$ linearly independent vectors $\mathbf{v}_1, \ldots,\mathbf{v}_k$ in $V$, we say that $W$ has dimension $k$. The set $\{ \mathbf{v}_1, \ldots,\mathbf{v}_k\}$ is called a basis of $W$.
\end{definition}

\begin{example}{}
In $\mathbb{R}^n$, the  $n$ standard unit vectors $\mathbf{e}_1$, $\ldots$, $\mathbf{e}_n$  
 are linearly independent and they span $\mathbb{R}^n$.  Hence, the dimension of $\mathbb{R}^n$ is $n$.
\end{example}

\begin{example}{}
In $\mathbb{R}^3$, the subspace spanned by the two linearly independent vectors $\mathbf{e}_1=(1,0,0)$ and $\mathbf{e}_3=(0,0,1)$ has dimension 2.
\end{example}

Next, we introduce the translate of a set.
\begin{definition}{Translate of a Set}
If $A$ is a subset of $\mb{R}^n$, $\mf{u}$ is a point in $\mb{R}^n$, the translate of the set $A$ by the vector $\mf{u}$ is the set 
\[A+\mf{u}=\left\{\mf{a}+\mf{u}\,|\,\mf{a}\in A\right\}.\]  
\end{definition}
\begin{example}{}
In $\mathbb{R}^3$, the translate of the set $A=\{(x,y,0)\,|\, x,y\in\mb{R}\}$ by the vector $\mf{u}=(0,0,-2)$ is the set
\[B=A+\mf{u}=\{(x,y,-2)\,|\, x,y\in\mb{R}\}.\]
\end{example}

In $\mathbb{R}^n$, the lines and the planes  are of particular interest. They are closely related to the concept of subspaces.
\begin{definition}{Lines   in $\pmb{\mathbb{R}^n}$}
A line $L$ in $\mathbb{R}^n$ is a translate of a subspace of $\mb{R}^n$ that has dimension 1. As a set, it contains all the points  $\mf{x}$ of the form
\[\mathbf{x}=\mathbf{x}_0+t\mathbf{v},\hspace{1cm}t\in\mathbb{R},\] where $\mathbf{x}_0$ is a fixed point in $\mathbb{R}^n$, and $\mathbf{v}$ is a nonzero vector in $\mathbb{R}^n$. The equation $\mathbf{x}=\mathbf{x}_0+t\mathbf{v}$, $t\in\mathbb{R}$, is known as the parametric equation of the line.

\end{definition}

A line is determined by two points.

\begin{example}{}
Given two distinct points $\mathbf{x}_1$ and $\mathbf{x}_2$ in $\mathbb{R}^n$, the line $L$ that passes through these two points have parametric equation given by 
\[\mathbf{x}=\mathbf{x}_1+t(\mathbf{x}_2-\mathbf{x}_1),\hspace{1cm}t\in\mathbb{R}.\]When $0\leq t\leq 1$, $\mathbf{x}=\mathbf{x}_1+t(\mathbf{x}_2-\mathbf{x}_1)$ describes all the points on the line segment with $\mathbf{x}_1$ and $\mathbf{x}_2$ as endpoints.
\end{example}

 \begin{figure}[ht]
\centering
\includegraphics[scale=0.2]{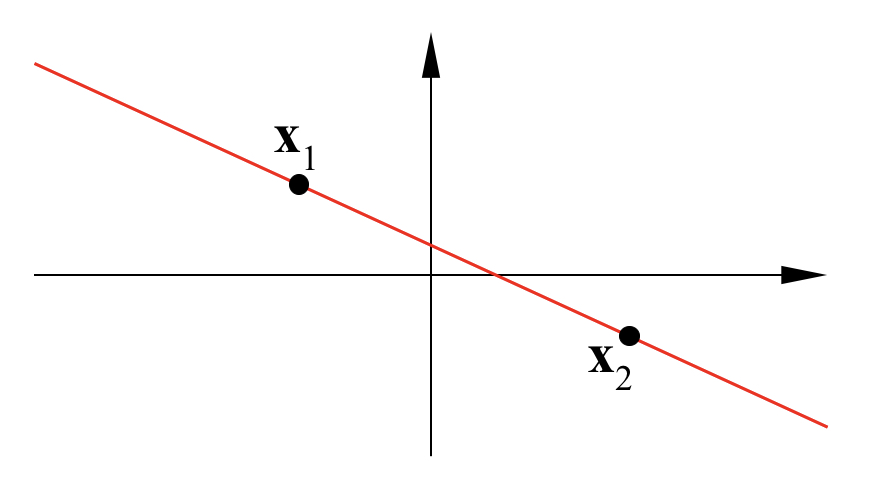}

\caption{A Line between two points.}\label{figure12}
\end{figure}

\begin{definition}{Planes   in $\pmb{\mathbb{R}^n}$}
A plane $W$ in $\mathbb{R}^n$  is a translate of a subspace of dimension 2. As a set, it contains all the points $\mf{x}$ of the form
\[\mathbf{x}=\mathbf{x}_0+t_1\mathbf{v}_1+t_2\mathbf{v}_2,\hspace{1cm}t_1, t_2\in\mathbb{R},\] where $\mathbf{x}_0$ is a fixed point in $\mathbb{R}^n$, and $\mathbf{v}_1$ and $\mathbf{v}_2$ are two linearly independent vectors in $\mathbb{R}^n$.

\end{definition}

Besides being a real vector space, $\mathbb{R}^n$ has  an additional structure. 
 Its definition is motivated as follows. Let $P(x_1, x_2, x_3)$ and $Q(y_1, y_2, y_3)$ be two points in $\mathbb{R}^3$. By Pythagoras theorem, the distance between $P$ and $Q$ is given by 
 \[PQ=\sqrt{(x_1-y_1)^2+(x_2-y_2)^2+(x_3-y_3)^2}.\]
 
 \begin{figure}[ht]
\centering
\includegraphics[scale=0.18]{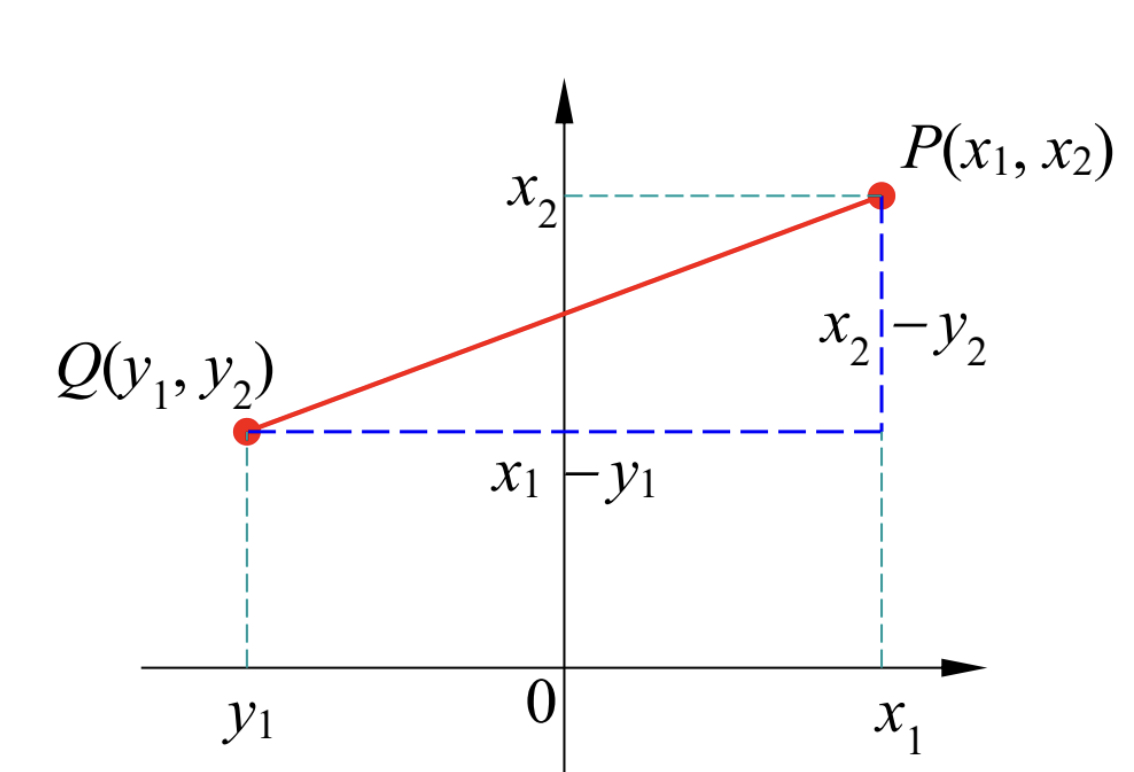}

\caption{Distance between two points in $\mb{R}^2$.}\label{figure2}
\end{figure}
  \begin{figure}[ht]
\centering
\includegraphics[scale=0.18]{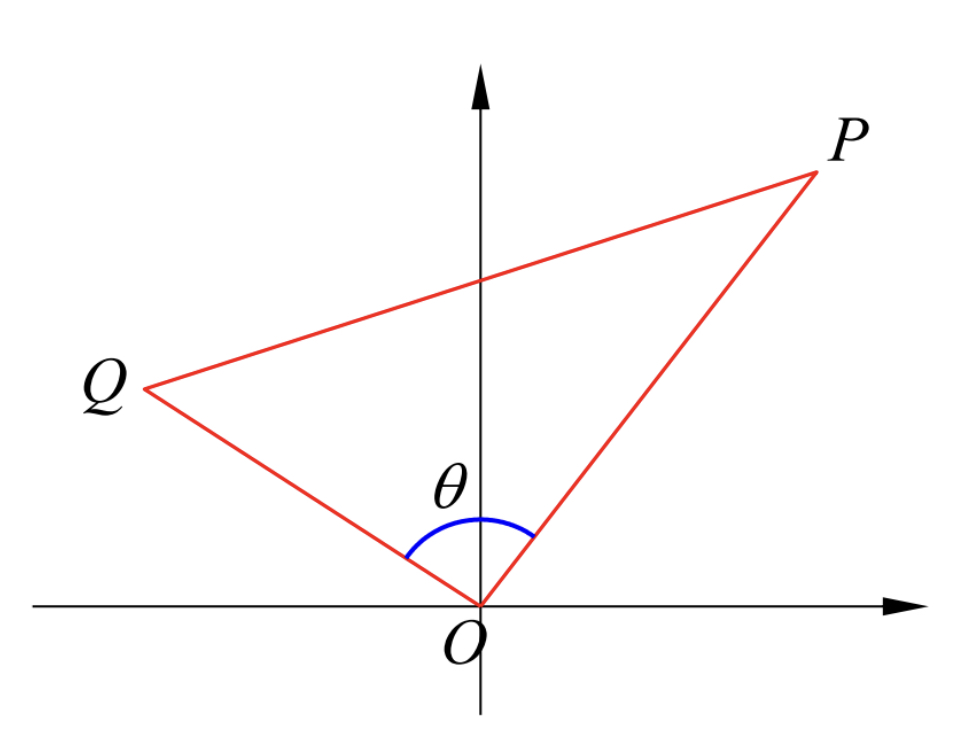}

\caption{Cosine rule.}\label{figure1}
\end{figure}
 
Consider the triangle $OPQ$ with vertices $O$, $P$, $Q$, where $O$ is the origin. Then    
 \begin{gather*}
 OP=\sqrt{x_1^2+x_2^2+x_3^2}, \hspace{1cm} OQ=\sqrt{y_1^2+y_2^2+y_3^2}.
 \end{gather*}Let $\theta$ be the minor angle between $OP$ and $OQ$. By cosine rule,
 \[PQ^2=OP^2+OQ^2-2\times OP\times OQ\times\cos\theta.\]
A straightforward computation gives
 \[OP^2+OQ^2-PQ^2=2(x_1y_1+x_2y_2+x_3y_3).\]
  Hence,
 \begin{equation}\label{eq230709_1}\cos\theta=\frac{x_1y_1+x_2y_2+x_3y_3}{\sqrt{x_1^2+x_2^2+x_3^2}\,\sqrt{y_1^2+y_2^2+y_3^2}}.\end{equation}It is a quotient of $x_1y_1+x_2y_2+x_3y_3$ by  the product of the lengths of $OP$ and $OQ$.

 Generalizing the expression $x_1y_1+x_2y_2+x_3y_3$ from $\mb{R}^3$ to $\mathbb{R}^n$ defines the dot product. For any two vectors $\mathbf{x}=(x_1, x_2, \ldots, x_n)$ and $\mathbf{y}=(y_1, y_2, \ldots, y_n)$ in $\mathbb{R}^n$, the dot product of $\mathbf{x}$ and $\mathbf{y}$ is defined as
 \[ \mathbf{x}\,\cdot\,\mathbf{y}=\sum_{i=1}^nx_iy_i=x_1y_1+x_2y_2+\cdots+x_ny_n.\]This is a special case of an inner product.
 
\begin{definition}{Inner Product Space}
A real vector space $V$ is an inner product space if for any two vectors $\mathbf{u}$ and $\mathbf{v}$ in $V$, an inner product $\langle\mathbf{u},\mathbf{v}\rangle$ of $\mathbf{u}$ and $\mathbf{v}$ is defined, and   the following conditions for any $\mathbf{u}, \mathbf{v}, \mathbf{w}$ in $V$ and $\alpha, \beta\in\mathbb{R}$ are satisfied. 
 \begin{enumerate}[1.]
 \item $\di \langle \mathbf{u}, \mathbf{v}\rangle =\langle \mathbf{v}, \mathbf{u}\rangle$.
 \item $\di \langle  \alpha\mathbf{u}+ \beta\mathbf{v},\mathbf{w}\rangle =\alpha\langle \mathbf{u}, \mathbf{w}\rangle+\beta\langle \mathbf{v}, \mathbf{w}\rangle$.
   \item $\di \langle\mathbf{v},\mathbf{v}\rangle\geq 0$ and 
 $\di \langle\mathbf{v},\mathbf{v}\rangle = 0$ if and only if $\mathbf{v}=\mathbf{0}$.
 \end{enumerate}
 \end{definition}

  \begin{proposition}{Euclidean Inner Product on $\pmb{\mathbb{R}^n}$}
On $\mathbb{R}^n$,  \[\langle\mathbf{x},\mathbf{y}\rangle=\mathbf{x}\,\cdot\,\mathbf{y}=\sum_{i=1}^nx_iy_i=x_1y_1+x_2y_2+\cdots+x_ny_n.\] defines an inner product, called the standard inner product or the Euclidean inner product. 
 
 \end{proposition}

 \begin{definition}{Euclidean Space}
 
 The vector  space $\mathbb{R}^n$ with the Euclidean inner product is called the Euclidean $n$-space. 
 \end{definition}
 In the future, when we do not specify, $\mb{R}^n$ always means the Euclidean $n$-space.
 
One can deduce some useful identities from the three axioms of an inner product space.
\begin{proposition}{} If $V$ is an inner product space, then the following holds.
\begin{enumerate}[(a)]
\item For any $\mathbf{v}\in V$,
$\di\langle\mathbf{0}, \mathbf{v}\rangle =0=\langle\mathbf{v},\mathbf{0}\rangle$.
 \item For any vectors $\mathbf{v}_1, \cdots,\mathbf{v}_k$, $\mathbf{w}_1, \cdots,\mathbf{w}_l$  in $V$,  and for any real numbers $\alpha_1, \cdots, \alpha_k$, $\beta_1, \cdots, \beta_l$,
 \[\left\langle\sum_{i=1}^k\alpha_i\mathbf{v}_i, \sum_{j=1}^l \beta_j\mathbf{w}_j\right\rangle =\sum_{i=1}^k\sum_{j=1}^l \alpha_i\beta_j\langle\mathbf{v}_i, \mathbf{w}_j\rangle.\]
 \end{enumerate}\end{proposition}

 Given that $V$ is an  inner product space, $\langle\mathbf{v}, \mathbf{v}\rangle\geq 0$ for any $\mathbf{v}$ in $V$. For example, for any  $\mathbf{x}=(x_1, x_2, \ldots, x_n)$ in $\mathbb{R}^n$, under the Euclidean inner product, 
 \[\langle\mathbf{x}, \mathbf{x}\rangle= \sum_{i=1}^nx_i^2 = x_1^2+x_2^2+\cdots+x_n^2\geq 0.\]
 When $n=3$,  the length of the vector $OP$ from the point $O(0,0,0)$ to the point $P(x_1, x_2, x_3)$ is
 \[OP=\sqrt{x_1^2+x_2^2+x_3^2}=\sqrt{\langle\mathbf{x}, \mathbf{x}\rangle},\hspace{1cm}\text{where}\;\mathbf{x}=(x_1, x_2, x_3).\]
 This motivates us to define to norm of a vector in an inner product space as follows.
 \begin{definition}{Norm of a Vector}
 Given that $V$ is an inner product space, the norm of a vector $\mathbf{v}$ is defined as  \[\Vert\mathbf{v}\Vert =\sqrt{\langle \mathbf{v},\mathbf{v}\rangle}.\]
 \end{definition}

 The norm of a vector in an inner product space satisfies some properties, which follow from the axioms for an inner product space.
 \begin{proposition}{}
 Let $V$ be an inner product space.
 \begin{enumerate}[1.]
 \item
 For any $\mathbf{v}$ in $V$, $\Vert \mathbf{v}\Vert\geq 0$ and $\Vert\mathbf{v}\Vert=0$ if and only if $\mathbf{v}=\mathbf{0}$.
 \item For any $\alpha\in\mathbb{R}$ and $\mathbf{v}\in V$, $\Vert\alpha\mathbf{v}\Vert\;=\;|\alpha|\,\Vert\mathbf{v}\Vert$.
 \end{enumerate}
 \end{proposition}
 Motivated by the distance between two points in $\mathbb{R}^3$, we   make the following definition.
  \begin{definition}{Distance Between Two Points}
 Given that $V$ is an inner product space, the distance between $\mathbf{u}$ and $\mathbf{v}$ in $V$ is defined as  \[d(\mathbf{u}, \mathbf{v}) =\Vert\mathbf{v}-\mathbf{u}\Vert=\sqrt{\langle \mathbf{v}-\mathbf{u},\mathbf{v}-\mathbf{u}\rangle}.\]
 \end{definition}
 For example, the distance between the points $\mathbf{x}=(x_1, \ldots, x_n)$ and $\mathbf{y}=(y_1, \ldots, y_n)$ in the Euclidean space $\mathbb{R}^n$ is
 \[d(\mathbf{x}, \mathbf{y})=\sqrt{\sum_{i=1}^n(x_i-y_i)^2}=\sqrt{(x_1-y_1)^2+\cdots+(x_n-y_n)^2}.\]
 
 For analysis in $\mathbb{R}$, an important inequality is the triangle inequality which says that $|x+y|\leq |x|+|y|$ for any $x$ and $y$ in $\mathbb{R}$. To generalize this inequality to $\mathbb{R}^n$, we need
 the celebrated Cauchy-Schwarz inequality. It holds on any inner product space.
 \begin{proposition}{Cauchy-Schwarz Inequality}
 Given that $V$ is an inner product space, for any $\mathbf{u}$ and $\mathbf{v}$ in $V$,
 \[\left|\langle\mathbf{u}, \mathbf{v}\rangle \right|\;\leq\; \Vert\mathbf{u}\Vert\,\Vert\mathbf{v}\Vert.\]The equality holds if and only if $\mathbf{u}$ and $\mathbf{v}$ are parallel.
 \end{proposition}
 \begin{myproof}{Proof}
It is obvious that if either $\mathbf{u}=\mathbf{0}$ or $\mathbf{v}=\mathbf{0}$, 
  \[\left|\langle\mathbf{u}, \mathbf{v}\rangle \right|\;=\;0\;=\; \Vert\mathbf{u}\Vert\,\Vert\mathbf{v}\Vert,\] and so the equality holds.
  
 Now assume that both $\mathbf{u}$ and $\mathbf{v}$ are nonzero vectors. Consider the quadratic function $f:\mathbb{R}\to\mathbb{R}$ defined by
 \[f(t)=\Vert t\mathbf{u}- \mathbf{v}\Vert^2=\langle t\mathbf{u}- \mathbf{v}, t\mathbf{u}-\mathbf{v}\rangle.\]
 Notice that $f(t)=at^2+bt+c$, where
 \[a=\langle \mathbf{u},\mathbf{u}\rangle=\Vert\mathbf{u}\Vert^2, \quad b=-2\langle\mathbf{u},\mathbf{v}\rangle,\quad c=\langle\mathbf{v},\mathbf{v}\rangle=\Vert\mathbf{v}\Vert^2.\] 
 The 3$^{\text{rd}}$ axiom of an inner product says that $f(t)\geq 0$ for all $t\in\mathbb{R}$. Hence, we must have $b^2-4ac\leq 0$. This gives
 \[\langle\mathbf{u}, \mathbf{v}\rangle^2\leq \Vert\mathbf{u}\Vert^2\Vert\mathbf{v}\Vert^2.\]
 Thus, we obtain the Cauchy-Schwarz inequality
  \[\left|\langle\mathbf{u}, \mathbf{v}\rangle \right|\;\leq\; \Vert\mathbf{u}\Vert\,\Vert\mathbf{v}\Vert.\]
  \bp
  
 The equality holds if and only if $b^2-4ac=0$. The latter means that $f(t)=0$ for some $t=\alpha$, which can happen if and only if \[\alpha\mathbf{u}-\mathbf{v}=\mathbf{0},\] or equivalently, $\mathbf{v}=\alpha\mathbf{u}$.
 
 \end{myproof}
 
 Now we can prove the triangle inequality.
 \begin{proposition}{Triangle Inequality}
 Let $V$ be an inner product space. For any vectors $\mathbf{v}_1, \mathbf{v}_2, \ldots, \mathbf{v}_{k}$   in $V$, 
 \[\Vert \mathbf{v}_1+\mathbf{v}_2+\cdots+\mathbf{v}_k\Vert\leq \Vert\mathbf{v}_1\Vert+\Vert\mathbf{v}_2\Vert+\cdots+\Vert\mathbf{v}_k\Vert.\]
 \end{proposition}
 \begin{myproof}{Proof}
 It is sufficient to prove the statement when $k=2$. The general case follows from induction.
 Given $\mathbf{v}_1$ and $\mathbf{v}_2$ in $V$, 
 \begin{align*}
 \Vert\mathbf{v}_1+\mathbf{v}_2\Vert^2&=\langle \mathbf{v}_1+\mathbf{v}_2, \mathbf{v}_1+\mathbf{v}_2\rangle\\
 &=\langle \mathbf{v}_1, \mathbf{v}_1\rangle +2\langle \mathbf{v}_1,\mathbf{v}_2\rangle +\langle\mathbf{v}_2,\mathbf{v}_2\rangle\\
 &\leq \Vert\mathbf{v}_1\Vert^2+2\Vert\mathbf{v}_1\Vert\Vert\mathbf{v}_2\Vert+\Vert\mathbf{v}_2\Vert^2\\
 &=\left(\Vert\mathbf{v}_1\Vert+\Vert\mathbf{v}_2\Vert\right)^2.
 \end{align*}
 
 This proves that
 \[\Vert \mathbf{v}_1+\mathbf{v}_2\Vert\leq \Vert\mathbf{v}_1\Vert+\Vert\mathbf{v}_2\Vert.\]
 \end{myproof}
 
From the triangle inequality, we can deduce the following.
 \begin{corollary}{}
 Let $V$ be an inner product space. For any vectors $\mathbf{u}$ and $\mathbf{v}$ in $V$,
 \[\bigl|\Vert\mathbf{u}\Vert -\Vert\mathbf{v}\Vert\bigr|\leq \Vert \mathbf{u}-\mathbf{v}\Vert.\]
 \end{corollary}
 
 Express in terms of distance, the triangle inequality takes the following form.
  \begin{proposition}{Triangle Inequality}
 Let $V$ be an inner product space. For any three points  $\mathbf{v}_1, \mathbf{v}_2,   \mathbf{v}_{3}$   in $V$, 
 \[d(\mathbf{v}_1,\mathbf{v}_2)\leq d(\mathbf{v}_1, \mathbf{v}_3)+d(\mathbf{v}_2, \mathbf{v}_3).\]
 \end{proposition}
 More generally, if $\mathbf{v}_1, \mathbf{v}_2, \ldots, \mathbf{v}_k$ are $k$ vectors in $V$, then
 \[d(\mathbf{v}_1, \mathbf{v}_k)\leq\sum_{i=2}^{k}d(\mathbf{v}_{i-1}, \mathbf{v}_i)= d(\mathbf{v}_1, \mathbf{v}_2)+\cdots +d(\mathbf{v}_{k-1}, \mathbf{v}_k).\]
 
 Since we can define the distance function on an inner product space, inner product space is a special case of metric spaces. 
 \begin{definition}{Metric Space}
Let $X$ be a set,  and let $d:X\times X\to\mathbb{R}$ be  a function defined on $X\times X$. We say that $d$ is a metric on $X$ provided that   the following conditions are satisfied.
 \begin{enumerate}[1.]
 \item For any $x$ and $y$ in $X$, $d(x, y)\geq 0$, and $d(x,y)=0$ if and only if $x=y$.
 \item $d(x, y)=d(y,x)$ for any $x$ and $y$ in $X$.
 \item For any $x$, $y$ and $z$   in $X$, $d(x,y)\leq d(x, z)+d(y,z)$.
 
 \end{enumerate}
If $d$ is a metric on $X$, we say that $(X, d)$ is a metric space.
  \end{definition}Metric spaces play important roles in advanced analysis. 
 If $V$ is an innner product space, it is a metric space with metric
 \[d(\mathbf{u}, \mathbf{v})=\Vert\mathbf{v}-\mathbf{u}\Vert.\]

 Using the Cauchy-Schwarz inequality, one can generalize the concept of angles to any two vectors in a real inner product space.  If $\mathbf{u}$ and $\mathbf{v}$ are two nonzero vectors in a real inner product space $V$, Cauchy-Schwarz inequality implies that
 \[\frac{\langle\mathbf{u}, \mathbf{v}\rangle}{\Vert\mathbf{u}\Vert\,\Vert\mathbf{v}\Vert}\] is a real number between $-1$ and 1. Generalizing the formula \eqref{eq230709_1}, we define the angle $\theta$ between $\mathbf{u}$ and $\mathbf{v}$ as
 \[\theta= \cos^{-1} \frac{\langle\mathbf{u}, \mathbf{v}\rangle}{\Vert\mathbf{u}\Vert\,\Vert\mathbf{v}\Vert}.\]This is an angle between $0^{\circ}$ and $180^{\circ}$. A necessary and sufficient condition for two vectors $\mathbf{u}$ and $\mathbf{v}$ to make a $90^{\circ}$ angle is $\langle \mathbf{u}, \mathbf{v}\rangle=0$. 
 \begin{definition}{Orthogonality}
 Let $V$ be a real inner product space. We say that the two vectors  $\mathbf{u}$ and $\mathbf{v}$ in $V$ are orthogonal if $\langle\mathbf{u}, \mathbf{v}\rangle=0$.
 \end{definition}

\begin{lemma}{Generalized Pythagoras Theorem}
Let $V$ be an inner product space. If $\mathbf{u}$ and $\mathbf{v}$ are orthogonal vectors in $V$, then
\[\Vert\mathbf{u}+\mathbf{v}\Vert^2=\Vert\mathbf{u}\Vert^2+\Vert\mathbf{v}\Vert^2.\]
\end{lemma}

Now we discuss the projection theorem.

\begin{theorem}{Projection Theorem}
Let $V$ be an inner product space, and let $\mathbf{w}$ be a nonzero vector in $V$. If $\mathbf{v}$ is a  vector in $V$, there is a unique way to write $\mathbf{v}$ as a sum of two vectors $\mathbf{v}_1$ and $\mathbf{v}_2$, such that $\mathbf{v}_1$ is parallel to $\mathbf{w}$ and $\mathbf{v}_2$ is orthogonal to $\mathbf{w}$.  Moreover, for any real number $\alpha$, 
\[\Vert \mathbf{v}-\alpha\mathbf{w}\Vert\geq \Vert\mathbf{v}-\mathbf{v_1}\Vert,\]
and the equality holds if and only if $\alpha$ is equal to the unique real number $\beta$ such that $\mathbf{v}_1=\beta\mathbf{w}$. 
\end{theorem}
  
 \begin{figure}[ht]
\centering
\includegraphics[scale=0.2]{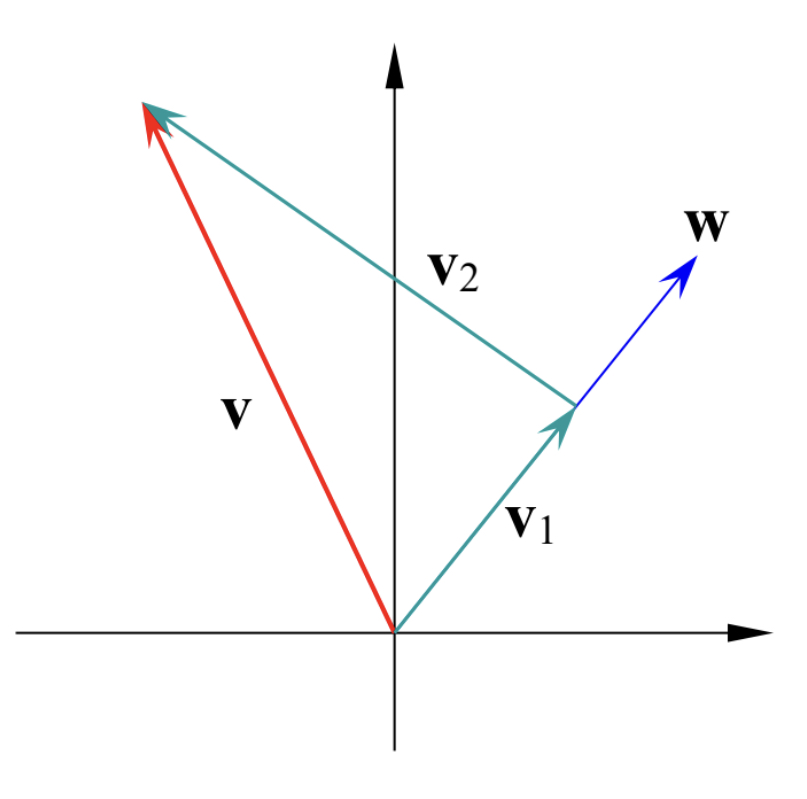}

\caption{The projection theorem.}\label{figure3}
\end{figure}
\begin{myproof}{Proof}Assume that  $\mathbf{v}$ can be written as a sum of two vectors $\mathbf{v}_1$ and $\mathbf{v}_2$, such that $\mathbf{v}_1$ is parallel to $\mathbf{w}$ and $\mathbf{v}_2$ is orthogonal to $\mathbf{w}$.
Since $\mathbf{w}$ is nonzero,   there is a real number $\beta$ such that $\mathbf{v}_1=\beta\mathbf{w}$. 
Since $\mathbf{v}_2=\mathbf{v}-\mathbf{v}_1=\mathbf{v}-\beta\mathbf{w}$ is orthogonal to $\mathbf{w}$, we have 
\[0=\langle \mathbf{v}-\beta\mathbf{w},\mathbf{w}\rangle=\langle\mathbf{v},\mathbf{w}\rangle-\beta\langle\mathbf{w},\mathbf{w}\rangle.\]
This implies that we must have
\[\beta=\frac{\langle\mathbf{v},\mathbf{w}\rangle}{\langle\mathbf{w},\mathbf{w}\rangle},\]
and
\[\mathbf{v}_1=\frac{\langle\mathbf{v},\mathbf{w}\rangle}{\langle\mathbf{w},\mathbf{w}\rangle}\mathbf{w}, \hspace{1cm}\mathbf{v}_2=\mathbf{v}-\frac{\langle\mathbf{v},\mathbf{w}\rangle}{\langle\mathbf{w},\mathbf{w}\rangle}\mathbf{w}.\]

It is easy to check that $\mathbf{v}_1$ and $\mathbf{v}_2$ given by these formulas indeed satisfy the requirements that $\mathbf{v}_1$ is parallel to $\mathbf{w}$ and $\mathbf{v}_2$ is orthogonal to $\mathbf{w}$. This establishes the existence and uniqueness of $\mathbf{v}_1$ and $\mathbf{v}_2$.
 
Now for any real number $\alpha$, 
\[\mathbf{v}-\alpha\mathbf{w}=\mathbf{v}-\mathbf{v}_1+(\beta-\alpha)\mathbf{w}.\] 
\bp
Since $\mathbf{v}-\mathbf{v}_1=\mathbf{v}_2$ is orthogonal to $(\beta-\alpha)\mathbf{w}$, the generalized Pythagoras theorem implies that
\[\Vert \mathbf{v}-\alpha\mathbf{w}\Vert^2=\Vert\mathbf{v}-\mathbf{v}_1\Vert^2+\Vert(\beta-\alpha)\mathbf{w}\Vert^2\geq\Vert\mathbf{v}-\mathbf{v}_1\Vert^2. \]
This proves that \[\Vert \mathbf{v}-\alpha\mathbf{w}\Vert\geq \Vert\mathbf{v}-\mathbf{v}_1\Vert.\] The equality holds if and only if
\[\Vert(\beta-\alpha)\mathbf{w}\Vert=|\alpha-\beta|\Vert\mathbf{w}\Vert=0.\]Since $\Vert\mathbf{w}\Vert\neq 0$, we must have $\alpha=\beta$.
\end{myproof}The vector $\mathbf{v}_1$ in this theorem is called the projection of $\mathbf{v}$ onto the subspace spanned by $\mathbf{w}$.

There is a more general projection theorem where the subspace $W$ spanned by $\mathbf{w}$ is replaced by a general subspace. We say that a vector $\mathbf{v}$ is orthogonal to the subspace $W$ if it is orthogonal to each vector $\mathbf{w}$ in $W$.
\begin{theorem}{General Projection Theorem}
Let $V$ be an inner product space, and let $W$ be a finite dimensional subspace of $V$. If $\mathbf{v}$ is a  vector in $V$, there is a unique way to write $\mathbf{v}$ as a sum of two vectors $\mathbf{v}_1$ and $\mathbf{v}_2$, such that $\mathbf{v}_1$ is in $W$ and $\mathbf{v}_2$ is orthogonal to $W$.  The vector $\mathbf{v}_1$   is denoted by $\text{proj}_W\mathbf{v}$. For any $\mathbf{w}\in W$,
\[\Vert \mathbf{v}- \mathbf{w}\Vert\geq \Vert\mathbf{v}-\text{proj}_W\mathbf{v}\Vert,\]
and the equality holds if and only if   $\mathbf{w} = \text{proj}_W\mathbf{v}$.  
\end{theorem}
\begin{myproof}{Sketch of Proof}
If $W$ is a  $k$- dimensional vector space, it has a basis consists of $k$ linearly independent vectors $\mathbf{w}_1, \ldots, \mathbf{w}_k$.  Since the vector $\mathbf{v}_1$ is in $W$, there are constants $c_1, \ldots, c_k$ such that
\[\mathbf{v}_1=c_1\mathbf{w}_1+\cdots+c_k\mathbf{w}_k.\]
\bp
The condition $\mathbf{v}_2=\mathbf{v}-\mathbf{v}_1$ is orthogonal to $W$ gives rise to $k$ equations
\begin{equation}\label{230710_1}\begin{split}c_1\langle \mathbf{w}_1, \mathbf{w}_1\rangle+ \cdots+c_k\langle \mathbf{w}_k,\mathbf{w}_1\rangle&=\langle\mathbf{v},\mathbf{w}_1\rangle,\\
 \vdots \hspace{2cm} &\\
c_1\langle \mathbf{w}_1, \mathbf{w}_k\rangle +\cdots+c_k\langle \mathbf{w}_k,\mathbf{w}_k\rangle&=\langle\mathbf{v},\mathbf{w}_k\rangle.
\end{split}\end{equation}
Using the fact that $\mathbf{w}_1, \ldots, \mathbf{w}_k$ are linearly independent, one can show that the $k\times k$ matrix 
\[A=\begin{bmatrix}\langle \mathbf{w}_1, \mathbf{w}_1\rangle&  \cdots & \langle \mathbf{w}_k,\mathbf{w}_1\rangle \\
  \vdots & \ddots &\vdots\\
 \langle \mathbf{w}_1, \mathbf{w}_k\rangle  &\cdots & \langle \mathbf{w}_k,\mathbf{w}_k\rangle\end{bmatrix}\]
 is invertible.  This shows that there is a unique $\mathbf{c}=(c_1, \ldots, c_k)$ satisfying the linear system \eqref{230710_1}.\end{myproof}

 If $V$ is an inner product space, a basis that consists of mutually orthogonal vectors are of special interest.

\begin{definition}{Orthogonal Set and Orthonormal Set}
Let $V$ be an inner product space. A subset of vectors $S=\{\mathbf{u}_1, \ldots, \mathbf{u}_k\}$ is called an orthogonal set if any two distinct vectors $\mathbf{u}_i$ and $\mathbf{u}_j$ in $S$ are orthogonal. Namely,
\[\langle\mathbf{u}_i, \mathbf{u}_j\rangle =0\hspace{1cm}\text{if}\; i\neq j.\]
$S$ is called an orthonormal set if it is an orthogonal set of unit vectors. Namely,
\[\langle\mathbf{u}_i, \mathbf{u}_j\rangle =\begin{cases}0\quad &\text{if}\; i\neq j, \\1\quad &\text{if}\; i= j\end{cases}.\]
\end{definition}
If $S=\{\mathbf{u}_1, \ldots, \mathbf{u}_k\}$ is an orthogonal set of nonzero vectors,   it is a linearly independent set of vectors. One can construct an orthonormal set by normalizing each vector  in the set. There is a standard algorithm, known as the Gram-Schmidt process, which can turn any linearly independent set of vectors $\{\mathbf{v}_1, \ldots, \mathbf{v}_k\}$ into an orthogonal set $\{\mathbf{u}_1, \ldots, \mathbf{u}_k\}$ of nonzero vectors. We start by the following lemma.

\begin{lemma}{}
Let $V$ be an inner  product space, and let $S=\{\mathbf{u}_1, \ldots, \mathbf{u}_k\}$ be an orthogonal set of nonzero vectors in $V$ that spans the subspace $W$. Given any vector $\mathbf{v}$ in $V$,
\[\text{proj}_W\mathbf{v}=\sum_{i=1}^k \frac{\langle\mathbf{v}, \mathbf{u}_i\rangle}{\langle\mathbf{u}_i, \mathbf{u}_i\rangle }\mathbf{u}_i.\]
\end{lemma}
\begin{myproof}{Proof}
By the general projection theorem, $\mathbf{v}=\mathbf{v}_1+\mathbf{v}_2$, where $\mathbf{v}_1=\text{proj}_W\mathbf{v}$ is in $W$ and $\mathbf{v}_2$ is orthogonal to $W$. 
Since $S$ is a basis for $W$, there exist scalars $c_1, c_2, \ldots, c_k$ such that
$\mathbf{v}_1=c_1\mathbf{u}_1+\cdots+c_k\mathbf{u}_k$. 
Therefore,
\[\mathbf{v}=c_1\mathbf{u}_1+\cdots+c_k\mathbf{u}_k+\mathbf{v}_2.\]Since $S$ is an orthogonal set of vectors and $\mathbf{v}_2$ is orthogonal to each $\mathbf{u}_i$, we find that for $1\leq i\leq k$,
\[\langle \mathbf{v}, \mathbf{u}_i\rangle =c_i\langle\mathbf{u}_i,\mathbf{u}_i\rangle.\]
This proves the lemma.
\end{myproof}
\begin{theorem}{Gram-Schmidt Process}
Let $V$ be an inner product space, and assume that $S=\{\mathbf{v}_1, \ldots, \mathbf{v}_k\}$ is a linearly independent set of vectors in $V$. Define the vectors $\mathbf{u}_1, \ldots, \mathbf{u}_k$ inductively by $\mathbf{u}_1=\mathbf{v}_1$, and for $2\leq j\leq k$, 
\[\mathbf{u}_j=\mathbf{v}_j-\sum_{i=1}^{j-1}\frac{\langle\mathbf{v}_j, \mathbf{u}_i\rangle}{\langle\mathbf{u}_i, \mathbf{u}_i\rangle }\mathbf{u}_i.\]
Then $S'=\left\{\mathbf{u}_1, \ldots, \mathbf{u}_k\right\}$ is a nonzero set of orthogonal vectors. Moreover, for each $1\leq j\leq k$, the set
$\{\mathbf{u}_i\,|\,1\leq i\leq j\}$ spans the same subspace as the set $\{\mathbf{v}_i\,|\,1\leq i\leq j\}$.
\end{theorem}
\begin{myproof}{Sketch of Proof}
For $1\leq j\leq k$, let $W_j$ be the subspace spanned by the set $\{\mathbf{v}_i\,|\,1\leq i\leq j\}$. The vectors $\mathbf{u}_1, \ldots, \mathbf{u}_k$  are constructed by letting $\mathbf{u}_1=\mathbf{v}_1$, and for $2\leq j\leq k$,
\[\mathbf{u}_j=\mathbf{v}_j-\text{proj}_{W_{j-1}}\mathbf{v}_j.\]
Since $\{\mathbf{v}_1, \ldots, \mathbf{v}_j\}$ is a linearly independent set, $\mathbf{u}_j\neq\mathbf{0}$. Using induction, one can show that $\text{span}\,\{\mathbf{u}_1, \ldots, \mathbf{u}_{j}\}=\text{span}\,\{\mathbf{v}_1, \ldots, \mathbf{v}_{j}\}$. By projection theorem, $\mathbf{u}_j$ is orthogonal to $W_{j-1}$. Hence, it is orthogonal to $\mathbf{u}_1, \ldots, \mathbf{u}_{j-1}$. This proves the theorem.
\end{myproof}

A mapping between two vector spaces that respect the linear structures is called a linear transformation.
\begin{definition}{Linear Transformation}
Let $V$ and $W$ be real vector spaces. A mapping $T:V\to W$ is called a linear transformation provided that for any $\mathbf{v}_1, \ldots,\mathbf{v}_k$ in $V$, for any real numbers $c_1, \ldots, c_k$,
\[T\left(c_1\mathbf{v}_1+\cdots+c_k\mathbf{v}_k\right)=c_1 T(\mathbf{v}_1)+\cdots+c_k T(\mathbf{v}_k).\]
\end{definition}
Linear transformations   play important roles in multivariable analysis. In the following, we first define a special class of linear transformations associated to special projections.

For $1\leq i\leq n$, let $\mathbb{L}_i$ be the subspace of $\mathbb{R}^n$ spanned by the unit vector $\mathbf{e}_i$. For the point  $\mathbf{x}=(x_1, \ldots, x_n)$, 
\[\text{proj}_{\mathbb{L}_i}\mathbf{x}=x_i\mathbf{e}_i.\]The number $x_i$ is the $i^{\text{th}}$-component of $\mathbf{x}$. It will play important roles later. The mapping from $\mathbf{x}$ to $x_i$ is a function from $\mathbb{R}^n$ to $\mathbb{R}$.
\begin{definition}{Projection Functions}
For $1\leq i\leq n$, the $i^{\text{th}}$-projection function on $\mathbb{R}^n$ is the function $\pi_i:\mathbb{R}^n\to\mathbb{R}$ defined by
\[\pi_i(\mathbf{x})=\pi_i(x_1, \ldots, x_n)=x_i.\]
\end{definition}

 \begin{figure}[ht]
\centering
\includegraphics[scale=0.2]{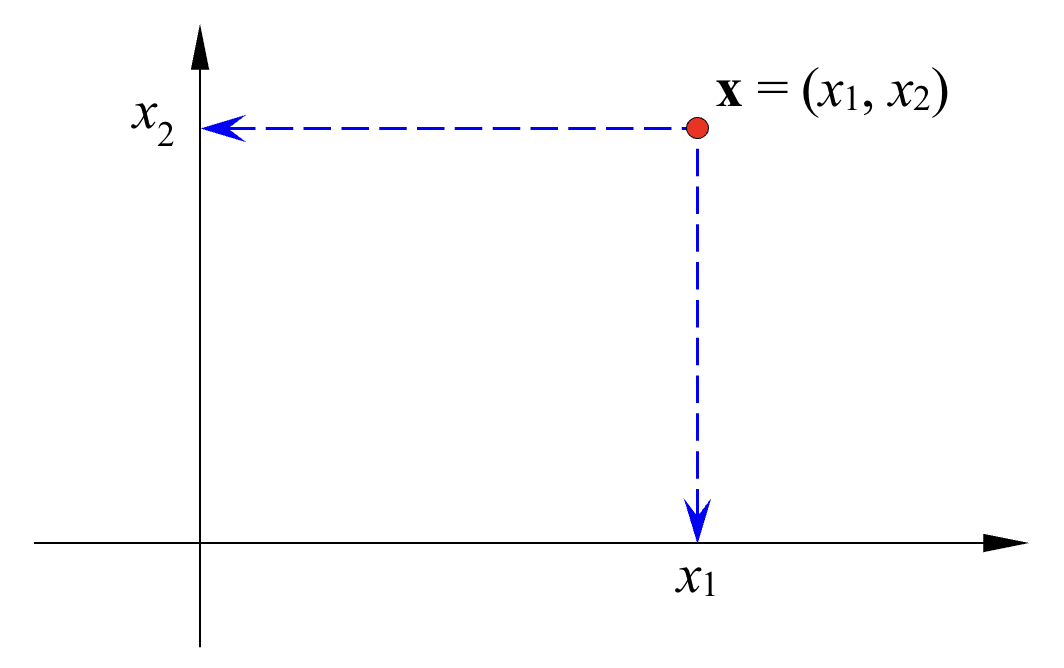}

\caption{The projection functions.}\label{figure47}
\end{figure}
The following is obvious.
\begin{proposition}{}
For $1\leq i\leq n$, the $i^{\text{th}}$-projection function on $\mathbb{R}^n$  is a linear transformation. Namely, for any $\mathbf{x}_1, \ldots, \mathbf{x}_k$ in $\mathbb{R}^n$, and any real numbers $c_1, \ldots, c_k$,
\[\pi_i \left(c_1\mathbf{x}_1+\cdots+c_k\mathbf{x}_k\right)=c_1 \pi_i(\mathbf{x}_1)+\cdots+c_k \pi_i(\mathbf{x}_k).\]
\end{proposition}

The following is a useful inequality.
\begin{proposition}{}
Let $\mathbf{x}$ be a vector in $\mathbb{R}^n$. Then 
\[\left|\pi_i(\mathbf{x})\right|\leq \Vert\mathbf{x}\Vert.\]
\end{proposition}

At the end of this section, let us introduce the concept of hyperplanes.
\begin{definition}{Hyperplanes}
In $\mb{R}^n$, a hyperplane is a translate of a subspace of dimension $n-1$. In other words, $\mb{H}$ is a hyperplane if there is a point $\mf{x}_0$ in $\mb{R}^n$, and $n-1$ linearly independent vectors $\mf{v}_1$, $\mf{v}_2$, $\ldots$, $\mf{v}_{n-1}$ such that $\mb{H}$ contains all points $\mf{x}$ of the form
\[\mf{x}=\mf{x}_0+t_1\mf{v}_1+\cdots+t_{n-1}\mf{v}_{n-1},\hspace{1cm}(t-1, \ldots, t_{n-1})\in\mb{R}^{n-1}.\]
\end{definition}
A hyperplane in $\mb{R}^1$ is a point. A hyperplane in $\mb{R}^2$ is a line. A hyperplane in $\mb{R}^3$ is a plane.

\begin{definition}{Normal Vectors}
Let $\mf{v}_1$, $\mf{v}_2$, $\ldots$, $\mf{v}_{n-1}$ be linearly independent vectors  in $\mb{R}^n$, and let $\mb{H}$ be the hyperplane
\[\mb{H}=\left\{\mf{x}_0+t_1\mf{v}_1+\cdots+t_{n-1}\mf{v}_{n-1}\,|\,(t_1, \ldots, t_{n-1})\in\mb{R}^{n-1}\right\}.\]
A nonzero vector $\mf{n}$ that is orthogonal to all the vectors $\mf{v}_1, \ldots, \mf{v}_{n-1}$ is called a  normal vector  to the hyperplane. If $\mf{x}_1$ and $\mf{x}_2$ are two points on $\mb{H}$, then $\mf{n}$ is orthogonal to the vector $\mf{v}=\mf{x}_2-\mf{x}_1$. Any two normal vectors of a hyperplane are scalar multiples of each other.
\end{definition}

\begin{proposition}{}
If $\mb{H}$ is a hyperplane with normal vector $\mf{n}=(a_1, a_2,\ldots, a_n)$, and $\mf{x}_0=(u_1, u_2, \ldots, u_n)$ is a point on $\mb{H}$, then the equation of $\mb{H}$ is given by
 \[a_1(x_1-u_1)+a_2(x_2-u_2)+\cdots+a_n(x_n-u_n)=\mf{n}\,\cdot\,(\mf{x}-\mf{x}_0)=0.\]
 Conversely, any equation of the form
 \[a_1x_1+a_2x_2+\cdots+a_nx_n=b\]
 is the equation of a hyperplane with normal vector $\mf{n}=(a_1, a_2,\ldots, a_n)$.
\end{proposition}
\begin{example}{}
Given  $1\leq i\leq n$, the equation $x_i=c$ is a hyperplane with normal vector $\mf{e}_i$. It is a hyperplane parallel to the coordinate plane $x_i=0$, and perpendicular to the $x_i$-axis.
\end{example}
\vp

\noindent
{\bf \large Exercises  \thesection}
\setcounter{myquestion}{1}

\begin{question}{\themyquestion}
  Let $V$ be an inner product space. If $\mathbf{u}$ and $\mathbf{v}$ are vectors in $V$, show that
 \[\bigl|\Vert\mathbf{u}\Vert -\Vert\mathbf{v}\Vert\bigr|\leq \Vert \mathbf{u}-\mathbf{v}\Vert.\]
\end{question}

\atc
\begin{question}{\themyquestion}
Let $V$ be an inner product space. If $\mathbf{u}$ and $\mathbf{v}$ are orthogonal vectors in $V$, show that
\[\Vert\mathbf{u}+\mathbf{v}\Vert^2=\Vert\mathbf{u}\Vert^2+\Vert\mathbf{v}\Vert^2.\]
\end{question}
 \atc
\begin{question}{\themyquestion}
Let $V$ be an inner product space, and let $\mathbf{u}$ and $\mathbf{v}$ be vectors in $V$. Show that
\[\langle\mathbf{u}, \mathbf{v}\rangle =\frac{\Vert\mathbf{u}+\mathbf{v}\Vert^2-\Vert\mathbf{u}-\mathbf{v}\Vert^2}{4}.\]
\end{question}

\atc
\begin{question}{\themyquestion}
Let $V$ be an inner product space, and let $\{\mathbf{u}_1, \ldots, \mathbf{u}_k\}$ be an orthonormal set of vectors in $V$. For any real numbers $\alpha_1, \ldots, \alpha_k$, show that
\[\Vert\alpha_1\mathbf{u}_1+\cdots+\alpha_k\mathbf{u}_k\Vert^2=\alpha_1^2+\cdots+\alpha_k^2.\]

\end{question}
 \atc
\begin{question}{\themyquestion}
Let $x_1, x_2, \ldots, x_n$ be real numbers. Show that
\begin{enumerate}[(a)]
\item $\di \sqrt{x_1^2+x_2^2\cdots+x_n^2}\leq |x_1|+|x_2|+\cdots+|x_n|$;
\item $\di |x_1+x_2+\cdots+x_n|\leq \sqrt{n}\sqrt{x_1^2+x_2^2\cdots+x_n^2}$.
\end{enumerate}
\end{question}

\section[Convergence of Sequences in $\mathbb{R}^n$]{Convergence of Sequences in $\pmb{\mathbb{R}^n}$} 

A point in the Euclidean space $\mathbb{R}^n$ is denoted by $\mathbf{x}=(x_1, x_2, \ldots, x_n)$. When $n=1$, we just denote it by $x$. When $n=2$ and $n=3$, it is customary to denote a point in $\mathbb{R}^2$ and $\mathbb{R}^3$ by $(x,y)$ and $(x,y,z)$ respectively.

The Euclidean inner product between the vectors $\mf{x}=(x_1, x_2, \ldots, x_n)$ and $\mf{y}=(y_1, y_2, \ldots, y_n)$ is
\[\langle\mf{x},\mf{y}\rangle=\mf{x}\,\cdot\,\mf{y}=\sum_{i=1}^n x_iy_i.\] The norm of $\mf{x}$ is
\[\Vert\mf{x}\Vert=\sqrt{\langle\mf{x},\mf{x}\rangle}=\sqrt{\sum_{i=1}^nx_i^2},\]while the distance between $\mf{x}$ and $\mf{y}$ is
\[d(\mf{x},\mf{y})=\Vert\mf{x}-\mf{y}\Vert=\sqrt{\sum_{i=1}^n(x_i-y_i)^2}.\]  

A sequence in $\mathbb{R}^n$ is a function $f:\mathbb{Z}^+\to\mathbb{R}^n$. For $k\in \mathbb{Z}^+$, let $\mathbf{a}_k=f(k)$. Then we can also denote the sequence by $\di\{\mathbf{a}_k\}_{k=1}^{\infty}$, or simply as $\{\mathbf{a}_k\}$. 

\begin{example}[label=230713_1]{}
The sequence $\di \left\{\left(\frac{k}{k+1}, \frac{2k+3}{k}\right)\right\}$ is a sequence in $\mathbb{R}^2$ with 
\[\mathbf{a}_k=\left(\frac{k}{k+1}, \frac{2k+3}{k}\right).\]
\end{example}

In volume I, we have seen that a sequence of real numbers $\{a_k\}_{k=1}^{\infty}$ is said to converge to a real number $a$ provided that for any $\varepsilon>0$, there is a positive integer $K$ such that 
\[|a_k-a|<\varepsilon\hspace{1cm}\text{for all} \;k\geq K.\] Notice that $|a_k-a|$ is the distance between $a_k$ and $a$. To define the convergence of a sequence in $\mathbb{R}^n$, we use the Euclidean distance.

\begin{definition}{Convergence of   Sequences}
A sequence $\{\mathbf{a}_k\}$ in $\mathbb{R}^n$ is said to converge to the point $\mathbf{a}$ in $\mathbb{R}^n$ provided that for any $\varepsilon>0$, there is a positive integer $K$ so that for all $k\geq K$,
\[\Vert\mathbf{a}_k-\mathbf{a}\Vert=d(\mathbf{a}_k,\mathbf{a})<\varepsilon.\]If $\{\mathbf{a}_k\}$ is a sequence   that converges to a point $\mathbf{a}$, we say that the sequence  $\{\mathbf{a}_k\}$ is convergent. A sequence that does not converge to any point in $\mathbb{R}^n$ is said to be divergent.
\end{definition}

\begin{figure}[ht]
\centering
\includegraphics[scale=0.2]{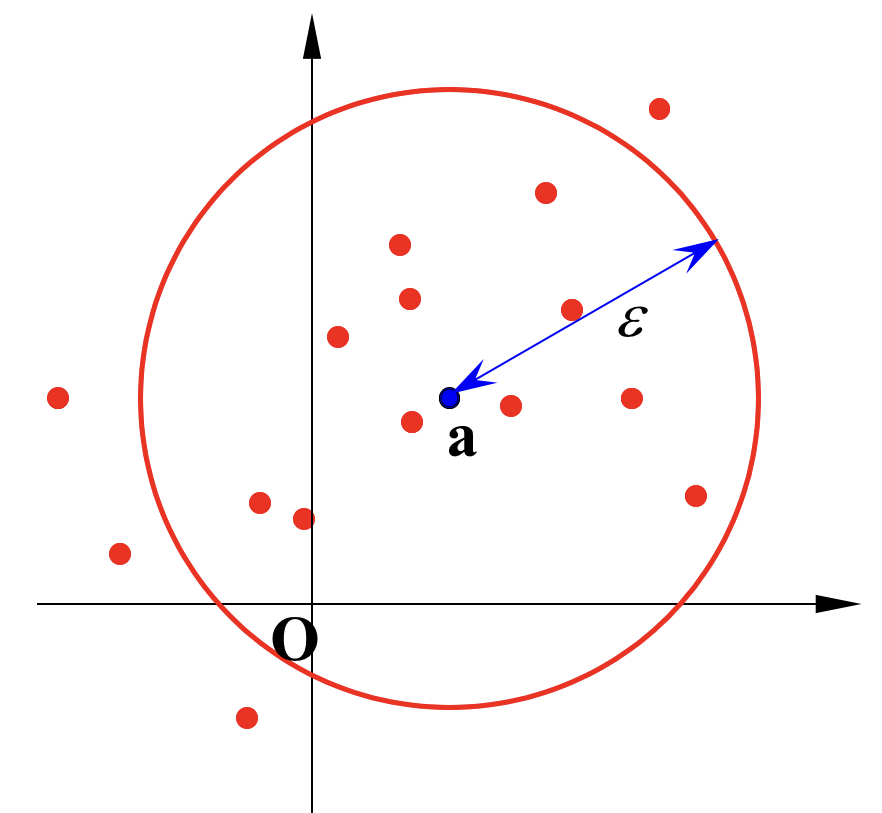}

\caption{The convergence of a sequence.}\label{figure4}
\end{figure}

As in the $n=1$ case, we have the following.
\begin{proposition}{}A sequence in $\mathbb{R}^n$ cannot converge to two different points.  \end{proposition}

\begin{definition}{Limit of a Sequence}
If $\{\mathbf{a}_k\}$ is a sequence in $\mathbb{R}^n$ that converges to the point $\mathbf{a}$, we call   $\mathbf{a}$ the limit of the sequence. This can   be expressed as
\[\lim_{k\to\infty}\mathbf{a}_k=\mathbf{a}.\] 
\end{definition}

The following is easy to establish.
\begin{proposition}[label=230720_4]{}
Let $\{\mathbf{a}_k\}$ be a sequence in $\mathbb{R}^n$. Then $\{\mathbf{a}_k\}$ converges to $\mf{a}$ if and only if
\[\lim_{k\to\infty}\Vert\mathbf{a}_k-\mathbf{a}\Vert=0.\]
\end{proposition}
\begin{myproof}{Proof}
 By definition, the sequence $\{\mathbf{a}_k\}$ is convergent if and only if for any $\varepsilon>0$, there is a positive integer $K$ so that for all $k\geq K$,
$\di \Vert\mathbf{a}_k-\mathbf{a}\Vert <\varepsilon$. This is the definition  of $\di \lim_{k\to\infty}\Vert\mathbf{a}_k-\mathbf{a}\Vert=0$.
\end{myproof}

As in the $n=1$ case,  $\{\mathbf{a}_{k_j}\}_{j=1}^{\infty}$ is a subsequence of $\{\mathbf{a}_k\} $ if $k_1, k_2, k_3, \ldots$ is a strictly increasing subsequence of positive integers.

\begin{corollary}{}
 If   $\{\mathbf{a}_k\}$ is a sequence in $\mathbb{R}^n$ that converges to the point $\mathbf{a}$, then any subsequence of $\{\mathbf{a}_k\}$ also converges to $\mathbf{a}$. 
\end{corollary}

\begin{example}[label=230713_2]{}
Let us investigate the convergence of the sequence     $\{\mathbf{a}_k\}$  in $\mathbb{R}^2$ with 
\[\mathbf{a}_k=\left(\frac{k}{k+1}, \frac{2k+3}{k}\right)\]that is defined in Example \ref{230713_1}. Notice that
\[\lim_{k\to\infty}\pi_1(\mathbf{a}_k)=\lim_{k\to\infty}\frac{k}{k+1}=1,\]
\[\lim_{k\to\infty}\pi_2(\mathbf{a}_k)=\lim_{k\to\infty}\frac{2k+3}{k}=2.\] 
It is natural for us to speculate that the sequence $\{\mathbf{a}_k\}$ converges to the point $\mathbf{a}=(1,2)$. \be
For $k\in\mathbb{Z}^+$,
\[\mathbf{a}_k-\mathbf{a}=\left(-\frac{1}{k+1}, \frac{3}{k}\right).\] 
Thus,
\[\Vert\mathbf{a}_k-\mathbf{a}\Vert=\sqrt{\frac{1}{(k+1)^2}+\frac{9}{k^2}}.\]

By squeeze theorem,
\[\lim_{k\to\infty} \Vert\mathbf{a}_k-\mathbf{a}\Vert=0.\]
This proves that  the sequence $\{\mathbf{a}_k\}$ indeed converges to the point $\mathbf{a}=(1,2)$.
\end{example2}
In the example above, we guess the limit of the sequence by looking at each components of the sequence. This in fact works for any sequences.   
\begin{theorem}[label=230713_3]{Componentwise Convergence of Sequences}
A sequence $\{\mathbf{a}_k\}$   in $\mathbb{R}^n$ converges to the point $\mathbf{a}$ if and only if for each $1\leq i\leq n$,
the sequence $\{\pi_i(\mathbf{a}_k)\}$ converges to the point $\{\pi_i(\mathbf{a})\}$.
\end{theorem}
\begin{myproof}{Proof}
Given $1\leq i\leq n$,
\[\pi_i(\mathbf{a}_k)-\pi_i(\mathbf{a})=\pi_i(\mathbf{a}_k-\mathbf{a}).\]
Thus, 
\[\left|\pi_i(\mathbf{a}_k)-\pi_i(\mathbf{a})\right|=\left|\pi_i(\mathbf{a}_k-\mathbf{a})\right|\leq \Vert\mathbf{a}_k-\mathbf{a}\Vert.\]
If the sequence $\{\mathbf{a}_k\}$    converges to the point $\mathbf{a}$, then
\[\lim_{k\to\infty}\Vert\mf{a}_k-\mf{a}\Vert=0.\]
By squeeze theorem,
\[\lim_{k\to\infty}\left|\pi_i(\mathbf{a}_k)-\pi_i(\mathbf{a})\right|=0.\]
This proves that the sequence $\{\pi_i(\mathbf{a}_k)\}$ converges to the point $\{\pi_i(\mathbf{a})\}$.
 \bp
Conversely, assume that for each $1\leq i\leq n$,
the sequence $\{\pi_i(\mathbf{a}_k)\}$ converges to the point $\{\pi_i(\mathbf{a})\}$. 
Then
\[\lim_{k\to\infty}\left|\pi_i(\mathbf{a}_k)-\pi_i(\mathbf{a})\right|=0\hspace{1cm}\text{for}\;1\leq i\leq n.\]
Since
\begin{align*}
\Vert\mathbf{a}_k-\mathbf{a}\Vert \leq\sum_{i=1}^n \left|\pi_i(\mathbf{a}_k-\mathbf{a})\right|,
\end{align*}
 squeeze theorem implies that
\[\lim_{k\to\infty} \Vert\mathbf{a}_k-\mathbf{a}\Vert=0.\]
This proves that the sequence $\{\mathbf{a}_k\}$  converges to the point $\mathbf{a}$.
\end{myproof}
Theorem \ref{230713_3} reduces the investigations of convergence of sequences in $\mathbb{R}^n$ to sequences in $\mathbb{R}$. Let us look at a few examples.

\begin{example}{}
Find the following limit.
\[\lim_{k\rightarrow\infty}\left(\frac{2^k+1}{3^k}, \left(1+\frac{1}{k}\right)^k, \frac{k}{\sqrt{k^2+1}}\right).\]
\end{example}
\begin{solution}{Solution}
We compute the limit componentwise.
\begin{gather*}
\lim_{k\to\infty}\frac{2^k+1}{3^k} =\lim_{k\to\infty}\left[\left(\frac{2}{3}\right)^k+\left(\frac{1}{3}\right)^k\right]=0+0=0,\\
\lim_{k\to\infty}\left(1+\frac{1}{k}\right)^k =e,\\
\lim_{k\to\infty}\frac{k}{\sqrt{k^2+1}} =\lim_{k\to\infty}\frac{k}{k\sqrt{1+\di\frac{1}{k^2}}}=1.
\end{gather*}
\bs Hence,
\[\lim_{k\rightarrow\infty}\left(\frac{2^k+1}{3^k}, \left(1+\frac{1}{k}\right)^k, \frac{k}{\sqrt{k^2+1}}\right)=(0,e,1).\]
\end{solution}

\begin{example}{}
Let $\{\mathbf{a}_k\}$ be the sequence with
\[\mathbf{a}_k=\left((-1)^k, \frac{(-1)^k}{k}\right).\]
Is the sequence convergent? Justify your answer.
\end{example}
\begin{solution}{Solution}
The sequence $\{\pi_1(\mathbf{a}_k)\}$ is the sequence $\{(-1)^k\}$, which is divergent. Hence, the sequence $\{\mathbf{a}_k\}$ is divergent.
\end{solution}

Using the componentwise convergence theorem, it is easy to establish the following.
\begin{proposition}{Linearity}
Let $\{\mathbf{a}_k\}$ and $\{\mathbf{b}_k\}$ be sequences in $\mathbb{R}^n$ that converges to $\mathbf{a}$ and $\mathbf{b}$ respectively. For any real numbers $\alpha$ and $\beta$, the sequence $\{\alpha\mathbf{a}_k+\beta\mathbf{b}_k\}$ converges to $\alpha\mathbf{a}+\beta\mathbf{b}$. Namely,
\[\lim_{k\to\infty}\left(\alpha\mathbf{a}_k+\beta\mathbf{b}_k\right)=\alpha\mathbf{a}+\beta\mathbf{b}.\]
\end{proposition}

\begin{example}[label=230717_1]{}
If $\{\mathbf{a}_k\}$ is a sequence in $\mathbb{R}^n$ that converges to $\mathbf{a}$, show that
\[\lim_{k\to\infty}\Vert\mathbf{a}_k\Vert =\Vert\mathbf{a}\Vert.\]
\end{example}
\begin{solution}{Solution}
Notice that
\[\Vert\mathbf{a}_k\Vert=\sqrt{\pi_1(\mathbf{a}_k)^2+\cdots+\pi_n(\mathbf{a}_k)^2}.\]
For $1\leq i\leq n$, 
\[\lim_{k\to\infty}\pi_i(\mathbf{a}_k)=\pi_i(\mathbf{a}).\]Using limit laws for sequences in $\mathbb{R}$, we have
\[\lim_{k\to\infty}\left(\pi_1(\mathbf{a}_k)^2+\cdots+\pi_n(\mathbf{a}_k)^2\right)=\pi_1(\mathbf{a})^2+\cdots+\pi_n(\mathbf{a})^2.\]
Using the fact that square root function is continuous, we find that
\begin{align*}
\lim_{k\to\infty}\Vert\mathbf{a}_k\Vert & =\lim_{k\to\infty}\sqrt{ \pi_1(\mathbf{a}_k)^2+\cdots+\pi_n(\mathbf{a}_k)^2}\\&=\sqrt{\pi_1(\mathbf{a})^2+\cdots+\pi_n(\mathbf{a})^2}=\Vert\mathbf{a}\Vert.\end{align*}

\end{solution}

There is also a Cauchy criterion for convergence of sequences in $\mathbb{R}^n$.
\begin{definition}{Cauchy Sequences}
A sequence $\{\mathbf{a}_k\}$ in $\mathbb{R}^n$ is a Cauchy sequence if for every $\varepsilon>0$, there is a positive integer $K$ such that for all $l\geq k\geq K$,
\[\Vert\mathbf{a}_l-\mathbf{a}_k\Vert<\varepsilon.\]
\end{definition}

\begin{theorem}{Cauchy Criterion}
A sequence $\{\mathbf{a}_k\}$ in $\mathbb{R}^n$ is convergent if and only if it is a Cauchy sequence.
\end{theorem}
Similar to the $n=1$ case, the Cauchy criterion allows us to determine whether a sequence in $\mathbb{R}^n$ is convergent without having to guess what is the limit first.
\begin{myproof}{Proof}
Assume that the sequence  $\{\mathbf{a}_k\}$  converges to $\mathbf{a}$. Given $\varepsilon>0$, there is a positive integer $K$ such that for all $k\geq K$, $\Vert\mathbf{a}_k-\mathbf{a}\Vert<\varepsilon/2$. Then for all $l\geq k\geq K$,
\[\Vert\mathbf{a}_l-\mathbf{a}_k\Vert \leq \Vert\mathbf{a}_l-\mathbf{a} \Vert +\Vert\mathbf{a}_k-\mathbf{a}\Vert <\varepsilon.\]
This proves that $\{\mathbf{a}_k\}$ is a Cauchy sequence.

Conversely, assume that   $\{\mathbf{a}_k\}$ is a Cauchy sequence.  Given $\varepsilon>0$, there is a positive integer $K$ such that for all $l\geq k\geq K$, 
\[\Vert\mathbf{a}_l-\mathbf{a}_k\Vert<\varepsilon.\]
 
For each $1\leq i\leq n$, 
\[\left|\pi_i(\mathbf{a}_l)-\pi_i(\mathbf{a}_k)\right|=\left|\pi_i
\left(\mathbf{a}_l-\mathbf{a}_k\right)\right|\leq \Vert\mathbf{a}_l-\mathbf{a}_k\Vert.\]

Hence, $\{\pi_i(\mathbf{a}_k)\}$ is a Cauchy sequence in $\mathbb{R}$. Therefore, it is convergent. By componentwise convergence theorem, the sequence $\{\mathbf{a}_k\}$ is convergent.

\end{myproof}
\vp
\noindent
{\bf \large Exercises  \thesection}
\setcounter{myquestion}{1}
\begin{question}{\themyquestion}
Show that a sequence in $\mathbb{R}^n$ cannot converge to two different points.
\end{question}
 \atc
 \begin{question}{\themyquestion}
Find the limit of the sequence $\{\mathbf{a}_k\}$, where
\[\mathbf{a}_k=\left(\frac{2k+1}{k+3}, \frac{\sqrt{2k^2+k}}{k}, \left(1+\frac{2}{k}\right)^k\right).\]
\end{question}

\atc
 \begin{question}{\themyquestion}
Let $\{\mathbf{a}_k\}$ be the sequence with 
\[\mathbf{a}_k=\left(\frac{1+(-1)^{k-1}k}{1+k}, \frac{1}{2^k}\right).\]Determine whether the sequence is convergent.
\end{question}

\atc
 \begin{question}{\themyquestion}
Let $\{\mathbf{a}_k\}$ be the sequence with 
\[\mathbf{a}_k=\left(\frac{ k}{1+k}, \frac{k}{\sqrt{k}+1}\right).\]Determine whether the sequence is convergent.
\end{question}

\atc
 \begin{question}{\themyquestion}
 Let $\{\mathbf{a}_k\}$ and $\{\mathbf{b}_k\}$ be sequences in $\mathbb{R}^n$ that converges to $\mathbf{a}$ and $\mathbf{b}$ respectively. Show that
 \[\lim_{k\to\infty}\langle\mathbf{a}_k, \mathbf{b}_k\rangle =\langle\mathbf{a},\mathbf{b}\rangle.\]Here $\langle\mathbf{x},\mathbf{y}\rangle=\mathbf{x}\,\cdot\,\mathbf{y}$ is the standard inner product on $\mathbb{R}^n$.
 \end{question}
 \atc
 \begin{question}{\themyquestion}
 Suppose that $\{\mathbf{a}_k\}$ is a sequence in $\mathbb{R}^n$ that converges to $\mathbf{a}$, and $\{c_k\}$ is a sequence of real numbers that converges to $c$, show that
 \[\lim_{k\to\infty}c_k\mathbf{a}_k=c\mathbf{a}.\]
 \end{question}
   \atc
 \begin{question}{\themyquestion}
 Suppose that $\{\mathbf{a}_k\}$ is a sequence of nonzero vectors in $\mathbb{R}^n$ that converges to $\mathbf{a}$ and $\mathbf{a}\neq\mathbf{0}$, show that
 \[\lim_{k\to\infty}\frac{\mathbf{a}_k}{\Vert\mathbf{a}_k\Vert}=\frac{\mathbf{a}}{\Vert\mathbf{a}\Vert}.\]
 \end{question}
 \atc
 \begin{question}{\themyquestion}
Let $\{\mathbf{a}_k\}$ and $\{\mathbf{b}_k\}$ be sequences in $\mathbb{R}^n$. If $\{\mathbf{a}_k\}$ is convergent and $\{\mathbf{b}_k\}$ is divergent, show that the sequence $\{\mathbf{a}_k+\mathbf{b}_k\}$ is divergent.
\end{question}
 \atc
 \begin{question}{\themyquestion}
 Suppose that $\{\mathbf{a}_k\}$ is a sequence in $\mathbb{R}^n$ that converges to $\mathbf{a}$. If $r=\Vert\mathbf{a}\Vert\neq 0$, show that there is a positive integer $K$ such that
 \[\Vert\mathbf{a}_k\Vert >\frac{r}{2}\hspace{1cm}\text{for all}\;k\geq K.\]
 \end{question}
 \atc
 \begin{question}{\themyquestion}
Let $\{\mathbf{a}_k\}$ be a sequence in $\mathbb{R}^n$ and let $\mathbf{b}$ be a point in $\mb{R}^n$. Assume that the sequence $\{\mathbf{a}_k\}$ does not converge to $\mf{b}$. Show that there is an $\varepsilon>0$ and a subsequence $\{\mathbf{a}_{k_j}\}$ of $\{\mf{a}_k\}$ such that \[\Vert\mf{a}_{k_j}-\mf{b}\Vert\geq \varepsilon\hspace{1cm}\text{for all}\; j\in\mathbb{Z}^+.\] \end{question}

\section{Open Sets and Closed Sets} 

In volume I, we call an interval of the form $(a, b)$ an {\it open} interval. Given a point $x$ in $\mathbb{R}$, a neighbourhood of $x$ is an open interval $(a,b)$ that contains $x$. Given a subset $S$ of $\mathbb{R}$, we say that $x$ is an interior point of $S$ if there is a neighboirhood of $x$ that is contained in $S$. We say that $S$ is closed in $\mathbb{R}$ provided that if $\{a_k\}$ is a sequence of points in $S$ that converges to $a$, then $a$ is also in $S$. These describe the topology of $\mathbb{R}$. It is relatively simple.

For $n\geq 2$, the topological features of $\mathbb{R}^n$ are much more complicated.

An open interval $(a,b)$ in $\mathbb{R}$ can be described as a set of the form
\[B=\left\{x\in\mathbb{R}\,|\,|x-x_0|<r\right\},\]
where $\di x_0=\frac{a+b}{2}$ and $r=\di\frac{b-a}{2}$.

\begin{figure}[ht]
\centering
\includegraphics[scale=0.2]{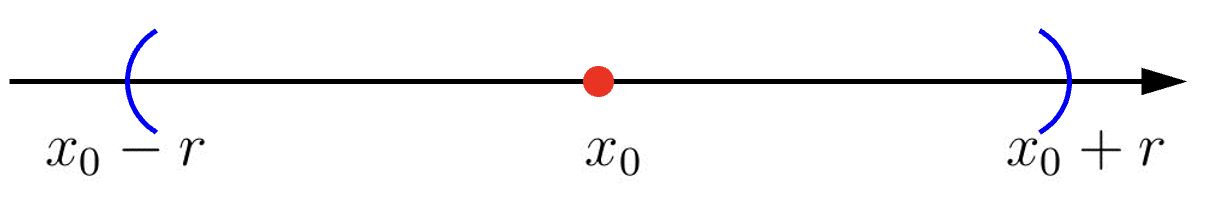}

\caption{An open interval.}\label{figure5}
\end{figure}

Generalizing this, we define open balls in $\mathbb{R}^n$.
\begin{definition}{Open Balls}
Given $\mathbf{x}_0$ in $\mathbb{R}^n$ and $r>0$, an open ball $B(\mathbf{x}_0,r)$ of radius $r$ with center at $\mathbf{x}_0$ is a subset of $\mathbb{R}^n$ of the form
\[B(\mathbf{x}_0,r)=\left\{\mathbf{x}\in\mathbb{R}^n\,|\, \Vert\mathbf{x}-\mathbf{x}_0\Vert<r\right\}.\]
It consists of all points of $\mathbb{R}^n$ whose distance to the center $\mathbf{x}_0$ is less than $r$.
\end{definition}
\begin{figure}[ht]
\centering
\includegraphics[scale=0.2]{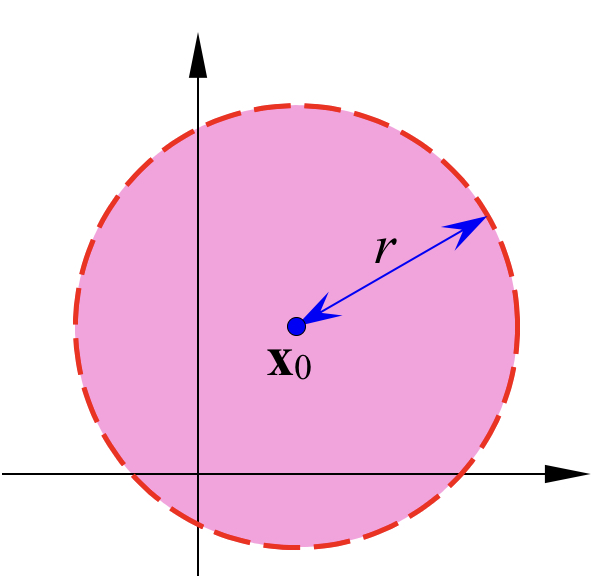}

\caption{An open ball.}\label{figure6}
\end{figure}
Obviously, it $0<r_1\leq r_2$, then $B(\mathbf{x}_0,r_1)\subset B(\mathbf{x}_0,r_2)$.  
The following is a useful lemma for balls with different centers.

\begin{lemma}[label=230714_1]{}
Let $\mathbf{x}_1$ be a point in the open ball $B(\mathbf{x}_0,r)$. Then $\Vert\mathbf{x}_1-\mathbf{x}_0\Vert<r$. If $r_1$ is a positive number satisfying
\[r_1\leq r-\Vert\mathbf{x}_1-\mathbf{x}_0\Vert,\]
then the open ball $B(\mathbf{x}_1,r_1)$ is contained in the open ball $B(\mathbf{x}_0,r)$.
\end{lemma}

\begin{figure}[ht]
\centering
\includegraphics[scale=0.2]{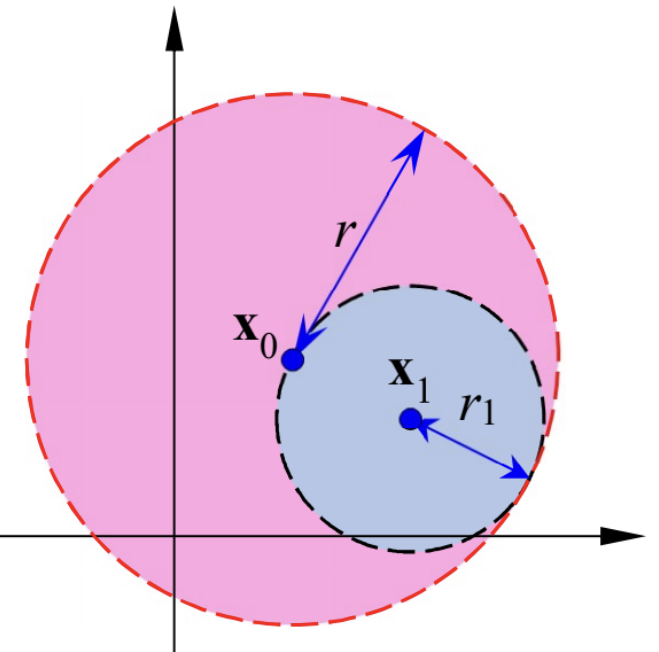}

\caption{An open ball containing another open ball with different center.}\label{figure7}
\end{figure}

\begin{myproof}{Proof}
Let $\mathbf{x}$ be a point in $B(\mathbf{x}_1,r_1)$. Then
\[\Vert\mathbf{x}-\mathbf{x}_1\Vert<r_1\leq r-\Vert\mathbf{x}_1-\mathbf{x}_0\Vert.\]\bp
By triangle inequality,
\[\Vert\mathbf{x}-\mathbf{x}_0\Vert\leq \Vert\mathbf{x}-\mathbf{x}_1\Vert+\Vert\mathbf{x}_1-\mathbf{x}_0\Vert<r.\]
Therefore, $\mathbf{x}$ is a point in $B(\mathbf{x}_0,r)$. This proves the assertion.
\end{myproof}

Now we define open sets   in $\mathbb{R}^n$.
\begin{definition}{Open Sets}
Let $S$ be a subset of $\mathbb{R}^n$. We say that $S$ is an open set if for each $\mathbf{x}\in S$, there is a ball $B(\mathbf{x}, r)$ centered at $\mathbf{x}$ that is contained in $S$. 
\end{definition}
The following example justifies that an open interval of the form $(a,b)$ is   an open set.
\begin{example}[label=230715_9]
{}
Let $S$ to be the open interval $S=(a,b)$ in $\mathbb{R}$. If $x\in S$, then $a<x<b$. Hence, $x-a$ and $b-x$ are positive. Let $r=\min\{x-a, b-x\}$. Then $r>0$,   $r\leq x-a$ and $r\leq b-x$. These imply that $a\leq x-r<x+r\leq b$. Hence, $B(x,r)=(x-r, x+r)\subset (a,b)=S$. This shows that the interval $(a,b)$ is   an open set.
\end{example}

\begin{figure}[ht]
\centering
\includegraphics[scale=0.2]{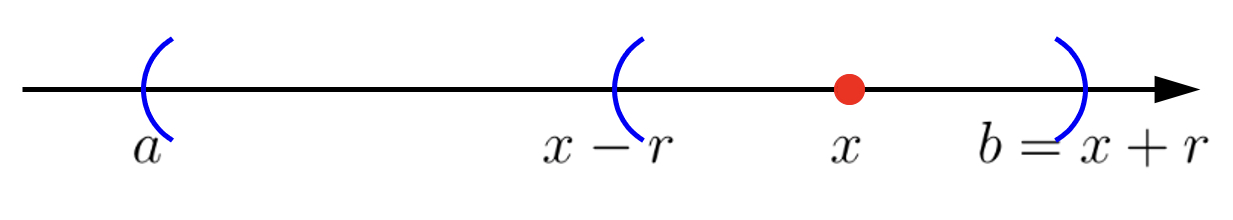}

\caption{The interval $(a,b)$ is an open set.}\label{figure13}
\end{figure}
The following example justifies that an open ball is indeed an open set.
\begin{example}{}
Let $S=B(\mathbf{x}_0, r)$ be the open ball with center at $\mathbf{x}_0$ and radius $r>0$ in $\mathbb{R}^n$. Show that $S$ is an open set.
\end{example}
\begin{solution}{Solution}
Given $\mathbf{x}\in S$, $d=\Vert\mathbf{x}-\mathbf{x}_0\Vert<r$. Let $r_1=r-d$. Then $r_1>0$. Lemma \ref{230714_1} implies that the ball $B(\mathbf{x}, r_1)$ is inside $S$. Hence, $S$ is an open set.

\end{solution}

\begin{example}{}
As   subsets of $\mathbb{R}^n$, $\emptyset$ and $\mathbb{R}^n$ are open sets.
\end{example}
\begin{example}{}
A one-point set $S=\{\mathbf{a}\}$ in $\mathbb{R}^n$  cannot be open, for there is no $r>0$ such that $B(\mathbf{a}, r)$ in contained in $S$.
\end{example}

Let us look at some other examples of open sets.
\begin{definition}{Open Rectangles}
 A set of the form 
\[U=\prod_{i=1}^n(a_i,b_i)=(a_1, b_1)\times\cdots\times (a_n, b_n) \]in $\mb{R}^n$, which is a cartesian product of  open  bounded intervals, in called an open rectangle.
\end{definition}

\begin{figure}[ht]
\centering
\includegraphics[scale=0.2]{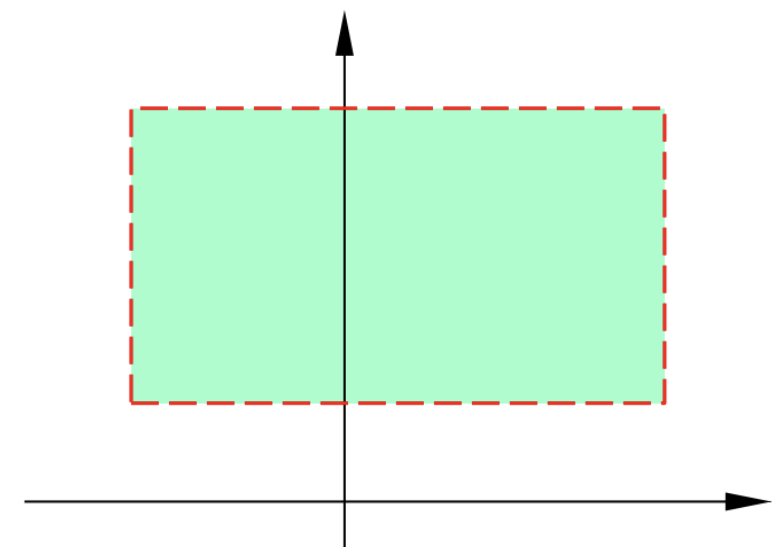}

\caption{A rectangle in $\mathbb{R}^2$.}\label{figure9}
\end{figure}
\begin{example}{}
Let $\di U=\prod_{i=1}^n(a_i,b_i)$ be an open rectangle in $\mathbb{R}^n$. Show that $U$ is an open set.
\end{example}
\begin{solution}{Solution}
Let $\mathbf{x}=(x_1, \ldots, x_n)$ be a point in $U$. Then for $1\leq i\leq n$, \[r_i=\min\{x_i-a_i, b_i-x_i\}>0\] and  \[(x_i-r_i, x_i+r_i)\subset (a_i, b_i).\] Let $r=\min\{r_1, \ldots, r_n\}$. Then $r>0$. We claim that $B(\mathbf{x}, r)$ is contained in $U$.  

If $\mathbf{y}\in B(\mathbf{x}, r)$, then $\Vert\mathbf{y}-\mathbf{x}\Vert<r$. This implies that
\[|y_i-x_i|\leq\Vert\mathbf{y}-\mathbf{x}\Vert<r\leq r_i\hspace{1cm}\text{for all}\;1\leq i\leq n.\]Hence,
\[y_i\in (x_i-r_i, x_i+r_i)\subset (a_i, b_i)\hspace{1cm}\text{for all}\;1\leq i\leq n.\] This proves  that
$\mathbf{y}\in U$, and thus, completes the proof that $B(\mathbf{x}, r)$ is contained in $U$. Therefore, $U$ is an open  set.
 
\end{solution}
\begin{figure}[ht]
\centering
\includegraphics[scale=0.2]{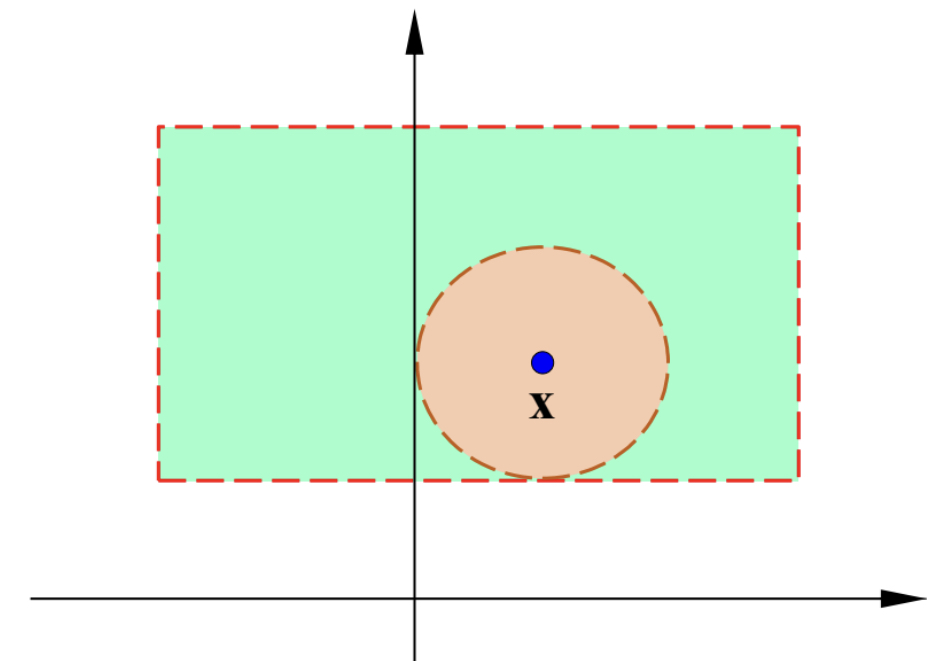}

\caption{An open rectangle is an open set.}\label{figure14}
\end{figure}
Next, we define closed sets. The definition is a straightforward generalization of the $n=1$ case.
\begin{definition}{Closed Sets}

Let $S$ be a subset of $\mathbb{R}^n$. We say that $S$ is closed in $\mathbb{R}^n$ provided that if $\{\mathbf{a}_k\}$
is a sequence of points in $S$ that converges to the point  $\mathbf{a}$, the point $\mathbf{a}$ is also in $S$.
\end{definition}
\begin{example}{}
As   subsets of $\mathbb{R}^n$, $\emptyset$ and $\mathbb{R}^n$ are closed sets. Since  $\emptyset$ and $\mathbb{R}^n$  are also open, a subset $S$ of $\mathbb{R}^n$ can be both open and closed.
\end{example}
\begin{example}{}
Let $S=\{\mathbf{a}\}$ be a one-point set in $\mathbb{R}^n$. A sequence $\{\mathbf{a}_k\}$ in $S$ is just the constant sequence where $\mathbf{a}_k=\mathbf{a}$ for all $k\in\mathbb{Z}^+$.  Hence, it converges to $\mathbf{a}$ which is in $S$. Thus, a one-point set $S$ is a closed set.
\end{example}

In volume I, we have proved the following. 
 
\begin{proposition}{}
Let $I$  be intervals of the form $(-\infty, a]$, $[a, \infty)$ or $[a, b]$. Then $I$ is a closed subset of $\mathbb{R}$.
\end{proposition}

\begin{definition}{Closed Rectangles}
A set of the form 
\[R=\prod_{i=1}^n[a_i,b_i]=[a_1, b_1]\times\cdots\times [a_n, b_n] \]in $\mb{R}^n$, which is a cartesian product of closed and bounded intervals, is called a closed rectangle.
\end{definition}
The following justifies that a closed rectangle is indeed a closed set.
\begin{example}{}
Let  
\[R=\prod_{i=1}^n[a_i,b_i]=[a_1, b_1]\times\cdots\times [a_n, b_n]  \] be a closed rectangle in $\mathbb{R}^n$. Show that $R$ is a closed set.
\end{example}
\begin{solution}{Solution}
Let $\{\mathbf{a}_k\}$ be a sequence in $R$ that converges to a point $\mathbf{a}$. For each $1\leq i\leq n$, $\{\pi_i(\mathbf{a}_k)\}$ is a sequence in $[a_i, b_i]$ that converges to $\pi_i(\mathbf{a})$. Since $[a_i, b_i]$ is a closed set in $\mathbb{R}$, $\pi_i(\mathbf{a})\in [a_i, b_i]$. Hence, $\mathbf{a}$ is in $R$. This proves that $R$ is a closed set.
 
\end{solution}

 It is not true that a set that is not open is closed. 
\begin{example}[label=230714_6]{}
Show that an interval of the form $I=(a, b]$ in $\mathbb{R}$ is neither open nor closed.
\end{example}
\begin{solution}{Solution}
If $I$ is open, since $b$ is in $I$, there is an $r>0$ such that $(b-r, b+r)=B(b, r)\subset I$. But then $b+r/2 $ is a point in $(b-r, b+r)$ but not in $I=(a,b]$, which gives a contradiction.  Hence, $I$ is not open.

For $k\in\mathbb{Z}^+$, let \[a_k=\di a+\frac{b-a}{k}.\] Then $\{a_k\}$ is a sequence in $I$ that converges to $a$, but $a$ is not in $I$. Hence, $I$ is not closed.
\end{solution}
Thus, we have seen that a subset $S$ of $\mathbb{R}^n$ can be both open and closed, and it can also be neither open nor closed.

Let us look at some other examples of closed sets.

 \begin{definition}{Closed Balls}
Given $\mathbf{x}_0$ in $\mathbb{R}^n$ and $r>0$, a closed ball  of radius $r$ with center at $\mathbf{x}_0$ is a subset of $\mathbb{R}^n$ of the form
\[CB(\mathbf{x}_0,r)=\left\{\mathbf{x}\in\mathbb{R}^n\,|\, \Vert\mathbf{x}-\mathbf{x}_0\Vert\leq r\right\}.\]
It consists of all points of $\mathbb{R}^n$ whose distance to the center $\mathbf{x}_0$ is less than or equal to $r$.
\end{definition}
The following justifies that a closed ball is indeed a closed set.
\begin{example}{}
Given $\mathbf{x}_0\in \mathbb{R}^n$ and $r>0$, show that the closed ball \[CB(\mathbf{x}_0,r)=\left\{\mathbf{x}\in\mathbb{R}^n\,|\, \Vert\mathbf{x}-\mathbf{x}_0\Vert\leq r\right\} \] is a closed set.
\end{example}
\begin{solution}{Solution}
Let $\{\mathbf{a}_k\}$ be a sequence in  $CB(\mathbf{x}_0,r)$ that converges to the point $\mathbf{a}$. Then
\[\lim_{k\to\infty}\Vert\mathbf{a}_k-\mathbf{a}\Vert=0.\] 
For each $k\in\mathbb{Z}^+$, $\Vert\mathbf{a}_k-\mathbf{x}_0\Vert\leq r$. By triangle inequality,
\[\Vert\mathbf{a} -\mathbf{x}_0\Vert\leq \Vert\mathbf{a}_k-\mathbf{x}_0\Vert+\Vert\mathbf{a}_k-\mathbf{a} \Vert\leq r+\Vert\mathbf{a}_k-\mathbf{a} \Vert.\]
Taking the $k\to\infty$ limit, we find that
\[ \Vert\mathbf{a} -\mathbf{x}_0\Vert\leq r.\]
Hence, $\mathbf{a}$ is in $CB(\mathbf{x}_0,r)$. This proves that $CB(\mathbf{x}_0,r)$ is a closed set.
\end{solution}

The following theorem gives the  relation between open and closed sets.
\begin{theorem}[label=230714_2]{}
Let $S$ be a subset of $\mathbb{R}^n$ and let $A=\mathbb{R}^n\setminus S$ be its complement in $\mathbb{R}^n$.  Then $S$ is open if and only if $A$ is closed.
\end{theorem}
\begin{myproof}{Proof}
 Assume that $S$ is open. Let $\{\mathbf{a}_k\}$ be a sequence in $A$ that converges to the point $\mathbf{a}$. We want to show that $\mathbf{a}$ is in $A$. Assume to the contrary that $\mathbf{a}$ is not in $A$. Then $\mathbf{a}$ is in $S$. Since $S$ is open, there is an $r>0$ such that $B(\mathbf{a}, r)$ is contained in $S$. Since the sequence $\{\mathbf{a}_k\}$ converges to $\mathbf{a}$, there is a positive integer $K$ such that for all $k\geq K$,
 $\di \Vert\mathbf{a}_k-\mathbf{a}\Vert<r$.
 But then this implies that $\mathbf{a}_K\in B(\mathbf{a}, r)\subset S$. This contradicts to $\mathbf{a}_K$ is in $A=\mathbb{R}^n\setminus S$. Hence, we must have $\mathbf{a}$ is in $A$, which proves that $A$ is closed.
 
 Conversely, assume that $A$ is closed. We want to show that $S$ is open. Assume to the contrary that $S$ is not open. Then there is a point $\mathbf{a}$ in $S$ such that for every $r>0$, $B(\mathbf{a}, r)$ is not contained in $S$. For every $k\in \mathbb{Z}^+$, since $B(\mathbf{a}, 1/k)$ is not contained in $S$, there is a point $\mathbf{a}_k$ in $B(\mathbf{a}, 1/k)$ such that $\mathbf{a}_k$ is not in $S$. Thus, $\{\mathbf{a}_k\}$ is a sequence in $A$ and 
 \[\Vert\mathbf{a}_k-\mathbf{a}\Vert<\frac{1}{k}.\]
 This shows that $\{\mathbf{a}_k\}$ converges to $\mathbf{a}$. Since $A$ is closed, $\mathbf{a}$ is in $A$, which contradicts to $\mathbf{a}$ is in $S$. Thus, we must have $S$ is open.
\end{myproof}
\begin{figure}[ht]
\centering
\includegraphics[scale=0.18]{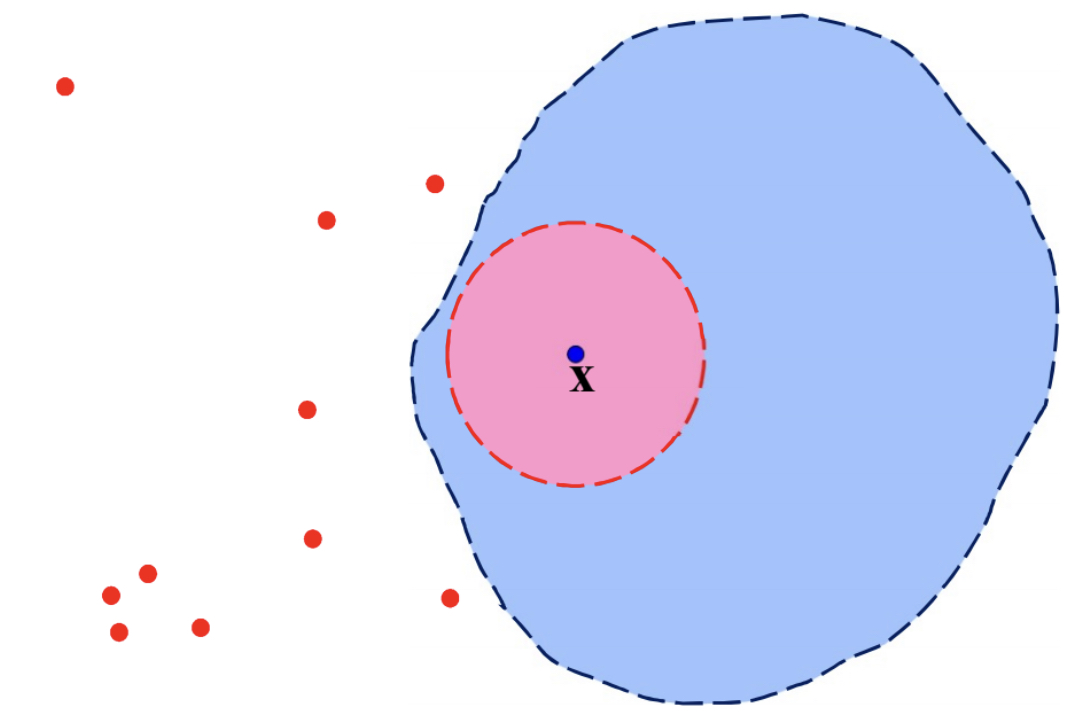}

\caption{A sequence outside an open set cannot converge to a point in the open set.}\label{figure15}
\end{figure}
Next, we consider unions and intersections of sets.

\begin{theorem}[label=230714_3]{}
\begin{enumerate}[1.]
\item Arbitrary union of open sets is open. Namely, if
  $\{U_{\alpha}\,|\, \alpha \in J\}$ is a collection of open sets in $\mathbb{R}^n$, then their union $\di U= \bigcup_{\alpha\in J}U_{\alpha}$ is also an open set.
\item Finite intersections of open sets is open. Namely, if $V_1, \ldots, V_k$ are open sets in $\mathbb{R}^n$, then their intersection $\di V=\bigcap_{i=1}^k V_i$ is also an open set.
\end{enumerate}
\end{theorem}
\begin{myproof}{Proof}
To prove the first statement, let $\mathbf{x}$ be a point in $\di U= \bigcup_{\alpha\in J}U_{\alpha}$. Then there is an $\alpha\in J$ such that $\mathbf{x}$ is in $U_{\alpha}$. Since $U_{\alpha}$ is open, there is an $r>0$ such that $B(\mathbf{x}, r)\subset U_{\alpha}\subset U$. Hence, $U$ is open.

For the second statement, let $\mathbf{x}$ be a point in $V=\di\bigcap_{i=1}^k V_i$. Then for each $1\leq i\leq k$, $\mathbf{x}$ is in the open set $V_i$. Hence, there is an $r_i>0$ such that $B(\mathbf{x}, r_i)\subset V_i$. Let $r=\min\{r_1, \ldots, r_k\}$. Then for $1\leq i\leq k$, $r\leq r_i$ and so $B(\mathbf{x}, r)\subset B(\mathbf{x}, r_i)\subset V_i$. Hence, $B(\mathbf{x}, r)\subset V$. This proves that $V$ is open.
\end{myproof}

As an application of this theorem, let us show that any open interval in $\mathbb{R}$ is indeed an open set.
\begin{proposition}{}
Let $I$  be an interval of the form $(-\infty, a)$, $(a, \infty)$ or $(a, b)$. Then $I$ is an open subset of $\mathbb{R}$.
\end{proposition}
\begin{myproof}{Proof}
We have shown in Example \ref{230715_9} that if $I$ is an interval of the form $(a,b)$, then $I$ is an open subset of $\mathbb{R}$. Now
\[(a, \infty)=\bigcup_{k=1}^{\infty}(a, a+k)\] is a union of open sets. Hence, $(a,\infty)$ is open. In the same way, one can show that an interval of the form $(-\infty, a)$ is open.
\end{myproof}

The next example shows that arbitrary   intersections of open sets is not necessary open.  
\begin{example}{}
For $k\in \mathbb{Z}^+$, let $U_k$ be the open set in $\mathbb{R}$ given by
\[U_k=\left(-\frac{1}{k}, \frac{1}{k}\right).\]Notice that the set
\[U=\bigcap_{k=1}^{\infty}U_k=\{0\}\] is a one-point set. Hence, it is not open in $\mathbb{R}$. 
\end{example}

De Morgan's law in set theory says that if $\{U_{\alpha}\,|\,\alpha\in J\}$ is a collection of sets in $\mathbb{R}^n$, then
\begin{gather*}
\mathbb{R}^n\setminus \bigcup_{\alpha\in J}U_{\alpha}=\bigcap_{\alpha\in J}\left(\mathbb{R}^n\setminus U_{\alpha}\right),\\
\mathbb{R}^n\setminus \bigcap_{\alpha\in J}U_{\alpha}=\bigcup_{\alpha\in J}\left(\mathbb{R}^n\setminus U_{\alpha}\right).
\end{gather*}
Thus, we obtain the counterpart of  Theorem \ref{230714_3} for closed sets.
\begin{theorem}[label=230714_4]{}
\begin{enumerate}[1.]
\item Arbitrary intersection of closed sets is closed. Namely, if
  $\{A_{\alpha}\,|\, \alpha \in J\}$ is a collection of closed sets in $\mathbb{R}^n$, then their intersection $\di A= \bigcap_{\alpha\in J}A_{\alpha}$ is also a closed set.
\item Finite  union of closed sets is closed. Namely, if $C_1, \ldots, C_k$ are closed sets in $\mathbb{R}^n$, then their union $\di C=\bigcup_{i=1}^k C_i$ is also a closed set.
\end{enumerate}
\end{theorem}
\begin{myproof}{Proof}
We prove the first statement. The proof of the second statement is similar. Given that
  $\{A_{\alpha}\,|\, \alpha \in J\}$ is a collection of closed sets in $\mathbb{R}^n$, for each $\alpha\in J$, let $U_{\alpha}=\mathbb{R}^n\setminus A_{\alpha}$. Then $\{U_{\alpha}\,|\, \alpha \in J\}$ is a collection of open sets in $\mathbb{R}^n$. By Theorem \ref{230714_3}, 
the set $\di \bigcup_{\alpha\in J}U_{\alpha}$ is open. By Theorem \ref{230714_2},  $\di \mathbb{R}^n\setminus\bigcup_{\alpha\in J}U_{\alpha}$ is closed. By De Morgan's law,
 \[\mathbb{R}^n\setminus\bigcup_{\alpha\in J}U_{\alpha}=\bigcap_{\alpha\in J}\left(\mathbb{R}^n\setminus U_{\alpha}\right)=\bigcap_{\alpha\in J}A_{\alpha}.\]
 This proves that $\di\bigcap_{\alpha\in J}A_{\alpha}$ is a closed set.
\end{myproof}
The following example says that any finite point set is a closed set.
\begin{example}{}
Let $S=\{\mathbf{x}_1, \ldots,\mathbf{x}_k\}$ be a finite point set in $\mathbb{R}^n$. Then $S=\di\bigcup_{i=1}^k\{\mathbf{x}_i\}$ is a finite union of one-point sets. Since one-point set is closed, $S$ is closed.
\end{example}

\vp
\noindent
{\bf \large Exercises  \thesection}
\setcounter{myquestion}{1}

\begin{question}{\themyquestion}
Let $A$ be the subset of $\mathbb{R}^2$ given by
\[A=\left\{(x,y)\,|\, x>0, y>0\right\}.\]Show that $A$ is an open set.
\end{question}
\atc
\begin{question}{\themyquestion}
Let $A$ be the subset of $\mathbb{R}^2$ given by
\[A=\left\{(x,y)\,|\, x\geq 0, y\geq 0\right\}.\]Show that $A$ is a closed set.
\end{question}

\atc
\begin{question}{\themyquestion}
Let $A$ be the subset of $\mathbb{R}^2$ given by
\[A=\left\{(x,y)\,|\, x>0, y\geq 0\right\}.\]
Is $A$ open?  Is $A$ closed? Justify your answers.
 
\end{question}
 \atc
 
\begin{question}{\themyquestion}
Let $C$ and $U$  be subsets of $\mathbb{R}^n$. Assume that $C$ is closed and $U$ is open, show that $U\setminus C$ is open and $C\setminus U$ is closed.
\end{question}

 \atc
 
\begin{question}{\themyquestion}
Let  $A$  be a subset  of $\mathbb{R}^n$, and let $B=A+\mf{u}$ be the translate of $A$ by the vector $\mf{u}$.
\begin{enumerate}[(a)]
\item Show that $A$ is open if and only if $B$ is open.
\item Show that $A$ is closed if and only if $B$ is closed.
\end{enumerate}
\end{question}

\section{Interior, Exterior, Boundary and Closure}
First, we introduce the interior of a set.
\begin{definition}{Interior}
Let $S$ be a subset of $\mathbb{R}^n$. We say that $\mathbf{x}\in \mathbb{R}^n$ is an interior point of $S$ if  there exists $r>0$ such that $B(\mathbf{x}, r)\subset S$.  The interior of $S$, denoted by $\text{int}\,S$, is defined to be the collection of all the interior points  of $S$. 
\end{definition}

\begin{figure}[ht]
\centering
\includegraphics[scale=0.2]{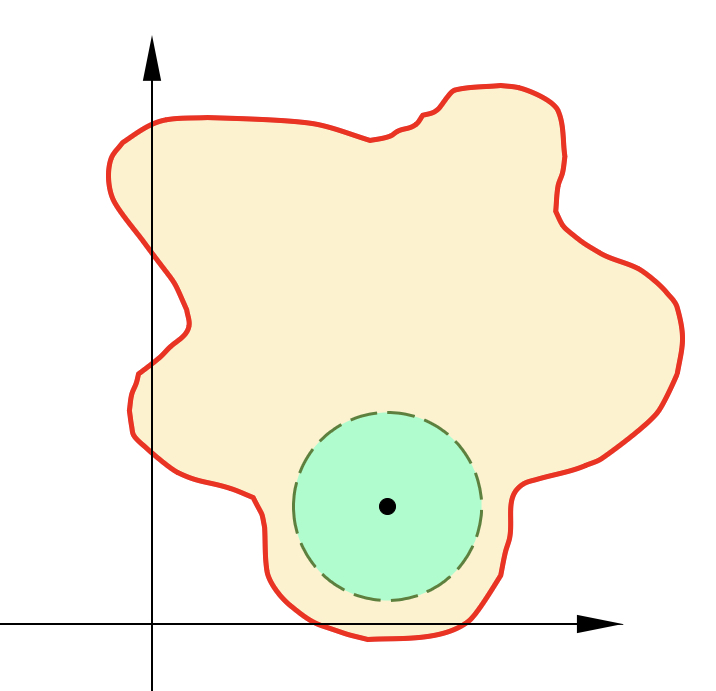}

\caption{The interior point of a set.}\label{figure10}
\end{figure}
The following gives a characterization of the interior of a set.
\begin{theorem}[label=230715_3]{}
Let $S$ be a subset of $\mathbb{R}^n$. Then we have the followings.
\begin{enumerate}[1.]
\item $\text{int}\,S$ is a subset of $S$.
\item $\text{int}\,S$ is an open set.
\item $S$ is an open set if and only if $S=\text{int}\, S$.
\item If $U$ is an open set that is contained in $S$, then $U\subset \text{int}\, S$.

\end{enumerate}These imply that $\text{int}\,S$ is the largest open set that is contained in $S$.
\end{theorem}
\begin{myproof}{Proof}
Let $\mathbf{x}$ be a point  in $\text{int}\,S$. By definition, there exists $r>0$ such that $B(\mathbf{x}, r)\subset S$. 
Since $\mathbf{x}\in B(\mathbf{x}, r)$ and $B(\mathbf{x}, r)\subset S$,  $\mathbf{x}$ is a point in $S$. Since we have shown that every point in $\text{int}\,S$ is in $S$,   $\text{int}\,S$ is a subset of $S$.

If $\mathbf{y}\in B(\mathbf{x}, r)$, Lemma \ref{230714_1} says that there is an $r'>0$ such that $B(\mathbf{y}, r')\subset B(\mathbf{x}, r)\subset S$. Hence, $\mathbf{y}$ is also in $\text{int}\, S$. This proves that $B(\mathbf{x}, r)$  is contained in $\text{int}\, S$. Since we have shown that for any $\mathbf{x}\in \text{int}\, S$, there is an $r>0$ such that $B(\mathbf{x}, r)$  is contained in $\text{int}\, S$, this shows that $\text{int}\, S$ is open.

If $S=\text{int}\, S$, $S$ is open. Conversely, if $S$ is open, for every $\mathbf{x}$ in $S$, there is an $r>0$ such that $B(\mathbf{x},r)\subset S$. Then $\mathbf{x}$ is in $\text{int}\,S$. Hence, $S\subset\text{int}\,S$. Since we have shown that $\text{int}\,S\subset S$ is always true, we conclude that if $S$ is open, $S=\text{int}\, S$.

If $U$ is a subset of $S$ and $U$ is open, for every $\mathbf{x}$ in $U$, there is an $r>0$ such that $B(\mathbf{x}, r)\subset U$. But then $B(\mathbf{x}, r)\subset S$. This shows that $\mathbf{x}$ is in $\text{int}\,S$. Since every point of $U$ is in $\text{int}\,S$, this proves that $U\subset \text{int}\,S$. 
\end{myproof}

\begin{example}[label=230715_1]{}
Find the interior of each of the following subsets of $\mathbb{R}$.
\begin{enumerate}[(a)]
\begin{tabular}{p{5cm}p{5cm}}
\item $A=(a, b)$ &
\item $B=(a, b]$\\
\item $C=[a, b]$ &
\item $\mathbb{Q}$
\end{tabular}

\end{enumerate}
\end{example}
\begin{solution}{Solution}
\begin{enumerate}[(a)]
\item
Since $A$ is an open set, $\text{int}\,A=A=(a,b)$.

\item Since $A$ is an open set that is contained in $B$, $A=(a,b)$ is contained in $\text{int}\, B$. Since $\text{int}\, B\subset B$, we only left to determine whether $b$ is in $\text{int}\, B$. The same argument as given in Example \ref{230714_6} shows that $b$ is not an interior point of $B$. Hence, $\text{int}\,B=A=(a,b)$.
\end{enumerate}\bs\begin{enumerate}[(a)]

\item[(c)] Similar arguments as given in (b) show  that $A\subset\text{int}\, C$, and both $a$ and $b$ are not interior points of $C$. Hence, $\text{int}\,C=A=(a,b)$.

\item[(d)] For any $x\in \mathbb{R}$ and any $r>0$, $B(x,r)=(x-r,x+r)$ contains an irrational number. Hence, $B(x,r)$ is not contained in $\mathbb{Q}$. This shows that $\mathbb{Q}$ does not have interior points. Hence, $\text{int}\,\mathbb{Q}=\emptyset$.
\end{enumerate}
\end{solution}

\begin{definition}{Neighbourhoods}
Let $\mathbf{x}$ be a point in $\mathbb{R}^n$ and let $U$ be a subset of $\mathbb{R}^n$. We say that $U$ is a  neighbourhood of $\mathbf{x}$ if $U$ is an {\it open} set that contains $\mathbf{x}$.
\end{definition}Notice that this definition is slightly different from  the one we use in volume I for the $n=1$ case. 

\begin{highlight}{Neighbourhoods}
By definition, if $U$ is a neighbourhood of $\mathbf{x}$, then $\mathbf{x}$ is an interior point of $U$, and there is an $r>0$ such that $B(\mathbf{x},r)\subset U$. 
\end{highlight}

\begin{example}{}
Consider the point $\mathbf{x}=(1,2)$ and the sets 
\begin{align*}
U&=\left\{(x_1,x_2)\,|\,x_1^2+x_2^2<9 \right\},\\
V&=\left\{(x_1, x_2)\,|\, 0<x_1<2, -1<x_2<3\right\}
\end{align*} in $\mathbb{R}^2$. The sets $U$ and $V$ are neighbourhoods of $\mathbf{x}$.
\end{example}
 \begin{figure}[ht]
\centering
\includegraphics[scale=0.2]{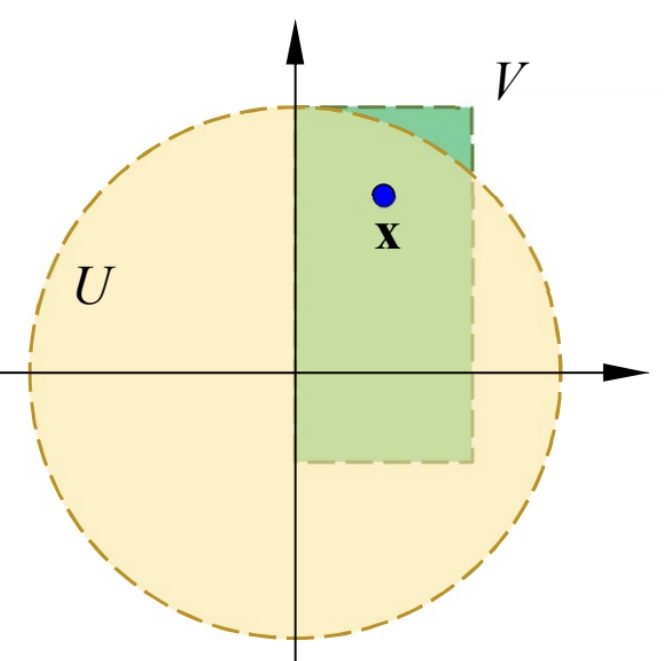}

\caption{The sets $U$ and $V$ are neighbourhoods of the point $\mathbf{x}$.}\label{figure16}
\end{figure}
Next, we introduce the exterior and boundary of a set.
\begin{definition}{Exterior}
Let $S$ be a subset of $\mathbb{R}^n$. We say that $\mathbf{x}\in \mathbb{R}^n$ is an exterior point of $S$ if  there exists $r>0$ such that $B(\mathbf{x}, r)\subset \mathbb{R}^n\setminus S$.  The exterior of $S$, denoted by $\text{ext}\,S$, is defined to be the collection of all the exterior points  of $S$. 
\end{definition}
\begin{definition}{Boundary}
Let $S$ be a subset of $\mathbb{R}^n$. 
We say that $\mathbf{x}\in \mathbb{R}^n$ 
  is a boundary point of $S$ if  for every $r>0$, the ball $B(\mathbf{x}, r)$   intersects both $S$ and    $\mathbb{R}^n\setminus S$.  The boundary of $S$, denoted by $\text{bd}\,S$ or $\pa S$, is defined to be the collection of all the  boundary points  of $S$. 
\end{definition}
 \begin{figure}[ht]
\centering
\includegraphics[scale=0.2]{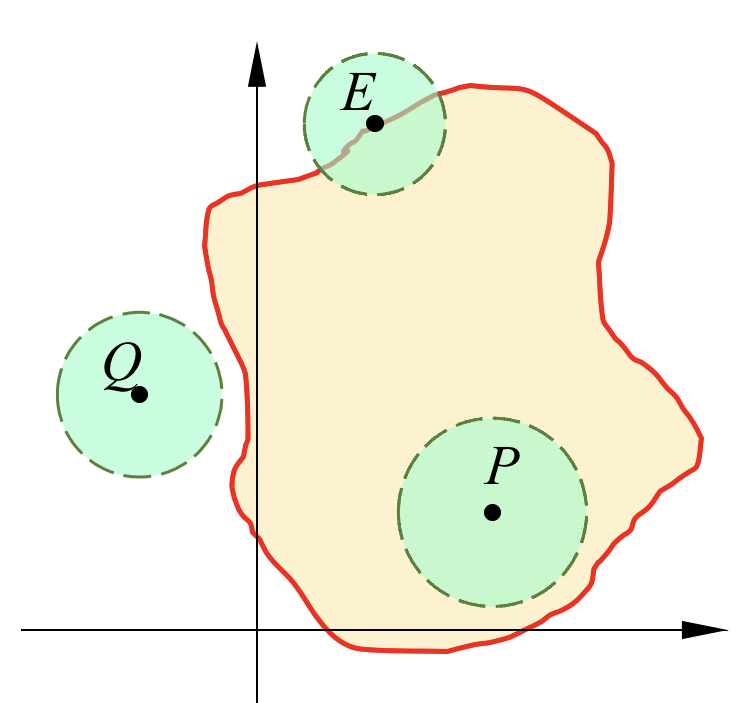}

\caption{$P$ is an interior point, $Q$ is an exterior point, $E$ is a boundary point.}\label{figure11}
\end{figure}
\begin{theorem}{}
Let $S$ be a subset of $\mathbb{R}^n$. We have the followings.
 \begin{enumerate}[(a)]
\item $\di \text{ext}\,(S)=\text{int}\,(\mathbb{R}^n\setminus S)$.
 \item $\di \text{bd}\,(S)=\text{bd}\,(\mathbb{R}^n\setminus S)$.
 \item $\text{int}\,S$, $\text{ext}\,S$ and $\text{bd}\,S$ are mutually disjoint sets.
 \item $\mathbb{R}^n =\text{int}\,S\,\cup\,\text{ext}\,S\,\cup\,\text{bd}\,S$.
\end{enumerate}
 \end{theorem}
 \begin{myproof}{Proof}
 (a) and (b) are obvious from definitions.
 
 For parts (c) and (d), we notice that for a point $\mathbf{x}\in \mathbb{R}^n$, exactly one of the following three statements holds.
 \begin{enumerate}[(i)]
 \item There exists $r>0$ such that $B(\mathbf{x}, r)\subset   S$. 
 \item There exists $r>0$ such that $B(\mathbf{x}, r)\subset \mathbb{R}^n\setminus S$. 
 \item For every $r>0$,  $B(\mathbf{x}, r)$   intersects both $S$ and    $\mathbb{R}^n\setminus S$.
 \end{enumerate}Thus, $\text{int}\,S$, $\text{ext}\,S$ and $\text{bd}\,S$ are mutually disjoint sets, and their union is $\mathbb{R}^n$.
 \end{myproof}

 \begin{example}[label=230715_10]{}
Find the exterior and boundary of each of the following subsets of $\mathbb{R}$.
\begin{enumerate}[(a)]
\begin{tabular}{p{5cm}p{5cm}}
\item $A=(a, b)$ &
\item $B=(a, b]$\\
\item $C=[a, b]$ &
\item $\mathbb{Q}$
\end{tabular}

\end{enumerate}
\end{example}
\begin{solution}{Solution}
We have seen in Example \ref{230715_1} that 
\[\text{int}\,A=\text{int}\,B=\text{int}\,C=(a,b).\]\bs
For any $r>0$, the ball $B(a,r)=(a-r,a+r)$ contains a point less than $a$, and a point larger than $a$. Hence, $a$ is a boundary point of the sets $A$, $B$ and $C$. Similarly, $b$ is a boundary point of the sets $A$, $B$ and $C$.

For every point $x$ which satisfies $x<a$, let $r=a-x$. Then $r>0$. Since  $x+r=a$, the ball $B(x,r)=(x-r, x+r)$ is contained in $(-\infty, a)$. Hence, $x$ is an exterior point of the sets $A$, $B$ and $C$.  Similarly every point $x$ such that $x>b$ is an exterior point of the sets $A$, $B$ and $C$. 

Since the interior, exterior and boundary of a set in $\mathbb{R}$  are three mutually disjoint sets whose union is $\mathbb{R}$, we conclude that
\begin{gather*}
\text{bd}\,A=\text{bd}\,B=\text{bd}\,C=\{a,b\},\\
\text{ext}\,A=\text{ext}\,B=\text{ext}\,C=(-\infty, a)\cup (b,\infty).
\end{gather*}
 For every $x\in \mathbb{R}$ and every $r>0$, the ball $B(x,r)=(x-r,x+r)$ contains a point in $\mathbb{Q}$ and a point not in $\mathbb{Q}$. Therefore, $x$ is a boundary point of $\mathbb{Q}$. This shows that $\text{bd}\,\mathbb{Q}=\mathbb{R}$, and thus, $\text{ext}\,\mathbb{Q}=\emptyset$.

\end{solution}

\begin{example}[label=230721_3]{}
Let $A=B(\mathbf{x}_0, r)$, where $\mathbf{x}_0$ is a point in $\mathbb{R}^n$, and $r$ is a positive number. Find the interior, exterior and boundary of $A$.
\end{example}
\begin{solution}{Solution}
We have shown that $A$ is open. Hence, $\text{int}\,A=A$. Let
\[U=\left\{\mathbf{x}\in\mathbb{R}^n\,|\, \Vert \mathbf{x}-\mathbf{x}_0\Vert>r\right\},\hspace{1cm} C=\left\{\mathbf{x}\in\mathbb{R}^n\,|\, \Vert \mathbf{x}-\mathbf{x}_0\Vert=r\right\}.\]Notice that $A$, $U$ and $C$ are mutually disjoint sets whose union is $\mathbb{R}^n$.

If $\mathbf{x}$ is in $U$, $d=\Vert \mathbf{x}-\mathbf{x}_0\Vert>r$. Let $r'=d-r$. Then $r'>0$. If $\mathbf{y}\in B(\mathbf{x}, r')$, then $\Vert\mathbf{y}-\mathbf{x}\Vert<r'$.
It follows that
\[\Vert\mathbf{y}-\mathbf{x}_0\Vert\geq\Vert\mathbf{x}-\mathbf{x}_0\Vert-\Vert\mathbf{y}-\mathbf{x}\Vert>d-r'=r.\]

This proves that $\mathbf{y}\in U$.  Hence, $B_d(\mathbf{x}, r')\subset U\subset \mathbb{R}^n\setminus A$, which shows that $\mathbf{x}$ is an exterior point of $A$. Thus, $U\subset\text{ext}\,A$.\bs

Now if $\mathbf{x}\in C$, $\Vert\mathbf{x}-\mathbf{x}_0\Vert=r$. For  every $r'>0$, let $a=\di\frac{1}{2}\min\{r'/r,1\}$. Then $\di a\leq\frac{1}{2}$ and $a\leq \di \frac{r'}{2r}$. Consider the point
\[\mathbf{v}=\mathbf{x}-a(\mathbf{x}-\mathbf{x}_0).\]
Notice that
\[\Vert\mathbf{v}-\mathbf{x}\Vert=ar\leq \frac{r'}{2}<r'.\]Thus, $\mathbf{v}$ is in $B(\mathbf{x}, r')$. On the other hand,
\[\Vert\mathbf{v}-\mathbf{x}_0\Vert=(1-a)r<r.\]Thus, $\mathbf{v}$ is in $A$. This shows that   $B(\mathbf{x}, r')$ intersects $A$. Since $\mathbf{x}$ is in $B(\mathbf{x}, r')$ but not in $A$, we find that $B(\mathbf{x}, r')$ intersects  $\mathbb{R}^n\setminus A$. Hence, $\mathbf{x}$ is a boundary point of $A$. This shows that $C\subset\text{bd}\,A$.

Since $\text{int}\,A$, $\text{ext}\,A$ and $\text{bd}\,A$ are mutually disjoint sets, we conclude that $\text{int}\,A=A$, $\text{ext}\,A=U$ and $\text{bd}\,A=C$.
\end{solution}

Now we introduce the closure of a set. 
\begin{definition}{Closure}
Let $S$ be a subset of $\mathbb{R}^n$. 
The closure of $S$, denoted by $\overline{S}$, is defined as 
\[\overline{S}=\text{int}\,S\cup \text{bd}\, S.\]
\end{definition}
\begin{example}{}
Example \ref{230721_3} shows that the closure of the open ball $B(\mf{x}_0, r)$ is the closed ball $CB(\mathbf{x}_0, r)$. 
\end{example}
\begin{example}{}
 Consider the sets $A=(a,b)$, $B=(a, b]$ and $C=[a, b]$  in Example \ref{230715_1} and Example \ref{230715_10}. We have shown that $\text{int}\,A=\text{int}\,B=\text{int}\,C=(a,b)$, and    $\text{bd}\,A=\text{bd}\,B=\text{bd}\,C=\{a,b\}$.  Therefore,  $\overline{A}=\overline{B}=\overline{C}=[a,b]$.
\end{example}

Since $\mathbb{R}^n$ is a disjoint union of $\text{int}\,S$, $\text{bd}\,S$ and $\text{ext}\,S$, 
we obtain the  following immediately from the definition.
\begin{theorem}{}
Let $S$ be a subset of $\mathbb{R}^n$. Then $\overline{S}$ and $\text{ext}\,S$ are complement of each other in $\mathbb{R}^n$.
\end{theorem}

The following theorem gives a characterization of the closure of a set.
\begin{theorem}[label=230715_5]{}
Let $S$ be a subset of $\mathbb{R}^n$, and let  $\mathbf{x}$ be a point in $\mathbb{R}^n$. The following statements are equivalent. 
\begin{enumerate}[(a)]
\item $\mathbf{x}\in \overline{S}$.
\item For every $r>0$, $B(\mathbf{x}, r)$ intersects $S$.
\item There is a sequence $\{\mathbf{x}_k\}$ in $S$ that converges to $\mathbf{x}$.
\end{enumerate}
\end{theorem}
\begin{myproof}{Proof}
If $\mathbf{x}$ is in $\overline{S}$, $\mathbf{x}$ is not in $\text{int}\,\left(\mathbb{R}^n\setminus S\right)$. Thus,
for every $r>0$, $B(\mathbf{x}, r)$ is not contained in $\mathbb{R}^n\setminus S$. Then it must  intersect  $S$. This proves (a) implies (b).

If (b) holds, for every $k\in\mathbb{Z}^+$, take $r=1/k$. The  ball $B(\mathbf{x}, 1/k)$ intersects $S$ at some point $\mathbf{x}_k$. This gives a sequence $\{\mathbf{x}_k\}$ satisfying
\[\Vert\mathbf{x}_k-\mathbf{x}\Vert<\frac{1}{k}.\]
Thus, $\{\mathbf{x}_k\}$ is a sequences in $S$ that converges to $\mathbf{x}$. This proves (b) implies (c).

If (c) holds, for every $r>0$, there is a positive integer $K$ such that for all $k\geq K$, $\Vert\mathbf{x}_k-\mathbf{x}\Vert<r$, and thus $\mathbf{x}_k\in B(\mathbf{x}, r)$. This shows that $B(\mathbf{x}, r)$ is not contained in $\mathbb{R}^n\setminus S$. Hence, $\mathbf{x}\notin \text{ext}\,S$, and thus we must have $\mathbf{x}\in\overline{S}$. This proves (c) implies (a).

\end{myproof}

The following theorem gives further properties of the closure of a set.
\begin{theorem}[label=230715_4]{}
Let $S$ be a subset of $\mathbb{R}^n$. 
\begin{enumerate}[1.]
\item $\overline{S}$ is a closed set that contains $S$.
\item $S$ is closed if and only if $S=\overline{S}$.
\item If $C$ is a closed subset of $\mathbb{R}^n$ and $S\subset C$, then $\overline{S}\subset C$.
\end{enumerate}
These imply that $\overline{S}$ is the smallest closed set that contains $S$.
\end{theorem}
\begin{myproof}{Proof}
These statements are counterparts of the statements in Theorem \ref{230715_3}. 

Since    $\text{ext}\,S=\text{int}\,(\mathbb{R}^n\setminus S)$, and the interior of a set is open, $\text{ext}\,S$ is open. Since
$\overline{S}=\mb{R}^n\setminus\text{ext}\,S$, $\overline{S}$ is a closed set. Since $\text{ext}\,S\subset \mathbb{R}^n\setminus S$, we find that 
\[\overline{S}=\mathbb{R}^n\setminus \text{ext}\,S\supset S.\] 

If $S=\overline{S}$, then $S$ must be closed since $\overline{S}$ is closed. Conversely, if $S$ is closed, $\mathbb{R}^n\setminus S$ is open, and so 
$\text{ext}\,S=\text{int}\,(\mathbb{R}^n\setminus S)=\mathbb{R}^n\setminus S$. It follows that $\overline{S}=\mathbb{R}^n\setminus\text{ext}\,S=S$.

If $C$ is a closed set that contains $S$, then $\mathbb{R}^n\setminus C$ is an open set that is contained in $\mathbb{R}^n\setminus S$. Thus, 
$\di \mathbb{R}^n\setminus C\subset \text{int}\,(\mathbb{R}^n\setminus S)=\text{ext}\,S$. This shows that $C\supset\mathbb{R}^n\setminus \text{ext}\,S=\overline{S}$.
\end{myproof}

\begin{corollary}{}
If $S$ be a subset of $\mathbb{R}^n$, 
$\di \overline{S}=S\cup\text{bd}\,S$. 
\end{corollary}
\begin{myproof}{Proof}
Since $\text{int}\,S\subset S$, $\overline{S}=\text{int}\,S\cup \text{bd}\,S\subset S\cup\text{bd}\,S$.
Since $S$ and $\text{bd}\,S$ are both subsets of $\overline{S}$, $S\cup\text{bd}\,S\subset \overline{S}$.  This proves that $\di \overline{S}=S\cup\text{bd}\,S$. 
\end{myproof}

\begin{example}{}
Let $U$ be the open rectangle 
$\di U=\prod_{i=1}^n(a_i, b_i)$
in $\mathbb{R}^n$. Show that the closure of  $U$ is the closed rectangle 
$\di R=\prod_{i=1}^n[a_i, b_i]$. 
\end{example}
\begin{solution}{Solution}Since $R$ is a closed set that contains $U$, $\overline{U}\subset R$.

If $\mathbf{x}=(x_1, \ldots, x_n)$ is a point in $R$, then $x_i\in [a_i, b_i]$ for each $1\leq i\leq n$. Since $[a_i, b_i]$ is the closure of $(a_i, b_i)$ in $\mathbb{R}$, there is a sequence $\{x_{i,k}\}_{k=1}^{\infty}$ in $(a_i, b_i)$ that converges to $x_i$. For $k\in\mathbb{Z}^+$, let
\[\mathbf{x}_k=(x_{1,k}, \ldots, x_{n,k}).\]Then $\{\mathbf{x}_k\}$ is a sequence in $U$ that converges to $\mathbf{x}$. This shows that $\mathbf{x}\in \overline{U}$, and thus completes the proof that $\overline{U}=R$.
\end{solution}

The proof of the following theorem shows the usefulness of the characterization of $\text{int}\,S$ as the largest open set that is contained in $S$, and $\overline{S}$ is the smallest closed set that contains $S$.
\begin{theorem}[label=230715_8]{}
If $A$ and $B$ are  subsets of $\mathbb{R}^n$ such that $A\subset B$, then
\begin{enumerate}[(a)]
\item
  $\text{int}\,A\subset \text{int}\,B$; and 
\item  $\overline{A}\subset\overline{B}$.
\end{enumerate}
\end{theorem}
\begin{myproof}{Proof}
Since $\text{int}\,A$ is an open set that is contained in $A$, it is an open set that is contained in $B$. By the fourth statement in Theorem \ref{230715_3}, $\text{int}\,A\subset\text{int}\,B$.

Since $\overline{B}$ is a closed set that contains $B$, it is a closed set that contains $A$. By the third statement in Theorem  \ref{230715_4}, $\overline{A}\subset\overline{B}$.
\end{myproof}

Notice that as subsets of $\mathbb{R}$, $(a, b)\subset (a, b]\subset [a,b]$. We have shown  that $\overline{(a,b)}=\overline{(a,b]}=\overline{[a,b]}$. In general, we have the following.
\begin{theorem}[label=230715_11]{}
If $A$ and $B$ are  subsets of $\mathbb{R}^n$ such that $A\subset B\subset \overline{A}$, then
   $\overline{A}=\overline{B}$.   
 
\end{theorem}

\begin{myproof}{Proof}
By Theorem \ref{230715_8},   $A\subset B$ implies that
  $\overline{A}\subset\overline{B}$, while $B\subset\overline{A}$  implies that $\overline{B}$ is contained in $\overline{\overline{A}}=\overline{A}$.
  Thus, we have 
  \[\overline{A}\subset \overline{B}\subset\overline{A},\]which proves that $\overline{B}=\overline{A}$.

\end{myproof}
\begin{example}{}
In general,   if $S$ is a subset of $\mathbb{R}^n$, it is not necessary true that $\text{int}\,S=\text{int}\,\overline{S}$, even when $S$ is an open set.  For example, take $S=(-1,0)\cup (0,1)$ in $\mathbb{R}$. Then $S$ is an open set and $\overline{S}=[-1,1]$. Notice that $\text{int}\,S=S=(-1,0)\cup(0,1)$, but $\text{int}\,\overline{S}=(-1,1)$. 
\end{example}

\vp
\noindent
{\bf \large Exercises  \thesection}
\setcounter{myquestion}{1}
\begin{question}{\themyquestion}
Let $S$ be a subset of $\mathbb{R}^n$. Show that $\text{bd}\,S$ is a closed set.
\end{question}
\atc

\begin{question}{\themyquestion}
Let $A$ be the subset of $\mathbb{R}^2$ given by \[
A=\left\{(x, y)\,|\, x<0, y\geq 0\right\}.\]
Find the interior, exterior, boundary and closure of $A$.

\end{question}

\atc
\begin{question}{\themyquestion}
Let $\mathbf{x}_0$ be a point in $\mathbb{R}^n$, and let $r$ be a positive number. Consider the subset of $\mathbb{R}^n$ given by
\[A=\left\{\mathbf{x}\in\mathbb{R}^n\,|\,0< \Vert\mathbf{x}-\mathbf{x}_0\Vert\leq r\right\}.\]  Find the interior, exterior, boundary and closure of $A$.

\end{question}
 \atc
 
 \begin{question}{\themyquestion}
Let $A$ be the subset of $\mathbb{R}^2$ given by
\[A=\left\{(x,y)\,|\, 1\leq x<3, -2<y\leq 5\right\}\cup\{(0,0), (2, -3)\}.\]
 Find the interior, exterior, boundary and closure of $A$.

\end{question}

\atc
\begin{question}{\themyquestion}
Let $S$ be a subset of $\mathbb{R}^n$. Show that \[\text{bd}\,S=\overline{S}\,\cap \,\overline{\mathbb{R}^n\setminus S}.\]
\end{question}

\atc
\begin{question}{\themyquestion}
Let $S$ be a subset of $\mathbb{R}^n$. Show that $\text{bd}\,\overline{S}\subset \text{bd}\,S$. Give an example where $\text{bd}\,\overline{S}\neq  \text{bd}\,S$.

\end{question}

\atc
\begin{question}{\themyquestion}
Let $S$ be a subset of $\mathbb{R}^n$.
\begin{enumerate}[(a)]
\item Show that $S$ is open if and only if $S$ does not contain any of its boundary points.
\item Show that $S$ is closed if and only if $S$ contains all its boundary points.

\end{enumerate}

\end{question}

\atc
\begin{question}{\themyquestion}
Let $S$ be a subset of $\mathbb{R}^n$, and let $\mathbf{x}$ be a point in $\mathbb{R}^n$. 
\begin{enumerate}[(a)]
\item Show that $\mathbf{x}$ is an interior point of $S$ if and only if there is a neighbourhood of $\mathbf{x}$ that is contained in $S$.

\item
Show that $\mathbf{x}\in \overline{S}$ if and only if every neighbourhood of $\mathbf{x}$ intersects $S$.
\item Show that $\mathbf{x}$ is a boundary point of $S$ if and only if every neighbourhood of $\mathbf{x}$ contains a point in $S$ and a point not in $S$.

\end{enumerate}

\end{question}

\atc
\begin{question}{\themyquestion}
Let $S$ be a subset of $\mathbb{R}^n$, and let $\mathbf{x}=(x_1, \ldots, x_n)$ be a point in the interior of $S$.
\begin{enumerate}[(a)]
\item
Show that there is an $r_1>0$ such that $CB(\mf{x}, r_1)\subset S$.
\item Show that there is an $r_2>0$ such that $\di \prod_{i=1}^n(x_i-r_2, x_i+r_2)\subset S$.

\item Show that there is an $r_3>0$ such that $\di \prod_{i=1}^n[x_i-r_3, x_i+r_3]\subset S$.
\end{enumerate}

\end{question}

\section{Limit Points and Isolated Points}
In this section, we generalize the concepts of limit points and isolated points to subsets of $\mathbb{R}^n$.

\begin{definition}{Limit Points}
Let $S$ be a subset of $\mathbb{R}^n$. A point $\mathbf{x}$ in $\mathbb{R}^n$ is a limit point of $S$ provided that there is a sequence $\{\mathbf{x}_k\}$ in $S\setminus\{\mathbf{x}\}$ that converges to $\mathbf{x}$. The set of limit points of $S$ is denoted by $S'$.

\end{definition}

By Theorem \ref{230715_5}, we obtain the following immediately.
\begin{theorem}[label=230715_7]{}
Let $S$ be a subset of $\mathbb{R}^n$, and let $\mathbf{x}$ be a point in $\mathbb{R}^n$. The following are equivalent.
\begin{enumerate}[(a)]
\item $\mathbf{x}$ is a limit point of $S$.
\item $\mathbf{x}$ is in $\overline{S\setminus\{\mathbf{x}\}}$.
\item For  every $r>0$, $B(\mathbf{x}, r)$ intersects $S$ at a point other than $\mathbf{x}$.\end{enumerate}
\end{theorem}
\begin{corollary}[label=230715_6]{}
If $S$ is a subset of $\mathbb{R}^n$, then $S'\subset \overline{S}$.
\end{corollary}
\begin{myproof}{Proof}
If $\mathbf{x}\in S'$, $\mathbf{x}\in\overline{S\setminus\{\mathbf{x}\}}$. Since $S\setminus\{\mathbf{x}\}\subset S$, we have $\overline{S\setminus\{\mathbf{x}\}}\subset \overline{S}$. Therefore, $\mathbf{x}\in\overline{S}$.
\end{myproof}

The following theorem says that the closure of a set is the union of the set with all its limit points.
\begin{theorem}{}
If $S$ is a subset of $\mathbb{R}^n$, then
$\di \overline{S}=S\cup S'$. 
\end{theorem}
\begin{myproof}{Proof}
By Corollary \ref{230715_6}, $S'\subset\overline{S}$. Since we also have $S\subset\overline{S}$, we find that $S\cup S'\subset \overline{S}$. 

Conversely, if $\mathbf{x}\in \overline{S}$, then by Theorem \ref{230715_5}, there is a sequence $\{\mathbf{x}_k\}$ in $S$ that converges to $\mathbf{x}$.  If $\mathbf{x}$ is not in $S$, then the sequence $\{\mathbf{x}_k\}$ is in $S\setminus\{\mathbf{x}\}$. In this case, $\mathbf{x}$ is a limit point of $S$. This shows that $\overline{S}\setminus S\subset S'$, and hence, $\overline{S}\subset S\cup  S'$.

\end{myproof}
In the proof above, we have shown the following.
\begin{corollary}
{}
Let $S$ be a subset of $\mathbb{R}^n$. Every point in $\overline{S}$ that is not in $S$ is a limit point of $S$. Namely,
\[\overline{S}\setminus S\subset S'.\]
\end{corollary}
Now we introduce the definition of isolated points. 
\begin{definition}{Isolated Points}
Let $S$ be a subset of $\mathbb{R}^n$. A point $\mathbf{x}$ in $\mathbb{R}^n$ is an isolated point of $S$ if
\begin{enumerate}[(a)]
\item
$\mathbf{x}$ is in $S$;
\item $\mathbf{x}$ is not a limit point of $S$.
\end{enumerate}
\end{definition}
\begin{remark}{}
By definition, a point $\mathbf{x}$ in $S$ is either an isolated point of $S$ or a limit point of $S$. 
\end{remark}
Theorem \ref{230715_7} gives the following immediately.
\begin{theorem}{}
Let $S$ be a subset of $\mathbb{R}^n$ and let $\mathbf{x}$ be a point in $S$. Then $\mathbf{x}$ is an isolated point of $S$ if and only if there is an $r>0$ such that  the ball $B(\mathbf{x}, r)$ does not contain other points of $S$ except the point $\mathbf{x}$. 
\end{theorem}
\begin{example}{}
Find the set of limit points and isolated points of the set $A=\mathbb{Z}^2$ as a subset of $\mathbb{R}^2$.
\end{example}
\begin{solution}{Solution}
If $\{\mathbf{x}_k\}$ is a sequence in $A$ that converges to a point $\mathbf{x}$, then there is a positive  integer $K$ such that for all $l\geq k\geq K$,
\[\Vert\mathbf{x}_l-\mathbf{x}_k\Vert<1.\] 
This implies that $\mathbf{x}_k=\mathbf{x}_K$ for all $k\geq K$. Hence, $\mathbf{x}=\mathbf{x}_K\in A$. This shows that $A$ is closed. Hence, $\overline{A}=A$. Therefore, $A'\subset A$.

For every $\mathbf{x}=(k,l)\in\mathbb{Z}^2$, $B(\mathbf{x}, 1)$ intersects $A$ only at the point $\mathbf{x}$ itself. Hence, $\mathbf{x}$ is an isolated point of $A$. This shows that every point of $A$ is an isolated point. Since $A'\subset A$, we must have $A'=\emptyset$.

\end{solution}

\begin{figure}[ht]
\centering
\includegraphics[scale=0.2]{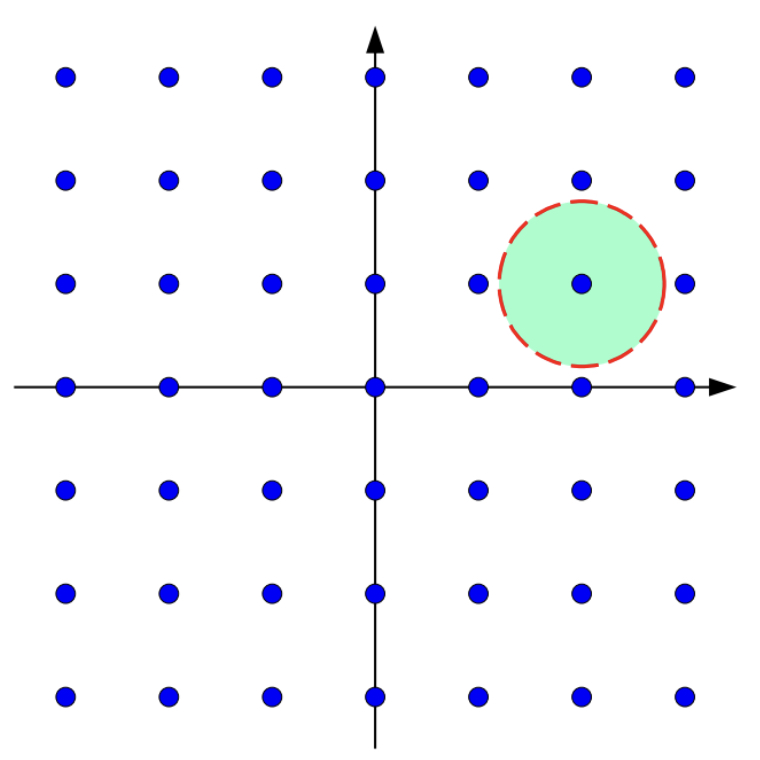}

\caption{The set $\mathbb{Z}^2$ does not have limit points.}\label{figure8}
\end{figure}

Let us prove the following useful fact.

\begin{theorem}{}
If $S$ is a subset of $\mathbb{R}^n$, every interior point of $S$ is a limit point of $S$. 
\end{theorem}
\begin{myproof}{Proof}
If $\mathbf{x}$ is an interior point of $S$, there exists $r_0>0$ such that $B(\mathbf{x}, r_0)\subset S$. Given $r>0$, let $r'=\di\frac{1}{2}\min\{r, r_0\}$. Then $r'>0$. Since $r'<r$ and $r'<r_0$, the point 
\[\mathbf{x}'=\mathbf{x}+r'\mathbf{e}_1\]   is in $B(\mathbf{x}, r)$ and $S$. Obviously, $\mathbf{x}'\neq \mathbf{x}$. Therefore, for every $r>0$, $B(\mathbf{x}, r)$ intersects $S$ at a point other than $\mathbf{x}$. This proves that $\mathbf{x}$ is a limit point of $S$.
\end{myproof}
Since $S\subset \text{int}\,S\cup\text{bd}\,S$, and $\text{int}\,S$ and $\text{bd}\,S$ are disjoint, we deduce the following.
\begin{corollary}{}
Let $S$ be a subset of $\mathbb{R}^n$. An isolated point of $S$ must be a boundary point.
\end{corollary}

Since every point in an open set $S$ is an interior point of $S$, we obtain the following.
\begin{corollary}{}
If $S$ is an open subset of $\mathbb{R}^n$, every point of $S$ is a limit point. Namely, $S\subset S'$.  
\end{corollary}

\begin{example}{}
If $I$ is an interval of the form $(a, b)$, $(a, b]$, $[a, b)$ or $[a,b]$ in $\mathbb{R}$, then $\text{bd}\,I=\{a, b\}$. It is easy to check that $a$ and $b$ are not isolated points of $I$. Hence, $I$ has no isolated points. Since $\overline{I}=I\cup I'$ and $I\subset I'$, we find that $I'=\overline{I}=[a,b]$. 
\end{example}

 In fact, we can prove a general theorem.
 \begin{theorem}{}
 Let $A$ and $B$ be subsets of $\mathbb{R}^n$ such that $A$ is open and $A\subset B\subset \overline{A}$. Then $B'=\overline{A}$. In particular, the set of limit points of $A$ is $\overline{A}$.
 \end{theorem}
 \begin{myproof}{Proof}By Theorem \ref{230715_11}, $\overline{A}=\overline{B}$. 
 Since $A$ is open, $A\subset A'$. 
 Since $\overline{A}=A\cup A'$, we find that $\overline{A}=A'$.
 
 In the exercises, one is asked to show that $A\subset B$ implies $A'\subset B'$. Therefore, $\overline{A}=A'\subset B'\subset \overline{B}$. Since $\overline{A}=\overline{B}$, we must have $B'=\overline{B}=\overline{A}$.
 \end{myproof}

\begin{example}{}Let $A$ be the subset of $\mathbb{R}^2$ given by
\[A=[-1,1]\times (-2,2]=\left\{(x,y)\,|\, -1\leq x\leq 1, -2<y\leq 2\right\}.\] Since
 $U=(-1,1)\times (-2,2)$ is open, $\overline{U}=[-1,1]\times [-2,2]$, and $U\subset A\subset \overline{U}$, the set of limit points of $A$ is $\overline{U}=[-1,1]\times [-2, 2]$. 
\end{example}
 
\vp
\noindent
{\bf \large Exercises  \thesection}
\setcounter{myquestion}{1}

\begin{question}{\themyquestion}
Let $A$ and $B$ be   subsets of $\mathbb{R}^n$ such that $A\subset B$. Show that $A'\subset B'$.

\end{question}
\atc
\begin{question}{\themyquestion}
Let $\mathbf{x}_0$ be a point in $\mathbb{R}^n$ and let $r$ be a positive number. Find the set of limit points of the open ball $B(\mathbf{x}_0, r)$.

\end{question}
\atc

\begin{question}{\themyquestion}
Let $A$ be the subset of $\mathbb{R}^2$ given by \[
A=\left\{(x, y)\,|\, x<0, y\geq 0\right\}.\]
Find the set of limit points of $A$.

\end{question}
\atc

\begin{question}{\themyquestion}
Let $\mathbf{x}_0$ be a point in $\mathbb{R}^n$, and let $r$ is a positive number. Consider the subset of $\mathbb{R}^n$ given by
\[A=\left\{\mathbf{x}\in\mathbb{R}^n\,|\,0< \Vert\mathbf{x}-\mathbf{x}_0\Vert\leq r\right\}.\]   \begin{enumerate}[(a)]
\item  Find the set of limit points of $A$.
\item Find the set of isolated points of the set $S=\mathbb{R}^n\setminus A$.
\end{enumerate}
\end{question}
 \atc
 
 \begin{question}{\themyquestion}
Let $A$ be the subset of $\mathbb{R}^2$ given by
\[A=\left\{(x,y)\,|\, 1\leq x<3, -2<y\leq 5\right\}\cup\{(0,0), (2, -3)\}.\]
Determine the set of isolated points and the set of limit points of $A$.

\end{question}

\atc
\begin{question}{\themyquestion}
Let $A=\mathbb{Q}^2$  as a subset of $\mathbb{R}^2$.
\begin{enumerate}[(a)]
\item Find the interior, exterior, boundary and closure of $A$.
\item Determine the set of isolated points   and  the set of limit points of $A$.
\end{enumerate}

\end{question}

\atc
\begin{question}{\themyquestion}
Let $S$ be a subset of $\mathbb{R}^n$. Show that $S$ is closed if and only if it contains all its limit points.
\end{question}

\atc
\begin{question}{\themyquestion}
Let $S$ be a subset of $\mathbb{R}^n$, and let $\mathbf{x}$ be a point in $\mathbb{R}^n$. 
Show that $\mathbf{x}$ is a limit point of $S$ if and only if every neighbourhood of $\mathbf{x}$ intersects $S$ at a point other than itself.

\end{question}
\atc
\begin{question}{\themyquestion}
Let $\mathbf{x}_1$, $\mf{x}_2$, $\ldots$, $\mf{x}_k$ be points in $\mathbb{R}^n$ and let $A=\mathbb{R}^n\setminus\left\{\mathbf{x}_1, \mf{x}_2, \ldots, \mf{x}_k\right\}$.
Find the set of limit points of $A$.

\end{question}

\chapter{Limits of Multivariable Functions and Continuity}\label{chapter2}
 We are interested in functions  $\mathbf{F}:\mathfrak{D}\to \mathbb{R}^m$ that are defined on   subsets  $\mathfrak{D}$ of $\mathbb{R}^n$, taking values in $\mb{R}^m$. When $n\geq 2$, these are called multivariable functions. When $m\geq 2$, they are called  vector-valued functions. When $m=1$, we usually write the function as  $f:\mathfrak{D}\to \mathbb{R}$.
 
\section{Multivariable Functions}\label{sec_mvf}

In this section, let us define some special classes of multivariable functions.

\subsection{Polynomials and Rational Functions}
A special class of functions is the set of polynomials in $n$ variables.
\begin{definition}{Polynomials}
Let $\mathbf{k}=(k_1, \ldots, k_n)$ be an $n$-tuple of nonnegative integers. 
Associated to this $n$-tuple $\mathbf{k}$, there is a monomial $p_{\mathbf{k}}:\mathbb{R}^n\to\mathbb{R}$ of degree $|\mf{k}|=k_1+\cdots+k_n$  of the form 
$\di p_{\mathbf{k}}(\mathbf{x})=x_1^{k_1}\cdots x_n^{k_n}$.

A polynomial in $n$ variables is a function $p:\mathbb{R}^n\to \mathbb{R}$ that is a {\it finite} linear combination of monomials in $n$ variables. It takes the form
\[p(\mathbf{x})=\sum_{j=1}^m c_{\mathbf{k}_j} p_{\mathbf{k}_j}(\mathbf{x}),\]
where  $\mf{k}_1, \mathbf{k}_2, \ldots, \mf{k}_m$ are  distinct $n$-tuples of nonnegative integers, and $c_{\mf{k}_1}, c_{\mathbf{k}_2}, \ldots, c_{\mf{k}_m}$  are nonzero real numbers. The degree of the polynomial $p(\mathbf{x})$ is $\max\{|\mf{k}_1|, |\mf{k}_2|, \ldots, |\mf{k}_m|\}$.
\end{definition}
\begin{example}{}
The following are examples of polynomials in three variables.
\begin{enumerate}[(a)]
\item $\di p(x_1, x_2, x_3)=x_1^2+x_2^2+x_3^2$
\item $\di p(x_1, x_2, x_3)=4x_1^2x_2-3x_1x_3+x_1x_2x_3$

\end{enumerate}
\end{example}

\begin{example}{}
The function $f:\mathbb{R}^n\to\mathbb{R}$,
\[f(\mathbf{x})=\Vert\mathbf{x}\Vert=\sqrt{x_1^2+\cdots+x_n^2}\] is not a polynomial.
\end{example}

When the domain of a function is not specified, we always assume that the domain is the largest set on which the function can be defined.
\begin{definition}{Rational Functions}
A rational function $f:\mathfrak{D}\to \mathbb{R}$ is the quotient of two polynomials $p:\mathbb{R}^n\to\mathbb{R}$ and $q:\mathbb{R}^n\to \mathbb{R}$. Namely,
\[f(\mathbf{x})=\frac{p(\mathbf{x})}{q(\mathbf{x})}.\]Its domain $\mathfrak{D}$ is the set
\[\mathfrak{D}=\left\{\mathbf{x}\in\mathbb{R}^n\,|\, q(\mathbf{x})\neq 0\right\}.\]

\end{definition}
\begin{example}{}
The function
\[f(x_1, x_2)=\frac{x_1x_2+3x_1^2}{x_1-x_2}\]
is a rational function defined on the set
\[\mathfrak{D}=\left\{(x_1, x_2)\in\mathbb{R}^2\,|\, x_1\neq x_2\right\}.\]
\end{example}

\subsection{Component Functions of a Mapping}
 If the codomain $\mathbb{R}^m$ of the function $\mathbf{F}:\mathfrak{D}\to \mathbb{R}^m$ has dimension $m\geq 2$, we usually call the function   a {\it mapping}. In this case, it would be good to consider the component functions. 

For $1\leq j\leq m$, the projection function $\pi_j:\mathbb{R}^m\to\mathbb{R}$ is the function 
\[\pi_j(x_1, \ldots, x_m)=x_j.\]
\begin{definition}{Component Functions}
Let $\mathbf{F}:\mathfrak{D}\to \mathbb{R}^m$ be a function defined on $\mathfrak{D}\subset\mathbb{R}^n$. For $1\leq j\leq m$, the $j^{\text{th}}$ component function of  $\mathbf{F}$ is the function $F_j:\mathfrak{D}\to\mathbb{R}$ defined as 
\[F_j=(\pi_j\circ \mathbf{F}):\mathfrak{D}\to\mathbb{R}.\]For each $\mathbf{x}\in\mathfrak{D}$, 
\[\mathbf{F}(\mathbf{x})=(F_1(\mathbf{x}), \ldots, F_m(\mathbf{x})).\]

\end{definition}
\begin{example}{}
For the function $\mathbf{F}:\mathbb{R}^3\to\mathbb{R}^3$, 
$\mathbf{F}(\mathbf{x})=-3\mathbf{x}$, 
the component functions are $F_1(x_1, x_2, x_3)=-3x_1$, $F_2(x_1, x_2, x_3)=-3x_2$, $F_3(x_1, x_2, x_3)=-3x_3$.
\end{example}

For convenience, we also define the notion of polynomial mappings.
 
\begin{definition}{Polynomial Mappings}
We call a function $\mf{F}:\mb{R}^n\to\mb{R}^m$ a polynomial mapping if each of its components  $F_j:\mb{R}^n\to\mb{R}$, $1\leq j\leq m$, is a polynomial function. The degree of the polynomial mapping $\mf{F}$ is the maximum of the degrees of the polynomials $F_1, F_2, \ldots, F_m$.

\end{definition}

\begin{example}{}
The mapping $\mf{F}:\mb{R}^3\to \mb{R}^2$,

\vspace{-0.4cm}
\[\mf{F}(x,y,z)=(x^2y+3xz, 8yz^3-7x)\] is a polynomial mapping of degree 4.
\end{example}
\subsection{Invertible Mappings}

The invertibility of a function $\mf{F}:\mk{D}\to \mb{R}^m$ is defined in the following way.
\begin{definition}{Inverse Functions}

Let $\mk{D}$ be a subset of $\mb{R}^n$, and let  $\mf{F}:\mk{D}\to \mb{R}^m$ be a function defined on $\mk{D}$. We say that $\mf{F}$ is invertible if $\mf{F}$ is one-to-one. In this case,  the inverse function $\mf{F}^{-1}:\mf{F}(\mk{D})\to \mk{D}$ is defined  so that for each $\mf{y}\in \mf{F}(\mk{D})$,
\[\mf{F}^{-1}(\mf{y})=\mf{x}\hspace{1cm}\text{if and only if}\hspace{1cm}\mf{F}(\mf{x})=\mf{y}.\]
\end{definition}

\begin{example}{}
Let $\mk{D}=\left\{(x,y)\,|\,x>0, y>0\right\}$ and let $\mf{F}:\mk{D}\to\mb{R}^2$ be the function defined as 
\[\mf{F}(x,y)=(x-y, x+y).\]Show that $\mf{F}$ is invertible and find its inverse.
\end{example}
\begin{solution}{Solution}
Let $u=x-y$ and $v=x+y$. Then
\[x=\frac{u+v}{2}, \hspace{1cm}y=\frac{v-u}{2}.\]
This shows that for any $(u, v)\in \mb{R}^2$, there is at most one pair of $(x,y)$ such that $\mf{F}(x,y)=(u,v)$. Thus, $\mf{F}$ is one-to-one, and hence, it is invertible. Observe that
\[\mf{F}(\mk{D})=\left\{(u,v)\,|\, v>0, -v<u<v.\right\}.\]

The inverse mapping is given by 
$\mf{F}^{-1}:\mf{F}(\mk{D})\to\mb{R}^2$,
\[\mf{F}^{-1}(u,v)=\left(\frac{u+v}{2}, \frac{v-u}{2}\right).\]
\end{solution}

\subsection{Linear Transformations}
Another special class of functions consists of linear transformations. A function $\mathbf{T}:\mathbb{R}^n\to\mathbb{R}^m$ is a  linear transformation if for any $\mathbf{x}_1, \ldots, \mathbf{x}_k$ in $\mathbb{R}^n$, and for any $c_1, \ldots, c_k$ in $\mathbb{R}$,
\[\mathbf{T}(c_1\mathbf{x}_1+\cdots+c_k\mathbf{x}_k)=c_1\mathbf{T}(\mathbf{x}_1)+\cdots+c_k\mathbf{T}(\mathbf{x}_k).\]Linear transformations are closely related to matrices.

An $m\times n$ matrix $A$ is an array with $m$ rows and $n$ columns of real numbers.  It has the form
\[A=[a_{ij}]=\begin{bmatrix} a_{11} & a_{12}& \cdots & a_{1n}\\a_{21}& a_{22} & \cdots & a_{2n}\\\vdots & \vdots & \ddots & \vdots\\
a_{m1} & a_{m2} & \cdots & a_{mn}\end{bmatrix}.\]
If $A=[a_{ij}]$ and $B=[b_{ij}]$ are $m\times n$ matrices, $\alpha$ and $\beta$ are real numbers,  $\alpha A+\beta B$ is defined  to be the $m\times n$ matrix $C=\alpha A+\beta B=[c_{ij}]$ with
\[c_{ij}=\alpha a_{ij}+\beta b_{ij}.\]
If $A=[a_{il}]$ is a $m\times k$ matrix, $B=[b_{lj}]$ is a $k\times n$ matrix, the product $AB$ is defined to be the $m\times n$ matrix $C=AB=[c_{ij}]$, where
\[c_{ij}=\sum_{l=1}^k a_{il}b_{lj}.\]It is easy to verify that matrix multiplications are associative.

Given $\mathbf{x}=(x_1, \ldots, x_n)$ in $\mathbb{R}^n$, we identify it with the column  vector
\[\mathbf{x}=\begin{bmatrix} x_1\\x_2\\\vdots\\x_n\end{bmatrix},\]
which is an $n\times 1$ matrix. If $A$ is an $m\times n$ matrix, 
   and $\mathbf{x}$ is a vector in $\mathbb{R}^n$, then $\mathbf{y}=A\mathbf{x}$ is  the vector in $\mathbb{R}^m$ given by
\[\mathbf{y}=A\mathbf{x}=\begin{bmatrix} a_{11} & a_{12}& \cdots & a_{1n}\\a_{21}& a_{22} & \cdots & a_{2n}\\\vdots & \vdots & \ddots & \vdots\\
a_{m1} & a_{m2} & \cdots & a_{mn}\end{bmatrix}\begin{bmatrix} x_1\\x_2\\\vdots\\x_n\end{bmatrix}=\begin{bmatrix}a_{11}x_1+a_{12}x_2+\cdots+a_{1n}x_n\\a_{21}x_1+a_{22}x_2+\cdots+a_{2n}x_n\\\vdots \\a_{m1}x_1+a_{m2}x_2+\cdots+a_{mn}x_n\end{bmatrix}.\]

The following is a standard result in linear algebra.
\begin{theorem}
{}A function $\mathbf{T}:\mathbb{R}^n\to\mathbb{R}^m$ is a  linear transformation if and only if there exists an $m\times n$ matrix $A=[a_{ij}]$ such that
\[\mathbf{T}(\mathbf{x})=A\mathbf{x}.\]
In this case, $A$ is called the matrix associated to the linear transformation  $\mathbf{T}:\mathbb{R}^n\to\mathbb{R}^m$.

\end{theorem}
\begin{myproof}{Sketch of Proof}
It is easy to verify that  the mapping  $\mathbf{T}:\mathbb{R}^n\to\mathbb{R}^m$, $\mathbf{T}(\mathbf{x})=A\mathbf{x}$ is a  linear transformation if $A$ is an $m\times n$ matrix.

Conversely, if $\mathbf{T}:\mathbb{R}^n\to\mathbb{R}^m$ is a  linear transformation, 
then for any $\mathbf{x}\in\mathbb{R}^n$, \[\mathbf{T}(\mathbf{x})=\mathbf{T}(x_1\mathbf{e}_1+x_2\mathbf{e}_2+\cdots+x_n\mathbf{e}_n)=x_1\mathbf{T}(\mathbf{e}_1)+x_2\mathbf{T}(\mathbf{e}_2)+\cdots +x_n\mathbf{T}(\mathbf{e}_n).\]
Define the vectors $\mathbf{a}_1$, $\mathbf{a}_2$, $\ldots$,  $\mathbf{a}_n$ in $\mathbb{R}^m$ by
\[\mathbf{a}_1=\mathbf{T}(\mathbf{e}_1), \;\mathbf{a}_2=\mathbf{T}(\mathbf{e}_2), \;\ldots, \;\mathbf{a}_n=\mathbf{T}(\mathbf{e}_n).\]
Let $A$ be the $m\times n$ matrix with column vectors  $\mathbf{a}_1$, $\mathbf{a}_2$, $\ldots$,  $\mathbf{a}_n$. Namely,
\[A=\begin{bmatrix}   \mathbf{a}_1 &\rvline& \mathbf{a}_2 &\rvline& \cdots &\rvline&\mathbf{a}_n\end{bmatrix}.\]
Then we have
$\di \mathbf{T}(\mathbf{x})=A\mathbf{x}$.
\end{myproof}

\begin{example}{}
Let $\mf{F}:\mb{R}^2\to\mb{R}^2$ be the function defined as

\[\mf{F}(x,y)=(x-y, x+y).\] Then $\mf{F}$ is a linear transformation with matrix
$\di A=\begin{bmatrix} 1 &-1\\ 1 & 1\end{bmatrix}$.
\end{example}
 
For the linear transformation $T:\mathbb{R}^n\to\mathbb{R}^m$, $ \mathbf{T}(\mathbf{x})=A\mathbf{x}$, the component functions are
\begin{align*}
T_1(\mathbf{x})&=a_{11}x_1+a_{12}x_2+\cdots+a_{1n}x_n,\\
T_2(\mathbf{x})&=a_{21}x_1+a_{22}x_2+\cdots+a_{2n}x_n,\\
&\hspace{2cm}\vdots\\
T_m(\mathbf{x})&=a_{m1}x_1+a_{m2}x_2+\cdots+a_{mn}x_n.
\end{align*}Each of them is a polynomial of degree at most one. Thus, a linear transformation is a polynomial mapping of degree at most one. It is easy to deduce the following.
 
\begin{corollary}{}
A mapping  $\mathbf{T}:\mathbb{R}^n\to\mathbb{R}^m$ is a linear transformation if and only if each component function is a linear transformation.
\end{corollary}

The followings are some standard results about linear transformations.
\begin{theorem}{}
If $\mathbf{S}:\mathbb{R}^n\to\mathbb{R}^m$ and $\mathbf{T}:\mathbb{R}^n\to\mathbb{R}^m$ are linear transformations with matrices $A$ and $B$ respectively, then for any real numbers $\alpha$ and $\beta$, $\alpha\mathbf{S}+\beta\mathbf{T}:\mathbb{R}^n\to\mathbb{R}^m$ is a linear transformation with matrix $\alpha A+\beta B$.
\end{theorem}
\begin{theorem}{}
If $\mathbf{S}:\mathbb{R}^n\to\mathbb{R}^m$ and $\mathbf{T}:\mathbb{R}^m\to\mathbb{R}^k$ are linear transformations with matrices $A$ and $B$, then $\mathbf{T}\circ\mathbf{S}:\mathbb{R}^n\to\mathbb{R}^k$ is a linear transformation with matrix $BA$.
\end{theorem}
\begin{myproof}{Sketch of Proof}
This follows from 
\[(\mathbf{T}\circ\mathbf{S})(\mathbf{x})=\mathbf{T}(\mathbf{S}(\mathbf{x}))=B(A\mathbf{x})=(BA)\mathbf{x}.\]
\end{myproof}

In the particular case when $m=n$, we have the following.
\begin{theorem}[label=230725_2]{}
Let   $\mathbf{T}:\mathbb{R}^n\to\mathbb{R}^n$ be a linear transformation represented by the matrix $A$.
The following are equivalent.
\begin{enumerate}[(a)]
\item The mapping $\mathbf{T}:\mathbb{R}^n\to\mathbb{R}^n$ is one-to-one.
\item The mapping $\mathbf{T}:\mathbb{R}^n\to\mathbb{R}^n$ is onto.

\item The matrix $A$ is invertible.
\item $\det A\neq 0$.

\end{enumerate}
\end{theorem}
In other words, if the linear transformation $\mathbf{T}:\mathbb{R}^n\to\mathbb{R}^n$ is one-to-one {\it or} onto, then it is bijective. In this case,   the linear transformation is invertible, and we can define the inverse function $\mathbf{T}^{-1}:\mathbb{R}^n\to\mathbb{R}^n$.

\begin{theorem}{}
Let   $\mathbf{T}:\mathbb{R}^n\to\mathbb{R}^n$ be an invertible  linear transformation represented by the matrix $A$.
Then the inverse mapping $\mf{T}^{-1}:\mathbb{R}^n\to\mathbb{R}^n$ is also a linear transformation and 
\[\mf{T}^{-1}(\mf{x})=A^{-1}\mf{x}.\]
\end{theorem}

\begin{example}{}
Let $\mf{T}:\mb{R}^2\to\mb{R}^2$ be the linear transformation   
\[ \mf{T}(x,y)=(x-y, x+y).\] The matrix associated with $\mf{T}$ is 
$\di A=\begin{bmatrix} 1 &-1\\ 1 & 1\end{bmatrix}$. Since $\det A=2\neq 0$, $\mf{T}$ is invertible. 
Since
$\di A^{-1}=\frac{1}{2}\begin{bmatrix} 1 & 1\\ -1 &1\end{bmatrix}$, we have
 \[\mf{T}^{-1}(x,y)=\left(\frac{x+y}{2}, \frac{-x+y}{2}\right).\]
\end{example}

\subsection{Quadratic Forms}\label{quadraticforms}

 Given an $m\times n$ matrix $A=[a_{ij}]$, its {\it transpose} is the $n\times m$ matrix $A^T=[b_{ij}]$, where 
  \[b_{ij}=a_{ji}\hspace{1cm}\text{for all}\; 1\leq i\leq n, 1\leq j\leq m.\]
   An $n\times n$  matrix $A$ is {\bf symmetric} if \[A=A^T.\]
   An $n\times n$ matrix $P$ is {\bf orthogonal} if  \[P^TP=PP^T=I.\] If the column vectors of $P$ are $\mf{v}_1$, $\mf{v}_2$, $\ldots$, $\mf{v}_n$, so that 
   \begin{equation}\label{230724_2}P=\begin{bmatrix}\mf{v}_1&\rvline &\mf{v}_2 &\rvline & \cdots & \rvline &\mf{v}_n\end{bmatrix},\end{equation}
   then $P$ is orthogonal if and only if $\{\mf{v}_1, \ldots, \mf{v}_n\}$ is an orthonormal set of vectors in $\mb{R}^n$.

   If $A$ is an $n\times n$ symmetric matrix, its characteristic polynomial
   \[p(\lambda)=\det(\lambda I_n -A)\]
   is a monic polynomial of degree $n$ with $n$ real roots $\lambda_1, \lambda_2, \ldots, \lambda_n$, counting with multiplicities. These roots are called the {\bf eigenvalues} of $A$. There is an orthonormal set of vectors  $\{\mf{v}_1, \ldots, \mf{v}_n\}$ in $\mb{R}^n$ such that
   \begin{equation}\label{230724_1}A\mf{v}_i=\lambda_i \mf{v}_i\hspace{1cm} \text{for all}\; 1\leq i\leq n.\end{equation}
  Let $D$ be the diagonal matrix 
  \begin{equation}\label{230724_7} D=\begin{bmatrix} \lambda_1 & 0 & \cdots & 0\\
  0 &\lambda_2 & \cdots & 0\\
  \vdots & \vdots & \ddots & \vdots\\
  0 & 0 & \cdots & \lambda_n\end{bmatrix},\end{equation}
  and let $P$ be the orthogonal matrix \eqref{230724_2}. Then \eqref{230724_1} is equivalent to $AP=PD$, or equivalently,
  \[A=PDP^T=PDP^{-1}.\]
  This is known as the orthogonal diagonalization of the real symmetric matrix $A$.

 A quadratic form in $\mb{R}^n$ is a polynomial function $Q:\mb{R}^n\to\mb{R}$ of the form
 \[Q(\mf{x})=\sum_{1\leq i<j\leq n}c_{ij}x_ix_j.\]
 An $n\times n $ symmetric matrix $A=[a_{ij}]$ defines a quadratic form  $Q_A:\mb{R}^n\to\mb{R}$ by
 \[Q_A(\mf{x})=\mathbf{x}^TA\mf{x}=\sum_{i=1}^n\sum_{j=1}^n a_{ij}x_ix_j.\]

\begin{example}{}The symmetric matrix $A=\di\begin{bmatrix} 1 & -2\\-2 & 5\end{bmatrix}$ defines the quadratic form
 \[Q_A(x,y)=x^2-4xy+5y^2.\]
 \end{example}
 
 Conversely,  given a quadratic form \[Q(\mf{x})=\sum_{1\leq i<j\leq n}c_{ij}x_ix_j,\]then $Q=Q_A$, where the entries of $A=[a_{ij}]$ are
 \[a_{ij}=\begin{cases} c_{ii}, \quad & \text{if}\; i=j,\\
 c_{ij}/2,\quad &\text{if}\; i<j,\\
 c_{ji}/2,\quad &\text{if}\;i>j.\end{cases}\]  Thus, there is a one-to-one correspondence between quadratic forms and symmetric matrices.
  
 If $A=PDP^T$ is  an orthogonal diagonalization of $A$, under the change of variables
 \[\mf{y}=P^T\mf{x},\quad \text{ or equivalently},\quad
 \mf{x}=P\mf{y}\]
 we find that
 \begin{equation}\label{230724_3}Q_A =\mf{y}^TD\mf{y}=\lambda_1y_1^2+\cdots+\lambda_ny_n^2.\end{equation}
  A  consequence of \eqref{230724_3} is the following.
 \begin{theorem}[label=230724_8]{}
 Let $A$ be an $n\times n$ symmetric matrix, and let $Q_A(\mf{x})=\mf{x}^TA\mf{x}$ be the associated quadratic form.    Let $\lambda_1, \lambda_2, \ldots, \lambda_n$ be the eigenvalues of $A$. Assume that
 \[  \lambda_n\leq\cdots\leq \lambda_{2} \leq \lambda_1.\]
 Then for any $\mf{x}\in\mb{R}^n$,  
 \[\lambda_n\Vert\mf{x}\Vert^2\leq Q_A(\mf{x})\leq \lambda_1\Vert \mf{x}\Vert^2.\]
 \end{theorem}
 \begin{myproof}{Sketch of Proof}
   Given $\mf{x}\in\mb{R}^n$, let $\mf{y}=P^T\mf{x}$. Then
    \[\Vert\mf{y}\Vert^2=\mf{y}^T\mf{y}=\mf{x}^TPP^T\mf{x}=\mf{x}^T\mf{x}=\Vert\mf{x}\Vert^2.\]  By \eqref{230724_3} ,
 \[Q_A(\mf{x})=\lambda_1y_1^2+\cdots+\lambda_ny_n^2.\]
Since 
 $  \lambda_n\leq\cdots\leq \lambda_{2} \leq \lambda_1$, we find that
 \[\lambda_n(y_1^2+\ldots+y_n^2)\leq Q_A(\mf{x})\leq \lambda_1(y_1^2+\ldots+y_n^2).\]
The assertion follows.

 \end{myproof}
 
 At the end of this section, let us recall the   classification of quadratic forms.
 
 \begin{highlight}{Definiteness of Symmetric Matrices}
 Given an $n\times n$ symmetric matrix $A=[a_{ij}]$, let $Q_A:\mb{R}^n\to\mb{R}$,
 \[Q_A(\mf{x})=\mathbf{x}^TA\mf{x}=\sum_{i=1}^n\sum_{j=1}^n a_{ij}x_ix_j \] be the associated quadratic form.
 \end{highlight}\begin{highlight}{}
 \begin{enumerate}[1.]
 \item We say that the matrix $A$ is positive definite, or the quadratic form $Q_A$ is positive definite, if $Q_A(\mf{x})>0$ for all $\mf{x}\neq \mf{0}$ in $\mb{R}^n$. 
  
  \item[2.] We say that the matrix $A$ is negative definite, or the quadratic form $Q_A$ is negative definite, if $Q_A(\mf{x})<0$ for all $\mf{x}\neq \mf{0}$ in $\mb{R}^n$. 
  
   \item[3.] We say that the matrix $A$ is indefinite, or the quadratic form $Q_A$ is indefinite, if there exist $\mf{u}$ and $\mf{v}$ in $\mb{R}^n$ such that $Q_A(\mf{u})>0$  and  $Q_A(\mf{v})<0$.

   \item[4.] We say that the matrix $A$ is positive semi-definite, or the quadratic form $Q_A$ is positive semi-definite, if $Q_A(\mf{x})\geq 0$ for all $\mf{x}$ in $\mb{R}^n$. 
     
  \item[5.] We say that the matrix $A$ is negative semi-definite, or the quadratic form $Q_A$ is negative semi-definite, if $Q_A(\mf{x})\leq 0$ for all $\mf{x}$ in $\mb{R}^n$. 
 \end{enumerate}
 \end{highlight}
 
 Obviously, a symmetric matrix $A$ is negative definite if and only if $-A$ is positive definite.
 
 The following is a standard result in linear algebra, which can be deduced  from \eqref{230724_3}.
 \begin{theorem}[label=230724_20]{}
 Let $A$ be an $n\times n$ symmetric matrix, and let $Q_A(\mf{x})=\mf{x}^TA\mf{x}$ be the associated quadratic form. Let $\{\lambda_1, \ldots, \lambda_n\}$ be the set of eigenvalues of $A$, repeated with multiplicities.
 \begin{enumerate}[(a)]
 \item $Q_A$ is positive definite if and only if $\lambda_i>0$ for all $1\leq i\leq n$.
  \item $Q_A$ is negative definite if and only if $\lambda_i<0$ for all $1\leq i\leq n$.
  \item $Q_A$ is indefinite if there exist $i$ and $j$ so that $\lambda_i>0$ and $\lambda_j<0$. 
  \item $Q_A$ is positive semi-definite if and only if $\lambda_i\geq 0$ for all $1\leq i\leq n$.
  \item $Q_A$ is negative semi-definite if and only if $\lambda_i\leq 0$ for all $1\leq i\leq n$.
 \end{enumerate}
 \end{theorem}
 From Theorem \ref{230724_8} and Theorem \ref{230724_20}, we obtain the following.
 \begin{corollary}[label=230725_5]{}
Let $Q:\mb{R}^n\to\mb{R}$ be a quadratic form. If $Q$ is positive definite, then there exists a positive constant $c$ such that
\[Q(\mf{x})\geq c\Vert\mf{x}\Vert^2\hspace{1cm}\text{for all}\;\mf{x}\in \mb{R}^n.\]
 \end{corollary}
In fact,   $c$ can be any positive number that is less than or equal to the smallest eigenvalue of the symmetric matrix $A$ associated to the quadratic form $Q$.

\section{Limits of Functions} 

 In this section, we study limits  of multivariable functions.
 
 \begin{definition} 
 {Limits of   Functions}
 Let $\mathfrak{D}$ be a subset of $\mathbb{R}^n$ and let $\mathbf{x}_0$ be a limit point of $\mathfrak{D}$. Given a function $\mathbf{F}:\mathfrak{D}\rightarrow \mathbb{R}^m$, we say that {\it the limit of $\mathbf{F}(\mathbf{x})$ as $\mathbf{x}$ approaches $\mathbf{x}_0$ is $\mathbf{v}$}, provided that whenever $\{\mathbf{x}_k\}$ is a sequence of points in $\mathfrak{D}\setminus\{\mathbf{x}_0\}$ that converges to $\mathbf{x}_0$, the sequence $\{\mathbf{F}(\mathbf{x}_k)\}$ of points in $\mathbb{R}^m$ converges to the point $\mathbf{v}$.
 
If the limit of $\mathbf{F}:\mathfrak{D}\rightarrow \mathbb{R}^m$ as $\mathbf{x}$ approaches $\mathbf{x}_0$ is $\mathbf{v}$, we write
 \[\lim_{\mathbf{x}\rightarrow \mathbf{x}_0}\mathbf{F}(\mathbf{x})=\mathbf{v}.\]
 \end{definition}
 \begin{example}[label=230719_4]{}
 For $1\leq i\leq n$, let $\pi_i:\mb{R}^n\to\mb{R}$ be the projection function $\pi_i(x_1, \ldots, x_n)=x_i$. By  the theorem on componentwise convergence of sequences, if $\{\mathbf{x}_k\}$ is a sequence in $\mb{R}^n\setminus\{\mf{x}_0\}$ that converges to the point $\mf{x}_0$, then 
 \[\lim_{k\to\infty}\pi_i(\mathbf{x}_k)=\pi_i(\mathbf{x}_0).\]
 This means that
 \[\lim_{\mf{x}\to\mf{x}_0}\pi_i(\mathbf{x})=\pi_i(\mf{x}_0).\]
 \end{example}
From the theorem on componentwise convergence of sequences, we also obtain the following immediately.
\begin{proposition}[label=230719_5]{}
Let $\mathfrak{D}$ be a subset of $\mathbb{R}^n$ and let $\mathbf{x}_0$ be a limit point of $\mathfrak{D}$. Given a function $\mathbf{F}:\mathfrak{D}\rightarrow \mathbb{R}^m$, \[\lim_{\mathbf{x}\rightarrow \mathbf{x}_0}\mathbf{F}(\mathbf{x})=\mathbf{v}\] if and only if for each $1\leq j\leq m$, \[\lim_{\mathbf{x}\rightarrow \mathbf{x}_0}F_j(\mathbf{x})=\pi_j(\mathbf{v}).\]
\end{proposition}

\begin{example}[label=230719_6]{}
Let $f:\mathbb{R}^n\to\mathbb{R}$ be the function defined as $f(\mathbf{x})=\Vert\mathbf{x}\Vert$. If $\mathbf{x}_0$ is a point in $\mathbb{R}^n$, find $\di\lim_{\mathbf{x}\to\mathbf{x}_0}f(\mathbf{x})$.
\end{example}
\begin{solution}{Solution}
We have shown in Example \ref{230717_1} that If $\{\mathbf{x}_k\}$ is a sequence in $\mathbb{R}^n\setminus\{\mathbf{x}_0\}$ that converges to $\mathbf{x}_0$, then
\[\lim_{k\to\infty}\Vert\mathbf{x}_k\Vert =\Vert\mathbf{x}_0\Vert.\]Therefore,
 $\di\lim_{\mathbf{x}\to\mathbf{x}_0}f(\mathbf{x})=\Vert\mathbf{x}_0\Vert$.
\end{solution}
By the limit laws for sequences, we also have the followings.
\begin{proposition}[label=230717_3]{}
Let $\mf{F}:\mk{D}\to \mb{R}^m$ and $\mf{G}:\mk{D}\to\mb{R}^m$ be functions defined on $\mk{D}\subset \mb{R}^n$. If $\mf{x}_0$ is a limit point of $\mk{D}$ and 
\[\lim_{\mf{x}\to\mf{x}_0}\mf{F}(\mf{x})=\mf{u},\hspace{1cm} \lim_{\mf{x}\to\mf{x}_0}\mf{G}(\mf{x})=\mf{v},\] then for any real numbers $\alpha$ and $\beta$,
\[\lim_{\mf{x}\to\mf{x}_0}(\alpha\mf{F}+\beta\mf{G})(\mf{x})=\alpha\mf{u}+\beta\mf{v}.\]
\end{proposition}
\begin{proposition}[label=230717_2]{}
Let $f:\mk{D}\to \mb{R}$ and $g:\mk{D}\to\mb{R}$ be functions defined on $\mk{D}\subset \mb{R}^n$. If $\mf{x}_0$ is a limit point of $\mk{D}$ and 
\[\lim_{\mf{x}\to\mf{x}_0}f(\mf{x})=u,\hspace{1cm} \lim_{\mf{x}\to\mf{x}_0}g(\mf{x})=v,\] then  
\[\lim_{\mf{x}\to\mf{x}_0}(fg)(\mf{x})=uv.\]If $g(\mathbf{x})\neq 0$ for all $\mf{x}\in\mk{D}$, and $v\neq 0$, then
\[\lim_{\mf{x}\to\mf{x}_0}\left(\frac{f}{g}\right)(\mathbf{x})=\frac{u}{v}.\]
\end{proposition}
\begin{example}[label=230719_9]{}
If $\mathbf{k}=(k_1, \ldots, k_n)$ is a $k$-tuple of nonnegative integers, the monomial $p_{\mf{k}}:\mb{R}^n\to\mb{R}$,
\[p_{\mathbf{k}}(\mf{x})=x_1^{k_1}\cdots x_n^{k_n}\] can be written as a product of the projection functions $\pi_i:\mb{R}^n\to\mb{R}$, $\pi_i(\mf{x})=x_i$, $1\leq i\leq n$. By Proposition \ref{230717_2}, 
\[\lim_{\mf{x}\to\mf{x}_0}p_{\mf{k}}(\mf{x})=p_{\mf{k}}(\mf{x}_0)\]for any $\mf{x}_0$ in $\mb{R}^n$. If $p:\mb{R}^n\to\mb{R}$ is a polynomial, it is  a finite linear combination of monomials.   Proposition \ref{230717_3} then implies that for any $\mathbf{x}_0$ in $\mathbb{R}^n$,
\[\lim_{\mf{x}\to\mf{x}_0}p(\mf{x})=p(\mf{x}_0).\]If $f:\mk{D}\to \mb{R}$, $f(\mf{x})=p(\mf{x})/q(\mf{x})$ is a rational function which is equal to the quotient of the polynomial $p(\mf{x})$  by the polynomial $q(\mf{x})$, then Proposition \ref{230717_2}  implies that
\[\lim_{\mf{x}\to\mf{x}_0}f(\mathbf{x})=f(\mf{x}_0)\] for any $\mf{x}_0\in \mk{D}=\left\{\mf{x}\in\mb{R}^n\,|\,q(\mf{x})\neq 0\right\}$.
\end{example}
\begin{example}[label=230806_1]{}
Find $\di\lim_{(x, y)\to (1, -1)}\frac{x^2+3xy+2y^2}{x^2+y^2}$.

\end{example}
\begin{solution}{Solution}
Since 
\begin{align*}\lim_{(x, y)\to (1, -1)}(x^2+3xy+2y^2)&=1-3+2=0,\\  \lim_{(x, y)\to (1, -1)}(x^2+y^2)&=1+1=2,\end{align*}we find that
\[\lim_{(x, y)\to (1, -1)}\frac{x^2+3xy+2y^2}{x^2+y^2}=\frac{0}{2}=0.\]
\end{solution}

It is easy to deduce the limit law for composite functions.
\begin{proposition}[label=230717_4]{}
Let $\mk{D}$ be a subset of $\mathbb{R}^n$, and let $\mathcal{U}$ be a subset of $\mathbb{R}^k$. Given the two functions $\mf{F}:\mk{D}\rightarrow \mathbb{R}^k$ and $\mf{G}: \mathcal{U}\rightarrow\mathbb{R}^m$, if $\mf{F}(\mk{D})\subset \mathcal{U}$, we can define the composite function $\mf{H}=\mf{G}\circ \mf{F}:\mk{D}\rightarrow \mathbb{R}^m$ by
$\mf{H}(\mf{x})=\mf{G}(\mf{F}(\mf{x}))$. If $\mf{x}_0$ is a limit point of $\mk{D}$, $\mf{y}_0$ is a limit point of $\mathcal{U}$, $\mf{F}(\mk{D}\setminus\{\mf{x}_0\})\subset \mathcal{U}\setminus\{\mf{y}_0\}$,
\[\lim_{\mf{x}\rightarrow \mf{x}_0}\mf{F}(\mf{x})=\mf{y}_0,\hspace{1cm}\lim_{\mf{y}\rightarrow \mf{y}_0}\mf{G}(\mf{y})=\mf{v},\] then
\[\lim_{\mf{x}\rightarrow \mf{x}_0}\mf{H}(\mf{x})=\lim_{\mf{x}\rightarrow \mf{x}_0}(\mf{G}\circ \mf{F})(\mf{x})=\mf{v}.\]
\end{proposition}The proof repeats verbatim the proof of the corresponding theorem for single variable functions.

\begin{figure}[ht]
\centering
\includegraphics[scale=0.18]{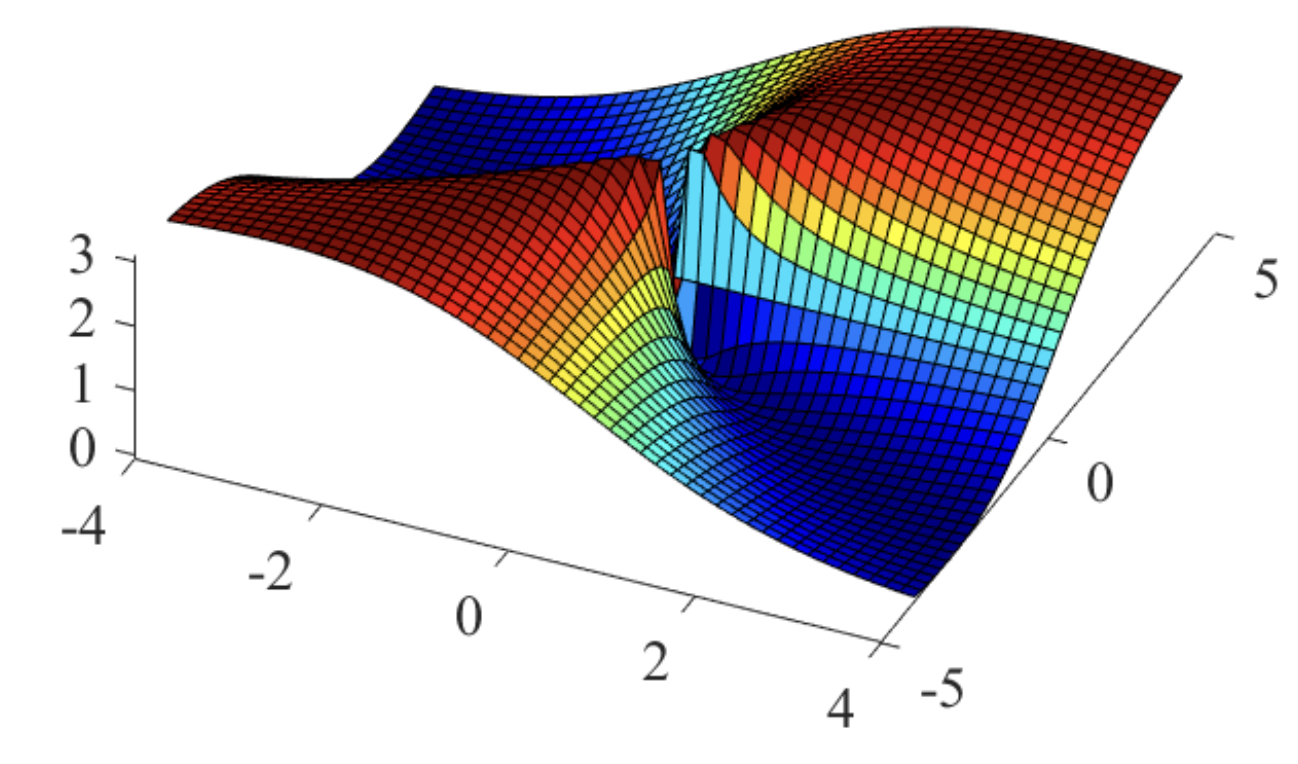}

\caption{The function  $f(x,y)=\di \frac{x^2+3xy+2y^2}{x^2+y^2}$ in Example \ref{230806_1}.}\label{figure48}
\end{figure}

\begin{example}[label=230806_2]{}
Find the limit $\di\lim_{(x,y)\to(0,0)}\frac{\sin (2x^2+3y^2)}{2x^2+3y^2}$.
\end{example}

\begin{solution}{Solution}
Since 
\[\lim_{(x,y)\to (0,0)}(2x^2+3y^2)=2\times 0+3\times 0=0,\hspace{1cm}\lim_{u\to 0}\frac{\sin u}{u}=1,\]
the limit law for composite functions implies that
\[\lim_{(x,y)\to(0,0)}\frac{\sin (2x^2+3y^2)}{2x^2+3y^2}=1.\]
\end{solution}
\begin{figure}[ht]
\centering
\includegraphics[scale=0.18]{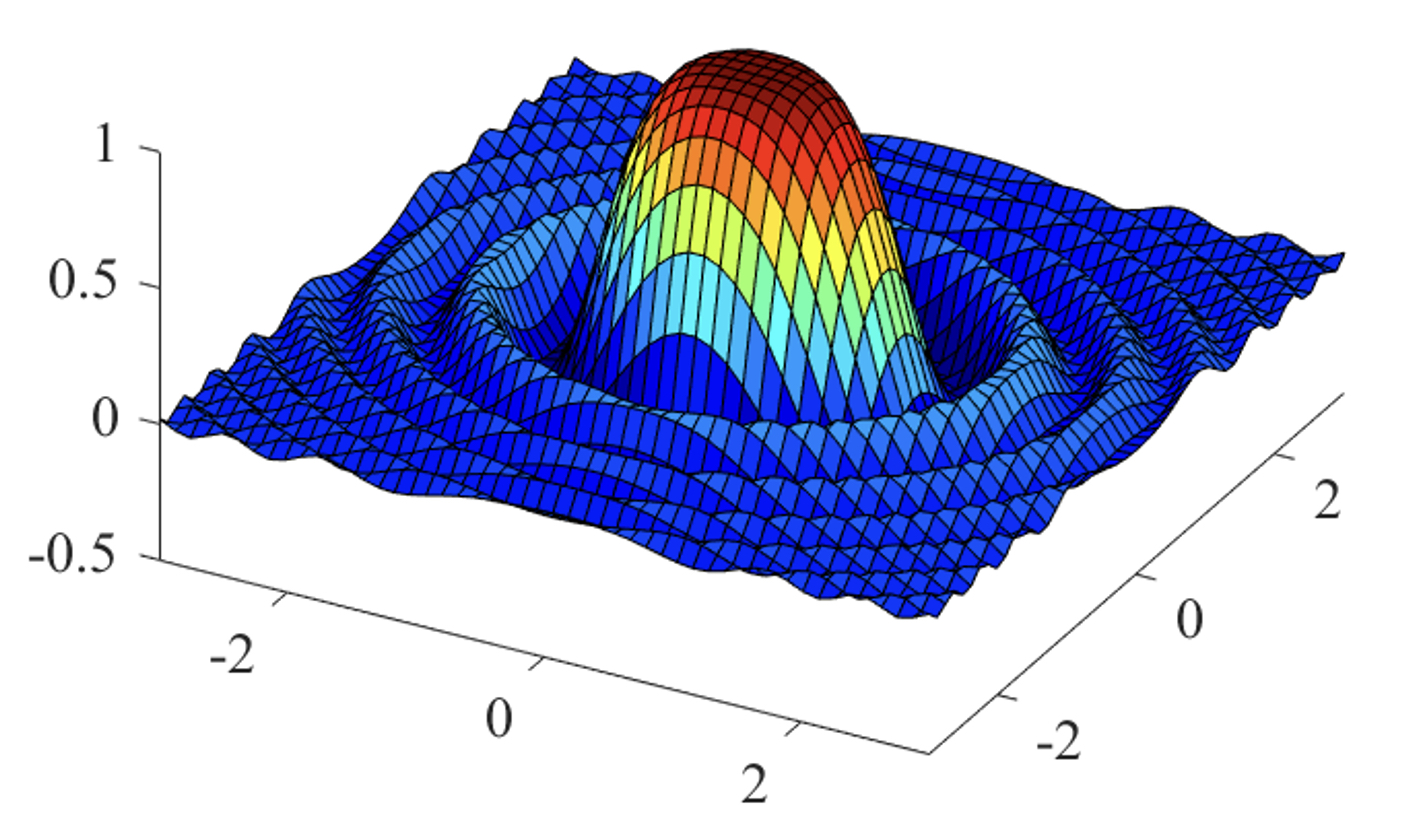}

\caption{The function  $f(x,y)=\di \frac{\sin (2x^2+3y^2)}{2x^2+3y^2}$ in Example \ref{230806_2}.}\label{figure49}
\end{figure}

Let us look at some examples where the rules we have studied cannot be applied.

\begin{example}[label=230806_3]{}
Determine whether the limit $\di\lim_{(x,y)\to (0,0)}\frac{x^2-2y^2}{x^2+y^2}$ exists.

\end{example}

\begin{solution}{Solution}
Let 
\[f(x,y)=\frac{x^2-2y^2}{x^2+y^2}=\frac{p(x,y)}{q(x,y)}.\]
When $(x,y)\to (0,0)$, $q(x,y)=x^2+y^2\to 0$. Hence, we cannot apply limit law  for quotients of functions. 

Consider the sequences of points $\{\mathbf{u}_k\}$ and $\{\mathbf{v}_k\}$ in $\mathbb{R}^2\setminus \{0,0\}$ given by
\[\mathbf{u}_k=\left(\frac{1}{k}, 0\right),\hspace{1cm}\mathbf{v}_k=\left(0, \frac{1}{k}\right).\]Notice that both the sequences   $\{\mathbf{u}_k\}$ and $\{\mathbf{v}_k\}$ converge to $(0,0)$. If $\di \lim_{(x,y)\to(0,0)}f(x,y)=a$, then both the sequences $\{f(\mathbf{u}_k)\}$ and $\{f(\mathbf{v}_k)\}$ should converge to $a$. Since
\[f(\mathbf{u}_k)=1,\hspace{1cm}f(\mathbf{v}_k)=-2\hspace{1cm}\text{for all}\;k\in\mathbb{Z}^+,\] 

 the sequence $\{f(\mathbf{u}_k)\}$ converges to 1, while the sequence $\{f(\mathbf{v}_k)\}$ converges to $-2$. These imply that $a=1$ and $a=-2$, which is a contradiction. Hence, the 
 limit $\di\lim_{(x,y)\to (0,0)}\frac{x^2-2y^2}{x^2+y^2}$ does not exist.
\end{solution}

\begin{figure}[ht]
\centering
\includegraphics[scale=0.18]{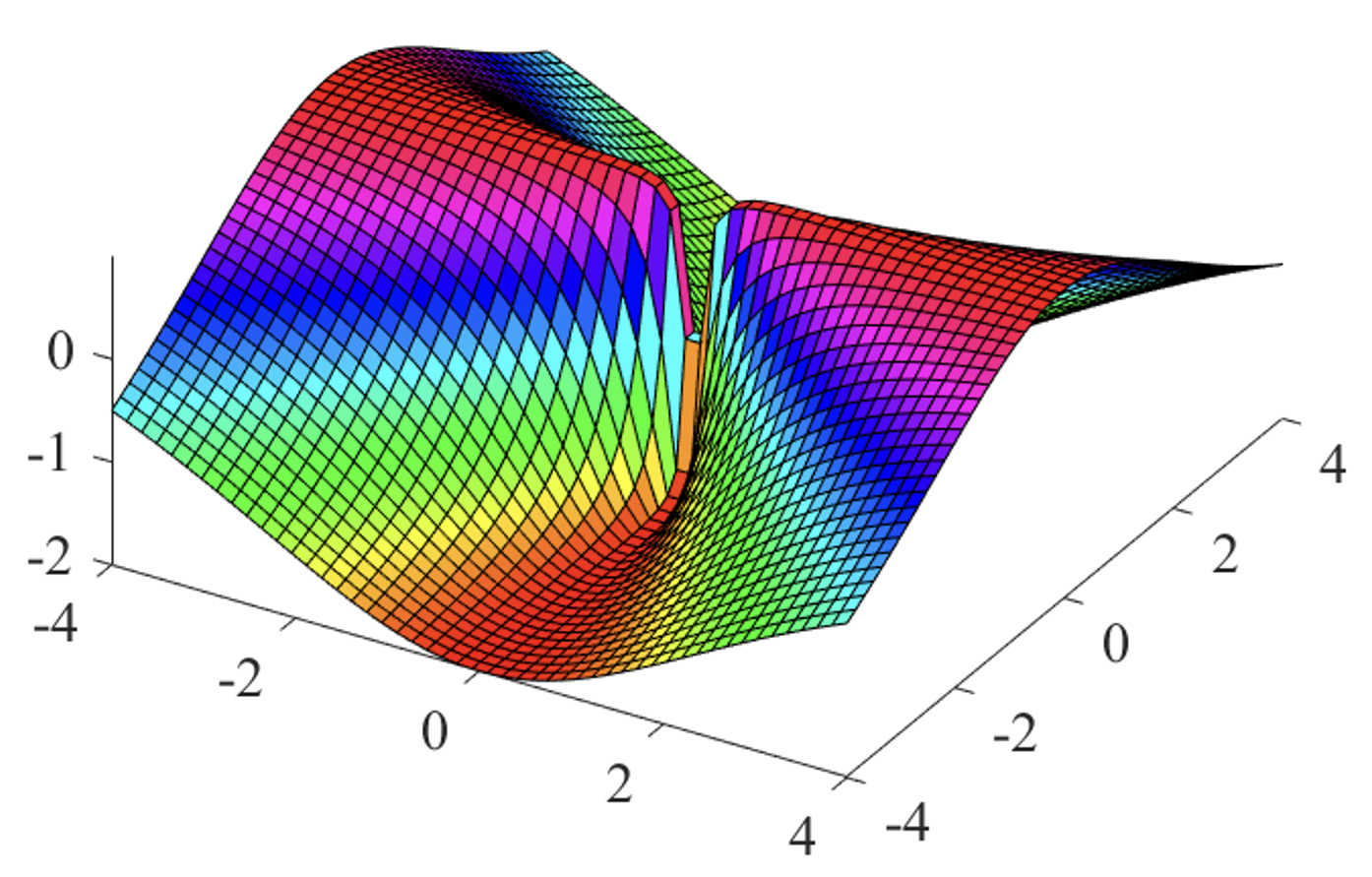}

\caption{The function  $f(x,y)=\di \frac{x^2-2y^2}{x^2+y^2}$ in Example \ref{230806_3}.}\label{figure50}
\end{figure}

\begin{example}[label=230806_4]{}
Determine whether the limit $\di\lim_{(x,y)\to (0,0)}\frac{xy}{x^2+2y^2}$ exists.

\end{example}

\begin{solution}{Solution}
Let 
\[f(x,y)=\frac{xy}{x^2+2y^2}.\]
Consider the sequences of points $\{\mathbf{u}_k\}$ and $\{\mathbf{v}_k\}$ in $\mathbb{R}^2\setminus \{0,0\}$ given by
\[\mathbf{u}_k=\left(\frac{1}{k}, 0\right),\hspace{1cm}\mathbf{v}_k=\left(\frac{1}{k}, \frac{1}{k}\right),\]Notice that both the sequences   $\{\mathbf{u}_k\}$ and $\{\mathbf{v}_k\}$ converge to $(0,0)$. If $\di \lim_{(x,y)\to(0,0)}f(x,y)=a$, then both the sequences $\{f(\mathbf{u}_k)\}$ and $\{f(\mathbf{v}_k)\}$ should converge to $a$. Since

\[f(\mathbf{u}_k)=0,\hspace{1cm}f(\mathbf{v}_k)=\frac{1}{3}\hspace{1cm}\text{for all}\;k\in\mathbb{Z}^+,\] 
 the sequence $\{f(\mathbf{u}_k)\}$ converges to 0, while the sequence $\{f(\mathbf{v}_k)\}$ converges to $1/3$. These imply that $a=0$ and $a=1/3$, which is a contradiction. Hence, the 
 limit $\di\lim_{(x,y)\to (0,0)}\frac{xy}{x^2+2y^2}$ does not exist.
\end{solution}
\begin{figure}[ht]
\centering
\includegraphics[scale=0.18]{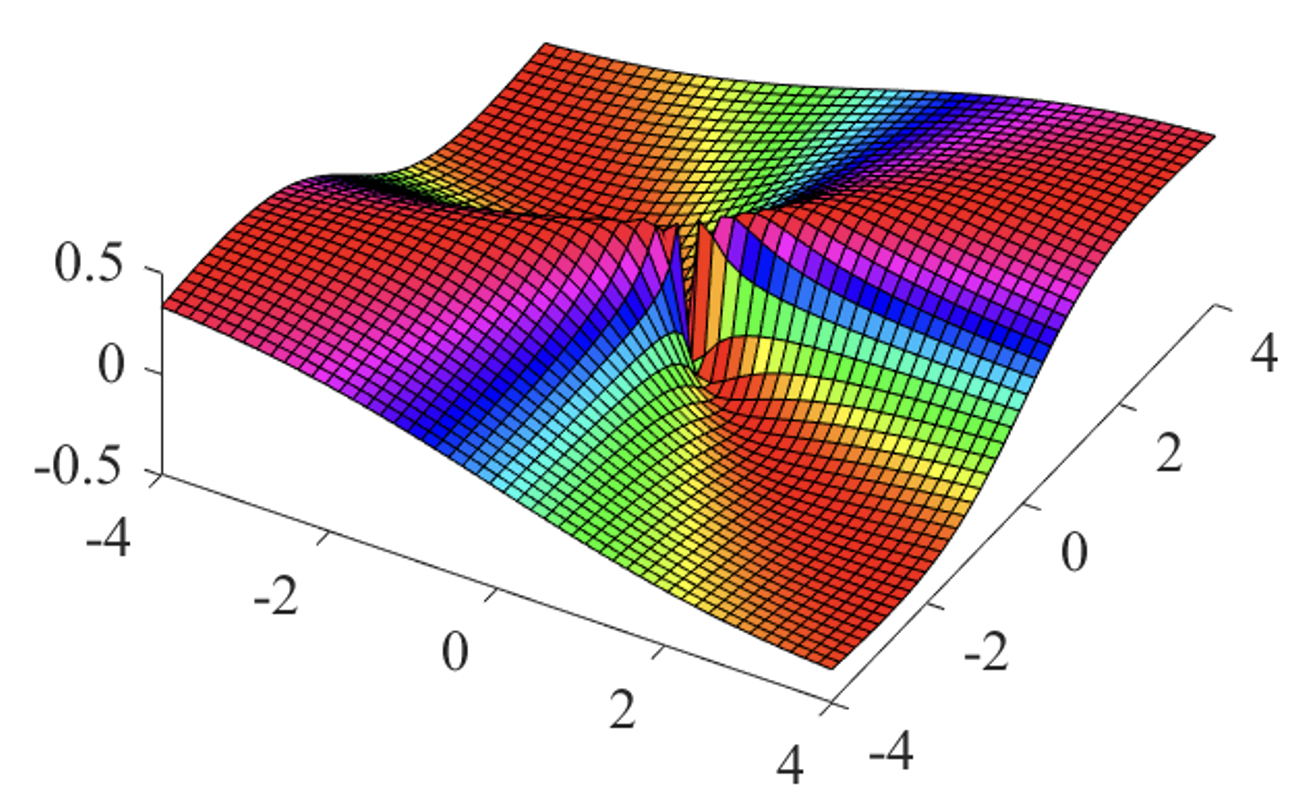}

\caption{The function  $f(x,y)=\di \frac{xy}{x^2+2y^2}$ in Example \ref{230806_4}.}\label{figure51}
\end{figure}

\begin{example}[label=230806_5]{}
Determine whether the limit $\di\lim_{(x,y)\to (0,0)}\frac{xy^2}{x^2+2y^4}$ exists.

\end{example}
\begin{solution}{Solution}
Let 
\[f(x,y)=\frac{xy^2}{x^2+2y^4}.\]
Consider the sequences of points $\{\mathbf{u}_k\}$ and $\{\mathbf{v}_k\}$ in $\mathbb{R}^2\setminus \{0,0\}$ given by
\[\mathbf{u}_k=\left(\frac{1}{k}, 0\right),\hspace{1cm}\mathbf{v}_k=\left(\frac{1}{k^2}, \frac{1}{k}\right),\]Notice that both the sequences   $\{\mathbf{u}_k\}$ and $\{\mathbf{v}_k\}$ converge to $(0,0)$. If $\di \lim_{(x,y)\to(0,0)}f(x,y)=a$, then both the sequences $\{f(\mathbf{u}_k)\}$ and $\{f(\mathbf{v}_k)\}$ should converge to $a$. Since
\[f(\mathbf{u}_k)=0,\hspace{1cm}f(\mathbf{v}_k)=\frac{1}{3}\hspace{1cm}\text{for all}\;k\in\mathbb{Z}^+,\]  the sequence $\{f(\mathbf{u}_k)\}$ converges to 0, while the sequence $\{f(\mathbf{v}_k)\}$ converges to $1/3$. These imply that $a=0$ and $a=1/3$, which is a contradiction. Hence, the 
 limit $\di\lim_{(x,y)\to (0,0)}\frac{xy^2}{x^2+2y^4}$ does not exist.
\end{solution}

\begin{figure}[ht]
\centering
\includegraphics[scale=0.18]{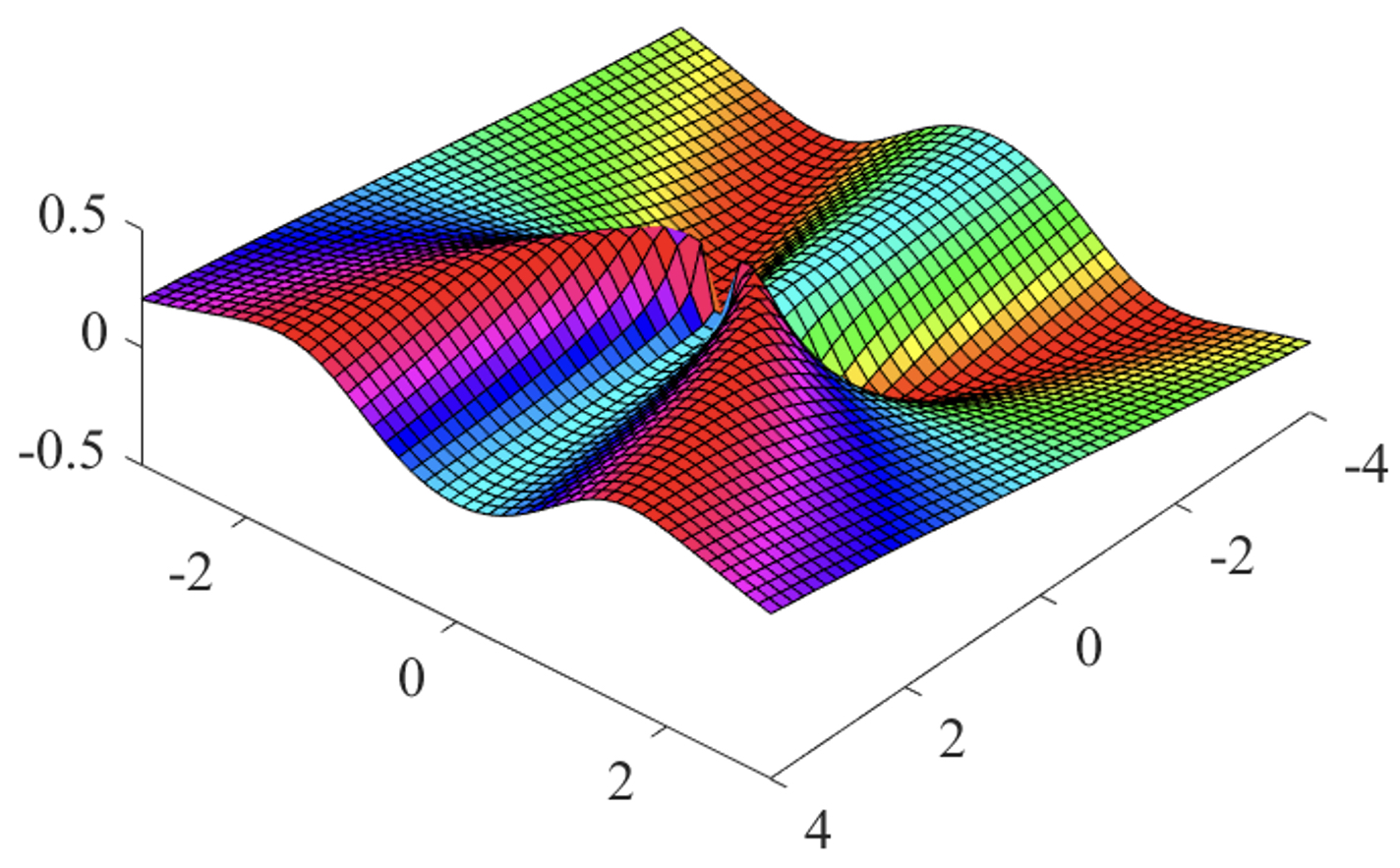}

\caption{The function  $f(x,y)=\di \frac{xy^2}{x^2+2y^4}$ in Example \ref{230806_5}.}\label{figure52}
\end{figure}

\begin{example}[label=230719_2]{}
Determine whether the limit $\di\lim_{(x,y)\to (0,0)}\frac{xy^2}{x^2+2y^2}$ exists.

\end{example}
\begin{solution}{Solution}
Let 
\[f(x,y)=\frac{xy^2}{x^2+2y^2}.\]
 If $\{(x_k, y_k)\}$ is a sequence of points in $\mathbb{R}^2\setminus\{0,0\}$ that converges to $(0,0)$, then
 \[\left|f(x_k,y_k)\right|=|x_k|\frac{y_k^2}{x_k^2+2y_k^2}\leq |x_k|.\]
 The sequence $\{x_k\}$ converges to 0. By squeeze theorem, the sequence $\{f(x_k, y_k)\}$ also converges to 0. This proves that 
 \[\lim_{(x,y)\to (0,0)}\frac{xy^2}{x^2+2y^2}=0.\]
\end{solution}

\begin{figure}[ht]
\centering
\includegraphics[scale=0.18]{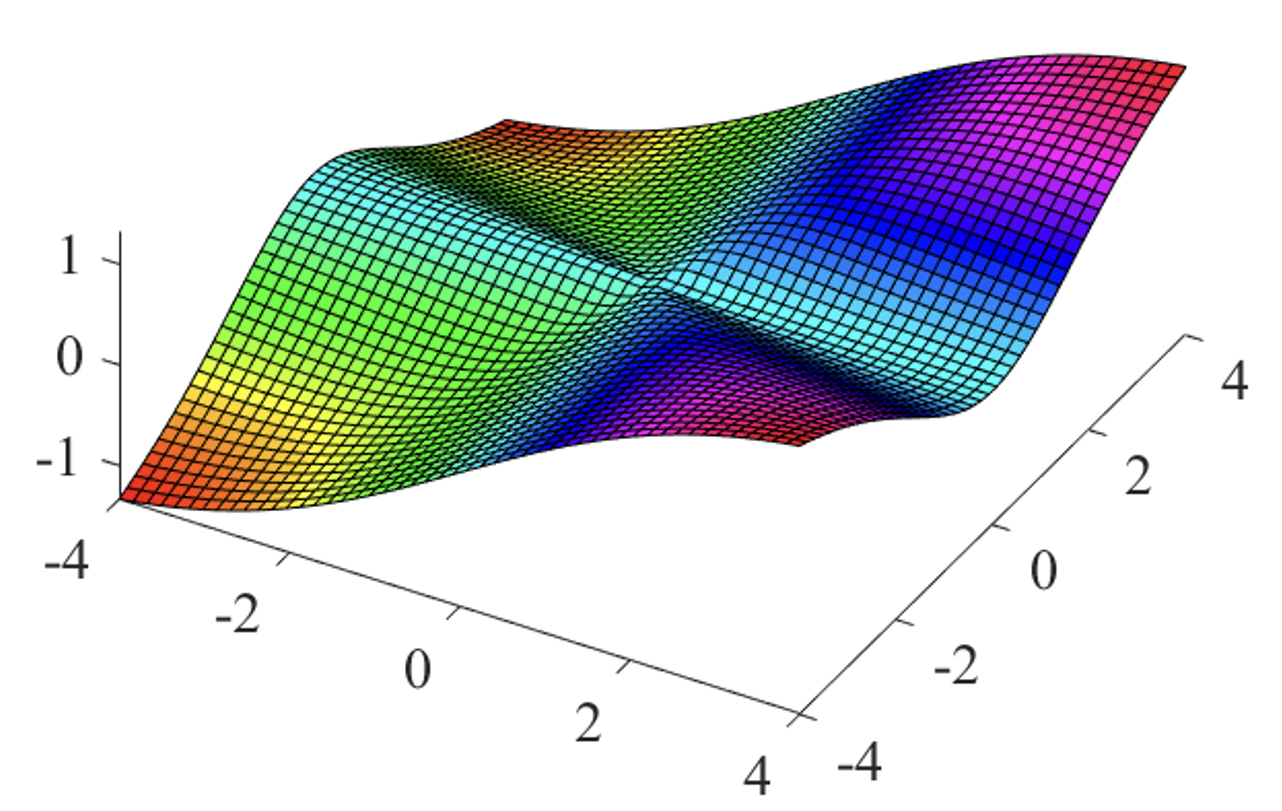}

\caption{The function  $f(x,y)=\di \frac{xy^2}{x^2+2y^2}$ in Example \ref{230719_2}.}\label{figure53}
\end{figure}

Similar to the single variable case, there is an equivalent definition of limits in terms of $\varepsilon$ and $\delta$. 

 \begin{theorem}[label=230719_1]{Equivalent Definitions for Limits}
  Let $\mk{D}$ be a subset of $\mb{R}^n$, and let $\mathbf{x}_0$ be a limit point of $\mk{D}$. Given a function $\mf{F}:\mk{D}\rightarrow \mathbb{R}^m$, 
  the following two definitions for 
  \[\lim_{\mf{x}\rightarrow \mf{x}_0}\mf{F}(\mf{x})=\mf{v}\] are equivalent.
  \begin{enumerate}[(i)]
  \item 
  Whenever $\{\mf{x}_k\}$ is a sequence of points in $\mk{D}\setminus\{\mf{x}_0\}$ that converges to $\mf{x}_0$, the sequence $\{\mf{F}(\mf{x}_k)\}$ converges to $\mf{v}$. 
  \item For any $\varepsilon>0$, there is a $\delta>0$ such that if the point $\mf{x}$ is in $\mk{D}$ and $0<\Vert\mf{x}-\mf{x}_0\Vert<\delta$, then $\Vert \mf{F}(\mf{x})-\mf{v}\Vert<\varepsilon$.
  \end{enumerate} 
 \end{theorem}
  \begin{myproof}{Proof}
We will prove that if (ii) holds, then (i) holds; and if (ii) does not hold, then (i) also does not hold.

First assume that (ii) holds. If $\{\mf{x}_k\}$ is a sequence   in $\mk{D}\setminus\{\mf{x}_0\}$ that converges to the point $\mf{x}_0$, we need to show that the sequence $\{\mf{F}(\mf{x}_k)\}$ converges to $\mf{v}$.  
  Given $\varepsilon>0$,  (ii)  implies that there is a $\delta>0$ such that for all $\mf{x}$ that is in $\mk{D}\setminus\{\mathbf{x}_0\}$ with $\Vert\mf{x}-\mf{x}_0\Vert<\delta$, we have $\Vert \mf{F}(\mf{x})-\mf{v}\Vert <\varepsilon$. 
  \bp
  Since $\{\mf{x}_k\}$ converges to $\mf{x}_0$, there is a positive integer $K$ such that for all $k\geq K$, $\Vert\mf{x}_k-\mf{x}_0\Vert<\delta$. Therefore, for all $k\geq K$,  
 $\di\Vert\mf{F}(\mf{x}_k)-\mf{v}\Vert<\varepsilon$. This shows that the sequence $\{\mf{F}(\mf{x}_k)\}$ indeed converges to $\mf{v}$.
 
Now assume that (ii) does not hold. Then there is an $\varepsilon>0$ such that for any $\delta>0$, there is a point $\mathbf{x}$ in $\mk{D}\setminus\{\mf{x}_0\}$ with $\Vert\mf{x}-\mf{x}_0\Vert<\delta$ but $\Vert\mf{F}(\mf{x}_k)-\mf{v}\Vert\geq \varepsilon$.
For this $\varepsilon>0$, we  construct a sequence $\{\mf{x}_k\}$ in $\mk{D}\setminus \{\mf{x}_0\}$ in the following way. For each positive integer $k$, there is a point $\mf{x}_k$ in $\mk{D}\setminus \{\mf{x}_0\}$ such that
$\Vert\mf{x}-\mf{x}_0\Vert<1/k$ but $\Vert\mf{F}(\mf{x}_k)-\mf{v}\Vert\geq\varepsilon$.
Then $\{\mf{x}_k\}$ is a sequence in $\mk{D}\setminus \{\mf{x}_0\}$ that satisfies
\[\Vert\mf{x}-\mf{x}_0\Vert<1/k\hspace{1cm}\text{for all}\;k\in\mathbb{Z}^+.\]
Hence, it converges to $\mf{x}_0$. Since $\Vert\mf{F}(\mf{x}_k)-\mf{v}\Vert\geq \varepsilon$ for all $k\in\mathbb{Z}^+$, the sequence $\{\mf{F}(\mf{x}_k)\}$ cannot converge to $\mf{v}$.  This proves that (i) does not hold.
 \end{myproof}
 
 We can give an alternative solution to Example \ref{230719_2} as follows.
 \begin{solution}{\linkt Alternative Solution to  Example \ref{230719_2}}
Let 
\[f(x,y)=\frac{xy^2}{x^2+2y^2}.\]
 Given $\varepsilon>0$, let $\delta=\varepsilon$. If $(x,y)$ is a point in $\mathbb{R}^2\setminus \{(0,0)\}$ such that 
 \[\sqrt{x^2+y^2}=\Vert(x,y)-(0,0)\Vert<\delta=\varepsilon,\] then $|x|<\varepsilon$. This implies that
 \[\left|f(x,y)-0\right|=|x|\frac{y^2}{x^2+2y^2}\leq |x|<\varepsilon.\]
 Hence,
 \[\lim_{(x,y)\to (0,0)}\frac{xy^2}{x^2+2y^2}=0.\]
\end{solution}

\vp
\noindent
{\bf \large Exercises  \thesection}
\setcounter{myquestion}{1}
\begin{question}{\themyquestion}
Determine whether the limit exists. If it exists, find the limit.
\begin{enumerate}[(a)]
\item
$\di \lim_{(x,y)\to (1, 2)}\frac{4x^2-y^2}{x^2+y^2}$
\item $\di \lim_{(x,y)\to (1, 2)}\sqrt{\frac{4x^2-y^2}{x^2+y^2}}$
\item $\di \lim_{(x,y)\to (1, 2)}\sqrt{\frac{4x^2+y^2}{x^2+y^2}}$

\end{enumerate}
\end{question}
\atc
\begin{question}{\themyquestion}
Determine whether the limit exists. If it exists, find the limit.
\begin{enumerate}[(a)]
\item $\di \lim_{(x,y)\to (0,0)}\frac{x^3+y^3}{x^2+y^2}$
\item $\di \lim_{(x,y)\to (0,0)}\frac{x^2+y^3}{x^2+y^2}$
\item $\di \lim_{(x,y)\to (0, 0)}\frac{e^{4x^2+y^2}-1}{4x^2+y^2 }$
\item $\di \lim_{(x,y)\to (0, 0)}\frac{e^{x^2+y^2}-1}{4x^2+y^2 }$
\end{enumerate}
\end{question}
\atc
\begin{question}{\themyquestion}
Determine whether the limit
\[\lim_{(x,y)\to (0,0)}\frac{x^2+4y^4}{4x^2+y^4}\] exists. If it exists, find the limit.

\end{question}
\atc
\begin{question}{\themyquestion}
Determine whether the limit
\[\lim_{(x,y)\to (1,1)}\frac{\cos(x^2+y^2-2)-1}{(x^2+y^2-2)^2}\] exists. If it exists, find the limit.

\end{question}
\atc
\begin{question}{\themyquestion}
Let $\mathbf{x}_0$ be a point in $\mathbb{R}^n$. Find the limit
$\di \lim_{\mf{x}\to\mathbf{x}_0} \frac{\mathbf{x}}{\Vert\mathbf{x}\Vert}$.
\end{question}

\atc
 \begin{question}{\themyquestion}
Let $\mk{D}$ be a subset of $\mathbb{R}^n$, and let 
$f:\mk{D}\to\mathbb{R}$ and $\mf{G}:\mk{D}\to\mb{R}^m$ be functions defined on $\mk{D}$. We can define the function $\mf{H}:\mk{D}\to\mb{R}^m$ by
\[\mf{H}(\mf{x})=f(\mf{x})\mf{G}(\mf{x})\hspace{1cm}\text{for all}\;\mathbf{x}\in\mk{D}.\]
If
$\mathbf{x}_0$ is a point in $\mk{D}$ and
\[\lim_{\mathbf{x}\to\mf{x}_0}f(\mathbf{x})=a,\hspace{1cm}\lim_{\mathbf{x}\to\mf{x}_0}\mf{G}(\mf{x})=\mf{v},\]show that
\[\lim_{\mf{x}\to\mf{x}_0}\mf{H}(\mf{x})=a\mf{v}.\]
\end{question}

\section{Continuity}\label{sec2.3}

The definition of continuity is  a direct generalization of the single variable case.
\begin{definition}{Continuity}
Let $\mk{D}$ be a subset of $\mb{R}^n$ that contains the point $\mf{x}_0$, and let $\mf{F}:\mk{D}\rightarrow\mathbb{R}^m$ be a function defined on $\mk{D}$. We say that the function $\mf{F}$ is {\bf continuous at } $\mf{x}_0$ provided that whenever $\{\mf{x}_k\}$ is a sequence of points in $\mk{D}$ that converges to $\mf{x}_0$, the sequence $\{\mf{F}(\mf{x}_k)\}$ converges to $\mf{F}(\mf{x}_0)$. 

We say that $\mf{F}:\mk{D}\rightarrow\mathbb{R}^m$ is a \textbf{continuous function} if it is continuous at every point of its domain $\mk{D}$.
\end{definition}
From the definition, we obtain the following immediately.
\begin{proposition}{Limits and Continuity}
Let $\mk{D}$ be a subset of $\mb{R}^n$ that contains the point $\mf{x}_0$, and let $\mf{F}:\mk{D}\rightarrow\mathbb{R}^m$ be a function defined on $\mk{D}$. 
\begin{enumerate}[1.]
\item
If $\mf{x}_0$ is an isolated point of $\mk{D}$, then $\mf{F}$ is continuous at $\mathbf{x}_0$.
\item If $\mf{x}_0$ is a limit point of $\mk{D}$,  then $\mf{F}$ is continuous at $\mathbf{x}_0$ if and only if
\[\lim_{\mf{x}\to\mf{x}_0}\mf{F}(\mf{x})=\mf{F}(\mf{x}_0).\]
\end{enumerate}
\end{proposition}
\begin{example}{}
Example \ref{230719_4} says that for each $1\leq i\leq n$,  the projection function $\pi_i:\mb{R}^n\to\mb{R}$, $\pi_i(\mf{x})=x_i$, is a continuous function.
\end{example}

\begin{example}{}
Example \ref{230719_6} says that the norm function    $f:\mb{R}^n\to\mb{R}$, $f(\mf{x})=\Vert\mathbf{x}\Vert$, is a continuous function.
\end{example}
From Proposition \ref{230719_5}, we have the following.
\begin{proposition}[label=230721_9]{}
Let $\mk{D}$ be a subset of $\mb{R}^n$ that contains the point $\mf{x}_0$, and let $\mf{F}:\mk{D}\rightarrow\mathbb{R}^m$ be a function defined on $\mk{D}$. The function $\mf{F}:\mk{D}\rightarrow\mathbb{R}^m$ is continuous at $\mathbf{x}_0$ if and only if each of the component functions $F_j=(\pi_j\circ\mf{F}):\mk{D}\rightarrow\mathbb{R}$, $1\leq j\leq m$, is continuous at $\mathbf{x}_0$.
\end{proposition}

\begin{example}{}
The function $\mf{F}:\mb{R}^3\to\mb{R}^2$,
\[\mf{F}(x,y,z)=(x, z),\] is a continuous function since each component function  is continuous.
\end{example}
Proposition \ref{230717_3} gives the following.
\begin{proposition}[label=230719_7]{}
Let $\mf{F}:\mk{D}\to \mb{R}^m$ and $\mf{G}:\mk{D}\to\mb{R}^m$ be functions defined on $\mk{D}\subset \mb{R}^n$, and let $\mf{x}_0$ be a point in $\mk{D}$. If $\mf{F}:\mk{D}\to \mb{R}^m$ and $\mf{G}:\mk{D}\to\mb{R}^m$ are continuous at $\mf{x}_0$, then for any real numbers $\alpha$ and $\beta$,
the function $(\alpha\mf{F}+\beta\mf{G}):\mk{D}\to\mb{R}^m$ is continuous at $\mf{x}_0$.
\end{proposition}

Proposition \ref{230717_2} gives the following.
\begin{proposition}[label=230719_8]{}
Let $f:\mk{D}\to \mb{R}$ and $g:\mk{D}\to\mb{R}$ be functions defined on $\mk{D}\subset \mb{R}^n$, and let $\mf{x}_0$ be a point in $\mk{D}$. Assume that the functions $f:\mk{D}\to \mb{R}$ and $g:\mk{D}\to\mb{R}$ are continuous at $\mf{x}_0$. 
\begin{enumerate}[1.]
\item
The function $(fg):\mk{D}\to\mb{R}$ is continuous at $\mf{x}_0$. 
\item If $g(\mf{x})\neq 0$ for all $\mf{x}\in\mk{D}$, then the function $(f/g):\mk{D}\to\mb{R}$ is continuous at $\mf{x}_0$. 
\end{enumerate}
\end{proposition}
Example \ref{230719_9} gives the following.
\begin{proposition}{}
Polynomials and rational functions are continuous functions.
\end{proposition}
Since each component of a linear transformation $\mf{T}:\mb{R}^n\to\mb{R}^m$ is a polynomial, we have the following.
\begin{proposition}{}
A linear transformation  $\mf{T}:\mb{R}^n\to\mb{R}^m$  is a continuous function.
\end{proposition}
Since a quadratic form $Q:\mb{R}^n\to\mb{R}$ is a polynomial, we have the following.
\begin{proposition}{}
A quadratic form $Q:\mb{R}^n\to\mb{R}$ given by \[Q(\mf{x})=\di\sum_{i=1}^n\sum_{j=1}^n a_{ij}x_ix_j \] is a continuous function.
\end{proposition}

The following is obvious from the definition of continuity.
\begin{proposition}[label=230719_10]{}
Let $\mathfrak{D}$ be a subset of $\mb{R}^n$, and let $\mf{F}:\mk{D}\to \mb{R}^m$ be a function that is continuous at the point $\mf{x}_0\in\mk{D}$. If $\mk{D}_1$ is a subset of $\mk{D}$ that contains $\mf{x}_0$, then the function $\mf{F}:\mk{D}_1\to\mb{R}^m$ is also continuous at $\mathbf{x}_0$.
\end{proposition}
\begin{example}[label=230806_6]{}
Let $\mk{D}$ be the set 
\[\mk{D}=\left\{(x,y)\,|\, x^2+y^2<1\right\},\] and let $f:\mk{D}\to\mb{R}$ be the function defined as
\[f(x, y)=\frac{xy}{1-x^2-y^2}.\]\be
Since $f_1(x,y)=xy$ and $f_2(x,y)=1-x^2-y^2$ are polynomials, they are continuous. Since $f_2(x,y)\neq 0$ for all $(x,y)\in\mk{D}$, $f:\mk{D}\to\mb{R}$ is a continuous function.
\end{example2}

\begin{figure}[ht]
\centering
\includegraphics[scale=0.18]{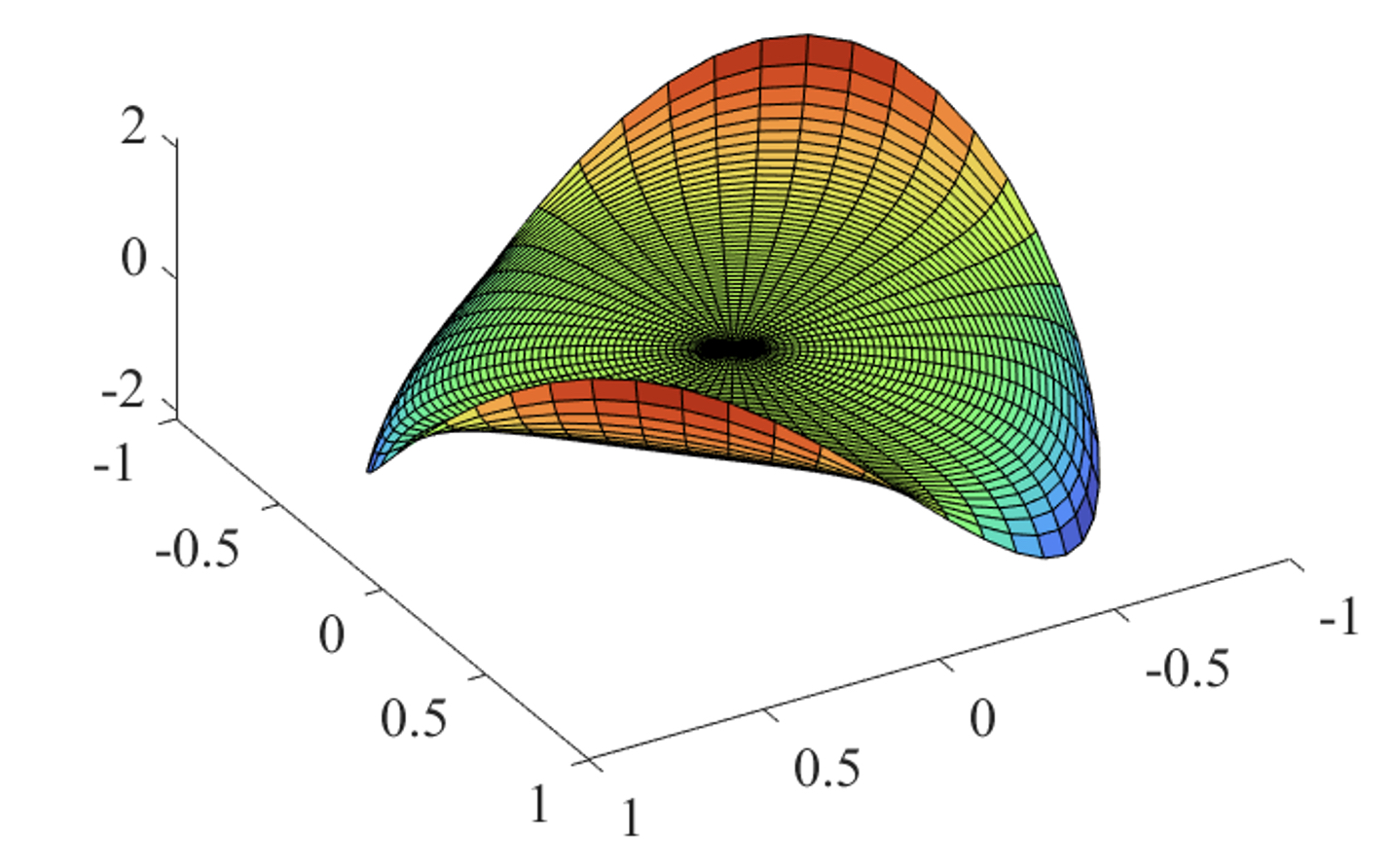}

\caption{The function  $f(x,y)=\di \frac{xy}{1-x^2-y^2}$ in Example \ref{230806_6}.}\label{figure54}
\end{figure}

Proposition \ref{230717_4} implies the following.
 \begin{proposition}[label=230719_11]{}
Let $\mk{D}$ be a subset of $\mathbb{R}^n$, and let $\mathcal{U}$ be a subset of $\mathbb{R}^k$. If $\mf{F}:\mk{D}\rightarrow \mathbb{R}^k$ and $\mf{G}: \mathcal{U}\rightarrow\mathbb{R}^m$ are functions such that $\mf{F}(\mk{D})\subset \mathcal{U}$,    $\mathbf{F}:\mk{D}\rightarrow \mathbb{R}^k$ is continuous at $\mf{x}_0$, $\mf{G}: \mathcal{U}\rightarrow\mathbb{R}^m$ is continuous at $\mathbf{y}_0$, then the composite function $\mf{H}=(\mf{G}\circ\mf{F}):\mk{D}\to\mb{R}^m$ is continuous at $\mf{x}_0$.
\end{proposition}

A direct proof of this theorem using the definition of continuity is actually much simpler. 
\begin{myproof}{Proof}
If $\{\mathbf{x}_k\}$ is a sequence of points in $\mk{D}$ that converges to $\mf{x}_0$, then since  $\mathbf{F}:\mk{D}\rightarrow \mathbb{R}^k$ is continuous at $\mf{x}_0$, $\{\mf{F}(\mf{x}_k)\}$ is a sequence of points in $\mathcal{U}$ that converges to $\mf{y}_0$. Since  $\mf{G}: \mathcal{U}\rightarrow\mathbb{R}^m$ is continuous at $\mathbf{y}_0$, $\{\mf{G}(\mf{F}(\mf{x}_k))\}$ is a sequence of points in $\mb{R}^m$ that converges to $\mf{G}(\mf{y}_0)=\mf{G}(\mf{F}(\mf{x}_0))$.
\bp
 In other words, the sequence $\{\mf{H}(\mf{x}_k)\}$ converges to $\mf{H}(\mf{x}_0)$. This shows that the function  $\mf{H}=(\mf{G}\circ\mf{F}):\mk{D}\to\mb{R}^m$ is continuous at $\mf{x}_0$.
\end{myproof}
\begin{figure}[ht]
\centering
\includegraphics[scale=0.18]{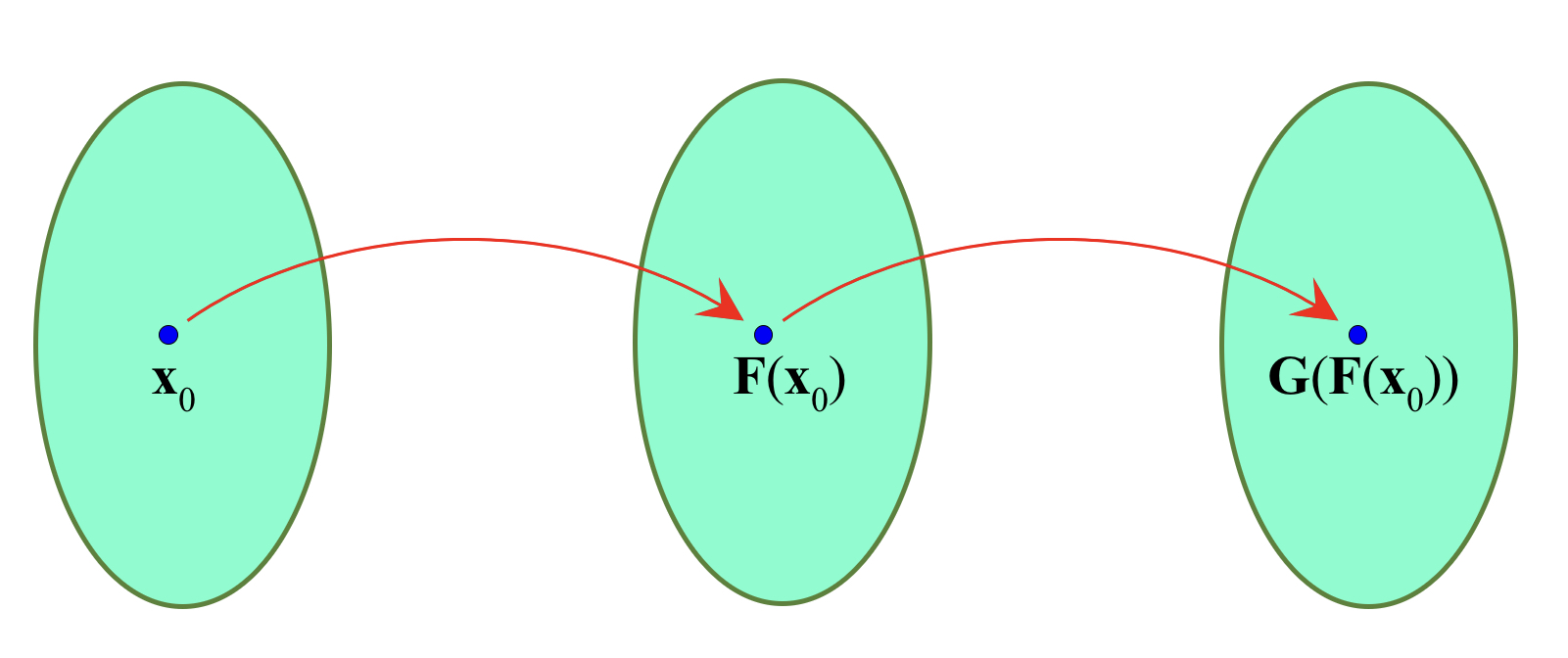}

\caption{Composition of functions.}\label{figure20}
\end{figure}
\begin{corollary}{}
Let $\mk{D}$ be a subset of $\mathbb{R}^n$, and let $\mathbf{x}_0$ be a point in $\mk{D}$. If the function $\mf{F}:\mk{D}\to\mb{R}^m$ is   continuous at $\mathbf{x}_0\in\mk{D}$, then the function $\Vert\mf{F}\Vert:\mk{D}\to\mb{R}$ is also continuous at $\mf{x}_0$. 
\end{corollary}

\begin{figure}[ht]
\centering
\includegraphics[scale=0.18]{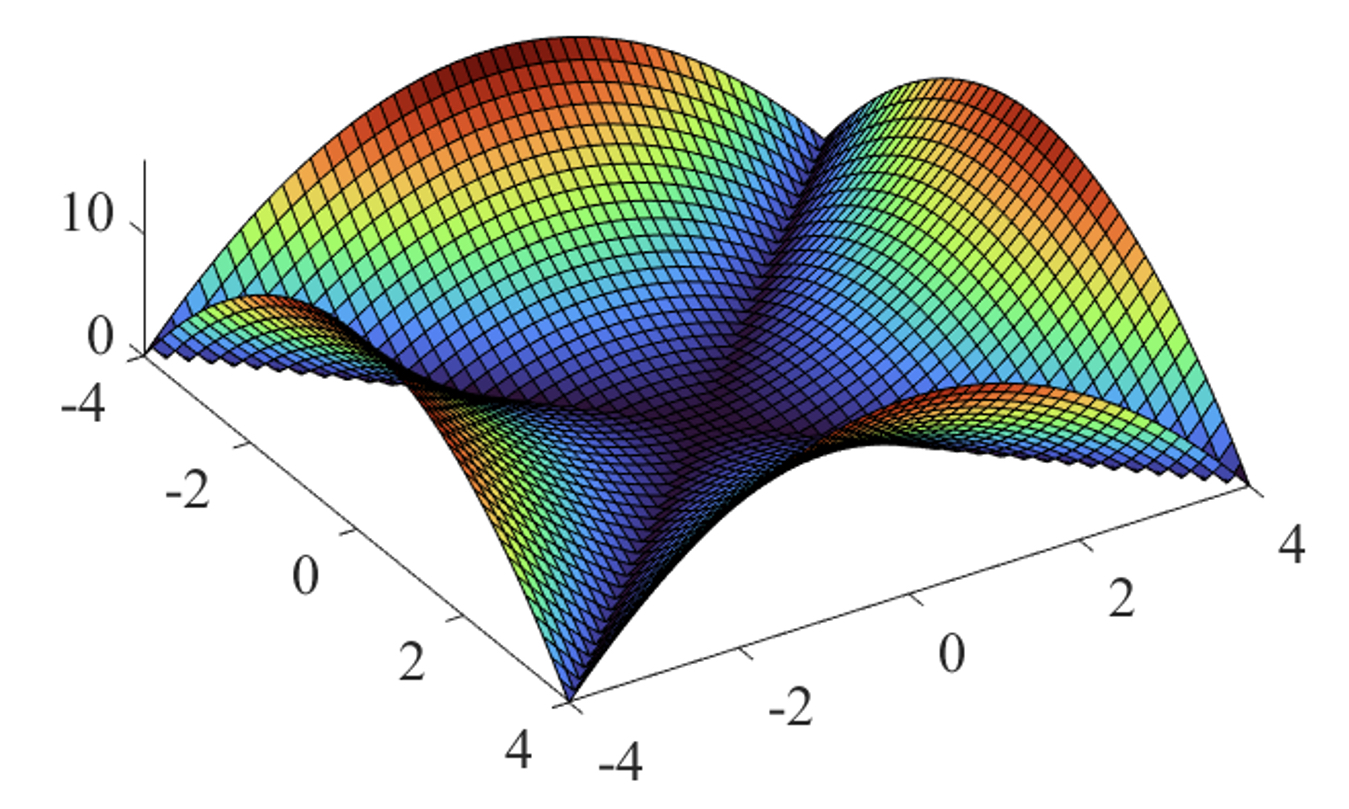}

\caption{The function $f(x,y)=|x^2-y^2|$.}\label{figure55}
\end{figure}

\begin{example}{}
The function $f:\mb{R}^2\to \mb{R}$, $f(x,y)=|x^2-y^2|$ is a continuous function since $f(x,y)=|p(x,y)|$, where $p(x,y)=x^2-y^2$ is a polynomial function, which is continuous.
\end{example}
\begin{example}{}
Consider the function $f:\mb{R}^2\to \mb{R}$, $f(x,y)=\di\sqrt{e^{2xy}+x^2+y^2}$.  Notice that $f(x,y)=\Vert\mf{F}(x,y)\Vert$, where $\mf{F}:\mb{R}^2\to\mb{R}^3$ is the function given by
\[\mf{F}(x,y)=\left(e^{xy},x,y\right).\] Since $g(x,y)=xy$ is a polynomial function, it is continuous. Being a composition of the continuous function $h(x)=e^x$ with the continuous function $g(x,y)=xy$, $F_1(x,y)=(h\circ g)(x,y)=e^{xy}$ is a continuous function. The functions $F_2(x,y)=x$ and $F_3(x,y)=y$ are continuous functions. Hence, $\mathbf{F}:\mb{R}^2\to\mb{R}^3$ is a continuous function. This implies that $f:\mb{R}^2\to \mb{R}$ is also a continuous function.
\end{example}

\begin{figure}[ht]
\centering
\includegraphics[scale=0.18]{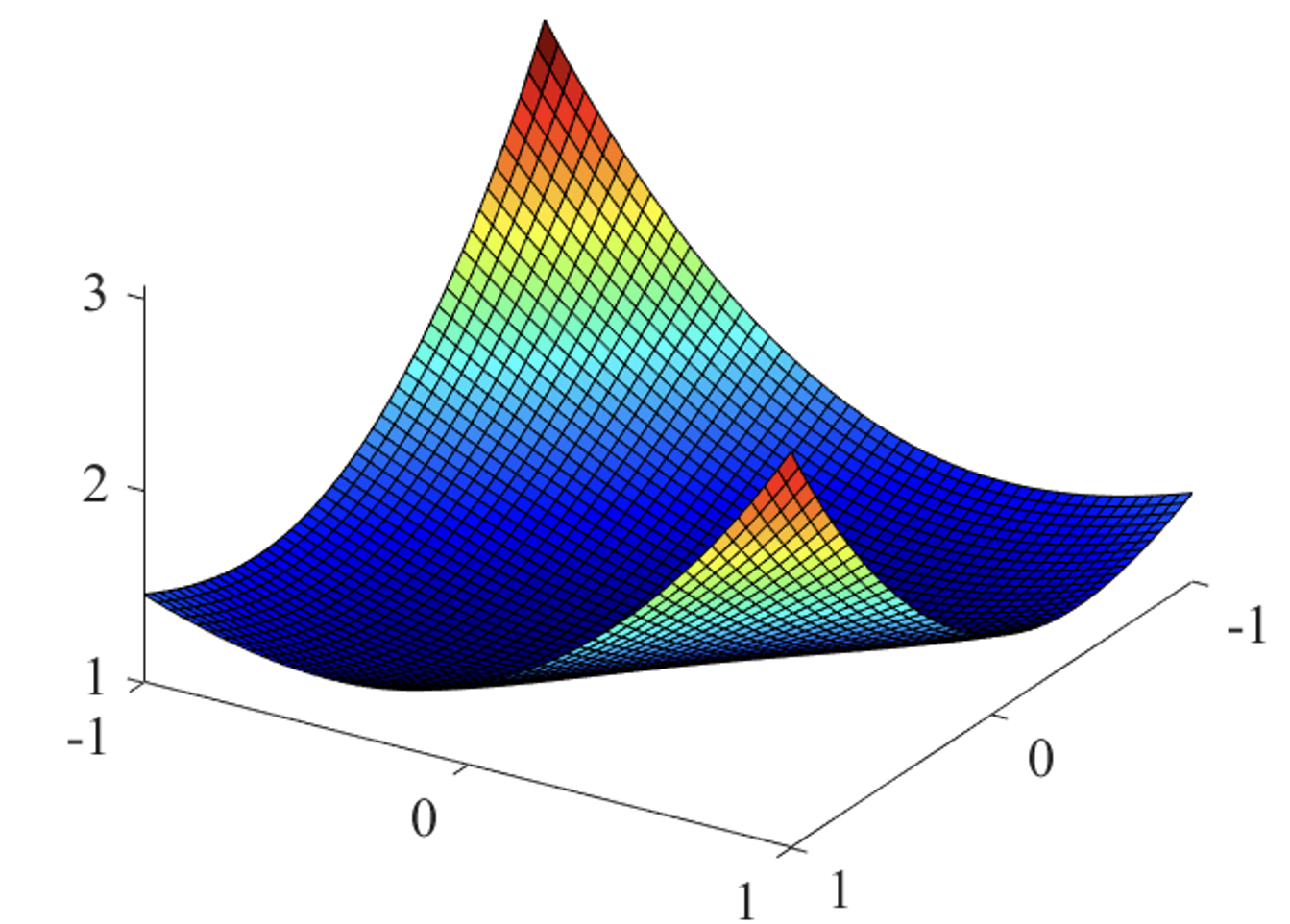}

\caption{The function $f(x,y)=\di\sqrt{e^{2xy}+x^2+y^2}$.}\label{figure56}
\end{figure}
\begin{example}{}
We have shown in volume I that the function $f:\mathbb{R}\to\mathbb{R}$,
\[f(x)=\begin{cases}\di \frac{\sin x}{x},\quad &\text{if}\;x\neq 0,\\1,\quad &\text{if}\;x=0,\end{cases}\] is a continuous function. Define the function $h:\mb{R}^3\to\mb{R}$ by
\[h(x,y,z)=\begin{cases}\di \frac{\sin (x^2+y^2+z^2)}{x^2+y^2+z^2},\quad &\text{if}\;(x,y,z)\neq (0,0,0),\\1,\quad &\text{if}\;(x,y,z)=(0,0,0).\end{cases}\] 
 
Since $h=f\circ g$, where $g:\mb{R}^3\to\mb{R}$ is the polynomial function $g(x,y,z)=x^2+y^2+z^2$, which is continuous, the function $h:\mb{R}^3\to\mb{R}$ is continuous.
\end{example}

The following gives an equivalent definition of continuity in terms of $\varepsilon$ and $\delta$.

 \begin{theorem}{Equivalent Definitions of Continuity}
  Let $\mk{D}$ be a subset of $\mb{R}^n$, and let $\mathbf{x}_0$ be a limit point of $\mk{D}$. Given a function $\mf{F}:\mk{D}\rightarrow \mathbb{R}^m$, 
  the following two definitions for the continuity of $\mathbf{F}$ at $\mathbf{x}_0$ are equivalent.
  \begin{enumerate}[(i)]
  \item 
  Whenever $\{\mf{x}_k\}$ is a sequence of points in $\mk{D}$ that converges to $\mf{x}_0$, the sequence $\{\mf{F}(\mf{x}_k)\}$ converges to $\mf{F}(\mf{x}_0)$. 
  \item For any $\varepsilon>0$, there is a $\delta>0$ such that if the point $\mf{x}$ is in $\mk{D}$ and $ \Vert\mf{x}-\mf{x}_0\Vert<\delta$, then $\Vert \mf{F}(\mf{x})-\mf{F}(\mf{x}_0)\Vert<\varepsilon$.
  \end{enumerate} 
 \end{theorem}  
 The proof is left as an exercise. Notice that   statement (ii) can be reformulated as follows. For any $\varepsilon>0$, there is a $\delta>0$ such that  if the point $\mf{x}$ is in $\mk{D}$ and $\mathbf{x}\in B(\mf{x}_0,\delta)$, then $\mf{F}(\mf{x})\in B(\mf{F}(\mf{x}_0), \varepsilon)$.

 \begin{figure}[ht]
\centering
\includegraphics[scale=0.2]{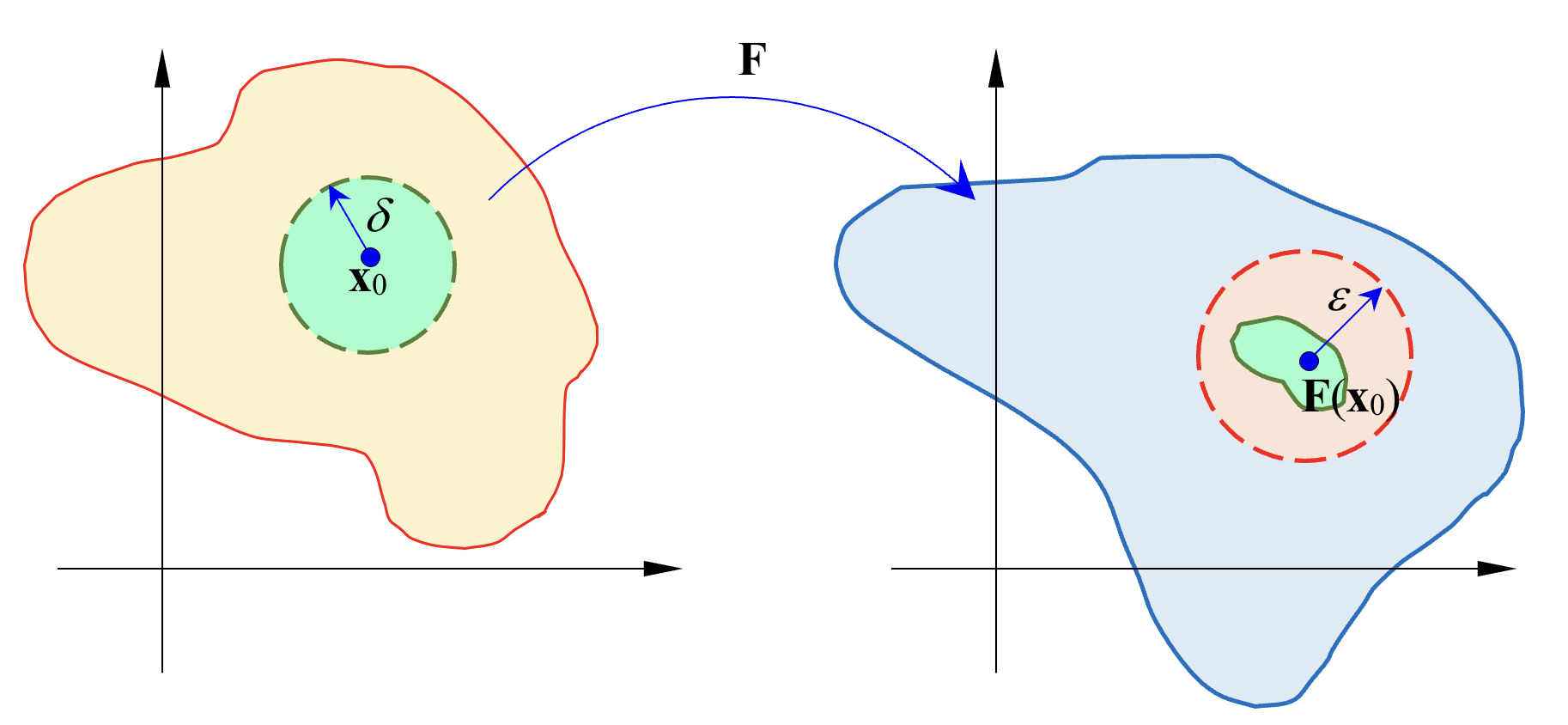}

\caption{The definition of continuity in terms of  $\varepsilon$ and $\delta$.}\label{figure17}
\end{figure}

Now we want to explore another important property of continuity.
\begin{theorem}[label=230720_5]{}
Let $\mathcal{O}$ be an open subset of $\mathbb{R}^n$, and let   $\mf{F}:\mathcal{O}\to\mathbb{R}^m$ be a function defined on $\mathcal{O}$. The following  are equivalent.
\begin{enumerate}[(a)]
\item  $\mf{F}:\mathcal{O}\to\mathbb{R}^m$ is continuous.
\item For every open subset $V$ of $\mb{R}^m$, $\mf{F}^{-1}(V)$ is an open subset of $\mb{R}^n$.
\end{enumerate}
\end{theorem}Note that for this theorem to hold, it is important that the domain of the function $\mf{F}$ is an open set.
\begin{myproof}{Proof}
Assume that (a) holds. Let $V$ be an open subset of $\mathbb{R}^m$, and let
\[U=\mf{F}^{-1}(V)=\left\{\mf{x}\in\mathcal{O}\,|\,\mf{F}(\mf{x})\in V\right\}.\]We need to show that $U$ is an open subset of $\mb{R}^n$. If $\mathbf{x}_0$ is in $U$, then it is in $\mathcal{O}$. Since $\mathcal{O}$ is open, there exists $r_0>0$ such that $B(\mf{x}_0, r_0)\subset\mathcal{O}$. Since $\mf{y}_0=\mf{F}(\mf{x}_0)$ is in $V$ and $V$ is open, there exists $\varepsilon>0$ such that $B(\mf{y}_0, \varepsilon)\subset V$. By (a), there exists $\delta>0$ such that for any $\mathbf{x}\in\mathcal{O}$, if $\Vert\mf{x}-\mf{x}_0\Vert<\delta$, then $\Vert \mf{F}(\mf{x})-\mf{F}(\mf{x}_0)\Vert<\varepsilon$. 
\bp
Take $r=\min\{\delta, r_0\}$. Then $r>0$, $r\leq r_0$ and $r\leq \delta$. If $\mathbf{x}$ is in $B(\mf{x}_0,r )$, then $\mf{x}\in \mathcal{O}$ and $\Vert\mf{x}-\mf{x}_0\Vert<r\leq \delta$. It follows that $\Vert \mf{F}(\mf{x})-\mf{F}(\mf{x}_0)\Vert<\varepsilon$. This implies that $\mf{F}(\mf{x})\subset B(\mf{y}_0, \varepsilon)\subset V$.
 Thus, $\mathbf{x}\in U$. In other words, we have shown that $B(\mf{x}_0,r )$ is contained in $U$. This proves that $U$ is open, which is the assertion of (b).

Conversely, assume that (b) holds. Let $\mathbf{x}_0$ be a point in $\mathcal{O}$, and let $\mathbf{y}_0=\mf{F}(\mf{x}_0)$. Given $\varepsilon>0$, the ball $V=B(\mathbf{y}_0,\varepsilon)$ is an open subset of $\mathbb{R}^m$. By (b), $U=\mf{F}^{-1}(V)$ is open in $\mb{R}^n$. By definition, $U$ is a subset of $\mathcal{O}$. Since $\mf{F}(\mf{x}_0)$ is in $V$, $\mf{x}_0$ is in $U$. Since  $U$ is open and it contains $\mf{x}_0$, there is an $r>0$ such that $B(\mf{x}_0, r)\subset U$. Take $\delta=r$. Then if $\mathbf{x}$ is a point in $\mathcal{O}$ and $\Vert\mf{x}-\mf{x}_0\Vert<r$, $\mf{x}\in B(\mf{x}_0, r)\subset U$. This implies that $\mf{F}(\mf{x})\in V=B(\mathbf{y}_0,\varepsilon)$. Namely, $\Vert\mf{F}(\mf{x})-\mf{F}(\mf{x}_0)\Vert<\varepsilon$. This proves that $\mf{F}:\mathcal{O}\to\mathbb{R}^m$ is continuous at $\mf{x}_0$. Since $\mf{x}_0$ is an arbitrary point in $\mathcal{O}$, $\mf{F}:\mathcal{O}\to\mathbb{R}^m$ is continuous.
\end{myproof}

Using the fact that a set is open if and only if its complement  is  closed, it is natural to expect the following.
\begin{theorem}[label=230720_6]{}
Let $\mathcal{A}$ be a closed subset of $\mathbb{R}^n$, and let   $\mf{F}:\mathcal{A}\to\mathbb{R}^m$ be a function defined on $\mathcal{A}$. The following are equivalent.
\begin{enumerate}[(a)]
\item  $\mf{F}:\mathcal{A}\to\mathbb{R}^m$ is continuous.
\item For every closed subset $C$ of $\mb{R}^m$, $\mf{F}^{-1}(C)$ is a closed subset of $\mb{R}^n$.
\end{enumerate}
\end{theorem} 
\begin{myproof}{Proof}
Assume that (a) holds. Let $C$ be a closed subset of $\mathbb{R}^m$, and let
\[D=\mf{F}^{-1}(C)=\left\{\mf{x}\in\mathcal{A}\,|\,\mf{F}(\mf{x})\in D\right\}.\]We need to show that $D$ is a closed subset of $\mb{R}^n$. If $\{\mf{x}_k\}$ is a sequence in $D$ that converges to the point $\mathbf{x}_0$ in $\mathbb{R}^n$, since $D\subset \mathcal{A}$ and $\mathcal{A}$ is closed, $\mathbf{x}_0$ is in $\mathcal{A}$. Since $\mf{F}$ is continuous at $\mf{x}_0$, the sequence $\{\mf{F}(\mathbf{x}_k)\}$ is a sequence in $C$ that converges to the point $\mf{F}(\mf{x}_0)$ in $\mb{R}^m$. Since $C$ is closed, $\mf{F}(\mf{x}_0)$ is in $C$. Therefore, $\mf{x}_0$ is in $D$. This proves that $D$ is closed.
 \bp
 
 Conversely, assume that (a) does not hold. Then $\mf{F}:\mathcal{A}\to\mb{R}^m$  is not continuous at some $\mf{x}_0\in\mathcal{A}$. Thus,  there exists $\varepsilon>0$ such that for any $\delta>0$, there exists a point $\mf{x}$ in $\mathcal{A}\cap B(\mf{x}_0, \delta)$ such that $\Vert\mf{F}(\mf{x})-\mf{F}(\mf{x}_0)\Vert\geq\varepsilon$. For $k\in\mb{Z}^+$, let $\mf{x}_k$ be a point in $\mathcal{A}\cap B(\mf{x}_0, 1/k)$ such that $\Vert\mf{F}(\mf{x}_k)-\mf{F}(\mf{x}_0)\Vert\geq\varepsilon$. Since 
 \[\Vert\mf{x}_k-\mf{x}_0\Vert<\frac{1}{k}\hspace{1cm}\text{for all}\;k\in\mb{Z}^+,\]
 the sequence $\{\mf{x}_k\}$ is a sequence in $\mathcal{A}$ that converges to $\mf{x}_0$.  
 Let
 \[C= \left\{\mf{y}\in\mb{R}^m \,|\, \Vert\mf{y}-\mf{F}(\mf{x}_0)\Vert\geq\varepsilon \right\}.\]
 Then $C$  is the complement of the open set  $B(\mf{F}(\mf{x}_0), \varepsilon)$. Hence, $C$ is closed. It contains   $\mf{F}(\mf{x}_k)$ for all $k\in\mb{Z}^+$, but it does not contain $\mf{F}(\mf{x}_0)$. Thus, the set $D=\mf{F}^{-1}(C)$ contains the sequence $\{\mf{x}_k\}$,   but does not contain its limit $\mf{x}_0$. This means $D$ is not closed. Therefore, (b) does not hold.

\end{myproof}
There is a much easier proof of Theorem \ref{230720_6} if $\mathcal{A}=\mathbb{R}^n$, using Theorem \ref{230720_5}, and the fact that a set is closed if and only if its complement is open. 

Theorem \ref{230720_5} and Theorem \ref{230720_6} provide useful tools to justfy that a set is open or closed in $\mathbb{R}^n$, using our known library of continuous functions.

\begin{example}[label=230720_7]{}
Let $A$ be the subset of $\mathbb{R}^2$ given by
\[A=\left\{(x,y)\,|\,x^2+y^2<20, y>x^2\right\}.\]
Show that $A$ is open.
\end{example}
\begin{solution}{Solution}
Let $\mathcal{O}=\left\{(x,y)\,|\,x^2+y^2<20\right\}$. This is a ball of radius $\sqrt{20}$ centered at the origin. Hence, $\mathcal{O}$ is open. Define the function $f:\mathcal{O}\to \mathbb{R}$ by 
$\di f(x,y)=y-x^2$.
Since $f$ is a polynomial, it is continuous. Notice that $y>x^2$ if and only if $f(x,y)>0$, if and only if $f(x,y)\in (0,\infty)$. This shows that    $A=f^{-1}((0,\infty))$. Since $(0,\infty)$ is open in $\mathbb{R}$, Theorem  \ref{230720_5} implies that $A$ is an open set.
\end{solution}
 \begin{figure}[ht]
\centering
\includegraphics[scale=0.2]{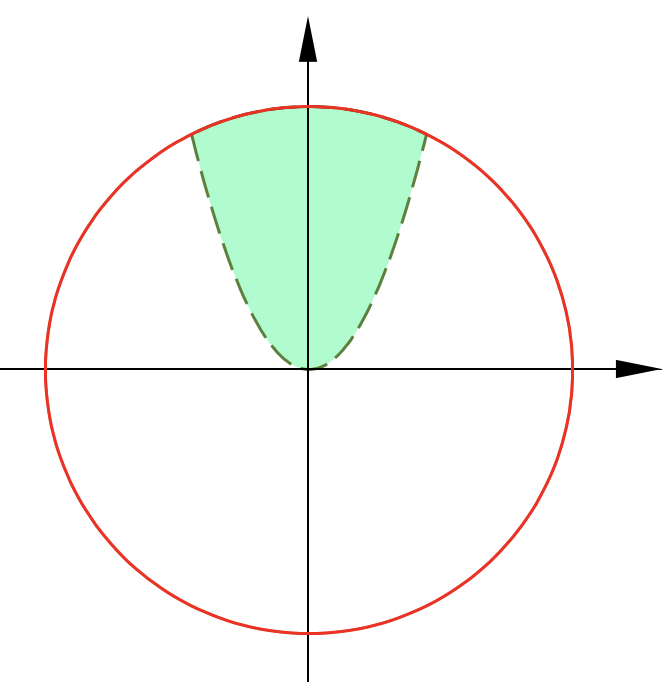}

\caption{The set $A$ in Example \ref{230720_7}.}\label{figure19}
\end{figure}
\begin{example}{} 
Let $C$ be the subset of $\mb{R}^3$ given by
\[C=\left\{(x,y,z)\,|\,x\geq 0, y\geq 0, y^2+z^2\leq 20.\right\}.\]
Show that $C$ is closed.
\end{example}
\begin{solution}{Solution}
Let $\pi_x:\mathbb{R}^3\to\mb{R}$ and $\pi_y:\mb{R}^3\to\mb{R}$ be the projection functions $\pi_x(x,y,z)=x$ and $\pi_y(x,y,z)=y$, and consider the function $g:\mb{R}^3\to\mb{R}$ defined as
\[g(x,y,z)=20-(y^2+z^2).\]Notice that  $y^2+z^2\leq 20$ if and only if $g(x,y,z)\geq 0$, if and only if $g(x,y,z)\in I=[0,\infty)$.  The projection functions $\pi_x$ and $\pi_y$ are continuous. Since $g$ is a polynomial, it is also continuous. The set $I=[0,\infty)$ is closed in $\mathbb{R}$. Therefore, the sets
$\pi_x^{-1}(I)$, $\pi_y^{-1}(I)$ and $g^{-1}(I)$  are closed in $\mathbb{R}^3$.  
Since
\[A= \pi_x^{-1}(I)\cap\pi_y^{-1}(I)\cap g^{-1}(I),\]  being an intersection of three closed sets, $A$ is closed in $\mathbb{R}^3$. 
\end{solution}

Using the same reasonings, we  obtain the following.
\begin{theorem}{}
Let $I_1, \ldots, I_n$ be intervals in $\mathbb{R}$.
\begin{enumerate}[1.]

\item If each of $I_1, \ldots, I_n$ are open intervals of the form $(a,b)$, $(a,\infty)$, $(-\infty, a)$ or $\mathbb{R}$, then $I_1\times \cdots\times I_n$ is an open subset of $\mathbb{R}^n$.
\item If each of $I_1, \ldots, I_n$ are closed intervals of the form $[a,b]$, $[a,\infty)$, $(-\infty, a]$ or $\mathbb{R}$, then $I_1\times \cdots\times I_n$ is a closed subset of $\mathbb{R}^n$.
\end{enumerate}
\end{theorem}
\begin{myproof}{Sketch of Proof}
Use the fact that 
\[I_1\times \cdots\times I_n=\bigcap_{i=1}^n\pi_i^{-1}(I_i),\]
where $\pi_i:\mathbb{R}^n\to\mathbb{R}$ is the projection function $\pi_i(x_1, \ldots, x_n)=x_i$.
\end{myproof}
\begin{example}{}
The set \[A=\left\{(x,y,z)\,|\,x<0, y>2, -10<z<-3\right\}\] is open in $\mathbb{R}^3$, since
\[A=(-\infty, 0)\times (2,\infty)\times (-10,-3).\] The set \[C=\left\{(x,y,z)\,|\,x\leq 0, y\geq 2, -10\leq z\leq -3\right\}\] is closed in $\mathbb{R}^3$, since
\[C=(-\infty, 0]\times [2,\infty)\times [-10,-3].\] 
\end{example}

We also have the following.
\begin{theorem}[label=230724_16]{}
Let $a$ and $b$ be real numbers, and assume that $f:\mathbb{R}^n\to\mathbb{R}$ is a continuous function.
Define the sets $A, B, C, D, E$ and $F$ as follows.

\begin{enumerate}[(a)]
\item  $\di A=\left\{\mf{x}\in\mb{R}^n\,|\,f(\mf{x})>a\right\}$  
\item $\di B=\left\{\mf{x}\in\mb{R}^n\,|\,f(\mf{x})\geq a\right\}$
\item $\di C=\left\{\mf{x}\in\mb{R}^n\,|\,f(\mf{x})< a\right\}$  
\item $\di D=\left\{\mf{x}\in\mb{R}^n\,|\,f(\mf{x})\leq a\right\}$
\item  $\di E=\left\{\mf{x}\in\mb{R}^n\,|\,a<f(\mf{x})<b\right\}$ 
\item $\di F=\left\{\mf{x}\in\mb{R}^n\,|\,a\leq f(\mf{x})\leq b\right\}$
\end{enumerate}
 
Then $A, C$ and $E$ are open sets, while $B$, $D$ and $F$ are closed sets.
\end{theorem}
The proof is left as an exercise.

\begin{example}[label=230724_19]{}

 Find the interior, exterior and boundary of each of the following sets.
\begin{enumerate}[(a)]
\item $\di A=\left\{(x,y)\,|\,0<x^2+4y^2<4\right\}$  
\item  $\di B=\left\{(x,y)\,|\,0< x^2+4y^2\leq 4\right\}$ 
\item  $\di C=\left\{(x,y) \,|\,x^2+4y^2\leq 4\right\}$ 
 \end{enumerate}

\end{example}

\begin{figure}[ht]
\centering
\includegraphics[scale=0.17]{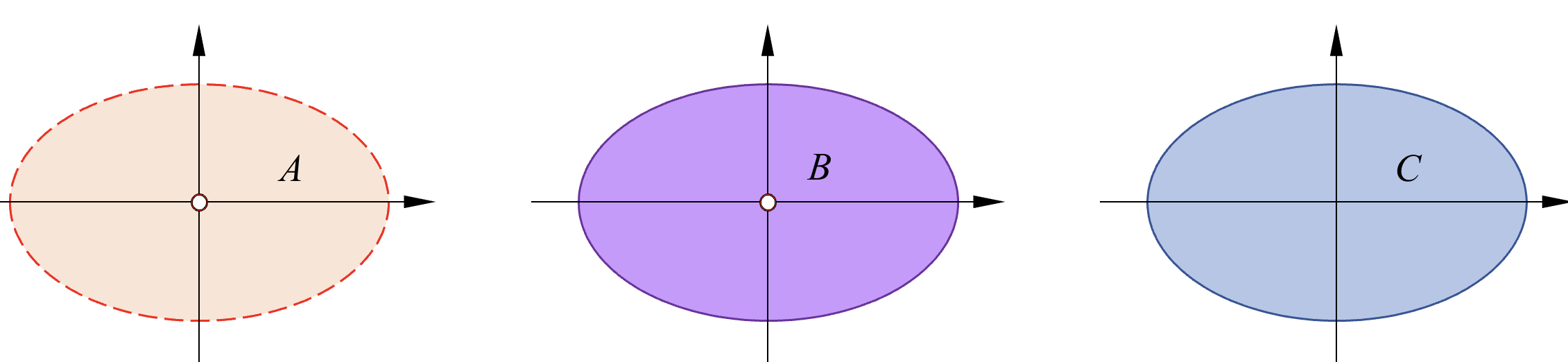}

\caption{The sets $A$, $B$ and $C$ defined  in Example \ref{230724_19}.}\label{figure39}
\end{figure}
\begin{solution}{Solution}
Let 
\[
D=\left\{(x,y) \,|\,x^2+4y^2< 4\right\},\hspace{1cm}
E=\left\{(x,y) \,|\,x^2+4y^2> 4\right\},\]
and let $f:\mathbb{R}^2\to\mathbb{R}$ be the function defined as
\[f(x, y)=x^2+4y^2.\]

Since $f$ is a polynomial, it is continuous. By Theorem \ref{230724_16}, $A, D$ and $E$ are open sets and $C$ is a closed set. Since $A\subset B$ and $D\subset C$, we have
\[A=\text{int}\,A\subset\text{int}\,B\subset B, \hspace{1cm} D\subset\text{int}\, C.\]
Since $\di E= \mb{R}^2\setminus C\subset \mb{R}^2\setminus B\subset \mb{R}^2\setminus A$, 
We have
\[E=\text{ext}\,C\subset \text{ext}\, B\subset\text{ext}\, A.\]
Let
\[F= \left\{(x,y) \,|\,x^2+4y^2= 4\right\}.\]
Then $\mb{R}^n$ is a disjoint union of $D$, $E$ and $F$. If $\mf{u}_0=(x_0, y_0)\in F$, either $x_0\neq 0$ or $y_0\neq 0$, but not both. If $x_0\neq 0$, define the sequences $\{\mf{u}_k\}$ and $\{\mf{v}_k\}$ by
\[\mf{u}_k=\left( \frac{k}{k+1}x_0, y_0\right),\hspace{1cm} \mf{v}_k=\left( \frac{k+1}{k}x_0, y_0\right).\] 
 
If $x_0=0$, then $y_0\neq 0$.
Define the sequences $\{\mf{u}_k\}$ and $\{\mf{v}_k\}$ by
\[\mf{u}_k=\left( x_0, \frac{k}{k+1}y_0\right),\hspace{1cm} \mf{v}_k=\left( x_0, \frac{k+1}{k} y_0\right).\]

In either case, $\{\mf{u}_k\}$ is a sequence of points in $A$ that converges to $\mf{u}_0$, while $\{\mf{v}_k\}$ is a sequence of points in $E$ that converges to $\mf{u}_0$. This proves that $\mf{u}_0$ is a boundary point of $A$, $B$ and $C$.  For the point $\mf{0}$, since it is not in $A$ and $B$, it is not an interior point of $A$ and $B$, but it is the limit of the sequence $\{(1/k, 0)\}$ that is in both $A$ and $B$. Hence, $\mf{0}$ is in the closure of $A$ and $B$, and hence, is a boundary point of $A$ and $B$. We conclude that
\bs
\[\text{int}\,A=\text{int}\,B=\left\{(x,y) \,|\,0<x^2+4y^2< 4\right\},\]
\[\text{int}\, C=\left\{(x,y) \,|\,x^2+4y^2< 4\right\},\]
\[\text{ext}\,A=\text{ext}\,B=\text{ext}\, C=\left\{(x,y) \,|\,x^2+4y^2> 4\right\},\]
\[\text{bd}\,A=\text{bd}\,B=\left\{(x,y) \,|\,x^2+4y^2= 4\right\}\cup\{\mf{0}\},\]
\[\text{bd}\, C=\left\{(x,y) \,|\,x^2+4y^2= 4\right\}.\]
\end{solution}

\begin{remark}{}Let $f:\mb{R}^n\to\mb{R}$ be a continuous function and let 
\[C=\left\{\mf{x}\in\mb{R}^n\,|\, a\leq f(\mf{x})\leq b\right\}.\]
One is tempting to say that \[\text{bd}\,C=\left\{\mf{x}\in\mb{R}^n\,|\,  f(\mf{x})=a\;\text{or}\; f(\mf{x})=b\right\}.\]This is not necessary true.
For example, consider the set $C$ in Example \ref{230724_19}. It can be written as
\[C=\left\{(x,y) \,|\,0\leq x^2+4y^2\leq 4\right\} \]However, the point where $f(x,y)=x^2+4y^2=0$ is not a boundary point of $C$. 
\end{remark}

Now we return to  continuous functions.  
\begin{theorem}[label=230721_5]{Pasting of Continuous Functions}
Let $A$ and $B$ be closed subsets of $\mathbb{R}^n$, and let $S=A\cup B$. If $\mf{F}:S\to\mb{R}^m$ is a function such that $\mf{F}_A=\mf{F}|_A:A\to\mb{R}^m$ and $\mf{F}_B=\mf{F}|_B:B\to\mb{R}^m$ are both continuous, then $\mf{F}:S\to\mb{R}^m$ is continuous.
\end{theorem}
\begin{myproof}{Proof}
Since $S$ is a union of two closed sets, it is closed. Applying Theorem \ref{230720_6}, it suffices to show that if $C$ is a closed subset of $\mb{R}^m$, then $\mf{F}^{-1}(C)$ is closed in $\mb{R}^n$. 
Notice that 
\begin{align*}
\mf{F}^{-1}(C)&=\left\{\mf{x}\in S\,|\, \mf{F}( \mf{x})\in C\right\}\\
&= \left\{\mf{x}\in A\,|\, \mf{F}( \mf{x})\in C\right\}\cup\left\{\mf{x}\in B\,|\, \mf{F}( \mf{x})\in C\right\}\\
&=\mf{F}_A^{-1}(C)\cup \mf{F}_B^{-1}(C).
\end{align*}Since $\mf{F}_A:A\to\mb{R}^m$ and $\mf{F}_B:B\to\mb{R}^m$ are both continuous functions, $\mf{F}_A^{-1}(C)$ and $\mf{F}_B^{-1}(C)$ are closed subsets of $\mathbb{R}^n$. Being a union of two closed subsets, $\mf{F}^{-1}(C)$ is closed. This completes the proof.
\end{myproof}
 
 \begin{example}{}
 Let $f:\mathbb{R}^2\to\mb{R}$ be the function defined as
 \[f(x,y)=\begin{cases}x^2+y^2,\quad &\text{if}\;x^2+y^2<1\\1,\quad &\text{if}\;x^2+y^2\geq 1.\end{cases}\]
 Show that $f$ is a continuous function.
 \end{example}
 \begin{solution}{Solution}
 Let $A=\left\{(x,y)\,|\, x^2+y^2\leq 1\right\}$ and $B=\left\{(x,y)\,|\, x^2+y^2\geq 1\right\}$. Then $A$ and $B$ are closed subsets of $\mb{R}^2$ and $\mb{R}^2=A\cup B$. Notice that $f|_A:A\to\mb{R}$ is the function $f(x,y)=x^2+y^2$, which is continuous since it is a polynomial. By definition, $f|_B:B\to\mb{R}$ is the constant function $f_B(x,y)=1$, which is also continuous. By Theorem \ref{230721_5}, the function $f:\mathbb{R}^2\to\mb{R}$ is continuous.
 \end{solution}
 
 Given positive integers $n$ and $m$, there is a natural bijective correspondence between $\mb{R}^n\times \mb{R}^m$ and $\mb{R}^{n+m}$ given by
 $\mf{T}:\mb{R}^n\times \mb{R}^m \to \mb{R}^{n+m}$,
 \[(\mathbf{x}, \mathbf{y})\mapsto (x_1, \ldots, x_n, y_1, \ldots, y_m),\]  where
 \[\mf{x}=(x_1, \ldots, x_n)\quad\text{and}\quad \mf{y}=(y_1, \ldots, y_m).\]
 Hence, sometimes we will denote a point in $\mb{R}^{n+m}$ as $(\mathbf{x}, \mathbf{y})$, where $\mathbf{x}\in\mb{R}^n$ and $\mf{y}\in\mb{R}^m$. By generalized Pythagoras theorem,
 \[\Vert (\mf{x}, \mf{y})\Vert^2=\Vert\mf{x}\Vert^2+\Vert\mf{y}\Vert^2.\]
 If $A$ is a subset of $\mathbb{R}^n$, $B$ is a subset of $\mathbb{R}^m$, $A\times B$ can be considered  as a subset of $\mb{R}^{n+m}$ given by
 \[A\times B=\left\{(\mathbf{x}, \mathbf{y})\,|\,\mathbf{x}\in A, \mathbf{y}\in B\right\}.\]
 
The following  is more general than Proposition \ref{230721_9}.
\begin{proposition}[label=230721_10]{}
Let $\mk{D}$   be a subset of $\mathbb{R}^n$, and let $\mf{F}:\mk{D}\to \mb{R}^k$ and $\mf{G}:\mk{D}\to \mb{R}^l$ be functions defined on $\mk{D}$. Define the function $\mf{H}:\mk{D}\to\mb{R}^{k+l}$ by
\[\mf{H}(\mf{x})=(\mf{F}(\mf{x}), \mf{G}(\mf{x})).\]
Then the function $\mf{H}:\mk{D} \to\mb{R}^{k+l}$ is continuous if and only if the functions $\mf{F}:\mk{D}\to \mb{R}^k$ and $\mf{G}:\mk{D}\to \mb{R}^l$  are continuous.
\end{proposition}
\begin{myproof}{Sketch of Proof}
This proposition follows immediately from Proposition \ref{230721_9}, since 
\[\mf{H}(\mf{x})=(F_1(\mf{x}), \ldots, F_k(\mf{x}), G_1(\mf{x}), \ldots, G_l(\mf{x})).\]
 
\end{myproof}
 
 For a function defined on a subset of $\mb{R}^n$, we can define its graph in the following way.
 \begin{definition}{The Graph of a Function}
 Let  $\mf{F}:\mk{D}\to\mb{R}^m$ be a function defined on $\mk{D}\subset\mb{R}^n$.
 The graph of $\mathbf{F}$, denoted by $G_{\mf{F}}$, is the subset of $\mb{R}^{n+m}$ defined as
 \[G_{\mf{F}}=\left\{(\mathbf{x}, \mathbf{y})\,|\, \mathbf{x}\in\mk{D}, \mf{y}=\mf{F}(\mf{x})\right\}.\]
 \end{definition}
 
 \begin{example}{}
 Let $\mk{D}=\left\{(x,y)\,|\,x^2+y^2\leq 1\right\}$, and let $f:\mk{D}\to\mb{R}$ be the function defined as
 \[ f(x,y)=\sqrt{1-x^2-y^2}.\]
 The graph of $f$ is 
 \[G_f=\left\{(x,y,z)\,|\, x^2+y^2\leq 1, z=\sqrt{1-x^2-y^2}\right\},\]
 which is the upper hemisphere.
 \end{example}
 \begin{figure}[ht]
\centering
\includegraphics[scale=0.2]{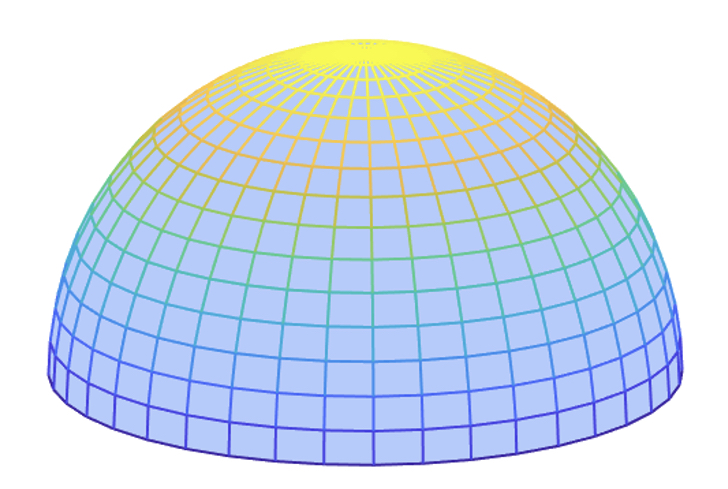}

\caption{The upper hemisphere is the graph of a function.}\label{figure24}
\end{figure}

Notice that if $\mk{D}$ is  a subset of $\mb{R}^n$, then the graph of the function $\mf{F}:\mk{D}\to\mb{R}^m$ is the image of the function $\mf{H}:\mk{D}\to \mb{R}^{n+m}$ defined as
\[\mf{H}(\mf{x})=\left(\mf{x}, \mf{F}(\mf{x})\right).\]
From Proposition \ref{230721_10}, we obtain the following.
\begin{corollary}[label=230721_11]{}
Let $\mk{D}$ be  a subset of $\mb{R}^n$, and let $\mf{F}:\mk{D}\to\mb{R}^m$ be a function defined on $\mk{D}$. The image of the function $\mf{H}:\mk{D}\to \mb{R}^{n+m}$,  
\[\mf{H}(\mf{x})=\left(\mf{x}, \mf{F}(\mf{x})\right),\] is the graph of $\mf{F}$. If the function  $\mf{F}:\mk{D}\to\mb{R}^m$ is continuous, then the function  $\mf{H}:\mk{D}\to \mb{R}^{n+m}$ is continuous.
\end{corollary}
 
 Now we consider 
  a special class of   functions called Lipschitz functions. 
\begin{definition}{}
Let $\mk{D}$ be a subset of $\mathbb{R}^n$. A function $\mf{F}:\mk{D}\to \mathbb{R}^m$ is Lipschitz provided that there exists a positive constant $c$ such that
\[\Vert\mf{F}(\mf{u})-\mf{F}(\mf{v})\Vert\leq c\Vert\mf{u}-\mf{v}\Vert\hspace{1cm}\text{for all}\;\mf{u}, \mf{v}\in\mk{D}.\]
The constant $c$ is called a Lipschitz constant of the function. If $c<1$, then $\mf{F}:\mk{D}\to \mathbb{R}^m$  is  called a contraction.
\end{definition}
 The following is easy to establish.

\begin{proposition}{}
Let $\mk{D}$ be a subset of $\mathbb{R}^n$, and let $\mf{F}:\mk{D}\to \mathbb{R}^m$ be a Lipschitz function. Then $\mf{F}:\mk{D}\to \mathbb{R}^m$ is continuous.
\end{proposition}
 
 \begin{example}{}
 A linear transformation of the form $\mf{T}:\mb{R}^n\to\mb{R}^n$, $\mf{T}(\mf{x})=a\mf{x}$,
 is a Lipschitz function with Lipschitz constant $|a|$.
 \end{example}
 In fact,  we have the following.
 \begin{theorem}[label=230724_9]{}
 A linear transformation $\mf{T}:\mb{R}^n\to\mb{R}^m$ is a Lipschitz function.
 \end{theorem}
 \begin{myproof}{Proof}
 Let $A$ be the $m\times n$ matrix such that $\mf{T}(\mf{x})=A\mf{x}$.
When $\mf{x}$ is in $\mb{R}^n$, 
 \[\Vert\mf{T}(\mf{x})\Vert^2=\left(A\mf{x}\right)^T(A\mf{x})=\mf{x}^T(A^TA)\mf{x}.\]
 The matrix $B=A^TA$ is a   positive semi-definite $n\times n$ symmetric matrix.  By Theorem \ref{230724_8},
 \[\mf{x}^T(A^TA)\mf{x}\leq \lambda_{\max}\Vert\mf{x}\Vert^2,\]
 where $\lambda_{\max}$ is the largest eigenvalue of $A^TA$. 
 \bp
 Therefore, for any $\mf{x}\in \mb{R}^n$, 
 \[\Vert \mf{T}(\mf{x})\Vert\leq \sqrt{\lambda_{\max}}\Vert\mf{x}\Vert.\]It follows that for any $\mf{u}$ and $\mf{v}$ in $\mb{R}^n$, 
  \[\Vert \mf{T}(\mf{u})-\mf{T}(\mf{v})\Vert=\Vert \mf{T}(\mf{u}-\mf{v})\Vert\leq \sqrt{\lambda_{\max}}\Vert\mf{u}-\mf{v}\Vert.\]
 Hence, $\mf{T}:\mb{R}^n\to\mb{R}^m$ is a Lipschitz mapping with Lipschitz constant $\sqrt{\lambda_{\max}}$.
 \end{myproof}
 
 \begin{example}{}
 Let $\mf{T}:\mb{R}^2\to\mb{R}^2$ be the mapping defined as
 \[\mf{T}(x,y)=(x-3y, 7x+4y).\]
 Find the smallest constant $c$ such that 
 \[\Vert\mf{T}(\mf{u})-\mf{T}(\mf{v})\Vert\leq c\Vert\mf{u}-\mf{v}\Vert\]for all $\mf{u}$ and $\mf{v}$ in $\mb{R}^2$.
 \end{example}
 \begin{solution}{Solution}
 Notice that
$\mf{T}(\mf{u})=A\mf{u}$, where $A$ is the $2\times 2$ matrix
$A=\di\begin{bmatrix} 1 & -3\\7 & 4\end{bmatrix}$. Hence,
\[\Vert\mf{T}(\mf{u})\Vert^2=\mf{u}^TA^TA\mf{u}=\mf{u}^TC\mf{u},\]where
\begin{align*}
C= \begin{bmatrix} 1 & 7 \\ -3 &4\end{bmatrix} \begin{bmatrix} 1 & -3\\7 & 4\end{bmatrix}=\begin{bmatrix} 50 & 25\\25 & 25\end{bmatrix}=25\begin{bmatrix} 2 & 1\\ 1 & 1\end{bmatrix}.
\end{align*}
 
For the matrix $G=\di\begin{bmatrix} 2 & 1 \\ 1 &1 \end{bmatrix}$, the eigenvalues are the solutions of 
\bs
\[\lambda^2-3\lambda+1=0,\]
which are
\[\lambda_1=\frac{3+\sqrt{5}}{2}\quad\text{and}\quad \lambda_2=\frac{3-\sqrt{5}}{2}.\]
Hence,
\[\Vert\mf{T}(\mf{u})\Vert^2\leq \frac{25(3+\sqrt{5})}{2}\Vert\mf{u}\Vert^2.\]
The smallest $c$ such that $\Vert\mf{T}(\mf{u})-\mf{T}(\mf{v})\Vert\leq c\Vert\mf{u}-\mf{v}\Vert$ for all $\mf{u}$ and $\mf{v}$ in $\mb{R}^2$ is
\[c=\sqrt{\frac{25(3+\sqrt{5})}{2}}=8.0902.\]
 \end{solution}
 
 \begin{remark}{}
 If $A$ is an $m\times n$ matrix, the matrix $B=A^TA$ is a positive semi-definite $n\times n$ symmetric matrix. Thus, all its eigenvalues are nonnegative. Let $\lambda_1, \ldots,\lambda_n$ be its eigenvalues with
  \[0= \lambda_n=\cdots=\lambda_{r+1}<\lambda_r\leq\lambda_{r-1}\leq \cdots \leq \lambda_1.\]Then $\lambda_1, \cdots, \lambda_r$ are the nonzero eigenvalues of $A^TA$. The singular values of $A$ are the numbers $\sigma_1, \ldots, \sigma_r$, where
  \[\sigma_i=\sqrt{\lambda_i},\quad 1\leq i\leq r.\]
  Theorem \ref{230724_9}  says that $\sigma_1$ is a Lipschitz constant of the linear transformation $\mf{T}(\mf{x})=A\mf{x}$.
 \end{remark}
 
 At the end of this section, we 
 want to discuss the vector space of $m\times n$ matrices $\mathcal{M}_{m,n}$. There is a natural vector space isomorphism between $\mathcal{M}_{m,n}$ and $\mb{R}^{mn}$, by mapping the matrix $A=[a_{ij}]$ to $\mf{x}=(x_k)$, where
 \[x_{(i-1)n+j}=a_{ij}\hspace{1cm}\text{for}\; 1\leq i\leq m, 1\leq j\leq n.\]
 In other words, if
 \begin{align*}\mathbf{a}_1&=(a_{11}, a_{12}, \ldots, a_{1n}), \\
 \mathbf{a}_2&=(a_{21}, a_{22}, \ldots, a_{2n}),\\
 &\hspace{1cm}\vdots\\
 \mf{a}_m&=(a_{m,1}, a_{m,2}, \ldots, a_{m,n})\end{align*} are the row vectors of $A$, then $A$ is mapped to
 the vector $(\mf{a}_1, \mf{a}_2, \ldots, \mf{a}_m)$ in $\mb{R}^{mn}$.
 Under this isomorphism, 
 the norm of a matrix $A=[a_{ij}]$ is
 \[\Vert  A\Vert=\sqrt{\sum_{i=1}^m\sum_{j=1}^n a_{ij}^2}=\sqrt{\sum_{i=1}^m\Vert\mf{a}_i\Vert^2},\]and the distance between two matrices $A=[a_{ij}]$ and $B=[b_{ij}]$ is
 \[d(A,B)=\Vert  A-B\Vert=\sqrt{\sum_{i=1}^n\sum_{j=1}^n (a_{ij}-b_{ij})^2}.\]
 The following proposition can be used to give an alternative proof of Theorem \ref{230724_9}.
 \begin{proposition}[label=230806_18]{}
 Let $A$ be an $m\times n$ matrix. If $\mf{x}$ is in $\mb{R}^n$, then 
 \[\Vert A\mf{x}\Vert\leq \Vert A\Vert\Vert \mf{x}\Vert.\]
 
 \end{proposition}
 \begin{myproof}{Proof}
 Let $\mf{a}_1$, $\ldots$, $\mf{a}_m$ be the row vectors of $A$, and let $\mf{w}=A\mf{x}$. Then 
 \[w_i=\langle\mathbf{a}_i, \mf{x}\rangle \hspace{1cm}\text{for}\;1\leq i\leq m.\]
 By Cauchy-Schwarz inequality, 
 \[|w_i|\leq \Vert\mf{a}_i\Vert \Vert\mf{x}\Vert \hspace{1cm}\text{for}\;1\leq i\leq m.\]
 \bp
 Thus,
 \begin{align*}\Vert\mf{w}\Vert &=\sqrt{w_1^2+w_2^2+\cdots+w_m^2}\\&\leq \Vert\mf{x}\Vert \sqrt{\Vert\mf{a}_1\Vert^2+ \Vert\mf{a}_2\Vert^2+\cdots +\Vert\mf{a}_m\Vert^2}=\Vert A\Vert \Vert\mf{x}\Vert.\end{align*}
 \end{myproof}
 
 The difference between the proofs of Theorem \ref{230724_9} and Proposition \ref{230806_18} is that, in the proof of Theorem \ref{230724_9}, we find that the smallest possible $c$ such that $\Vert A\mf{x}\Vert\leq c\Vert\mf{x}\Vert$ for all $\mf{x}$ in $\mb{R}^n$ is the largest singular value of the matrix $A$. In Proposition \ref{230806_18}, we find a candidate for $c$, which is the norm of the matrix $A$, but this  is usually not the optimal one.
 
When $m=n$, we  denote the space  of $n\times n$ matrices $\mathcal{M}_{n,n}$ simply as
   $\mathcal{M}_n$  .  
 The determinant of the matrix $A=[a_{ij}]\in \mathcal{M}_n$ is given by
 \[\det A=\sum_{\sigma }\text{sgn}(\sigma)a_{1\sigma(1)}a_{2\sigma(2)}\cdots a_{n\sigma(n)}.\]Here the summation is over all the $n!$ permutations $\sigma$ of the set $S_n=\{1, 2, \ldots, n\}$, and $\text{sgn}(\sigma)$ is the sign of the permutation $\sigma$, which is equal to 1 or $-1$, depending on whether $\sigma$ can be written as the product of an even number or an odd number of transpositions. 
For example, when $n=1$, $\det [a]=a$. When $n=2$,
\[\det\begin{bmatrix}a_{11} & a_{12}\\a_{21} &a_{22}\end{bmatrix}=a_{11}a_{22}-a_{12}a_{21}.\]
When $n=3$, 
\begin{align*}
\det\begin{bmatrix}a_{11} & a_{12} &a_{13}\\a_{21} &a_{22}& a_{23}\\a_{31} & a_{32} & a_{33}\end{bmatrix}&=a_{11}a_{22}a_{33}+a_{12}a_{23}a_{31}+a_{13}a_{21}a_{32}\\&\quad -a_{11}a_{23}a_{32}-a_{13}a_{22}a_{31}-a_{12}a_{21}a_{33}.
\end{align*}

The determinant function $\det: \mathcal{M}_n\to\mb{R}$ is a polynomial function on the variables $(a_{ij})$.  
 Hence,  it is a continuous function. Recall that a matrix $A\in\mathcal{M}_n$ is invertible if and only if $\det A\neq 0$. Let
 \[\text{GL}\,(n, \mb{R})=\left\{A\in \mathcal{M}_n\,|\, \det A\neq 0\right\}\] be the subset of $\mathcal{M}_n$ that consist of invertible $n\times n$ matrices. It is a group under matrix multiplication, called the general linear group. By definition, 
 \[\text{GL}\,(n,\mb{R})=\det\!^{-1}(\mb{R}\setminus\{0\}).\]
 Since  $\mb{R}\setminus\{0\}$ is an open subset of $\mb{R}$, $\text{GL}\,(n,\mb{R})$ is an open subset of $\mathcal{M}_n$. This gives the following.
 \begin{proposition}{}
 Given that $A$ is an invertible $n\times n$ matrix, there exists $r>0$ such that if $B$ is an $n\times n$ matrix with $\Vert B-A\Vert <r$, then $B$ is also invertible.
 \end{proposition}
\begin{myproof}{Sketch of Proof}
This is simply a rephrase of the statement that if $A$ is a point in the open set $\text{GL}\,(n,\mb{R})$, then there is a ball $B(A, r)$ with center at $A$ that is contained in $\text{GL}\,(n,\mb{R})$.  
\end{myproof}

Let $A$ be an $n\times n$ matrix. For $1\leq i, j\leq n$, the $(i,j)$-minor of $A$, denoted by $M_{i,j}$, is the determinant of the $(n-1)\times (n-1)$ matrix obtained by deleting the $i^{\text{th}}$-row and $j^{\text{th}}$- column of $A$. 
Using the same reasoning as above, we find that the function $M_{i,j}:\mathcal{M}_n\to \mb{R}$ is a continuous function. 
The $(i,j)$ cofactor  $C_{i,j}$ of $A$ is given by $C_{i,j}=(-1)^{i+j}M_{i,j}$.
The cofactor matrix of $A$ is $C_A=[C_{ij}]$. Since each of the components is continuous, the function $C:\mathcal{M}_n\to \mathcal{M}_n$ taking $A$ to $C_A$ is  a continuous function.

If $A$ is invertible,
\[A^{-1}=\frac{1}{\det A}C_A^T.\]Since both $C:\mathcal{M}_n\to \mathcal{M}_n$ and $\det:\mathcal{M}_n\to\mb{R}$ are continuous functions, and $\det:\text{GL}\,(n,\mb{R})\to\mb{R}$ is a function that is never equal to 0, we obtain the following.

\begin{theorem}{}
The map $\mathscr{I}:\text{GL}\,(n,\mb{R})\to \text{GL}\,(n,\mb{R})$ that takes $A$ to $A^{-1}$ is continuous.
\end{theorem}
\vp
\noindent
{\bf \large Exercises  \thesection}
\setcounter{myquestion}{1}

\begin{question}{\themyquestion}
Let $\mathbf{x}_0$ be a point in $\mb{R}^n$. Define the function $f:\mb{R}^n\to\mb{R}$ by
\[f(\mf{x})=\Vert \mf{x}-\mf{x}_0\Vert.\]
Show that $f$ is a continuous function.
\end{question}
\atc
\begin{question}{\themyquestion}
Let $\mathcal{O}=\mathbb{R}^3\setminus\{(0,0,0)\}$ and define the function $\mf{F}:\mathcal{O}\to \mathbb{R}^2$ by
\[\mf{F}(x,y,z)=\left(\frac{y}{x^2+y^2+z^2}, \frac{z}{x^2+y^2+z^2}\right).\]
Show that $\mf{F}$ is a continuous function.
\end{question}
\atc
\begin{question}{\themyquestion}
Let   $f:\mb{R}^n\to\mb{R}$ be the function defined as
\[f(\mathbf{x})=\begin{cases}1,\quad &\text{if at least one of the $x_i$ is rational},\\
0,\quad &\text{otherwise}.\end{cases}\] At which point of $\mathbb{R}^n$ is the function $f$ continuous?
\end{question}
\atc
\begin{question}{\themyquestion}
Let   $f:\mb{R}^n\to\mb{R}$ be the function defined as
\[f(\mathbf{x})=\begin{cases}x_1^2+\cdots+x_n^2,\quad &\text{if at least one of the $x_i$ is rational},\\
0,\quad &\text{otherwise}.\end{cases}\] At which point of $\mathbb{R}^n$ is the function $f$ continuous?
\end{question}
\atc
\begin{question}{\themyquestion}Let $f:\mb{R}^3\to\mb{R}$ be the function defined by
\[f(x,y,z)=\begin{cases}\di \frac{\sin(x^2+4y^2+z^2)}{x^2+4y^2+z^2},\quad &\text{if}\;(x,y,z)\neq (0,0,0),\\
a,\quad &\text{if}\;(x,y,z)=(0,0,0).\end{cases}\]
Show that there exists a value $a$ such that $f$ is a continuous function, and find this value of $a$.
\end{question}
\atc
\begin{question}{\themyquestion}
Let $a$ and $b$ be positive numbers, and let $\mathcal{O}$ be the subset of $\mb{R}^n$ defined as
\[\mathcal{O}=\left\{\mf{x}\in\mb{R}^n\,|\, a<\Vert\mf{x}\Vert<b\right\}.\]
Show that $\mathcal{O}$ is open.
\end{question}
\atc
\begin{question}{\themyquestion}
Let $A$ be the subset of $\mathbb{R}^2$ given by
\[A=\left\{(x,y)\,|\, \sin(x+y)+xy>1\right\}.\]
Show that $A$ is an open set.
\end{question}
 \atc
 \begin{question}{\themyquestion}
Let $A$ be the subset of $\mathbb{R}^3$ given by
\[A=\left\{(x,y,z)\,|\,x\geq 0, y\leq 1, e^{xy}\leq z\right\}.\]
Show that $A$ is a closed set.
\end{question}
 \atc
\begin{question}{\themyquestion}
A plane in $\mathbb{R}^3$ is the set of all points $(x,y,z)$ satisfying an equation of the form 
\[ax+by+cz=d,\]where $(a,b,c)\neq (0,0,0)$. Show that a plane is a closed subset of $\mb{R}^3$.
\end{question}
\atc
\begin{question}{\themyquestion}
Define the sets $A, B, C$ and $D$ as follows.

\begin{enumerate}[(a)]
\item $\di A=\left\{(x,y, z)\,|\, x^2+4y^2+9z^2<36\right\}$  
\item  $\di B=\left\{(x,y,z) \,|\,x^2+4y^2+9z^2\leq 36\right\}$ 

\item $\di C=\left\{(x,y, z)\,|\,0<x^2+4y^2+9z^2<36\right\}$  
\item  $\di D=\left\{(x,y,z)\,|\,0< x^2+4y^2+9z^2\leq 36\right\}$ 
 \end{enumerate}
For each of these sets, find its interior, exterior and boundary.
\end{question}
  \atc
  
\begin{question}{\themyquestion}
Let $a$ and $b$ be real numbers, and assume that $f:\mathbb{R}^n\to\mathbb{R}$ is a continuous function.
Consider the following subsets of $\mathbb{R}^n$.

\begin{enumerate}[(a)]
\item  $\di A=\left\{\mf{x}\in\mb{R}^n\,|\,f(\mf{x})>a\right\}$  
\item $\di B=\left\{\mf{x}\in\mb{R}^n\,|\,f(\mf{x})\geq a\right\}$
\item $\di C=\left\{\mf{x}\in\mb{R}^n\,|\,f(\mf{x})< a\right\}$  
\item $\di D=\left\{\mf{x}\in\mb{R}^n\,|\,f(\mf{x})\leq a\right\}$
\item  $\di E=\left\{\mf{x}\in\mb{R}^n\,|\,a<f(\mf{x})<b\right\}$ 
\item $\di F=\left\{\mf{x}\in\mb{R}^n\,|\,a\leq f(\mf{x})\leq b\right\}$
\end{enumerate}
Show that $A, C$ and $E$ are open sets, while $B$, $D$ and $F$ are closed sets.

\end{question}

\atc
\begin{question}{\themyquestion}
Let  $f:\mathbb{R}^2\to\mb{R}$ be the function defined as
 \[f(x,y)=\begin{cases}x^2+y^2,\quad &\text{if}\;x^2+y^2<4\\8-x^2-y^2,\quad &\text{if}\;x^2+y^2\geq 4.\end{cases}\]
 Show that $f$ is a continuous function.
\end{question}

\atc
\begin{question}{\themyquestion}
Show that the distance function on $\mathbb{R}^n$, $d:\mathbb{R}^n\times\mathbb{R}^n\to\mathbb{R}$,
\[d(\mf{u},\mf{v})=\Vert\mf{u}-\mf{v}\Vert,\] is continuous in the following sense. If $\{\mf{u}_k\}$ and $\{\mf{v}_k\}$ are sequences in $\mb{R}^n$ that converges to $\mf{u}$ and $\mf{v}$ respectively, then the sequence $\{d(\mf{u}_k,\mf{v}_k)\}$ converges to $d(\mf{u}, \mf{v})$.
 
\end{question}

\atc
 
\begin{question}{\themyquestion}
Let $\mf{T}:\mb{R}^2\to\mb{R}^3$ be the mapping 
\[\mf{T}(x,y)=(x+y, 3x-y, 6x+5y).\]
Show that $\mf{T}:\mb{R}^2\to\mb{R}^3$ is a Lipschitz mapping, and find the smallest Lipschitz constant for this mapping.
 \end{question}
\atc
\begin{question}{\themyquestion}
Given that $A$ is a subset of $\mb{R}^m$  and   $B$ is a subset of $\mb{R}^n$, let $C=A\times B$. Then $C$ is a subset of $\mb{R}^{m+n}$.
\begin{enumerate}[(a)]
\item
If $A$ is open in $\mb{R}^m$ and $B$ is open in $\mb{R}^n$, show that $A\times B$ is open in $\mb{R}^{m+n}$.
\item If $A$ is closed in $\mb{R}^m$ and $B$ is closed in $\mb{R}^n$, show that $A\times B$ is closed in $\mb{R}^{m+n}$.
\end{enumerate}
\end{question}

\atc
 
\begin{question}{\themyquestion}
Let $\mk{D}$ be a   subset of $\mb{R}^n$, and let $f:\mk{D}\to\mb{R}$ be a continuous function defined on $\mk{D}$. 
   Let $A=\mk{D}\times \mb{R}$ and define the function $g: A\to\mb{R}$ by
\[g(\mf{x}, y)= y-f(\mf{x}).\]
Show that $g: A\to\mb{R}$ is continuous.

\end{question}
\atc
 
\begin{question}{\themyquestion}
Let $U$ be an open subset of $\mb{R}^n$, and let $f:U\to\mb{R}$ be a continuous function defined on $U$. 
Show that the sets
\begin{align*}
\mathcal{O}_1&=\left\{(\mf{x}, y)\,|\, \mf{x}\in U, y< f(\mf{x})\right\},\quad
\mathcal{O}_2=\left\{(\mf{x}, y)\,|\, \mf{x}\in U, y>f(\mf{x})\right\}
\end{align*}
are open subsets of $\mb{R}^{n+1}$.

\end{question}

 \atc
 
\begin{question}{\themyquestion}
Let $C$ be a closed subset of $\mb{R}^n$, and let $f:C\to\mb{R}$ be a continuous function defined on $C$. 
Show that the sets
\begin{align*}
\mathcal{A}_1&=\left\{(\mf{x}, y)\,|\, \mf{x}\in C, y\leq  f(\mf{x})\right\},\quad
\mathcal{A}_2=\left\{(\mf{x}, y)\,|\, \mf{x}\in C, y\geq f(\mf{x})\right\}
\end{align*}
are closed subsets of $\mb{R}^{n+1}$.

\end{question}

\section{Uniform Continuity}\label{sec1_9}

  In volume I, we have seen that uniform continuity plays important role in single variable analysis. In this section, we extend this concept to multivariable functions.
 
 \begin{definition}{Continuity}
Let $\mk{D}$ be a subset of $\mb{R}^n$, and let $\mf{F}:\mk{D}\rightarrow\mathbb{R}^m$ be a function defined on $\mk{D}$. We say that the function $\mf{F}$ is {\bf uniformly continuous   }   provided that for any $\varepsilon>0$, there exists $\delta>0$ such that if $\mf{u}$ and $\mf{v}$ are points in $\mk{D}$ and $\Vert\mf{u}-\mf{v}\Vert<\delta$, then
\[\Vert\mf{F}(\mf{u})-\mf{F}(\mf{v})\Vert<\varepsilon.\]
\end{definition}
The following two propositions are obvious.
\begin{proposition}{}
A   uniformly continuous function is continuous. 
\end{proposition}
\begin{proposition}{}
Given that $\mk{D}$ is a subset of $\mb{R}^n$, and $\mk{D}'$ is a subset of $\mk{D}$, if the function $\mf{F}:\mk{D}\to\mb{R}^m$ is uniformly continuous, then  the function $\mf{F}:\mk{D}'\to\mb{R}^m$ is also uniformly continuous.
\end{proposition}
A special class of uniformly continuous functions is the class of Lipschitz functions.
\begin{theorem}[label=230724_11]{}
Let $\mk{D}$ be a subset of $\mb{R}^n$, and let $\mf{F}:\mk{D}\rightarrow\mathbb{R}^m$ be a function defined on $\mk{D}$. If $\mf{F}:\mk{D}\rightarrow\mathbb{R}^m$ is Lipschitz, then it is uniformly continuous.
\end{theorem}
The proof is straightforward. 

\begin{remark}{}Theorem \ref{230724_9} and Theorem
\ref{230724_11} imply  that a linear transformation is uniformly continuous.
\end{remark}

There is an equivalent definition for uniform continuity in terms of sequences.

\begin{theorem}[label=230724_12]{}
Let $\mk{D}$ be a subset of $\mb{R}^n$, and let $\mf{F}:\mk{D}\rightarrow\mathbb{R}^m$ be a function defined on $\mk{D}$. Then  the following are equivalent.
\begin{enumerate}[(i)]
\item $\mf{F}:\mk{D}\rightarrow\mathbb{R}^m$  is uniformly continuous. Namely, given $\varepsilon>0$, there exists $\delta>0$ such that if $\mf{u}$ and $\mf{v}$ are points in $\mk{D}$ and $\Vert\mf{u}-\mf{v}\Vert<\delta$, then
\[\Vert\mf{F}(\mf{u})-\mf{F}(\mf{v})\Vert<\varepsilon.\]
\item If $\{\mf{u}_k\}$ and $\{\mf{v}_k\}$ are two sequences in $\mk{D}$ such that 
\[\lim_{k\to\infty} \left(\mf{u}_k-\mf{v}_k \right)=\mf{0},\]
then
\[\lim_{k\to\infty}\left(\mf{F}(\mf{u}_k)-\mf{F}(\mf{v}_k)\right)=\mf{0}.\]
\end{enumerate}
\end{theorem}
Let us give a proof of this theorem here.
\begin{myproof}{Proof}
Assume that (i) holds, and $\{\mf{u}_k\}$ and $\{\mf{v}_k\}$ are two sequences in $\mk{D}$ such that 
\[\lim_{k\to\infty} \left(\mf{u}_k-\mf{v}_k \right)=\mf{0}.\] Given $\varepsilon>0$, (i) implies that there exists $\delta>0$ such that if $\mf{u}$ and $\mf{v}$ are points in $\mk{D}$ and $\Vert\mf{u}-\mf{v}\Vert<\delta$, then
\[\Vert\mf{F}(\mf{u})-\mf{F}(\mf{v})\Vert<\varepsilon.\] Since $\di \lim_{k\to\infty} \left(\mf{u}_k-\mf{v}_k \right)=\mf{0}$, there is a positive integer $K$ such that for all $k\geq K$, $\Vert \mf{u}_k-\mf{v}_k \Vert<\delta$. It follows that
\[\Vert\mf{F}(\mf{u}_k)-\mf{F}(\mf{v}_k)\Vert<\varepsilon
\hspace{1cm}\text{for all}\;k\geq K.\] 
\bp
This shows that \[\lim_{k\to\infty}\left(\mf{F}(\mf{u}_k)-\mf{F}(\mf{v}_k)\right)=\mf{0},\]and thus completes the proof of (i) implies (ii).

Conversely, assume that (i) does not hold. This means  there exists an $\varepsilon>0$, for all $\delta>0$, there exist points $\mf{u}$ and $\mf{v}$ in $\mk{D}$ such that $\Vert\mf{u}-\mf{v}\Vert<\delta$ and $\Vert\mf{F}(\mf{u})-\mf{F}(\mf{v})\Vert\geq \varepsilon$. Thus, for every $k\in\mb{Z}^+$, there exists $\mf{u}_k$ and $\mf{v}_k$ in $\mk{D}$ such that
\begin{equation}\label{230724_13}\Vert\mf{u}_k-\mf{v}_k\Vert<\frac{1}{k},\end{equation}
and $\Vert\mf{F}(\mf{u}_k)-\mf{F}(\mf{v}_k)\Vert\geq \varepsilon$. Notice that $\{\mf{u}_k\}$ and $\{\mf{v}_k\}$ are sequences in $\mk{D}$. Eq. \eqref{230724_13} implies that $\di \lim_{k\to\infty} \left(\mf{u}_k-\mf{v}_k \right)=\mf{0}$. Since  $\Vert\mf{F}(\mf{u}_k)-\mf{F}(\mf{v}_k)\Vert\geq \varepsilon$, \[\lim_{k\to\infty}\left(\mf{F}(\mf{u}_k)-\mf{F}(\mf{v}_k)\right)\neq \mf{0}.\]This shows that if (i) does not hold, then (ii) does not hold.
\end{myproof}

From Theorem \ref{230724_12}, we can deduce the following.
\begin{proposition}{}
Let $\mk{D}$ be a subset of $\mb{R}^n$, and let $\mf{F}:\mk{D}\rightarrow\mathbb{R}^m$ be a function defined on $\mk{D}$. Then  $\mf{F}:\mk{D}\rightarrow\mathbb{R}^m$ is uniformly continuous if and only if each of the component functions $F_j=(\pi_j\circ\mf{F}):\mk{D}\rightarrow\mathbb{R}$, $1\leq j\leq m$, is uniformly continuous.
\end{proposition}

Let us look at some more  examples.
\begin{example}[label=230724_14]{}Let $\mk{D}$ be the open rectangle $\mk{D}=(0, 5)\times (0, 7)$, and
consider the function $f:\mk{D}\to\mb{R}$ defined by
\[f(x,y)= xy.\]
Determine whether $f:\mk{D}\to\mb{R}$ is uniformly continuous. 
\end{example}
\begin{solution}{Solution}
For any two points $\mf{u}_1=(x_1,y_1)$ and $\mf{u}_2=(x_2,y_2)$ in $\mk{D}$, $0<x_1, x_2<5$ and $0<y_1, y_2<7$. Since
\[f(\mf{u}_1)-f(\mf{u}_2)= x_1y_1-x_2y_2
=x_1(y_1-y_2)+y_2(x_1-x_2),\]
we find that
\begin{align*}
\left|f(\mf{u}_1)-f(\mf{u}_2)\right|&\leq |x_1||y_1-y_2|+|y_2||x_1-x_2|\\&\leq 5\Vert\mf{u}_1-\mf{u}_2\Vert+7\Vert\mf{u}_1-\mf{u}_2\Vert=12\Vert\mf{u}_1-\mf{u}_2\Vert.\end{align*}
This shows that $f:\mk{D}\to\mb{R}$ is a Lipschitz function. Hence, it is uniformly continuous. 
\end{solution}

\begin{example}[label=230724_15]{}Consider the function $f:\mb{R}^2\to\mb{R}$ defined by
\[f(x,y)= xy.\]
Determine whether $f:\mb{R}^2\to\mb{R}$ is uniformly continuous. 
\end{example}
\begin{solution}{Solution}
 For $k\in\mb{Z}^+$, let 
 \[\mf{u}_k=\left(k+\frac{1}{k}, k\right), \hspace{1cm}\mf{v}_k=(k,k).\]
 Then $\{\mf{u}_k\}$ and $\{\mf{v}_k\}$ are sequences of points in $\mb{R}^2$ and 
 \[\lim_{k\to \infty}\left(\mf{u}_k-\mf{v}_k\right)=\lim_{k\to\infty}\left(\frac{1}{k}, 0\right)=(0,0).\] 
 However,
 \[f(\mf{u}_k)-f(\mf{v}_k)=k\left(k+\frac{1}{k}\right) -k^2=1.\]
 \bs
 Thus,
 \[\lim_{k\to\infty}\left(f(\mf{u}_k)-f(\mf{v}_k)\right)=1\neq 0.\]
 Therefore, the function $f:\mb{R}^2\to\mb{R}$ is not uniformly continuous. 
\end{solution}

Example \ref{230724_14} and \ref{230724_15} show that whether a function is uniformly continuous depends on the domain of the function.

\vp
\noindent
{\bf \large Exercises  \thesection}
\setcounter{myquestion}{1}
\begin{question}{\themyquestion}
Let $\mf{F}:\mb{R}^3\to\mb{R}^2$ be the function defined as
\[\mf{F}(x,y,z)=(3x-2z+7, x+y+z-4).\]
Show that $\mf{F}:\mb{R}^3\to\mb{R}^2$ is uniformly continuous.
\end{question}
\atc
 \begin{question}{\themyquestion}
Let $\mk{D}=(0,1)\times (0,2)$. Consider the function $f:\mk{D}\to\mb{R}$ defined as
\[f(x,y)=x^2+3y.\]
Determine whether $f$ is uniformly continuous. 
\end{question}
 
\atc
  \begin{question}{\themyquestion}
Let $\mk{D}=(1,\infty)\times (1,\infty)$. Consider the function $f:\mk{D}\to\mb{R}$ defined as
\[f(x,y)=\sqrt{x+y}.\]
Determine whether $f$ is uniformly continuous. 
\end{question}
\atc
 \begin{question}{\themyquestion}
Let $\mk{D}=(0,1)\times (0,2)$. Consider the function $f:\mk{D}\to\mb{R}$ defined as
\[f(x,y)=\frac{1}{\sqrt{x+y}}.\]
Determine whether $f$ is uniformly continuous. 
\end{question}

\section{Contraction Mapping Theorem}

Among the Lipschitz functions, there is a subset called contractions.
\begin{definition}{Contractions}
Let $\mk{D}$ be a subset of $\mathbb{R}^n$. A function $\mf{F}:\mk{D}\to \mathbb{R}^m$ is called a contraction if there exists a  constant $0\leq c<1 $ such that
\[\Vert\mf{F}(\mf{u})-\mf{F}(\mf{v})\Vert\leq c\Vert\mf{u}-\mf{v}\Vert\hspace{1cm}\text{for all}\;\mf{u}, \mf{v}\in\mk{D}.\]
In other words, a contraction is a Lipschitz function which has a Lipschitz constant that is less than 1.
\end{definition}
  
\begin{example}{}
Let $\mf{b}$ be a point in $\mb{R}^n$, and let $\mf{F}:\mb{R}^n\to\mb{R}^n$ be the function defined as
\[\mf{F}(\mf{x})=c\mf{x}+\mf{b}.\]
The mapping $\mf{F}$ is a contraction if and only if $|c|<1$.
\end{example}

The contraction mapping theorem is an important result in analysis. Extended  to metric spaces, it is an important tool to prove the existence and uniqueness of solutions of ordinary differential equations.

\begin{theorem}{Contraction Mapping Theorem}
Let $\mk{D}$ be a closed   subset of $\mb{R}^n$, and let $\mf{F}:\mk{D}\to \mk{D}$ be a contraction. Then $\mf{F}$ has a unique fixed point. Namely, there is a unique $\mf{u}$ in $\mk{D}$ such that $\mf{F}(\mf{u})=\mf{u}$.  
\end{theorem}
\begin{myproof}{Proof}By definition, there is a constant $c\in [0,1)$ such that 
\[\Vert\mf{F}(\mf{u})-\mf{F}(\mf{v})\Vert\leq c\Vert\mf{u}-\mf{v}\Vert\hspace{1cm}\text{for all}\;\mf{u}, \mf{v}\in \mk{D}.\]
\bp
We start with any point $\mf{x}_0$ in $\mk{D}$ and construct the sequence $\{\mf{x}_k\}$ inductively by
\[\mf{x}_{k+1}=\mf{F}(\mf{x}_k)\hspace{1cm} \text{for all}\;k\geq 0.\]
  
Notice that for all $k\in \mb{Z}^+$,
\[\Vert\mf{x}_{k+1}-\mf{x}_k\Vert=\Vert\mf{F}(\mf{x}_{k})-\mf{F}(\mf{x}_{k-1})\Vert\leq c\Vert \mf{x}_k-\mf{x}_{k-1}\Vert.\]
By iterating, we find that
\[\Vert\mf{x}_{k+1}-\mf{x}_k\Vert\leq c^k\Vert\mf{x}_1-\mf{x}_0\Vert.\]
Therefore, if $l> k\geq 0$,  triangle inequality implies that
\begin{align*}
\Vert\mf{x}_l-\mf{x}_k\Vert &\leq\Vert\mf{x}_l-\mf{x}_{l-1}\Vert+\cdots +\Vert\mf{x}_{k+1}-\mf{x}_k\Vert\\
&\leq (c^{l-1}+\ldots+c^k)\Vert\mf{x}_1-\mf{x}_0\Vert.\end{align*}Since $c\in [0,1)$,
\[c^{l-1}+\ldots+c^k=c^k(1+c+\cdots+c^{l-k-1})<\frac{c^k}{1-c}.\]
Therefore, for all $l> k\geq 0$,
\[\Vert\mf{x}_l-\mf{x}_k\Vert <\frac{c^k}{1-c}\Vert\mf{x}_1-\mf{x}_0\Vert.\]
 
 Given $\varepsilon>0$, there exists a positive integer $K$ such that for all $k\geq K$, 
\[\frac{c^k}{1-c}\Vert\mf{x}_1-\mf{x}_0\Vert<\varepsilon.\]
 
This implies that for all $l>k\geq K$,
\[\Vert\mf{x}_l-\mf{x}_k\Vert<\varepsilon.\]
  \bp
In other words, we have shown that $\{\mf{x}_k\}$ is  a Cauchy sequence. Therefore, it converges to a point $\mf{u}$ in $\mb{R}^n$. Since $\mk{D}$ is closed, $\mf{u}$ is in $\mk{D}$.

Since $\mf{F}$ is continuous, the sequence $\{\mf{F}(\mf{x}_k)\}$ converges to $\mf{F}(\mf{u})$. But $\mf{F}(\mf{x}_{k})=\mf{x}_{k+1}$. Being a subsequence of $\{\mf{x}_k\}$, the sequence $\{\mf{x}_{k+1}\}$   converges to $\mf{u}$ as well. This shows that
\[\mf{F}(\mf{u})=\mf{u},\]which says that  $\mf{u}$ is a fixed point of $\mf{F}$.
Now if $\mf{v}$ is another point in $\mk{D}$ such that 
$\mf{F}(\mf{v})=\mf{v}$, then
\[\Vert \mf{u}-\mf{v}\Vert=\Vert\mf{F}(\mf{u})-\mf{F}(\mf{v})\Vert\leq c\Vert\mf{u}-\mf{v}\Vert.\]
Since $c\in [0,1)$, this can only be true if $\Vert\mf{u}-\mf{v}\Vert=0$, which implies that $\mf{v}=\mf{u}$. Hence, the fixed point of $\mf{F}$ is unique.
\end{myproof}

As an application of the contraction mapping theorem, we  prove the following.
\begin{theorem}[label=230808_2]{}
Let $r$ be a positive number and let $\mf{G}:B(\mf{0}, r)\to \mb{R}^n$ be a mapping such that $\mf{G}(\mf{0})=\mf{0}$, and 
\[\Vert \mf{G}(\mf{u})-\mf{G}(\mf{v})\Vert\leq\frac{1}{2}\Vert\mf{u}-\mf{v}\Vert\hspace{1cm}\text{for all}\;\mf{u}, \mf{v}\in B(\mf{0}, r).\]If $\mf{F}:B(\mf{0}, r)\to \mb{R}^n$ is the function defined as
\[ \mf{F}(\mf{x})=\mf{x}+\mf{G}(\mf{x}),\] 
then $\mf{F}$ is a one-to-one continuous mapping whose image  contains the open ball $B(\mf{0}, r/2)$. 
\end{theorem}
\begin{myproof}{Proof}
By definition, $\mf{G}$ is a contraction. Hence, it is continuous. Therefore, $\mf{F}:B(\mf{0}, r)\to \mb{R}^n$ is also continuous. If $\mf{F}(\mf{u})=\mf{F}(\mf{v})$, then
\[\mf{u}-\mf{v}=\mf{G}(\mf{v})-\mf{G}(\mf{u}).\]
 
Therefore,
\[\Vert\mf{u}-\mf{v}\Vert =\Vert\mf{G}(\mf{v})-\mf{G}(\mf{u})\Vert\leq \frac{1}{2}\Vert\mf{u}-\mf{v}\Vert.\] This implies that $\Vert\mf{u}-\mf{v}\Vert=0$, and thus, $\mf{u}=\mf{v}$. Hence, $\mf{F}$ is one-to-one.

Given $\mf{y}\in B(\mf{0}, r/2)$, 
let $r_1=2\Vert\mf{y}\Vert$. Then $r_1<r$. 
Consider the map $\mf{H}: CB(\mf{0}, r_1)\to \mb{R}^n$  defined as
\[\mf{H}(\mf{x})=\mf{y}-\mf{G}(\mf{x}).\]
For any $\mf{u}$ and $\mf{v}$ in $CB(\mf{0}, r_1)$,
\[\Vert\mf{H}(\mf{u})-\mf{H}(\mf{v})\Vert=\Vert\mf{G}(\mf{u})-\mf{G}(\mf{v})\Vert\leq \frac{1}{2}\Vert\mf{u}-\mf{v}\Vert.\]

Therefore, $\mf{H}$ is also a contraction. Notice that if $\mf{x}\in CB(\mf{0}, r_1)$,
\[\Vert\mf{H}(\mf{x})\Vert \leq \Vert\mf{y}\Vert+\Vert\mf{G}(\mf{x})-\mf{G}(\mf{0})\Vert\leq \frac{r_1}{2}+\frac{1}{2}\Vert\mf{x}\Vert\leq\frac{r_1}{2}+\frac{r_1}{2}=r_1.\]

Therefore, $\mf{H}$ is a contraction that maps the closed set $CB(\mf{0}, r_1)$ into itself.  By the contraction mapping theorem, there exists $\mf{u}$ in $CB(\mf{0}, r_1)$  such that $\mf{H}(\mf{u})=\mf{u}$. This gives
\[\mf{y}-\mf{G}(\mf{u}) =\mf{u},\]or equivalently,
\[\mf{y}=\mf{u}+\mf{G}(\mf{u})=\mf{F}(\mf{u}).\]
 In other words, we have shown that there exists $\mf{u}\in CB(\mf{0}, r_1)\subset B(\mf{0}, r)$ such that $\mf{F}(\mf{u})=\mathbf{y}$. This proves that the image of the map   $\mf{F}:B(\mf{0}, r)\to \mb{R}^n$    contains the open ball $B(\mf{0}, r/2)$. 
\end{myproof}

\vp
\noindent
{\bf \large Exercises  \thesection}
\setcounter{myquestion}{1}

 \begin{question}{\themyquestion}
 Let \[S^n=\left\{(x_1, \ldots, x_n, x_{n+1})\in\mb{R}^{n+1}\,|\,x_1^2+\cdots+x_n^2+x_{n+1}^2=1\right\}\] be the $n$-sphere, and
 let $\mf{F}:S^n\to S^n$ be a mapping such that
 \[\Vert\mf{F}(\mf{u})-\mf{F}(\mf{v})\Vert\leq \frac{2}{3}\Vert\mf{u}-\mf{v}\Vert\hspace{1cm}\text{for all}\;\mf{u}, \mf{v} \in S^n.\] Show that there is a unique $\mf{w}\in S^n$ such that $\mf{F}(\mf{w})=\mf{w}$.
\end{question}
\atc
\begin{question}{\themyquestion}
Let $r$ be a positive number, and let $c$ be a positive number less than 1. Assume that $\mf{G}:B(\mf{0}, r)\to \mb{R}^n$ is a mapping such that $\mf{G}(\mf{0})=\mf{0}$, and 
\[\Vert \mf{G}(\mf{u})-\mf{G}(\mf{v})\Vert\leq c\Vert\mf{u}-\mf{v}\Vert\hspace{1cm}\text{for all}\;\mf{u}, \mf{v}\in B(\mf{0}, r).\]If $\mf{F}:B(\mf{0}, r)\to \mb{R}^n$ is the function defined as
\[ \mf{F}(\mf{x})=\mf{x}+\mf{G}(\mf{x}),\] show that
  $\mf{F}$ is a one-to-one continuous mapping whose image  contains the open ball $B(\mf{0}, ar)$, where $a=1-c$. 
\end{question}

 \chapter{Continuous Functions on Connected Sets  and Compact Sets}\label{chapter3}
In volume I, we have seen that intermediate value theorem and extreme value theorem play important roles in analysis. In order to extend these two theorems to multivariable functions, we need to consider two  topological properties of sets -- the connectedness and compactness.

\section{Path-Connectedness and Intermediate Value Theorem} 

We want to extend the intermediate value theorem to multivariable functions. For this, we need to consider a topological property called connectedness. In this section, we will discuss the topological property called path-connectedness first, which is a more natural concept.

\begin{definition}{Path}
Let $S$ be a subset of $\mb{R}^n$, and let $\mf{u}$ and $\mf{v}$ be two points in $S$. A {\it path in $S$ joining $\mf{u}$ to $\mf{v}$} is a {\bf continuos} function $\boldsymbol{\boldsymbol{\gamma}}:[a,b]\to S$ such that $\boldsymbol{\boldsymbol{\gamma}}(a)=\mf{u}$ and $\boldsymbol{\boldsymbol{\gamma}}(b)=\mf{v}$. 
\end{definition}
For any real numbers $a$ and $b$ with $a<b$, the map $u:[0,1]\to [a,b]$ defined by
\[u(t)=a+t(b-a)\] is a continuous bijection. Its inverse $u^{-1}:[a,b]\to [0,1]$ is
\[u^{-1}(t)=\frac{t-a}{b-a},\] which is also continuous. Hence, in the definition of a path, we can let the domain be any $[a, b]$ with $a<b$. 
 \begin{figure}[ht]
\centering
\includegraphics[scale=0.18]{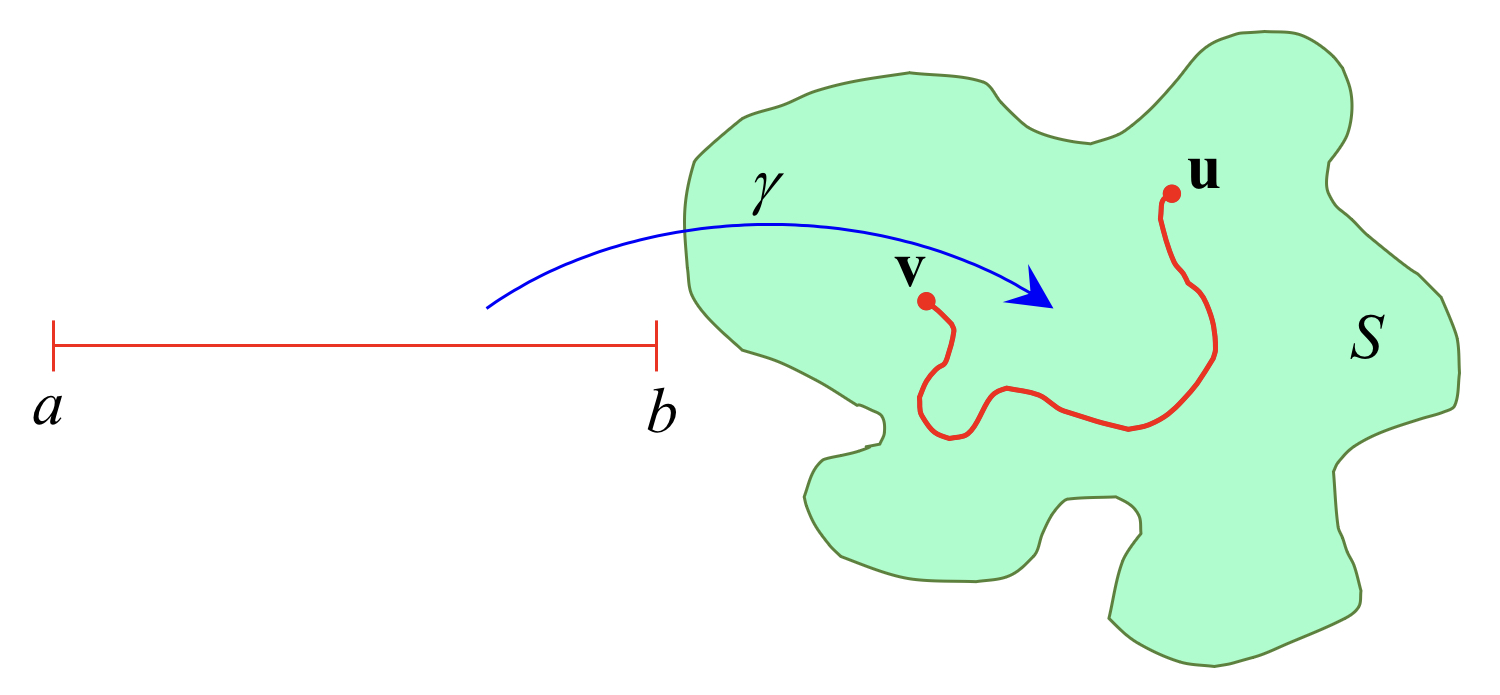}

\caption{A path in $S$ joining $\mf{u}$ to $\mf{v}$.}\label{figure21}
\end{figure}
\begin{example}{}
Given a set $S$ and a point $\mathbf{x}_0$ in $S$, the constant function $\boldsymbol{\boldsymbol{\gamma}}:[a,b]\to S$, $\boldsymbol{\boldsymbol{\gamma}}(t)=\mathbf{x}_0$, is a path in $S$.
\end{example}
If $\boldsymbol{\boldsymbol{\gamma}}:[a,b]\to S$ is a path in $S\subset\mb{R}^n$, and $S'$ is any other subset of $\mb{R}^n$ that contains the image of $\boldsymbol{\boldsymbol{\gamma}}$, then $\boldsymbol{\boldsymbol{\gamma}}$ is also a path in $S'$.
\begin{example}[label=230721_1]{}Let $R$ be the rectangle $R=[-2, 2]\times [-2, 2]$. 
The function $\boldsymbol{\boldsymbol{\gamma}}:[0,1]\to\mb{R}^2$, $\boldsymbol{\boldsymbol{\gamma}}(t)=(\cos(\pi t), \sin(\pi t))$ is a path in $R$ joining $\mathbf{u}=(1,0)$ to $\mf{v}=(-1,0)$.
\end{example}

 \begin{figure}[ht]
\centering
\includegraphics[scale=0.2]{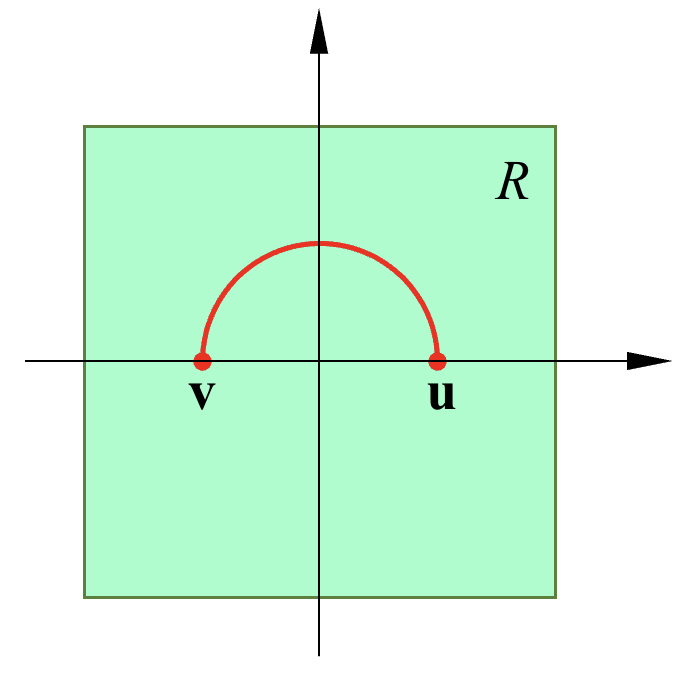}

\caption{The path in Example \ref{230721_1}.}\label{figure22}
\end{figure}
\begin{example}{}
Let $S$ be a subset of $\mb{R}^n$. If $\boldsymbol{\boldsymbol{\gamma}}:[a,b]\to S$ is a path in $S$ joining $\mf{u}$ to $\mf{v}$, then $\widetilde{\boldsymbol{\boldsymbol{\gamma}}}:[-b,-a]\to S$, $\widetilde{\boldsymbol{\boldsymbol{\gamma}}}(t)=\boldsymbol{\boldsymbol{\gamma}}(-t)$, is a path in $S$ joining $\mathbf{v}$ to $\mf{u}$.
\end{example}

Now we define path-connectedness.
\begin{definition}{Path-Connected}
Let $S$ be a subset of $\mb{R}^n$. We say that $S$ is path-connected if any two points $\mf{u}$ and $\mf{v}$ in $S$ can be joined by a path in $S$.
\end{definition}
It is easy to characterize a path-connected subset of $\mathbb{R}$. In volume I, we have defined the concept of convex sets. A subset $S$ of $\mathbb{R}$ is a convex set provided that for any $u$ and $v$ in $S$ 
and any $t\in [0,1]$, $(1-t)u+tv$ is also in $S$. This is equivalent to if $u$ and $v$ are points in $S$
with $u<v$, all the points $w$ satisfying $u<w<v$ is also in $S$. We have shown that a subset $S$ of $\mb{R}$ is a convex set if and only if it is an interval. 

The following theorem characterize a path-connected subset of $\mathbb{R}$.
\begin{theorem}[label=230721_8]{}
Let $S$ be a  subset of $\mathbb{R}$. Then $S$ is path-connected if and only if $S$ is an interval. 
\end{theorem}
\begin{myproof}{Proof}
If $S$ is an interval, then for any $u$ and $v$ in $S$, and for any $t\in [0,1]$, $(1-t)u+tv$ is in $S$. Hence, the function $ \gamma:[0,1]\to S$, $\gamma(t)=(1-t)u+tv$ is a path in $S$ that joins $u$ to $v$.

Conversely, assume that $S$ is a path-connected subset of $\mathbb{R}$. To show that $S$ is an interval, we need to show that for any $u$ and $v$   in $S$ with $u<v$,   any $w$ that is in the interval $[u,v]$ is also in $S$. Since $S$ is path-connected, there is a path $\gamma:[0,1]\to S$ such that $\gamma(0)=u$ and $\gamma(1)=v$. Since $\gamma$ is continuous,  and $w$ is in between $\gamma(0)$ and $\gamma(1)$, intermediate value theorem implies that there is a $c\in [0,1]$ so that $\gamma(c)=w$.  Thus, $w$ is in $S$.  
\end{myproof}

To explore path-connected subsets of $\mathbb{R}^n$ with $n\geq 2$, we first extend the concept of convex sets to $\mathbb{R}^n$. Given two points $\mf{u}$ and $\mf{v}$ in $\mb{R}^n$, when $t$ runs through all the points in the interval $[0,1]$, $(1-t)\mathbf{u}+t\mathbf{v}$ describes all the points on the line segment between $\mathbf{u}$ and $\mf{v}$. 

\begin{definition}{Convex Sets}
Let $S$ be a subset of $\mb{R}^n$. We say that $S$ is convex if   for any  two points $\mf{u}$ and $\mf{v}$ in $S$, the line segment between $\mf{u}$ and $\mf{v}$ lies entirely in $S$. Equivalently, $S$ is convex provided that for any two points $\mf{u}$ and $\mf{v}$ in $S$, the point $(1-t)\mathbf{u}+t\mathbf{v}$ is in $S$ for any $t\in [0,1]$.
\end{definition} 

\begin{figure}[ht]
\centering
\includegraphics[scale=0.2]{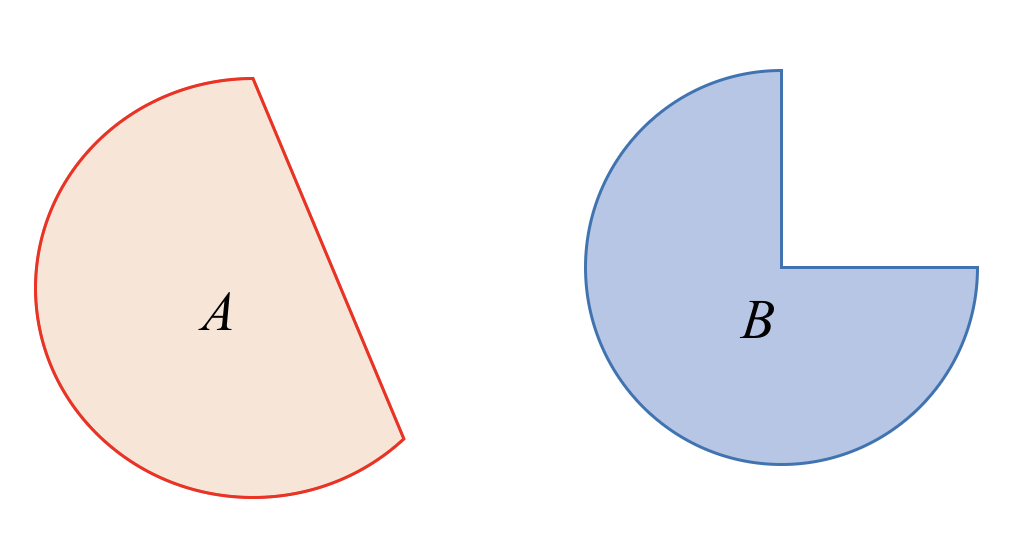}

\caption{$A$ is a convex set, $B$ is not.}\label{figure26}
\end{figure}

If $\mathbf{u}=(u_1, \ldots, u_n)$ and $\mathbf{v}=(v_1, \ldots, v_n)$ are two points in $\mb{R}^n$, the map $\boldsymbol{\gamma}:[0,1]\to\mathbb{R}^n$,
\[\boldsymbol{\gamma}(t)=(1-t)\mf{u}+t\mf{v}=((1-t)u_1+tv_1, \ldots, (1-t)u_n+tv_n)\] is a continuous functions, since each of its components is continuous. Thus, we have the following.
\begin{theorem}[label=230721_4]
{}
Let $S$ be a   subset of $\mathbb{R}^n$. If $S$ is convex, then it is path-connected.
\end{theorem}

Let us look at some examples of convex sets.
\begin{example}{}
Let $I_1$, $\ldots$, $I_n$ be intervals in $\mathbb{R}$. Show that the set $S=I_1\times \cdots\times I_n$ is path-connected.
\end{example}
\begin{solution}{Solution}
We claim that $S$ is convex. Then Theorem \ref{230721_4} implies that $S$ is path-connected.

Given that $\mathbf{u}=(u_1, \ldots, u_n)$ and $\mathbf{v}=(v_1, \ldots, v_n)$ are two points in $S$, for each $1\leq i\leq n$, $u_i$ and $v_i$ are in $I_i$. Since $I_i$ is an interval, for any $t\in [0,1]$, $(1-t)u_i+tv_i$ is in $I_i$. Hence,
\[(1-t)\mf{u}+t\mf{v}=((1-t)u_1+tv_1, \ldots, (1-t)u_n+tv_n)\]  is in $S$. This shows that $S$ is convex.
\end{solution}

Special cases of sets of the form $S=I_1\times \cdots\times I_n$ are open and closed rectangles.
\begin{example}{}
An open rectangle \[U=(a_1, b_1)\times \cdots\times (a_n, b_n)\] and its closure  \[R=[a_1, b_1]\times \cdots\times [a_n, b_n]\] are convex sets. Hence, they are path-connected. 
\end{example}

\begin{example}{}
Let $\mf{x}_0$ be a point in $\mb{R}^n$, and let $r$ be a positive number. Show that the open ball $B(\mathbf{x}_0, r)$ and the closed ball $CB(\mathbf{x}_0, r)$ are path-connected sets.
\end{example}
\begin{solution}{Solution}
Let $\mathbf{u}$ and $\mf{v}$ be two points in $B(\mathbf{x}_0, r)$. Then $\Vert\mf{u}-\mf{x}_0\Vert<r$ and $\Vert\mf{v}-\mf{x}_0\Vert<r$.
For any $t\in [0,1]$, $t\geq 0$ and $1-t\geq 0$. By triangle inequality,
\begin{align*}
\Vert(1-t)\mathbf{u}+t\mf{v}-\mf{x}_0\Vert&\leq \Vert(1-t)(\mf{u}-\mf{x}_0)\Vert+\Vert t(\mf{v}-\mf{x}_0)\Vert\\
&=(1-t)\Vert\mf{u}-\mf{x}_0\Vert+t\Vert\mf{v}-\mf{x}_0\Vert\\&<(1-t)r+tr=r.
\end{align*}
\bs
 
This shows that $(1-t)\mathbf{u}+t\mf{v}$ is in  $B(\mathbf{x}_0, r)$. Hence,  $B(\mathbf{x}_0, r)$ is convex. 

Replacing $<$ by $\leq$, one can show that $CB(\mathbf{x}_0, r)$ is convex.

By Theorem \ref{230721_4}, the open ball $B(\mathbf{x}_0, r)$ and the closed ball $CB(\mathbf{x}_0, r)$ are path-connected sets.
\end{solution}

Not all the path-connected sets  are   convex. Before we give an example, let us first prove the following useful lemma.
\begin{lemma}[label=230721_13]{}
Let $A$ and $B$ be path-connected subsets of $\mathbb{R}^n$. If $A\cap B$ is nonempty, then $S=A\cup B$ is path-connected.
\end{lemma}
\begin{myproof}{Proof}
Let  $\mf{u}$ and $\mf{v}$ be two points in $S$. If both  $\mf{u}$ and $\mf{v}$ are in the set $A$, then they can be joined by a path in $A$, which is also in $S$. Similarly, if both $\mf{u}$ and $\mf{v}$ are in the set $B$, then they can be joined by a path  in $S$. If $\mf{u}$ is in $A$ and $\mf{v}$ is in $B$, let $\mathbf{x}_0$ be any point in $A\cap B$. Then $\mf{u}$ and $\mf{x}_0$ are both in the path-connected set $A$, and $\mf{v}$ and $\mf{x}_0$ are both in the path-connected set $B$. Therefore, there exist continuous functions $\boldsymbol{\gamma}_1:[0,1]\to A$ and $\boldsymbol{\gamma}_2:[1,2]\to B$ such that
$\boldsymbol{\gamma}_1(0)=\mf{u}$, $\boldsymbol{\gamma}_1(1)=\mf{x}_0$, $\boldsymbol{\gamma}_2(1)=\mf{x}_0$ and $\boldsymbol{\gamma}_2(2)=\mf{v}$. Define the function $\boldsymbol{\gamma}:[0,2]\to A\cup B$ by
\[\boldsymbol{\gamma}(t)=\begin{cases}\boldsymbol{\gamma}_1(t),\quad &\text{if}\; 0\leq t\leq 1,\\\boldsymbol{\gamma}_2(t),\quad &\text{if}\; 1\leq t\leq 2.\end{cases}\]

Since $[0,1]$ and $[1,2]$ are closed subsets of $\mb{R}$, the function $\boldsymbol{\gamma}:[0,2]\to S$ is continuous. Thus, $\boldsymbol{\gamma}$ is a path in $S$ from $\mf{u}$ to $\mf{v}$. This proves that $S$ is path-connected.
\end{myproof}

\begin{figure}[ht]
\centering
\includegraphics[scale=0.18]{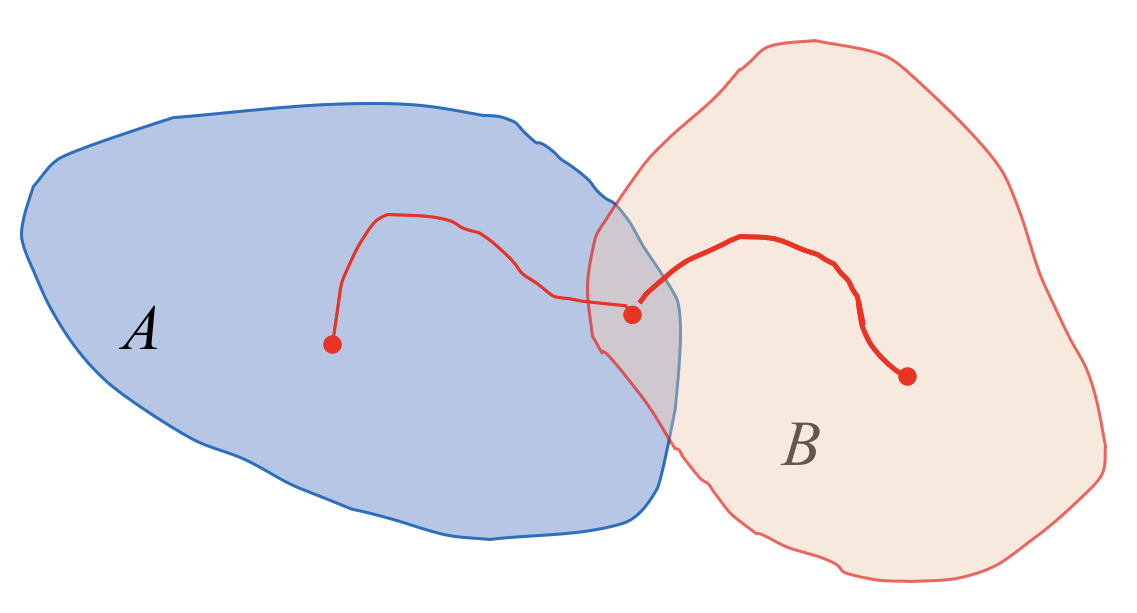}

\caption{If two sets $A$ and $B$ are path-connected and $A\cap B$ is nonempty, then $A\cup B$ is also path-connected.}\label{figure27}
\end{figure}

Now we can give an example of a path-connected set that is not convex.
\begin{example}{}
Show that the set 
\[S=\left\{(x,y)\,|\, 0\leq x\leq 1, -2\leq y\leq 2\right\}\cup
\left\{(x,y)\,|\, (x-2)^2+y^2\leq 1\right\}\] is path-connected, but not convex.

\end{example}
\begin{solution}{Solution}
The set \[A=\left\{(x,y)\,|\, 0\leq x\leq 1, -2\leq y\leq 2\right\}=[0,1]\times [-2,2]\]is a closed rectangle. Therefore, it is path- connected. The set
\[B=
\left\{(x,y)\,|\, (x-2)^2+y^2\leq 1\right\}\]
 
is a closed ball with center at $(2,0)$ and radius 1. Hence, it is also path-connected. Since the point $\mf{x}_0=(1,0)$ is in both $A$ and $B$, $S=A\cup B$ is path-connected.

The points $\mathbf{u}=(1,2)$ and $\mf{v}=(2,1)$ are in $S$. Consider the point
\[\mathbf{w}=\frac{1}{2}\mf{u}+\frac{1}{2}\mf{v}=\left(\frac{3}{2},\frac{3}{2}\right).\]
It is not    in $S$. This shows that $S$ is not convex.
\end{solution}

 \begin{figure}[ht]
\centering
\includegraphics[scale=0.2]{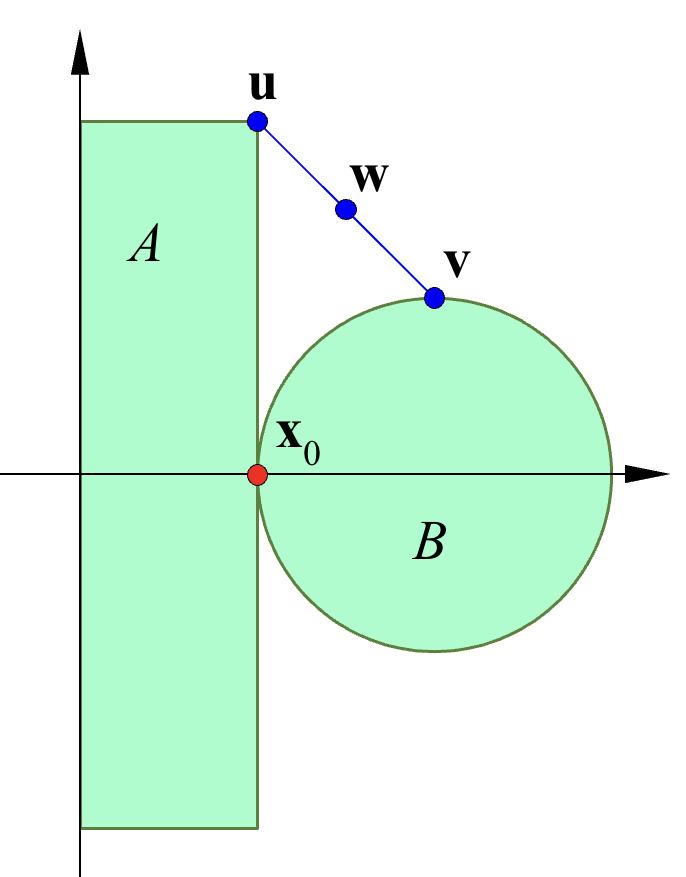}

\caption{The set $A\cup B$ is path-connected but not convex.}\label{figure23}
\end{figure}

Let us  now prove the following important theorem which says that continuous functions preserve path-connectedness.
\begin{theorem}[label=230721_7]{}
Let $\mk{D}$ be a path-connected subset of $\mathbb{R}^n$. If $\mf{F}:\mk{D}\to\mb{R}^m$ is a continuous function, then $\mf{F}(\mk{D})$ is path-connected.
\end{theorem}
\begin{myproof}{Proof}
Let $\mathbf{v}_1$ and $\mathbf{v}_2$ be two points in $\mf{F}(\mk{D})$. Then there exist $\mf{u}_1$ and $\mf{u}_2$ in $\mk{D}$ such that $\mf{F}(\mf{u}_1)=\mf{v}_1$ and $\mf{F}(\mf{u}_2)=\mf{v}_2$. Since $\mk{D}$ is path-connected, there is a continuous function $\boldsymbol{\gamma}:[0,1]\to\mk{D}$ such that $\boldsymbol{\gamma}(0)=\mf{u}_1$ and $\boldsymbol{\gamma}(1)=\mf{u}_2$. The map $\boldsymbol{\alpha}=(\mf{F}\circ\boldsymbol{\gamma}):[0,1]\to \mf{F}(\mk{D})$ is then a conitnuous map with $\boldsymbol{\alpha}(0)=\mathbf{v}_1$ and $\boldsymbol{\alpha}(1)=\mf{v}_2$. This shows that $\mf{F}(\mk{D})$ is path-connected.
\end{myproof}

From Theorem \ref{230721_7}, we obtain the following.
\begin{theorem}[label=230721_14]{Intermediate Value Theorem for Path-Connected Sets}
Let $\mk{D}$ be a path-connected subset of $\mathbb{R}^n$, and let $f:\mk{D}\to \mb{R}$ be a function defined on $\mk{D}$. If $f$ is continuous, then $f(\mk{D})$ is an interval. 
\end{theorem}
\begin{myproof}{Proof}
By Theorem \ref{230721_7}, $f(\mk{D})$ is a path-connected subset of $\mathbb{R}$. By Theorem \ref{230721_8}, $f(\mk{D})$ is an interval.
\end{myproof}

We can also use Theorem \ref{230721_7} to establish more examples of path-connected sets. Let us first look at an example.

\begin{example}{}
Show that the circle
\[S^1=\left\{(x,y)\,|\,x^2+y^2=1\right\}\] is path-connected.
\end{example}
\begin{solution}{Solution}
Define the function $f:[0,2\pi]\to\mb{R}^2$ by
\[f(t)=(\cos t, \sin t).\]Notice that $S^1=f([0,2\pi)]$.
Since each component of $f$ is a continuous function, $f$ is a continuous function. Since $[0,2\pi]$ is an interval, it is path-connected. By  Theorem \ref{230721_7}, $S^1=f([0,2\pi])$ is path-connected.
\end{solution}

A more general theorem is as follows.
\begin{theorem}[label=230721_12]{}
Let $\mk{D}$ be a path-connected subset of $\mathbb{R}^n$, and let $\mf{F}:\mk{D}\to\mb{R}^m$ be a function defined on $\mk{D}$. If $\mf{F}:\mk{D}\to\mb{R}^m$ is continuous, then the graph of $\mf{F}$, 
\[G_{\mf{F}}= \left\{(\mathbf{x}, \mathbf{y})\,|\, \mathbf{x}\in\mk{D}, \mf{y}=\mf{F}(\mf{x})\right\} \] is  a path-connected subset of $\mb{R}^{n+m}$.
\end{theorem}
\begin{myproof}{Proof}
By Corollary \ref{230721_11}, the function $\mf{H}:\mk{D}\to \mb{R}^{n+m}$, $\mf{H}(\mf{x})=\left(\mf{x}, \mf{F}(\mf{x})\right)$, is continuous. Since $\mf{H}(\mk{D})=G_{\mf{F}}$, 
  Theorem \ref{230721_7} implies that $G_{\mf{F}}$ is  a path-connected subset of $\mb{R}^{n+m}$.
\end{myproof}

Now let us consider spheres, which are boundary of balls.
\begin{definition}{The Standard Unit $\pmb{n}$-Sphere $\pmb{S^n}$}
A standard unit $n$-sphere $S^n$ is a subset of $\mb{R}^{n+1}$ consists of all points $\mathbf{x}=(x_1,   \ldots, x_n, x_{n+1})$ in $\mb{R}^{n+1}$ satisfying the equation $\Vert\mathbf{x}\Vert=1$, namely,
\[x_1^2+\cdots+x_n^2+x_{n+1}^2=1.\]
\end{definition}
The $n$-sphere $S^n$ is the boundary of the $(n+1)$ open ball $B^{n+1}=B(\mathbf{0}, 1)$ with center at the origin and radius 1.
\begin{figure}[ht]
\centering
\includegraphics[scale=0.2]{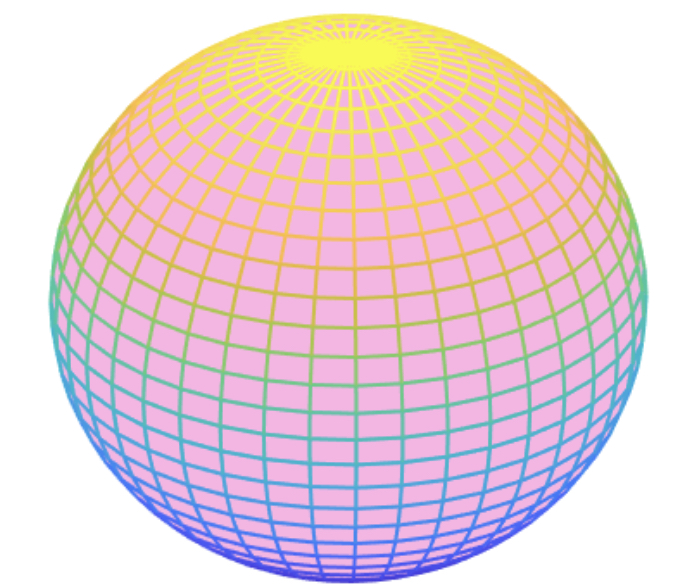}

\caption{A sphere.}\label{figure25}
\end{figure}

\begin{example}{}
Show that the standard unit $n$-sphere $S^n$ is path-connected.
\end{example}
\begin{solution}{Solution}
Notice that $S^n=S^n_+\cup S^n_-$, where $S^n_+$ and $S^n_-$ are respectively the upper and lower hemispheres with $x_{n+1}\geq 0$ and $x_{n+1}\leq 0$ respectively. 

\bs
If $\mathbf{x}\in S^n_+$, then 
\[x_{n+1}=\sqrt{1-x_1^2-\ldots-x_n^2};\]whereas if $\mathbf{x}\in S^n_-$,
\[x_{n+1}=-\sqrt{1-x_1^2-\ldots-x_n^2}.\]
Let 
\[CB^n=\left\{(x_1, \ldots, x_n)\,|\, x_1^2+\cdots+x_n^2\leq 1\right\}\]
 
be the closed ball in $\mathbb{R}^n$ with center at the origin and radius 1. Define the functions $f_{\pm}:CB^n\to\mb{R}$ by
\[f_{\pm}(x_1, \ldots, x_n)=\pm\sqrt{1-x_1^2-\ldots-x_n^2}.\]

Notice that $S^n_+$ and $S^n_-$ are respectively the graphs of $f_+$ and $f_-$. 
Since they are compositions of the square root function and a polynomial function, which are both continuous, $f_{+}$ and $f_-$ are continuous functions. 
The closed ball $CB^n$ is path-connected. Theorem \ref{230721_12} then implies that $S_+^n$ and $S_-^n$ are path-connected.

 Since both $S^n_+$ and $ S_-^n$ contain the unit vector $\mathbf{e}_1$ in $\mathbb{R}^{n+1}$, the set $S^n_+\cap S_-^n$ is nonempty. By Lemma \ref{230721_13}, $S^n=S^n_+\cup S_-^n$ is path-connected.
\end{solution}

\begin{remark}{}
There is an alternative way to prove that the $n$-sphere $S^n$ is path-connected. Given two distinct points $\mathbf{u}$ and $\mathbf{v}$ in $S^n$, they are unit vectors in $\mb{R}^{n+1}$. We want to show that there is a path in $S^n$ joining $\mathbf{u}$ to $\mathbf{v}$.

Notice that the line segment  $L=\left\{(1-t)\mathbf{u}+t\mathbf{v}\,|\,0\leq t\leq 1\right\}$  in $\mb{R}^{n+1}$ contains the origin if and only if $\mathbf{u}$ and $\mathbf{v}$ are parallel, if and only if $\mathbf{v}= -\mf{u}$. Thus, we discuss two cases.
\end{remark}
\begin{highlight}{}
\noindent \textbf{Case I:} $\mathbf{v}\neq -\mf{u}$. \\
In this case, let $\boldsymbol{\gamma}:[0,1]\to\mathbb{R}^{n+1}$ be the function defined as
\[\boldsymbol{\gamma}(t)=\frac{(1-t)\mathbf{u}+t\mathbf{v}}{\Vert (1-t)\mathbf{u}+t\mathbf{v}\Vert}.\]
Since $(1-t)\mathbf{u}+t\mathbf{v}\neq\mathbf{0}$ for all $0\leq t\leq 1$, $\boldsymbol{\gamma}$ is a continuous function. It is easy to check that its image lies in $S^n$. Hence, $\boldsymbol{\gamma}$ is a path in $S^n$ joining $\mathbf{u}$ to $\mathbf{v}$.

\noindent \textbf{Case 2:} $\mathbf{v}=-\mf{u}$.\\
In this case, let $\mathbf{w}$ be a unit vector orthogonal to $\mathbf{u}$, and let $\boldsymbol{\gamma}:[0,\pi]\to \mathbb{R}^{n+1}$ be the function defined as
\[\boldsymbol{\gamma}(t)=(\cos t) \mathbf{u}+(\sin t)\mathbf{w}.\]
Since $\sin t$ and $\cos t$ are continuous functions, $\boldsymbol{\gamma}$ is a continuous function.
Since $\mathbf{u}$ and $\mathbf{w}$ are orthogonal, the generalized Pythagoras theorem implies that
\[\Vert\boldsymbol{\gamma}(t)\Vert^2=\cos^2t\Vert\mathbf{u}\Vert^2+\sin^2 t\Vert\mathbf{w}\Vert^2=\cos^2t+\sin^2t=1.\]
Therefore, the image of $\boldsymbol{\gamma}$ lies in $S^n$. It is easy to see that $\boldsymbol{\gamma}(0)=\mf{u}$ and $\boldsymbol{\gamma}(\pi)=-\mf{u}=\mf{v}$. Hence, $\boldsymbol{\gamma}$ is a path in $S^n$ joining $\mathbf{u}$ to $\mathbf{v}$.
\end{highlight}

\begin{example}{}
Let $f:S^n\to \mb{R}$ be a continuous function.
Show that there is a point $\mf{u}_0$ on $S^n$ such that $f(\mf{u}_0)=f(-\mf{u}_0)$.
\end{example}
\begin{solution}{Solution}
The function $\mf{g}:\mb{R}^{n+1}\to\mb{R}^{n+1}$, $\mf{g}(\mf{u})=-\mf{u}$ is a linear transformation. Hence, it is continuous. Restricted to $S^n$, $\mf{g}(S^n)=S^n$. Thus, the function $f_1:S^n\to \mb{R}$, $f_1(\mf{u})=f(-\mf{u})$, is also continuous. 
\bs
It follows that the function $h:S^n\to\mb{R}$ defined by 
\[h(\mf{u})=f(\mf{u})-f(-\mf{u})\] is continuous.
Notice that 
\[h(-\mf{u})=f(-\mf{u})-f(\mf{u})=-h(\mf{u}).\]
This implies that if the number $a$ is in the range of $h$, so does the number $-a$. Since the number 0 is in between $a$ and $-a$ for any $a$, and $S^n$ is path-connected, intermediate value theorem implies that the number 0 is also in the range of $h$. This means that there is an $\mf{u}_0$ on $S^n$ such that $h(\mf{u}_0)=0$. Equivalently, $f(\mf{u}_0)=f(-\mf{u}_0)$.
\end{solution}

Theorem \ref{230721_14} says that a continuous function defined on a path-connected set satisfies the intermediate value theorem. We make the following definition.
\begin{definition}{Intermediate Value Property}
Let $S$ be a subset of $\mathbb{R}^n$. We say that $S$ has intermediate value property provided that whenever $f:S\to\mathbb{R}$ is a continuous function, then $f(S)$ is an interval.
\end{definition}
  Theorem \ref{230721_14} says that if $S$ is a path-connected set, then it has intermediate value property. It is natural to ask whether it is true that any set $S$ that has the intermediate value property must be path-connected. Unfortunately, it turns out that the answer is yes only when $S$ is  a subset of $\mathbb{R}$. If $S$ is a subset of $\mb{R}^n$ with $n\geq 2$, this is not true. This leads us to define a new property of sets called connectedness in the next section.

\vp
\noindent
{\bf \large Exercises  \thesection}
\setcounter{myquestion}{1}
 \begin{question}{\themyquestion}
 Is the set $A=(-1, 2)\cup(2, 5]$ path-connected? Justify your answer.
 \end{question}
 \atc
 \begin{question}{\themyquestion}
 Let $a$ and $b$ be positive numbers, and let $A$ be the subset of $\mathbb{R}^2$ given by
 \[A=\left\{(x,y)\,\left|\, \frac{x^2}{a^2}+\frac{y^2}{b^2}\leq 1\right.\right\}.\]
 Show that $A$ is convex, and deduce that it is path-connected.
 \end{question}
  \atc
 \begin{question}{\themyquestion}
 Let $(a,b,c)$ be a nonzero vector, and let $\mb{P}$ be the plane in $\mb{R}^3$ given by
 \[\mb{P}=\left\{(x,y,z)\,|\, ax+by+cz=d\right\},\]
 where $d$ is a constant.  Show that $\mb{P}$ is convex, and deduce that it is path-connected.
 \end{question}
 \atc
  \begin{question}{\themyquestion}
Let $S$ be the subset of $\mathbb{R}^3$ given by
 \[S=\left\{(x,y, z)\,|\, x>0, y\leq 1, 2\leq z<7\right\}.\]
 Show that $S$ is path-connected.
 \end{question}
  \atc
 \begin{question}{\themyquestion}
 Let $a$, $b$ and $c$ be positive numbers, and let $S$ be the subset of $\mathbb{R}^3$ given by
 \[S=\left\{(x,y, z)\,\left|\, \frac{x^2}{a^2}+\frac{y^2}{b^2}+\frac{z^2}{c^2}=1\right.\right\}.\]
 Show that $S$ is path-connected.
 \end{question}
  \atc

 \begin{question}{\themyquestion}
 Let $\mf{u}=(3, 0)$ and let $A$ be the subset of $\mb{R}^2$ given by
 \[A=\left\{(x,y)\,|\,x^2+y^2\leq 1\right\}.\]
 Define the function $f:A\to\mb{R}$ by
 $f(\mf{x})=d(\mf{x},\mf{u})$.
 \begin{enumerate}[(a)]
 \item Find $f(\mf{x}_1)$ and $f(\mf{x}_2)$, where $\mf{x}_1=(1,0)$ and $\mf{x}_2=(-1,0)$.
 \item Use intermediate value theorem to justify that there is a point $\mathbf{x}_0$ in $A$ such that $d(\mf{x}_0, \mf{u})=\pi$.
 \end{enumerate}
 \end{question}
 \atc
 \begin{question}{\themyquestion}
 Let $A$ and $B$ be subsets of $\mb{R}^n$. If $A$ and $B$ are convex, show that $A\cap B$ is also convex.
 \end{question}

\section{Connectedness and Intermediate Value Property}

In this section, we  study a property of sets  which is   known as connectedness. Let us first look at the path-connected subsets of $\mb{R}$ from a different perpective. We have shown in the previous section that a subset of $\mb{R}$ is path-connected if and only if it is an interval. A set of the form 
\[A=(-2,2]\setminus \{0\}=(-2,0)\cup (0,2]\] is not path-connected, since it contains the points $-1$ and $1$, but it does not contain the point 0 that is in between. Intuitively, there is no way to go from the point $-1$ to 1 {\it continuously} without leaving the set $A$. 

Let $U=(-\infty, 0)$ and $V=(0,\infty)$. Notice that $U$ and $V$ are open subsets of $\mb{R}$ which both intersect the set $A$. Moreover,
\[A=(A\cap U)\cup (A\cap V),\]
or equivalently,
\[A\subset U\cup V.\]
We  say that $A$ is {\it separated} by the open sets $U$ and $V$.

\begin{definition}{Separation of a Set}
Let $A$ be a subset of $\mb{R}^n$. A {\it separation} of $A$ is a pair $(U,V)$ of subsets of $\mb{R}^n$ which satisfies the following conditions.
\begin{enumerate}[(a)]
\item $U$ and $V$ are open sets.
\item $A\cap U\neq \emptyset$ and $A\cap V\neq \emptyset$.
\item $A\subset U\cup V$, or equivalently, $A$ is the union of $A\cap U$ and $A\cap V$.
\item $A$ is disjoint from $U\cap V$, or equivalently, $A\cap U$ and $A\cap V$ are disjoint.

\end{enumerate}
If $(U, V)$ is a separation of  $A$, we say that $A$ is separated by the open sets $U$ and $V$, or the open sets $U$ and $V$ separate $A$.
\end{definition}
\begin{example}{}
Let $A=(-2, 0)\cup (0, 2]$, and let $U=(-\infty, 0)$ and $V=(0, \infty)$. Then the open sets $U$ and $V$ separate $A$.

Let $U_1=(-3,0)$ and $V_1=(0, 3)$. The open sets $U_1$ and $V_1$ also separate $A$.
\end{example}

Now we define connectedness.
\begin{definition}{Connected Sets}
Let $A$ be a subset of $\mb{R}^n$. We say  that $A$ is connected if there does not exist a pair of open sets $U$ and $V$ that separate $A$.
\end{definition}

\begin{example}[label=230722_1]{}
Determine whether the set
\[A=\left\{(x,y)\,|\, y=0\right\}\cup \left\{(x,y)\,\left|\,y=\frac{2}{1+x^2}\right.\right\}\] is connected.
\end{example}
\begin{solution}{Solution}
Let $f:\mathbb{R}^2\to\mb{R}$ be the function defined as
\[f(x,y)=y(x^2+1).\]
Since $f$ is a polynomial function, it is continuous. The intervals $V_1=(-1,1)$ and $V_2=(1,3)$ are open sets in $\mb{R}$. Hence, the sets $U_1=f^{-1}(V_1)$ and $U_2=f^{-1}(V_2)$ are disjoint and they are open in $\mathbb{R}^2$. Notice that
\[A\cap U_1=\left\{(x,y)\,|\, y=0\right\},\hspace{1cm} A\cap U_2=\left\{(x,y)\,\left|\,y=\frac{2}{1+x^2}\right.\right\}. \]
Thus, $A\cap U_1$ and $A\cap U_2$ are nonempty, $A\cap U_1$ and $A\cap U_2$ are disjoint, and $A$ is a union of $A\cap U_1$ and $A\cap U_2$. This shows that the open sets $U_1$ and $U_2$ separate $A$. Hence, $A$ is not connected.
\end{solution}
 \begin{figure}[ht]
\centering
\includegraphics[scale=0.2]{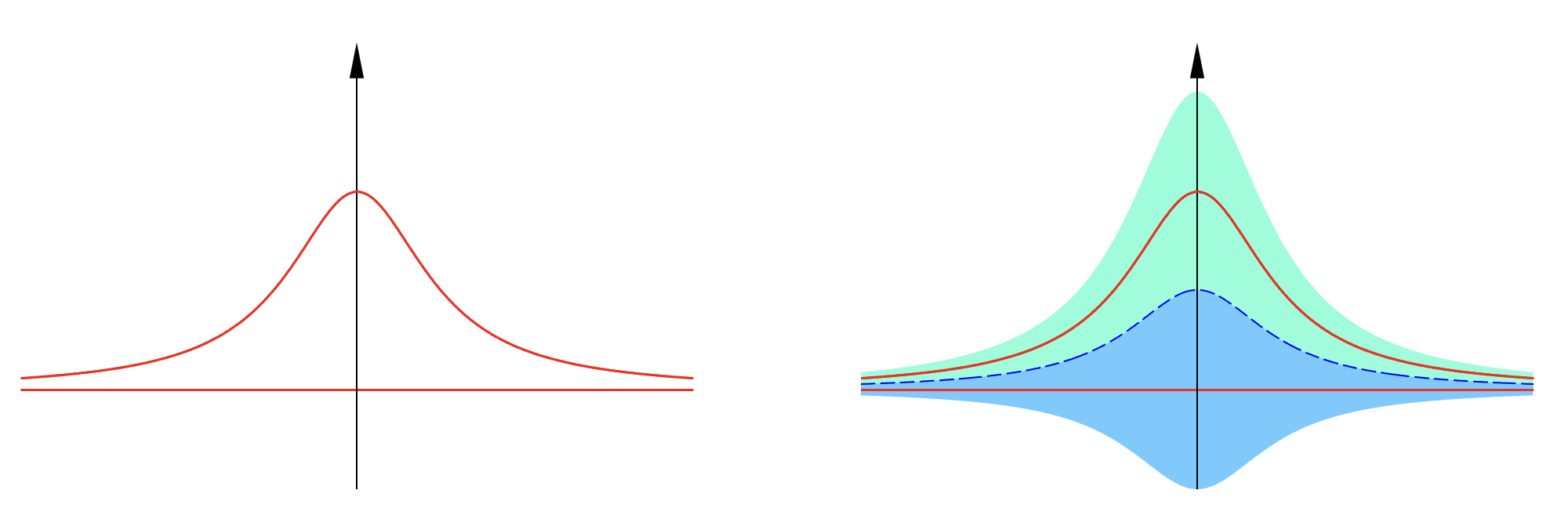}

\caption{The set $A$ defined in Example \ref{230722_1} is not connected.}\label{figure28}
\end{figure}

Now let us explore the relation between path-connected and connected. We first prove the following.
\begin{theorem}[label=230722_2]{}
Let $A$ be a subset of $\mb{R}^n$, and assume that the open sets $U$ and $V$ separate $A$. Define the function $f:A\to \mb{R}$ by
\[f(\mf{x})=\begin{cases}0,\quad &\text{if}\;\mf{x}\in A\cap U,\\1,\quad &\text{if}\;\mf{x}\in A\cap V.\end{cases}\]
Then $f$ is continuous.
\end{theorem}
Notice that the function $f$ is well defined since $A\cap U$ and $A\cap V$ are disjoint.
\begin{myproof}{Proof}
Let $\mf{x}_0$ be  a point in $A$. We want to prove that $f$ is continuous at $\mf{x}_0$. Since $A$ is contained in $U\cup V$,   $\mathbf{x}_0$ is in $U$ or in $V$. It suffices to consider the case where $\mf{x}_0$ is in $U$. The case where $\mf{x}_0$ is in $V$ is similar.

 If $\mathbf{x}_0$ is in $U$, since $U$ is open, there is an $r>0$ such that $B(\mf{x}_0, r)\subset U$. If $\{\mf{x}_k\}$ is a sequence in $A$ that converges $\mf{x}_0$, there exists a positive integer $K$ such that for all $k\geq K$, $\Vert\mf{x}_k-\mf{x}_0\Vert<r$. Thus, for all $k\geq K$, $\mf{x}_k\in B(\mf{x}_0, r)\subset U$, and hence, $f(\mf{x}_k)=0$. This proves that the sequence $\{f(\mf{x}_k)\}$ converges to $0$, which is  $f(\mf{x}_0)$. Therefore, $f$ is continuous at $\mf{x}_0$.
\end{myproof}

Now we can prove the theorem which says that a path-connected set is connected.
\begin{theorem}[label=230722_4]{}
Let $A$ be a subset of $\mb{R}^n$. If $A$ is path-connected, then it is connected.
\end{theorem}
\begin{myproof}{Proof}
We prove the contrapositive, which says that if $A$ is not connected, then it is not path-connected. 

If $A$ is not connected, there is a pair of open sets $U$ and $V$ that separate $A$. By Theorem \ref{230722_2}, the function $f:A\to \mb{R}$ defined by
\[f(\mf{x})=\begin{cases}0,\quad &\text{if}\;\mf{x}\in A\cap U,\\1,\quad &\text{if}\;\mf{x}\in A\cap V \end{cases}\]
is continuous. Since $f(A)=\{0,1\}$ is not an interval, by the contrapositive of the intermediate value theorem for path-connected sets, $A$ is not path-connected.
\end{myproof}

Theorem \ref{230722_4} provides us a large library of connected sets.
\begin{example}{}
The following sets are path-connected. Hence, they are also connected.
\begin{enumerate}[1.]
\item A set $S$ in $\mathbb{R}^n$ of the form $S=I_1\times \cdots\times I_n$, where $I_1, \ldots, I_n$ are intervals in $\mathbb{R}$.
\item Open rectangles and closed rectangles.
\item Open balls and closed balls.
\item The $n$-sphere $S^n$.

\end{enumerate}
\end{example}

The following theorem says that path-connectedness and connectedness are equivalent in $\mathbb{R}$.
\begin{theorem}[label=230722_8]{}
Let $S$ be a subset of $\mathbb{R}$. Then the following are equivalent.
\begin{enumerate}[(a)]
\item $S$ is an interval.
\item $S$ is path-connected.
\item $S$ is connected.
\end{enumerate}
\end{theorem}
\begin{myproof}{Proof}
We have proved (a) $\iff$ (b) in the previous section. In particular, (a) implies (b). Theorem \ref{230722_4} says that (b) implies (c). Now we only need to prove that (c) implies (a).

Assume that (a) is not true. Namely, $S$ is not an interval. Then there are points $u$ and $v$ in $S$ with $u<v$, such that there is a $w\in (u,v)$ that is not in $S$. Let $U=(-\infty, w)$ and $V=(w, \infty)$. Then $U$ and $V$ are disjoint open subsets of $\mb{R}$. Since $w\notin S$, $S\subset U\cup V$. Since $u\in S\cap U$ and $v\in S\cap V$, $S\cap U$ and $S\cap V$ are nonempty. Hence, $U$ and $V$ are open sets that separate $S$. This shows that $S$ is not connected. Thus, we have proved that if (a) is not true, then (c) is not true. This is equivalent to (c) implies (a).
\end{myproof}

Connectedness is also preserved by continuous functions. 
\begin{theorem}[label=230722_7]{}
Let $\mk{D}$ be a  connected subset of $\mathbb{R}^n$. If $\mf{F}:\mk{D}\to\mb{R}^m$ is a continuous function, then $\mf{F}(\mk{D})$ is   connected.
\end{theorem}
\begin{myproof}{Proof}
 We prove the contra-positive.
 Assume that $\mf{F}(\mk{D})$ is not connected. Then there are open sets $V_1$ and $V_2$ in $\mathbb{R}^m$ that separate $\mf{F}(\mk{D})$.   Let
 \begin{align*}
 \mk{D}_1&=\left\{\mf{x}\in\mk{D}\,|\,\mf{F}(\mf{x})\in V_1\right\},\\
 \mk{D}_2&=\left\{\mf{x}\in\mk{D}\,|\,\mf{F}(\mf{x})\in V_2\right\}.
 \end{align*}
 \bp
Since $\mf{F}(\mk{D})\cap V_1$ and $\mf{F}(\mk{D})\cap V_2$ are nonempty, $\mk{D}_1$ and $\mk{D}_2$ are nonempty. 
Since $\mf{F}(\mk{D})\subset V_1\cup V_2$, $\mk{D}=\mk{D}_1\cup\mk{D}_2$. Since $V_1\cap V_2$ is disjoint from $\mf{F}(\mk{D})$, $\mk{D}_1$ and $\mk{D}_2$ are disjoint. However, $\mk{D}_1$  and $\mk{D}_2$ are not necessary open sets. We will define two open sets $U_1$ and $U_2$ in $\mathbb{R}^n$ such that $\mk{D}_1=\mk{D}\cap U_1$ and $\mk{D}_2=\mk{D}\cap U_2$. Then $U_1$ and $U_2$ are open sets that separate $\mk{D}$.  
 
 For each $\mf{x}_0$ in $\mk{D}_1$, $\mf{F}(\mf{x}_0)\in  V_1$. Since $V_1$ is open, there exists $\varepsilon_{\mf{x}_0}>0$ such that the ball $B( \mf{F}(\mf{x}_0), \varepsilon_{\mf{x}_0})$ is contained in $V_1$. By the continuity of $\mathbf{F}$ at $\mathbf{x}_0$, there exists $\delta_{\mf{x}_0}>0$ such that for all $\mathbf{x}$ in $\mk{D}$, if $\mf{x}\in B(\mf{x}_0, \delta_{\mf{x}_0})$, then $\mf{F}(\mf{x})\in B( \mf{F}(\mf{x}_0), \varepsilon_{\mf{x}_0})\subset V_1$. In other words, 
 \[\mk{D}\cap B(\mf{x}_0, \delta_{\mf{x}_0})\subset \mf{F}^{-1}(V_1)=\mk{D}_1.\]
 Notice that  $B(\mf{x}_0, \delta_{\mf{x}_0})$ is an open set. 
 Define
 \[U_1=\bigcup_{\mathbf{x}_0\in \mk{D}_1}B(\mf{x}_0, \delta_{\mf{x}_0}).\]
 Being a union of open sets, $U_1$ is open. Since
 \[\mk{D}\cap U_1 =\bigcup_{\mathbf{x}_0\in \mk{D}_1}\left(\mk{D}\cap B(\mf{x}_0, \delta_{\mf{x}_0})\right)\subset \mk{D}_1,\] and \[\mk{D}_1=\bigcup_{\mathbf{x}_0\in \mk{D}_1}\{\mf{x}_0\}\subset \bigcup_{\mathbf{x}_0\in \mk{D}_1}\left(\mk{D}\cap B(\mf{x}_0, \delta_{\mf{x}_0})\right)=\mk{D}\cap U_1,\]
 we find that
 $\mk{D}\cap U_1=\mk{D}_1$. Similarly,
 define  \[U_2=\bigcup_{\mathbf{x}_0\in \mk{D}_2}B(\mf{x}_0, \delta_{\mf{x}_0}).\]
 Then $U_2$ is an open set and  $\mk{D}\cap U_2=\mk{D}_2$.  This completes the construction of the open sets $U_1$ and $U_2$ that separate $\mk{D}$. Thus, $\mk{D}$ is not connected.
\end{myproof}

From Theorem \ref{230722_8} and Theorem \ref{230722_7}, we also have an intermediate value theorem for connected sets. 
\begin{theorem}[label=230722_9]{Intermediate Value Theorem for Connected Sets}
Let $\mk{D}$ be a   connected subset of $\mathbb{R}^n$, and let $f:\mk{D}\to \mb{R}$ be a function defined on $\mk{D}$. If $f$ is continuous, then $f(\mk{D})$ is an interval. 
\end{theorem}
\begin{myproof}{Proof}
By Theorem \ref{230722_7}, $f(\mk{D})$ is a  connected subset of $\mathbb{R}$. By Theorem \ref{230722_8}, $f(\mk{D})$ is an interval.
\end{myproof}

Now we can prove the following.
\begin{theorem}{}
Let $S$ be a subset of $\mathbb{R}^n$. Then $S$ is connected if and only if it has the intermediate value property. 
\end{theorem}
\begin{myproof}{Proof}
If $S$ is connected and $f:S\to\mb{R}$ is continuous,  Theorem \ref{230722_9} implies that $f(S)$ is an interval. Hence, $S$ has the intermediate value property.

If $S$ is not connected, Theorem \ref{230722_2} gives a continuous function $f:S\to\mb{R}$ such that $f(S)=\{0,1\}$ is not an interval. Thus, $S$ does not  have the intermediate value property.
\end{myproof}

To give an example of a connected set that is not path-connected, we need    a lemma.

\begin{lemma}[label=230722_11]{}
Let $A$ be a subset of $\mathbb{R}^n$  that is separated by the open sets $U$ and $V$. If   $C$ is a connected subset of $A$, then $C\cap U=\emptyset$ or $C\cap V=\emptyset$.
\end{lemma}
\begin{myproof}{Proof}
Since $C\subset A$, $C\subset U\cup V$, and $C$ is disjoint from $U\cap V$.
If $C\cap U\neq \emptyset$ and $C\cap V\neq \emptyset$, then the open sets $U$ and $V$ also separate $C$. This contradicts to $C$ is connected. Thus, we must have $C\cap U=\emptyset$ or $C\cap V=\emptyset$.
\end{myproof}

\begin{theorem}[label=230722_12]{}
Let $A$ be a connected subset of $\mathbb{R}^n$. If $B$ is a subset of $\mb{R}^n$ such that \[A\subset B\subset \overline{A},\] then $B$ is also connected.
\end{theorem}
\begin{myproof}{Proof}
If $B$ is not connected, there exist open sets $U$ and $V$ in $\mathbb{R}^n$ that separate $A$. Since $A$ is connected, Lemma \ref{230722_11}  says that $A\cap U=\emptyset$ or $A\cap V=\emptyset$. Without loss of generality, assume that $A\cap V=\emptyset$. Then $A\subset\mb{R}^n\setminus V$. Thus, $\mb{R}^n\setminus V$ is a closed set that contains $A$. This implies that $\overline{A}\subset \mb{R}^n\setminus V$. Hence, we also have $B\subset   \mb{R}^n\setminus V$, which contradicts to the fact that the set $B\cap V$ is not empty.
\end{myproof}
\begin{example}[label=230722_13]{The Topologist's Sine Curve}
Let $S$ be the subset of $\mathbb{R}^2$ given by $S=A\cup L$, where
\[A=\left\{(x,y)\,\left|\, 0<x\leq 1, y=\sin\left(\frac{1}{x}\right)\right.\right\},\]and 
\[L=\left\{(x,y)\,|\, x=0, -1\leq y\leq 1\right\}.\]
 
\begin{enumerate}[(a)]
 \item Show that $S\subset \overline{A}$.
\item Show that $S$ is connected.
\item Show that $S$ is not path-connected.
\end{enumerate} 
\end{example}
\begin{solution}{Solution}
\begin{enumerate}[(a)]
\item Since $A\subset \overline{A}$, it suffices to show that $L\subset \overline{A}$.  Given $(0,u)\in L$, $-1\leq u\leq 1$. Thus,  $a=\sin^{-1}u\in [-\pi/2, \pi/2]$. Let
\[x_k=\frac{1}{a+2\pi k}\hspace{1cm}\text{for}\;k\in\mb{Z}^+.\]\end{enumerate}\bs
\begin{enumerate}[(a)]
\item[]
Notice that $x_k\in (0,1]$ and \[\sin \di\frac{1}{x_k}=\sin a=u.\] Thus, $\{(x_k, \sin (1/x_k))\}$ is a sequence of points in $A$ that converges to $(0, u)$. This proves that $(0, u)\in \overline{A}$. Hence, $  L\subset \overline{A}$.  
\item[(b)]
 The interval $(0, 1]$ is path-connected and the function $f:(0, 1]\to\mb{R}$, $\di f(x)=\sin\left(\frac{1}{x}\right)$ is continuous. Thus, $A=G_f$ is path-connected, and hence it is connected. Since $A\subset S\subset\overline{A}$, Theorem \ref{230722_12} implies that $S$ is connected.
\item[(c)] If $S$ is path connected, there is a path $\boldsymbol{\gamma}:[0,1]\to S$ such that $\boldsymbol{\gamma}(0)=(0,0)$ and $\boldsymbol{\gamma}(1)=(1,\sin 1)$. Let $\boldsymbol{\gamma}(t)=(\gamma_1(t), \gamma_2(t))$. Then $\gamma_1:[0,1]\to\mb{R}$ and $\gamma_2:[0,1]\to\mb{R}$ are continuous functions. 
Consider the sequence $\{x_k\}$ with
\[x_k=\frac{1}{\di\frac{\pi}{2}+\pi k},\hspace{1cm}k\in\mb{Z}^+.\] Notice that $\{x_k\}$ is a decreasing sequence of points in $[0,1]$ that converges to $0$.  For each $k\in\mb{Z}^+$, 
  $(x_k, y_k)\in S$ if and onlly if $y_k=\sin(1/x_k)$.  
  
  Since $\gamma_1:[0,1]\to\mb{R}$ is continuous, $\gamma_1(0)=0$ and $\gamma_1(1)=1$, intermediate value theorem implies that there exists $t_1\in [0,1]$ such that $\gamma_1(t_1)=x_1$. Similarly, there exists $t_2\in [0, t_1]$ such that $\gamma_1(t_2)=x_2$. Continue the argument gives a decreasing sequence $\{t_k\}$ in $[0,1]$ such that $\di \gamma_1(t_k)=x_k$  for all $k\in\mb{Z}^+$. Since the sequence $\{t_k\}$ is bounded below, it converges to some $t_0$ in $[0,1]$. Since $\gamma_2:[0,1]\to \mathbb{R}$ is also continuous, the sequence $\{\gamma_2(t_k)\}$ should converge to $\gamma_2(t_0)$. 

Since  $\boldsymbol{\gamma}(t_k)\in S$ and $\gamma_1(t_k)=x_k$, we must have  $\gamma_2(t_k)=y_k=(-1)^k$. But then the sequence $\{\gamma_2(t_k)\}$ is not convergent. This gives a contradiction.
Hence, there does not exist a path in $S$ that joins the point $(0,0)$ to the point $(1,\sin 1)$. This proves that $S$ is not path-connected.
\end{enumerate}
\end{solution}

 \begin{figure}[ht]
\centering
\includegraphics[scale=0.2]{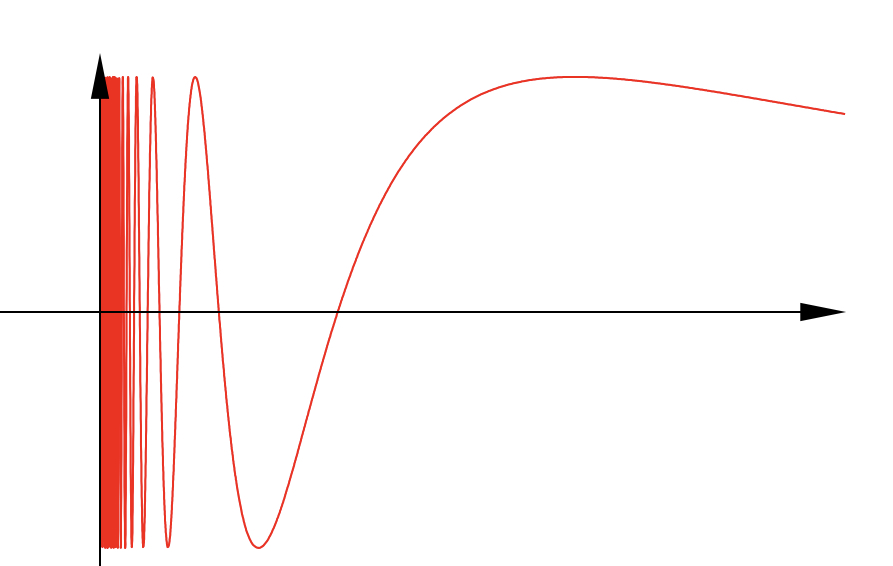}

\caption{The topologist's sine curve.}\label{figure29}
\end{figure}

\begin{remark}{}
Example \ref{230722_13} gives a set that is connected but not path-connected.
\begin{enumerate}[1.]
\item
One can in fact show that $S=\overline{A}$.
\item To show that $A$ is connected, we can also use the fact that if $\mk{D}$ is a connected subset of $\mathbb{R}^n$, and   $\mf{F}:\mk{D}\to\mb{R}^m$ is a continuous function, then the graph of $\mf{F}$ 
 is    connected. The proof of this fact is left as an exercise.
\end{enumerate}
\end{remark}

At the end of this section, we want to give a sufficient condition for a connected subset of $\mb{R}^n$ to be path-connected.

First we define the meaning of a polygonal path.
\begin{definition}{Polygonal Path}
Let $S$ be a subset of $\mb{R}^n$, and let $\mf{u}$ and $\mf{v}$ be two points in $S$. A path  $\boldsymbol{\gamma}:[a,b]\to S$ in $S$ that joins $\mf{u}$ to $\mf{v}$ is a polygonal path provided that there is a partition $P=\{t_0, t_1, \ldots, t_k\}$ of $[a,b]$ such that for $1\leq i\leq k$,
\[\boldsymbol{\gamma}(t)=\mf{x}_{i-1}+\frac{t-t_{i-1}}{t_i-t_{i-1}}\left(\mf{x}_i-\mf{x}_{i-1}\right), \quad \text{when}\;t_{i-1}\leq t\leq t_i.\]
\end{definition}
\begin{figure}[ht]
\centering
\includegraphics[scale=0.2]{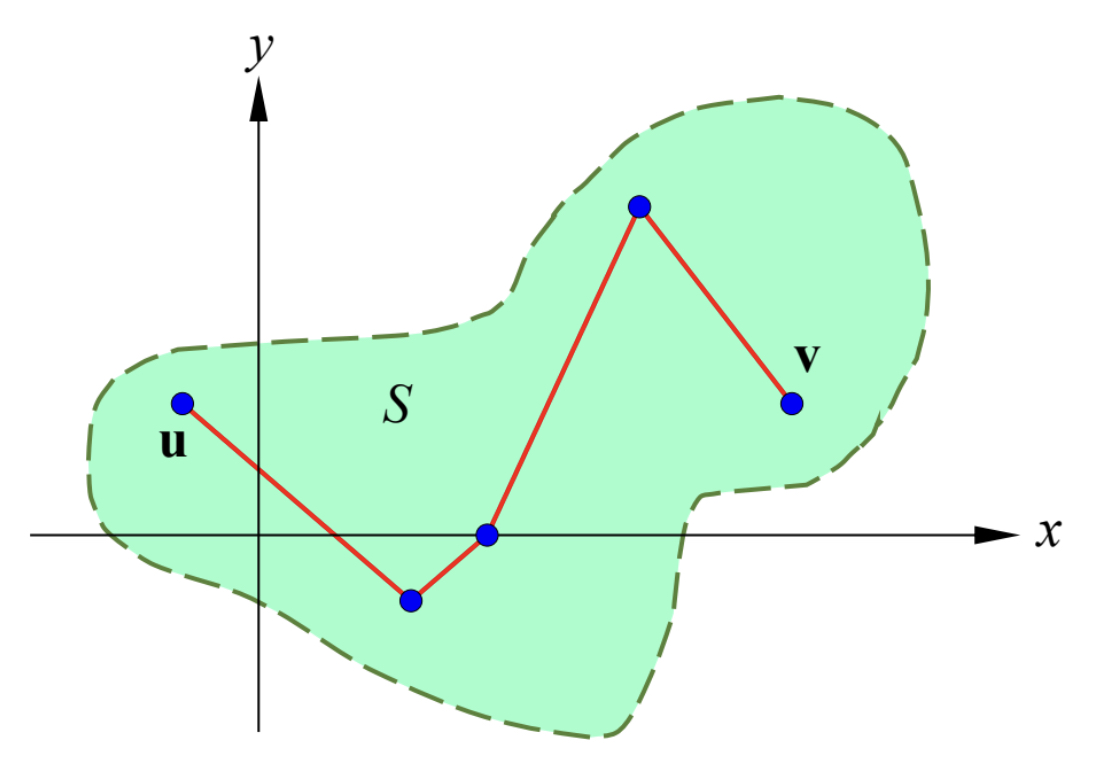}

\caption{A polygonal path.}\label{figure43}
\end{figure}

Obviously, we have the following.
\begin{proposition}{}
If $S$ is a convex subset of $\mb{R}^n$, then any two points in $S$ can be joined by a polygonal path in $\mb{R}^n$.

\end{proposition}
If $\boldsymbol{\gamma}_1:[a,c]\to A$ is a polygonal path in $A$ that joins $\mf{u}$ to $\mf{w}$, $\boldsymbol{\gamma}_2:[c,b]\to B$ is a polygonal path in $B$ that joins $\mf{w}$ to $\mf{v}$, then
the path $\boldsymbol{\gamma}:[a,b]\to A\cup B$,
\[\boldsymbol{\gamma}(t)=\begin{cases}\boldsymbol{\gamma}_1(t),\quad &\text{if}\;a\leq t\leq c,\\\boldsymbol{\gamma}_2(t),\quad &\text{if}\;c\leq t\leq b,\end{cases}\] is a polygonal path in $A\cup B$ that joins $\mf{u}$ to $\mf{v}$.
Using this, we can prove the following useful theorem.
\begin{theorem}[label=230802_2]{}
Let $S$ be a  connected subset of $\mb{R}^n$. If $S$ is an open set, then any two points in $S$ can be joined by a polygonal path in $S$. In particular, $S$ is path connected.
\end{theorem}
\begin{myproof}{Proof}
We use proof by contradiction. Supposed that $S$ is open but there are two points $\mf{u}$ and $\mf{v}$ in $S$ that cannot be joined by a polygonal path in $S$. Consider the sets
\[U=\left\{\mf{x}\in S\,|\, \text{there is a polygonal path in $S$ that joins $\mf{u}$ to $\mf{x}$}\right\},\] 
\[V=\left\{\mf{x}\in S\,|\, \text{there is no polygonal path in $S$ that joins $\mf{u}$ to $\mf{x}$}\right\}.\] Obviously $\mf{u}$ is in $U$ and $\mf{v}$ is in $V$, and $S=U\cup V$. We claim that both $U$ and $V$ are open sets. 
\bp

If $\mf{x}$ is in the open set $S$, there is an $r>0$ such that $B(\mf{x}, r)\subset S$.    Since $B(\mf{x},r)$ is convex, any   point $\mf{w}$ in $B(\mf{x}, r)$ can be joined by a polygonal path in $B(\mf{x}, r)$ to $\mf{x}$. Hence, if $\mf{x}$ is in $U$, $\mf{w}$ is in $U$. If $\mf{x}$ is in $V$, $\mf{w}$ is in $V$. This shows that if $\mf{x}$ is in $U$, then $B(\mf{x}, r)\subset U$. If $\mf{x}$ is in $V$, then $B(\mf{x}, r)\subset V$.
Hence, $U$ and $V$ are open sets.

Since $U$ and $V$ are nonempty open sets and 
$S=U\cup V$, they form a separation of $S$. This contradicts to $S$ is connected. Hence, any two points in $S$ can be joined by  a polygonal path in $S$.
\end{myproof}

\vp
\noindent
{\bf \large Exercises  \thesection}
\setcounter{myquestion}{1}

\begin{question}{\themyquestion}
Determine whether the set
\[A=\left\{(x,y)\,|\, y=0\right\}\cup \left\{(x,y)\,\left|\,x>0, y=\frac{2}{x}\right.\right\}\] is connected.
\end{question}

\atc
\begin{question}{\themyquestion}
Let $\mk{D}$ be a connected subset of $\mathbb{R}^n$, and let $\mf{F}:\mk{D}\to\mb{R}^m$ be a function defined on $\mk{D}$. If $\mf{F}:\mk{D}\to\mb{R}^m$ is continuous, show that the graph of $\mf{F}$, 
\[G_{\mf{F}}= \left\{(\mathbf{x}, \mathbf{y})\,|\, \mathbf{x}\in\mk{D}, \mf{y}=\mf{F}(\mf{x})\right\} \]is also  connected.
 \end{question}

 \atc
 \begin{question}{\themyquestion}
Determine whether the set
\[A=\left\{(x,y)\,|\,0\leq x<1, -1<y\leq 1\right\}\cup  \{(1,0), (1,1)\} \] is connected.
\end{question}
\atc
 \begin{question}{\themyquestion}
Assume that $A$ is a connected subset of $\mb{R}^3$ that contains the points $\mathbf{u}=(0,2,0)$ and $\mathbf{v}=(2,-6, 3)$.
\begin{enumerate}[(a)]
\item Show that there is a point $\mathbf{x}=(x,y,z)$ in $A$ that lies in the plane $y=0$.
\item Show that there exists a point $\mathbf{x}=(x,y,z)$ in $A$ that lies on the sphere $x^2+y^2+z^2=25$.

\end{enumerate}
\end{question}
 \atc
\begin{question}{\themyquestion}
Let $A$ and $B$ be  connected subsets of $\mathbb{R}^n$. If $A\cap B$ is nonempty, show that $S=A\cup B$ is  connected.  
 \end{question}

\section{Sequential Compactness and Compactness}

In volume I, we have seen that sequential compactness plays important role in extreme value theorem. 
In this section, we extend the definition of sequential compactness to subsets of $\mb{R}^n$. We will also consider another concept called compactness. 

Let us start with the definition of bounded sets.
\begin{definition}{Bounded Sets}
Let $S$ be a subset of $\mb{R}^n$. We say that $S$ is bounded if there exists a positive number $M$ such that
\[\Vert\mf{x}\Vert\leq M\hspace{1cm}\text{for all}\; \mf{x}\in S.\]
\end{definition}
\begin{remark}{}
Let $S$ be a subset of $\mb{R}^n$. If $S$ is bounded and $S'$ is a subset of $S$, then it is obvious  that $S'$ is also bounded.
\end{remark}

\begin{example}{}
Show that a ball $B(\mathbf{x}_0, r)$ in $\mathbb{R}^n$ is bounded.
\end{example}
\begin{solution}{Solution}
Given $\mf{x}\in B(\mathbf{x}_0, r)$, $\Vert\mf{x}-\mf{x}_0\Vert <r$. Thus,
\[\Vert\mf{x}\Vert\leq \Vert\mf{x}_0\Vert + \Vert\mf{x}-\mf{x}_0\Vert <\Vert\mf{x}_0\Vert +r.\]
Since $M=\Vert\mf{x}_0\Vert +r$ is a constant independent of the points in the  ball $B(\mathbf{x}_0, r)$, the ball $B(\mathbf{x}_0, r)$ is bounded.
\end{solution}

Notice that if $\mf{x}_1$ and $\mf{x}_2$ are points in $\mb{R}^n$, and $S$ is a set in $\mb{R}^n$ such that
\[\Vert\mf{x}-\mf{x}_1\Vert <r_1\hspace{1cm}\text{for all}\;\mf{x}\in S,\]
then
\[\Vert\mf{x}-\mf{x}_2\Vert <r_1+\Vert\mf{x}_2-\mf{x}_1\Vert\hspace{1cm}\text{for all}\;\mf{x}\in S.\]
Thus, we have the following.
\begin{proposition}{}
Let $S$ be a subset in $\mb{R}^n$. The following are equivalent.
\begin{enumerate}[(a)]
\item
$S$ is bounded.
\item There is a point $\mf{x}_0$ in $\mb{R}^n$ and a positive constant $M$ such that
\[\Vert\mf{x}-\mf{x}_0\Vert \leq M\hspace{1cm}\text{for all}\;\mf{x}\in S.\]

\item For any $\mf{x}_0$ in $\mb{R}^n$, there is a positive constant $M$ such that
\[\Vert\mf{x}-\mf{x}_0\Vert \leq M\hspace{1cm}\text{for all}\;\mf{x}\in S.\]
\end{enumerate}
\end{proposition}

 \begin{figure}[ht]
\centering
\includegraphics[scale=0.2]{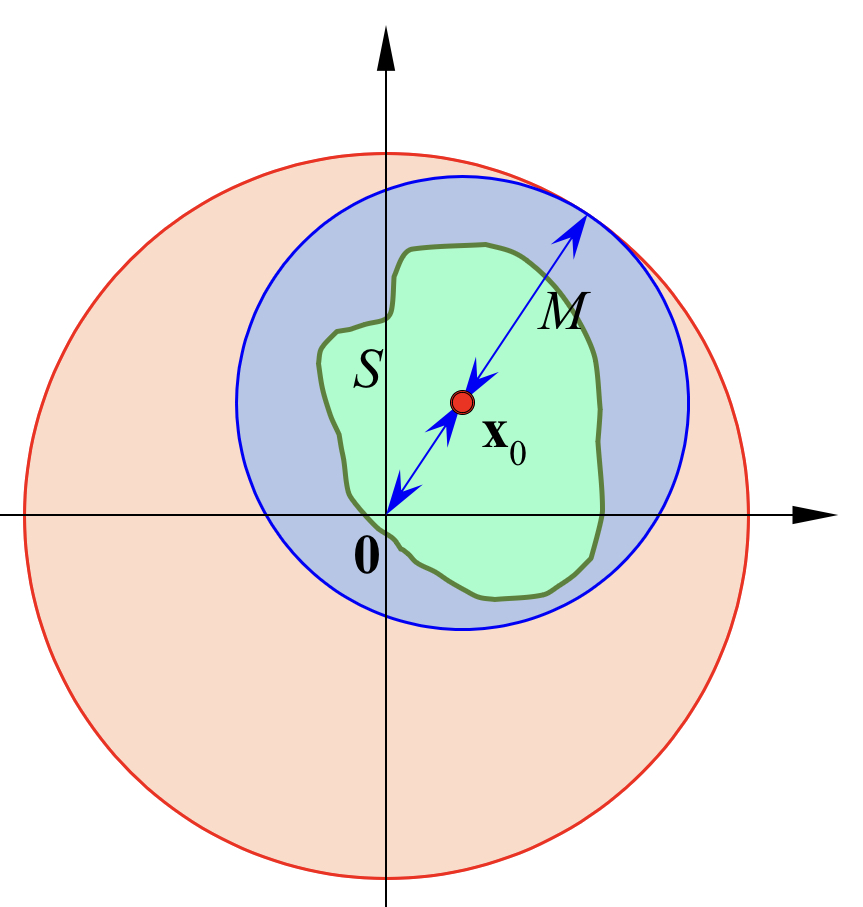}

\caption{The set $S$ is bounded.}\label{figure30}
\end{figure}
We say that a sequence $\{\mf{x}_k\}$ is bounded if the set $\{\mf{x}_k\,|\,k\in\mb{Z}^+\}$ is bounded. The following is a standard theorem about convergent sequences.

\begin{proposition}{}
If $\{\mathbf{x}_k\}$ is a sequence in $\mb{R}^n$ that is convergent, then it is bounded.
\end{proposition}
\begin{myproof}{Proof}
Assume that the sequence $\{\mf{x}_k\}$ converges to the point $\mf{x}_0$. Then there is a positive integer $K$ such that 
\[\Vert\mf{x}_k-\mf{x}_0\Vert<1\hspace{1cm}\text{for all}\;k\geq K.\]
Let
\[M=\max\{\Vert\mathbf{x}_k-\mf{x}_0\Vert\,|\, 1\leq k\leq K-1\}+1.\]
Then   $M$ is finite and 
\[\Vert\mf{x}_k-\mf{x}_0\Vert\leq M\hspace{1cm}\text{for all}\;k\in\mb{Z}^+.\]
Hence, the sequence $\{x_k\}$ is bounded.
\end{myproof}

\begin{figure}[ht]
\centering
\includegraphics[scale=0.2]{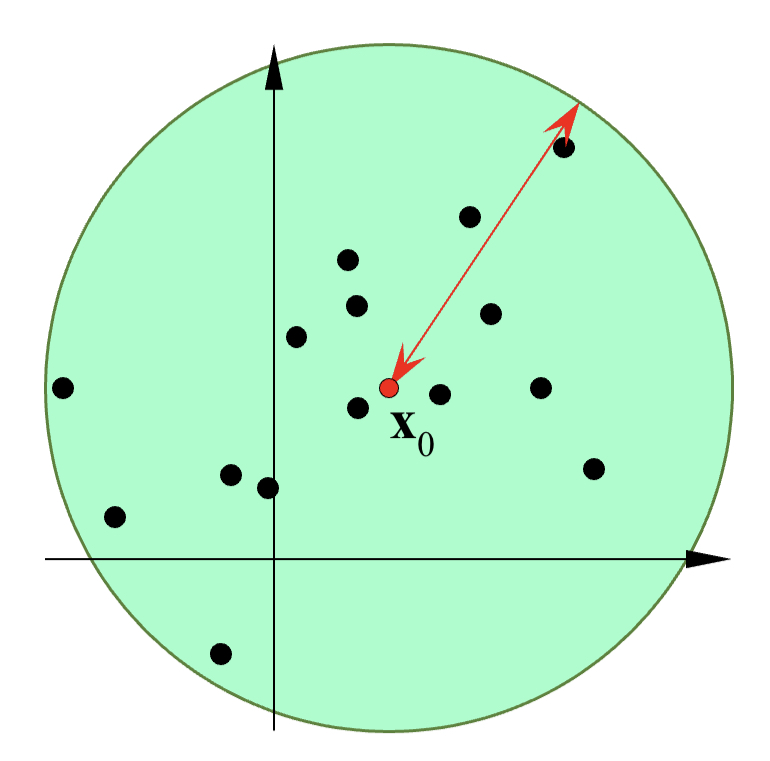}

\caption{A convergent sequence is bounded.}\label{figure31}
\end{figure}

Let us now define the {\it diameter} of a bounded set.
If $S$ is a subset of $\mathbb{R}^n$ that is bounded, there is a positive number $M$ such that
\[\Vert\mf{x}\Vert\leq M\hspace{1cm}\text{for all}\; \mf{x}\in S.\]It follows from triangle inequality that for any $\mf{u}$ and $\mf{v}$ in $S$,
\[\Vert\mf{u}-\mf{v}\Vert\leq\Vert\mf{u}\Vert+\Vert\mf{v}\Vert\leq 2M.\]
Thus, the set
\begin{equation}\label{230723_7}D_S=\left\{d(\mf{u},\mf{v})\,|\, \mf{u}, \mf{v}\in S\right\}=\left\{\Vert\mf{u}-\mf{v}\Vert\,|\, \mf{u}, \mf{v}\in S\right\}
\end{equation} is a set of  nonnegative real numbers that is bounded above. In fact, for any subset $S$ of $\mb{R}^n$, one can define the set of real numbers $D_S$ by \eqref{230723_7}.   Then $S$ is a bounded set if and only if the set $D_S$ is bounded above.

\begin{definition}{Diameter of a Bounded Set}
Let $S$ be a bounded subset of $\mb{R}^n$. The diameter of $S$, denoted by $\text{diam}\,S$, is defined as
\[\text{diam}\,S=\sup\left\{d(\mf{u},\mf{v})\,|\, \mf{u}, \mf{v}\in S\right\}=\sup\left\{\Vert\mf{u}-\mf{v}\Vert\,|\, \mf{u}, \mf{v}\in S\right\}.\]  
\end{definition}
\begin{example}[label=230723_8]{}
Consider the rectangle $R=[a_1, b_1]\times \cdots\times [a_n,b_n]$. If $\mf{u}$ and $\mf{v}$ are two points in $R$, for each $1\leq i\leq n$,
$u_i, v_i\in [a_i, b_i]$. Thus,
\[|u_i-v_i|\leq b_i-a_i.\]
It follows that
\[\Vert\mf{u}-\mf{v}\Vert \leq \sqrt{(b_1-a_1)^2+\cdots+(b_n-a_n)^2}.\]
If $\mf{u}_0=\mf{a}=(a_1, \ldots, a_n)$ and $\mf{v}_0=\mf{b}=(b_1, \ldots, b_n)$, then $\mf{u}_0$ and $\mf{v}_0$ are in $R$, and 
\[\Vert\mf{u}_0-\mf{v}_0\Vert =\sqrt{(b_1-a_1)^2+\cdots+(b_n-a_n)^2}.\] This shows that the diameter of $R$ is
\[\text{diam}\,R=\Vert\mf{b}-\mf{a}\Vert=\sqrt{(b_1-a_1)^2+\cdots+(b_n-a_n)^2}.\]
\end{example}

\begin{figure}[ht]
\centering
\includegraphics[scale=0.2]{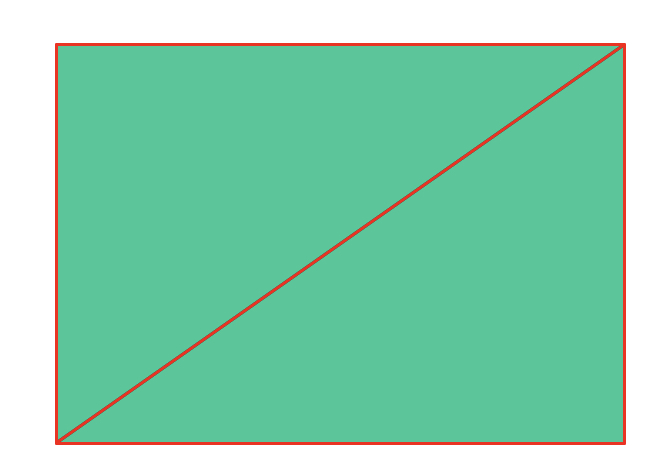}

\caption{The diameter of a rectangle.}\label{figure33}
\end{figure}

Intuitively, the diameter of the open rectangle $U=(a_1, b_1)\times\cdots\times (a_n, b_n)$ is also equal to \[d=\sqrt{(b_1-a_1)^2+\cdots+(b_n-a_n)^2}.\]
However, the points $\mf{a}=(a_1, \ldots, a_n)$ and $\mf{b}=(b_1, \ldots, b_n)$ are not in $U$. There does not exist two points in $U$ whose distance is  $d$, but there are sequences of points $\{\mf{u}_k\}$ and $\{\mf{v}_k\}$ such that their distances $\{\Vert\mf{u}_k-\mf{v}_k\Vert\}$ approach $d$ as $k\to\infty$. We will formulate this as a more general theorem.

\begin{theorem}[label=230723_4]{}
Let $S$ be a subset of $\mb{R}^n$. If $S$ is bounded, then its closure $\overline{S}$ is also bounded. Moreover, $\text{diam}\,\overline{S}=\text{diam}\,S$.
\end{theorem}
\begin{myproof}{Proof}

If $\mf{u}$ and $\mf{v}$ are two points in $\overline{S}$, there exist sequences $\{\mf{u}_k\}$ and $\{\mf{v}_k\}$ in $S$ that converge respectively to $\mf{u}$ and $\mf{v}$.
Then
\begin{equation}\label{231012_1}d(\mf{u},\mf{v})=\lim_{k\to\infty}d(\mf{u}_k,\mf{v}_k).\end{equation}
For each $k\in\mb{Z}^+$, since $\mf{u}_k$ and $\mf{v}_k$ are in $S$, 
\[d(\mf{u}_k,\mf{v}_k)\leq\, \text{diam}\,S.\]
Eq. \eqref{231012_1} implies that 
\[d(\mf{u},\mf{v})\leq \, \text{diam}\,S.\]Since this is true for any $\mf{u}$ and $\mf{v}$ in $\overline{S}$,
$\overline{S}$ is bounded and 
\[\text{diam}\,\overline{S}\;\leq\;\text{diam}\,S.\]
Since $S\subset \overline{S}$, we also have
  $\text{diam}\,S\leq\text{diam}\,\overline{S}$. We conclude  that $\text{diam}\,\overline{S} =\text{diam}\,S$.
\end{myproof}

The following example justifies that the diameter of a ball of radius $r$ is indeed $2r$.
\begin{example}{}
Find the diameter of the open ball $B(\mf{x}_0, r)$ in $\mb{R}^n$.
\end{example}
\begin{solution}{Solution}
By Theorem  \ref{230723_4}, the diameter of the open ball $B(\mf{x}_0, r)$ is the same as the diameter of its closure, the closed ball $CB(\mf{x}_0, r)$. Given $\mf{u}$ and $\mf{v}$ in $CB(\mf{x}_0, r)$, $\Vert\mf{u}-\mf{x}_0\Vert\leq r$ and $\Vert\mf{v}-\mf{x}_0\Vert\leq r$. Therefore,
\[\Vert\mf{u}-\mf{v}\Vert\leq \Vert\mf{u}-\mf{x}_0\Vert+\Vert\mf{v}-\mf{x}_0\Vert\leq 2r.\]
This shows that $\text{diam}\,CB(\mf{x}_0, r)\leq 2r$. The points $\mf{u}_0=\mf{x}_0+r\mf{e}_1$ and    $\mf{v}_0=\mf{x}_0-r\mf{e}_1$ are in the closed ball $CB(\mf{x}_0, r)$. Since
\[\Vert\mf{u}_0-\mf{v}_0\Vert=\Vert 2r\mf{e}_1\Vert=2r,\] $\text{diam}\,CB(\mf{x}_0, r)\geq 2r$. 
Therefore, the diameter of the closed ball $CB(\mf{x}_0, r)$ is exactly $2r$. By Theorem  \ref{230723_4}, the diameter of the open ball $B(\mf{x}_0, r)$ is also $2r$.
\end{solution}

\begin{figure}[ht]
\centering
\includegraphics[scale=0.2]{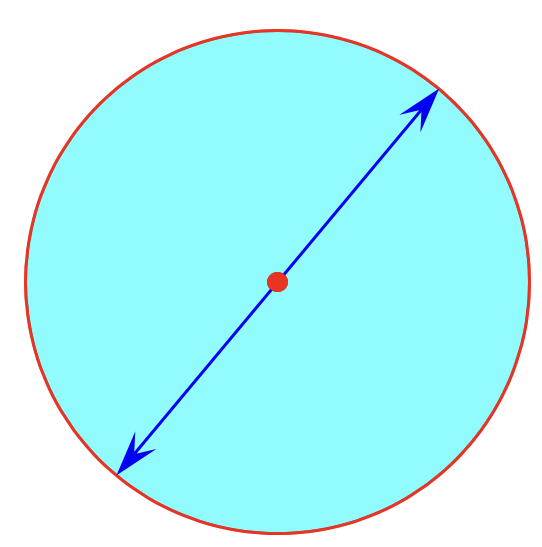}

\caption{The diameter of a ball.}\label{figure34}
\end{figure}

In volume I, we have shown that a bounded sequence in $\mb{R}$ has a convergent subsequence. This is achieved by using the monotone convergence theorem, which says that a bounded monotone sequence in $\mb{R}$ is convergent. 
For points in $\mb{R}^n$ with $n\geq 2$, we cannot apply monotone convergence theorem, as we cannot define a simple order on the points in $\mb{R}^n$ when $n\geq 2$. Nevertheless, we can use the result of $n=1$ and the componentwise convergence theorem to show that a bounded sequence in $\mb{R}^n$ has a convergent subsequence.

\begin{theorem}[label=230722_16]{}
Let $\{\mf{u}_k\}$ be a sequence in $\mb{R}^n$. If $\{\mf{u}_k\}$ is bounded, then there is a subsequence that is convergent.
\end{theorem}
\begin{myproof}{Sketch of Proof}
The $n=1$ case is already established in volume I. Here we prove the $n=2$ case. The $n\geq 3$ case can be proved by induction using the same reasoning.

For $k\in\mb{Z}^+$, let $\mf{u}_k=(x_k,y_k)$. Since
\[|x_k|\leq \Vert\mf{u}_k\Vert\quad\text{and}\quad  |y_k|\leq \Vert\mf{u}_k\Vert,\]the sequences $\{x_k\}$ and $\{y_k\}$ are bounded sequences. Thus, there is a subsequence $\{x_{k_j}\}_{j=1}^{\infty}$ of $\{x_k\}_{k=1}^{\infty}$ that converges to a point $x_0$ in $\mb{R}$. Consider the subsequence $\{y_{k_j}\}_{j=1}^{\infty}$ of the sequence $\{y_k\}_{k=1}^{\infty}$. It is also bounded. Hence, there is a subsequence $\{y_{k_{j_l}}\}_{l=1}^{\infty}$ that converges to a point $y_0$ in $\mb{R}$. Notice that the subsequence $\{x_{k_{j_l}}\}_{l=1}^{\infty}$ of $\{x_k\}_{k=1}^{\infty}$ is also a subsequence of $\{x_{k_j}\}_{j=1}^{\infty}$. Hence, it also converges to $x_0$. By componentwise convergence theorem,  $\{\mf{u}_{k_{j_l}}\}_{l=1}^{\infty}$ is a subsequence of $\{\mf{u}_k\}_{k=1}^{\infty}$ that converges to $(x_0,y_0)$. This proves the theorem when $n=2$.

\end{myproof}

Now we study the concept of sequential compactness.
It is the same as the $n=1$ case.

\begin{definition}{Sequentially Compact}
Let $S$ be a subset of $\mb{R}^n$. We say that $S$ is sequentially compact
provided that every sequence in $S$ has a subsequence that converges to a
point in $S$.
\end{definition}

In volume I, we proved the Bolzano-Weierstrass theorem, which says that a subset of $\mb{R}$ is sequentially compact if and only if it is closed and bounded. In fact, the same is true for the $n\geq 2$ case. Let us first look at some examples.
\begin{example}[label=230722_14]{}
Show that the set $A=\left\{(x,y)\,|\,x^2+y^2<1\right\}$ is not sequentially compact.
\end{example}
\begin{solution}{Solution}
For $ k\in\mb{Z}^+$, let
\[\mf{u}_k=\left(\frac{k}{k+1}, 0\right).\]Then $\{\mf{u}_k\}$ is a sequence in $A$ that converges to the point $\mf{u}_0=(1,0)$ that is not in $A$. Thus, every subsequence of $\{\mf{u}_k\}$ converges to the point $\mf{u}_0$, which is not in $A$. This means the sequence $\{\mf{u}_k\}$ in $A$ does not have a subsequence that converges to a point in $A$. Hence, $A$ is not sequentially compact.

\end{solution}
Note that the set $A$ in Example  \ref{230722_14} is not closed.

 \begin{example}[label=230722_15]{}
Show that the set $C=\left\{(x,y)\,|\,1\leq x\leq 3, y\geq 0\right\}$ is not sequentially compact.
\end{example}
\begin{solution}{Solution}
For $ k\in\mb{Z}^+$, let
$\mf{u}_k=\left(2, k\right)$. Then $\{\mf{u}_k\}$ is a sequence in $C$. If  $\{\mf{u}_{k_j}\}_{j=1}^{\infty}$ is  a subsequence of $\{\mf{u}_k\}$, then $k_1,k_2, k_3, \ldots$ is a strictly increasing sequence of positive integers. Therefore $k_j\geq j$ for all $j\in\mb{Z}^+$. It follows that
\[\Vert \mathbf{u}_{k_j}\Vert=\Vert (2, k_j) \Vert\geq k_j\geq j\hspace{1cm}\text{for all}\;j\in\mb{Z}^+.\]
Hence, the subsequence $\{\mf{u}_{k_j}\}$ is not bounded. Therefore, it is not convergent. This means that the sequence $\{\mf{u}_k\}$ in $C$ does not have a convergent subsequence. Therefore, $C$ is not sequentially compact.

\end{solution}
Note that the set $C$ in Example  \ref{230722_15} is not bounded.

Now we prove the main theorem.
\begin{theorem}{Bolzano-Weierstrass Theorem}
Let $S$ be a subset of $\mb{R}^n$. The following are equivalent.
\begin{enumerate}[(a)]
\item $S$ is closed and bounded.
\item $S$ is sequentially compact.

\end{enumerate}\end{theorem}

\begin{myproof}{Proof}
First assume that $S$ is closed and bounded. Let $\{\mf{x}_k\}$ be a sequence in $S$. Then $\{\mf{x}_k\}$ is also bounded. By Theorem \ref{230722_16}, there is subsequence $\{\mf{x}_{k_j}\}$ that converges to some $\mf{x}_0$. Since $S$ is closed, we must have $\mf{x}_0$ is in $S$. This proves that every sequence in $S$ has a subsequence that converges to a point in $S$. Hence, $S$ is sequentially compact. This completes the proof of (a) implies (b).

To prove that (b) implies (a), it suffices to show that if $S$ is not closed or $S$ is not bounded,  then $S$ is not sequentially compact. 

If $S$ is not closed, there is  a sequence $\{\mf{x}_k\}$ in $S$ that converges to a point $\mf{x}_0$, but $\mf{x}_0$ is not in $S$. Then every subsequence of   $\{\mf{x}_k\}$ converges to the point $\mf{x}_0$, which is not in $S$. Thus, $\{\mf{x}_k\}$ is a sequence in $S$  that does not have any subsequence that converges to a point in $S$. This shows that $S$ is not sequentially compact.

If $S$ is not bounded, for each positive integer $k$, there  is a point $\mf{x}_k$ in $S$ such that $\Vert\mf{x}_k\Vert\geq k$. If $\{\mf{x}_{k_j}\}_{j=1}^{\infty}$ is  a subsequence of $\{\mf{x}_k\}$, then $k_1,k_2, k_3, \ldots$ is a strictly increasing sequence of positive integers. Therefore $k_j\geq j$ for all $j\in\mb{Z}^+$. It follows that
$\Vert \mathbf{x}_{k_j}\Vert \geq k_j\geq j$ for al $j\in\mb{Z}^+$.
Hence, the subsequence $\{\mf{x}_{k_j}\}$ is not bounded. Therefore, it is not convergent. This means that the sequence $\{\mf{x}_k\}$ in $S$ does not have a convergent subsequence. Therefore, $S$ is not sequentially compact.
\end{myproof}

\begin{corollary}{}
A closed rectangle $\di R=[a_1, b_1]\times\cdots\times [a_n, b_n]$ in $\mb{R}^n$   is sequentially compact.
\end{corollary}
\begin{myproof}{Proof}
We have shown in Chapter \ref{chapter1} that $R$ is closed. Example \ref{230723_8} shows that $R$ is bounded. Thus, $R$ is sequentially compact.
\end{myproof}

An interesting consequence of Theorem \ref{230723_4} is the following.
\begin{corollary}{}
If $S$ be a bounded subset of $\mb{R}^n$, then its closure $\overline{S}$ is sequentially compact.
\end{corollary}
\begin{example}[label=230723_1]{}
Determine whether the following subsets of $\mb{R}^3$ is sequentially compact.
\begin{enumerate}[(a)]
\item $A=\left\{(x,y,z)\,|\, xyz=1\right\}$.
\item $B=\left\{(x,y,z)\,|\, x^2+4y^2+9z^2\leq 36\right\}$.
\item $C=\left\{(x,y,z)\,|\, 1\leq x\leq 2, 1\leq y\leq 3, 0<xyz\leq 4\right\}$.

\end{enumerate}
\end{example}
\begin{solution}{Solution}
\begin{enumerate}[(a)]
\item
For any $k\in\mb{Z}^+$, let
\[\mf{u}_k=\left(k, \frac{1}{k}, 1\right).\]
Then $\{\mf{u}_k\}$ is a sequence in $A$, and $\Vert\mf{u}_k\Vert\geq k$. Therefore, $A$ is not bounded. Hence, $A$ is not sequentially compact.
\item For any $\mf{u}=(x,y,z)\in B$,
\[\Vert\mf{u}\Vert^2=x^2+y^2+z^2\leq x^2+4y^2+9z^2\leq 36.\]
Hence, $B$ is bounded. The function $f:\mb{R}^3\to\mb{R}$, $f(x,y,z)=x^2+4y^2+9z^2$ is a polynomial. Hence, it is continuous. Since the set $I=(-\infty, 36]$ is closed in $\mb{R}$, and $B=f^{-1}(I)$, $B$ is closed in $\mb{R}^3$. Since $B$ is closed and bounded, it is sequentially compact.
\end{enumerate}\bs
\begin{enumerate}[(a)]
 \item[(c)] For any $k\in\mb{Z}^+$, let
\[\mf{u}_k=\left(1,1,\frac{1}{k}\right).\]
Then $\{\mf{u}_k\}$ is a sequence of points in $C$ that converges to the point $\mf{u}_0=(1,1,0)$, which is not in $C$. Thus, $C$ is not closed, and so $C$ is not sequentially compact.
\end{enumerate}
\end{solution}

The following theorem asserts that continuous functions preserve sequential compctness.
\begin{theorem}[label=230722_17]{}
Let $\mk{D}$ be a sequentially compact subset of $\mb{R}^n$. If the function $\mf{F}:\mk{D}\to\mb{R}^m$ is continuous, then $\mf{F}(\mk{D})$ is a sequentially compact subset of $\mb{R}^m$.
\end{theorem}
The proof of this theorem is identical to the $n=1$ case.
\begin{myproof}{Proof}
Let $\{\mf{y}_k\}$ be a sequence in $\mf{F}(\mk{D})$. For each $k\in\mb{Z}^+$, there exists $\mf{x}_k\in\mk{D}$ such that $\mf{F}(\mf{x}_k)=\mf{y}_k$. Since $\mk{D}$ is sequentially compact, there is a subsequence $\{\mf{x}_{k_j}\}$ of $\{\mf{x}_k\}$ that converges to a point $\mf{x}_0$ in $\mk{D}$. Since $\mf{F}$ is continuous, the sequence $\{\mf{F}(\mf{x}_{k_j})\}$ converges to $\mf{F}(\mf{x}_0)$. Note that $\mf{F}(\mf{x}_0)$ is in $\mf{F}(\mk{D})$. In other words,  $\{\mf{y}_{k_j}\}$ is a subsequence  of the sequence $\{\mf{y}_k\}$ that converges to $\mf{F}(\mf{x}_0)$ in $\mf{F}(\mk{D})$. This shows that every sequence in $\mf{F}(\mk{D})$ has a subsequence that converges to a point in $\mf{F}(\mk{D})$. Thus, $\mf{F}(\mk{D})$ is a sequentially compact subset of $\mb{R}^m$.
\end{myproof}

We are going to  discuss important consequences of   Theorem \ref{230722_17} in the coming section. For the rest of this section, we  introduce the concept of compactness, which plays a central role in modern analysis.
We start with the definition of an open covering.
\begin{definition}{Open Covering}
Let $S$ be a subset of $\mb{R}^n$, and let $\mathscr{A}=\left\{U_{\alpha}\,|\,\alpha\in J\right\}$ be a collection of open sets in $\mb{R}^n$ indexed by the set $J$. We say that $\mathscr{A}$ is an open covering of $S$ provided that
\[S\subset\bigcup_{\alpha\in J}U_{\alpha}.\]
\end{definition}

\begin{example}[label=230722_18]{}
For each $k\in\mb{Z}^+$, let $U_k=(1/k, 1)$. Then $U_k$ is an open set in $\mb{R}$ and 
\[\bigcup_{k=1}^{\infty}U_k=(0,1).\]
Hence, $\mathscr{A}=\left\{U_k\,|\,k\in\mb{Z}^+\right\}$ is an open covering of the set $S=(0,1)$.
\end{example}

\begin{remark}{}
If  $\mathscr{A}=\left\{U_{\alpha}\,|\,\alpha\in J\right\}$ is an open covering of $S$ and $S'$ is a subset of $S$, then   $\mathscr{A}=\left\{U_{\alpha}\,|\,\alpha\in J\right\}$ is also an open covering of $S'$.\end{remark}

\begin{example}[label=230722_19]{}
For each $k\in\mb{Z}^+$, let $U_k=B(\mf{0}, k)$ be the ball in $\mb{R}^n$ centered at the origin and having radius $k$. Then
\[\bigcup_{k=1}^{\infty}U_k=\mb{R}^n.\]
Thus,  $\mathscr{A}=\left\{U_k\,|\,k\in\mb{Z}^+\right\}$ is an open covering of any subset $S$ of $\mb{R}^n$.
\end{example}

\begin{definition}{Subcover}
Let $S$ be a subset of $\mb{R}^n$, and let $\mathscr{A}=\left\{U_{\alpha}\,|\,\alpha\in J\right\}$ be  an open covering of $S$. A subcover is a subcollection of $\mathscr{A}$ which is also a covering of  $S$. A {\it finite subcover} is a subcover that contains only finitely many elements.
\end{definition}

\begin{example}[label=230723_29]{}
For each $k\in\mb{Z}$, let $U_k=(k, k+2)$. Then
$\di \bigcup_{k=-\infty}^{\infty}U_k=\mb{R}$. Thus,  $\mathscr{A}=\left\{U_k\,|\,k\in\mb{Z}\right\}$ is an open covering of the set $S=[-3, 4)$. There is a finite subcover of $S$ given by \[\mathscr{B}=\{U_{-4}, U_{-3}, U_{-2}, U_{-1}, U_0, U_1, U_2\}.\]
\end{example}
\begin{definition}{Compact Sets}
Let $S$ be a subset of $\mb{R}^n$. We say that $S$ is compact provided that every open covering of $S$ has a finite subcover. Namely, if $\mathscr{A}=\left\{U_{\alpha}\,|\,\alpha\in J\right\}$ is an open covering of $S$, then there exist $\alpha_1, \ldots, \alpha_k\in J$ such that
\[S\subset \bigcup_{j=1}^kU_{\alpha_j}.\]
\end{definition}

\begin{example}{}
The subset $S=(0,1)$ of $\mb{R}$ is not compact. For 
 $k\in\mb{Z}^+$, let $U_k=(1/k, 1)$.  Example \ref{230722_18} says that $\mathscr{A}=\left\{U_k\,|\,k\in\mb{Z}^+\right\}$ is an open covering of the set $S$.  We claim that there is no finite subcollection of $\mathscr{A}$ that covers $S$.
 
 Assume to the contrary that there exists a finite subcollection of $\mathscr{A}$ that covers $S$. Then there are positive integers $k_1, \ldots, k_m$ such that
 \[(0,1)\subset\bigcup_{j=1}^m U_{k_j}=\bigcup_{j=1}^m\left(\frac{1}{k_j}, 1\right).\]
 Notice that if $k_i\leq k_j$, then $U_{k_i}\subset U_{k_j}$. Thus, if $K=\max\{k_1, \ldots, k_m\}$, then
 \[\bigcup_{j=1}^m U_{k_j}=U_K=\left(\frac{1}{K}, 1\right),\]
 and so $S=(0,1)$ is not contained in $U_K$. This  gives a contradiction. Hence, $S$ is not compact.
\end{example}

\begin{example}{}
As a subset of itself, $\mathbb{R}^n$ is not compact. For  $k\in\mb{Z}^+$, let $U_k=B(\mf{0}, k)$ be the ball in $\mb{R}^n$ centered at the origin and having radius $k$. Example \ref{230722_19} says that $\mathscr{A}=\left\{U_k\,|\,k\in\mb{Z}^+\right\}$ is an open covering of  $\mb{R}^n$. We claim that there is no finite subcover. 

Assume to the contrary that there is a finite subcover. Then there exist  positive integers $k_1, \ldots, k_m$ such that
\[\mb{R}^n=\bigcup_{j=1}^m U_k.\]
Notice that if $k_i\leq k_j$, then $U_{k_i}\subset U_{k_j}$. Thus, if $K=\max\{k_1, \ldots, k_m\}$, then
 \[\bigcup_{j=1}^m U_{k_j}=U_K=B(\mf{0}, K).\]Obviously, $B(\mf{0}, K)$ is not equal to $\mb{R}^n$. This gives a contradiction. Hence, $\mb{R}^n$ is not compact.
\end{example}

Our goal is to prove the Heine-Borel theorem, which says that a subset of $\mb{R}^n$ is compact if and only if it is closed and bounded. We first prove the easier direction.
\begin{theorem}[label=230723_5]{}
Let $S$ be a subset of $\mb{R}^n$. If $S$ is compact, then it is closed and bounded.
\end{theorem}
\begin{myproof}{Proof}
We show that if $S$ is compact, then it is bounded; and if $S$ is compact, then it is closed. 

First we prove that if $S$ is compact, then it is bounded.  For  $k\in\mb{Z}^+$, let $U_k=B(\mf{0}, k)$ be the ball in $\mb{R}^n$ centered at the origin and having radius $k$. Example \ref{230722_19} says that $\mathscr{A}=\left\{U_k\,|\,k\in\mb{Z}^+\right\}$ is an open covering of  $S$. Since $S$ is compact, there exist   positive integers $k_1, \ldots, k_m$ such that
\[S\subset\bigcup_{j=1}^m U_{k_j}=U_K=B(\mf{0}, K),\]
\bp
where $K=\max\{k_1, \ldots, k_m\}$. This shows that 
\[\Vert\mf{x}\Vert\leq K\hspace{1cm}\text{for all}\;\mf{x}\in S.\]
Hence, $S$ is bounded.

Now we prove that if $S$ is compact, then it is closed. For this, it suffices to show that $\overline{S}\subset S$, or equivalently, any point that is not in $S$ is not in $\overline{S}$.  Assume that    $\mf{x}_0$  is not in $S$. For each $k\in\mb{Z}^+$, let \[V_k=\text{ext}\, B(\mf{x}_0, 1/k)=\left\{\mathbf{x}\in\mb{R}^n\,\left|\,\Vert\mf{x}-\mf{x}_0\Vert>\frac{1}{k}\right.\right\}.\] Then $V_k$ is open in $\mb{R}^n$. If $\mf{x}$ is a point in $\mb{R}^n$ and $\mf{x}\neq\mf{x}_0$, then $r=\Vert\mf{x}-\mf{x}_0\Vert>0$. There is a $k\in\mb{Z}^+$ such that $1/k<r$. Then $\mf{x}$ is in $ V_k$.  This shows that
\[\bigcup_{k=1}^{\infty}V_k=\mb{R}^n\setminus\{\mf{x}_0\}.\] 
Therefore, $\mathscr{A}=\left\{V_k\,|\,k\in\mb{Z}^+\right\}$ is an open covering of $S$. Since $S$ is compact, there is a finite subcover. Namely, there exist positive integers $k_1, \ldots, k_m$ such that
\[S\subset\bigcup_{j=1}^m V_{k_j}=V_K,\]
where $K=\max\{k_1, \ldots, k_m\}$.  Since $B(\mf{x}_0,1/K)$ is disjoint from $V_K$, it does not contain any point of $S$. This shows that $\mathbf{x}_0$ is not in $\overline{S}$, and thus   the proof is completed.

\end{myproof}

\begin{example}[label=230723_2]{}
 The set
  \[A=\left\{(x,y,z)\,|\, xyz=1\right\}\] in Example \ref{230723_1} is not compact because it is not bounded. The set \[C=\left\{(x,y,z)\,|\, 1\leq x\leq 2, 1\leq y\leq 3, 0<xyz\leq 4\right\}\] is not compact because it is not closed.
 
\end{example}

We are now left to show that a closed and bounded subset of $\mb{R}^n$ is compact.  We start by proving a special case.
\begin{theorem}[label=230723_6]{}
A closed rectangle $R=[a_1, b_1]\times \cdots\times [a_n, b_n]$ in $\mb{R}^n$ is compact.
\end{theorem}
\begin{myproof}{Proof}
We will prove by contradiction. Assume that $R$ is not compact, and we show that this  will lead to a contradiction.
The idea   is to use the bisection method.

If $R$ is not compact, there is an open covering
  $\mathscr{A}=\left\{U_{\alpha}\,|\,\alpha\in J\right\}$ of $R$ which does not have a finite subcover. 
  
  Let $R_1=R$, and let $d_1=\text{diam}\,R_1$.
 For $1\leq i\leq n$, let $a_{i,1}=a_i$ and $b_{i,1}=b_i$, and let
$m_{i,1}$ to be the midpoint of the interval $[a_{i,1}, b_{i,1}]$. The hyperplanes $x_i=m_{i,1}$, $1\leq i\leq n$, divides the rectangle $R_{1}$ into $2^n$ subrectangles. Notice that $\mathscr{A}$ is also an open covering of each of these subrectangles. If each of these subrectangles can be covered by a finite subcollection of open sets in $\mathscr{A}$, then $R$ also can be covered by a finite subcollection of open sets in $\mathscr{A}$. Since we assume $R$ cannot be covered by any finite subcollection of open sets in $\mathscr{A}$, there is at least one of the   $2^n$ subrectangles which cannot be covered by any finite subcollection of open sets in $\mathscr{A}$. Choose one of these, and denote it by $R_2$. 
 
 Define $a_{i,2}, b_{i,2}$ for $1\leq i\leq n$ so that
 \[R_2=[a_{1,2}, b_{1,2}]\times \cdots\times [a_{n,2}, b_{n,2}].\] Note that
 \[b_{i,2}-a_{i,2}=\frac{b_{i,1}-a_{i,1}}{2}\hspace{1cm}\text{for}\;1\leq i\leq n.\] Therefore,  $d_2=\text{diam}\,R_2=d_1/2$. 
 
 We continue this bisection process to obtain the rectangles $R_1, R_2, \cdots$, so that $R_{k+1}\subset R_k$ for all $k\in \mb{Z}^+$,  and $R_k$ cannot be covered by any finite subcollections of $\mathscr{A}$.
  \end{myproof}
 \begin{figure}[ht]
\centering
\includegraphics[scale=0.2]{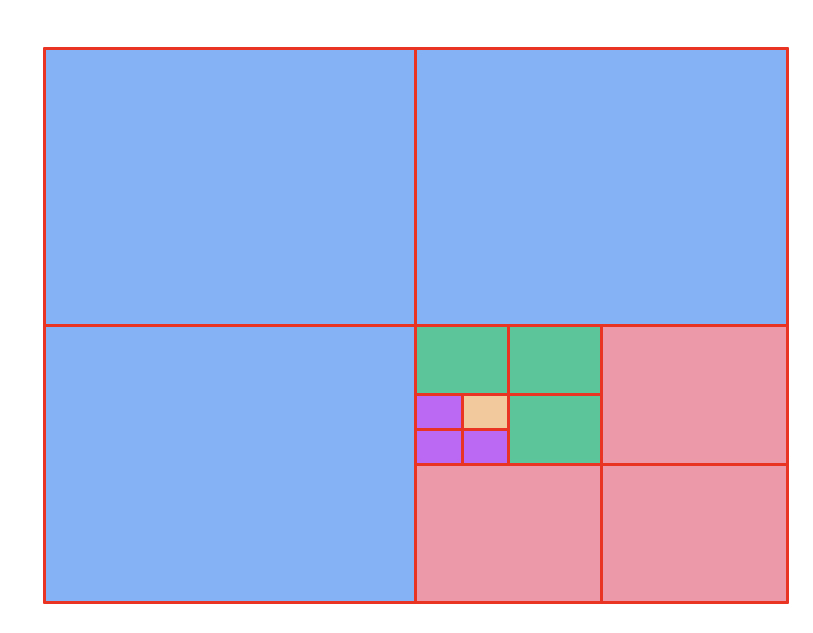}

\caption{Bisection method.}\label{figure32}
\end{figure}

 \begin{myproof}{}
 
Define $a_{i,k}, b_{i,k}$ for $1\leq i\leq n$ so that
 \[R_k=[a_{1,k}, b_{1,k}]\times \cdots\times [a_{n,k}, b_{n,k}].\]Then for all $k\in\mb{Z}^+$, 
 \[b_{i,k+1}-a_{i,k+1}=\frac{b_{i,k}-a_{i,k}}{2}\hspace{1cm}\text{for}\;1\leq i\leq n.\] It follows that $d_{k+1}=\text{diam}\,R_{k+1}=d_k/2$. 
 
 For any $1\leq i\leq n$, $\{a_{i,k}\}_{k=1}^{\infty}$ is an increasing sequence that is bounded above by $b_i$, and  $\{b_{i,k}\}_{k=1}^{\infty}$ is a decreasing sequence that is bounded below by $a_i$. By monotone convergence theorem, the sequence  $\{a_{i,k}\}_{k=1}^{\infty}$ converges to $a_{i,0}=\di\sup_{k\in\mb{Z}^+}a_{i,k}$; while the sequence   $\{b_{i,k}\}_{k=1}^{\infty}$ converges to $b_{i,0}=\di\inf_{k\in\mb{Z}^+}b_{i,k}$.
 Since
 \[ b_{i,k}-a_{i,k}=\frac{b_i-a_i}{2^{k-1}}\hspace{1cm}\text{for all}\;k\in\mb{Z}^+,\]we find that $a_{i,0}=b_{i,0}$. Let $c_i=a_{i,0}=b_{i,0}$. Then $a_{i,k}\leq c_i\leq b_{i,k}$ for all $1\leq i\leq n$ and all $k\in\mathbb{Z}^+$. Thus, $\mf{c}=(c_1, \ldots, c_n)$ is a point in $R_k$ for all $k\in\mb{Z}^+$. By assumption that $\mathscr{A}$ is an open covering of $R=R_1$, there exists $\beta\in J$ such that $\mf{c}\in U_{\beta}$. Since $U_{\beta}$ is an open set, there is an $r>0$ such that $B(\mf{c}, r)\subset U_{\beta}$. Since 
 \[d_k=\text{diam}\,R_k=\frac{d_1}{2^{k-1}}\hspace{1cm}\text{for all}\;k\in\mb{Z}^+,\]
 \bp
 we find that $\di\lim_{k\to\infty}d_k=0$. Hence, there is a positive integer $K$ such that $d_K<r$. If $\mf{x}\in R_K$, then 
 \[\Vert\mf{x}-\mf{c}\Vert\leq \text{diam}\,R_K=d_K<r.\]
 This implies that $\mf{x}$ is in   $B(\mf{c}, r)$.  Thus, we have shown that $R_K\subset B(\mf{c}, r)$. Therefore, $R_K$ is contained in the single element $ U_{\beta}$ of $\mathscr{A}$, which contradicts to $R_K$ cannot be covered by any finite subcollection of open sets in $\mathscr{A}$.
 
 We conclude that $R$ must be compact.
\end{myproof}

Now we can prove the Heine-Borel theorem.
\begin{theorem}{Heine-Borel Theorem}
Let $S$ be a subset of $\mb{R}^n$. Then $S$ is compact if and only if it is closed and bounded.
\end{theorem}
\begin{myproof}{Proof}
We have shown in Theorem \ref{230723_5} that if $S$ is compact, then it must be closed and bounded.

Now assume that $S$ is closed and bounded, and let $\mathscr{A}=\left\{U_{\alpha}\,|\,\alpha\in J\right\}$ be an open covering of $S$. Since $S$ is bounded, there exists a positive number $M$ such that
\[\Vert\mf{x}\Vert\leq M\hspace{1cm}\text{for all}\;\mathbf{x}\in S.\]
Thus, if $\mf{x}=(x_1, \ldots, x_n)$ is in $S$, then for all $1\leq i\leq n$,
$|x_i|\leq\Vert\mf{x}\Vert \leq M$. This implies that $S$ is contained in the closed rectangle 
\[R=[-M,M]\times\cdots\times [-M, M].\] Let $V=\mathbb{R}^n\setminus S$. Since $S$ is closed, $V$ is an open set. Then $\widetilde{\mathscr{A}}=\mathscr{A}\cup\{V\}$ is an open covering of $\mathbb{R}^n$, and hence it is an open covering of $R$. By Theorem \ref{230723_6}, $R$ is compact.  Thus, there exists  $\widetilde{\mathscr{B}}\subset   \widetilde{\mathscr{A}}$ which is a finite subcover of $R$. Then $\mathscr{B}=\widetilde{\mathscr{B}}\setminus\{V\}$ is a finite subcollection of $\mathscr{A}$ that covers $S$. This proves that $S$ is compact.
\end{myproof}

\begin{example}{}
We have shown in Example \ref{230723_1} that the set
  \[B=\left\{(x,y,z)\,|\, x^2+4y^2+9z^2\leq 36\right\}\] is closed and bounded. Hence, it is compact.
\end{example}
 
 Now we can conclude our main theorem from the Bolzano-Weierstrass theorem and the Heine-Borel theorem.
 \begin{theorem}{}
 Let $S$ be a subset of $\mb{R}^n$. Then the following are equivalent.
 \begin{enumerate}[(a)]

 \item $S$ is sequentially compact.
 \item $S$ is closed and bounded.
  \item
 $S$ is compact.
 \end{enumerate}
 \end{theorem}
 
\begin{remark}{}
Henceforth, when we say a subset $S$  of $\mb{R}^n$ is compact, we mean it is a closed and bounded set, and it is sequentially compact. By Theorem \ref{230723_4}, a subset $S$ of $\mb{R}^n$ has compact closure if and only if it is a bounded set.
\end{remark}

 Finally, we can conclude the following, which says that continuous functions preserve compactness.
\begin{theorem}[label=230723_28]{}
Let $\mk{D}$ be a  compact subset of $\mb{R}^n$. If the function $\mf{F}:\mk{D}\to\mb{R}^m$ is continuous, then $\mf{F}(\mk{D})$ is a   compact subset of $\mb{R}^m$.
\end{theorem}
\begin{myproof}{Proof}
Since $\mk{D}$ is compact, it is sequentially compact. By Theorem \ref{230722_17}, $\mf{F}(\mk{D})$ is a sequentially compact subset of $\mb{R}^m$. Hence, $\mf{F}(\mk{D})$ is a   compact subset of $\mb{R}^m$.
\end{myproof}

\vp
\noindent
{\bf \large Exercises  \thesection}
\setcounter{myquestion}{1}

 \begin{question}{\themyquestion}
Determine whether the following subsets of $\mb{R}^2$ is sequentially compact.
\begin{enumerate}[(a)]
\item $A=\left\{(x,y)\,|\, x^2+y^2=9\right\}$.
\item $B=\left\{(x,y)\,|\, 0<x^2+4y^2 \leq 36\right\}$.
\item $C=\left\{(x,y)\,|\, x\geq 0,  0\leq y\leq x^2\right\}$.

\end{enumerate}
\end{question}

\atc
 \begin{question}{\themyquestion}
Determine whether the following subsets of $\mb{R}^3$ is   compact.
\begin{enumerate}[(a)]
\item $A=\left\{(x,y, z)\,|\, 1\leq x\leq 2\right\}$.
\item $B=\left\{(x,y, z)\,|\,|x|+|y|+|z|\leq 10\right\}$.
\item $C=\left\{(x,y, z)\,|\, 4\leq x^2+y^2+z^2\leq 9\right\}$.

\end{enumerate}
\end{question}

\atc
\begin{question}{\themyquestion}
Given that $A$ is a compact subset of $\mb{R}^n$ and $B$ is a subset of $A$, show that $B$ is compact if and only if it is closed. 
\end{question}

\atc
\begin{question}{\themyquestion}
 If $S_1, \ldots, S_k$ are compact subsets of $\mb{R}^n$, show that $S=S_1\cup \cdots\cup S_n$ is also compact.
\end{question}

\atc
\begin{question}{\themyquestion}
If $A$ is a compact subset of $\mb{R}^m$, $B$ is a compact subset of $\mb{R}^n$, show that $A\times B$ is a compact subset of $\mb{R}^{m+n}$.
\end{question}

\section{Applications of Compactness}
In this section, we consider the applications of compactness. We are going to use repeatedly the fact that a subset $S$ of $\mb{R}^n$ is compact if and only if it is closed and bounded, if and only if it is sequentially compact.

\subsection{The Extreme Value Theorem}
 First we define bounded functions.
 
 \begin{definition}{Bounded Functions}
 Let $\mk{D}$ be a subset of $\mb{R}^n$, and let $\mf{F}:\mk{D}\to\mb{R}^m$ be a function defined on $\mk{D}$. We say that the function $\mf{F}$ is bounded if the set $\mf{F}(\mk{D})$ is a bounded subset of $\mb{R}^m$. In other words, the function $\mf{F}:\mk{D}\to\mb{R}^m$ is bounded if there is positive number $M$ such that
 \[\Vert\mf{F}(\mf{x})\Vert\leq M\hspace{1cm}\text{for all}\;\mf{x}\in \mk{D}.\]
 \end{definition}
 
 \begin{example}{}
 Let $\mk{D}=\left\{(x,y, z)\,|\, 0<x^2+y^2+z^2<4\right\}$, and let $\mf{F}:\mk{D}\to \mb{R}^2$ be the function defined as
 \[\mf{F}(x,y,z)=\left(\frac{1}{x^2+y^2+z^2}, x+y+z\right).\]For $k\in \mb{Z}^+$, the point $\mf{u}_k=\di (1/k, 0, 0)$ is in $\mk{D}$ and 
 \[\mf{F}(\mf{u}_k)=\left(k^2, \frac{1}{k}\right).\]
 Thus, $\Vert \mf{F}(\mf{u}_k)\Vert\geq k^2$. This shows that $\mf{F}$ is not bounded, even though $\mk{D}$ is a bounded set.
 \end{example}
 
 Theorem \ref{230722_17} gives the following.
 \begin{theorem}[label=230723_9]{}
 Let $\mk{D}$ be a compact subset of $\mb{R}^n$. If the function $\mf{F}:\mk{D}\to\mb{R}^m$ is continuous, then it is bounded.
 \end{theorem}
 \begin{myproof}{Proof}
 By Theorem  \ref{230723_28}, $\mf{F}(\mk{D})$ is  compact. Hence, it is bounded.
 \end{myproof}
 
  \begin{example}{}
 Let $\mk{D}=\left\{(x,y, z)\,|\,  1<x^2+y^2+z^2<4\right\}$, and let $\mf{F}:\mk{D}\to \mb{R}^2$ be the function defined as
 \[\mf{F}(x,y,z)=\left(\frac{1}{x^2+y^2+z^2}, x+y+z\right).\]Show that $\mf{F}:\mk{D}\to \mb{R}^2$ is a bounded function.
 \end{example}
 \begin{solution}{Solution}
 Notice that the set $\mk{D}$ is not closed. Therefore, we cannot apply Theorem \ref{230723_9} directly. Consider the set $\mathcal{U}=\left\{(x,y, z)\,|\,  1\leq x^2+y^2+z^2\leq 4\right\}$. For any $\mf{u}=(x,y,z)$ in $\mathcal{U}$, $\Vert\mf{u}\Vert\leq 2$. Hence, $\mathcal{U}$ is bounded. The function $f:\mb{R}^3\to\mb{R}$ defined as $f(x,y,z)=x^2+y^2+z^2$ is continuous, and $\mathcal{U}=f^{-1}([1, 4])$. Since $[1,4]$ is closed in $\mb{R}$, $\mathcal{U}$ is closed in $\mb{R}^3$. 
 Since $f(x,y,z)\neq 0$ on $\mathcal{U}$, \[F_1(x,y,z)=\frac{1}{x^2+y^2+z^2}\] is continuous on $\mathcal{U}$. Being a polynomial function, $F_2(x,y,z)= x+y+z$ is continuous. Thus, $\mathbf{F}:\mathcal{U} \to\mb{R}^2$ is continuous.
 Since $\mathcal{U}$ is  closed and bounded, Theorem \ref{230723_9} implies that  $\mathbf{F}:\mathcal{U} \to\mb{R}^2$ is bounded. Since $\mk{D}\subset \mathcal{U}$, $\mf{F}:\mk{D}\to \mb{R}^2$ is  also a bounded function.
 \end{solution}

 Recall that if $S$ is a subset of $\mb{R}$, $S$ has maximum value if and only if $S$ is bounded above and $\sup S$ is in $S$; while $S$ has minimum value if and only if $S$ is bounded below and $\inf S$ is in $S$.
 
 \begin{definition}{Extremizer and Extreme Values}
 Let $\mk{D}$  be a subset of $\mb{R}^n$, and let $f:\mk{D}\to\mb{R}$ be a function defined on $\mk{D}$. 
 \begin{enumerate}[1.]
 \item
 The function $f$ has  maximum value if there is a point $\mf{x}_0$ in $\mk{D}$ such that
 \[f(\mf{x}_0)\geq f(\mf{x})\hspace{1cm}\text{for all}\;\mf{x}\in\mk{D}.\]
The point  $\mf{x}_0$ is called a maximizer of $f$; and $f(\mf{x}_0)$ is the maximum value of $f$.
 \item
  The function $f$ has  minimum value if there is a point $\mf{x}_0$ in $\mk{D}$ such that
 \[f(\mf{x}_0)\leq f(\mf{x})\hspace{1cm}\text{for all}\;\mf{x}\in\mk{D}.\]
The point  $\mf{x}_0$ is called a minimizer of $f$; and $f(\mf{x}_0)$ is the minimum value of $f$.
 
 \end{enumerate}
 \end{definition}

  We have proved in volume I that a  sequentially compact subset of $\mb{R}$ has a maximum value and a minimum value.  This gives us the extreme value theorem.
 
 \begin{theorem}{Extreme Value Theorem}
 Let $\mk{D}$ be a compact subset of $\mb{R}^n$. If the function $f:\mk{D}\to\mb{R}$ is continuous, then it has a maximum value and  a minimum value.  
 \end{theorem}
  \begin{myproof}{Proof}
 By Theorem  \ref{230722_17}, $f(\mk{D})$ is a sequentially compact  subset of $\mb{R}$. Therefore, $f$ has a maximum value and  a  minimum value.  
 \end{myproof}
 
 \begin{example}{}
 Let $\mk{D}=\left\{(x,y)\,|\, x^2+2x+y^2\leq 3\right\}$, and let $f:\mk{D}\to \mb{R}$ be the function defined by
 \[f(x,y)=x^2+xy^3+e^{x-y}.\]
 Show that $f$ has a maximum value and a minimum value.
 
 \end{example}
 
 \begin{solution}{Solution}
 Notice that  \[\mk{D}=\left\{(x,y)\,|\, x^2+2x+y^2\leq 3\right\}=\left\{(x,y)\,|\, (x+1)^2+y^2\leq 4\right\}\]is a closed ball. Thus, it is closed and bounded. The function $f_1(x,y)=x^2+xy^3$ and the function $g(x,y)=x-y$ are polynomial functions. Hence, they are continuous. The exponential function $h(x)=e^x$ is continuous. Hence, the function $f_2(x,y)=(h\circ g)(x,y)=e^{x-y}$ is continuous. Since $f=f_1+f_2$, the function $f:\mk{D}\to \mb{R}$ is continuous. Since $\mk{D}$ is compact, the function $f:\mk{D}\to \mb{R}$ has a maximum value and a minimum value.
 \end{solution}
 
 \begin{remark}{Extreme Value Property}
Let $S$ be a subset of $\mathbb{R}^n$. We say that $S$ has  {\it extreme value property} provided that whenever $f:S\to\mathbb{R}$ is a continuous function, then $f$ has maximum and minimum values.

The extreme value theorem says that if $S$ is compact, then it has extreme value property. Now let us show the converse. Namely, if $S$ has extreme value property, then it is compact, or equivalently, it is closed and bounded.

If $S$ is not bounded, the function $f:S\to\mb{R}$, $f(\mf{x})=\Vert \mf{x}\Vert$ is continuous, but it does not have maximum value. If $S$ is not closed, there is  a sequence $\{\mf{x}_k\}$ in $S$ that converges to a point $\mf{x}_0$ that is not in $S$. The function $g:S\to \mb{R}$, 
$g(\mf{x})= \Vert\mf{x}-\mf{x}_0\Vert$ is continuous and $g(\mf{x})\geq 0$ for all $\mf{x}\in S$. Since $\di\lim_{k\to\infty}g(\mf{x}_k)=0$, we find that $\inf g(S)=0$. Since $\mf{x}_0$ is not in $S$, there is no point $\mf{x}$ in $S$ such that $g(\mf{x})=0$. Hence, $g$ does not have minimum value. This shows that for $S$ to have extreme value property, it is necessary that $S$ is closed and bounded.

Therefore, a subset $S$ of $\mb{R}^n$ has extreme value property if and only if it is compact.
\end{remark}

\subsection{Distance Between Sets}
 The distance between two sets is defined in the following way.
 \begin{definition}{Distance Between Two Sets}
 Let $A$ and $B$ be two subsets of $\mb{R}^n$. The distance between $A$ and $B$ is defined as
 \[d(A,B)=\inf\left\{d(\mf{a},\mf{b})\,|\, \mf{a}\in A, \mf{b}\in B\right\}.\]
 \end{definition}The distance between two sets is always well-defined and nonnegative. If $A$ and $B$ are not disjoint, then their distance is 0.
 
 \begin{example}[label=230723_10]{}
 Let $A=\left\{(x,y)\,|\,x^2+y^2<1\right\}$ and let $B=[1, 3]\times [-1,1]$. Find the distance between the two sets $A$ and $B$.
 \end{example}
 \begin{solution}{Solution}
 For $k\in\mb{Z}^+$, let $\mf{a}_k$ be the point in $A$ given by
 \[\mf{a}_k=\left(1-\frac{1}{k}, 0\right).\]
 Let $\mf{b}=(1,0)$. Then $\mf{b}$ is in $B$. Notice that
 \[d(\mf{a}_k,\mf{b})=\Vert\mf{a}_k-\mf{b}\Vert =\frac{1}{k}.\]
 Hence, $d(A,B)\leq \di\frac{1}{k}$ for all $k\in \mb{Z}^+$. This shows that the distance between $A$ and $B$ is 0.
 \end{solution}

  \begin{figure}[ht]
\centering
\includegraphics[scale=0.2]{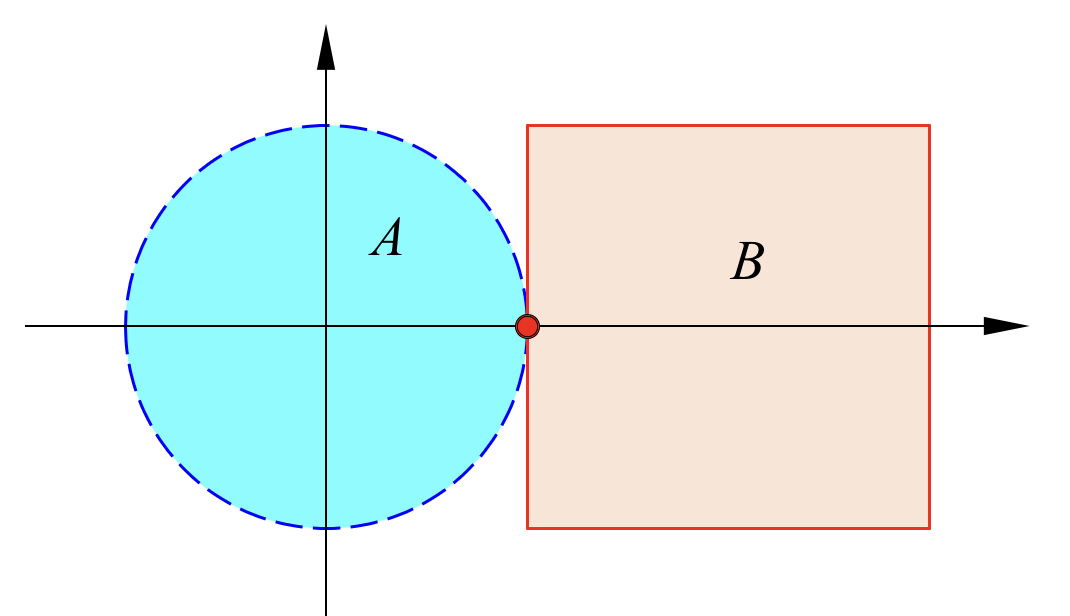}

\caption{The sets $A$ and $B$ in Example \ref{230723_10}.}\label{figure35}
\end{figure}
 In Example \ref{230723_10}, we find that the distance between two disjoint sets can be 0, even though they are both bounded.

 \begin{example}[label=230723_11]{}
 Let $A=\left\{(x,y)\,|\,y=0\right\}$ and let $B=\di \left\{(x,y)\,|\, xy=1\right\}$. Find the distance between the two sets $A$ and $B$.
 \end{example}
 \begin{solution}{Solution}
 For $k\in\mb{Z}^+$, let  $\mf{a}_k=\left(k, 0\right)$ and $\mf{b}_k=\left(k, 1/k\right)$. Then $\mf{a}_k$ is in $A$ and $\mf{b}_k$ is in $B$.
  Notice that
 \[d(\mf{a}_k,\mf{b}_k)= \Vert\mf{a}_k-\mf{b}_k\Vert =\frac{1}{k}.\]
 Hence, $d(A,B)\leq \di\frac{1}{k}$ for all $k\in \mb{Z}^+$. This shows that the distance between $A$ and $B$ is 0.
 \end{solution}
 
 \begin{figure}[ht]
\centering
\includegraphics[scale=0.2]{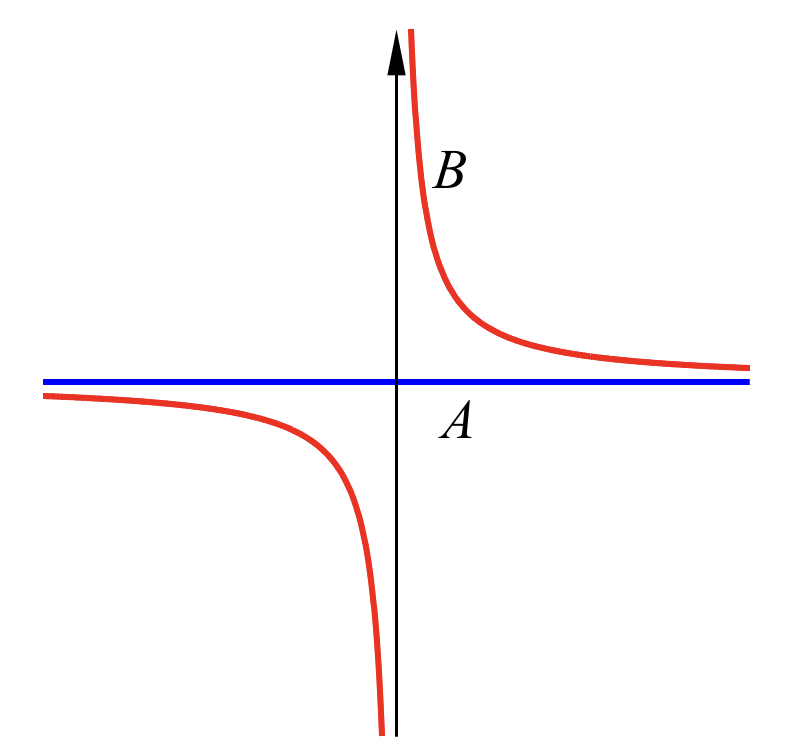}

\caption{The sets $A$ and $B$ in Example \ref{230723_11}.}\label{figure36}
\end{figure}

  In Example \ref{230723_11}, we find that the distance between two disjoint sets can be 0, even though both of them are   closed.
 
 When $B$ is the one-point set $\{\mf{x}_0\}$, the distance between $A$ and $B$ is the distance from the point $\mf{x}_0$ to the set $A$. We  denote this distance as $d(\mf{x}_0, A)$. In other words,
 \[d(\mf{x}_0, A)=\inf\left\{d(\mf{a},\mf{x}_0)\,|\,\mf{a}\in A\right\}.\]
 If $\mathbf{x}_0$ is a point in $A$, then  $d(\mf{x}_0, A)=0$. 
 However, the distance from a point $\mathbf{x}_0$ to a set $A$ can be 0 even though $\mathbf{x}_0$ is not in $A$. For example, the distance between the point $\mathbf{x}_0=(1,0)$ and the set $A=\left\{(x,y)\,|\, x^2+y^2<1\right\}$ is 0, even thought $\mathbf{x}_0$ is not in $A$. The following proposition says that this cannot happen if $A$ is closed.
 
 \begin{proposition}[label=230723_19]{}
 Let $A$ be a closed subset of $\mb{R}^n$ and let $\mf{x}_0$ be a point in $\mb{R}^n$. Then $d(\mf{x}_0, A)=0$ if and only if $\mathbf{x}_0$ is in $A$.
 \end{proposition}
 \begin{myproof}{Proof}
 If $ \mathbf{x}_0$ is in $A$, it is obvious that $d(\mf{x}_0, A)=0$.
 
 Conversely, if $ \mathbf{x}_0$ is   not in $A$, $\mathbf{x}_0$ is in the open set $\mathbb{R}^n\setminus A$. Therefore, there is an $r>0$ such that $B(\mf{x}_0, r)\subset \mathbb{R}^n\setminus A$. For any $\mf{a}\in A$, $\mf{a}\notin B(\mf{x}_0, r)$. Therefore, $\Vert\mf{x}_0- \mf{a}\Vert \geq r$. Taking infimum over $\mf{a}\in A$, we find that 
 $d(\mf{x}_0, A)\geq r$. Hence, $d(\mf{x}_0, A)\neq 0$.
 \end{myproof}
 
 \begin{figure}[ht]
\centering
\includegraphics[scale=0.2]{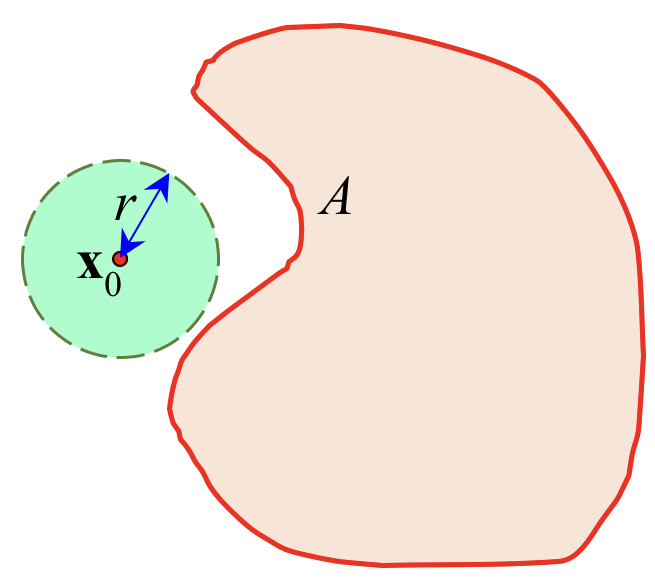}

\caption{A point outside a closed set has positive distance from the set.}\label{figure38}
\end{figure}
 \begin{proposition}[label=230723_12]{}
  Given a subset $A$ of $\mb{R}^n$, define the function $f:\mb{R}^n\to \mb{R}$ by
 \[f(\mf{x})=d(\mf{x}, A).\] Then $f$ is a continuous function.
 
 \end{proposition}
 \begin{myproof}{Proof}
 We prove something stronger. For any $\mathbf{u}$ and $\mf{v}$ in $\mb{R}^n$, we claim that
 \[|f(\mf{u})-f(\mf{v})|\leq\Vert \mf{u}-\mf{v}\Vert.\]This means that $f$ is a Lipschitz function with Lipschitz constant 1, which   implies that it is continuous.
 
 Given $\mathbf{u}$ and $\mf{v}$ in $\mb{R}^n$, if $\mf{a}$ is in $A$, then
 \[d(\mf{u},A)\leq \Vert \mf{u}-\mf{a}\Vert\leq \Vert\mf{v}-\mf{a}\Vert  +\Vert\mf{u}-\mf{v}\Vert.\]
 
 This shows that
 \[\Vert\mf{v}-\mf{a}\Vert\geq d(\mf{u},A)-\Vert\mf{u}-\mf{v}\Vert.\]
 Taking infimum over $\mf{a}\in A$, we find that 
  \[d(\mf{v},A)\geq d(\mf{u},A)-\Vert\mf{u}-\mf{v}\Vert.\]

  Therefore,
  \[f(\mf{u})-f(\mf{v})\leq \Vert\mf{u}-\mf{v}\Vert.\]

  Interchanging $\mf{u}$ and $\mf{v}$, we obtain 
    \[f(\mf{v})-f(\mf{u})\leq \Vert\mf{u}-\mf{v}\Vert.\] This proves that  \[|f(\mf{u})-f(\mf{v})|\leq\Vert \mf{u}-\mf{v}\Vert.\]
 \end{myproof}
 
 Now we can prove the following.
 \begin{theorem}[label=230723_13]{}
 Let $A$ and $C$ be disjoint subsets of $\mb{R}^n$. If $A$ is compact and $C$ is closed, then the distance between $A$ and $C$  is positive.
 \end{theorem}

 \begin{myproof}{Proof}
 By Proposition \ref{230723_12}, the function $f: A\to \mb{R}$, $f(\mf{a})=d(\mf{a}, C)$ is continuous. Since $A$ is compact,  $f$ has a minimum value. Namely, there is a point $\mf{a}_0$ in $A$ such that
 \[d(\mf{a}_0, C)\leq d(\mf{a}, C)\hspace{1cm}\text{for all}\;\mf{a}\in A.\]
 For any $\mf{a}$ in $A$ and $\mf{c}\in C$,
\[d(\mf{a}, \mf{c})\geq d(\mf{a}, C)\geq d(\mf{a}_0, C).\]
Taking infimum over all $\mf{a}\in A$ and $\mf{c}\in C$, we find that
\[d(A, C)\geq d(\mf{a}_0, C).\] By definition, we also have $ d(A,C)\leq d(\mf{a}_0, C)$. Thus, $ d(A,C)= d(\mf{a}_0, C)$. 
Since $A$ and $C$ are disjoint and $C$ is closed,  Proposition \ref{230723_19} implies that $d(A,C)=d(\mf{a}_0, C)>0$. 
 
 \end{myproof}

 An equivalent form of Theorem \ref{230723_13} is the following important theorem.
 \begin{theorem}[label=230906_1]{}
 Let $A$ be a compact subset of $\mb{R}^n$, and let $U$ be an open subset of $\mb{R}^n$ that contains $A$. Then there is a positive number $\delta$ such that if $\mf{x}$ is a point in $\mb{R}^n$ that has a distance less than $\delta$ from the set $A$, then $\mf{x}$ is in $U$.
 \end{theorem}
 
  \begin{figure}[ht]
\centering
\includegraphics[scale=0.2]{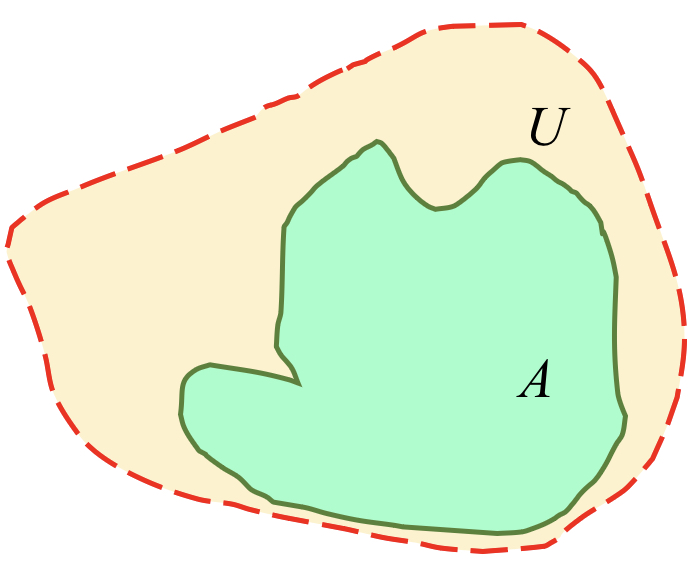}

\caption{A compact set has a positive distance from the boundary of the open set that contains it.}\label{figure37}
\end{figure}
 \begin{myproof}{Proof}
 Let $C=\mb{R}^n\setminus U$. Then $C$ is a closed subset of $\mb{R}^n$ that is disjoint from $A$. By Theorem \ref{230723_13}, $\delta=d(A, C)>0$. If $\mf{x}$ is in $\mb{R}^n$ and $d(\mf{x}, A)<\delta$, then $\mf{x}$ cannot be in $C$. Therefore, $\mf{x}$ is in $U$.
 \end{myproof}
 
As a corollary, we have the following.
\begin{corollary}[label=230905_12]{}
 Let $A$ be a compact subset of $\mb{R}^n$, and let $U$ be an open subset of $\mb{R}^n$ that contains $A$. Then there is a positive number $r$ and a compact set $K$   such that $A\subset K\subset U$, and if $\mf{x}$ is a point in $\mb{R}^n$ that has a distance less than $r$ from the set $A$, then $\mf{x}$ is in $K$.
 \end{corollary}
 \begin{myproof}{Proof}
 By Theorem \ref{230906_1}, there is a positive number $\delta$ such that  if $\mf{x}$ is a point in $\mb{R}^n$ that has a distance less than $\delta$ from the set $A$, then $\mf{x}$ is in $U$. Take $r=\delta/2$, and let
 \[K=\overline{V}, \hspace{1cm}\text{where}\;\;V=\bigcup_{\mf{u}\in A}B(\mf{u}, r).\]
 Since $A$ is compact, it is bounded. There is a positive number $M$ such that $\Vert\mf{u}\Vert\leq M$ for all $\mf{u}\in A$. If $\mf{x}\in V$, then there is an $\mf{u}\in A$ such that $\Vert\mf{x}-\mf{u}\Vert<r$. This implies that $\Vert\mf{x}\Vert\leq M+r$. Hence, the set $V$ is also bounded. Since $K$ is the closure of a bounded set, $K$ is compact. 
 Since $A\subset V$, $A\subset K$. If $\mf{w}\in K$, since $K$ is the closure of $V$, there is a point $\mf{v}$ in $V$ that lies in $B(\mf{w}, r)$. By the definition of $V$, there is a point $\mf{u}$ in $A$ such that $\mf{v}\in B(\mf{u}, r)$.
 Thus,
 \[\Vert \mf{w}-\mf{u}\Vert\leq \Vert\mf{w}-\mf{v}\Vert+\Vert\mf{v}-\mf{u}\Vert<r+r=\delta.\]
 This implies that $\mf{w}$ has a distance less than $\delta$ from $A$. Hence, $\mf{w}$ is in $U$. This shows that $K\subset U$.
 
 Now if $\mf{x}$ is a point that has distance $d$ less than $r$ from the set $A$,  there is a point $\mf{u}$ is $A$ such that $\Vert \mf{x}-\mf{u}\Vert<r$. This implies that $\mf{x}\in B(\mf{u},r)\in V\subset K$. 
 
 \end{myproof}
 
\subsection{Uniform Continuity}
In Section \ref{sec1_9}, we have discussed uniform continuity.
Let $\mk{D}$ be a subset of $\mb{R}^n$ and let $\mf{F}:\mk{D}\to\mb{R}^m$ be a  function defined on $\mk{D}$. We say that $\mf{F}:\mk{D}\to\mb{R}^m$ is uniformly continuous provided that for any $\varepsilon>0$, there exists $\delta>0$ such that for any points $\mf{u}$ and $\mf{v}$ in $\mk{D}$, if $\Vert\mf{u}-\mf{v}\Vert<\delta$, then $\Vert\mf{F}(\mf{u})-\mf{F}(\mf{v})\Vert<\varepsilon$. If a function is uniformly continuous, it is continuous. The converse is not true. 
However, a continuous function that is defined on a compact subset of $\mb{R}^n$ is  uniformly continuous.  This is an important theorem in analysis.

\begin{theorem}[label=230725_1]{}
Let $\mk{D}$ be a subset of $\mb{R}^n$, and let $\mf{F}:\mk{D}\to\mb{R}^m$ be a continuous function defined on $\mk{D}$. If $\mk{D}$ is compact, then $\mf{F}:\mk{D}\to\mb{R}^m$ is uniformly continuous.
\end{theorem}
\begin{myproof}{Proof}
Assume to the contrary that $\mf{F}:\mk{D}\to\mb{R}^m$ is not uniformly continuous. Then there exists an $\varepsilon>0$, for any $\delta>0$, there exist points $\mf{u}$ and $\mf{v}$ in $\mk{D}$ such that $\Vert\mf{u}-\mf{v}\Vert<\delta$ and $\Vert\mf{F}(\mf{u})-\mf{F}(\mf{v})\Vert\geq \varepsilon$. 
This implies that for any $k\in \mb{Z}^+$, there exist $\mf{u}_k$ and $\mf{v}_k$ in $\mk{D}$ such that $\Vert\mf{u}_k-\mf{v}_k\Vert<1/k$ and $\Vert\mf{F}(\mf{u}_k)-\mf{F}(\mf{v}_k)\Vert\geq \varepsilon$. Since $\mk{D}$ is sequentially compact, there is a subsequence $\{\mf{u}_{k_j}\}$ of $\{\mf{u}_k\}$ that converges to a point $\mf{u}_0$ in $\mk{D}$. Consider the sequence $\{\mf{v}_{k_j}\}$ in $\mk{D}$.   It has  a subsequence $\{\mf{v}_{k_{j_l}}\}$ that converges to a point $\mf{v}_0$ in $\mk{D}$. Being a subsequence of $\{\mf{u}_{k_j}\}$, the sequence $\{\mf{u}_{k_{j_l}}\}$ also converges to $\mf{u}_0$.

Since $\mf{F}:\mk{D}\to\mb{R}^m$ is continuous, the sequences $\{\mf{F}(\mf{u}_{k_{j_l}})\}$ and $\{\mf{F}(\mf{v}_{k_{j_l}})\}$ converge to $\mf{F}(\mf{u}_0)$ and $\mf{F}(\mf{v}_0)$ respectively. Notice that by construction, 
\[\Vert \mf{F}(\mf{u}_{k_{j_l}})-\mf{F}(\mf{v}_{k_{j_l}})\Vert\geq \varepsilon\hspace{1cm}\text{for all}\;l\in\mb{Z}^+.\]

Thus, 
$\Vert\mf{F}(\mf{u}_0)-\mf{F}(\mf{v}_0)\Vert\geq \varepsilon$. This implies that $\mf{F}(\mf{u}_0)\neq \mf{F}(\mf{v}_0)$, and so $\mf{u}_0\neq\mf{v}_0$.

  Since $k_{j_1}, k_{j_2}, \ldots$ is a strictly increasing sequence of positive integers, $k_{j_l}\geq l$. Thus,
\[\Vert\mf{u}_{k_{j_l}}-\mf{v}_{k_{j_l}}\Vert<\frac{1}{k_{j_l}}\leq\frac{1}{l}.\]
\bp
Taking $l\to\infty$ implies that $\mf{u}_0=\mf{v_0}$. This gives a contradiction.
Thus, $\mf{F}:\mk{D}\to\mb{R}^m$ must be uniformly continuous.
\end{myproof}

\begin{example}{}
Let $\mk{D}=(-1, 4)\times (-7, 5]$ and let $\mf{F}:\mk{D}\to\mb{R}^3$ be the function defined as
\[\mf{F}(x,y)=\left(\sin (x+y), \sqrt{x+y+8}, e^{xy}\right).\]
Show that $\mf{F}$ is uniformly continuous.
\end{example}
\begin{solution}
{Solution}
Let $\mathcal{U}=[-1, 4]\times [-7, 5]$. Then $\mathcal{U}$ is a closed and bounded subset of $\mb{R}^2$ that contains $\mk{D}$. The functions $f_1(x,y)=x+y$, $f_2(x,y)=x+y+8$ and $f_3(x,y)=xy$ are polynomial functions. Hence, they are continuous. If $(x,y)\in \mathcal{U}$, $x\geq -1$, $y\geq -7$ and so $f_2(x,y)=x+y+8\geq 0$. Thus, $f_2(\mathcal{U})$ is contained in the domain of the square root function. Since the square root function, the sine function and the exponential function are continuous on their domains, we find that the functions
\[F_1(x,y)=\sin (x+y),\quad F_2(x,y)=\sqrt{x+y+8}, \quad F_3(x,y)=e^{xy}\]are continuous on $\mathcal{U}$. Since $\mathcal{U}$ is closed and bounded, $\mf{F}:\mathcal{U}\to\mb{R}^3$ is uniformly continuous. Since $\mk{D}\subset\mathcal{U}$, $\mf{F}:\mk{D}\to\mb{R}^3$ is uniformly continuous.
\end{solution}

 \subsection{Linear Transformations and Quadratic Forms}

  In Chapter \ref{chapter2}, we have seen that a linear transformation $\mf{T}:\mb{R}^n\to\mb{R}^m$ is a matrix transformation. Namely, there exists an $m\times n$ matrix such that 
  \[\mf{T}(\mf{x})=A\mf{x}\hspace{1cm}\text{for all}\;\mf{x}\in\mb{R}^n.\]A linear transformation is continuous. 
    Theorem \ref{230724_9} says that a linear transformation is Lipschitz. More precisely, there exists a positive constant $c>0$ such that
    \[\Vert\mf{T}(\mf{x})\Vert\leq c\Vert\mf{x}\Vert\hspace{1cm}\text{for all}\;\mf{x}\in\mb{R}^n.\]

  Theorem \ref{230725_2} says that when $m=n$, a linear transformation $\mf{T}:\mb{R}^n\to\mb{R}^n$  is invertible if and only if it is one-to-one, if and only if the matrix $A$ is invertible, if and only if $\det A\neq 0$. Here we want to give a stronger characterization of a linear transformation $\mf{T}:\mb{R}^n\to\mb{R}^n$ that is invertible.
  
  Recall that to show that a linear transformation $\mf{T}:\mb{R}^n\to\mb{R}^m$  is one-to-one, it is sufficient to show that $\mf{T}(\mf{x})=\mf{0}$ implies that $\mf{x}=\mf{0}$.  
  
  \begin{theorem}[label=230725_6]{}
  Let $\mf{T}:\mb{R}^n\to\mb{R}^n$ be a linear transformation. The following are equivalent.
  \begin{enumerate}[(a)]
  \item
  $\mf{T}$ is invertible.
  \item There is a positive constant $a$ such that
    \[\Vert\mf{T}(\mf{x})\Vert\geq a\Vert\mf{x}\Vert\hspace{1cm}\text{for all}\;\mf{x}\in\mb{R}^n.\]
  \end{enumerate}
  \end{theorem}
  \begin{myproof}{Proof}
  (b) implies (a) is easy. Notice that (b) says that
  \begin{equation}\label{230725_4}\Vert\mf{x}\Vert\leq \frac{1}{a}\Vert\mf{T}(\mf{x})\Vert\hspace{1cm}\text{for all}\;\mf{x}\in\mb{R}^n.\end{equation} If $\mf{T}(\bf{x})=\mf{0}$, then  $\Vert\mf{T}(\mf{x})\Vert=0$. Eq. \eqref{230725_4} implies that $\Vert\mf{x} \Vert=0$. Thus, $\mf{x}=\mf{0}$. This proves that $\mf{T}$ is one-to-one. Hence, it is invertible.

  Conversely, assume that $\mf{T}:\mb{R}^n\to\mb{R}^n$ is invertible. Let \[S^{n-1}=\left\{(x_1, \ldots, x_n)\,|\,x_1^2+\cdots+x_n^2=1\right\}\]
   be the standard unit $(n-1)$-sphere in $\mb{R}^n$.  We have seen that $S^{n-1}$ is compact. 
  For any $\mf{u}\in S^{n-1}$, $\mf{u}\neq \mf{0}$. Therefore, $\mf{T}(\mf{u})\neq \mf{0}$ and so $\Vert\mf{T}(\mf{u})\Vert>0$.
 The function $f:S^{n-1}\to\mb{R}^n$,
 $\di f(\mf{u})= \Vert \mf{T}(\mf{u})\Vert $ is continuous. Hence, it has a mimimum value at some $\mf{u}_0$ on $S^{n-1}$. Let $a=\Vert \mf{T}(\mf{u}_0)\Vert$. Then $a>0$. Since $a$ is the minimum value of $f$, 
  \[\Vert \mf{T}(\mf{u})\Vert\geq a\hspace{1cm} \text{for all}\;\mf{u}\in S^{n-1}.\]
  \bp
Notice that if $\mf{x}=\mf{0}$, $\Vert\mf{T}(\mf{x})\Vert\geq a\Vert\mf{x}\Vert$ holds trivially.  If $\mf{x}$ is in $\mb{R}^n$ and $\mf{x}\neq 0$, let $\mf{u}=\alpha\mf{x}$, where  $\alpha=1/\Vert\mf{x}\Vert$. Then $\mf{u}$ is in $S^{n-1}$. Therefore,
   $\di \Vert \mf{T}(\mf{u})\Vert\geq a$. 
   Since $\mf{T}(\mf{u})=\alpha \mf{T}(\mf{x})$, and $\alpha>0$, we find that
$\Vert   \mf{T}(\mf{u})\Vert=\alpha\Vert\mf{T}(\mf{x})\Vert$. Hence, $\di \alpha\Vert\mf{T}(\mf{x})\Vert\geq a$. This gives
\[\Vert\mf{T}(\mf{x})\Vert \geq \frac{a}{\alpha}=a\Vert\mf{x}\Vert.\]
 
  \end{myproof}

  In Section \ref{quadraticforms}, we have reviewed some theories of quadratic forms from linear algebra. In Theorem \ref{230724_8}, we state for a quadratic form $Q_A:\mb{R}^n\to \mb{R}$, $Q_A(\mf{x})=\mf{x}^TA\mf{x}$ defined by the symmetric matrix $A$, we have
 \[\lambda_n\Vert\mf{x}\Vert^2\leq Q_A(\mf{x})\leq \lambda_1\Vert \mf{x}\Vert^2 \hspace{1cm}\text{for all}\; \mf{x}\in\mb{R}^n.\]Here $\lambda_n$ is the smallest eigenvalue of $A$, and $\lambda_1$ is the largest eigenvalue of $A$.
 
  We have used Theorem \ref{230724_8} to prove that a linear transformation is Lipschitz in Theorem \ref{230724_9}. It boils down to the fact that
  if $\mf{T}(\mf{x})=A\mf{x}$, then $\Vert\mf{T}(\mf{x})\Vert^2=\mf{x}^T(A^TA)\mf{x}$, and $A^TA$ is a positive semi-definite quadractic form. In fact, we can also use Theorem \ref{230724_8} to prove Theorem \ref{230725_6}, using the fact that if $A$ is invertible, then $A^TA$ is positive definite.

 Let us prove a weaker version of  Theorem \ref{230724_8} here, which is sufficient to establish Theorem \ref{230725_6} and the theorem which says that a linear transformation is Lipschitz.
 
 \begin{theorem}[label=230725_7]{}
 Let $A$ be an $n\times n$ symmetric matrix, and let $Q_A:\mb{R}^n\to\mb{R}$ be the quadratic form $Q_A(\mf{x})=\mf{x}^TA\mf{x}$ defined by $A$. There exists constants $a$ and $b$ such that
 \[a\Vert \mf{x}\Vert^2 \leq Q_A(\mf{x})\leq b\Vert \mf{x}\Vert^2\hspace{1cm}\text{for all}\;\mf{x}\in\mb{R}^n,\]
 $ Q_A(\mf{u})=a\Vert\mf{u}\Vert^2$ and $Q_A(\mf{v})=b\Vert\mf{v}\Vert^2$ for some $\mf{u}$ and $\mf{v}$ in $\mb{R}^n$. Therefore,
 \begin{enumerate}[(i)]
 \item
 if $A$ is positive semi-definite, $b\geq a\geq 0$;
 \item if $A$ is positive definite, $b\geq a>0$.
 \end{enumerate}
 \end{theorem}
 \begin{myproof}{Proof}
 As in the proof of  Theorem \ref{230725_6}, consider the continuous function  $Q_A: S^{n-1}\to\mb{R}$.  Since $S^{n-1}$ is compact, there exsits $\mf{u}$ and $\mf{v}$ in $S^{n-1}$ such that
 \[Q_A(\mf{u})\leq Q_A(\mf{w})\leq Q_A(\mf{v})\hspace{1cm}\text{for all}\;\mf{w}\in S^{n-1}.\]
 Let $a=Q_A(\mf{u})$ and $b=Q_A(\mf{v})$.  
 If $\mf{x}=\mf{0}$, $a\Vert \mf{x}\Vert^2 \leq Q_A(\mf{x})\leq b\Vert \mf{x}\Vert^2$ holds trivially.
 Now if $\mf{x}$ is in $\mb{R}^n$ and $\mf{x}\neq \mf{0}$, let $\mf{w}=\alpha\mf{x}$, where $\alpha=1/\Vert\mf{x}\Vert$. Then $\mf{w}$ in on $S^{n-1}$. Notice that 
 \[Q_A(\mf{w})=\alpha^2Q_A(\mf{x}).\]
 Hence,
 \[a\leq \frac{1}{\Vert\mf{x}\Vert^2}Q_A(\mf{x})\leq b.\]
 This proves that
 \[a\Vert \mf{x}\Vert^2 \leq Q_A(\mf{x})\leq b\Vert \mf{x}\Vert^2.\]
 \end{myproof}

 \subsection{Lebesgue Number Lemma}
Now let us prove the following important theorem.
\begin{theorem}{Lebesgue Number Lemma}
Let $A$ be a subset of $\mb{R}^n$, and let $\mathscr{A}=\left\{U_{\alpha}\,|\,\alpha\in J\right\}$ be an open covering of $A$. If $A$ is compact, there exists a positive number $\delta$ such that if $S$ is a subset of $A$ and $\text{diam}\,S<\delta$, then $S$ is contained in one of the elements of $\mathscr{A}$. Such a positive number $\delta$ is called the Lebesgue number of the covering $\mathscr{A}$.
\end{theorem}
We give two proofs of this theorem.

\begin{myproof}{First  Proof of the Lebesgue Number Lemma}
We use proof by contradiction. Assume that there does not exist a positive number $\delta$ such that any subset $S$ of $A$ that has diameter less than $\delta$ lies inside an open set in $\mathscr{A}$. Then for any $k\in \mathbb{Z}^+$, there is a subset $S_k$ of $A$ whose diameter is less than $1/k$, but $S_k$ is not contained in any element of $\mathscr{A}$. 

For each $k\in\mb{Z}^+$, 
the set $S_k$ cannot be empty. Let $\mf{x}_k$ be any point in $S_k$. Then $\{\mf{x}_k\}$ is a sequence of points in $A$. Since $A$ is  sequentially compact, there is a subsequence $\{\mathbf{x}_{k_m}\}$ that converges to a point $\mf{x}_0$ in $A$. Since $\mathscr{A}$ is an open covering of $A$, there exists $\beta\in J$ such that $\mf{x}_0\in U_{\beta}$. Since $U_{\beta}$ is open, there exists $r>0$ such that $B(\mf{x}_0, r)\subset U_{\beta}$. Since the sequence $\{\mathbf{x}_{k_m}\}$   converges   $\mf{x}_0$, there is a positive integer $M$ such that for all $m\geq M$, $\mf{x}_{k_m}\in B(\mf{x}_0, r/2)$. There exists an integer $j\geq M$ such that $1/k_j <r/2$. If $\mf{x}\in A_{k_j}$, then
\[\Vert\mf{x}-\mf{x}_{k_j}\Vert\leq\text{diam}\, A_{k_j}<\frac{1}{k_j}<\frac{r}{2}.\]
Since $\mf{x}_{k_j}\in B(\mf{x}_0, r/2)$, $\Vert\mf{x}_{k_j}-\mf{x}_0\Vert<r/2$. Therefore,
$\Vert\mf{x}-\mf{x}_0\Vert<r$. This proves that $\mf{x}\in B(\mf{x}_0, r)\subset U_{\beta}$. Thus, we have shown that $A_{k_j}\subset U_{\beta}$. But this contradicts to $A_{k_j}$ does not lie in any element of $\mathscr{A}$.

\end{myproof}
\begin{myproof}{Second Proof of the Lebesgue Number Lemma}
Since $A$ is compact, there are finitely many indices $\alpha_1, \ldots, \alpha_m$ in $J$ such that
\[A\subset\bigcup_{j=1}^m U_{\alpha_j}.\]

For $1\leq j\leq m$, let $C_{j}=\mathbb{R}^n\setminus U_{\alpha_j}$. Then $C_{j}$ is a closed set and
$\di\bigcap_{j=1}^m C_j$ is disjoint from $A$. By Theorem \ref{230723_12}, the function $f_j:A\to \mb{R}$, $f_j(\mf{x})=d(\mf{x}, C_j)$ is continuous. Define $f:A\to \mb{R}$ by
\[f(\mf{x})=\frac{f_1(\mf{x})+\cdots +f_{m}(\mf{x})}{m}.\]
 
 \bp 
Then $f$ is also a continuous function. Since $A$ is compact, there is a point $\mf{a}_0$ in $A$ such that
\[f(\mf{a}_0)\leq f(\mf{a})\hspace{1cm}\text{for all}\;\mf{a}\in A.\]

Notice that $f_j(\mf{a}_0)\geq 0$ for all $1\leq j\leq m$.
Since $\di \bigcap_{j=1}^m C_j$ is disjoint from $A$, there is an $1\leq k\leq m$ such that $\mf{a}_0\notin C_k$. Proposition 
\ref{230723_19} says that $f_k(\mf{a}_0)=d(\mf{a}_0, C_k)>0$. Hence, $ f(\mf{a}_0)>0$. Let $\delta= f(\mf{a}_0)$. It is the minimum value of the function $f:A\to \mb{R}$.

Now let $S$ be a nonempty subset of $A$ such that $\text{diam}\,S<\delta$. Take a point $\mf{x}_0$   in $S$. Let $1\leq l\leq m$ be an integer such that
\[f_l(\mf{x}_0)\geq f_j(\mf{x}_0)\hspace{1cm}\text{for all}\;1\leq j\leq m.\]
Then 
\[\delta\leq f(\mf{x}_0)\leq f_l(\mf{x}_0)=d(\mf{x}_0, C_l).\]
For any $\mf{u}\in C_l$,
\[d(\mf{x}_0,\mf{u})\geq d(\mf{x}_0, C_l)\geq \delta.\]
 If $\mf{x} \in S$, then $d(\mf{x},\mf{x}_0)\leq\text{diam}\,S<\delta$. This implies that $\mf{x}$ is not in $C_l$. Hence, it must  be in $U_{\alpha_l}$. This shows that $S$ is contained in $U_{\alpha_l}$, which is an element of $\mathscr{A}$. This completes the proof of the theorem.
\end{myproof}

The Lebesgue number lemma can be used to give an alternative proof of Theorem \ref{230725_1}, which says that a continuous function defined on a compact subset of $\mb{R}^n$ is uniformly continuous.

\begin{myproof}{Alternative Proof of \linkt Theorem \ref{230725_1}}
Fixed $\varepsilon>0$. We want to show that there exists $\delta>0$ such that if $\mf{u}$ and $\mf{v}$ are in $\mk{D}$ and $\Vert\mf{u}-\mf{v}\Vert<\delta$, then $\Vert \mf{F}(\mf{u})-\mf{F}(\mf{v})\Vert<\varepsilon$.

We will construct an open covering of $\mk{D}$ indexed by $J=\mk{D}$. Since $\mf{F}:\mk{D}\to\mb{R}^m$ is continuous, for each $\mf{x}\in \mk{D}$, there is a positive number $\delta_{\mf{x}}$ (depending on $\mf{x}$), such that if $\mf{u}$ is in $\mk{D}$ and $\Vert\mf{u}-\mf{x}\Vert<\delta_{\mf{x}}$, then $\Vert \mf{F}(\mf{u})-\mf{F}(\mf{x})\Vert<\varepsilon/2$. Let $U_{\mf{x}}=B(\mf{x}, \delta_{\mf{x}})$. Then $U_{\mf{x}}$ is an open set. If $\mf{u}$ and $\mf{v}$ are in $U_{\mf{x}}$, $\Vert \mf{F}(\mf{u})-\mf{F}(\mf{x})\Vert<\varepsilon/2$ and $\Vert \mf{F}(\mf{v})-\mf{F}(\mf{x})\Vert<\varepsilon/2$. Thus,
 $\Vert \mf{F}(\mf{u})-\mf{F}(\mf{v})\Vert<\varepsilon$. 
 
 Now $\mathscr{A}=\left\{U_{\mf{x}}\,|\,\mf{x}\in\mk{D}\right\}$ is an open covering of $\mk{D}$. Since $\mk{D}$ is compact, the Lebesgue number lemma implies that there exists a number $\delta>0$ such that if $S$ is a subset of $\mk{D}$ that has diameter less than $\delta$, then $S$ is contained in one of the $U_{\mf{x}}$ for some $\mf{x}\in\mk{D}$. We claim that this is the $\delta$ that we need. 
 
 If $\mf{u}$ and $\mf{v}$ are two points in $\mk{D}$ and $\Vert\mf{u}-\mf{v}\Vert<\delta$, then $S=\{\mf{u}, \mf{v}\}$ is a set with diameter less than $\delta$. Hence, there is an $\mf{x}\in \mk{D}$ such that  $S\subset U_{\mf{x}}$. This implies that
 $\mf{u}$ and $\mf{v}$ are in $U_{\mf{x}}$. Hence, $\Vert \mf{F}(\mf{u})-\mf{F}(\mf{v})\Vert<\varepsilon$. This completes the proof.
\end{myproof}

\vp
\noindent
{\bf \large Exercises  \thesection}
\setcounter{myquestion}{1}
\begin{question}{\themyquestion}
Let $\mk{D}=\left\{(x,y)\,|\, 2<x^2+4y^2<10\right\}$, and let $\mf{F}:\mk{D}\to\mb{R}^3$ be the function defined as
\[\mf{F}(x,y)=\left(\frac{x}{x^2+y^2}, \frac{y}{x^2+y^2}, \frac{x^2-y^2}{x^2+y^2}\right).\]
Show that the function $\mf{F}:\mk{D}\to\mb{R}^3$  is bounded.
\end{question}
\atc
\begin{question}{\themyquestion}
Let $\mk{D}=\left\{(x,y,z)\,|\,1\leq x^2+4y^2\leq 10, 0\leq z\leq 5\right\}$, and let $f:\mk{D}\to\mb{R}$ be the function defined as
\[f(x,y,z)= \frac{x^2-y^2}{x^2+y^2+z^2}.\]
Show that the function $f:\mk{D}\to\mb{R}$   has a maximum value and a minimum value.
\end{question}

\atc
\begin{question}{\themyquestion}
Let $A=\left\{(x,y)\,|\,  x^2+4y^2\leq 16\right\}$  and $B=\left\{(x,y)\,|\, x+y\geq 10\right\}$. Show that the distance between the sets $A$ and $B$ is positive.
\end{question}
 \atc
\begin{question}{\themyquestion}
Let $\mk{D}=\left\{(x,y,z)\,|\, x^2+y^2+z^2\leq 20\right\}$ and let $f:\mk{D}\to \mb{R}$ be the function defined as
\[f(x,y,z)=e^{x^2+4z^2}.\]
Show that $f:\mk{D}\to \mb{R}$ is uniformly continuous.
\end{question}
 \atc
\begin{question}{\themyquestion}
Let $\mk{D}=(-1, 2)\times (-6, 0)$ and let $f:\mk{D}\to \mb{R}$ be the function defined as
\[f(x,y)=\sqrt{x+y+7}+\ln(x^2+y^2+1).\]
Show that $f:\mk{D}\to \mb{R}$ is uniformly continuous.
\end{question}

\chapter{Differentiating Functions of Several Variables}\label{chapter4}

In this chapter, we study differential calculus of functions of several variables. 
 
\section{Partial Derivatives} 
When $f:(a,b)\to\mb{R}$ is a function defined on an open interval $(a,b)$, the derivative of the function at a point $x_0$ in $(a,b)$ is defined as
\[f'(x_0)=\lim_{h\to 0}\frac{f(x_0)+h)-f(x_0)}{h},\]
provided that the limit exists. The derivative gives the instantaneous rate of change of the function at the point $x_0$. Geometrically, it is the slope of the tangent line to the graph of  the function $f:(a,b)\to\mb{R}$  at the point $(x_0, f(x_0))$. 

\begin{figure}[ht]
\centering
\includegraphics[scale=0.2]{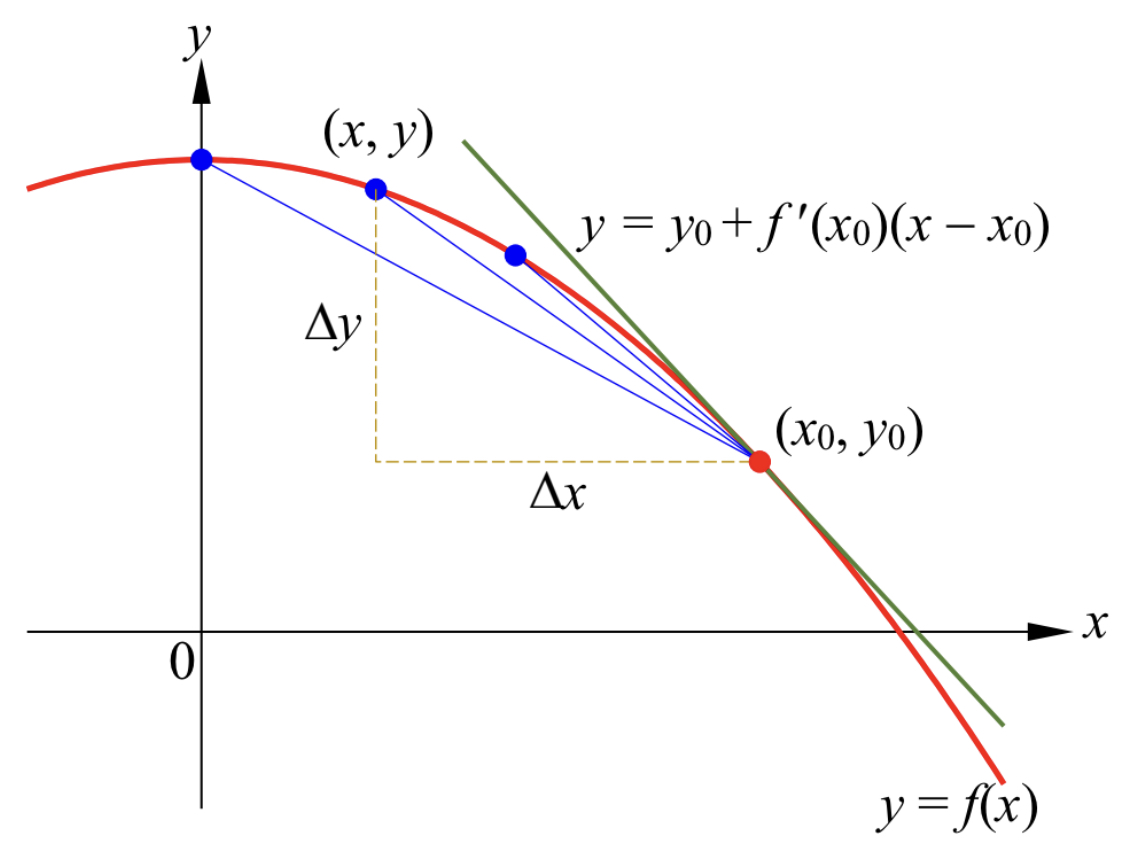}

\caption{Derivative as slope of tangent line.}\label{figure40}
\end{figure}

Now   consider a function  $f:\mathcal{O}\to\mb{R}$ that is  defined on an open subset  $\mathcal{O}$ of $\mb{R}^n$, where $n\geq 2$. What is the  natural way to extend the concept of derivatives to this function?

From the perspective  of rate of change, we   need to consider the change of $f$ in various  {\it different directions}. This   leads us to consider directional derivatives. Another perspective is to regard existence of derivatives as {\it differentiability} and {\it first-order approximation}. Later we will see that all these are closely related.

First let us consider the rates of change of the function  $f:\mathcal{O}\to\mb{R}$ at a point $\mathbf{x}_0$ in $\mathcal{O}$ along the directions of the coordinate axes. These are called partial derivatives.

\begin{definition}{Partial Derivatives}
Let $\mathcal{O}$ be an open subset of $\mb{R}^n$ that contains the point $\mf{x}_0$, and let $f:\mathcal{O}\to \mb{R}$ be a function defined on $\mathcal{O}$. For $1\leq i\leq n$, we say that the function $f:\mathcal{O}\to \mb{R}$ has a partial derivative with respect to its $i^{\text{th}}$ component at the point $\mf{x}_0$ if the limit
\[\lim_{h\to 0}\frac{f(\mf{x}_0+h\mf{e}_i)-f(\mf{x}_0)}{h}\] exists. In this case, we denote the limit by $\di \frac{\pa f}{\pa x_i}(\mf{x}_0)$, and call it the partial derivative of $f:\mathcal{O}\to \mb{R}$  with respect to $x_i$ at $\mf{x}_0$.

We say that the function $f:\mathcal{O}\to \mb{R}$ has partial derivatives at $\mf{x}_0$ if $\di \frac{\pa f}{\pa x_i}(\mf{x}_0)$ exists for all $1\leq i\leq n$.
\end{definition}

\begin{remark}{}
When we consider partial derivatives of a function, we always assume   that   the domain of the function is an open set $\mathcal{O}$, so that each point $\mf{x}_0$ in the domain is an interior point of $\mathcal{O}$, and   a limit point of $\mathcal{O}\setminus\{\mf{x}_0\}$. By definition of open sets, there exists $r>0$ such that $B(\mf{x}_0, r)$ is contained in $\mathcal{O}$. This allows us to compare the function values of $f$ in a neighbourhood of $\mf{x}_0$ from various different directions.

By definition, $\di\frac{\pa f}{\pa x_i}(\mf{x}_0)$ measures the rate of change of $f$ at $\mf{x}_0$ in the direction of $\mf{e}_i$.  It can also be interpreted as the slope of a curve at the point $(\mf{x}_0, f(\mf{x}_0))$ on the surface $x_{n+1}=f(\mf{x})$, as shown in Figure \ref{figure41}\end{remark}

\begin{highlight}{Notations for Partial Derivatives}
An alternative notation for $\di\frac{\pa f}{\pa x_i}(\mf{x}_0)$ is $f_{x_i}(\mf{x}_0)$.
\end{highlight}

\begin{figure}[ht]
\centering
\includegraphics[scale=0.2]{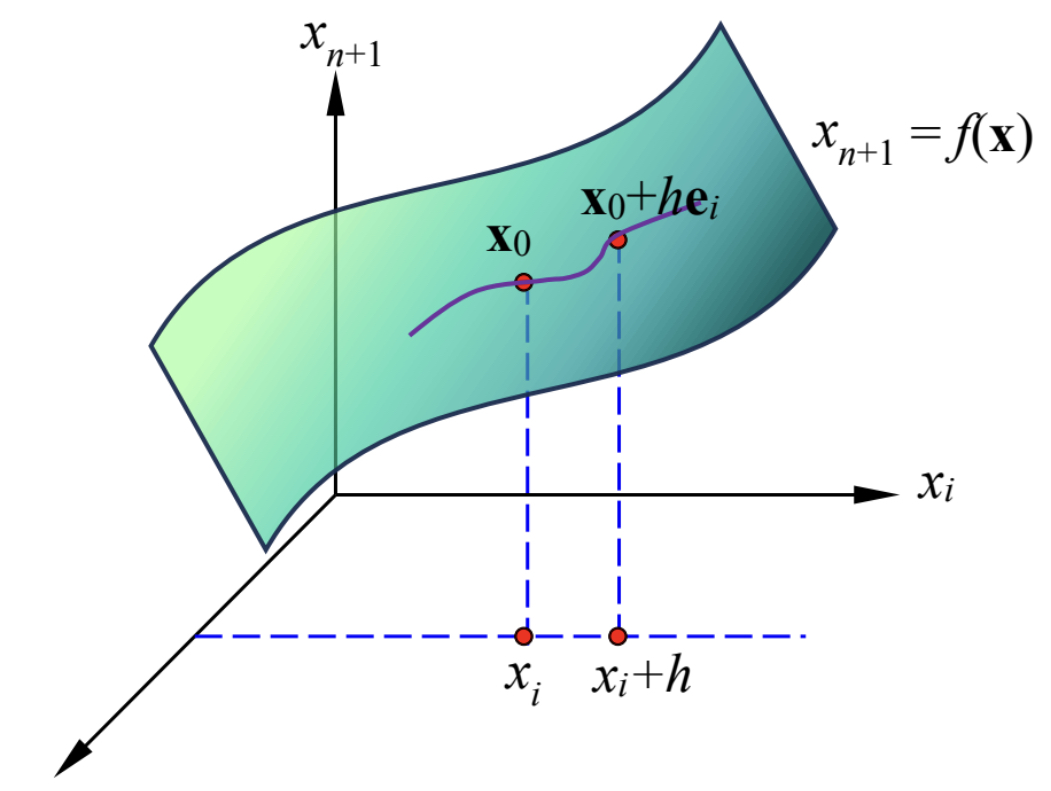}

\caption{Partial derivative.}\label{figure41}
\end{figure}

\begin{remark}{Partial Derivatives}
Let $\mf{x}_0=(a_1, a_2, \ldots, a_n)$ and
define the function $g:(-r, r)\to \mb{R}$ by
\[g(h)=f(\mathbf{x}_0+h\mf{e}_i)=f(a_1, \ldots, a_{i-1}, a_i+h, a_{i+1}, \ldots, a_n).\]
Then
\[\lim_{h\to 0}\frac{f(\mf{x}_0+h\mf{e}_i)-f(\mf{x}_0)}{h}=\lim_{h\to 0}\frac{g(h)-g(0)}{h}=g'(0).\]
Thus, 
$\di f_{x_i}(\mf{x}_0)$ exists if and only if $g(h)$ is differentiable at $h=0$. Moreover, 
to find $\di  f_{x_i}(\mf{x}_0)$, we regard the variables $x_1, \ldots, x_{i-1}, x_{i+1}, \ldots, x_n$ as constants, and differentiate with respect to $x_i$. Hence, the derivative rules such as sum rule, product rule and quotient rule   still work for partial derivatives, as long as one is clear which variable to take derivative, which variable to be regarded as constant.

\end{remark}

\begin{example}{}
Let $f:\mathbb{R}^2\to\mb{R}$ be the function defined as $f(x,y)=x^2y$.
Find $\di f_x(1,2)$ and $\di f_y(1,2)$.
\end{example}
\begin{solution}{Solution}
\[\frac{\pa f}{\pa x}=2xy,\hspace{1cm}\frac{\pa f}{\pa y}=x^2.\]
Therefore,
\[f_x(1,2)=4, \hspace{1cm}f_y(1,2)=1.\]
\end{solution}

\begin{example}{}
Let $f:\mathbb{R}^2\to\mb{R}$ be the function defined as $f(x,y)=|x+y|$.
Determine whether $\di f_x(0,0)$ exists.
\end{example}
\begin{solution}{Solution}
By definition, $f_x(0,0)$ is given by the limit
\[\lim_{h\to 0}\frac{f(h,0)-f(0,0)}{h}\]
if it exists. Since\[ \lim_{h\to 0}\frac{f(h,0)-f(0,0)}{h}=\lim_{h\to 0}\frac{|h|}{h}, \]
and
\[\lim_{h\to 0^-}\frac{|h|}{h}=-1\quad \text{and}\quad \lim_{h\to 0^+}\frac{|h|}{h}=1,\]
the limit \[ \lim_{h\to 0}\frac{f(h,0)-f(0,0)}{h}\] does not exist. Hence, $f_x(0,0)$ does not exist.
\end{solution}

\begin{definition}{}
Let $\mathcal{O}$ be an open subset of $\mb{R}^n$, and let $f:\mathcal{O}\to \mb{R}$ be a function defined on $\mathcal{O}$. 
If the function $f:\mathcal{O}\to \mb{R}$ has   partial derivative with respect to $x_i$ at every point   of $\mathcal{O}$, this defines the function $\di f_{x_i}:\mathcal{O}\to \mb{R}$. In this case, we say that the partial derivative of $f$ with respect to $x_i$ exists.  

If $\di f_{x_i}:\mathcal{O}\to \mb{R}$ exists for all $1\leq i\leq n$, we say that the function $f:\mathcal{O}\to \mb{R}$ has partial derivatives.
\end{definition}

\begin{example}[label=230729_1]{}
Find the partial derivatives of the function $f:\mb{R}^3\to \mb{R}$ defined as
\[f(x,y,z)=\sin (xy+z)+\frac{3x}{y^2+z^2+1}.\]
\end{example}
\begin{solution}{Solution}
\begin{align*}
\frac{\pa f}{\pa x}(x,y,z)&=y\cos(xy+z)+\frac{3}{y^2+z^2+1},
\\
\frac{\pa f}{\pa y}(x,y,z)&=x\cos(xy+z)-\frac{6xy}{(y^2+z^2+1)^2},
\\
\frac{\pa f}{\pa z}(x,y,z)&= \cos(xy+z)-\frac{6xz}{(y^2+z^2+1)^2}.\hspace{4cm}
\end{align*}
\end{solution}

For a function defined on an open subset of $\mb{R}^n$, there are $n$ partial derivatives with respect to the $n$ directions defined by the coordinate axes. These   define a vector in $\mb{R}^n$.
\begin{definition}{Gradient}
Let $\mathcal{O}$ be an open subset of $\mathbb{R}^n$, and let $\mf{x}_0$ be a point in $\mathcal{O}$.
If the function $f:\mathcal{O}\to\mb{R}$  has partial derivatives at $\mf{x}_0$, we define the gradient of the function $f$ at $\mf{x}_0$ as the vector in $\mb{R}^n$ given by
\[\nabla f(\mathbf{x}_0)=\left(\frac{\pa f}{\pa x_1}(\mf{x}_0), \frac{\pa f}{\pa x_2}(\mf{x}_0), \cdots, \frac{\pa f}{\pa x_n}(\mf{x}_0)\right).\]
 
\end{definition} 
Let us revisit Example \ref{230729_1}.
\begin{example}{}
The gradient of the function  $f:\mb{R}^3\to \mb{R}$ defined as
\[f(x,y,z)=\sin (xy+z)+\frac{3x}{y^2+z^2+1} \] in Example \ref{230729_1} is the function $\nabla f:\mb{R}^3\to\mb{R}^3$, 
\begin{align*}\nabla f(x,y,z)&=\begin{bmatrix} \di y\cos(xy+z)+\frac{3}{y^2+z^2+1}\\\di x\cos(xy+z)-\frac{6xy}{(y^2+z^2+1)^2}\\\di  \cos(xy+z)-\frac{6xz}{(y^2+z^2+1)^2}\end{bmatrix}.\end{align*}
In particular,
\[\nabla f(1,-1,1)=\left(0, \frac{5}{3}, \frac{1}{3}\right).\]
\end{example} 

It is straightforward to extend the definition of partial derivative to a function  $\mf{F}:\mathcal{O}\to \mb{R}^m$ whose codomain is $\mb{R}^m$ with $m\geq 2$.
\begin{definition}{}
Let $\mathcal{O}$ be an open subset of $\mb{R}^n$, and let $\mf{F}:\mathcal{O}\to \mb{R}^m$ be a function defined on $\mathcal{O}$. Given $\mf{x}_0$ in $\mathcal{O}$ and $1\leq i\leq n$, we say that $\mf{F}:\mathcal{O}\to\mb{R}^m$ has partial derivative with respect to $x_i$ at the point $\mf{x}_0$ if the limit 
\[\frac{\pa\mf{F}}{\pa x_i}(\mf{x}_0)=\lim_{h\to 0}\frac{\mf{F}(\mf{x}_0+h\mf{e}_i)-\mf{F}(\mf{x}_0)}{h}\] exists.
We say that  $\mf{F}:\mathcal{O}\to\mb{R}^m$ has partial derivative  at the point $\mf{x}_0$ if $\di \frac{\pa\mf{F}}{\pa x_i}(\mf{x}_0)$ exists for each $1\leq i\leq n$. We say that $\mf{F}:\mathcal{O}\to \mb{R}^m$ has partial derivative if it has partial derivative at each point of $\mathcal{O}$.
\end{definition}
Since the limit of a function $\mf{G}:(-r,r)\to\mb{R}^m$ when $h\to 0$ exists if and only if the limit of each component function $G_j:(-r,r)\to \mb{R}$, $1\leq j\leq m$ when $h\to 0$ exists, we have the following.
\begin{proposition}{}
Let $\mathcal{O}$ be an open subset of $\mb{R}^n$, and let $\mf{F}:\mathcal{O}\to \mb{R}^m$ be a function defined on $\mathcal{O}$. Given $\mf{x}_0$ in $\mathcal{O}$ and $1\leq i\leq n$,  $\mf{F}:\mathcal{O}\to\mb{R}^m$ has partial derivative with respect to $x_i$ at the point $\mf{x}_0$ if and only if
 if each component function $F_j:\mathcal{O}\to \mb{R}$, $1\leq j\leq m$ has partial derivative with respect to $x_i$ at the point $\mf{x}_0$. In this case, we have
\[\frac{\pa\mf{F}}{\pa x_i}(\mf{x}_0)=\left(\frac{\pa F_1}{\pa x_i}(\mf{x}_0), \ldots, \frac{\pa F_m}{\pa x_i}(\mf{x}_0)\right).\]
\end{proposition}

To capture all the partial derivatives, we define a derivative matrix.
\begin{definition}{The Derivative Matrix}
Let $\mathcal{O}$ be an open subset of $\mb{R}^n$ that contains the point $\mf{x}_0$, and let $\mf{F}:\mathcal{O}\to \mb{R}^m$ be a function defined on $\mathcal{O}$.   If $\mf{F}:\mathcal{O}\to\mb{R}^m$ has partial derivative  at the point $\mf{x}_0$, the derivative matrix of $\mf{F}:\mathcal{O}\to\mb{R}^m$ at $\mf{x}_0$ is the $m\times n$ matrix
\[\mf{DF}(\mf{x}_0)=\begin{bmatrix} \nabla F_1(\mf{x}_0)\\\nabla F_2(\mf{x}_0)\\ \vdots \\\nabla F_m(\mf{x}_0)\end{bmatrix}
=\begin{bmatrix}\vspace{0.3cm} \di \frac{\pa F_1}{\pa x_1}(\mf{x}_0) & \di \frac{\pa F_1}{\pa x_2}(\mf{x}_0) &\cdots & \di \frac{\pa F_1}{\pa x_n}(\mf{x}_0)\\
\vspace{0.3cm}
\di \frac{\pa F_2}{\pa x_1}(\mf{x}_0)& \di \frac{\pa F_2}{\pa x_2}(\mf{x}_0) &\cdots & \di \frac{\pa F_2}{\pa x_n}(\mf{x}_0)\\
\vspace{0.3cm}\vdots &\vdots &\ddots & \vdots \\
\di \frac{\pa F_m}{\pa x_1}(\mf{x}_0) & \di \frac{\pa F_m}{\pa x_2}(\mf{x}_0) &\cdots & \di \frac{\pa F_m}{\pa x_n}(\mf{x}_0)
\end{bmatrix}.\]
\end{definition} When $m=1$, the derivative matrix is just the gradient of the function as a row matrix.

\begin{example}{}
Let $\mf{F}:\mb{R}^3\to\mb{R}^2$ be the function defined as
\[\mf{F}(x,y,z)=\left(xy^2z^3, x+3y-7z\right).\]
Find the derivative matrix of $\mf{F}$ at the point $(1, -1, 2)$.
\end{example}
\begin{solution}{Solution}
\[
\mf{DF}(x,y,z)=\begin{bmatrix} y^2z^3 & 2xyz^3 & 3xy^2z^2\\
1 & 3 & -7\end{bmatrix}.\]
Thus,  the derivative matrix of $\mf{F}$ at the point $(1, -1, 2)$ is
\[\mf{DF}(1, -1, 2)=\begin{bmatrix} 8 & -16 & 12\\
1 & 3 & -7\end{bmatrix}.\]
\end{solution}

Since the partial derivatives of a function is defined componentwise, we can focus on functions $f:\mathcal{O}\to \mb{R}$ whose codomain is $\mb{R}$. One might wonder why we have not mentioned the word ''differentiable'' so far. For single variable functions, we have seen in volume I that if a function is differentiable at a point, then it is continuous  at that point. For multivariable functions, the existence of partial derivatives is not enough to guarantee continuity, as is shown in the next example.

\begin{example}[label=230729_2]{}
Let $f:\mb{R}^2\to\mb{R}$ be the function defined as
\[f(x,y)=\begin{cases}\di \frac{xy}{x^2+y^2},\quad &\text{if}\;(x,y)\neq (0,0),\\
0,\quad &\text{if}\;(x,y)=(0,0).\end{cases}\]
Show that $f$ is not continuous at $(0,0)$, but it has partial derivatives at $(0,0)$.
\end{example}
\begin{figure}[ht]
\centering
\includegraphics[scale=0.18]{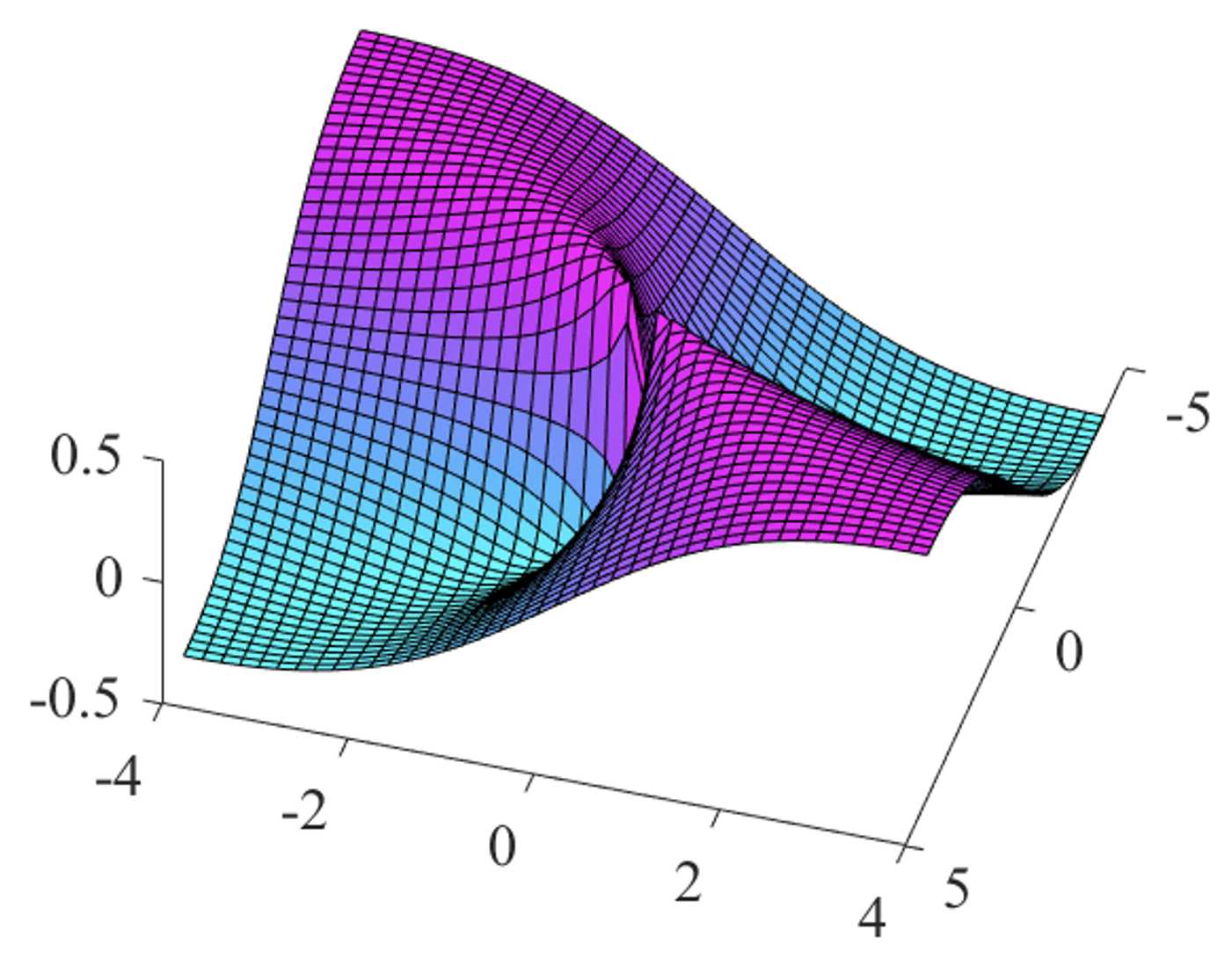}

\caption{The function $\di f(x,y)$ defined in Example \ref{230729_2}.}\label{figure46}
\end{figure}
\begin{solution}{Solution}
Consider the sequence $\{\mf{u}_k\}$ with
\[\mf{u}_k=\left(\frac{1}{k}, \frac{1}{k}\right).\]
It is a sequence in $\mb{R}^2$ that converges to $(0,0)$. Since
\[f(\mf{u}_k)=\frac{1}{2}\hspace{1cm}\text{for all}\;k\in \mb{Z}^+,\]
\bs
the sequence $\{f(\mf{u}_k)\}$ converges to $1/2$. But $f(0,0)=0\neq 1/2$. Since there is a sequence $\{\mf{u}_k\}$ that converges to $(0,0)$, but the sequence $\{f(\mf{u}_k)\}$ does not converge to $f(0,0)$, $f$ is not continuous at $(0,0)$.

To find partial derivatives at $(0,0)$, we use definitions.
\begin{align*}
f_x(0,0)&=\lim_{h\to 0}\frac{f(h,0)-f(0,0)}{h}=\lim_{h\to 0}\frac{0-0}{h}=0,\\
f_y(0,0)&=\lim_{h\to 0}\frac{f(0,h)-f(0,0)}{h}=\lim_{h\to 0}\frac{0-0}{h}=0.
\end{align*}These show that $f$ has partial derivatives at $(0,0)$, and $f_x(0,0)=f_y(0,0)=0$.
\end{solution}
For the function defined in Example  \ref{230729_2}, it has partial derivatives at all points. In fact, when $(x,y)\neq (0,0)$, we can apply derivative rules directly and find that
\[\frac{\pa f}{\pa x}(x,y)=\frac{(x^2+y^2)y-2x^2y}{(x^2+y^2)^2}=\frac{y(y^2-x^2)}{(x^2+y^2)^2}.\]Similarly,
\[\frac{\pa f}{\pa y}(x,y)=\frac{x(x^2-y^2)}{(x^2+y^2)^2}.\]

Let us highlight again our conclusion.
\begin{highlight}{Partial Derivative vs Continuity}
The existence of partial derivatives does not imply continuity.
\end{highlight}
This prompts us to find a better definition of {\it differentiability}, which can imply continuity. This will be considered in a latter section.

When the function $f:\mathcal{O}\to \mb{R}$ has partial derivative with respect to $x_i$, we obtain the function $f_{x_i}:\mathcal{O}\to \mb{R}$. Then we can discuss whether the function $f_{x_i}$ has partial derivative at a point in $\mathcal{O}$.

\begin{definition}{Second Order Partial Derivatives}

Let $\mathcal{O}$ be an open subset of $\mathbb{R}^n$ that contains the point $\mf{x}_0$, and let $f:\mathcal{O}\to\mathbb{R}$ be a function defined on $\mathcal{O}$. Given that $1\leq i\leq n$, $1\leq j\leq n$, we say that the second order partial derivative $\di \frac{\pa^2 f}{\pa x_j\pa x_i}$ exists at $\mf{x}_0$ provided that there exists an open ball $B(\mf{x}_0, r)$ that is contained in $\mathcal{O}$ such that $\di \frac{\pa f}{\pa x_i}: B(\mf{x}_0,r)\to \mb{R}$ exists, and it has partial derivative with respect to $x_j$ at the point $\mf{x}_0$. In this case, we define the second order partial derivative  $\di \frac{\pa^2 f}{\pa x_j\pa x_i}(\mf{x}_0)$ of $f$ at $\mf{x}_0$ as
\[\frac{\pa^2 f}{\pa x_j\pa x_i}(\mf{x}_0)=\frac{\pa f_{x_i}}{\pa x_j}(\mf{x}_0)=\lim_{h\to 0}\frac{f_{x_i}(\mf{x}_0+h\mf{e}_j)-f_{x_i}(\mf{x}_0)}{h}.\]

We say that the function $f:\mathcal{O}\to\mathbb{R}$ has   second order partial derivatives at $\mf{x}_0$ provided that $\di \frac{\pa^2 f}{\pa x_j\pa x_i}(\mf{x}_0)$ exists for all $1\leq i\leq n, 1\leq j\leq n$.
\end{definition}

In the same way, one can also define second order partial derivatives  for a function $\mf{F}:\mathcal{O}\to \mb{R}^m$ with codomain $\mb{R}^m$ when $m\geq 2$.

\begin{remark}{}
In the definition of the second order partial derivative $\di\frac{\pa^2 f}{\pa x_j\pa x_i}(\mf{x}_0)$, instead of assuming $f_{x_i}(\mf{x})$ exists for all $\mf{x}$ in a ball of radius $r$ centered at $\mf{x}_0$, 
it is sufficient to assume that there exists $r>0$ such that $f_{x_i}(\mf{x}_0+h\mf{e}_j)$ exists for all $|h|<r$.
\end{remark}

\begin{definition}{}
Given $1\leq i\leq n, 1\leq j\leq n$,
we say that the function $f:\mathcal{O}\to\mathbb{R}$ has the second order partial derivative $\di \frac{\pa^2f}{\pa x_j\pa x_i}$ provided that $\di \frac{\pa^2f}{\pa x_j\pa x_i}(\mf{x}_0)$ exists for all $\mf{x}_0$ in $\mathcal{O}$. 

We say that the function $f:\mathcal{O}\to\mathbb{R}$ has   second order partial derivatives provided that $\di \frac{\pa^2 f}{\pa x_j\pa x_i}$ exists for all $1\leq i\leq n, 1\leq j\leq n$.
\end{definition}
\begin{highlight}{Notations for Second Order Partial Derivatives}
Alternative notations for second order partial derivatives are
\[\frac{\pa^2f}{\pa x_j\pa x_i}=(f_{x_i})_{x_j}=f_{x_ix_j}.\]
Notice that the orders of $x_i$ and $x_j$ are different in different notations.
\end{highlight}
\begin{remark}{}
Given $1\leq i\leq n, 1\leq j\leq n$, the function $f:\mathcal{O}\to\mathbb{R}$ has the second order partial derivative $\di \frac{\pa^2f}{\pa x_j\pa x_i}$ provided that $f_{x_i}:\mathcal{O}\to\mb{R}$ exists, 
  and $f_{x_i}$ has partial derivative with respect to $x_j$.

\end{remark}

\begin{example}[label=230729_3]{}
Find the second order partial derivatives of the function 
$\di f:\mb{R}^2\to \mb{R}$ defined as
\[ f(x,y)=xe^{2x+3y}.\]

\end{example}
\begin{solution}{Solution}
We find the first order partial derivatives first.
\begin{align*}
\frac{\pa f}{\pa x}(x,y)&=e^{2x+3y}+2xe^{2x+3y}=(1+2x)e^{2x+3y},\\
\frac{\pa f}{\pa y}(x,y)&=3xe^{2x+3y}.\end{align*}
Then we compute the second order partial derivatives.
\begin{align*}
\frac{\pa^2 f}{\pa x^2}(x,y)&=2e^{2x+3y}+2(1+2x)e^{2x+3y}=(4+4x)e^{2x+3y},\\
\frac{\pa^2 f}{\pa y\pa x}(x,y)&=3(1+2x)e^{2x+3y}=(3+6x)e^{2x+3y},\\
\frac{\pa^2 f}{\pa x\pa y}(x,y)&=3e^{2x+3y}+6xe^{2x+3y}=(3+6x)e^{2x+3y},\\
\frac{\pa^2 f}{\pa y^2}(x,y)&=9xe^{2x+3y}.
\end{align*}
\end{solution}

\begin{definition}{The Hessian Matrix}
Let $\mathcal{O}$ be an open subset of $\mb{R}^n$ that contains the point $\mf{x}_0$. If $f:\mathcal{O}\to \mb{R}$ is a function that has second order partial derivatives at $\mf{x}_0$, the Hessian matrix of $f$ at $\mf{x}_0$ is the $n\times n$ matrix defined as
\begin{align*}H_f(\mf{x}_0)&=\begin{bmatrix}\di\frac{\pa^2f}{\pa x_i\pa x_j}(\mf{x}_0)\end{bmatrix}\\&=\begin{bmatrix}\vspace{0.3cm} \di \frac{\pa^2 f}{\pa x_1^2}(\mf{x}_0) &\di \frac{\pa^2 f}{\pa x_1\pa x_2}(\mf{x}_0) & \cdots &\di  \frac{\pa^2 f}{\pa x_1\pa x_n}(\mf{x}_0)\\
\vspace{0.3cm} \di \frac{\pa^2 f}{\pa x_2\pa x_1}(\mf{x}_0) &\di \frac{\pa^2 f}{ \pa x_2^2}(\mf{x}_0) & \cdots &\di  \frac{\pa^2 f}{\pa x_2\pa x_n}(\mf{x}_0)\\
\vdots & \vdots & \ddots& \vdots\\
\vspace{0.3cm} \di \frac{\pa^2 f}{\pa x_n\pa x_1}(\mf{x}_0) &\di \frac{\pa^2 f}{\pa x_n\pa x_2}(\mf{x}_0) & \cdots &\di  \frac{\pa^2 f}{ \pa x_n^2}(\mf{x}_0)\end{bmatrix}.\end{align*}
\end{definition}
We do not define Hessian matrix for a function $\mf{F}:\mathcal{O}\to\mb{R}^m$ with codomain $\mb{R}^m$ when $m\geq 2$.

\begin{example}{}
For the function $\di f:\mb{R}^2\to \mb{R}$ defined as
$\di f(x,y)=xe^{2x+3y}$ in Example \ref{230729_3},
\[H_f(x,y)=\begin{bmatrix}(4+4x)e^{2x+3y} & (3+6x)e^{2x+3y}\\
(3+6x)e^{2x+3y} & 9xe^{2x+3y} \end{bmatrix}.\]
\end{example}
In Example \ref{230729_3}, we notice that 
\[\frac{\pa^2 f}{\pa y\pa x}(x,y)=\frac{\pa^2 f}{\pa x\pa y}(x,y)\]for all $(x,y)\in \mb{R}^2$. The following example  shows that {\it this is not always true}.

\begin{example}[label=230729_4]{}
Consider the function $f:\mb{R}^2\to \mb{R}$ defined as
\[f(x,y)=\begin{cases}\di \frac{xy(x^2-y^2)}{x^2+y^2},\quad &\text{if}\; (x,y)\neq (0,0),\\
0,\quad &\text{if}\;(x,y)=(0,0).\end{cases}\]
Find $f_{xy}(0,0)$ and $f_{yx}(0,0)$.
\end{example}

\begin{figure}[ht]
\centering
\includegraphics[scale=0.18]{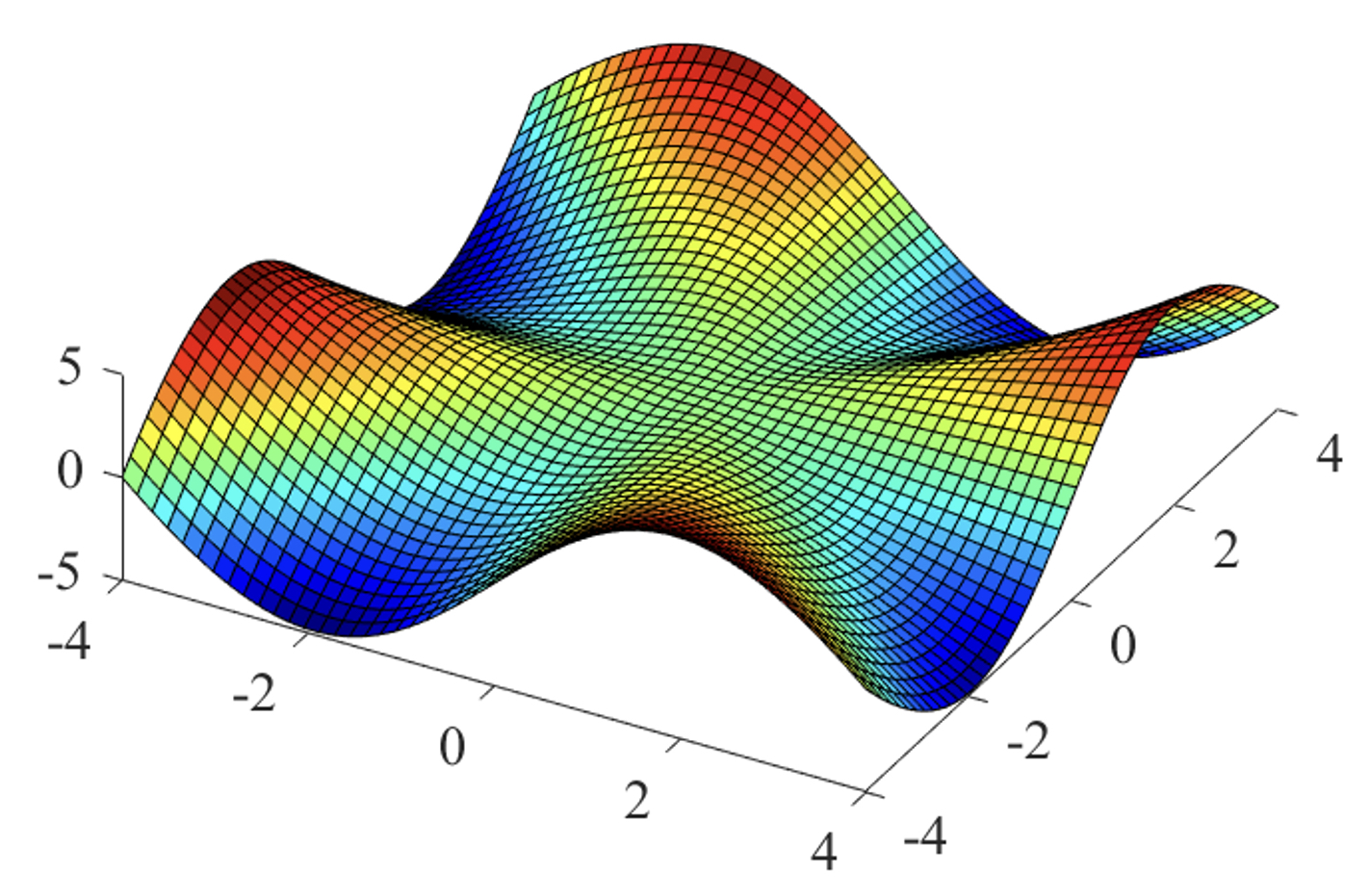}

\caption{The function $\di f(x,y)$ defined in Example \ref{230729_4}.}\label{figure57}
\end{figure}
\begin{solution}{Solution}
To compute $f_{xy}(0,0)$, we need to compute $f_x(0, h)$ for all $h$ in a neighbourhood of 0. To compute $f_{yx}(0,0)$, we need to compute $f_y(h,0)$ for all $h$ in a neighbourhood of $0$. Notice that for any $h\in \mb{R}$, $f(0,h)=f(h,0)=0$.  By considering $h=0$ and $h\neq 0$ separately, we find that
\begin{align*}
f_x(0, h)&=\lim_{t\to 0}\frac{f(t, h)-f(0,h)}{t}=\lim_{t\to 0} \frac{h(t^2-h^2)  }{t^2+h^2}= -h,  \\
f_y(h,0)&=\lim_{t\to 0}\frac{f(h,t)-f(h,0)}{t}=\lim_{t\to 0}\frac{h(h^2-t^2)}{h^2+t^2}=h.
\end{align*}
It follows that
\begin{align*}
f_{xy}(0,0)&=\lim_{h\to 0}\frac{f_x(0,h)-f_x(0,0)}{h}=\lim_{h\to 0}\frac{-h}{h}=-1,\\
f_{yx}(0,0)&=\lim_{h\to 0}\frac{f_y(h,0)-f_y(0,0)}{h}=\lim_{h\to 0}\frac{h}{h}=1.
\end{align*}
\end{solution}

\begin{highlight}{}
Example \ref{230729_4} shows that there exists a function $f:\mathbb{R}^2\to \mb{R}$ which has second order partial derivatives at $(0,0)$ but
\[\frac{\pa^2 f}{\pa  x\pa y}(0,0)\neq \frac{\pa^2 f}{\pa  y\pa x}(0,0).\]
\end{highlight}

\begin{remark}[label=230729_6]{}
If $\mathcal{O}$ is an open subset of $\mb{R}^n$ that contains the point $\mf{x}_0$, there exists $r>0$ such that $B(\mf{x}_0, r)\subset \mathcal{O}$. Given that $f:\mathcal{O}\to \mb{R}$ is a function defined on $\mathcal{O}$, and $1\leq i<j\leq n$, let $\mk{D}$ be the ball with center at $(0,0)$ and radius $r$ in $\mb{R}^2$.  Define the function $g:\mk{D}\to \mb{R}$ by
\[g(u, v)=f(\mf{x}_0+u\mf{e}_i+v\mf{e}_j).\]
Then $\di \frac{\pa^2 f}{\pa x_j\pa x_i}(\mf{x}_0)$ exists if and only if $\di \frac{\pa^2 g}{\pa v\pa u}(0,0)$ exists. In such case, we have
\[ \frac{\pa^2 f}{\pa x_j\pa x_i}(\mf{x}_0)=\frac{\pa^2 g}{\pa v\pa u}(0,0).\]
\end{remark}

The following gives a sufficient condition to interchange the order of taking partial derivatives. 

\begin{theorem}{Clairaut's Theorem or Schwarz's Theorem}
Let $\mathcal{O}$ be an open subset of $\mb{R}^n$ that contains the point $\mathbf{x}_0$, and let $f:\mathcal{O}\to\mb{R}$ be a function defined on $\mathcal{O}$. Assume that $1\leq i<j\leq n$, and the second order partial derivatives  $\di \frac{\pa^2 f}{\pa x_j\pa x_i}:\mathcal{O}\to \mb{R}$ and $\di \frac{\pa^2 f}{\pa x_i\pa x_j}:\mathcal{O}\to \mb{R}$ exist. If the functions $\di \frac{\pa^2 f}{\pa x_j\pa x_i}$ and $\di \frac{\pa^2 f}{\pa x_i\pa x_j}:\mathcal{O}\to \mb{R}$ are continuous at $\mf{x}_0$, then
\[\frac{\pa^2 f}{\pa x_j\pa x_i}(\mf{x}_0)=\frac{\pa^2 f}{\pa x_i\pa x_j}(\mf{x}_0).\]
\end{theorem}
\begin{myproof}{Proof}
Since $\mathcal{O}$ is an open set that contains the point $\mf{x}_0$, there exists $r>0$ such that $B(\mf{x}_0, r)\subset \mathcal{O}$. Let
\[\mk{D}=\left\{(u,v)|\, u^2+v^2<r^2\right\},\]
and define the function $g:\mk{D}\to \mb{R}$ by
\[g(u, v)=f(\mf{x}_0+u\mf{e}_i+v\mf{e}_j).\]

By Remark \ref{230729_6}, 
$g$ has second order partial derivatives, and $\di \frac{\pa^2 g}{\pa v\pa u}$ and $\di \frac{\pa^2 g}{\pa u\pa v}$ are continuous  at $(0,0)$. We need to show that
\[\frac{\pa^2 g}{\pa v\pa u}(0,0)=\frac{\pa^2 g}{\pa u\pa v}(0,0).\]
Consider the function
\[G(u,v)=g(u,v)-g(u,0)-g(0,v)+g(0,0).\] 
 
Notice that
\[G(u,v)=H_v(u)-H_v(0)=S_u(v)-S_u(0),\]

\bp
where

\[H_v(u)=g(u,v)-g(u,0),\hspace{1cm} S_u(v)=g(u,v)-g(0,v).\]

For fixed $v$ with $|v|<r$, the function $H_v(u)$ is defined for those $u$ with $|u|<\sqrt{r^2-v^2}$, such that $(u,v)$ is in $\mk{D}$. It is differentiable with 
\[H_v'(u)=\frac{\pa g}{\pa u}(u,v)-\frac{\pa g}{\pa u}(u,0).\]
Hence, if $(u,v)$  is in $\mk{D}$,
 mean value theorem for single variable functions implies that there exists $c_{u,v}\in (0,1)$ such that
\begin{align*}
G(u,v)&=H_v(u)-H_v(0)\\&=uH_v'(c_{u,v}u)\\
&=u\left(\frac{\pa g}{\pa u}(c_{u,v}u,v)-\frac{\pa g}{\pa u}(c_{u,v}u,0)\right).\end{align*}
 
Regard this now as a function of $v$, the mean value theorem for single variable functions implies that there exists $d_{u,v}\in (0,1)$ such that
\begin{equation}\label{230729_7}G(u,v)=uv \frac{\pa^2g}{\pa v\pa u}(c_{u,v}u, d_{u,v}v).\end{equation}

Using the same reasoning, we find that for $(u,v)\in\mk{D}$, there exists
$\widetilde{d}_{u,v}\in (0,1)$ such that 
\[G(u,v)= vS_u'(\widetilde{d}_{u,v}v)=v\left(\frac{\pa g}{\pa v}(u,\widetilde{d}_{u,v}v)-\frac{\pa g}{\pa v}(0, \widetilde{d}_{u,v}v)\right).\]
Regard this as a function of $u$, mean value theorem implies that there exists $\widetilde{c}_{u,v}\in (0,1)$ such that
\begin{equation}\label{230729_8}G(u,v)=uv \frac{\pa^2g}{\pa u\pa v}(\widetilde{c}_{u,v}u, \widetilde{d}_{u,v}v).\end{equation}

Comparing \eqref{230729_7} and \eqref{230729_8}, we find that
\[\frac{\pa^2g}{\pa v\pa u}(c_{u,v}u, d_{u,v}v)=\frac{\pa^2g}{\pa u\pa v}(\widetilde{c}_{u,v}u, \widetilde{d}_{u,v}v).\]
 \bp
When $(u,v)\to (0,0)$, $(c_{u,v}u, d_{u,v}v)\to (0,0)$ and $(\widetilde{c}_{u,v}u, \widetilde{d}_{u,v}v)\to (0,0)$. The continuities of $g_{uv}$ and $g_{vu}$ at $(0,0)$ then imply that 
\[\frac{\pa^2 g}{\pa v\pa u}(0,0)=\frac{\pa^2 g}{\pa u\pa v}(0,0).\]
 
This completes the proof.
\end{myproof}

\begin{example}{}
Consider the function $f:\mb{R}^2\to \mb{R}$ in Example \ref{230729_4} defined as
\[f(x,y)=\begin{cases}\di \frac{xy(x^2-y^2)}{x^2+y^2},\quad &\text{if}\; (x,y)\neq (0,0),\\
0,\quad &\text{if}\;(x,y)=(0,0).\end{cases}\]

  When $(x,y)\neq (0,0)$, we find that
\begin{align*}
\frac{\pa f}{\pa x}(x,y)&= \frac{y(x^4+4x^2y^2-y^4)}{(x^2+y^2)^2},\\
\frac{\pa f}{\pa y}(x,y)&= \frac{x(x^4-4x^2y^2-y^4)}{(x^2+y^2)^2}.
\end{align*}

It follows that
\[
\frac{\pa^2 f}{\pa y\pa x}(x,y) = \frac{x^6+9x^4y^2-9x^2y^4-y^6   }{(x^2+y^2)^3} =\frac{\pa^2 f}{\pa x\pa y}(x,y).
\]
Indeed, both $f_{xy}$ and $f_{yx}$ are continuous on $\mathbb{R}^2\setminus\{(0,0)\}$.
\end{example}

\begin{corollary}{}
Let $\mathcal{O}$ be an open subset of $\mb{R}^n$ that contains the point $\mathbf{x}_0$, and let $f:\mathcal{O}\to\mb{R}$ be a function defined on $\mathcal{O}$. If all the second order partial derivatives of the function $f:\mathcal{O}\to\mb{R}$ at $\mathbf{x}_0$ are continuous, then the Hessian matrix $H_f(\mf{x}_0)$ of $f$ at $\mf{x}_0$ is a symmetric matrix.
\end{corollary}

\begin{remark}{}
One can define partial derivatives of higher orders following the same rationale as we define the second order partial derivatives. Extension of 
Clairaut's theorem  to higher order partial derivatives is straightforward. The key point is the continuity of the partial derivatives involved.
\end{remark}

\vp
\noindent
{\bf \large Exercises  \thesection}
\setcounter{myquestion}{1}

\begin{question}{\themyquestion}
Let $f:\mb{R}^3\to \mb{R}$ be the function defined as 
\[f(x,y,z)=\frac{x z}{e^y+1}.\]
Find $\nabla f(1, 0, -1)$, the gradient of $f$ at the point $(1,0,-1)$.
\end{question}
\atc

\begin{question}{\themyquestion}
Let $\mf{F}:\mb{R}^2\to \mb{R}^3$ be the function defined as
\[\mf{F}(x,y)= \left( x^2y, xy^2, 3x^2+4y^2\right).\]
Find $\mf{DF}(2, -1)$, the derivative matrix of $\mf{F}$ at the point $(2, -1)$.
\end{question}
\atc

\begin{question}{\themyquestion}
Let $f:\mb{R}^3\to \mb{R}$ be the function defined as 
\[f(x,y,z)=x^2+3xyz+2y^2z^3.\]
Find $H_f(1, -1, 2)$, the Hessian matrix of $f$ at the point $(1,-1, 2)$.
\end{question}
\atc

\begin{question}{\themyquestion}
Let $f:\mb{R}^2\to\mb{R}$ be the function defined as
\[f(x,y)=\begin{cases}\di \frac{3xy}{x^2+4y^2},\quad &\text{if}\;(x,y)\neq (0,0),\\
0,\quad &\text{if}\;(x,y)=(0,0).\end{cases}\]
Show that $f$ is not continuous at $(0,0)$, but it has partial derivatives at $(0,0)$.
\end{question}
\atc
\begin{question}{\themyquestion}
Let $f:\mathbb{R}^2\to\mb{R}$ be the function defined as $f(x,y)=|x^2+y|$.
Determine whether $\di f_y(1,-1)$ exists.
\end{question}
\atc

\begin{question}{\themyquestion}
Let $f:\mb{R}^2\to\mb{R}$ be the function defined as
\[f(x,y)=\begin{cases}\di \frac{x^2y}{x^2+y^2},\quad &\text{if}\;(x,y)\neq (0,0),\\
0,\quad &\text{if}\;(x,y)=(0,0).\end{cases}\]
Show that $f$ is   continuous, it has   partial derivatives, but the partial derivatives   are not continuous.
\end{question}
\atc

\begin{question}{\themyquestion}
Consider the function $f:\mb{R}^2\to \mb{R}$ defined as
\[f(x,y)=\begin{cases}\di \frac{xy(x^2+9y^2)}{4x^2+y^2},\quad &\text{if}\; (x,y)\neq (0,0),\\
0,\quad &\text{if}\;(x,y)=(0,0).\end{cases}\]
Find the Hessian matrix $H_f(0,0)$ of $f$ at $(0,0)$. 
\end{question}

\section{Differentiability and First Order Approximation} 
Let $\mathcal{O}$ be an open subset of $\mb{R}^n$ that contains the point $\mf{x}_0$, and let $\mf{F}:\mathcal{O}\to \mb{R}^m$ be a function defined on $\mathcal{O}$. 
As we have seen in the previous section, even if $\mf{F}$ has partial derivatives at $\mf{x}_0$, it does not imply that $\mf{F}$ is continuous at $\mf{x}_0$. Heuristically, this is because the partial derivatives only consider the change of the function  along the $n$ directions defined by the coordinate axes, while continuity of $\mf{F}$ requires us to consider the change of  $\mf{F}$ along {\it all} directions.

\subsection{Differentiability}
In this section, we will give a suitable definition of {\it differentiability} to ensure that we can capture the change of $\mf{F}$ in all   directions. 
Let us first revisit an alternative perpective of {\it differentiability} for a single variable function $f:(a,b)\to \mb{R}$, which we have discussed in volume I.
If $x_0$ is a point in $(a,b)$, then the function $f:(a,b)\to \mb{R}$ is differentiable at $x_0$ if and only if there is a number $c$ such that
\begin{equation}\label{230730_1}\lim_{h\to 0}\frac{f(x_0+h)-f(x_0)-ch}{h}=0.\end{equation}
In fact, if $f$ is differentiable at $x_0$, then this number $c$ has to equal to $f'(x_0)$. 

Now for a function $\mf{F}:\mathcal{O}\to \mb{R}^m$ defined on an open subset $\mathcal{O}$ of $\mb{R}^n$, to consider the differentiability of $\mf{F}$ at $\mf{x}_0\in\mathcal{O}$, we should compare $\mf{F}(\mf{x}_0)$ to $\mf{F}(\mf{x}_0+\mf{h})$ for all $\mf{h}$ in a neighbourhood of $\mf{0}$. But then a reasonable substitute of the number $c$ should be a linear transformation $\mf{T}:\mb{R}^n \to \mb{R}^m$, so that for each $\mf{h}$ in a neighbourhood of $\mf{0}$, it gives a vector $\mf{T}(\mf{h})$ in $\mb{R}^m$. As now $\mf{h}$ is a vector in $\mb{R}^n$, we cannot divide by $\mf{h}$ in \eqref{230730_1}. It should be replaced  with $\Vert\mf{h}\Vert$, the norm of $\mf{h}$.

\begin{definition}{Differentiability}
Let $\mathcal{O}$ be an open subset of $\mb{R}^n$ that contains the point $\mf{x}_0$, and let $\mf{F}:\mathcal{O}\to \mb{R}^m$ be a function defined on $\mathcal{O}$. The function $\mf{F}:\mathcal{O}\to \mb{R}^m$ is {\it differentiable} at $\mf{x}_0$ provided that there exists a linear transformation $\mf{T}:\mb{R}^n \to \mb{R}^m$  so that 
\[\lim_{\mf{h}\to\mf{0}}\frac{\mf{F}(\mf{x}_0+\mf{h})-\mf{F}(\mf{x}_0)-\mf{T}(\mf{h})}{\Vert\mf{h}\Vert}=\mf{0}.\]  $\mf{F}:\mathcal{O}\to \mb{R}^m$ is differentiable if it is differentiable at each point of $\mathcal{O}$.
\end{definition}
 
\begin{remark}{}
 
The differentiability of $\mf{F}:\mathcal{O}\to\mb{R}^m$ at $\mf{x}_0$ amounts to the existence of a linear transformation $\mf{T}:\mb{R}^n \to \mb{R}^m$  so that \[\mf{F}(\mf{x}_0+\mf{h})=\mf{F}(\mf{x}_0)+\mf{T}(\mf{h})+ \boldsymbol{\varepsilon}(\mf{h})\Vert\mf{h}\Vert,\] where $ \boldsymbol{\varepsilon}(\mf{h})\to\mf{0}$ as $\mf{h}\to \mf{0}$.
\end{remark}

The following is obvious from the definition.
\begin{proposition}[label=230730_5]{}
Let $\mathcal{O}$ be an open subset of $\mb{R}^n$ that contains the point $\mf{x}_0$, and let $\mf{F}:\mathcal{O}\to \mb{R}^m$ be a function defined on $\mathcal{O}$. The function $\mf{F}:\mathcal{O}\to \mb{R}^m$ is differentiable at $\mf{x}_0$ if and only if each of its component functions  $F_j:\mathcal{O}\to \mb{R}$, $1\leq j\leq m$ is differentiable at $\mf{x}_0$. 
\end{proposition}

\begin{myproof}{Proof}
Let the components of the function \[ \boldsymbol{\varepsilon}(\mf{h})=\frac{\mf{F}(\mf{x}_0+\mf{h})-\mf{F}(\mf{x}_0)-\mf{T}(\mf{h})}{\Vert\mf{h}\Vert} \] be $ \varepsilon_1(\mf{h}), \varepsilon_2(\mf{h}), \ldots, \varepsilon_m(\mf{h})$. Then  for $1\leq j\leq m$,
\[ \varepsilon_j(\mf{h})=\frac{F_j(\mf{x}_0+\mf{h})-F_j(\mf{x}_0)-T_j(\mf{h})}{\Vert\mf{h}\Vert}.\]
The assertion of the proposition follows  from the fact that \[\di\lim_{\mf{h}\to \mf{0}}\boldsymbol{\varepsilon}(\mf{h})=\mf{0}\;\;\text{ if and only if }\;\;   \lim_{\mf{h}\to \mf{0}}\varepsilon_j(\mf{h})=0\;\; \text{for all }\;1\leq j\leq m,\] while $\di \lim_{\mf{h}\to \mf{0}}\varepsilon_j(\mf{h})=0$ if and only if $F_j:\mathcal{O}\to \mb{R}$ is differentiable at $\mf{x}_0$.
\end{myproof}

Let us look at a simple example of differentiable functions.
\begin{example}[label=230730_4]{}
Let $A$ be an $m\times n$ matrix, and let $\mf{b}$ be a point in $\mb{R}^m$. Define the function $\mf{F}:\mb{R}^n\to\mb{R}^m$ by
\[\mf{F}(\mf{x})=A\mf{x}+\mf{b}.\]
Show that $\mf{F}:\mb{R}^n\to\mb{R}^m$ is differentiable.
\end{example}
\begin{solution}{Solution}
Given $\mf{x}_0$ and $\mf{h}$ in $\mb{R}^n$, notice that
\begin{equation}\label{230730_3}\mf{F}(\mf{x}_0+\mf{h})-\mf{F}(\mf{x}_0)=A(\mf{x}_0+\mf{h})+\mf{b}-A\mf{x}_0-\mf{b}=A\mf{h}.\end{equation}
The map $\mf{T}:\mb{R}^n\to\mb{R}^m$ defined as $\mf{T}(\mf{h})=A\mf{h}$  is a linear transformation. Eq. \eqref{230730_3} says that
\[\mf{F}(\mf{x}_0+\mf{h})-\mf{F}(\mf{x}_0)-\mf{T}(\mf{h})=\mf{0}.\]
Thus,
\[\lim_{\mf{h}\to\mf{0}}\frac{\mf{F}(\mf{x}_0+\mf{h})-\mf{F}(\mf{x}_0)-\mf{T}(\mf{h})}{\Vert\mf{h}\Vert}=\mf{0}.\]Therefore, $\mf{F}$ is differentiable at $\mf{x}_0$. Since the point $\mf{x}_0$ is arbitrary, the function $\mf{F}:\mb{R}^n\to\mb{R}^m$ is differentiable.
\end{solution}

The next theorem says that differentiability implies continuity.
\begin{theorem}{Differentiability Implies Continuity}
Let $\mathcal{O}$ be an open subset of $\mb{R}^n$ that contains the point $\mf{x}_0$, and let $\mf{F}:\mathcal{O}\to \mb{R}^m$ be a function defined on $\mathcal{O}$. If the function $\mf{F}:\mathcal{O}\to \mb{R}^m$ is differentiable at $\mf{x}_0$, then it is continuous at $\mf{x}_0$.
\end{theorem}
\begin{myproof}{Proof}
Since $\mf{F}:\mathcal{O}\to \mb{R}^m$ is differentiable at $\mf{x}_0$, there exists a linear transformation $\mf{T}:\mb{R}^n\to\mb{R}^m$ such that  
\[\boldsymbol{\varepsilon}(\mf{h})=\frac{\mf{F}(\mf{x}_0+\mf{h})-\mf{F}(\mf{x}_0)-\mf{T}(\mf{h})}{\Vert\mf{h}\Vert}\;\;\xrightarrow{\;\;\mf{h}\to\mf{0}\;\;} \;\;\mf{0}.\]
\bp

By Theorem \ref{230724_9}, there is a positive constant $c$ such that
\[\Vert \mf{T}(\mf{h})\Vert\leq c\Vert \mf{h}\Vert\hspace{1cm}\text{for all}\;\mf{h}\in\mb{R}^n.\]
Therefore,
\[\Vert \mf{F}(\mf{x}_0+\mf{h})-\mf{F}(\mf{x}_0)\Vert \leq 
\Vert\mf{T}(\mf{h})\Vert +\Vert\mf{h}\Vert \Vert\boldsymbol{\varepsilon}(\mf{h})\Vert\leq \Vert\mf{h}\Vert\left(c+\Vert\boldsymbol{\varepsilon}(\mf{h})\Vert\right).\]
This implies that 
\[\lim_{\mf{h}\to\mf{0}}\mf{F}(\mf{x}_0+\mf{h})=\mf{F}(\mf{x}_0).\]
Thus, $\mf{F}:\mathcal{O}\to \mb{R}^m$ is  continuous at $\mf{x}_0$.
\end{myproof}

\begin{example}[label=230730_10]{}The function $f:\mb{R}^2\to\mb{R}$   defined as
\[f(x,y)=\begin{cases}\di \frac{xy}{x^2+y^2},\quad &\text{if}\;(x,y)\neq (0,0),\\
0,\quad &\text{if}\;(x,y)=(0,0)\end{cases}\] in Example \ref{230729_2} is not differentiable at $(0,0)$ since
it is not continuous at $(0,0)$. However, we have shown that it has partial derivatives at $(0,0)$.
\end{example}

Let us study the function  $\mf{F}:\mb{R}^n\to\mb{R}^m$, 
$\di \mf{F}(\mf{x})=A\mf{x}+\mf{b}$  that is defined in Example \ref{230730_4}.  
The component functions of $\mf{F}$ are
\begin{align*}
F_1(x_1, x_2, \ldots, x_n)&=a_{11}x_1+a_{12}x_2+\cdots+a_{1n}x_n+b_1,\\
F_2(x_1, x_2, \ldots, x_n)&=a_{21}x_1+a_{22}x_2+\cdots+a_{2n}x_n+b_2,\\
&\hspace{2cm}\vdots\\
F_m(x_1, x_2, \ldots, x_n)&=a_{m1}x_1+a_{m2}x_2+\cdots+a_{mn}x_n+b_m.
\end{align*}
Notice that 
\begin{align*}
\nabla F_1(\mathbf{x})&=\mf{a}_1=\left(a_{11}, a_{12}, \ldots, a_{1n}\right),\\
\nabla F_2(\mathbf{x})&=\mf{a}_2=\left(a_{21}, a_{22}, \ldots, a_{2n}\right),\\
&\hspace{1.5cm}\vdots\\
\nabla F_m(\mathbf{x})&=\mf{a}_m=\left(a_{m1}, a_{m2}, \ldots, a_{mn}\right)
\end{align*}are the row vectors of $A$. Hence,
the derivative matrix of $\mf{F}$ is a given by
\[\mf{DF}(\mf{x})=\begin{bmatrix}  \nabla F_1(\mathbf{x}) \\\nabla F_2(\mathbf{x})\\ \vdots \\ \nabla F_m(\mathbf{x}) \end{bmatrix}=\begin{bmatrix} a_{11} & a_{12} & \cdots & a_{1n}\\
a_{21} & a_{22} & \cdots & a_{2n}\\
\vdots & \vdots & \ddots& \vdots\\
a_{m1} & a_{m2} & \cdots & a_{mn}\end{bmatrix},\] which is the matrix $A$ itself.
Observe that
\[\text{DF}(\mf{x})\mf{h}=\begin{bmatrix}a_{11}h_1+a_{12}h_2+\cdots+a_{1n}h_n\\a_{21}h_1+a_{22}h_2+\cdots+a_{2n}h_n\\
\vdots\\a_{m1}h_1+a_{m2}h_2+\cdots+a_{mn}h_n\end{bmatrix}=
\begin{bmatrix}\langle \nabla F_1(\mf{x}), \mf{h}\rangle \\\langle \nabla F_2(\mf{x}), \mf{h}\rangle\\ \vdots \\ \langle \nabla F_m(\mf{x}), \mf{h}\rangle\end{bmatrix}.\]

 From  Example \ref{230730_4}, we suspect that the linear transformation $\mf{T}:\mb{R}^n\to \mb{R}^m$ that appears in the definition of differentiability of a function should be the linear transformation defined by the derivative matrix. In fact, this is the case.

\begin{theorem}[label=230730_11]{}
Let $\mathcal{O}$ be an open subset of $\mb{R}^n$ that contains the point $\mf{x}_0$, and let $\mf{F}:\mathcal{O}\to \mb{R}^m$ be a function defined on $\mathcal{O}$. The following are equivalent.
\begin{enumerate}[(a)]
\item The function $\mf{F}:\mathcal{O}\to \mb{R}^m$ is differentiable at $\mf{x}_0$.

\item  The function $\mf{F}:\mathcal{O}\to \mb{R}^m$ has partial derivatives at $\mf{x}_0$, and \begin{equation}\label{230730_7}\lim_{\mf{h}\to\mf{0}}\frac{\mf{F}(\mf{x}_0+\mf{h})-\mf{F}(\mf{x}_0)-\mf{DF}(\mf{x}_0) \mf{h}}{\Vert\mf{h}\Vert}=\mf{0}.\end{equation}

\item For each $1\leq j\leq m$, the component function $F_j:\mathcal{O}\to \mb{R}$ has partial derivatives at $\mf{x}_0$, and \[ \lim_{\mf{h}\to\mf{0}}\frac{F_j(\mf{x}_0+\mf{h})-F_j(\mf{x}_0)-\langle \nabla F_j(\mf{x}_0),\mf{h}\rangle}{\Vert\mf{h}\Vert}=0.\]
\end{enumerate}
\end{theorem}
\begin{myproof}{Proof}
The equivalence of (b) and (c) is Proposition \ref{230730_5}, the componentwise differentiability. Thus, we are left to prove the equivalence of (a) and (b).

First, we prove (b) implies (a).
If (b) holds, let $\mf{T}:\mb{R}^n\to\mb{R}^m$ be the linear transformation defined by the derivative matrix $\mf{DF}(\mf{x}_0) $. Then \eqref{230730_7} says that $\mf{F}:\mathcal{O}\to\mb{R}^m$ is differentiable at $\mf{x}_0$.

Conversely,  assume that $\mf{F}:\mathcal{O}\to \mb{R}^m$ is differentiable at $\mf{x}_0$. Then there exists a linear transformation $\mf{T}:\mb{R}^n\to\mb{R}^m$ such that
\begin{equation}\label{230730_6}\lim_{\mf{h}\to\mf{0}}\frac{\mf{F}(\mf{x}_0+\mf{h})-\mf{F}(\mf{x}_0)-\mf{T}( \mf{h})}{\Vert\mf{h}\Vert}=\mf{0}.\end{equation}
Let $A$ be a $m\times n$ matrix so that $\mf{T}( \mf{h})=A\mf{h}$. 
For $1\leq i\leq n$, eq. \eqref{230730_6} implies that
\[\lim_{h\to 0}\frac{\mf{F}(\mf{x}_0+h\mf{e}_i)-\mf{F}(\mf{x}_0)-A( h\mf{e}_i)}{h}=\mf{0}.\]

This gives 
\[A\mf{e}_i=\lim_{h\to 0}\frac{\mf{F}(\mf{x}_0+h\mf{e}_i)-\mf{F}(\mf{x}_0) }{h}.\]
This shows that $\di\frac{\pa \mf{F}}{\pa x_i}(\mf{x}_0)$ exists and
\[ \frac{\pa \mf{F}}{\pa x_i}(\mf{x}_0)=A\mf{e}_i.\]
 
Therefore,  $\mf{F}:\mathcal{O}\to \mb{R}^m$ has partial derivatives at $\mf{x}_0$. Since
\begin{align*}
A&=\begin{bmatrix}A\mf{e}_1 &\rvline & A\mf{e}_2 & \rvline & \cdots & \rvline & A\mf{e}_n\end{bmatrix}\\&=\begin{bmatrix} \di\frac{\pa \mf{F}}{\pa x_1}(\mf{x}_0)&\rvline &\di \frac{\pa \mf{F}}{\pa x_2}(\mf{x}_0)&\rvline &\cdots &\rvline & \di\frac{\pa \mf{F}}{\pa x_n}(\mf{x}_0)\end{bmatrix}=\mf{DF}(\mf{x}_0),\end{align*}eq. \eqref{230730_6} says that  \[\lim_{\mf{h}\to\mf{0}}\frac{\mf{F}(\mf{x}_0+\mf{h})-\mf{F}(\mf{x}_0)-\mf{DF}(\mf{x}_0) \mf{h}}{\Vert\mf{h}\Vert}=\mf{0}.\]
This proves (a) implies (b).
\end{myproof}

\begin{corollary}[label=230731_2]{}
Let $\mathcal{O}$ be an open subset of $\mb{R}^n$ that contains the point $\mf{x}_0$, and let $\mf{F}:\mathcal{O}\to \mb{R}^m$ be a function defined on $\mathcal{O}$. If the partial derivatives of  $\mf{F}:\mathcal{O}\to \mb{R}^m$ exist at $\mf{x}_0$, but
\[\lim_{\mf{h}\to\mf{0}}\frac{\mf{F}(\mf{x}_0+\mf{h})-\mf{F}(\mf{x}_0)-\mf{DF}(\mf{x}_0) \mf{h}}{\Vert\mf{h}\Vert}\neq \mf{0},\]then $\mf{F}$ is not differentiable at $\mf{x}_0$.
\end{corollary}
\begin{myproof}{Proof}
If  $\mf{F}$ is   differentiable at $\mf{x}_0$, Theorem \ref{230730_11} says that we must have \[\lim_{\mf{h}\to\mf{0}}\frac{\mf{F}(\mf{x}_0+\mf{h})-\mf{F}(\mf{x}_0)-\mf{DF}(\mf{x}_0) \mf{h}}{\Vert\mf{h}\Vert}=\mf{0}.\]
By contrapositive, since \[\lim_{\mf{h}\to\mf{0}}\frac{\mf{F}(\mf{x}_0+\mf{h})-\mf{F}(\mf{x}_0)-\mf{DF}(\mf{x}_0) \mf{h}}{\Vert\mf{h}\Vert}\neq \mf{0},\]we find that $\mf{F}$ is not differentiable at $\mf{x}_0$.
\end{myproof}

\begin{example}[label=230730_13]{}
Let $f:\mb{R}^2\to\mb{R}$ be the function defined as
\[f(x,y)=\begin{cases}\di \frac{x^3}{x^2+y^2},\quad &\text{if}\;(x,y)\neq (0,0),\\
0,\quad &\text{if}\;(x,y)=(0,0).\end{cases}\]Determine whether $f$ is differentiable at $(0,0)$. 
\end{example}
\begin{figure}[ht]
\centering
\includegraphics[scale=0.18]{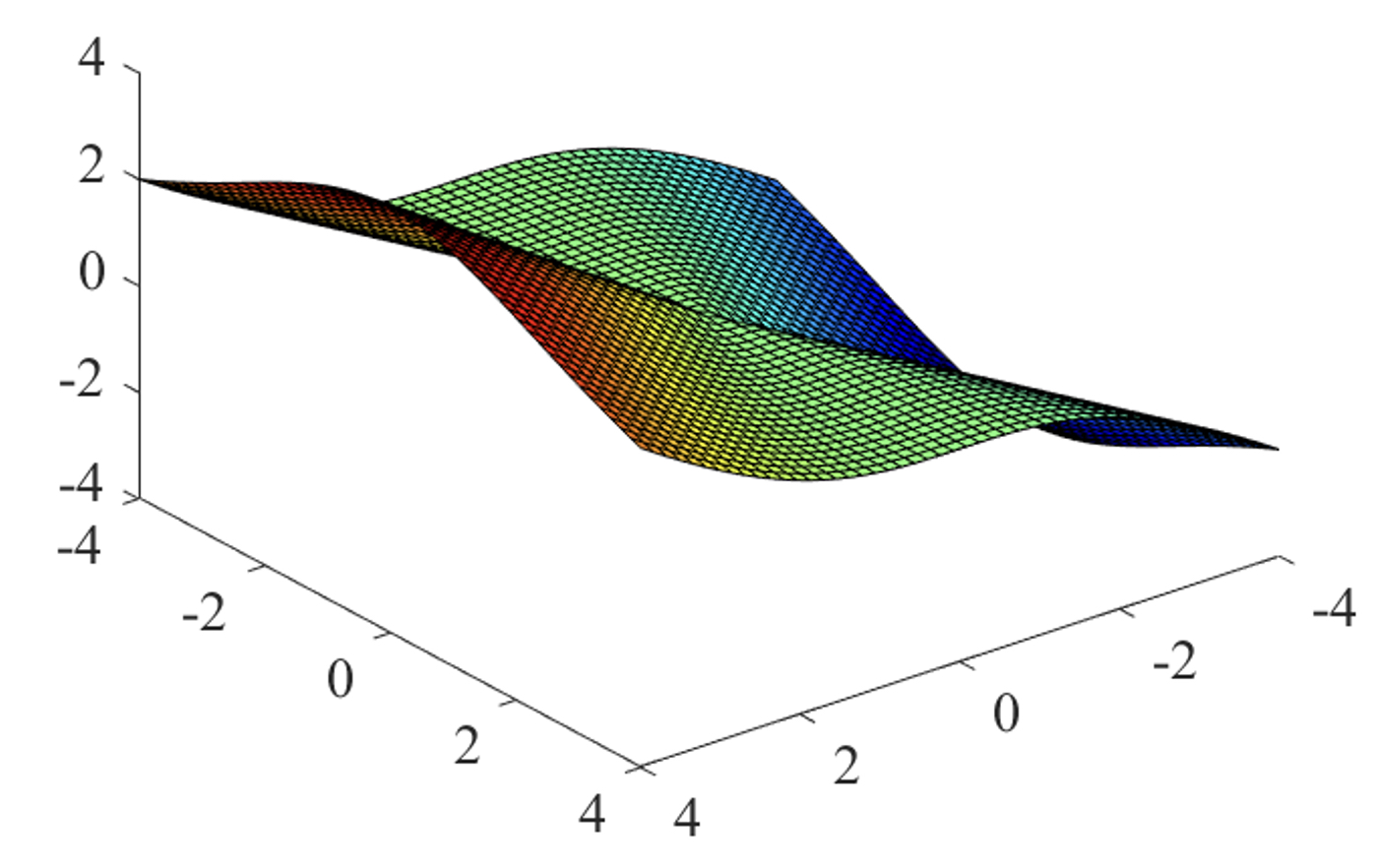}

\caption{The function $\di f(x,y)$ defined in Example \ref{230730_13}.}\label{figure58}
\end{figure}
\begin{solution}{Solution}
One can show that $f$  is continuous at $\mf{0}=(0,0)$.  Hence, we cannot use continuity to determine whether $f$ is differentiable at $\mf{x}_0$. Notice that
\[f_x(0,0)=\lim_{h\to 0}\frac{f(h,0)-f(0,0)}{h}=\lim_{h\to 0}\frac{h-0}{h}=1,\]
\bs
\[f_y(0,0)=\lim_{h\to 0}\frac{f(0,h)-f(0,0)}{h}=\lim_{h\to 0}\frac{0-0}{h}=0.\]
Therefore, $f$ has partial derivatives at $\mf{0}$, and $\nabla f(\mf{0})=( 1,  0)$. 
Now we consider the function
\[\varepsilon(\mf{h})=\frac{f(\mf{h})-f(\mf{0})-\langle \nabla f (\mf{0}), \mf{h}\rangle }{\Vert \mf{h}\Vert}=-\frac{h_1h_2^2}{(h_1^2+h_2^2)^{3/2}}.\]Let $\{\mf{h}_k\}$ be the sequence  with 
$\di \mf{h}_k=\left(\frac{1}{k}, \frac{1}{k}\right)$. 
It converges to $\mf{0}$. Since
\[\varepsilon(\mf{h}_k)=-\frac{1}{2\sqrt{2}}\hspace{1cm}\text{for all}\;k\in\mb{Z}^+,\]
The sequence $\{\varepsilon(\mf{h}_k)\}$ does not converge to 0.
Hence,
\[\lim_{\mf{h}\to\mf{0}}\frac{f(\mf{h})-f(\mf{0})-\langle \nabla f (\mf{0}), \mf{h}\rangle }{\Vert \mf{h}\Vert}\neq 0.\]
Therefore, $f$ is not differentiable at $(0,0)$.
\end{solution}
Example \ref{230730_13} gives a function which is continuous and has partial derivatives at a point, yet it fails to be differentiable at that point. In the following,  we are going to give a sufficient condition for differentiability.  
We begin with a lemma.
\begin{lemma}[label=230730_9]{}
Let $\mf{x}_0$ be a point in $\mb{R}^n$ and let $f:B(\mf{x}_0, r)\to \mb{R}$ be a function defined on an open ball centered at $\mf{x}_0$. Assume that  $f:B(\mf{x}_0, r)\to \mb{R}$  has first order partial derivatives. For each $\mf{h}$ in $\mb{R}^n$ with $\Vert \mf{h}\Vert<r$, there exists $\mf{z}_1, \ldots, \mf{z}_n$ in $B(\mf{x}_0, r)$ such that 
\[f(\mf{x}_0+\mf{h})-f(\mf{x}_0)=\sum_{i=1}^nh_i \frac{\pa f}{\pa x_i}(\mf{z}_i),\]
and
\[\Vert\mf{z}_i-\mf{x}_0\Vert<\Vert\mf{h}\Vert\hspace{1cm}\text{for all}\;1\leq i\leq n.\]
\end{lemma}
\begin{myproof}{Proof}We will take a zigzag path from $\mf{x}_0$ to $\mf{x}_0+\mf{h}$, which is a union of paths parallel to the coordinate axes.
For $1\leq i\leq n$, let
\[\mf{x}_i=\mf{x}_0+\sum_{k=1}^ih_k\mf{e}_k=\mathbf{x}_0+h_1\mf{e}_1+\cdots+h_i\mf{e}_i.\]Then  $\mf{x}_i$
is in $B(\mf{x}_0,r)$. 
Notice that $B(\mf{x}_0, r)$ is a convex set. Therefore, for any $1\leq i\leq n$, the line segment between $\mf{x}_{i-1}$ and $\mf{x}_i=\mf{x}_{i-1}+h_i\mf{e}_i$  lies entirely inside $B(\mf{x}_0, r)$. Since $f:B(\mf{x}_0, r)\to \mb{R}$  has first order partial derivative with respect to $x_i$, the function $g_i:[0,1]\to \mb{R}$,
\[g_i(t)=f(\mf{x}_{i-1}+th_i\mf{e}_i)\]is differentiable   and 
\[g_i'(t)=h_i\frac{\pa f}{\pa x_i} (\mf{x}_{i-1}+th_i\mf{e}_i).\]
By mean value theorem, there exists $c_i\in (0,1)$ such that
\[f(\mf{x}_i)-f(\mf{x}_{i-1})=g_i(1)-g_i(0)=g'_i(c_i)=h_i\frac{\pa f}{\pa x_i} (\mf{x}_{i-1}+c_ih_i\mf{e}_i).\]
\bp
Let \[\mf{z}_i=\mf{x}_{i-1}+c_ih_i\mf{e}_i=\mf{x}_0+\sum_{k=1}^{i-1}h_k\mf{e}_k+c_ih_i\mf{e}_i.\] 

Then $\mf{z}_i$ is a point in $B(\mf{x}_0, r)$. Moreover,
\[f(\mf{x}_0+\mf{h})-f(\mf{x}_0)=\sum_{i=1}^n\left(f(\mf{x}_i)-f(\mf{x}_{i-1})\right)=\sum_{i=1}^nh_i \frac{\pa f}{\pa x_i}(\mf{z}_i).\]

For $1\leq i\leq n$, since $c_i\in (0, 1)$, we   have
\[\Vert\mf{z}_i-\mf{x}_0\Vert=\sqrt{h_1^2+\cdots+h_{i-1}^2+c_i^2h_i^2}<\sqrt{h_1^2+\cdots+h_{i-1}^2+h_i^2}\leq \Vert\mf{h}\Vert.\]
This completes the proof.
\end{myproof}
 
\begin{figure}[ht]
\centering
\includegraphics[scale=0.2]{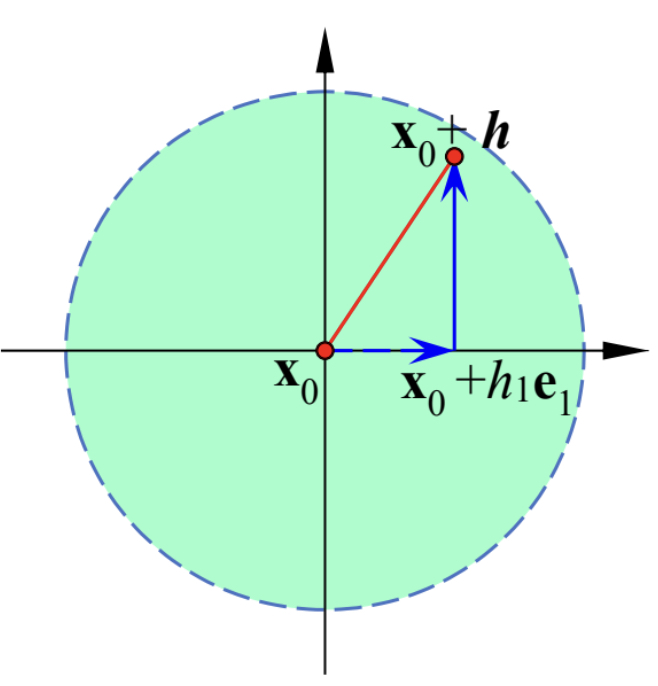}

\caption{A zigzag path from $\mf{x}_0$ to $\mf{x}_0+\mf{h}$.}\label{figure42}
\end{figure}

\begin{theorem}[label=230730_14]{}
Let $\mathcal{O}$ be an open subset of $\mb{R}^n$ that contains the point $\mf{x}_0$, and let $\mf{F}:\mathcal{O}\to \mb{R}^m$ be a function defined on $\mathcal{O}$. If the partial derivatives of  $\mf{F}:\mathcal{O}\to \mb{R}^m$ exists and are continuous at $\mf{x}_0$, then $\mf{F}$ is differentiable at $\mf{x}_0$.
\end{theorem}
\begin{myproof}{Proof}
By Proposition \ref{230730_5}, it suffices to prove the theorem for a function $ f:\mathcal{O}\to \mb{R}$ with codomain $\mb{R}$. Since $\mathcal{O}$ is an open set that contains the point $\mf{x}_0$, there exists $r>0$ such that $B(\mf{x}_0, r)\subset \mathcal{O}$. By Lemma \ref{230730_9}, for each $\mf{h}$ that satisfies $0<\Vert \mf{h}\Vert<r$, there exists $\mf{z}_1, \mf{z}_2, \ldots, \mf{z}_n$ such that
\[f(\mf{x}_0+\mf{h})-f(\mf{x}_0)=\sum_{i=1}^nh_i \frac{\pa f}{\pa x_i}(\mf{z}_i),\]
and
\[\Vert\mf{z}_i-\mf{x}_0\Vert<\Vert\mf{h}\Vert\hspace{1cm}\text{for all}\;1\leq i\leq n.\]Therefore, 
\begin{align*}
 \frac{f(\mf{x}_0+\mf{h})-f(\mf{x}_0)-\langle \nabla f(\mf{x}_0),\mf{h}\rangle}{\Vert\mf{h}\Vert}  =\sum_{i=1}^n\frac{h_i}{\Vert\mf{h}\Vert} \left(\frac{\pa f}{\pa x_i}(\mf{z}_i)-\frac{\pa f}{\pa x_i}(\mf{x}_0)\right).\end{align*}
 
Fixed $\varepsilon>0$. For $1\leq i\leq n$, since $f_{x_i}:B(\mf{x}_0, r)\to\mb{R}$ is continuous at $\mf{x}_0$, there exists $0<\delta_i\leq r$ such that if $0<\Vert \mf{z} -\mf{x}_0\Vert <\delta_i$, then
\[|f_{x_i}(\mf{z})-f_{x_i}(\mf{x}_0)|<\frac{\varepsilon}{n}.\]

Take $\delta=\min\{\delta_1, \ldots, \delta_n\}$. Then $\delta>0$. If $\Vert\mf{h}\Vert <\delta$, then for $1\leq i\leq n$, $\Vert \mf{z}_i-\mf{x}_0\Vert<\Vert\mf{h}\Vert<\delta\leq\delta_i$. Thus,
\[|f_{x_i}(\mf{z_i})-f_{x_i}(\mf{x}_0)|<\frac{\varepsilon}{n}.\]
This implies that
\begin{align*}
\left| \frac{f(\mf{x}_0+\mf{h})-f(\mf{x}_0)-\langle \nabla f(\mf{x}_0),\mf{h}\rangle}{\Vert\mf{h}\Vert}\right|& \leq \sum_{i=1}^n\frac{|h_i|}{\Vert\mf{h}\Vert} \left|\frac{\pa f}{\pa x_i}(\mf{z}_i)-\frac{\pa f}{\pa x_i}(\mf{x}_0)\right|\\
&<\varepsilon.
\end{align*}
Hence,
\[\lim_{\mf{h}\to \mf{0}}\frac{f(\mf{x}_0+\mf{h})-f(\mf{x}_0)-\langle \nabla f(\mf{x}_0),\mf{h}\rangle}{\Vert\mf{h}\Vert}=0.\]
This proves that $f$ is differentiable at $\mf{x}_0$. 
\end{myproof}

Theorem \ref{230730_14} says that a function which has continuous partial derivatives is differentiable.  This prompts us to make the following definition.
\begin{definition}{Continuously Differentiable}
Let $\mathcal{O}$ be an open subset of $\mb{R}^n$, and let $\mf{F}:\mathcal{O}\to\mb{R}^m$ be a function defined on $\mathcal{O}$. We say that  $\mf{F}:\mathcal{O}\to\mb{R}^m$ is continuously differentiable, or $C^1$, provided that it has partial derivatives that are continuous. 
\end{definition}
 Theorem \ref{230730_14} says that a continuously differentiable function  is differentiable.
 
 Analogously, we define $C^k$ for any $k\geq 1$.
 \begin{definition}{$\pmb{C^k}$ Functions}
Let $\mathcal{O}$ be an open subset of $\mb{R}^n$, and let $\mf{F}:\mathcal{O}\to\mb{R}^m$ be a function defined on $\mathcal{O}$. We say that  $\mf{F}:\mathcal{O}\to\mb{R}^m$ is $k$-times continuously differentiable, or $C^k$, provided that it has all partial derivatives of order $k$, and each of them is continuous. 
\end{definition}
\begin{definition}{$\pmb{C^{\infty}}$ Functions}
Let $\mathcal{O}$ be an open subset of $\mb{R}^n$, and let $\mf{F}:\mathcal{O}\to\mb{R}^m$ be a function defined on $\mathcal{O}$. We say that  $\mf{F}:\mathcal{O}\to\mb{R}^m$ is infinitely differentiable, or $C^{\infty}$, provided that it is $C^k$ for all positive integers $k$.
\end{definition}
\begin{proposition}{}
Polynomials and rational functions are infinitely differentiable functions.
\end{proposition}
\begin{myproof}{Sketch of Proof}
A partial derivative of a rational function is still a rational function, which is continuous.
\end{myproof}
Obviously, for any $k\in\mb{Z}^+$, a $C^{k+1}$ function is $C^k$.  

\begin{remark}{Higher Order Differentiability}
We can define second order differentiability in the following way. We say that a function $\mf{F}:\mathcal{O}\to \mb{R}$ is twice differentiable at a point $\mf{x}_0$ in $\mathcal{O}$ if there is a neighbourhood of $\mf{x}_0$ which $\mf{F}$ has first order partial derivatives, and each of them is differentiable at the point $\mathbf{x}_0$. Theorem \ref{230730_14} says that a $C^2$ function is twice differentiable.

Similarly, we can define higher order differentiability.
\end{remark}
\subsection{First Order Approximations}
First we extend the concept of order of approximation to multivariable functions.
\begin{definition}{Order of Approximation}
Let $\mathcal{O}$ be an open subset of $\mb{R}^n$ that contains the point $\mf{x}_0$, and let $k$ be a positive integer. We say that the two functions $\mf{F}:\mathcal{O}\to \mb{R}^m$ and $\mf{G}:\mathcal{O}\to \mb{R}^m$ are  $k^{\text{th}}$-order of approximations of each other at $\mathbf{x}_0$ provided that
\[\lim_{\mf{h}\to \mf{0}}\frac{\mf{F}(\mf{x}_0+\mf{h})-\mf{G}(\mf{x}_0+\mf{h})}{\Vert\mf{h}\Vert^k}=\mf{0}.\]
\end{definition}

Recall that a mapping  $\mf{G}:\mathcal{O}\to\mb{R}^m$ is     a   polynomial  mapping of degree at most one if it has the form
\[\mathbf{G}(\mf{x})=\begin{bmatrix} a_{11}x_1+a_{12}x_2+\cdots+a_{1n}x_n+b_1\\
 a_{21}x_1+a_{22}x_2+\cdots+a_{2n}x_n+b_2\\
  \vdots\\
 a_{m1}x_1+a_{m2}x_2+\cdots+a_{mn}x_n+b_m\end{bmatrix}=A\mathbf{x}+\mf{b},\]
 where $A=[a_{ij}]$ and $\mf{b}=(b_1, \ldots, b_m)$. The mapping $\mf{G}$ is a linear transformation if and only if  $\mf{b}=\mf{0}$. 

The following theorem shows that first order approximation is closely related to differentiability. It is a consequence of Theorem \ref{230730_11}.
\begin{theorem}{First Order Approximation Theorem}
Let $\mathcal{O}$ be an open subset of $\mb{R}^n$ that contains the point $\mf{x}_0$, and let $\mf{F}:\mathcal{O}\to\mb{R}^m$ be a function defined on $\mathcal{O}$.
\begin{enumerate}[(a)]
\item If $\mf{F}:\mathcal{O}\to\mb{R}^m$ is continuous at $\mf{x}_0$, and there is a   polynomial mapping $\mf{G}:\mathcal{O}\to \mb{R}^m$ of degree at most one which is a first order approximation of   $\mf{F}:\mathcal{O}\to \mb{R}^m$ at the point $\mf{x}_0$, then $\mf{F}:\mathcal{O}\to\mb{R}^m$ is differentiable at $\mf{x}_0$.
\item If  
$\mf{F}:\mathcal{O}\to\mb{R}^m$ is differentiable at $\mf{x}_0$, then there is a unique polynomial mapping $\mf{G}:\mathcal{O}\to \mb{R}^m$ of degree at most one which is a first order approximation of $\mf{F}$ at $\mf{x}_0$. It is   given by
\[\mf{G}(\mf{x})=\mf{F}(\mf{x}_0)+\mf{DF}(\mf{x}_0)(\mf{x}-\mf{x}_0).\]
\end{enumerate}
\end{theorem}
\begin{myproof}{Proof}
First we prove (a). Assume that $\mf{G}:\mathcal{O}\to \mb{R}^m$ is a   polynomial mapping  of degree at most one which is a first order approximation of   $\mf{F}:\mathcal{O}\to \mb{R}^m$ at the point $\mf{x}_0$. There exists an $m\times n$ matrix $A$ and a vector $\mf{b}$ in $\mb{R}^m$ such that
\[\mf{G}(\mf{x})=A\mathbf{x}+\mf{b}.\]
By assumption,
\begin{equation}\label{230730_15}\lim_{\mf{h}\to \mf{0}}\frac{\mf{F}(\mf{x}_0+\mf{h})-A(\mf{x}_0+\mf{h})-\mf{b}}{\Vert\mf{h}\Vert}=\mf{0}.\end{equation}
 
This implies that
\[\lim_{\mf{h}\to\mf{0}}\left(\mf{F}(\mf{x}_0+\mf{h})-A(\mf{x}_0+\mf{h})-\mf{b}\right)=\mf{0},\]

which gives
\[A\mf{x}_0+\mf{b}=\lim_{\mf{h}\to\mf{0}}\mf{F}(\mf{x}_0+\mf{h})=\mf{F}(\mf{x}_0).\]
Substitute back into \eqref{230730_15}, we find that
\[\lim_{\mf{h}\to \mf{0}}\frac{\mf{F}(\mf{x}_0+\mf{h})-\mf{F}(\mf{x}_0)-A\mf{h}}{\Vert\mf{h}\Vert}=\mf{0}.\]
\bp
Since $\mf{T}(\mf{h})=A\mf{h}$ is a linear transformation, this shows that $\mf{F}:\mathcal{O}\to\mb{R}^m$ is differentiable at $\mf{x}_0$.

Next, we prove (b). If $\mf{F}:\mathcal{O}\to \mb{R}^m$ is differentiable at $\mf{x}_0$, Theorem \ref{230730_11} says that
\[ \lim_{\mf{h}\to\mf{0}}\frac{\mf{F}(\mf{x}_0+\mf{h})-\mf{F}(\mf{x}_0)-\mf{DF}(\mf{x}_0) \mf{h}}{\Vert\mf{h}\Vert}=\mf{0}.\]

This precisely means that the polynomial mapping  $\mf{G}:\mathcal{O}\to \mb{R}^m$,
 \[\mf{G}(\mf{x})=\mf{F}(\mf{x}_0)+\mf{DF}(\mf{x}_0)(\mf{x}-\mf{x}_0), \] is a first order approximation of $\mf{F}:\mathcal{O}\to \mb{R}^m$  at $\mf{x}_0$. By definition, the polynomial mapping $\mf{G}$ has degree at most one. The uniqueness of $\mf{G}$ is also asserted in Theorem \ref{230730_11}.
\end{myproof}

\begin{remark}{}
The first order approximation theorem says that if the function $\mf{F}:\mathcal{O}\to\mb{R}^m$ is differentiable at the point $\mf{u}$, then there is a unique polynomial mapping $\mf{G }:\mathcal{O}\to\mb{R}^m$ of degree at most one which is a first order approximation of $\mf{F}:\mathcal{O}\to\mb{R}^m$ at the point $\mf{u}$. The components of the mapping $\mf{G }:\mathcal{O}\to\mb{R}^m$ are given by
\[G_j(x_1, \ldots, x_n)=F_j(u_1, \ldots, u_n)+\sum_{i=1}^n\frac{\pa F_j}{\pa x_i}(u_1, \ldots, u_n)(x_i-u_i).\]
Notice that this is a (generalization) of Taylor polynomial of order 1.
\end{remark}
\begin{example}{}
Let $\mf{F}:\mb{R}^3\to \mb{R}^2$ be the function defined as
\[\mf{F}(x,y,z)=(xyz^2, x+2y+3z),\] and let $\mf{x}_0=(1, -1, 1)$.
Find a vector $\mf{b}$ in $\mb{R}^2$ and a $2\times 3$ matrix $A$ such that
\[\lim_{\mf{h}\to \mf{0}}\frac{ \mf{F}(\mf{x}_0+\mf{h})-A\mf{h}-\mf{b}}{\Vert\mf{h}\Vert}=\mf{0}.\]
\end{example}
\begin{solution}{Solution}
The function $\mf{F}:\mb{R}^3\to \mb{R}^2$ is itself a polynomial mapping. Hence, it is differentiable. The derivative matrix is given by
\[\mf{DF}(\mf{x})=\begin{bmatrix} yz^2 & xz^2 & 2xyz\\ 1 & 2 & 3\end{bmatrix}.\] By the first order approximation theorem, 
$\mf{b}=\mf{F}(\mf{x}_0)=(-1, 2)$ and
\[A=\mf{DF}(1, -1, 1)=\begin{bmatrix} -1 & 1 & -2\\ 1 & 2 & 3\end{bmatrix}.\]
\end{solution}

\begin{example}{}
 Determine whether the limit $\di\lim_{(x,y)\to (0,0)}\frac{e^{x+2y}-1-x-2y}{\sqrt{x^2+y^2}}$ exists.
 \end{example}
 \begin{solution}{Solution}
 Let $f(x,y)=e^{x+2y}$. Then
 \[\frac{\pa f}{\pa x}(x,y)=e^{x+2y},\hspace{1cm}\frac{\pa f}{\pa y}(x,y)=2e^{x+2y}.\]It follows that
 \[f(0,0)=1,\;\;\frac{\pa f}{\pa x}(0,0)=1,\;\;\frac{\pa f}{\pa y}(0,0)=2.\]
 Since the function $g(x,y)=x+2y$ is continuous and the exponential function is also continuous, $f$ has continuous first order partial derivatives. Hence, $f$ is differentiable. By first order approximation theorem,
 \[\lim_{(x,y)\to (0,0)}\frac{f(x,y)-f(0,0)-x\di\frac{\pa f}{\pa x}(0,0)-y\frac{\pa f}{\pa y}(0,0)}{\sqrt{x^2+y^2}}=0.\]
 \bs
Since 
 \[f(x,y)-f(0,0)-x\frac{\pa f}{\pa x}(0,0)-y\frac{\pa f}{\pa y}(0,0)=e^{x+2y}-1-x-2y,\]we find that
 \[\lim_{(x,y)\to (0,0)}\frac{e^{x+2y}-1-x-2y}{\sqrt{x^2+y^2}}=0.\]
 \end{solution}
 
 \subsection{Tangent Planes}
 
The tangent plane to a graph is  closely related to the concept of differentiability and first order approximations. Recall that the graph of a function $f:\mathcal{O}\to \mb{R}$ defined on a subset of $\mb{R}^n$ is the subset of $\mb{R}^{n+1}$ consists of all the points of the form $(\mf{x}, f(\mf{x}))$ where $\mf{x}\in\mathcal{O}$.
\begin{definition}{Tangent Planes}
Let $\mathcal{O}$ be an open subset of $\mb{R}^n$ that contains the point $\mf{x}_0$, and let $f:\mathcal{O}\to \mb{R}$ be a function defined on $\mathcal{O}$. The graph of $f$ has a tangent plane at $\mf{x}_0$ if it is differentiable at $\mf{x}_0$. In this case, the tangent plane is the hyperplane of  $\mb{R}^{n+1}$ that satisfies the equation
\[x_{n+1}=f(\mf{x}_0)+ \langle \nabla f(\mf{x}_0), \mf{x}-\mf{x}_0\rangle, \hspace{1cm}\text{where}\;\mf{x}=(x_1, \ldots, x_n).\]

\end{definition}
The tangent plane is the graph of the polynomial function  of degree at most one which is the first order approximation of the function $f$ at the point $\mf{x}_0$.

\begin{figure}[ht]
\centering
\includegraphics[scale=0.18]{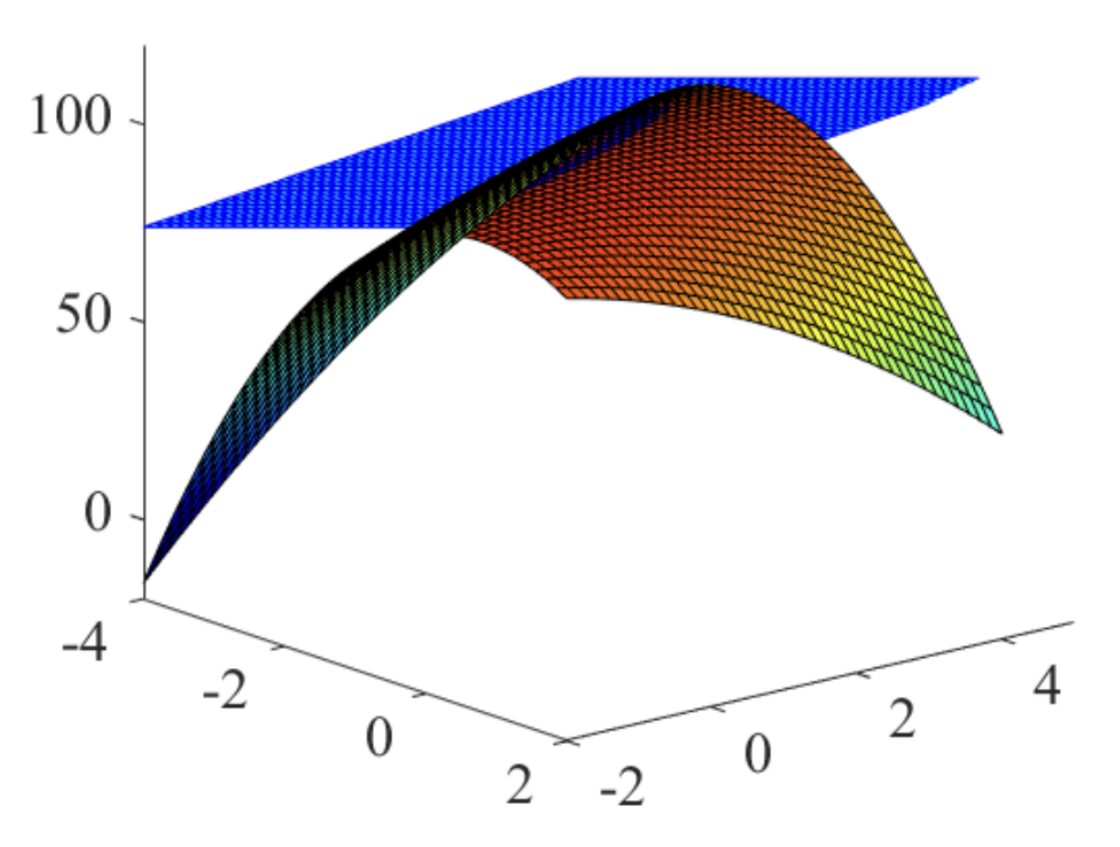}

\caption{The tangent plane to the graph of a function.}\label{figure45}
\end{figure}
\begin{example}{}
Find the equation of the tangent plane to the graph of the  function $f:\mb{R}^2\to \mb{R}$, $f(x,y)=x^2+4xy+5y^2$ at the point where $(x, y)=(1, -1)$. 
\end{example}
\begin{solution}
{Solution}The function $f$ is a polynomial. Hence, it is a differentiable function with \[\nabla f(x,y)=(2x+4y, 4x+10y).\]
\bs
 From this, we find that $\nabla f(1, -1)=(-2, -6)$. Together with
$f(1, -1)=2$, we find that the equation of the tangent plane to the graph of $f$ at the point where $(x,y)=(1,-1)$ is
\[z=2-2(x-1)-6(y+1)=-2x-6y-2.\]
\end{solution}

\subsection{Directional Derivatives}
As we mentioned before, the partial derivatives measure the rate of change of the function when it varies along the directions of the coordinate axes. To capture the rate of change of a function along other directions, we define the concept of directional derivatives.
Notice that a direction in $\mb{R}^n$ is specified by a {\it unit} vector.
\begin{definition}{Directional Derivatives}
Let $\mathcal{O}$ be an open subset of $\mb{R}^n$ that contains the point $\mf{x}_0$, and let $\mf{F}:\mathcal{O}\to \mb{R}^m$ be a function defined on $\mathcal{O}$.
Given a {\bf unit} vector $\mf{u}$ in $\mb{R}^n$, we say that $\mf{F}$ has directional derivative in the direction of $\mf{u}$ at the point $\mf{x}_0$ provided that the limit
\[\lim_{h\to 0}\frac{\mf{F}(\mf{x}_0+h\mf{u})-\mf{F}(\mf{x}_0)}{h}\] exists.  This limit, denoted as $\mf{D}_{\mf{u}}\mf{F}(\mf{x}_0)$, is called the  directional derivative of $\mf{F}$ in the direction of $\mf{u}$ at the point $\mf{x}_0$. 
\end{definition}
When $m=1$, it is customary to denote the directional derivative of $f:\mathcal{O}\to\mb{R}$ in the direction of $\mf{u}$ at the point $\mf{x}_0$  as $D_{\mf{u}}f(\mf{x}_0)$.

\begin{remark}{}
For any nonzero vector $\mf{v}$ in $\mb{R}^n$, we can also define $\mf{D}_{\mf{v}}\mf{F}(\mf{x}_0)$ as
\[\mf{D}_{\mf{v}}\mf{F}(\mf{x}_0)=\lim_{h\to 0}\frac{\mf{F}(\mf{x}_0+h\mf{v})-\mf{F}(\mf{x}_0)}{h}.\]However, we will not call it a directional derivative unless $\mf{v}$ is a unit vector.
\end{remark}
\begin{remark}{}
From the definition, it is obvious that when $\mf{u}$ is one of  the standard unit vectors $\mf{e}_1$, $\ldots$, $\mf{e}_n$, then the directional derivative in the direction of $\mf{u}$ is a partial derivative.  More precisely,
\[\mf{D}_{\mf{e}_i}\mf{F}(\mf{x}_0)=\frac{\pa\mf{F}}{\pa x_i}(\mf{x}_0),\hspace{1cm}1\leq i\leq n.\]
\end{remark}

The following is obvious.
\begin{proposition}{}
Let $\mathcal{O}$ be an open subset of $\mb{R}^n$ that contains the point $\mf{x}_0$, and let $\mf{F}:\mathcal{O}\to \mb{R}^m$ be a function defined on $\mathcal{O}$.
Given a nonzero  vector $\mf{v}$ in $\mb{R}^n$, $\mf{D }_{\mf{v}}\mf{F}(\mf{x}_0)$ exists if only if $D_{\mf{v}}F_j(\mf{x}_0)$ exists for all $1\leq j\leq m$. Moreover,
\[\mf{D }_{\mf{v}}\mf{F}(\mf{x}_0)=\left(D_{\mf{v}}F_1(\mf{x}_0), D_{\mf{v}}F_2(\mf{x}_0), \ldots, D_{\mf{v}}F_m(\mf{x}_0)\right).\]
\end{proposition}

\begin{example}{}
Let $f:\mb{R}^2\to \mb{R}$ be the function defined as
\[f(x,y)=x^2y.\]
Given that $\mf{v}=(v_1, v_2)$ is a nonzero vector in $\mb{R}^2$, find $D_{\mf{v}}f(3,2)$. 
\end{example}
\begin{solution}{Solution}
By definition,
\[
D_{\mf{v}}f(3,2)=\lim_{h\to 0}\frac{f(3+hv_1, 2+hv_2)-f(3,2)}{h}
=g'(0),
\]

where 
\[g(h)=f(3+hv_1, 2+hv_2)=(3+hv_1)^2(2+hv_2).\]
 
Since
\[g'(h)=2v_1(3+hv_1)(2+hv_2)+v_2(3+hv_1)^2,\]

we find that \[D_{\mf{v}}f(3,2)=g'(0)=12v_1+9v_2.\]

Take $\mf{v}=\mf{e}_1=(1,0)$ and $\mf{v}=\mf{e}_2=(0,1)$ respectively, we find that 
$f_{x}(3,2)=12$ and $f_y(3,2)=9$. For general $\mf{v}=(v_1, v_2)$, we notice that
\[D_{\mf{v}}f(3,2)=\langle \nabla f(3,2), \mf{v}\rangle.\]
\end{solution}

\begin{example}[label=230806_7]{}
Consider the function $f:\mb{R}^2\to \mb{R}$  defined as
\[f(x,y)=\begin{cases}\di \frac{xy}{x^2+y^2},\quad &\text{if}\;(x,y)\neq (0,0),\\
0,\quad &\text{if}\;(x,y)=(0,0)\end{cases}\] in Example \ref{230729_2}.
 Find all the nonzero vectors $\mf{v}$  for which $D_{\mf{v}}f(0,0)$ exists.
\end{example}
\begin{solution}{Solution}
Given a nonzero vector $\mf{v}=(v_1, v_2)$, $v_1^2+v_2^2\neq 0$. By definition,
\begin{align*}
D_{\mf{v}}f(0,0)=\lim_{h\to 0}\frac{f(hv_1, hv_2)-f(0,0)}{h}=\lim_{h\to 0}\frac{1}{h} \frac{v_1v_2}{v_1^2+v_2^2}.
\end{align*}This limit exists if and only if $v_1v_2=0$, which is the case if $v_1=0$ or $v_2=0$. 
\end{solution}
\begin{figure}[ht]
\centering
\includegraphics[scale=0.17]{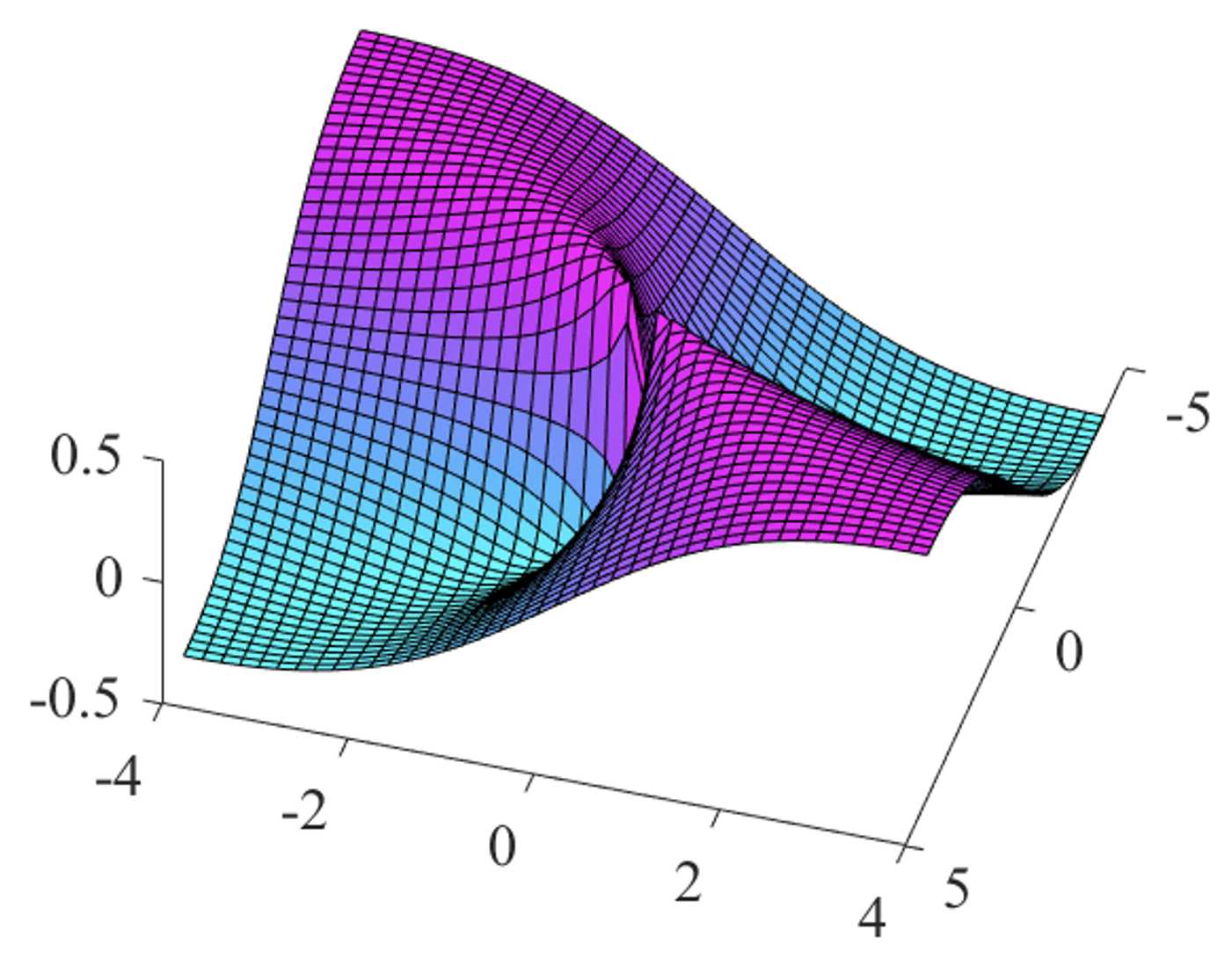}

\caption{The function $f(x,y)$ in Example \ref{230806_7}.}\label{figure46_2}
\end{figure}
\begin{example}[label=230731_1]{}
Let $f:\mb{R}^2\to \mb{R}$ be the function defined as
\[f(x,y)=\begin{cases}\di \frac{y\sqrt{x^2+y^2}}{|x|},\quad &\text{if}\; x\neq 0,\\
0,\quad &\text{if} \;x=0.\end{cases}\]  
 Find all the nonzero vectors $\mf{v}$  for which $D_{\mf{v}}f(0,0)$ exists.
\end{example}

\begin{figure}[ht]
\centering
\includegraphics[scale=0.18]{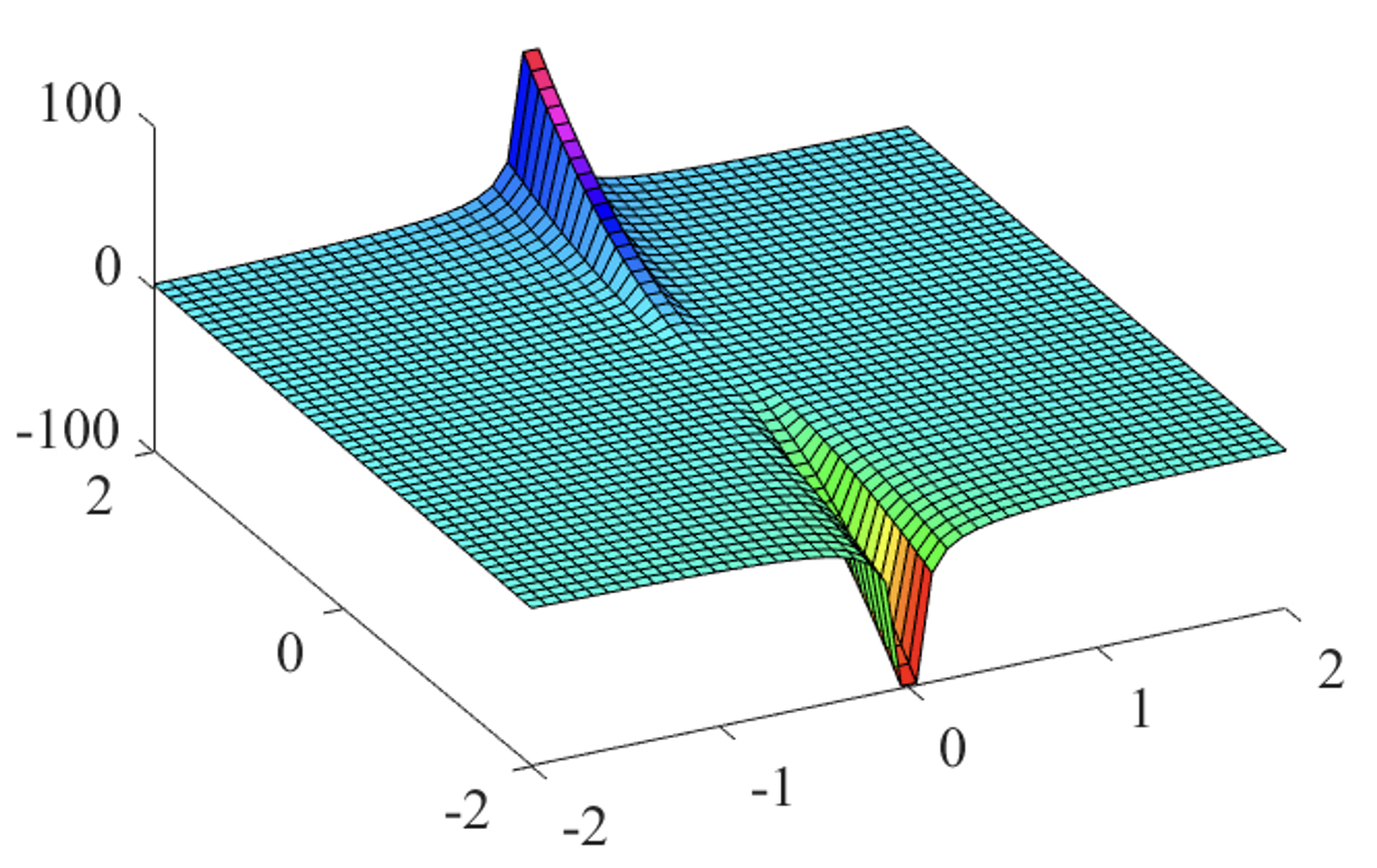}

\caption{The function $f(x,y)$ in Example \ref{230731_1}.}\label{figure59}
\end{figure}
\begin{solution}{Solution}
Given a nonzero vector $\mf{v}=(v_1, v_2)$, we consider two cases. 

\textbf{Case I:} $v_1=0$. \\
Then $\mf{v}=(0,v_2)$. In this case,
\[D_{\mf{v}}f(0,0)=\lim_{h\to 0}\frac{f(0, hv_2)-f(0,0)}{h}=\lim_{h\to 0}\frac{0-0}{h}=0.\]

\textbf{Case 2:} $v_1\neq 0$.  
\begin{align*}
D_{\mf{v}}f(0,0)&=\lim_{h\to 0}\frac{f(hv_1, hv_2)-f(0,0)}{h}\\
&=\lim_{h\to 0}\frac{1}{h}\frac{hv_2}{|hv_1|}\sqrt{h^2(v_1^2+v_2^2)}\\
&=\frac{v_2\sqrt{v_1^2+v_2^2}}{|v_1|}.
\end{align*} 
We conclude that $D_{\mf{v}}f(0,0)$ exists for all nonzero vectors $\mf{v}$.
\end{solution}

\begin{remark}{}
For the function considered in Example \ref{230731_1}, by taking $\mf{v}$ to be $(1,0)$ and $(0,1)$ respectively, we find that $f_{x}(0,0)=0$ and $f_y(0,0)=0$. 
Notice that
\[\lim_{\mf{h}\to\mf{0}}\frac{f(\mf{h})-f(\mf{0})-\langle \nabla f(\mf{0}), \mf{h}\rangle}{\Vert\mf{h}\Vert}=\lim_{\mf{h}\to\mf{0}}\frac{h_2}{|h_1|}.\]
This limit does not exist. By Corollary \ref{230731_2}, $f$ is not differentiable at $(0,0)$. This gives an example of a function which is not differentiable at $(0,0)$ but has directional derivatives at $(0,0)$ in all directions. 
In fact, one can show that $f$ is not continuous at $(0,0)$. 
\end{remark}

The following theorem says that differentiability of a function implies existence of directional derivatives.
\begin{theorem}[label=230802_1]{}
Let $\mathcal{O}$ be an open subset of $\mb{R}^n$ that contains the point $\mf{x}_0$, and let $\mf{F}:\mathcal{O}\to \mb{R}^m$ be a function defined on $\mathcal{O}$. If $\mf{F}$ is differentiable at $\mf{x}_0$, then for any nonzero vector $\mf{v}$, $\mf{D}_{\mf{v}}\mf{F}(\mf{x}_0)$ exists and
\[\mf{D}_{\mf{v}}\mf{F}(\mf{x}_0)=\mf{DF}(\mf{x}_0)\mf{v}=\begin{bmatrix} \langle \nabla F_1(\mf{x}_0), \mf{v}\rangle\\  \langle \nabla F_2(\mf{x}_0), \mf{v}\rangle\\\vdots \\  \langle \nabla F_m(\mf{x}_0), \mf{v}\rangle\end{bmatrix}.\]
\end{theorem}
\begin{myproof}{Proof}
Again, it is sufficient to consider a function $f:\mathcal{O}\to \mb{R}$ with codomain $\mb{R}$.
By definition, $D_{\mf{v}}f(\mf{x}_0)$ is given by the limit
\[\lim_{h\to 0}\frac{f(\mf{x}_0+h\mf{v})-f(\mf{x}_0)}{h}\]if it exists. Since $f$ is differentiable at $\mf{x}_0$, it has partial derivatives at $\mf{x}_0$ and 
\[\lim_{\mf{h}\to \mf{0}}\frac{f(\mf{x}_0+\mf{h})-f(\mf{x}_0)-\langle \nabla f(\mf{x}_0),\mf{h}\rangle }{\Vert\mf{h}\Vert}=0.\]
 
As $h\to 0$, $h\mf{v}\to\mf{0}$. By limit law for composite functions, we find that
\[\lim_{h\to 0}\frac{f(\mf{x}_0+h\mf{v})-f(\mf{x}_0)-\langle \nabla f(\mf{x}_0),h\mf{v}\rangle }{|h|\Vert  \mf{v}\Vert}=0.\]

This implies that
\[\lim_{h\to 0}\frac{f(\mf{x}_0+h\mf{v})-f(\mf{x}_0)-h\langle \nabla f(\mf{x}_0), \mf{v}\rangle }{h}=0.\]
Thus, \[D_{\mf{v}}f(\mf{x}_0)=\lim_{h\to 0}\frac{f(\mf{x}_0+h\mf{v})-f(\mf{x}_0)}{h}=\langle \nabla f(\mf{x}_0), \mf{v}\rangle.\]
\end{myproof}

\begin{example}{}
Consider the function $\mf{F}:\mb{R}^2\to \mb{R}^2$ defined as $\mf{F}(x,y)=(x^2y, xy^2)$. Find $\mf{D}_{\mf{v}}\mf{F}(2, 3)$ when $\mf{v}=(-1, 2)$.

\end{example}
\begin{solution}{Solution}
Since $\mf{F}$ is a polynomial mapping, it is differentiable. The derivative matrix is
$\di \mf{DF}(x,y)=\begin{bmatrix} 2xy & x^2\\ y^2 & 2xy\end{bmatrix}$. Therefore,
\[\mf{D}_{\mf{v}}\mf{F}(2, 3)=\mf{DF}(2,3)\begin{bmatrix} -1\\2\end{bmatrix}=\begin{bmatrix} 12 & 4\\ 9 & 12\end{bmatrix}\begin{bmatrix} -1\\2\end{bmatrix}=\begin{bmatrix} -4\\15\end{bmatrix}.\]

\end{solution}

Theorem \ref{230802_1} can be used to determine the direction which a differentiable function increase fastest at a point.

\begin{corollary}{}
Let $\mathcal{O}$ be an open subset of $\mb{R}^n$ that contains the point $\mf{x}_0$, and let $f:\mathcal{O}\to \mb{R}$ be a function defined on $\mathcal{O}$. If $f$ is differentiable at $\mf{x}_0$ and $\nabla f(\mf{x}_0)\neq\mf{0}$, then at the point $\mf{x}_0$, the function $f$ increases fastest in the direction of $\nabla f(\mf{x}_0)$.
\end{corollary}
\begin{myproof}{Proof}
Let $\mf{u}$ be a unit vector. Then the rate of change of the function $f$ at the point $\mf{x}_0$ in the direction of $\mf{u}$ is given by
\[D_{\mf{u}}f(\mf{x}_0)=\langle\nabla f(\mf{x}_0),\mf{u}\rangle.\]By Cauchy-Schwarz inequality, 
\[\langle\nabla f(\mf{x}_0),\mf{u}\rangle\leq\Vert \nabla f(\mf{x}_0)\Vert\Vert\mf{u}\Vert=\Vert \nabla f(\mf{x}_0)\Vert,\]
and the equality holds if and only if $\mf{u}$ has the same direction as $ \nabla f(\mf{x}_0)$.
\end{myproof}

\vp
\noindent
{\bf \large Exercises  \thesection}
\setcounter{myquestion}{1}

\begin{question}{\themyquestion}
Let $f:\mb{R}^3\to\mb{R}$ be the function defined as
\[ f(x, y, z)=xe^{y^2+4z}.\]
  Find a  vector $\mf{c}$ in $\mb{R}^3$ and a constant $b$ such that
\[\lim_{\mf{h}\to \mf{0}}\frac{f(\mf{x}_0+\mf{h})-\langle \mf{c},\mf{h}\rangle -b}{\Vert \mf{h}\Vert}=0,\]where $\mf{x}_0=(3,2,-1)$.
 
\end{question}
\atc
\begin{question}{\themyquestion}
Let $\mf{F}:\mb{R}^2\to\mb{R}^3$ be the function defined as
\[\mf{F}(x,y)=(x^2+4y^2, 7xy, 2x+y).\]
Find a polynomial mapping $\mf{G}:\mb{R}^2\to\mb{R}^3$ of degree at most one which is a first order approximation of $\mf{F}:\mb{R}^2\to\mb{R}^3$ at the point $(1, -1)$.
\end{question}

\atc
\begin{question}{\themyquestion}
 Let $\mf{x}_0= (1,2, 0, -1)$, and let $\mf{F}:\mb{R}^4\to\mb{R}^3$ be the function defined as
\[\mf{F}(x_1, x_2, x_3, x_4)=\left( x_2x_3^2, \,x_3x_4^3+x_2, \,x_4+2x_1+1\right).\]
  
 Find a $3\times 4$ matrix $A$ and a vector $\mf{b}$ in $\mb{R}^3$ such that
\[\lim_{\mf{x}\to\mf{x}_0}\frac{\mf{F}(\mf{x})-A\mf{x}-\mf{b}}{\Vert \mf{x}-\mf{x}_0\Vert}=\mf{0}.\]
 
\end{question}

\atc
\begin{question}{\themyquestion}
Let $f:\mb{R}^2\to\mb{R}$ be the function defined as
\[f(x,y)=\sin(x^2+y)+5xy^2.\]
 Find $D_{\mf{v}}f(1,-1)$ for any nonzero vector $\mf{v}=(v_1, v_2)$.
 
\end{question}
\atc
\begin{question}{\themyquestion}
Let $f:\mb{R}^2\to\mb{R}$ be the function defined as
\[f(x,y)=\begin{cases}\di\frac{x^2y^2}{x^2+y^2},\quad &\text{if}\;(x,y)\neq (0,0)\\
0,\quad &\text{if}\;(x,y)=(0,0).\end{cases}\]
  Show that $f:\mb{R}^2\to\mb{R}$ is  continuously differentiable.
\end{question}
\atc
\begin{question}{\themyquestion}
Find the equation of the tangent plane to the graph of the  function $f:\mb{R}^2\to \mb{R}$, $f(x,y)=4x^2+3xy-y^2$ at the point where $(x, y)=(2, -1)$.
\end{question}
\atc
\begin{question}{\themyquestion}
Let $f:\mb{R}^2\to\mb{R}$ be the function defined as
\[f(x,y)=\begin{cases}\di\frac{x^2y}{x^2+y^2},\quad &\text{if}\;(x,y)\neq (0,0),\\
0,\quad &\text{if}\;(x,y)=(0,0).\end{cases}\]
\begin{enumerate}[(a)]
\item Show that $f:\mb{R}^2\to\mb{R}$ is continuous.
\item Show that $f:\mb{R}^2\to\mb{R}$ has partial derivatives.
\item Show that $f:\mb{R}^2\setminus\{(0,0)\}\to\mb{R}$ is differentiable.
\item Show that $f:\mb{R}^2\to\mb{R}$ is not differentiable at $(0,0)$.
\item Find all the nonzero vectors $\mf{v}=(v_1, v_2)$ for which $D_{\mf{v}}f(0,0)$ exists.
\end{enumerate}
\end{question}

\atc
\begin{question}{\themyquestion}
Let $f:\mb{R}^2\to\mb{R}$ be the function defined as
\[f(x,y)=\begin{cases}\di\frac{|x| \sqrt{x^2+y^2}}{y},\quad &\text{if}\;y\neq 0,\\
0,\quad &\text{if}\;y=0.\end{cases}\]
\begin{enumerate}[(a)]
\item Show that $f:\mb{R}^2\to\mb{R}$ is not continuous at $(0,0)$.
\item Show that $D_{\mf{v}}f(0,0)$ exists for all nonzero vectors $\mf{v}$.
\end{enumerate}
\end{question}
 \atc
\begin{question}{\themyquestion}
Let $f:\mb{R}^2\to\mb{R}$ be the function defined as
\[f(x,y)=\begin{cases}\di (x^2+y^2)\sin\left( \frac{1}{\sqrt{x^2+y^2}}\right),\quad &\text{if}\;(x,y)\neq (0,0),\\
0,\quad &\text{if}\;(x,y)=(0,0).\end{cases}\]
\begin{enumerate}[(a)]
\item Show that $f:\mb{R}^2\to\mb{R}$ is differentiable at $(0,0)$.
\item Show that  $f:\mb{R}^2\to\mb{R}$ is not continuously differentiable at $(0,0)$.
\end{enumerate}
\end{question}

\section{The Chain Rule and the Mean Value Theorem} 

In volume I, we have seen that the chain rule plays an important role in calculating the derivative of a composite function. Given that $f:(a,b)\to\mb{R}$ and $g:(c,d)\to\mb{R}$ are functions such that $f((a,b))\subset (c,d)$, the chain rule says that if $f$ is differentiable at $x_0$, $g$ is differentiable at $y_0=f(x_0)$, then the composite function $(g\circ f):(a,b)\to\mb{R}$ is differentiable at $x_0$, and
\[(g\circ f)'(x_0)=g'(f(x_0))f'(x_0).\]

For multivariable functions, the chain rule takes the following form.
\begin{theorem}{The Chain Rule}
Let $\mathcal{O}$ be an open subset of $\mb{R}^n$, and let $\mathcal{U}$ be an open subset of $\mb{R}^k$. Assume that $\mf{F}:\mathcal{O}\to \mb{R}^k$ and $\mf{G}:\mathcal{U}\to \mb{R}^m$ are functions such that $\mf{F}(\mathcal{O})\subset \mathcal{U}$. If $\mf{F}$ is differentiable at $\mf{x}_0$, $\mf{G}$ is differentiable at $\mf{y}_0=\mf{F}(\mf{x}_0)$, then the composite function $\mf{H}=(\mf{G}\circ\mf{F}):\mathcal{O}\to \mb{R}^m$ is differentiable at $\mf{x}_0$ and
\[\mf{DH}(\mf{x}_0)=\mf{D}(\mf{G}\circ \mf{F})(\mf{x}_0)=\mf{DG}(\mf{F}(\mf{x}_0))\mf{DF}(\mf{x}_0).\]
\end{theorem}
Notice that on the right hand side, $\mf{DG}(\mf{F}(\mf{x}_0))$ is an $m\times k$ matrix, $\mf{DF}(\mf{x}_0)$ is an $k\times n$ matrix. Hence, the product $\mf{DG}(\mf{F}(\mf{x}_0))\mf{DF}(\mf{x}_0)$ makes sense, and it is an $m\times n$ matrix, which is the correct size for the derivative matrix $\mf{DH}(\mf{x}_0)$.

Let us spell out more explicitly. 
Assume that 
\begin{align*}&\mf{F}(x_1, x_2, \ldots, x_n)\\&=(F_1(x_1, x_2, \ldots, x_n), F_2(x_1, x_2, \ldots, x_n), \ldots, F_k(x_1, x_2, \ldots, x_n)),\\
&\mf{G}(y_1, y_2, \ldots, y_k)\\&=(G_1(y_1, y_2, \ldots, y_k), G_2(y_1, y_2, \ldots, y_k), \ldots, G_m(y_1, y_2, \ldots, y_k)),\\
&\mf{H}(x_1, x_2, \ldots, x_n)\\&=(H_1(x_1, x_2, \ldots, x_n), H_2(x_1, x_2, \ldots, x_n), \ldots, H_m(x_1, x_2, \ldots, x_n)).
\end{align*} Then for $1\leq j\leq m$,
\begin{align*}
&H_j(x_1, x_2, \ldots, x_n)\\&= G_j\left(F_1(x_1, x_2, \ldots, x_n), F_2(x_1,x_2,  \ldots, x_n), \ldots, F_k(x_1, x_2, \ldots, x_n)\right).
\end{align*}For $1\leq l\leq k$, let 
\[y_l= F_l\left(x_1, x_2, \ldots, x_n\right).\] The chain rule says that if $1\leq q\leq n$,
\begin{align*}
   \frac{\pa H_j}{\pa x_q}(x_1, x_2, \ldots, x_n)   &=\sum_{l=1}^k \frac{\pa G_j}{\pa y_l}(y_1, y_2,  \ldots, y_k)\frac{\pa F_l}{\pa x_q}(x_1, x_2, \ldots, x_n)\\
   &= \frac{\pa G_j}{\pa y_1}(y_1, y_2,  \ldots, y_k)\frac{\pa F_1}{\pa x_q}(x_1, x_2, \ldots, x_n)\\
  & \quad  + \frac{\pa G_j}{\pa y_2}(y_1, y_2,  \ldots, y_k)\frac{\pa F_2}{\pa x_q}(x_1, x_2, \ldots, x_n)\\
  &\hspace{3cm}\vdots
  \\
  &\quad+ \frac{\pa G_j}{\pa y_k}(y_1, y_2,  \ldots, y_k)\frac{\pa F_k}{\pa x_q}(x_1, x_2, \ldots, x_n).
\end{align*}
Namely, to differentiate $H_j= G_j\circ \mf{F}$ with respect to $x_q$, we differentiate $G_j$ with respect to each of the variables $y_1, \ldots, y_k$, multiply each by the partial derivatives of $F_1, \ldots, F_k$ with respect to $x_q$, then take the sum.

Let us illustrate this with a simple example. 
\begin{example}{}Consider the function $h:\mb{R}^2\to\mb{R}$ defined as
\[h(x,y)=\sin(2x+3y)+e^{xy}.\]
It is straightforward to find that
\[\frac{\pa h}{\pa x}=2\cos(2x+3y)+ye^{xy},\hspace{1cm}\frac{\pa h}{\pa y}=3\cos(2x+3y)+xe^{xy}.\]
Notice that we can write $h=g\circ \mf{F}$, where $\mf{F}:\mb{R}^2\to \mb{R}^2$ is the function 
\[\mf{F}(x,y)=(2x+3y, xy),\]
and $g:\mb{R}^2\to \mb{R}$ is the function 
\[g(u, v)=\sin u+ e^v.\]
\be
Obviously, $\mf{F}$ and $g$ are continuously differentiable functions. 
\[\mf{DF}(x,y)=\begin{bmatrix} 2 & 3\\ y& x\end{bmatrix},\hspace{1cm}
Dg(u,v)=\begin{bmatrix}\cos u & e^v\end{bmatrix}.\]
Taking $u=2x+3y$ and $v=xy$, we find that
\begin{align*}Dg(u,v)\mf{DF}(x,y)&=\begin{bmatrix}\cos(2x+3y) & e^{xy}\end{bmatrix}
\begin{bmatrix} 2 & 3\\ y& x\end{bmatrix}\\&=\begin{bmatrix} 2\cos(2x+3y)+ye^{xy} & 3\cos(2x+3y)+xe^{xy}\end{bmatrix}\\&=Dh(x,y).
\end{align*}
\end{example2}

Now let us prove the chain rule.
\begin{myproof}{\linkt Proof of the Chain Rule}
Since $\mf{F}$  is differentiable at $\mf{x}_0$ and $\mf{G}$ is differentiable at $\mf{y}_0=\mf{F}(\mf{x}_0)$, $\mf{DF}(\mf{x}_0)$ and $\mf{DG}(\mf{y}_0)$ exist. There exists positive numbers $r_1$ and $r_2$ such that $B(\mf{x}_0, r_1)\subset\mathcal{O}$ and $B(\mf{y}_0, r_2)\subset \mathcal{U}$. 
Let
\begin{align*}
\boldsymbol{\varepsilon}_1(\mf{h})&= \frac{\mf{F}(\mf{x}_0+\mf{h})-\mf{F}(\mf{x}_0)-\mf{DF}(\mf{x}_0)\mf{h}}{\Vert\mf{h}\Vert},\hspace{1cm}\mf{h}\in B(\mf{0}, r_1),
\\
\boldsymbol{\varepsilon}_2(\mf{v})&= \frac{\mf{G}(\mf{y}_0+\mf{v})-\mf{G}(\mf{y}_0)-\mf{DG}(\mf{y}_0)\mf{v}}{\Vert\mf{v}\Vert},\hspace{1cm}\mf{v}\in B(\mf{0}, r_2).
\end{align*}Since $\mf{F}$ is differentiable at $\mf{x}_0$ and $\mf{G}$ is differentiable at $\mf{y}_0$, 
\[\lim_{\mf{h}\to \mf{0}}\boldsymbol{\varepsilon}_1(\mf{h})=\mf{0},\hspace{1cm} \lim_{\mf{v}\to \mf{0}}\boldsymbol{\varepsilon}_2(\mf{v})=\mf{0}.\]
There exist positive constants $c_1$ and $c_2$ such that
\[\Vert \mf{DF}(\mf{x}_0)\mf{h}\Vert\leq c_1\Vert\mf{h}\Vert\hspace{1cm}\text{for all}\;\mf{h}\in\mb{R}^n,\]
\[\Vert \mf{DG}(\mf{y}_0)\mf{v}\Vert\leq c_2\Vert\mf{v}\Vert\hspace{1cm}\text{for all}\;\mf{v}\in\mb{R}^k.\]

Now since $\mf{F}$ is differentiable at $\mf{x}_0$, it is continuous at $\mf{x}_0$. Hence, there exists a positive number $r$ such that $r\leq r_1$ and $\mf{F}(B(\mf{x}_0, r))\subset B(\mf{y}_0, r_2)$. 
\bp
For $\mf{h}\in B(\mf{0}, r)$, let
\[\mf{v}=\mf{F}(\mf{x}_0+\mf{h})-\mf{F}(\mf{x}_0).\]
Then $\mf{v}\in B(\mf{0}, r_2)$ and 
\[\mf{v}=\mf{DF}(\mf{x}_0)\mf{h}+\Vert\mf{h}\Vert \boldsymbol{\varepsilon}_1(\mf{h}).\]
It follows that
\[\Vert\mf{v}\Vert \leq \Vert\mf{DF}(\mf{x}_0)\mf{h}\Vert+\Vert\mf{h}\Vert \Vert\boldsymbol{\varepsilon}_1(\mf{h})\Vert\leq \Vert\mf{h}\Vert\left(c_1+  \Vert\boldsymbol{\varepsilon}_1(\mf{h})\Vert\right).\]
In particular, we find that when $\mf{h}\to \mf{0}$, $\mf{v}\to \mf{0}$.
Now,
\begin{align*}
&\mf{H}(\mf{x}_0+\mf{h})-\mf{H}(\mf{x}_0)\\&=\mf{G}(\mf{F}(\mf{x}_0+\mf{h}))-\mf{G}(\mf{F}(\mf{x}_0))\\
&=\mf{G}(\mf{y}_0+\mf{v})-\mf{G}(\mf{y}_0)\\
&=\mf{DG}(\mf{y}_0)\mf{v}+\Vert\mf{v}\Vert \boldsymbol{\varepsilon}_2(\mf{v})\\
&=\mf{DG}(\mf{y}_0)\mf{DF}(\mf{x}_0)\mf{h}+\Vert\mf{h}\Vert\mf{DG}(\mf{y}_0)\boldsymbol{\varepsilon}_1(\mf{h})+\Vert\mf{v}\Vert \boldsymbol{\varepsilon}_2(\mf{v}).
\end{align*}Therefore, for $\mf{h}\in B(\mf{0},r)\setminus\{\mf{0}\}$, 
\[\frac{\mf{H}(\mf{x}_0+\mf{h})-\mf{H}(\mf{x}_0)-\mf{DG}(\mf{y}_0)\mf{DF}(\mf{x}_0)\mf{h}}{\Vert\mf{h}\Vert} =\mf{DG}(\mf{y}_0)\boldsymbol{\varepsilon}_1(\mf{h})+\frac{\Vert\mf{v}\Vert}{\Vert\mf{h}\Vert} \boldsymbol{\varepsilon}_2(\mf{v})
.\]This implies that
\begin{align*}
&\left\Vert \frac{\mf{H}(\mf{x}_0+\mf{h})-\mf{H}(\mf{x}_0)-\mf{DG}(\mf{y}_0)\mf{DF}(\mf{x}_0)\mf{h}}{\Vert\mf{h}\Vert}\right\Vert\\ &\leq \Vert\mf{DG}(\mf{y}_0)\boldsymbol{\varepsilon}_1(\mf{h})\Vert+\frac{\Vert\mf{v}\Vert}{\Vert\mf{h}\Vert} \Vert\boldsymbol{\varepsilon}_2(\mf{v})\Vert\\
&\leq c_2\Vert \boldsymbol{\varepsilon}_1(\mf{h})\Vert+\left(c_1+  \Vert\boldsymbol{\varepsilon}_1(\mf{h})\Vert\right)\Vert\boldsymbol{\varepsilon}_2(\mf{v})\Vert.
\end{align*}Since $\mf{v}\to\mf{0}$ when $\mf{h}\to \mf{0}$, we find that $\boldsymbol{\varepsilon}_2(\mf{v})\to\mf{0}$ when $\mf{h}\to\mf{0}$. Thus, we find that
\[\lim_{\mf{h}\to \mf{0}} \frac{\mf{H}(\mf{x}_0+\mf{h})-\mf{H}(\mf{x}_0)-\mf{DG}(\mf{y}_0)\mf{DF}(\mf{x}_0)\mf{h}}{\Vert\mf{h}\Vert}=\mf{0}.\]
\bp
This concludes that $\mf{H}$ is differentiable at $\mf{x}_0$ and
\[\mf{DH}(\mf{x}_0)=\mf{DG}(\mf{y}_0)\mf{DF}(\mf{x}_0).\]
\end{myproof}

\begin{example}{}
Let  $\mf{F}:\mb{R}^3\to\mb{R}^2$ be the function defined as
\[\mf{F}(x,y,z)=(x^2+4y^2+9z^2, xyz).\] Find a vector $\mf{b}$ in $\mb{R}^2$ and a $2\times 3$ matrix $A$ such that
\[\lim_{(u, v, w)\to (1, -1, 0)}\frac{\mf{F}(2u+v, v+w, u+w)-\mf{b}-A\mf{p}}{\sqrt{(u-1)^2+(v+1)^2+w^2}}=\mf{0},\hspace{0.2cm}\text{where}\;\mf{p}=\di\begin{bmatrix} u\\v\\w\end{bmatrix}.\]
\end{example}
\begin{solution}{Solution}
Let $\mf{p}_0=(1,-1,0)$, and let $\mf{G}:\mb{R}^3\to\mb{R}^3$ be the mapping \[\mf{G}(u,v,w)=(2u+v, v+w, u+w).\]
Then $\mf{H}(\mf{p})=\mf{H}(u,v,w)=\mf{F}(2u+v, v+w, u+w)=(\mf{F}\circ \mf{G})(u,v,w)$. Notice that $\mf{F}$ and $\mf{G}$ are polynomial mappings. Hence, they are infinitely differentiable. To have
\begin{align*}
&\lim_{\mf{p}\to \mf{p}_0}\frac{\mf{H}(\mf{v})-\mf{b}-A\mf{p}}{\Vert\mf{p}-\mf{p}_0\Vert}\\&=\lim_{(u, v, w)\to (1, -1, 0)}\frac{\mf{F}(2u+v, v+w, u+w)-\mf{b}-A\mf{p}}{\sqrt{(u-1)^2+(v+1)^2+w^2}}=\mf{0},
\end{align*}
the first order approximation theorem says that
 
\[\mf{b}+A \mf{p}=\mf{H}(\mf{p}_0)+\mf{DH}(\mf{p}_0)\left(\mf{p}-\mf{p}_0\right).\]Therefore,
\[A=\mf{DH}(\mf{p}_0)\quad\text{and} \quad \mf{b}=\mf{H}(\mf{p}_0)-A\mf{p}_0.\]
\bs
Notice that $\mf{G}(\mf{p}_0)=\mf{G}(1,-1,0)=(1,-1,1)$, \[\mf{H}(\mf{p}_0)=\mf{H}(1,-1,0)=\mf{F}(1,-1,1)=(14, -1),\]
\[\mf{DG}(u,v,w)=\begin{bmatrix} 2 & 1 & 0\\0 & 1 & 1\\1 & 0 & 1\end{bmatrix},\hspace{1cm}\mf{DF}(x,y,z)=\begin{bmatrix}2x&8y&18z\\yz& xz & xy \end{bmatrix}.\]
By chain rule,
\begin{align*}
A&=\mf{DF}(1,-1,1)\mf{DG}(1,-1,0)\\&=\begin{bmatrix}2 &-8&18 \\-1& 1 & -1 \end{bmatrix}\begin{bmatrix} 2 & 1 & 0\\0 & 1 & 1\\1 & 0 & 1\end{bmatrix}=\begin{bmatrix} 22 & -6 & 10\\-3 & 0 & 0 \end{bmatrix}. \end{align*}
It follows that
\[\mf{b}=\begin{bmatrix} 14\\-1\end{bmatrix}-\begin{bmatrix} 22 & -6 & 10\\-3 & 0 & 0 \end{bmatrix}\begin{bmatrix}  1\\-1\\0\end{bmatrix}=\begin{bmatrix}  -14\\2\end{bmatrix}.\]
\end{solution}

\begin{example}{}
Let $\alpha$ be a positive number, and let $f:\mb{R}^n\to\mb{R}$ be the function defined as
\[f(\mathbf{x})=\Vert \mf{x}\Vert^{\alpha}.\]Find the values of $\alpha$ so that $f$ is differentiable.
\end{example}
\begin{solution}{Solution}
Let $g:\mb{R}^n\to\mb{R}$ be the function \[g(\mf{x})=\Vert\mf{x}\Vert^2= x_1^2+x_2^2+\cdots+x_n^2.\]Then $g(\mb{R}^n)= [0, \infty)$, and $g(\mf{x})=0$ if and only if $\mf{x}=\mf{0}$.
\bs
Since $g$ is a polynomial, it is infinitely differentiable. Let $h:[0,\infty)\to \mb{R}$ be the function $h(u)=u^{\alpha/2}$. Then $h$ is differentiable on $(0,\infty)$. Since $f(\mf{x})=(h\circ g)(\mf{x})$, chain rule implies that for all $\mf{x}_0\in \mb{R}^n\setminus\{\mf{0}\}$, $f$ is differentiable at $\mf{x}_0$.

Now consider the point $\mf{x}=\mf{0}$. Notice that for $1\leq i\leq n$, $f_{x_i}(\mf{0})$ exists provided that the limit
\[\lim_{h\to 0}\frac{f(h\mf{e}_i)-f(\mf{0})}{h}=\lim_{h\to 0} \frac{|h|^{\alpha }}{h} \] exists. This is the case if $\alpha>1$. Therefore, $f$ is not differentiable at $\mf{x}=\mf{0}$ if $\alpha\leq 1$. If $\alpha>1$, we find that $f_{x_i}(\mf{0})=0$ for all $1\leq i\leq n$. Hence, $\nabla f(\mf{0})=\mf{0}$. Since
\[\lim_{\mf{h}\to\mf{0}}\frac{f(\mf{h})-f(\mf{0})-\langle \nabla f(\mf{0}), \mf{h}\rangle}{\Vert\mf{h}\Vert}
=\lim_{\mf{h}\to\mf{0}}\Vert\mf{h}\Vert^{\alpha-1}=0,\]
we conclude that when $\alpha>1$, $f$ is differentiable at $\mf{x}=\mf{0}$.

Therefore, $f$ is differentiable if and only if $\alpha>1$.
\end{solution}

\begin{example}[label=230801_1]{}
Let $f:\mathbb{R}^2\to \mb{R}$ be a twice continuously differentiable function, and let
$g:\mb{R}^2\to \mb{R}$ be the function defined as
\[g(r,\theta)=f(r\cos\theta, r\sin\theta).\]
Show that 
\[\frac{\pa^2 g}{\pa r^2}+\frac{1}{r}\frac{\pa g}{\pa r}+\frac{1}{r^2}\frac{\pa^2 g}{\pa\theta^2}=\frac{\pa^2 f}{\pa x^2}+\frac{\pa^2 f}{\pa y^2}.\]
\end{example}
\begin{solution}{Solution}
Let $\mf{H}:\mb{R}^2\to\mb{R}^2$ be the mapping defined by
\[\mf{H}(r,\theta)=(r\cos\theta, r\sin \theta).\]
\bs
Then $\mf{H}$ is infinitely differentiable, and $g=f\circ\mf{H}$. Let 
$x=H_1(r,\theta)=r\cos\theta$ and $y=H_2(r,\theta)=r\sin\theta$. By chain rule,
\begin{align*}
\frac{\pa g}{\pa r}&=\frac{\pa f}{\pa x}\frac{\pa x}{\pa r}+\frac{\pa f}{\pa y}\frac{\pa y}{\pa r}=\cos\theta\frac{\pa f}{\pa x}+\sin\theta\frac{\pa f}{\pa y},\\
\frac{\pa g}{\pa \theta}&=\frac{\pa f}{\pa x}\frac{\pa x}{\pa \theta}+\frac{\pa f}{\pa y}\frac{\pa y}{\pa \theta}=-r\sin\theta\frac{\pa f}{\pa x}+r\cos\theta\frac{\pa f}{\pa y}.
\end{align*}Using product rule and chain rule, we then have
\begin{align*}
\frac{\pa^2 g}{\pa r^2}&=\cos\theta \left(\frac{\pa^2 f}{\pa x^2}\frac{\pa x}{\pa r}+\frac{\pa^2f}{\pa y\pa x}\frac{\pa y}{\pa r}\right)+\sin\theta \left(\frac{\pa^2 f}{\pa x\pa y}\frac{\pa x}{\pa r}+\frac{\pa^2f}{\pa y^2}\frac{\pa y}{\pa r}\right).\end{align*}
Since $f$ has continuous second order partial derivatives, $f_{xy}=f_{yx}$. Therefore,
\begin{align*}
\frac{\pa^2 g}{\pa r^2}&=\cos^2\theta\frac{\pa^2f}{\pa x^2}+2\sin\theta\cos\theta\frac{\pa^2f}{\pa x\pa y}+\sin^2\theta\frac{\pa^2f}{\pa  y^2}.
\end{align*}
Similarly, we have
\begin{align*}
\frac{\pa^2 g}{\pa \theta^2}&=-r\sin\theta \left(\frac{\pa^2 f}{\pa x^2}\frac{\pa x}{\pa \theta}+\frac{\pa^2f}{\pa y\pa x}\frac{\pa y}{\pa \theta}\right)+r\cos\theta \left(\frac{\pa^2 f}{\pa x\pa y}\frac{\pa x}{\pa \theta}+\frac{\pa^2f}{\pa y^2}\frac{\pa y}{\pa \theta}\right)\\
&\quad -r\cos\theta\frac{\pa f}{\pa x}-r\sin\theta\frac{\pa f}{\pa y}\\
&=r^2\sin^2\theta\frac{\pa^2f}{\pa x^2}-2r^2\sin\theta\cos\theta\frac{\pa^2f}{\pa x\pa y}+r^2\cos^2\theta\frac{\pa^2f}{\pa  y^2}-r\frac{\pa g}{\pa r}.
\end{align*}From these, we obtain
\[\frac{\pa^2 g}{\pa r^2}+\frac{1}{r}\frac{\pa g}{\pa r}+\frac{1}{r^2}\frac{\pa^2 g}{\pa\theta^2}=\frac{\pa^2 f}{\pa x^2}+\frac{\pa^2 f}{\pa y^2}.\]
\end{solution}
 
Example \ref{230801_1}  gives the  {\it Laplacian} \[ \Delta f=\frac{\pa^2 f}{\pa x^2}+\frac{\pa^2 f}{\pa y^2}\] of $f$ in polar coordinates. It is  customary that one would abuse notation and write $g=f$, so that the formula takes the form
\[\frac{\pa^2 f}{\pa x^2}+\frac{\pa^2 f}{\pa y^2}=\frac{\pa^2 f}{\pa r^2}+\frac{1}{r}\frac{\pa f}{\pa r}+\frac{1}{r^2}\frac{\pa^2 f}{\pa\theta^2}.\]

\begin{remark}{}
We can use the chain rule to prove Theorem \ref{230802_1}. Given that  $\mathcal{O}$ is an open subset of $\mb{R}^n$ that contains the point $\mf{x}_0$, and   $\mf{F}:\mathcal{O}\to \mb{R}^m$ is a function that is differentiable at $\mf{x}_0$, we want to show that  $\mf{D}_{\mf{v}}\mf{F}(\mf{x}_0)$ exists for any nonzero vector $\mf{v}$, and 
\[\mf{D}_{\mf{v}}\mf{F}(\mf{x}_0)=\mf{DF}(\mf{x}_0)\mf{v}.\]Since $\mathcal{O}$ is an open set that contains the point $\mf{x}_0$, there is an $r>0$ such that $B(\mf{x}_0,r)\subset\mathcal{O}$. By definition,
\[\mf{D}_{\mf{v}}\mf{F}(\mf{x}_0)=\lim_{h\to 0}\frac{\mf{F}(\mf{x}_0+h\mf{v})-\mf{F}(\mf{x}_0)}{h}=\mf{g}'(0),\]
where  $\mf{g}:(-r,r)\to\mb{R}^m$ is the function $\mf{g}(h)=\mf{F}(\mf{x}_0+h\mf{v})$. Let $\boldsymbol{\gamma} :(-r,r)\to\mb{R}^n$ be the function defined as $\boldsymbol{\gamma} (h)=\mf{x}_0+h\mf{v}$. Then $\boldsymbol{\gamma}$ is a differentiable function with $\boldsymbol{\gamma}'(h)=\mf{v}$. Since $\mf{g}=\mf{F}\circ \boldsymbol{\gamma}$, and $\boldsymbol{\gamma}(0)=\mf{x}_0$, the chain rule implies that $\mf{g}$ is differentiable at $h=0$ and 
\[\mf{g}'(0)=\mf{DF}(\mf{x}_0)\boldsymbol{\gamma}'(0)=\mf{DF}(\mf{x}_0)\mf{v}.\]
This completes the proof.
\end{remark}

\begin{definition}{Tangent Line to a Curve}
A curve in $\mb{R}^{n}$ is  a continuous function $\boldsymbol{\gamma}:[a,b]\to \mb{R}^n$.  Let $c_0$   be a point in $(a, b)$. If  the curve $\boldsymbol{\gamma}$ is differentiable at $c_0$, the tangent vector to the curve $\boldsymbol{\gamma}$ at the point $\boldsymbol{\gamma}(c_0)$ is the vector $\boldsymbol{\gamma}'(c_0)$ in $\mb{R}^n$, while the tangent line to the curve $\boldsymbol{\gamma}$ at the point $\boldsymbol{\gamma}(c_0)$ is the line in $\mb{R}^n$ given by $\mf{x}:\mb{R}\to\mb{R}^n$,
\[\mathbf{x}(t)=\boldsymbol{\gamma}(c_0)+t\boldsymbol{\gamma}'(c_0).\]
\end{definition}

\begin{remark}{Tangent Lines and Tangent Planes}

Let $\mathcal{O}$ be an open subset of $\mb{R}^n$ that contains the point $\mf{x}_0$, and let $f:\mathcal{O}\to\mb{R}$ be a function that is differentiable at $\mf{x}_0$. We have seen that the tangent plane to the graph of $f$ at the point $(\mf{x}_0, f(\mf{x}_0))$ has equation
\end{remark}
\begin{highlight}{}
\[x_{n+1}=f(\mf{x}_0)+\langle \nabla f(\mf{x}_0), \mf{x}-\mf{x}_0\rangle.\]

Now assume that $r>0$ and $\boldsymbol{\gamma}:(-r, r)\to \mb{R}^{n+1}$ is a differentiable curve in $\mb{R}^{n+1}$ that lies on the graph of $f$, and $\gamma(0)=(\mf{x}_0, f(\mf{x}_0))$. For all $t\in (-r, r)$,
\[\gamma_{n+1}(t)=f(\gamma_1(t), \ldots, \gamma_n(t)).\]  
By chain rule, we find that
\[\gamma_{n+1}'(0)=\langle \nabla f(\mf{x}_0), \mf{v}\rangle, \hspace{1cm}\text{where}\;\mf{v}=(\gamma_1'(0), \ldots, \gamma_n'(0)).\]
The vector $\mf{w}=(\mf{v}, \gamma_{n+1}'(0))$ is the tangent vector to the curve $\boldsymbol{\gamma}$  at the point $(\mf{x}_0, f(\mf{x}_0))$. The equation of the tangent line is 
\[(x_1(t), \ldots, x_n(t), x_{n+1}(t))=(\mf{x}_0, f(\mf{x}_0))+t(\gamma_1'(0), \ldots, \gamma_n'(0), \gamma_{n+1}'(0)).\]
Thus, we find that 
\[(x_1(t), \ldots, x_n(t))=\mf{x}(t)=\mf{x}_0+t\mf{v},\]and 
\[x_{n+1}(t) =f(\mf{x}_0)+t\gamma_{n+1}'(0).\]These imply that
\begin{align*}x_{n+1}(t) 
&=f(\mf{x}_0)+t\langle \nabla f(\mf{x}_0), \mf{v}\rangle\\&=f(\mf{x}_0)+\langle \nabla f(\mf{x}_0), \mf{x}(t)-\mf{x}_0\rangle.\end{align*}
Thus, the tangent line to the curve $\boldsymbol{\gamma}$ lies in the tangent plane.

In fact, the tangent plane to the graph of a function $f$ at a point   can be characterized as the unique plane that contains all the tangent lines to the differentiable curves  that lie on the graph and passing through that point.
\end{highlight}

Now we turn to the mean value theorem. 
For a single variable function, the mean value theorem says that given  that $f:I\to\mb{R}$ is a differentiable function defined on the open interval $I$,    if $x_0$ and $x_0+h$ are two points in $I$,    there exists $c\in (0,1)$ such that
\[f(x_0+h)-f(x_0)=hf'(x_0+ch).\] 
Notice that the point $\di x_0+ch$ is a point strictly in between $x_0$ and $x_0+h$.
To generalize this theorem to multivariable functions, one natural question to ask is the following. If $\mf{F}:\mathcal{O}\to\mb{R}^m$ is a differentiable function defined on the open subset $\mathcal{O}$ of $\mb{R}^n$, $\mf{x}_0$ and $\mf{x}_0+\mf{h}$ are points in $\mathcal{O}$ such that the line segment between them lies entirely in $\mathcal{O}$, does there exist a constant $c\in (0,1)$ such that 
\[\mf{F}(\mf{x}_0+\mf{h})-\mf{F}(\mf{x}_0)=\mf{DF}(\mf{x}_0+c \mf{h})\mf{h}?\]
When $m\geq 2$,
the answer is no in general.
Let us  look at the following example.
\begin{example}{}
Consider the function $\mf{F}:\mb{R}^2\to \mb{R}^2$ defined as
\[\mf{F}(x,y)=(x^2y, xy).\]
Show that there does not exist  a contant $c\in (0,1)$ such that
\[\mf{F}(\mf{x}_0+\mf{h})-\mf{F}(\mf{x}_0)=\mf{DF}(\mf{x}_0+c \mf{h})\mf{h},\]when $\mf{x}_0=(0,0)$ and  $\mf{h}=(1,1)$.
\end{example}
\begin{solution}{Solution}
Notice that
\[\mf{DF}(x,y)=\begin{bmatrix} 2xy & x^2\\y& x\end{bmatrix}.\]When $\mf{x}_0=(0,0)$ and $\mf{h}=(1,1)$, $\mf{x}_0+c\mf{h}=(c,c)$.
If there exists a constant $c\in (0,1)$ such that 
\[\mf{F}(\mf{x}_0+\mf{h})-\mf{F}(\mf{x}_0)=\mf{DF}(\mf{x}_0+c \mf{h})\mf{h},\]
\bs
then
\[
\begin{bmatrix}1\\1\end{bmatrix}=\begin{bmatrix}2c^2 & c^2\\c & c\end{bmatrix}\begin{bmatrix}1\\1\end{bmatrix}.
\]This gives
\[3c^2=1\quad\text{and}\quad 2c=1.\]
But $2c=1$ gives $c=1/2$. When $c=1/2$, $3c^2=3/4\neq 1$. Hence, no such $c$ can exist.
\end{solution}

However, when $m=1$, we indeed have a mean value theorem.
\begin{theorem}{The Mean Value Theorem}
Let $\mathcal{O}$ be an open subset of $\mb{R}^n$, and let $\mf{x}_0$ and $\mf{x}_0+\mf{h}$ be two points in $\mathcal{O}$ such that the line segment between them lies entirely in $\mathcal{O}$. If 
 $f:\mathcal{O}\to\mb{R}$ is a differentiable function, there exist a constant $c\in (0,1)$ such that 
\[f(\mf{x}_0+\mf{h})-f(\mf{x}_0)=\langle \nabla f(\mf{x}_0+c\mf{h}), \mf{h}\rangle =\sum_{i=1}^n h_i\frac{\pa f}{\pa x_i}(\mf{x}_0+c\mf{h}).\]
\end{theorem}
\begin{myproof}{Proof}
Define the function $\gamma:[0,1]\to \mb{R}$ by  $\gamma(t)=\mf{x}_0+t\mf{h}$. Then $\gamma$ is a differentiable function with $\gamma'(t)=\mf{h}$. Let $g=(f\circ \gamma):[0,1]\to \mb{R}$. Then
\[g(t)=(f\circ \gamma)(t)=f(\mf{x}_0+t\mf{h}).\]Since $f$ and $\gamma$ are differentiable, the chain rule implies that $g$ is also differentiable and
\[g'(t)=\langle\nabla f(\mf{x}_0+t\mf{h}), \gamma'(t)\rangle=\langle\nabla f(\mf{x}_0+t\mf{h}),\mf{h}\rangle.\]
By mean value theorem for single variable functions, we find that there exists $c\in (0,1)$ such that
\[g(1)-g(0)=g'(c).\]
\bp
In other words, there exists $c\in (0,1)$ such that
\[f(\mf{x}_0+\mf{h})-f(\mf{x}_0)=\langle \nabla f(\mf{x}_0+c\mf{h}), \mf{h}\rangle.\]This completes the proof.
\end{myproof}
As in the single variable case, the mean value theorem has the following application.
\begin{corollary}{}
Let $\mathcal{O}$ be an open connected subset of $\mb{R}^n$, and let $f:\mathcal{O}\to\mb{R}$ be a function defined on $\mathcal{O}$. If $f$ is differentiable and $\nabla f(\mf{x})=\mf{0}$ for all $\mf{x}\in\mathcal{O}$, then $f$ is a constant function.
\end{corollary}
\begin{myproof}{Proof}
If $\mf{u}$ and $\mf{v}$ are two points in $\mathcal{O}$ such that the line segment between them lies entirely in $\mathcal{O}$, then the mean value theorem implies that $f(\mf{u})=f(\mf{v})$. 

Since $\mathcal{O}$ is an open connected subset of $\mb{R}^n$, Theorem \ref{230802_2} says that any two points $\mf{u}$ and $\mf{v}$  in $\mathcal{O}$ can be joined by a polygonal path in $\mathcal{O}$. In other words, there are points $\mf{x}_0, \mf{x}_1, \ldots, \mf{x}_k$ in $\mathcal{O}$ such that $\mf{x}_0=\mf{u}$, $\mf{x}_k=\mf{v}$, and for $1\leq i\leq k$, the line segment between $\mf{x}_{i-1}$ and $\mf{x}_i$ lies entirely in $\mathcal{O}$. Therefore,
\[f(\mf{x}_{i-1})=f(\mf{x}_i)\hspace{1cm}\text{for all}\; 1\leq i\leq k.\]
This proves that $f(\mf{u})=f(\mf{v})$.
Hence, $f$ is a constant function.
\end{myproof}

\vp
\noindent
{\bf \large Exercises  \thesection}
\setcounter{myquestion}{1}
\begin{question}{\themyquestion}
Let  $\mf{F}:\mb{R}^2\to\mb{R}^3$ be the function defined as
\[\mf{F}(x,y)=(x^2+ y^2, xy, x+y).\] Find a vector $\mf{b}$ in $\mb{R}^3$ and a $3\times 2$ matrix $A$ such that
\[\lim_{(u, v)\to (1, -1)}\frac{\mf{F}(5u+3v,  u-2v)-\mf{b}-A\mf{w}}{\sqrt{(u-1)^2+(v+1)^2 }}=\mf{0},\hspace{0.2cm}\text{where}\;\mf{w}=\di\begin{bmatrix} u\\v \end{bmatrix}.\]
\end{question}
 \atc
 \begin{question}{\themyquestion}
Let  $\phi:\mb{R} \to\mb{R}$ and $\psi:\mb{R} \to\mb{R}$ be functions that have continuous second order  derivatives, and let $c$ be a constant. Define the function $f:\mb{R}^2\to\mb{R}$ by 
\[f(t,x)=\phi(x+ct)+\psi(x-ct).\]
Show that
\[\frac{\pa^2 f}{\pa t^2}-c^2\frac{\pa^2f}{\pa x^2}=0.\]
\end{question}

 \atc
 \begin{question}{\themyquestion}
Let $\alpha$ be a constant, and let $f:\mb{R}^n\setminus \{\mf{0}\} \to\mb{R}$ be the function defined by
\[f(\mf{x})=\Vert\mf{x}\Vert^{\alpha}.\]
Find the value(s) of $\alpha$ such that
\[\Delta f(\mf{x})=\sum_{i=1}^n\frac{\pa^2 f}{\pa x_i^2}(\mf{x})=\frac{\pa^2 f}{\pa x_1^2}(\mf{x})+\frac{\pa^2 f}{\pa x_2^2}(\mf{x})+\cdots+\frac{\pa^2f}{\pa x_n^2}(\mf{x})=0.\]
\end{question}
 \atc
 \begin{question}{\themyquestion}
Let $f:\mb{R}^2\to\mb{R}$ be a function such that $f(0,0)=2$ and 
\[\frac{\pa f}{\pa x}(x,y)=11\quad\text{and}\quad\frac{\pa f}{\pa y}=-7\hspace{1cm}\text{for all}\;(x,y)\in\mb{R}^2.\]
Show that
\[f(x,y)=2+11x-7y\hspace{1cm}\text{for all}\;(x,y)\in\mb{R}^2.\]
\end{question}

 \atc
 \begin{question}{\themyquestion}
Let $\mathcal{O}$ be an open subset of $\mb{R}^2$, and let $u:\mathcal{O}\to\mb{R}$ and $v:\mathcal{O}\to\mb{R}$ be twice continuously differentiable functions. Define the function $\mf{F}:\mathcal{O}\to \mb{R}^2$ by
\[\mf{F}(x,y)=(u(x,y), v(x,y)).\]
Let $\mathcal{U}$ be an open subset of $\mb{R}^2$ that contains $\mf{F}(\mathcal{O})$, and let $f:\mathcal{U}\to\mb{R}$ be a twice continuously differentiable function. Define the function $g:\mathcal{O}\to\mb{R}$ by
\[g(x,y)=(f\circ\mf{F})(x,y)=f(u(x,y), v(x,y)).\]
Find $g_{xx}, g_{xy}$ and $g_{yy}$ in terms of the first and second order partial derivatives of $u, v$ and $f$.

\end{question}
 
\section{Second Order Approximations} 
In this section, we turn to consider second order approximations. We  only consider a function $f:\mathcal{O}\to\mb{R}$ defined on an open subset $\mathcal{O}$ of $\mb{R}^n$ and whose codomain is $\mb{R}$. The function is said to be twice differentiable if it has first order partial derivatives, and each $f_{x_i}:\mathcal{O}\to\mb{R}$, $1\leq i\leq n$, is   a differentiable function. Notice that a twice differentiable function has continuous first order partial derivatives. Hence, it is  differentiable. The differentiability of each $f_{x_i}$, $1\leq i\leq n$ also implies that $f$ has second order partial derivatives.

\begin{lemma}[label=230803_1]{}
Let $\mathcal{O}$ be an open subset of $\mb{R}^n$, and let $f:\mathcal{O}\to\mb{R}$ be a twice differentiable function defined on $\mathcal{O}$.  If  $\mf{x}_0$ and $\mf{x}_0+\mf{h}$ are two points in $\mathcal{O}$ such that the line segment between them lies entirely in $\mathcal{O}$, then there is a $c\in (0,1)$ such that
\begin{align*}
f(\mf{x}_0+\mf{h})-f(\mf{x}_0)-\langle\nabla f(\mf{x}_0), \mf{h}\rangle &=\frac{1}{2}\mf{h}^TH_f(\mf{x}_0+c\mf{h})\mf{h}\\&=\frac{1}{2}\sum_{i=1}^n\sum_{j=1}^n h_ih_j\frac{\pa^2 f}{\pa x_j\pa x_i}(\mf{x}_0+c\mf{h}).\end{align*}
\end{lemma}

\begin{myproof}{Proof}
Given $\mf{x}_0\in\mathcal{O}$, let $r$ be a positive number such that  $B(\mf{x}_0, r)\subset \mathcal{O}$. Define the function $g:(-r,r)\to\mb{R}$ by
\[g(t)=f(\mf{x}_0+t\mf{h}).\]

 Since $f:\mathcal{O}\to\mb{R}$ is differentiable, chain rule implies that $g:(-r,r)\to\mb{R}$ is differentiable and 
\[g'(t)=\sum_{i=1}^n h_i\frac{\pa f}{\pa x_i}(\mf{x}_0+t\mf{h})=\langle\nabla f(\mf{x}_0+t\mf{h}), \mf{h}\rangle.\]
 
Since each $f_{x_i}:\mathcal{O}\to\mb{R}$, $1\leq i\leq n$  is differentiable,
chain rule again implies that $g'$ is differentiable and
\bp
\[g''(t)=\sum_{i=1}^n \sum_{j=1}^nh_ih_j\frac{\pa f}{\pa x_j\pa x_i}(\mf{x}_0+t\mf{h})=\mf{h}^TH_f(\mf{x}_0+t\mf{h})\mf{h}.\]

By Lagrange's remainder theorem, there is a $c\in (0,1)$ such that
\[g(1)-g(0)-g'(0)(1-0)=\frac{g''(c)}{2}(1-0)^2.\]
This gives 
\[
f(\mf{x}_0+\mf{h})-f(\mf{x}_0)-\langle\nabla f(\mf{x}_0), \mf{h}\rangle = \frac{1}{2}\sum_{i=1}^n\sum_{j=1}^n h_ih_j\frac{\pa^2 f}{\pa x_j\pa x_i}(\mf{x}_0+c\mf{h}).\]
\end{myproof}

If a function has continuous second order partial derivatives, then it is twice differentiable, and  Clairaut's theorem implies that its Hessian matrix is symmetric.  For such a function, we can prove the second order approximation theorem.
\begin{theorem}{Second Order Approximation Theorem}
Let $\mathcal{O}$ be an open subset of $\mb{R}^n$ that contains the point $\mf{x}_0$, and let $f:\mathcal{O}\to\mb{R}$ be a twice continuously differentiable function defined on $\mathcal{O}$. We have the followings.
\begin{enumerate}[(a)]
\item
$\di \lim_{\mf{h}\to\mf{0}}\frac{f(\mf{x}_0+\mf{h})-f(\mf{x}_0)-\di \langle\nabla f(\mf{x}_0), \mf{h}\rangle -\frac{1}{2}\mf{h}^TH_f(\mf{x}_0 )\mf{h}}{\Vert\mf{h}\Vert^2}=0$.
\item
If $Q(\mf{x})$ is a polynomial of degree at most two such that 
\[\lim_{\mf{h}\to\mf{0}}\frac{f(\mf{x}_0+\mf{h})-Q(\mf{x}_0+\mf{h}) }{\Vert\mf{h}\Vert^2}=0,\]then
\begin{equation}\label{230805_1}Q(\mf{x})= f(\mf{x}_0)+\di \langle\nabla f(\mf{x}_0), \mf{x}-\mf{x}_0\rangle +\frac{1}{2}(\mf{x}-\mf{x}_0)^TH_f(\mf{x}_0)(\mf{x}-\mf{x}_0).\end{equation}
\end{enumerate}
\end{theorem}Combining (a) and (b), the second order approximation theorem says that for a twice continuously differentiable function, there exists a  unique polynomial of degree at most 2 which is a second order approximation of the function.

\begin{myproof}{Proof}Let us prove part (a) first. Since $\mathcal{O}$ is open,
there is an $r>0$ such that $B(\mf{x}_0, r)\subset\mathcal{O}$. For each $\mf{h}$ in $\mb{R}^n$ with $\Vert\mf{h}\Vert<r$, Lemma \ref{230803_1} says that there is a $c_{\mf{h}}\in (0,1)$ such that
\[f(\mf{x}_0+\mf{h})-f(\mf{x}_0)-\di \langle\nabla f(\mf{x}_0), \mf{h}\rangle =\frac{1}{2} \mf{h}^TH_f(\mf{x}_0+c\mf{h})\mf{h}.\]
 
Therefore, if $0<\Vert\mf{h}\Vert<r$, 
\begin{align*}&\left|\frac{f(\mf{x}_0+\mf{h})-f(\mf{x}_0)-\di \langle\nabla f(\mf{x}_0), \mf{h}\rangle -\frac{1}{2}\mf{h}^TH_f(\mf{x}_0)\mf{h}}{\Vert\mf{h}\Vert^2}\right|\\
&=\frac{1}{2}\left|\sum_{i=1}^n\sum_{j=1}^n\frac{ h_ih_j}{\Vert\mf{h}\Vert^2}\left(\frac{\pa^2 f}{\pa x_j\pa x_i}(\mf{x}_0 +c_{\mf{h}}\mf{h})-\frac{\pa^2 f}{\pa x_j\pa x_i}(\mf{x}_0) \right)\right|\\
&\leq \frac{1}{2}\sum_{i=1}^n\sum_{j=1}^n\frac{ |h_i||h_j|}{\Vert\mf{h}\Vert^2}\left|\frac{\pa^2 f}{\pa x_j\pa x_i}(\mf{x}_0 +c_{\mf{h}}\mf{h})-\frac{\pa^2 f}{\pa x_j\pa x_i}(\mf{x}_0)  \right|
\\
&\leq \frac{1}{2}\sum_{i=1}^n\sum_{j=1}^n\left|\frac{\pa^2 f}{\pa x_j\pa x_i}(\mf{x}_0 +c_{\mf{h}}\mf{h})-\frac{\pa^2 f}{\pa x_j\pa x_i}(\mf{x}_0)  \right|.
\end{align*}Since $c_{\mf{h}}\in (0,1)$,  $\di \lim_{\mf{h}\to\mf{0}}\left(\mf{x}_0+c_{\mf{h}}\mf{h}\right)=\mf{x}_0$. For all $1\leq i\leq n$, $1\leq j\leq n$, $f_{x_j x_i}$ is continuous. Hence, 
\[\lim_{\mf{h}\to\mf{0}}\frac{\pa^2 f}{\pa x_j\pa x_i}(\mf{x}_0 +c_{\mf{h}}\mf{h})=\frac{\pa^2 f}{\pa x_j\pa x_i}(\mf{x}_0).\]This proves that
\[\lim_{\mf{h}\to\mf{0}}\frac{f(\mf{x}_0+\mf{h})-f(\mf{x}_0)-\di \langle\nabla f(\mf{x}_0), \mf{h}\rangle-\frac{1}{2}\mf{h}^TH_f(\mf{x}_0)\mf{h}}{\Vert\mf{h}\Vert^2}=0.\]
To prove part (b), let
\[P(\mf{x})=f(\mf{x}_0)+\di \langle\nabla f(\mf{x}_0), \mf{x}-\mf{x}_0\rangle +\frac{1}{2}(\mf{x}-\mf{x}_0)^TH_f(\mf{x}_0)(\mf{x}-\mf{x}_0).\]
Part (a) says that
\begin{equation}\label{230805_2}\lim_{\mf{h}\to\mf{0}}\frac{f(\mf{x}_0+\mf{h})-P(\mf{x}_0+\mf{h})}{\Vert \mf{h}\Vert^2}=0.\end{equation}
\bp
Since $Q(\mf{x})$ is a polynomial of degree at most two in $\mf{x}$, $Q(\mf{x}_0+\mf{h})$ is a  polynomial of degree at most two in $\mf{h}$. Therefore, we can write $Q(\mf{x}_0+\mf{h})$ as
\[Q(\mf{x}_0+\mf{h})=c+\sum_{i=1}^n b_ih_i+\frac{1}{2}\sum_{i=1}^na_{ii}h_i^2+ \sum_{1\leq i< j\leq n}  a_{ij}h_ih_j.\]
Since
\[\lim_{\mf{h}\to\mf{0}}\frac{f(\mf{x}_0+\mf{h})-Q(\mf{x}_0+\mf{h}) }{\Vert\mf{h}\Vert^2}=0,\]subtracting \eqref{230805_2} gives
\begin{equation}\label{230803_5}\lim_{\mf{h}\to\mf{0}}\frac{P(\mf{x}_0+\mf{h})-Q(\mf{x}_0+\mf{h}) }{\Vert\mf{h}\Vert^2}=0.\end{equation}
It follows that
\begin{equation}\label{230803_3}
\lim_{\mf{h}\to\mf{0}}\left(P(\mf{x}_0+\mf{h})-Q(\mf{x}_0+\mf{h}) \right)=0,
\end{equation}and
\begin{equation}\label{230803_4}
\lim_{\mf{h}\to\mf{0}}\frac{P(\mf{x}_0+\mf{h})-Q(\mf{x}_0+\mf{h}) }{\Vert\mf{h}\Vert}=0.
\end{equation}Since $f$ has continuous second order partial derivatives, $f_{x_jx_i}(\mf{x}_0)=f_{x_ix_j}(\mf{x}_0)$. Thus,
\begin{align*}
&P(\mf{x}_0+\mf{h})-Q(\mf{x}_0+\mf{h})\\&=(f(\mf{x}_0)-c)+\sum_{i=1}^n  h_i\left(\frac{\pa f}{\pa x_i}(\mf{x}_0)-b_i\right)\\&\quad +\frac{1}{2}\sum_{i=1}^nh_i^2\left(\frac{\pa^2f}{\pa x_i^2}(\mf{x}_0)-a_{ii}\right)+\sum_{1\leq i< j\leq n}  h_ih_j\left(\frac{\pa^2 f}{\pa x_j\pa x_i}(\mf{x}_0)-a_{ij}\right).\end{align*}
Eq. \eqref{230803_3} implies that $c=f(\mf{x}_0)$. Then eq. \eqref{230803_4} implies that 
\[b_i=\frac{\pa f}{\pa x_i}(\mf{x}_0)\hspace{1cm}\text{for all}\;1\leq i\leq n.\]
Finally, \eqref{230803_5} implies that for any $1\leq i\leq   j\leq n$,
\[a_{ij} =\frac{\pa^2 f}{\pa x_i\pa x_j}(\mf{x}_0).\]
This completes the proof that
$Q(\mf{x})=P(\mf{x})$.
\end{myproof}

\begin{example}[label=230803_6]{}
Find a polynomial $Q(x,y)$ of degree at most 2 such that 
\[\lim_{(x,y)\to (1,2)}\frac{\sin(4x^2-y^2)- Q(x,y)}{(x-1)^2+(y-2)^2}=0.\]
\end{example}
\begin{solution}{Solution}
Since $g(x,y)=4x^2-y^2$ is a polynomial function, it is infinitely differentiable. Since the sine function is also infinitely differentiable, the function $f(x,y)=\sin(4x^2-y^2)$ is  infinitely differentiable.
\begin{gather*}
f_x(x,y)=8x\cos(4x^2-y^2), \hspace{1cm}f_y(x,y)=-2y\cos(4x^2-y^2),\\
f_{xx}(x,y)=8\cos(4x^2-y^2)-64x^2\sin(4x^2-y^2),\\
f_{xy}(x,y)=f_{yx}(x,y)=16xy\sin(4x^2-y^2),\\
f_{yy}(x,y)=-2\cos(4x^2-y^2)-4y^2\sin(4x^2-y^2).
\end{gather*}Hence,
\[f(1,2)=0,\;\;f_{x}(1,2)=8, \;\;f_y(1,2)=-4,\]
\[f_{xx}(1,2)=8,\;\;f_{xy}(1,2)=0,\;\;f_{yy}(1,2)=-2.\]
By the second order approximation theorem,
\begin{align*}
Q(x,y)&=f(1,2)+f_x(1,2)(x-1)+f_y(1,2)(y-2)+\frac{1}{2}f_{xx}(1,2)(x-1)^2\\&\quad+f_{xy}(1,2)(x-1)(y-2)+\frac{1}{2}f_{yy}(1,2)(y-2)^2\\
&=8(x-1)-4(y-2)+4(x-1)^2-(y-2)^2\\
&=4x^2-y^2.
\end{align*}

\end{solution}
 \begin{example}{}
 Determine whether the limit $\di\lim_{(x,y)\to (0,0)}\frac{e^{x+y}-1-x-y}{x^2+y^2}$ exists. If yes, find the limit.
 \end{example}
 \begin{solution}{Solution}
 Since the exponential funtion  and the function $g(x,y)=x+y$ are infinitely differentiable,   the function $f(x,y)=e^{x+y}$ is  infinitely differentiable. 
 By the second order approximation theorem,
\[\lim_{(x,y)\to (0,0)}\frac{f(x,y)-Q(x,y)}{x^2+y^2}=0,\]where
  \begin{align*}Q(x,y)&=f(0,0)+x\frac{\pa f}{\pa x}(0,0)+y\frac{\pa f}{\pa y}(0,0)\\&\quad +\frac{1}{2} x^2 \frac{\pa^2 f}{\pa x^2}(0,0)
 +xy\frac{\pa^2 f}{\pa x\pa y}(0,0)+\frac{1}{2}y^2\frac{\pa^2 f}{\pa y^2}(0,0).
 \end{align*}
 Now
 \[
\frac{\pa f}{\pa x}(x,y)= \frac{\pa f}{\pa y}(x,y)=
\frac{\pa^2f}{\pa x^2}(x,y)=\frac{\pa^2 f}{\pa x\pa y}(x,y)=\frac{\pa^2 f}{\pa y^2}(x,y)=e^{x+y}.
 \]
 Thus,
 \[f(0,0)=\frac{\pa f}{\pa x}(0,0)=\frac{\pa f}{\pa y}(0,0)=\frac{\pa^2f}{\pa x^2}(0,0)=\frac{\pa^2 f}{\pa x\pa y}(0,0)=\frac{\pa^2 f}{\pa y^2}(0,0)=1.\] 
  It follows that
 \[Q(x,y) =1+x+y+\frac{1}{2}x^2+xy+\frac{1}{2}y^2.
 \]
 Hence,
\begin{equation}\label{230805_3}\lim_{(x,y)\to(0,0)}\frac{\di e^{x+y}-1-x-y-\frac{1}{2}x^2-xy-\frac{1}{2}y^2}{x^2+y^2}=0.\end{equation}
 If \[\lim_{(x,y)\to (0,0)}\frac{e^{x+y}-1-x-y}{x^2+y^2}=a\] exists, subtracting \eqref{230805_3} shows that  
 \[a=\lim_{(x,y)\to (0,0)}h(x,y),\hspace{1cm}\text{where}\;\;h(x,y)=\frac{\di\frac{1}{2}x^2+xy+\frac{1}{2}y^2}{x^2+y^2}.\]    
 \bs
 This implies that  if $\{\mf{w}_k\}$ is a sequence in $\mb{R}^2\setminus\{\mf{0}\}$ that converges to $(0,0)$, then the sequence $\{h(\mf{w}_k)\}$   converges to $a$. 
 For $k\in\mb{Z}^+$, let
 \[\mf{u}_k=\left(\frac{1}{k}, 0\right),\quad \mf{v}_k=\left(\frac{1}{k}, \frac{1}{k}\right).\]
 Then   $\{\mf{u}_k\}$ and $\{\mf{v}_k\}$  are sequences in $\mb{R}^2\setminus\{\mf{0}\}$ that converge to $(0,0)$. Hence, the sequences $\{h(\mf{u}_k)\}$ and $\{h(\mf{v}_k)\}$ both converge to $a$. Since
 \[h(\mf{u}_k)=\frac{1}{2},\quad h(\mf{v}_k)=1\hspace{1cm}\text{for all}\;k\in\mb{Z}^+,\]
 the sequence $\{h(\mf{u}_k)\}$ converges to $\frac{1}{2}$, while the sequence $\{h(\mf{v}_k)\}$ converges to 1. This gives a contradiction. Hence, the limit \[\lim_{(x,y)\to (0,0)}\frac{e^{x+y}-1-x-y}{x^2+y^2}\] does not exist.
 \end{solution}

\vp
\noindent
{\bf \large Exercises  \thesection}
\setcounter{myquestion}{1}
\begin{question}{\themyquestion}
Let $f:\mb{R}^2\to \mb{R}$ be the function 
\[f(x,y)=x^2y+4xy^2.\]
Find  a polynomial $Q(x,y)$ of degree at most 2 such that
\[\lim_{(x,y)\to (1, -1)}\frac{f(x,y)-Q(x,y)}{(x-1)^2+(y+1)^2}=0.\]
\end{question}
 \atc
 \begin{question}{\themyquestion}
  Determine whether the limit $\di\lim_{(x,y)\to (0,0)}\frac{\sin(x+y)-x-y}{x^2+y^2}$ exists.  If yes, find the limit.
\end{question}
  \atc
 \begin{question}{\themyquestion}
 Determine whether the limit $\di\lim_{(x,y)\to (0,0)}\frac{\cos(x+y)-1}{x^2+y^2}$ exists.  If yes, find the limit.
\end{question}
\section{Local Extrema} 

In this section, we use differential calculus to study local extrema of a  function $f:\mathcal{O}\to \mb{R}$ that is defined on an open subset $\mathcal{O}$ of $\mb{R}^n$. The definition of local extrema that we give here is only restricted to such functions. 

\begin{definition}{Local Maximum and Local Minimum}
Let $\mathcal{O}$ be an open subset of $\mb{R}^n$ that contains the point $\mf{x}_0$, and let $f:\mathcal{O}\to \mb{R}$ be a function defined on $\mathcal{O}$.

\begin{enumerate}[1.]
\item
The point $\mf{x}_0$ is called a {\it local maximizer} of $f$ provided that there is a $\delta>0$ such that $B(\mf{x}_0, \delta)\subset \mathcal{O}$ and for all $\mf{x}\in B(\mf{x}_0, \delta)$, 
\[f(\mf{x})\leq f(\mf{x}_0).\]
The value $f(\mf{x}_0)$ is called a local maximum value of $f$.

\item
The point $\mf{x}_0$ is called a {\it local minimizer} of $f$ provided that there is a $\delta>0$ such that $B(\mf{x}_0, \delta)\subset \mathcal{O}$ and for all $\mf{x}\in B(\mf{x}_0, \delta)$, 
\[f(\mf{x})\geq f(\mf{x}_0).\]
The value $f(\mf{x}_0)$ is called a local minimum value of $f$.
\item The point $\mf{x}_0$ is called a local extremizer if it is either a local maximizer or a local minimizer. The value $f(\mf{x}_0)$ is called a local extreme value if it is either a local maximum value or a local minimum value.
\end{enumerate}
 
\end{definition}
From the definition, it is obvious that  $\mf{x}_0$ is a local minimizer of the function $f:\mathcal{O}\to\mb{R}$ if and only if it is a local maximizer of the function $-f:\mathcal{O}\to \mb{R}$.

\begin{example}[label=230805_6]{}
\begin{enumerate}[(a)]
\item
For the function $f:\mb{R}^2\to \mb{R}$, $f(x,y)=x^2+y^2$, $(0,0)$ is a local minimizer.
\item
For the function $g:\mb{R}^2\to \mb{R}$, $g(x,y)=-x^2-y^2$, $(0,0)$ is a local maximizer.\end{enumerate}
\be
\begin{enumerate}[(a)]
\item[(c)]
For the function $h:\mb{R}^2\to \mb{R}$, $h(x,y)=x^2-y^2$, $\mf{0}=(0,0)$ is neither a local maximizer nor a local minimizer.  For any $\delta>0$,  let $r=\delta/2$. The points $\mf{u}=(r, 0)$ and $\mf{v}=(0, r)$ are in $B(\mf{0}, \delta)$, but 
\[h(\mf{u})=r^2>0=h(\mf{0}), \hspace{1cm} h(\mf{v})=-r^2<0=h(\mf{0}).\]
\end{enumerate}
\end{example2}
 
\begin{figure}[ht]
\centering
\includegraphics[scale=0.18]{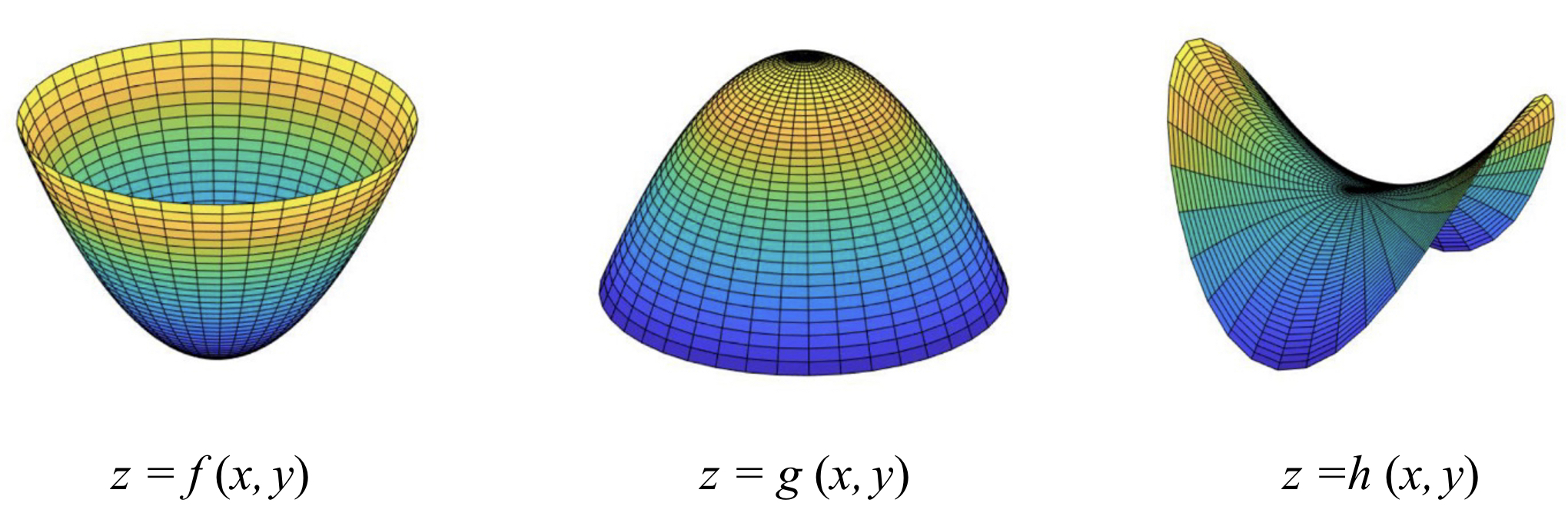}

\caption{The functions $f(x,y)$, $g(x,y)$ and $h(x,y)$ defined in Example \ref{230805_6}.}\label{figure44}
\end{figure}

The following theorem gives a necessary condition for a point to be a local extremum if the function has partial derivatives at that point.
\begin{theorem}[label=230805_8]{}
Let $\mathcal{O}$ be an open subset of $\mb{R}^n$ that contains the point $\mf{x}_0$, and let $f:\mathcal{O}\to \mb{R}$ be a function defined on $\mathcal{O}$.  If $\mf{x}_0$ is a local extremizer and $f$ has partial derivatives at $\mf{x}_0$, then the gradient of $f$ at $\mf{x}_0$ is the zero vector, namely, $\nabla f(\mf{x}_0)=\mf{0}$.
\end{theorem}
\begin{myproof}{Proof}
Without loss of generality, assume that $\mathbf{x}_0$ is a local minimizer. Then there is a $\delta>0$ such that $B(\mf{x}_0, \delta)\subset \mathcal{O}$ and\begin{equation}\label{230805_7} f(\mf{x})\geq f(\mf{x}_0)
\hspace{1cm}\text{ for all}\;\mf{x}\in B(\mf{x}_0, \delta).\end{equation} For $1\leq i\leq n$, consider the function $g_i:(-\delta, \delta)\to \mb{R}$ defined by
 $g_i(t)=f(\mf{x}_0+t\mf{e}_i)$. 
 By the definition of partial derivatives, $g_i$ is differentiable at $t=0$ and 
  \bp
 \[g_i'(0)=\frac{\pa f}{\pa x_i}(\mf{x}_0).\]

 Eq. \eqref{230805_7} implies that
 \[g_i(t)\geq g_i(0)\hspace{1cm}\text{for all}\;t\in (-\delta, \delta).\]
 In other words, $t=0$ is a local minimizer of the function $g_i:(-\delta, \delta)\to \mb{R}$.  From the theory of single variable analysis, we must have $g_i'(0)=0$. Hence,  $f_{x_i}(\mf{x}_0)=0$ for all $1\leq i\leq n$. This  proves that $\nabla f(\mf{x}_0)=\mf{0}$.
 
\end{myproof}

Theorem \ref{230805_8} prompts us to make the following definition.
\begin{definition}{Stationary Points}
Let $\mathcal{O}$ be an open subset of $\mb{R}^n$ that contains the point $\mf{x}_0$, and let $f:\mathcal{O}\to \mb{R}$ be a function defined on $\mathcal{O}$. If $f$ has partial derivatives at $\mf{x}_0$ and $\nabla f(\mf{x}_0)=\mf{0}$, we call $\mf{x}_0$ a stationary point of $f$. 
\end{definition}

Theorem \ref{230805_8} says that if $f:\mathcal{O}\to\mb{R}$  has partial derivatives at $\mf{x}_0$, a necessary condition for $\mathbf{x}_0$ to be a local extremizer is that it is a stationary point.

\begin{example}{}
For all the three functions $f$, $g$ and $h$ defined in  Example \ref{230805_6}, the point $\mf{0}=(0,0)$ is a stationary point. However, $\mf{0}$ is local minimizer of  $f$, a local maximizer of $g$, but neither a local maximizer nor a local minimizer of $h$. 
\end{example}

The behavior of the function $h(x,y)=x^2-y^2$ in Example \ref{230805_6} prompts us to make the following definition.
\begin{definition}{Saddle Points}
Let $\mathcal{O}$ be an open subset of $\mb{R}^n$ that contains the point $\mf{x}_0$, and let $f:\mathcal{O}\to \mb{R}$ be a function defined on $\mathcal{O}$. The point $\mf{x}_0$  is a saddle point of the function $f$ if it is a stationary point of $f$, but it is not a local extremizer. In other words, $\nabla f(\mf{x}_0)=\mf{0}$, but for any $\delta>0$, there exist $\mf{x}_1$ and $\mf{x}_2$ in $B(\mf{x}_0, \delta)\cap \mathcal{O}$ such that
\[f(\mf{x}_1)>f(\mf{x}_0)\quad\text{and}\quad f(\mf{x}_2)<f(\mf{x}_0).\]

\end{definition}
\begin{example}{}
 
$(0,0)$ is a saddle point of the function $h:\mb{R}^2\to\mb{R}$, $h(x,y)=x^2-y^2$.
\end{example}

By definition, if $\mf{x}_0$ is a stationary point of the function $f:\mathcal{O}\to\mb{R}$, then it is either a local maximizer, a local minimizer, or a saddle point. If $f:\mathcal{O}\to\mb{R}$ has continuous second order partial derivatives at $\mf{x}_0$, we can use the second derivative test to partially determine whether $\mf{x}_0$ is a local maximizer, a local minimizer, or a saddle point. 
When $n=1$, we have seen that a stationary point $x_0$ of a function $f$ is a local minimum if $f''(x_0)>0$. It is a local maximum if $f''(x_0)<0$. For multivariable functions, it is natural to expect that whether $\mf{x}_0$ is a local extremizer 
  depends on the definiteness of the Hessian matrix $H_f(\mf{x}_0)$. 
  
In Section \ref{sec_mvf}, we have discussed the classification of a symmetric matrix. It is either positive semi-definite, negative semi-definite or indefinite. Among the positive semi-definite ones, there are those that are positive definite.  Among the negative semi-definite matrices, there are those which are  negative definite.
  
\begin{theorem}{Second Derivative Test}
Let $\mathcal{O}$ be an open subset of $\mathbb{R}^n$, and let $f:\mathcal{O}\to\mb{R}$ be a twice continuously differentiable function defined on $\mathcal{O}$.  Assume that  $\mathbf{x}_0$ is a stationary point of $f:\mathcal{O}\to\mb{R}$.
\begin{enumerate}[(i)]
\item If $H_f(\mf{x}_0)$ is positive definite, then $\mf{x}_0$ is a local minimizer of $f$.
\item If $H_f(\mf{x}_0)$ is negative definite, then $\mf{x}_0$ is a local maximizer of $f$.
\item If $H_f(\mf{x}_0)$ is indefinite, then $\mf{x}_0$ is a saddle point.
\end{enumerate}
\end{theorem}
The cases that are not covered in the second derivative test are the cases where $H_f(\mf{x}_0)$ is positive semi-definite but not positive definite, or $H_f(\mf{x}_0)$ is negative semi-definite but not negative definite. These are the inconclusive cases.

\begin{myproof}{Proof of the Second Derivative Test}
Notice that (i) and (ii) are equivalent since $\mf{x}_0$ is a local minimizer of $f$ if and only if it is a local maximizer of $-f$, and $H_{-f}=-H_f$.  A symmetric matrix $A$ is positive definite if and  only if $-A$ is negative definite. 
Thus, we only need to prove  (i) and  (iii).

Since $\mf{x}_0$ is a stationary point, $\nabla f(\mf{x}_0)=\mf{0}$.
It follows from the second order approximation theorem that
\begin{equation}
\label{230805_9}\lim_{\mf{h}\to\mf{0}}\frac{f(\mf{x}_0+\mf{h})-f(\mf{x}_0)-\frac{1}{2}\mf{h}^TH_f(\mf{x}_0 )\mf{h}}{\Vert\mf{h}\Vert^2}=0.\end{equation}

To prove (i), asume that $H_f(\mf{x}_0)$ is positive definite. By Theorem \ref{230725_5}, there is a positive number $c$ such that
\[\mf{h}^TH_f(\mf{x}_0 )\mf{h}\geq c\Vert\mf{h}\Vert^2\hspace{1cm}\text{for all}\;\mf{h}\in\mb{R}^n.\]
Eq. \ref{230805_9} implies that there is a $\delta>0$ such that $B(\mf{x}_0, \delta)\subset \mathcal{O}$ and for all $\mf{h}$ with $0<\Vert\mf{h}\Vert<\delta$,
\[\left|\frac{f(\mf{x}_0+\mf{h})-f(\mf{x}_0)-\frac{1}{2}\mf{h}^TH_f(\mf{x}_0 )\mf{h}}{\Vert\mf{h}\Vert^2}\right|<\frac{c}{3}.\]
Therefore,
\[\left|f(\mf{x}_0+\mf{h})-f(\mf{x}_0)-\frac{1}{2}\mf{h}^TH_f(\mf{x}_0 )\mf{h}\right|\leq \frac{c}{3}\Vert\mf{h}\Vert^2\hspace{1cm} \text{for all}\;\Vert\mf{h}\Vert<\delta.\]
This implies that for all $\mf{h}$ with $\Vert\mf{h}\Vert<\delta$,
\[ f(\mf{x}_0+\mf{h})- f(\mf{x}_0)\geq \frac{1}{2}\mf{h}^TH_f(\mf{x}_0 )\mf{h}- \frac{c}{3}\Vert\mf{h}\Vert^2\geq \frac{c}{6}\Vert\mf{h}\Vert^2\geq 0.\]
Thus, $f(\mf{x})\geq f(\mf{x}_0)$ for all $\mf{x}\in B(\mf{x}_0, \delta)$. This shows that $\mf{x}_0$ is a local minimizer of $f$.
\bp
Now to prove (iii), assume that $H_f(\mf{x}_0)$ is indefinite. Then there exist  unit vectors $\mf{u}_1$ and $\mf{u}_2$ so that
\[\varepsilon_1=\mathbf{u}_1^TH_f(\mf{x}_0)\mf{u}_1<0,\hspace{1cm} \varepsilon_2=\mathbf{u}_2^TH_f(\mf{x}_0)\mf{u}_2>0.\]
Let $\varepsilon=\frac{1}{2}\min\{|\varepsilon_1|, \varepsilon_2\}$.
 Eq. \eqref{230805_9} implies that there is a $\delta_0>0$ such that $B(\mf{x}_0, \delta_0)\subset \mathcal{O}$ and for all $\mf{h}$ with $0< \Vert\mf{h}\Vert<\delta_0$,
\begin{equation}\label{230805_10}\left| f(\mf{x}_0+\mf{h})-f(\mf{x}_0)-\frac{1}{2}\mf{h}^TH_f(\mf{x}_0 )\mf{h} \right|<\varepsilon\Vert\mf{h}\Vert^2.\end{equation}
For any $\delta>0$, let $r=\frac{1}{2}\min\{\delta, \delta_0\}$. Then the points $\mf{x}_1=\mf{x}_0+r\mf{u}_1$ and $\mf{x}_2=\mf{x}_0+r\mf{u}_2$ are in the ball $B(\mf{x}_0, \delta)$ and the ball $B(\mf{x}_0, \delta_0)$. Eq. \eqref{230805_10} implies that for $i=1,2$, 
\[  -r^2\varepsilon\leq f(\mf{x}_0+r\mf{u}_i)-f(\mf{x}_0)-\frac{r^2}{2}\mf{u}_i^TH_f(\mf{x}_0 )\mf{u}_i  < r^2\varepsilon.\]
Therefore,
\[f(\mf{x}_0+r\mf{u}_1)-f(\mf{x}_0)<r^2\left(\frac{1}{2}\mf{u}_1^TH_f(\mf{x}_0 )\mf{u}_1+\varepsilon \right)= r^2\left(\frac{1}{2}\varepsilon_1+\varepsilon\right)\leq 0\]since $\varepsilon\leq -\frac{1}{2}\varepsilon_1$; while
\[f(\mf{x}_0+r\mf{u}_2)-f(\mf{x}_0)> r^2\left(\frac{1}{2}\mf{u}_2^TH_f(\mf{x}_0 )\mf{u}_2-\varepsilon \right)= r^2\left(\frac{1}{2}\varepsilon_2-\varepsilon\right)\geq 0\]since $\varepsilon\leq  \frac{1}{2}\varepsilon_2$. Thus, $\mf{x}_1$ and $\mf{x}_2$ are points in $B(\mf{x}_0,\delta)$, but $f(\mf{x}_1)<f(\mf{x}_0)$ while $f(\mf{x}_2)>f(\mf{x}_0)$. These show that $\mf{x}_0$ is a saddle point.
\end{myproof}

A symmetric matrix is positive definite if and only if all its eigenvalues are positive. It is negative definite if and only if all its eigenvalues are negative. It is indefinite if it has at least one positive eigenvalue, and at least one negative eigenvalue. For a diagonal matrix, its eigenvalues are the entries on the diagonal.

Let us revisit Example \ref{230805_6}.

\begin{example}{}
For the functions considered in Example \ref{230805_6}, we have seen that $(0,0)$ is a stationary point of each of them.
Notice that $H_f(0,0)=\di\begin{bmatrix}2 & 0\\0 & 2\end{bmatrix}$ is positive definite, $H_g(0,0)=\di\begin{bmatrix}-2 & 0\\0 & -2\end{bmatrix}$ is negative definite, $H_h(0,0)=\di\begin{bmatrix}2 & 0\\0 & -2\end{bmatrix}$ is indefinite. Therefore, $(0,0)$ is a local minimizer of $f$, a local maximizer of $g$, and a saddle point of $h$. 
\end{example}

Now let us look at an example which shows that when the Hessian matrix is positive semi-definite but not positive definite,   we cannot  make any conclusion about the nature of a stationary point.
\begin{example}{}
Consider the functions $f:\mb{R}^2\to\mb{R}$ and $g:\mb{R}^2\to\mb{R}$   given respectively by
\[f(x,y)=x^2+y^4, \quad g(x,y)=x^2-y^4.\]These are infinitely differentiable functions. It is easy to check that $(0,0)$ is a stationary point of both of them. Now,
\[H_f(0,0)=H_g(0,0)=\begin{bmatrix} 2 & 0\\0 & 0\end{bmatrix} \]
is a positive semi-definite matrix. However, $(0,0)$ is a local minimizer of $f$, but a saddle point of $g$.
\end{example}
To determine the definiteness of an $n\times n$ symmetric matrix by  looking at the sign of its eigenvalues is ineffective when $n\geq 3$. There is an easier way to determine whether  a symmetric matrix is positive definite.  Let us first introduce the definition of principal submatrices.

\begin{definition}
{Principal Submatrices}Let $A$ be an $n\times n$ matrix. For $1\leq k\leq n$, the $k^{\text{th}}$-principal submatrix $M_k$ of $A$ is the $k\times k$ matrix consists of the first $k$ rows and first $k$ columns of $A$. 
\end{definition}

\begin{example}{}
For the matrix $A=\begin{bmatrix} 1 & 2 & 3\\ 4 & 5 & 6\\ 7 & 8 & 9\end{bmatrix}$, the first, second and third principal submatrices are
\[M_1=\begin{bmatrix} 1\end{bmatrix}, \;M_2=\begin{bmatrix} 1 & 2\\ 4 & 5\end{bmatrix}, \; M_3=\begin{bmatrix} 1 & 2 & 3\\ 4 & 5 & 6\\ 7 & 8 & 9\end{bmatrix}\]respectively.

\end{example}

\begin{theorem}{Sylvester's Criterion for Positive Definiteness}
An $n\times n$ symmetric matrix $A$  is positive definite if and only if $\det M_k>0$ for all $1\leq k\leq n$, where $M_k$ is its $k^{\text{th}}$ principal submatrix.
\end{theorem}
The proof of this theorem is given in Appendix \ref{appA}. 
Using the fact that a symmetric matrix $A$ is negative definite if and only if $-A$ is positive definite, it is easy to obtain a criterion for a symmetric matrix to be negative definite in terms of the determinants of its principal submatrices. 

\begin{theorem}{Sylvester's Criterion for Negative Definiteness}
An $n\times n$ symmetric matrix $A$  is negative  definite if and only if   $(-1)^k\det M_k>0$ for all $1\leq k\leq n$,   where $M_k$ is its $k^{\text{th}}$ principal submatrix.
\end{theorem}

\begin{example}{}
Consider the matrix 
\[A=\begin{bmatrix}  1  &   2  &  -3\\
    -1   &  4   &  2\\
    -3   &  5 &    8\end{bmatrix}.\]
    Since
    \[\det M_1=1, \;\det M_2=6, \;\det M_3=\det A=5\] are all positive, $A$ is positive definite.
\end{example}

For a function $f:\mathcal{O}\to\mb{R}$ defined on an open subset $\mathcal{O}$ of $\mb{R}^2$, we have the following.
\begin{theorem}{}
Let $\mathcal{O}$ be an open subset of $\mb{R}^2$. Suppose that $(x_0, y_0)$ is a stationary point of the twice continuously differentiable function $f:\mathcal{O}\to\mb{R}$. Let
\[D(x_0, y_0)=\frac{\pa^2 f}{\pa x^2}(x_0, y_0)\frac{\pa^2 f}{\pa y^2}(x_0, y_0)-\left[\frac{\pa^2 f}{\pa x\pa y}(x_0, y_0)\right]^2.\]
\begin{enumerate}[(i)]
\item If $\di\frac{\pa^2 f}{\pa x^2}(x_0,y_0)>0$ and $D(x_0, y_0)>0$, then
the point $(x_0, y_0)$ is  a local minimizer of $f$.
\item  If $\di\frac{\pa^2 f}{\pa x^2}(x_0,y_0)<0$ and $D(x_0, y_0)>0$, then
the point $(x_0, y_0)$ is  a local maximizer of $f$.
\item If $D(x_0, y_0)<0$, the point $(x_0, y_0)$ is a saddle point of $f$.
\end{enumerate}
\end{theorem}
\begin{myproof}{Proof}
We notice that
\[H_f(x_0,y_0)=\begin{bmatrix}\vspace{0.3cm}\di \frac{\pa^2 f}{\pa x^2}(x_0, y_0) &\di \frac{\pa^2 f}{\pa x\pa y}(x_0, y_0)\\
\di \frac{\pa^2 f}{\pa x\pa y}(x_0, y_0) & \di \frac{\pa^2 f}{\pa y^2}(x_0, y_0)\end{bmatrix}.\]\bp
Hence, $\di \di\frac{\pa^2 f}{\pa x^2}(x_0,y_0)$ is the determinant of the first principal submatrix of $H_f(x_0, y_0)$, while $D(x_0, y_0)$ is the determinant of $H_f(x_0, y_0)$, the second principal submatrix of $H_f(x_0,y_0)$. Thus, (i) and (ii) follow  from the Sylvester criteria as well as the second derivative test.

For (iii), we notice that the $2\times 2$ matrix $H_f(x_0, y_0)$ is indefinite if and only if it has one positive eigenvalue and one negative eigenvalue, if and only if $D(x_0, y_0)=\det H_f(x_0, y_0)<0$.
\end{myproof}

Now we look at some examples of the applications of the second derivative test.

\begin{example}{}
Let $f:\mb{R}^2\to\mb{R}$ be the function defined as 
\[f(x,y)=x^4+y^4+4xy.\]
Find the stationary points of $f$ and classify them.
\end{example}
\begin{solution}{Solution}
 Since $f$ is a polynomial function, it is infinitely differentiable. 
\[\nabla f(x,y)=(4x^3+4y, 4y^3+4x).\]
To find the stationary points, we need to solve the system of equations
\[\begin{cases} x^3+y=0\\
y^3+x=0\end{cases}\hspace{-0.4cm}.\]From the first equation, we have $y=-x^3$. Substitute into the second equation gives
\[-x^9+x=0,\]
or equivalently, 
\[x(x^8-1)=0.\]
\bs
Thus,  $x=0$ or $x=\pm 1$. When $x=0$, $y=0$. When $x=\pm 1$, $y=\mp 1$. Therefore, the stationary points of $f$ are $\mf{u}_1=(0,0)$, $\mf{u}_2=(1,-1)$ and $\mf{u}_3=(-1,1)$.
Now,
\[H_f(x,y)=\begin{bmatrix} 12x^2 & 4 \\ 4 & 12y^2\end{bmatrix}.\]Therefore, \[H_f(\mf{u}_1)=\begin{bmatrix} 0 & 4 \\ 4 & 0\end{bmatrix},\quad H_f(\mf{u}_2)=H_f(\mf{u}_3)=\begin{bmatrix} 12 & 4 \\ 4 & 12\end{bmatrix}.\]
It follows that
\[D(\mf{u}_1)=-16<0,\quad D(\mf{u}_2)=D(\mf{u}_3)=128>0.\]
Since $f_{xx}(\mf{u}_2)=f_{xx}(\mf{u}_3)=12>0$, we conclude that
$\mf{u}_1$ is a saddle point, $\mathbf{u}_2$ and $\mf{u}_3$ are local minimizers.
\end{solution}

\begin{figure}[ht]
\centering
\includegraphics[scale=0.18]{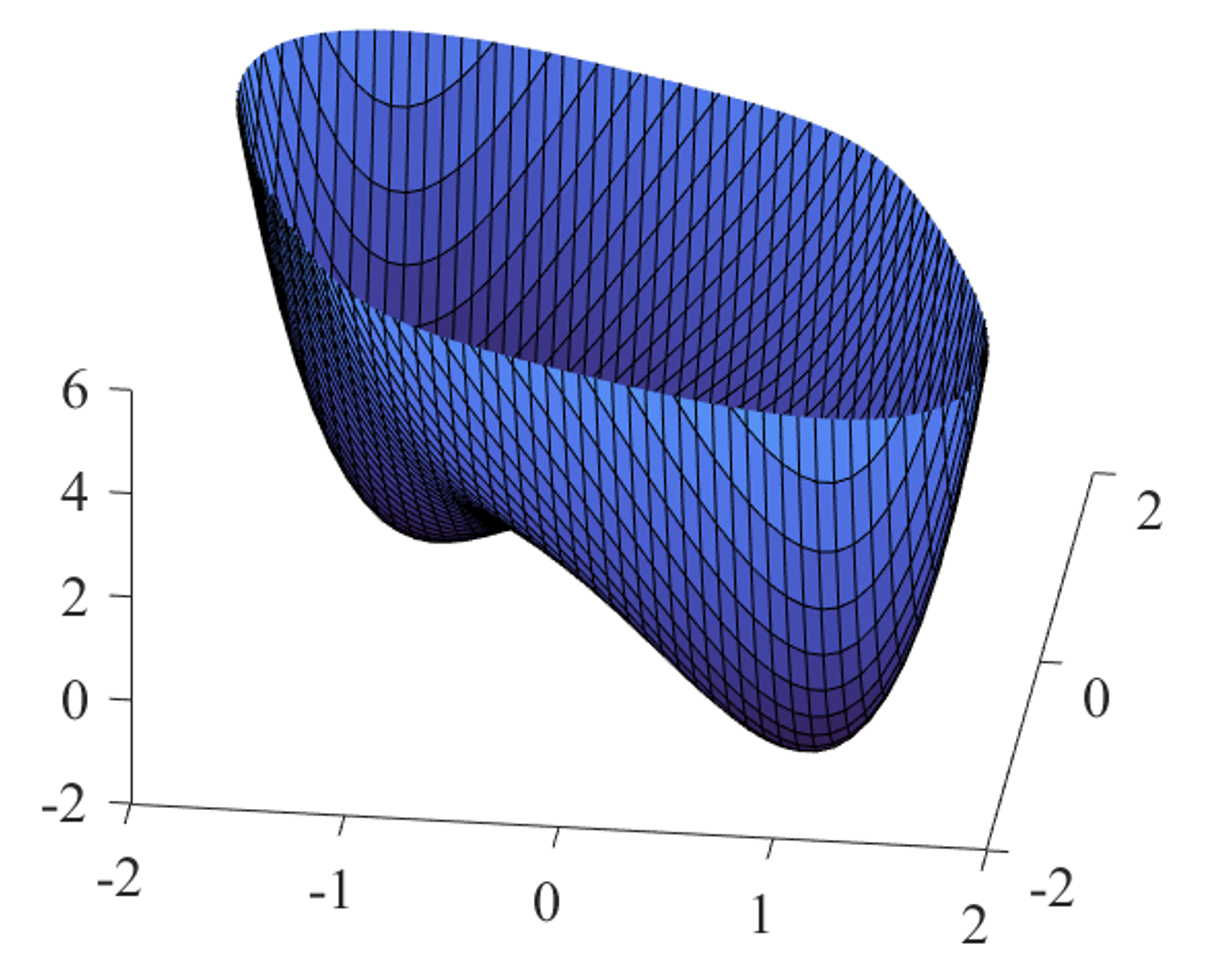}

\caption{The function  $f(x,y)=x^4+y^4+4xy$.}\label{figure60}
\end{figure}

\begin{example}{}
 Consider the function $f:\mathbb{R}^3\rightarrow\mathbb{R}$ defined as
\begin{align*}
f(x,y,z)=x^3-xy^2+5x^2 - 4xy - 2xz +y^2 + 6yz + 37z^2.
\end{align*}Show that $(0,0,0)$ is a local minimizer of $f$.
\end{example}
\begin{solution}{Solution}
 Since $f$ is a polynomial function, it is infinitely differentiable. Since
\[\nabla f(x,y,z)=(3x^2-y^2+10x-4y-2z, -2xy-4x+2y+6z, -2x+6y+74z),\]we find that
\[\nabla f(0,0,0)=(0,0,0).\]
Hence, $(0,0,0)$ is a stationary point.

Now, 
\[H_f(x,y, z)=\begin{bmatrix}  6x+10 & -2y-4 & -2\\ -2y-4 & -2x+2 & 6\\
-2 & 6 & 74\end{bmatrix}.\] 
Therefore,
\[ H_f(0,0,0)=\begin{bmatrix}   10 &  -4 & -2\\  -4 & 2 & 6\\
-2 & 6 & 74\end{bmatrix}.\] 

The determinants of the three principal submatrices of $H_f(0,0,0)$ are
\begin{gather*}
\det M_1 =10,\hspace{1cm}
\det M_2 =\left|\begin{matrix} 10 & -4\\ -4 & 2\end{matrix}\right|=4,\\
\det M_3 =\left|\begin{matrix}   10 &  -4 & -2\\  -4 & 2 & 6\\
-2 & 6 & 74\end{matrix}\right|=24.
\end{gather*}This shows that $H_f(0,0,0)$ is positive definite. Hence, $(0,0,0)$ is a local minimizer of $f$.
\end{solution}
\vp
\noindent
{\bf \large Exercises  \thesection}
\setcounter{myquestion}{1}
\begin{question}{\themyquestion}
Let $f:\mb{R}^2\to\mb{R}$ be the function defined as 
\[f(x,y)=x^2+4y^2+5xy-8x-11y+7.\]
Find the stationary points of $f$ and classify them.
\end{question}
\atc
 \begin{question}{\themyquestion}
Let $f:\mb{R}^2\to\mb{R}$ be the function defined as 
\[f(x,y)=x^2+4y^2+3xy-5x-18y+1.\]
Find the stationary points of $f$ and classify them.
\end{question}

\atc
 \begin{question}{\themyquestion}
Let $f:\mb{R}^2\to\mb{R}$ be the function defined as 
\[f(x,y)=x^3+y^3+12xy.\]
Find the stationary points of $f$ and classify them.
\end{question}

\atc
 \begin{question}{\themyquestion}
Consider the function $f:\mathbb{R}^3\rightarrow\mathbb{R}$ defined as
\begin{align*}
f(x,y,z)= z^3-2z^2-x^2-y^2-xy+x-y.
\end{align*}Show that $(1,-1,0)$ is a stationary point of $f$ and determine the nature of this stationary point.
\end{question}

\atc
 \begin{question}{\themyquestion}
Consider the function $f:\mathbb{R}^3\rightarrow\mathbb{R}$ defined as
\begin{align*}
f(x,y,z)= z^3+2z^2-x^2-y^2-xy+x-y.
\end{align*}Show that $(1,-1,0)$ is a stationary point of $f$ and determine the nature of this stationary point.
\end{question}

\chapter{The Inverse  and   Implicit Function Theorems}

In this chapter, we discuss the inverse function theorem and implicit function theorem, which are two important theorems in multivariable analysis. Given a function that maps  a subset of $\mb{R}^n$ to $\mb{R}^n$, the inverse function theorem gives   sufficient conditions for the existence of a local inverse  and its differentiability. Given a system of $m$ equations with $n+m$ variables, the implicit function theorem gives sufficient conditions to   solve $m$ of the  variables   in terms of the other $n$ variables locally   such that the solutions are differentiable functions. We want to emphasize that these theorems are {\it local}, in the sense that each of them   asserts the existence of a function defined in a neighbourhood  of a point.

 In some sense, the two theorems are equivalent, which means one can deduce one from the other. In this book, we will prove the inverse function theorem first, and use it to deduce the implicit function theorem.

\section{The Inverse Function Theorem} 

Let $\mathfrak{D}$ be a  subset of $\mb{R}^n$. If the function $\mf{F}:\mk{D}\to\mb{R}^n$ is one-to-one,  we can define the inverse function $\mf{F}^{-1}:\mf{F}(\mk{D})\to \mb{R}^n$. The question we want to study here is the following. If $\mk{D}$ is an open set and $\mf{F}$ is differentiable at the point $\mf{x}_0$ in $\mk{D}$, is the inverse function $\mf{F}^{-1}$ differentiable at $\mf{y}_0=\mf{F}(\mf{x}_0)$? For this, we  also want the point $\mf{y}_0$ to be an interior point of $\mf{F}(\mk{D})$. More precisely, is there a neighbourhood $U$ of $\mf{x}_0$ that is mapped bijectively by $\mf{F}$ to a neighbourhood $V$ of $\mf{y}_0$? If the answer is yes, and $\mf{F}^{-1}$ is differentiable at $\mf{y}_0$, then the chain rule would imply that
\[\mf{DF}^{-1}(\mf{y}_0)\mf{DF}(\mf{x}_0)=I_n.\]
Hence, a necessary condition for $\mf{F}^{-1}$ to be differentiable at $\mf{y}_0$ is that the derivative matrix $\text{DF}(\mf{x}_0)$ has to be invertible.

Let us study the map $f:\mb{R}\to\mb{R}$ given by $f(x)=x^2$. The range of the function is $[0, \infty)$. Notice that if $x_0>0$, then $I=(0,\infty)$ is a neighbourhood of $x_0$ that is mapped bijectively by $f$ to the neighbourhood $J=(0,\infty)$ of $f(x_0)$. If $x_0<0$, then $I=(-\infty, 0)$ is a neighbourhood of $x_0$ that is mapped bijectively by $f$ to the neighbourhood $J=(0,\infty)$ of $f(x_0)$. However, if $x_0=0$, the point $f(x_0)=0$ is not an interior point of $f(\mb{R})=[0,\infty)$. Notice that $f'(x)=2x$. Therefore, $x=0$ is the point which $f'(x)=0$. 

If $x_0>0$, take $I=(0,\infty)$ and $J=(0,\infty)$. Then $f:I\to J$ has an inverse given by $f^{-1}:J\to I$, $f^{-1}(x)=\sqrt{x}$. It is a differentiable function with
\[(f^{-1})'(x)=\frac{1}{2\sqrt{x}}.\]
In particular, at $y_0=f(x_0)=x_0^2$,
\[(f^{-1})'(y_0)=\frac{1}{2\sqrt{y_0}}=\frac{1}{2x_0}=\frac{1}{f'(x_0)}.\]
Similarly, if $x_0<0$, take $I=(-\infty, 0)$ and $J=(0,\infty)$. Then $f:I\to J$ has an inverse given by $f^{-1}:J\to I$, $f^{-1}(x)=-\sqrt{x}$. It is a differentiable function with
\[(f^{-1})'(x)=-\frac{1}{2\sqrt{x}}.\]
In particular, at $y_0=f(x_0)=x_0^2$,
\[(f^{-1})'(y_0)=-\frac{1}{2\sqrt{y_0}}=\frac{1}{2x_0}=\frac{1}{f'(x_0)}.\]

\begin{figure}[ht]
\centering
\includegraphics[scale=0.2]{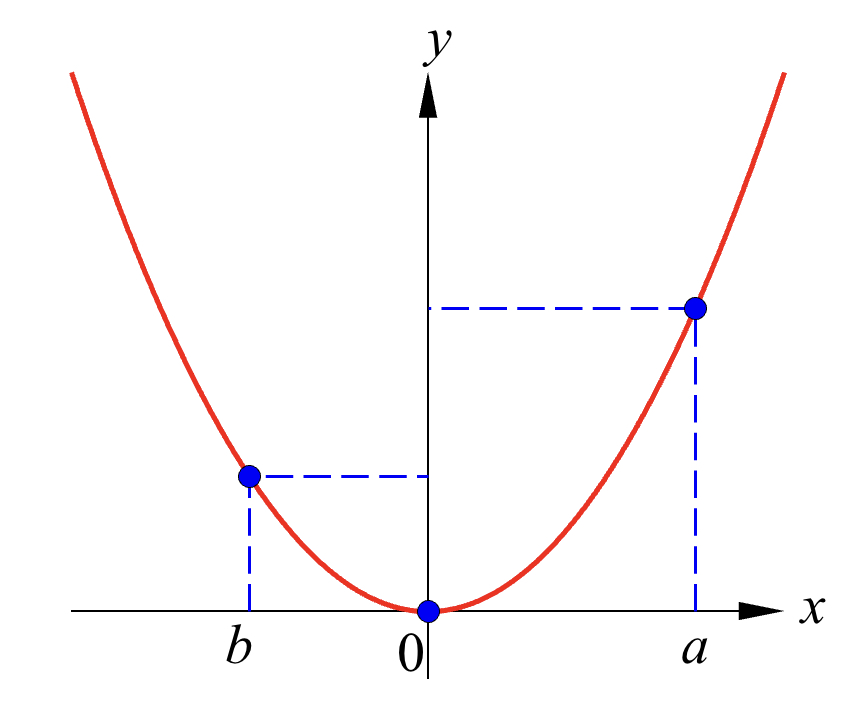}

\caption{The function $f:\mb{R}\to\mb{R}$, $f(x)=x^2$.}\label{figure61}
\end{figure}
For a single variable function, the inverse function theorem takes the following form.
\begin{theorem}
{ (Single Variable) Inverse Function Theorem}
Let $\mathcal{O}$ be an open subset  of $\mb{R}$ that contains the point $x_0$, and let $f:\mathcal{O}\to\mb{R}$ be a continuously differentiable function defined on $\mathcal{O}$. Suppose that $f'(x_0)\neq 0$. Then there exists an open interval $I$ containing $x_0$ such that $f$ maps $I$ bijectively onto the open interval $J=f(I)$. The inverse function $f^{-1}:J\to I$ is continuously differentiable. For any $y\in J$, if $x$ is the point in $I$ such that $f(x)=y$, then 
\[(f^{-1})'(y)=\frac{1}{f'(x)}.\]
\end{theorem}

\begin{myproof}{Proof}
Without loss of generality, assume that $f'(x_0)>0$. Since $\mathcal{O}$ is an open set and $f'$ is continuous at $x_0$, there is an $r_1>0$ such that $(x_0-r_1, x_0+r_1)\subset \mathcal{O}$ and for all $x\in (x_0-r_1, x_0+r_1)$, 
\[|f'(x)-f'(x_0)|<\frac{f'(x_0)}{2}.\]
This implies that
\[f'(x)>\frac{f'(x_0)}{2}>0\hspace{1cm}\text{for all}\;x\in  (x_0-r_1, x_0+r_1).\]
Therefore, $f$ is strictly increasing on $ (x_0-r_1, x_0+r_1)$. Take any $r>0$ that is less that $r_1$. Then $[x-r, x+r]\subset  (x_0-r_1, x_0+r_1)$. By intermediate value theorem, the function $f$ maps $[x-r, x+r]$ bijectively onto  $[f(x-r), f(x+r)]$. Let $I=(x-r, x+r)$ and $J=(f(x-r), f(x+r))$. Then $f:I\to J$ is a bijection and $f^{-1}:J\to I$ exists. In volume I, we have proved that    $f^{-1}$  is differentiable, and 
\[(f^{-1})'(y)=\frac{1}{f'(f^{-1}(y))}\hspace{1cm}\text{for all}\;y\in J.\] This formula shows that $(f^{-1})':J\to\mb{R}$ is continuous.
\end{myproof}

\begin{remark}{}

 In the inverse function theorem, we   determine the   invertibility of the function in a neighbourhood of a point $x_0$. The theorem says that if $f$ is continuously differentiable and  $f'(x_0)\neq 0$, then $f$ is locally invertible at $x_0$. Here the assumption that $f'$ is continuous   is essential.
In volume I, we have seen that for a continuous function $f:I\to \mb{R}$ defined on an open interval $I$ to be one-to-one, it is necessary that it is strictly monotonic. The function $f:\mb{R}\to\mb{R}$,
\[f(x)=\begin{cases}\di  x+x^2\sin\left(\frac{1}{x}\right),\quad &\text{if}\;x\neq 0,\\0,\quad & \text{if}\;x=0,\end{cases}\] is an example of a differentiable function where $f'(0)=1\neq 0$, but $f$ fails to be strictly monotonic in any neighbourhood of the point $x=0$. 

 This annoying behavior can be removed if we assume that $f'$ is continuous. If $f'(x_0)\neq 0$ and $f'$ is continuous, there is a neighbourhood $I$ of $x_0$ such that $f'(x)$ has the same sign as $f'(x_0)$ for all $x\in I$. This implies that $f$ is strictly monotonic on $I$.

\end{remark}

\begin{example}{}
Let $f:\mb{R}\to\mb{R}$ be the function defined as
\[f(x)= 2x+4\cos x.\]
 Show that there is an open interval $I$ containing $0$ such that $f:I\to \mb{R}$ is one-to-one, and $f^{-1}:f(I)\to \mb{R}$ is continuously differentiable. Determine $(f^{-1})'(f(0))$.
\end{example}
\begin{solution}
{Solution}
The function $f$ is infinitely differentiable and $f'(x)=2-4\sin x$. Since $f'(0)=2\neq 0$, the inverse function theorem says that there is an open interval $I$ containing $0$ such that $f:I\to \mb{R}$ is one-to-one, and $f^{-1}:f(I)\to \mb{R}$ is continuously differentiable. Moreover,
\[(f^{-1})'(f(0))=\frac{1}{f'(0)}=\frac{1}{2}.\]
\end{solution}

Now let us consider functions defined on open subsets of $\mb{R}^n$, where $n\geq 2$. We first consider a linear transformation $\mf{T}:\mb{R}^n\to\mb{R}^n$. There is an $n\times n$ matrix $A$ such that 
\[\mf{T}(\mf{x})=A\mf{x}.\]
The mapping $\mf{T}:\mb{R}^n\to\mb{R}^n$ is one-to-one if and only if $A$ is invertible, if and only if $\det A\neq 0$. In this case, $\mf{T}$ is a bijection and $\mf{T}^{-1}:\mb{R}^n\to\mb{R}^n$ is the linear transformation  given by
\[\mf{T}^{-1}(\mf{x})=A^{-1}\mf{x}.\] 
Notice that for any $\mf{x}$ and $\mf{y}$ in $\mb{R}^n$,
\[\mf{DT}(\mf{x})=A,\hspace{1cm}\mf{DT}^{-1}(\mf{y})=A^{-1}.\] 
The content of the inverse function theorem is to extend this to nonlinear mappings.

\begin{theorem}[label=230808_1]
{ Inverse Function Theorem}
Let $\mathcal{O}$ be an open subset  of $\mb{R}^n$ that contains the point $\mf{x}_0$, and let $\mf{F}:\mathcal{O}\to\mb{R}^n$ be a continuously differentiable function defined on $\mathcal{O}$. If $\det \mf{DF}(\mf{x}_0)\neq 0$, then we have the followings.
\begin{enumerate}[(i)]
\item There exists a  neighbourhood  $U$  of $\mf{x}_0$ such that $\mf{F}$ maps $U$ bijectively onto the {\it open} set $V=\mf{F}(U)$.
\item  The inverse function $\mf{F}^{-1}:V\to U$ is continuously differentiable.

\item  For any $\mf{y}\in V$, if $\mf{x}$ is the point in $U$ such that $\mf{F}(\mf{x})=\mf{y}$, then 
\[\mf{DF}^{-1}(\mf{y})=\mf{DF}(\mf{F}^{-1}(\mf{y}))^{-1}=\mf{DF}(\mf{x})^{-1}.\]

\end{enumerate}
\end{theorem}

\begin{figure}[ht]
\centering
\includegraphics[scale=0.2]{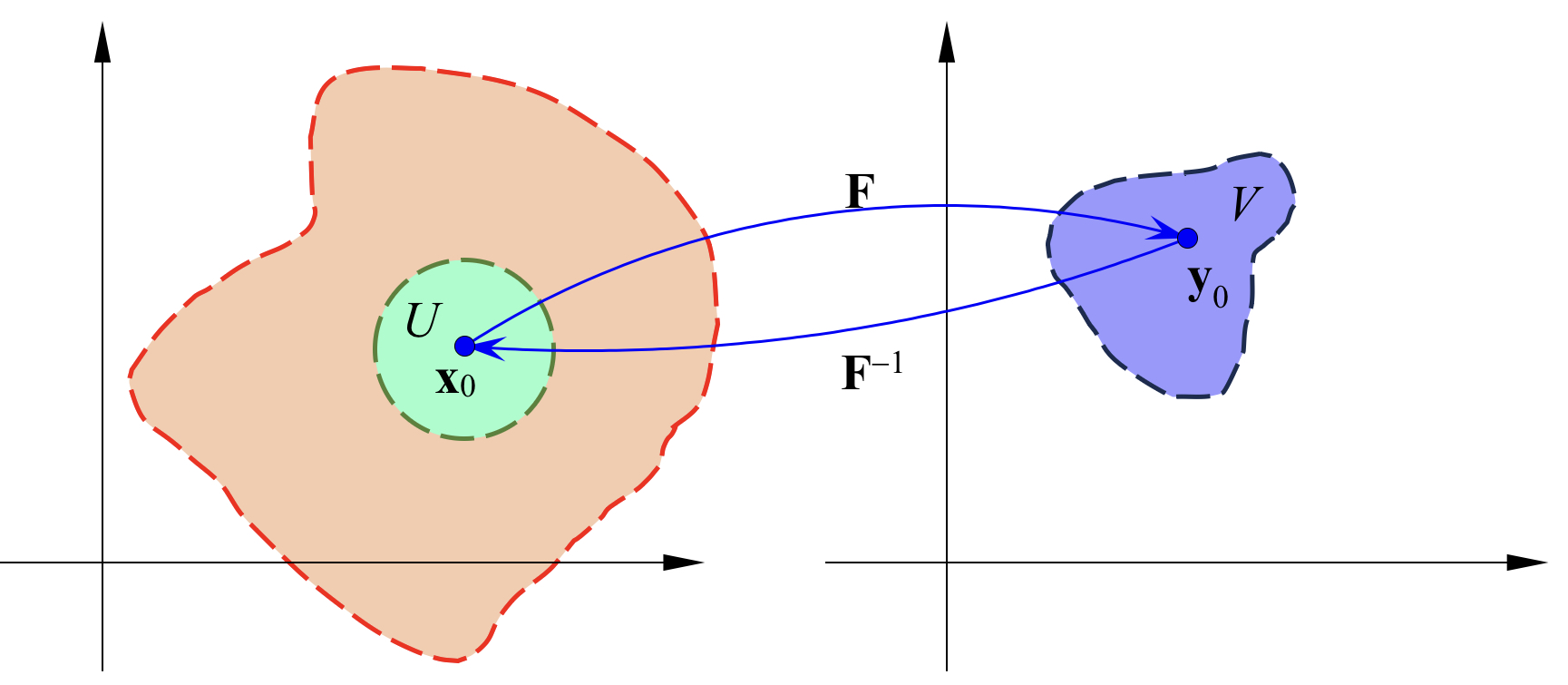}

\caption{The inverse function theorem.}\label{figure62}
\end{figure}
For a linear transformation which is a degree one polynomial mapping, the inverse function theorem holds {\it globally}. 
For a general continuously differentiable mapping, the inverse function theorem says that the first order approximation of the function at a point can determine the local invertibility of the function at that point.

When $n\geq 2$,
the proof of the inverse function theorem   is substantially more complicated than the $n=1$ case, as we do not have the monotonicity argument used in the $n=1$ case. The proof will be presented in Section \ref{sec5.2}.  We will   discuss the examples and applications in this section.

\begin{example}[label=230808_4]{}
 Let $\mf{F}:\mb{R}^2\to\mb{R}^2$ be the mapping defined by
 \[\mf{F}(x,y)=(3x-2y+7, 4x+5y-2).\]
 Show that $\mf{F}$ is a bijection, and find $\mf{F}^{-1}(x,y)$ and $\mf{DF}^{-1}(x,y)$.
\end{example}

\begin{solution}{Solution}
The mapping $\mf{F}:\mb{R}^2\to \mb{R}^2$
can be written as $\mf{F}(\mf{x})=\mf{T}(\mf{x})+\mf{b}$, where $\mf{T}:\mb{R}^2\to\mb{R}^2$ is the linear transformation
\[\mf{T}(x,y)=(3x-2y, 4x+5y),\]
\bs 
and $\mf{b}=(7, -2)$. For $\mf{u}=(x,y)$, $\mf{T}(\mf{u})=A\mf{u}$, where $A=\di\begin{bmatrix} 3 & -2\\4 & 5\end{bmatrix}$. Since $\det A=23\neq 0$, the linear transformation $\mf{T}:\mb{R}^2\to \mb{R}^2$ is  one-to-one. Hence, $\mf{F}:\mb{R}^2\to\mb{R}^2$ is also one-to-one. Given $\mf{v}\in\mb{R}^2$, let $\mf{u}=A^{-1}(\mf{v}-\mf{b})$. Then $\mf{F}(\mf{u})=\mf{v}$. Hence, $\mf{F}$ is also onto. The inverse $\mf{F}^{-1}:\mb{R}^2\to\mb{R}^2$ is given by
\[\mf{F}^{-1}(\mf{v})=A^{-1}(\mf{v}-\mf{b}).\]
Since 
\[A^{-1}=\frac{1}{23}\begin{bmatrix} 5 & 2\\-4 & 3\end{bmatrix},\]
we find that
\begin{align*}\mf{F}^{-1}(x,y)&=\left(\frac{5(x-7)+2(y+2)}{23}, \frac{-4(x-7)+3(y+2)}{23}\right)\\
&=\left(\frac{5x+2y-31}{23}, \frac{-4x+3y+34}{23}\right),
\end{align*}
and
\[\mf{DF}^{-1}(x,y)=\frac{1}{23}\begin{bmatrix} 5 & 2\\-4 & 3\end{bmatrix}.\]
\end{solution}
\begin{example}[label=230808_5]{}
Determine the values of $a$ such that the mapping $\mf{F}:\mb{R}^3\to \mb{R}^3$ defined by
\[\mf{F}(x,y,z)=(2x+y+az, x-y+3z, 3x+2y+z+7)\]is invertible.
\end{example}

\begin{solution}{Solution}
The mapping $\mf{F}:\mb{R}^3\to \mb{R}^3$
can be written as $\mf{F}(\mf{x})=\mf{T}(\mf{x})+\mf{b}$, where $\mf{T}:\mb{R}^3\to\mb{R}^3$ is the linear transformation
\[\mf{T}(x,y,z)=(2x+y+az, x-y+3z, 3x+2y+z),\] 
\bs
and $\mf{b}=(0,0,7)$.
 Thus, $\mf{F}$ is a degree one polynomial mapping with
\[\mf{DF}(\mf{x})=\begin{bmatrix} 2  & 1  & a\\ 1 & -1 & 3\\3 & 2 & 1\end{bmatrix}.\]
The mapping $\mf{F}$ is invertible if and only if it is one-to-one, if and only if $\mf{T}$ is one-to-one, if and only if $\det \mf{DF}(\mf{x})\neq 0$.  Since
\[\det \mf{DF}(\mf{x})=5a-6,\] the mapping $\mf{F}$ is invertible if and only if $a\neq 6/5$. 
\end{solution}

\begin{example}{}
 Let $\mf{\Phi}:\mb{R}^2\to\mb{R}^2$ be the mapping defined as
 \[\mf{\Phi}(r,\theta)=(r\cos\theta, r\sin\theta).\]
 Determine the points $(r,\theta )\in \mb{R}^2$ where the inverse function theorem can be applied to this mapping. Explain the significance of this result.
\end{example}
\begin{solution}{Solution}
Since $\sin\theta$ and $\cos\theta$ are infinitely differentiable functions, the mapping $\mf{\Phi}$ is infinitely differentiable with
\[\mf{D\Phi}(r,\theta)=\begin{bmatrix} \cos\theta & -r\sin\theta\\\sin\theta & r\cos\theta\end{bmatrix}.\]
Since \[\det \mf{D\Phi}(r,\theta)=r\cos^2\theta+r\sin^2\theta=r,\]

the inverse function theorem is not applicable at the point $(r,\theta)$ if $r=0$. 

The mapping $\mf{\Phi}$ is a change from polar coordinates to rectangular coordinates. The result above shows that the change of coordinates is locally one-to-one away from the origin of the $xy$-plane.
\end{solution}
\begin{example}{}
 Consider the mapping $\mf{F}:\mb{R}^2\to\mb{R}^2$ given by
 \[\mf{F}(x,y)=(x^2-y^2, xy).\]
 Show that there is a neighbourhood $U$ of the point $\mf{u}_0=(1,1)$ such that $\mf{F}:U\to\mb{R}^2$ is one-to-one, $V=\mf{F}(U)$ is an open set, and $\mf{G}=\mf{F}^{-1}:V\to U$ is continuously differentiable. Then find $\di \frac{\pa G_1}{\pa y}(0,1)$.
\end{example}

\begin{solution}{Solution}
The mapping $\mf{F}$ is a polynomial mapping. Thus, it is continuously differentiable. Notice that $\mf{F}(\mf{u}_0)=(0,1)$ and 
\[\mf{DF}(x,y)=\begin{bmatrix} 2x &  -2y\\ y & x\end{bmatrix},\quad \mf{DF}(\mf{u}_0)=\begin{bmatrix} 2  &  -2 \\ 1 & 1\end{bmatrix}.\]Since 
$\det \mf{DF}(\mf{u}_0)=4\neq 0$, the inverse function theorem implies that there is a neighbourhood $U$ of the point $\mf{u}_0$ such that $\mf{F}:U\to\mb{R}^2$ is one-to-one, $V=\mf{F}(U)$ is an open set, and $\mf{G}=\mf{F}^{-1}:V\to U$ is continuously differentiable. 
Moreover,
\[\mf{DG}(0,1)=\mf{DF}(1,1)^{-1}=\frac{1}{4}\begin{bmatrix} 1  &   2 \\ -1 & 2\end{bmatrix}.\]
 
From here, we find that
\[ \frac{\pa G_1}{\pa y}(0,1)=\frac{2}{4}=\frac{1}{2}.\]
\end{solution}

\begin{example}{}
 Consider the system of equations
 \begin{align*}
 \sin(x+y)+x^2y+3xy^2&=2,\\
 2xy+5x^2-2y^2&=1.
 \end{align*}
 \be
 Observe that $(x,y)=(1,-1)$ is a solution of this system. Show that there is a neighbourhood $U$ of $\mf{u}_0=(1, -1)$ and  an $r>0$ such that for all $(a, b)$ satisfying $(a-2)^2+(b-1)^2<r^2$, the system \begin{align*}
 \sin(x+y)+x^2y+3xy^2&=a,\\
 2xy+5x^2-2y^2&=b 
 \end{align*}has a unique solution $(x,y)$ that lies in $U$.
\end{example2}

\begin{solution}{Solution}
Let $\mf{F}:\mb{R}^2\to \mb{R}^2$ be the function defined by
\begin{align*}
\mf{F}(x,y)=\left(\sin(x+y)+x^2y+3xy^2, 
 2xy+5x^2-2y^2\right).
\end{align*}Since the sine function is infinitely differentiable, $\sin(x+y)$ is infinitely differentiable. The functions $g(x,y)=x^2y+3xy^2$ and $F_2(x,y)= 2xy+5x^2-2y^2$ are polynomial functions. Hence, they are also infinitely differentiable. This shows that $\mf{F}$ is infinitely differentiable. Since
\[\mf{DF}(x,y)=\begin{bmatrix} \cos(x+y)+2xy+3y^2 & \cos(x+y)+x^2+6xy\\
2y+10x & 2x-4y\end{bmatrix},\]

we find that
\[\mf{DF}(1,-1)=\begin{bmatrix} 2 & -4\\8 & 6\end{bmatrix}.\]
It follows that $\det \text{DF}(1,-1)=44\neq 0$. 

By the inverse function theorem, there exists  a neighbourhood $U_1$ of $\mf{u}_0$ such that $\mf{F}:U_1\to\mb{R}^2$ is one-to-one and $V=\mf{F}(U_1)$ is an open set. Since $\mf{F}(\mf{u}_0)=(2,1)$, the point $\mf{v}_0=(2,1)$ is a point in the open set $V$. Hence, there exists $r>0$ such that $B(\mf{v}_0, r)\subset V$. Since $B(\mf{v}_0, r)$ is open and $\mf{F}$ is continuous, $U=\mf{F}^{-1}\left(B(\mf{v}_0, r)\right)$ is an open subset of $\mb{R}^2$.
 The map $\mf{F}:U\to B(\mf{v}_0, r)$ is a bijection. For all $(a, b)$ satisfying $(a-2)^2+(b-1)^2<r^2$, $(a,b)$ is in $B(\mf{v}_0, r)$. Hence, there is a unique $(x,y)$ in $U$ such that $\mf{F}(x,y)=(a,b)$. This means that the system
 \bs
 \vspace{-0.5cm}
  \begin{align*}
 \sin(x+y)+x^2y+3xy^2&=a,\\
 2xy+5x^2-2y^2&=b 
 \end{align*}has a unique solution $(x,y)$ that lies in $U$.
\end{solution}

At the end of this section, let us prove the following theorem.

\begin{theorem}[label=230808_8]{}
Let $A$ be an $n\times n$ matrix, and let $\mf{x}_0$ and $\mf{y}_0$ be two points in $\mb{R}^n$. Define the mapping $\mf{F}:\mb{R}^n\to \mb{R}^n$ by 
\[\mf{F}(\mf{x})=\mf{y}_0+A\left(\mf{x}-\mf{x}_0\right).\]
Then $\mf{F}$ is infinitely differentiable with $\mf{DF}(\mf{x})=A$. It is one-to-one and onto if and only if $\det A\neq 0$. In this case, 
\[\mf{F}^{-1}(\mf{y})=\mf{x}_0+A^{-1}\left(\mf{y}-\mf{y}_0\right),\quad\text{and}\quad \mf{DF}^{-1}(\mf{y})=A^{-1}.\]In particular, $\mf{F}^{-1}$ is also infinitely differentiable.
\end{theorem}
\begin{myproof}{Proof}
Obviously, $\mf{F}$ is a polynomial mapping. Hence, $\mf{F}$ is infinitely differentiable. By a straightforward computation, we find  that $\mf{DF}=A$.

Notice that $\mf{F}=\mf{F}_2\circ \mf{T}\circ \mf{F}_1$, where $\mf{F}_1:\mb{R}^n\to\mb{R}^n$ is the translation $\mf{F}_1(\mf{x})=\mf{x}-\mf{x}_0$, $\mf{T}:\mb{R}^n\to\mb{R}^n$ is the linear transformation $\mf{T}(\mf{x})=A\mf{x}$, and $\mf{F}_2:\mb{R}^n\to\mb{R}^n$ is the translation $\mf{F}_2(\mf{y})=\mf{y}+\mf{y}_0$. 
 Since translations are bijective mappings, $\mf{F}$ is one-to-one and onto if and only if $\mf{T}:\mb{R}^n\to\mb{R}^n$ is one-to-one and onto, if and only if $\det A\neq 0$. 
 
If  \[\mathbf{y}=\mf{y}_0+A\left(\mf{x}-\mf{x}_0\right),\]
then \[\mf{x}=\mf{x}_0+A^{-1}\left(\mf{y}-\mf{y}_0\right).\] This gives the formula for $\mf{ F}^{-1}(\mf{y})$. The formula for $\mf{DF}^{-1}(\mf{y})$ follows.
\end{myproof}
 
\vp
\noindent
{\bf \large Exercises  \thesection}
\setcounter{myquestion}{1}
\begin{question}{\themyquestion}
Let $f:\mb{R}\to\mb{R}$ be the function defined as
\[f(x)=e^{2x}+4x\sin x+2\cos x.\]
 Show that there is an open interval $I$ containing $0$ such that $f:I\to \mb{R}$ is one-to-one, and $f^{-1}:f(I)\to \mb{R}$ is continuously differentiable. Determine $(f^{-1})'(f(0))$.
\end{question}

\atc
\begin{question}{\themyquestion}
 Let $\mf{F}:\mb{R}^2\to\mb{R}^2$ be the mapping defined by
 \[\mf{F}(x,y)=(3x+2y-5, 7x+4y-3).\]
 Show that $\mf{F}$ is a bijection, and find $\mf{F}^{-1}(x,y)$ and $\mf{DF}^{-1}(x,y)$.
\end{question}
\atc

 \begin{question}{\themyquestion}
 Consider the mapping $\mf{F}:\mb{R}^2\to\mb{R}^2$ given by
 \[\mf{F}(x,y)=(x^2+y^2, xy).\]
 Show that there is a neighbourhood $U$ of the point $\mf{u}_0=(2,1)$ such that $\mf{F}:U\to\mb{R}^2$ is one-to-one, $V=\mf{F}(U)$ is an open set, and $\mf{G}=\mf{F}^{-1}:V\to U$ is continuously differentiable. Then find $\di \frac{\pa G_2}{\pa x}(5,2)$.
\end{question}

\atc

 \begin{question}{\themyquestion}
Let $\mf{\Phi}:\mb{R}^3\to\mb{R}^3$ be the mapping defined as
 \[\mf{\Phi}(\rho, \phi,\theta)=(\rho\sin\phi\cos\theta, \rho\sin\phi\sin\theta, \rho\cos\phi).\]
 Determine the points $(\rho,\phi ,\theta)\in \mb{R}^3$ where the inverse function theorem can be applied to this mapping. Explain the significance of this result.
\end{question}
\atc

 \begin{question}{\themyquestion}
Consider the system of equations
 \begin{align*}
4x+y-5xy&=2,\\
 x^2+y^2-3xy^2&=5.
 \end{align*}Observe that $(x,y)=(-1,1)$ is a solution of this system. Show that there is a neighbourhood $U$ of $\mf{u}_0=(-1, 1)$ and  an $r>0$ such that for all $(a, b)$ satisfying $(a-2)^2+(b-5)^2<r^2$, the system
  \begin{align*}
4x+y-5xy&=a,\\
 x^2+y^2-3xy^2&=b
 \end{align*}has a unique solution $(x,y)$ that lies in $U$.
\end{question}

\section{The Proof of the Inverse Function Theorem} \label{sec5.2}

In this section, we prove the inverse function theorem stated in Theorem \ref{230808_1}.
 The hardest part of the proof is the first statement, which asserts that there is a neighbourhood $U$ of $\mf{x}_0$ such that restricted to $U$, $\mf{F}$ is one-to-one, and the image of $U$ under $\mf{F}$ is open in $\mb{R}^n$.
 
In the statement of the inverse function theorem, we assume that the derivative matrix of the continuously differentiable mapping $\mf{F}:\mathcal{O}\to\mb{R}^n$ is  invertible at the point $\mf{x}_0$. The continuities of the partial derivatives of $\mf{F}$ then implies that there is a neighbourhood $\mathcal{N}$ of $\mf{x}_0$ such that the derivative matrix of $\mf{F}$ at any $\mf{x}$ in $\mathcal{N}$ is also invertible. 

Theorem \ref{230725_6} asserts that a linear transformation $\mf{T}:\mb{R}^n\to\mb{R}^n$ is invertible if and only if there is a positive constant $c$ such that
\[\Vert\mf{T}(\mf{u})-\mf{T}(\mf{v})\Vert\geq c\Vert\mf{u}-\mf{v}\Vert \hspace{1cm}\text{for all}\;\mf{u}, \mf{v}\in\mb{R}^n.\]
\begin{definition}{Stable Mappings}{}
A mapping $\mf{F}:\mk{D}\to \mb{R}^n$ is {\it stable} if there is a positive constant $c$ such that 
\[\Vert\mf{F}(\mf{u})-\mf{F}(\mf{v})\Vert\geq c\Vert\mf{u}-\mf{v}\Vert \hspace{1cm}\text{for all}\;\mf{u}, \mf{v}\in\mathfrak{D}.\]
\end{definition}
In other words, a linear transformation $\mf{T}:\mb{R}^n\to\mb{R}^n$ is invertible if and only if it is stable.

\begin{remark}[label=230808_13]{Stable Mappings vs Lipschitz Mappings}
Let $\mk{D}$ be a subset of $\mb{R}^n$. Observe that if $\mf{F}:\mk{D}\to\mb{R}^n$ is a stable mapping, there is a constant $c>0$ such that
\[\Vert\mf{F}(\mf{u}_1)-\mf{F}(\mf{u}_2)\Vert\geq c\Vert\mf{u}_1-\mf{u}_2\Vert \hspace{1cm}\text{for all}\;\mf{u}_1, \mf{u}_2\in\mathfrak{D}.\]
This implies that $\mf{F}$ is one-to-one, and thus the inverse $\mf{F}^{-1}: \mf{F}(\mk{D})\to\mb{R}^n$ exists. Notice that for any $\mf{v}_1$ and $\mf{v}_2$ in $\mf{F}(\mk{D})$, 
\[\Vert\mf{F}^{-1}(\mf{v}_1)-\mf{F}^{-1}(\mf{v}_2)\Vert\leq \frac{1}{c}\Vert\mf{v}_1-\mf{v}_2\Vert.\] This means that  $\mf{F}^{-1}: \mf{F}(\mk{D})\to\mb{R}^n$ is a Lipschitz mapping.
\end{remark}
For a mapping  $\mf{F}:\mk{D}\to \mb{R}^n$ that satisfies the assumptions in the statement of the inverse function theorem,   it is stable in a neighbourhood of $\mf{x}_0$. 
 \begin{theorem}[label=230808_10]{}
 Let $\mathcal{O}$ be an open subset  of $\mb{R}^n$ that contains the point $\mf{x}_0$, and let $\mf{F}:\mathcal{O}\to\mb{R}^n$ be a continuously differentiable function defined on $\mathcal{O}$. If $\det \mf{DF}(\mf{x}_0)\neq 0$, then there exists a  neighbourhood  $U$  of $\mf{x}_0$ such that $\mf{DF}(\mf{x})$ is invertible for all $\mf{x}\in U$,  $\mf{F}$ maps $U$ bijectively onto the   open set $V=\mf{F}(U)$, and the map $\mf{F}:U\to V$ is stable.  
 \end{theorem}
Recall that when $A$ is a subset of $\mb{R}^n$, $\mf{u}$ is a point in $\mb{R}^n$,  
\[ A+\mf{u}=\left\{\mf{a}+\mf{u}\,|\,\mf{a}\in A\right\} \] is the translate of the set $A$ by the vector $\mf{u}$. The set $A$ is open if and only if  $A+\mf{u}$ is open, $A$ is closed if and only if $A+\mf{u}$ is closed. 
 
\begin{lemma}[label=230808_9]{}
It is sufficient to prove Theorem \ref{230808_10}  when  $\mf{x}_0=\mf{0}$, $\mf{F}(\mf{x}_0)=\mf{0}$ and $\mf{DF}(\mf{x}_0)=I_n$.
\end{lemma}

\begin{myproof}{\linkt Proof of Lemma \ref{230808_9}}
Assume that   Theorem \ref{230808_10} holds when  $\mf{x}_0=\mf{0}$, $\mf{F}(\mf{x}_0)=\mf{0}$ and $\mf{DF}(\mf{x}_0)=I_n$.

Now given that $\mf{F}:\mathcal{O}\to\mb{R}^n$  is a continuously differentiable mapping with $\det\mf{DF}(\mf{x}_0)\neq 0$,
let $\mf{y}_0=\mf{F}(\mf{x}_0)$ and $A=\mf{DF}(\mf{x}_0)$. Then $A$ is invertible.  Define the open set $\mathcal{D}$ as
$\mathcal{D}=\mathcal{O}-\mf{x}_0$. It is a neighbourhood of the point $\mf{0}$.
 Let $\mf{G}:\mathcal{D}\to\mb{R}^n$ be the mapping 
\[\mf{G}(\mf{x})= A^{-1}\left(\mf{F}(\mf{x}+\mf{x}_0)-\mf{y}_0\right).\]
 
 Then $\mf{G}(\mf{0})=\mf{0}$.
Using the same reasoning as the proof of  Theorem \ref{230808_8}, we find that $\mf{G}$ is continuously differentiable and

\[\mf{DG} (\mf{x})= A^{-1} \mf{DF}(\mf{x}+\mf{x}_0).\] 

\bp
 This gives \[\mf{DG}(\mf{0})=A^{-1}\mf{DF}(\mf{x}_0)=I_n.\] 
By assumption, Theorem \ref{230808_10} holds for the mapping $\mf{G}$. Namely,  there exist  neighbourhoods $\mathcal{U}$ and $\mathcal{V}$ of $\mf{0}$ such that $\mf{G}:\mathcal{U}\to\mathcal{V}$ is a bijection and $\mf{DG}(\mf{x})$ is invertible for all $\mf{x}\in\mathcal{U}$. Moreover, there is a positive constant $a$ such that
  \[\Vert \mf{G}(\mf{u}_1)-\mf{G}(\mf{u}_2)\Vert\geq a\Vert\mf{u}_1-\mf{u}_2\Vert\hspace{1cm}\text{for all}\;\mf{u}_1, \mf{u}_2\in \mathcal{U}.\]
 
 Let $U$ be the neighbourhood of $\mf{x}_0$ given by
$U=\mathcal{U}+\mf{x}_0$. 
 By Theorem \ref{230808_8}, the mapping $\mf{H}:\mb{R}^n\to\mb{R}^n$,
 \[\mf{H}(\mf{y})=A^{-1}(\mf{y}-\mf{y}_0)\]

 is a continuous   bijection. Therefore, $V=\mf{H}^{-1}(\mathcal{V})$ is an open subset of $\mb{R}^n$ that contains $\mf{y}_0$. By definition,  $\mf{F}$ maps $U$ bijectively to $V$.  
 Since 
 \[\mf{F}(\mf{x})=\mathbf{y}_0+A\mf{G}(\mf{x}-\mf{x}_0),\]
 we find that 
 \[\mf{DF}(\mf{x})=A\left(\mf{DG}(\mf{x}-\mf{x}_0)\right).\]
 Since $A$ is invertible, $\mf{DF}(\mf{x})$ is invertible for all $\mf{x}\in U$. 
   Theorem \ref{230725_6} says that there is a positive constant $\alpha$ such that 
   
 \[\Vert A\mf{x}\Vert \geq \alpha\Vert\mf{x}\Vert\hspace{1cm}\text{for all}\;\mf{x}\in \mb{R}^n.\]
 Therefore,
 for any $\mf{u}_1$ and $\mf{u}_2$ in $U$,
 \begin{align*}
 \Vert\mf{F}(\mf{u}_1)-\mf{F}(\mf{u}_2)\Vert &= \Vert A\left(\mf{G}(\mf{u}_1-\mf{x}_0)-\mf{G}(\mf{u}_2-\mf{x}_0)\right)\Vert\\
 &\geq\alpha \Vert\mf{G}(\mf{u}_1-\mf{x}_0)-\mf{G}(\mf{u}_2-\mf{x}_0)\Vert\\
 &\geq a\alpha \Vert\mf{u}_1-\mf{u}_2\Vert.
 \end{align*}This shows that $\mf{F}:U\to V$ is stable, and thus completes the proof of the lemma.
\end{myproof}

 Now we prove Theorem \ref{230808_10}.

 \begin{myproof}{\linkt Proof of Theorem  \ref{230808_10}}
 By Lemma \ref{230808_9}, we only need to consider the case where $\mf{x}_0=\mf{0}$, $\mf{F}(\mf{x}_0)=\bf{0}$ and $\mf{DF}(\mf{x}_0)=I_n$.
 
 Since $\mf{F}:\mathcal{O}\to\mb{R}^n$ is continuously differentiable, the map $\mf{DF}:\mathcal{O}\to \mathcal{M}_n$ is continuous. Since $\det :\mathcal{M}_n\to\mb{R}$ is also continuous, and $\det\mf{DF}(\mf{0})=1$, there is an $r_0>0$ such that $B(\mf{0}, r_0)\subset \mathcal{O}$ and for all $\mf{x}\in B(\mf{0}, r_0)$, $\det \mf{DF}(\mf{x})>\frac{1}{2}$. In particular, $\mf{DF}(\mf{x})$ is invertible for all $\mf{x}\in B(\mf{0}, r_0)$.
 
 Let  $\mf{G}:\mathcal{O}\to \mb{R}^n$ be the mapping defined as
 \[\mf{G}(\mf{x})=\mf{F}(\mf{x})-\mf{x},\]
 so that $\mf{F}(\mf{x})=\mf{x}+\mf{G}(\mf{x})$.  
 The mapping  $\mf{G}$ is   continuosly differentiable. It satisfies $\mf{G}(\mf{0})=\mf{0}$ and \[\mf{DG}(\mf{0})=\mf{DF}(\mf{0})-I_n=\mf{0}.\] Since $\mf{G}$ is continuously differentiable, for any $1\leq i\leq n$, $1\leq j\leq n$, there exists $r_{i,j}>0$ such that $B(\mf{0}, r_{i,j})\subset\mathcal{O}$ and for all $\mf{x}\in B(\mf{0}, r_{i,j})$, 
 
 \[\left|\frac{\pa G_i}{\pa x_j}(\mf{x})\right|=\left|\frac{\pa G_i}{\pa x_j}(\mf{x})-\frac{\pa G_i}{\pa x_j}(\mf{0})\right|<\frac{1}{2n}.\]
  
 Let \[r=\min\left(\{r_{i, j}\,|\, 1\leq i\leq n, 1\leq j\leq n\}\cup\{r_0\}\right).\] Then $r>0$, $B(\mf{0}, r)\subset B(\mf{0}, r_{0})$ and $B(\mf{0}, r)\subset B(\mf{0}, r_{i,j})$ for all $1\leq i\leq n$, $1\leq j\leq n$.
 The ball $B(\mf{0},r)$ is a convex set. If $\mf{u}$ and $\mf{v}$ are two points in  $B(\mf{0},r)$,  mean value theorem implies that for $1\leq i\leq n$, there exists $\mf{z}_i\in B(\mf{0}, r)$ such that
 \[ G_i(\mf{u})-G_i(\mf{v})=\sum_{j=1}^n(u_j-v_j)\frac{\pa G_i}{\pa x_j}(\mf{z}_i).\]

It follows that
 \begin{align*}
 \left|G_i(\mf{u})-G_i(\mf{v})\right| &\leq \sum_{j=1}^n |u_j-v_j|\left|\frac{\pa G_i}{\pa x_j}(\mf{z}_i)\right|\\&\leq\frac{1}{2n}\sum_{j=1}^n |u_j-v_j|\leq \frac{1}{2\sqrt{n}}\Vert\mf{u}-\mf{v}\Vert.\end{align*}
  \bp
Therefore,
 \[\Vert\mf{G}(\mf{u})-\mf{G}(\mf{v})\Vert=\sqrt{\sum_{i=1}^n \left(G_i(\mf{u})-G_i(\mf{v})\right)^2}\leq \frac{1}{2}\Vert\mf{u}-\mf{v}\Vert.\]
  
 This shows that $\mf{G}:B(\mf{0},r)\to \mb{R}^n$ is a map satisfying $\mf{G}(\mf{0})=\mf{0}$, and
 \[\Vert \mf{G}(\mf{u})-\mf{G}(\mf{v})\Vert\leq\frac{1}{2}\Vert\mf{u}-\mf{v}\Vert\hspace{1cm}\text{for all}\;\mf{u}, \mf{v}\in B(\mf{0}, r).\]

 By Theorem \ref{230808_2}, the map $\mf{F}:B(\mf{0}, r)\to \mb{R}^n$ is one-to-one, and its image contains the open ball $B(\mf{0}, r/2)$. Let $V=B(\mf{0},r/2)$. Then $V$ is an open subset of $\mb{R}^n$ that is contained in the image of $\mf{F}$. Since $\mf{F}:B(\mf{0}, r)\to \mb{R}^n$ is continuous, $U=\mf{F}|_{B(\mf{0},r)}^{-1}(V)$ is an open set.   By definition,   $\mf{F}:U\to V$ is a bijection. Since $U$ is contained in $B(\mf{0}, r_0)$, $\mf{DF}(\mf{x})$ is invertible for all $\mf{x}$ in $U$.
Finally,   for any $\mf{u}$ and $\mf{v}$ in $U$,
 \[\Vert\mf{F}(\mf{u})-\mf{F}(\mf{v})\Vert \geq \Vert\mf{u}-\mf{v}\Vert-\Vert\mf{G}(\mf{u})-\mf{G}(\mf{v})\Vert\geq\frac{1}{2} \Vert\mf{u}-\mf{v}\Vert.\]

 This completes the proof of the theorem.

 \end{myproof} 

To complete the proof of  the inverse function theorem, it remains to prove that $\mf{F}^{-1}:V\to  U$ is continuously differentiable, and \[\mf{DF}^{-1}(\mf{y})=\mf{DF}(\mf{F}^{-1}(\mf{y}))^{-1}.\]

\begin{theorem}{}
Let $\mathcal{O}$ be an open subset  of $\mb{R}^n$ that contains the point $\mf{x}_0$, and let $\mf{F}:\mathcal{O}\to\mb{R}^n$ be a continuously differentiable function defined on $\mathcal{O}$. If $\det \mf{DF}(\mf{x}_0)\neq 0$, then there exists a  neighbourhood  $U$  of $\mf{x}_0$ such that $\mf{F}$ maps $U$ bijectively onto the open set $V=\mf{F}(U)$, the inverse function $\mf{F}^{-1}:V\to U$ is continuously differentiable, and for any $\mf{y}\in V$, if $\mf{x}$ is the point in $U$ such that $\mf{F}(\mf{x})=\mf{y}$, then 
\[\mf{DF}^{-1}(\mf{y})=\mf{DF}(\mf{x})^{-1}.\]

\end{theorem}
\begin{myproof}{Proof}
Theorem \ref{230808_10} asserts that there exists a  neighbourhood  $U$  of $\mf{x}_0$ such that $\mf{F}$ maps $U$ bijectively onto the open set $V=\mf{F}(U)$, $\mf{DF}(\mf{x})$ is invertible for all $\mf{x}$ in $U$, and there is a positive constant $c$ such that
 \begin{equation}\label{230808_14}\Vert \mf{F}(\mf{u}_1)-\mf{F}(\mf{u}_2)\Vert\geq c\Vert\mf{u}_1-\mf{u}_2\Vert\hspace{1cm}\text{for all}\;\mf{u}_1, \mf{u}_2\in U.\end{equation}
 
 Now given $\mf{y}$ in $V$, we want to show that $\mf{F}^{-1}$ is differentiable at $\mf{y}$ and $\mf{DF}^{-1}(\mf{y})=\mf{DF}(\mf{x})^{-1}$, where $\mf{x}=\mf{F}^{-1}(\mf{y})$. Since $V$ is open, there is an $r>0$ such that $B(\mf{y}, r)\subset V$. For $\mf{k}\in \mb{R}^n$ such that $\Vert\mf{k}\Vert<r$, let
 \[\mf{h}(\mf{k})=\mf{F}^{-1}(\mf{y}+\mf{k})-\mf{F}^{-1}(\mf{y}).\]
 
 Then
 \[\mf{F}(\mf{x})=\mf{y}\quad\text{and}\quad\mf{F}(\mf{x}+\mf{h})=\mf{y}+\mf{k}.\]
 Eq. \eqref{230808_14} implies that
 \begin{equation}\label{230808_15}\Vert\mf{h}\Vert\leq \frac{1}{c}\Vert\mf{k}\Vert.\end{equation}
 
 Let $A=\mf{DF}(\mf{x})$. By assumption, $A$ is invertible. Notice that
 \begin{align*}\quad \mf{F}^{-1}(\mf{y}+\mf{k})-\mf{F}^{-1}(\mf{y})-A^{-1}\mf{k}
 & =-A^{-1}\left(\mf{k}-A\mf{h}\right)\\
&=-A^{-1}\left(\mf{F}(\mf{x}+\mf{h})-\mf{F}(\mf{x})-A\mf{h}\right).\end{align*}

There is a positive constant $\beta$ such that
\[\Vert A^{-1}\mf{y}\Vert \leq \beta\Vert \mf{y}\Vert\hspace{1cm}\text{for all}\;\mf{y}\in\mb{R}^n.\]
Therefore,
\begin{equation}\label{230808_16}\begin{split}
&\left\Vert \frac{\mf{F}^{-1}(\mf{y}+\mf{k})-\mf{F}^{-1}(\mf{y})-A^{-1}\mf{k}}{\Vert\mf{k}\Vert}\right\Vert
\\& \leq\frac{ \beta }{\Vert\mf{k}\Vert}\Vert \left(\mf{F}(\mf{x}+\mf{h})-\mf{F}(\mf{x})-A\mf{h}\right)\Vert\\
&\leq \frac{ \beta }{c}\left\Vert \frac{\mf{F}(\mf{x}+\mf{h})-\mf{F}(\mf{x})-A\mf{h}}{\Vert\mf{h}\Vert}\right\Vert.\end{split}
\end{equation}Since $\mf{F}$ is differentiable at $\mf{x}$, 
\[\lim_{\mf{h}\to\mf{0}}\frac{\mf{F}(\mf{x}+\mf{h})-\mf{F}(\mf{x})-A\mf{h}}{\Vert\mf{h}\Vert}=\mf{0}.\]
\bp
Eq. \eqref{230808_15} implies that $\di\lim_{\mf{k}\to\mf{0}}\mf{h}=\mf{0}$. 
 Eq. \eqref{230808_16} then implies that
\[\lim_{\mf{k}\to\mf{0}}\frac{\mf{F}^{-1}(\mf{y}+\mf{k})-\mf{F}^{-1}(\mf{y})-A^{-1}\mf{k}}{\Vert\mf{k}\Vert}=\mf{0}.\]
This proves that $\mf{F}^{-1}$ is differentiable at $\mf{y}$ and
\[\mf{DF}^{-1}(\mf{y})=A^{-1}=\mf{DF}(\mf{x})^{-1}.\]
Now the map $\mf{DF}^{-1}:V\to\text{GL}\,(n,\mb{R})$ is the compositions of the maps 
$\mf{F}^{-1}:V\to U$, $\mf{DF}:U\to\text{GL}\,(n,\mb{R})$  and $\mathscr{I}:\text{GL}\,(n,\mb{R})\to\text{GL}\,(n,\mb{R})$ which takes $A$ to $A^{-1}$. Since each of these maps is continuous,   the map $\mf{DF}^{-1}:V\to\text{GL}\,(n,\mb{R})$ is continuous. This completes the proof that $\mf{F}^{-1}:V\to U$ is continuously differentiable.
\end{myproof}

At the end of this section, let us give a brief discussion about the concept of homeomorphism and diffeomorphism.
\begin{definition}{Homeomorphism}
Let $A$ be a subset of $\mb{R}^m$ and let $B$ be a subset of $\mb{R}^n$. We say that $A$ and $B$ are homeomorphic if there exists a continuous bijective function $\mf{F}:A\to B$ whose inverse $\mf{F}^{-1}:B\to A$ is also continuous. Such a function $\mf{F}$ is called a  homeomorphism between $A$ and $B$.
\end{definition}
\begin{definition}{Diffeomorphism}
Let $\mathcal{O}$ and $\mathcal{U}$ be open subsets of $\mb{R}^n$.  We say that $\mathcal{U}$ and $\mathcal{O}$ are diffeomorphic if there exists a homeomorphism $\mf{F}:\mathcal{O}\to \mathcal{U}$ between $\mathcal{O}$ and $\mathcal{U}$ such that $\mf{F}$ and $\mf{F}^{-1}$ are differentiable.
\end{definition}
\begin{example}[label=230811_1]{}
Let $A=\left\{(x,y)\,|\, x^2+y^2<1\right\}$ and $B=\left\{(x,y)\,|\, 4x^2+9y^2<36\right\}$. Define the map $\mf{F}:\mb{R}^2\to\mb{R}^2$ by
\[\mf{F}(x,y)=(3x, 2y).\]
\be
Then $\mf{F}$ is an invertible linear transformation with 
\[\mf{F}^{-1}(x,y)=\left(\frac{x}{3}, \frac{y}{2}\right).\]The mappings $\mf{F}$ and $\mf{F}^{-1}$ are continuously differentiable.
It is easy to show that $\mf{F}$ maps $A$ bijectively onto $B$. Hence, $\mf{F}:A\to B$ is a diffeomorphism between $A$ and $B$.

\end{example2}

\begin{figure}[ht]
\centering
\includegraphics[scale=0.2]{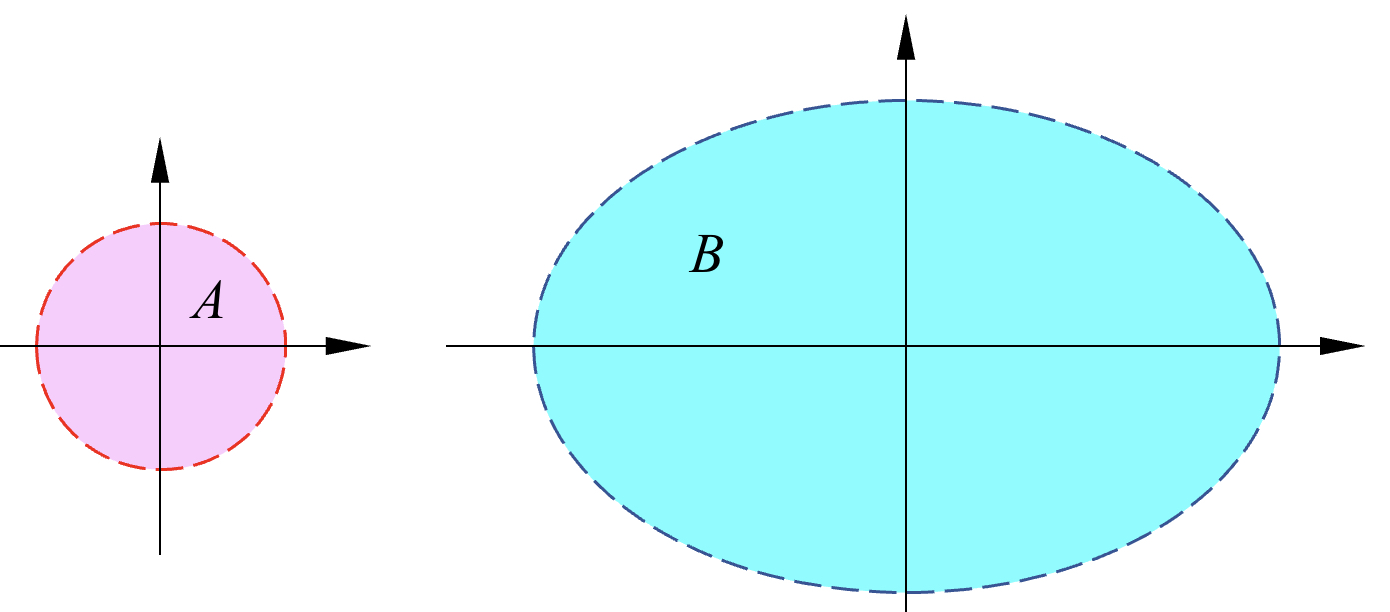}

\caption{  $A=\left\{(x,y)\,|\, x^2+y^2<1\right\}$ and $B=\left\{(x,y)\,|, 4x^2+9y^2<36\right\}$ are diffeomorphic.}\label{figure65}
\end{figure}

Theorem \ref{230808_8} gives the following.
\begin{theorem}{}
Let $A$ be an invertible $n\times n$ matrix, and let $\mf{x}_0$ and $\mf{y}_0$ be two points in $\mb{R}^n$. Define the mapping $\mf{F}:\mb{R}^n\to \mb{R}^n$ by 
\[\mf{F}(\mf{x})=\mf{y}_0+A\left(\mf{x}-\mf{x}_0\right).\]
If $\mathcal{O}$ is an open subset of $\mb{R}^n$, then $\mf{F}:\mathcal{O}\to\mf{F}(\mathcal{O})$ is a diffeomorphism. 
\end{theorem}

The inverse function theorem gives the following.
\begin{theorem}{}
Let $\mathcal{O}$ be an open subset of $\mb{R}^n$, and let $\mf{F}:\mathcal{O}\to\mb{R}^n$ be a continuously differentiable mapping such that $  \mathbf{DF}(\mf{x})   $ is invertible for all $\mf{x}\in\mathcal{O}$. If $\mathcal{U}$ is an open subset contained in $\mathcal{O}$ such that $\mf{F}:\mathcal{U}\to\mb{R}^n$ is one-to-one, then $\mf{F}:\mathcal{U}\to \mf{F}(\mathcal{U})$ is a diffeomorphism.
\end{theorem}
The proof of this theorem is left as an exercise.  
 
\vp
\noindent
{\bf \large Exercises  \thesection}
\setcounter{myquestion}{1}

 \begin{question}{\themyquestion}
Let $\mf{F}:\mb{R}^2\to\mb{R}^2$ be the mapping given by
\[\mf{F}(x,y)=(xe^y+xy, 2x^2+3y^2).\]
Show that there is a neighbourhood $U$ of $(-1, 0)$ such that the mapping  $\mf{F}:U\to\mb{R}^2$ is stable.
\end{question}
 \atc
 \begin{question}{\themyquestion}
Let $\mathcal{O}$ be an open subset of $\mb{R}^n$, and let $\mf{F}:\mathcal{O}\to\mb{R}^n$ be a continuously differentiable mapping such that $\det \mathbf{DF}(\mf{x})\neq 0$ for all $\mf{x}\in\mathcal{O}$. Show that $\mf{F}(\mathcal{O})$ is an open set.
\end{question}

\atc
 \begin{question}{\themyquestion}
Let $\mathcal{O}$ be an open subset of $\mb{R}^n$, and let $\mf{F}:\mathcal{O}\to\mb{R}^n$ be a continuously differentiable mapping such that $  \mathbf{DF}(\mf{x})   $ is invertible for all $\mf{x}\in\mathcal{O}$. If $\mathcal{U}$ is an open subset contained in $\mathcal{O}$ such that $\mf{F}:\mathcal{U}\to\mb{R}^n$ is one-to-one, then $\mf{F}:\mathcal{U}\to \mf{F}(\mathcal{U})$ is a diffeomorphism.

\end{question}

\atc
 \begin{question}{\themyquestion}
Let $\mathcal{O}$ be an open subset of $\mb{R}^n$, and let $\mf{F}:\mathcal{O}\to\mb{R}^n$ be a   differentiable mapping. Assume that there is a positive constant $c$ such that
\[\Vert\mf{F}(\mf{u})-\mf{F}(\mf{v})\Vert\geq c\Vert\mf{u}-\mf{v}\Vert\hspace{1cm}\text{for all}\;\mf{u},\mf{v}\in \mathcal{O}.\]
Use first order approximation theorem to show that for any $\mf{x}\in \mathcal{O}$ and any $\mf{h}\in \mb{R}^n$,
\[\Vert\mf{DF}(\mf{x})\mf{h}\Vert\geq c\Vert\mf{h}\Vert.\]
\end{question}
 
 \atc
 \begin{question}{\themyquestion}
Let $\mathcal{O}$ be an open subset of $\mb{R}^n$, and let $\mf{F}:\mathcal{O}\to\mb{R}^n$ be a continuously   differentiable mapping. 
\begin{enumerate}[(a)]
\item If $\mf{F}:\mathcal{O}\to\mb{R}^n$ is stable, show that the derivative matrix $\mf{DF}(\mf{x})$ is invertible at every $\mf{x}$ in $\mathcal{O}$.
\item Assume that the derivative matrix $\mf{DF}(\mf{x})$ is invertible at every $\mf{x}$ in $\mathcal{O}$. If $C$ is a compact subset of $\mathcal{O}$, show that the mapping 
$\mf{F}:C\to\mb{R}^n$ is stable.
\end{enumerate}
\end{question}
\section{The Implicit Function Theorem} \label{sec5.3}

The implicit function theorem is about the possibility of solving $m$ variables  from a system of $m$ equations with $n+m$ variables. Let us study some special cases.

Consider the function $f:\mb{R}^2\to\mb{R}$ given by $f(x,y)=x^2+y^2-1$. For a point $(x_0, y_0)$ that satisfies $f(x_0, y_0)=0$, we want to ask whether there is a neighbourhood $I$ of $x_0$, a neighbourhood $J$ of $y_0$, and a function $g:I\to\mb{R}$ such that for  $(x,y)\in I\times J$, 
$f(x,y)=0$  if and only if $y=g(x)$.

\begin{figure}[ht]
\centering
\includegraphics[scale=0.2]{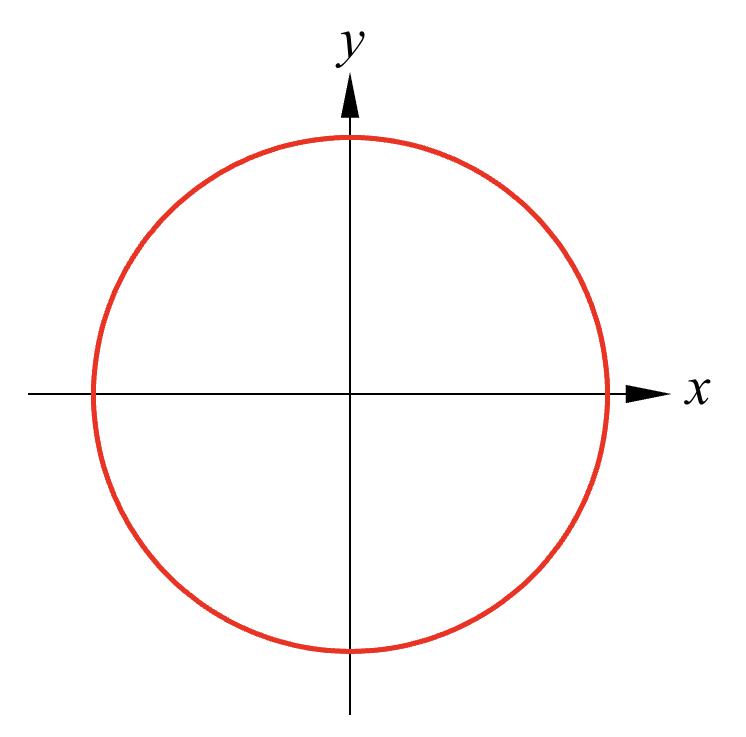}

\caption{The points in the $(x,y)$ plane satisfying $x^2+y^2-1=0$.}\label{figure63}
\end{figure}

If $(x_0, y_0)$ is a point with $y_0>0$ and $f(x_0, y_0)=0$, then we can take the neighbourhoods $I=(-1,1)$ and $J=(0,\infty)$ of $x_0$ and $y_0$ respectively, and define the function $g:I\to\mb{R}$ by
\[g(x)=\sqrt{1-x^2}.\]
We then find that for $(x,y)\in I\times J$, $f(x,y)=0$ if and only if $y=\sqrt{1-x^2}=g(x)$.

If $(x_0, y_0)$ is a point with $y_0<0$ and $f(x_0, y_0)=0$, then we can take the neighbourhoods $I=(-1,1)$ and $J=(-\infty,0)$ of $x_0$ and $y_0$ respectively, and define the function $g:I\to\mb{R}$ by
\[g(x)=-\sqrt{1-x^2}.\]
We then find that for $(x,y)\in I\times J$, $f(x,y)=0$ if and only if $y=-\sqrt{1-x^2}=g(x)$.

However, if $(x_0, y_0)=(1,0)$, any neighbourhood $J$ of $y_0$ must contain an interval of the form $(-r,r)$. If $I$ is a neighbourhood of $1$, $(x,y)$ is a point in $I\times (-r,r)$ such that $f(x,y)=0$, then $(x,-y)$ is another point in $I\times (-r,r)$  satisfying $f(x, -y)=0$. This shows that there does not exist any function $g:I\to\mb{R}$ such that when $(x,y)\in I\times J$, $f(x,y)=0$ if and only if $y=g(x)$. We say that we cannot solve $y$ as a function of $x$ in a neighbourhood of the point $(1,0)$.

Similarly, we cannot solve $y$ as a  function of $x$ in a neighbourhood of the point $(-1,0)$.

However, in a neighbourhood of the points $(1,0)$ and $(-1,0)$, we can solve $x$ as a function of $y$.

For a function $f:\mathcal{O}\to \mb{R}$ defined on an open subset $\mathcal{O}$ of $\mb{R}^2$, the implicit function theorem takes the following form.

\begin{theorem}{Dini's Theorem}
Let $\mathcal{O}$ be an open subset of $\mb{R}^2$ that contains the point $(x_0, y_0)$, and let $f:\mathcal{O}\to\mb{R}$ be a continuously differentiable function defined on $\mathcal{O}$ such that $f(x_0, y_0)=0$. If $\di \frac{\pa f}{\pa y}(x_0, y_0)\neq 0$, then there is a neighbourhood $I$ of $x_0$, a neighbourhood $J$ of $y_0$, and a continuously differentiable function $g:I\to J$ such that for any $(x,y)\in I\times J$, $f(x,y)=0$ if and only if $y=g(x)$. Moreover, for any $x\in I$, 
\[\frac{\pa f}{\pa x}(x, g(x))+\frac{\pa f}{\pa y}(x, g(x))g'(x)=0.\]
\end{theorem}
Dini's theorem says that to be able to solve $y$ as a function of $x$, a sufficient condition is that the function $f$ has continuous partial derivatives, and $f_y$ does not vanish. By interchanging the roles of $x$ and $y$, we see that if $f_x$ does not vanish, we can solve $x$ as a function of $y$.

For the function $f:\mb{R}^2\to\mb{R}$, $f(x,y)=x^2+y^2-1$, the points on the set $x^2+y^2=1$ which $f_y(x,y)=2y$ vanishes are the points $(1,0)$ and $(-1,0)$. In fact, we have seen that we cannot solve $y$ as functions of $x$ in neighbourhoods of these two points.

\begin{myproof}{Proof of Dini's Theorem}
 Without loss of generality, assume that $f_y(x_0, y_0)>0$. Let $\mf{u}_0=(x_0,y_0)$.
Since $f_y:\mathcal{O}\to\mb{R}$ is continuous, there is an $r_1>0$ such that the closed rectangle $R=[x_0-r_1, x_0+r_1]\times [y_0-r_1, y_0+r_1]$ lies in $\mathcal{O}$, and for all $(x,y)\in R$, $f_y(x,y)>f_y(x_0, y_0)/2>0$. For any $x\in [x_0-r_1, x_0+r_1]$, the function $h_x:[y_0-r_1, y_0+r_1]\to\mb{R}$ has derivative $h_x'(y)=\di f_y(x,y)$ that is positive. Hence, $h_x(y)=g(x,y)$ is strictly increasing in $y$. This implies that
\[f(x, y_0-r_1)<f(x, y_0)<f(x, y_0+r_1).\]
When $x=x_0$, we find that
\[f(x_0, y_0-r_1)<0<f(x_0, y_0+r_1).\]
Since $f$ is continuously differentiable, it is continuous. Hence, there is an $r_2>0$ such that $r_2\leq r_1$, and for all $x\in [x_0-r_2, x_0+r_2]$,
\[f(x,y_0-r_1)<0\quad\text{and}\quad f(x, y_0+r_1)>0.\] Let $I=(x_0-r_2, x_0+r_2)$. For $x\in I$, since $h_x:[y_0-r_1, y_0+r_1]\to\mb{R}$ is continuous, and
\[h_x(y_0-r_1)<0<h_x(y_0+r_1),\]
intermediate value theorem implies that there is a   $y \in (y_0-r_1, y_0+r_1)$ such that $h_x(y)=0$. Since $h_x$ is strictly increasing, this $y$ is unique, and  we denote  it by $g(x)$.  This defines the function $g:I\to\mb{R}$. Let $J=(y_0-r_1, y_0+r_1)$. By our argument, 
for each $x\in I$, $y=g(x)$  is a unique $y\in J$ such that $f(x,y)=0$. Thus, for any $(x,y)\in I\times J$, $f(x,y)=0$ if and only if $y=g(x)$. 

It remains to prove that $g:I\to\mb{R}$ is continuosly differentiable. By our convention above, there is a positive constant $c$ such that
\[\frac{\pa f}{\pa y}(x,y)\geq c\hspace{1cm}\text{for all}\;(x,y)\in I\times J.\]
\bp
Fixed $x\in I$. There exists an $r>0$ such that $(x-r,x+r)\subset I$. For $h$ satisfying $0<|h|<r$, $x+h$ is in $I$. 
By mean value theorem, there is a $c_h\in (0,1)$ such that  
\[f(x+h, g(x+h))-f(x,g(x))=h\frac{\pa f}{\pa x}(\mf{u}_h)+(g(x+h)-g(x))\frac{\pa f}{\pa y}(\mf{u}_h),\]
where
\begin{equation}\label{230808_20}\mf{u}_h=(x, g(x))+c_h(h, g(x+h)-g(x)).\end{equation}
Since
\[f(x+h, g(x+h))=0=f(x, g(x)),\]
we find that
\begin{equation}\label{230808_19} \frac{g(x+h)-g(x)}{h}=-\frac{f_x(\mf{u}_h)}{f_y(\mf{u}_h)}.\end{equation}
Since $f_x$ is continuous on the compact set $R$, it is bounded. Namely, there exists a constant $M$ such that
\[|f_x(x,y)|\leq M\hspace{1cm}\text{for all}\;(x,y)\in R.\]
Eq. \eqref{230808_19} then implies that
\[|g(x+h)-g(x)|\leq\frac{M}{c}|h|.\]
Taking $h\to 0$ proves that $g$ is continuous at $x$. From \eqref{230808_20}, we find that
\[\lim_{h\to 0}\mf{u}_h=(x, g(x)).\]
Since $f_x$ and $f_y$ are continuous at $(x, g(x))$, eq. \eqref{230808_19} gives
\[\lim_{h\to 0}\frac{g(x+h)-g(x)}{h}=-\lim_{h\to 0}\frac{f_x(\mf{u}_h)}{f_y(\mf{u}_h)}=-\frac{f_x(x,g(x))}{f_y(x,g(x))}.\]
This proves that $g$ is differentiable at $x$ and 
\[\frac{\pa f}{\pa x}(x, g(x))+\frac{\pa f}{\pa y}(x, g(x))g'(x)=0.\]
\end{myproof}

\begin{figure}[ht]
\centering
\includegraphics[scale=0.2]{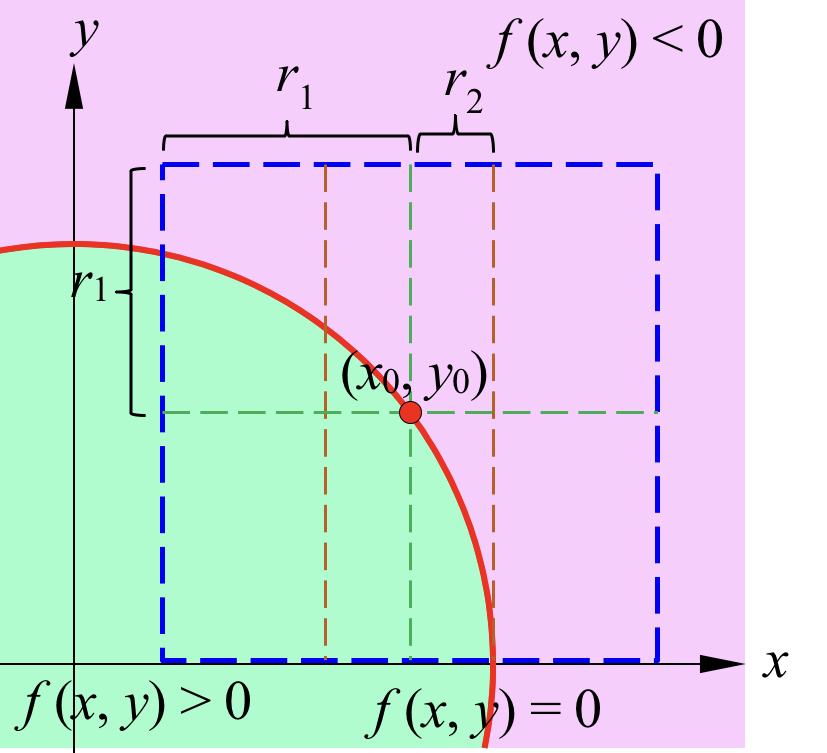}

\caption{Proof of Dini's Theorem.}\label{figure67}
\end{figure}

\begin{example}{}
Consider the equation

\vspace{-0.4cm}
\[xy^3+\sin(x+y)+4x^2y=3.\]
Show that in a neighbourhood of $(-1,1)$, this equation defines $y$ as a function of $x$. If this function is denoted as $y=g(x)$, find $g'(-1)$. 
\end{example}
\begin{solution}{Solution}
Let $f:\mb{R}^2\to\mb{R}$ be the function defined as

\vspace{-0.4cm}
\[f(x,y)=xy^3+\sin(x+y)+4x^2y-3.\]
Since sine function and polynomial functions are infinitely differentiable, $f$ is infinitely differentiable. 
\[\frac{\pa f}{\pa y}(x,y)=3xy^2+\cos(x+y)+4x^2,\hspace{1cm}
\frac{\pa f}{\pa y}(-1,1)=2\neq 0.\]
By Dini's theorem, there is a neighbourhood of $(-1,1)$ such that $y$ can be solved as a function of $x$. Now, 
\[\frac{\pa f}{\pa x}(x,y)=y^3+\cos(x+y)+8xy,\hspace{1cm}\frac{\pa f}{\pa x}(-1,1)=-6.\]
Hence,
$\di g'(0)=-\frac{-6}{2}=3$.
\end{solution}

Now we turn to the general case. First we consider polynomial mappings of degree at most one. Let $A=[a_{ij}]$ be an $m\times n$ matrix, and let $B=[b_{ij}]$ be an $m\times m$ matrix. Given $\mf{x}\in \mb{R}^n$, $\mf{y}\in \mb{R}^m$, $\mf{c}\in \mb{R}^m$, the system of equations
\[A\mf{x}+B\mf{y}=\mf{c}\] is 
 the following $m$ equations in $m+n$ variables $x_1, \ldots, x_n, y_1, \ldots, y_m$.
\begin{align*}
a_{11}x_1+a_{12}x_2+\cdots+a_{1n}x_n+b_{11}y_1+b_{12}y_2+\cdots+b_{1m}y_m&=c_1,\\
a_{21}x_1+a_{22}x_2+\cdots+a_{2n}x_n+b_{21}y_1+b_{22}y_2+\cdots+b_{2m}y_m&=c_2,\\
 \vdots\hspace{8cm} &\\
a_{m1}x_1+a_{m2}x_2+\cdots+a_{mn}x_n+b_{m1}y_1+b_{m2}y_2+\cdots+b_{mm}y_m&=c_m.
\end{align*}

Let us look at an example.
\begin{example}{}
Consider the linear system
\begin{align*}
2x_1+3x_2-5x_3+2y_1-y_2&=1\\
3x_1-x_2+2x_3-3y_1+ y_2&=0
\end{align*}Show that $\mf{y}=(y_1, y_2)$ can be solved as a function of $\mf{x}=(x_1, x_2, x_3)$. Write down the function $\mf{G}:\mb{R}^3\to\mb{R}^2$ such that the solution is given by $\mf{y}=\mf{G}(\mf{x})$, and find $\mf{DG}(\mf{x})$.
\end{example}
\begin{solution}{Solution}
Let 
\[A=\begin{bmatrix} 2 & 3 & -5\\3 & -1 & 2\end{bmatrix},\quad B=\begin{bmatrix} 2 & -1\\-3 & 1\end{bmatrix}.\]
Then the system can be written as
\[A\mf{x}+B\mf{y}=\mf{c},\hspace{1cm}\text{where}\;\mf{c}=\begin{bmatrix} 1\\ 0\end{bmatrix}.\]
This implies that
\begin{equation}\label{230809_5}B\mf{y}=\mf{c}-A\mf{x}.\end{equation}
\bs
For every $\mf{x}\in \mb{R}^3$, $\mf{c}-A\mf{x}$ is a vector in $\mb{R}^2$. Since $\det B=-1\neq 0$, $B$ is invertible. Therefore, there is a unique $\mf{y}$ satisfying \eqref{230809_5}. It is given by
 
\begin{align*}
\mf{G}(\mf{x})&=\mf{y}=B^{-1}\left(\mf{c}-A\mf{x}\right)\\
&=-\begin{bmatrix} 1 & 1\\3 & 2\end{bmatrix}\begin{bmatrix} 1\\ 0\end{bmatrix}+\begin{bmatrix} 1 & 1\\3 & 2\end{bmatrix}\begin{bmatrix} 2 & 3 & -5\\3 & -1 & 2\end{bmatrix}\mf{x}\\
&=- \begin{bmatrix} 1\\ 3\end{bmatrix}+ \begin{bmatrix} 5 & 2 & -3\\12 & 7 & -11\end{bmatrix}\mf{x}\\
&=\begin{bmatrix}5x_1+2x_2-3x_3-1\\12x_1+7x_2-11x_3-3 \end{bmatrix}.\end{align*}
 
It follows that $\mf{DG}=\begin{bmatrix}5 & 2& -3 \\12 & 7 &-11  \end{bmatrix}$.
\end{solution}

The following theorem gives a general scenario.
\begin{theorem}[label=230810_1]{}
Let $A=[a_{ij}]$ be an $m\times n$ matrix, and let $B=[b_{ij}]$ be an $m\times m$ matrix. 
Define the function $\mf{F}:\mb{R}^{m+n}\to\mb{R}^m$ by
\[\mf{F}(\mf{x}, \mf{y})=A\mf{x}+B\mf{y}-\mf{c},\]
where $\mf{c}$ is a constant vector in $\mb{R}^m$. 
The equation $\mf{F}(\mf{x},\mf{y})=\mf{0}$ defines the variable $\mf{y}=(y_1, \ldots, y_m)$ as a function of $\mf{x}=(x_1, \ldots, x_n)$ if and only if the matrix $B$ is invertible. If we denote this function as  $\mf{G}:\mb{R}^n\to \mb{R}^m$, then
\[\mf{G}(\mf{x})= B^{-1}\left(\mf{c}-A\mf{x}\right),\] and
\[\mf{DG}(\mf{x})=-B^{-1}A.\]
\end{theorem}
\begin{myproof}{Proof}
The equation  $\mf{F}(\mf{x},\mf{y})=\mf{0}$ defines the variables $\mf{y}$ as a function of $\mf{x}$ if and only for for each $\mf{x}\in \mb{R}^n$, there is a unique $\mf{y}\in \mb{R}^m$ satisfying 
\[B\mf{y} = \mf{c}-A\mf{x}.\]
This is a linear system for the variable $\mf{y}$. By the theory of linear algebra, a unique solution $\mf{y}$  exists if and only if $B$ is invertible. In this case, the solution is given by
\[\mf{y}= B^{-1}\left(\mf{c}-A\mf{x}\right).\]The rest of the assertion follows. 
\end{myproof}

 Write a point in $\mb{R}^{m+n}$ as $(\mf{x}, \mf{y})$, where $\mf{x}\in \mb{R}^n$ and $\mf{y}\in\mb{R}^m$. If  $\mf{F}:\mb{R}^{m+n}\to\mb{R}^m$ is a  function that is differentiable at the point   $(\mf{x}, \mf{y})$,   the $m\times (m+n)$ derivative matrix $\mf{DF}(\mf{x}, \mf{y})$ can be written as
 \[\mf{DF}(\mf{x}, \mf{y})=\begin{bmatrix} \mf{D}_{\mf{x}}\mf{F}(\mf{x}, \mf{y}) &\rvline & \mf{D}_{\mf{y}}\mf{F}(\mf{x}, \mf{y})\end{bmatrix},\]
 where
\begin{align*}
\mf{D}_{\mf{x}}\mf{F}(\mf{x}, \mf{y})&=\begin{bmatrix}\vspace{0.3cm} \di \frac{\pa F_1}{\pa x_1} (\mf{x}, \mf{y})& \di \frac{\pa F_1}{\pa x_2} (\mf{x}, \mf{y})  &\cdots & \di \frac{\pa F_1}{\pa x_n} (\mf{x}, \mf{y})\\\vspace{0.3cm} \di \frac{\pa F_2}{\pa x_1} (\mf{x}, \mf{y})& \di \frac{\pa F_2}{\pa x_2} (\mf{x}, \mf{y})  &\cdots & \di \frac{\pa F_2}{\pa x_n} (\mf{x}, \mf{y})\\
\vspace{0.3cm} \vdots & \vdots & \ddots & \vdots\\
 \di \frac{\pa F_m}{\pa x_1} (\mf{x}, \mf{y})& \di \frac{\pa F_m}{\pa x_2} (\mf{x}, \mf{y})  &\cdots & \di \frac{\pa F_m}{\pa x_n} (\mf{x}, \mf{y})\end{bmatrix},
\end{align*}
\begin{align*}
\mf{D}_{\mf{y}}\mf{F}(\mf{x}, \mf{y})&=\begin{bmatrix}\vspace{0.3cm} \di \frac{\pa F_1}{\pa y_1} (\mf{x}, \mf{y})& \di \frac{\pa F_1}{\pa y_2} (\mf{x}, \mf{y})  &\cdots & \di \frac{\pa F_1}{\pa y_m} (\mf{x}, \mf{y})\\\vspace{0.3cm} \di \frac{\pa F_2}{\pa y_1} (\mf{x}, \mf{y})& \di \frac{\pa F_2}{\pa y_2} (\mf{x}, \mf{y})  &\cdots & \di \frac{\pa F_2}{\pa y_m} (\mf{x}, \mf{y})\\
\vspace{0.3cm} \vdots & \vdots & \ddots & \vdots\\
 \di \frac{\pa F_m}{\pa y_1} (\mf{x}, \mf{y})& \di \frac{\pa F_m}{\pa y_2} (\mf{x}, \mf{y})  &\cdots & \di \frac{\pa F_m}{\pa y_m} (\mf{x}, \mf{y})\end{bmatrix}.
\end{align*}Notice that $\mf{D}_{\mf{y}}\mf{F}(\mf{x}, \mf{y})$ is a square matrix. 
 
When $A=[a_{ij}]$ is an $m\times n$ matrix,  $B=[b_{ij}]$ is an $m\times m$ matrix, $\mf{c}$ is a vector in $\mb{R}^m$, 
and $\mf{F}:\mb{R}^{m+n}\to\mb{R}^m$ is the function defined as
\[\mf{F}(\mf{x}, \mf{y})=A\mf{x}+B\mf{y}-\mf{c},\]
it is easy to compute that 
\[\mf{D}_{\mf{x}}\mf{F}(\mf{x}, \mf{y})=A,\hspace{1cm} \mf{D}_{\mf{y}}\mf{F}(\mf{x}, \mf{y})=B.\]
Theorem \ref{230810_1} says that we can solve $\mf{y}$ as a function of $\mf{x}$ from the system of $m$ equations
\[\mf{F}(\mf{x}, \mf{y})=\mf{0}\] if and only if 
\[B=\mf{D}_{\mf{y}}\mf{F}(\mf{x}, \mf{y})\] is invertible. In this case, if   $\mf{G}:\mb{R}^n\to\mb{R}^m$ is the function so that $\mf{y}=\mf{G}(\mf{x})$ is the solution, then
\[\mf{DG}(\mf{x})=-B^{-1}A=-\mf{D}_{\mf{y}}\mf{F}(\mf{x}, \mf{y})^{-1}\mf{D}_{\mf{x}}\mf{F}(\mf{x}, \mf{y}).\]
In fact, this latter follows  from $\di\mf{F}(\mf{x}, \mf{G}(\bf{x}))=\mf{0}$ and the chain rule. 

The special case of degree one polynomial mappings   gives us sufficient insight into the general implicit function theorem. However, for  nonlinear mappings, the conclusions can only be made {\it locally}. 

 \begin{theorem}{Implicit Function Theorem}
Let $\mathcal{O}$ be an open subset of $\mb{R}^{m+n}$, and let $\mf{F}:\mathcal{O}\to\mb{R}^m$ be a continuously differentiable function defined on $\mathcal{O}$.  Assume that $\mf{x}_0$ is a point in $ \mb{R}^n$ and $\mf{y}_0$ is a point in $\mb{R}^m$ such that the point $(\mf{x}_0, \mf{y}_0)$ is in $\mathcal{O}$ and  $\mf{F}(\mf{x}_0, \mf{y}_0)=\mf{0}$. If $\di \det \mf{D}_{\mf{y}}\mf{F} (\mf{x}_0, \mf{y}_0)\neq 0$, then we have the followings.
\begin{enumerate}[(i)]
\item There is a neighbourhood $U$ of $\mf{x}_0$, a neighbourhood $V$ of $\mf{y}_0$, and a continuously differentiable function $\mf{G}:U\to \mb{R}^m$ such that for any $(\mf{x},\mf{y})\in U\times V$, $\mf{F}(\mf{x}, \mf{y})=\mf{0}$ if and only if $\mf{y}=\mf{G}(\mf{x})$. 
\item For any $\mf{x}\in U$, 
\[\mf{D}_{\mf{x}}\mf{F}(\mf{x}, \mf{G}(\mf{x}))+\mf{D}_{\mf{y}}\mf{F}(\mf{x}, \mf{G}(\mf{x}))\mf{DG}(\mf{x})=\mf{0}.\]
\end{enumerate}
\end{theorem}

Here we will give a proof of the implicit function theorem using the inverse function theorem. The idea of the proof is to construct a mapping which one can apply the inverse function theorem. Let us look at an example first.
\begin{example}{}
Let $\mf{F}:\mb{R}^5\to\mb{R}^2$ be the function defined as
\[\mf{F}(x_1, x_3, x_3, y_1, y_2) =(x_1y_2^2, x_2x_3y_1^2+x_1y_2).\]
Define the mapping $\mf{H}:\mb{R}^5\to \mb{R}^{5}$ as
\[\mf{H}(\mf{x}, \mf{y})=\left(\mf{x}, \mf{F}(\mf{x}, \mf{y})\right)=(x_1, x_2, x_3, x_1y_2^2, x_2x_3y_1^2+x_1y_2).\]

Then we find that
\begin{align*}
\mf{DH}(\mf{x}, \mf{y})=\begin{bmatrix}\; \begin{matrix}1 & 0 & 0 &\rvline &   0 & 0\\
0 & 1 & 0 &\rvline &   0 & 0\\
0 & 0 & 1  &\rvline &   0 & 0\\
\hline
y_2^2 & 0 & 0&\rvline & 0 & 2x_1y_2\\y_2 & x_3y_1^2 & x_2y_1^2 &\rvline &
  2x_2x_3y_1 & x_1\end{matrix}\;\end{bmatrix}.
\end{align*}
Notice that
\[\mf{DH}(\mf{x}, \mf{y})=\begin{bmatrix}\;\begin{matrix} I_3 &\rvline & \mf{0}\\
\hline
  \mf{D}_{\mf{x}}\mf{F}(\mf{x},\mf{y}) &\rvline &  \mf{D}_{\mf{y}}\mf{F}(\mf{x},\mf{y})\end{matrix}\;\end{bmatrix}.\]
\end{example}

\begin{myproof}
{Proof of the Implicit Function Theorem}
  Let $\mf{H}:\mathcal{O}\to \mb{R}^{m+n}$ be the mapping defined as
\[\mf{H}(\mf{x}, \mf{y})=\left(\mf{x}, \mf{F}(\mf{x}, \mf{y})\right).\]
Notice that $\mf{F}(\mf{x}, \mf{y})=\mf{0}$ if and only if $\mf{H}(\mf{x}, \mf{y})=(\mf{x}, \mf{0})$.
Since the first $n$ components of $\mf{H}$ are infinitely differentiable functions, the mapping $\mf{H}:\mathcal{O}\to \mb{R}^{m+n}$ is continuously differentiable. 
\bp
Now, 
\begin{align*}
\mf{DH}(\mf{x}, \mf{y})=\begin{bmatrix}\;\begin{matrix} I_n &\rvline & \mf{0}\\
\hline
  \mf{D}_{\mf{x}}\mf{F}(\mf{x},\mf{y}) &\rvline &  \mf{D}_{\mf{y}}\mf{F}(\mf{x},\mf{y})\end{matrix}\;\end{bmatrix}. 
\end{align*}

Therefore,
\[\det \mf{DH}(\mf{x}_0, \mf{y}_0)=\det \mf{D}_{\mf{y}}\mf{F}(\mf{x}_0,\mf{y}_0)\neq 0.\]

By the inverse function theorem, there is a neighbourhood $W$ of $(\mf{x}_0, \mf{y}_0)$ and a neighbourhood $Z$ of $\mf{H}(\mf{x}_0, \mf{y}_0)=(\mf{x}_0, \mf{0})$ such that $\mf{H}:W\to Z$ is a bijection and $\mf{H}^{-1}:Z\to W$ is continuously differentiable. 
For $\mf{u}\in\mb{R}^n$, $\mf{v}\in \mb{R}^m$ so that $(\mf{u}, \mf{v})\in Z$, let
\[\mf{H}^{-1}(\mf{u}, \mf{v})=(\mf{\Phi}(\mf{u}, \mf{v}), \mf{\Psi}(\mf{u}, \mf{v})),\]where $\mf{\Phi}$ is a map from $Z$ to $\mb{R}^n$ and $\mf{\Psi}$ is a map from $Z$ to $\mb{R}^m$. Since $\mf{H}^{-1}$ is continuously differentiable, $\mf{\Phi}$ and $\mf{\Psi}$ are continuously differentiable.

Given $r>0$, let $D_r$ be the open cube $\di D_r=\prod_{i=1}^{m+n}(-r,r)$. 
Since $W$ and $Z$ are  open sets that contain $(\mf{x}_0,\mf{y}_0)$ and $(\mf{x}_0, \mf{0})$ respectively, there exists $r>0$ such that
\[(\mf{x}_0,\mf{y}_0)+D_r\subset W,\hspace{1cm}(\mf{x}_0,\mf{0})+D_r\subset Z.\]
 If $\di A_r=\prod_{i=1}^{n}(-r,r)$, $\di B_r=\prod_{i=1}^{m}(-r,r)$, $U=\mf{x}_0+A_r$, $V=\mf{y}_0+B_r$, then
\[(\mf{x}_0,\mf{y}_0)+D_r=U\times V,\hspace{1cm}(\mf{x}_0,\mf{0})+D_r=U\times B_r.\]
 Hence, $U\times V\subset W$ and $U\times B_r\subset Z$.
Define $\mf{G}:U\to \mb{R}^m$ by
\[\mf{G}(\mf{x})=\mf{\Psi}(\mf{x},\mf{0}).\]Since $\mf{\Psi}$ is continuously differentiable, $\mf{G}$ is continuously differentiable. If $\mf{x}\in U$, $\mf{y}\in V$, then $(\mf{x},\mf{y})\in W$. For such $(\mf{x}, \mf{y})$,   $\mf{F}(\mf{x}, \mf{y})=\mf{0}$ implies $\mf{H}(\mf{x}, \mf{y})=(\mf{x}, \mf{0})$. Since $\mf{H}:W\to Z$ is a bijection, $(\mf{x}, \mf{0})\in Z$ and $\mf{H}^{-1}(\mf{x}, \mf{0})=(\mf{x}, \mf{y})$. Comparing the last $m$ components give
\[\mf{y}=\mf{\Psi}(\mf{x},\mf{0})=\mf{G}(\mf{x}).\]
 
 \bp
Conversely,  
since $\mf{H}(\mf{H}^{-1}(\mf{u},\mf{v}))=(\mf{u},\mf{v})$ for all $(\mf{u}, \mf{v})\in Z$, we find that
\[\left(\mf{\Phi}(\mf{u},\mf{v}), \mf{F}(\mf{\Phi}(\mf{u},\mf{v}), \mf{\Psi}(\mf{u},\mf{v}))\right)=(\mf{u}, \mf{v})\] for all $(\mf{u}, \mf{v})\in Z$. For all $\mf{u}\in U$, $(\mf{u}, \mf{0})$ is in $Z$. Therefore, 
\[\mf{\Phi}(\mf{u},\mf{0})=\mf{u}, \quad \mf{F}(\mf{\Phi}(\mf{u},\mf{0}), \mf{\Psi}(\mf{u},\mf{0})) =\mf{0}.\]
This implies that if $\mf{x}\in U$, then
$\di \mf{F}(\mf{u}, \mf{G}(\mf{u}))=\mf{0}$.
In other words, if $(\mf{x}, \mf{y})$ is in $U\times V$ and $\mf{y}=\mf{G}(\mf{x})$, we must have $\mf{F}(\mf{x},\mf{y})=\mf{0}$.
Since we have shown that $\mf{G}:U\to \mb{R}^m$ is continuously differentiable, the formula \[\mf{D}_{\mf{x}}\mf{F}(\mf{x}, \mf{G}(\mf{x}))+\mf{D}_{\mf{y}}\mf{F}(\mf{x}, \mf{G}(\mf{x}))\mf{DG}(\mf{x})=\mf{0} \]follows from $\mf{F}(\mf{x}, \mf{G}(\mf{x}))=\mf{0}$ and the chain rule.
\end{myproof}

\begin{example}[label=230810_7]{}
Consider the system of equations
\begin{equation}\label{231029_1}\begin{split}
2x^2y+3xy^2u+xyv+uv&=7\\
4xu-5yv+u^2y+v^2x&=1\end{split}
\end{equation}Notice that when $(x,y)=(1,1)$, $(u,v)=(1,1)$ is a solution of this system. Show that there  are neighbourhoods $U$ and $V$ of $(1,1)$, and a continuously differentiable function $\mf{G}:U\to\mb{R}^2$ such that if $(x,y, u,v)\in U\times V$, then $(x,y,u,v)$ is a solution of the system of equations above if and only if $u=G_1(x,y)$ and $v=G_2(x,y)$. Also, find the values of $\di \frac{\pa G_1}{\pa x}(1,1)$, $\di \frac{\pa G_1}{\pa y}(1,1)$, $\di \frac{\pa G_2}{\pa x}(1,1)$ and $\di \frac{\pa G_2}{\pa y}(1,1)$.
\end{example}
\begin{solution}{Solution}
Define the function $\mf{F}:\mb{R}^4\to \mb{R}^2$ by
\[\mf{F}(x,y,u,v)=(2x^2y+3xy^2u+xyv+uv-7, 4xu-5yv+u^2y+v^2x-1).\]
\bs
This is a polynomial mapping. Hence, it is continuously differentiable. It is easy to check that $\mf{F}(1,1,1,1)=\mf{0}$.
Now,
\[\mf{D}_{(u,v)}\mf{F}(x,y,u,v)=\begin{bmatrix} 3xy^2+v & xy+u\\
4x+2uy & -5y+2vx \end{bmatrix}.\]
Thus,
\[\det \mf{D}_{(u,v)}\mf{F}(1,1,1,1)=\begin{bmatrix}4 & 2\\
6 & -3\end{bmatrix}=-24\neq 0.\]
By implicit function theorem, there  are neighbourhoods $U$ and $V$ of $(1,1)$, and a continuously differentiable function $\mf{G}:U\to\mb{R}^2$ such that, if $(x,y, u,v)\in U\times V$, then $(x,y,u,v)$ is a solution of the system of equations \eqref{231029_1}  if and only if $u=G_1(x,y)$ and $v=G_2(x,y)$.

Finally, 
\[\mf{D}_{(x,y)}\mf{F}(x,y,u,v)=\begin{bmatrix} 4xy+3y^2u+yv & 2x^2+6xyu+xv \\
4 u+v^2 &  -5 v+u^2 \end{bmatrix}, \]\[ \mf{D}_{(x,y)}\mf{F}(1,1,1,1)=\begin{bmatrix} 8 & 9 \\
5 &  -4 \end{bmatrix}.\]Chain rule gives
\begin{align*}
\mf{DG}(1,1)&=-\mf{D}_{(u,v)}\mf{F}(1,1,1,1)^{-1}\mf{D}_{(x,y)}\mf{F}(1,1,1,1)\\
&=\frac{1}{24}\begin{bmatrix}-3 & -2\\
-6 & 4\end{bmatrix}\begin{bmatrix} 8 & 9 \\
5 &  -4 \end{bmatrix}\\
&=\frac{1}{24}\begin{bmatrix}-34 & -19\\
-28 & -70\end{bmatrix}.
\end{align*}
Therefore,
\[\frac{\pa G_1}{\pa x}(1,1)=-\frac{17}{12},\hspace{1cm} \frac{\pa G_1}{\pa y}(1,1)=-\frac{19}{24}, \]\[\frac{\pa G_2}{\pa x}(1,1)=-\frac{7}{6},\hspace{1cm} \frac{\pa G_2}{\pa y}(1,1)=-\frac{35}{12}.\]
\end{solution}

\begin{remark}{The Rank of a Matrix}
In the formulation of the implicit function theorem, the assumption that $\det\mf{D}_{\mf{y}}\mf{F}(\mf{x}_0,\mf{y}_0)\neq 0$ can be replaced by the assumption that there are $m$ variables $u_1, \ldots, u_m$ among the $n+m$ variables $x_1, \ldots, x_n, y_1, \ldots, y_m$ such that $\det\mf{D}_{(u_1, \ldots, u_m)}\mf{F}(\mf{x}_0,\mf{y}_0)\neq 0$.

Recall that the rank $r$ of an $m\times k$ matrix $A$ is the dimension of its row space or the dimension of its column space. Thus, the rank $r$ of a $m\times k$ matrix $A$  is the maximum number of column vectors of $A$ which are linearly independent, or the maximum number of row vectors of $A$ that are linearly independent. Hence, the maximum possible value of $r$ is $\max\{m, k\}$. If $r=\max\{m,k\}$, we say that the matrix $A$ has maximal rank. For a $m\times k$ matrix where $m\leq k$, it has maximal rank if $r=m$. In this case, there is a $m\times m$ submatrix of $A$ consists of $m$ linearly independent vectors in $\mb{R}^m$. The determinant of this submatrix is nonzero. 

Thus, the condition  $\det\mf{D}_{\mf{y}}\mf{F}(\mf{x}_0,\mf{y}_0)\neq 0$ in the formulation of the implicit function theorem can be replaced by the condition that the  $m\times (m+n)$ matrix $ \mf{DF}(\mf{x}_0,\mf{y}_0)$ has maximal rank.
\end{remark}
\begin{example}[label=230810_8]{}
Consider the system
\begin{equation}\label{231029_2}\begin{split}
2x^2y+3xy^2u+xyv+uv&=7\\
4xu-5yv+u^2y+v^2x&=1
\end{split}\end{equation} defined in Example  \ref{230810_7}.  Show that there  are neighbourhoods $U$ and $V$ of $(1,1)$, and a continuously differentiable function $\mf{H}:V\to\mb{R}^2$ such that if $(x,y, u,v)\in U\times V$, then $(x,y,u,v)$ is a solution of the system of equations  if and only if $x=H_1(u,v)$ and $y=H_2(u,v)$. Find $\mf{DH}(1,1)$.
\end{example}
\begin{solution}{Solution}
Define the function $\mf{F}:\mb{R}^4\to\mb{R}^2$ as in the solution of Example  \ref{230810_7}. Since
\[\det \mf{D}_{(x,y)}\mf{F}(1,1,1,1)=\begin{bmatrix} 8 & 9 \\
5 &  -4\end{bmatrix}=-77\neq 0,\]the implicit function theorem implies there  are neighbourhoods $U$ and $V$ of $(1,1)$, and a continuously differentiable function $\mf{H}:V\to\mb{R}^2$ such that if $(x,y, u,v)\in U\times V$, then $(x,y,u,v)$ is a solution of the system of equations \eqref{231029_2} if and only if $x=H_1(u,v)$ and $y=H_2(u,v)$. Moreover,
\begin{align*}
\mf{DH}(1,1)&=-\mf{D}_{(x,y)}\mf{F}(1,1,1,1)^{-1}\mf{D}_{(u,v)}\mf{F}(1,1,1,1)\\
&=\frac{1}{77}\begin{bmatrix}-4 & -9\\
-5 & 8\end{bmatrix}\begin{bmatrix}4 & 2\\
6 & -3\end{bmatrix}\\
&=\frac{1}{77}\begin{bmatrix}-70 & 19\\
28 & -34\end{bmatrix}.
\end{align*}
\end{solution}

\begin{remark}{}
The function $\mf{G}:U\to \mb{R}^2$ in Example \ref{230810_7} and the function $\mf{H}:V\to\mb{R}^2$ in Example \ref{230810_8} are in fact inverses of each other. 

Notice that $\mf{DG}(1,1)$ is invertible. By the inverse function theorem,  there is a neighbourhood $U'$ of $(1,1)$ such that $V'=\mf{G}(U)$ is open, and $\mf{G}:U'\to V'$ is a bijection with continuously differentiable inverse. By shrinking down the sets $U$ and $V$, we can assume that $U=U'$, and $V=V'$. If $(x,y)\in U$ and $(u,v)\in V$, $\mf{F}(x,y,u,v)=0$ if and only if $(u,v)=\mf{G}(x,y)$, if and only if $(x,y)=\mf{H}(u,v)$. This implies that $\mf{G}:U\to V$ and $\mf{H}:V\to U$ are inverses of each other.
\end{remark}

At the end of this section, let us consider a geometric  application of the implicit function theorem.
First let us revisit the example where $f(x,y)=x^2+y^2-1$. At each point $(x_0, y_0)$ such that $f(x_0, y_0)=0$, 
\[x_0^2+y_0^2=1.\]
Hence,   $\nabla f(x_0, y_0)= (2x_0, 2y_0)\neq \mf{0}$. Notice that the vector  $\nabla f(x_0, y_0)=(2x_0, 2y_0)$ is normal to the circle $x^2+y^2=1$ at the point $(x_0, y_0)$.

\begin{figure}[ht]
\centering
\includegraphics[scale=0.2]{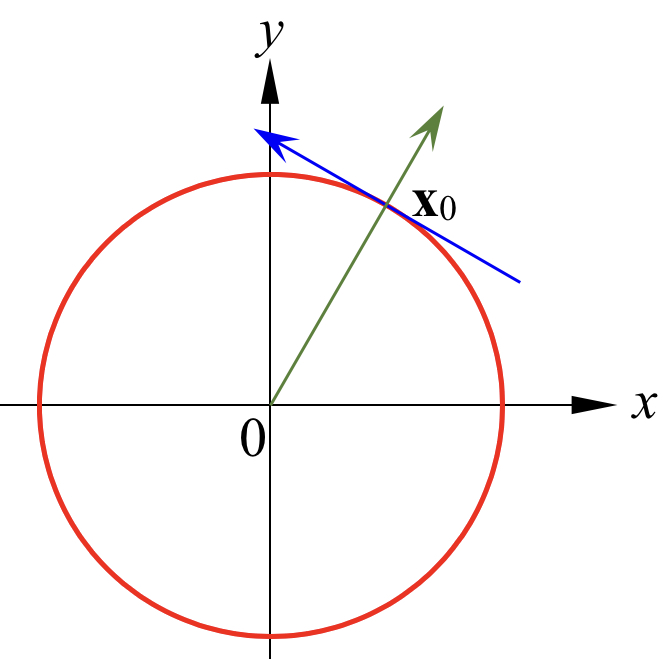}

\caption{The tangent vector and normal vector at a point on the circle $x^2+y^2-1=0$.}\label{figure64}
\end{figure}

If $y_0>0$, let   $U=(-1,1)\times (0,\infty)$. Restricted to $U$,  the points where $f(x,y)=0$ is the graph of the function $g:(-1,1)\to \mb{R}$, $g(x)=\sqrt{1-x^2}$.

If $y_0<0$, let   $U=(-1,1)\times (-\infty,0)$. Restricted to $U$,  the points where $f(x,y)=0$ is the graph of the function $g:(-1,1)\to \mb{R}$, $g(x)=-\sqrt{1-x^2}$. 

If $y_0=0$, then $x_0=1$ or $-1$. In fact, we can consider more generally the cases where $x_0>0$ and $x_0<0$.

If $x_0>0$, let   $U= (0,\infty)\times (-1,1)$. Restricted to $U$,  the points where $f(x,y)=0$ is the graph of the function $g:(-1,1)\to \mb{R}$, $g(y)=\sqrt{1-y^2}$.

If $x_0<0$, let   $U= (-\infty,0)\times (-1,1)$. Restricted to $U$,  the points where $f(x,y)=0$ is the graph of the function $g:(-1,1)\to \mb{R}$, $g(y)=-\sqrt{1-y^2}$.

\begin{definition}{Surfaces}
Let $S$ be a subset of $\mb{R}^k$ for some positive integer $k$. We say that $S$ is a $n$-dimensional surface if for each $\mf{x}_0$ on $S$, there is an open subset $\mathcal{D}$ of $\mb{R}^{n}$, an open neighbourhood $\mathcal{U}$ of $\mf{x}_0$ in $\mb{R}^k$, and a one-to-one differentiable mapping $\mf{G}:\mathcal{D}\to \mb{R}^k$ such that $\mf{G}(\mathcal{D})\subset S$, $\mf{G}(\mathcal{D})\cap\mathcal{U}=S\cap \mathcal{U}$, and $\mf{DG}(\mf{u})$ has  rank $n$ at each $\mf{u}\in\mathcal{D}$.

\end{definition} 
 
\begin{example}[label=230811_3]{}
We claim that the $n$-sphere
\[S^n=\{(x_1, \ldots, x_n, x_{n+1})\,|\, x_1^2+\cdots+x_n^2+x_{n+1}^2=1\}\] is an $n$-dimensional surface. Let  $(a_1, \ldots, a_n, a_{n+1})$ be a point on $S^n$. Then at least one of the components $a_1, \ldots, a_n, a_{n+1}$ is nonzero. Without loss of generality, assume that $  a_{n+1}>0$. Let 
\[\mathcal{D}=\left\{(x_1, \ldots, x_n)\,|\, x_1^2+\cdots+x_n^2<1\right\},\]
\[\mathcal{U}=\left\{(x_1, \ldots, x_n, x_{n+1})\,|\, x_{n+1}>0\right\},\]
and define the mapping $\mf{G}:\mathcal{D}\to \mathcal{U}$ by
\[\mf{G}(x_1, \ldots, x_n)=\left(x_1, \ldots, x_n, \sqrt{1-x_1^2-\cdots-x_n^2}\right).\]
Then $\mf{G}$ is a differentiable mapping, $\mf{G}(\mathcal{D})\subset S^n$ and $\mf{G}(\mathcal{D})\cap\mathcal{U}=S^n\cap \mathcal{U}$. 
Now,
\[\mf{DG}(x_1, \ldots, x_n)=\begin{bmatrix} \;\begin{matrix} I_n\\\hline \mf{v}\end{matrix}\;\end{bmatrix},\]
where $\mf{v}=\nabla G_{n+1}(x_1,\ldots, x_n)$. Since the first $n$-rows of $\mf{DG}(x_1, \ldots, x_n)$ is the $n\times n$ identity matrix, it has rank $n$.
Thus, $S^n$ is an $n$-dimensional surface.
\end{example}

Generalizing Example \ref{230811_3}, we find that a large class of surfaces is provided by graphs of differentiable functions.
\begin{theorem}
{}
Let $\mathcal{D}$ be an open subset of $\mb{R}^n$, and let $g:\mathcal{D}\to\mb{R}$ be a differentiable mapping. Then the graph of $g$ given by
\[G_g=\left\{(x_1, \ldots, x_n, x_{n+1})\,|\, (x_1, \ldots, x_n)\in \mathcal{D}, x_{n+1}=g(x_1, \ldots, x_n)\right\},\]
is an $n$-dimensional surface.
\end{theorem}

A hyperplane in $\mb{R}^{n+1}$ is the set of points in $\mb{R}^{n+1}$ which satisfies an equation of the form
\[a_1x_1+\cdots+a_nx_n+a_{n+1}x_{n+1}=b,\]
where $\mf{a}=(a_1, \ldots, a_n, a_{n+1})$ is a nonzero vector in $\mb{R}^{n+1}$. By definition, if $\mf{u}$ and $\mf{v}$ are two points on the plane, then
\[\langle\mf{a}, \mf{u}-\mf{v}\rangle =\mf{0}.\]
This shows that $\mf{a}$ is a vector normal to the plane.

When $\mathcal{D}$ is an open subset of $\mb{R}^n$, and  $g:\mathcal{D}\to\mb{R}$ is a differentiable mapping, the graph $G_g$ of $g$ 
is an $n$-dimensional surface. If $\mf{u}=(u_1, \ldots, u_n)$ is a point on $\mathcal{D}$, $ (\mf{u}, g(\mf{u}))$ is a point on $G_g$, we have seen that the equation of the tangent plane at the point $ (\mf{u}, g(\mf{u}))$ is given by
\[x_{n+1}=f(\mf{u})+\sum_{i=1}^n\frac{\pa g}{\pa x_i}(\mf{u},  g(\mf{u}))(x_i-u_i).\]

Implicit function theorem gives the following.
\begin{theorem}[label=230810_9]{}
Let $\mathcal{O}$ be an open subset of $\mb{R}^{n+1}$, and let $f:\mathcal{O}\to\mb{R}$ be a continuously differentiable function. If $\mf{x}_0$ is a point in $\mathcal{O}$  such that $f(\mf{x}_0)=0$ and  $\nabla f(\mf{x}_0)\neq\mf{0}$, then there is neighbourhood $U$ of $\mf{x}_0$ contained in $\mathcal{O}$ such that restricted to $U$, $f(\mf{x})=0$ is the graph of a continuously differentiable function $g:\mathcal{D}\to\mb{R}$, and  $\nabla f(\mf{x})$ is a vector normal to the tangent plane of the graph at the point $\mf{x}$.
\end{theorem}

\begin{myproof}{Proof}
Assume that $\mf{x}_0=(a_1, \ldots, a_n, a_{n+1})$.
Since $\nabla f(\mf{x}_0)\neq 0$, there is a $1\leq k\leq n+1$ such that $\di \frac{\pa f}{\pa x_k}(\mf{x}_0)\neq 0$. Without loss of generality, assume that $k=n+1$. 

\bp
Given a point $\mf{x}=(x_1, \ldots, x_n, x_{n+1})$ in $\mb{R}^{n+1}$, let $\mf{u}=(x_1, \ldots, x_n)$ so that $\mf{x}=(\mf{u}, x_{n+1})$. By the implicit function theorem, there is a neighbourhood  $\mathcal{D}$ of $\mf{u}_0=(a_1, \ldots, a_{n})$, an $r>0$, and a continuously differentiable function $g:\mathcal{D}\to\mb{R}$ such that if $U=\mathcal{D}\times (a_{n+1}-r, a_{n+1}+r)$, $(\mf{u}, u_{n+1})\in U$, then   $f(\mf{u}, u_{n+1})=0$ if and only if $u_{n+1}=g(\mf{u})$. In other words, in the neighbourhood $U$ of $\mf{x}_0=(\mf{u}_0, a_{n+1})$, $f(\mf{u}, u_{n+1})=0$ if and only if $(\mf{u}, u_{n+1})$ is a point on the graph of the function $g$. The equation of the tangent plane at the point $ (\mf{u}, u_{n+1})$ is
\[x_{n+1}-u_{n+1}=\sum_{i=1}^{n}\frac{\pa g}{\pa x_i}(\mf{u})(x_i-u_i).\]

By chain rule,
\[\frac{\pa g}{\pa x_i}(\mf{u})=-\frac{\di \frac{\pa f}{\pa x_i}(\mf{u}, u_{n+1})}{\di\frac{\pa f}{\pa x_{n+1}}(\mf{u}, u_{n+1})}.\]
Hence, the equation of the tangent plane can be rewritten as
\[\sum_{i=1}^{n+1} (x_i-u_i)\frac{\pa f}{\pa x_i}(\mf{u}, u_{n+1})=0.\]
This shows that $\nabla f(\mf{u}, u_{n+1})$ is a vector normal to the tangent plane.
\end{myproof}

 \begin{example}{}
 Find the equation of the tangent plane to the surface
$x^2+4y^2+9z^2=36$ at the point $(6,1,-1)$.
 \end{example}
 \begin{solution}{Solution}
 Let $f(x,y,z)=x^2+4y^2+9z^2$. Then
$\nabla f(x,y,z)=(2x, 8y, 18z)$.
 It follows that $\nabla f(6,1,-1) =2(6,4,-9)$.
 Hence, the equation of the tangent plane to the surface at 
 $(6,1,-1)$ is
 \[6x+4y-9z=36+4+9=49.\]
 \end{solution}
\vp
\noindent
{\bf \large Exercises  \thesection}
\setcounter{myquestion}{1}

\begin{question}{\themyquestion}
Consider the equation
\[ 4yz^2+3xz^3-11xyz=14.\] 
Show that in a neighbourhood of $(-1,1, 2)$, this equation defines $z$ as a function of $(x, y)$. If this function is denoted as $z=g(x, y)$, find $\nabla g(-1, 1)$. 
\end{question}
\atc

\begin{question}{\themyquestion}
Consider the system of equations
\begin{align*}
2xu^2+vyz+3uv&=2\\
5x+7yzu-v^2&=1
\end{align*}
\begin{enumerate}[(a)]
\item Show that when $(x,y,z)=(-1,1,1)$, $(u,v)=(1,1)$ is a solution of this system. 
\item  Show that there  are neighbourhoods $U$ and $V$ of $(-1,1,1)$ and $(1,1)$, and a continuously differentiable function $\mf{G}:U\to\mb{R}^2$ such that, if $(x,y, z, u,v)\in U\times V$, then $(x,y,z,u,v)$ is a solution of the system of equations above if and only if $u=G_1(x,y, z)$ and $v=G_2(x,y, z)$. 
\item Find the values of $\di \frac{\pa G_1}{\pa x}(-1,1,1)$, $\di \frac{\pa G_2}{\pa x}(-1,1,1)$  and $\di \frac{\pa G_2}{\pa z}(-1, 1,1)$.
\end{enumerate}
\end{question}

\atc

\begin{question}{\themyquestion}
Let $\mathcal{O}$ be an open subset of $\mb{R}^{2n}$, and let $\mf{F}:\mathcal{O}\to \mb{R}^n$ be a continuously differentiable function. Assume that $\mf{x}_0$ and $\mf{y}_0$ are points in $\mb{R}^n$ such that $(\mf{x}_0, \mf{y}_0)$ is a point in $\mathcal{O}$, $\mf{F}(\mf{x}_0, \mf{y}_0)=\mf{0}$, and $\mf{D}_{\mf{x}}\mf{F}(\mf{x}_0, \mf{y}_0)$ and $\mf{D}_{\mf{y}}\mf{F}(\mf{x}_0, \mf{y}_0)$ are invertible. Show that there exist neighbourhoods $U$ and $V$ of $\mf{x}_0$ and $\mf{y}_0$, and a continuously differentiable bijective function $\mf{G}:U\to V$ such that, if $(\mf{x}, \mf{y})$ is in $U\times V$, $\mf{F}(\mf{x}, \mf{y})=\mf{0}$ if and only if $\mf{y}=\mf{G}(\mf{x})$.
\end{question}
 
\section{Extrema Problems and the Method of Lagrange Multipliers} 

Optimization problems are very important in our daily life and in mathematical sciences. Given a function $f:\mk{D}\to\mb{R}$, we would like to know whether it has a maximum value or a minimum value. In Chapter \ref{chapter3}, we have dicusssed the extreme value theorem, which asserts that a continuous function that is defined on a compact set must have   maximum and   minimum values. In Chapter \ref{chapter4}, we showed that if a function $f:\mk{D}\to\mb{R}$ has (local) extremum at an interior point $\mf{x}_0$ of its domain $\mk{D}$ and it is differentiable at $\mf{x}_0$, then $\mf{x}_0$ must be a stationary point. Namely, $\nabla f(\mf{x}_0)=\mf{0}$.

Combining these various results, we can formulate a strategy for solving a special type of optimization problems. 
Let us first consider the following example.
\begin{example}[label=230811_5]{}Let
\[K=\left\{(x,y)\,|\, x^2+4y^2\leq 100\right\},\] and let $f:K\to\mb{R}$ be the function defined as
\[f(x,y)=x^2+y^2.\]
Find the maximum and minimum values of $f:K\to\mb{R}$, and the points where these values appear.
\end{example}
\begin{solution}{Solution}
Let $g:\mb{R}^2\to \mb{R}$ be the function defined as $g(x,y)=x^2+4y^2-100$. It is a polynomial function. Hence, it is continuous. Since $K=g^{-1}((-\infty, 0])$ and $(-\infty, 0]$ is closed in $\mb{R}$, $K$ is a closed set. By a previous exercise, \[\mathcal{O}=\text{int}\,K=\left\{(x,y)\,|\, x^2+4y^2< 100\right\}\]and \[\mathcal{C}=\text{bd}\,K=\left\{(x,y)\,|\, x^2+4y^2= 100\right\}.\]

\bs
For any $(x,y)\in K$, $\Vert (x,y)\Vert^2=x^2+y^2\leq x^2+4y^2\leq 100$. Therefore, $K$ is bounded. Since $K$ is closed and bounded, and the function $f:K\to\mb{R}$, $f(x)=x^2+ y^2$ is continuous, extreme value theorem says that $f$ has maximum and minimum values.
These values  appear either in $\mathcal{O}$ or on $\mathcal{C}$.  

Since $f:\mathcal{O}\to\mb{R}$ is differentiable, if $(x_0, y_0) $ is an extremizer of $f:\mathcal{O}\to\mb{R}$, we must have $\nabla f(x_0, y_0)=(0,0)$, which gives $(x_0, y_0)=(0,0)$. 

The other candidates of extremizers are on $\mathcal{C}$. Therefore, we need to find the maximum and minimum values of $f(x,y)=x^2+y^2$ subject to the constraint $x^2+4y^2=100$. From $x^2+4y^2=100$, we find that $x^2=100-4y^2$, and $y$ can only take values in the interval $[-5, 5]$. Hence, we want to find the maximum and minimum values of $h:[-5,5]\to \mb{R}$,
\[h(y)=100-4y^2+y^2=100-3y^2.\]
When $y=0$, $h$ has maximum value $100$, and when $y=\pm 5$, it has minimum value $100-3\times 25=25$. Notice that when $y=0$, $x =\pm 10$; while when $y=\pm 5$, $x=0$.

Hence, we have five candidates for the extremizers of $f$. Namely, $\mf{u}_1=(0,0)$, $\mf{u}_2=(10,0)$, $\mf{u}_3=(-10,0)$, $\mf{u}_4=(0, 5)$ and $\mf{u}_5=(0,-5)$. The function values at these 5 points are
\[f(\mf{u}_1)=0,\quad f(\mf{u}_2)=f(\mf{u}_3)=100, \quad
f(\mf{u}_4)=f(\mf{u}_5)=25.\]
Therefore, the minimum value of $f:K\to\mb{R}$ is 0, and the maximum value is 100. The minimum value appears at the point $(0,0)\in\text{int}\,K$, while the maximum value appears at $(\pm 10,0)\in \text{bd}\,K$.
\end{solution}

\begin{figure}[ht]
\centering
\includegraphics[scale=0.2]{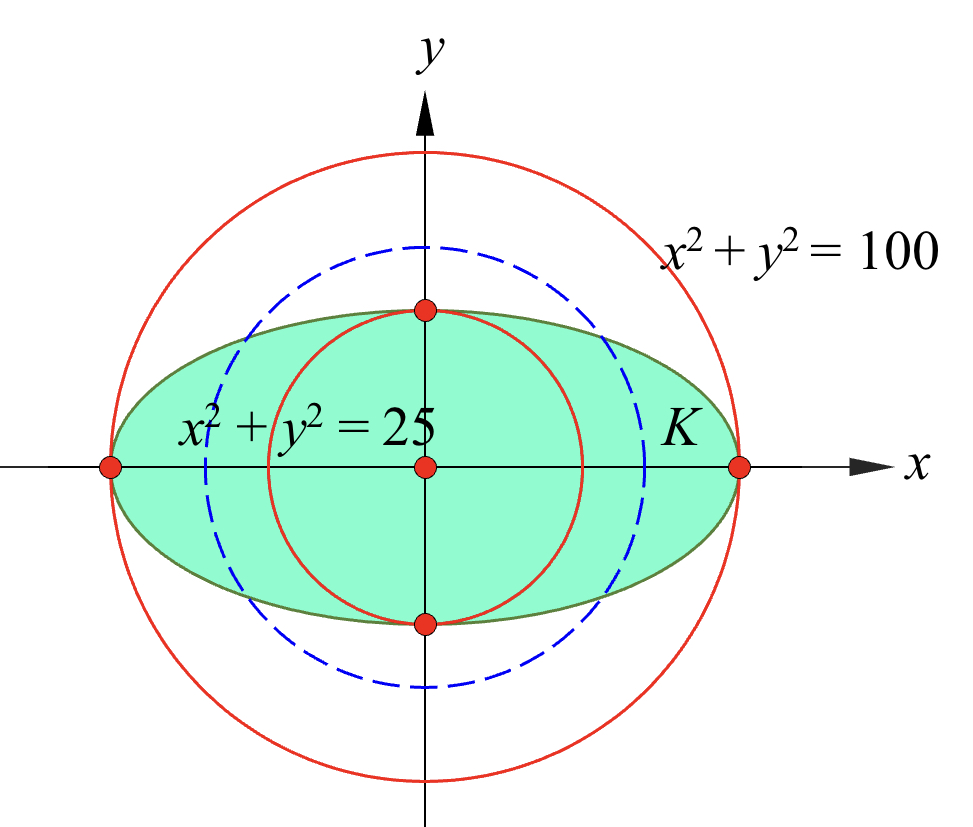}

\caption{The extreme values of $f(x,y)=x^2+y^2$ on the sets $K=\left\{(x,y)\,|\, x^2+4y^2\leq  100\right\}$ and $\mathcal{C}=\left\{(x,y)\,|\, x^2+4y^2= 100\right\}$.}\label{figure66}
\end{figure}

Example \ref{230811_5} gives a typical scenario of the optimization problems that we want to study in this section.

\begin{highlight}{Optimization Problem}
Let $K$ be a compact subset  of $\mb{R}^n$ with interior $\mathcal{O}$, and  let $f:K\to\mb{R}$ be a function continuous on $K$, differentiable on $\mathcal{O}$.
We want to  find the maximum and minimum values  of $f:K\to\mb{R}$.

\begin{enumerate}[(i)]
\item By the extreme value theorem, $f:K\to\mb{R}$ has maximum and minimum values.
\item Since $K$ is closed,   $K$ is a disjoint union of its interior $\mathcal{O}$ and its boundary $\mathcal{C}$. Since $\mathcal{C}$ is a subset of $K$, it is bounded. On the other hand, being the boundary of a set, $\mathcal{C}$ is   closed. Therefore, $\mathcal{C}$ is  compact.
\item The extreme values of $f$ can appear in $\mathcal{O}$ or on $\mathcal{C}$.
\item If $\mathbf{x}_0$ is an extremizer of $f:K\to\mb{R}$ and it is in $\mathcal{O}$, we must have $\nabla f(\mf{x}_0)=\mf{0}$. Namely, $\mf{x}_0$ is a stationary point of $f:\mathcal{O}\to\mb{R}$.
 \item If $\mf{x}_0$ is an extremizer of $f:K\to\mb{R}$ and it is not in $\mathcal{O}$, it is an extremizer of $f:\mathcal{C}\to\mb{R}$.
 \item  Since $\mathcal{C}$ is compact, $f:\mathcal{C}\to\mb{R}$ has maximum and minimum values.
\end{enumerate}
\end{highlight}
\begin{highlight}{}
Therefore, the steps to find the maximum and minimum values of $f:K\to\mb{R}$ are as follows.
\begin{enumerate} 
\item[~\hspace{0.5cm}\textbf{Step 1}] Find the stationary points of $f:\mathcal{O}\to\mb{R}$.
\item[\textbf{Step 2}]  Find the extremizers   of $f:\mathcal{C}\to\mb{R}$.
\item[\textbf{Step 3}]  Compare the values of $f$ at the stationary points of $f:\mathcal{O}\to\mb{R}$ and the extremizers of $f:\mathcal{C}\to\mb{R}$ to determine the extreme values of $f:K\to\mb{R}$.

\end{enumerate}
\end{highlight}
Of particular interest is when the boundary of $K$ can be expressed as $g(\mf{x})=0$, where $g:\mathcal{D}\to\mb{R}$ is a continuously differentiable function defined on an open subset $\mathcal{D}$ of $\mb{R}^n$. If $f$ is also defined and differentiable on $\mathcal{D}$,  the problem of finding the extreme values of $f:\mathcal{C}\to\mb{R}$ becomes finding the extreme values of $f:\mathcal{D}\to\mb{R}$ subject to the constraint $g(\mf{x})=0$. In Example \ref{230811_5}, we have used $g(\mf{x})=0$ to solve one of the variables in terms of the others and substitute into $f$ to transform the optimization problem to a problem with fewer variables. However, this strategy can be quite complicated because it is often not possible to solve one variable in terms of the others explicitly from the constraint $g(\mf{x})=0$. The method of Lagrange  multipliers   provides a way to solve constraint optimization problems without having to explicitly solve some variables in terms of the others. The validity of this method is justified by the implicit function theorem.

\begin{theorem}[label=230812_1]{The Method of Lagrange Multiplier (One Constraint)}
Let $\mathcal{O}$ be an open subset of $\mb{R}^{n+1}$ and let $f:\mathcal{O}\to\mb{R}$ and $g:\mathcal{O}\to\mb{R}$ be continuously differentiable functions defined on $\mathcal{O}$. Consider the subset of $\mathcal{O}$ defined as
\[\mathcal{C}=\left\{\mf{x}\in \mathcal{O}\,|\, g(\mf{x})=0\right\}.\]
If $\mf{x}_0$ is an extremizer of the function $f:\mathcal{C}\to\mb{R}$ and $\nabla g(\mf{x}_0)\neq\mf{0}$, then there is a constant $\lambda$, known as the Lagrange multiplier, such that
\[\nabla f(\mf{x}_0)=\lambda \nabla g(\mf{x}_0).\]
 
\end{theorem}
\begin{myproof}{Proof}
 Without loss of generality, assume that $\mf{x}_0$ is a maximizer of $f:\mathcal{C}\to\mb{R}$. Namely,
 \begin{equation}\label{230811_7}f(\mf{x})\leq f(\mf{x}_0)\hspace{1cm}\text{for all}\;\mf{x}\in\mathcal{C}.\end{equation}
 Given that $\nabla g(\mf{x}_0)\neq 0$, there exists a $1\leq k\leq n+1$ such that $\di\frac{\pa g}{\pa x_k}(\mf{x}_0)\neq 0$. Without loss of generality, assume that $k=n+1$. Let $\mf{x}_0=(a_1, \ldots, a_{n}, a_{n+1})$. Given a point $\mf{x}=(x_1, \ldots, x_{n}, x_{n+1})$ in $\mb{R}^{n+1}$,   let $\mf{u}=(x_1, \ldots, x_{n})$ so that $\mf{x}=(\mf{u}, x_{n+1})$.  By implicit function theorem, there is a neighbourhood $\mathcal{D}$ of $\mf{u}_0=(a_1, \ldots, a_{n})$, an $r>0$, and a continuously differentiable function $h:\mathcal{D}\to\mb{R}$  such that for $(\mf{u}, x_{n+1})\in \mathcal{D}\times (a_{n+1}-r,a_{n+1}+r)$, $g(\mf{u}, x_{n+1})=0$  if and only if $x_{n+1}=h(\mf{u})$. Consider the function $F:\mathcal{D}\to \mb{R}$ defined as
 \[F(\mf{u})=f(\mf{u}, h(\mf{u})).\]
 By \eqref{230811_7},  we find that
 \[F(\mf{u}_0) \geq F(\mf{u})\hspace{1cm}\text{for all}\; \mf{u}\in\mathcal{D}.\] In other words, $\mf{u}_0$ is a maximizer of the function $F:\mathcal{D}\to\mb{R}$. Since $\mf{u}_0$ is an interior point of $\mathcal{D}$ and $F:\mathcal{D}\to\mb{R}$ is continuously differentiable, $\nabla F(\mathbf{u}_0)=0$. Since $F(\mf{u})=f(\mf{u}, h(\mf{u}))$, we find that for $1\leq i\leq n$,
\begin{equation}\label{230811_8}\frac{\pa F}{\pa x_i}(\mf{u}_0)=\frac{\pa f}{\pa x_i}(\mf{u}_0, a_{n+1})+\frac{\pa f}{\pa x_{n+1}}(\mf{u}_0, a_{n+1})\frac{\pa h}{\pa x_i}(\mf{u}_0) =0.\end{equation}
On the other hand,  applying chain rule to $g(\mf{u}, h(\mf{u}))=0$ and set $\mf{u}=\mf{u}_0$, we find that
\begin{equation}\label{230811_9}\frac{\pa g}{\pa x_i}(\mf{u}_0, a_{n+1})+\frac{\pa g}{\pa x_{n+1}}(\mf{u}_0, a_{n+1})\frac{\pa h}{\pa x_i}(\mf{u}_0) =0\hspace{1cm}\text{for}\; 1\leq i\leq n.\end{equation}
By assumption, $\di \frac{\pa g}{\pa x_{n+1}}(\mf{x}_0)\neq 0$. Let
\[\lambda =\frac{\di \frac{\pa f}{\pa x_{n+1}}(\mf{x}_0)}{\di \frac{\pa g}{\pa x_{n+1}}(\mf{x}_0)}.\]
\bp
Then
\begin{equation}\label{230811_10} \frac{\pa f}{\pa x_{n+1}}(\mf{x}_0)=\lambda  \frac{\pa g}{\pa x_{n+1}}(\mf{x}_0).\end{equation}
Eqs. \eqref{230811_8} and \eqref{230811_9} show that for $1\leq i\leq n$,
\begin{equation}\label{230811_11}\frac{\pa f}{\pa x_i}(\mf{x}_0)=-\lambda \frac{\pa g}{\pa x_{n+1}}(\mf{x}_0)\frac{\pa h}{\pa x_i}(\mf{u}_0) =\lambda \frac{\pa g}{\pa x_i}(\mf{x}_0).\end{equation}
Eqs. \eqref{230811_10} and \eqref{230811_11} together imply that 
\[\nabla f(\mf{x}_0)=\lambda \nabla g(\mf{x}_0).\]
This completes the proof of the theorem.
\end{myproof}

\begin{remark}{}
Theorem \ref{230812_1} says that if $\mf{x}_0$ is an extremizer of the constraint optimization problem $\max/\min f(\mf{x})$ subject to $g(\mf{x})=0$, then the gradient of $f$ at $\mf{x}_0$ should be parallel to the gradient of $g$ at $\mf{x}_0$ if the latter is nonzero. One can refer to Figure \ref{figure66} for an illustration. Recall that the gradient of $f$ gives the direction where $f$ changes most rapidly, while the gradient of $g$ here represents the normal vector to the curve $g(\mf{x})=0$.

Using the method of Lagrange multiplier, there are $n+2$ variables $x_1, \ldots, x_{n+1}$ and $\lambda$ to be solved. The equation $\nabla f(\mf{x})=\lambda \nabla g(\mf{x})$ gives $n+1$ equations, while the equation $g(\mf{x})=0$ gives one. Therefore, we need to solve $n+2$ variables from $n+2$ equations.
\end{remark}

\begin{example}{}
Let us solve the constraint optimization problem that appears in Example \ref{230811_5} using the Lagrange multiplier method. Let $f:\mb{R}^2\to\mb{R}$ and $g:\mb{R}^2\to \mb{R}$ be respectively the functions $f(x,y)=x^2+y^2$ and $g(x,y)=x^2+4y^2-100$. They are both continuously differentiable. We want to find the maximum and minimum values of the function $f(x,y)$ subject to the constraint $g(x,y)=0$. Notice that
$\nabla g(x,y)=(2x, 8y)$ is the zero vector if and only if $(x,y)=(0,0)$, but $(0,0)$ is not on the curve $g(x,y)=0$. Hence, for any $(x,y)$ satisfying $g(x,y)=0$, $\nabla g(x, y)\neq \mf{0}$.

\be By the method of Lagrange multiplier, we need to find $(x,y)$ satisfying

\[\nabla f(x,y)=\lambda \nabla g(x,y)\quad\text{and}\quad g(x,y)=0.\]

Therefore,
\[2x=2\lambda x,\hspace{1cm} 2y=8\lambda y.\]
This gives
\[x(1-\lambda)=0,\quad y(1-4\lambda)=0.\]
The first equation says that either $x=0$ or $\lambda=1$.

If $x=0$, from $x^2+4y^2=100$, we must have $y=\pm 5$. 

If $\lambda=1$, then $y(1-4\lambda)=0$ implies that $y=0$. From $x^2+4y^2=100$, we then obtain $x=\pm 10$.

Hence, we find that the candidates for the extremizers are $(\pm 10, 0)$ and $(0, \pm 5)$. Since $f(\pm  10, 0)=100$ and $f(0, \pm 5)=25$, we conclude that subject to $x^2+4y^2=100$, the maximum value of $f(x,y)=x^2+y^2$ is 100, and the minimum value of $f(x,y)=x^2+y^2$ is 25.
\end{example2}

\begin{example}{}
Use the Lagrange multiplier method to find the maximum and minimum values of the function  $f(x,y,z)=8x+24y+27z$ on the set
\[S=\left\{(x,y,z)\,|\, x^2+4y^2+9z^2=289\right\},\]and the points where each of them appears.

\end{example}
\begin{solution}{Solution}
Let $g:\mb{R}^3\to\mb{R}$ be the function
\[g(x,y,z)=x^2+4y^2+9z^2-289.\]
The functions $f:\mb{R}^3\to\mb{R}$, $f(x,y,z)=8x+24y+27z$ and $g:\mb{R}^3\to\mb{R}$ are both continuously differentiable.
\bs
Notice that $\nabla g(x,y,z)=(2x, 8y, 18z)=\mf{0}$ if and only if $(x,y,z)=\mf{0}$, and $\mf{0}$ does not lie on $S$. By Lagrange multiplier method, to find the maximum and minimum values of $f:S\to\mb{R}$, we need to solve the equations
\[\nabla f(x,y,z)=\lambda \nabla g(x,y,z)\quad\text{and}\quad g(x,y,z)=0.\]

These give
\begin{gather*}
8=2\lambda x,\quad 
24=8\lambda y,\quad
27=18\lambda z\\
x^2+4y^2+9z^2=289.
\end{gather*}To satisfy the first three equations, none of the $\lambda$, $x$, $y$ and $z$ can be zero. We find that
\[x=\frac{4}{\lambda},\quad y=\frac{3}{\lambda},\quad z=\frac{3}{2\lambda}.\]
Substitute into the last equation, we have
\[\frac{64+144+81}{4\lambda^2}=289.\]
This gives $4\lambda^2=1$. Hence, $\lambda=\di \pm \frac{1}{2}$.
When $\lambda=\di \frac{1}{2}$, $(x,y,z)=(8, 6, 3)$. When $\lambda=-\di \frac{1}{2}$, $(x,y,z)=(-8, -6, -3)$. These are the two candidates for the extremizers of $f:S\to\mb{R}$.

Since $f(8,6,3)=289$ and $f(-8,-6,-3)=-289$, we find that the maximum and minimum values of $f:S\to\mb{R}$ are 289 and $-289$ respectively, and the maximum value appear at $(8,6,3)$, the minimum value appear at $(-8,-6,-3)$.
\end{solution}

Now we  consider   more general constraint optimization problems which can have more than one constraints.
\begin{theorem}[label=230812_7]{The Method of Lagrange Multiplier (General)}
Let $\mathcal{O}$ be an open subset of $\mb{R}^{m+n}$ and let $f:\mathcal{O}\to\mb{R}$ and $\mf{G}:\mathcal{O}\to\mb{R}^m$ be continuously differentiable functions defined on $\mathcal{O}$. Consider the subset of $\mathcal{O}$ defined as
\[\mathcal{C}=\left\{\mf{x}\in \mathcal{O}\,|\, \mf{G}(\mf{x})=\mf{0}\right\}.\]
If $\mf{x}_0$ is an extremizer of the function $f:\mathcal{C}\to\mb{R}$ and the matrix $\text{DG}(\mf{x}_0)$ has (maximal) rank $m$, then there are constants  $\lambda_1$, $\ldots$, $\lambda_m$, known as the Lagrange multipliers, such that  
\[\nabla f(\mf{x}_0)=\sum_{i=1}^m\lambda_i \nabla G_i(\mf{x}_0).\]

\end{theorem}
\begin{myproof}{Proof}
 Without loss of generality, assume that $\mf{x}_0$ is a maximizer of $f:\mathcal{C}\to\mb{R}$. Namely,
 \begin{equation}\label{230811_12}f(\mf{x})\leq f(\mf{x}_0)\hspace{1cm}\text{for all}\;\mf{x}\in\mathcal{C}.\end{equation}
 Given that the matrix $\mf{DG}(\mf{x}_0)$ has rank $m$, $m$ of the column vectors are linearly independent. Without loss of generality, assume that  the column vectors in the last $m$ columns are linearly independent.   Write a point $\mf{x}$ in $\mb{R}^{m+n}$ as  $\mf{x}=(\mf{u}, \mf{v})$, where $\mf{u}=(u_1, \ldots, u_{n})$ is in $\mb{R}^n$ and $\mf{v}=(v_{1}, \ldots, v_{m})$ is in $\mb{R}^m$. 
 By our assumption, $\mf{D}_{\mf{v}}\mf{G}(\mf{u}_0, \mf{v}_0)$ is invertible. 
  By implicit function theorem, there is a neighbourhood $\mathcal{D}$ of $\mf{u}_0$, a neighbourhood $\mathcal{V}$ of $\mf{v}_0$, and a continuously differentiable function $\mf{H}:\mathcal{D}\to\mb{R}^m$  such that for $(\mf{u}, \mf{v})\in \mathcal{D}\times \mathcal{V}$, $\mf{G}(\mf{u}, \mf{v})=\mf{0}$  if and only if $\mf{v}=\mf{H}(\mf{u})$. Consider the function $F:\mathcal{D}\to \mb{R}$ defined as
 \[F(\mf{u})=f(\mf{u}, \mf{H}(\mf{u})).\]
 By \eqref{230811_12},  we find that
 \[F(\mf{u}_0) \geq F(\mf{u})\hspace{1cm}\text{for all}\; \mf{u}\in\mathcal{D}.\] 
 
 \bp
 
 In other words, $\mf{u}_0$ is a maximizer of the function $F:\mathcal{D}\to\mb{R}$. Since $\mf{u}_0$ is an interior point of $\mathcal{D}$ and $F:\mathcal{D}\to\mb{R}$ is continuously differentiable, $\nabla F(\mathbf{u}_0)=0$. Since $F(\mf{u})=f(\mf{u}, \mf{H}(\mf{u}))$, we find that 
\begin{equation}\label{230811_13}\nabla F(\mf{u}_0)=D_{\mf{u}}f(\mf{u}_0,\mf{v}_0)+D_{\mf{v}}f(\mf{u}_0, \mf{v}_0)\mf{DH}(\mf{u}_0)=\mf{0}.\end{equation}

On the other hand,  applying chain rule to $\mf{G}(\mf{u}, \mf{H}(\mf{u}))=\mf{0}$ and set $\mf{u}=\mf{u}_0$, we find that
\begin{equation}\label{230811_14}\mf{D}_{\mf{u}}\mf{G}(\mf{u}_0, \mf{v}_0)+\mf{D}_{\mf{v}}\mf{G}(\mf{u}_0, \mf{v}_0)\mf{DH}(\mf{u}_0)=\mf{0}.\end{equation} Take 
\[\begin{bmatrix}\lambda_1 & \lambda_2 &\cdots &\lambda_m\end{bmatrix}=\boldsymbol{\lambda} =D_{\mf{v}}f(\mf{x}_0)\mf{D}_{\mf{v}}\mf{G}(\mf{x}_0)^{-1}.\]

Then
\begin{equation}\label{230811_15} D_{\mf{v}}f(\mf{x}_0)=\boldsymbol{\lambda} {D}_{\mf{v}}\mf{G}(\mf{x}_0).\end{equation}

Eqs. \eqref{230811_13} and \eqref{230811_14} show that 
\begin{equation}\label{230811_16} D_{\mf{u}}f(\mf{x}_0)=-\boldsymbol{\lambda} \mf{D}_{\mf{v}}\mf{G}(\mf{x}_0)\mf{DH}(\mf{u}_0)=\boldsymbol{\lambda}\mf{D}_{\mf{u}}\mf{G}(\mf{x}_0).\end{equation}
Eqs. \eqref{230811_15} and \eqref{230811_16} together imply that 
\[\nabla f(\mf{x}_0)=\boldsymbol{\lambda} \mf{DG}(\mf{x}_0)=\sum_{i=1}^m\lambda_i \nabla G_i(\mf{x}_0).\] 
This completes the proof of the theorem.
\end{myproof}

In the general constraint optimization problem proposed in Theorem \ref{230812_7}, there are $n+2m$ variables $u_1, \ldots, u_n$, $v_1, \ldots, v_m$ and $\lambda_1, \ldots, \lambda_m$ to be solved. 
The components of
\[\nabla f(\mf{x}) =\sum_{i=1}^m\lambda_i \nabla G_i(\mf{x}) \] give $n+m$ equations, while the components of $\mf{G}(\mf{x})=\mf{0}$ give $m$ equations. Hence, we have to solve $n+2m$ variables from $n+2m$ equations.
Let us look at  an example.
\begin{example}{}
Let $K$ be the subset of $\mb{R}^3$ given by
\[K=\left\{(x,y,z)\,|\, x^2+y^2\leq 4,  x+y+z=1\right\}.\]
Find the maximum and minimum values of the function $f:K\to\mb{R}$, $f(x,y,z)=x+3y+z$.

\end{example}
\begin{solution}{Solution}
Notice that $K$ is the intersection of the two closed sets $K_1=\left\{(x,y,z)\,|\, x^2+y^2 \leq 4\right\}$ and $K_2=\left\{(x,y,z)\,|\,   x+y+z=1\right\}$. Hence, $K$ is a closed set. If $(x,y,z)$ is in $K$, $x^2+y^2\leq 4$. Thus, $|x|\leq 2$, $|y|\leq 2$ and hence $|z|\leq 1+|x|+|y|\leq 5$. 
 This shows that $K$ is bounded. Since $K$ is closed and bounded, $f:K\to\mb{R}$ is continuous, $f:K\to\mb{R}$ has maximum and minimum values.

 Let 
\[D=\left\{(x,y,z)\,|\, x^2+y^2 < 4,  x+y+z=1\right\},\]
\[C=\left\{(x,y,z)\,|\, x^2+y^2=4,  x+y+z=1\right\}.\]
Then $K=C\cup D$. We can consider the extremizers of $f:D\to\mb{R}$ and $f:C\to \mb{R}$ separately.

To find the extremizers of $f:D\to\mb{R}$, we can regard this as a constraint optimization problem where we want to find the extreme values of $f:\mathcal{O}\to\mb{R}$, $f(x,y,z)=x+3y+z$ on  \[\mathcal{O}=\left\{(x,y,z)\,|\, x^2+y^2 < 4 \right\},\]

  subject to the constraint $g(x,y,z)=0$, where $g:\mathcal{O}\to\mb{R}$ is the function $g(x,y,z)=x+y+z-1$. Now $\nabla g(x,y,z)=(1,1,1)\neq\mf{0}$. Hence, at an extremizer, we must have
$\nabla f(x,y,z)=\lambda g(x,y,z)$, which gives
\[(1,3,1)=\lambda (1, 1, 1).\]
\bs
This says that the two vectors $(1,3,1)$ and $(1,1,1)$ must be parallel, which is a contradiction. Hence, $f:\mathcal{O}\to\mb{R}$ does not have extremizers.

Now, to find the extremizers of $f:C\to\mb{R}$, we can consider it as finding the extreme values of $f:\mb{R}^3\to\mb{R}$, $f(x,y,z)=x+3y+z$, subject to $\mf{G}(x,y,z)=0$, where
\[\mf{G}(x,y,z)=(x^2+y^2-4, x+y+z-1).\]
Now
\[\mf{DG}(x,y,z)=\begin{bmatrix} 2x & 2y & 0\\1 & 1 & 1\end{bmatrix}.\]
This matrix has rank less than 2 if and only if $(2x, 2y, 0)$ is parallel to $(1, 1, 1)$, which gives $x=y=z=0$. But the point $(x,y,z)=(0,0,0)$ is not on $C$. Therefore, $\mf{DG}(x,y,z)$ has maximal rank for every $(x,y,z)\in C$. Using the Lagrange multiplier method, to solve for the extremizer of $f:C\to\mb{R}$, we need to solve the system
\[\nabla f(x,y,z)=\lambda \nabla G_1(x,y,z)+\mu G_2(x,y,z),\quad \mf{G}(x,y,z)=\mf{0}.\]

These gives
\begin{gather*}
1 =2\lambda x+\mu,\quad 3=2\lambda y+\mu,\quad 1= \mu,\\
x^2+y^2 =4,\quad x+y+z=1.
\end{gather*}
 
From $\mu=1$, we have
$2\lambda x=0$ and $2\lambda y=2$. The latter implies that $\lambda\neq 0$. Hence, we must have $x=0$. Then $x^2+y^2=4$ gives   $y=\pm 2$. When $(x,y)=(0,2)$, $z=-1$. When $(x,y)=(0,-2)$, $z=3$. Hence, we only have two candidates for extremizers, which are $(0, 2, -1)$ and $(0,-2, 3)$.
Since \[f(0, 2, -1)=5,\quad f(0, -2, 3)=-3,\] we find that $f:K\to\mb{R}$ has maximum value $5$ at the point $(0, 2, -1)$, and minimum value $-3$ at the point $(0, -2, 3)$.

\end{solution}
\vp
\noindent
{\bf \large Exercises  \thesection}
\setcounter{myquestion}{1}
\begin{question}{\themyquestion}
Find the extreme values of the function $f(x,y,z)= 4x^2+y^2+yz+z^2 $ on the set
\[S=\left\{(x,y,z)\,|\, 2x^2+y^2 +z^2\leq 8\right\}.\]
\end{question}
\atc

\begin{question}{\themyquestion}
Find the point in the set
\[S=\left\{(x,y)\,|\, 4x^2+y^2\leq 36, x^2+4y^2\geq 4\right\}\] that is closest to and farthest from the point  $(1,0)$.
\end{question}
\atc

\begin{question}{\themyquestion}
Use the Lagrange multiplier method to find the maximum and minimum values of the function $f(x,y,z)=x+2y-z $ on the set
\[S=\left\{(x,y,z)\,|\, x^2+y^2+4z^2\leq 84\right\},\]and the points where each of them appears.

\end{question}
\atc

\begin{question}{\themyquestion}
Find the extreme values of the function $f(x,y,z)=x$ on the set
\[S=\left\{(x,y,z)\,|\, x^2=y^2+z^2, 7x+3y+4z=60\right\}.\]
\end{question}
\atc
 
 \begin{question}{\themyquestion}
Let $K$ be the subset of $\mb{R}^3$ given by
\[K=\left\{(x,y,z)\,|\, 4x^2+z^2\leq 68,  y+z=12\right\}.\]
Find the maximum and minimum values of the function $f:K\to\mb{R}$, $f(x,y,z)=x+2y$.

\end{question}
 \atc
 
 \begin{question}{\themyquestion}
Let $A$ be an $n\times n$ symmetric matrix, and let $Q_A:\mb{R}^n\to\mb{R}$ be the quadratic form $Q_A(\mf{x})=\mf{x}^TA\mf{x}$ defined by $A$. Show that the minimum and maximum values of $Q_A:S^{n-1}\to\mb{R}$ on the unit sphere $S^{n-1}$ are the smallest and largest eigenvalues of $A$.

\end{question}

\chapter{Multiple Integrals}
For a single variable functions, we have discussed the Riemann integrability of a function $f:[a,b]\to\mb{R}$ defined on a compact interval $[a,b]$.
 In this chapter, we consider the theory of Riemann integrals for multivariable functions. For a function $\mf{F}:\mk{D}\to\mb{R}^m$ that takes values in $\mb{R}^m$ with $m\geq 2$, we define the integral componentwise. Namely, we say that the function $\mf{F}:\mk{D}\to\mb{R}^m$ is Riemann integrable if and only if each of the component functions $F_j:\mk{D}\to\mb{R}$, $1\leq j\leq m$ is Riemann integrable, and we  define 
 \[\int_{\mk{D}}\mf{F}=\left(\int_{\mk{D}}F_1, \int_{\mk{D}}F_2, \ldots, \int_{\mk{D}}F_m\right).\]
 Thus, in this chapter, we will only discuss the theory of integration for functions $f:\mk{D}\to\mb{R}$ that take values in $\mb{R}$.
 
  A direct generalization of a compact interval $[a,b]$ to $\mb{R}^n$ is a product of compact intervals $\di\mf{I}=\prod_{i=1}^n [a_i, b_i]$, which is a closed rectangle.  In this chapter, when we say $\mf{I}$ is a rectangle, it means     $\di \mf{I}$ can be written as $\di\prod_{i=1}^n [a_i, b_i]$ with  $a_i<b_i$ for all $1\leq i\leq n$. The  edges of $\mf{I}=\di\prod_{i=1}^n[a_i,b_i]$ are $[a_1, b_1]$, $[a_2, b_2]$, $\ldots$, $[a_n, b_n]$.

 We first discuss the integration theory of functions defined on closed rectangles of the form $\di\prod_{i=1}^n[a_i,b_i]$. For applications, we need to consider functions defined on other subsets $\mk{D}$ of $\mb{R}^n$. 
 
One of the most useful theoretical tools for evaluating single integrals is the fundamental theorem of calculus. To apply this tool for multiple integrals, we need to consider iterated integrals. Another useful tool is the change of variables formula. For multivariable functions, the change of variables theorem is much more complicated. Nevertheless, we will discuss these in this chapter.
 
\section{Riemann Integrals} \label{sec6.1}
In this section, we define the Riemann integral of a function $f:\mk{D}\to\mb{R}$ defined on a  subset $\mk{D}$ of $\mb{R}^n$.  We first consider the case where $\mk{D}=\di\prod_{i-1}^n[a_i,b_i]$.

Let us first consider partitions. We say that  $P=\{x_0, x_1, \ldots, x_k\}$ is a partition of the interval $[a,b]$  if
$\di a=x_0<x_1<\cdots<x_{k-1}<x_k=b$.
It divides $[a,b]$ into $k$ subintervals $J_1, \ldots, J_k$, where $J_i=[x_{i-1}, x_i]$. 

\begin{definition}{Partitions}
  A partition $\mf{P}$ of  a closed rectangle $\mf{I}=\di\prod_{i=1}^n [a_i, b_i]$ is achieved by having a partition $P_i$ of $[a_i, b_i]$ for each $1\leq i\leq n$. We write $\mf{P}=(P_1, P_2, \ldots, P_n)$ for such a partition. The partition $\mf{P}$ divides the rectangle $\mf{I}$ into a collection $\mathcal{J}_{\mf{P}}$ of rectangles, any two of which have disjoint interiors. A closed rectangle $\mf{J}$ in  $\mathcal{J}_{\mf{P}}$ can be written as \[\mf{J}=J_1\times J_2\times \cdots\times J_n,\] where $J_i$, $1\leq i\leq n$ is a subinterval in the partition $P_i$. 
\end{definition}If the partition $P_i$ divides $[a_i,b_i]$ into $k_i$ subintervals, then the partition $\mf{P}=(P_1, \ldots, P_n)$ divides the rectangle $\mf{I}=\di\prod_{i=1}^n [a_i,b_i]$ into $|\mathcal{J}_{\mf{P}}|=k_1k_2\cdots k_n$ rectangles.

\begin{example}[label=230821_1]{}
Consider the rectangle $\mf{I}=[-2, 9]\times [1, 6]$. Let $P_1=\{-2, 0, 4, 9\}$ and $P_2=\{1, 3, 6\}$. The partition $P_1$ divides the interval $I_1=[-2,9]$ into the three subintervals $[-2,0]$, $[0,4]$ and $[4,9]$. The partition $P_2$ divides the interval $I_2=[1,6]$ into the two subintervals $[1,3]$ and $[3,6]$. Therefore, the partition $\mf{P}=(P_1, P_2)$ divides the rectangle $\mf{I}$ into the following six rectangles.

\vspace{-1cm}
\begin{gather*}
[-2,0]\times [1,3],\quad [0,4]\times [1,3],\quad[4,9]\times [1,3],\\ [-2,0]\times [3,6],\quad
[0,4]\times [3,6],  \quad [4,9]\times [3,6].
\end{gather*}
\end{example}

 \begin{figure}[ht]
\centering
\includegraphics[scale=0.2]{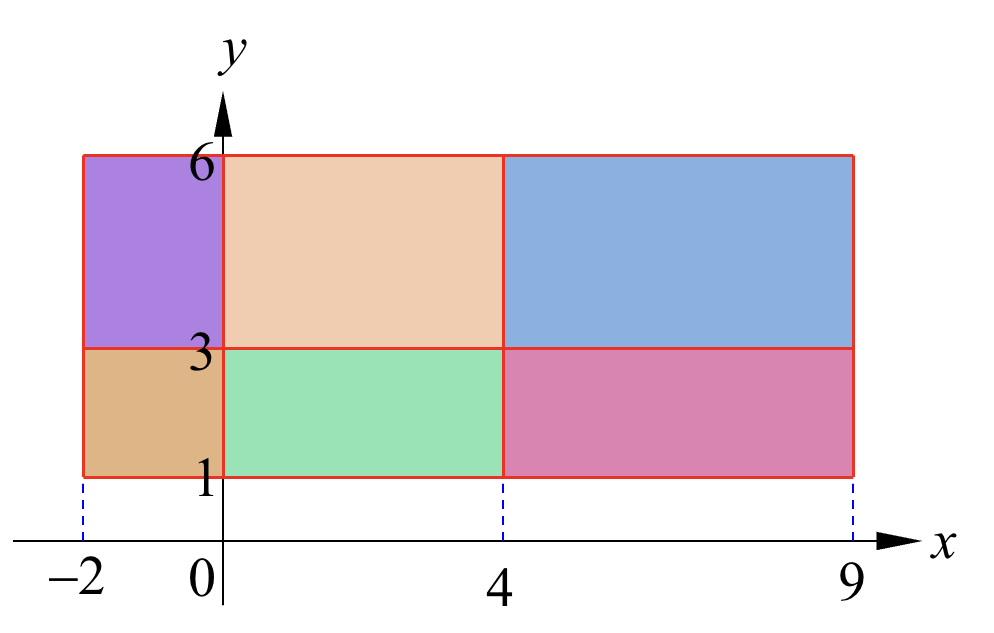}
\caption{A partition of the rectangle $[-2,9]\times [1,6]$ given in Example \ref{230821_1}.}\label{figure84}
\end{figure}

\begin{definition}{Regular and Uniformly Regular Partitions}
Let $\mf{I}=\di\prod_{i=1}^n[a_i, b_i]$ be a rectangle in $\mb{R}^n$. We say that $\mf{P}=(P_1, \ldots, P_n)$ is a {\it regular partition} of $\mf{I}$ if for each $1\leq i\leq n$, $P_i$ is a regular partition of $[a_i, b_i]$ into $k_i$ intervals. We say that $\mf{P}$ is a {\it uniformly regular partition} of $\mf{P}$ into $k^n$ rectangles if for each $1\leq i\leq n$, $P_i$ is a regular partition of $[a_i, b_i]$ into $k$ intervals.
\end{definition}

\begin{example}[label=230824_4]{}
Consider the rectangle $\mf{I}=[-2,7]\times[-4, 8]$. 
\begin{enumerate}[(a)]
\item The partition $\mf{P}=(P_1, P_2)$ where $P_1=\{-2, 1, 4, 7\}$ and $P_2= \{-4,-1,2,5,8\}$ is a regular partition of $\mf{I}$. 
\item
The partition $\mf{P}=(P_1, P_2)$ where $P_1=\{-2, 1, 4, 7\}$ and $P_2= \{-4,0,4,8\}$ is a uniformly regular partition of $\mf{I}$ into $3^2=9$ rectangles. 
\end{enumerate}
\end{example}

 \begin{figure}[ht]
\centering
\includegraphics[scale=0.2]{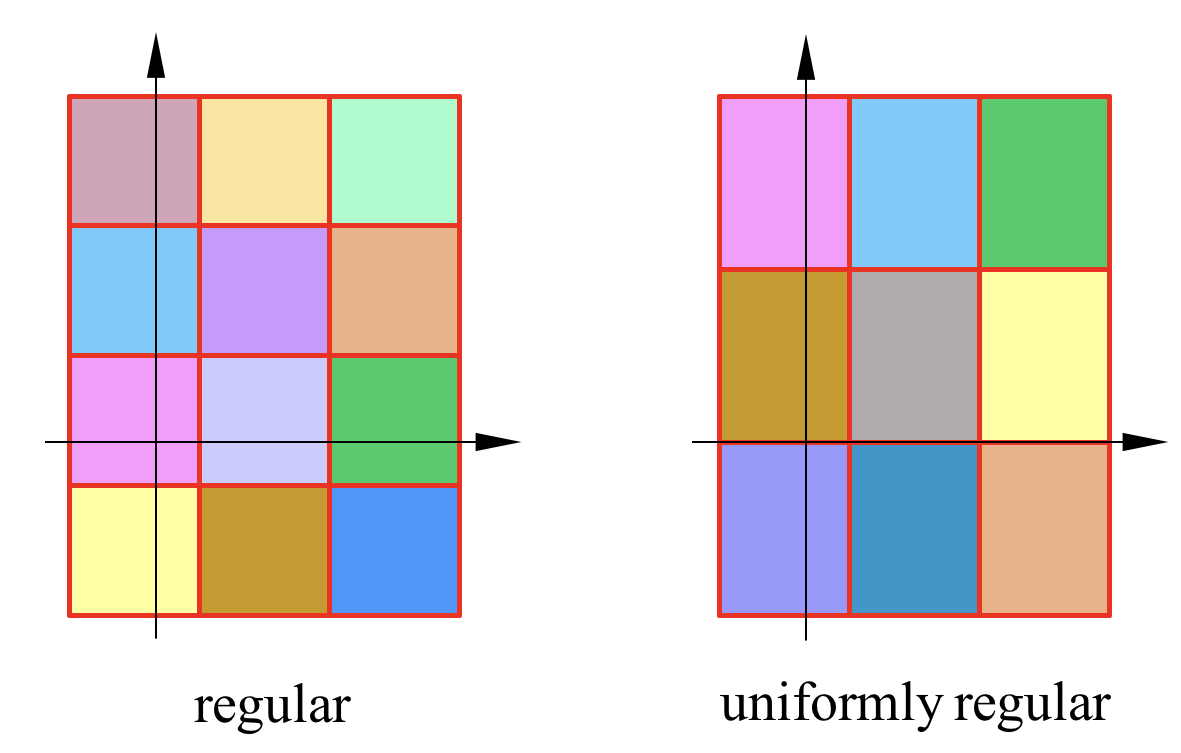}
\caption{A  regular and a uniformly regular partition of $[-2,7]\times [-4,8]$ discussed in Example \ref{230824_4}.}\label{figure87}
\end{figure}

The length of an interval $[a,b]$ is $b-a$. The area of a rectangle $[a,b]\times [c,d]$ is $(b-a)\times (d-c)$. In general, we define the volume of a closed rectangle of the form $\mf{I}=\di\prod_{i=1}^n[a_i,b_i]$ in $\mb{R}^n$ as follows.
\begin{definition}{Volume of a Rectangle}
The volume of the closed rectangle $\mf{I}=\di\prod_{i=1}^n [a_i, b_i]$ is defined as the product of the lengths of all its edges. Namely,
\[\text{vol}\,(\mf{I})=\prod_{i=1}^n(b_i-a_i).\]
\end{definition}
\begin{example}{}
The volume of the rectangle $\mf{I}=[-2,9]\times [1,6]$ is 
\[\text{vol}\,(\mf{I})=11\times 5=55.\]
\end{example}
When $P=\{x_0, x_1, \ldots, x_k\}$ is a partition of $[a,b]$, it divides $[a,b]$ into $k$ subintervals $J_1, \ldots, J_k$, where $J_i=[x_{i-1}, x_i]$. Notice that
\[\sum_{i=1}^k\text{vol}\,(J_i)=\sum_{i=1}^k(x_i-x_{i-1})=b-a.\]
Assume that $\mf{P}=(P_1, \cdots, P_n)$ is a partition of the rectangle $\mf{I}=\di\prod_{i=1}^n[a_i,b_i]$ in $\mb{R}^n$. Then for $1\leq i\leq n$, $P_i$ is a partition of $[a_i,b_i]$. If $P_i$ divides $[a_i,b_i]$ into the $k_i$ subintervals $J_{i,1}, J_{i,2}, \ldots, J_{i,k_i}$, then the collection of rectangles in the partition $\mf{P}$ is
\[\mathcal{J}_{\mf{P}}=\left\{J_{1,m_1}\times \cdots\times J_{n, m_n}\,|\,1\leq m_i \leq k_i\;\text{for}\;1\leq i\leq n\right\}.\]
Notice that
\[\text{vol}\,\left(J_{1,m_1}\times \cdots\times J_{n, m_n}\right)=\text{vol}\,(J_{1,m_1})\times \cdots\times \text{vol}\,(J_{n, m_n}).\]
From this, we obtain the sum of volumes formula:
\begin{align*}
\sum_{\mf{J}\in\mathcal{J}_{\mf{P}}}\text{vol}\,(\mf{J})&=
\sum_{m_n=1}^{k_n}\cdots\sum_{m_1=1}^{k_1}\text{vol}\, (J_{1,m_1})\times \cdots\times \text{vol}\,(J_{n, m_n})\\
&=\left[\sum_{m_1=1}^{k_1}\text{vol}\,(J_{1, m_1})\right]\times \cdots\times\left[\sum_{m_n=1}^{k_n}\text{vol}\,(J_{n, m_n})\right] \\
&=(b_1-a_1)\times\cdots\times (b_n-a_n)\\
&=\text{vol}\,(\mf{I}).
\end{align*}

\begin{proposition}{}
Let $\mf{P}$ be a partition of   $\mf{I}=\di\prod_{i=1}^n[a_i,b_i]$. Then the sum of the volumes of the rectangles $\mf{J}$ in the partition $\mf{P}$ is equal to the volume of the rectangle $\mf{I}$. 
\end{proposition}

One of the motivations to define the integral $\di\int_{\mf{I}}f$ for a nonnegative function $f:\mf{I}\to\mb{R}$ is to find the volume bounded between  the graph of $f$ and the rectangle $\mf{I}$ in $\mb{R}^{n+1}$. To find the volume, we partition $\mf{I}$ into small rectangles,  pick a point $\boldsymbol{\xi}_{\mf{J}}$ in each of these rectangles $\mf{J}$, and approximate the function on $\mf{J}$ as a constant given by the value $f(\boldsymbol{\xi}_{\mf{J}})$. The volume between the rectangle $\mf{J}$ and the graph of $f$ over $\mf{J}$ is then approximated by $ f(\boldsymbol{\xi}_{\mf{J}}) \,\text{vol}\,(\mf{J})$.
This leads us to the concept of  Riemann sums. 

If  $\mf{P}$ is a partition of $\mf{I} =\di\prod_{i=1}^n[a_i, b_i]$, we say that $A$ is a  set of {\it intermediate points} for the partition $\mf{P}$ if $A=\left\{\boldsymbol{\xi}_{\mf{J}}\,|\,\mf{J}\in \mathcal{J}_{\mf{P}}\right\}$ is a subset of $\mf{I}$ indexed by $\mathcal{J}_{\mf{P}}$, such that $\boldsymbol{\xi}_{\mf{J}}\in\mf{J}$ for each $\mf{J}\in \mathcal{J}_{\mf{P}}$. 
 
\begin{definition}{Riemann Sums}
Let $\mf{I} =\di\prod_{i=1}^n[a_i, b_i]$, and let $f:\mf{I}\to\mb{R}$ be a function defined on $\mf{I}$. Given a partition $\mf{P}$ of $\mf{I}$,  a set  $A=\left\{\boldsymbol{\xi}_{\mf{J}}\,|\,\mf{J}\in \mathcal{J}_{\mf{P}}\right\}$ of intermediate points for the partition $\mf{P}$, the Riemann sum of $f$ with respect to the partition $\mf{P}$ and the set of intermediate points $A=\{\boldsymbol{\xi}_{\mf{J}}\}$ is the sum
\[R(f,\mf{P}, A)=\sum_{\mf{J}\in\mathcal{J}_{\mf{P}}}f(\boldsymbol{\xi}_{\mf{J}})\text{vol}\,(\mf{J}).\]
\end{definition}

\begin{example}[label=230821_2]{}
Let $\mf{I}=[-2,9]\times [1,6]$, and let $\mf{P}=(P_1, P_2)$ be the partition of $\mf{I}$ with $P_1=\{-2,0,4,9\}$ and $P_2=\{1,3,6\}$. Let $f:\mf{I}\to\mb{R}$ be the function defined as $f(x,y)=x^2+y$. Consider a set of intermediate points $A$ as follows.

\vspace{0.2cm} ~\hspace{2cm}
\begin{tabular}{||c||c|c|c||}
\hline
\hline
$\mf{J}$ &$\boldsymbol{\xi}_{\mf{J}}$  & $f(\boldsymbol{\xi}_{\mf{J}})$ & $\text{vol}\,(\mf{J})$\\
\hline
\hline
 $[-2,0]\times [1,3] $ & $(-1,1)$ & 2 & 4\\
 \hline
 $[-2,0]\times [3,6] $ & $(0,3)$ & 3 & 6\\
\hline
 $[0,4]\times [1,3] $ & $(1,1)$ & 2 & 8\\
 \hline
 $[0,4]\times [3,6] $ & $(2,4)$ & 8 & 12\\
\hline
 $[4,9]\times [1,3] $ & $(4,2)$ & 18 &10\\
 \hline
 $[4,9]\times [3,6] $ & $(9,3)$ & 84 &15\\

\hline
\hline
\end{tabular}

\vspace{0.3cm}
The Riemann sum $R(f,\mf{P},A)$ is equal to

\vspace{-0.4cm}
\[2\times 4+3\times 6+2\times 8+8\times 12+18\times 10+84\times 15=1578.\]
\end{example}

\begin{example}{}
If $f:\mf{I}\to\mb{R}$ is the constant function $f(\mf{x})=c$, then for any partition $\mf{P}$ of $\mf{I}$ and any set of intermediate points $A=\{\boldsymbol{\xi}_{\mf{J}}\}$,
\[R(f,\mf{P}, A)=c\,\text{vol}\,(\mf{I}).\]
When $c>0$, this is  the volume of the rectangle $\mf{I}\times [0,c]$ in $\mb{R}^{n+1}$. 
\end{example}

As in the single variable case,    Darboux sums provide  bounds for Riemann sums.
\begin{definition}{Darboux Sums}
Let $\mf{I} =\di\prod_{i=1}^n[a_i, b_i]$ , and let $f:\mf{I}\to\mb{R}$ be a  {\it bounded} function defined on $\mf{I}$. Given a partition $\mf{P}$ of $\mf{I}$,  let $\mathcal{J}_{\mf{P}}$ be the collection of rectangles in the partition $\mf{P}$. For each $\mf{J}$ in $\mathcal{J}_{\mf{P}}$, 
let
\[m_{\mf{J}}=\inf\left\{f(\mf{x})\,|\,\mf{x}\in\mf{J}\right\}\quad\text{and}\quad M_{\mf{J}}=\sup\left\{f(\mf{x})\,|\,\mf{x}\in\mf{J}\right\}.\]The Darboux lower sum $L(f,\mf{P})$ and the Darboux upper sum $U(f,\mf{P})$ are defined as
\[L(f,\mf{P})=\sum_{\mf{J}\in\mathcal{J}_{\mf{P}}} m_{\mf{J}}\; \text{vol}\,(\mf{J})\quad \text{and}\quad U(f,\mf{P})=\sum_{\mf{J}\in\mathcal{J}_{\mf{P}}} M_{\mf{J}}\; \text{vol}\,(\mf{J}).\]
\end{definition}

\begin{example}{}
If $f:\mf{I}\to\mb{R}$ is the constant function $f(\mf{x})=c$, then  
\[L(f,\mf{P})=c\,\text{vol}\,(\mf{I})=U(f,\mf{P})\hspace{1cm}\text{for any partition $\mf{P}$ of $\mf{I}$}.\]
\end{example}

\begin{example}[label=230821_3]{}
Consider the function $f:\mf{I}\to\mb{R}$, $f(x,y)=x^2+y$ defined in Example \ref{230821_2}, where 
  $\mf{I}=[-2,9]\times [1,6]$. For the partition $\mf{P}=(P_1, P_2)$  with $P_1=\{-2,0,4,9\}$ and $P_2=\{1,3,6\}$, we have the followings.

\vspace{0.5cm} ~\hspace{0.5cm}
\begin{tabular}{||c||c|c|c||}
\hline
\hline
$\mf{J}$ &$m_{\mf{J}}$  &$M_{\mf{J}}$     & $\text{vol}\,(\mf{J})$\\
\hline
\hline
 $[-2,0]\times [1,3] $ &  $0^2+1=1$ & $(-2)^2+3=7$& 4\\
 \hline
 $[-2,0]\times [3,6] $ &  $0^2+3=3 $ & $(-2)^2+6=10$ & 6\\
\hline
 $[0,4]\times [1,3] $ &  $0^2+1=1$  & $4^2+3=19$& 8\\
 \hline
 $[0,4]\times [3,6] $ &  $0^2+3=3$  & $4^2+6=22$& 12\\
\hline
 $[4,9]\times [1,3] $ &  $4^2+1=17$  & $9^2+3=84$ &10\\
 \hline
 $[4,9]\times [3,6] $ &  $4^2+3=19$  & $9^2+6=87$ &15\\

\hline
\hline
\end{tabular}

\be
Therefore, the Darboux lower sum is
\[L(f,\mf{P})=1\times 4+3\times 6+1\times 8+3\times 12+17\times 10+19\times 15=521;\]
while the Darboux upper sum is
\[U(f,\mf{P})=7\times 4+10\times 6+19\times 8+22\times 12+84\times 10+87\times 15=2649.\]
\end{example2}

Notice that we can only define Darboux sums if the function $f:\mf{I}\to\mb{R}$ is bounded. This means that there are constants $m$ and $M$ such that
\[m\leq f(x)\leq M\hspace{1cm}\text{for all}\;\mf{x}\in\mf{I}.\]
If $\mf{P}$ is a partition of the rectangle $\mf{I}$, and $\mf{J}$ is a rectangle in the partition $\mf{P}$, $\boldsymbol{\xi}_J$ is a point in $\mf{J}$, then
\[m\leq m_{\mf{J}}\leq f(\boldsymbol{\xi}_{\mf{J}}) \leq M_{\mf{J}}\leq M.\]
Multipluying throughout by $\text{vol}\,(\mf{J})$ and summing over $\mf{J}\in\mathcal{J}_{\mf{P}}$, we obtain the following.

\begin{proposition}[label=230822_3]{}
Let $\mf{I}=\di\prod_{i=1}^n[a_i,b_i]$, and let $f:\mf{I}\to\mb{R}$ be a bounded function defined on $\mf{I}$. If
\[m\leq f(\mf{x})\leq M\hspace{1cm}\text{for all}\;\mf{x}\in\mf{I},\]then for any partition $\mf{P}$ of $\mf{I}$, and for any choice of intermediate points $A=\left\{\boldsymbol{\xi}_{\mf{J}}\right\}$ for the partition $\mf{P}$, we have
\[m\,\text{vol}\,(\mf{I})\;\leq\; L(f,\mf{P})\;\leq \;R(f,\mf{P},A)\;\leq\;U(f,\mf{P})\;\leq\;M\,\text{vol}\,(\mf{I}).\]
\end{proposition}

To study the behaviour of the Darboux sums when we modify the  partitions, we first extend the concept of refinement of a partition to rectangles in $\mb{R}^n$. Recall that if $P$ and $P^*$ are  partitions of the interval $[a,b]$, $P^*$ is a refinement of $P$ if each partition point of $P$ is also a partition point of $P^*$.

\begin{definition}{Refinement of a Partition}
Let $\mf{I}=\di\prod_{i=1}^n[a_i,b_i]$, and let $\mf{P}=(P_1, \ldots, P_n)$ and $\mf{P}^*=(P_1^*, \ldots, P_n^*)$ be partitions of $\mf{I}$. We say that $\mf{P}^*$ is a refinement of $\mf{P}$ if for each $1\leq i\leq n$, $P_i^*$ is a refinement of $P_i$. 
\end{definition}

 \begin{figure}[ht]
\centering
\includegraphics[scale=0.2]{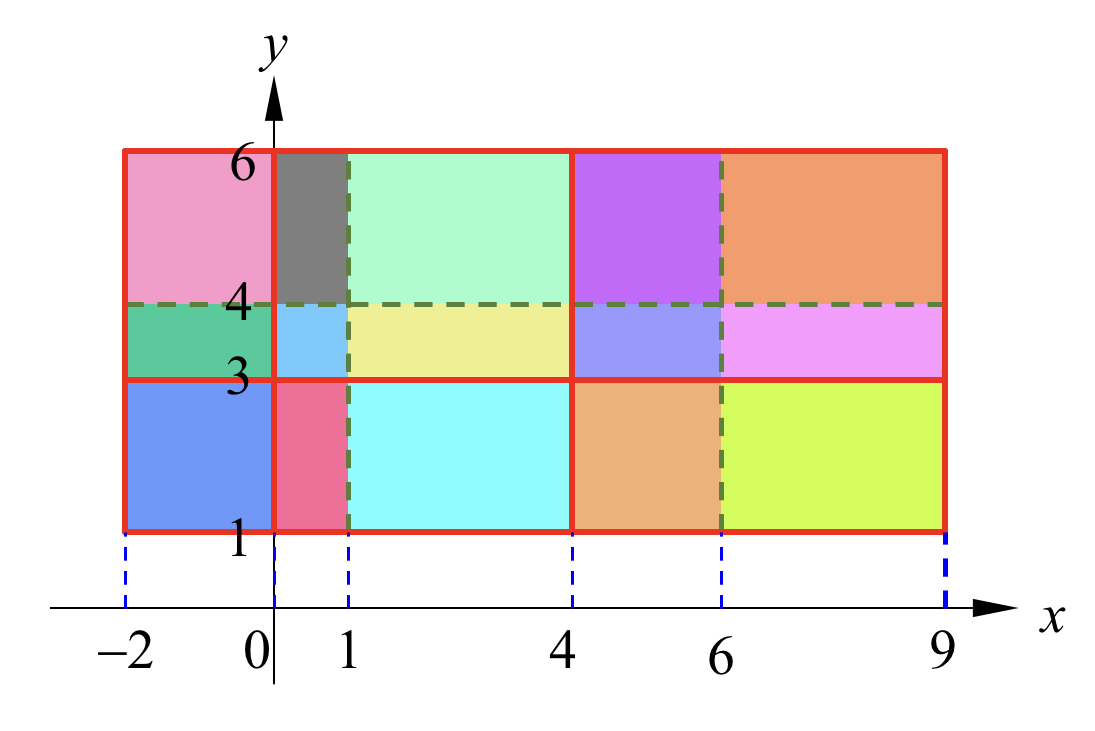}
\caption{A refinement of the partition of the rectangle $[-2,9]\times [1,6]$ given in Figure \ref{figure84}.}\label{figure85}
\end{figure}

\begin{example}[label=230822_1]{}
Let us consider the partition $\mf{P}=(P_1, P_2)$ of the rectangle $\mf{I}=[-2,9]\times [1,6]$ given in Example \ref{230821_1}, with $P_1=\{-2,0,4,9\}$ and $P_2=\{1, 3, 6\}$. Let $P_1^*=\{-2,0,1,4,6,9\}$ and $P_2^*=\{1,3,4,6\}$. Then $\mf{P}^*=(P_1^*,P_2^*)$ is a refinement of $\mf{P}$.
\end{example}

If the partition $\mf{P}^*$ is a refinement of the partition $\mf{P}$, then for each $\mf{J}$ in $\mathcal{J}_{\mf{P}}$, $\mf{P}^*$ induces a partition of $\mf{J}$, which we denote by $\mathbf{P}^*(\mf{J})$.  

\begin{example}{}The partition $\mf{P}^*$ in Example \ref{230822_1} induces the partition $\mf{P}^*(\mf{J})=(P_1^*(\mf{J}), P_2^*(\mf{J}))$ of the rectangle $\mf{J}=[0,4]\times [3,6]$, where $P_1^*(\mf{J})=\{0,1,4\}$ and $P_2^*(\mf{J})=\{3,4,6\}$. The partition $\mf{P}^*(\mf{J})$ divides the rectangle $\mf{J}$ into 4 rectangles, as shown in Figure \ref{figure85}.
\end{example}

If the partition $\mf{P}^*$ is a refinement of the partition $\mf{P}$, 
then the collection of rectangles in $\mf{P}^*$ 
is the union of the collection of rectangles in $\mf{P}^*(\mf{J})$ when $\mf{J}$ ranges over the collection of  rectangles in $\mf{P}$. 
Namely,
\[\mathcal{J}_{\mf{P}^*}=\bigcup_{\mf{J}\in\mathcal{J}_{\mf{P}}}\mathcal{J}_{\mf{P}^*(\mf{J})}.\]
Using this, we can deduce the following.
\begin{proposition}[label=230822_4]{}
Let $\mf{I}=\di\prod_{i=1}^n[a_i,b_i]$, and let $f:\mf{I}\to\mb{R}$ be a bounded function defined on $\mf{I}$. If $\mf{P}$ and $\mf{P}^*$ are partitions of $\mf{I}$ and $\mf{P}^*$ is a refinement of $\mf{P}$, then 
\[L(f,\mf{P}^*)=\sum_{\mf{J}\in\mathcal{J}_{\mf{P}}}L(f, \mf{P}^*(\mf{J})),\hspace{1cm} U(f,\mf{P}^*)=\sum_{\mf{J}\in\mathcal{J}_{\mf{P}}}U(f, \mf{P}^*(\mf{J})).\]
\end{proposition}

From this, we can show that a refinement improves the Darboux sums, in the sense that a lower sum  increases, and an upper sum decreases.
\begin{theorem}[label=230824_5]{}
Let $\mf{I}=\di\prod_{i=1}^n[a_i,b_i]$, and let $f:\mf{I}\to\mb{R}$ be a bounded function defined on $\mf{I}$. If $\mf{P}$ and $\mf{P}^*$ are partitions of $\mf{I}$ and $\mf{P}^*$ is a refinement of $\mf{P}$, then 
\[L(f,\mf{P})\leq L(f,\mf{P}^*)\leq U(f,\mf{P}^*)\leq U(f,\mf{P}).\]
\end{theorem}
\begin{myproof}{Proof}
For each rectangle $\mf{J}$ in the partition $\mf{P}$,
\[m_{\mf{J}}\leq f(\mf{x})\leq M_{\mf{J}}\hspace{1cm}\text{for all}\;\mf{x}\in\mf{J}.\]
Applying  Proposition \ref{230822_3} to the function $f:\mf{J}\to\mb{R}$ and the partition $\mf{P}^*(\mf{J})$, we find that
\[m_{\mf{J}}\,\text{vol}\;(\mf{J})\leq L(f,\mf{P}^*(\mf{J}))\leq U(f,\mf{P}^*(\mf{J}))\leq M_{\mf{J}}\,\text{vol}\;(\mf{J}).\]
\bp
Summing over $\mf{J}\in\mathcal{J}_{\mf{P}}$, we find that
\[L(f,\mf{P})\leq \sum_{\mf{J}\in\mathcal{J}_{\mf{P}}}L(f, \mf{P}^*(\mf{J}))\leq \sum_{\mf{J}\in\mathcal{J}_{\mf{P}}}U(f, \mf{P}^*(\mf{J}))\leq U(f,\mf{P}).\]
The assertion follows from Proposition \ref{230822_4}.
\end{myproof}
It is difficult to visualize the Darboux sums with a multivariable functions. Hence, we illustrate  refinements improve Darboux sums using  single variable functions, as shown in Figure \ref{figure88} and Figure \ref{figure89}.

 \begin{figure}[ht]
\centering
\includegraphics[scale=0.2]{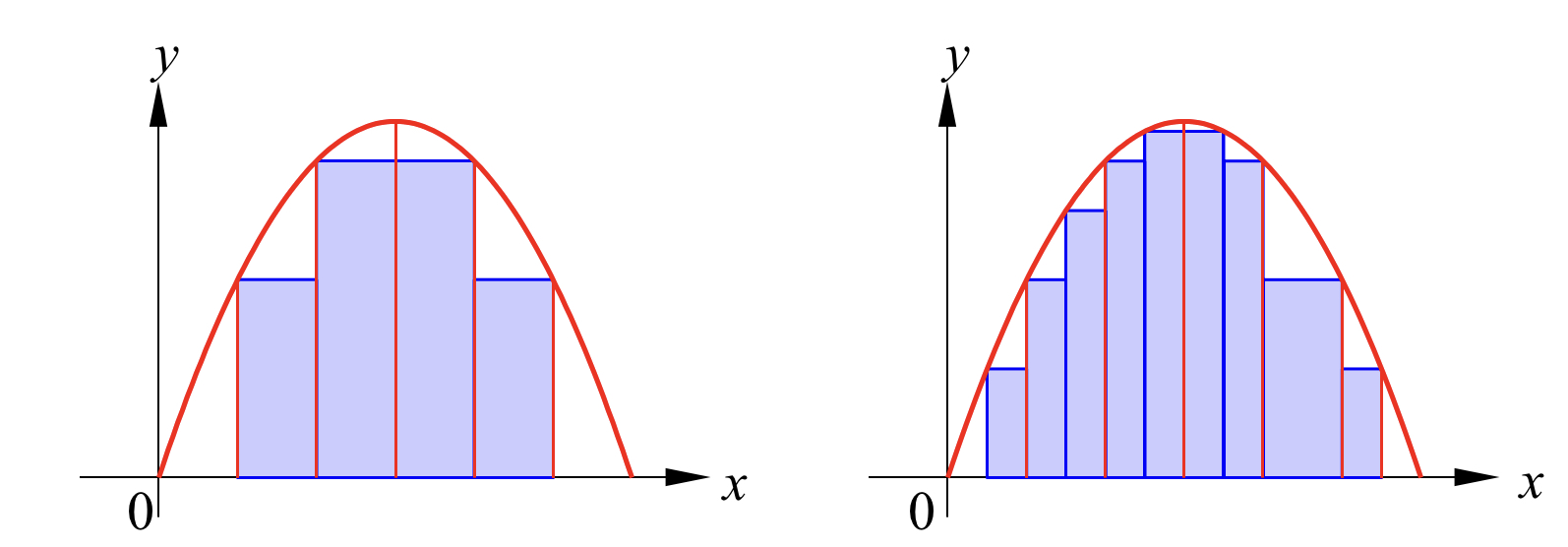}
\caption{A refinement of the partition increases the Darboux lower sum.}\label{figure88}
\end{figure}

\begin{figure}[ht]
\centering
\includegraphics[scale=0.2]{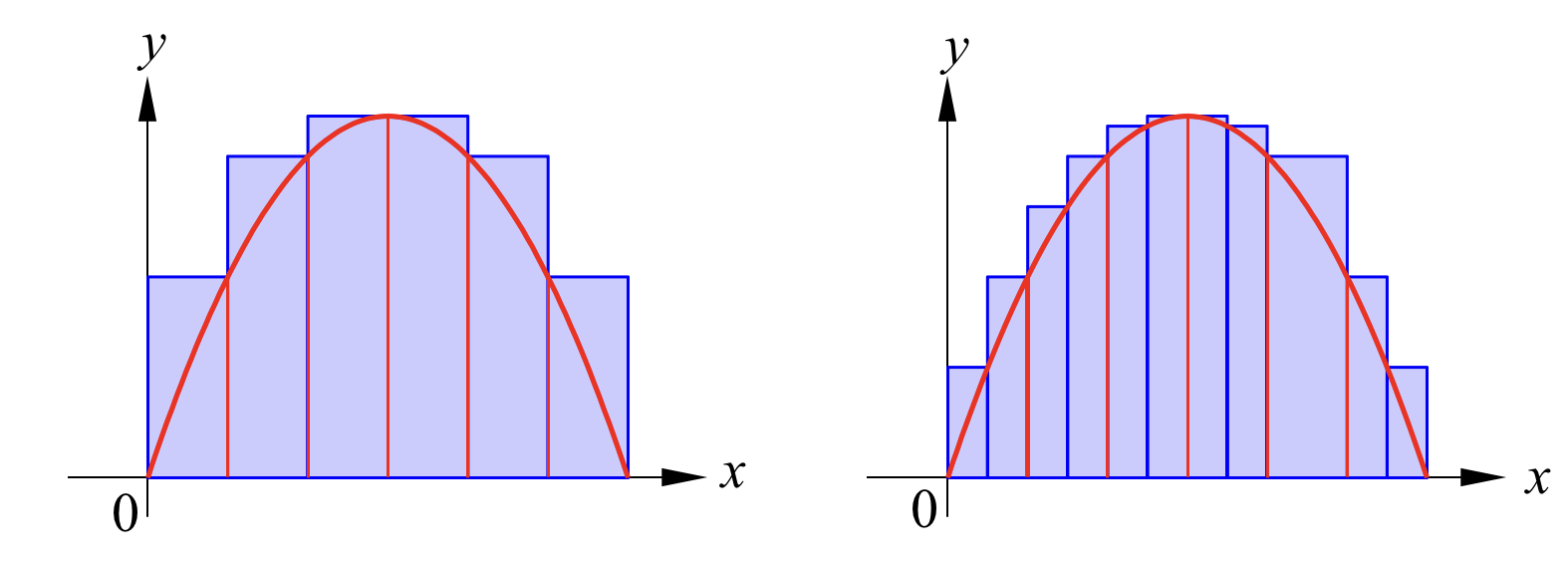}
\caption{A refinement of the partition decreases the Darboux upper sum.}\label{figure89}
\end{figure}
As a consequence of Theorem \ref{230824_5}, we can prove the following.
\begin{corollary}[label=230822_5]{}
Let $\mf{I}=\di\prod_{i=1}^n[a_i,b_i]$, and let $f:\mf{I}\to\mb{R}$ be a bounded function defined on $\mf{I}$. For any two partitions $\mf{P}_1$ and $\mf{P}_2$   of $\mf{I}$,  
\[L(f,\mf{P}_1)\leq U(f,\mf{P}_2).\]
\end{corollary}
\begin{myproof}{Proof}Let  $\mf{P}_1=(P_{1,1}, P_{1,2},\ldots, P_{1,n})$ and $\mf{P}_2=(P_{2,1},P_{2,2}, \ldots, P_{2,n})$.
 For $1\leq i\leq n$, let $P^*_i$ be the common refinement of $P_{1,i}$ and $P_{2,i}$ obtained by taking the union of the partition points in $P_{1,i}$ and $P_{2,i}$. Then $\mf{P}^*=(P_1^*, \ldots, P_n^*)$ is a common refinement of the partitions $\mf{P}_1$ and $\mf{P}_2$. By Theorem \ref{230824_5},
\[L(f,\mf{P}_1)\leq L(f,\mf{P}^*)\leq U(f,\mf{P}^*)\leq U(f,\mf{P}_2).\]
\end{myproof}

Now we define lower and upper integrals of a bounded function $f:\mf{I}\to\mb{R}$ .

\begin{definition}{Lower Integrals and Upper Integrals}
 Let $\mf{I}=\di\prod_{i=1}^n[a_i,b_i]$, and let $f:\mf{I}\to\mb{R}$ be a bounded function defined on $\mf{I}$. Let $S_L(f)$ be the set of Darboux lower sums of $f$, and let $S_U(f)$ be the set of Darboux upper sums of $f$.  
 \begin{enumerate}[1.]
 \item The lower integral of $f$, denoted by $\di\underline{\int_{\mf{I}}}f$, is defined as the least upper bound of the Darboux lower sums.
 \[\underline{\int_{\mf{I}}}f=\sup S_L(f)=\sup \left\{L(f,\mf{P})\,|\, \text{$\mf{P}$ is a partition of $\mf{I}$}\right\}.\]
 \item The upper integral of $f$, denoted by $\di\overline{\int_{\mf{I}}}f$, is defined as the greatest lower bound of the Darboux upper sums.
 \[\overline{\int_{\mf{I}}}f=\inf S_U(f)=\inf \left\{U(f,\mf{P})\,|\, \text{$\mf{P}$ is a partition of $\mf{I}$}\right\}.\]
 
 \end{enumerate}
\end{definition}

\begin{example}[label=230823_1]{}
If $f:\mf{I}\to\mb{R}$ is the constant function $f(\mf{x})=c$, then for any partition $\mf{P}$ of $\mf{I}$,
\[L(f,\mf{P})=c\,\text{vol}\,(\mf{I})=U(f,\mf{P}).\]Therefore, both $S_L(f)$ and $S_U(f)$ are the one-element set $\{c\,\text{vol}\,(\mf{I})\}$. This shows that
\[\underline{\int_{\mf{I}}}f=\overline{\int_{\mf{I}}}f=c\,\text{vol}\,(\mf{I}).\]
\end{example}

For a constant function, the lower integral and the upper integral are the same. For a general bounded funtion, we have the following.

\begin{theorem}{}
Let $\mf{I}=\di\prod_{i=1}^n[a_i,b_i]$, and let $f:\mf{I}\to\mb{R}$ be a bounded function defined on $\mf{I}$. Then we have
\[\underline{\int_{\mf{I}}}f\;\leq \; \overline{\int_{\mf{I}}}f.\]
\end{theorem}

\begin{myproof}{Proof}
By Corollary \ref{230822_5}, every element of $S_L(f)$ is less than or equal to any element of $S_U(f)$. This implies that 
\[\underline{\int_{\mf{I}}}f\;= \; \sup S_L(f)\;\leq \inf S_U(f)\;=\;\overline{\int_{\mf{I}}}f.\]
\end{myproof}

\begin{example}[label=230823_5]{The Dirichlet's Function}
Let $\di\mf{I}=\di\prod_{i=1}^n[a_i,b_i]$, and let $f:\mf{I}\to\mb{R}$ be the function defined as 

\vspace{-0.4cm}
\[f(\mf{x})=\begin{cases}1,\quad &\text{if all components of $\mf{x}$ are rational},\\0,\quad &\text{otherwise}.\end{cases}\]
This is known as the Dirichlet's function. Find the lower inegral and the upper integral of $f:\mf{I}\to\mb{R}$.
\end{example}
\begin{solution}{Solution}
Let $\mf{P}=(P_1, \ldots, P_n)$ be a partition of $\mf{I}$. A rectangle $\mf{J}$ in the partition $\mf{P}$  can be written in the form $\mf{J}=\di\prod_{i=1}^n[u_i, v_i]$. By denseness of rational numbers and irrational numbers, there exist a rational number $\alpha_i$ and an irrational number $\beta_i$ in $(u_i, v_i)$. Let $\boldsymbol{\alpha}=(\alpha_1, \ldots, \alpha_n)$ and $\boldsymbol{\beta}=(\beta_1, \ldots, \beta_n)$. Then $\boldsymbol{\alpha}$ and $\boldsymbol{\beta}$ are points in $\mf{J}$, and 
\[0=f(\boldsymbol{\beta})\leq f(\mf{x})\leq f(\boldsymbol{\alpha})=1\hspace{1cm}\text{for all}\;\mf{x}\in\mf{J}.\]
Therefore,
\[m_{\mf{J}}=\inf_{\mf{x}\in\mf{J}}f(\mf{x})=0,\hspace{1cm}M_{\mf{J}}=\sup_{\mf{x}\in\mf{J}}f(\mf{x})=1.\]
It follows that
\[L(f,\mf{P})=\sum_{\mf{J}\in\mathcal{J}_{\mf{P}}}m_{\mf{J}}\,\text{vol}\,(\mf{J})=0,\]\[
U(f,\mf{P})=\sum_{\mf{J}\in\mathcal{J}_{\mf{P}}}M_{\mf{J}}\,\text{vol}\,(\mf{J})=\sum_{\mf{J}\in\mathcal{J}_{\mf{P}}} \text{vol}\,(\mf{J})=\text{vol}\,(\mf{I}).\]
Therefore, \[S_L(f)=\{0\},\quad \text{while}\quad S_U(f)=\{\text{vol}\,(\mf{I})\}.\]
This shows that   the lower inegral and the upper integral of $f:\mf{I}\to\mb{R}$ are given respectively by
\[\underline{\int_{\mf{I}}}f=0\quad\text{and}\quad \overline{\int_{\mf{I}}}f=\text{vol}\,(\mf{I}).\]
\end{solution}

As we mentioned before, one of the motivations to define the integral $f:\mf{I}\to\mb{R}$ is to calculate volumes. Given that $f:\mf{I}\to\mb{R}$ is a nonnegative continuous function defined on the rectangle $\mf{I}$ in $\mb{R}^n$, let  
\[S=\left\{(\mf{x}, y)\,|\, \mf{x}\in\mf{I}, 0\leq y\leq f(\mf{x})\right\},\] which is the solid bounded between $\mf{I}$ and the graph of $f$. It is reasonable to expect that $S$ has a volume, which we denote by $\text{vol}\,(S)$.
We want to define the integral $\di\int_{\mf{I}}f$ so that it gives  $\text{vol}\,(S)$. Notice that if $\mf{P}$ is a partition of $\mf{I}$, then the Darboux lower sum 
\[L(f,\mf{P})=\sum_{\mf{J}\in\mathcal{J}_{\mf{P}}}m_{\mf{J}}\,\text{vol}\,(\mf{J})\]
is the sum of volumes of the collection of rectangles 
\[\left\{\mf{J}\times [0,m_{\mf{J}}]\,|\,\mf{J}\in\mathcal{J}_{\mf{P}}\right\}\]
 in $\mb{R}^{n+1}$, each of which is contained in $S$. Since any two of these rectangles can only intersect on the boundaries, it is reasonable to expect that \[L(f,\mf{P})\leq \text{vol}\,(S).\] Similarly, the Darboux upper sum 
\[U(f,\mf{P})=\sum_{\mf{J}\in\mathcal{J}_{\mf{P}}}M_{\mf{J}}\,\text{vol}\,(\mf{J})\]
is the sum of volumes of the collection of rectangles 
\[\left\{\mf{J}\times [0,M_{\mf{J}}]\,|\,\mf{J}\in\mathcal{J}_{\mf{P}}\right\}\]
 in $\mb{R}^{n+1}$, the union of which contains $S$. Therefore, it is reasonable to expect that
 \[\text{vol}\,(S)\leq U(f,\mf{P}).\]
 Hence, the volume of $S$ should be a number between $L(f,\mf{P})$ and $U(f,\mf{P})$ for any  partition $\mf{P}$. To make the volume well-defined, there should be only one number between $L(f,\mf{P})$ and $U(f,\mf{P})$ for all  partitions $\mf{P}$.  By definition, any number between the lower integral and the upper integral is in between $L(f,\mf{P})$ and $U(f,\mf{P})$ for any   partition $\mf{P}$.  Hence, to have the volume well-defined, we must require the lower integral and the upper integral to be the same. This motivates the following definition of integrability for a general bounded function.
 
 \begin{definition}{Riemann integrability}
 Let $\mf{I}=\di\prod_{i=1}^n[a_i,b_i]$, and let $f:\mf{I}\to\mb{R}$ be a bounded function defined on $\mf{I}$. We say that $f:\mf{I}\to\mb{R}$ is Riemann integrable, or simply integrable, if
 \[\underline{\int_{\mf{I}}}f\;= \; \overline{\int_{\mf{I}}}f.\]
 In this case, we define the integral of $f$ over the rectangle $\mf{I}$ as
 \[\int_{\mf{I}}f\;=\;\underline{\int_{\mf{I}}}f\;= \; \overline{\int_{\mf{I}}}f.\]
 It is the unique number larger than or equal to all Darboux lower sums, and smaller than or equal to all Darboux upper sums.
 
 \end{definition}
 
 \begin{example}{}
 Example \ref{230823_1} says that a constant function $f:\mf{I}\to\mb{R}$, $f(\mf{x})=c$ is integrable and 
 \[\int_{\mf{I}}f=c\,\text{vol}\,(\mf{I}).\]
 
 \end{example}
 \begin{example}[label=230827_1]{}
 The Dirichlet's function defined in Example \ref{230823_5} is not Riemann integrable since the lower integral and the upper integral are not equal.
 
 \end{example}
 
 \begin{highlight}{Leibniz Notation for Riemann Integrals}
 The Leibniz notation of the Riemann integral of $f:\mf{I}\to\mb{R}$ is
 \[\int_{\mf{I}}f(\mf{x})d\mf{x},\quad\text{or equivalently},  \quad\int_{\mf{I}}f(x_1, \ldots, x_n)dx_1\cdots dx_n.\]
 \end{highlight}
 
 As in the single variable case, there are some criteria for Riemann integrability which follows directly from the criteria that the lower integral and the upper integral are the same.
 
 \begin{theorem}[label=230824_7]{}
 Let $\mf{I}=\di\prod_{i=1}^n[a_i,b_i]$, and let $f:\mf{I}\to\mb{R}$ be a bounded function defined on $\mf{I}$. The following   are equivalent.
 \begin{enumerate}[(a)]
 \item The function $f:\mf{I}\to\mb{R}$ is Riemann integrable.
 \item For every $\varepsilon>0$, there is a partition $\mf{P}$ of the rectangle $\mf{I}$ such that 
 \[U(f,\mf{P})-L(f,\mf{P})<\varepsilon.\]
 \end{enumerate}
 \end{theorem}
 
 We define an Archimedes sequence of partitions exactly the same as in the single variable case.
 \begin{definition}{Archimedes Sequence of Partitions}
  Let $\mf{I}=\di\prod_{i=1}^n[a_i,b_i]$, and let $f:\mf{I}\to\mb{R}$ be a bounded function defined on $\mf{I}$.  If $\{\mf{P}_k\}$ is a sequence of partitions of the rectangle $\mf{I}$ such that
  \[\lim_{k\to\infty}\left(U(f,\mf{P}_k)-L(f,\mf{P}_k)\right)=0,\]
  we call $\{\mf{P}_k\}$ an Archimedes sequence of partitions for the function $f$.
 
 \end{definition}
 Then we have the following theorem.
 \begin{theorem}{The Archimedes-Riemann Theorem}
 Let $\mf{I}=\di\prod_{i=1}^n[a_i,b_i]$, and let $f:\mf{I}\to\mb{R}$ be a bounded function defined on $\mf{I}$. The function $f:\mf{I}\to\mb{R}$ is Riemann integrable if and only if   $f$ has an Archimedes sequence of partitions $\{\mf{P}_k\}$. In this case, the integral $\di\int_{\mf{I}}f$ can be computed by
 \[\int_{\mf{I}}f=\lim_{k\to\infty}L(f,\mf{P}_k)=\lim_{k\to\infty}U(f,\mf{P}_k).\]
 
 \end{theorem}
A candidate for an Archimedes sequence of partitions  is the sequence $\{\mf{P}_k\}$, where $\mf{P}_k$ is the uniformly regular partition of $\mf{I}$ into $k^n$ rectangles.  

\begin{example}[label=230823_6]{}Let $\mf{I}=[0,1]\times[0,1]$.
Consider the function $f:\mf{I}\to\mb{R}$ defined as
\[f(x,y)=\begin{cases}1,\quad &\text{if}\;x\geq y,\\0,\quad &\text{if}\;x<y.\end{cases}\]
For $k\in\mb{Z}^+$, let  $\mf{P}_k$ be the uniformly regular partition of $\mf{I}$ into $k^2$ rectangles.
\begin{enumerate}[(a)]
\item For each $k\in\mb{Z}^+$, compute the Darboux lower sum $L(f,\mf{P}_k)$ and the Darboux upper sum $U(f,\mf{P}_k)$.
\item Show that $f:\mf{I}\to\mb{R}$ is Riemann integrable and find the integral $\di\int_{\mf{I}}f$.

\end{enumerate}
\end{example}
\begin{solution}{Solution}
Fixed $k\in\mb{Z}^+$, let $P_k=\left\{u_0, u_1, \ldots, u_k\right\}$, where $u_i=\di \frac{i}{k}$ for $0\leq i\leq k$. Then   $\mf{P}_k=(P_k,P_k)$, and it divides $\mf{I}=[0,1]\times [0,1]$ into the $k^2$ rectangles $\mf{J}_{i,j}$, $1\leq i\leq k$, $1\leq j\leq k$, where
$\di \mf{J}_{i,j}=\left[u_{i-1}, u_i\right]\times \left[u_{j-1}, u_j\right]$. We have \[\text{vol}\,(\mf{J}_{i,j})=\di\frac{1}{k^2}.\]

Let
\[m_{i,j}=\inf_{(x,y)\in \mf{J}_{i,j}}f(x,y)\quad\text{and}\quad M_{i,j}=\sup_{(x,y)\in \mf{J}_{i,j}}f(x,y).\]

Notice that if $i<j-1$, then  
\[x\leq u_i<u_{j-1}\leq y\hspace{1cm}\text{for all}\; (x,y)\in\mf{J}_{i,j}.\]
Hence,
\[f(x,y)=0\hspace{1cm}\text{for all}\; (x,y)\in\mf{J}_{i,j}.\]
This implies that 
\[m_{i,j}=M_{i,j}=0\hspace{1cm}\text{when}\; i<j-1.\]
\bs
If $i\geq j+1$, then  
\[x\geq u_{i-1}\geq u_{j}\geq y\hspace{1cm}\text{for all}\; (x,y)\in\mf{J}_{i,j}.\]
Hence,
\[f(x,y)=1\hspace{1cm}\text{for all}\; (x,y)\in\mf{J}_{i,j}.\]
This implies that 
\[m_{i,j}= M_{i,j}=1\hspace{1cm}\text{when}\;i\geq j+1.\]
When $i=j-1$,  if $(x,y)$ is in $\mf{J}_{i,j}$, \[x\leq u_i=u_{j-1}\leq y,\] and $x=y$ if and only if $(x,y)$ is the point $(u_i, u_{j-1})$. Hence, $f(x,y)=0$ for all $(x,y)\in\mf{J}_{i,j}$, except for $(x,y)=(u_i, u_{j-1})$, where $f(u_i, u_{j-1})=1$.  Hence,
\[m_{i,j}=0,\quad M_{i,j}=1\hspace{1cm}\text{when}\;i =j-1.\]
When $i=j$,  $0\leq f(x,y)\leq 1$ for all $(x,y)\in\mf{J}_{i,j}$. Since $(u_{i-1}, u_j)$ and $(u_i, u_j)$ are in $\mf{J}_{i,j}$, and $f(u_{i-1}, u_j)=0$ while $f(u_i,u_j)=1$, we find that
 \[m_{i,j}=0,\quad M_{i,j}=1\hspace{1cm}\text{when}\;i =j.\]
 
It follows that
\begin{align*}L(f,\mf{P}_k)&=\sum_{i=1}^k\sum_{j=1}^k m_{i,j}\,\text{vol}\,(\mf{J}_{i,j})
 =\sum_{i=2}^k\sum_{j=1}^{i-1} \frac{1}{k^2}\\
 &=\frac{1}{k^2}\sum_{i=2}^k(i-1) =\frac{1}{k^2}\sum_{i=1}^{k-1}i =\frac{k(k-1) }{2k^2}.
 \\U(f,\mf{P}_k)&=\sum_{i=1}^k\sum_{j=1}^k M_{i,j}\,\text{vol}\,(\mf{J}_{i,j})
 =\sum_{i=1}^{k-1}\sum_{j=1}^{i+1}\frac{1}{k^2}+\sum_{j=1}^k\frac{1}{k^2}\\
 &=\frac{1}{k}+\frac{1}{k^2}\sum_{i=1}^{k-1}(i+1)=\frac{1}{ k^2}\left(\frac{k(k+1)}{2}-1+k\right)=\frac{k^2+3k-2}{2k^2}.
 \end{align*}
 \bs
 Since
 \[U(f,\mf{P}_k)-L(f,\mf{P}_k)=\frac{2k-1}{k^2}\hspace{1cm}\text{for all}\;k\in\mb{Z}^+,\]
 we find that
 \[\lim_{k\to\infty}\left(U(f,\mf{P}_k)-L(f,\mf{P}_k)\right)=0.\]
 Hence, $\{\mf{P}_k\}$ is an Archimedes sequence of partitions for $f$. By the Arichimedes-Riemann theorem, $f:\mf{I}\to\mb{R}$ is Riemann integrable, and
 \[\int_{\mf{I}}f=\lim_{k\to\infty}L(f,\mf{P}_k)=\lim_{k\to\infty} \frac{k(k-1)}{2k^2}=\frac{1}{2}.\]
\end{solution}

 \begin{figure}[ht]
\centering
\includegraphics[scale=0.2]{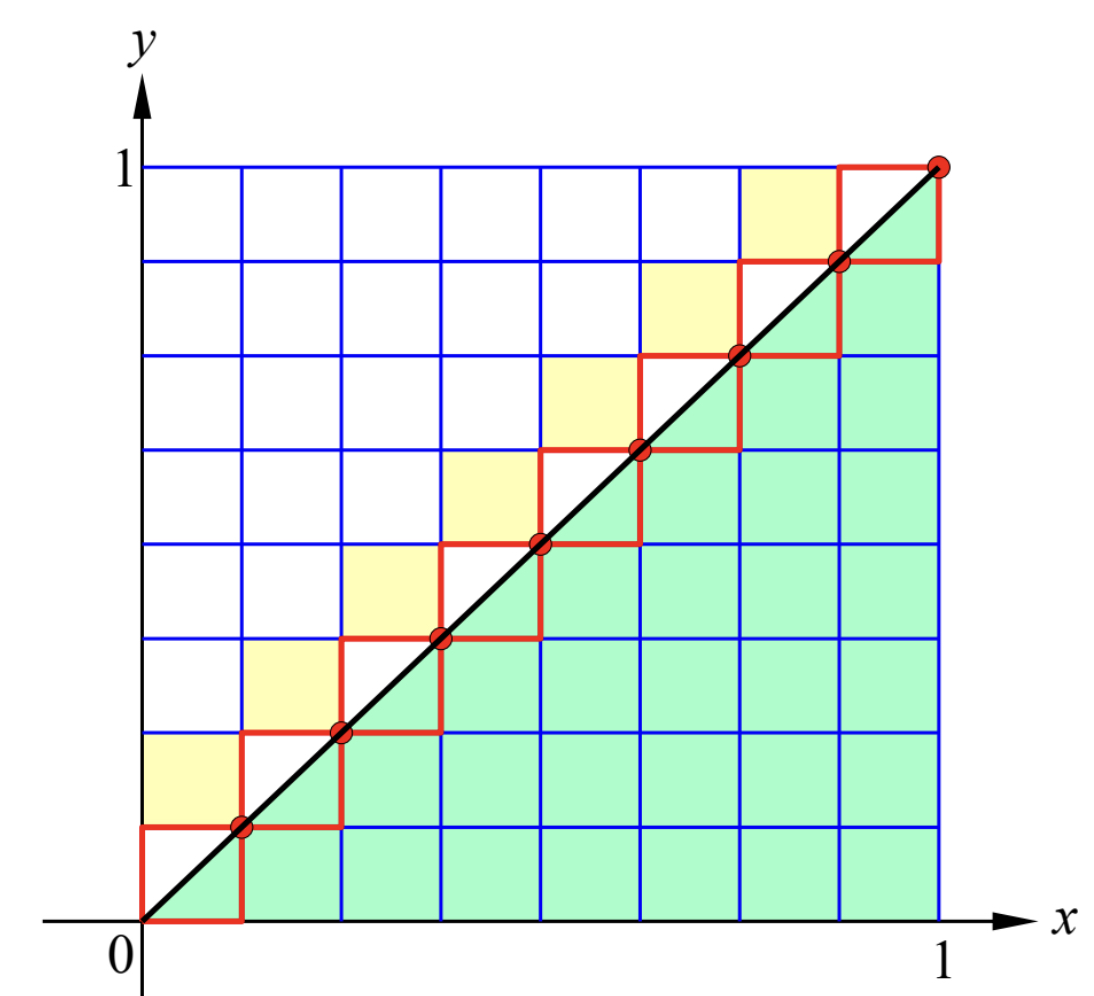}
\caption{This figure illustrates the different cases considered in  Example \ref{230823_6} when $k=8$.}\label{figure86}
\end{figure}

 As in the single variable case, there is an equivalent definition  for Riemann integrability using Riemann sums. 

For a partition $P=\{x_0, x_1, \ldots, x_k\}$ of an interval $[a,b]$, we define the gap of the partition $P$ as
\[|P|=\max\left\{x_i-x_{i-1}\,|\,1\leq i\leq k\right\}.\]
For a closed rectangle $\mf{I}=\di\prod_{i=1}^n[a_i,b_i]$, we  replace the length $x_i-x_{i-1}$ of an interval in the partition by the diameter of a rectangle in the partition. Recall that the diameter of a rectangle $\di \mf{J}=\prod_{i=1}^n[u_i,v_i]$ is
\[\text{diam}\, \mf{J}=\sqrt{(v_1-u_1)^2+\cdots+(v_n-u_n)^2}.\]

\begin{definition}{Gap of a Partition}
Let
$\mf{P}$ be a partition of the rectangle $\mf{I}=\di\prod_{i=1}^n[a_i,b_i]$. Then the gap of the partition $\mf{P}$ is defined as
\[|\mf{P}|=\max\left\{\text{diam}\, \mf{J}\,|\,\mf{J}\in\mathcal{J}_{\mf{P}}\right\}.\]

\end{definition}

\begin{example}{}
Find the gap of the partition $\mf{P}=(P_1, P_2)$ of the rectangle $\mf{I}=[-2,9]\times [1,6]$ defined in Example \ref{230821_1}, where $P_1=\{-2,0,4,9\}$ and $P_2=\{1,3,6\}$.\end{example}
\begin{solution}{Solution} The length of the three invervals in the partition $P_1=\{-2,0,4,9\}$ of the interval $[-2, 9]$ are $2, 4$ and $5$ respectively. The lengths of the two intervals in the partition $P_2=\{1,3,6\}$ of the interval $[1,6]$ are $2$ and $3$ respectively.
Therefore, the diameters of the 6 rectangles in the partition $\mf{P}$ are
\begin{gather*}
\sqrt{2^2+2^2},\quad\sqrt{4^2+2^2},\quad\sqrt{5^2+2^2},\\
\sqrt{2^2+3^2},\quad
\sqrt{4^2+3^2},\quad \sqrt{5^2+3^2}.
\end{gather*}From this, we see that the gap of $\mf{P}$ is $\sqrt{5^2+3^2}=\sqrt{34}$.
\end{solution}

In the example above, notice that $|P_1|=5$ and $|P_2|=3$. 
In general, it is not difficult to see the following.
\begin{proposition}{}
Let $\mf{P}=(P_1, \ldots, P_n)$ be a partition of the closed rectangle $\mf{I}=\di\prod_{i=1}^n[a_i,b_i]$. Then 
\[|\mf{P}|=\sqrt{|P_1|^2+\cdots+|P_n|^2}.\]
\end{proposition}

The following theorem gives equivalent definitions of Riemann integrability of a bounded function. 
\begin{theorem}[label=230824_1]{Equivalent Definitions for Riemann Integrability}
Let $\mf{I}=\di\prod_{i=1}^n[a_i,b_i]$, and let $f:\mf{I}\to\mb{R}$ be a bounded function defined on $\mf{I}$. The following three statements are equivalent for  saying that $f:\mf{I}\to\mb{R}$ is Riemann integrable.
\begin{enumerate}[(a)]
\item
The lower integral and the upper integral are the same. Namely, \[\underline{\int_{\mf{I}}}f\;=\;\overline{\int_{\mf{I}}}f.\]
\item There exists a number $I$ that satisfies the following. For any $\varepsilon>0$, there exists a $\delta>0$ such that   if $\mf{P}$ is a partition of the rectangle $\mf{I}$ with $|\mf{P}|<\delta$, then  
\[\left|R(f,\mf{P}, A)-I\right|<\varepsilon \]for any choice of intermediate points $A=\left\{\boldsymbol{\xi}_{\mf{J}}\right\}$ for the partition $\mf{P}$.
\item For any $\varepsilon>0$, there exists a $\delta>0$  such that  if $\mf{P}$ is a partition of the rectangle $\mf{I}$ with $|\mf{P}|<\delta$, then  
\[U(f,\mf{P})-L(f,\mf{P})<\varepsilon.\]
\end{enumerate}

\end{theorem}
The most   useful definition is in fact the second one in terms of Riemann sums. It says that a bounded function $f:\mf{I}\to\mb{R}$ is Riemann integrable if the limit 
\[\lim_{|\mf{P}|\to 0}R(f, \mf{P}, A)\] exists. 
As a  consequence of Theorem \ref{230824_1}, we have the following.
\begin{theorem}[label=230824_2]{}
Let $\mf{I}=\di\prod_{i=1}^n[a_i,b_i]$, and let $f:\mf{I}\to\mb{R}$ be a bounded function defined on $\mf{I}$. If $f:\mf{I}\to\mb{R}$ is Riemann integrable, then for any sequence $\{\mf{P}_k\}$  of partitions of $\mf{I}$ satisfying 
\[\lim_{k\to\infty}|\mf{P}_k|=0,\]
we have
\begin{enumerate}[(i)]
\item $\di \int_{\mf{I}}f=\lim_{k\to\infty}L(f,\mf{P}_k)=\lim_{k\to\infty}U(f,\mf{P}_k)$.
\item $\di\int_{\mf{I}}f=\lim_{k\to\infty}R(f,\mf{P}_k,A_k)$, where for each $k\in\mb{Z}^+$, $A_k$ is a choice of intermediate points for the partition $\mf{P}_k$.
\end{enumerate}
\end{theorem}The proof is exactly the same as the single variable case.
The contrapositive of Theorem \ref{230824_2} gives the following.
\begin{theorem}[label=230824_11]{}
Let $\mf{I}=\di\prod_{i=1}^n[a_i,b_i]$, and let $f:\mf{I}\to\mb{R}$ be a bounded function defined on $\mf{I}$.  Assume that    $\{\mf{P}_k\}$   is a sequence of partitions of $\mf{I}$ such that
\[\lim_{k\to\infty}|\mf{P}_k|=0.\] 
\begin{enumerate}[(a)]
\item If for each $k\in\mb{Z}^+$, there exists a choice of intermediate points $A_k$ for the partition $\mf{P}_k$ such that the limit 
$\di \lim_{k\to\infty}R(f,\mf{P}_k,A_k)$  does not exist, then $f:\mf{I}\to\mb{R}$ is not Riemann integrable.
\item If for each $k\in\mb{Z}^+$, there exist two  choices of intermediate points $A_k$ and $B_k$ for the partition $\mf{P}_k$ so that the two limits $\di\lim_{k\to\infty}R(f,\mf{P}_k,A_k)$ and $\di\lim_{k\to\infty}R(f,\mf{P}_k, B_k)$ are not the same, then $f:\mf{I}\to\mb{R}$ is not Riemann integrable.
\end{enumerate}
\end{theorem}
Theorem \ref{230824_11} is useful for justifying that a bounded function is not Riemann integrable, without having to compute the lower integral or the upper integral.
To apply this theorem, we usually consider the sequence of partitions $\{\mf{P}_k\}$, where $\mf{P}_k$ is the uniformly regular partition of $\mf{I}$ into $k^n$ rectangles.
\begin{example}{}
Let $\mf{I}=[0,1]\times [0,1]$, and let $f:\mf{I}\to\mb{R}$ be the function defined as
\[f(x,y)=\begin{cases}0,\quad  & \text{if $x$ is rational},\\y,\quad &\text{if $x$ is irrational}.\end{cases}\]
Show that $f:\mf{I}\to\mb{R}$ is not Riemann integrable.
\end{example}
\begin{solution}{Solution}
For $k\in\mb{Z}^+$, let $\mf{P}_k$ be the  uniformly regular partition of $\mf{I}$ into $k^2$ rectangles. Then $\mf{P}_k=(P_k, P_k)$, where $P_k=\left\{u_0, u_1, \ldots, u_k\right\}$ with $u_i=\di \frac{i}{k}$ when $0\leq i\leq k$. Notice that $|\mf{P}_k|=\di\frac{\sqrt{2}}{k}$, and so
$\di \lim_{k\to\infty}\mf{P}_k=0$.

The partition $\mf{P}_k$ divides the square $\mf{I}$ into $k^2$ squares $\mf{J}_{i,j}$, $1\leq i\leq k$, $1\leq j\le k$, where $\mf{J}_{i,j}=[u_{i-1}, u_i]\times [u_{j-1},u_j]$. For $1\leq i\leq k$, since irrational numbers are dense, there is an irrational number $c_i$ in the interval $(u_{i-1}, u_i)$. For $1\leq i\leq k$, $1\leq j\leq k$, let $\boldsymbol{\alpha}_{i,j}$ and $\boldsymbol{\beta}_{i,j}$ be the points in $\mf{J}_{i,j}$ given respectively by
\[\boldsymbol{\alpha}_{i,j}= (u_{i}, u_{j}),\hspace{1cm}\boldsymbol{\beta}_{i,j}= (c_{i}, u_{j}).\]
Then
\[f(\boldsymbol{\alpha}_{i,j})=0,\hspace{1cm}f(\boldsymbol{\beta}_{i,j})=u_j.\]
Let $\di A_k=\left\{\boldsymbol{\alpha}_{i,j} \right\}$  and $B_k =\left\{\boldsymbol{\beta}_{i,j} \right\}$. 
Then the Riemann sums $R(f,\mf{P}_k, A_k)$ and $R(f,\mf{P}_k, B_k)$ are given respectively by
\[R(f,\mf{P}_k, A_k)=\sum_{i=1}^{k}\sum_{j=1}^k f(\boldsymbol{\alpha}_{i,j})\;\text{vol}\,(\mf{J}_{i,j})=0,\]
\bs
and
\begin{align*}
R(f,\mf{P}_k, B_k)&=\sum_{i=1}^{k}\sum_{j=1}^k f(\boldsymbol{\beta}_{i,j})\;\text{vol}\,(\mf{J}_{i,j})=\sum_{i=1}^k\sum_{j=1}^k 
\frac{j}{k}\times \frac{1}{k^2}\\
&=\frac{k\times k(k+1)}{2k^3}=\frac{k+1}{2k}.
\end{align*}Therefore, we find that
\[\lim_{k\to\infty} R(f,\mf{P}_k, A_k)=0,\hspace{1cm} \lim_{k\to\infty} R(f,\mf{P}_k, B_k)=\frac{1}{2}.\]
Since the two limits are not the same, we conclude that $f:\mf{I}\to\mb{R}$ is not Riemann integrable.
\end{solution}

Now we return to the proof of Theorem \ref{230824_1}. 
 To prove this theorem, it is easier to   show that (a) is equivalent to (c), and (b) is equivalent to (c).  
We will prove the equivalence of (a) and (c). The proof of the equivalence of (b) and (c) is left to the exercises.  It  is a consequence of the inequality
\[L(f,\mf{P})\leq R(f, \mf{P}, A)\leq U(f,\mf{P}),\] which holds for any partition $\mf{P}$ of the rectangle $\mf{I}$, and any choice of intermediate points $A$ for the partition $\mf{P}$.  

By Theorem \ref{230824_7}, (a) is equivalent to

\vspace{-0.2cm}
\begin{enumerate}
\item[(a$^{\prime}$)] For every $\varepsilon>0$, there is a partition $\mf{P}$ of $\mf{I}$ such that
\[U(f,\mf{P})-L(f,\mf{P})<\varepsilon.\]

\end{enumerate}
Thus, to prove the equivalence of (a) and (c), it is sufficient to show the equivalence of (a$^{\prime}$) and (c). But then
(c) implies (a$^{\prime}$) is  obvious. Hence, we are left with the most technical part, which is the proof of (a$^{\prime}$) implies (c). 

We formulate this as a standalone theorem.
\begin{theorem}[label=230824_8]{}
Let $\mf{I}=\di\prod_{i=1}^n[a_i,b_i]$, and let $\mf{P}_0$ be a fixed a partition of $\mf{I}$.  Given that $f:\mf{I}\to\mb{R}$ is a bounded   function defined on $\mf{I}$, for any $\varepsilon>0$, there is a $\delta>0$ such that for all partitions $\mf{P}$ of $\mf{I}$, if $|\mf{P}|<\delta$, then 
\begin{equation}\label{230824_10}U(f,\mf{P})-L(f,\mf{P})<U(f,\mf{P}_0)-L(f,\mf{P}_0)+\varepsilon.\end{equation}
\end{theorem}

\begin{highlight}{}
If Theorem \ref{230824_8} is proved, we can show that (a$^{\prime}$) implies (c) in Theorem \ref{230824_1} as follows. Given $\varepsilon>0$,  (a$^{\prime}$)  implies that we can choose a $\mf{P}_0$ such that 
\[U(f,\mf{P}_0)-L(f,\mf{P}_0)<\frac{\varepsilon}{2}.\]By Theorem \ref{230824_8}, there is a $\delta>0$ such that for all partitions $\mf{P}$ of $\mf{I}$, if $|\mf{P}|<\delta$, then 
\[U(f,\mf{P})-L(f,\mf{P})<U(f,\mf{P}_0)-L(f,\mf{P}_0)+\frac{\varepsilon}{2}<\varepsilon.\] This proves that  (a$^{\prime}$) implies (c). 
\end{highlight}

Hence, it remains for us to prove theorem \ref{230824_8}. Let us introduce some additional notations. Given the rectangle $\mf{I}=\di\prod_{i=1}^n[a_i,b_i]$, for $1\leq i\leq n$, let
\begin{equation}\label{230824_9}\begin{split}\mathcal{S}_i&=\frac{\text{vol}\,(\mf{I})}{b_i-a_i}\\&=(b_1-a_1)\times \cdots\times (b_{i-1}-a_{i-1})(b_{i+1}-a_{i+1})\times \cdots \times (b_n-a_n).\end{split}\end{equation}
This is the area of   the bounday of $\mf{I}$ that is contained in the hyperplane $x_i=a_i$ or $x_i=b_i$. For example, when $n=2$, $\mf{I}=[a_1, b_1]\times [a_2, b_2]$,  
$\mathcal{S}_1=b_2-a_2$ is the length of the vertical side, while $\mathcal{S}_2=b_1-a_1$ is the length of the horizontal side of the rectangle $\mf{I}$.

\begin{myproof}{\linkt Proof of Theorem \ref{230824_8}}
Since $f:\mf{I}\to\mb{R}$ is bounded, there is a positive number $M$ such that
\[|f(\mf{x})|\leq M\hspace{1cm}\text{for all}\;\mf{x}\in\mf{I}.\]
Assume that $\mf{P}_0=(\widetilde{P}_1, \ldots,\widetilde{P}_n)$.   For $1\leq i\leq n$, let $k_i$ be the number of intervals in the partition $\widetilde{P}_i$.  Let \[K=\max\{k_1, \ldots, k_n\},\] and 
\[\mathcal{S}=\mathcal{S}_1+\cdots+\mathcal{S}_n,\] where $\mathcal{S}_i$, $1\leq i\leq n$ are defined by \eqref{230824_9}. Given $\varepsilon>0$, let
\[\delta= \frac{\varepsilon}{4MK\mathcal{S}}.\]
Then $\delta>0$. If $\mf{P}=(P_1, \ldots, P_n)$ is a partition of $\mf{I}$ with $|\mf{P}|<\delta$, we want to show that \eqref{230824_10} holds. 
 Let $\mf{P}^*=(P_1^*, \ldots, P_n^*)$ be the common refinement of $\mf{P}_0$ and $\mf{P}$ such that $P_i^*$ is the partition of $[a_i,b_i]$ that contains  all the partition points of $\widetilde{P}_i$ and $P_i$. 
For $1\leq i\leq n$,  let $U_i$ be the collection of intervals in $P_i$ which contain  partition points of $\widetilde{P}_i$, and let $V_i$ be the collection of the intervals of $P_i$ that is not in $U_i$. 
  Each interval in $V_i$ must be in the interior of one of the intervals in $\widetilde{P}_i$. Thus, each interval in $V_i$ is an interval in the partition $P_i^*$.
Since each partition point of $\widetilde{P}_i$ can be contained in at most two intervals of $P_i$, but the first and last partition points of $P_i$ and $\widetilde{P}_i$ are the same, we find that $|U_i|\leq 2k_i$. 

 Since $|P_i|\leq |\mf{P}|<\delta$, each interval in $P_i$ has length less than $\delta$. Therefore, the sum of the lengths of the intervals in $U_i$ is less than $2k_i\delta$. Let
\[\mathcal{Q}_i=\left\{\mf{J}\in\mathcal{J}_{\mf{P}}\,|\, \text{the $i^{\text{th}}$-edge of $\mf{J}$ is from $U_i$}\right\}.\]
Then 
\[\sum_{\mf{J}\in\mathcal{Q}_i}\text{vol}\,(\mf{J})<2k_i\delta \mathcal{S}_i\leq 2K\delta \mathcal{S}_i.\]
 \end{myproof}
 \begin{figure}[ht]
\centering
\includegraphics[scale=0.2]{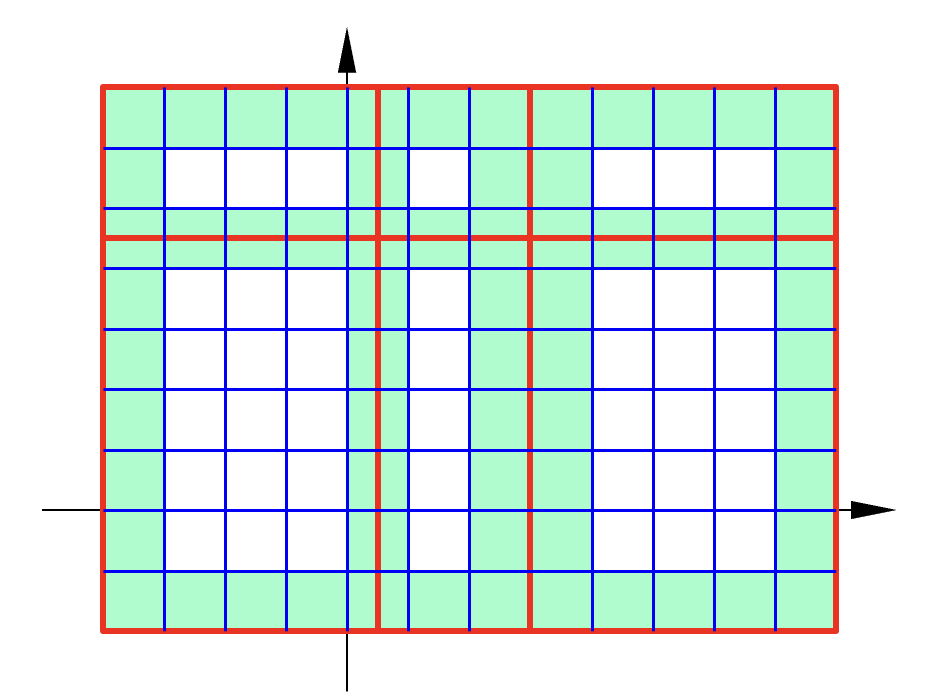}
\caption{The   partitions $\mf{P}_0$ and $\mf{P}$ in the proof of Theorem \ref{230824_8}, $\mf{P}_0$ is the partition with red grids, while $\mf{P}$ is the partition with blue grids. Those shaded rectangles are rectangles in $\mf{P}$ that contain  partition points of $\mf{P}_0$.}\label{figure90}
\end{figure}
 \begin{myproof}{} 
 Now let \[\mathcal{Q}=\di\bigcup_{i=1}^n\mathcal{Q}_i.\]

 Then \[\sum_{\mf{J}\in\mathcal{Q}}\text{vol}\,(\mf{J})<2K\delta \sum_{i=1}^n\mathcal{S}_i=2K\delta\mathcal{S}.\]
 For each of the rectangles $\mf{J}$ that is in $\mathcal{Q}$, we do a simple estimate
 \[M_{\mf{J}}-m_{\mf{J}}\leq 2M.\]
 Therefore,
 \[\sum_{\mf{J}\in\mathcal{Q}}\left(M_{\mf{J}}-m_{\mf{J}}\right)\text{vol}\,(\mf{J})<4MK\delta\mathcal{S}\leq \varepsilon.\]
For the rectangles $\mf{J}$ that are in $\mathcal{J}_{\mf{P}}\setminus\mathcal{Q}$, each of them is a rectangle in the partition $\mf{P}^*$.
 Therefore,
 \[\sum_{\mf{J}\in\mathcal{J}_{\mf{P}}\setminus\mathcal{Q}}\left(M_{\mf{J}}-m_{\mf{J}}\right)\text{vol}\,(\mf{J})\leq U(f,\mf{P}^*)-L(f,\mf{P}^*)\leq U(f,\mf{P}_0)-L(f,\mf{P}_0).\]
 \bp
 Hence,
 \begin{align*}
 U(f,\mf{P})-L(f,\mf{P}) &=\sum_{\mf{J}\in\mathcal{J}_{\mf{P}} }\left(M_{\mf{J}}-m_{\mf{J}}\right)\text{vol}\,(\mf{J})\\
 &=\sum_{\mf{J}\in\mathcal{J}_{\mf{P}}\setminus\mathcal{Q}}\left(M_{\mf{J}}-m_{\mf{J}}\right)\text{vol}\,(\mf{J})+\sum_{\mf{J}\in\mathcal{Q}}\left(M_{\mf{J}}-m_{\mf{J}}\right)\text{vol}\,(\mf{J})\\
 &<U(f,\mf{P}_0)-L(f,\mf{P}_0)+\varepsilon.
 \end{align*}
 This completes the proof.
\end{myproof}

Finally we  extend Riemann integrals to functions $f:\mk{D}\to\mb{R}$ that are defined on bounded subsets $\mk{D}$ of $\mb{R}^n$.  If $\mk{D}$ is bounded, there is a positive number $L$ such that \[\Vert\mf{x}\Vert\leq L\hspace{1cm} \text{for all}\;\mf{x}  \in  \mk{D}. \]This implies that $\mk{D}$ is contained in the closed rectangle $\mf{I}_L=\di\prod_{i=1}^n[-L,L]$. To define the Riemann integral of $f:\mk{D}\to\mb{R}$, we need to extend the domain of $f$ from $\mk{D}$ to $\mf{I}_L$. To avoid affecting the integral, we should extend by zero.

\begin{definition}{Zero Extension}
Let $\mk{D}$ be a subset of $\mb{R}^n$, and let $f:\mk{D}\to\mb{R}$ be a function defined on $\mk{D}$. The zero extension of $f:\mk{D}\to\mb{R}$ is the function $  \check{f}:\mb{R}^n\to\mb{R}$ which is defined as
\[      \check{f}(\mf{x})=\begin{cases} f(\mf{x}),\quad &\text{if}\;\mf{x}\in \mk{D},\\
0,\quad &\text{if}\;\mf{x}\notin\mk{D}.\end{cases}\]
If $\mathcal{U}$ is any subset of $\mb{R}^n$ that contains $\mk{D}$, then the zero extension of $f$ to $\mathcal{U}$ is the function $      \check{f}:\mathcal{U}\to\mb{R}$.
\end{definition}
Obviously, if $f:\mk{D}\to\mb{R}$ is a bounded function, its zero extension  $  \check{f}:\mb{R}^n\to\mb{R}$ is also bounded.
Since we have defined Riemann integrability for a bounded function $g:\mf{I}\to\mb{R}$ that is defined on a closed rectangle $\mf{I}$, it is natural to say that a function $f:\mk{D}\to\mb{R}$ is Riemann integrable if its zero extension $  \check{f}:\mf{I}\to\mb{R}$  to a closed rectangle $\mf{I}$ is Riemann integrable, and define
\[\int_{\mk{D}}f=\int_{\mf{I}}  \check{f}.\] For this to be unambiguous, we have to check that if $\mf{I}_1$ and $\mf{I}_2$ are closed rectangles that contain the bounded set $\mk{D}$, the zero extension $  \check{f}:\mf{I}_1\to\mb{R}$ is Riemann integrable if and only if the zero extension $  \check{f}:\mf{I}_2\to\mb{R}$ is Riemann integrable. Moreover, \[\int_{\mf{I}_1}  \check{f}=\int_{\mf{I}_2}  \check{f}.\] This small technicality would be proved in Section \ref{sec6.2}. Assuming this, we can give the following formal definition for Riemann integrality of a bounded function defined on a bounded domain.

\begin{definition}{Riemann Integrals of General Functions}
Let $\mk{D}$ be a bounded subset of $\mb{R}^n$, and let $\mf{I}=\di\prod_{i=1}^n[a_i,b_i]$ be a closed rectangle in $\mb{R}^n$ that contains $\mk{D}$. Given that $f:\mk{D}\to\mb{R}$ is a bounded function defined on $\mk{D}$, we say that $f:\mk{D}\to\mb{R}$ is Riemann integrable if its zero extension $  \check{f}:\mf{I}\to\mb{R}$  is Riemann integrable. If this is the case, we define the integral of $f$ over $\mk{D}$ as
\[\int_{\mk{D}}f=\int_{\mf{I}}  \check{f}.\]
\end{definition}

\begin{example}[label=230825_1]{} 
Let  $\mf{I}=[0,1]\times [0,1]$, and let $f:\mf{I}\to\mb{R}$ be the function defined as
\[f(x,y)=\begin{cases}1,\quad &\text{if}\;x\geq y,\\0,\quad &\text{if}\;x<y.\end{cases}\]
which is considered in Example \ref{230823_6}. 
Let 
\[\mk{D}=\left\{(x,y)\in\mf{I}\,|\,   y\leq x \right\},\]and let $g:\mk{D}\to\mb{R}$ be the constant function $g(\mf{x})=1$. Then $f:\mf{I}\to\mb{R}$ is the zero extension of $g$ to the square $\mf{I}$ that contains $\mk{D}$. 

\be
In Example \ref{230823_6}, we have shown that $f:\mf{I}\to\mb{R}$ is Riemann integrable and \[\int_{\mf{I}}f(\mf{x})d\mf{x}=\frac{1}{2}.\] Therefore, $g:\mk{D}\to\mb{R}$ is Riemann integrable and 
\[\int_{\mk{D}}g(\mf{x})d\mf{x}=\frac{1}{2}.\]
\end{example2}
\begin{remark}{}
Here we  make two remarks about the Riemann integrals.
\begin{enumerate}[1.]
\item
When $f:\mk{D}\to\mb{R}$ is the constant function, we should expect that it is Riemann integrable if and only if $\mk{D}$ has a volume, which should be defined as
\[\text{vol}\,(\mk{D})=\int_{\mk{D} } d\mf{x}.\]
\item If $f:\mk{D}\to\mb{R}$ is a nonnegative continuous function defined on the bounded set $\mk{D}$ that has a volume, we would expect that $f:\mk{D}\to\mb{R}$ is Riemann integrable, and the integral $\di\int_{\mk{D}}f(\mf{x})d\mf{x}$ gives the volume of the solid bounded between $\mk{D}$ and the graph of $f$.
\end{enumerate}
\end{remark}

In Section \ref{sec6.3}, we will give a characterization of sets $\mk{D}$ that have volumes. We will also prove that if $f:\mk{D}\to\mb{R}$ is a continuous function defined on a set $\mk{D}$ that has volume, then $f:\mk{D}\to\mb{R}$ is Riemann integrable.

\vp
\noindent
{\bf \large Exercises  \thesection}
\setcounter{myquestion}{1}

\begin{question}{\themyquestion}
 Let $\mf{I}=[-5, 8]\times [2, 5]$, and let  $\mf{P}=(P_1, P_2)$ be the partition  of $\mf{I}$ with $P_1=\{-5, -1,2,7,8\}$ and $P_2=\{2, 4, 5\}$. 
 Find gap of the partition $\mf{P}$.
 
\end{question}
\atc

\begin{question}{\themyquestion}
  Let $\mf{I}=[-5, 8]\times [2, 5]$, and let $f:\mf{I}\to\mb{R}$ be the function defined as $f(x,y)=x^2+2y$. Consider the partition $\mf{P}=(P_1, P_2)$   of $\mf{I}$ with $P_1=\{-5, -1,2,7,8\}$ and $P_2=\{2, 4, 5\}$. 
 Find the Darboux lower sum $L(f,\mf{P})$ and the Darboux upper sum $U(f, \mf{P})$.
\end{question}
\atc

\begin{question}{\themyquestion}
   Let $\mf{I}=[-5, 8]\times [2, 5]$, and let $f:\mf{I}\to\mb{R}$ be the function defined as $f(x,y)=x^2+2y$. Consider the partition $\mf{P}=(P_1, P_2)$   of $\mf{I}$ with $P_1=\{-5, -1,2,7,8\}$ and $P_2=\{2, 4, 5\}$. 
 For each rectangle $\mf{J}=[a,b]\times [c,d]$ in the partition $\mf{P}$, let $\boldsymbol{\alpha}_{\mf{J}}=(a,c)$ and $\boldsymbol{\beta}_{\mf{J}}=(b, d)$. Find the Riemann sums $R(f,\mf{P},A)$ and $R(f,\mf{P}, B)$, where  $A=\{\boldsymbol{\alpha}_{\mf{J}}\}$ and $B= \{\boldsymbol{\beta}_{\mf{J}}\}$.
\end{question}
\atc

\begin{question}{\themyquestion}
 Let $\di\mf{I}=[-1,1]\times [2, 5]$, and let $f:\mf{I}\to\mb{R}$ be the function defined as 
\[f(x,y)=\begin{cases}1,\quad &\text{if $x$ and $y$ are rational},\\0,\quad &\text{otherwise}.\end{cases}\]
\begin{enumerate}[(a)]
\item Given that $\mf{P}$ is a partition of $\mf{I}$, find the Darboux lower sum $L(f,\mf{P})$ and the Darboux upper sum 
$U(f,\mf{P})$.
\item  Find the lower integral $\di\underline{\int_{\mf{I}}} f$ and the upper integral  $\di\overline{\int_{\mf{I}}} f$.
\item Explain why $f:\mf{I}\to\mb{R}$ is not Riemann integrable.

\end{enumerate}
\end{question}
 
\atc

\begin{question}{\themyquestion}
 Let $\mf{I}=[0,4]\times [0,2]$.
 Consider the function $f:\mf{I}\to\mb{R}$ defined as
\[f(x,y)=2x+3y+1.\]
For $k\in\mb{Z}^+$, let   $\mf{P}_k$ be the uniformly regular partition of $\mf{I}=[0,4]\times [0,2]$ into $k^2$ rectangles.
\begin{enumerate}[(a)]
\item For each $k\in\mb{Z}^+$, compute the Darboux lower sum $L(f,\mf{P}_k)$ and the Darboux upper sum $U(f,\mf{P}_k)$.
\item Show that $f:\mf{I}\to\mb{R}$ is Riemann integrable and find the integral $\di\int_{\mf{I}}f$.

\end{enumerate}
\end{question}

\atc
\begin{question}{\themyquestion}
 Let $\mf{I}=\di\prod_{i=1}^n[a_i,b_i]$, and let $f:\mf{I}\to\mb{R}$ be a function defined on $\mf{I}$. Show that the following are equivalent.
\begin{enumerate}[(a)]
 \item There exists a number $I$ that satisfies the following. For any $\varepsilon>0$, there exists a $\delta>0$ such that   if $\mf{P}$ is a partition of the rectangle $\mf{I}$ with $|\mf{P}|<\delta$, then  
\[\left|R(f,\mf{P}, A)-I\right|<\varepsilon \]for any choice of intermediate points $A=\left\{\boldsymbol{\xi}_{\mf{J}}\right\}$ for the partition $\mf{P}$.
\item For any $\varepsilon>0$, there exists a $\delta>0$  such that  if $\mf{P}$ is a partition of the rectangle $\mf{I}$ with $|\mf{P}|<\delta$, then  
\[U(f,\mf{P})-L(f,\mf{P})<\varepsilon.\]
\end{enumerate}
\end{question}
 
\section{Properties of Riemann Integrals}\label{sec6.2}

 In this section, we discuss  properties of Riemann integrals.
 Let us first consider Riemann integrals of functions $f:\mf{I}\to\mb{R}$ defined on closed rectangles of the form $\mf{I}=\di\prod_{i=1}^n[a_i,b_i]$. Using some of these properties, we  prove that the definition of Riemann integrabililty for functions $f:\mk{D}\to\mb{R}$ defined on general bounded  sets, as given in Section \ref{sec6.1}, is unambiguous. Finally, we will extend the properties of Riemann integrals to functions $f:\mk{D}\to\mb{R}$ defined on bounded sets.

Linearity   is one of the most important properties. For functions defined on closed rectangles of the form $\mf{I}=\di\prod_{i=1}^n[a_i,b_i]$, the proof is straightforward using the Riemann sum definition of Riemann integrability, as in the single variable case.
\begin{theorem}[label=230826_3]{Linearity}
Let $\mf{I}=\di\prod_{i=1}^n[a_i,b_i]$, and let $f:\mf{I}\to\mb{R}$ and $g:\mf{I}\to\mb{R}$ be Riemann integrable functions. For any real numbers $\alpha$ and $\beta$, $(\alpha f+\beta g):\mf{I}\to\mb{R}$ is also Riemann integrable, and 
\[\int_{\mf{I}}(\alpha f+\beta g)=\alpha\int_{\mf{I}}f+\beta\int_{\mf{I}}g.\]
\end{theorem}
\begin{myproof}{Sketch of Proof}
If $\mf{P}$ is a partition of $\mf{I}$ and $A$ is a set of intermediate points for $\mf{P}$, then 
\[R(\alpha f+\beta g, \mf{P}, A)=\alpha R(f, \mf{P}, A)+\beta R(g, \mf{P}, A).\] The results follows by taking the $|\mf{P}|\to 0$ limit.
 
\end{myproof}

\begin{example}{}
Let $\mf{I}=[0,2]\times [0,2]$, and let $f:\mf{I}\to\mb{R}$ and $g:\mf{I}\to\mb{R}$  be Riemann integrable functions. Find  the integrals $\di \int_{\mf{I}} f$ and $\di\int_{\mf{I}}g$ if  

\vspace{-0.4cm} 
\[f(x,y)=g(y,x)\quad\text{and}\quad (f+g)(x,y)=6\hspace{1cm}\text{for all}\;(x,y)\in \mf{I}.\]

\end{example}
\begin{solution}{Solution}
Since $\mf{I}$ is symmetric with respect to the line $y=x$ and $f(x,y)=g(y,x)$ for all $(x,y)\in \mf{I}$, we have $\di \int_{\mf{I}} f=\di\int_{\mf{I}}g$. By linearity,
\[\int_{\mf{I}} f+\di\int_{\mf{I}}g=\int_{\mf{I}}(f+g)=6\times \text{vol}\,(\mf{I})=24.\]
Hence, \[\int_{\mf{I}} f=\di\int_{\mf{I}}g=12.\]
\end{solution}
 
The following theorem is about the integral of a nonnegative function.
\begin{theorem}[label=230825_2]{}
Let $\mf{I}=\di\prod_{i=1}^n[a_i,b_i]$, and let $f:\mf{I}\to\mb{R}$ be a bounded function defined on $\mf{I}$. Assume that
 $f(\mf{x})\geq 0$ for all $\mf{x}$ in $\mf{I}$. If $f:\mf{I}\to\mb{R}$ is Riemann integrable, then
 
 \vspace{-0.4cm}
 \[\int_{\mf{I}}f \;\geq \;0.\]
 
\end{theorem}
\begin{myproof}{Proof}
  For any partition $\mf{P}$ of $\mf{I}$, $L(f,\mf{P})\geq 0$. Therefore,
\[\int_{\mf{I}}f \;=\;\underline{\int_{\mf{I}}}\,f \;\geq \;L(f,\mf{P})\;\geq \;0.\]
\end{myproof}

The monotonicity theorem then follows from linearity and Theorem \ref{230825_2}.
\begin{theorem}[label=230825_3]{Monotonicity}
Let $\mf{I}=\di\prod_{i=1}^n[a_i,b_i]$, and let $f:\mf{I}\to\mb{R}$ and $g:\mf{I}\to\mb{R}$ be Riemann integrable functions. If
 $f(\mf{x})\geq g(\mf{x})$ for all $\mf{x}$ in $\mf{I}$, then
 
 \vspace{-0.4cm}
 \[\int_{\mf{I}}f \;\geq \;\int_{\mf{I}}g.\]
 
\end{theorem}
\begin{myproof}{Proof}
By linearity, the function $(f-g):\mf{I}\to\mb{R}$ is integrable, and 
\[\int_{\mf{I}}(f-g)\,=\,\int_{\mf{I}}f\,-\,\int_{\mf{I}}g.\]
By Theorem \ref{230825_2}, $\di \int_{\mf{I}}(f-g)\geq 0$, and the assertion follows.
\end{myproof}

The next important property is the additivity of the Riemann integrals.
\begin{theorem}{Additivity}
Let $\mf{I}=\di\prod_{i=1}^n[a_i,b_i]$, and let $\mf{P}_0$ be a partition of $\mf{I}$. If $f:\mf{I}\to\mb{R}$ is a bounded function defined on $\mf{I}$, then $f:\mf{I}\to\mb{R}$ is Riemann integrable if and only if for each $\mf{J}\in\mathcal{J}_{\mf{P}_0}$, $f:\mf{J}\to\mb{R}$ is Riemann integrable. In such case, we also have
\[\int_{\mf{I}}f\;=\;\sum_{\mf{J}\in\mathcal{J}_{\mf{P}_0}}\int_{\mf{J}}f.\]
\end{theorem}
\begin{myproof}{Proof}
It is sufficient to consider the case that $\mf{P}_0=(P_1, \ldots, P_n)$ divides $\mf{I}$ into two rectangles $\mf{I}_1$ and $\mf{I}_2$ by having a partition point $c$ inside the $j^{\text{th}}$-edge $[a_j, b_j]$ for some $1\leq j\leq n$. Namely,  $P_j=\{a_j, c, b_j\}$, and for $i\neq j$, $P_i=\{a_i, b_i\}$. The general case can be proved by induction, adding one partition point at a time. 

Assume that $f:\mf{I}\to\mb{R}$ is Riemann integrable. Given $\varepsilon>0$, there is a partition $\mf{P}$ of $\mf{I}$ such that
\[U(f,\mf{P})-L(f,\mf{P})<\varepsilon.\]
 
Let $\mf{P}^*$ be a common refinement of $\mf{P}$ and $\mf{P}_0$. Then
\[U(f,\mf{P}^*)-L(f,\mf{P}^*)\leq U(f,\mf{P})-L(f,\mf{P})<\varepsilon.\]

\bp
But $\mf{P}^*$ induces a partition $\mf{P}^*(\mf{I}_1)$ and $\mf{P}^*(\mf{I}_2)$ of $\mf{I}_1$ and $\mf{I}_2$, and we have
\[U(f,\mf{P}^*)=U(f,\mf{P}^*(\mf{I}_1))+U(f,\mf{P}^*(\mf{I}_2)),\]\[L(f,\mf{P}^*)=L(f,\mf{P}^*(\mf{I}_1))+L(f,\mf{P}^*(\mf{I}_2)).\]
Therefore,
\[U(f,\mf{P}^*(\mf{I}_1))-L(f,\mf{P}^*(\mf{I}_1))+U(f,\mf{P}^*(\mf{I}_2))-L(f,\mf{P}^*(\mf{I}_2))<\varepsilon.\]
 
This implies that
\[U(f,\mf{P}^*(\mf{I}_j))-L(f,\mf{P}^*(\mf{I}_j))<\varepsilon\hspace{1cm}\text{for}\;j=1,2.\]
Hence, $f:\mf{I}_1\to\mb{R}$ and $f:\mf{I}_2\to\mb{R}$ are Riemann integrable.

Conversely, assume that $f:\mf{I}_1\to\mb{R}$ and $f:\mf{I}_2\to\mb{R}$ are Riemann integrable. Let $\{\mf{P}_{1,k}\}$ and $\{\mf{P}_{2,k}\}$ be Archimedes sequences of partitions for  $f:\mf{I}_1\to\mb{R}$ and $f:\mf{I}_2\to\mb{R}$ respectively. Then
\[\int_{\mf{I}_j}f=\lim_{k\to\infty} U(f, \mf{P}_{j,k})=\lim_{k\to\infty} L(f, \mf{P}_{j,k}) \hspace{1cm} \text{for}\;j=1,2.\]
For $k\in\mb{Z}^+$, let $\mf{P}^*_k$ be the partition  of $\mf{I}$ obtained by taking unions of partition points in $\mf{P}_{1, k}$ and $\mf{P}_{2,k}$. Then $\mf{P}_{1,k}=\mf{P}_k^*(\mf{I}_1)$ and $\mf{P}_{2,k}=\mf{P}_k^*(\mf{I}_2)$. It follows   that
\[U(f,\mf{P}_k^*) = U(f,\mf{P}_{1,k})+ U(f,\mf{P}_{2,k}), \quad  L(f,\mf{P}_k^*) = L(f,\mf{P}_{1,k})+ L(f,\mf{P}_{2,k}).\]
Therefore,
\[\lim_{k\to\infty}\left(U(f,\mf{P}_k^*) -L(f,\mf{P}_k^*) \right)=0.\]
Hence, $\{\mf{P}_k^*\}$ is an Archimedes sequence of partitions for $f:\mf{I}\to\mb{R}$. This shows that $f:\mf{I}\to\mb{R}$ is Riemann integrable, and
\[\int_{\mf{I}}f=\lim_{k\to\infty}U(f,\mf{P}_k^*) = \lim_{k\to\infty} \left(U(f,\mf{P}_{1,k})+ U(f,\mf{P}_{2,k})\right)=\int_{\mf{I}_1}f+\int_{\mf{I}_2}f.\]

\end{myproof}

 \begin{figure}[ht]
\centering
\includegraphics[scale=0.2]{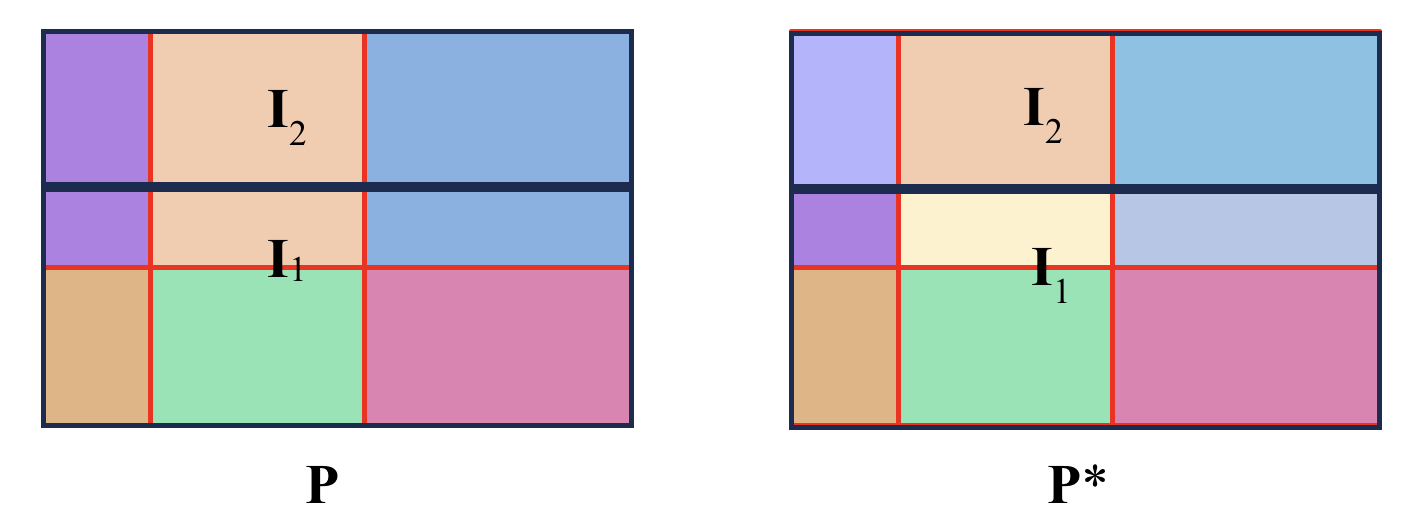}
\caption{A partition $\mf{P}$ of $\mf{I}$ and the refined partition $\mf{P}^*$ that induces partitions on $\mf{I}_1$ and $\mf{I}_2$.}\label{figure91}
\end{figure}
 
Next we state a lemma which is useful.
\begin{lemma}[label=230825_4]{}
Given that $\mf{I}=\di\prod_{i=1}^n[a_i,b_i]$ is a closed rectangle in $\mb{R}^n$,   let 

 \[\omega=\frac{1}{2}\min\{b_i-a_i\,|\,1\leq i\leq n\}.\]
\begin{enumerate}[(i)]
\item Given $\eta>0$, let  $\mf{I}_{\eta}$ be the closed rectangle
$\di\mf{I}_{\eta}=\prod_{i=1}^n[a_i-\eta, b_i+\eta]$.
For any $\varepsilon>0$, there exists $\delta>0$ such that for any $0<\eta <\delta$, 
$0<\text{vol}\,(\mf{I}_{\eta})-\text{vol}\,(\mf{I})<\varepsilon$.
\item Given $0<\kappa<\omega$, let  $  \check{\mf{I}}_{\kappa}$ be the closed rectangle
$\di   \check{\mf{I}}_{\kappa}=\prod_{i=1}^n[a_i+\kappa, b_i-\kappa]$.
For any $\varepsilon>0$, there exists $0<\delta\leq \omega$ such that for any $0<\kappa<\delta$, 
$0<\text{vol}\,(\mf{I})-\text{vol}\,(  \check{\mf{I}}_{\kappa})<\varepsilon$.
\end{enumerate}
\end{lemma}

 \begin{figure}[ht]
\centering
\includegraphics[scale=0.2]{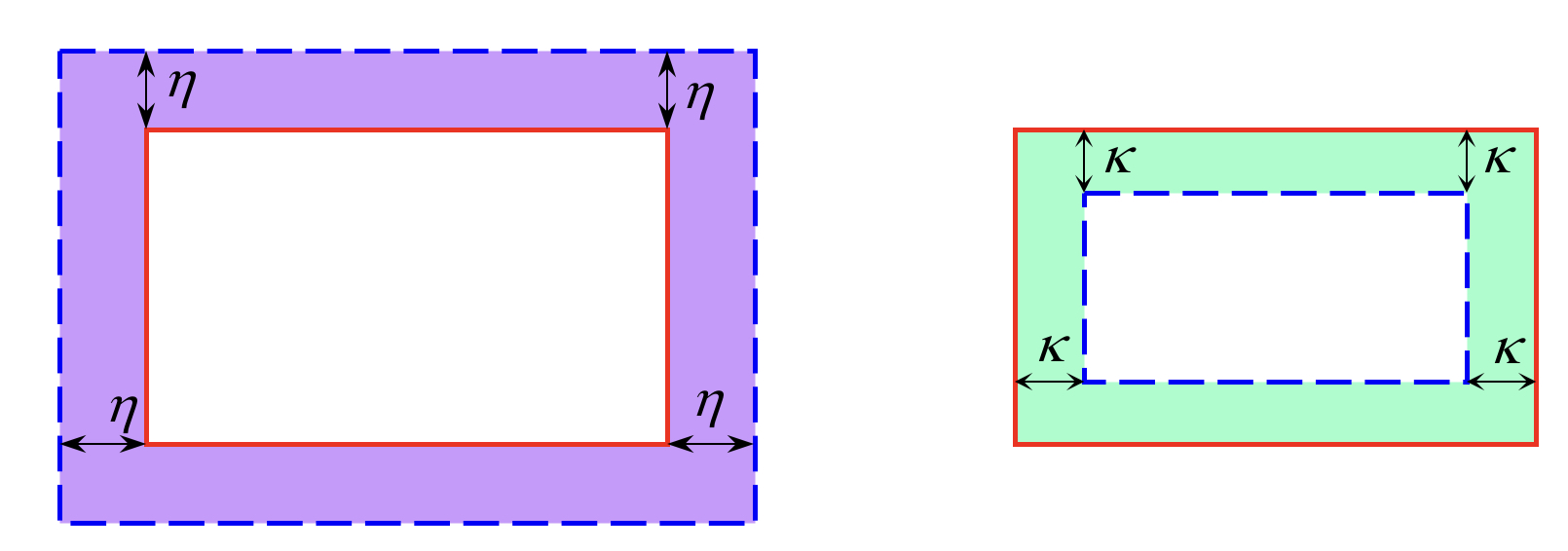}
\caption{Enlarging or shrinking a rectangle by an arbitrary amount.}\label{figure92}
\end{figure}
This lemma says that one can enlarge or shrink a rectangle by an arbitrarily small amount. It can be proved by elementary means. But here we  use some analysis technique to prove it.
\begin{myproof}{Proof}
We prove part (i). The argument for part (ii) is the same.
Consider the function $h:[0,\infty)\to\mb{R}$ defined by
\[h(\eta)=\text{vol}\,(\mf{I}_{\eta})=\prod_{i=1}^n(b_i-a_i+2\eta).\]
As a function of $\eta$, $h(\eta)$ is a polynomial, and it is a strictly increasing continuous function. The assertion is basically the definition of the limit $\di\lim_{\eta\to 0^+}h(\eta)=h(0)$.
\end{myproof}

The following theorem says that a bounded function $f:\mf{I}\to\mb{R}$ which is identically zero on the interior of $\mf{I}$ is Riemann integrable with integral $0$. This is something we would have expected.
\begin{theorem}[label=230826_1]{}
Let $\mf{I}=\di\prod_{i=1}^n[a_i,b_i]$, and let $f:\mf{I}\to\mb{R}$ be a bounded function such that
\[f(\mf{x})=0\hspace{1cm}\text{for all}\;\mf{x}\in \text{int}\,(\mf{I}).\]
Then $f:\mf{I}\to\mb{R}$ is Riemann integrable and
\[\int_{\mf{I}}f=0.\]
\end{theorem}
\begin{myproof}{Proof}
Let 
 $\di \omega=\frac{1}{2}\min\{b_i-a_i\,|\,1\leq i\leq n\}$.
Since $f:\mf{I}\to\mb{R}$ is bounded, there is a positive number $M$ such that
\[|f(\mf{x})|\leq M\hspace{1cm}\text{for all }\;\mf{x}\in \mf{I}.\]
\bp
By Lemma \ref{230825_4}, there is a $\kappa\in (0,\omega)$ such that 
\[\text{vol}\,(\mf{I})-\text{vol}\,(\mf{I}_{\kappa})<\frac{\varepsilon}{M},\]
where $\di \mf{I}_{\kappa}= \prod_{i=1}^n[a_i+\kappa, b_i-\kappa]$. It is a rectangle that is contained in $\text{int}\,(\mf{I})$.

Let $\mf{P}=(P_1, \ldots, P_n)$ be the partition of $\mf{I}$ with $P_i=\{a_i, a_i+\kappa, b_i-\kappa, b_i\}$. Then $\mf{I}_{\kappa}$ is one of the rectangles in the partition $\mf{P}$. On $\mf{J}=\mf{I}_{\kappa}$, $f(\mf{x})=0$, and so $M_{\mf{J}}=m_{\mf{J}}=0$. For all other rectangles $\mf{J}$ in $\mathcal{J}_{\mf{P}}$, we use the crude estimate
\[-M\leq m_{\mf{J}}\leq M_{\mf{J}}\leq M.\]
Then
\begin{align*}
L(f,\mf{P}) &=\sum_{\mf{J}\in\mathcal{J}_{\mf{P}}}m_{\mf{J}} \,\text{vol}\,(\mf{J})=\sum_{\mf{J}\in\mathcal{J}_{\mf{P}}\setminus\{\mf{I}_{\kappa}\}} m_{\mf{J}} \,\text{vol}\,(\mf{J})\\
&\geq  -M\left(\text{vol}\,(\mf{I})-\text{vol}\,(\mf{I}_{\kappa})\right)>-\varepsilon.
\end{align*}
In the same way, we find that
$\di U(f,\mf{P})<\varepsilon$. 
Since $\varepsilon>0$ is arbitrary, we find that
\[\underline{\int_{\mf{I}}}f \geq 0\quad\text{and}\quad \overline{\int_{\mf{I}}}f\leq 0.\]
Since $\di \underline{\int_{\mf{I}}}f \leq \overline{\int_{\mf{I}}}f$, we conclude that \[ \underline{\int_{\mf{I}}}f = \overline{\int_{\mf{I}}}f=0.\]
This proves that $f:\mf{I}\to\mb{R}$ is Riemann integrable and
\[ \int_{\mf{I}}f=0.\]
\end{myproof}

Now let us give a proof that the definition given in Section \ref{sec6.1}  for a bounded function $f:\mk{D}\to\mb{R}$ defined on a bounded subset of $\mb{R}^n$ to be Riemann integrable is unambiguous.
The crucial point is the following.
\begin{lemma}[label=230826_2]{}
Let $\mf{I}=\di\prod_{i=1}^n [a_i, b_i]$ and $  \check{\mf{I}}=\di\prod_{i=1}^n [\check{a}_i, \check{b}_i]$ be closed rectangles  in $\mb{R}^n$ such that $\mf{I} \subset\check{\mf{I}}$, and let $f:\mf{I}\to\mb{R}$ be a bounded function defined on $\mf{I}$. Then $f:\mf{I}\to\mb{R}$ is Riemann integrable if and only if  its zero extension  $  \check{f}:  \check{\mf{I}}\to\mb{R}$ is Riemann integrable. In such case, we also have
\[\int_{\mf{I}}f=\int_{  \check{\mf{I}}}  \check{f}.\]
\end{lemma}
\begin{myproof}{Proof}
 Let  $\check{\mf{P}}=\{\check{P}_1, \ldots, \check{P}_n\}$ be the partition  of $  \check{\mf{I}}$ such that the set $\check{P}_i$ is the set that contains  
$\check{a}_i, a_i, b_i, \check{b}_i$.  

For each rectangle  $\mf{J}$ in $\mathcal{J}_{\check{\mf{P}}}\setminus\{\mf{I}\}$, it is disjoint from the interior of $\mf{I}$. Hence,  $\check{f}$ vanishes in the interior of $\mf{J}$. By Theorem \ref{230826_1},  $\check{f}:\mf{J}\to\mb{R}$ is Riemann integrable and $\di\int_{\mf{J}}\check{f}=0$. 
   It follows from the additivity theorem that $\check{f}:\check{\mf{I}}\to\mb{R}$ is Riemann integrable if and only if $\check{f}: \mf{I} \to\mb{R}$
 is Riemann integrable, and
 \[\int_{  \check{\mf{I}}}  \check{f}=\int_{\mf{I}}\check{f}.\]
 However, restricted to $\mf{I}$, $\check{f}(\mf{x})=f(\mf{x})$. Hence, $\check{f}:\check{\mf{I}}\to\mb{R}$ is Riemann integrable if and only if $f: \mf{I} \to\mb{R}$
 is Riemann integrable. In such case, we have
 \[\int_{  \check{\mf{I}}}  \check{f}=\int_{\mf{I}}f.\]
\end{myproof}

 \begin{figure}[ht]
\centering
\includegraphics[scale=0.2]{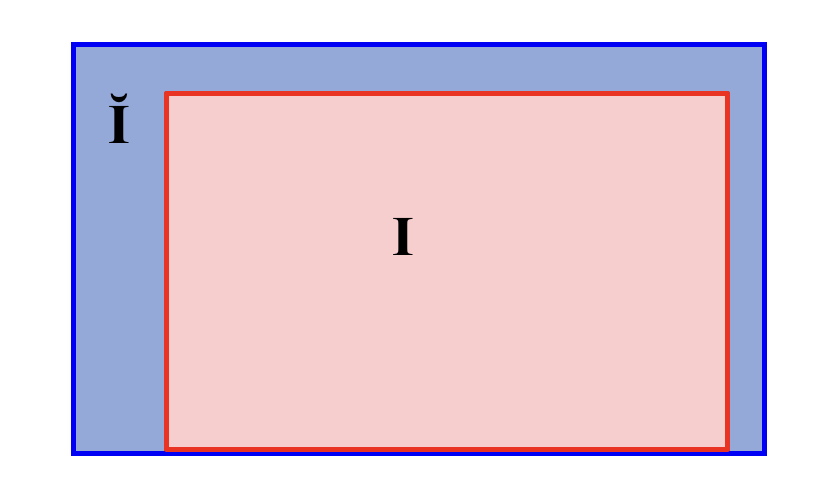}
\caption{The rectangle $\mf{I}$ is contained in the rectangle $\check{\mf{I}}$.}\label{figure93}
\end{figure}
Finally we can prove the main result.
\begin{theorem}{}
Let $\mk{D}$ be a bounded set in $\mb{R}^n$, and let $f:\mk{D}\to\mb{R}$ be  a bounded function defined on $\mk{D}$. The definition for Riemann integrability of $f:\mk{D}\to\mb{R}$ is unambiguous. Namely, if  $\di\mf{I}_1=\di\prod_{i=1}^n[a_i',b_i']$ and    $\di\mf{I}_2= \di\prod_{i=1}^n[a_i'',b_i'']$  contain $\mk{D}$, the zero extension $\check{f}:\mf{I}_1\to\mb{R}$ is Riemann integrable if and only if the zero extension  $\check{f}:\mf{I}_2\to\mb{R}$ is Riemann integrable. In the latter case, 

\vspace{-0.2cm}
\[\int_{\mf{I}_1}\check{f}=\int_{\mf{I}_2}\check{f},\] and so we can define unambiguously 

\vspace{-0.2cm}
\[\int_{\mk{D}}f=\int_{\mf{I}}\check{f},\] where $\mf{I}$  is any rectangle of the form $ \di\prod_{i=1}^n[a_i,b_i]$ that contains $\mk{D}$.
\end{theorem}

\begin{myproof}{Proof}
Let $\mf{I}=\mf{I}_1\cap \mf{I}_2$. Then $\mf{I}$ is a rectangle that is contained in $\mf{I}_1$ and $\mf{I}_2$. Lemma \ref{230826_2} then says that $\check{f}:\mf{I}_1\to\mb{R}$ is Riemann integrable if and only if  $\check{f}:\mf{I}\to\mb{R}$ is Riemann integrable, if and only if $\check{f}:\mf{I}_2\to\mb{R}$ is Riemann integrable. In latter case,
\[\int_{\mf{I}_1}\check{f}=\int_{\mf{I}}\check{f}=\int_{\mf{I}_2}\check{f}.\]
\end{myproof}
\begin{figure}[ht]
\centering
\includegraphics[scale=0.2]{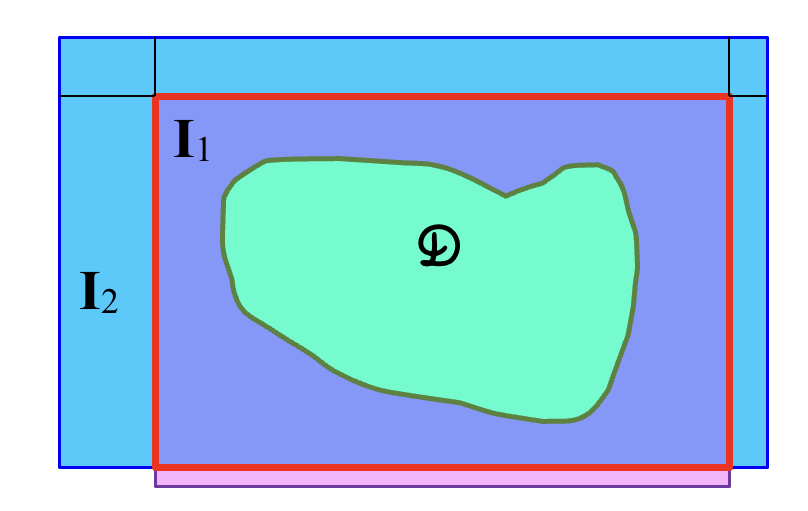}
\caption{The set $\mk{D}$ is contained in the rectangles $ \mf{I}_1$ and $\mf{I}_2$.}\label{figure94}
\end{figure}

Now we can extend the linearity and monotonicity to Riemann integrals over any bounded domains.

\begin{theorem}{Linearity}
Let $\mk{D}$ be bounded subset of $\mb{R}^n$, and let $f:\mk{D}\to\mb{R}$ and $g:\mk{D}\to\mb{R}$ be Riemann integrable functions. For any real numbers $\alpha$ and $\beta$, $(\alpha f+\beta g):\mk{D}\to\mb{R}$ is also Riemann integrable, and 
\[\int_{\mk{D}}(\alpha f+\beta g)=\alpha\int_{\mk{D}}f+\beta\int_{\mk{D}}g.\]
\end{theorem}
\begin{myproof}{Proof}
Let $\mf{I}=\di\prod_{i=1}^n[a_i,b_i]$ be a closed rectangle that contains $\mk{D}$, and let $\check{f}:\mf{I}\to\mb{R}$ and  $\check{g}:\mf{I}\to\mb{R}$ be the zero extensions of $f:\mk{D}\to\mb{R}$  and $g:\mk{D}\to\mb{R}$ to $\mf{I}$. It is easy to check that   $(\alpha \check{f}+\beta \check{g}):\mf{I}\to\mb{R}$ is the zero extension of $(\alpha f+\beta g):\mk{D}\to\mb{R}$ to $\mf{I}$. 
Since $f:\mk{D}\to\mb{R}$ and $g:\mk{D}\to\mb{R}$ are Riemann integrable,  $\check{f}:\mf{I}\to\mb{R}$ and  $\check{g}:\mf{I}\to\mb{R}$  are Riemann integrable and
\[\int_{\mf{I}}\check{f}=\int_{\mk{D}}f,\hspace{1cm}\int_{\mf{I}}\check{g}=\int_{\mk{D}}g.\]
 By Theorem \ref{230826_3},   $(\alpha \check{f}+\beta \check{g}):\mf{I}\to\mb{R}$ is Riemann integrable, and 
\[  \int_{\mf{I}}(\alpha \check{f}+\beta \check{g})=\alpha\int_{\mf{I}}\check{f}+\beta\int_{\mf{I}}\check{g}=\alpha\int_{\mk{D}}f+\beta\int_{\mk{D}}g.\]
It follows that  $(\alpha f+\beta g):\mk{D}\to\mb{R}$ is also Riemann integrable, and 
\[\int_{\mk{D}}(\alpha f+\beta g)= \int_{\mf{I}}(\alpha \check{f}+\beta \check{g})=\alpha\int_{\mk{D}}f+\beta\int_{\mk{D}}g.\]
\end{myproof}

\begin{theorem}[label=230826_4]{}
Let $\mk{D}$ be a bounded subset of $\mb{R}^n$, and let $f:\mk{D}\to\mb{R}$ be a bounded function defined on $\mk{D}$. Assume that
 $f(\mf{x})\geq 0$ for all $\mf{x}$ in $\mk{D}$. If $f:\mk{D}\to\mb{R}$ is Riemann integrable, then
 \[\int_{\mk{D}}f \;\geq \;0.\]
 
\end{theorem}
\begin{myproof}{Proof}
Let $\mf{I}=\di\prod_{i=1}^n[a_i,b_i]$ be a closed rectangle that contains $\mk{D}$, and let $\check{f}:\mf{I}\to\mb{R}$ be the zero extension of $f:\mk{D}\to\mb{R}$ to $\mf{I}$.  Since  $f:\mk{D}\to\mb{R}$ is Riemann integrable, $\check{f}:\mf{I}\to\mb{R}$ is also Riemann integrable. It is easy to check that $\check{f}(\mf{x})\geq 0$ for all $\mf{x}$ in $\mf{I}$.  Therefore,
\[\int_{\mk{D}}f=\int_{\mf{I}}\check{f}\geq 0.\]
\end{myproof}

As before, monotonicity is a consequence of linearity and Theorem \ref{230826_4}.
\begin{theorem}[label=230825_5]{Monotonicity}
Let $\mk{D}$ be a bounded subset of  $\mb{R}^n$, and let $f:\mk{D}\to\mb{R}$ and $g:\mk{D}\to\mb{R}$ be Riemann integrable functions. If
 $f(\mf{x})\geq g(\mf{x})$ for all $\mf{x}$ in $\mk{D}$, then
 \[\int_{\mk{D}}f \;\geq \;\int_{\mk{D}}g.\]
 
\end{theorem}

At the end of this section, we want to present two theorems whose proofs are almost verbatim those for the $n=1$ case.  The first theorem says that if a function is Riemann integrable, so is its absolute value.
\begin{theorem}[label=230829_1]{Absolute Value of Riemann Integrable Functions}
Let $\mk{D}$ be a bounded subset of  $\mb{R}^n$, and let $f:\mk{D}\to\mb{R}$ be a bounded function defined on $\mk{D}$. If the function $f:\mk{D}\to\mb{R}$ is Riemann integrable, then  the function $|f|:\mk{D}\to\mb{R}$ is also Riemann integrable.
\end{theorem}
\begin{myproof}{Sketch of Proof}
If $\check{f}:\mf{I}\to\mb{R}$ is the zero extension of $f:\mk{D}\to\mb{R}$ to the closed rectangle $\mf{I}=\di\prod_{i=1}^n[a_i,b_i]$ that contains $\mk{D}$, then $|\check{f}|:\mf{I}\to\mb{R}$ is the zero extension of $|f|:\mk{D}\to\mb{R}$. Hence,
it is sufficient to consider the case where $\mk{D}$ is a closed rectangle of the form $\mf{I}=\di\prod_{i=1}^n[a_i,b_i]$. The proof is almost the same as the $n=1$ case.
The key of the proof is the fact that for any subset $A$ of $\mf{I}$, 
\[\sup_{\mf{x}\in A}|f(\mf{x})|-\inf_{\mf{x}\in A}|f(\mf{x})|
\leq \sup_{\mf{x}\in A}f(\mf{x})-\inf_{\mf{x}\in A}f(\mf{x}).\]
\end{myproof}

The second theorem says  that products of Riemann integrable functions are Riemann integrable.
\begin{theorem}[label=230829_2]{Products of Riemann Integrable Functions}
Let $\mk{D}$ be a bounded subset of  $\mb{R}^n$, and let $f:\mk{D}\to\mb{R}$ and $g:\mk{D}\to\mb{R}$ be bounded functions defined on $\mk{D}$. If the functions $f:\mk{D}\to\mb{R}$ and $g:\mk{D}\to\mb{R}$ are Riemann integrable, then the function $(fg):\mk{D}\to\mb{R}$ is also Riemann integrable.
 
\end{theorem}
 \begin{myproof}{Sketch of Proof}
 
It is sufficient to consider the case where $\mk{D}$ is a closed rectangle of the form $\mf{I}=\di\prod_{i=1}^n[a_i,b_i]$.  The proof is almost the same as the $n=1$ case. 
The key of the proof is the fact that if $M$ is positive number such that

\vspace{-0.4cm}
\[|f(\mf{x})|\leq M\quad\text{and}\quad |g(\mf{x})|\leq M \hspace{1cm}\text{for all}\; \mf{x}\in\mf{I},\] then
for any subset $A$ of $\mf{I}$, 

\vspace{-0.8cm}
\begin{align*}
&\sup_{\mf{x}\in A}(fg)(\mf{x}) -\inf_{\mf{x}\in A}(fg)(\mf{x})
\\&\leq M\left(  \sup_{\mf{x}\in A}f(\mf{x})-\inf_{\mf{x}\in A}f(\mf{x})+ \sup_{\mf{x}\in A}g(\mf{x})-\inf_{\mf{x}\in A}g(\mf{x})\right).\end{align*}
\end{myproof}
\vp
\noindent
{\bf \large Exercises  \thesection}
\setcounter{myquestion}{1}

\begin{question}{\themyquestion}
Let $\mf{I}=[0,3]\times [0,3]$, and let $f:\mf{I}\to\mb{R}$ and $g:\mf{I}\to\mb{R}$  be Riemann integrable functions. Suppose that 
\[f(x,y)=g(y,x)\quad\text{and}\quad (3f+2g)(x,y)=10\hspace{1cm}\text{for all}\;(x,y)\in \mf{I},\]find $\di \int_{\mf{I}} f$ and $\di\int_{\mf{I}}g$.
\end{question}

\atc

\begin{question}{\themyquestion}
 Complete the details in the proof of Theorem \ref{230829_1}.
\end{question}
 \atc

\begin{question}{\themyquestion}
  Complete the details in the proof of Theorem \ref{230829_2}.
\end{question}
 
\section{Jordan Measurable Sets and Riemann Integrable Functions} \label{sec6.3}

In this section, we will give some sufficient conditions for a bounded function $f:\mk{D}\to\mb{R}$ to be Riemann integrable.
We start with the following theorem.
\begin{theorem}[label=230827_7]{}
Let $\mf{I}=\di\prod_{i=1}^n[a_i,b_i]$, and let $f:\mf{I}\to\mb{R}$ be a continuous function defined on $\mf{I}$. Then $f:\mf{I}\to\mb{R}$ is Riemann integrable.
\end{theorem}
\begin{myproof}{Proof}
  Since $f:\mf{I}\to\mb{R}$ is continuous and $\mf{I}$ is compact,   $f:\mf{I}\to\mb{R}$ is uniformly continuous.  Given $\varepsilon>0$, there exists $\delta>0$ such that if $\mf{u}$ and $\mf{v}$ are points in $\mf{I}$ and $\Vert\mf{u}-\mf{v}\Vert<\delta$, then 

\vspace{-0.5cm}
\[|f(\mf{u})-f(\mf{v})|<\frac{\varepsilon}{\text{vol}\,(\mf{I})}.\]
Let $\mf{P}$ be any partition of $\mf{I}$ with $|\mf{P}|<\delta$. A rectangle $\mf{J}$ in $\mathcal{J}_{\mf{P}}$  is a compact set. Since $f:\mf{J}\to\mb{R}$ is continuous, the extreme value theorem says that there exist points $\mf{u}_{\mf{J}}$ and $\mf{v}_{\mf{J}}$ in $\mf{J}$ such that
\[f(\mf{u}_{\mf{J}})\leq f(\mf{x})\leq f(\mf{v}_{\mf{J}})\hspace{1cm}\text{for all}\;\mf{x}\in\mf{J}.\]
Therefore, 

\vspace{-0.5cm}
\[m_{\mf{J}}=\inf_{\mf{x}\in\mf{J}}f(\mf{x})=f(\mf{u}_{\mf{J}})\quad\text{and}\quad  M_{\mf{J}}=\sup_{\mf{x}\in\mf{J}}f(\mf{x})=f(\mf{v}_{\mf{J}}).\]Since $|\mf{P}|<\delta$, 
\[\Vert\mf{u}_{\mf{J}}-\mf{v}_{\mf{J}}\Vert\leq\text{diam}\,\mf{J}\leq |\mf{P}|<\delta.\]
Therefore,
\[M_{\mf{J}}-m_{\mf{J}}=f(\mf{v}_{\mf{J}})-f(\mf{u}_{\mf{J}})<\frac{\varepsilon}{\text{vol}\,(\mf{I})}.\]
This implies that

\vspace{-0.4cm}
\[
U(f,\mf{P})-L(f,\mf{P}) =\sum_{\mf{J}\in\mathcal{J}_{\mf{P}}}\left(M_{\mf{J}}-m_{\mf{J}}\right)\,\text{vol}\,(\mf{J}) <\frac{\varepsilon}{\text{vol}\,(\mf{I})}\sum_{\mf{J}\in\mathcal{J}_{\mf{P}}}\text{vol}\,(\mf{J})=\varepsilon.\]
Hence, $f:\mf{I}\to\mb{R}$ is Riemann integrable.
\end{myproof}

\begin{example}{}
Let $f:[0,1]\times [0,1]\to\mb{R}$ be the function defined as 
\[f(x,y)=\sin(xy).\]
This is a composition of the sine function and a polynomial, both of which are continuous functions. Hence, $f:[0,1]\times [0,1]\to\mb{R}$ is  a continuous function. Therefore, $f:[0,1]\times [0,1]\to\mb{R}$ is Riemann integrable.
\end{example}

In Section \ref{sec6.1}, we have seen that  a constant function $f:\mf{I}\to\mb{R}$, $f(\mf{x})=c$ defined on   $\mf{I}=\di\prod_{i=1}^n[a_i,b_i]$ is Riemann integrable and  its integral is \[\int_{\mf{I}}f=c\,\text{vol}\,(\mf{I}).\]
Since constant functions are the simplest bounded functions, it is natural to ask whether a constant function  $f:\mk{D}\to\mb{R}$, $f(\mf{x})=c$ on a bounded set $\mk{D}$ is always Riemann integrable. By linearity, it is sufficient to consider the case when $c=1$. When $\mk{D} $ is a closed rectangle of the form $\mf{I}=\di\prod_{i=1}^n[a_i,b_i]$, the answer is affirmative and we have 
\[\int_{\mf{I}}\,d\mf{x}=\text{vol}\,(\mf{I})\]

To consider a general set $\mk{D}$, let us first define the characteristic function of a set.
\begin{definition}{Characteristic Functions}
Let $A$ be a subset of $\mb{R}^n$. The characteristic function or indicator function of the set $A$  is the function $\chi_A:\mb{R}^n\to\mb{R}$ defined as
\[\chi_A(\mf{x})=\begin{cases} 1,\quad &\text{if}\;\mf{x}\in A,\\0,\quad &\text{if}\;\mf{x}\notin A.\end{cases}\]
\end{definition}

\begin{example}{}
Let 
$\di A=\di\left\{(x,y)\,|\, x>0\right\}$. Notice that the function $\chi_A:\mb{R}^2\to\mb{R}$ is
\[\chi_A(x,y)=\begin{cases} 1,\quad &\text{if}\;x>0,\\0,\quad &\text{if}\;x\leq 0.\end{cases}\] It is continuous at $(x,y)$ if and only if $x>0$ or $x<0$. 
\end{example}

 \begin{figure}[ht]
\centering
\includegraphics[scale=0.2]{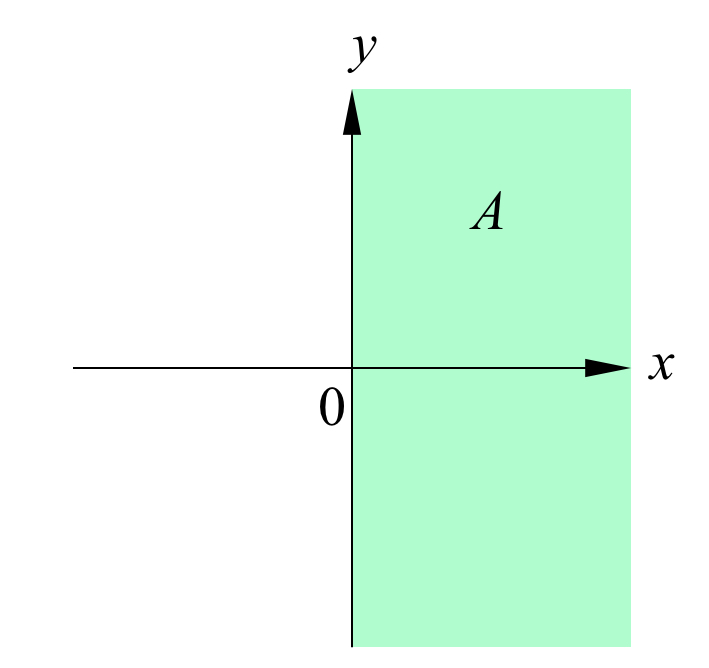}
\caption{The  set $A=\di\left\{(x,y)\,|\, x>0\right\}$.}\label{figure97}
\end{figure}

\begin{highlight}{Interior, Exterior and Boundary of a Set}
In Chapter \ref{chapter1}, we have seen that if $A$ is a subset of $\mb{R}^n$,  then $\mb{R}^n$ is a disjoint union of $\text{int}\,A$, $\text{ext}\, A$ and $\pa A$.

If $\mf{x}_0$ is a point  in $\mb{R}^n$, $\mf{x}_0\in \text{int}\,A$ if and only if there is an $r>0$ such that $B(\mf{x}_0, r)\subset A$;  $\mf{x}_0\in \text{ext}\,A$ if and only if there is an $r>0$ such that $B(\mf{x}_0, r)\subset\mb{R}^n\setminus A$; and $\mf{x}_0\in \pa A$ if for every $r>0$, $B(\mf{x}_0,r)$ contains a point in $A$ and a point not in $A$.

\end{highlight}

\begin{theorem}[label=230827_10]
{}
Let $A$ be a subset of $\mb{R}^n$, and let $\chi_A:\mb{R}^n\to\mb{R}$ be the characteristic function   of   $A$. Then the set of discontinuities of the function $\chi_A$ is the set $\pa A$.
\end{theorem}
\begin{myproof}{Proof}
Since $\mb{R}^n$ is a disjoint union of $\text{int}\,A$, $\text{ext}\, A$ and $\pa A$, we will show that $\chi_A$ is continuous on $\text{int}\,A$ and $\text{ext}\,A$, and discontinuous at every point in $\pa A$.

The sets $\text{int}\,A$ and $\text{ext}\, A$ are open sets, and $f$ is equal to 1 on $\text{int}\,A$ and 0 on $\text{ext}\,A$. For every $\mf{x}_0$  in $ \text{int}\,A$, there is an $r>0$ such that $B(\mf{x}_0, r)\subset A$. Therefore, for any $\varepsilon>0$, if $\mf{x}$ is such that $\Vert\mf{x}-\mf{x}_0\Vert<r$, then \[|f(\mf{x})-f(\mf{x}_0)|=0<\varepsilon.\] This shows that $f$ is continuous at $\mf{x}_0$. Similarly, if $\mf{x}_0$ is in $\text{ext}\, A$, there is an $r>0$ such that $B(\mf{x}_0, r)\subset \mb{R}^n\setminus A$. The same reasoning shows that $f$ is continuous at $\mf{x}_0$.

Now consider a point $\mf{x}_0$ that is in $\pa A$. For any $k\in\mb{Z}^+$, there is a point $\mf{u}_k\in A$ and a point $\mf{v}_k\notin A$ such that $\mf{u}_k$ and $\mf{v}_k$ are in the neighbourhood $B(\mf{x}_0, 1/k)$ of $\mf{x}_0$. The two sequences $\{\mf{u}_k\}$ and $\{\mf{v}_k\}$ both converge to $\mf{x}_0$, but the sequence $\{f(\mf{u}_k)\}$ converges to 1, the sequence $\{f(\mf{v}_k)\}$ converges to 0. This shows that $f$ is not continuous at $\mf{x}_0$.
\end{myproof}

 \begin{figure}[ht]
\centering
\includegraphics[scale=0.2]{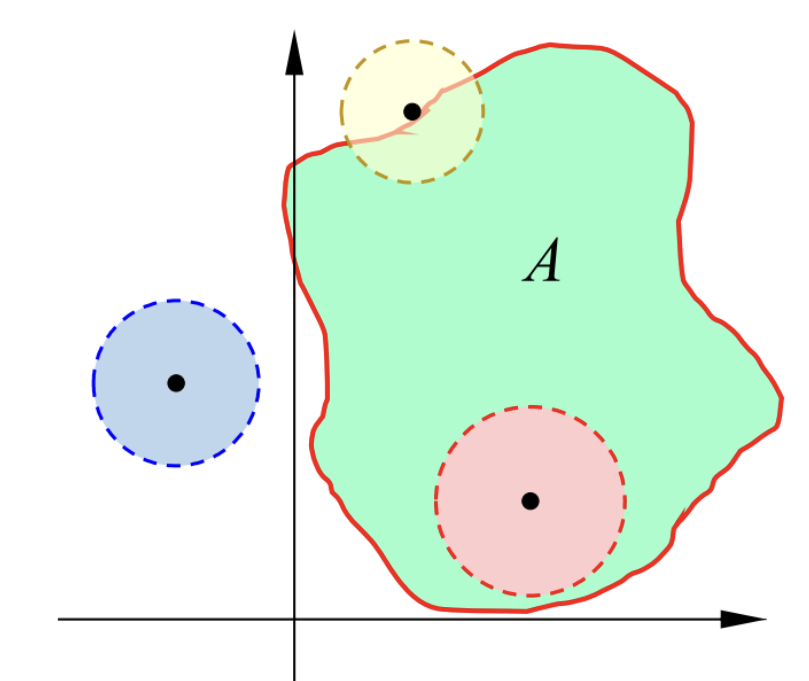}
\caption{The characteristic function of a set $A$ is not continuous at $\mf{x}_0$ if and only if $\mf{x}_0\in\pa A$.}\label{figure95}
\end{figure}

By definition, restricted to the set $A$, $\chi_A:A\to\mb{R}$ is the constant function $\chi_A(A)=1$.
Now we define Jordan measurable sets and its volume.
\begin{definition}{Jordan Measurable Sets and Volume}
Let $\mk{D}$ be a bounded subset of $\mb{R}^n$. We say that $\mk{D}$ is Jordan measurable if the constant function $\chi_D:\mk{D}\to\mb{R}$ is Riemann integrable. In this case, we define the volume of $\mk{D}$ as
\[\text{vol}\,(\mk{D})=\int_{\mk{D}}\chi_{\mk{D}}=\int_{\mk{D}}\,d\mf{x}.\]
\end{definition}

\begin{example}{}
The closed rectangle $\mf{I}=\di\prod_{i=1}^n[a_i,b_i]$ is Jordan measurable, and its volume is
\[\text{vol}\,(\mf{I})=\int_{\mf{I}}\;d\mf{x}=\prod_{i=1}^n(b_i-a_i),\] as what we have  defined earlier.
\end{example}
\begin{example}{}
Example \ref{230823_6} says that the set
\[\mk{D}=\left\{(x,y)\,|\, 0\leq y\leq x\leq 1\right\}\] is Jordan measurable and $\di\text{vol}\,(\mk{D})=\frac{1}{2}$.
\end{example}

 \begin{figure}[ht]
\centering
\includegraphics[scale=0.2]{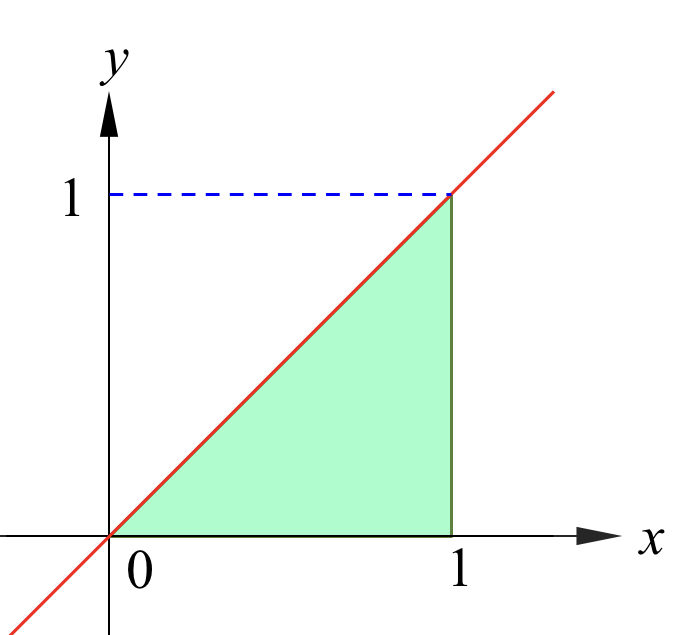}
\caption{The  set $\mk{D}=\left\{(x,y)\,|\, 0\leq y\leq x\leq 1\right\} $ is Jordan measurable.}\label{figure98}
\end{figure}

One might think that all bounded subsets of $\mb{R}^n$ has volumes. This is not true. An example is given below.
\begin{example}{}
Let $\mf{I}=[0,1]^n$ and let
\[\mk{D}=\left\{\mf{x}\in \mf{I}\,|\,\mf{x}\in \mb{Q}^n\right\}.\]
Notice that $\mk{D}$ is a subset of the rectangle $\mf{I}$, and the zero extension of $\chi_{\mk{D}}:\mk{D}\to\mb{R}$ is the function $\chi_{\mk{D}}:\mf{I}\to\mb{R}$,
\[\chi_{\mk{D}}(\mf{x})=\begin{cases}1,\quad &\text{if all components of $\mf{x}$ are rational},\\0,\quad &\text{otherwise},\end{cases}\]
which is the Dirichlet's function. We have seen in Example \ref{230827_1} that the function $\chi_{\mk{D}}:\mf{I}\to\mb{R}$ is not Riemann integrable. Hence, $\chi_{\mk{D}}:\mk{D}\to\mb{R}$  is not Riemann integrable. This means the set $\mk{D}$ is not Jordan measurable and so it does not have a volume.
\end{example}
This example also shows that if  $B$ is  a subset of $A$, and the function $f:A\to\mb{R}$ is Riemann integrable, the function $f:B\to\mb{R}$ is not necessary Riemann integrable.

The next example says that the boundary of a rectangle has volume 0.
\begin{example}[label=230827_6]{}
Let $\mf{I}=\di\prod_{i=1}^n[a_i,b_i]$, and let $\mk{D}=\pa\, \mf{I}$. Notice that $\mk{D}$ is contained in $\mf{I}$.  The zero extension of $\chi_{\mk{D}}:\mk{D}\to\mb{R}$ is the function $\chi_{\mk{D}}:\mf{I}\to\mb{R}$ which vanishes on the interior of $\mf{I}$. By Theorem \ref{230826_1},  $\chi_{\mk{D}}:\mf{I}\to\mb{R}$ is Riemann integrable and $\di \int_{\mf{I}}\chi_{\mk{D}}=0$. Therefore, $\mk{D}=\pa\,\mf{I}$ has zero volume.
\end{example}

\begin{remark}
{  Darboux Sums for a Characteristic Function}
Given a bounded set $\mk{D}$ that is contained in the rectangle $\mf{I}$, if $\mf{P}$ is a partition of $\mf{I}$, 
$L(\chi_{\mk{D}},\mf{P})$ is the sum of the volumes of the rectangles in $\mf{P}$ that is contained in $\mk{D}$; while $U(\chi_{\mk{D}},\mf{P})$ is the sum of the volumes of the rectangles in $\mf{P}$ that intersect $\mk{D}$. See Figure \ref{figure101}.

Thus, for $\mk{D}$ to have volume, the two numbers $L(\chi_{\mk{D}},\mf{P})$ and $U(\chi_{\mk{D}},\mf{P})$ should get closer and closer when the partitions $\mf{P}$ gets finer. 
\end{remark}

 \begin{figure}[ht]
\centering
\includegraphics[scale=0.2]{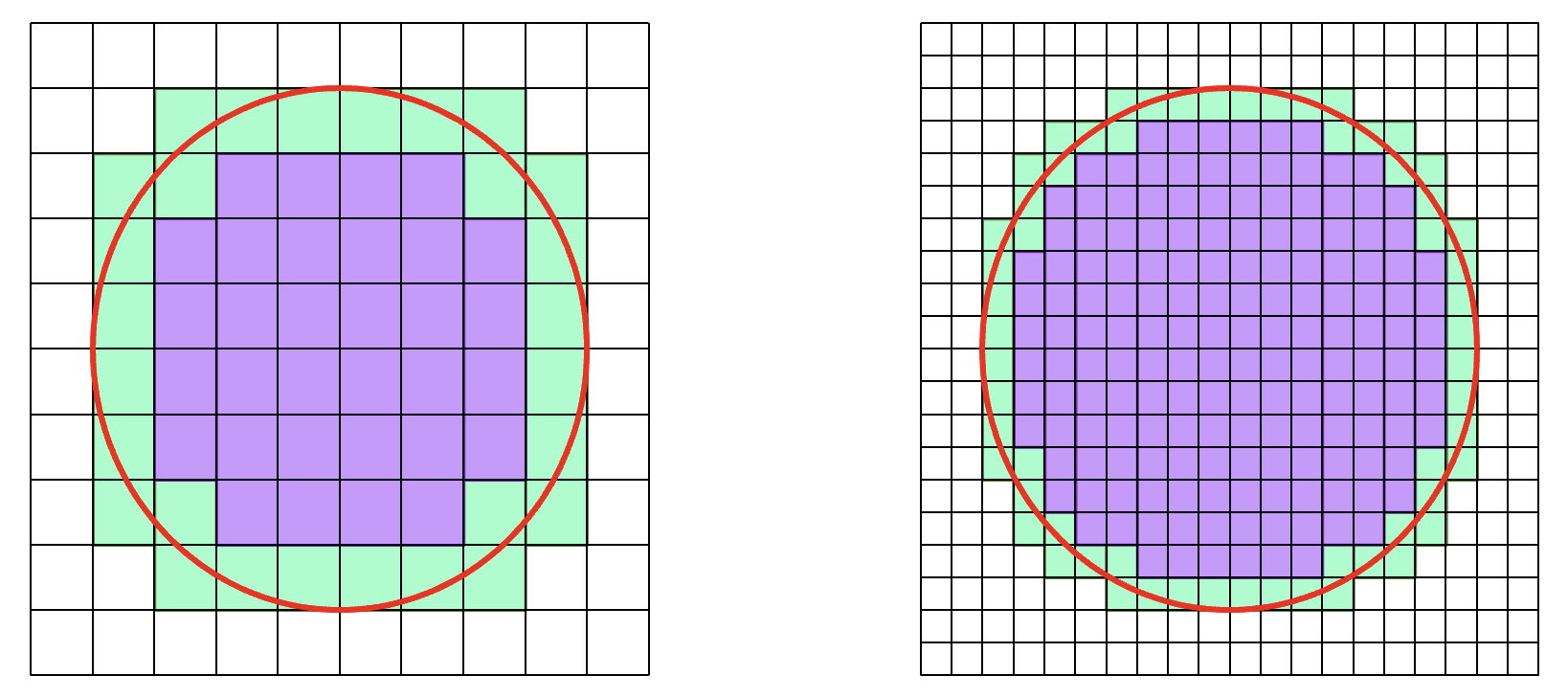}
\caption{The geometric quantities represented by $L(\chi_{\mk{D}},\mf{P})$ and $U(\chi_{\mk{D}},\mf{P})$ when $\mk{D}$ is the region bounded inside the circle.}\label{figure101}
\end{figure}
 
Our goal is to give characterization of sets that are Jordan measurable. 
We will consider those that have zero volumes first. 
The following is a  useful lemma.
\begin{lemma}[label=230827_2]{}
Let $\mf{I}$ be a closed rectangle in $\mb{R}^n$ that contains the closed rectangles $\mf{I}_1, \ldots, \mf{I}_k$. There is a partition $\mf{P}$ of $\mf{I}$ such that if $\mf{J}$ is a rectangle in the partition $\mf{P}$, then $\mf{J}$ is either contained in an $\mf{I}_j$ for some $1\leq j\leq k$, or $\mf{J}$ is disjoint from the interiors of $\mf{I}_j$ for all $1\leq j\leq k$.
\end{lemma}
\begin{myproof}{Sketch of Proof}
We construct the partition $\mf{P} =(P_1, \ldots, P_n)$ in the following way. For each $1\leq i\leq n$, the partition points in $P_i$ is the set of end points of the $i^{\text{th}}$-edge of $\mf{I}$, $\mf{I}_1$, $\ldots$, $\mf{I}_{k}$. One can check that this partition satisfies the requirement. See Figure \ref{figure96} for an illustration.
\end{myproof}

 \begin{figure}[ht]
\centering
\includegraphics[scale=0.2]{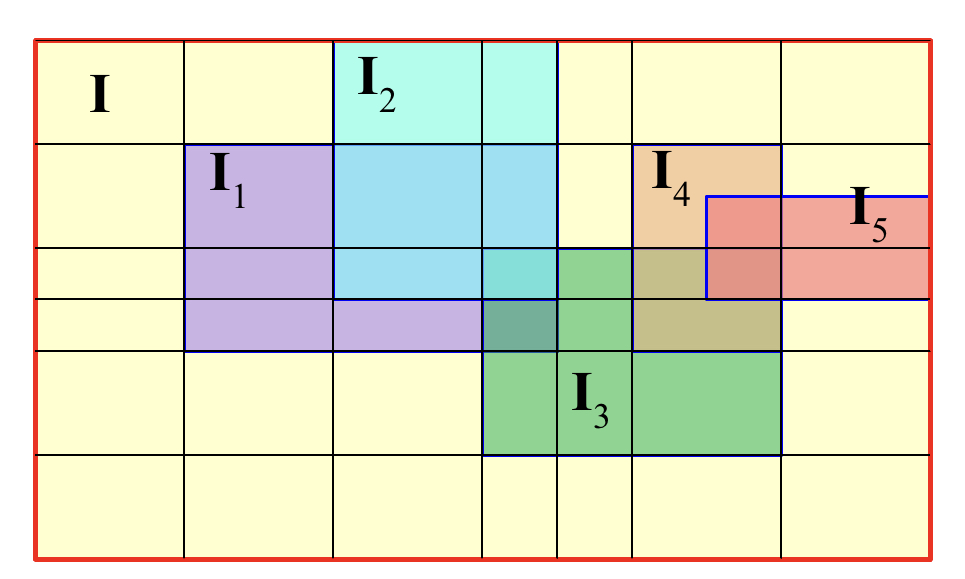}
\caption{A partition of the rectangle $\mf{I}$ that satisfies the conditions in Lemma \ref{230827_2}. }\label{figure96}
\end{figure}

Let us introduce the definition of a cube.
\begin{definition}{Cubes}
A rectangle of the form $\di\prod_{i=1}^n[a_i,b_i]$ such that 
\[b_1-a_1=b_2-a_2=\cdots=b_n-a_n=\ell=2r\] is called a (closed) cube with side length $\ell=2r$. The center of the cube is
\[\mf{c}=(c_1, c_2, \ldots, c_n)=\left(\frac{a_1+b_1}{2}, \frac{a_2+b_2}{2}, \ldots, \frac{a_n+b_n}{2}\right).\]
We will denote such a cube by $Q_{\mf{c},r}$. 
\end{definition}There are also cubes whose edges are not parallel to the coordinate axes. In this chapter, when we say a cube, we always mean a cube defined above.

Now we can give a characterization of sets with zero volume.
\begin{theorem}[label=230827_5]{}
Let $\mk{D}$ be a bounded subset of $\mb{R}^n$. The following are equivalent.
\begin{enumerate}[(a)]
\item The set $\mk{D}$ is Jordan measurable and it has zero volume.
\item For any $\varepsilon>0$, there are finitely many closed cubes $Q_1, \ldots, Q_k$  such that
\[\mk{D}\subset \di \bigcup_{j=1}^kQ_j\quad\text{and}\quad \sum_{j=1}^k\text{vol}\,(Q_j)<\varepsilon.\]
\item For any $\varepsilon>0$, there are finitely many closed rectangles $\mf{I}_1, \ldots, \mf{I}_k$  such that
\[\mk{D}\subset \di \bigcup_{j=1}^k\mf{I}_j\quad\text{and}\quad \sum_{j=1}^k\text{vol}\,(\mf{I}_j)<\varepsilon.\]
\end{enumerate}
\end{theorem}

\begin{myproof}{Proof}

First assume that $\mk{D}$ is a Jordan measurable set with zero volume. There is a positive number $R$ such that the closed cube $Q_{\mf{0},R}=[-R,R]^n$ contains the set $\mk{D}$. Let $\mf{I}=Q_{\mf{0},R}$. Then the function $\chi_{\mk{D}}:\mf{I}\to\mb{R}$ is Riemann integrable and 
$\di \int_{\mf{I}}\chi_{\mk{D}}=0$.
Given $m\in\mb{Z}^+$, let $\mf{P}_m$ be the uniformly regular partition of $\mf{I}$ into $m^n$ rectangles. Notice that each rectangle in the  partition $\mf{P}_m$ is a cube. Since  $\di\lim_{m\to\infty}|\mf{P}_m|=0$, we have
\[\lim_{m\to\infty}U(\chi_{\mk{D}},\mf{P}_m)=\int_{\mf{I}}\chi_{\mk{D}}=0.\]
Given $\varepsilon>0$, there is a positive integer $M$ such that for all $m\geq M$, 
\[U(\chi_{\mk{D}},\mf{P}_m)<\varepsilon.\]
Consider the partition $\mf{P}_M$. 
Notice that for $\mf{J}\in\mathcal{J}_{\mf{P}_{M}}$, 
\[M_{\mf{J}}(\chi_{\mk{D}})=\begin{cases}1,\quad &\text{if}\;\mf{J}\cap\mk{D}\neq\emptyset,\\
0,\quad &\text{if}\;\mf{J}\cap\mk{D}=\emptyset.\end{cases}\]
 \bp
Let \[\mathscr{A}=\left\{\mf{J}\in\mathcal{J}_{\mf{P}_{M}}\,|\, \mf{J}\cap\mk{D}\neq\emptyset\right\}.\]

Then
\[U(\chi_{\mk{D}},\mf{P}_{M})=\sum_{\mf{J}\in \mathscr{A}}\text{vol}\,(\mf{J}).\]
$\mathscr{A}$ is a finite collection of cubess. Hence, we can named the cubes in $\mathscr{A}$ as $Q_1, \ldots, Q_k$. By construction, \[\mk{D}\subset \di \bigcup_{j=1}^kQ_j\quad\text{and }\quad
  \sum_{j=1}^k\text{vol}\,(Q_j)<\varepsilon.\]
This proves that (a) implies (b).

(b) implies (c) is obvious since a cube is a rectangle.

Now assume that (c) holds. Given $\varepsilon>0$,  (c) says that  there are  closed rectangles $\mf{I}_1, \ldots, \mf{I}_k$ such that
\[\mk{D}\subset \di \bigcup_{j=1}^k\mf{I}_j \quad \text{and} 
\quad \sum_{j=1}^k\text{vol}\,(\mf{I}_j)<\frac{\varepsilon}{2}.\]
By Lemma \ref{230825_4}, for each $1\leq j\leq k$, there is a closed rectangle $\check{\mf{I}}_j$ such that $\mf{I}_j\subset \text{int}\, \check{\mf{I}}_j$ and
\[\text{vol}\,(\check{\mf{I}}_j)-\text{vol}\,(\mf{I}_j)<\frac{\varepsilon}{2k}.\]
 
It follows that
\[\mk{D}\subset\bigcup_{j=1}^k\text{int}\, \check{\mf{I}}_j\quad\text{and}\quad \sum_{j=1}^k\text{vol}\,(\check{\mf{I}}_j)<\varepsilon.\]
Let $\mf{I}$ be a closed rectangle whose interior contains the bounded set $\di \bigcup_{j=1}^k \check{\mf{I}}_j$.

\bp
  By Lemma \ref{230827_2}, there is a partition $\mf{P}$ of $\mf{I}$ such that  each rectangle $\mf{J}$ in the partition $\mf{P}$  is either contained in an $\check{\mf{I}}_j$ for some $1\leq j\leq k$, or  is disjoint from the interiors of $ \check{\mf{I}}_j$ for all $1\leq j\leq k$.
Let
\[\mathscr{B}=\left\{\mf{J}\in\mathcal{J}_{\mf{P}}\,|\, \mf{J}\subset\check{\mf{I}}_j\;\text{for some}\;1\leq j\leq k\right\}.\]
If $\mf{J}\notin \mathscr{B}$, then $\mf{J}\cap \text{int}\,\check{\mf{I}}_j=\emptyset$ for all $1\leq j\leq k$. Therefore, $\mf{J}\cap\mk{D}=\emptyset$. For these $\mf{J}$, $M_{\mf{J}}(\chi_{\mk{D}})=m_{\mf{J}}(\chi_{\mk{D}})=0$.
 If $\mf{J}$ is in $\mathscr{B}$, we use the simple estimate $M_{\mf{J}}\leq 1$. Thus,
\[U(\chi_{\mk{D}}, \mf{P}) =\sum_{\mf{J}\in \mathscr{B}} M_{\mf{J}} \text{vol}\,(\mf{J})\leq
  \sum_{\mf{J}\in \mathscr{B}} \text{vol}\,(\mf{J})\leq \sum_{j=1}^k\text{vol}\,(\check{\mf{I}}_j)<\varepsilon.
\]
Since $L(\chi_{\mk{D}},\mf{P})\geq 0$, we   find that 
\[U(\chi_{\mk{D}}, \mf{P}) -L(\chi_{\mk{D}}, \mf{P}) <\varepsilon.\]
This proves that $\chi_{\mk{D}}:\mf{I}\to\mb{R}$ is Riemann integrable.  Since we have shown that there exists a partition $\mf{P}$ such that $
\di U(\chi_{\mk{D}},\mf{P})<\varepsilon$,
we  have
\[\text{vol}\,(\mk{D})=\int_{\mf{I}}\chi_{\mk{D}}\leq U(\chi_{\mk{D}},\mf{P}) <\varepsilon.\]Since $\varepsilon>0$ is arbitrary, we find that
$\text{vol}\,(\mk{D})=0$. This completes the proof of (c) implies (a).

\end{myproof}

Motivated by Theorem \ref{230827_5}, we make the following definition.
\begin{definition}{Jordan Content Zero}
Let $\mk{D}$ be a bounded subset of $\mb{R}^n$. We say that $\mk{D}$ has Jordan content zero provided that 
for any $\varepsilon>0$, there are finitely many closed rectangles $\mf{I}_1, \ldots, \mf{I}_k$ such that
\[\mk{D}\subset \di \bigcup_{j=1}^k\mf{I}_j\quad\text{and}\quad \sum_{j=1}^k\text{vol}\,(\mf{I}_j)<\varepsilon.\]
\end{definition}

\begin{highlight}{Sets that have Jordan Content Zero}
Let $\mk{D}$ be a bounded subset of $\mb{R}^n$. Theorem \ref{230827_5} says that $\mk{D}$ is Jordan measurable with volume zero if and only if it has Jordan content zero.
\end{highlight}

The characterization of sets with zero volume given in Theorem \ref{230827_5} facilitates the proofs of properties of such sets.
\begin{theorem}[label=230827_9]{}
 Let $\mk{D}_1$ and $\mk{D}_2$ be bounded subsets of $\mb{R}^n$. If $\mk{D}_1$ has Jordan content zero and $\mk{D}_2\subset \mk{D}_1$, then $\mk{D}_2$ also has Jordan content zero.
 
\end{theorem}
\begin{myproof}{Proof}
Given $\varepsilon>0$, since $\mk{D}_1$ has Jordan content zero,  there are closed rectangles $\mf{I}_1, \ldots, \mf{I}_k$ such that
\[\mk{D}_1\subset \di \bigcup_{j=1}^k\mf{I}_j\quad\text{and}\quad \sum_{j=1}^k\text{vol}\,(\mf{I}_j)<\varepsilon.\]
Since $\mk{D}_2\subset \mk{D}_1$, we find that
\[\mk{D}_2\subset \di \bigcup_{j=1}^k\mf{I}_j\quad\text{and}\quad \sum_{j=1}^k\text{vol}\,(\mf{I}_j)<\varepsilon.\]Therefore, $\mk{D}_2$ also has Jordan content zero.
\end{myproof}

\begin{example}{}
Let $\mk{D}$ be the subset of $\mb{R}^3$ given by
\[\mk{D}=\left\{(x,y,2)\,|\, -2\leq x\leq 3, -5\leq y\leq 7\right\}.\]
Show that $\mk{D}$ is a Jordan measurable set with zero volume.
\end{example}
\begin{solution}{Solution}
Let $\mf{I}=[-2, 3]\times [-5, 7]\times [2, 3]$. Then $\mf{I}$ is a closed rectangle in $\mb{R}^3$. Example \ref{230827_6} says that $\pa\, \mf{I}$ has Jordan content zero. Since $\mk{D}\subset\pa\,\mf{I}$, Theorem \ref{230827_9} says that $\mk{D}$ has Jordan content zero. Hence, $\mk{D}$ is a Jordan measurable set with zero volume.
\end{solution}

The next theorem concerns unions and intersections of sets of Jordan content zero.
\begin{theorem}[label=230827_13]{}
\begin{enumerate}[(a)]
\item
If $\mathscr{A}=\left\{ \mk{D}_{\alpha}\,|\,\alpha\in J\right\}$ is a collection of sets that have Jordan content zero, then their intersection $\mathcal{U}=\di\bigcap_{\alpha\in J}\mk{D}_{\alpha}$ also has Jordan content zero.
\item If $\mk{D}_1$, $\ldots$, $\mk{D}_m$ are finitely many sets that have   Jordan content zero, then their union $\mk{D}=\di\bigcup_{j=1}^m \mk{D}_j$ is also a set that has Jordan content zero.
\end{enumerate}
\end{theorem}

\begin{myproof}{Proof}
(a) is obvious since $\mathcal{U}\subset \mk{D}_{\alpha}$ for any $\alpha\in J$.

(b) is basically a consequence of  the fact that finite union of finite sets is finite. Given $\varepsilon>0$, for each $1\leq j\leq m$, since $\mk{D}_j$ has Jordan content zero, there is a finite collection $\mathscr{B}_j=\{\mf{I}_{\beta_j}\,|\,\beta_j\in J_j\}$ of closed rectangles such that
\[\mk{D}_j\subset \bigcup_{\beta_j\in J_j}\mf{I}_{\beta_j},\hspace{1cm}\sum_{\beta_j\in J_j}\text{vol}\,\left(\mf{I}_{\beta_j}\right)<\frac{\varepsilon}{m}.\]

Let
\[\mathscr{B}=\bigcup_{j=1}^m \mathscr{B}_j.\]

Since each $\mathscr{B}_j$, $1\leq j\leq m $ is finite, $\mathscr{B}$ is also a finite collection of closed rectangles. Moreover,
\[\mk{D}=\bigcup_{j=1}^m\mk{D}_j\subset \bigcup_{j=1}^m \bigcup_{\beta_j\in J_j}\mf{I}_{\beta_j}=\bigcup_{\mf{I}_{\beta}\in \mathscr{B}}\mf{I}_{\beta},\]
\bp
and
\[\sum_{\mf{I}_{\beta}\in \mathscr{B}}\text{vol}\,(\mf{I}_{\beta})\leq\sum_{j=1}^m \sum_{\beta_j\in J_j}\text{vol}\,\left(\mf{I}_{\beta_j}\right)<\varepsilon.\]
This shows that $\mk{D}$ has Jordan content zero.
\end{myproof}

\begin{example}{}
It is obvious that a one-point subset of $\mb{R}^n$ has Jordan content zero. It follows that any finite subset of $\mb{R}^n$ has Jordan content zero.
\end{example}

Now we want to consider general Jordan measurable sets. We first prove the following two theorems, giving more examples of Riemann integrable functions. The first one is a special case of the second one, but we need to prove it first to prove the second theorem.
\begin{theorem}[label=230827_8]{}
Let $\mf{I}=\di\prod_{i=1}^n[a_i, b_i]$, and let $f:\mf{I}\to\mb{R}$ be a bounded function defined on $\mf{I}$. If $f:\mf{I}\to\mb{R}$ is continuous on the interior of $\mf{I}$, then $f:\mf{I}\to\mb{R}$ is Riemann integrable.
\end{theorem}
\begin{myproof}{Proof}

We will show that  for any $\varepsilon>0$, there is a partition $\mf{P}$ of $\mf{I}$ such that $U(f,\mf{P})-L(f,\mf{P})<\varepsilon$. 

 Since $f:\mf{I}\to\mb{R}$ is a bounded function, there is a positive number $M$ such that
\[|f(\mf{x})|\leq M\hspace{1cm}\text{for all}\;\mf{x}\in \mf{I}.\]

By Lemma \ref{230825_4}, there is a closed rectangle $\check{\mf{I}}=\di\prod_{i=1}^n[u_i, v_i]$ contained in the interior of $\mf{I}$, such that
\[\text{vol}\,(\mf{I})-\text{vol}\,(\check{\mf{I}})<\frac{\varepsilon}{4M}.\]

\bp
Let $\mf{P}_0=(P_1, \ldots, P_n)$ be the partition of $\mf{I}$ given by $P_i=\{a_i, u_i, v_i, b_i\}$ for $1\leq i\leq n$. Then $\check{\mf{I}}$ is a rectangle in the partition $\mf{P}_0$. 

Since $f:\check{\mf{I}}\to\mb{R}$ is continuous, Theorem \ref{230827_7} implies that there is a partition $\mf{P}_1$ of $\check{\mf{I}}$ such that
\[U(f,\mf{P}_1)-L(f,\mf{P}_1)<\frac{\varepsilon}{2}.\]
Let $\mf{P}$ be the partition of $\mf{I}$ so that it contains all the partition points in $\mf{P}_0$ and $\mf{P}_1$. Then $\mf{P}$ is a refinement of $\mf{P}_0$ and the partition that $\mf{P}$ induces on $\check{\mf{I}}$ is $\mf{P}(\check{\mf{I}})=\mf{P}_1$. 
By Proposition \ref{230822_4},
\begin{align*}
&U(f,\mf{P})-L(f,\mf{P}) =\sum_{\mf{J}\in\mathcal{J}_{\mf{P}_0}}\left(U(f,\mf{P}(\mf{J}))-L(f,\mf{P}(\mf{J}))\right)\\
&=U(f,\mf{P}_1)-L(f,\mf{P}_1)+\sum_{\mf{J}\in\mathcal{J}_{\mf{P}_0}\setminus\{\check{\mf{I}}\}}\left(U(f,\mf{P}(\mf{J}))-L(f,\mf{P}(\mf{J}))\right).\end{align*}

For each $\mf{J}$ in $\mathcal{J}_{\mf{P}_0}\setminus\{\check{\mf{I}}\}$,  we use the crude estimate  
\[U(f, \mf{P}(\mf{J})-L(f, \mf{P}(\mf{J}))\leq  2M \,\text{vol}\,(\mf{J}).\]

These imply that
\begin{align*}
U(f,\mf{P})-L(f,\mf{P})
&<\frac{\varepsilon}{2}+2M\sum_{\mf{J}\in\mathcal{J}_{\mf{P}_0}\setminus\{\check{\mf{I}}\}}\text{vol}\,(\mf{J})\\
&=\frac{\varepsilon}{2}+2M\left(\text{vol}\,(\mf{I})-\text{vol}\,(\check{\mf{I}})\right)<\varepsilon.
\end{align*}
\end{myproof}

\begin{highlight}{Set of Discontinuities of a Function}

Given a function $f:A\to\mb{R}$ defined on the set $A$, the set of discontinuities of $f$ is the set of all points $\mf{x}_0$ in $A$ such that $f$ is not continuous at $\mf{x}_0$. 

If $B$ is a subset of $A$, and $\mf{x}_0$ is a point  of $B$, $f:A\to\mb{R}$ is continuous at $\mf{x}_0$ implies that $f:B\to\mb{R}$ is continuous at $\mf{x}_0$. Hence, the set of discontinuities of the function $f:B\to\mb{R}$ is a subset of the set of discontinuities of the function $f:A\to\mb{R}$.  
\end{highlight}

\begin{theorem}[label=230827_11]{}
Let $\mf{I}=\di\prod_{i=1}^n[a_i, b_i]$. Given that $f:\mf{I}\to\mb{R}$ is a bounded function defined on $\mf{I}$, let $\mathcal{N}_f$ be the set of discontinuities of $f:\mf{I}\to\mb{R}$. If $\mathcal{N}_f$ is a set that has Jordan content zero, then $f:\mf{I}\to\mb{R}$ is Riemann integrable.
\end{theorem}
\begin{myproof}{Proof}
We will  show that for any $\varepsilon>0$,  there is a partition $\mf{P}$ of $\mf{I}$ such that $U(f,\mf{P})-L(f,\mf{P})<\varepsilon$. 
 
Since $f:\mf{I}\to\mb{R}$ is a bounded function, there is a positive number $M$ such that
\[|f(\mf{x})|\leq M\hspace{1cm}\text{for all}\;\mf{x}\in \mf{I}.\]

Since $\mathcal{N}_f$ is a set of Jordan content zero that is contained in $\mf{I}$, there are closed rectangles $\mf{I}_1, \ldots, \mf{I}_k$ such that
\[\mathcal{N}_f\subset \bigcup_{j=1}^k\mf{I}_j\quad\text{and}\quad\sum_{j=1}^k\text{vol}\,(\mf{I}_j)<\frac{\varepsilon}{4M}.\]

By Lemma \ref{230827_2}, there is a partition $\mf{P}_0$ of $\mf{I}$ such that  each rectangle $\mf{J}$ in the partition $\mf{P}_0$  is either contained in an $ \mf{I}_j$ for some $1\leq j\leq k$, or  is disjoint from the interiors of $\mf{I}_j$ for all $1\leq j\leq k$.

Let
\[\mathscr{A}=\left\{\mf{J}\in\mathcal{J}_{\mf{P}_0}\,|\, \mf{J}\subset \mf{I}_j\;\text{for some}\;1\leq j\leq k\right\},\]
and
\[\mathscr{B}=\left\{\mf{J}\in\mathcal{J}_{\mf{P}_0}\,|\, \mf{J}\cap \text{int}\, (\mf{I}_j)=\emptyset\;\text{for all}\;1\leq j\leq k\right\}.\]

Assume that $\mathscr{B}$ contains $N$ rectangles. If $\mf{J}\in\mathscr{B}$, $f:\mf{J}\to\mb{R}$ is continuous on the interior of $\mf{J}$. 
By Theorem  \ref{230827_8}, $f:\mf{J}\to\mb{R}$ is Riemann integrable. 
Therefore, there is a partition $\mf{P}_{\mf{J}}$ of $\mf{J}$ such that
\[U(f,\mf{P}_{\mf{J}})-L(f, \mf{P}_{\mf{J}})<\frac{\varepsilon}{2N}.\]

\bp
The rest of the proof is similar to the proof of Theorem \ref{230827_8}. Let $\mf{P}$ be the partition of $\mf{I}$ which contains all the partition points in $\mf{P}_0$ and $\mf{P}_{\mf{J}}$ for all $\mf{J}\in\mathscr{B}$.  
Then 
\begin{align*}
&U(f,\mf{P})-L(f,\mf{P}) =\sum_{\mf{J}\in\mathcal{J}_{\mf{P}_0}}\left(U(f,\mf{P}(\mf{J}))-L(f,\mf{P}(\mf{J}))\right)\\
&=\sum_{\mf{J}\in\mathscr{A}}\left(U(f,\mf{P}(\mf{J}))-L(f,\mf{P}(\mf{J}))\right)+\sum_{\mf{J}\in\mathscr{B}}\left(U(f,\mf{P}(\mf{J}))-L(f,\mf{P}(\mf{J}))\right).\end{align*}

For each $\mf{J}$ in $\mathscr{A}$,  we use the crude estimate  
\[U(f, \mf{P}(\mf{J})-L(f, \mf{P}(\mf{J}))\leq  2M \,\text{vol}\,(\mf{J}).\]

Using the fact that \[ \bigcup_{\mf{J}\in\mathscr{A}} \mf{J}\;\subset\;\bigcup_{j=1}^k\mf{I}_j,\] we have
\[ \sum_{\mf{J}\in \mathscr{A}}\left(U(f,\mf{P}(\mf{J}))-L(f, \mf{P}(\mf{J)})\right)\leq 2M\sum_{j=1}^k\text{vol}\,(\mf{I}_j)<\frac{\varepsilon}{2}.\]

For each $\mf{J}\in\mathscr{B}$, $\mf{P}(\mf{J})$ is a refinement of $\mf{P}_{\mf{J}}$, and thus
\[U(f,\mf{P}(\mf{J}))-L(f, \mf{P}(\mf{J)})<\frac{\varepsilon}{2N}.\]This implies that
\[\sum_{\mf{J}\in \mathscr{B}}\left(U(f,\mf{P}(\mf{J}))-L(f, \mf{P}(\mf{J)})\right)<\frac{\varepsilon}{2}.\]
These give us
$\di 
U(f,\mf{P})-L(f,\mf{P})<\varepsilon$, as desired.
\end{myproof}

Now we can prove the following characterization of Jordan measurable sets.
\begin{theorem}[label=230827_15]{}
Let $\mk{D}$ be a bounded subset of $\mb{R}^n$. The following are equivalent.
\begin{enumerate}[(a)]
\item
$\mk{D}$ is a Jordan measurable set.
\item The boundary of $\mk{D}$ has Jordan content zero.
\end{enumerate}

\end{theorem}
\begin{myproof}{Proof}
Let $\mf{I}$ be a closed rectangle that contains $\mk{D}$. By definition, $\mk{D}$ is Jordan measurable if and only if the function $\chi_{\mk{D}}:\mf{I}\to\mb{R}$ is Riemann integrable.

By Theorem \ref{230827_10}, the set of discontinuities of the function $\chi_{\mk{D}}:\mf{I}\to\mb{R}$ is the set $\pa\mk{D}$. 
If the boundary of $\mk{D}$ has Jordan content zero, Theorem \ref{230827_11} implies that $\chi_{\mk{D}}:\mf{I}\to\mb{R}$ is Riemann integrable. This proves (b) implies (a).

Conversely, if $\mk{D}$  is Jordan measurable, given $\varepsilon>0$, there is a partition $\mf{P}$ such that
\[U(\chi_{\mk{D}},\mf{P})-L(\chi_{\mk{D}},\mf{P})<\frac{\varepsilon}{2}.\]

For each $\mf{J}$ in $\mathcal{J}_{\mf{P}}$, there are only three possibilities for the pair $(m_{\mf{J}}, M_{\mf{J}})$. Namely, $(1,1)$, $(0,0)$ or $(0,1)$.
Let
\begin{align*}
\mathscr{A}&=\left\{\mf{J}\in\mathcal{J}_{\mf{P}}\,|\, m_{\mf{J}}=M_{\mf{J}}=1\right\},\\
\mathscr{B}&=\left\{\mf{J}\in\mathcal{J}_{\mf{P}}\,|\, m_{\mf{J}}=M_{\mf{J}}=0\right\},\\
\mathscr{C}&=\left\{\mf{J}\in\mathcal{J}_{\mf{P}}\,|\, m_{\mf{J}}=0, \,M_{\mf{J}}=1\right\}.
\end{align*}

Then $\mathcal{J}_{\mf{P}}=\mathscr{A}\cup\mathscr{B}\cup\mathscr{C}$, and we have
\[U(\chi_{\mk{D}},\mf{P})-L(\chi_{\mk{D}},\mf{P})=\sum_{\mf{J}\in\mathcal{J}_{\mf{P}}}(M_{\mf{J}}-m_{\mf{J}})\,\text{vol}\,(\mf{J})=
\sum_{\mf{J}\in\mathscr{C}}\text{vol}\,(\mf{J}).\]
This implies that
\[\sum_{\mf{J}\in\mathscr{C}}\text{vol}\,(\mf{J})<\frac{\varepsilon}{2}.\]

Notice that $\mf{J}$ is in $\mathscr{A}$ if and only if  $f(\mf{x})=1$ for all $\mf{x}\in\mf{J}$, if and only if $\mf{J}\subset \mk{D}$. This implies that \[\text{int}\,\mf{J}\subset\text{int}\,\mk{D}\hspace{1cm}\text{for all}\;\mf{J}\in\mathscr{A}.\]
\bp
 Similarly, $\mf{J}$ is in $\mathscr{B}$ if and only if  $f(\mf{x})=0$ for all $\mf{x}\in\mf{J}$, if and only if $\mf{J}\subset \mb{R}^n\setminus \mk{D}$. This implies that  \[\text{int}\,\mf{J}\subset\text{ext}\,\mk{D}\hspace{1cm}\text{for all}\;\mf{J}\in\mathscr{B}.\]
 
Let
\[\mathcal{S}=\bigcup_{\mf{J}\in\mathscr{A}\cup\mathscr{B}}\pa\mf{J}.\]

Since $\mb{R}^n$ is a disjoint union of $\text{int}\,\mk{D}$, $\text{ext}\,\mk{D}$ and $\pa\,\mk{D}$, we must have
\[\pa\mk{D}\subset\left(\bigcup_{\mf{J}\in\mathscr{C}}\mf{J}\right)\;\cup\;\mathcal{S}.\]

Since the boundary of a closed rectangle has Jordan content zero, and    $\mathscr{A}\cup\mathscr{B}$ is a finite set, Theorem \ref{230827_13} implies that $\mathcal{S}$ has Jordan content zero. Hence, there is a finite collection of rectangles $\mathscr{D}=\left\{\mf{I}_j\,|\, 1\leq j\leq k\right\}$ such that
\[\mathcal{S}\subset \bigcup_{j=1}^k \mf{I}_j\quad\text{and}\quad \sum_{j=1}^k\text{vol}\,(\mf{I}_j)<\frac{\varepsilon}{2}.\]
Let $\mathscr{E}=\mathscr{C}\cup\mathscr{D}$. Then $\mathscr{E}$ is a finite collection of closed rectangles,
\[\pa\,\mk{D}\subset \bigcup_{\mf{J}\in \mathscr{E}}\mf{J}\quad\text{and}\quad \sum_{\mf{J}\in \mathscr{E}} \text{vol}\,(\mf{J})< \varepsilon.\]This shows that $\mk{D}$ has Jordan content zero.

\end{myproof}
Using Theorem \ref{230827_15}, we can obtain more examples of Jordan measurable sets. First we prove  the following.
\begin{lemma}[label=230827_16]{}
Let $A$ and $B$ be subsets of $\mb{R}^n$. Then
\[\pa (A\cup B)\subset \pa A\cup \pa B,\hspace{1cm}\pa(A\cap B)\subset \pa A\cup \pa B.\]
\end{lemma}
\begin{myproof}{Proof}
If  $\mf{x}_0$ be a point in $\pa (A\cup B)$,   there is sequence    of points $\{\mf{u}_k\}$ in $A\cup B$ that converges to $\mf{x}_0$. Each point in this sequence is either in $  A$ or in $  B$. Therefore, there is a subsequence $\{\mf{u}_{k_j}\}$ that is in $ A$ or in $ B$. There is also a sequence $\{\mf{v}_k\}$ in $\mb{R}^n\setminus (A\cup B)$ that converges to $\mf{x}_0$. This sequence is in both $\mb{R}^n\setminus A$ and in $\mb{R}^n\setminus B$. Therefore, $\mf{x}_0$ is in $\pa  A$ or in  $\pa B$.

If  $\mf{x}_0$ is a point in $\pa (A\cap B)$, there is sequence    of points $\{\mf{u}_k\}$ in $\mb{R}^n\setminus (A\cap B)$ that converges to $\mf{x}_0$. Each point in this sequence is either in $\mb{R}^n\setminus A$ or in $\mb{R}^n\setminus B$. Therefore, there is a subsequence $\{\mf{u}_{k_j}\}$ that is in $\mb{R}^n\setminus A$ or in $\mb{R}^n\setminus B$. There is also a sequence $\{\mf{v}_k\}$ in $A\cap B$ that converges to $\mf{x}_0$. This sequence is in both $A$ and   $B$. Therefore, $\mf{x}_0$ is in $\pa  A$ or in  $\pa B$.

\end{myproof}

One is tempted to think that $\pa(A\cap B)\subset \pa A \cap \pa B$. But this is not true, as shown in the following example.
\begin{example}{}
Let $A=[0,2]\times [0,2]$ and $B=[1,3]\times [1,3]$. We find that 
$A\cap B=[1,2]\times [1,2]$. 
As shown in Figure \ref{figure99}, $\pa A\cap \pa B$ is a set with 4 points,   $\pa (A\cap B)\neq \pa A\cap \pa B$, but $\pa (A\cap B)\subset \pa A\cup\pa B$.
\end{example}

 \begin{figure}[ht]
\centering
\includegraphics[scale=0.2]{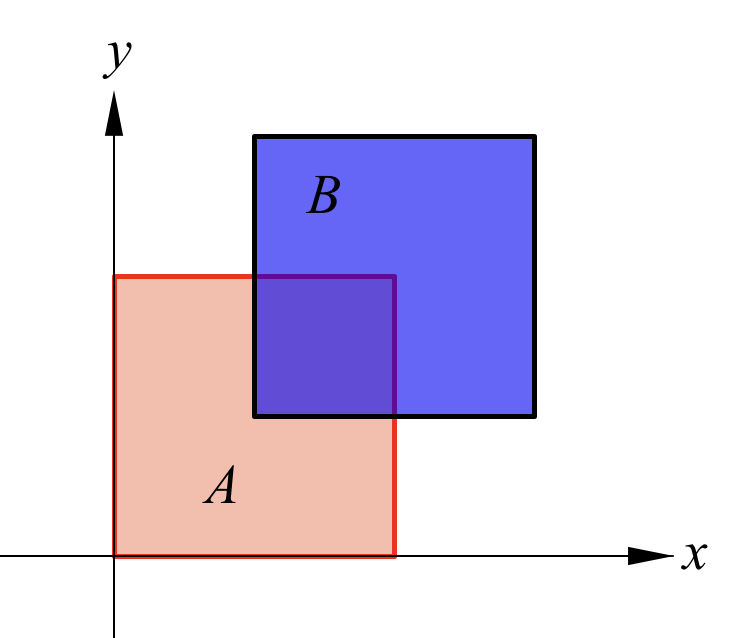}
\caption{$\pa (A\cap B)\neq \pa A\cap \pa B$, but $\pa (A\cap B)\subset \pa A\cup\pa B$. }\label{figure99}
\end{figure}

Using Lemma \ref{230827_16}, we obtain the following.
\begin{theorem}
{}
If $\mk{D}_1$, $\mk{D}_2$,  $\ldots$, $\mk{D}_m$ are Jordan measurable sets, then the set  $\mk{D}_1\cap \mk{D}_2\cap\cdots\cap\mk{D}_m$  and the set 
$\mk{D}_1\cup \mk{D}_2\cup\cdots\cup\mk{D}_m$ are also Jordan measurable.
\end{theorem}
\begin{myproof}{Proof}
It suffices to prove the case where $m=2$. The general case follows by induction. 

If $\mk{D}_1$ and $\mk{D}_2$ are Jordan measurable, Theorem \ref{230827_15} says that $\pa\mk{D}_1$ and $\pa\mk{D}_2$ have Jordan content zero. Theorem \ref{230827_13} says that $\pa\mk{D}_1\cup \pa\mk{D}_2$ has Jordan content zero. Lemma \ref{230827_16} and Theorem \ref{230827_9} imply that $\pa(\mk{D}_1\cap\mk{D}_2)$ and $\pa (\mk{D}_1\cup\mk{D}_2)$ have Jordan content zero. Theorem \ref{230827_15}  again implies that $\mk{D}_1\cap\mk{D}_2$ and $\mk{D}_1\cup\mk{D}_2$ are Jordan measurable sets.
\end{myproof}

Observe that the concept of Jordan measurable sets and Riemann integrable functions are closely related. In  a nutshell, a set $\mk{D}$ is a Jordan measurable set if and only if all the constant functions $f:\mk{D}\to\mb{R}$ are Riemann integrable. In fact,  it is also if and only if all the continuous functions  $f:\mk{D}\to\mb{R}$ are Riemann integrable. We can even allow discontinuities on a set that has Jordan content zero. We will prove this after a few preparatory remarks and lemmas.

\begin{remark}{}
If $\mk{D}$ is a bounded subset of $\mb{R}^n$ and $\mk{D}$ is contained in the closed rectangle $\mf{I}$, then $\overline{\mk{D}}$ is also contained in $\mf{I}$. This implies that $\pa \mk{D}$ is also contained in $\mf{I}$.
\end{remark}

The following example depicts the relation between the set of dicontinuities of a function $f:\mk{D}\to\mb{R}$ and its zero extension.
\begin{example}[label=230827_17]{}
Let $\mk{D}=\left\{(x,y)\,|\, 0\leq y\leq x\leq 2\right\}$, and let $f:\mk{D}\to\mb{R}$ be the function defined as
\[f(x,y)=\begin{cases} 1,\quad &\text{if}\; (x,y)\in\mk{D}\;\text{and}\; 0 \leq x<1,\\2,\quad &\text{if}\; (x,y)\in\mk{D}\;\text{and}\;1\leq x\leq 2.\end{cases}\]

\be
The set of discontinuities of the function $f:\mk{D}\to\mb{R}$ is
\[\mathcal{N}_f=\{(1,y)\,|\, 0\leq y\leq x\}.\] 
The set $\mk{D}$ is contained in the rectangle $\mf{I}=[0,2]\times [0,2]$. 
Let $\check{f}:\mf{I}\to\mb{R}$ be the zero extension of $f:\mk{D}\to\mb{R}$. 
Then the set of discontinuities of $\check{f}:\mf{I}\to\mb{R}$ is
\[\mathcal{N}_{\check{f}}=\left\{(1,y)\,|\, 0\leq y\leq 1\right\}\cup \{(x,x)\,|\, 0\leq x\leq 2\},\]
which is a subset of $\mathcal{N}_f\cup \pa \mk{D}$. 
\end{example2}

 \begin{figure}[ht]
\centering
\includegraphics[scale=0.2]{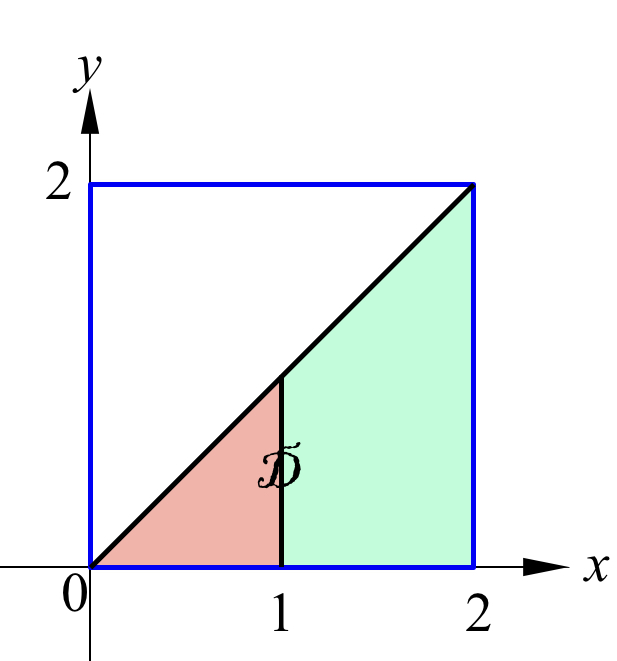}
\caption{The set of discontinuities of the functions discussed in Example \ref{230827_17}. }\label{figure100}
\end{figure}
The following lemma gives the general case.

\begin{lemma}[label=230827_18]{}
Let $\mk{D}$ be a bounded subset of $\mb{R}^n$. Given that $f:\mk{D}\to\mb{R}$ is a bounded function defined on $\mk{D}$, let $\check{f}:\mb{R}^n\to\mb{R}$ be its zero extension. Let
\begin{align*}
\mathcal{N}_f&=\left\{\mf{x}_0\in\mk{D}\,|\, f:\mk{D}\to\mb{R} \;\text{is not continuous at $\mf{x}_0$}\right\},\\
\mathcal{N}_{\check{f}}&=\left\{\mf{x}_0\in\mb{R}^n\,|\, \check{f}:\mb{R}^n\to\mb{R} \;\text{is not continuous at $\mf{x}_0$}\right\}.
\end{align*} Then 
\[\mathcal{N}_{\check{f}}\subset \pa \mk{D}\;\cup\; \mathcal{N}_f.\]
 
\end{lemma}
\begin{myproof}{Proof}
We will    show that $\check{f}:\mb{R}^n\to\mb{R}$ is continuous at $\mf{x}_0$ if $\mf{x}_0$ is not in $\pa \mk{D}\;\cup\; \mathcal{N}_f$. If $\mf{x}_0\notin \pa \mk{D}\;\cup\; \mathcal{N}_f$, there are two possibilities. 
\begin{enumerate}[$\bullet$\;\;]
\item  $\mf{x}_0$ is in $\text{ext}\,\mk{D}$.
\item $\mf{x}_0$ is in $\text{int}\,\mk{D}$ and $f:\mk{D}\to\mb{R}$ is continuous at $\mf{x}_0$.

\end{enumerate}

If $\mf{x}_0$ is in $\text{ext}\,\mk{D}$, there is an $r>0$ such that $B(\mf{x}_0, r)\subset \mb{R}^n\setminus\mk{D}$. Since $\check{f}(\mf{x})=0$ for all $\mf{x}\in B(\mf{x}_0, r)$, $\check{f}:\mb{R}^n\to\mb{R}$ is continuous at $\mf{x}_0$.

If $\mf{x}_0$ is in $\text{int}\,\mk{D}$ and $f:\mk{D}\to\mb{R}$ is continuous at $\mf{x}_0$, we want to show  that $\check{f}:\mb{R}^n\to\mb{R}$ is continuous at $\mf{x}_0$. Given that $\{\mf{x}_k\}$ is a sequence in $\mb{R}^n$ that converges to $\mf{x}_0$,   there is a positive integer $K$ such that $\mf{x}_k\in \text{int}\,\mk{D}$ for all $k\geq K$. This implies that $\check{f}(\mf{x}_k)=f(\mf{x}_k)$ for all $k\geq K$. Since $f:\mk{D}\to\mb{R}$ is continuous at $\mf{x}_0$, 
\[\lim_{k\to\infty}f(\mf{x}_{K+k})=f(\mf{x}_0).\]

This implies that
\[\lim_{k\to\infty} \check{f}(\mf{x}_k)=\lim_{k\to\infty}\check{f}(\mf{x}_{K+k})=\lim_{k\to\infty}f(\mf{x}_{K+k})=f(\mf{x}_0)=\check{f}(\mf{x}_0).\]
Hence, $\check{f}:\mb{R}^n\to\mb{R}$ is continuous at $\mf{x}_0$. This completes the proof.

\end{myproof}

Now we can prove the main theorem.
\begin{theorem}[label=230830_1]{}
 Let $\mk{D}$ be a bounded subset of $\mb{R}^n$. Given that $f:\mk{D}\to\mb{R}$ is a bounded function defined on $\mk{D}$, let $\mathcal{N}_f$ be the set of discontinuities of $f:\mk{D}\to\mb{R}$. If $\mk{D}$ is a Jordan measurable set, and $\mathcal{N}_f$ is a set that has Jordan content zero, then $f:\mk{D}\to\mb{R}$ is Riemann integrable.
\end{theorem}
\begin{myproof}{Proof}
Let $\mf{I}$ be a closed rectangle in $\mb{R}^n$ that contains $\mk{D}$, and let $\check{f}:\mf{I}\to\mb{R}$ be the zero extension of $f:\mk{D}\to\mb{R}$ to $\mf{I}$. 
We want to show that $\check{f}:\mf{I}\to\mb{R}$ is Riemann integrable.
\bp
By Lemma \ref{230827_18}, the set of discontinuities of  $\check{f}:\mf{I}\to\mb{R}$ is contained in  $\pa \mk{D}\cup\mathcal{N}_f$. Since $\mk{D}$ is Jordan measurable, Theorem \ref{230827_15} says that $\pa \mk{D}$ has Jordan content zero. Theorem \ref{230827_13} and Theorem \ref{230827_9} then imply that the set $\mathcal{N}_{\check{f}}$ has Jordan content zero. By Theorem \ref{230827_11}, $\check{f}:\mf{I}\to\mb{R}$ is Riemann integrable. This completes the proof.
\end{myproof}
In particular, we have the following.

\begin{corollary}[label=230828_1]{}
Let $\mk{D}$ be a subset of $\mb{R}^n$ that is Jordan measurable, and let $f:\mk{D}\to\mb{R}$ be a bounded function defined on $\mk{D}$. If $f:\mk{D}\to\mb{R}$ is continuous, then it is Riemann integrable.
\end{corollary}

Let us emphasize the result in Corollary \ref{230828_1}.
\begin{highlight}{Riemann Integrability of Continuous Functions}
Any continuous function defined on  a Jordan measurable set is Riemann integrable.
\end{highlight}

Another interesting corollary of Theorem \ref{230830_1} is the following.
\begin{corollary}[label=230828_2]{}
Let $\mk{D}$ be a subset of $\mb{R}^n$ that has Jordan content zero. If $f:\mk{D}\to\mb{R}$ is a bounded function defined on $\mk{D}$, then $f:\mk{D}\to\mb{R}$ is Riemann integrable and 
\[\int_{\mk{D}}f\;=\;0.\]
\end{corollary}
\begin{myproof}{Proof}
Since $\mk{D}$ has Jordan content zero, and the set $\mathcal{N}_f$ of discontinuities of $f:\mk{D}\to\mb{R}$ is a subset of $\mk{D}$, $\mathcal{N}_f$ has Jordan content zero.  Theorem \ref{230830_1}  implies  that $f:\mk{D}\to\mb{R}$ is Riemann integrable. Since $f:\mk{D}\to \mb{R}$ is bounded, there is a positive number $M$ such that 
$\di -M\leq f(\mf{x})\leq M$ for all $\mf{x}\in\mk{D}$.
\bp
By monotonicity theorem,
\[\int_{\mk{D}}-M\,d\mf{x}\;\leq\;\int_{\mk{D}}f(\mf{x})\,d\mf{x}\;\leq\; \int_{\mk{D}}M\,d\mf{x}.\]
For any constant $c$, linearity implies that
$\di \int_{\mk{D}}c\,d\mf{x}\;=\;c\,\text{vol}\,(\mk{D})=0$.
This proves that 
\[\int_{\mk{D}}f\;=\;0.\]
\end{myproof}

Let us highlight this important result. 
\begin{highlight}{Bounded Functions Defined on Sets that has Jordan Content Zero}
Any  bounded function defined on  a set that has  Jordan content zero is Riemann integrable with integral zero.
\end{highlight}

Corollary \ref{230828_2} also gives the following.
\begin{lemma}[label=230828_3]{}
Let $\mk{D}$ be a bounded subset of $\mb{R}^n$, and let $f:\mk{D}\to\mb{R}$ be a  bounded function defined on $\mk{D}$. If there is a subset $\mathcal{A}$ of $\mk{D}$ with Jordan content zero such that  $f(\mf{x})=0$ for all $\mf{x}\in\mk{D}\setminus\mathcal{A}$, then $f:\mk{D}\to\mb{R}$ is Riemann integrable and $\di \int_{\mk{D}}f=0$.
\end{lemma}
\begin{myproof}{Proof}
Let $\mf{I}$ be a closed rectangle in $\mb{R}^n$ that contains $\mk{D}$, and let $\check{f}:\mf{I}\to\mb{R}$ be the zero extension of the function $f:\mk{D}\to\mb{R}$. Notice that it is also the zero extension of the function $f|_{\mathcal{A}}:\mathcal{A}\to\mb{R}$. By Corollary \ref{230828_2}, $f|_{\mathcal{A}}:\mathcal{A}\to\mb{R}$ is Riemann integrable and $\di \int_{\mathcal{A}}f|_{\mathcal{A}}=0$. Hence,  $\check{f}:\mf{I}\to\mb{R}$ is Riemann integrable, and so is $f:\mk{D}\to\mb{R}$. Moreover,
\[\int_{\mk{D}}f=\int_{\mf{I}}\check{f}=\int_{\mathcal{A}}f|_{\mathcal{A}}=0.\]
\end{myproof}
Using Lemma \ref{230828_3}, we obtain the following important result, which says that Riemann integrability is not affected by the definition of the function on a set that has  Jordan content zero.
\begin{theorem}[label=230828_4]{}
Let $\mk{D}$ be a bounded subset of $\mb{R}^n$, and let $f:\mk{D}\to\mb{R}$ and $g:\mk{D}\to\mb{R}$ be bounded functions defined on $\mk{D}$. If $f:\mk{D}\to\mb{R}$ is Riemann integrable and there is a subset $\mathcal{A}$ of $\mk{D}$ which has Jordan content zero such that 
\[g(\mf{x})=f(\mf{x})\hspace{1cm}\text{for all}\;\mf{x}\in \mk{D}\setminus \mathcal{A},\]
then $g:\mk{D}\to\mb{R}$ is Riemann integrable and
\[\int_{\mk{D}}g\;=\;\int_{\mk{D}}f.\]
\end{theorem}
\begin{myproof}{Sketch of Proof}
Let $h:\mk{D}\to\mb{R}$ be the function $h(\mf{x})=f(\mf{x})-g(\mf{x})$. Then $h(\mf{x})=0$ for all $\mf{x}\in\mk{D}\setminus\mathcal{A}$. By Lemma \ref{230828_3}, $h:\mk{D}\to\mb{R}$ is Riemann integrable and $\di\int_{\mk{D}}h=0$. The assertion follows from linearity. 
\end{myproof}

Using Theorem \ref{230828_4}, we can generalize additivity to arbitrary sets.

\begin{theorem}{Additivity}
Given that $\mk{D}_1$ and $\mk{D}_2$ are bounded subsets of $\mb{R}^n$ such that $\mk{D}_1\cap\mk{D}_2$ is a set that has Jordan content zero, let $\mk{D}=\mk{D}_1\cup\mk{D}_2$. Assume that $f:\mk{D}\to\mb{R}$ is a bounded function defined on $\mk{D}$.
  If the functions $f:\mk{D}_1\to\mb{R}$ and $f:\mk{D}_2\to \mb{R}$ are Riemann integrable, then the function  $f:\mk{D}\to\mb{R}$ is Riemann integrable and
\[\int_{\mk{D}}f=\int_{\mk{D}_1\cup\mk{D}_2}f=\int_{\mk{D}_1}f+\int_{\mk{D}_2}f.\]
 
\end{theorem}
\begin{myproof}{Proof}
Let $\mf{I}$ be a closed rectangle in $\mb{R}^n$ that contains $\mk{D}$, and let $\mk{D}_0=\mk{D}_1\cap \mk{D}_2$. We are given that $\mk{D}_0$ has   Jordan content zero. Let $\check{f}:\mf{I}\to\mb{R}$,  $\check{f}_1:\mf{I}\to\mb{R}$,   $\check{f}_2:\mf{I}\to\mb{R}$ and  $\check{f}_0:\mf{I}\to\mb{R}$ be resepctively the zero extensions of $f:\mk{D}\to\mb{R}$, $f:\mk{D}_1\to\mb{R}$,   $f:\mk{D}_2\to \mb{R}$ and $f:\mk{D}_0\to \mb{R}$. It is easy to see that
\begin{equation}\label{230828_5}\check{f}(\mf{x})=\check{f}_1(\mf{x})+\check{f}_2(\mf{x})-\check{f}_0(\mf{x})\hspace{1cm}\text{for all}\;\mf{x}\in\mf{I}.\end{equation}

By Corollary \ref{230828_2}, $\di f:\mk{D}_0\to\mb{R}$ is Riemann integrable and $\di \int_{\mk{D}_0}f=0$. Hence, $\di  \check{f}_0:\mf{I}\to\mb{R}$ is Riemann integrable.

Since $f:\mk{D}_1\to\mb{R}$ and $f:\mk{D}_2\to \mb{R}$ are Riemann integrable, $\check{f}_1:\mf{I}\to\mb{R}$,   $\check{f}_2:\mf{I}\to\mb{R}$ are Riemann integrable. Linearly and \eqref{230828_5} imply that  $\check{f}:\mf{I}\to\mb{R}$ is Riemann integrable and
\[\int_{\mk{D}}f=\int_{\mk{D}_1}f+\int_{\mk{D}_2}f-\int_{\mk{D}_0}f=\int_{\mk{D}_1}f+\int_{\mk{D}_2}f.\]

\end{myproof}

By induction and Theorem \ref{230827_13}, we obtain the following.
\begin{theorem}[label=230828_6]{}
Given that $\mk{D}_1$, $\mk{D}_2$, $\ldots$, $\mk{D}_m$  are bounded subsets of $\mb{R}^n$ such that for any pairs of $(i,j)$ with $i\neq j$, $\mk{D}_i\cap \mk{D}_j$ has Jordan content zero, let $\mk{D}=\mk{D}_1\cup\mk{D}_2\cup \cdots\cup\mk{D}_m$. Assume that  $f:\mk{D}\to\mb{R}$ is a bounded function defined on $\mk{D}$.  If the functions $f:\mk{D}_j\to\mb{R}$, $1\leq j\leq m$, are Riemann integrable, then the function $f:\mk{D}\to\mb{R}$ is Riemann integrable and
\[\int_{\mk{D}}f=\int_{\mk{D}_1\cup\mk{D}_2\cup\cdots\cup\mk{D}_m}f=\int_{\mk{D}_1}f+\int_{\mk{D}_2}f+\cdots+\int_{\mk{D}_m}f.\]
\end{theorem}

\begin{remark}{}
If a function $f:\mk{D}\to\mb{R}$ is Riemann integrable and $\mk{D}_1$ is a subset of $\mk{D}$, $f:\mk{D}_1\to\mb{R}$ is not necessarily Riemann integrable. For example, consider the constant function $f$ on $\mk{D}=[0,1]^n$ which takes value  1. Its restriction to $\mk{D}_1=\mk{D}\cap\mb{Q}^n$ is not Riemann integrable.

Thus for the general additivity theorem, we do not have if and only if.
\end{remark}

Using the fact that a set $\mk{D}$ is Jordan measurable if the function $\chi_{\mk{D}}:\mk{D}\to\mb{R}$ is measurable, Theorem \ref{230828_6} gives the following.
\begin{corollary}{}
Let $\mk{D}_1$, $\mk{D}_2$, $\ldots$, $\mk{D}_m$ be Jordan measurable subsets of $\mb{R}^n$ such that for any pairs of $(i,j)$ with $i\neq j$, $\mk{D}_i\cap \mk{D}_j$ has Jordan content zero. Then the set   $\mk{D}=\mk{D}_1\cup\mk{D}_2\cup \cdots\cup\mk{D}_m$ is also Jordan measurable. Moreover, 
\[\text{vol}\,(\mk{D})=\text{vol}\,(\mk{D}_1)+\text{vol}\,(\mk{D}_2)+\cdots+\text{vol}\,(\mk{D}_m).\]
\end{corollary}

Let us return to explore more on Jordan measurable sets. So far we only know explicitly that a closed rectangle is Jordan measurable, and some examples of sets that have Jordan content  zero.
Since a bounded subset $\mk{D}$ is Jordan measurable  if and only if its boundary has Jordan content zero, we will first explore sets that have Jordan content zero. The following theorem will give us a lots of examples of sets that have Jordan content zero. 

\begin{theorem}[label=230830_2]{}
Let $\mk{D}$ be a Jordan measurable subset of $\mb{R}^n$, and let $f:\mk{D}\to\mb{R}$ be a Riemann integrable function. Then the graph of $f$ defined by
\[G_f=\left\{(\mf{x}, y)\in\mb{R}^{n+1}\,|\,\mf{x}\in\mk{D}, y=f(\mf{x})\right\}\]
is a subset of $\mb{R}^{n+1}$ that has Jordan content zero.
\end{theorem}
\begin{myproof}{Proof}
Let $\mf{I}$ be a closed rectangle in $\mb{R}^n$ that contains $\mk{D}$, and let $\check{f}:\mf{I}\to\mb{R}$ be the zero extension of $f:\mk{D}\to\mb{R}$. 
Fixed $\varepsilon>0$. Since $f:\mk{D}\to\mb{R}$ is Riemann integrable, there is a partition $\mf{P} $ of $\mf{I}$ such that
\[U(\check{f}, \mf{I})-L(\check{f}, \mf{I})<\frac{\varepsilon}{2}.\]

 Let 
 \[\eta=\frac{\varepsilon}{4\text{vol}\,(\mf{I})}.\]
 \bp
 Then $\eta>0$.
Let 
\[\mathscr{A}=\left\{\mf{J}\times [m_{\mf{J}}-\eta, M_{\mf{J}}+\eta]\,|\, \mf{J}\in\mathcal{J}_{\mf{P}}\right\}.\]
Then $\mathscr{A}$ is a finite collection of closed rectangles in $\mb{R}^{n+1}$. If $(\mf{x}, f(\mf{x}))$ is in $G_f$, there is a $\mf{J}\in\mathcal{J}_{\mf{P}}$ such that $\mf{x}\in \mf{J}$. Then $m_{\mf{J}}\leq f(\mf{x})\leq M_{\mf{J}}$ implies that $ (\mf{x}, f(\mf{x}))$ is in $\mf{J}\times [m_{\mf{J}}-\eta, M_{\mf{J}}+\eta]$. This proves that 
\[G_f\subset\;\bigcup_{\mf{K}\in\mathscr{A}}\;\mf{K}.\]
    
Now,
\begin{align*}
&\sum_{\mf{K}\in\mathscr{A}}\text{vol}\,(\mf{K}) 
 =\sum_{\mf{J}\in\mathcal{J}_{\mf{P}}}\left(M_{\mf{J}}-m_{\mf{J}}+2\eta\right)\,\text{vol}\,(\mf{J})\\
&=U(\check{f},\mf{P})-L(\check{f},\mf{P})+2\eta \sum_{\mf{J}\in\mathcal{J}_{\mf{P}}}\text{vol}\,(\mf{J})<\frac{\varepsilon}{2}+2\eta \,\text{vol}\,(\mf{I})<\varepsilon.
\end{align*}This proves that $G_f$ has Jordan content zero.
\end{myproof}

Specialize to continuous functions, we have the following.
\begin{corollary}[label=230829_5]{}
Let $\mk{D}$ be a Jordan measurable set in $\mb{R}^n$, and let $f:\mk{D}\to\mb{R}$ be a continuous function. Then the graph of $f$ defined by
\[G_f=\left\{(\mf{x}, y)\in\mb{R}^{n+1}\,|\,\mf{x}\in\mk{D}, y=f(\mf{x})\right\}\]
is a subset of $\mb{R}^{n+1}$ that has Jordan content zero.
\end{corollary}

\begin{example}{}
Any line segment $L$ between two  points $(x_1, y_1)$ and $(x_2, y_2)$   in $\mb{R}^2$ has Jordan content zero. 

If $x_1=x_2$, the line segment $L$ is vertical. It is a subset of the boundary of the closed rectangle $[x_1, x_1+1]\times [y_1, y_2]$. Hence, $L$ has Jordan content zero.  
\be
If $x_1\neq x_2$, $L$ is the graph of the continuous function $f:[x_1, x_2]\to\mb{R}$,
\[f(x)=y_1+\frac{y_2-y_1}{x_2-x_1}(x-x_1).\]
Therefore $L$ also has Jordan content zero.
\end{example2}

 \begin{figure}[ht]
\centering
\includegraphics[scale=0.2]{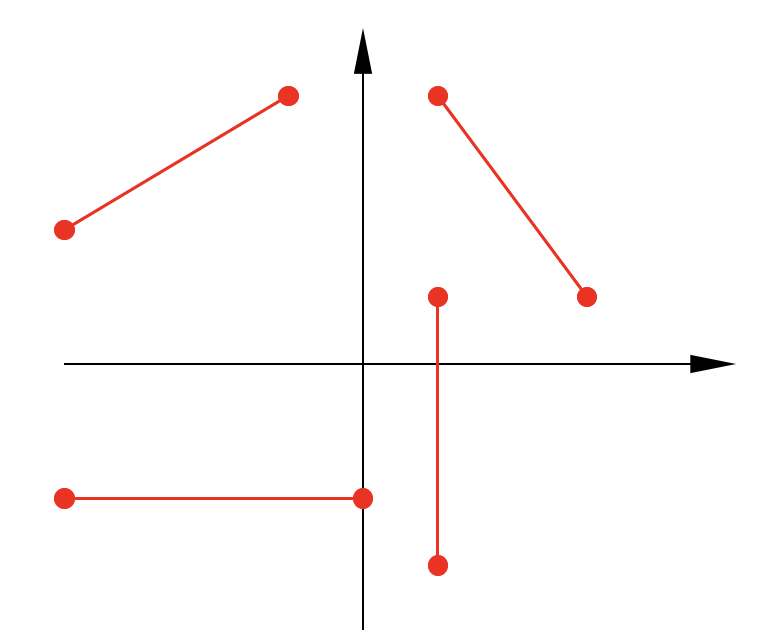}
\caption{Line segments in the plane have Jordan content zero. }\label{figure102}
\end{figure}

\begin{example}{}
Since the boundary of a polygon in $\mb{R}^2$ is a finite union of line segments, a polygon is Jordan measurable. The interior of the polygon is also Jordan measurable as it has the same boundary.
\end{example}

\begin{figure}[ht]
\centering
\includegraphics[scale=0.2]{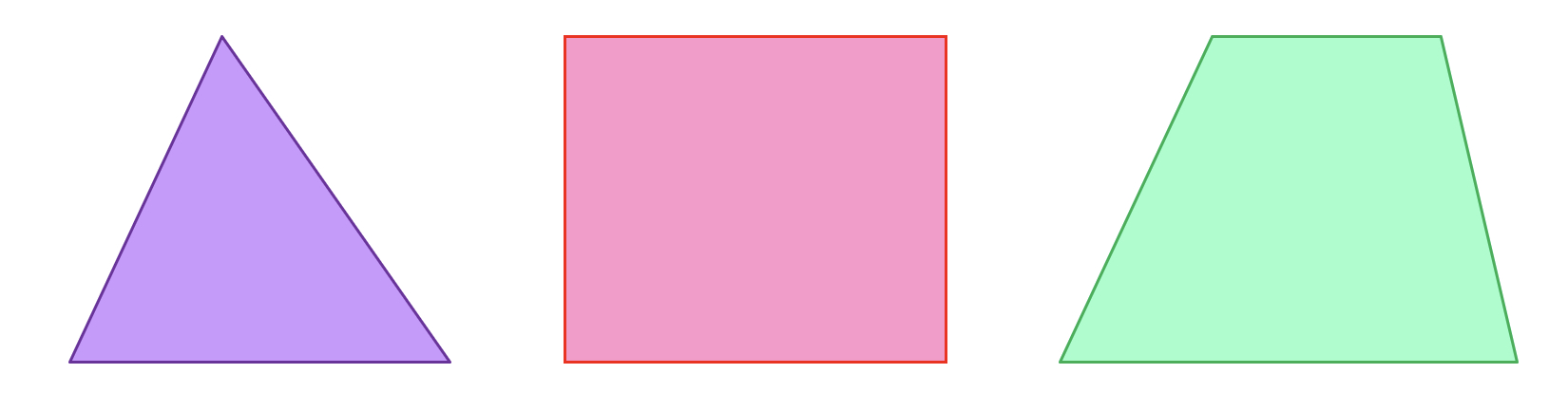}
\caption{Polygons in the plane are Jordan measurable. }\label{figure103}
\end{figure}

\begin{example}{}
We can argue that a line segment or the part of a plane in $\mb{R}^3$ contained inside a bounded  set has Jordan content zero. A plane in $\mb{R}^3$ satisfies an equation of the form 
\[ax+by+cz=d,\] where $(a,b,c)\neq \mf{0}$. Therefore, one can always solve one of the variables as a function of the other two. For example, if $c\neq 0$, then
\[z=f(x,y)=\frac{d-ax-by}{c}.\]
Hence, a plane is the graph of a continuous function. If we consider the part of the plane contained within  a bounded set, then it must have Jordan content zero. For example, 
\[ S=\left\{(x,y,z)\,|\, x+y+z=3, x\geq 0, y\geq 0, z\geq 0\right\}\] is the part of the plane $ x+y+z=3$ bounded inside the rectangle $[0,3]\times [0,3]\times [0,3]$. Hence, $S$ has Jordan content zero. 

A line segment in $\mb{R}^3$ can always be regarded as a subset of a part of a plane that is contained in a bounded set. Hence, it also has Jordan content zero.
\end{example}

\begin{example}[label=230829_3]{}
The boundary of the open rectangle $U=\di\prod_{i=1}^n (a_i, b_i)$  is the same as the boundary of its closure $R=\di\prod_{i=1}^n[a_i,b_i]$. Hence, $U$ is also a Jordan measurable set. Since $U$ and $\pa U$ are disjoint, and their union is $R$,
\[\text{vol}\,U+\text{vol}\,\pa U=\text{vol}\, R.\]
Since $\text{vol}\,(\pa U)=0$, we have
\[\text{vol}\,U=\text{vol}\, R.\]
\end{example}

Motivated by Example \ref{230829_3}, an interesting question to ask is if the subset $\mk{D}$ of $\mb{R}^n$ is Jordan measurable, is its closure $\overline{\mk{D}}$ Jordan measurable? This is answered in the following theorem.
\begin{theorem}[label=230829_4]{}
 If $\mk{D}$  is a subset of $\mb{R}^n$ that is Jordan measurable, so is $\overline{\mk{D}}$. Moreover,
\[\text{vol}\,\overline{\mk{D}}=\text{vol}\,\mk{D}.\]
\end{theorem}
\begin{myproof}{Proof} First we claim that
$\di \pa \overline{\mk{D}}\subset \pa \mk{D}$. As the closure of $\mk{D}$,
$\overline{\mk{D}}$ is a disjoint union of $\text{int}\, \mk{D}$ and $\pa \mk{D}$. As the closure of $\overline{\mk{D}}$, $\overline{\mk{D}}$  is  a disjoint union of $\text{int}\, \overline{\mk{D}}$ and $\pa \overline{\mk{D}}$. Since $\mk{D}\subset \overline{\mk{D}}$, we have $\text{int}\,\mk{D}\subset\text{int}\,\overline{\mk{D}}$. Hence, we must have $\di \pa \overline{\mk{D}}\subset \pa \mk{D}$. 

If $\mk{D}$ is Jordan measurable, $\pa \mk{D}$ has Jordan content zero. Since $\di \pa \overline{\mk{D}}\subset \pa \mk{D}$, $\pa \overline{\mk{D}}$ also has Jordan content zero. Hence, $\overline{\mk{D}}$ is Jordan measurable.

For the last statement, we use the fact that $\overline{\mk{D}}=\mk{D}\cup\pa\mk{D}$. Notice that $\mk{D}\cap \pa\mk{D}\subset\pa\mk{D}$. Hence, $\mk{D}\cap\pa\mk{D}$ has Jordan content zero.  By the additivity theorem, 
\[\text{vol}\,\mk{D}+\text{vol}\,\pa\mk{D} =\text{vol}\,\overline{\mk{D}}.\]
Since $\text{vol}\,\pa\mk{D}=0$, we conclude that $\text{vol}\,\overline{\mk{D}}=\text{vol}\,\mk{D}$.
\end{myproof}

\begin{example}{}
Consider the set $A=(-1,0)\cup (0,1)$. Its closure is $\overline{A}=[-1,1]$. Hence, $\pa\overline{A}=\{-1,1\}$ is not equal to $\pa A=\{-1, 0, 1\}$. This also shows that even for an open set $A$, we does not have $\pa A=\pa\overline{A}$.
\end{example}
\begin{remark}{}
If $\mk{D}$ is a bounded subset of $\mb{R}^n$ such that $\overline{\mk{D}}$ is Jordan measurable, one cannot deduce that $\mk{D}$ is Jordan measurable. An example is given by $\mk{D}=[0,1]^n\cap \mb{Q}^n$, which is not Jordan measurable, but $\overline{\mk{D}}=[0,1]^n$ is Jordan measurable.
\end{remark}

\begin{example}[label=230829_8]{}
Let $r$ be a positive number. We claim that the disc \[\mk{D}=\di\left\{(x,y)\,|\,x^2+y^2<r^2\right\}\] and its closure are Jordan measurable sets. By Theorem \ref{230829_4}, it is sufficient to show that $\mk{D}$ is Jordan measurable.  Notice that
\[\pa\mk{D}= \left\{(x,y)\,|\,x^2+y^2=r^2\right\}=S_+\cup S_-,\]
where
\[S_{\pm}=\left\{(x,y)\,\left|\, -1\leq x\leq 1, y=\pm \sqrt{r^2-x^2}\right.\right\}.\]

$S_{\pm}$ are the graphs of the functions $f_{\pm}:[-1,1]\to\mb{R}$,
\[f_{\pm}(x)=\pm\sqrt{r^2-x^2}.\]  Since $[-1,1]$ is a Jordan measurable set in $\mb{R}$, and $f_{\pm}:[-1,1]\to\mb{R}$ are continuous functions, Corollary \ref{230829_5} implies that $S_{\pm}=G_{f_{\pm}}$ have Jordan content zero. Therefore, $\mk{D}$ is Jordan measurable.
\end{example}

\begin{figure}[ht]
\centering
\includegraphics[scale=0.2]{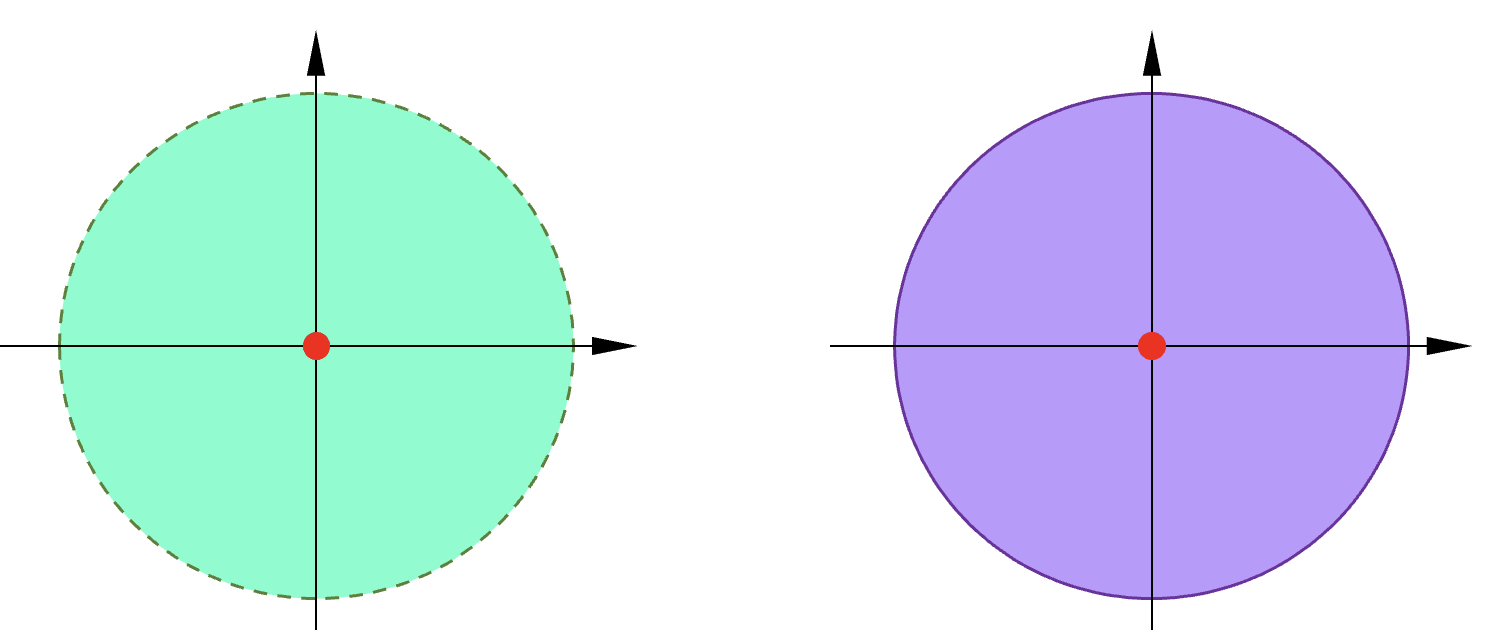}
\caption{An open ball and its closure are Jordan measurable. }\label{figure104}
\end{figure}

More generally than the open balls, we have the following.
\begin{example}[label=230829_7]{}
Let $[a,b]$ be a closed interval in $\mb{R}$, and let $f:[a,b]\to\mb{R}$ and $g:[a,b]\to\mb{R}$ be continuous functions satisfying $f(x)\leq g(x)$ for all $x\in [a,b]$. Define $\mk{D}_1$ and $\mk{D}_2$ to be the sets
\begin{align*}\mk{D}_1=&\left\{(x,y)\,|\, a<x<b, f(x)<y<g(x)\right\},\\
\mk{D}_2=&\left\{(x,y)\,|\, a\leq x\leq b, f(x)\leq y\leq g(x)\right\}
\end{align*}
Show that $\mk{D}_1$ and $\mk{D}_2$ are Jordan measurable sets.
\end{example}
\begin{solution}{Solution}
Since $f:[a,b]\to\mb{R}$ and $g:[a,b]\to\mb{R}$ are continuous,  $\mk{D}_1$ is open and $\mk{D}_2$ is closed.  Since $[a,b]$ is compact, they are bounded. There is a positive number $M$ such that
\[|f(x)|\leq M\quad\hspace{1cm}\quad |g(x)|\leq M\hspace{1cm}\text{for all}\;x\in [a,b].\]
Let
\begin{align*}
S_1&=\left\{(a,y)\,|\, -M\leq y\leq M\right\},\\
S_2&=\left\{(b,y)\,|\, -M\leq y\leq  M\right\},\\
S_3&=\left\{(x,y)\,|\, a\leq x\leq b, y=f(x)\right\},\\
S_4&=\left\{(x,y)\,|\, a\leq x\leq b, y=g(x)\right\},
\end{align*}and let  $S=S_1\cup S_2\cup S_3\cup S_4$.
 Then 
 \[\mk{D}_2\setminus\mk{D}_1\subset S.\]
 Since $\mk{D}_1$ is open and $\mk{D}_1\subset \mk{D}_2$,
 \[\mk{D}_1=\text{int}\,\mk{D}_1\subset\text{int}\,\mk{D}_2.\]
 Since $\mk{D}_2$ is closed and  $\mk{D}_1\subset \mk{D}_2$,
 \[\overline{\mk{D}_1}\subset\overline{\mk{D}_2}=\mk{D}_2.\]
 \bs
 Therefore,
 \[\pa\mk{D}_1=\overline{\mk{D}_1}\setminus\text{int}\,\mk{D}_1\subset \mk{D}_2\setminus \mk{D}_1\subset S,\]
 \[\pa\mk{D}_2=\overline{\mk{D}_2}\setminus\text{int}\,\mk{D}_2\subset \mk{D}_2\setminus \mk{D}_1\subset S.\]
  Since $S_1$ and $S_2$ are line segments, they have Jordan content zero. Since $S_3$ and $S_4$ are graphs of continuous functions defined on the Jordan measurable set $[a,b]$, $S_3$ and $S_4$ also have Jordan content zero. These imply  that $S$ has Jordan content zero. Thus, $\pa\mk{D}_1$ and $\pa\mk{D}_2$ also have Jordan content zero, which imply that $\mk{D}_1$ and $\mk{D}_2$ are Jordan measurable sets.
\end{solution}

\begin{figure}[ht]
\centering
\includegraphics[scale=0.2]{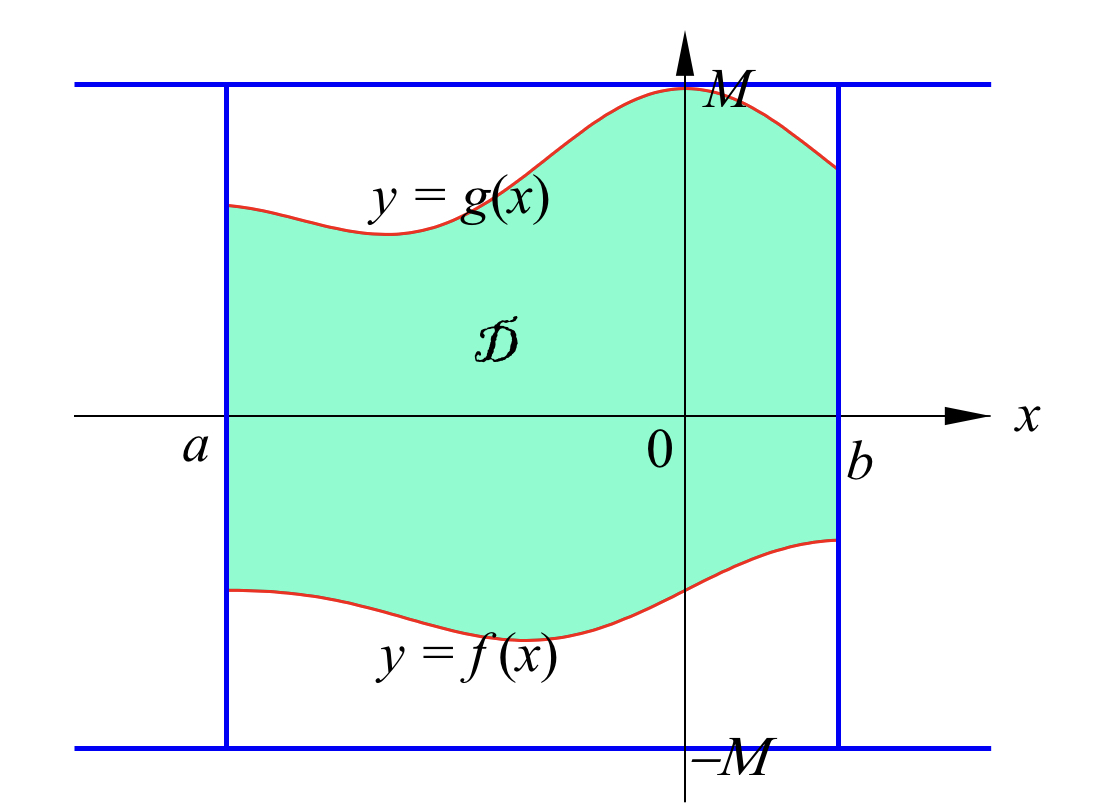}
\caption{The set $\mk{D}=\left\{(x,y)\,|\, a\leq x\leq b, f(x)\leq y\leq g(x)\right\}$ is Jordan measurable. }\label{figure105}
\end{figure}

More generally than Example \ref{230829_7}, we can prove the following.
\begin{theorem}[label=230829_9]
{}
Let $\mathcal{U}$ be a Jordan measurable set in $\mb{R}^n$, and let $f:\mathcal{U}\to\mb{R}$ and $g:\mathcal{U}\to\mb{R}$ be bounded continuous functions on $\mathcal{U}$ satisfying $f(\mf{x})\leq g(\mf{x})$ for all $\mf{x}\in\mathcal{U}$. Then the subsets

\vspace{-0.8cm}
\begin{align*}
\mk{D}_1&=\left\{(\mf{x},y)\,|\, \mf{x}\in\mathcal{U}, f(\mf{x})<y<g(\mf{x})\right\},\\
\mk{D}_2&=\left\{(\mf{x},y)\,|\, \mf{x}\in\mathcal{U}, f(\mf{x})\leq y\leq g(\mf{x})\right\}
\end{align*} of $\mb{R}^{n+1}$ are Jordan measurable.
\end{theorem}
\begin{myproof}{Sketch  of Proof}
Let $M$ be a positive number such that 
\[|f(\mf{x})|\leq M\quad\text{and}\quad |g(\mf{x})|\leq M\hspace{1cm}\; \text{for all} \;\mf{x}\in \mathcal{U}.\] The sets $\pa\mk{D}_1$ and $\pa\mk{D}_2$ are contained in the set $S=S_1\cup S_2\cup S_3$, where 

\vspace{-0.8cm}
\begin{align*}
S_1&=\left\{(\mf{x}, y)\,|\,\mf{x}\in\pa\,\mathcal{U}, -M\leq y\leq M\right\},\\
S_2&=\left\{(\mf{x}, y)\,|\,\mf{x}\in \mathcal{U}, y=f(\mf{x})\right\},\\
S_3&=\left\{(\mf{x}, y)\,|\,\mf{x}\in  \mathcal{U}, y=g(\mf{x})\right\}.
\end{align*}The sets $S_1$, $S_2$ and $S_3$   have  Jordan content zero. 
\end{myproof}

\begin{example}{}
We claim that an open ball $B(\mf{x}_0, r)$ in $\mb{R}^n$ and its closure are Jordan measurable sets. It is sufficient to consider the case where $\mf{x}_0=\mf{0}$ and $r=1$. Let $B^n=B(\mf{0}, 1)$. We will show that $B^n$ is  Jordan measurable by induction on $n$. Then $\overline{B^n}$ is also Jordan measurable.

When $n=1$, $B^1=(-1,1)$ is an interval whose boundary is the two point set $\{-1,1\}$ which has Jordan content zero.  For $n\geq 1$, assume that $B^{n}$ is a Jordan measurable subset of $\mb{R}^n$. Notice that 
\[B^{n+1}=\left\{(\mf{x}, y)\,|\, \mf{x}\in B^n, f_-(\mf{x})<y<f_+(\mf{x})\right\},\]
where  $f_{\pm}:B^{n}\to\mb{R}$ are bounded continuous functions defined by
\[f_{\pm}(x_1, \ldots, x_n)=\pm   \sqrt{1-x_1^2-\cdots-x_{n}^2}.\]
By inductive hypothesis, 
  $ B^{n}$ is Jordan measurable. By Theorem \ref{230829_9}, $B^{n+1}$ is also Jordan measurable.

\end{example}

Now we give some examples of Riemann integrable functions.
\begin{example}[label=230830_3]{}
Let $\mk{D}=\left\{(x,y,z)\,|\,x^2+y^2<4, -3\leq z\leq 3\right\}$, and let $f:\mk{D}\to\mb{R}$ be the function defined as
\[f(x,y,z)= x^2+4y^2+9z^2.\]
Explain why $f:\mk{D}\to\mb{R}$ is Riemann integrable. 
\end{example}
\begin{solution}{Solution}
The set $\mathcal{U}=\left\{(x,y,z)\,|\,x^2+y^2<4\right\}$ is an open ball. Hence, it is Jordan measurable. The functions $g_{\pm}:\mathcal{U}\to\mb{R}$, $g_{\pm}(x,y)=\pm 3$ are continuous functions. Theorem \ref{230829_9} implies that $\mk{D}$ is Jordan measurable. The function $f(x,y,z)=x^2+4y^2+9z^2$ is a polynomial. Hence, it is continuous on $\overline{\mk{D}}$, and hence, it is bounded and continuous on $\mk{D}$. Therefore,  $f:\mk{D}\to\mb{R}$ is Riemann integrable. 
\end{solution}

\begin{figure}[ht]
\centering
\includegraphics[scale=0.16]{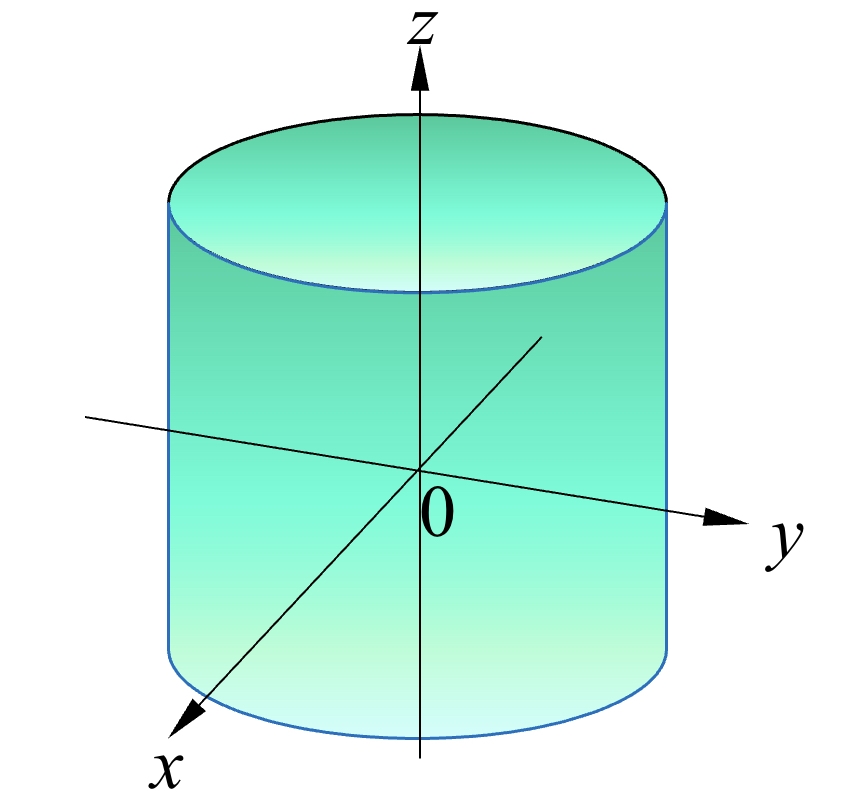}
\caption{The Jordan measurable set in Example \ref{230830_3}. }\label{figure107}
\end{figure}

\begin{example}[label=230830_4]{}
Let $\mk{D}=\left\{(x,y)\,|\,x^2+y^2<1\right\}$, and let $f:\mk{D}\to\mb{R}$ be the function  
\[f(x,y)=\begin{cases} x^2,\quad &\text{if}\; x<y,\\
y^2+1,\quad &\text{if}\;  x\geq y.\end{cases}\]
Explain why $f:\mk{D}\to\mb{R}$ is Riemann integrable. 
\end{example}
\begin{solution}{Solution}
The set $\mk{D}$ is an open ball. Hence, it is Jordan measurable. The set of discontinuities of the function $f:\mk{D}\to\mb{R}$ is contained in the line segment $L$ from the point $(-1,-1)$ to the point $(1,1)$. Since $L$ has Jordan content zero, $f:\mk{D}\to\mb{R}$ is Riemann integrable.
\end{solution}

\begin{figure}[ht]
\centering
\includegraphics[scale=0.18]{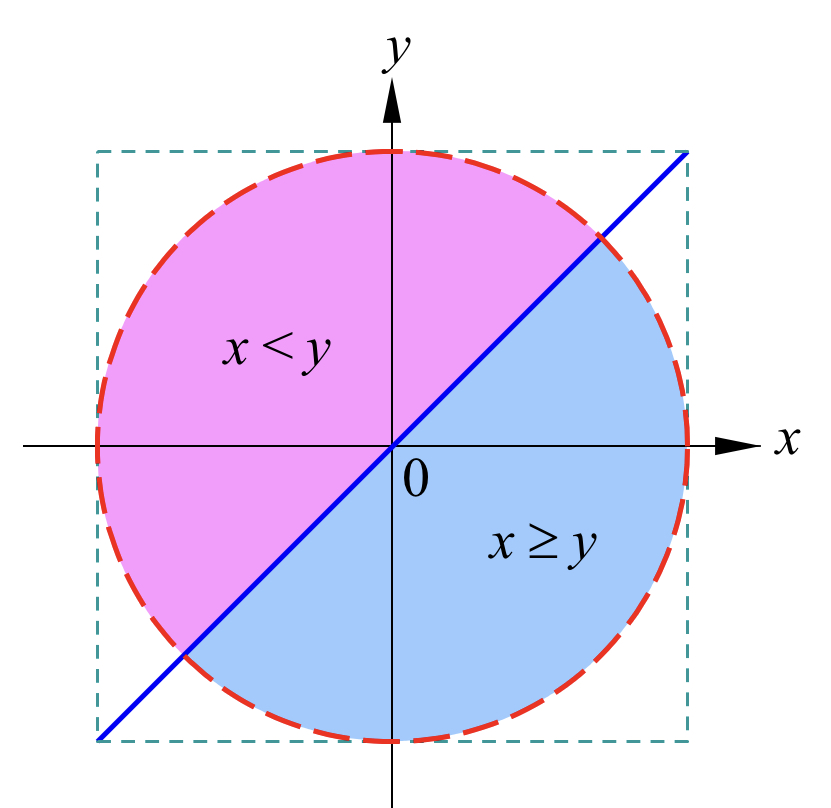}
\caption{The Jordan measurable set  in Example \ref{230830_4}. }\label{figure106}
\end{figure}

At the end of this section, let us  prove a few interesting theorems. The next theorem shows that a set with Jordan content zero cannot have nonempty interior.
\begin{theorem}{}
Let $\mk{D}$ be a  subset of $\mb{R}^n$ that has Jordan content zero.  Then $\text{int}\,\mk{D}=\emptyset$.

\end{theorem}
\begin{myproof}{Proof}
  
We use proof by contradiction. 
If  $\text{int}\, \mk{D}\neq\emptyset$,
it is an open set that contains at least one point $\mf{u}=(u_1, \ldots, u_n)$. By definition of interior points, there exists $r>0$ such that $B(\mf{u}, r)\subset \mk{D}$. There exists $\delta>0$ such that the rectangle $\mf{I}_{\delta}=\di\prod_{i=1}^n[u_i-\delta, u_i+\delta]$ is   contained in $B(\mf{u}, r)$, and hence in $\mk{D}$. 

Since $\mk{D}$ has Jordan content zero, we find that $\mf{I}_{\delta}$ also has Jordan content zero. But 
\[\text{vol}\,(\mf{I}_{\delta})=(2\delta)^n>0.\]
This gives a contradiction.

\end{myproof}

The next one is the mean value theorem for integrals.
\begin{theorem}{Mean Value Theorem for Integrals}
Let $\mk{D}$ be a closed and bounded Jordan measurable set in $\mb{R}^n$, and let $f:\mk{D}\to\mb{R}$ be a continuous function. If $\mk{D}$ is connected or path-connected, then there is a point $\mf{x}_0$ in $\mk{D}$ such that 
\[\int_{\mk{D}}f(\mf{x})d\mf{x}=f(\mf{x}_0)\,\text{vol}\,(\mk{D}).\]
\end{theorem}
\begin{myproof}{Proof}

Since $\mk{D}$ is compact and $f:\mk{D}\to\mb{R}$ is continuous, extreme value theorem asserts that there exist points $\mf{u}$ and $\mf{v}$ in $\mk{D}$ such that
\[f(\mf{u})\leq f(\mf{x})\leq f(\mf{v})\hspace{1cm}\text{for all}\;\mf{x}\in\mk{D}.\]

Since $\mk{D}$ is a Jordan measurable set and $f:\mk{D}\to\mb{R}$ is continuous, $f:\mk{D}\to\mb{R}$ is Riemann integrable. The monotonicity theorem implies that
\[f(\mf{u})\,\text{vol}\,(\mk{D})\leq \int_{\mk{D}}f(\mf{x})d\mf{x}\leq f(\mf{v})\,\text{vol}\,(\mk{D}).\]
 
If $\text{vol}\,(\mk{D})=0$, we can take $\mf{x}_0$ to be any point in $\mk{D}$. If $\text{vol}\,(\mk{D})\neq 0$, notice that
\[c=\frac{1}{\text{vol}\,(\mk{D})}\int_{\mk{D}}f(\mf{x})d\mf{x}\]
satisfies
\[f(\mf{u})\leq c\leq f(\mf{v}).\]
Since $\mk{D}$ is connected or path-connected, and  $f:\mk{D}\to\mb{R}$ is continuous, intermediate value theorem asserts that $f(\mk{D})$ must be an interval. Since $f(\mf{u})$ and $f(\mf{v})$ are in $f(\mk{D})$ and $c$ is in between them, $c$ must also be in $f(\mk{D})$. This means that there is an $\mf{x}_0$ in $\mk{D}$ such that
\[\frac{1}{\text{vol}\,(\mk{D})}\int_{\mk{D}}f(\mf{x})d\mf{x}=c=f(\mf{x}_0).\]

\end{myproof}

\vp
\noindent
{\bf \large Exercises  \thesection}
\setcounter{myquestion}{1}

\begin{question}{\themyquestion}
Explain why the set
\[\mk{D}=\left\{(x,y,z)\,|\, 4x^2+y^2+9z^2<36\right\}\] is Jordan measurable.
\end{question}
\atc
 \begin{question}{\themyquestion}
Explain why the set
\[\mk{D}=\left\{(x,y,z)\,|\, x\geq 0, y\geq 0, z\geq 0, 3x+4y+6z\leq 12\right\}\] is Jordan measurable.
\end{question}
\atc

\begin{question}{\themyquestion}
Let 
$\di \mk{D}=\left\{(x,y,z)\,|\, x^2+y^2+z^2= 25\right\}$, and let $f:\mk{D}\to\mb{R}$ be the function defined as \[f(x,y,z)=\begin{cases} 1,\quad &\text{if $x$, $y$ and $z$ are rational},\\0,\quad &\text{otherwise}.\end{cases}\] Explain why $f:\mk{D}\to\mb{R}$ is Riemann integrable and find $\di \int_{\mk{D}} f$.
\end{question}
\atc

\begin{question}{\themyquestion}
Let 
$\di \mk{D}=\left\{(x,y,z)\,|\, x^2+y^2+z^2\leq 25\right\}$, and let $f:\mk{D}\to\mb{R}$ be the function defined as \[f(x,y,z)=ze^{|xy|}.\] Explain why $f:\mk{D}\to\mb{R}$ is Riemann integrable.
\end{question}
\atc
\begin{question}{\themyquestion}
Let 
$\di \mk{D}=[0,2]\times (-2, 5)\times (1, 7]$, and let $f:\mk{D}\to\mb{R}$ be the function defined as \[f(x,y,z)=\begin{cases} x+y,\quad &\text{if}\;x<y+z,\\
2x-y,\quad &\text{if}\; x\geq y+z.\end{cases}\] Explain why $f:\mk{D}\to\mb{R}$ is Riemann integrable.
\end{question}
\atc

\begin{question}{\themyquestion}
Let 
$\di \mk{D}=\left\{(x,y,z)\,|\, 4x^2+9y^2\leq 36, 0 \leq z\leq x^2+y^2\right\} $, and let $f:\mk{D}\to\mb{R}$ be the function defined as \[f(x,y,z)=\begin{cases} x,\quad &\text{if}\;x<y+z,\\
y+z,\quad &\text{if}\; x\geq y+z.\end{cases}\] Explain why $f:\mk{D}\to\mb{R}$ is Riemann integrable.
\end{question}

\atc
\begin{question}{\themyquestion}
If $\mk{D}$ is a Jordan measurable set that is contained in the closed rectangle $\mf{I}$, show that $\mf{I}\setminus\mk{D}$ is also Jordan measurable. Moreover, 
\[\text{vol}\,(\mf{I}\setminus\mk{D})=\text{vol}\,(\mf{I})-\text{vol}\,(\mk{D}).\]
\end{question}
\atc

\begin{question}{\themyquestion}
If $\mk{D}_1$ and $\mk{D}_2$ are Jordan measurable sets  and $\mk{D}_2$  is contained in $\mk{D}_1$, show that $\mk{D}_1\setminus\mk{D}_2$ is also Jordan measurable.
Moreover, 
\[\text{vol}\,(\mk{D}_1\setminus\mk{D}_2)=\text{vol}\,(\mk{D}_1)-\text{vol}\,(\mk{D}_2).\]
\end{question}

\atc

\begin{question}{\themyquestion}
If $\mk{D}$ is a Jordan measurable set, show that $\text{int}\,\mk{D}$ is also Jordan measurable.
Moreover, 
\[\text{vol}\,(\text{int}\,\mk{D})=\text{vol}\,(\mk{D}).\]
\end{question}
\atc

\begin{question}{\themyquestion}
Let $\mk{D}_1$ and $\mk{D}_2$ be  Jordan measurable sets in $\mb{R}^m$ and $\mb{R}^n$ respectively. Assume that $\mk{D}_1$ has Jordan content zero. Show that the set $\mk{D}=\mk{D}_1\times\mk{D}_2$  in $\mb{R}^{m+n}$ also has Jordan content zero.
\end{question}

\atc

\begin{question}{\themyquestion}
 Let $\di\mk{D}=\left\{(x,y)\,|\, x^2+y^2\leq 9\right\}$. Show that the integral $\di \int_{\mk{D}} xdxdy$ exist and is equal to 0.
\end{question}
\atc

\begin{question}{\themyquestion}
let $\mk{D}$ be a Jordan measurable set, and let $f:\mk{D}\to\mb{R}$ be a Riemann integrable function. If $g:\overline{\mk{D}}\to\mb{R}$ is a bounded function such that $g(\mf{x})=f(\mf{x})$ for all $\mf{x}$ in $\mk{D}$, show that $g:\overline{\mk{D}}\to\mb{R}$ is Riemann integrable and
\[\int_{\overline{\mk{D}}}g=\int_{\mk{D}}f.\]
\end{question}

\section{Iterated Integrals and Fubini's Theorem} 

In Section \ref{sec6.3},  we have given a sufficient condition for a function $f:\mk{D}\to\mb{R}$ to be Riemann integrable. 
\begin{highlight}{Riemann Integrable Functions}
If $\mk{D}$ is a subset of $\mb{R}^n$ such that  a constant function on $\mk{D}$ is Riemann integrable, then any bounded function on $\mk{D}$ whose set of discontinuities is a set that has Jordan content zero is Riemann integrable. 
\end{highlight}

However, we have not discussed any   strategy to  compute a Riemann integral, except by using a sequence of partitions $\{\mf{P}_k\}$ with \[\lim_{k\to\infty}|\mf{P}_k|=0.\]  This is a practical approach if one has a computer, but it is not feasible for hand calculations. Besides, it might also be difficult for us to understand the dependence of the integral on the parameters in the integrand. When $n=1$, we have seen that the fundamental theorem of calculus gives us a powerful tool to calculate a Riemann integral when the integrand is a continuous function that has explicit antiderivatives. 
 To be able to apply this powerful tool   in the  multivariable context, we need to relate multiple integrals with iterated integrals. This is the topic that is studied in this section.
 
 As a motivation, consider a continuous function $f:[a,b]\times [c,d]\to\mb{R}$ defined on the closed rectangle $\mf{I}=[a,b]\times [c,d]$ in $\mb{R}^2$. If $\mf{P}=(P_1, P_2)$ is a partition of $\mf{I}$ with 
 \[P_1=\{x_0, x_1, \ldots, x_k\}\quad \text{and}\quad P_2=\{y_0, y_1, \ldots, y_l\},\] there are $kl$ rectangles in the partition $\mf{P}$. Denote the rectangles by    $\mf{J}_{i,j}$ with $1\leq i\leq k$ and $1\leq j\leq l$, where 
 \[\mf{J}_{i,j}=[x_{i-1}, x_i]\times [y_{j-1}, y_j].\] 
 Choose a set of intermediate points $A=\{\alpha_i\}$   for  the partition $P_1$,  and a set of intermediate points  $B=\{\beta_j\}$  for the partition $P_2$. Let 
 \[\boldsymbol{\xi}_{i,j}=(\alpha_i,\beta_j)\hspace{1cm}\text{for}\;1\leq i\leq k, 1\leq j\leq l.\]
 Then $C=\left\{\boldsymbol{\xi}_{i,j}\,|\, 1\leq i\leq k, 1\leq j\leq l\right\}$ is a choice of intermediate points for the partition $\mf{P}$. The Riemann sum $R(f,\mf{P}, C)$ is given by
 \begin{equation}\label{230830_6}
 R(f, \mf{P}, C)=\sum_{i=1}^k\sum_{j=1}^l f(\alpha_i,\beta_j)(x_i-x_{i-1})(y_j-y_{j-1}).\end{equation}
 Since it is a finite sum, it does not matter which order we perform the summation. 
   For fixed $x\in [a,b]$, let $g_x:[c,d]\to\mb{R}$ be the function \[g_x(y)=f(x,y),\hspace{1cm}y\in [c,d].\]
 If we perform the sum over $j$ in \eqref{230830_6} first, we find that
 \begin{align*}
 R(f, \mf{P}, C)&=\sum_{i=1}^k \left(\sum_{j=1}^lg_{\alpha_i}(\beta_j)(y_j-y_{j-1})\right)(x_i-x_{i-1}) \\&=\sum_{i=1}^k R(g_{\alpha_i}, P_2, B)(x_i-x_{i-1}).\end{align*}
 Since $g_{\alpha_i}:[c,d]\to\mb{R}$ is continuous, it is Riemann integrable. Therefore,
 \[\lim_{|P_2|\to 0}R(g_{\alpha_i}, P_2, B)=\int_c^d g_{\alpha_i}(y) dy.\]This prompts us to define the function $F:[a,b]\to\mb{R}$ by
 \[F(x)=\int_c^d g_x(y)dy=\int_c^df(x,y)dy.\]
 Then  
 \begin{align*}
 \lim_{|P_2|\to 0}R(f, \mf{P}, C)&= \sum_{i=1}^k\lim_{|P_2|\to 0}R(g_{\alpha_i}, P_2, B)(x_i-x_{i-1}) \\&=\sum_{i=1}^kF(\alpha_i)(x_i-x_{i-1})=R(F,P_1, A).\end{align*}
If $F:[a,b]\to\mb{R}$ is also Riemann integrable, we would have
 \begin{align*}
   \lim_{|P_1|\to 0}\lim_{|P_2|\to 0}R(f, \mf{P}, C) &=\lim_{|P_1|\to 0}R(F, P_1, A)\\&=\int_{a}^b F(x)dx=\int_a^b\left(\int_c^d f(x,y)dy\right)dx.\end{align*}
 Interchanging the roles of  of $x$ and $y$, or equivalently, summing over $i$ first in \eqref{230830_6}, we find that 
 \[  \lim_{|P_2|\to 0}\lim_{|P_1|\to 0}R(f, \mf{P}, C)  =\int_c^d\left(\int_a^b f(x,y)dx\right)dy.\]
 \begin{figure}[ht]
\centering
\includegraphics[scale=0.2]{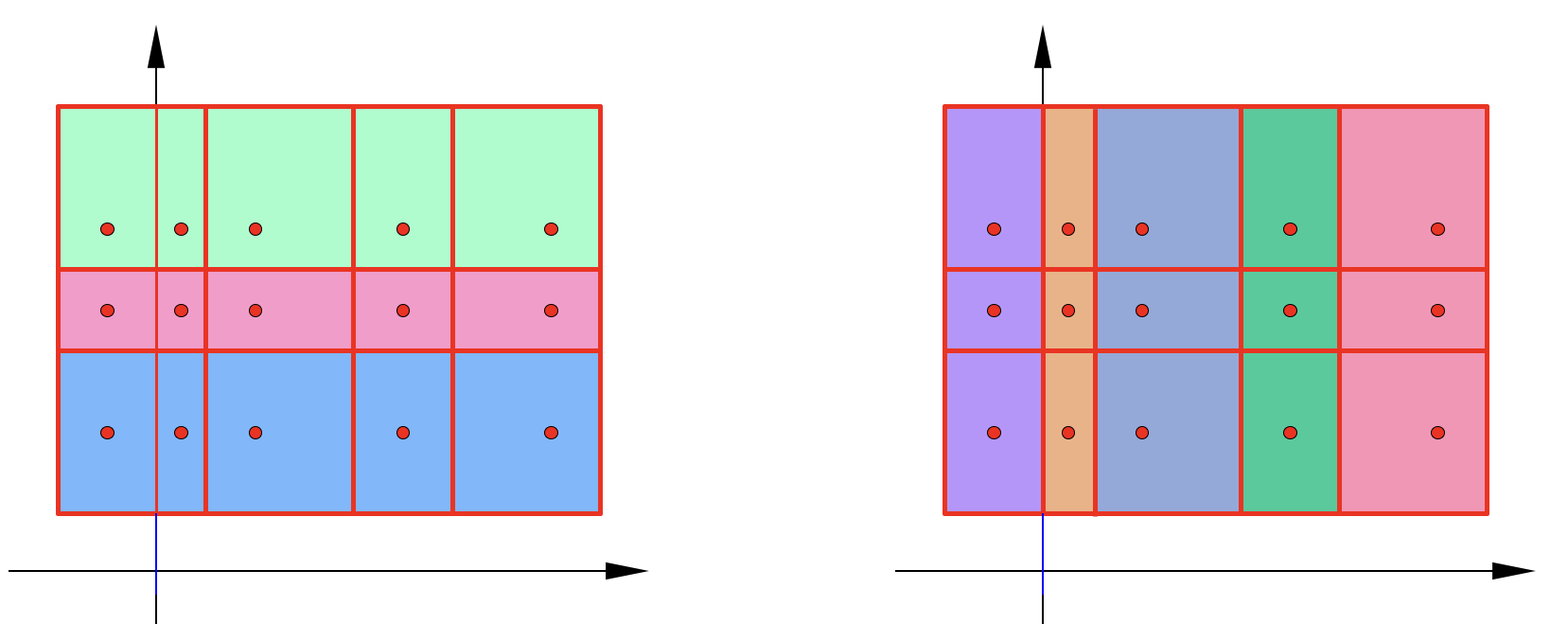}
\caption{Given a partition $\mf{P}$ of a rectangle, one can sum over the rectangles row by row, or column by column. }\label{figure108}
\end{figure}

Since 
\[|\mf{P}|=\sqrt{|P_1|^2+|P_2|^2},\]
$|\mf{P}|\to 0$ if and only if $(|P_1|,|P_2|)\to (0,0)$. The question becomes whether  the two limits
\[\lim_{|P_1|\to 0}\lim_{|P_2|\to 0}R(f, \mf{P}, C) \quad  \text{and} \quad   \lim_{|P_2|\to 0}\lim_{|P_1|\to 0}R(f, \mf{P}, C)  \] are equal; and whether they are equal to the limit \[  \lim_{(|P_1|,|P_2|)\to (0,0)}R(f, \mf{P}, C).\]

\begin{remark}{}
Let $f:\mb{R}^2\to\mb{R}$ be the function defined as
\[f(x,y)=\frac{x^2}{x^2+y^2},\hspace{1cm}(x,y)\in\mb{R}^2\setminus\{(0,0)\}.\]
\end{remark}\begin{highlight}{}
We find that
\begin{align*}
\lim_{y\to 0}\lim_{x\to 0}f(x,y)&=\lim_{y\to 0}0=0,\\
\lim_{x\to 0}\lim_{y\to 0}f(x,y)&=\lim_{x\to 0}1=1.
\end{align*}

Hence,
\[\lim_{y\to 0}\lim_{x\to 0}f(x,y)\neq \lim_{x\to 0}\lim_{y\to 0}f(x,y).\]
This example shows that we cannot simply interchange the order of limits. In fact, the limit  
\[\lim_{(x,y)\to (0,0)}f(x,y)\]
does not exist.
\end{highlight}

The integrals 
\[\int_c^d\left(\int_a^b f(x,y)dx\right)dy\quad \text{and}\quad \int_a^b\left(\int_c^d f(x,y)dy\right)dx\]
are called {\it iterated integrals}.

\begin{definition}{Iterated Integrals}
Let $n$ be a positive integer larger than 1, and let $k$ be a positive integer less than $n$. Denote a point in $\mb{R}^n$ by $(\mf{x}, \mf{y})$, where $\mf{x}\in\mb{R}^k$ and $\mf{y}\in\mb{R}^{n-k}$. Given that $\mf{I}=\di\prod_{i=1}^n[a_i,b_i]$ is a closed rectangle in $\mb{R}^n$, 
let \[\mf{I}_{\mf{x}}=\prod_{i=1}^k[a_i,b_i]\quad   \text{and}\quad \di\mf{I}_{\mf{y}}=\prod_{i=k+1}^n[a_i,b_i].\] 
If $f:\mf{I}\to\mb{R}$ is a bounded function defined on $\mf{I}$, an interated integral is an integral of the form 
\[\int_{\mf{I}_{\mf{y}}}\int_{\mf{I}_{\mf{x}}}f(\mf{x},\mf{y})d\mf{x} d\mf{y}\quad\text{or}\quad \int_{\mf{I}_{\mf{x}}}\int_{\mf{I}_{\mf{y}}}f(\mf{x},\mf{y})d\mf{y} d\mf{x},\]whenever they exist. 
\end{definition}
Let us consider the following example.

\begin{example}[label=230831_3]{}
Let $g:[a,b]\to\mb{R}$ and $h:[a,b]\to\mb{R}$ be continuous functions defined on $[a,b]$ such that $g(x)\leq h(x)$ for all $x\in [a,b]$. Consider the set $\mk{D}$ defined as
\[\mk{D}=\left\{(x,y)\,|\, a\leq x\leq b, g(x)\leq y\leq h(x)\right\}.\]
Let $\chi_{\mk{D}}:\mb{R}^2\to\mb{R}$ be the corresponding characteristic function. If 
\[c\leq g(x)\leq h(x)\leq d\hspace{1cm} \text{for all }\;x\in [a,b],\]
 then $\mf{I}=[a,b]\times [c,d]$ is a closed rectangle that contains $\mk{D}$. We have seen that $\mk{D}$ is a Jordan measurable set. Hence, $\chi_{\mk{D}}:\mf{I}\to\mb{R}$ is a Riemann integrable function. 
 For any $x\in [a,b]$,
\[\int_c^d \chi_{\mk{D}}(x,y) dy=\int_{g(x)}^{h(x)}\,dy= h(x)-g(x).\]
 
Therefore, the iterated integral $\di \int_a^b\left(\int_c^d \chi_{\mk{D}}(x,y)dy\right)dx$ is equal to
\[\int_a^b\left(\int_c^d \chi_{\mk{D}}(x,y)dy\right)dx=\int_a^b \left(h(x)-g(x)\right) dx.\]
In single variable calculus, we have learned that the integral $\di \int_a^b \left(h(x)-g(x)\right) dx$  gives the area of $\mk{D}$. Thus, in this case, we have
\[\int_a^b\left(\int_c^d \chi_{\mk{D}}(x,y)dy\right)dx=\text{vol}\,(\mk{D})=\int_{[a,b]\times [c,d]}\chi_{\mk{D}}(x,y)dxdy.\]
Namely, the iterated integral is equal to the double integral.
\end{example}

\begin{figure}[ht]
\centering
\includegraphics[scale=0.2]{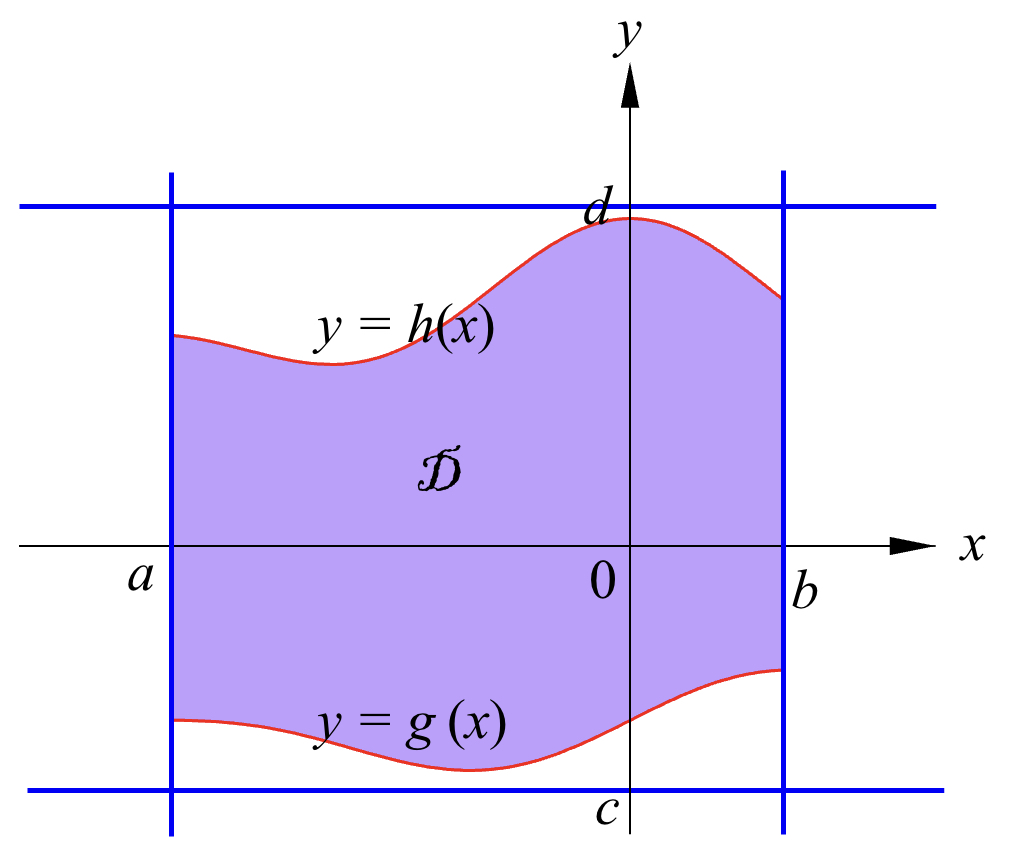}
\caption{The domain  $\mk{D}=\left\{(x,y)\,|\, a\leq x\leq b, g(x)\leq y\leq h(x)\right\}$.}\label{figure109}
\end{figure}

 The following theorem is the general case 
when $n=2$.
\begin{theorem}[label=230831_1]{Fubini's Theorem in the Plane}
Let $\mf{I}=[a,b]\times [c,d]$, and let $f:\mf{I}\to\mb{R}$ be a Riemann integrable function. For each $x\in [a,b]$, define the function $g_x:[c,d]\to\mb{R}$ by 
\[g_x(y)=f(x,y),\quad y\in [c,d].\]
If $g_x:[c,d]\to\mb{R}$ is Riemann integrable for each $x\in [a,b]$, let $F:[a,b]\to\mb{R}$ be the function defined as
\[F(x)=\int_c^d g_x(y)dy=\int_c^df(x,y)dy.\]
Then we have the following.
\begin{enumerate}[(a)]
\item  The function $F:[a,b]\to\mb{R}$ is Riemann integrable.
\item The integral of $f:[a,b]\times [c,d]\to\mb{R}$ is equal to the integral of $F:[a,b]\to\mb{R}$. Namely, 
\[\int_{[a,b]\times [c,d]}f=\int_a^bF.\]
Equivalently,
\[\int_{[a,b]\times [c,d]}f(x,y)dxdy=\int_a^b \int_c^df(x,y)dydx.\]
\end{enumerate}\end{theorem}

\begin{myproof}{Proof}
Let
\[I=\int_{[a,b]\times [c,d]}f(x,y)dxdy.\]
We will show that for any $\varepsilon>0$, there exists $\delta>0$ such that if $P$ is a partition of $[a,b]$ with $|P|<\delta$, and  $A=\{\alpha_i\}$ is any set of intermediate points for $P$, then 
\[\left|R(F, P, A)-I\right|<\varepsilon.\]This will prove  both (a) and (b).

Fixed $\varepsilon>0$. Since $f:[a,b]\times [c,d]\to\mb{R}$ is Riemann integrable with integral $I$, there exists $\delta_0>0$ such that if $\mf{P}=(P_1, P_2)$  is a partition of $[a,b]\times [c,d]$ with $|\mf{P}|<\delta_0$,  then
\[U(f,\mf{P})-L(f,\mf{P})<\varepsilon.\]
Take $\delta=\delta_0/2$. Let $P=\{x_0, x_1, \ldots, x_k\}$ be a partition of $[a,b]$ with $|P|<\delta$. Take any partition $P_2=\{y_0, y_1, \ldots, y_l\}$ of $[c,d]$ such that $|P_2|<\delta$. Let $\mf{P}=(P_1, P_2)$, where $P_1=P$. Then $|\mf{P}|<\sqrt{2}\delta<\delta_0$.  
For $1\leq i\leq k, 1\leq j\leq l$, let 
\begin{align*}m_{i,j}&=\inf_{(x, y)\in [x_{i-1}, x_i]\times [y_{j-1}, y_j]}f(x,y),\\
 M_{i,j}&=\sup_{(x, y)\in [x_{i-1}, x_i]\times [y_{j-1}, y_j]}f(x,y).\end{align*}
Then
\begin{align*}
L(f,\mf{P})&=\sum_{i=1}^k\sum_{j=1}^l m_{i,j}(x_i-x_{i-1})(y_j-y_{j-1}),\\
U(f,\mf{P})&=\sum_{i=1}^k\sum_{j=1}^l M_{i,j}(x_i-x_{i-1})(y_j-y_{j-1}).
\end{align*}
Now let $A=\{\alpha_i\}$ be any choice of intermediate points for the partition $P=P_1$. Notice that
for any  $1\leq i\leq k$, additivity theorem says that
\[  F(\alpha_i)=\sum_{j=1}^l\int_{y_{j-1}}^{y_j}g_{\alpha_i}(y)dy=\sum_{j=1}^l\int_{y_{j-1}}^{y_j}f(\alpha_i,y)dy.\]
\bp
Since
\[m_{i,j}\leq f(\alpha_i, y)\leq M_{i,j}\hspace{1cm}\text{for all}\, y\in [y_{j-1}, y_j],\]
we find that for  $1\leq j\leq l$,  
\[m_{i,j}(y_i-y_{j-1})\leq \int_{y_{j-1}}^{y_j}f(\alpha_i,y)dy\leq M_{i,j}(y_i-y_{j-1}).\]
It follows that
\[\sum_{j=1}^l m_{i,j}(y_j-y_{j-1})\leq F(\alpha_i)\leq \sum_{j=1}^l M_{i,j}(y_j-y_{j-1}).\]
Multiply by $(x_i-x_{i-1})$, and sum over $i$ from 1 to $k$, we find that
\[
L(f,\mf{P})\leq \sum_{i=1}^kF(\alpha_i)(x_i-x_{i-1})\leq U(f,\mf{P}).
\]
In other words,
\[L(f,\mf{P})\leq R(F,P, A)\leq U(f,\mf{P}).\]
Since we also have
\[L(f,\mf{P})\leq  I\leq U(f,\mf{P}),\] we find that
\[\left|R(F, P, A)-I\right|\leq U(f,\mf{P})-L(f,\mf{P})<\varepsilon.\]This completes the proof.
\end{myproof}
\begin{example}{}
Evaluate the integral $\di\int_{[0,1]\times [0,1]}x\sin(xy)dxdy$.
\end{example}
\begin{solution}{Solution}
The function $f:[0,1]\times [0,1]\to\mb{R}$, $f(x,y)=x\sin(xy)$ is a continuous function. Hence, it is Riemann integrable. 

\bs
For each $x\in [0,1]$, the function $g_x:[0,1]\to\mb{R}$, $g_x(y)=x\sin (xy)$ is also continuous. Hence, $g_x:[0,1]\to\mb{R}$ is Riemann integrable.
By Fubini's theorem,
\begin{align*}\int_{[0,1]\times [0,1]}x\sin(xy)dxdy&=\int_0^1\int_0^1x\sin(xy)dydx\\&=\int_0^1\left[-\cos(xy)\right]_{y=0}^{y=1}dx
\\&=\int_0^1 \left(1-\cos x\right)dx\\&=1-\left[\sin x\right]_0^1=1-\sin 1.
\end{align*}
\end{solution}

The roles of $x$ and $y$ in Fubini's theorem can be interchanged, and we obain the following.
\begin{corollary}[label=230831_2]{}
Assume that $f:[a,b]\times [c,d]\to\mb{R}$ is a Riemann integrable function such that for each $x\in [a,b]$, the function $g_x:[c,d]\to\mb{R}$, $g_x(y)=f(x,y)$ is Riemann integrable; and for each $y\in [c,d]$, the function $h_y:[a,b]\to\mb{R}$, $h_y(x)=f(x,y)$ is Riemann integrable. Then we can interchange the order of integration. Namely,
\[ \int_{c}^d\int_a^b f(x,y)dxdy=\int_a^b \int_c^df(x,y)dydx.\]
\end{corollary}
\begin{example}{}
If we evaluate the iterated integral $\di \int_0^1\int_0^1x\sin(xy)dx dy$ directly, it would be quite tedious as we need to apply integration by parts to evaluate the integral $\di\int_0^1 x\sin(xy)dx$. Using Corollary \ref{230831_2}, we can interchange the order of integration and obtain
\[ \int_0^1\int_0^1x\sin(xy)dx dy= \int_0^1\int_0^1x\sin(xy)dy dx=1-\sin 1.\]
\end{example}

\begin{remark}{}
The assumption that $f:[a,b]\times [c,d]\to\mb{R}$ is Riemann integrable is essential in Fubini's theorem. It does not follow from the fact that for each $x\in [a,b]$, the function $g_x:[c,d]\to\mb{R}$ is Riemann integrable, and the function $F:[a,b]\to\mb{R}$, 
\[F(x)=\int_c^dg_x(y)dy\] is Riemann integrable. For example, let $g:[-1,1]\to\mb{R}$ and $h:[-1,1]\to\mb{R}$ be the functions defined as
\[g(x)=\begin{cases} 1,\quad &\text{if $x$ is rational}, \\-1,\quad &\text{if $x$ is irrational},\end{cases}\hspace{1cm}
 h(y)=\begin{cases} 1,\quad &\text{if }\;y\geq 0, \\-1,\quad &\text{if }\;y<0.\end{cases}.\]
Then define the function $f:[-1,1]\times [-1,1]\to \mb{R}$ by
\[f(x,y)=g(x)h(y).\]Since $h:[-1,1]\to\mb{R}$ is a step function, it is Riemann integrable and 
\[\int_{-1}^1h(y)dy=\int_{-1}^0h(y)dy+\int_0^1h(y)dy=0.\] Hence,
for fixed $x\in [-1,1]$,  
\[\int_{-1}^1f(x,y)dy=0.\]
Thus, the function $F:[-1,1]\to\mb{R}$, 
\[F(x)=\int_{-1}^1f(x,y)dy,\]being a function that is always zero, is  Riemann integrable with integral 0. It follows that
\[\int_{-1}^1\int_{-1}^1 f(x,y)dydx=0.\]However, one can prove  that the function $f:[-1,1]\times [-1,1]\to \mb{R}$ is not Riemann integrable, using the same way that we show  that a Dirichlet's function is not Riemann integrable.
\end{remark}
\begin{highlight}{}
The fact that for each $x\in [a,b]$, the function $g_x:[c,d]\to\mb{R}$, $g_x(y)=f(x,y)$ is Riemann integrable also does not follow from the fact that $f:[a,b]\times [c,d]\to\mb{R}$ is Riemann integrable. Consider for example the function 
$f:[-1,1]\times [0,1]\to\mb{R}$,
\[f(x,y)=\begin{cases} 1,\quad &\text{if $x=0$ and $y$ is rational},\\0,\quad &\text{otherwise}.\end{cases}\]
Then the set of discontinuities $\mathcal{N}$ of $f:[-1,1]\times [0,1]\to\mb{R}$ is the line segment between the point $(0,0)$ and the point $(0,1)$.  Hence, $\mathcal{N}$ has Jordan content 0. Therefore, $f:[-1,1]\times [0,1]\to\mb{R}$ is Riemann integrable. For $x=0$, $g_0:[0,1]\to\mb{R}$ is the Dirichlet's function. Hence, $g_0:[0,1]\to\mb{R}$ is not Riemann integrable.

\end{highlight}

Now we consider the  case depicted in Example \ref{230831_3} for more general functions. 
\begin{theorem}[label=230831_5]{}
Let $g:[a,b]\to\mb{R}$ and $h:[a,b]\to\mb{R}$ be continuous functions defined on $[a,b]$ such that $g(x)\leq h(x)$ for all $x\in [a,b]$, and let $\mk{D}$ be the set  
\[\mk{D}=\left\{(x,y)\,|\, a\leq x\leq b, g(x)\leq y\leq h(x)\right\}.\] If $f:\mk{D}\to\mb{R}$ is a continuous function, then
\[\int_{\mk{D}}f(x,y)dxdy=\int_a^b\int_{g(x)}^{h(x)}f(x,y)dy dx.\]
\end{theorem}
\begin{myproof}{Proof}

Since  $g$ and $h$ are continuous functions, there exist numbers $c$ and $d$ such that
\[c\leq g(x)\leq h(x)\leq d\hspace{1cm}\text{for all}\; x\in [a,b].\] Then $\mf{I}=[a,b]\times [c,d]$ be a closed rectangle that contains $\mk{D}$.
We have shown before that $\mk{D}$ is a Jordan measurable set and  $f:\mk{D}\to\mb{R}$ is Riemann integrable.  Therefore, $\check{f}:\mf{I}\to\mb{R}$ is  Riemann integrable. 

\bp
On the other hand, for each $x\in [a,b]$, 
the function $g_x:[c,d]\to\mb{R}$ is a piecewise continuous function given by
\[g_x(y)=\begin{cases}0,\quad &\text{if}\; c\leq y<g(x),\\f(x,y),\quad &\text{if}\;g(x)\leq y\leq h(x),\\0,\quad &\text{if}\; h(x)<y\leq d.\end{cases}\]
Hence, $g_x:[c,d]\to\mb{R}$ is Riemann integrable and for $x\in [a,b]$,
\[F(x)=\int_c^dg_x(y)dy=\int_{g(x)}^{h(x)}f(x,y)dy.\]
 By Fubini's theorem in the plane, the function $F:[a,b]\to\mb{R}$ is Riemann integrable, and 
\[\int_a^b\int_{g(x)}^{h(x)}f(x,y)dydx=\int_a^b F(x)dx=\int_{\mf{I}}\check{f}(x,y)dxdy=\int_{\mk{D}}f(x,y)dxdy.\]
\end{myproof}

Again, the roles of $x$ and $y$ in Theorem \ref{230831_5} can be interchanged. Let us look at the following example.

\begin{example}[label=230831_6]{}
Let
$\mk{D}$ be the region in the plane bounded between the curve $y^2=x$ and the line $L$ between the points $(1,1)$ and $(4, -2)$. Evaluate the integral $\di\int_{\mk{D}}ydxdy$. 
\end{example}
\begin{solution}{Solution}
The  equation of the line $L$ is $x+y=2$. Hence, 
\[\mk{D}=\left\{(x,y)\,|\, -2\leq y\leq 1, y^2\leq x\leq 2-y\right\}.\]
\bs
Using Fibini's theorem, we find that
\begin{align*}
\int_{\mk{D}}ydxdy&=\int_{-2}^1 \int_{y^2}^{2-y} ydxdy=\int_{-2}^1 y(2-y-y^2)dy\\&=\int_{-2}^1(2y-y^2-y^3)dy
=\left[y^2-\frac{y^3}{3}-\frac{y^4}{4}\right]_{-2}^1=-\frac{9}{4}.
\end{align*}
\end{solution}

\begin{figure}[ht]
\centering
\includegraphics[scale=0.2]{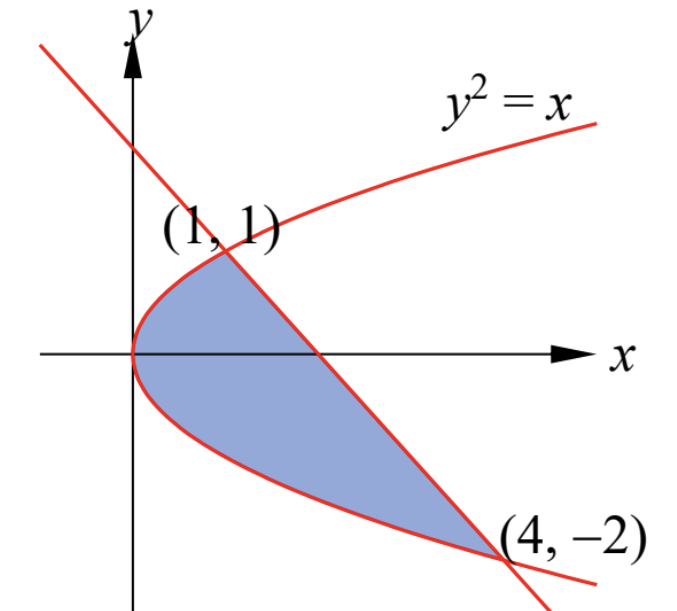}
\caption{The domain $\mk{D}=\left\{(x,y)\,|\, -2\leq y\leq 1, y^2\leq x\leq 2-y\right\}$.}\label{figure110}
\end{figure}
In Example \ref{230831_6}, it will be harder if one prefers to integrate over $y$ first. 

\begin{example}{}
Let
\[\mk{D}=\left\{(x,y)\,|\, x\geq 0, y\geq 0, 4x^2+9y^2\leq 36\right\}.\]Evaluate the integral $\di\int_{\mk{D}}xdxdy$. 
\end{example}
\begin{solution}{Solution}
The set $\mk{D}$ can be expressed in two different ways.
\begin{align*}
\mk{D}&=\left\{(x,y)\,\left|\, 0\leq x\leq 3, 0\leq y\leq \frac{1}{3}\sqrt{36-4x^2}\right.\right\},\\
\mk{D}&=\left\{(x,y)\,\left|\, 0\leq y\leq 2, 0\leq x\leq \frac{1}{2}\sqrt{36-9y^2}\right.\right\}.
\end{align*}
\bs
The function $f:\mk{D}\to\mb{R}$, $f(x,y)=x$ is continuous. Hence, the integral $\di\int_{\mk{D}}xdxdy$ is equal to iterated integrals, which we can integrate with respect to $x$ first, or with respect to $y$ first. If we integrate with respect to $y$ first, we find that
\begin{align*}
 \int_{\mk{D}}xdxdy&=\int_0^3\int_0^{\frac{1}{3}\sqrt{36-4x^2}}xdydx=\frac{1}{3}\int_0^3x \sqrt{36-4x^2}dx.
\end{align*}This integral needs to be computed using integration by substitution.
If we integrate over $x$ first, we find that
\begin{align*}
 \int_{\mk{D}}xdxdy&=\int_0^2\int_0^{\frac{1}{2}\sqrt{36-9y^2}}xdxdy\\
 &=\frac{1}{2}\left[x^2\right]_{y=0}^{y=\frac{1}{2}\sqrt{36-9y^2}}=\frac{1}{8}\int_0^2(36-9y^2) dy.
\end{align*}This integral can be easily evaluated to give
\[ \int_{\mk{D}}xdxdy=\frac{1}{8}\left[36y-3y^3\right]_0^2=6.\]
\end{solution}
\begin{figure}[ht]
\centering
\includegraphics[scale=0.19]{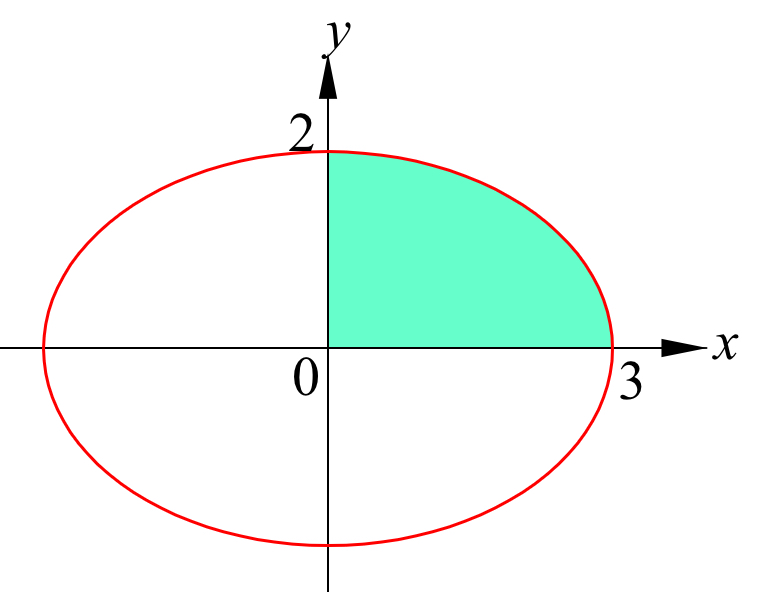}
\caption{The domain $ \mk{D}=\left\{(x,y)\,|\, x\geq 0, y\geq 0, 4x^2+9y^2\leq 36\right\}$.}\label{figure111}
\end{figure}

Now let us generalize the Fubini's theorem to arbitrary positive integer $n$ that is larger than 1.

\begin{theorem}[label=230831_7]{Fubini's Theorem}
Let $n$ be a positive integer larger than 1, and let $k$ be a positive integer less than $n$. Given that $\mf{I}=\di\prod_{i=1}^n[a_i, b_i]$ is a closed rectangle in $\mb{R}^n$, let $\mf{I}_{\mf{x}}=\di\prod_{i=1}^k[a_i, b_i]$ and $\mf{I}_{\mf{y}}=\di\prod_{i=k+1}^n[a_i, b_i]$. Assume that $f:\mf{I}\to\mb{R}$ is a Riemann integrable function such that for each $\mf{x}$ in $\mf{I}_{\mf{x}}$,   the function $g_{\mf{x}}:\mf{I}_{\mf{y}}\to\mb{R}$ defined  by 
\[g_{\mf{x}}(\mf{y})=f(\mf{x},\mf{y}),\quad \mf{y}\in\mf{I}_{\mf{y}} \]
  is Riemann integrable.  Let $F:\mf{I}_{\mf{x}}\to\mb{R}$ be the function defined as
\[F(\mf{x})=\int_{\mf{I}_{\mf{y}}}g_{\mf{x}}(\mf{y})d\mf{y}=\int_{\mf{I}_{\mf{y}}}f(\mf{x},\mf{y})d\mf{y}.\]
Then we have the followings.
\begin{enumerate}[(a)]
\item  The function $F:\mf{I}_{\mf{x}}\to\mb{R}$ is Riemann integrable.
\item The integral of $f:\mf{I}\to\mb{R}$ is equal to the integral of $F:\mf{I}_{\mf{x}}\to\mb{R}$. Namely, 
\[\int_{\mf{I}}f(\mf{x},\mf{y})d\mf{x}d\mf{y}=\int_{\mf{I}_{\mf{x}}}\int_{\mf{I}_{\mf{y}}}f(\mf{x},\mf{y})d\mf{y}d\mf{x}.\]
\item For each $\mf{y}$ in $\mf{I}_{\mf{y}}$, define the function $h_{\mf{y}}:\mf{I}_{\mf{x}}\to\mb{R}$ by 
\[h_{\mf{y}}(\mf{x})=f(\mf{x},\mf{y}),\quad \mf{x}\in\mf{I}_{\mf{x}}.\] If the function $h_{\mf{y}}:\mf{I}_{\mf{x}}\to\mb{R}$ is Riemann integrable for each $\mf{y}\in\mf{I}_{\mf{y}}$, then we can interchange the order of integration. Namely,
\[\int_{\mf{I}_{\mf{x}}}\int_{\mf{I}_{\mf{y}}}f(\mf{x},\mf{y})d\mf{y}d\mf{x}=\int_{\mf{I}_{\mf{y}}}\int_{\mf{I}_{\mf{x}}}f(\mf{x},\mf{y})d\mf{x}d\mf{y}.\]
\end{enumerate}\end{theorem}

 The proof is similar to the $n=2$ case and we leave it to the readers.

A useful case is the following which generalizes Theorem \ref{230831_5}. 
\begin{theorem}[label=230901_1]{}
Let  $\mathcal{U}$ be a Jordan measurable set in $\mb{R}^{n-1}$, and let $g:\mathcal{U}\to\mb{R}$ and $h:\mathcal{U}\to\mb{R}$ be bounded continuous functions on $\mathcal{U}$ satisfying $g(\mf{x})\leq h(\mf{x})$ for all $\mf{x}\in\mathcal{U}$. Consider the subset $\mk{D}$ of $\mb{R}^n$ defined as
\[
\mk{D} =\left\{(\mf{x},y)\,|\, \mf{x}\in\mathcal{U}, g(\mf{x})\leq y\leq h(\mf{x})\right\}.\]
If $f:\mk{D}\to\mb{R}$ is a bounded continuous function, then it is Riemann integrable, and
\[\int_{\mk{D}}f(\mf{x}, y)d\mf{x}dy=\int_{\mathcal{U}}\int_{g(\mf{x})}^{h(\mf{x})}f(\mf{x}, y)dyd\mf{x},\]
\end{theorem}

Let us look at an example.
\begin{example}{}
Evaluate the integral $\di\int_{\mathcal{S}}xdxdydz$, where $\mathcal{S}$ is the solid bounded between the plane $x+y+z=1$ and the three coordinate planes. Then find the integral $\di\int_{\mathcal{S}}(x+5y+3z)dxdydz$.
\end{example}

\begin{solution}{Solution}
The solid $\mathcal{S}$  can be expressed as
\[\mathcal{S}=\left\{(x,y,z)\,|\, (x,y)\in\mk{D}, 0\leq z\leq 1-x-y\right\},\]where\[\mk{D}=\left\{(x,y)\,|\,0\leq x\leq 1, 0\leq y\leq 1-x\right\}.\]Since $\mk{D}$ is a triangle, it is a Jordan measurable set. The function $f(x,y,z)=x$ is continuous. 
\bs
Hence, we can apply Theorem \ref{230901_1}.
\begin{align*}\int_{\mathcal{S}}xdxdydz
&=\int_{\mk{D}}\left(\int_{0}^{1-x-y}xdz \right)dxdy\\
&=\int_0^1\int_{0}^{1-x}x(1-x-y)dydx\\
&=\int_0^1 x\left[(1-x)y-\frac{y^2}{2}\right]_0^{1-x} dx\\
&=\frac{1}{2}\int_0^1x(1-x)^2dx=\frac{1}{2}\int_0^1x^2(1-x)dx\\
&=\frac{1}{2}\left(\frac{1}{3}-\frac{1}{4}\right)=\frac{1}{24}.
\end{align*}
Since the solid $\mathcal{S}$ is symmetric in $x$, $y$ and $z$, we have
\[\int_{\mathcal{S}}xdxdydz=\int_{\mathcal{S}}ydxdydz=\int_{\mathcal{S}}zdxdydz.\]
Therefore,
\[\int_{\mathcal{S}}(x+5y+3z)dxdydz=9\int_{\mathcal{S}}xdxdydz=\frac{3}{8}.\]
\end{solution}

\begin{figure}[ht]
\centering
\includegraphics[scale=0.18]{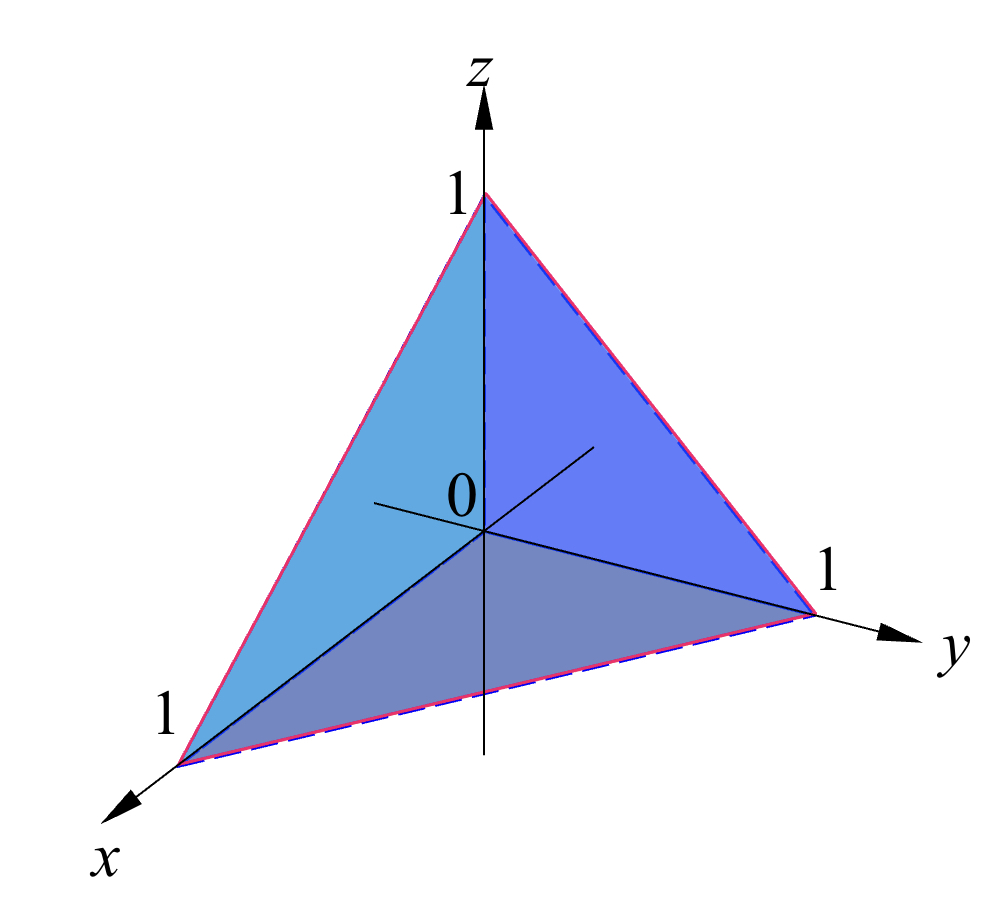}
\caption{The solid $\mathcal{S}$ bounded between the plane $x+y+z=1$ and the three coordinate planes.}\label{figure112}
\end{figure}
\vp
\noindent
{\bf \large Exercises  \thesection}
\setcounter{myquestion}{1}
\begin{question}{\themyquestion}
 Let $\mf{I}=[0,2]\times [0,2]$, and let $f:\mf{I}\to \mb{R}$ be the function defined as $f(x,y)= x^7y^3$. For a positive integer $k$, let $\mf{P}_k$ be the uniformly regular partition of $\mf{I}$ into $k^2$ rectangles. Write down the summation formula for the Darboux upper sum $U(f,\mf{P}_k)$. Show that the limit $\di\lim_{k\to\infty}U(f,\mf{P}_k)$ exists and find the limit.
\end{question}

\atc
\begin{question}{\themyquestion}
 Let $\mk{D}$ be the triangle with vertices $(0,0)$, $(1,0)$ and $(1,1)$. Evaluate the integral $\di \int_{\mk{D}} e^{x^2}dxdy$.
\end{question}

\atc

\begin{question}{\themyquestion}
 Let
$\mk{D}$ be the region in the plane bounded between the cuve $y=x^2$ and the line $y=2x+3$. Evaluate the integral $\di \int_{\mk{D}}(x+2y)dxdy$.
\end{question}

\atc

\begin{question}{\themyquestion}
 Let
$\mk{D}$ be the region in the plane bounded between the cuve $y^2=4x$ and the line $y=2x-4$. Evaluate the integral $\di \int_{\mk{D}}(x+y)dxdy$.
\end{question}

\atc

\begin{question}{\themyquestion}
 Evaluate the integral $\di \int_0^1\int_y^1 \sqrt{9x^2+16}\,dx dy$.
\end{question}

\atc

\begin{question}{\themyquestion}
 Evaluate the integral $\di\int_{\mathcal{S}}xydxdydz$, where $\mathcal{S}$ is the solid bounded between the plane $x+y+z=4$ and the three coordinate planes. Then find the integral $\di\int_{\mathcal{S}}(4xy+5yz+6xz)dxdydz$.
\end{question}

\atc

\begin{question}{\themyquestion}
Let $f:[a,b]\to\mb{R}$ be a continuous function, and let $G$ be the set
\[G=\left\{(x,y,0)\,|\, a\leq x\leq b, 0\leq y\leq |f(x)|\right\}\]
in $\mb{R}^3$ that lies in the plane $z=0$. Rotate the set $G$ about the $x$-axis, we obtain a solid of revolution $S$, which can be described as
\[S=\left\{(x,y,z)\,|\, a\leq x\leq b, y^2+z^2\leq f(x)^2\right\}.\]
Show that the volume of $S$ is 
\[\text{vol}\,(S)=\pi \int_{a}^bf(x)^2dx.\]
\end{question}
 
\atc

\begin{question}{\themyquestion}
Let $\mk{D}_1$ and $\mk{D}_2$ be  Jordan measurable sets in $\mb{R}^m$ and $\mb{R}^n$ respectively. Show that the set $\mk{D}=\mk{D}_1\times\mk{D}_2$ is a Jordan measurable set in $\mb{R}^{m+n}$ and
\[\text{vol}\,(\mk{D}_1\times\mk{D}_2)=\text{vol}\,(\mk{D}_1)\times\text{vol}\,(\mk{D}_2).\]
\end{question}

 \section{Change of Variables Theorem}
 Consider the problem of evaluating an  integral of the form $\di\int_{\mk{D}}f(x,y)dxdy$ when $\mk{D}$ is the disc $\mk{D}=\left\{(x,y)\,|\,x^2+y^2\leq r^2\right\}$. When $f:\mk{D}\to\mb{R}$ is a continuous function, Fubini's theorem says that we can write the integral as
 \[\int_{\mk{D}}f(x,y)dxdy=\int_{-r}^r \int_{-\sqrt{r^2-x^2}}^{\sqrt{r^2-x^2}}f(x,y)dydx.\]
 However, it is usually quite complicated to evaluate this integral due to the square roots.  In some sense, we have not fully utilized the circular symmetry of the region of integration $\mk{D}$. For regions that have circular symmetry, it might be easier if we use polar coordinates $(r,\theta)$ instead of rectangular coordinates $(x,y)$. The goal of this section is to discuss the change of variables formula for multiple integrals.

 For single variable functions,  the change of variable formula is usually known as integration by substitution. We have proved the following theorem in volume I.
 \begin{theorem}[label=230223_5]{Integration by Substitution}
Let  $\psi:[a,b]\to\mathbb{R}$ be a function that satisfies the following conditions:
\begin{enumerate}[(i)]
\item $\psi$ is continuous and one-to-one on $[a,b]$;
\item $\psi$  is continuously differentiable on $(a,b)$;
\item $\psi'(x)$ is bounded on $(a,b)$.
\end{enumerate}  If $\psi([a,b])=[c,d]$, and    $f:[c,d]\to\mathbb{R}$ is a  bounded function that is   continuous on $(c,d)$, then the function $h:[a,b]\to \mathbb{R}$,
\[h(x)=f(\psi(x))|\psi'(x)|\] is  Riemann integrable and
\[\int_c^d f(u)du=\int_a^bf(\psi(x))|\psi'(x)|dx.\]
\end{theorem}
 The function $\psi:[a,b]\to\mathbb{R}$ that satisfies all the three   conditions (i)--(iii) in Theorem \ref{230223_5} defines  a \emph{smooth} change of variables $u=\psi(x)$ from $x$ to $u$.  
 
\begin{definition}{Smooth Change of Variables}
Let $\mathcal{O}$ be an open subset of $\mb{R}^n$. A mapping $\mathbf{\Psi}:\mathcal{O}\to\mb{R}^n$ from $\mathcal{O}$ to $\mb{R}^n$ is called a smooth change of variables provided that it satisfies the following conditions.
\begin{enumerate}[(i)]
\item $\mathbf{\Psi}:\mathcal{O}\to\mb{R}^n$ is one-to-one.
\item $\mathbf{\Psi}:\mathcal{O}\to\mb{R}^n$ is continuously differentiable.
\item For each $\mf{x}\in\mathcal{O}$, the derivative matrix $\mathbf{D\Psi}(\mf{x})$ is invertible.

\end{enumerate}
\end{definition}
\begin{remark}{}If the mapping $\mathbf{\Psi}:\mathcal{O}\to\mb{R}^n$ is   continuously differentiable, and the derivative matrix $\mathbf{D\Psi}(\mf{x})$ is invertible for each $\mf{x}\in\mathcal{O}$, the inverse function theorem  implies that the mapping
 $\mathbf{\Psi}:\mathcal{O}\to\mb{R}^n$ is locally one-to-one. However, it might not be {\it globally one-to-one}. For $\mathbf{\Psi}:\mathcal{O}\to\mb{R}^n$ to be a smooth change of variables, we need to impose the additional condition that it is globally one-to-one.

\end{remark}
\begin{example}[label=230902_2]{}
Let $\mf{x}_0$ be a point in $\mb{R}^n$. The mapping $\mf{\Psi}:\mb{R}^n\to\mb{R}^n$, $\mf{\Psi}(\mf{x})=\mf{x}+\mf{x}_0$ is a smooth change of variables. It is one-to-one, continuously differentiable, and the derivative matrix is $\mf{D\Psi}(\mf{x})=I_n$, which is invertible.
\end{example} 

\begin{example}[label=230902_1]{}
Let $\mf{x}_0$ and $\mf{y}_0$ be points in $\mb{R}^n$, and let $A$ be an invertible $n\times n$ matrix. The mapping $\mf{\Psi}:\mb{R}^n\to\mb{R}^n$ defined by

\vspace{-0.4cm}
\[\mf{\Psi}(\mf{x})=\mf{y}_0+A(\mf{x}-\mf{x}_0)\]
  is a one-to-one continuously differentiable mapping.
Its derivative matrix is $\mf{D\Psi}(\mf{x})=A$, which  is invertible for all $\mf{x}$ in $\mb{R}^n$. This shows that $\mf{\Psi}$ is a smooth change of variables. 
\end{example}
The mapping   $\mf{\Psi}:\mb{R}^n\to\mb{R}^n$ in Example \ref{230902_1}  is a composition of translations and an invertible linear transformation.

\begin{example}{}
Let $\mathcal{O}=\left\{(x,y)\,|\, x>0, y>0\right\}$, and let $\mf{\Psi}:\mathcal{O}\to\mb{R}^2$ be the mapping defined as
\[\mf{\Psi}(x,y)=(x^2-y^2, 2xy).\]Show that $\mf{\Psi}:\mathcal{O}\to\mb{R}^2$ is a smooth change of variables.
\end{example}
\begin{solution}{Solution}
First we show that $\mf{\Psi}$ is one-to-one. If $\mf{\Psi}(x_1, y_1)=\mf{\Psi}(x_2, y_2)$, then
\[x_1^2-y_1^2=x_2^2-y_2^2,\quad 2x_1y_1=2x_2y_2.\]
Let $z_1=x_1+iy_1$ and $z_2=x_2+iy_2$. Then we find that
\[z_1^2=x_1^2-y_1^2+2ix_1y_1=x_2^2-y_2^2+2ix_2y_2=z_2^2.\]
Hence, we must have $z_2=\pm z_1$. Restricted to $\mathcal{O}$, $x_1, x_2, y_1, y_2>0$. Hence, we must have $z_1=z_2$, or equivalently, $(x_1, y_1)=(x_2, y_2)$. This shows that $\mf{\Psi}:\mathcal{O}\to\mb{R}^2$ is one-to-one.
Now 
\[\mf{D\Psi}(x,y)=\begin{bmatrix} 2x & -2y\\2y & 2x\end{bmatrix}\]
is continuous, and $\det \mf{D\Psi}(x,y)=4x^2+4y^2\neq 0$ for all $(x,y)\in\mathcal{O}$. This proves that $\mf{\Psi}$ is continuously differentiable and the derivative matrix $\mf{D\Psi}(x,y)$ is invertible for all $(x,y)\in\mathcal{O}$.
Hence, $\mf{\Psi}:\mathcal{O}\to\mb{R}^2$ is a smooth change of variables.
\end{solution}

In this section, we will state  the change of variables theorem, and give some discusssions about why this theorem holds. We will also look at examples of how this theorem is applied, especially for polar and spherical coordinates. The proof of the theorem is quite technical and will be given in next section.

\begin{theorem}[label=230903_4]{The Change of Variables Theorem}
Let $\mathcal{O}$ be an open subset of $\mb{R}^n$, and let $\mf{\Psi}:\mathcal{O}\to\mb{R}^n$ be a smooth change of variables. If $\mk{D}$ is a Jordan measurable set such that its closure $\overline{\mk{D}}$ is contained in $\mathcal{O}$, then $\mf{\Psi}(\mk{D})$ is also Jordan measurable. If $f:\mf{\Psi}(\mk{D})\to\mb{R}$ is a bounded continuous function, then the function $g:\mk{D}\to\mb{R}$ defined as
\[g(\mf{x})=f\left(\mf{\Psi}(\mf{x})\right)\left|\det \mf{D\Psi}(\mf{x})\right|\] is Riemann integrable, and  
\begin{equation}\label{230902_3}\int_{\mf{\Psi}(\mk{D})}f(\mf{x})d\mf{x}=\int_{\mk{D}}g(\mf{x})d\mf{x}=\int_{\mk{D}}f\left(\mf{\Psi}(\mf{x})\right)\left|\det \mf{D\Psi}(\mf{x})\right|d\mf{x}.\end{equation}
\end{theorem}
Notice that the two vertical lines on $\det \mf{D\Psi}(\mf{x})$ in \eqref{230902_3} means the absolute value, not the determinant.

\begin{remark}{Jacobian}
For a mapping $\mf{\Psi}:\mathcal{O}\to\mb{R}^n$ from a subset of $\mb{R}^n$ to $\mb{R}^n$, the derivative matrix $\mf{D\Psi}(\mf{x})$ is also called  the {\it Jacobian matrix} of the mapping $\mf{\Psi}$.
The determinant of the Jacobian matrix is denoted by
\[ \frac{\pa(\Psi_1, \ldots, \Psi_n)}{\pa (x_1, \ldots, x_n)}.\]It 
is known as the {\it Jacobian determinant}, or simply as {\it Jacobian}.
In practice, we will often denote a change of variables $\mf{\Psi}:\mathcal{O}\to\mb{R}^n$ by $\mf{u}=\bf{\Psi}(\mf{x})$. Then the Jacobian can be written as
\[\frac{\pa (u_1, \ldots, u_n)}{\pa (x_1, \ldots, x_n)}.\]Using this notation, the change of variables formula \eqref{230902_3} reads as
\begin{align*}
&\int_{\mf{\Psi}(\mk{D})}f(u_1, \ldots, u_n)du_1\cdots du_n\\&=\int_{\mk{D}}f(u_1(\mf{x}), \ldots, u_n(\mf{x}))\left|\frac{\pa (u_1, \ldots, u_n)}{\pa (x_1, \ldots, x_n)}\right|dx_1 \ldots dx_n.\end{align*}
\end{remark}

  \subsection{Translations and  Linear Transformations} 
In the single variable case, a translation is a map $T:\mb{R}\to\mb{R}$, $T(x)=x+c$. If $f:[a,b]\to\mb{R}$ is a Riemann integrable function, then
 \[\int_a^b f(x)dx=\int_{a-c}^{b-c}f(x+c)dx.\]
For $n\geq 2$,  we have the following theorem, which is a stronger version of Theorem \ref{230903_4}.

\begin{theorem}[label=230903_8]{}Let $\mf{x}_0$ be a fixed point in $\mb{R}^n$, and let $\mf{\Psi}:\mb{R}^n\to\mb{R}^n$ be the translation $\mf{\Psi}(\mf{x})=\mf{x}+\mf{x}_0$.
If $\mk{D}$ is a Jordan measurable subset of $\mb{R}^n$, then $\mf{\Psi}(\mk{D})$ is Jordan measurable. If $f:\mf{\Psi}(\mk{D})\to\mb{R}$ is a Riemann integrable function,  then  $g=(f\circ \mf{\Psi}):\mk{D}\to\mb{R}$ is also Riemann integrable, and 
  \begin{equation}\label{230902_6}\int_{ \mk{D}+\mf{x}_0}f(\mf{x})d\mf{x}=\int_{\mf{\Psi}(\mk{D})}f(\mf{x})d\mf{x}=\int_{\mk{D}}g(\mf{x})d\mf{x}=\int_{\mk{D} }f( \mf{x}+\mf{x}_0)d\mf{x}.
  \end{equation}
\end{theorem}

\begin{myproof}{Proof}
Obviously, translation maps a rectangle to a rectangle with the same volume. Hence, it maps sets that have Jordan content zero to sets that have Jordan content zero.  It is also obvious that $\mf{\Psi}$ maps the boundary of $\mk{D}$ to the boundary of $\mf{\Psi}(\mk{D})$. This shows that $\mf{\Psi}(\mk{D})$ is Jordan measurable.

 If $\mf{I}$ is a closed rectangle that contains $\mk{D}$, then $\mf{I}'=\mf{I}+\mf{x}_0$ is a closed rectangle that contains $\mf{\Psi}(\mk{D})$. Let $\check{f}: \mf{I}'\to\mb{R}$ and $\check{g}:\mf{I}\to\mb{R}$ be the zero extensions of $f:\mf{\Psi}(\mk{D})\to\mb{R}$ and $g=(f\circ \mf{\Psi}):\mk{D}\to\mb{R}$ respectively. Then $\check{g}=\check{f}\circ\mf{\Psi}$. Since $f:\mf{\Psi}(\mk{D})\to\mb{R}$ is a Riemann integrable,
 \[\underline{\int_{\mf{I}'}}\check{f}=\overline{\int_{\mf{I}'}}\check{f}=\int_{\mf{I}'}\check{f}=\int_{\mf{\Psi}(\mk{D})}f.\]
 
  Given   a partition $\mf{P}'=(P_1', \ldots, P_n')$ of $\mf{I}'$, let $\mf{P}=(P_1, \ldots, P_n)$ be the partition of $\mf{I}$ induced by the translation $\mf{\Psi}$. Namely, $\mf{P}$ is a partition such that the rectangle $\mf{J}$ is in $\mathcal{J}_{\mf{P}}$ if and only if $\mf{J}'=\mf{\Psi}(\mf{J})=\mf{J}+\mf{x}_0$ is in $\mf{P}'$.
  \bp
  
   Then
 \begin{align*}m_{\mf{J}}(\check{g})&=\inf\left\{g(\mf{x})\,|\, \mf{x}\in\mf{J}\right\}= \inf\left\{f(\mf{x}+\mf{x}_0)\,|\, \mf{x}\in\mf{J}\right\}\\
 &=   \inf\left\{f(\mf{x})\,|\, \mf{x}'\in\mf{J}+\mf{x}_0\right\}=m_{\mf{J}'}(\check{f}).
 \end{align*}
 
  It follows that
 \[L(\check{g},\mf{P})=\sum_{\mf{J}\in\mathcal{J}_{\mf{P}}}m_{\mf{J}}(\check{g})\,\text{vol}\,(\mf{J})=\sum_{\mf{J}'\in\mathcal{J}_{\mf{P}'}}m_{\mf{J}'}(\check{f})\,\text{vol}\,(\mf{J}')=L(\check{f}, \mf{P}').\]
 
 Similarly, we have
 \[U(\check{g},\mf{P})=U(\check{f}, \mf{P}').\]

Thus, the sets $S_L(\check{g})$ and $S_L(\check{f})$ of lower sums of $\check{g}$ and $\check{f}$ are the same, and the sets $S_U(\check{g})$ and $S_U(\check{f})$ of upper sums of $\check{g}$ and $\check{f}$ are also the same. These imply that
\begin{align*}\underline{\int_{\mf{I}}}\check{g}&=\sup S_L(\check{g})=\sup S_L(\check{f})=\underline{\int_{\mf{I}'}}\check{f}=\int_{\mf{\Psi}(\mk{D})}f,\\
\overline{\int_{\mf{I}}}\check{g}&=\inf S_U(\check{g})=\inf S_U(\check{f})=\overline{\int_{\mf{I}'}}\check{f}=\int_{\mf{\Psi}(\mk{D})}f.
\end{align*}
Hence, $g:\mk{D}\to\mf{I}$ is Riemann integrable and
\[\int_{\mk{D}}g=\int_{\mf{I}}\check{g}=\int_{\mf{\Psi}(\mk{D})} f.\]

\end{myproof}

\begin{figure}[ht]
\centering
\includegraphics[scale=0.18]{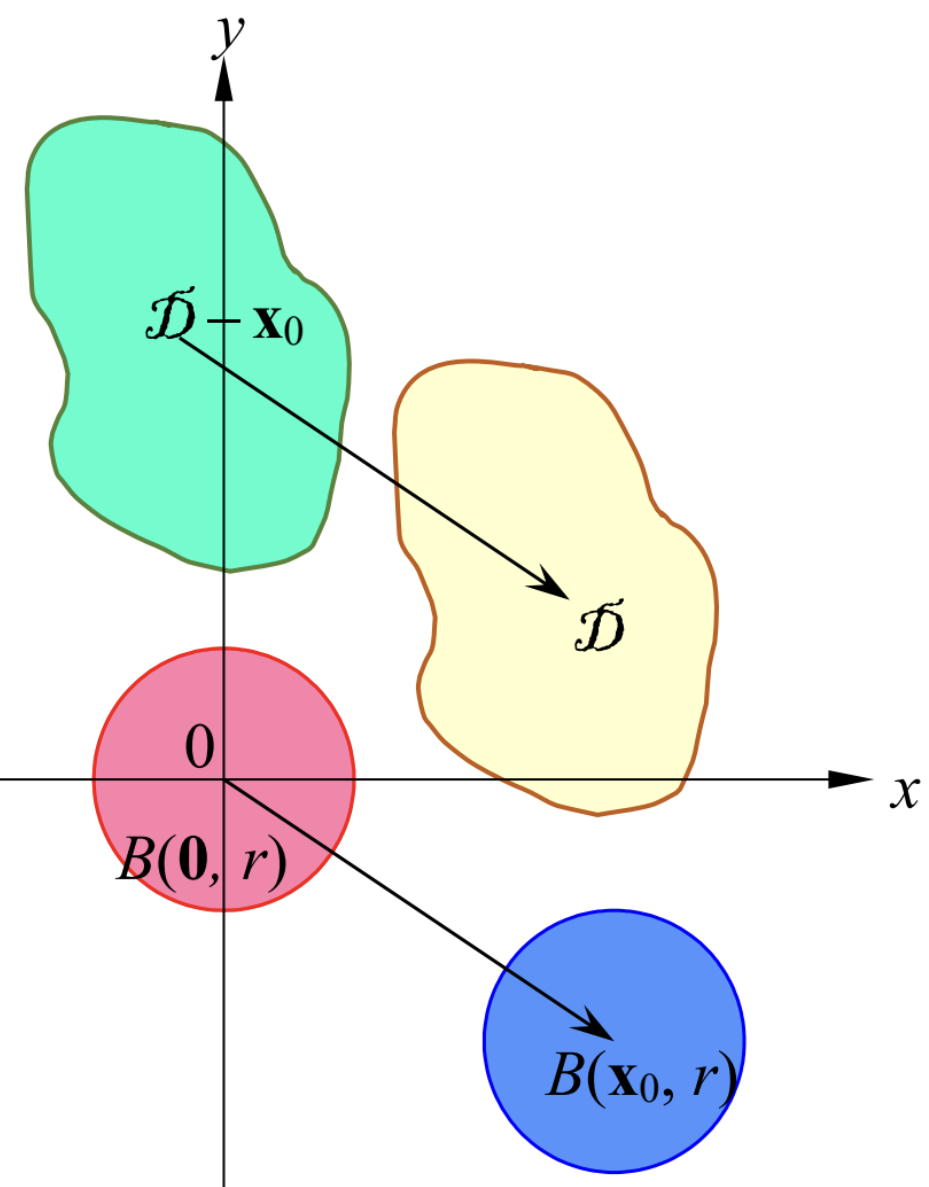}
\caption{A translation in the plane. }\label{figure113}
\end{figure}

\begin{remark}{}
It is easy to check that for the translation $\mf{\Psi}:\mb{R}^n\to\mb{R}^n$, $\mf{\Psi}(\mf{x})=\mf{x}+\mf{x}_0$, the change of variables formula \eqref{230902_6} is precisely the formula \eqref{230902_3}, since $\mf{D\Psi}(\mf{x})=I_n$ in this case.
\end{remark}
\begin{example}[label=230902_4]{}
 Let 
\[ \mk{D}=\left\{(x,y)\,|\,(x-2)^2+(y+3)^2\leq 16\right\}.\] Evaluate the integral $\di \int_{\mk{D}} (4x+y)dxdy$.
  \end{example}
  \begin{solution}{Solution}
  Make the change of variables $u=x-2$, $v=y+3$, which is a translation. Then $x=u+2$, $y=v-3$, and we have
  \begin{align*}
  \int_{\mk{D}} (4x+y)dxdy&=\int_{u^2+v^2\leq 16}(4u+8+v-3)dudv\\&=\int_{u^2+v^2\leq 16} (4u+v+5)dudv.\end{align*}
  Since the disc $\mathcal{B}=\left\{(u,v)\,|\, u^2+v^2\leq 16\right\}$ is invariant when we change $u$ to $-u$, or change $v$ to $-v$, 
  the integrals
  \bs
  \[\int_{u^2+v^2\leq 16} ududv\quad\text{and}\quad \int_{u^2+v^2\leq 16} vdudv\]
   are equal to 0. Therefore,
  \[\int_{\mk{D}} (4x+y)dxdy=5\int_{u^2+v^2\leq 16} dudv.\]

  In single variable analysis, we have shown that the area of a disc of radius $r$ is $\pi r^2$. Hence,
  \[\int_{\mk{D}} (4x+y)dxdy=5\times\text{area}\,(\mathcal{B})=5\times 16\pi=80\pi.\]
  \end{solution}

Now we consider a linear transformation $\mf{T}:\mb{R}^n\to\mb{R}^n$, $\mf{T}(\mf{x})=A\mf{x}$ defined by an invertible matrix $A$.  
Since $\mf{DT}(\mf{x})=A$, the change of variables theorem says that for  any function $f:\mk{D}\to\mb{R}$ that is bounded and continuous on $\mk{D}$, 
\begin{equation}\label{230902_10}\int_{\mf{T}(\mk{D})}f(\mf{x})d\mf{x}=\int_{\mk{D}}f(\mf{T}(\mf{x}))|\det A|d\mf{x}=|\det A|\int_{\mk{D}}f(\mf{T}(\mf{x}))d\mf{x}.\end{equation}
In the special case where $f$ is a constant function, we have
\[\text{vol}\,(\mf{T}(\mk{D}))=|\det A|\,\text{vol}\,(\mk{D}).\]
A very crucial fact to the proof of the change of variables theorem is a special case of this formula when $\mk{D}$ is a rectangle.
\begin{theorem}[label=230903_1]{}
Let $\mf{I}=\di\prod_{i=1}^n[a_i,b_i]$, and let $\mf{T}:\mb{R}^n\to\mb{R}^n$, $\mf{T}(\mf{x})=A\mf{x}$ be an invertible linear transformation. Then
\begin{equation}\label{230905_6}\text{vol}\,(\mf{T}(\mf{I}))=|\det A|\,\text{vol}\,(\mf{I}).\end{equation}
\end{theorem}

Linear transformations  map linear objects to linear objects. However, the image of a rectangle under a linear transformation is not necessary a rectangle, but  is always a parallelepiped. 
If $\mf{I}$ is the closed rectangle $ \di\prod_{i=1}^n[a_i, b_i]$, a point  $\mf{x}$ in $\mf{I}$ can be written as
\[\mf{x}=\mf{a}+t_1 (b_1-a_1)\mf{e}_1+\cdots+t_n(b_n-a_n)\mf{e}_n,\]
where $ \mf{a}=(a_1, \ldots, a_n)$, and $\mf{t}=(t_1, \ldots, t_n)\in [0,1]^n$. Hence, we say that the rectangle $\mf{I}$ is a parallelepiped based at the point $\mf{a}$ and spanned by the $n$-linearly independent vectors $\mf{v}_i=(b_i-a_i)\mf{e}_i$, $1\leq i\leq n$.

\begin{definition}{Parallelepipeds}
A (closed) parallelepiped in $\mb{R}^n$ is a solid $\mathscr{P}$ in $\mb{R}^n$ based at a point $\mf{a}$ and spanned by $n$-linearly independent vectors $\mf{v}_1, \ldots, \mf{v}_n$. It can be described as 
\[\mathscr{P}=\left\{\mathbf{a}+t_1\mf{v}_1+\cdots+t_n\mf{v}_n\,|\, \mf{t}=(t_1, \ldots, t_n)\in [0,1]^n\right\}.\]
\end{definition}

\begin{figure}[ht]
\centering
\includegraphics[scale=0.18]{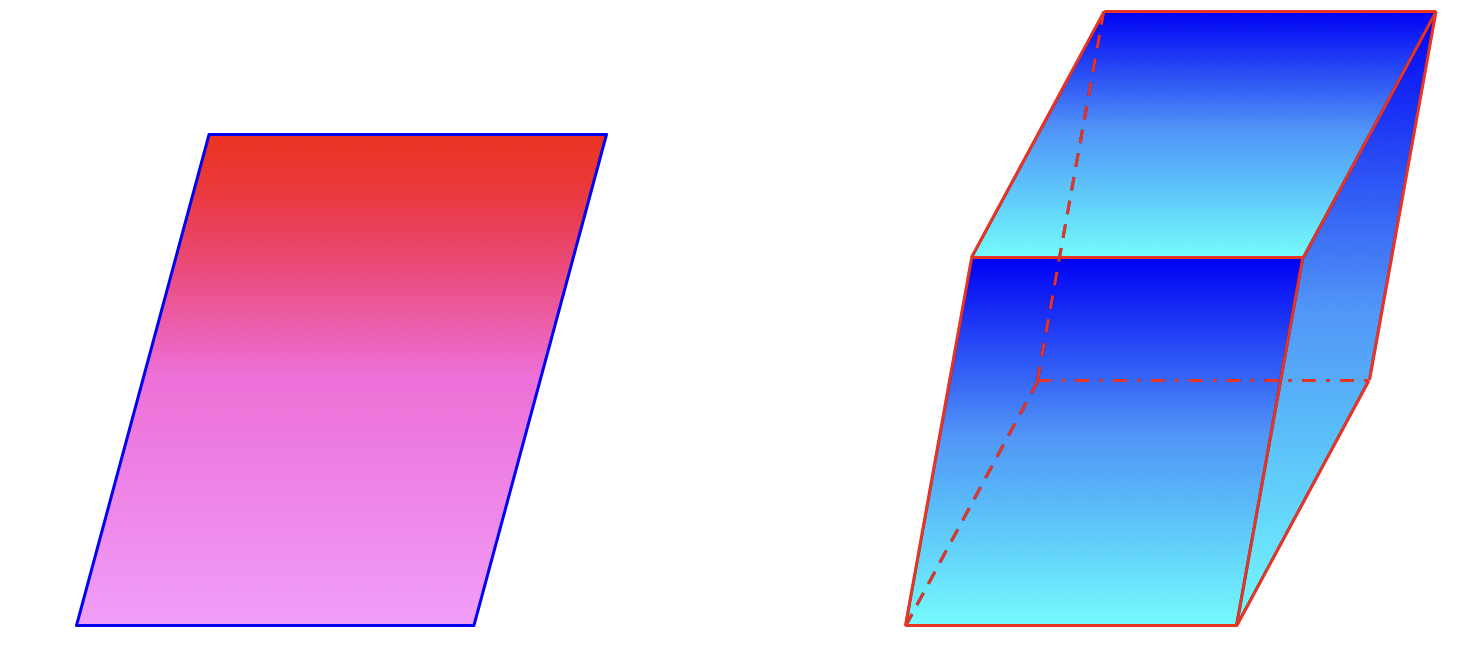}
\caption{Parallelepipeds in $\mb{R}^2$ and $\mb{R}^3$. }\label{figure114}
\end{figure}
 The boundary of a parallelepiped    is a union of $2n$ bounded subsets, each of them is contained in a hyperplane. Thus, the boundary of a parallelepiped has Jordan content zero. Therefore, a parallelepiped is a Jordan measurable set. 

If  $\mathscr{P}$ be a parallelepiped in $\mb{R}^n$   based at the point $\mf{a}$ and spanned by the vectors   $\mf{v}_1, \ldots, \mf{v}_n$, and  $\mf{T}:\mb{R}^n\to\mb{R}^n$ is an invertible linear transformation, then $\mf{T}(\mathscr{P})$ is the parallelepiped in $\mb{R}^n$   based at the point $\mf{T}(\mf{a})$ and spanned by the vectors   $\mf{T}(\mf{v}_1), \ldots, \mf{T}(\mf{v}_n)$.

The cube $[0,1]^n$ is called the standard unit cube and it is often denoted by $Q_n$. 
 If $\mathscr{P}$ is a parallelepiped in $\mb{R}^n$   based at the point $\mf{a}$ and spanned by the vectors   $\mf{v}_1, \ldots, \mf{v}_n$, then
$\mathscr{P}=\mf{\Psi}(Q_n)$, where   $\mf{\Psi}:\mb{R}^n\to\mb{R}^n$ is the mapping
\[\mf{\Psi}(\mf{x})=A\mf{x}+\mf{a},\hspace{1cm}A=\begin{bmatrix} \mf{v}_1 & \rvline & \cdots & \rvline &\mf{v}_n\end{bmatrix},\]
which is a composition of an invertible linear transformation and a translation.  Theorem \ref{230903_1} says that
\begin{equation}\label{230903_2}\text{vol}\,(\mathscr{P})=|\det A|,\end{equation} where $A$ is the matrix whose column vectors are $\mf{v}_1, \ldots, \mf{v}_n$. For example, for a parallelogram in $\mb{R}^2$ which is spanned by the vectors
\[\mf{v}_1=\begin{bmatrix} a_1\\b_1\end{bmatrix}\quad \text{and}\quad \mf{v}_2=\begin{bmatrix} a_2\\b_2\end{bmatrix},\]
the area of the parallogram is
\[\det \begin{bmatrix} a_1 & a_2\\b_1 & b_2\end{bmatrix}.\]
For a parallelepiped in $\mb{R}^3$ which is spanned by the vectors
\[\mf{v}_1=\begin{bmatrix} a_1\\b_1\\c_1\end{bmatrix},\quad   \mf{v}_2=\begin{bmatrix} a_2\\b_2\\c_2\end{bmatrix}\quad \text{and}\quad \mf{v}_3=\begin{bmatrix} a_3\\b_3\\c_3\end{bmatrix},\]
the volume of the parallelepiped is 
\[\det \begin{bmatrix} a_1 & a_1 &a_3\\b_1 & b_2 &b_3\\c_1 &c_2 & c_3\end{bmatrix}.\]
These formulas have been derived in an elementary course. For general $n$, we will prove \eqref{230903_2} in Appendix \ref{appB} using geometric arguments. This will then imply Theorem \ref{230903_1}.

From the theory of linear algebra, we know that an invertible matrix is a product of elementary matrices. Hence,
an invertible linear transformation $\mf{T}:\mb{R}^n\to\mb{R}^n$ can be written as
\[\mf{T}=\mf{T}_m\circ \cdots\circ \mf{T}_2\circ\mf{T}_1,\]
where $\mf{T}_1, \mf{T}_2, \ldots, \mf{T}_m$ is one of the three types of elementary transformations, corresponding to the three types of elementary matrices.
\begin{enumerate}[I.]
\item  When $E$ is the elementary matrix obtained from the identity matrix $I_n$ by interchanging two distinct  rows  $i$ and $j$, the linear tranformation $\mf{T}:\mb{R}^n\to\mb{R}^n$, $\mf{T}(\mf{x})=E\mf{x}$ interchanges $x_i$ and $x_j$, and fixes the other variables. In this case, $\det E=-1$ and $|\det E|=1$.  
\item When $E$ is the elementary matrix obtained from the identity matrix  $I_n$ by multiplying   row $i$ by a nonzero constant $c$, the linear   transformation  $\mf{T}:\mb{R}^n\to\mb{R}^n$, $\mf{T}(\mf{x})=E\mf{x}$  maps the point $\mf{x}=(x_1, \ldots, x_n)$ to \[\mf{T}(\mf{x})=(x_1, \ldots, x_{i-1}, cx_i, x_{i+1},\ldots, x_n).\]  In this case, $\det E=c$, and $|\det E|=|c|$.
\item When $E$ is the elementary matrix obtained from the identity matrix $I_n$ by adding a constant  $c$ times of row $j$ to another row $i$, the linear   transformation  $\mf{T}:\mb{R}^n\to\mb{R}^n$, $\mf{T}(\mf{x})=E\mf{x}$  maps the point $\mf{x}=(x_1, \ldots, x_n)$ to \[\mf{T}(\mf{x})=(x_1, \ldots, x_{i-1}, x_i+cx_j, x_{i+1},\ldots, x_n).\]  In this case, $\det E=1$, and $|\det E|=1$.
\end{enumerate}
Since each of the elementary transformations involves changes in at most two variables, it is sufficient to consider these transformations when $n=2$.
\begin{example}[label=230903_5]{}
Let $\mf{T}:\mb{R}^2\to \mb{R}^2 $ be the linear transformation
\[\mf{T}(x,y)=(y,  x).\]
The matrix $E$ corresponding to this transformation is
\[E=\begin{bmatrix}  0 & 1\\   1 & 0\end{bmatrix}.\]
Under this transformation, the rectangle $\mf{I}=[a, b]\times [c,d]$ is mapped to the rectangle $\mf{I}'=[c,d]\times [a,b]$. It is easy to see that
\[\text{vol}\,(\mf{I}')= \text{vol}\,(\mf{I}).\]
\end{example}

\begin{figure}[ht]
\centering
\includegraphics[scale=0.2]{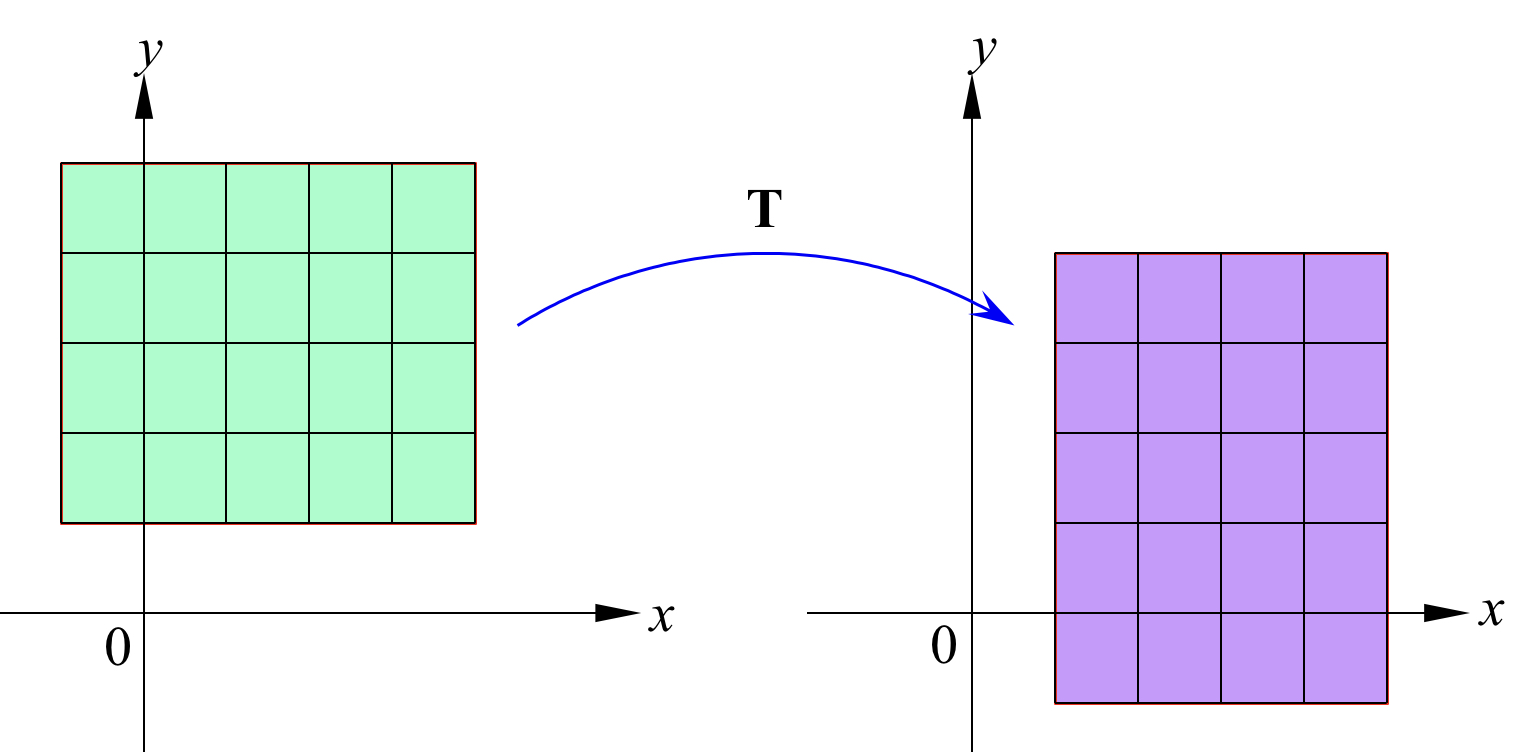}
\caption{The linear transformation $\mf{T}(x,y)=(y,x)$. }\label{figure115}
 \end{figure}

 \begin{figure}[ht]
\centering
\includegraphics[scale=0.2]{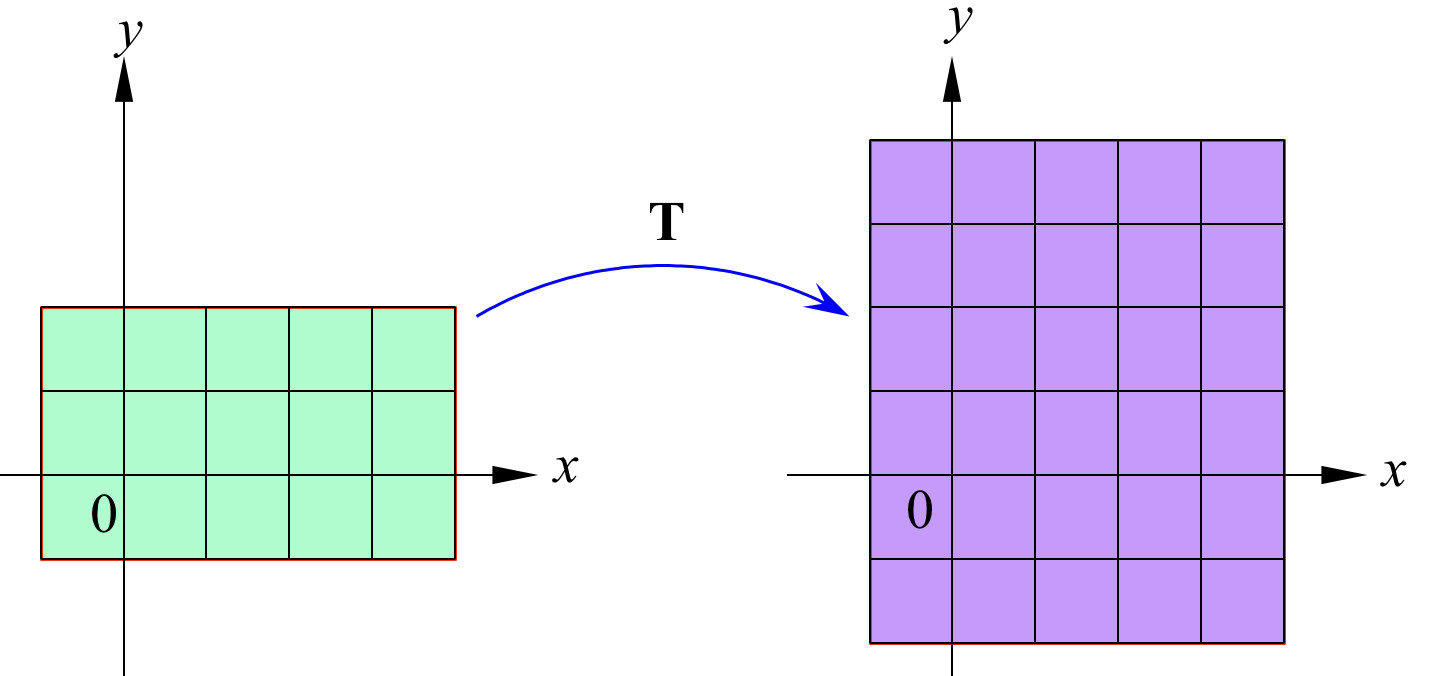}
\caption{The linear transformation $\mf{T}(x,y)=(x,2y)$. }\label{figure116}
\end{figure}

\begin{example}[label=230903_6]{}
Let $\mf{T}:\mb{R}^2\to \mb{R}^2 $ be the linear transformation
\[\mf{T}(x,y)=(x,  ky),\hspace{1cm}k\neq 0.\]
The matrix $E$ corresponding to this transformation is
\[E=\begin{bmatrix}  1 & 0\\   0 &k\end{bmatrix}.\]
Under this transformation, the rectangle $\mf{I}=[a, b]\times [c,d]$ is mapped to the rectangle $\mf{I}'=[a,b]\times [kc,kd]$ if $k>0$; and to the rectangle $\mf{I}'=[a,b]\times [kd,kc]$ if $k<0$. In any case, we find that 
\[\text{vol}\,(\mf{I}')=|k|\,\text{vol}\,(\mf{I}).\]
\end{example}

\begin{example}[label=230903_7]{}
Let $\mf{T}:\mb{R}^2\to \mb{R}^2 $ be the linear transformation
\[\mf{T}(x,y)=(x+ky,  y).\]
The matrix $E$ corresponding to this transformation is
\[E=\begin{bmatrix}  1 & k\\   0 &1\end{bmatrix}.\]
Under this transformation, the rectangle $\mf{I}=[a, b]\times [c,d]$ is mapped to     the parallepiped $\mathscr{P}$  with vertices $(a+kc, c)$, $(a+kd, d)$, $(b+kc, c)$ and $(b+kd, d)$. Using elementary geometric argument, one can show that
\[\text{vol}\,(\mathscr{P})= \text{vol}\,(\mf{I}).\]
\end{example}

\begin{figure}[ht]
\centering
\includegraphics[scale=0.2]{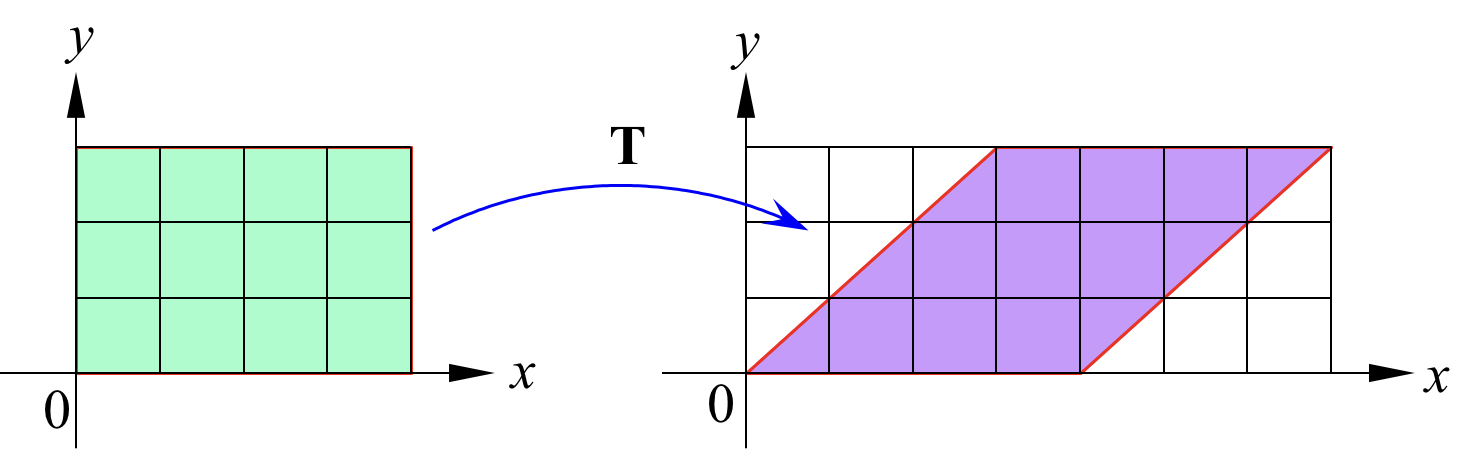}
\caption{The linear transformation $\mf{T}(x,y)=(x+y,y)$. }\label{figure117}

\includegraphics[scale=0.2]{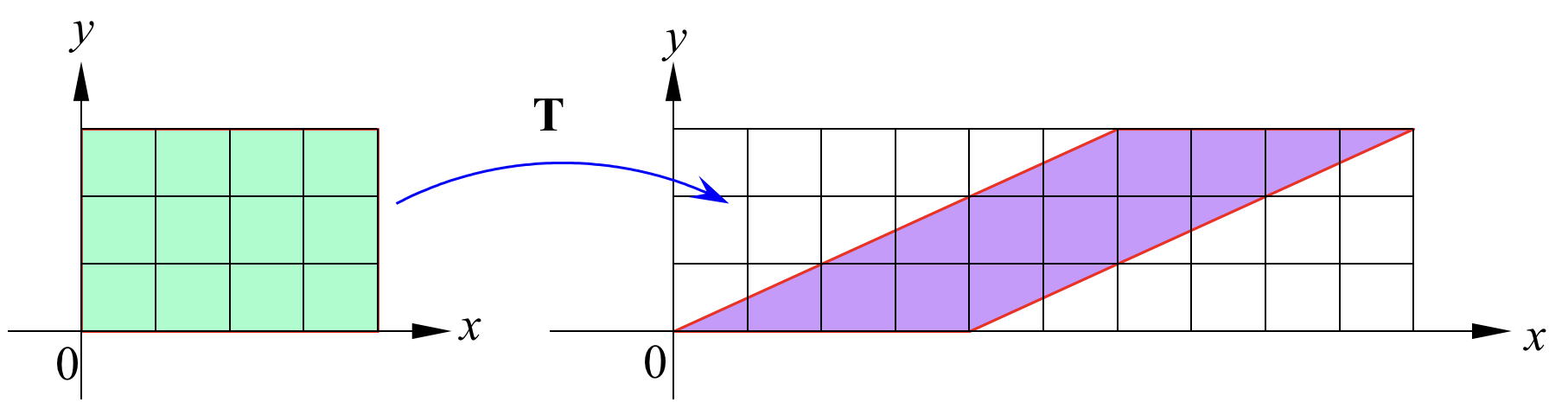}
\caption{The linear transformation $\mf{T}(x,y)=(x+2y,y)$. }\label{figure118}

\includegraphics[scale=0.2]{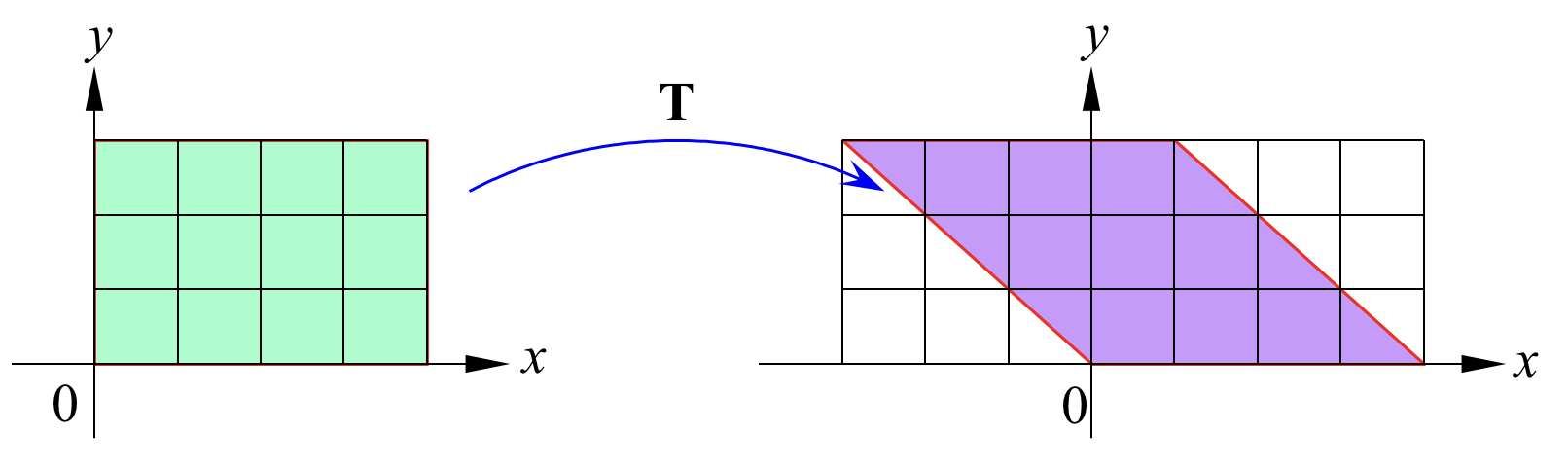}
\caption{The linear transformation $\mf{T}(x,y)=(x-y,y)$. }\label{figure119}
\end{figure}
Combining Example \ref{230903_5}, Example \ref{230903_6} and Example \ref{230903_7}, we conclude that \eqref{230905_6} holds when $\mf{T}:\mb{R}^n\to\mathbb{R}^n$ is an elementary transformation. 

The type II elementary transformations maps rectangles to rectangles, so do their compositions. Therefore, \eqref{230905_6} also holds if the linear transformation $\mf{T}:\mb{R}^n\to\mathbb{R}^n$ is a composition of   type II elementary transformations. This  gives the following.

\begin{theorem}[label=230903_9]{}Let $\mf{x}_0=(u_1, \ldots, u_n)$ be a fixed point in $\mb{R}^n$, and let $\mf{\Psi}:\mb{R}^n\to\mb{R}^n$ be the mapping 
\[\Psi_i(\mf{x})=\alpha_i x_i+u_i.\] Equivalently,
$\mf{\Psi}(\mf{x})=A\mf{x}+\mf{x}_0$, where 
$A$ is a diagonal matrix with diagonal entries $\alpha_1, \ldots, \alpha_n$.
If $\mk{D}$ is a Jordan measurable subset of $\mb{R}^n$, then $\mf{\Psi}(\mk{D})$ is Jordan measurable. If $f:\mf{\Psi}(\mk{D})\to\mb{R}$ is a Riemann integrable function,  then  $h=(f\circ \mf{\Psi}):\mk{D}\to\mb{R}$ is also Riemann integrable, and 
  \begin{equation}\label{230902_9} \int_{\mf{\Psi}(\mk{D})}f(\mf{x})d\mf{x}=|\det A|\int_{\mk{D}}h(\mf{x})d\mf{x}=|\det A|\int_{\mk{D} }f(A \mf{x}+\mf{x}_0)d\mf{x}.
  \end{equation}
\end{theorem}
Notice that $\det A=\alpha_1\cdots \alpha_n$. If $\mf{y}=\mf{\Psi}(\mf{x})$, then
\[y_i=\alpha_i x_i+u_i,\hspace{1cm} 1\leq i\leq n.\]

 The proof of Theorem \ref{230903_9} is similar to the proof Theorem \ref{230903_8}, by establishing one-to-one correspondence between the partitions, and using the fact that for any rectangles $\mf{J}$, 
\[\text{vol}\,(\mf{\Psi}(\mf{J}))=|\alpha_1\cdots \alpha_n|\,\text{vol}\,(\mf{J}).\]

\begin{example}{}
Find the area of the ellipse 
\[\mathcal{E}=\left\{(x,y)\,|\, 4(x+1)^2+9(y-5)^2\leq 49\right\}.\]
\end{example}
\begin{solution}{Solution}
Make a change of variables $u=2(x+1)$ and $v=3(y-5)$. The Jacobian   is
\[\frac{\pa (u,v)}{\pa (x,y)}=\det \begin{bmatrix} 2 & 0 \\0 & 3\end{bmatrix} =6,\]
and so
\[\frac{\pa (x,y)}{\pa (u,v)}=\frac{1}{6}.\]
Therefore,
\begin{align*}
\text{area}\,(\mathcal{E})&=\int_{ 4(x+1)^2+9(y-5)^2\leq 49}dxdy=\int_{u^2+v^2\leq 49}\left|\frac{\pa (x,y)}{\pa (u,v)}\right|dudv\\
&=\frac{1}{6}\int_{u^2+v^2\leq 49}dudv=\frac{49}{6}\pi.
\end{align*}
\end{solution}

Finally, let us consider an example of applying a general linear transformation. 
 \begin{example}{}
 Evaluate the integral
 \[\int_{\mk{D}}\frac{2x+3y+3}{2x-3y+8}dxdy,\]
 where \[\mk{D}=\di\left\{(x,y)\,|\, 2|x|+3|y|\leq 6\right\}.\]
 \end{example}
 \begin{solution}{Solution}
 Notice that for any $(x,y)\in \mk{D}$, \[|2x-3y+8|\geq 8-2|x|-3|y|\geq 2. \]

 \end{solution}
 
\begin{figure}[ht]
\centering
\includegraphics[scale=0.2]{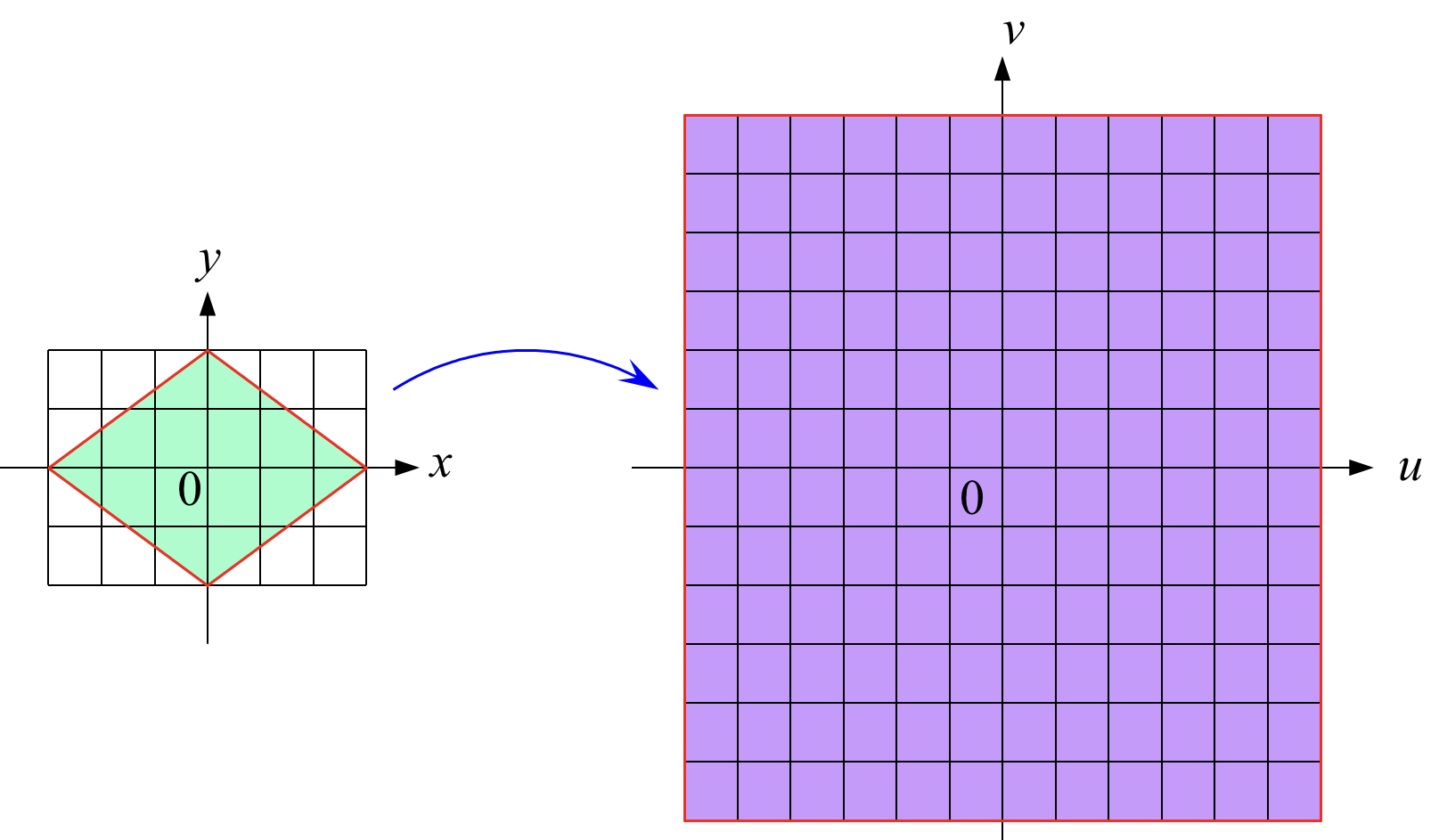}
\caption{The transformation  $u=2x-3y$ and $v=2x+3y$. }\label{figure120}
\end{figure}
 \begin{solution}{}
  Hence, the function \[h(x,y)= \frac{2x+3y+3}{2x-3y+8}\] is continuous on $\mk{D}$.
 The region $\mk{D}$ is enclosed by the 4 lines $2x+3y=6$, $2x+3y=-6$, $2x-3y=6$ and $2x-3y=-6$. This prompts us to define a change of variables by
 $u=2x-3y$ and $v=2x+3y$. This is a linear transformation with Jacobian 
 \[\frac{\pa (u,v)}{\pa (x,y)}=\det \begin{bmatrix} 2 & -3\\2 & 3\end{bmatrix} =12.\]
  Therefore,
 \[\frac{\pa (x,y)}{\pa (u,v)}=\frac{1}{12}.\]

 The region $\mk{D}$ in the $(x,y)$-plane is mapped to the rectangle
 \[\mathcal{R}=\left\{(u,v)\,|\, -6\leq u\leq 6, -6\leq v\leq 6\right\}\] in the $(u,v)$-plane.
 \bs
 Thus,
 \begin{align*}
 \int_{\mk{D}}\frac{2x+3y+3}{2x-3y+8}dxdy&=\int_{\mathcal{R}}\frac{v}{u+8}\left|\frac{\pa (x,y)}{\pa (u,v)}\right|dudv\\
 &=\frac{1}{12}\int_{-6}^6\int_{-6}^6 \frac{v+3}{u+8}dudv\\
 &=\frac{1}{12}\left[\frac{v^2}{2}+3v\right]_{-6}^6\left[\ln (u+8)\right]_{-6}^6\\
 &=3\ln 7.
 \end{align*} 
 \end{solution}

\subsection{Polar Coordinates}
Given a point $(x,y)$ in the plane $\mb{R}^2$, if $r$ is a nonnegative number  and $\theta$ is a real number such that \[x=r\cos\theta,\quad y=r\sin\theta,\]
then $(r,\theta)$ are called the polar coordinates of the point $(x,y)$. Notice that
\[r=\sqrt{x^2+y^2}.\]
Restricted to 
\[V=\left\{(r,\theta)\,|\,r>0, 0\leq \theta<2\pi\right\},\]
the map $\mf{\Phi}:V\to\mb{R}^2$,
\[\mf{\Phi}(r,\theta)=(r\cos\theta, r\sin\theta)\] is one-to-one, and its range is $\mb{R}^2\setminus\{(0,0)\}$. However, the inverse of $\mf{\Phi}$ fails to be continuous. We can extend the map $\mf{\Phi}$ to $\mb{R}^2$ continuously. Namely, 
given $(r,\theta)\in \mb{R}^2$, let $(x,y)=\mf{\Phi}(r,\theta)$, where
\[x=r\cos\theta, \quad  y=r\sin\theta.\]
Then $\mf{\Phi}:\mb{R}^2\to\mb{R}^2$ is continuously differentiable, but it fails to be one-to-one. Nevertheless, for any real number $\alpha$, the map is continuous and one-to-one on the open set  
\[\mathcal{O}_{\alpha}=\left\{(r,\theta)\,|\,r>0, \alpha<\theta<\alpha+2\pi\right\}.\]  The derivative 
matrix of $\mf{\Phi}:\mb{R}^2\to\mb{R}^2$  is given by
\[\mf{D\Phi}(r,\theta)=\begin{bmatrix} \cos\theta & -r\sin\theta\\\sin\theta & r\cos\theta\end{bmatrix}.\]
Since 
\[\det \mf{D\Phi}(r,\theta)=r\cos^2\theta+r\sin^2\theta=r,\]
we find that for any $(r,\theta)\in \mathcal{O}_{\alpha}$, $ \mf{D\Phi}(r,\theta)$ is invertible. Hence, $\mf{\Phi}:\mathcal{O}_{\alpha}\to\mb{R}^2$ is a smooth change of variables.

\begin{figure}[!ht]
\centering
\includegraphics[scale=0.2]{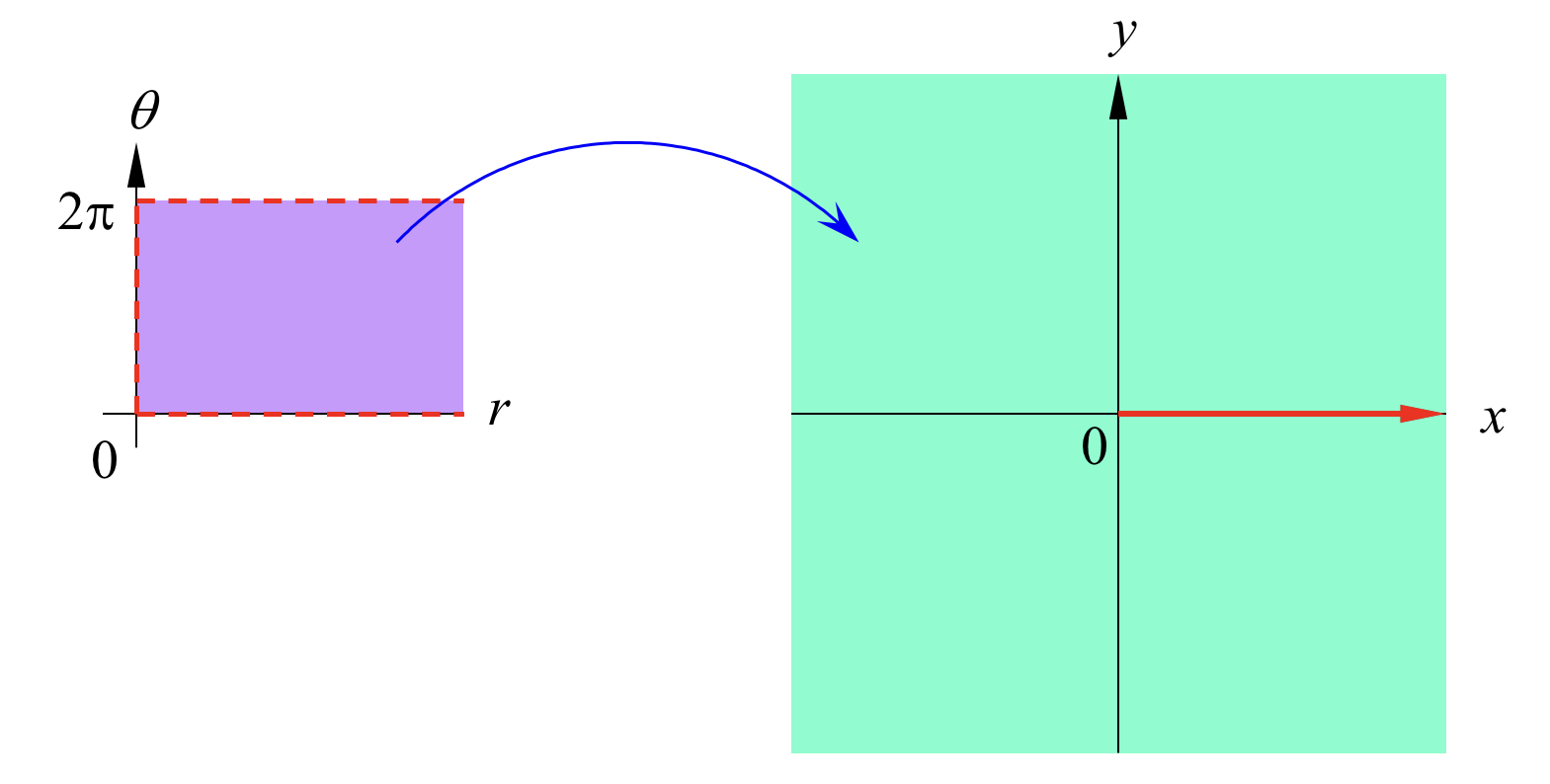}
\caption{The mapping $\mf{\Phi}:\mathcal{O}\to\mb{R}^2$, $\mf{\Phi}(r,\theta)=(r\cos\theta, r\sin\theta)$ maps  $\mathcal{O}=\left\{(r,\theta)\,|\,r>0, 0<\theta<2\pi\right\}$ to $\mb{R}^2\setminus L$, where $L$ is the positive $x$-axis.}\label{figure122}
 \end{figure}
 
Let us consider the special case where $\alpha=0$. In this case, let $\mathcal{O}=\mathcal{O}_0$. 
 The map $\mf{\Phi}:\mathcal{O}\to \mb{R}^2$, 
\[\mf{\Phi}(r,\theta)=(r\cos\theta, r\sin\theta)\] is a smooth change of variables from  polar coordinates to rectangular coordinates. 
Under this change of variables, 
\[\mf{\Phi}(\mathcal{O})=\mb{R}^2\setminus\left\{(x,0)\,|\, x\geq 0\right\}\]is an open set in $\mb{R}^2$.
If $\mathcal{D}$ is the open rectangle $(r_1, r_2)\times (\theta_1, \theta_2)$, with
\[0<r_1<r_2\quad\text{and}\quad 0<\theta_1<\theta_2<2\pi,\]
then 
$\mf{\Phi}(\mathcal{D})$ is the open set
bounded between the two circles $x_1^2+y_1^2=r_1^2$ and $x_2^2+y_2^2=r_2^2$, and the two rays $y=x\tan\theta_1$, $x\geq 0$ and $y=x\tan\theta_2$, $x\geq 0$.

 \begin{figure}[!ht]
\centering
\includegraphics[scale=0.2]{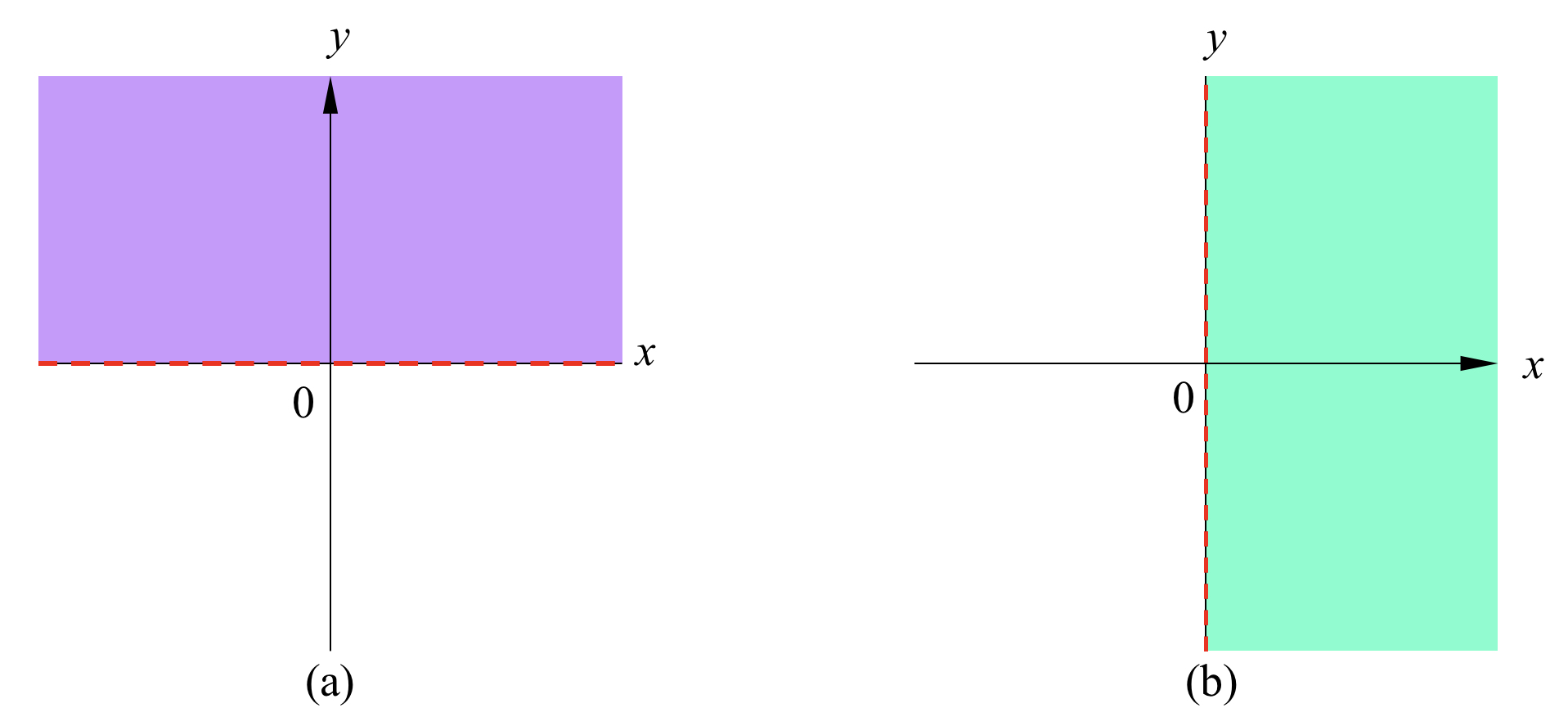}
\caption{The region $\mk{D}=\left\{(r\cos\theta, r\sin\theta)\,|\,r_1<r<r_2, \theta_1<\theta<\theta_2\right\}$  in the $(x,y)$-plane. (a) $r_1=0$, $r_2=\infty$, $\theta_1=0$, $\theta_2=\pi$. (b) $r_1=0$, $r_2=\infty$, $\theta_1=-\frac{\pi}{2}$, $\theta_2=\frac{\pi}{2}$.  }\label{figure123}
 
\includegraphics[scale=0.2]{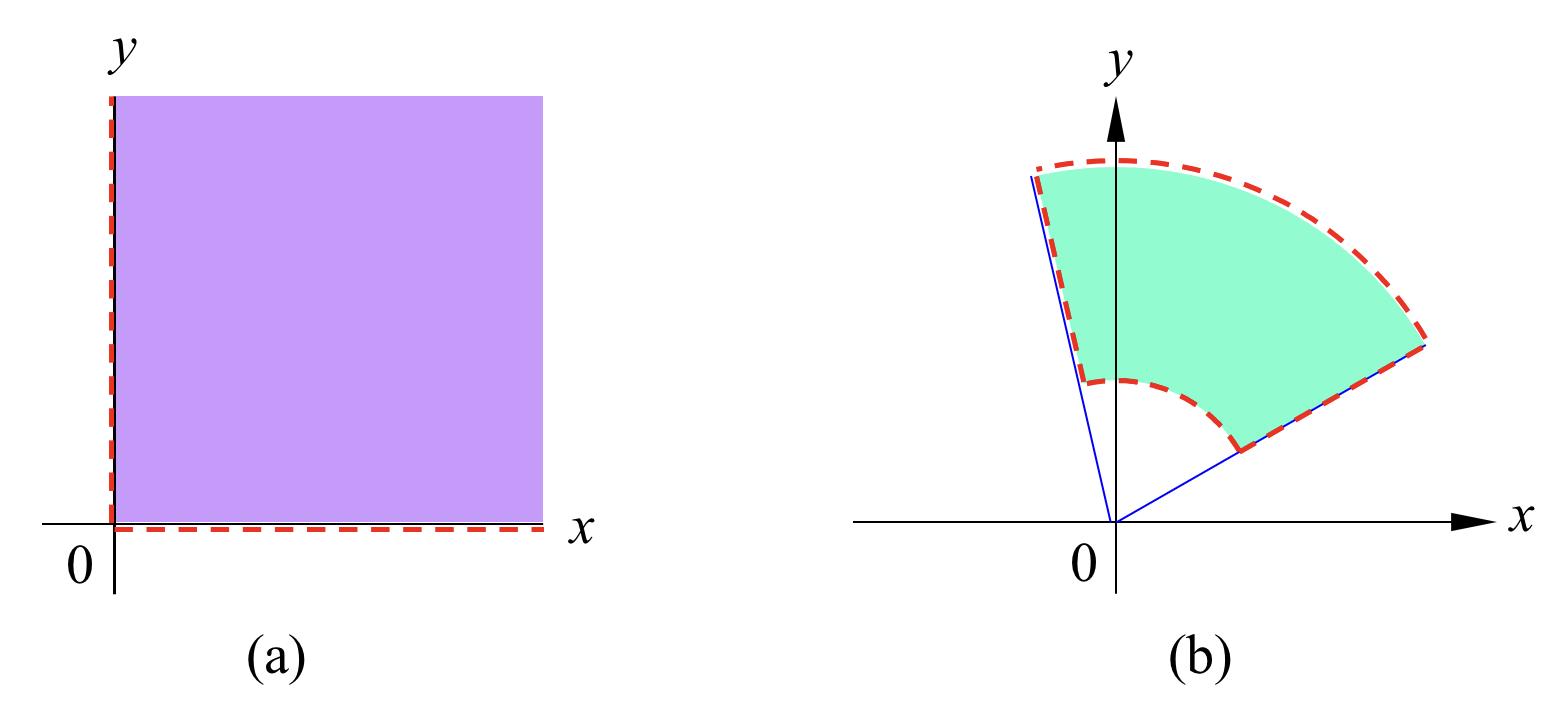}
\caption{The region $\mk{D}=\left\{(r\cos\theta, r\sin\theta)\,|\,r_1<r<r_2, \theta_1<\theta<\theta_2\right\}$  in the $(x,y)$-plane. (a) $r_1=0$, $r_2=\infty$, $\theta_1=0$, $\theta_2=\frac{\pi}{2}$. (b) $r_1=2$, $r_2=5$, $\theta_1= \frac{\pi}{6}$, $\theta_2=\frac{4\pi}{7}$.  }\label{figure124}
\end{figure}

To apply the change of variables theorem, we notice that the Jacobian is
\[\frac{\pa (x,y)}{\pa (r,\theta)}=r.\]
Thus, 
\[dxdy =\frac{\pa (x,y)}{\pa (r,\theta)}drd\theta=rdrd\theta.\]
The change of variables theorem says the following. 
\begin{theorem}{}

Let $[\alpha, \beta]$ be a closed interval such that $\beta\leq \alpha+2\pi$. Assume that $g:[\alpha,\beta]\to\mb{R}$ and $h:[\alpha, \beta]\to\mb{R}$ are continuous functions satisfying
\[0\leq g(\theta)\leq h(\theta)\hspace{1cm}\text{for all}\; \alpha\leq \theta\leq\beta.\]
Let $\mk{D}$ be the region in the $(x,y)$-plane given by
\[\mk{D}=\left\{(r\cos\theta, r\sin\theta)\,|\, \alpha\leq\theta\leq\beta,  g(\theta)\leq r\leq h(\theta)\right\}.\]
If $f:\mk{D}\to\mb{R}$ is a continuous function, then
\begin{equation}\label{230908_1}\int_{\mk{D}}f(x,y)dxdy=\int_{\alpha}^{\beta}\int_{g(\theta)}^{h(\theta)}f(r\cos\theta, r\sin\theta)rdrd\theta.\end{equation}

\end{theorem}
\begin{myproof}{Proof}
Let
\[\mathcal{U}=\left\{(r,\theta)\,|\, \alpha\leq\theta\leq\beta, g(\theta)\leq r\leq h(\theta)\right\}.\]
Then $\mathcal{U}$ is a compact Jordan measurable set in $\mb{R}^2$.

 If 
$\di\beta<\alpha+2\pi$, take any $\alpha_0$ such that 
\[\alpha_0<\alpha<\beta<\alpha_0+2\pi.\]
For example, we can take \[\alpha_0=\alpha-\di \frac{2\pi -(\beta-\alpha)}{2}.\] 
If we also have
\[ g(\theta)>0\hspace{1cm}\text{for all}\;\theta\in [a,b],\]
then $\mathcal{U}$ is contained in the set
\[\mathcal{O}_{\alpha_0}=\left\{(r,\theta)\,|\,r>0, \alpha_0<\theta<\alpha_0+2\pi\right\},\] 
 
and $\mf{\Phi}(\mathcal{U})=\mk{D}$.
Applying the change of variables theorem to the mapping $\mf{\Phi}:\mathcal{O}_{\alpha_0}\to\mb{R}^2$ gives the desired formula \eqref{230908_1} immediately. \bp

If $g(\theta)=0$ for some $\theta\in [\alpha, \beta]$, then we consider the set
\[\mathcal{U}_{\varepsilon}=\left\{(r,\theta)\,|\, \alpha\leq\theta\leq\beta, g(\theta)+\varepsilon\leq r\leq h(\theta)+\varepsilon\right\},\quad\text{where}\;\varepsilon>0.\]
It is contained in $\mathcal{O}_{\alpha_0}$. Using boundedness of the continuous functions $g:[\alpha,\beta]\to\mb{R}$, $h:[\alpha,\beta]\to\mb{R}$ and $f:\mk{D}\to\mb{R}$, it is easy to show that
\[\int_{\mk{D}}f(x,y)dxdy=\lim_{\varepsilon\to 0^+}\int_{\mf{\Phi}(\mathcal{U}_{\epsilon})}f(x,y)dxdy.\]
By the change of variables formula, we have
\[\int_{\mf{\Phi}(\mathcal{U}_{\epsilon})}f(x,y)dxdy=\int_{\alpha}^{\beta}\int_{g(\theta)+\varepsilon}^{h(\theta)+\varepsilon}f(r\cos\theta, r\sin\theta)rdrd\theta.\] Taking the $\varepsilon\to 0^+$ limit yields again the desired formula \eqref{230908_1}.

The last case we have to consider is when $\beta=\alpha+2\pi$.  The technicality is that $\mathcal{U}$ is not contained in any of the $\mathcal{O}_{\alpha_0}$ restricted to which the mapping $\mf{\Phi}(r,\theta)=(r\cos\theta, r\sin\theta)$ is a smooth change of variables. Instead of taking limits, there is an  alternative way to resolve the problem. We write $\mathcal{U}=\mathcal{U}_1\cup\mathcal{U}_2$, where
\begin{align*}
\mathcal{U}_1&=\left\{(r,\theta)\,|\, \alpha\leq\theta\leq\alpha+\pi, g(\theta)\leq r\leq h(\theta)\right\},\\
\mathcal{U}_2&=\left\{(r,\theta)\,|\, \alpha+\pi\leq\theta\leq\alpha+2\pi, g(\theta)\leq r\leq h(\theta)\right\},
\end{align*}We have shown that the change of variables formula \eqref{230908_1} is valid for $\mathcal{U}_1$ and $\mathcal{U}_2$. Apply the additivity theorem, we find that
\begin{align*}
\int_{\mk{D}}f(x,y)dxdy&=\int_{\mf{\Psi}(\mathcal{U}_1)}f(x,y)dxdy+\int_{\mf{\Psi}(\mathcal{U}_2)}f(x,y)dxdy\\
&=\int_{\alpha}^{\alpha+\pi}\int_{g(\theta)}^{h(\theta)}f(r\cos\theta, r\sin\theta)rdrd\theta\\&\quad +\int_{\alpha+\pi}^{\alpha+2\pi}\int_{g(\theta)}^{h(\theta)}f(r\cos\theta, r\sin\theta)rdrd\theta\\
&=\int_{\alpha}^{\beta}\int_{g(\theta)}^{h(\theta)}f(r\cos\theta, r\sin\theta)rdrd\theta.
\end{align*}Namely, the formula  \eqref{230908_1} is still valid when $\beta=\alpha+2\pi$.  
\end{myproof}

Let us give a
geometric explanation for the Jacobian
\[\frac{\pa(x,y)}{\pa (r,\theta)}=r,\quad\text{where}\;\;x=r\cos\theta, y=r\sin\theta.\]
Assume that \[\theta_1<\theta_2<\theta_1+2\pi\quad\text{and}\quad  0<r_1<r_2.\] Let
\[\mathcal{D}=\left\{(r\cos\theta,r\sin\theta)\,|\,  r_1\leq r\leq r_2, \theta_1\leq \theta\leq \theta_2\right\}.\]
The area bounded between the circles $x^2+y^2=r_1^2$ and $x^2+y^2=r_2^2$ is $\pi ( r_2^2-r_1^2)$. By rotational symmetry of the circle, the area of $\mathcal{D}$ is
\[\Delta A=\pi(r_2^2-r_1^2)\times\frac{\theta_2-\theta_1}{2\pi}=\overline{r}\Delta r\Delta\theta,\]where
\[\overline{r}=\frac{r_1+r_2}{2},\quad \Delta r=r_2-r_1\quad\text{and}\quad \Delta\theta=\theta_2-\theta_1.\]
When $\Delta r\to 0$, then 
\[\Delta A\sim r\Delta r\Delta \theta,\]where $r\sim r_1\sim r_2$.

 \begin{figure}[!ht]
\centering
\includegraphics[scale=0.2]{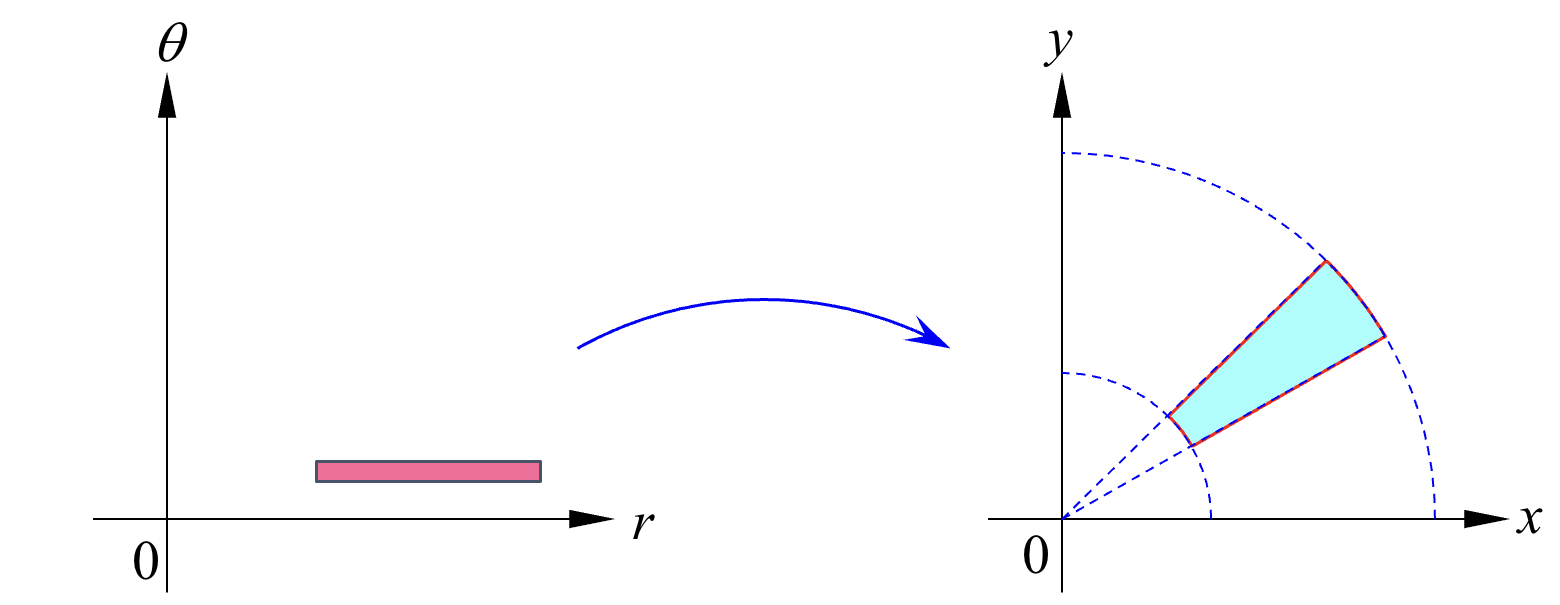}
\caption{The rectangle $[r_1,r_2]\times [\theta_1, \theta_2]$ in the $(r,\theta)$-plane is mapped to the region $\left\{(r\cos\theta, r\sin\theta)\,|\,r_1\leq r\leq r_2, \theta_1\leq \theta\leq \theta_2\right\}$  in the $(x,y)$-plane.  }\label{figure125}
\end{figure}
Now let us  look at some examples.

\begin{example}{}
 Evaluate the integral $\di \int_{\mk{D}}(2x^2+y^2)dxdy$,  where 
 \[\mk{D}=\left\{(x,y)\,|\, x\geq 0, y\geq 0, 4x^2+9y^2\leq 36\right\}.\]
\end{example}
\begin{solution}{Solution}
Making a change of variables $x=3u$, $y=2v$, the Jacobian is
\[\frac{\pa(x,y)}{\pa(u,v)}=\det\begin{bmatrix} 3 & 0\\0 & 2\end{bmatrix} =6.\]
Then
\begin{align*}
\int_{\mk{D}}(2x^2+y^2)dxdy&=\int_{\mathcal{U}} \left(18u^2+4v^2\right)\left|\frac{\pa(x,y)}{\pa(u,v)}\right|dudv\\&=6\int_{\mathcal{U}} \left(18u^2+4v^2\right)dudv,\end{align*} where $\mathcal{U}$ is 
 the region \[\mathcal{U}=\left\{(u,v)\,|\,u\geq 0, v\geq 0, u^2+v^2\leq 1\right\}.\] Since $\mathcal{U}$ is symmetric if we interchange $u$ and $v$, we find that
\[\int_{\mathcal{U}}u^2dudv=\int_{\mathcal{U}}v^2dudv=\frac{1}{2}\int_{\mathcal{U}}(u^2+v^2)dudv.\]
Using polar coordinates $u=r\cos\theta$, $v=r\sin\theta$, we find that
\begin{align*}
 \int_{\mathcal{U}}(u^2+v^2)dudv&=\int_0^{\frac{\pi}{2}}\int_0^1 r^2\times rdrd\theta= \frac{\pi}{2} \int_0^1 r^3dr=\frac{\pi}{8}.
\end{align*}
Therefore,
\[\int_{\mk{D}}(2x^2+y^2)dxdy=\frac{6\times (18+4)}{2} \int_{\mathcal{U}}(u^2+v^2)dudv=\frac{33\pi}{4}.\]
\end{solution}

 \begin{figure}[!ht]
\centering
\includegraphics[scale=0.2]{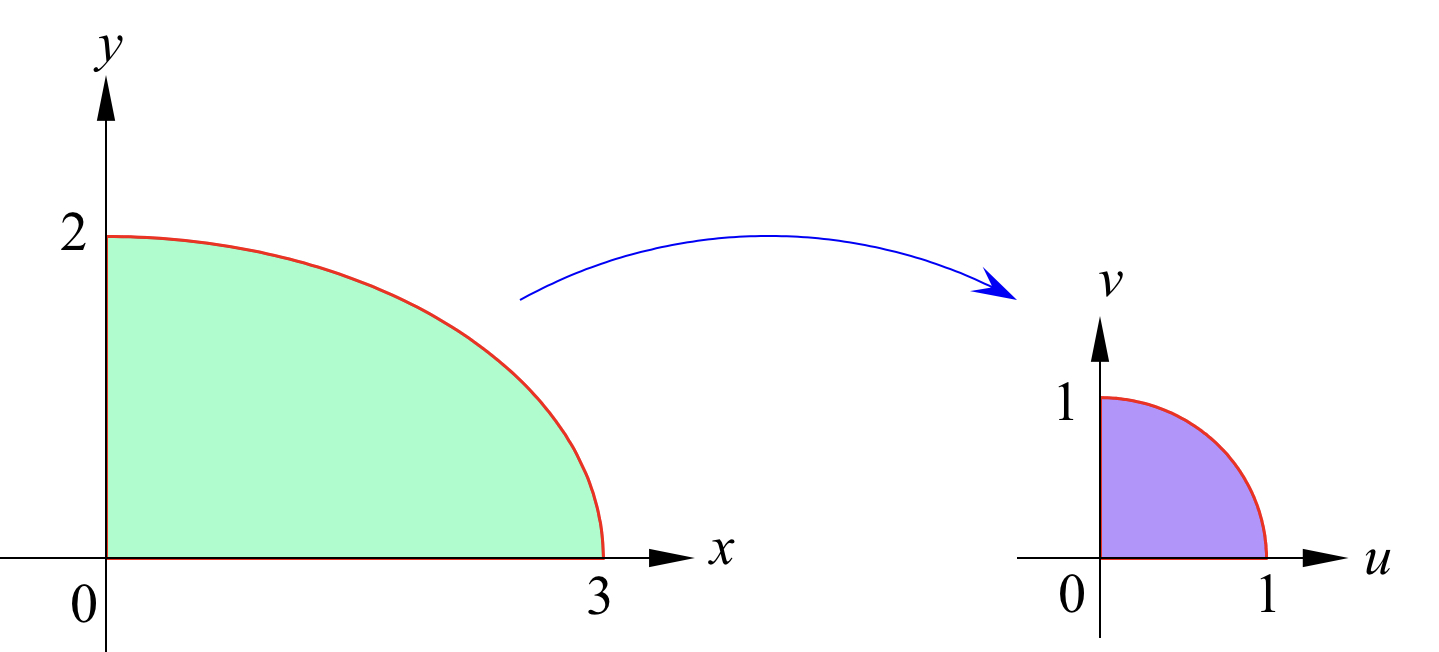}
\caption{ The   regions $\mk{D}=\left\{(x,y)\,|\, x\geq 0, y\geq 0, 4x^2+9y^2\leq 36\right\}$ and $\mathcal{U}=\left\{(u,v)\,|\,u\geq 0, v\geq 0, u^2+v^2\leq 1\right\}$.  }\label{figure127}
\end{figure}

\begin{example}{}
Find the volume of the solid bounded between the surface $z=x^2+y^2$ and the plane $z=9$.

\end{example}
\begin{solution}{Solution}
The solid $\mathcal{S}$ bounded between the surface $z=x^2+y^2$ and the plane $z=9$ can be expressed as
\[\mathcal{S}=\left\{(x,y,z)\,|\, (x,y)\in\mk{D}, x^2+y^2\leq z\leq 9\right\},\]where\[\mk{D}=\left\{(x,y)\,|\,x^2+y^2\leq 9\right\}.\]Since $\mk{D}$ is a closed ball, it is a Jordan measurable set. The volume of $\mathcal{S}$ is the integral of the constant function $\chi_{\mathcal{S}}:\mathcal{S}\to\mb{R}$. It is a continuous function. By Fubini's theorem,
\begin{align*}\text{vol}\,(\mathcal{S})
&=\int_{\mk{D}}\left(\int_{x^2+y^2}^9 dz \right)dxdy=\int_{\mk{D}}(9-x^2-y^2)dxdy\end{align*}
Using polar coordinatets, we have
\begin{align*}\text{vol}\,(\mathcal{S})&=\int_0^{2\pi}\int_0^3(9-r^2)rdrd\theta=2\pi \left[\frac{9r^2}{2}-\frac{r^4}{4}\right]_0^3=\frac{81\pi}{2}.
\end{align*}
\end{solution}

 \begin{figure}[!ht]
\centering
\includegraphics[scale=0.2]{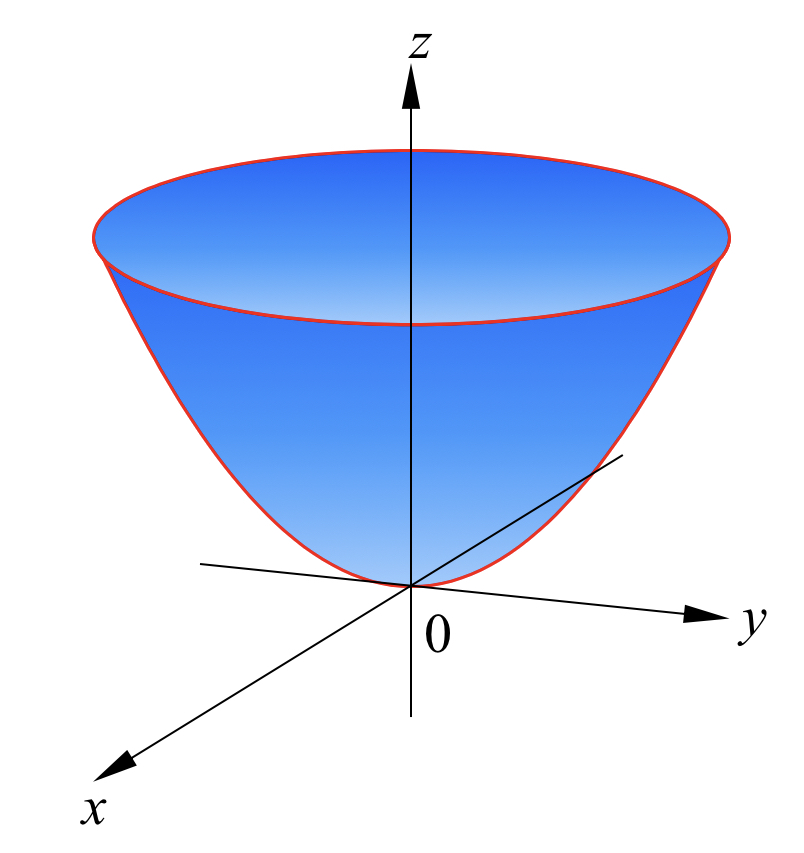}
\caption{The solid bounded between the surface $z=x^2+y^2$ and the plane $z=9$. }\label{figure126}
\end{figure}

\begin{example}{}
Let $a$ be a positive number.
Find the volume of ball \[B=\left\{(x,y,z)\,|\, x^2+y^2+z^2\leq a^2\right\}.\]
 
\end{example}
\begin{solution}{Solution}
Let \[\mk{D}=\left\{(x,y)\,|\,x^2+y^2\leq a^2\right\}.\] 
Then the ball $B$ can be decsribed as
\[B=\left\{(x,y,z)\,|\, (x,y)\in\mk{D}, -\sqrt{a^2-x^2-y^2}\leq z\leq \sqrt{a^2-x^2-y^2}\right\}.\]
Thus, by Fubini's theorem, its volume is
\begin{align*}
\text{vol}\,(B)&=\int_{B}dxdydz=\int_{\mk{D}}\left(\int_{-\sqrt{a^2-x^2-y^2}}^{\sqrt{a^2-x^2-y^2}}dz\right)dxdy\\&=\int_{\mk{D}} 2\sqrt{a^2-x^2-y^2} dxdy\end{align*}
\bs
Using polar coordinates, \[\text{vol}\,(B)
  =2\int_0^{2\pi}\int_0^a \sqrt{a^2-r^2}\,rdrd\theta=4\pi\int_0^ar\sqrt{a^2-r^2}\,dr.\]
 
Let $u=a^2-r^2$. Then $du=-2rdr$. When $r=0$, $u=a^2$. When $r=a$, $u=0$. Therefore,
\begin{align*}
\text{vol}\,(B)=2\pi \int_0^{a^2} u^{\frac{1}{2}}du=2\pi \left[ \frac{2}{3} u^{\frac{3}{2}}\right]_0^{a^2}=\frac{4\pi}{3}a^3.
\end{align*}
\end{solution}

\begin{example}[label=230909_1]{}
Let $a$ be a positive number, and let $\alpha$ be a number in the interval $(0, \frac{\pi}{2})$.
 Find the volume of the solid $E$ bounded between the sphere \[S=\left\{(x,y,z)\,|\, x^2+y^2+z^2=a^2\right\}\] and the cone
 \[C=\left\{(x,y,z)\,|\, z=\cot\alpha \sqrt{x^2+y^2}\right\}.\]
\end{example}
\begin{solution}{Solution}
The surfaces $S$ and $C$ intersect at the points $(x,y,z)$ satisfying
\[(x^2+y^2)(1+\cot^2\alpha)=a^2.\]
Namely,
\[x^2+y^2=a^2\sin^2\alpha.\]Therefore,
\[E=\left\{(x,y,z)\,|\, (x,y)\in\mk{D}, \cot\alpha \sqrt{x^2+y^2}\leq z\leq \sqrt{a^2-x^2-y^2}\right\},\]
where
\[\mk{D}=\left\{(x,y)\,|\,x^2+y^2\leq a^2\sin^2\alpha\right\}.\]
\bs
Using Fubini's theorem and polar coordinates, we find that
\begin{align*}
\text{vol}\,(E)&=\int_Edxdydz=\int_{\mk{D}}\left(\int_{\cot \alpha \sqrt{x^2+y^2}}^{\sqrt{a^2-x^2-y^2}}dz\right)dxdy\\
&=\int_0^{2\pi}\int_0^{a\sin\alpha}\left(\sqrt{a^2-r^2}-r\cot\alpha \right)rdrd\theta\end{align*}
 
Using a change of variables $u=a^2-r^2$, we find that
\begin{align*}
\int_0^{a\sin\alpha} r\sqrt{a^2-r^2}dr &=\frac{1}{2}\int_{a^2\cos^2\alpha}^{a^2}u^{\frac{1}{2}} du\\&=\frac{1}{3}\left[u^{\frac{3}{2}}\right]_{a^2\cos^2\alpha}^{a^2}=\frac{a^3}{3}\left(1-\cos^3\alpha\right).
\end{align*}
On the other hand
\begin{align*}
\int_0^{a\sin\alpha}r^2\cot\alpha \,dr&= \cot\alpha\left[\frac{r^3}{3}\right]_0^{a\sin\alpha} \\&=\frac{a^3}{3}\frac{\cos\alpha}{\sin\alpha}\sin^3\alpha =
\frac{a^3}{3}(\cos\alpha-\cos^3\alpha).
\end{align*}
Therefore, the volume of $E$ is
\[\text{vol}\,(E)=\frac{2\pi a^3}{3}(1-\cos\alpha).\]
\end{solution}

 \begin{figure}[!ht]
\centering
\includegraphics[scale=0.2]{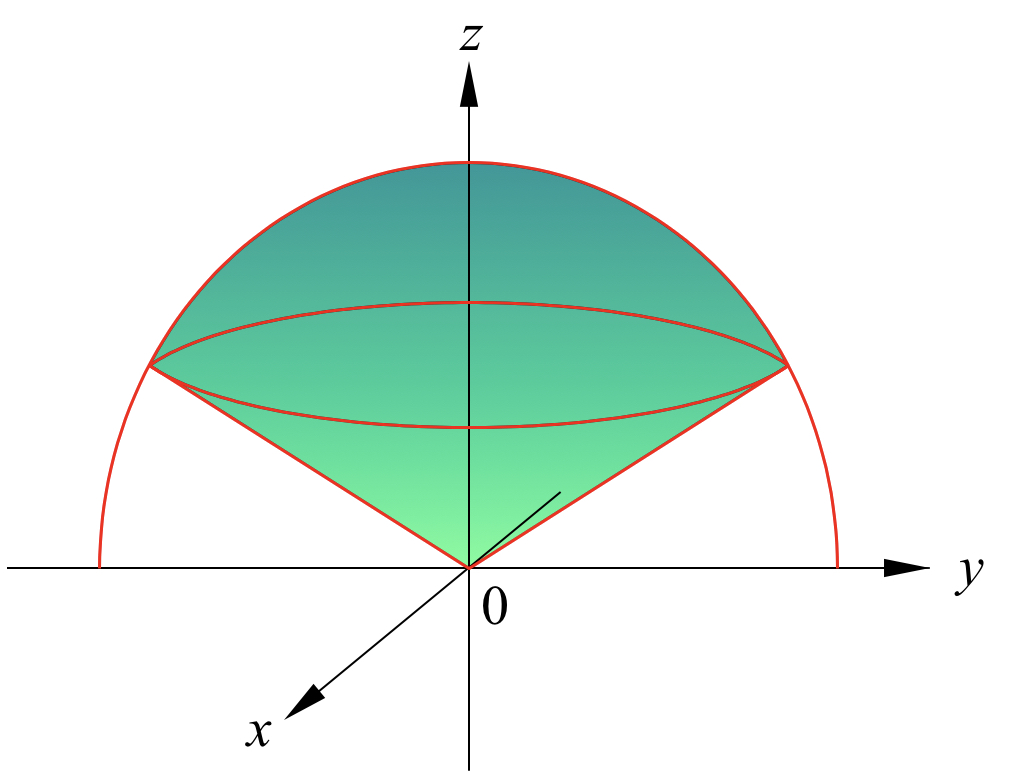}
\caption{ The   solid bounded between a sphere and a cone.  }\label{figure128}
\end{figure}
\subsection{Spherical Coordinates} 
Now we consider the spherical coordinates, which is an alternative coordinate system for $\mb{R}^3$. Consider the mapping
$\mf{\Psi}:\mb{R}^3\to\mb{R}^3$ given by 
\[\mf{\Psi} (\rho, \phi, \theta)=(x,y,z)=( \rho\sin\phi\cos\theta, \rho\sin\phi\sin\theta, \rho\cos\phi).\]
Namely,
\begin{align*}
x&=\rho\sin\phi \cos\theta,\\
y&=\rho\sin\phi\sin\theta,\\
z&=\rho\cos\phi.
\end{align*}
Let $V$ be the set
\[V=\left\{(\rho, \phi, \theta)\,|\,\rho>0, 0\leq \phi\leq\pi, 0\leq \theta<2\pi\right\}.\]  Given $\mf{u}=(x,y,z)\in\mb{R}^3\setminus\{(0,0,0)\}$, we claim that there is a unique $(\rho, \phi ,\theta)\in V$ such that $\mf{\Psi}(\rho,  \phi, \theta)=(x,y,z)$. This triple $(\rho,  \phi,\theta)$ is called a spherical coordinates of the point $\mf{u}=(x,y,z)$. It is easy to see that
\[\rho=\sqrt{x^2+y^2+z^2}=\Vert\mf{u}\Vert\] is the distance from the point $\mf{u}=(x,y,z)$ to the origin. If we let 
\[\phi=\cos^{-1}\frac{z}{\rho},\]
$\phi$ satisfies $0\leq \phi\leq\pi$, and 
\[\langle \mf{u}, \mf{e}_3\rangle =z=\rho \cos\phi=\Vert\mf{u}\Vert\cos\phi.\] Thus, geometrically, 
$\phi$ is the angle the vector  from $\mf{0}$ to $\mf{u}$ makes with the positive $z$-axis. Let $W$ be the $(x,y)$-plane in $\mb{R}^3$. Then 
\[(x,y,0)=\text{proj}_{W}\mf{u}.\]
Let
\[r=\rho\sin\phi.\]
Then 
\[r=\sqrt{\rho^2-\rho^2\cos^2\phi}=\sqrt{\rho^2-z^2}=\sqrt{x^2+y^2}.\]

 \begin{figure}[!ht]
\centering
\includegraphics[scale=0.2]{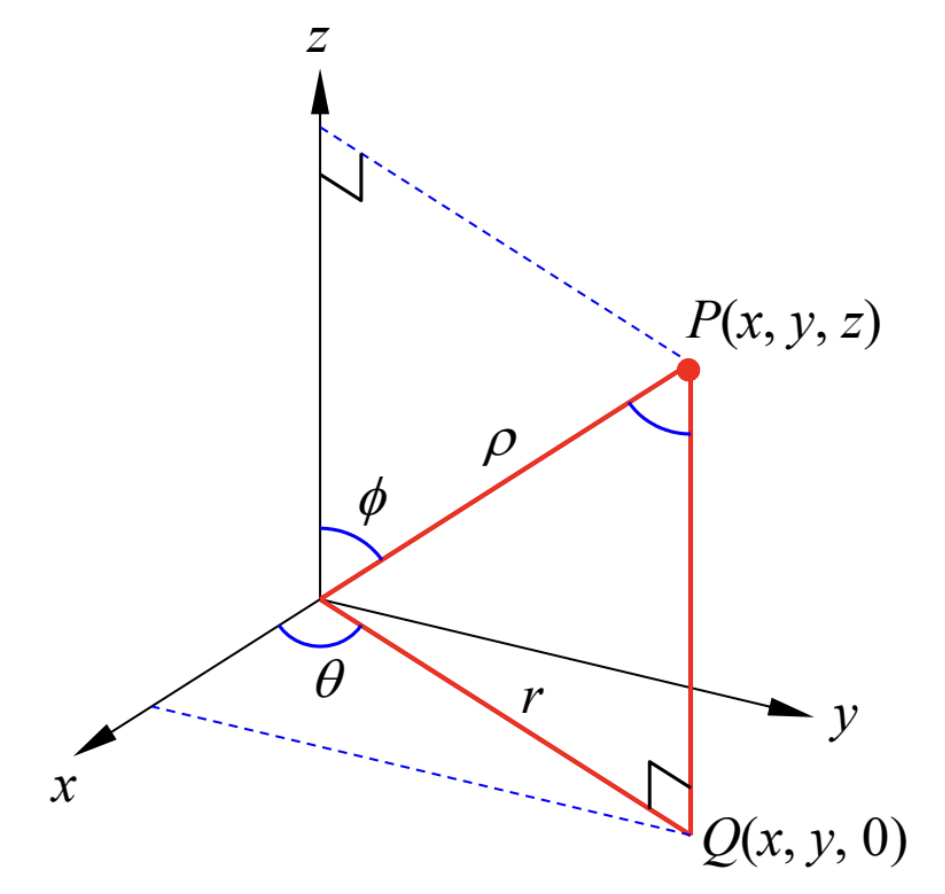}
\caption{Spherical Coordinates.  }\label{figure129}
\end{figure}
\noindent
Thus, $\theta\in [0,2\pi)$ is uniquely determined so that
\[x=r\cos\theta,\quad y=r\sin\theta.\] Equivalently, $(r,\theta)$ is the polar coordinates of the point $(x,y)$ in $\mb{R}^2\setminus\{(0,0)\}$. 
Hence, we find that   the map $\mf{\Psi}:V\to \mb{R}^3$,
\[\mf{\Psi} (\rho, \phi,\theta)=(x,y,z)=( \rho\sin\phi\cos\theta, \rho\sin\phi\sin\theta, \rho\cos\phi)\] is one-to-one on the set $V$, and the  range is $\mb{R}^3\setminus\{(0,0,0)\}$.
 However, the inverse map is not continuous. 

Let us calculate the derivative matrix of $\mf{\Psi}:\mb{R}^3\to\mb{R}^3$.
We find that
\begin{align*}
\mf{D\Psi}(\rho, \phi, \theta)&=\begin{bmatrix} \sin\phi\cos\theta  & \rho\cos\phi\cos\theta &-\rho\sin\phi\sin\theta\\ \sin\phi\sin\theta
& \rho\cos\phi\sin\theta& \rho\sin\phi\cos\theta \\ \cos\phi  & -\rho\sin\phi&  0\end{bmatrix}.
\end{align*}
Therefore, the Jacobian $\det \mf{D\Psi}(\rho, \phi, \theta)$  is
\begin{align*}
 \frac{\pa(x,y,z)}{\pa (\rho,\phi, \theta)} &=\cos\phi\times \det \begin{bmatrix}   \rho\cos\phi\cos\theta &-\rho\sin\phi\sin\theta\\   \rho\cos\phi\sin\theta& \rho\sin\phi\cos\theta  \end{bmatrix}
  \\&\quad +\rho\sin\phi \times \det \begin{bmatrix} \sin\phi\cos\theta    &-\rho\sin\phi\sin\theta\\ \sin\phi\sin\theta
 & \rho\sin\phi\cos\theta  \end{bmatrix}\\
 &=\rho^2\cos^2\phi\sin\phi+\rho^2\sin^3\phi=\rho^2\sin\phi.
\end{align*}
This shows that $\mf{D\Psi}(\rho, \phi, \theta)$ is invertible if and only if $\rho\neq 0$ and $\sin\phi\neq 0$, if and only if $(x,y,z)$ does not lie on the $z$ axis. 
Thus, for any real number $\alpha$, if $\mathcal{O}_{\alpha}$ is the open set
\[\mathcal{O}_{\alpha}=\left\{(\rho,\phi,\theta)\,|\,\rho>0, 0<\phi<\pi, \alpha<\theta<\alpha+2\pi\right\},\] then
$\mf{\Psi}:\mathcal{O}_{\alpha}\to\mb{R}^3$ is a smooth change of variables.

The change of variables theorem gives the following.
\begin{theorem}{}
Let $[\alpha,\beta]$ and $[\delta, \eta]$ be a closed intervals such that \[\beta\leq\alpha+2\pi\quad \text{and}\quad 0\leq\delta<\eta\leq \pi,\] and let $\mf{I}= [\delta, \eta]\times[\alpha,\beta]$. Assume that the functions $g:\mf{I}\to\mb{R}$ and $h:\mf{I}\to\mb{R}$ satisfy $0\leq g(\phi,\theta)\leq h(\phi,\theta)$ for all $(\phi,\theta)\in\mf{I}$, let $\mk{D}$ be the region in $\mb{R}^3$ defined by
\[\mk{D}=\left\{(\rho\sin\phi\cos\theta,\rho\sin\phi\sin\theta, \rho\cos\phi)\,|\, (\phi,\theta)\in\mf{I}, g(\phi,\theta)\leq\rho\leq h(\phi,\theta)\right\}.\]
If $f:\mk{D}\to\mb{R}$ is a continuous function, then
\begin{align*}&\int_{\mk{D}}f(x,y,z)dxdydz\\&=\int_{\alpha}^{\beta}\int_{\delta}^{\eta}\int_{g(\phi,\theta)}^{h(\phi,\theta)}f(\rho\sin\phi\cos\theta,\rho\sin\phi\sin\theta, \rho\cos\phi)\rho^2\sin\phi \,d\rho d\phi d\theta.
\end{align*}
\end{theorem}

Again, if $\beta<\alpha+2\pi$, $\delta>0$, $\eta<\pi$ and $g(\phi,\theta)>0$ for all $(\phi,\theta)\in\mf{I}$, this is just a direct consequence of the general change of variables theorem.  The rest can be argued by taking limits. 

The results of Example \ref{230909_1} can be used to give a hindsight about the Jacobian 
\[\frac{\pa (x,y,z)}{\pa(\rho,\phi,\theta)}=\rho^2\sin\phi\] that appears in the change from spherical coordinates to rectangular coordinates.
Consider the rectangle  $\mf{I}=[\rho_1,\rho_2]\times [\phi_1,\phi_2]\times [\theta_1, \theta_2]$ in the $(\rho,\phi,\theta)$ space, where $\theta_2<\theta_1+2\pi$, and for simplicity, assume that $0<\phi_1<\phi_2<\frac{\pi}{2}$. Under the mapping \[\mf{\Psi}(\rho,\phi,\theta)= (\rho\sin\phi\cos\theta,\rho\sin\phi\sin\theta, \rho\cos\phi),\]
$\mf{\Psi}(\mf{I})$ is a wedge in the solid $E$ in $\mb{R}^3$ bounded between the spheres $x^2+y^2+z^2=\rho_1^2$, $x^2+y^2+z^2=\rho_2^2$, and the cones
$z=\cot\phi_1\sqrt{x^2+y^2}$, $z=\cot\phi_2\sqrt{x^2+y^2}$. 
Since $E$ has a rotational symmetry with respect to $\theta$,
\[\Delta V=\text{vol}\,(\mf{\Psi}(\mf{I}))=\frac{\theta_2-\theta_1}{2\pi}\,\text{vol}\,(E).\]
Using inclusion and exclusion principle, the result of Example \ref{230909_1} gives
\begin{align*}
\text{vol}\,(E)&=\frac{2\pi}{3}\rho_2^3(1-\cos\phi_2)-\frac{2\pi}{3}\rho_1^3(1-\cos\phi_2)\\&-\frac{2\pi}{3}\rho_2^3(1-\cos\phi_1)+\frac{2\pi}{3}\rho_1^3(1-\cos\phi_1)\\
&=\frac{2\pi}{3}(\rho_2^3-\rho_1^3)(\cos\phi_1-\cos\phi_2).
\end{align*}

\begin{figure}[ht]
\centering
\includegraphics[scale=0.2]{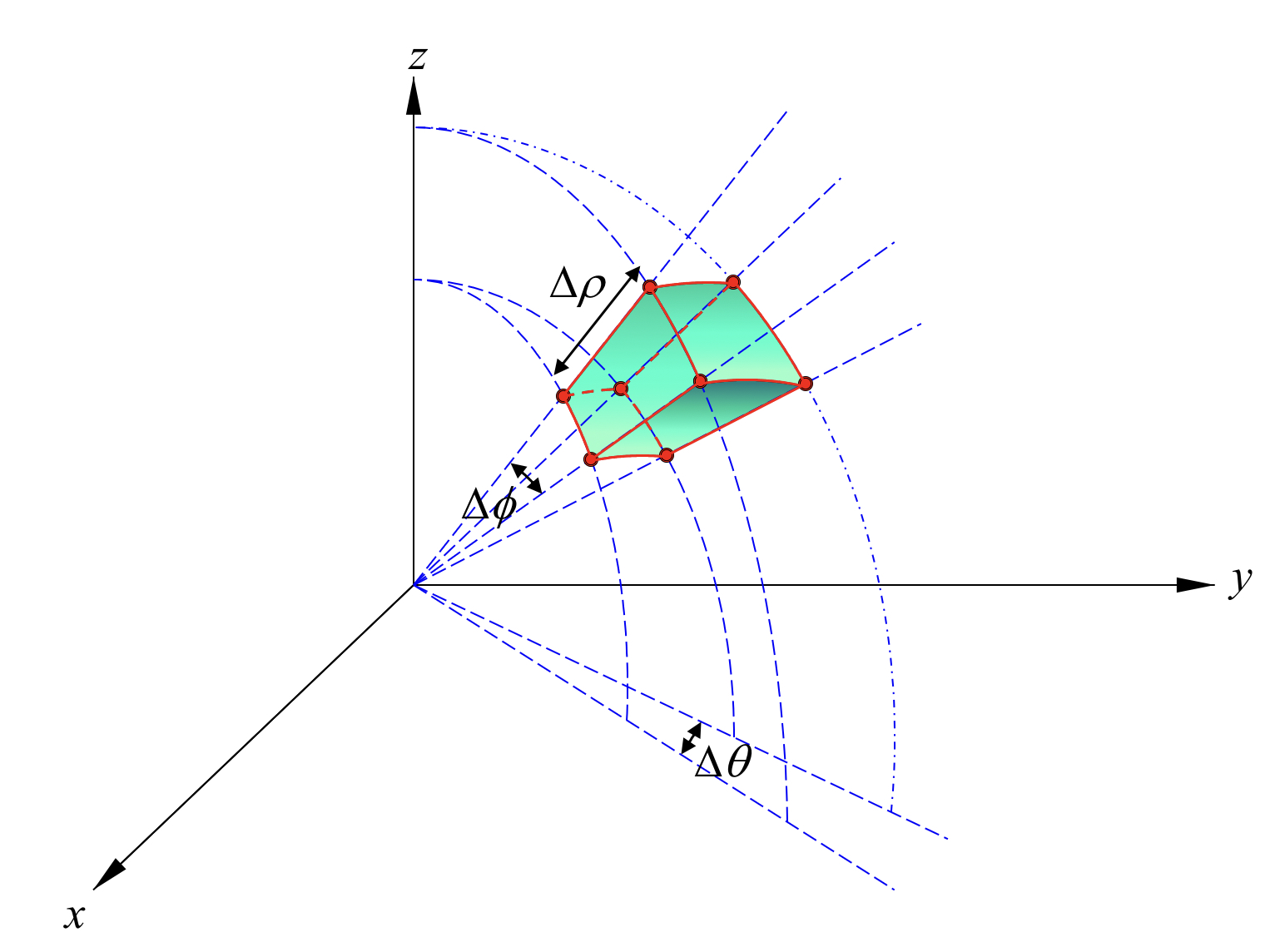}
\caption{Volume change under spherical coordinates. }\label{figure131}
\end{figure}

By mean value theorem there is a  $\overline{\rho}\in (\rho_1,\rho_2)$ and a $\overline{\phi}\in (\phi_1, \phi_2)$ such that
\[\rho_2^3-\rho_1^3=3\overline{\rho}^2\Delta\rho\quad\text{and}\quad \cos\phi_1-\cos\phi_2=\sin\overline{\phi}\Delta\phi,\]
where
\[\Delta\rho=\rho_2-\rho_1,\hspace{1cm}\Delta\phi=\phi_2-\phi_1.\]
Let $\Delta\theta=\theta_2-\theta_1$. Then we find that
\[\Delta V=\overline{\rho}^2\sin\overline{\phi}\,\Delta\rho\Delta \phi\Delta\theta.\]
This gives an interpretation of the Jacobian $\rho^2\sin\phi$.

Let us look at an example of applying spherical coordinates.
\begin{example}{}
Compute the integral $\di \int_{E}(x^2+4z)dxdydz$, where 
\[E=\left\{(x,y,z)\,|\,   x^2+y^2+z^2\leq 9, z\geq 0\right\}.\]

\end{example}
\begin{solution}{Solution}
Let $B$ be the sphere 
\[B=\left\{(x,y,z)\,|\,   x^2+y^2+z^2\leq 9\right\}.\]
By symmetry,
\[\int_Ex^2dxdydz=\frac{1}{2}\int_B x^2dxdydz=\frac{1}{6}\int_B(x^2+y^2+z^2)dxdydz.\]
Using spherical coordinates, we have 
\begin{align*}
\int_Ex^2dxdydz&=\frac{1}{6}\int_0^{2\pi}\int_0^{\pi}\int_0^3 \rho^2\times\rho^2\sin\phi \,d\rho d\phi d\theta\\
&=\frac{\pi}{3}\left[-\cos\phi\right]_0^{\pi}\left[\frac{\rho^5}{5}\right]_0^3=\frac{162\pi}{5}.
\end{align*}
On the other hand,
\begin{align*}
\int_Ezdxdydz&=\int_0^{2\pi}\int_0^{\frac{\pi}{2}}\int_0^3\rho\cos\phi\times\rho^2\sin\phi \,d\rho d\phi d\theta\\
&=2\pi \left[-\frac{\cos^2 \phi}{2}\right]_0^{\frac{\pi}{2}}\left[\frac{\rho^4}{4}\right]_0^3=\frac{81\pi}{4}. 
\end{align*}
\bs
Therefore,
\[\int_{E}(x^2+4z)dxdydz=\frac{162\pi}{5}+81\pi =\frac{567\pi}{5}.\]
\end{solution}
In the example above, we have used the symmetry of the region $E$ to avoid some complicated computations. Another example is the following.
\begin{example}{}Let $a$ be a positive number.
Evaluate the integral $\di\int_E x^4 dxdydz$, where
\[E=\left\{(x,y,z)\,|\,x\geq 0, y\geq 0, z\geq 0, x^2+y^2+z^2\leq a^2\right\}.\]

\end{example}
\begin{solution}{Solution}
The expression of the $z$ variable in terms of the spherical coordinates  is considerably simpler than the $x$ and $y$ variables.
By symmetry, we have
\[\int_E x^4 dxdydz=\int_E z^4 dxdydz.\]
Thus, using spherical coordinates, we find that
\begin{align*}
\int_E x^4 dxdydz=&\int_0^{\frac{\pi}{2}}\int_0^{\frac{\pi}{2}}\int_0^a\rho^4\cos^4\phi\times \rho^2\sin\phi \,d\rho d\phi d\theta\\
&=\frac{\pi}{2}\left[-\frac{\cos^5\theta}{5}\right]_0^{\frac{\pi}{2}}\left[\frac{\rho^7}{7}\right]_0^a=\frac{\pi a^7}{70}.
\end{align*}
\end{solution}

\subsection{Other Examples}
 
\begin{example}{}
Let $\mk{D}$ be the region

\vspace{-0.4cm}
\[\mk{D}=\left\{(x,y)\,|\,y>0, 4\leq x^2-y^2\leq 9, 3\leq xy\leq 7 \right\}.\] Compute the integral
\[\int_{\mk{D}} (x^3y+xy^3)dxdy.\]
\end{example}

\begin{solution}{Solution}
Let $\mathcal{O}=\left\{(x,y)\,|\,y>0\right\}$, and let
$\mf{\Psi}:\mathcal{O}\to\mb{R}^2$ be the mapping
\[\mf{\Psi}(x,y)=(x^2-y^2, xy).\]
If $(x_1, y_1)$ and $(x_2, y_2)$ are points in $\mathcal{O}$ such that $\mf{\Psi}(x_1, y_1)=\mf{\Psi}(x_2, y_2)$, then \[x_1^2-y_1^2=x_2^2-y_2^2\quad\text{and}\quad x_1y_1=x_2y_2.\]
Let $z_1=x_1+iy_1$ and $z_2=x_2+iy_2$. Then
\[z_1^2=(x_1+iy_1)^2=(x_2+iy_2)^2=z_2^2.\]
This implies that  $z_2=\pm z_1$. Thus, $y_2=\pm y_1$. Since $y_1$ and $y_2$ are positive, we find that $y_1=y_2$. Since $x_1y_1=x_2y_2$, we then deduce that $x_1=x_2$. Hence, $\mf{\Psi}:\mathcal{O}\to\mb{R}^2$ is one-to-one. Since it is a polynomial mapping, it is continuously differentiable. Since
\[\mf{D\Psi}(x,y)=\begin{bmatrix} 2x & -2y\\y & x\end{bmatrix},\hspace{1cm}\det \mf{D\Psi}(x,y)=2(x^2+y^2),\]
we find that $\det \mf{D\Psi}(x,y)\neq 0$ for all $(x,y)\in\mathcal{O}$. This implies that $\mf{\Psi}:\mathcal{O}\to\mb{R}^2$ is a smooth change of variables. 
Let $u=x^2-y^2$, $v=xy$. The Jacobian is
\[\frac{\pa(u,v)}{\pa (x,y)}=2(x^2+y^2).\]
Notice that
\[\mf{\Psi}(\mk{D})=\left\{(u,v)\,|\, 4\leq u\leq 9, 3\leq v\leq 7\right\}.\]
Therefore,
\begin{align*}
\int_{\mk{D}}xy(x^2+y^2)dxdy&=\frac{1}{2}\int_{\mk{D}}xy
\frac{\pa (u,v)}{\pa(x,y)}dxdy\\
&=\frac{1}{2}\int_{\mf{\Psi}(\mk{D})}vdudv=\frac{1}{2}\int_3^7\int_4^9 vdudv=50.
\end{align*}
\end{solution}

\begin{figure}[ht]
\centering
\includegraphics[scale=0.2]{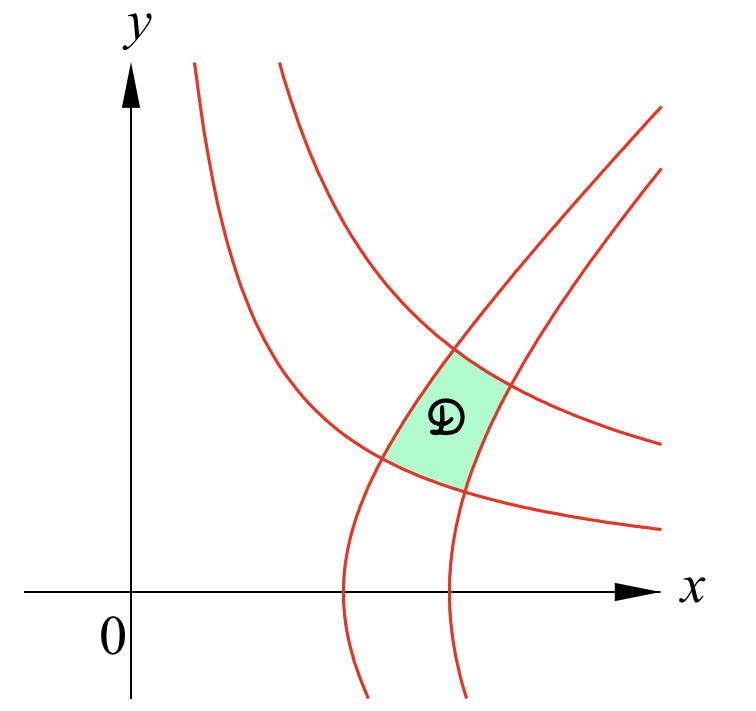}
\caption{The region $\di\mk{D}=\left\{(x,y)\,|\,y>0, 4\leq x^2-y^2\leq 9, 3\leq xy\leq 7 \right\}$. }\label{figure121}
\end{figure}

\begin{remark}{Hyperspherical Coordinates}
For any $n\geq 4$, the hyperspherical coorfinates in $\mb{R}^n$ are the coordinates $(r,\theta_1,\ldots,\theta_{n-1})$ such that
\begin{equation}\label{230909_2}
\begin{split}
x_1&=r\cos\theta_1,\\
x_2&=r\sin\theta_1\cos\theta_2,\\
x_3&=r\sin\theta_1\sin\theta_2\cos\theta_3,\\
&\hspace{1cm}\vdots\\
x_{n-1}&=r\sin\theta_1\cdots\sin\theta_{n-2}\cos\theta_{n-1}\\
x_{n}&=r\sin\theta_1\cdots\sin\theta_{n-2}\sin\theta_{n-1}.
\end{split}\end{equation}

Here
\[r=\sqrt{x_1^2+\cdots+x_n^2}.\]

If 
\[V=(0,\infty)\times [0,\pi)^{n-2}\times [0,2\pi),\]there is a one-to-one correspondence between $(r,\theta_1, \ldots,\theta_{n-1})$ in $V$ and $(x_1, \ldots, x_n)$ in $\mb{R}^n\setminus\{\mf{0}\}$ given by \eqref{230909_2}. One can show that the Jacobian of this transformation is
\[\frac{\pa(x_1, x_2, \ldots, x_n)}{\pa (r, \theta_1,\ldots,\theta_{n-1})}=r^{n-1}\sin^{n-2}\theta_1\cdots \sin\theta_{n-2}.\]

\end{remark}

\vp
\noindent
{\bf \large Exercises  \thesection}
\setcounter{myquestion}{1}

\begin{question}{\themyquestion}
 Let 
\[ \mk{D}=\left\{(x,y)\,|\,4(x+1)^2+9(y-2)^2\leq 144\right\}.\] Evaluate the integral $\di \int_{\mk{D}} (2x+3y)dxdy$.
\end{question}
\atc
\begin{question}{\themyquestion}
 Evaluate the integral
 \[\int_{\mk{D}}\frac{x+y}{(x-2y+8)^2}dxdy,\]
 where \[\mk{D}=\di\left\{(x,y)\,|\, |x|+2|y|\leq 7\right\}.\]
\end{question}
\atc

\begin{question}{\themyquestion}
Evaluate the integral $\di \int_{\mk{D}}(x^2-xy+y^2)dxdy$,  where 
 \[\mk{D}=\left\{(x,y)\,|\, x\geq 0,  x^2+9y^2\leq 36\right\}.\]
\end{question}
\atc

\begin{question}{\themyquestion}
Find the volume of the solid bounded between the surface $z=x^2+y^2$ and the surface $x^2+y^2+z^2=20$.
\end{question}
\atc

\begin{question}{\themyquestion}
Let $a$, $b$ and $c$ be positive numbers, and let $E$ be the solid
\[E=\left\{(x,y,z)\,\left|\,\frac{x^2}{a^2}+\frac{y^2}{b^2}+\frac{z^2}{c^2}\leq 1\right.\right\}.\]
Evaluate $\di\int_Ex^2dxdydz$.
\end{question}
\atc

\begin{question}{\themyquestion}
Let $a$, $b$ and $c$ be positive numbers, and let $E$ be the solid
\[E=\left\{(x,y,z)\,\left|\, x\geq 0, \frac{x^2}{a^2}+\frac{y^2}{b^2}+\frac{z^2}{c^2}\leq 1\right.\right\}.\]
Evaluate $\di\int_Ex^6dxdydz$.
\end{question}
 
\atc
\begin{question}{\themyquestion}
Let $\mk{D}$ be the region
\[\mk{D}=\left\{(x,y)\,|\,x>0, 1\leq x^2-y^2\leq 25, 1\leq xy\leq 6 \right\}.\] Compute the integral
\[\int_{\mk{D}} \frac{x^4-y^4}{xy}dxdy.\]
\end{question}
 
 \atc
\begin{question}{\themyquestion}
Let $\mk{D}$ be the region
\[\mk{D}=\left\{(x,y)\,|\,5x^2-2xy+10y^2\leq 9 \right\},\] and let
 $\mf{\Psi}:\mb{R}^2\to\mb{R}^2$ be the mapping defined by
\[\mf{\Psi}(x,y)=(2x+y, x-3y).\] 
\begin{enumerate}[(a)]
\item Explain why   $\mf{\Psi}:\mb{R}^2\to\mb{R}^2$   is a smooth change of variables.  
\item 
Find 
  $\mf{\Psi}(\mk{D})$.

  \item Compute the integral $\di \int_{\mk{D}}\frac{8}{5x^2-2xy+10y^2+16}\,dxdy$.

\end{enumerate} 
\end{question}
 
\section{Proof of the Change of Variables Theorem} 
In this section, we give a complete proof of the change of variables theorem, which we restate  here.

\begin{theorem}[label=230905_5]{The Change of Variables Theorem}
Let $\mathcal{O}$ be an open subset of $\mb{R}^n$, and let $\mf{\Psi}:\mathcal{O}\to\mb{R}^n$ be a smooth change of variables. If $\mk{D}$ is a Jordan measurable set such that its closure $\overline{\mk{D}}$ is contained in $\mathcal{O}$, then $\mf{\Psi}(\mk{D})$ is Jordan measurable, and for any function $f:\mf{\Psi}(\mk{D})\to\mb{R}$ that is  bounded and continuous, we have
\[\int_{\mf{\Psi}(\mk{D})}f(\mf{x})d\mf{x}=\int_{\mk{D}}f\left(\mf{\Psi}(\mf{x})\right)\left|\det \mf{D\Psi}(\mf{x})\right|d\mf{x}.\]
\end{theorem}
Among the assertions in the theorem, we will first establish the following.
\begin{theorem}[label=230905_4] {}
Let $\mathcal{O}$ be an open subset of $\mb{R}^n$, and let $\mf{\Psi}:\mathcal{O}\to\mb{R}^n$ be a smooth change of variables. If $\mk{D}$ is a Jordan measurable set such that its closure $\overline{\mk{D}}$ is contained in $\mathcal{O}$, then $\mf{\Psi}(\mk{D})$ is also Jordan measurable.  
\end{theorem}
A special case of the change of variables theorem is when  $f:\mf{\Psi}(\mk{D})\to\mb{R}$ is the characteristic function of $\mf{\Psi}(\mk{D})$. This gives the change of volume theorem.
\begin{theorem}[label=230906_2]{The Change of Volume Theorem}
Let $\mathcal{O}$ be an open subset of $\mb{R}^n$, and let $\mf{\Psi}:\mathcal{O}\to\mb{R}^n$ be a smooth change of variables. If $\mk{D}$ is a Jordan measurable set such that its closure $\overline{\mk{D}}$ is contained in $\mathcal{O}$, then  
\[\text{vol}\left(\mf{\Psi}(\mk{D})\right) =\int_{\mk{D}} \left|\det \mf{D\Psi}(\mf{x})\right|d\mf{x}.\]
\end{theorem}
In the following, let us   give some remarks about the statements in the theorem, and outline   the plan of the  proof.
\begin{highlight}{The Change of Variables Theorem}
 
\begin{enumerate}[1.]
\item 
The first step is to prove Theorem \ref{230905_4} which asserts that $\mf{\Psi}(\mk{D})$ is Jordan measurable. To do this, we first show that a smooth change of variables sets up a one-to-one correspondence between the open sets in the domain and the range. This basically follows from inverse function theorem.
\item  Since $f:\mf{\Psi}(\mk{D})\to\mb{R}$ is continuous and bounded, if $\mf{\Psi}(\mk{D})$ is Jordan measurable, $f:\mf{\Psi}(\mk{D})\to\mb{R}$ is Riemann integrable.
\item Let $g:\mathcal{O}\to\mb{R}$  be the function
\[g(\mf{x})=|\det  \mf{D\Psi}(\mf{x})|.\]
Since $\mf{\Psi}:\mathcal{O}\to\mb{R}^n$ is continuously differentiable, the function $\mf{D\Psi}:\mathcal{O}\to\mb{R}^{n^2}$ is continuous.
Since  determinant  and absolute value are   continuous functions,  $g:\mathcal{O}\to\mb{R}$ is a  continuous function. Since $\overline{\mk{D}}$ is a compact set contained in $\mathcal{O}$, $g:\overline{\mk{D}}\to\mb{R}$ is bounded. 
\item Since $\mf{\Psi}:\mathcal{O}\to\mb{R}^n$ is continuous, and the functions $f:\mf{\Psi}(\mk{D})\to\mb{R}$ and $g:\mk{D}\to\mb{R}$ are continuous and bounded, the function $h:\mk{D}\to\mb{R}$,
\[h(\mf{x})=f\left(\mf{\Psi}(\mf{x})\right)g(\mf{x})=f\left(\mf{\Psi}(\mf{x})\right)\left|\det \mf{D\Psi}(\mf{x})\right|\] is continuous and bounded. Hence, it is Riemann integrable.
\item To prove the change of variables theorem, we will first prove the change of volume theorem. This is the most technical part of  the proof. 
\item  To prove the change of volume theorem, we first consider the case where $\mf{\Psi}:\mb{R}^n\to\mb{R}^n$ is an invertible linear transformation. In this case, the theorem says that if $\mk{D}$ is a Jordan measurable set, and $\mf{T}:\mb{R}^n\to\mb{R}^n$, $\mf{T}(\mf{x})=A\mf{x}$ is an invertible linear transformation, then 
\begin{equation}\label{230907_2}\text{vol}\,\left(\mf{T}(\mk{D})\right)=|\det A|\,\text{vol}\,(\mk{D}).\end{equation}
\end{enumerate}
\end{highlight}
\begin{highlight}{}
\begin{enumerate}[1.]\item[7.] To prove \eqref{230907_2}, we first consider the case where $\mk{D}=\mf{I}$ is a closed rectangle. This is an easy consequence of the fact that the volume of a parallelepiped spanned by the vectors $\mf{v}_1, \ldots,\mf{v}_n$ is equal to $|\det A|$, where $A$ is the matrix with $\mf{v}_1, \ldots,\mf{v}_n$ as column vectors. This was proved in Appendix \ref{appB}.

\item[8.] After proving the change of volume theorem, we will prove the change of variables theorem for the special case where $\mk{D}=\mf{I}$ is a closed rectangle first. The general theorem then follows by some simple analysis argument.
\end{enumerate}
\end{highlight}

 We begin by the following
  proposition which says that a smooth change of variables maps open sets to open sets.
\begin{proposition}[label=230905_3]{}
Let $\mathcal{O}$ be an open subset of $\mb{R}^n$, and let $\mf{\Psi}:\mathcal{O}\to\mb{R}^n$ be a smooth change of variables.
Then for any open set $\mathcal{D}$ that  is contained in $\mathcal{O}$, $\mf{\Psi}(\mathcal{D})$ is open in $\mb{R}^n$.
In particular,  $\mf{\Psi}(\mathcal{O})$ is an open subset of $\mb{R}^n$.
 
\end{proposition}
\begin{myproof}{Proof}
Given that $\mathcal{D}$ is an open subset of $\mb{R}^n$, let $\mathcal{W}=\mf{\Psi}(\mathcal{D})$. We want to show that $\mathcal{W}$ is an open set. If  $\mf{y}_0$ is a point in $\mathcal{W}$, there is an $\mf{x}_0$ in $\mathcal{D}$ such that $\mf{y}_0=\mf{\Psi}(\mf{x}_0)$. 
Since $\mf{\Psi}:\mathcal{O}\to\mb{R}^n$ is continuously differentiable and $\mf{D\Psi}(\mf{x}_0)$ is invertible, we can apply inverse function theorem to conclude that there is an open set $\mathcal{U}_0$ containing $\mf{x}_0$ such that $\mathbf{\Psi}(\mathcal{U}_0)$ is also open, and $\mf{\Psi}^{-1}:\mf{\Psi}(\mathcal{U}_0)\to\mathcal{U}_0$ is continuously differentiable. Let $\mathcal{U}=\mathcal{U}_0\cap\mathcal{D}$. Then $\mathcal{U}$ is an open subset of $\mathcal{D}$ and $\mathcal{U}_0$. It follows that $\mathcal{V}=\mf{\Psi}(\mathcal{U})=(\mf{\Psi}^{-1})^{-1}(\mathcal{U})$ is an open subset of $\mb{R}^n$ that is contained in $\mathcal{W}=\mf{\Psi}(\mathcal{D})$. Notice that $\mathcal{V}$ is an open set that contains $\mf{y}_0$. Thus, we have shown that every point in $\mathcal{W}$ has a neighbourhood that lies in $\mathcal{W}$. This proves that $\mathcal{W}$ is an open set.
\end{myproof}

The following proposition says that the inverse of a smooth change of variables is  also a smooth change of variables.
\begin{proposition}[label=230905_7]
{}
Let $\mathcal{O}$ be an open subset of $\mb{R}^n$, and let $\mf{\Psi}:\mathcal{O}\to\mb{R}^n$ be a smooth change of variables. Then $\mf{\Psi}^{-1}:\mf{\Psi}(\mathcal{O})\to\mb{R}^n$ is also a smooth change of variables.
\end{proposition}
\begin{myproof}{Proof}
By Proposition \ref{230905_3}, $\mf{\Psi}(\mathcal{O})$ is an open set.
By default, $\mf{\Psi}^{-1}:\mf{\Psi}(\mathcal{O})\to\mb{R}^n$ is one-to-one.
As in the proof of Proposition \ref{230905_3}, the inverse function theorem implies that it is  continuously differentiable. If $\mf{x}_0=\mf{\Psi}^{-1}(\mf{y}_0)$, inverse function theorem says that
\[\mf{D\Psi}^{-1}(\mf{y}_0)=\mf{D\Psi}(\mf{x}_0)^{-1}.\]
The inverse of an invertible matrix is invertible. Hence, for any $\mf{y}_0$ in $\mf{\Psi}(\mathcal{O})$, $\mf{D\Psi}^{-1}(\mf{y}_0)$ is invertible. These prove that $\mf{\Psi}^{-1}:\mf{\Psi}(\mathcal{O})\to\mb{R}^n$ is a smooth change of variables.
\end{myproof}

\begin{remark}[label=230905_8]{Homeomorphisms and Diffeomorphisms}
Let $\mathcal{O}$ be an open subset of $\mb{R}^n$, and let $\mf{\Psi}:\mathcal{O}\to \mb{R}^n$ be a continuous injective map such that $\mf{\Psi}(\mathcal{O})$ is open, and the inverse map $\mf{\Psi}^{-1}:\mf{\Psi}(\mathcal{O})\to\mathcal{O}$ is   continuous. Then we say that  $\mf{\Psi}:\mathcal{O}\to \mf{\Psi}(\mathcal{O})$  is a {\it homeomorphism}. A homeomorphism sets up a one-to-one correspondence between open sets in $\mathcal{O}$ and open sets in $\mf{\Psi}(\mathcal{O})$. 

If $\mf{\Psi}:\mathcal{O}\to \mf{\Psi}(\mathcal{O})$ is a homeomorphism and both the maps $\mf{\Psi}:\mathcal{O}\to \mf{\Psi}(\mathcal{O})$ and $\mf{\Psi}^{-1} :\mf{\Psi}(\mathcal{O})\to\mathcal{O}$ are continuously differentiable, then we say that $\mf{\Psi}:\mathcal{O}\to \mf{\Psi}(\mathcal{O})$ is a {\it diffeomorphism}.  Proposition  \ref{230905_3} and Proposition \ref{230905_7} imply that a continuous change of variables is a diffeomorphism.

A map of the form $\mf{\Psi}:\mb{R}^n\to\mb{R}^n$, \[\mf{\Psi}(\mf{x})=\mf{y}_0+A(\mf{x}-\mf{x}_0),\] where $\mf{x}_0$ and $\mf{y}_0$ are points in $\mb{R}^n$ and $A$ is an invertible matrix, is a diffeomorphism.
\end{remark}

Now we can prove the following which is essential for the proof of Theorem \ref{230905_4}.
\begin{theorem}{}
Assume that $\mathcal{O}$ and $\mathcal{U}$ are open subsets of $\mb{R}^n$, and  $\mf{\Psi}:\mathcal{O}\to\mathcal{U}$ is a homeomorphism. If $\mk{D}$ is  a subset of $\mathcal{O}$ such that $\overline{\mk{D}}$ is also contained in $\mathcal{O}$, then
\[\text{int}\,\mf{\Psi}(\mk{D})=\mf{\Psi}(\text{int}\,\mk{D}),\hspace{1cm}\overline{\mf{\Psi}(\mk{D})}=\mf{\Psi}(\overline{\mk{D}}).\]
Thus,
\[\pa\mf{\Psi}(\mk{D})=\mf{\Psi}(\pa\mk{D}).\]
\end{theorem}
\begin{myproof}{Proof}
The interior of a set $A$ is an open set that contains all the open set that is contained in $A$.  By  Remark \ref{230905_8},  there is a one-to-one correspondence between the open sets that are contained in $\mk{D}$ and the open sets that are contained in $\mf{\Psi}(\mk{D})$. Therefore, 
\[\text{int}\,\mf{\Psi}(\mk{D})=\mf{\Psi}(\text{int}\,\mk{D}).\]
Since $\overline{\mk{D}}$ is a compact set and $\mf{\Psi}:\mathcal{O}\to\mathcal{U}$ is continuous, $\mf{\Psi}(\overline{\mk{D}})$ is a compact set. Therefore,  $\mf{\Psi}(\overline{\mk{D}})$  is a closed set that contains $\mf{\Psi}(\mk{D})$. This implies that
\begin{equation}\label{230905_9}\overline{\mf{\Psi}(\mk{D})}\subset \mf{\Psi}(\overline{\mk{D}}).\end{equation}
Since  $\mf{\Psi}^{-1}: \mathcal{U}\to\mathcal{O}$ is also continuous, the same argument gives
\[\overline{ \mk{D} }=\overline{\mf{\Psi}^{-1}\left(\mf{\Psi}(\mk{D})\right)}\subset \mf{\Psi}^{-1}(\overline{\mf{\Psi}(\mk{D})}).\]
This implies that
\begin{equation}\label{230905_10}\mf{\Psi}\left(\overline{ \mk{D} }\right)\subset  \overline{\mf{\Psi}(\mk{D})}.\end{equation}
Eq. \eqref{230905_9} and \eqref{230905_10} give

\vspace{-0.4cm}
\[\overline{\mf{\Psi}(\mk{D})}=\mf{\Psi}(\overline{\mk{D}}).\]
The last assertion follows from the fact that for any set $A$, $\overline{A}$ is a disjoint union of $\text{int}\,A$ and $\pa A$.
\end{myproof}

Recall that a set $\mk{D}$ in $\mb{R}^n$ has Jordan content zero if and only if for every $\varepsilon>0$, $\mk{D}$ can be covered by finitely many cubes $Q_1$, $\ldots$, $Q_k$, such that
\[\sum_{j=1}^k\text{vol}\,(Q_j)<\varepsilon.\]

The next proposition gives a control of the size of the cube under a smooth change of variables.
\begin{proposition}[label=230906_3]{}
Let $\mathcal{O}$ be an open subset of $\mb{R}^n$, and let $\mf{\Psi}:\mathcal{O}\to\mb{R}^n$ be a smooth change of variables. If $ Q_{\mf{c},r}$ is a cube with center at $\mf{c}$ and side length $2r$, then $\mf{\Psi}(Q_{\mf{c},r})$ is contained in the cube $Q_{\mf{\Psi}(\mf{c}), \lambda r}$, where
\[ \lambda=\max_{1\leq i\leq n}\max_{\mf{x}\in Q_{\mf{c},r}}\sum_{j=1}^n\left|\frac{\pa\Psi_i}{\pa x_j}(\mf{x})\right|.\]
Therefore,
\[\text{vol}\,\left(\mf{\Psi}(Q_{\mf{c},r})\right)\leq\lambda^n \text{vol}\,(Q_{\mf{c},r}).\]
\end{proposition}
\begin{remark}{}
Note that since $Q_{\mf{c},r}$ is a compact set and $\di \frac{\pa\Psi_i}{\pa x_j}(\mf{x})$ is continuous for all $1\leq i,j\leq n$, 
\[\max_{\mf{x}\in Q_{\mf{c},r}}\sum_{j=1}^n\left|\frac{\pa\Psi_i}{\pa x_j}(\mf{x})\right|\]exists.
\end{remark}
\begin{myproof}{\linkt Proof of Proposition \ref{230906_3}}
Notice that  $\mf{u}\in Q_{\mf{c}, r}$ if and only if \[|u_i-c_i|\leq  r \hspace{1cm}\text{for each}\;1\leq i\leq n.\] Let $\mf{d}=\mf{\Psi}(\mf{c})$.  
Given $\mf{v}=\mf{\Psi}(\mf{u})$ with $\mf{u}\in Q_{\mf{c}, r}$, we want to show that $\mf{v}$ is in $Q_{\mf{d},\lambda r}$, or equivalently, \[|v_i-d_i|\leq \lambda r \hspace{1cm}\text{for each}\;1\leq i\leq n.\]  

\bp 
This is basically an application of mean value theorem. 
The set $Q_{\mf{c},r}$ is convex and the map $\Psi_i:\mathcal{O}\to\mb{R}$ is continuously differentiable. Mean value theorem says that there is a point $\mf{x}$ in $Q_{\mf{c},r}$ such that
\[
v_i-d_i =\Psi_i(\mf{u})-\Psi_i(\mf{c})=\sum_{j=1}^n \frac{\pa \Psi_i}{\pa x_j}(\mf{x})(u_j-c_j).\]

Therefore,
\begin{align*}
|v_i-d_i|&\leq \sum_{j=1}^n \left|\frac{\pa \Psi_i}{\pa x_j}(\mf{x})\right|\left|u_j-c_j\right|\leq r \sum_{j=1}^n \left|\frac{\pa \Psi_i}{\pa x_j}(\mf{x})\right|\\
&\leq r\max_{\mf{x}\in Q_{\mf{c},r}}\sum_{j=1}^n \left|\frac{\pa \Psi_i}{\pa x_j}(\mf{u})\right|\leq \lambda r.
\end{align*}This proves that $\mf{\Psi}(Q_{\mf{c},r})$ is contained in $Q_{\mf{\Psi}(\mf{c}),\lambda r}$.
The last assertion in the proposition about the volumes is obvious.
\end{myproof}

Now we prove Theorem \ref{230905_4}.
\begin{myproof}{\linkt Proof of Theorem \ref{230905_4}}
Since $\overline{\mk{D}}$ is a compact set that is contain in the open set $\mathcal{O}$, Theorem \ref{230905_12} says that there is a positive number $d$ and a compact set $C$ such that $\mk{D}\subset C\subset \mathcal{O}$, and any point in $\mb{R}^n$ that has a distance less than $d$ from a point in $\overline{\mk{D}}$ lies in $C$.

Since $\mf{\Psi}:\mathcal{O}\to\mb{R}^n$ is continuously differentiable, for all $1\leq i,j\leq n$, $\di \frac{\pa \Psi_i}{\pa x_j}:C\to\mb{R}$ is a continuous function. Since $C$ is a compact set, for each $1\leq i\leq n$, the function
\[\sum_{j=1}^n\left|\frac{\pa\Psi_i}{\pa x_j}\right|(\mf{x})\] has a maximum   on $C$. Hence,

\[\lambda=\max_{1\leq i\leq n}\max_{\mf{x}\in C}\sum_{j=1}^n\left|\frac{\pa\Psi_i}{\pa x_j}\right|(\mf{x})\] 
 exists.
 \bp
Since $\mk{D}$ is Jordan measurable, $\pa\mk{D}$ has Jordan content zero. Since $C$   contains $\overline{\mk{D}}$, it contains $\pa\mk{D}$. 
  Given $\varepsilon>0$, there exist cubes $Q_1, Q_2, \ldots, Q_k$, each of which intersects $\mk{D}$, and such that 
\[\pa\mk{D}\subset \bigcup_{j=1}^n Q_j \quad\text{and}\quad\sum_{j=1}^k \text{vol}\,(Q_j)<\frac{\varepsilon}{\lambda^n}.\]

Since a uniformly regular partition of a cube will divide the cube into cubes, we can also assume that each of the cubes $Q_j$, $1\leq j\leq k$ has diameter less than $d$. This implies that each $Q_j$, $1\leq j\leq k$ is contained in $C$. For $1\leq j\leq k$, let $l_j$ be the side length of $Q_j$. Proposition \ref{230906_3} says that $\mf{\Psi}(Q_j)$ is contained in a cube $\widetilde{Q}_j$ with side length $\lambda l_j$. 
Therefore,

\[\pa \mf{\Psi}(\mk{D})=\mf{\Psi}(\pa\mk{D})\subset\bigcup_{j=1}^k \widetilde{Q}_j,\]

\[\text{and}\hspace{1cm}
\sum_{j=1}^k\text{vol}\,(\widetilde{Q}_j) \leq \lambda^n \sum_{j=1}^k\text{vol}\,(Q_j)<\varepsilon.
\]This shows that $\pa \mf{\Psi}(\mk{D})$ has Jordan content zero. Hence, $\mf{\Psi}(\mk{D})$ is a Jordan measurable set.
\end{myproof}

To prove the change of volume formula, the crucial thing is to first prove the special case where $\mk{D}=\mf{I}$ is a rectangle, and $\mf{\Psi}$ is an invertible linear transformation. In Appendix \ref{appB}, we prove the following theorem which gives the volume of a parallelepiped.

\begin{theorem}[label=230906_4]{}
Let $\mathscr{P}$ be a parallelepiped in $\mb{R}^n$ spanned by the linearly independent vectors $\mf{v}_1, \ldots, \mf{v}_n$. Then the volume of $\mathscr{P}$ is equal to $|\det A|$, where $A$ is the matrix whose column vectors are $\mf{v}_1, \ldots, \mf{v}_n$.
\end{theorem}
We then use this to deduce the following special case of the change of volume formula.

 \begin{theorem}[label=230906_5]{}
Let $\mf{I}$ be a closed rectangle in $\mb{R}^n$, and let $\mf{T}:\mb{R}^n\to\mb{R}^n$, $\mf{T}(\mf{x})=A\mf{x}$ be an invertible linear transformation. Then
\begin{equation*}\text{vol}\,(\mf{T}(\mf{I}))=|\det A|\,\text{vol}\,(\mf{I}).\end{equation*}
\end{theorem}

Using this, we can prove the more general change of volume formula for a Jordan measurable set under an invertible linear transformation.
\begin{theorem}[label=230906_11]{}
If $\mf{T}:\mb{R}^n\to\mb{R}^n$, $\mf{T}(\mf{x})=A\mf{x}$ is an invertible linear transformation, and $\mk{D}$ is a Jordan measurable set, then $\mf{T}(\mk{D})$ is also Jordan measurable and
\[\text{vol}\,(\mf{T}(\mk{D}))=|\det A|\,\text{vol}\,(\mk{D}).\]
\end{theorem}  
\begin{myproof}{Proof}
 Since $\mf{T}:\mb{R}^n\to\mb{R}^n$ is invertible, $\det A\neq 0$.
The fact that the set $\mf{T}(\mk{D})$ is Jordan measurable follows from Theorem \ref{230905_4}. Let $\mf{I}$ be a closed rectangle that contains $\mk{D}$. Since $\mk{D}$ is Jordan measurable, the characteristic function $\chi_{\mk{D}}:\mf{I}\to\mb{R}$ is Riemann integrable.

Given $\varepsilon>0$, there is a partition $\mf{P}$ of $\mf{I}$ such that
\[U(\chi_{\mk{D}},\mf{P})-\int_{\mf{I}}\chi_{\mk{D}} <\frac{\varepsilon}{|\det A|}.\]
Hence, 
\[U(\chi_{\mk{D}},\mf{P})<\text{vol}\,(\mk{D})+\frac{\varepsilon}{|\det A|}.\]
Let
\[\mathscr{A}=\left\{\mf{J}\in \mathcal{J}_{\mf{P}}\,|\, \mf{J}\cap\mk{D}\neq \emptyset\right\}.\]
Then
\[\sum_{\mf{J}\in\mathscr{A}}\,\text{vol}\,(\mf{J})=U(\chi_{\mk{D}},\mf{P})<\text{vol}\,(\mk{D})+\frac{\varepsilon}{|\det A|}.\]

Notice that
\[\mk{D}\subset \bigcup_{\mf{J}\in\mathscr{A}}\mf{J}.\]
\bp

Therefore,
\[\mf{T}(\mk{D})\subset \bigcup_{\mf{J}\in\mathscr{A}}\mf{T}(\mf{J}).\]
For each rectangle $\mf{J}$, $\mf{T}(\mf{J})$ is a parallelepiped. For any two distinct rectangles in $\mathscr{A}$, they are disjoint or intersect at a set that has Jordan content zero. Therefore, 
 additivity theorem implies that the set $K$ defined as
\[K=\bigcup_{\mf{J}\in\mathscr{A}}\mf{T}(\mf{J})\] 

is Jordan measurable, and 
\[
\text{vol}\,(K)=\sum_{\mf{J}\in\mathscr{A}}\text{vol}\,(\mf{T}(\mf{J}))=|\det A|\,\sum_{\mf{J}\in\mathscr{A}}\text{vol}\,(\ \mf{J})
<|\det A|\,\text{vol}\,(\mk{D})+\varepsilon.\]

Since $\mf{T}(\mk{D})\subset K$, we find that
\[\text{vol}\,(\mf{T}(\mk{D}))\leq \text{vol}\,(K)<|\det A|\,\text{vol}\,(\mk{D})+\varepsilon.\]
Since $\varepsilon>0$ is arbitrary, we conclude that
\begin{equation}\label{230906_8}\text{vol}\,(\mf{T}(\mk{D}))\leq |\det A|\,\text{vol}\,(\mk{D}).\end{equation}

Since $\mf{T}^{-1}:\mb{R}^n\to\mb{R}^n$ is also an invertible linear transformation, we find that
\begin{equation}\label{230906_9}
\begin{split}
\text{vol}\,(\mk{D})&=\text{vol}\,\left(\mf{T}^{-1}\left(\mf{T}(\mk{D})\right)\right)\\&\leq |\det A^{-1}|\,\text{vol}\,(\mf{T}(\mk{D}))=\frac{1}{|\det A|}\text{vol}\,(\mf{T}(\mk{D})).\end{split}
\end{equation}
Eq. \eqref{230906_8} and \eqref{230906_9} together give
\[\text{vol}\,(\mf{T}(\mk{D}))=|\det A|\,\text{vol}\,(\mk{D}).\]
\end{myproof}

Recall that by identitying an $n\times n$ matrix $A=[a_{ij}]$ as a point in $\mb{R}^{n^2}$, we have defined the norm of $A$ as
\[\Vert A\Vert=\sqrt{\sum_{i=1}^n\sum_{j=1}^na_{ij}^2}.\]
Besides the triangle inequality, this norm also satisfies the following identity.
\begin{lemma}[label=230906_14]
{}If $A=[a_{ij}]$ and $B=[b_{ij}]$ are $n\times n$ matrices, then 
\[\Vert AB\Vert\leq \Vert A\Vert\,\Vert B\Vert.\]
\end{lemma}
\begin{myproof}{Proof}
Let $[c_{ij}]=C=AB$. Then for any $1\leq i,j\leq n$,
\[c_{ij}=\sum_{k=1}^na_{ik}b_{kj}.\]
By Cauchy-Schwarz inequality, 
\[c_{ij}^2\leq \left( \sum_{k=1}^n a_{ik}^2\right)\left(\sum_{l=1}^nb_{lj}^2\right).\]
Therefore,
\[\Vert C\Vert^2=\sum_{i=1}^n\sum_{j=1}^n c_{ij}^2\leq \sum_{i=1}^n\sum_{j=1}^n \left( \sum_{k=1}^n a_{ik}^2\right)\left(\sum_{l=1}^nb_{lj}^2\right)=\Vert A\Vert^2\Vert B\Vert^2.\]
This proves that
\[\Vert AB\Vert=\Vert C\Vert\leq \Vert A\Vert\,\Vert B\Vert.\]
\end{myproof}

Now we prove the change of volume formula, which is the most technical part.
\begin{myproof}{\linkt Proof of Theorem \ref{230906_2}}
Given the smooth change of variables $\mf{\Psi}:\mathcal{O}\to\mb{R}^n$,
let $g:\mathcal{O}\to\mb{R}$  be the continuous function
\[g(\mf{x})=|\det  \mf{D\Psi}(\mf{x})|.\] 
We want to show that if $\mk{D}$ is a Jordan measurable set such that its closure $\overline{\mk{D}}$ is contained in $\mathcal{O}$, then
\[\text{vol}\,\left(\mf{\Psi}(\mk{D})\right)=\int_{\mk{D}}|\det  \mf{D\Psi}(\mf{x})|d\mf{x}=\int_{\mk{D}}g(\mf{x})d\mf{x}=\mathcal{I}.\]
\bp

We will first prove that 
\[\text{vol}\,(\mf{\Psi}(\mk{D}))\leq \mathcal{I}.\]
By Theorem \ref{230905_4}, $\mf{\Psi}(\mk{D})$ is Jordan measurable. By Theorem \ref{230829_4},   its closure $\overline{\mf{\Psi}(\mk{D})}=\mf{\Psi}(\overline{\mk{D}})$ is also Jordan measurable, and 
\[\text{vol} \left(\mf{\Psi}(\mk{D})\right)=\text{vol} \left( \overline{\mf{\Psi}(\mk{D} )}\right)=\text{vol} \left( \mf{\Psi}(\overline{\mk{D}})\right).\]

On the other hand, since $\overline{\mk{D}}\setminus \mk{D}$ has Jordan content zero, 
\[\int_{\mk{D}}g(\mf{x})d\mf{x}=\int_{\overline{\mk{D}}}g(\mf{x})d\mf{x}.\]
 Hence, we can assume from the beginning that $\mk{D}=\overline{\mk{D}}$, or equivalently, $\mk{D}$ is closed.

As in the proof of Theorem \ref{230905_4}, Theorem \ref{230905_12} says that there is a positive number $d$ and a compact set $C$ such that $\mk{D}\subset C\subset \mathcal{O}$, and any point in $\mb{R}^n$ that has a distance less than $d$ from a point in $ \mk{D}$ lies in $C$. On the compact set $C$, the function  $g:C\to\mb{R}$  is continuous. By extreme value theorem, there are points  $\mf{u}$ and $\mf{v}$  in $C$ such that
$g(\mf{u}) \leq g(\mf{x})\leq g(\mf{v}) $ for all $\mf{x}\in C$.
Let $m_g=g(\mf{u})$ and $M_g=g(\mf{v})$. Then $m_g>0$ and 
\[m_g\leq g(\mf{x})\leq M_g  \quad \hspace{1cm}\text{for all}\;\mf{x}\in C.\]

On the other hand, the function $\mf{D\Psi}^{-1}:C\to \mb{R}^{n^2}$ is continuous on the compact set $C$. Hence, it is bounded. Namely, there is a positive number $M_h$ such that
\[  \Vert \mf{D\Psi}^{-1}(\mf{x})\Vert \leq M_h\hspace{1cm}\text{for all}\;\mf{x}\in C.\]

Let $L$ be a positive number  such that $C$ is contained in the cube $\mf{I}=[-L, L]^n$. Let $\check{g}:\mf{I}\to\mb{R}$ be the zero extension of $g:\mk{D}\to\mb{R}$. For each positive integer $k$, let $\mf{P}_k$ be the uniformly regular partition of $\mf{I}$ into $k^n$ rectangles. Then $\di\lim_{k\to\infty}|\mf{P}_k|=0$. Therefore,  
\[\lim_{k\to\infty}U(\check{g}, \mf{P}_k) =\int_{\mf{I}}\check{g} (\mf{x})d\mf{x}=\int_{\mk{D}}g(\mf{x})d\mf{x}=\mathcal{I},\]
\bp 
and
\[\lim_{k\to\infty}\left(U(\chi_{\mk{D}}, \mf{P}_k)-L(\chi_{\mk{D}}, \mf{P}_k)\right)=0.\]
 
The compactness of $C$  implies that  the continuous functions   $\di \mf{D\Psi}:C\to\mb{R}^{n^2}$ and $ g:C\to\mb{R}$  are uniformly continuous. 
  Given $\varepsilon>0$,  there exists a $\delta_1>0$ such that if $\mf{u}$ and $\mf{v}$ are points in $C$ with $\Vert\mf{u}-\mf{v}\Vert <\delta_1$, then 
\[\left\Vert \mf{D\Psi}(\mf{u})-\mf{D\Psi}(\mf{v})\right\Vert<\frac{\varepsilon}{M_hn}\quad\text{and}\quad |g(\mf{u})-g(\mf{v})|< m_g\varepsilon.\]
 
Let $\delta=\min\{d, \delta_1\}$.
There is a positive integer $K$ such that for all $k\geq K$, $|\mf{P}_k|<\delta$,
\[   U(\check{g}, \mf{P}_k)<\mathcal{I}+ \varepsilon \quad\text{and}\quad U(\chi_{\mk{D}},\mf{P}_k)-L(\chi_{\mk{D}},\mf{P}_k)<\frac{\varepsilon}{M_g}.\]
Consider the partition $\mf{P}=\mf{P}_K$. Since $\mf{I}$ is a cube and $\mf{P}$ is a uniformly regular partition of $\mf{I}$, each rectangle in the partition $\mf{P}$ is also a cube, all with the same side length $2r$.  
Denote by $\mathscr{A}$ and $\mathscr{B}$ the sets
\[\mathscr{A}=\left\{\mf{J}\in\mathcal{J}_{\mf{P}}\,|\,\mf{J}\cap\mk{D}\neq\emptyset\right\},\hspace{1cm} \mathscr{B}=\left\{\mf{J}\in\mathcal{J}_{\mf{P}}\,|\,\mf{J}\subset\mk{D}\right\}.\]
 
$\mathscr{A}$  is a finite collection of cubes, and it contains the collection $\mathscr{B}$. 
By definition,
\[L(\chi_{\mk{D}},\mf{P})=\sum_{\mf{J}\in\mathscr{B}}\text{vol}\,(\mf{J}),\hspace{1cm}U(\chi_{\mk{D}},\mf{P})=\sum_{\mf{J}\in\mathscr{A}}\text{vol}\,(\mf{J}).\]
Therefore,
\[ \sum_{\mf{J}\in\mathscr{A}\setminus \mathscr{B}}\text{vol}\,(\mf{J})<\frac{\varepsilon}{M_g}.\]

After renaming, we can assume that
\[\mathscr{A}=\left\{Q_{\beta}\,|\, 1\leq \beta\leq s\right\},\]where $Q_{\beta}=Q_{\mf{c}_{\beta}, r}$ is a cube with center at $\mf{c}_{\beta}$ and side length $2r$.
By the definition of $\mathscr{A}$, 
\[\mk{D}\subset \bigcup_{{\beta}=1}^s Q_{\beta}.\]
\bp
Therefore,
\[\mf{\Psi}(\mk{D})\subset  \bigcup_{{\beta}=1}^s \mf{\Psi}\left(Q_{ \beta}\right).\]

Fixed a $1\leq {\beta}\leq s$. Since $Q_{\beta}$ intersects $\mk{D}$ and \[\text{diam}\,  Q_{\beta}= |\mf{P}|<\delta\leq d,\] we find that $Q_{  \beta}$ is contained in $ C$. 
 Define the invertible linear transformation $\mf{T}_{\beta}:\mb{R}^n\to\mb{R}^n$ by $\mf{T}_{\beta}(\mf{x})=A_{\beta}\mf{x}$, where $A_{\beta}=\mf{D\Psi}(\mf{c}_{\beta})$. Then let 
$\mf{\Phi}_{\beta}: Q_{\beta}\to\mb{R}^n$ be the map $\mf{\Phi}_{\beta}=\mf{T}_{\beta}^{-1}\circ\mf{\Psi}:Q_{ \beta}\to\mb{R}^n$. 
By chain rule, 
\[\mf{D\Phi}_{\beta}(\mf{x})=\mf{DT}_{\beta}^{-1}(\mf{\Psi}(\mf{x})) \,\mf{D\Psi}(\mf{x})=A_{\beta}^{-1}\,\mf{D\Psi}(\mf{x})\quad\text{for all}\;\mf{x}\in Q_{ \beta}.\]
 
Therefore, for $\mf{x}\in Q_{\beta}$, 
\[ \mf{D\Phi}_{\beta}(\mf{x})-I_n =\mf{D\Psi}^{-1}(\mf{\Psi}(\mf{c}_{\beta})) \left(\mf{D\Psi}(\mf{x})-\mf{D\Psi}(\mf{c}_{\beta})\right).\]
Since $\Vert \mf{x}-\mf{c}_{\beta}\Vert\leq\text{diam}\,Q_{\beta}<\delta\leq \delta_1$, Lemma \ref{230906_14} implies that
\begin{align*}
\left\Vert \mf{D\Phi}_{\beta}(\mf{x})-I_n\right\Vert\leq \left\Vert \mf{D\Psi}^{-1}(\mf{\Psi}(\mf{c}_{\beta})) \right\Vert\,\left\Vert \mf{D\Psi}(\mf{x})-\mf{D\Psi}(\mf{c}_{\beta})\right\Vert<\frac{\varepsilon}{n}.
\end{align*}
 
This implies that if $i\neq j$,
\[\left|\frac{\pa (\Phi_{\beta})_i}{\pa x_j}(\mf{x})\right|<\frac{\varepsilon}{n}\hspace{1cm} \text{for all}\;\mf{x}\in Q_{\beta};\]
while
\[\left|\frac{\pa (\Phi_{\beta})_i}{\pa x_i}(\mf{x})\right|<1+\frac{\varepsilon}{n}\hspace{1cm} \text{for all}\;\mf{x}\in Q_{\beta}.\]
Hence, 
\[\lambda_{\beta}=\max_{1\leq i\leq n}\max_{\mf{x}\in Q_{\beta}}\sum_{j=1}^n \left|\frac{\pa (\Phi_{\beta})_i}{\pa x_j}(\mf{x})\right|\leq 1+\varepsilon.\]
Since $\mf{\Psi}=\mf{T}_{\beta}\circ \mf{\Phi}_{\beta}$, by Theorem \ref{230906_11} and Proposition \ref{230906_3}, we have
\[\text{vol}\,(\mf{\Psi}(Q_{\beta})=|\det A_{\beta}|\,\text{vol}\,(\mf{\Phi}_{\beta}(Q_{\beta}))\leq |\det A_{\beta}|\,\lambda_{\beta}^n \,\text{vol}\left(Q_{\beta}\right).\]
 
\bp
Summing over $\beta$, we find that
\[
\text{vol}\left(\mf{\Psi}(\mk{D})\right)\leq \sum_{\beta=1}^s \text{vol}\left(\mf{\Psi}\left(Q_{\beta}\right)\right) \leq  (1+\varepsilon)^n \sum_{\beta=1}^s|\det \mf{D\Psi}(\mf{c}_{\beta})|\,\text{vol}\left(Q_{\beta}\right).
\]

We divide the sum into a sum over those $Q_{\beta}$ in $\mathscr{B}$ and a sum over those $Q_{\beta}$ in $\mathscr{A}\setminus\mathscr{B}$.
For the sum over those in $\mathscr{B}$, we find that
\begin{align*}\sum_{Q_{\beta}\in \mathscr{B}}|\det \mf{D\Psi}(\mf{c}_{\beta})|\,\text{vol}\left(Q_{\beta}\right)
\leq U(\check{g}, \mf{P})<\mathcal{I}+\varepsilon.
\end{align*}
  For the sum over those $Q_{\beta}$ in $\mathscr{A}\setminus\mathscr{B}$,
\begin{align*}\sum_{Q_{\beta}\in  \mathscr{A}\setminus\mathscr{B}} |\det \mf{D\Psi}(\mf{c}_{\beta})|\,\text{vol}\left(Q_{\beta}\right)
&\leq   M_g\sum_{Q_{\beta}\in\mathscr{A}\setminus\mathscr{B}}\text{vol}\left(Q_{\beta}\right)< \varepsilon.
\end{align*}
Hence,
\[
\text{vol}\left(\mf{\Psi}(\mk{D})\right)\leq (1+\varepsilon)^n\left(\mathcal{I}+2\varepsilon\right).\]
 
Since $\varepsilon>0$ is arbitrary, taking the limit $\varepsilon\to 0^+$, we find that
\[\text{vol}\,\left(\mf{\Psi}(\mk{D})\right)\leq\mathcal{I}=\int_{\mk{D}} \left|\det \mf{D\Psi}(\mf{x})\right|d\mf{x}.\]
This is true for any smooth change of variables $\mf{\Psi}:\mathcal{O}\to\mb{R}$ and any Jordan measurable closed subet $\mk{D}$ that is contained in $\mathcal{O}$.

Now we want to prove  the opposite inequality. First note that the same inequality applied to the smooth change of  variables $\mf{\Psi}^{-1}:\mf{\Psi}(\mathcal{O})\to\mb{R}^n$. Thus, if $\mathcal{F}$ is a closed Jordan measurable subset of $\mf{\Psi}(\mathcal{O})$, then
\begin{equation}\label{230906_13}\text{vol}\,\left(\mf{\Psi}^{-1}(\mathcal{F})\right)\leq \int_{\mathcal{F}}\left|\det \mf{D\Psi}^{-1}(\mf{u})\right|d\mf{u}.\end{equation}Using the same $\varepsilon$ and partition $\mf{P}$ as above,   
for $1\leq \beta\leq s$, let 
\[\mathcal{F}_{\beta}=\mf{\Psi}\left(\mk{D}\cap Q_{\beta}\right).\]
Since $\mk{D}$ and $Q_{\beta}$ are closed Jordan measurable sets, $\mk{D}\cap Q_{\beta}$ is also a closed Jordan measurable set, and so is $\mathcal{F}_{\beta}$. 
Additivity theorem implies that

\[\mathcal{I}=\int_{\mk{D}}g(\mf{x})d\mf{x}=\sum_{\beta=1}^s\int_{\mk{D}\cap Q_{\beta}}g(\mf{x})d\mf{x}.\]
\bp
For each $1\leq\beta\leq s$, since $\mk{D}\cap Q_{\beta}$ is compact, there is a point $\mf{v}_{\beta}\in \mk{D}\cap Q_{\beta}$ such that 
\[g(\mf{x})\leq g(\mf{v}_{\beta})\hspace{1cm}\text{for all}\;\mf{x}\in\mk{D}\cap Q_{\beta}.\]

This gives
\[\int_{\mk{D}\cap Q_{\beta}}g(\mf{x})d\mf{x}\leq g(\mf{v}_{\beta})\,\text{vol}\,(\mk{D}\cap Q_{\beta})=g(\mf{v}_{\beta})\,\text{vol}\,(\mf{\Psi}^{-1}(\mathcal{F}_{\beta})).\]

By \eqref{230906_13}, we find that
\begin{align*}
\int_{\mk{D}\cap Q_{\beta}}g(\mf{x})d\mf{x}\leq g(\mf{v}_{\beta})\, \int_{\mathcal{F}_{\beta}}\left|\det \mf{D\Psi}^{-1}(\mf{u})\right|d\mf{u}.
\end{align*}

Again, there is a point $\mf{w}_{\beta}\in \mk{D}\cap Q_{\beta}$ such that
\[\left|\det \mf{D\Psi}^{-1}(\mf{u})\right|\leq \left|\det \mf{D\Psi}^{-1}(\mf{\Psi}(\mf{w}_{\beta}))\right| \hspace{1cm}\text{for all}\;\mf{u}\in \mathcal{F}_{\beta}.\]
 
This implies that
\[\int_{\mathcal{F}_{\beta}}\left|\det \mf{D\Psi}^{-1}(\mf{u})\right|d\mf{u}\leq \left|\det \mf{D\Psi }(\mf{w}_{\beta}))\right|^{-1}\,\text{vol}\,(\mathcal{F}_{\beta}).\] Hence,
\[\int_{\mk{D}\cap Q_{\beta}}g(\mf{x})d\mf{x}\leq \frac{g(\mf{v}_{\beta})}{g(\mf{w}_{\beta})}\text{vol}\,(\mathcal{F}_{\beta})=\text{vol}\,(\mathcal{F}_{\beta})\left(1+\frac{g(\mf{v}_{\beta})-g(\mf{w}_{\beta})}{g(\mf{w}_{\beta})}\right).\] 
Now since $\mf{v}_{\beta}$ and $\mf{w}_{\beta}$ are in $Q_{\beta}$, $\Vert\mf{v}_{\beta}-\mf{w}_{\beta}\Vert<\delta_1$. Thus,
\[\left|\frac{g(\mf{v}_{\beta})-g(\mf{w}_{\beta})}{g(\mf{w}_{\beta})}\right|\leq \frac{1}{m_g}\left|g(\mf{v}_{\beta})-g(\mf{w}_{\beta})\right|<\varepsilon.\]
This gives
\[\int_{\mk{D}\cap Q_{\beta}}g(\mf{x})d\mf{x}\leq (1+\varepsilon)\,\text{vol}\left(\mf{\Psi}\left(\mk{D}\cap Q_{\beta}\right)\right).\]
Summing over $\beta$ and using additivity theorem, we find that
\[\mathcal{I}=\int_{\mk{D}}g(\mf{x})d\mf{x}\leq  (1+\varepsilon)\sum_{\beta=1}^s\text{vol}\left(\mf{\Psi}\left(\mk{D}\cap Q_{\beta}\right)\right)
=(1+\varepsilon)\,\text{vol}\left(\mf{\Psi}(\mk{D})\right).\]
\bp
Taking $\varepsilon\to 0^+$ gives the desired inequality
\[\mathcal{I}\leq \text{vol}\,\left(\mf{\Psi}(\mk{D})\right).\]
This completes the proof of the change of volume theorem.
\end{myproof}

To conclude the proof of the change of variables theorem, we need the following generalization of the mean value theorem for integrals.
 
\begin{theorem}{Generalized Mean Value Theorem for Integrals}
Let $\mk{D}$ be a compact Jordan measurable set, and let $f:\mk{D}\to\mb{R}$ and $g:\mk{D}\to\mb{R}$ be   continuous functions. If $\mk{D}$ is connected or path-connected, and $g(\mf{x})\geq 0$ for all $\mf{x}\in\mk{D}$, then there is a point $\mf{x}_0$ in $\mk{D}$ such that 
\[\int_{\mk{D}}f(\mf{x})g(\mf{x})d\mf{x}=f(\mf{x}_0) \int_{\mk{D}} g(\mf{x})d\mf{x}.\]
\end{theorem}
The proof of this generalized mean value theorem is almost the same as the mean value theorem. The latter  can be considered as the special case where $g(\mf{x})=1$ for all $\mf{x}$ in $\mk{D}$. 
\begin{myproof}{Proof}

Since $\mk{D}$ is compact and $f:\mk{D}\to\mb{R}$ is continuous, extreme value theorem asserts that there exist points $\mf{u}$ and $\mf{v}$ in $\mk{D}$ such that
\[f(\mf{u})\leq f(\mf{x})\leq f(\mf{v})\hspace{1cm}\text{for all}\;\mf{x}\in\mk{D}.\]
Since $g(\mf{x})\geq 0$ for all $\mf{x}\in\mk{D}$, we find that
\[f(\mf{u}) g(\mf{x})\leq f(\mf{x})g(\mf{x})\leq f(\mf{v})g(\mf{x})\hspace{1cm}\text{for all}\;\mf{x}\in\mk{D}.\]
 The monotonicity theorem implies that
\[f(\mf{u})\int_{\mk{D}} g(\mf{x})d\mf{x}\leq \int_{\mk{D}}f(\mf{x})g(\mf{x})d\mf{x}\leq f(\mf{v})\int_{\mk{D}} g(\mf{x})d\mf{x}.\]
\bp
Let\[U= \int_{\mk{D}} g(\mf{x})d\mf{x}.\]

If $U= 0$, we can take $\mf{x}_0$ to be any point in $\mk{D}$. 
If $U\neq 0$, notice that
\[c=\frac{1}{U}\int_{\mk{D}}f(\mf{x}) g(\mf{x})d\mf{x}\]
satisfies
\[f(\mf{u})\leq c\leq f(\mf{v}).\]
As in the proof of the mean value theorem, 
  $\mk{D}$ is connected or path-connected allows us to conclude that there is an $\mf{x}_0$ in $\mk{D}$ such that
\[\frac{1}{U}\int_{\mk{D}}f(\mf{x}) g(\mf{x})d\mf{x}=c=f(\mf{x}_0).\]
 This gives 
\[\int_{\mk{D}}f(\mf{x})g(\mf{x})d\mf{x}=f(\mf{x}_0) \int_{\mk{D}} g(\mf{x})d\mf{x}.\]
\end{myproof}
Next, we prove the special case of the change of variables theorem when $\mk{D}$ is a closed rectangle.
\begin{theorem}[label=230907_1]{}
Let $\mathcal{O}$ be an open subset of $\mb{R}^n$, and let $\mf{\Psi}:\mathcal{O}\to\mb{R}^n$ be a smooth change of variables. If $\mf{I}$ is a closed rectangle contained in $\mathcal{O}$,  and $f:\mf{\Psi}(\mf{I})\to\mb{R}$  is    continuous, then
\[ \int_{\mf{\Psi}(\mf{I})}f(\mf{x})d\mf{x}=\int_{\mf{I}}f\left(\mf{\Psi}(\mf{x})\right)\left|\det \mf{D\Psi}(\mf{x})\right|d\mf{x}.\]
\end{theorem}
\begin{myproof}{ Proof}
It is sufficient to show that for any $\varepsilon>0$,
\[\left| \int_{\mf{\Psi}(\mf{I})}f(\mf{x})d\mf{x}-\int_{\mf{I}}f\left(\mf{\Psi}(\mf{x})\right)\left|\det \mf{D\Psi}(\mf{x})\right|d\mf{x}\right|<\varepsilon.\]
\bp
Since $f:\mf{\Psi}(\mf{I})\to\mb{R}$  and $\mf{\Psi}:\mf{I}\to\mb{R}^n$ are continuous, $(f\circ \mf{\Psi}):\mf{I}\to\mb{R}$  is    continuous. Since $\mf{I}$ is compact, $(f\circ \mf{\Psi}):\mf{I}\to\mb{R}$   is uniformly continuous.  
Given $\varepsilon>0$, there is a $\delta>0$ such that if $\mf{u}$ and $\mf{v}$ are points in $\mf{I}$,
\[\left|f(\mf{\Psi}(\mf{u}))-f(\mf{\Psi}(\mf{v}))\right|<\frac{\varepsilon}{\text{vol}\,(\mf{\Psi}(\mf{I}))}.\]

Let $\mf{P}$ be a partition of $\mf{I}$ such that $|\mf{P}|<\delta$.
By additivity theorem,
\[\int_{\mf{\Psi}(\mf{I})}f(\mf{x})d\mf{x}=\sum_{\mf{J}\in\mathcal{J}_{\mf{P}}}\int_{\mf{\Psi}(\mf{J})}f(\mf{x})d\mf{x},\]
and 
\[\int_{\mf{I}}f\left(\mf{\Psi}(\mf{x})\right)\left|\det \mf{D\Psi}(\mf{x})\right|d\mf{x}=\sum_{\mf{J}\in\mathcal{J}_{\mf{P}}} \int_{\mf{J}}f\left(\mf{\Psi}(\mf{x})\right)\left|\det \mf{D\Psi}(\mf{x})\right|d\mf{x}.\]
Each $\mf{J}\in\mathcal{J}_{\mf{P}}$  is a closed rectangle. Hence, it is compact and path-connected. Since $\mf{\Psi}:\mf{I}\to\mb{R}^n$ is continuous, $\mf{\Psi}(\mf{I})$ is also compact and path-conneted.  By the generalized mean value theorem, for each $\mf{J}\in \mathcal{J}_{\mf{P}}$, there exist $\mf{u}_{\mf{J}}$ and $\mf{v}_{\mf{J}}$ in $\mf{J}$ such that
\[ \int_{\mf{\Psi}(\mf{J})}f(\mf{x})d\mf{x}=f(\mf{\Psi}(\mf{u}_{\mf{J}}))\,\text{vol}\,\left(\mf{\Psi}(\mf{J})\right)\]and\[
\int_{\mf{J}}f\left(\mf{\Psi}(\mf{x})\right)\left|\det \mf{D\Psi}(\mf{x})\right|d\mf{x}=f(\mf{\Psi}(\mf{v}_{\mf{J}}))\,\int_{\mf{J}} \left|\det \mf{D\Psi}(\mf{x})\right|d\mf{x}.\]By the change of volume theorem, 
\[\int_{\mf{J}} \left|\det \mf{D\Psi}(\mf{x})\right|d\mf{x}=\text{vol}\,\left(\mf{\Psi}(\mf{J})\right).\]
Since $\Vert\mf{u}_{\mf{J}}-\mf{v}_{\mf{J}}\Vert\leq\text{diam}\,\mf{J}<\delta$, we find that
\[\left|f(\mf{\Psi}(\mf{u}_{\mf{J}}))-f(\mf{\Psi}(\mf{v}_{\mf{J}}))\right|<\frac{\varepsilon}{\text{vol}\,(\mf{\Psi}(\mf{I}))}.\]
\bp
It follows that
\begin{align*}
&\left| \int_{\mf{\Psi}(\mf{I})}f(\mf{x})d\mf{x}-\int_{\mf{I}}f\left(\mf{\Psi}(\mf{x})\right)\left|\det \mf{D\Psi}(\mf{x})\right|d\mf{x}\right|\\
&\leq\sum_{\mf{J}\in\mathcal{J}_{\mf{P}}}\left| \int_{\mf{\Psi}(\mf{J})}f(\mf{x})d\mf{x}-\int_{\mf{J}}f\left(\mf{\Psi}(\mf{x})\right)\left|\det \mf{D\Psi}(\mf{x})\right|d\mf{x}\right|\\
&=\sum_{\mf{J}\in\mathcal{J}_{\mf{P}}}\left|f(\mf{\Psi}(\mf{u}_{\mf{J}}))-f(\mf{\Psi}(\mf{v}_{\mf{J}}))\right|\,\text{vol}\,\left(\mf{\Psi}(\mf{J})\right)\\
&<\frac{\varepsilon}{\text{vol}\,(\mf{\Psi}(\mf{I}))}\sum_{\mf{J}\in\mathcal{J}_{\mf{P}}}\text{vol}\,\left(\mf{\Psi}(\mf{J})\right)=\frac{\varepsilon}{\text{vol}\,(\mf{\Psi}(\mf{I}))}\times \text{vol}\,(\mf{\Psi}(\mf{I}))=\varepsilon.
\end{align*}This completes the proof.
\end{myproof}

Finally, we conclude the general case.
\begin{myproof}{ 
Conclusion of the Proof of the Change of Variables Theorem}
As in the special case, it is sufficient to show that for any $\varepsilon>0$,
\[\left| \int_{\mf{\Psi}(\mk{D})}f(\mf{x})d\mf{x}-\int_{\mk{D}}f\left(\mf{\Psi}(\mf{x})\right)\left|\det \mf{D\Psi}(\mf{x})\right|d\mf{x}\right|<\varepsilon.\]
Let us first proceed as in the proof of change of volume theorem.   There is a positive number $d$ and a compact set $C$ such that $\overline{\mk{D}}\subset C\subset \mathcal{O}$, and any point in $\mb{R}^n$ that has a distance less than $d$ from a point in $ \overline{\mk{D}}$ lies in $C$. On the compact set $C$, the function  $g:C\to\mb{R}$, $g(\mf{x})=|\det\mf{D\Psi}(\mf{x})|$  is continuous. Therefore, it is bounded. Namely, there is a positive number $M_g$ so that
\[\left|\det\mf{D\Psi}(\mf{x})\right|\leq M_g\quad \hspace{1cm}\text{for all}\;\mf{x}\in C.\]
Since we assume that the function $f:\mf{\Psi}(\mk{D})\to\mb{R}$ is bounded, there is a positive number $M_f$ such that
\[\left|f(\mf{\Psi}(\mf{x}))\right|\leq M_f\quad \hspace{1cm}\text{for all}\;\mf{x}\in \mk{D}.\]

\bp
Let $M=M_fM_g$, and let $\mf{I}$ be a rectangle that contains $\mk{D}$. Given $\varepsilon>0$, let $\mf{P}$ be a partition of $\mf{I}$ such that $|\mf{P}|<d$ and 
\[ L(\chi_{\mk{D}},\mf{P})>\text{vol}\,(\mk{D})-\frac{\varepsilon}{2M}.\]

Let 
\begin{align*}\mathscr{A} =\left\{\mf{J}\in\mathcal{J}_{\mf{P}}\,|\,\mf{J}\cap\mk{D}\neq\emptyset\right\},\hspace{1cm}
\mathscr{B} =\left\{\mf{J}\in\mathcal{J}_{\mf{P}}\,|\,\mf{J}\subset\mk{D}\right\}.\end{align*}

Since $|\mf{P}|<d$, each $\mf{J}$ in $\mathscr{A}$ is contained in $C$. 
Moreover, 
\[L(\chi_{\mk{D}},\mf{P})=\sum_{\mf{J}\in \mathscr{B}}\text{vol}\,(\mf{J}),\]
Denote by $\mathcal{Q}$ the set
\[\mathcal{Q}=\bigcup_{\mf{J}\in\mathscr{B}}\mf{J}.\]
Then $\mathcal{Q}$ is a compact subset of $\mk{D}$, $\mathcal{S}=\mk{D}\setminus \mathcal{Q}$ is Jordan measurable, 
and
\begin{align*}
\text{vol} (\mk{D}\setminus\mathcal{Q})=\text{vol}\,(\mk{D})-\sum_{\mf{J}\in\mathscr{B}}\text{vol}\,(\mf{J})<\frac{\varepsilon}{2M}.
\end{align*}
By additivitity theorem,
\begin{align*}
\int_{\mf{\Psi}(\mk{D})}f(\mf{x})d\mf{x}=\sum_{\mf{J}\in\mathscr{B}}\int_{\mf{\Psi}(\mf{J})}f(\mf{x})d\mf{x}+\int_{\mf{\Psi}(\mk{D}\setminus\mathcal{Q})}f(\mf{x})d\mf{x},
\end{align*}
\begin{align*}
\int_{\mk{D}}f\left(\mf{\Psi}(\mf{x})\right)\left|\det \mf{D\Psi}(\mf{x})\right|d\mf{x}&=\sum_{\mf{J}\in\mathscr{B}}\int_{\mf{J}}f\left(\mf{\Psi}(\mf{x})\right)\left|\det \mf{D\Psi}(\mf{x})\right|d\mf{x} \\&\quad +\int_{\mk{D}\setminus\mathcal{Q}}f\left(\mf{\Psi}(\mf{x})\right)\left|\det \mf{D\Psi}(\mf{x})\right|d\mf{x}.
\end{align*}
Theorem \ref{230907_1} says that for each $\mf{J}$ in $\mathscr{B}$, 
\[\int_{\mf{\Psi}(\mf{J})}f(\mf{x})d\mf{x}=\int_{\mf{J}}f\left(\mf{\Psi}(\mf{x})\right)\left|\det \mf{D\Psi}(\mf{x})\right|d\mf{x}.\]
\bp
Therefore,
\begin{align*}
&\left| \int_{\mf{\Psi}(\mk{D})}f(\mf{x})d\mf{x}-\int_{\mk{D}}f\left(\mf{\Psi}(\mf{x})\right)\left|\det \mf{D\Psi}(\mf{x})\right|d\mf{x}\right|\\
&\leq \left|\int_{\mf{\Psi}(\mk{D}\setminus\mathcal{Q})}f(\mf{x})d\mf{x}\right|+\left|\int_{\mk{D}\setminus\mathcal{Q}}f\left(\mf{\Psi}(\mf{x})\right)\left|\det \mf{D\Psi}(\mf{x})\right|d\mf{x}\right|.
\end{align*}
For the term $\di \int_{\mf{\Psi}(\mk{D}\setminus\mathcal{Q})}f(\mf{x})d\mf{x}$, we have
\[\left|\int_{\mf{\Psi}(\mk{D}\setminus\mathcal{Q})}f(\mf{x})d\mf{x}\right|\leq M_f \,\text{vol} \left(\mf{\Psi}(\mk{D}\setminus\mathcal{Q})\right).\]
By the change of volume theorem,
\[\text{vol} \left(\mf{\Psi}(\mk{D}\setminus\mathcal{Q})\right)=\int_{\mk{D}\setminus\mathcal{Q}} \left|\det \mf{D\Psi}(\mf{x})\right|d\mf{x}
\leq M_g \text{vol}\left( \mk{D}\setminus\mathcal{Q}\right).\]
Therefore,
\[\left|\int_{\mf{\Psi}(\mk{D}\setminus\mathcal{Q})}f(\mf{x})d\mf{x}\right|\leq M_fM_g \text{vol}\left( \mk{D}\setminus\mathcal{Q}\right)=M\text{vol}\left( \mk{D}\setminus\mathcal{Q}\right)<\frac{\varepsilon}{2}.\]
Similarly, 
\[\left|\int_{\mk{D}\setminus\mathcal{Q}}f\left(\mf{\Psi}(\mf{x})\right)\left|\det \mf{D\Psi}(\mf{x})\right|d\mf{x}\right|\leq M\text{vol}\left( \mk{D}\setminus\mathcal{Q}\right)<\frac{\varepsilon}{2}.\]
This gives
\[\left| \int_{\mf{\Psi}(\mk{D})}f(\mf{x})d\mf{x}-\int_{\mk{D}}f\left(\mf{\Psi}(\mf{x})\right)\left|\det \mf{D\Psi}(\mf{x})\right|d\mf{x}\right|<\varepsilon,\]
which completes the proof.
\end{myproof}

\section{Some Important Integrals and Their Applications} 

Up to now we have only discussed multiple integrals for bounded functions $f:\mk{D}\to\mb{R}$ defined on bounded domains. For practical applications, we need to consider improper integrals where the function is not bounded or the domain is not bounded. As in the single variable case, we need to take limits. In the  multi-variable case, things become considerably more complicated. Interested readers can read the corresponding sections in the book \cite{Zorich_2}. 
 In this section, we  use theories learned in multiple integrals to derive some  explicit formulas of improper integrals of single-variable functions, without introducing the definition of improper multiple integrals. We then give some applications of these formulas.
 \begin{proposition}[label=230911_5]{}
 For any positive number $a$,
 \[\int_{-\infty}^{\infty} e^{-ax^2}dx=\sqrt{\frac{\pi}{a}}.\]
 \end{proposition}
 \begin{myproof}{Proof}
 Since the function $f:\mb{R}\to\mb{R}$, $f(x)=e^{-ax^2}$ is positive for all $x\in\mb{R}$,
 \[\int_{-\infty}^{\infty} e^{-ax^2}dx=\lim_{L\to \infty} \int_{-L}^{L} e^{-ax^2}dx.\]
 Given a positive number $R$, we consider the double integral
 \[I_R=\int_{\overline{B(\mf{0}, R)}}e^{-a(x^2+y^2)}dxdy.\]
 For any positive number $L$,
\[\overline{B(\mf{0},L)}\;\subset\; [-L,L]\times [-L,L]\;\subset\; \overline{B(\mf{0},\sqrt{2}L)}.\]
Since the function $g:\mb{R}^2\to\mb{R}$, $g(x)=e^{-a(x^2+y^2)}$ is positive,
\begin{equation}\label{230911_7}I_L\;\leq \;
\int_{[-L,L]\times [-L,L]}e^{-a(x^2+y^2)}dxdy\;\leq\; I_{\sqrt{2}L}.\end{equation}
\bp
 Using polar coordinates, we find that
 \[I_R=\int_0^{2\pi}\int_0^R e^{-ar^2}rdrd\theta= 2\pi \left[-\frac{e^{-ar^2}}{2a}\right]_0^R=\frac{\pi}{a}\left(1-e^{-aR^2}\right).\]
 Thus,
 \[\lim_{R\to\infty} I_R=\frac{\pi}{a}.\]
 Eq. \eqref{230911_7} then implies that 
  \[\lim_{L\to\infty} \int_{[-L,L]\times [-L,L]}e^{-a(x^2+y^2)}dxdy=\frac{\pi}{a}.\]
   
 By Fubini's theorem,
 \[ \int_{[-L,L]\times [-L,L]}e^{-a(x^2+y^2)}dxdy=\left(\int_{-L}^Le^{-x^2}dx\right)^2.\]
Thus, we conclude that 
 \[\int_{-\infty}^{\infty} e^{-ax^2}dx=\lim_{L\to \infty} \int_{-L}^{L} e^{-ax^2}dx=\sqrt{\frac{\pi}{a}}.\]

 \end{myproof}

\begin{figure}[ht]
\centering
\includegraphics[scale=0.18]{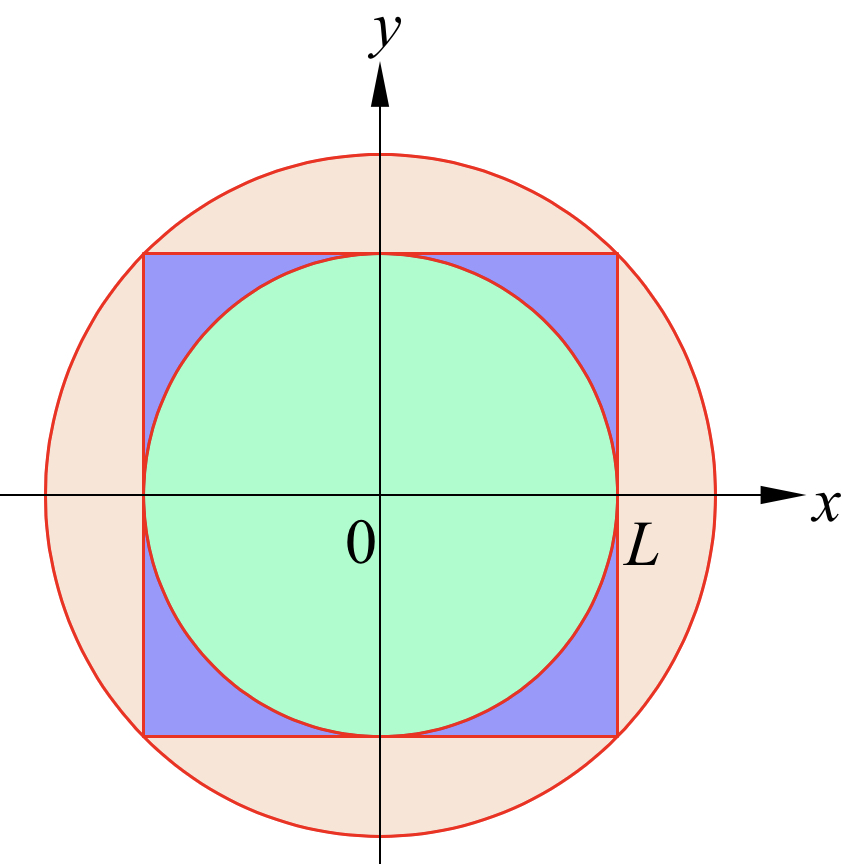}
\caption{$\overline{B(\mf{0},L)}\;\subset\; [-L,L]\times [-L,L]\;\subset\; \overline{B(\mf{0},\sqrt{2}L)}$. }\label{figure132}
\end{figure}
 The improper integral $\di\int_{-\infty}^{\infty} e^{-ax^2}dx$ with $a>0$  plays an important role in various areas of mathematics. For example, in probability theorem, the probability density function of a normal random variable with mean $\mu$ and standard deviation $\sigma$ is given by
 \[f(x)=\frac{1}{\sqrt{2\pi}\,\sigma}\exp\left(-\frac{(x-\mu)^2}{2\sigma^2}\right).\]The normalization factor $1/(\sqrt{2\pi}\sigma)$ is required such that
 \[\int_{-\infty}^{\infty}f(x)dx=1,\]which ensures that total probability is 1.
 
 Recall that we have defined the gamma function $\Gamma(s)$ for a real number $s>0$ by the improper integral 
 \[\Gamma(s)=\int_0^{\infty} t^{s-1}e^{-t}dt.\]
 The value of $\Gamma(1)$ is easy to compute.
 \[\Gamma(1)=\int_0^{\infty}e^{-t}dt=1.\]
 Using integration by parts, one can show that
 \begin{equation}\label{230911_8}
 \Gamma(s+1)=s\Gamma(s)\hspace{1cm}\text{when}\,s>0.
 \end{equation}
 From this, we find that
 \[\Gamma(n+1)=n!\hspace{1cm}\text{for all}\; n\in\mb{Z}^+.\]
 The value of $\Gamma(s)$ when $s=1/2$ is also of particular interest. 
 \begin{theorem}{}
 The value of the gamma function $\Gamma(s)$ at $s=1/2$ is
 \[\Gamma\left(\frac{1}{2}\right)=\sqrt{\pi}.\]
 
 \end{theorem}
 \begin{myproof}{Proof}
 By Proposition \ref{230911_5}, 
 \[\sqrt{\pi}=2\int_0^{\infty}e^{-x^2}dx= 2\lim_{a\to 0^+}\lim_{L\to\infty}\int_a^Le^{-x^2}dx.\]
 \bp
 Making a change of variables $t=x^2$, we find that
 \[2\int_a^L e^{-x^2}dx=\int_{a^2}^{L^2}t^{-\frac{1}{2}}e^{-t}dt.\]
 Therefore,
 \begin{align*}
 \Gamma\left(\frac{1}{2}\right)&=\int_0^{\infty}t^{-\frac{1}{2}}e^{-t}dt=\lim_{a\to 0^+}\lim_{L\to\infty} \int_{a^2}^{L^2}t^{-\frac{1}{2}}e^{-t}dt\\
 &=2\lim_{a\to 0^+}\int_a^Le^{-x^2}dx=\sqrt{\pi}.
 \end{align*}

 \end{myproof}

 Another useful formula we have mentioned in volume I is the formula for the beta function $B(\alpha,\beta)$ defined as
 \[B(\alpha,\beta)=\int_0^1 t^{\alpha-1}(1-t)^{\beta-1}dt\hspace{1cm}\text{when}\;\alpha>0, \beta>0.\]
 It is easy to show that the integral  is indeed convergent when $\alpha$ and $\beta$ are positive. 
 We have the following recursive formula.
 \begin{lemma}[label=230911_11]{}
 For $\alpha>0$ and $\beta>0$, we have
 \[B(\alpha,\beta)=\frac{(\alpha+\beta+1)(\alpha+\beta)}{\alpha\beta}B(\alpha+1,\beta+1).\]
 \end{lemma}
 
 \begin{myproof}{Proof}

 First notice that for $\alpha>0$ and $\beta>0$,
 \begin{align*}
 \int_0^1t^{\alpha-1}(1-t)^{\beta-1}dt&=\int_0^1t^{\alpha-1}(1-t)^{\beta-1}(1-t+t)dt\\
 &= \int_0^1t^{\alpha-1}(1-t)^{\beta}dt+ \int_0^1t^{\alpha}(1-t)^{\beta-1}dt.
 \end{align*}
This gives
 \[B(\alpha,\beta)=B(\alpha+1,\beta)+B(\alpha,\beta+1).\]
 \bp
 Apply this formula again to the two terms on the right, we find that
 \[B(\alpha,\beta)=B(\alpha+2, \beta)+2B(\alpha+1,\beta+1)+B(\alpha,\beta+2).\]
 
 Using integration by parts, one can show that when $\alpha>0$ and $ \beta>0$,
 \[ \int_0^1t^{\alpha}(1-t)^{\beta-1}dt=\frac{\alpha}{\beta}\int_0^1t^{\alpha-1}(1-t)^{\beta}dt.\]This gives
 \[B(\alpha+1,\beta) =\frac{\alpha}{\beta}B(\alpha,\beta+1),\]
 Therefore,
 \begin{align*}
 B(\alpha,\beta)&=\left(\frac{\alpha+1}{\beta}+2+\frac{\beta+1}{\alpha}\right)B(\alpha+1,\beta+1)\\
 &=\frac{(\alpha+\beta+1)(\alpha+\beta)}{\alpha\beta}B(\alpha+1,\beta+1).
 \end{align*}
 
 \end{myproof}
 Now we can derive the explicit formula for the beta function.
 \begin{theorem}{The Beta Function}
 For any positive real numbers $\alpha$  and $\beta$,
 \[B(\alpha,\beta)=\int_0^1 t^{\alpha-1}(1-t)^{\beta-1}dt=\frac{\Gamma(\alpha)\Gamma(\beta)}{\Gamma(\alpha+\beta)}.\]
 \end{theorem}
 \begin{myproof}{Proof}
 We first consider the case where $\alpha>1$ and $\beta>1$. Let $g:[0,\infty)\times [0,\infty)\to\mb{R}$ be the function defined as
 \[g(u,v)=u^{\alpha-1}v^{\beta-1} e^{-u-v}.\]
 This is a continuous function.
 For $L>0$, let 
 \begin{align*}\mathcal{U}_L&=\left\{(t,w)\,|\, 0\leq t\leq 1, 0\leq w\leq L\right\},\\\mk{D}_L&=\left\{(u,v)\,|\,  u\geq 0, v\geq 0, u+v\leq L\right\}.\end{align*}
\bp
Consider the map
 $\mf{\Psi}:\mb{R}^2\to \mb{R}^2$ defined as
 \[(u,v)=\mf{\Psi}(t,w)=(tw, (1-t)w).\]
 
 Notice that $\mf{\Psi}$ maps the interior of $\mathcal{U}_L$ one-to-one onto the interior of $\mk{D}_L$. The Jacobian of this map is
 \[\frac{\pa (u,v)}{\pa (t, w)}=\det \begin{bmatrix} w & t\\ -w & 1-t\end{bmatrix}=w.\]
 Thus, $\mf{\Psi}:(0,1)\times (0,\infty)\to\mb{R}^2$ is a smooth change of variables. By taking  limits, the change of variables theorem implies that
 \begin{equation}\label{230911_10}\int_{\mathcal{U}_L} (g\circ \mf{\Psi})(t,w)\frac{\pa (u,v)}{\pa (t,w)}dtdw=\int_{\mk{D}_L}g(u,v)dudv.\end{equation}
 Now Fubini's theorem says that
 \begin{align*}
 \int_{\mathcal{U}_L} (g\circ \mf{\Psi})(t,w)\frac{\pa (u,v)}{\pa (t,w)}dtdw=\int_0^1 t^{\alpha-1}(1-t)^{\beta-1}dt\int_0^L w^{\alpha+\beta-1}e^{-w}dw.
 \end{align*}
 
 Thus,
 \[\lim_{L\to\infty} \int_{\mathcal{U}_L} (g\circ \mf{\Psi})(t,w)\frac{\pa (u,v)}{\pa (t,w)}dtdw=\Gamma(\alpha+\beta) B(\alpha,\beta).\]

 On the other hand, we notice that
 \[\left[0, \frac{L}{2}\right]^2\;\subset\; \mk{D}_L\;\subset \;[0,L]^2.\]
 Since the function $g:[0,\infty)\times [0,\infty)\to\mb{R}$ is nonnegative, this implies that
 \begin{equation}\label{230911_9}I_{\frac{L}{2}}\leq \int_{\mk{D}_L}g(u,v)dudv\leq I_L,\end{equation}
 where
 \[I_L=\int_{[0,L]^2}g(u,v)dudv=\int_0^Lu^{\alpha-1}e^{-u}du\int_0^Lv^{\beta-1}e^{-v}dv.\]
 \bp
 By the definition of the gamma function,
 \[\lim_{L\to\infty} I_L=\Gamma(\alpha)\Gamma(\beta).\]
 
 Eq. \eqref{230911_9} then implies that
 \[\lim_{L\to \infty} \int_{\mk{D}_L}g(u,v)dudv=\Gamma(\alpha)\Gamma(\beta).\]

 It follows from \eqref{230911_10}  that
 \[\Gamma(\alpha+\beta) B(\alpha,\beta)=\Gamma(\alpha)\Gamma(\beta),\]
 which gives the desired formula when $\alpha>1$ and $\beta>1$.
 
 For the general case where $\alpha>0$ and $\beta>0$,  Lemma \ref{230911_11} and \eqref{230911_8} give
 \begin{align*}
 B(\alpha,\beta)&=\frac{(\alpha+\beta+1)(\alpha+\beta)}{\alpha\beta}\frac{\Gamma(\alpha+1)\Gamma(\beta+1)}{\Gamma(\alpha+\beta+2)}=
 \frac{\Gamma(\alpha)\Gamma(\beta)}{\Gamma(\alpha+\beta)}.
 \end{align*}
 This completes the proof.
 \end{myproof}
 
 \begin{figure}[ht]
\centering
\includegraphics[scale=0.18]{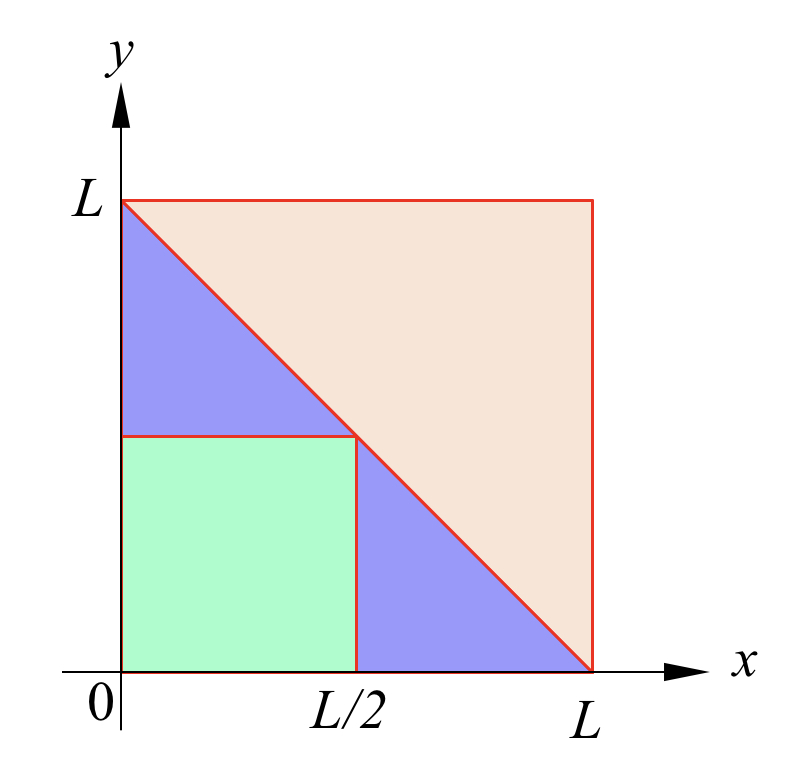}
\caption{$\di \left[0, \frac{L}{2}\right]^2\;\subset\; \left\{(u,v)\,|\,  u\geq 0, v\geq 0, u+v\leq L\right\}\;\subset \;[0,L]^2$. }\label{figure133}
\end{figure}
 
 Now we give an interesting application of the formula of the beta function.
 \begin{theorem}{}
 For $n\geq 1$, the volume of the $n$-ball of radius $a$,
 \[B_n(a)=\left\{(x_1, \ldots, x_n)\in \mb{R}^n\,|\, x_1^2+\cdots+x_n^2\leq  a^2\right\},\] is equal to
 \[V_n(a)=\frac{\pi^{\frac{n}{2}}a^n}{\Gamma\left(\frac{n+2}{2}\right)}.\]
 \end{theorem}
 \begin{myproof}{Proof}
 It is easy to see that for any $a>0$, 
 \[V_n(a)=V_na^n,\hspace{1cm}\text{where}\;V_n=V_n(1).\]
 Since $B_1(1)=[-1,1]$, we find that $V_1=2$.
 For $n\geq 2$, 
 notice that for fixed   $-1\leq y\leq 1$, the ball $B_n(1)$ intersects the plane $x_n=y$ on the set
 \[S_n(y)=\left\{(x_1, \ldots, x_{n-1},y)\,|\, x_1^2+\ldots+x_{n-1}^2\leq 1-y^2\right\}.\]
   Fubini's theorem implies that
 \begin{align*}
 V_n&=\int_{-1}^1 \int_{B_{n-1}\left(\sqrt{1-y^2}\right)}dx_1\ldots dx_{n-1}dy=\int_{-1}^1(1-y^2)^{\frac{n-1}{2}}V_{n-1}dy\\
 &=2V_{n-1}\int_0^1(1-y^2)^{\frac{n-1}{2}}dy= V_{n-1}\int_0^{1}t^{-\frac{1}{2}}(1-t)^{\frac{n-1}{2}}dt\\
 &=V_{n-1}B\left(\frac{1}{2}, \frac{n+1}{2}\right)=\sqrt{\pi} V_{n-1} \frac{\Gamma\left(\frac{n+1}{2}\right)}{\Gamma\left(\frac{n+2}{2}\right)}.
 \end{align*}
 This formula is still correct when $n=1$ if we define $V_0=1$. Therefore,
 \begin{align*}
 V_n =\prod_{k=1}^n\frac{V_k}{V_{k-1}}= \prod_{k=1}^n\left( \sqrt{\pi}\frac{\Gamma\left(\frac{k+1}{2}\right)}{\Gamma\left(\frac{k+2}{2}\right)}\right)=\frac{\pi^{\frac{n}{2}}}{\Gamma\left(\frac{n+2}{2}\right)}.
 \end{align*}
 \end{myproof}

\chapter{Fourier Series and Fourier Transforms} 
 
 In this chapter, we shift our attention to the theory of Fourier series and Fourier transforms. In volume I, we have considered  expansions  of functions as   power series, which are limits of polynomials.  In this chapter, we consider expansions of functions in another class of infinitely differentiable functions -- the trigonometric functions $\sin x$ and $\cos x$. The reason to consider $\sin x$ and $\cos x$ is that they are representative of periodic functions.
 
Recall that
a  function $f:\mathbb{R}\to\mathbb{R}$ is said to be periodic if there is a positive number $p$ so that
\[f(x+p)=f(x)\hspace{1cm}\text{for all}\; x\in \mathbb{R}.\]
Such a number $p$ is called a period of the function $f$. If $p$ is a period of $f$, then for any positive integer $n$, $np$ is also a period of $f$.

The functions $\sin x$ and $\cos x$ are periodic functions of period $2\pi$. If $f:\mb{R}\to \mb{R}$ is a periodic function of period $p=2L$, then the function $g:\mb{R}\to\mb{R}$ defined as 
\[g(x)=f\left(\frac{\pi x}{L}\right)\] is periodic of period $2\pi$. Hence, we can concentrate on functions that are periodic of period $2\pi$.
 
The celebrated Euler formula
\[e^{ix}=\cos x+i\sin x\] connects the trigonometric functions $\sin x$, $\cos x$ with the exponential function  with imaginary arguments. Hence, in this chapter, we shift our paradigm and consider complex-valued functions $f:D\to\mb{C}$ defined on a subset $D$ of $\mb{R}$. Since a complex number $z=x+iy$ with real part $x$ and imaginary part $y$ can be identified with the point $(x,y)$ in $\mb{R}^2$, such a function can be regarded as a function $f:D\to\mb{R}^2$, so that derivative and integrals are defined componentwise. More precisely, given $x\in D$, we write $f(x)=u(x)+iv(x)$, where $u(x)$ and $v(x)$ are respectively the real part and imaginary part of $f(x)$. If $x_0$ is an interior point of $D$, we say that $f$ is differentiable at $x_0$ if
the limit
\[\lim_{x\to x_0}\frac{f(x)-f(x_0)}{x-x_0} \]
exists. This is if and only if both $u:D\to\mb{R}$ and $v:D\to\mb{R}$ are differentiable at $x_0$, and we have
\[f'(x_0)=u'(x_0)+iv'(x_0).\]
Similarly, if $[a,b]$ is a closed interval that is contained in $D$, we say that $f$ is Riemann integrable over $[a,b]$ if and only if both $u$ and $v$ are Riemann integrable over $[a,b]$, and we have
\[\int_a^b f(x)dx=\int_a^b u(x)dx+ i\int_a^b v(x)dx.\]
If $F:[a,b]\to\mb{C}$ is a continuously differentiable function, the fundamental theorem of calculus implies that
\[F(b)-F(a)=\int_a^b F'(x)dx.\]
 
\section{Orthogonal Systems of Functions and Fourier Series } \label{sec7.1}
In the following, let $I=[a,b]$ be a compact interval in $\mb{R}$   unless otherwise specified. Denote by $\mathcal{R}(I, \mb{C})$ the set of all complex-valued functions $
f:I\to \mb{C}$ that are Riemann integrable. Given two functions $f$ and $g$ in $\mathcal{R}(I, \mb{C})$,   their sum  $f+g$ is the function $(f+g):I\to\mb{C}$,
\[(f+g)(x)=f(x)+g(x).\] If $\alpha$ is a complex number,   the scalar product of $\alpha$ with $f$ is the function $(\alpha f):I\to\mb{C}$, where
\[(\alpha f)(x)=\alpha f(x).\]
With the addition and scalar multiplication thus defined, $\mathcal{R}(I,\mb{C})$  is a complex vector space. From the theory of integration, we know that  the set of complex-valued continuous functions on $I$, denoted by $C(I, \mb{C})$, is a subspace of  $\mathcal{R}(I,\mb{C})$.

If $f:I\to\mb{C}$ is Riemann  integrable, so does its complex conjugate $\overline{f}:I\to\mb{C}$ defined as
\[\overline{f}(x)=\overline{f(x)}.\]
In volume I, we have proved that if two real-valued functions $f:I\to\mb{R}$ and $g:I\to\mb{R}$ are Riemann integrable, so is their product $(fg):I\to\mb{R}$. Using this, it is easy to check that if $f:I\to\mb{C}$ and $g:I\to\mb{C}$ are Riemann integrable complex-valued functions, $(f\overline{g}):I\to\mb{C}$ is also Riemann integrable. 
\begin{proposition}{}
Given $f$ and $g$ in  $\mathcal{R}(I, \mb{C})$,   define
\[\langle f, g\rangle =\int_a^b f(x)\overline{g(x)}dx.\]
For any $f, g, h$ in  $\mathcal{R}(I, \mb{C})$, and any complex numbers $\alpha$ and $\beta$, we have the followings.
\begin{enumerate}[(a)]
 \item $\di \langle g, f\rangle =\overline{\langle f,  g\rangle}$.
 \item $\di \langle  \alpha f+ \beta g, h\rangle =\alpha\langle f, h\rangle+\beta\langle g, h\rangle$.

   \item $\di \langle f, f\rangle\geq 0$.
 \end{enumerate}
 We call $\langle \;\cdot\;,\;\cdot\;\rangle$ a positive {\it semi-definite} inner product on $\mathcal{R}(I, \mb{C})$.
\end{proposition}

It follows from (a) and (b) that
\[  \langle f, \alpha g+ \beta  h\rangle =\overline{\alpha}\langle f, g\rangle+\overline{\beta}\langle f, h\rangle.\]
More generally, we have
\begin{equation}\label{230813_3}\left\langle \sum_{i=1}^m\alpha_i f_i ,\;\sum_{j=1}^n\beta_jg_j\right\rangle
=\sum_{i=1}^m\sum_{j=1}^n\alpha_i\overline{\beta_j}\langle f_i ,g_j\rangle.\end{equation}
If $f(x)=u(x)+iv(x)$, where $u(x)=\text{Re}\,f(x)$ and $v(x)=\text{Im}\,f(x)$, then
\[\langle f, f\rangle =\int_a^b \left(u(x)^2+v(x)^2\right)dx.\]

Notice that $\langle f, f\rangle=0$ does not imply that $f=0$.
For example, take any $c$ in $ [a,b]$, and define the function
  $f:I\to\mb{C}$ by
\[f(x)=\begin{cases}1,\quad &\text{if}\;x=c,\\0,\quad &\text{otherwise}.\end{cases}\]
Then 
\[\langle f, f\rangle =\int_a^b f(x)\overline{f(x)}dx=0\] even though $f$ is not the zero function.
This is why we call 
 $\langle \;\cdot\;,\;\cdot\;\rangle$ a positive semi-definite inner product. Restricted to the subspace of continuous functions $C(I, \mb{C})$, $\langle \;\cdot\;,\;\cdot\;\rangle$ is a  positive definite inner product, or simply an inner product in the usual sense.
 
 Using the positive semi-definite inner product, we can define a semi-norm on $\mathcal{R}(I, \mb{C})$.
 \begin{definition}
 {The $\pmb{L^2}$ Semi-Norm}
 Given $f:I\to\mb{C}$ in $\mathcal{R}(I, \mb{C})$, the semi-norm of $f$ is defined as
 \[\Vert f\Vert =\sqrt{\langle  f,f\rangle }=\sqrt{\int_a^b |f(x)|^2 dx}.\]
 \end{definition}
 
 It has the following properties.
 \begin{proposition}[label=230814_5]{}
 Given $f$ in $\mathcal{R}(I, \mb{C})$ and $\alpha\in\mb{C}$, we have the followings.
 \begin{enumerate}[(a)]
 \item $\Vert f\Vert\geq 0$.
 \item $\Vert\alpha f\Vert=|\alpha|\Vert f\Vert$.
 
 \end{enumerate}
 \end{proposition}
 The Cauchy-Schwarz inequality still holds for the positive semi-definite inner product on $\mathcal{R}(I,\mb{C})$.
 
 \begin{proposition}{Cauchy-Schwarz Inequality}
 Given $f$ and $g$ in $\mathcal{R}(I,\mb{C})$,
 \[\left|\langle f, g\rangle\right|\leq \Vert f\Vert \Vert g\Vert.\]
 \end{proposition}
 The proof is exactly the same as for an inner product on a real vector space.
 An immediate consequence of the Cauchy-Schwarz inequality is the triangle inequality.
 \begin{proposition}[label=230814_6]{Triangle Inequality}
 Let  $f_1, \ldots, f_n$ be functions in $\mathcal{R}(I, \mb{C})$ and let $\alpha_1, \ldots, \alpha_n$ be complex numbers. We have
 \[\Vert \alpha_1 f_1+\cdots+\alpha_nf_n\Vert\leq |\alpha_1|\Vert f_1\Vert+\cdots +|\alpha_n|\Vert f_n\Vert.\]
 \end{proposition}
  The proof is also the same as for a real inner product. One consider the case $n=2$ first and then prove the general case by induction on $n$.

 We can define orthogonality on $\mathcal{R}(I, \mb{C})$ the same way as for a real inner product space.
 \begin{definition}{Orthogonality}
 Given two functions $f$ and $g$ in $\mathcal{R}(I, \mb{C})$, we say that they are orthogonal if $\langle f, g\rangle=0$.
 \end{definition}

 \begin{example}[label=230813_1]{}
 Let $I=[0, 2\pi]$. For $n\in \mb{Z}$, define $\phi_n:I\to\mb{C}$ by
 $\di \phi_n(x)=e^{inx}$.
 Show that if $m$ and $n$ are distinct integers, then $\phi_m$ and $\phi_n$ are orthogonal.
 \end{example}
\begin{solution}{Solution}
Notice that
\[
\langle \phi_m, \phi_n\rangle =\int_0^{2\pi }\phi_m(x)\overline{\phi_n(x)}dx
=\int_0^{2\pi} e^{i(m-n)x}dx.\]
\bs
Since $m\neq n$, and
\[\frac{d}{dx}e^{i(m-n)x}=i(m-n)e^{i(m-n)x},\]
fundamental theorem of calculus implies that
\[\langle \phi_m, \phi_n\rangle=\left[\frac{e^{i(m-n)x}}{i(m-n)}\right]_0^{2\pi}=0.\]
 
Hence, $\phi_m$ and $\phi_n$ are orthogonal.\end{solution}

\begin{definition}{Orthogonal System and Orthonormal System}
Let $\mathcal{S}=\left\{\phi_{\alpha}\,|\, \alpha \in J\right\}$ be a subset of functions in $\mathcal{R}(I, \mb{C})$ indexed by the set $J$. We say that $\mathcal{S}$ is an orthogonal system of functions if 
\[\langle \phi_{\alpha}, \phi_{\beta}\rangle=0\hspace{1cm}\text{whenever}\;\alpha\neq \beta,\] and \[\Vert \phi_{\alpha}\Vert\neq 0\hspace{1cm}\text{for all}\;\alpha\in J.\]
We say that $\mathcal{S}$ is an orthonormal  system of functions if it is an orthogonal system and 
\[\Vert \phi_{\alpha}\Vert=1\hspace{1cm}\text{for all}\;\alpha\in J.\]
\end{definition}
Notice that in our definition of orthogonal system, we have an additional condition that each element in the set $\mathcal{S}$ cannot have zero norm.
By definition, it is obvious that if $\mathcal{S}$ is an orthogonal system, then any subset of $\mathcal{S}$ is also an orthogonal system. The same holds for orthonormal systems.

\begin{example}[label=230813_4]{}
Let $I=[0, 2\pi]$. For $n\in \mb{Z}$, define $\phi_n:I\to\mb{C}$ by
 $\di\phi_n(x)=e^{inx}$. Then\[\Vert \phi_n\Vert^2 =\int_0^{2\pi} e^{inx} e^{-inx}dx=2\pi.\]
Example \ref{230813_1} implies that  $\mathcal{S}=\{\phi_n \,|\, n\in\mb{Z}\}$ is an orthogonal system.
\be
If  we let $\varphi_n:I\to\mb{R}$, $n\in\mathbb{Z}$ be the function
\[\varphi_n(x)=\frac{\phi_n(x)}{\Vert\phi_n \Vert}=\frac{e^{inx}}{\sqrt{2\pi}},\]
then $\widetilde{\mathcal{S}}=\left\{\varphi_n\,|\,n\in\mb{Z}\right\}$ is an orthonormal system.
\end{example2}

Using the semi-norm, we can define a relation $\sim$ on $\mathcal{R}(I, \mb{C})$ in the following way. We say that $f\sim g$ if and only if $\Vert f-g\Vert=0$. It is easy to check that this is an equivalence relation. Reflexivity and symmetry are obvious. For transitivity, we note that if  $f\sim g$ and $g\sim h$, then $\Vert f-g\Vert=0$ and $\Vert g-h\Vert=0$. It follows from triangle inequality that
\[\Vert f-h\Vert\leq \Vert f-g\Vert+\Vert g-h\Vert=0.\]This implies that $\Vert f-h\Vert=0$, and thus $f\sim h$. Hence, $\sim$ is an equivalence relation on   $\mathcal{R}(I, \mb{C})$, which we call $L^2$-equivalent.

\begin{definition}{$\pmb{L^2}$ Equivalent Functions}
Two Riemann integrable functions $f:I\to\mb{C}$ and $g:I\to\mb{C}$ are $L^2$-equivalent if
\[\Vert f-g\Vert=0.\]
\end{definition}

\begin{example}{}
Let $I=[a,b]$, and let $S$ be a finite subset of $I$. If $f:I\to\mb{C}$ and $g:I\to\mb{C}$ are two Riemann integrable functions and
\[f(x)=g(x)\hspace{1cm}\text{for all}\;x\in [a,b]\setminus S,\]
then $f$ and $g$ are $L^2$-equivalent.
\end{example}
Regarding  $\mathcal{R}(I, \mb{C})$ as an additive group, the subset $\mathcal{K}(I, \mb{C})$ that contains all the functions in $  \mathcal{R}(I, \mb{C})$ that have zero norm is a normal subgroup. They are functions that are $L^2$-equivalent to the zero function.
Denote by \[\widehat{\mathcal{R}}(I, \mb{C})= \mathcal{R}(I, \mb{C})/\mathcal{K}(I, \mb{C})=\mathcal{R}(I, \mb{C})/\sim \]
  the quotient group. Then each element of $\widehat{\mathcal{R}}(I, \mb{C})$ is an $L^2$ equivalent class of functions. 

If $u$ is in $\mathcal{K}(I, \mb{C})$, $g$ is in $\mathcal{R}(I, \mb{C})$, then
Cauchy-Schwarz inequality implies that
\[\left|\langle u, g\rangle \right|\leq \Vert u\Vert \Vert g\Vert=0.\]
Thus, $\langle u, g\rangle=0$.
If  $f$ is $L^2$ equivalent to $f_1$, $g$ is $L^2$ equivalent to $g_1$, there exists $u$ and $v$ in $\mathcal{K}(I,\mb{C})$ such that $f_1=f+u$ and $g_1=g+v$. 
Therefore,
\[\langle f_1,g_1\rangle=\langle f+u, g+v\rangle =\langle f, g\rangle+\langle u, g\rangle +\langle f, v\rangle+\langle u, v\rangle=\langle f,g\rangle.\]
Hence, the positive semi-definite inner product $\langle \;\cdot\;,\;\cdot\;\rangle$ on $\mathcal{R}(I, \mb{C})$ induces an infinite product on $\widehat{\mathcal{R}}(I, \mb{C})$ by
\[\langle [f], [g]\rangle=\langle f, g\rangle.\]  If $[f]\in \widehat{\mathcal{R}}(I, \mb{C})$ is such that 
\[\langle [f], [f]\rangle =\langle f, f\rangle=0,\] then $f$ is in $\mathcal{K}(I, \mb{C})$, and thus, $[f]=[0]$. This says that the infinite product $\langle \;\cdot\;,\;\cdot\;\rangle$ on $\widehat{\mathcal{R}}(I, \mb{C})$ is positive definite. The additional condition we impose on a subset $\mathcal{S}$ of $\mathcal{R}(I, \mb{C})$ to be an orthogonal system just means that none of the elements in $\mathcal{S}$ is $L^2$-equivalent to the zero function.

For an orthogonal system, we have the following from \eqref{230813_3}.
\begin{theorem}{Generalized Pythagoras Theorem}
Let   $\mathcal{S}=\left\{\phi_{k}\,|\, 1\leq k\leq n\right\}$  be an orthogonal system of functions in $\mathcal{R}(I, \mb{C})$. For any complex numbers $\alpha_1, \ldots, \alpha_n$,
\[\Vert  \alpha_1 \phi_1+\cdots+\alpha_n\phi_n\Vert^2= |\alpha_1|^2\Vert \phi_1\Vert^2+\cdots +|\alpha_n|^2\Vert \phi_n\Vert^2.\]
\end{theorem}

The functions $\phi_n:[0,2\pi]\to\mb{C}$, $f_n(x)=e^{inx}$, $n\in\mb{Z}$ are easy to deal with because of $\di \frac{d}{dx}e^{ax}=ae^{ax}$ for any complex numbers $a$. The drawback is they are complex-valued functions. Since
\[e^{inx}=\cos nx+i\sin nx,\] 
if one wants to work with real-valued functions, one should consider the functions $\cos nx$ and $\sin nx$. 
\begin{proposition}[label=230813_5]{}
Let $I=[0,2\pi]$, and define the functions $C_n:I\to \mb{R}$, $n\geq 0$, and $S_n:I\to \mb{R}$, $n\geq 1$ by
\[C_n(x)=\cos nx, \hspace{1cm} S_n(x)=\sin nx.\]
Then $\mathscr{B}=\left\{C_n\,|\, n\geq 0\right\}\cup \left\{S_n\,|\, n\geq 1\right\}$ is an orthogonal system, and $\Vert C_0\Vert =\sqrt{2\pi}$,
\[\Vert C_n\Vert=\Vert S_n\Vert=\sqrt{\pi}
\hspace{1cm}\text{when}\;n\geq 1.\]
\end{proposition}
\begin{myproof}{Proof}
For $n\in\mb{Z}$, let $\phi_n:[0,2\pi]\to\mb{C}$ be the function $\phi_n(x)=e^{inx}$. Then  $C_0=\phi_0$, and when $n\in\mb{Z}^+$, 
\[C_n=\frac{\phi_n+\phi_{-n}}{2},\hspace{1cm} S_n=\frac{\phi_n-\phi_{-n}}{2i}.\]
 
Since $\{\phi_n\,|\, n\in\mathbb{Z}\}$ is an orthogonal system, we find that for $n\in\mb{Z}^+$,
\[\langle C_0, C_n\rangle =\frac{1}{2}\langle \phi_0, \phi_n\rangle+\frac{1}{2}\langle \phi_0, \phi_{-n}\rangle=0,\]
\[\langle C_0, S_n\rangle =\frac{i}{2}\langle \phi_0, \phi_n\rangle-\frac{i}{2}\langle \phi_0, \phi_{-n}\rangle=0.\]
For $m, n\in \mb{Z}^+$ such that $m\neq n$,
\[\langle C_m, C_n\rangle=\frac{1}{4 }\left( \langle \phi_m, \phi_n\rangle+\langle \phi_m, \phi_{-n}\rangle+\langle \phi_{-m},\phi_n\rangle+\langle \phi_{-m}, \phi_{-n}\rangle\right)=0,\]
\[\langle S_m, S_n\rangle=\frac{1}{4}\left( \langle \phi_m, \phi_n\rangle-\langle \phi_m, \phi_{-n}\rangle-\langle \phi_{-m}, \phi_n\rangle+\langle \phi_{-m},\phi_{-n}\rangle\right)=0.\]

For $m, n\in \mb{Z}^+$, considering the cases $m=n$ and $m\neq n$ separately, we find that
\[\langle S_m, C_n\rangle=\frac{1}{4i}\left( \langle \phi_m, \phi_n\rangle+\langle \phi_m, \phi_{-n}\rangle-\langle \phi_{-m}, \phi_n\rangle-\langle \phi_{-m}, \phi_{-n}\rangle\right)=0.\]
\bp
These show that $\mathscr{B}$ is an orthogonal system.
For $n\in \mb{Z}^+$, since $\Vert \phi_n\Vert=\Vert \phi_{-n}\Vert=\sqrt{2\pi}$, and $\phi_n$ and $\phi_{-n}$ are orthogonal, we have
\[\Vert C_n\Vert^2=\left\Vert\frac{1}{2}\phi_n+\frac{1}{2}\phi_{-n}\right\Vert^2= \frac{1}{4 }\Vert\phi_n\Vert^2+\frac{1}{4}\Vert\phi_{-n}\Vert^2=\pi,\]
\[\Vert S_n\Vert^2= \left\Vert\frac{1}{2i}\phi_n-\frac{1}{2i}\phi_{-n}\right\Vert^2= \frac{1}{4 }\Vert\phi_n\Vert^2+\frac{1}{4}\Vert\phi_{-n}\Vert^2=\pi.\]
These complete the proof.
\end{myproof}

Given a finite subset $\mathcal{S}=\{f_1, \ldots, f_n\}$  of  $\mathcal{R}(I, \mb{C})$,   let
\[W_{\mathcal{S}}=\text{span}\, \mathcal{S}=\left\{c_1 f_1+\cdots+c_nf_n\,|\, c_1, \ldots, c_n\in\mb{C}\right\}\]
be the subspace of $\mathcal{R}(I, \mb{C})$ spanned by $\mathcal{S}$. 
We say that an element $g$ of $\mathcal{R}(I, \mb{C})$ is orthogonal to $W_{\mathcal{S}}$ if it is orthogonal to each $f\in W_{\mathcal{S}}$. This is if and only if $g$ is orthogonal to $f_k$ for all $1\leq k\leq n$.
The projection theorem says the following.

\begin{theorem}{Projection Theorem}
Let $\mathcal{S}=\{\phi_1, \ldots, \phi_n\}$ be an orthogonal system of functions in $\mathcal{R}(I, \mb{C})$, and let $W_{\mathcal{S}}$ be the subspace of $\mathcal{R}(I, \mb{C})$ spanned by $\mathcal{S}$. Given $f$ in $\mathcal{R}(I, \mb{C})$, there is a unique $g\in W_{\mathcal{S}}$ such that $f-g$ is orthogonal to $W_{\mathcal{S}}$. It is called the projection of the function $f$ onto the subspace $W_{\mathcal{S}}$, denoted by $\text{proj}_{W_{\mathcal{S}}}f$, and it is given by
\[\text{proj}_{W_{\mathcal{S}}}f=\sum_{k=1}^n\frac{\langle f, \phi_k\rangle }{\langle \phi_k, \phi_k\rangle}\phi_k=\frac{\langle f, \phi_1\rangle }{\langle \phi_1, \phi_1\rangle}\phi_1+\cdots+
\frac{\langle f, \phi_n\rangle }{\langle \phi_n, \phi_n\rangle}\phi_n.\]For any $h\in W_{\mathcal{S}}$, 
\[ \Vert f-h\Vert\geq \Vert f-\text{proj}_{W_{\mathcal{S}}}f\Vert.\]
 \end{theorem}
\begin{myproof}{Proof}
Assume that $g$ is a function in $W_{\mathcal{S}}$ such that
$f-g$ is orthogonal to $W_{\mathcal{S}}$. Then there exist complex numbers $\alpha_1, \ldots, \alpha_n$ such that
\[g=\alpha_1\phi_1+\cdots+\alpha_n\phi_n.\]
\bp
Since $\langle \phi_k, \phi_l\rangle =0$ if $k\neq l$, we find that 
\[\langle g,\phi_k\rangle =\alpha_k\langle \phi_k, \phi_k\rangle\hspace{1cm}\text{for}\;1\leq k\leq n.\]
Since $f-g$ is orthogonal to $W_{\mathcal{S}}$, $\langle f-g, \phi_k\rangle =0$ for all $1\leq k\leq n$. This gives $\langle f,\phi_k\rangle=\langle g,\phi_k\rangle$, and thus
\[\alpha_k\langle \phi_k, \phi_k\rangle =\langle f, \phi_k\rangle\hspace{1cm}\text{for}\;1\leq k\leq n.\]
Hence, we must have \[\alpha_k=\frac{\langle f, \phi_k\rangle}{\langle \phi_k, \phi_k\rangle}.\]
This implies the uniqueness of $g$ if it exists.
It is easy to check that the function

\[g=\sum_{k=1}^n\frac{\langle f, \phi_k\rangle }{\langle \phi_k, \phi_k\rangle}\phi_k\]
 
 is indeed a function in $W_{\mathcal{S}}$ such that $f-g$ is orthogonal to $W_{\mathcal{S}}$.
 
Finally, for any $h$ in $W_{\mathcal{S}}$, $g-h$ is also in $W_{\mathcal{S}}$. Hence, $g-h$ is orthogonal to $f-g$. By the generalized Pythogoras theorem,
\[\Vert f-h\Vert^2=\Vert(f-g)+(g-h)\Vert^2=\Vert f-g\Vert^2+\Vert g-h\Vert^2\geq \Vert f-g\Vert^2.\]
This proves that \[\Vert f-h\Vert\geq \Vert f-g\Vert\hspace{1cm}\text{for all}\;h\in W_{\mathcal{S}}.\]
\end{myproof}

Now we restrict our consideration to functions $f$ that are periodic of period $2\pi$. In this case, the function is uniquely determined by its values on an interval $[a, b]$ of length $2\pi$. We often take $I=[0,2\pi]$ or $I=[-\pi, \pi]$. Notice that if $f:\mb{R}\to\mb{C}$ is a function of period $2\pi$, then for any $\alpha\in\mb{R}$,
\[\int_{\alpha}^{\alpha+2\pi} f(x)dx=\int_0^{2\pi}f(x)dx=\int_{-\pi}^{\pi}f(x)dx.\]
Any function  $f:[\alpha, \alpha+2\pi]\to \mb{C}$ defined on an interval of length $2\pi$ can be extended to be a $2\pi$-periodic function. 
\begin{definition}{Extension of Functions}
Let $I=[\alpha, \alpha+2\pi]$ be an inverval of length $2\pi$, and let $f:I\to\mb{C}$ be a function defined on $I$. We can extend $f$ to be a $2\pi$-periodic function $\widetilde{f}:\mb{R}\to\mb{C}$ in the following way.
\begin{enumerate}[(i)]
\item For $x\in (\alpha, \alpha+2\pi)$, define
\[\widetilde{f}(x+2n\pi)=f(x)
\hspace{1cm}\text{for all}\;n\in \mathbb{Z}.\]
\item For $x=\alpha$,  define
\[\widetilde{f}(\alpha+2n\pi)=\di\frac{f(\alpha)+f(\alpha+2\pi)}{2}\hspace{1cm}\text{for all}\;n\in \mathbb{Z}.\]
\end{enumerate}
\end{definition}
Examples are shown in Figures \ref{figure68} and \ref{figure70}.

\begin{figure}[ht]
\centering
\includegraphics[scale=0.2]{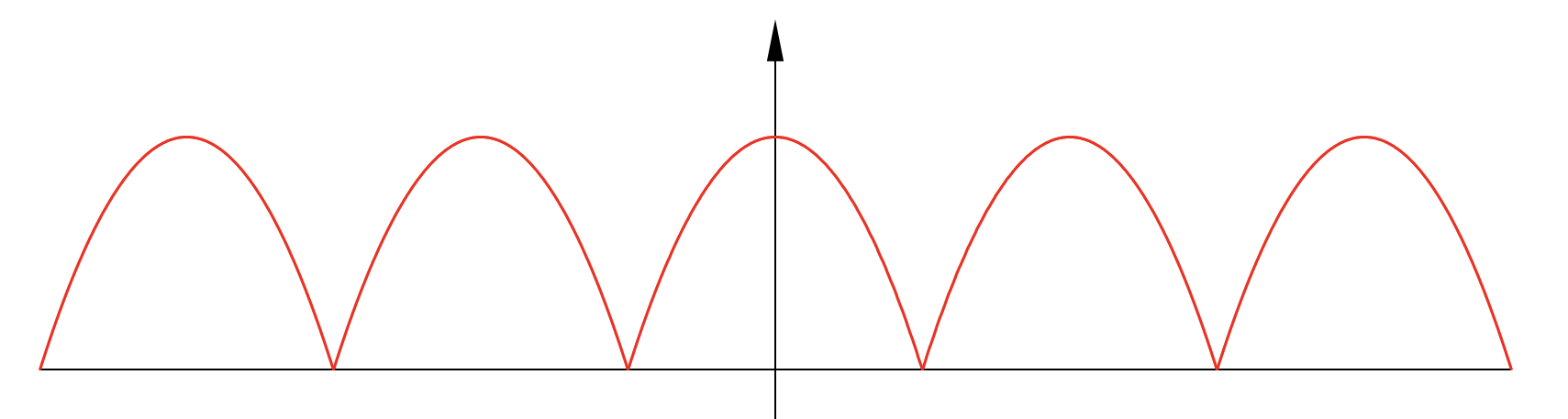}

\caption{Extending a function  $f:[-\pi,\pi]\to\mb{R}$ periodically.}\label{figure68}
\end{figure}
\begin{figure}[ht]
\centering
\includegraphics[scale=0.2]{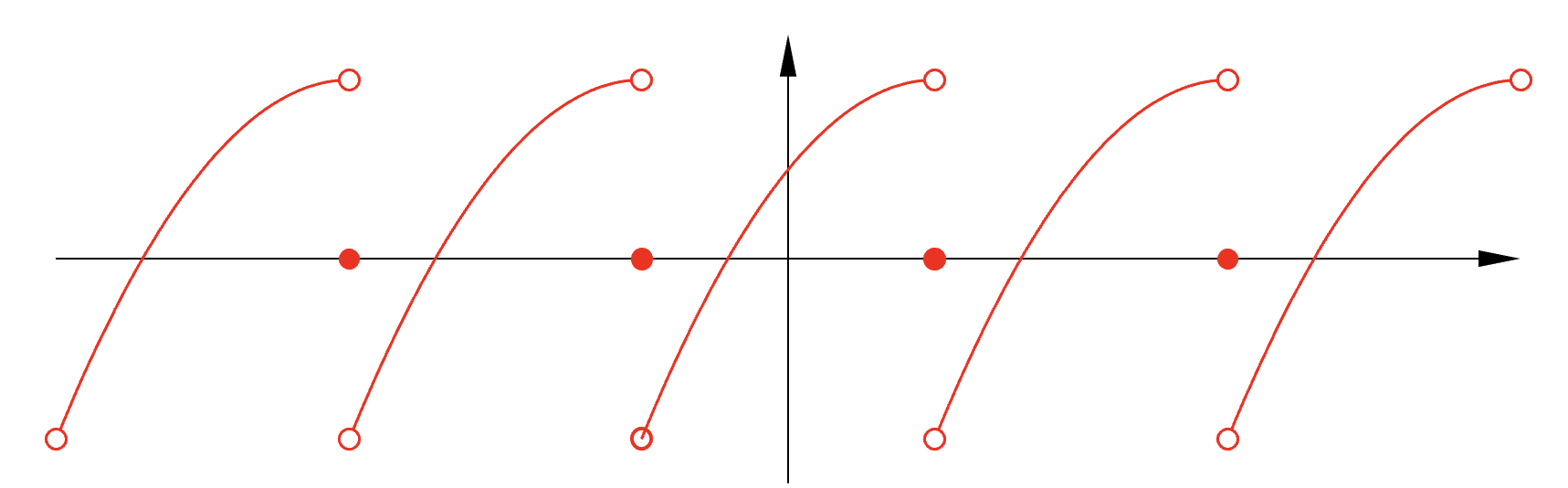}

\caption{Extending a function $f:[-\pi,\pi]\to\mb{R}$ periodically.}\label{figure70}
\end{figure}

Now let us define Fourier series.
Example \ref{230813_4} asserts that the set \[\mathcal{S}=\di\left\{\phi_n\,|\, n\in\mathbb{Z}\right\},\quad  \text{where}\; \phi_n(x)=e^{inx},\] is an orthogonal system of functions in $\mathcal{R}(I, \mb{C})$, where $I=[-\pi, \pi]$. For $n\geq 0$, let $W_n$ be the subspace of $\mathcal{R}(I, \mb{C})$ spanned by $\mathcal{S}_n=\left\{e^{ikx}\,|\,-n\leq k\leq n\right\}$. It is a vector space of dimension $2n+1$  with basis $\mathcal{S}_n$. Moreover, 
\[W_0\subset W_1\subset W_2\subset \cdots.\] A real basis of $W_n$ is given by
\[\mathscr{B}_n=\left\{\sin kx\,|\, k=1, \ldots, n\right\}\cup \left\{\cos kx \,|\, k=0, 1, \ldots, n\right\}.\]
Given $f\in\mathcal{R}(I, \mb{C})$, let $s_n=\text{proj}_{W_n} f$ be the projection of $f$ onto $W_n$.  The projection theorem says that
\[s_n(x)=\left(\text{proj}_{W_n} f\right)(x)= \sum_{k=-n}^nc_ke^{ikx},\]
where
\[c_k=\frac{\langle f, \phi_k\rangle}{\langle\phi_k,\phi_k\rangle}=\frac{1}{2\pi} \int_{-\pi}^{\pi}f(x)e^{-ikx}dx.\]
By Proposition \ref{230813_5} and the projection theorem, $s_n(x)$ can also be written as
\[s_n(x)=\left(\text{proj}_{W_n} f\right)(x)= \frac{a_0}{2}+\sum_{k=1}^n\left(a_k\cos kx+b_k\sin kx\right),\]
where  
\[a_k=\frac{1}{\pi}\int_{-\pi}^{\pi}f(x)\cos kx dx,\hspace{1cm} \text{for}\;0\leq k\leq n,\]\[b_k=\frac{1}{\pi}\int_{-\pi}^{\pi}f(x)\sin kx dx\hspace{1cm}\text{for}\;1\leq k\leq n.\]
 By definition, we find that \[c_0=\frac{a_0}{2}, \]and when $k\geq 1$,
 \[  c_{-k}=\frac{a_k+ib_k}{2}, \quad c_k=\frac{a_k-ib_k}{2}.\]
 If $f$ is a real-valued function, $a_k$ and $b_k$ are real and $c_{-k}=\overline{c_{k}}$.
  \begin{definition}{Trigonometric Series}
A trigonometric series is a series of the form 
 \[\frac{a_0}{2}+\sum_{k=1}^{\infty}(a_k\cos kx+b_k\sin kx).\]
 \end{definition}Since a trigonometric series can   be expressed in the form $\di \sum_{k=-\infty}^{\infty}c_ke^{ikx}$, we also call a series of the form $\di \sum_{k=-\infty}^{\infty}c_ke^{ikx}$ a trigonometric series. Fourier series of a function is a trigonometric series.
 \begin{definition}{Fourier Series and its $\pmb{n^{\text{th}}}$ Partial Sums}
 Let $I=[-\pi, \pi]$. The Fourier series of a function $f$ in $\mathcal{R}(I, \mb{C})$ is the infinite series
 \[\sum_{k=-\infty}^{\infty} c_ke^{ikx}\quad \text{or}\quad \frac{a_0}{2}+\sum_{k=1}^{\infty}\left(a_k\cos kx+b_k\sin kx\right),\]
 where
 \begin{gather*}
c_k=\frac{1}{2\pi}\int_{-\pi}^{\pi}  f(x)e^{-ikx}dx,\hspace{1cm}k\in \mathbb{Z},\\
a_k=\frac{1}{\pi}\int_{-\pi}^{\pi} f(x)\cos kx dx,\hspace{1cm}k\geq 0,\\
b_k=\frac{1}{\pi}\int_{-\pi}^{\pi} f(x)\sin kxdx, \hspace{1cm}k\geq 1.
\end{gather*}
The $n^{\text{th}}$-partial sum of the Fourier series is
\[s_n(x)=\sum_{k=-n}^{n} c_ke^{ikx}=\frac{a_0}{2}+\sum_{k=1}^{n}\left(a_k\cos kx+b_k\sin kx\right).\]
It is the projection of $f$ onto the subspace of $\mathcal{R}(I,\mb{C})$ spanned by $\mathcal{S}_n=\left\{e^{ikx}\,|\,-n\leq k\leq n\right\}$. 
 \end{definition}
 \begin{remark}{}
 If $f:\mb{R}\to\mb{C}$ is a $2\pi$-periodic function which is Riemann integrable over a closed interval of length $2\pi$, the Fourier series of $f$ is the Fourier series of $f:[-\pi,\pi]\to\mb{C}$. 
 \end{remark}

 \begin{remark}{}
 If $I=[-L, L]$, the Fourier series of a function $f\in \mathcal{R}(I, \mb{C})$ is the series
 \[\sum_{k=-\infty}^{\infty} c_k\exp\left(\frac{i\pi kx}{L}\right),\]
 where\[c_k=\frac{1}{2L}\int_{-L}^L f(x)\exp\left(-\frac{i\pi kx}{L}\right)dx.\]
 \end{remark}
 Henceforth, we only consider the case where $I$ is a closed interval of length $2\pi$.
Let us look at some examples.
 
 \begin{example}[label=230815_5]{}
 Find the Fourier series of the function $f:[-\pi, \pi]\to\mb{R}$ defined as \[f(x)=x.\]
 
 \end{example}
 
 \begin{solution}{Solution}
  When $k=0$,
 \begin{align*}
 c_0=\frac{1}{2\pi}\int_{-\pi}^{\pi}f(x)dx=\frac{1}{2\pi}\int_{-\pi}^{\pi}xdx=0.
 \end{align*}
When $k\neq 0$,
 \begin{align*}
 c_k&=\frac{1}{2\pi}\int_{-\pi}^{\pi}xe^{-ikx}dx\\
 &=\frac{1}{2\pi}\left\{\left[-\frac{1}{ik}xe^{-ikx}\right]_{-\pi}^{\pi}+\frac{1}{ik}\int_{-\pi}^{\pi}e^{-ikx}dx\right\}\\
 &=\frac{(-1)^{k-1}}{ik}.\end{align*} 
 
 Therefore, the Fourier series of $f$ is
 \[ \sum_{k=1}^{\infty}(-1)^{k-1}\frac{e^{ikx}-e^{-ikx}}{ik}= \sum_{k=1}^{\infty}2(-1)^{k-1} \frac{\sin kx}{k}.\]
 \end{solution}
 
\begin{figure}[ht]
\centering
\includegraphics[scale=0.2]{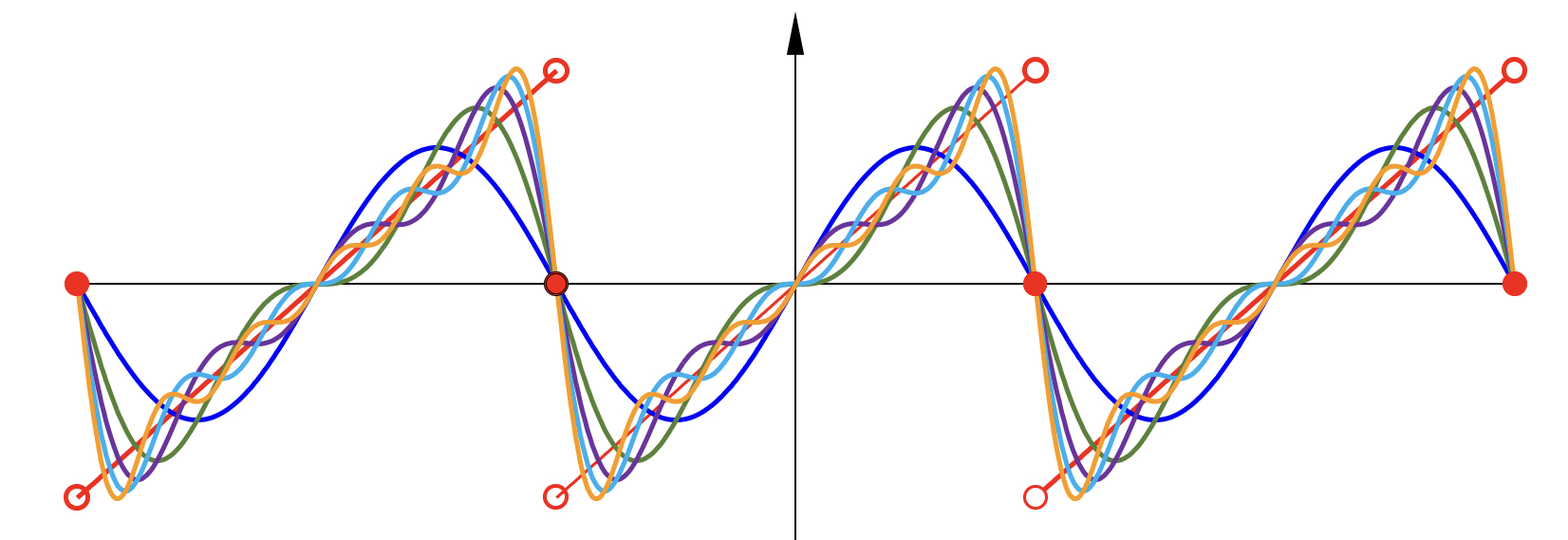}

\caption{The function $f(x)=x, -\pi < x<\pi$ and $s_n(x)$, $1\leq n\leq 5$.}\label{figure69}
\end{figure}

\begin{remark}{}
Let $I=[-\pi,\pi]$. Given $f$ in $\mathcal{R}(I, \mb{C})$, we call  each \[  c_k=\int_If(x)e^{-ikx}dx,\hspace{1cm} k\in\mathbb{Z}\] a Fourier coefficient  of $f$. The mapping $\mathfrak{F}_k$ from  $\mathcal{R}(I, \mb{C})$ to $\mb{C}$ which takes a function $f$ to $c_k$ is a linear transformation between vector spaces. 

When $f:I\to \mb{R}$ is a real-valued function, we usually prefer to work with the Fourier coefficients $a_k$ and $b_k$. One can show that if $f:[-\pi,\pi]\to\mb{R}$ is an odd function, then $a_k=0$ for all $k\geq 0$, so that the Fourier series of $f$ only has  sine terms. If $f:[-\pi,\pi]\to\mb{R}$ is an even function, then $b_k=0$ for all $k\geq 1$, so that the Fourier series of $f$ only has the constant and the cosine terms.
\end{remark}
 
  \begin{remark}{}
 If $f:I\to\mb{C}$ is a function of the form
 \[f(x)=\sum_{k \in J}c_ke^{ikx},\]
 where $J$ is a finite subset of integers, then the Fourier series of $f$ is equal to itself. 
 
 \end{remark}
  \begin{example}[label=230815_6]{}
 The Fourier series of a constant function $f:I\to \mb{C}$, $f(x)=c$ is just $c$ itself.
 \end{example}
 
 \begin{example}{}
 Let $f:[0,2\pi]\to \mb{R}$ be the function defined as
 \[f(x)=x(2\pi -x).\]
 Find its Fourier series,  and express it in terms of trigonometric functions.
 \end{example}
 \begin{figure}[ht]
\centering
\includegraphics[scale=0.2]{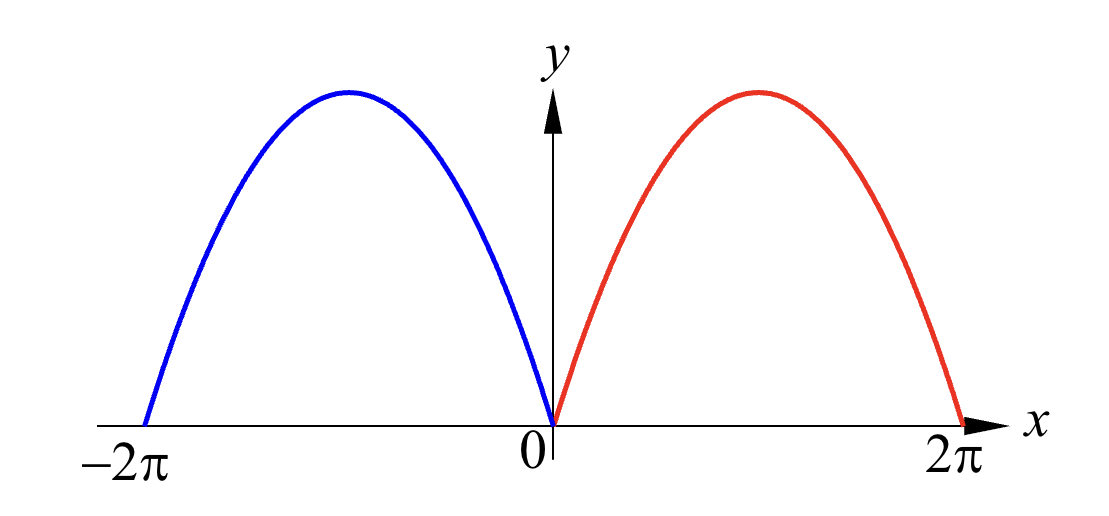}

\caption{The function $f:[0,2\pi]\to\mb{R}$, $f(x)=x(2\pi-x)$ and its entension.}\label{figure71}
\end{figure}
 \begin{solution}{Solution}
 When we extend $f$ periodically to the function $\widetilde{f}:\mb{R}\to\mb{R}$, we find that when $x\in [0, 2\pi]$, 
 \[\widetilde{f}(-x)=f(-x+2\pi)=(2\pi -x)x=f(x)=\widetilde{f}(x).\]
 Hence, $\widetilde{f}(x)$ is an even function. This implies that the Fourier series of $f:[0,2\pi]\to\mb{R}$ only has cosine terms.
 
 \[a_0=\frac{1}{\pi}\int_0^{2\pi}f(x)dx=\frac{2}{\pi}
 \int_0^{\pi}(2\pi x-x^2)dx=\frac{2}{\pi}\left[\pi x^2-\frac{x^3}{3}\right]_0^{\pi}=\frac{4}{3}\pi^2.\]
 \bs
 For $k\geq 1$,
 \begin{align*}
 a_k&=\frac{2}{\pi}\int_0^{\pi}f(x)\cos xdx=\frac{2}{\pi}
 \int_0^{\pi}(2\pi x-x^2)\cos x dx\\
 &=\frac{2}{\pi}\left(\left[\frac{(2\pi x -x^2)\sin kx}{k}\right]_0^{\pi}-\frac{1}{k}\int_0^{\pi} (2\pi -2x)\sin kx dx\right)\\
 &=-\frac{2}{\pi k}\left(\left[-\frac{(2\pi -2x)\cos kx}{k}\right]_0^{\pi}-\frac{2}{k}\int_0^{\pi}\cos kxdx\right)\\
 &=-\frac{4}{k^2}+\frac{4}{\pi k^3}\left[\sin kx\right]_0^{\pi} =-\frac{4}{k^2}.
 \end{align*}
 Therefore, the Fourier series of  $f$ is
 \[ \frac{2}{3}\pi^2-4\sum_{k=1}^{\infty}\frac{\cos kx}{k^2}.\]
 \end{solution}

  \begin{example}[label=230814_1]{}
 Let $a$ and $b$ be two numbers satisfying $-\pi\leq a<b\leq \pi$, and let $g:[-\pi, \pi]\to\mb{R}$ be the function defined as
 \[g(x)=\begin{cases} 1,\quad &\text{if}\; a\leq x\leq b,\\
 0,\quad &\text{otherwise}.\end{cases}\]
 Find the Fourier series of $g$ in exponential form.
 \end{example}
 
 \begin{figure}[ht]
\centering
\includegraphics[scale=0.2]{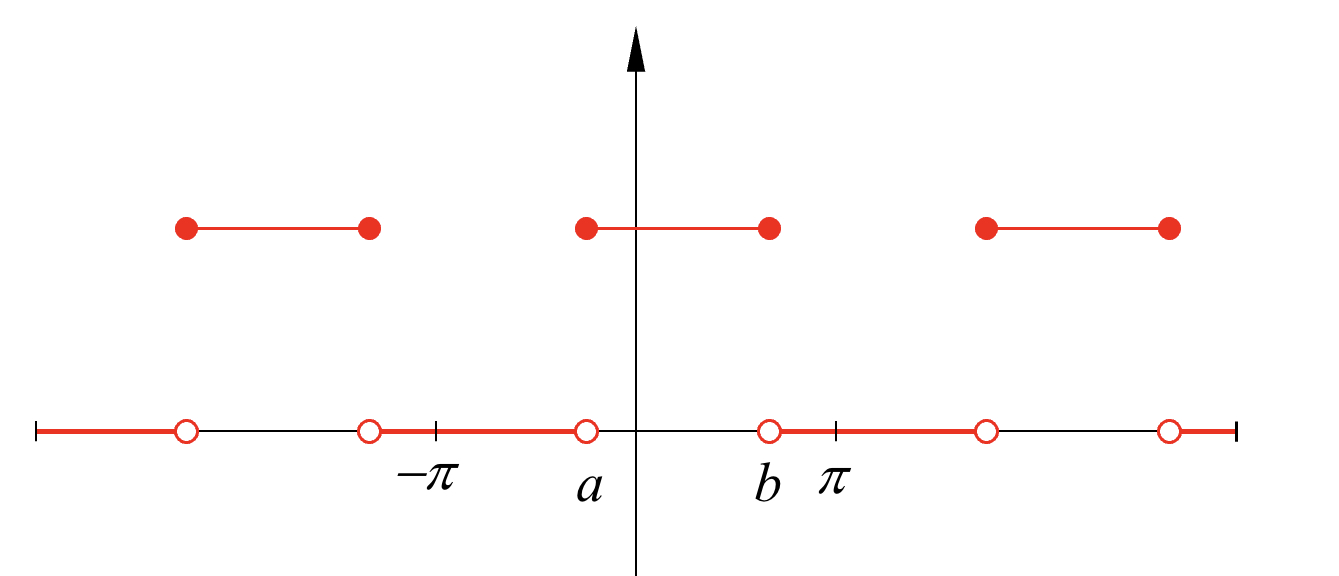}

\caption{The function $g:[0,2\pi]\to\mb{R}$ defined in Example \ref{230814_1} and its entension.}\label{figure72}
\end{figure}
 \begin{solution}{Solution}
Since $g$ is piecewise continuous, it is Riemann integrable. 
 \[c_0=\frac{1}{2\pi}\int_{-\pi}^{\pi}g(x)dx=\frac{1}{2\pi}\int_a^b dx=\frac{b-a}{2\pi}.\]
   For $ k\geq 1$,
   \[
  c_k
   =\frac{1}{2\pi}\int_{-\pi}^{\pi}g(x)e^{-ikx}dx =\frac{1}{2\pi}\int_a^b e^{-ikx}dx =\frac{e^{-ikb}-e^{-ika}}{-2  \pi i k},
   \]
   \[c_{-k}=\overline{c_k}=\frac{e^{ikb}-e^{ika}}{2  \pi i k}.\]
   Therefore, the Fourier series of $g$ is
 \[\frac{b-a}{2\pi}+\frac{i}{2\pi }\sum_{k=1}^{\infty} \frac{(e^{-ikb}-e^{-ika})e^{ikx}-(e^{ikb}-e^{ika})e^{-ikx}}{  k}.
 \]
 
 \end{solution}
  \begin{example}{}
Find the Fourier series of the function $f:[-\pi, \pi]\to\mb{R}$, $f(x)=x\sin x$. 

\end{example}

\begin{solution}{Solution}
Since $f$ is a real-valued even function, we only need to compute the Fourier coefficients
\[a_k(f)=\frac{1}{\pi}\int_{-\pi}^{\pi} x\sin x\cos kxdx\hspace{1cm}\text{when}\;k\geq 0.\]

Let $g:[-\pi, \pi]\to\mb{R}$  be the function $g(x)=x$. We have seen in Example \ref{230815_5} that the Fourier series of $g$ is given by
\[G(x)=  \sum_{k=1}^{\infty}2(-1)^{k-1} \frac{\sin kx}{k}.\]
 
Therefore, for $k\geq 1$,
\[b_k(g)=\frac{1}{\pi}\int_{-\pi}^{\pi} x\sin kx=\frac{2(-1)^{k-1}}{k}.\]
 \bs
From this, we find that 
\begin{gather*}
a_0(f)=b_1(g)=2,\\
a_1(f)=\frac{1}{2\pi}\int_{-\pi}^{\pi}x\sin 2xdx=\frac{1}{2}b_2(g)=-\frac{1}{2};
\end{gather*}
and when $k\geq 2$,
\begin{align*}
a_k(f)&=\frac{1}{2\pi}\int_{-\pi}^{\pi}x\left(\sin(k+1)x-\sin (k-1)x\right)dx\\
&=\frac{1}{2}\left(b_{k+1}(g)-b_{k-1}(g)\right)\\
&=(-1)^k\left(\frac{1}{k+1}-\frac{1}{k-1}\right)\\&=\frac{2(-1)^{k-1}}{k^2-1}.
\end{align*}
Hence, the Fourier series of the function $f:[-\pi,\pi]\to\mb{R}$, $f(x)=x\sin x$ is
\[1-\frac{1}{2}\cos x+\sum_{k=2}^{\infty}\frac{2(-1)^{k-1}}{k^2-1}\cos kx.\]
\end{solution}
 \begin{figure}[ht]
\centering
\includegraphics[scale=0.2]{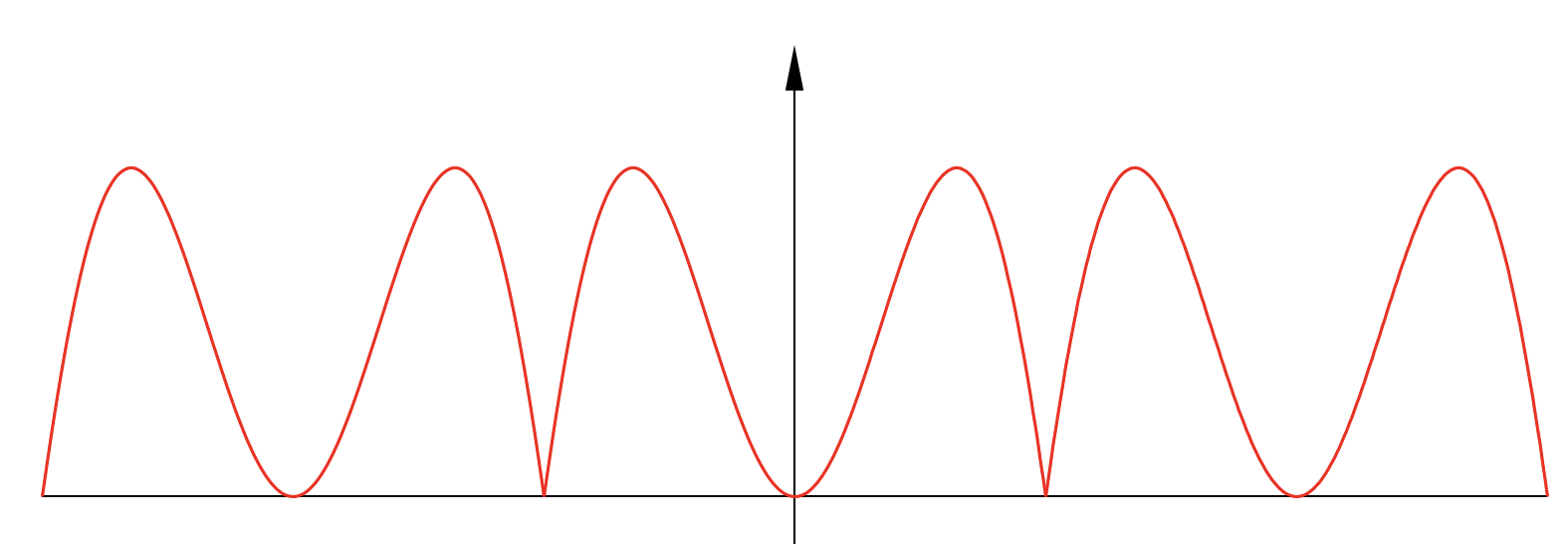}

\caption{The function $f:[-\pi,\pi]\to\mb{R}$, $f(x)=x\sin x$ and its periodic entension.}\label{figure76}
\end{figure}
 At the end of this section, let us make an additional remark.

 \begin{remark}{Semi-Norms}
A semi-norm on a complex vector space $V$ is a function $\Vert \;\cdot\;\Vert :V\to\mb{R}$ which defines the norm $\Vert\mf{v}\Vert$  for each $\mf{v}$ in $V$ such that the following hold.
\begin{enumerate}[(a)]
\item For any $\mf{v}\in V$, $\Vert\mf{v}\Vert \geq 0$.
\item For any $\alpha\in\mathbb{C}$, and any $\mf{v}\in V$, $\Vert \alpha \mf{v}\Vert=|\alpha|\Vert\mf{v}\Vert$.
\item For any $\mf{u}$ and $\mf{v}$ in $V$, $\Vert\mf{u}+\mf{v}\Vert\leq\Vert\mf{u}\Vert+\Vert\mf{v}\Vert$.

 \end{enumerate}
If in addition, we have
\begin{enumerate}[(a)]
\item[(d)] $\Vert\mf{v}\Vert = 0$ if and only if $\mf{v}=0$,
\end{enumerate}then $\Vert \; \cdot\;\Vert $ is called a norm on the vector space $V$.

Proposition \ref{230814_5} and Proposition \ref{230814_6} justify that the $L^2$-norm
\[\Vert f\Vert_2 =\sqrt{\int_I |f(x)|^2 dx}\]
is indeed a semi-norm on the vector space $\mathcal{R}(I,\mb{C})$. 

There are other semi-norms on $\mathcal{R}(I,\mb{C})$. One of them which will also be useful later is the $L^1$-norm defined as
\[\Vert f\Vert_1=\int_I|f(x)|dx.\]
The fact that this is a semi-norm is quite easy to establish.
\end{remark}

\vp
\noindent
{\bf \large Exercises  \thesection}
\setcounter{myquestion}{1}

\begin{question}{\themyquestion}
Let $f:[-\pi, \pi]\to\mb{R}$ be a real-valued Riemann integrable function.
\begin{enumerate}[(a)]
\item
 If $f:[-\pi,\pi]\to\mb{R}$ is an odd function, show that   the Fourier series of $f$ has the form
 \[\sum_{k=1}^{\infty}b_k\sin kx,\hspace{1cm}\text{where}\; b_k=\frac{2}{\pi}\int_0^{\pi}f(x)\sin kx.\]
 \item If $f:[-\pi,\pi]\to\mb{R}$ is an even function, show that   the Fourier series of $f$ has the form
 \[\frac{a_0}{2}+\sum_{k=1}^{\infty}a_k\cos kx,\hspace{1cm}\text{where}\; a_k=\frac{2}{\pi}\int_0^{\pi}f(x)\cos kx.\]
 \end{enumerate}
\end{question}
\atc
 
 \begin{question}{\themyquestion}
  Find the Fourier series of the function $f:[-\pi, \pi]\to\mb{R}$, $f(x)=|x|$,  and express it in terms of trigonometric functions.
\end{question}
 \atc
 
 \begin{question}{\themyquestion}
  Find the Fourier series of the function $f:[-\pi, \pi]\to\mb{R}$, $f(x)=x^2$,  and express it in terms of trigonometric functions.
\end{question}

 \atc
 
 \begin{question}{\themyquestion}
  Find the Fourier series of the function $f:[0, 2\pi]\to\mb{R}$, $f(x)=x^2$,  and express it in terms of trigonometric functions.
\end{question}

 \atc
 
 \begin{question}{\themyquestion}
  Find the Fourier series of the function $f:[-\pi , \pi]\to\mb{R}$, $f(x)=\sin 2x$,  and express it in terms of trigonometric functions.
\end{question}
 \atc
 
 \begin{question}{\themyquestion}
  Find the Fourier series of the function $f:[-\pi , \pi]\to\mb{R}$, 
  \[f(x)=\begin{cases} 0,\quad &\text{if}\; -\pi\leq x<0,\\
  \sin x,\quad &\text{if}\quad 0\leq x\leq\pi,\end{cases}\] and express it in terms of trigonometric functions.
\end{question}
\atc

\begin{question}{\themyquestion}
Find the Fourier series of the function $f:[-\pi,\pi]\to\mb{R}$, $f(x)=x\cos x$ from the Fourier series of the function $g:[-\pi,\pi]\to\mb{R}$, $g(x)=x$.
\end{question}

 \atc

\begin{question}{\themyquestion}
Let $x_0$ be a point in the interval $[a,b]$, and let $f:[a,b]\to\mb{C}$ and $g:[a,b]\to\mb{C}$ be $L^2$-equivalent Riemann integrable functions. Assume that both $f$ and $g$ are continuous at the point $x_0$, show that $f(x_0)=g(x_0)$.
\end{question}
 
\section{The Pointwise Convergence of a Fourier Series} \label{sec7.2}

Let $I=[-\pi, \pi]$ and let $\mathcal{R}(I, \mb{C})$ be the vector space that consists of all Riemann integrable functions $f:I\to\mb{C}$ that are defined on $I$. Each of these functions can be extended to a periodic function $\widetilde{f}:\mb{R}\to\mb{C}$ so that $\widetilde{f}(x)=f(x)$ for all $x$ in the interior of $I$.

Given $f\in \mathcal{R}(I,\mb{C})$, we define the Fourier series of $f$ as the infinite series
\[\sum_{k=-\infty}^{\infty}c_ke^{ikx}=\frac{a_0}{2}+\sum_{k=1}^{\infty}\left(a_k\cos kx+b_k\sin kx\right),\]
where
\begin{gather*}
c_k=\frac{1}{2\pi}\int_{I} f(x)e^{-ikx}dx,\hspace{1cm}k\in \mathbb{Z},\\
a_k=\frac{1}{\pi}\int_If(x)\cos kx dx,\hspace{1cm}k\geq 0,\\
b_k=\frac{1}{\pi}\int_If(x)\sin kxdx, \hspace{1cm}k\geq 1.
\end{gather*}

The problem of interest to us is the convergence of the Fourier series. 
Given $n\geq 0$, the $n^{\text{th}}$-partial sum of the Fourier series of $f$ is
 \[s_n(x)=\sum_{k=-n}^{n} c_ke^{ikx}.\]
 We say that the Fourier series 
 converges pointwise if the sequence of partial sum functions $\{s_n:I\to\mb{C}\}$ converges pointwise. Let us first give an integral expression for the partial sums $s_n(x)$ of the Fourier series. 
By definition,
\begin{align*}
s_n(x)&=\frac{1}{2\pi} \sum_{k=-n}^{n}\left( \int_{-\pi}^{\pi}
f(t)e^{-ikt}dt\right)e^{ikx}\\
&=\frac{1}{2\pi}\int_{-\pi}^{\pi}f(t)\sum_{k=-n}^ne^{ik(x-t)}dt.
\end{align*} 
If $x\in 2\pi\mb{Z}$, $e^{ikx}=1$ for all $-n\leq k\leq n$. Therefore,  
\[\sum_{k=-n}^ne^{ikx} =2n+1.\] If $x\notin 2\pi \mb{Z}$,  $e^{ix}\neq 1$. Using the sum formula for a geometric sequence, we have
\begin{align*}
\sum_{k=-n}^ne^{ikx}&=e^{-inx}\left(1+e^{ix}+\cdots+e^{2inx}\right) 
 =\frac{e^{-i\left(n+\frac{1}{2}\right)x}}{e^{-\frac{ix}{2}}}
\times\frac{ e^{i(2n+1)x}-1}{ e^{ix}-1}\\
&=\frac{e^{i(n+\frac{1}{2})x}-e^{-i(n+\frac{1}{2})x}}{e^{\frac{ix}{2}}-e^{-\frac{ix}{2}}} 
 =\frac{\di \sin \left(n+\frac{1}{2}\right)x}{\di \sin\frac{x}{2}}.
\end{align*} 

\begin{definition}{Dirichlet Kernel}
Given a nonnegative integer $n$, the Dirichlet kernel   $D_n:\mb{R}\to\mb{R}$ is 
\[D_n(x)=\sum_{k=-n}^ne^{ikx} =\begin{cases} 2n+1,\quad &\text{if}\; x\in 2\pi\mb{Z},
\\\di
\frac{\di \sin \left(n+\frac{1}{2}\right)x}{\di \sin\frac{x}{2}},\quad &\text{otherwise}.\end{cases}\]
\end{definition}
Our derivation above gives the following.
\begin{proposition}[label=230815_2]{}
Let $I=[-\pi,\pi]$, and let  $f:I\to\mb{C}$ be a Riemann integrable function. The $n^{\text{th}}$-partial sum $s_n(x)$ of the Fourier series of $f$ has an integral representation given by
\[s_n(x)=\frac{1}{2\pi}\int_{-\pi}^{\pi}f(t)D_n(x-t)dt,\]
where $D_n:\mb{R}\to\mb{R}$ is the Dirichlet kernel given by
\[D_n(x) =\begin{cases} 2n+1,\quad &\text{if}\; x\in 2\pi\mb{Z},
\\\di
\frac{\di \sin \left(n+\frac{1}{2}\right)x}{\di \sin\frac{x}{2}},\quad &\text{otherwise}.\end{cases}\]
\end{proposition}
 
  \begin{figure}[ht]
\centering
\includegraphics[scale=0.2]{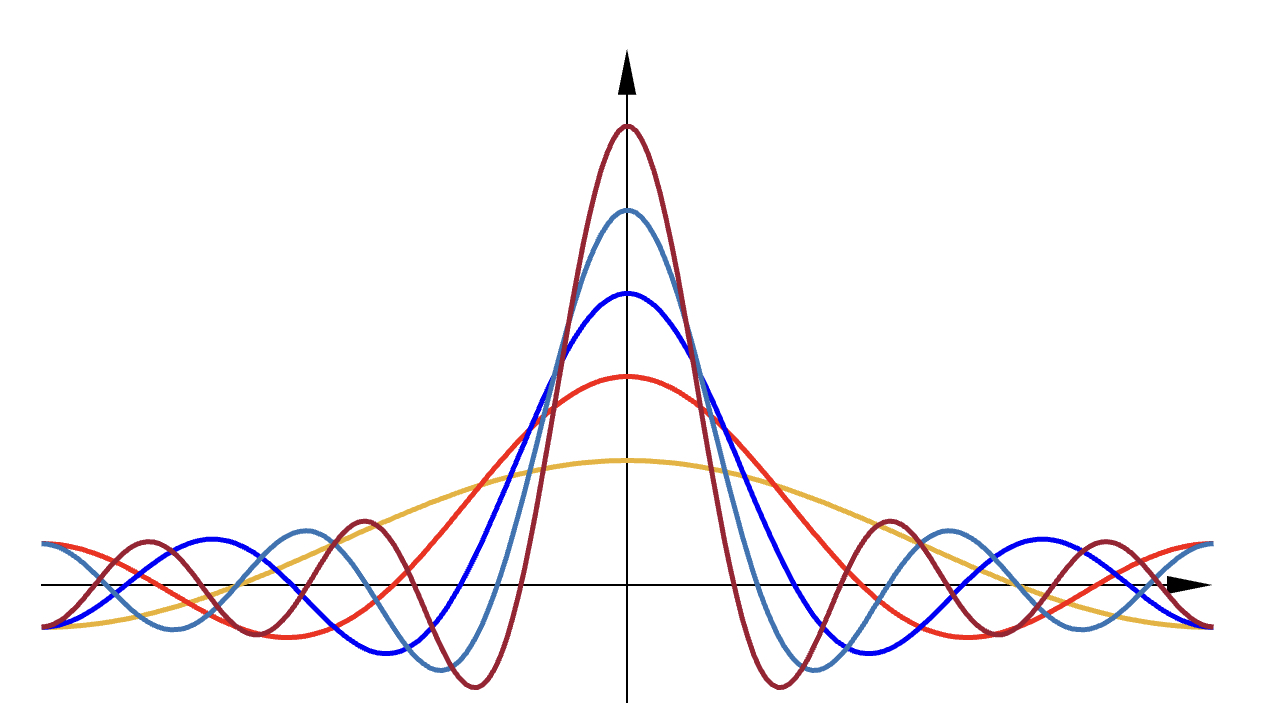}

\caption{The Dirichlet kernels $D_n:[-\pi, \pi]\to\mb{R}$ for $1\leq n\leq 5$.}\label{figure74}
\end{figure}

\begin{remark}{}
By definition,  the Dirichlet kernel $D_n(x)$ is equal to
\[D_n(x)=\sum_{k=-n}^{n}e^{ikx}.\]From this, one can see that $D_n(x)$ is an infinitely differentiable $2\pi$-periodic function, and it is an even function.
\end{remark}

 Recall that $g:[a,b]\to\mb{C}$ is a step function if there is a partition $P=\{x_0, x_1, \ldots, x_l\}$ of $[a,b]$ such that for each $1\leq j\leq l$, $g:(x_{j-1}, x_j)\to\mb{C}$ is a constant function. It is easy to see that  $g:[a,b]\to\mb{C}$ is a step function if and only if both its real and imaginary parts are step functions.
The following theorem asserts  that a Riemann integrable function  $f:[a,b]\to\mb{R}$ can be approximated in $L^1$ by   step functions.
 \begin{theorem}[label=230814_7]{}
 Let $f:[a,b]\to\mb{C}$ be a Riemann integrable function. 
   For every $\varepsilon>0$, there is a step function $g:[a,b]\to\mb{C}$ such that
 \[\int_a^b |f(x)-g(x)|dx<\varepsilon.\]
  
 \end{theorem}
 \begin{myproof}{Proof}
 Let $f(x)=u(x)+iv(x)$, where $u:[a,b]\to\mb{R}$ and $v:[a,b]\to\mb{R}$ are the real and imaginary parts of $f$.
 Assume that $u_1:[a,b]\to\mb{R}$ and $v_1:[a,b]\to\mb{R}$ are step functions such that
 \[\int_a^b|u(x)-u_1(x)|dx<\frac{\varepsilon}{2}, \hspace{1cm}
 \int_a^b|v(x)-v_1(x)|dx<\frac{\varepsilon}{2}.\]
 
 Let $g:[a,b]\to\mb{R}$ be the function $g=u_1+iv_1$. 
Then  $g$ is  a step function.
By triangle inequality, we have
\[\int_a^b|f(x)-g(x)|dx\leq \int_a^b |u(x)-u_1(x)|dx+\int_a^b|v(x)-v_1(x)|dx<\varepsilon.\]
Therefore,
it is sufficient to prove  the theorem when $f$ is a real-valued function.

Given $\varepsilon>0$, since $f:[a,b]\to\mb{R}$ is Riemann integrable, there is a partition $P=\{x_0, x_1, \ldots, x_l\}$ of $[a,b]$ such that 
\[U(f,P)-L(f,P)< \varepsilon,\]
where 
\vspace{-0.3cm}
\[
U(f,P)=\di\sum_{j=1}^l M_j(x_j-x_{j-1}),\hspace{1cm}M_j=\sup_{x_{j-1}\leq x\leq x_j}f(x),\]and
\[  L(f,P)=\di\sum_{j=1}^l m_j(x_j-x_{j-1}),\hspace{1cm} m_j=\inf_{x_{j-1}\leq x\leq x_j}f(x),
\]

are respectively the Darboux upper sum and Darboux lower sum of $f$ with respect to the partition $P$.
Define the function $g:[a,b]\to\mb{R}$ by  
\[g(x)=m_j \hspace{1cm}\text{when}\; x_{j-1}\leq x<x_j,\]and $g(b)=f(b)$. Then $g$ is a step function, and 
\[f(x)\geq g(x)\hspace{1cm}\text{for all}\;x\in [a,b].\]
  \bp
 It follows that
 \begin{align*}
 \int_a^b |f(x)-g(x)|dx&=\int_a^b (f(x)-g(x))dx\\
 &= \int_a^b f(x)dx-\sum_{j=1}^l\int_{x_{j-1}}^{x_j}g(x)dx\\
 &\leq U(f,P)-\sum_{j=1}^lm_j(x_j-x_{j-1})\\
 &=U(f,P)-L(f,P)<\varepsilon.
 \end{align*}
This completes the proof.

 \end{myproof}
 
   \begin{figure}[ht]
\centering
\includegraphics[scale=0.2]{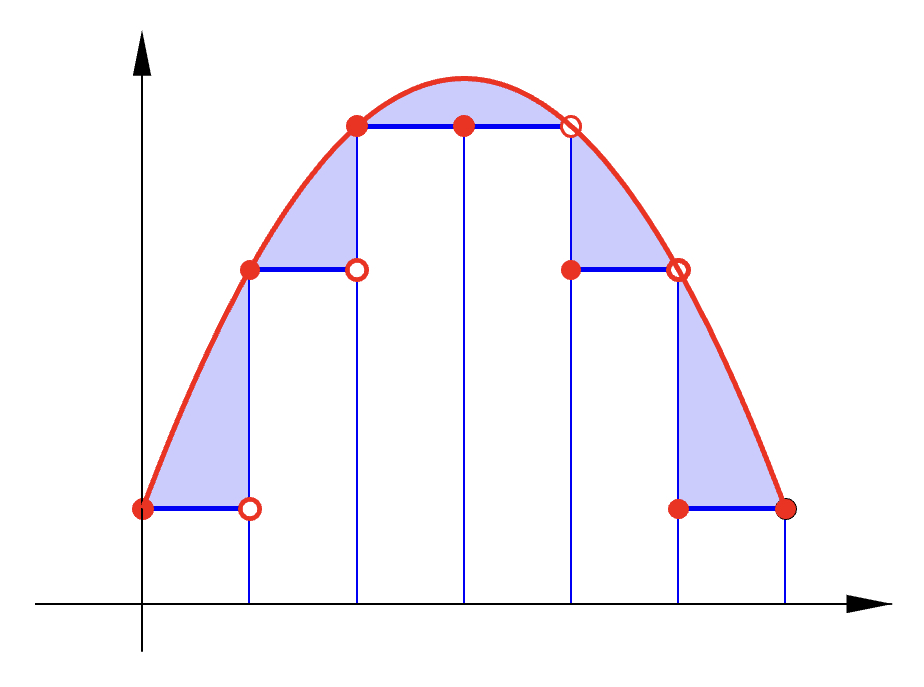}

\caption{Approximating a Riemann integrable function by a step function.}\label{figure75}
\end{figure}
An important tool in the proof of pointwise convergence of Fourier series is the Riemann-Lebesgue lemma. This lemma is important in its own right. Hence, we state it in the most general setting.
\begin{theorem}{The Riemann-Lebesgue Lemma}
Let $I=[a,b]$. If $f:I\to\mb{C}$ is a Riemann integrable function, then
\[\lim_{\beta\to\infty}\int_a^b f(x)e^{i\beta x}dx=0.\]
\end{theorem}
\begin{myproof}{Proof}
  Given $\varepsilon>0$,     Theorem \ref{230814_7} says  that there is a step function $g:[a,b]\to\mb{R}$
  such that  \[\int_a^b |f(x)-g(x)|dx<\frac{\varepsilon}{2}.\]
  This implies that
  \[  \left|\int_a^b (f(x)-g(x))e^{i\beta x}dx\right|\leq \int_a^b |f(x)-g(x)|dx<\frac{\varepsilon}{2}.\]
  Let
    $P=\{x_0, x_1, \ldots, x_l\}$ be  a partition of $[a,b]$ such that for $1\leq j\leq l$, $g(x)=m_j$   for all $x$ in  $(x_{j-1}, x_j)$.  
     Let $M=\max\{|m_1|, \ldots, |m_l|\}$. 
  Then
 \[\left|\int_{x_{j-1}}^{x_j} g(x)e^{i\beta x}dx\right|=\left|\frac{m_j}{i\beta}\left(e^{i\beta x_j}-e^{i\beta x_{j-1}}\right)\right| \leq \frac{2M}{\beta}.\]
 
 It follows that if  $\di\beta>\frac{4Ml}{\varepsilon}$,
 \[\left|\int_a^b g(x)e^{i\beta x}dx\right|\leq \sum_{j=1}^l\left|\int_{x_{j-1}}^{x_j} g(x)e^{i\beta x}dx\right|\leq  \frac{2Ml}{\beta}<\frac{\varepsilon}{2}.\]
 
 Therefore,
 \[
 \left|\int_a^b f(x)e^{i\beta x}dx\right|\leq  \left|\int_a^b (f(x)-g(x))e^{i\beta x}dx\right|+ \left|\int_a^b g(x)e^{i\beta x}dx\right|<  \varepsilon.\]
  This proves the assertion.
\end{myproof}
 Since
\[\sin\beta x=\frac{ e^{i\beta x}-e^{-i\beta x}}{2i},\] we obtain the following.
\begin{corollary}{}
Let $I=[a,b]$. If $f:I\to\mb{C}$ is a Riemann integrable function, then
\[\lim_{\beta\to\infty}\int_a^b f(x)\sin \beta x dx=0.\]
\end{corollary}

  Recall that a function $f:[a,b]\to\mb{C}$ is piecewise continuous if there is a partition $P=\{x_0, x_1, \ldots, x_l\}$ of $[a,b]$ such that for each $1\leq j\leq l$, $f:(x_{j-1}, x_j)\to\mb{C}$ is a continuous functions. It is piecewise differentiable if there is a partition $P=\{x_0, x_1, \ldots, x_l\}$ of $[a,b]$ such that for each $1\leq j\leq l$, $f:(x_{j-1}, x_j)\to\mb{C}$ is a differentiable functions. Obviously, if $f:[a,b]\to\mb{C}$ is piecewise differentiable, it is piecewise continuous. The piecewise continuity and piecewise differentiability do not impose any conditions on the partition points. In the following, we introduce  a class of functions which satisfy stronger conditions on the  partition points.
\begin{definition}{Strongly   Piecewise Differentiable Functions}
A function $f:[a,b]\rightarrow\mathbb{C}$ is  strongly piecewise continuous   if there is a partition $P=\{x_0, x_1, \ldots, x_l\}$ of $[a,b]$ such that for each $1\leq j\leq l$,  the limits 
\[
f_+(x_{j-1})=\lim_{x\rightarrow x_{j-1}^+}f(x)\quad \text{and}\quad f_-(x_j)=\lim_{x\rightarrow x_j^-}f(x)
\]exist, and the function $g_j:[x_{j-1},x_j]\rightarrow\mathbb{C}$   defined as
\[g_j(x)=\begin{cases} f_+(x_{j-1}),\quad & \text{if}\; x=x_{j-1}\\f(x),\quad & \text{if}\; x_{j-1}<x<x_j\\f_-(x_{j}),\quad &  \text{if}\;x=x_j\end{cases}\] is continuous. 
 If   $f$ is also differentiable on $(x_{j-1}, x_j)$, and the limits
\[f_+'(x_{j-1})=\lim_{h\rightarrow 0^+}\frac{f(x_{j-1}+h)-f_+(x_{j-1})}{h}\] and \[ f_-'(x_j)=\lim_{h\rightarrow 0=^+}\frac{f_-(x_j)-f(x_{j}-h) }{h}\] exist, we say that $f:[a,b]\rightarrow\mathbb{C}$ is  strongly piecewise differentiable.
\end{definition}
 Notice that a strongly piecewise differentiable function is strongly piecewise continuous and bounded. Therefore, it is Riemann integrable. 
 
 We have abused notation above and denote the limit 
\[\lim_{h\rightarrow 0^+}\frac{f(c+h)-f_+(c)}{h}\]as $f_+'(c)$. Strictly speaking, $f_+'(c)$ is the right derivative of $f$ at $c$ which is defined as
\[\lim_{h\rightarrow 0^+}\frac{f(c+h)-f(c)}{h}.\]The two expressions are equivalent if $f(c)=f_+(c)$, meaning that $f$ is right continuous at $c$. In fact, for the limit \[\lim_{h\rightarrow 0^+}\frac{f(c+h)-f(c)}{h}\] to exist, a necessary condition is $f_+(c)$ exists and is equal to $f(c)$. However, here we do not require the function $f$ to be continuous at the partition points $x_j$, $0\leq j\leq l$. We only require the function to have left and right limits at these points. Since $f(c)=f_+(c)=f_-(c)$ only when $f$ is continuous at $c$, we modify the definitions of $f_+'(c)$ and $f_-'(c)$ for functions that can have discontinuity at the point $c$.

If $x$ is an interior point of $I$ and $f:I\to\mb{C}$ is differentiable at $x$,
\[\lim_{h\to 0^+}\frac{f(x+h)-f(x)}{h}=f_+'(x)=f'(x)=f_-'(x)=\lim_{h\to 0^+}\frac{f(x)-f(x-h)}{h}.\]
Thus, if the function $f:[a,b]\to\mb{C}$ is strongly piecewise differentiable, then for any $x\in (a,b)$,  
\begin{equation}\label{230817_1}\lim_{h\to 0^+}\frac{f(x+h)+f(x-h)-f_+(x)-f_-(x)}{h}=f_+'(x)-f_-'(x).\end{equation}

\begin{example}[label=230815_1]{}
The function $f:[-\pi, \pi]\to\mb{R}$   defined as
\[f(x)=\begin{cases} \pi-x, \quad &\text{if}\; -\pi\leq x<0,\\
x^2,\quad &\text{if}\; 0\leq x\leq \pi,\end{cases}\]
is strongly piecewise differentiable.
\end{example}

\begin{figure}[ht]
\centering
\includegraphics[scale=0.2]{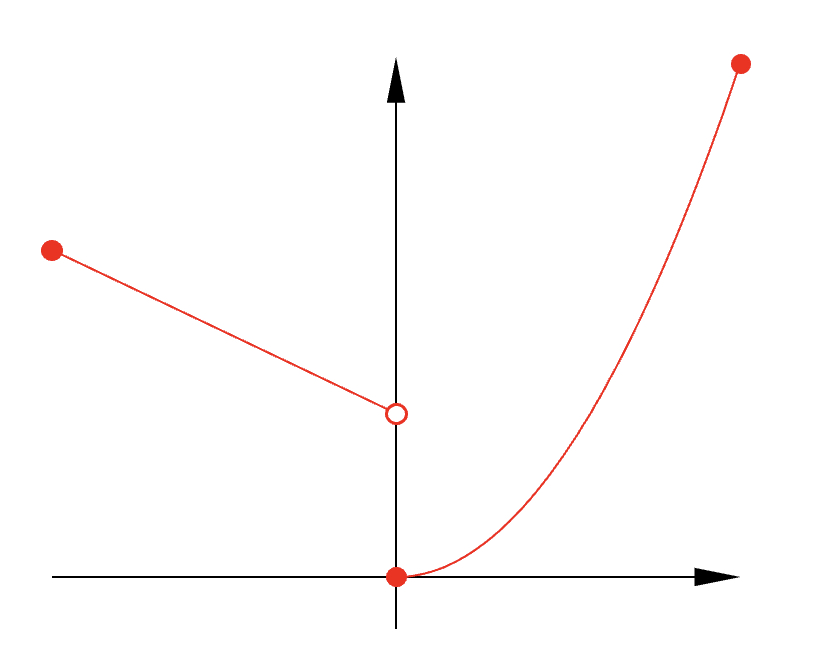}

\caption{The strongly piecewise differentiable function defined in Example \ref{230815_1}.}\label{figure73}
\end{figure}

\begin{lemma}[label=230820_4]{}
Let $f:[-\pi, \pi]\to \mb{C}$ be a strongly piecewise differentiable function, and let $\widetilde{f}:\mb{R}\to \mb{C}$ be the $2\pi$-periodic extension of $f$. Given $x\in \mb{R}$, define the function $h:[0,\pi]\to\mb{C}$ by   
\begin{equation}\label{230820_5}h(t)=\frac{\widetilde{f}(x+t)+\widetilde{f}(x-t)-\widetilde{f}_+(x)-\widetilde{f}_-(x)}{\di \sin\frac{t}{2}}, \quad t\in (0,\pi],\end{equation}
and $h(0)$ can be any value.  Then $h:[0,\pi]\to\mb{C}$ is a piecewise continuous bounded function. Hence, $h:[0,\pi]\to\mb{C}$ is Riemann integrable.
\end{lemma}
\begin{myproof}{Proof}
When $t\in [0, \pi]$, $\sin\di \frac{t}{2}=0$ only when $t=0$. Notice that $\widetilde{f}$ is a bounded piecewise continuous  function. Hence,   $h:[0,\pi]\to\mb{C}$ is  piecewise continuous. For any positive number $r$ that is less than $\pi$, $h$ is bounded on $[r,\pi]$.  To show that $h:[0,\pi]\to\mb{C}$ is bounded, it is sufficient to show that  $h(t)$ has a limit  when $t\to 0^+$. 
Now
\[\lim_{t\to 0^+}h(t) =\lim_{t\to 0^+}\frac{\widetilde{f}(x+t)+\widetilde{f}(x-t)-\widetilde{f}_+(x)-\widetilde{f}_-(x)}{t} \lim_{t\to 0^+}\frac{t}{\di \sin\frac{t}{2}}.\]
\bp
Since \[ \lim_{t\to 0^+}\frac{t}{\di \sin\frac{t}{2}}=2,\]
eq. \eqref{230817_1} gives
\[\lim_{t\to 0^+}h(t) =2\left(\widetilde{f}_+'(x)-\widetilde{f}_-'(x)\right).\]
This completes the proof.
\end{myproof}

Now we can prove the Dirichlet's theorem.
\begin{theorem}[label=230820_1]{Dirichlet's Theorem}
Let $f:[-\pi, \pi]\to \mb{C}$ be a strongly piecewise differentiable function, and let $\widetilde{f}:\mb{R}\to \mb{C}$ be the $2\pi$-periodic extension of $f$.
For every $x\in \mb{R}$,  
the Fourier series of $f$ converges at the point $x$ to 
\[\frac{\widetilde{f}_-(x)+\widetilde{f}_+(x)}{2}.\]
\end{theorem}
\begin{myproof}{Proof} Notice that the function $\widetilde{f}$ is also strongly piecewise differentiable on any compact interval.
The $n^{\text{th}}$-partial sum of the Fourier series of $f$ is
\[s_n(x)= \sum_{k=-n}^{n}c_ke^{ikx},\hspace{1cm}\text{where}\;c_k=\frac{1}{2\pi}\int_{-\pi}^{\pi} f(x)e^{-ikx}dx.\]
For a fixed real number $x$, we want to show that $s_n(x)$
converges to the number
\[u=\frac{\widetilde{f}_-(x)+\widetilde{f}_+(x)}{2}.\]
By Proposition \ref{230815_2},
\[s_n(x)=\frac{1}{2\pi}\int_{-\pi}^{\pi}f(t)D_n(x-t)dt=\frac{1}{2\pi}\int_{x-\pi}^{x+\pi}\widetilde{f}(x-t)D_n(t)dt.\]

 \bp
Since $\widetilde{f}$ and $D_n$ are $2\pi$-periodic functions, we find that
\[s_n(x)=\frac{1}{2\pi}\int_{-\pi}^{\pi}\widetilde{f}(x-t)D_n(t)dt.\]

Notice that
\[\frac{1}{2\pi}\int_{-\pi}^{\pi}D_n(t)dt=\frac{1}{2\pi}\sum_{k=-n}^n\int_{-\pi}^{\pi}e^{ikt}dt=1.\]
Therefore,
\[s_n(x)-u=\frac{1}{2\pi}\int_{-\pi}^{\pi}\left(\widetilde{f}(x-t)-u\right)D_n(t)dt.\]

Using the fact that $D_n:\mb{R}\to\mb{R}$ is an even function, we find that
\begin{align*}
s_n(x)-u&=\frac{1}{2\pi}\int_{0}^{\pi}\left(\widetilde{f}(x+t)+\widetilde{f}(x-t)-2u\right)D_n(t)dt\\&=\frac{1}{2\pi}\int_{0}^{\pi}\left(\widetilde{f}(x+t)+\widetilde{f}(x-t)-2u\right)\frac{\di\sin \left(n+\frac{1}{2}\right)t}{\di \sin\frac{t}{2}}dt\\
&=\frac{1}{2\pi}\int_{0}^{\pi}h(t) \sin \left[\left(n+\frac{1}{2}\right)t\right]dt, \end{align*}
 where $h:[0,\pi]\to\mb{C}$ is the function defined by \eqref{230820_5}.
 By Lemma \ref{230820_4}, $h:[0,\pi]\to\mb{C}$ is  Riemann integrable.
By the Riemann-Lebesgue lemma,
\[\lim_{n\to \infty}\int_{0}^{\pi}h(t) \sin \left[\left(n+\frac{1}{2}\right)t\right]dt=0.\]
This proves that
\[\lim_{n\to \infty}s_n(x)=u.\]

\end{myproof}

From the Dirichlet's theorem, we can spell out  the following explicitly.
\begin{corollary}{}
If $f:[-\pi, \pi]\to \mb{C}$ is a strongly piecewise differentiable function, then its Fourier series converges pointwise. Denote by $F:\mb{R}\to\mb{C}$ the Fourier series of $f$. Then 
for any $x\in\mb{R}$, 
\[F(x)=\frac{\widetilde{f}_+(x)+\widetilde{f}_-(x)}{2},\]
where  $\widetilde{f}:\mb{R}\to \mb{C}$ is the $2\pi$-periodic extension of $f$. This implies that
\begin{enumerate}[(a)]
\item If $x\in (-\pi, \pi)$ and $f$ is continuous at $x$, then $F(x)=f(x)$.
\item If $f$ is right continuous at $-\pi$, left continuous at $\pi$, and $f(-\pi)=f(\pi)$, then $F(-\pi)=F(\pi)=f(-\pi)=f(\pi)$.

\end{enumerate}
\end{corollary}
Let us look at a few examples.
\begin{example}[label=230816_3]{}
Let $f:[-\pi,\pi]\to\mb{R}$ be the function $f(x)=x$ considered 
in Example \ref{230815_5}. Since $f$ is a strongly differentiable function, the Fourier series of $f$ converges pointwise.  We have shown that the Fourier series   is given by
 \[ F(x)= \sum_{k=1}^{\infty}2(-1)^{k-1} \frac{\sin kx}{k}.\]
 
For any $x\in (-\pi , \pi)$, this series converges to $x$. When $x=\pi$, it converges to
 \[0=\frac{f_+(-\pi)+f_-(\pi)}{2}.\]
  
 When $x=\di\frac{\pi}{2}$, since
$\di\sin\frac{\pi k}{2}$ is 0 when $k=2n$, it is equal to 
 1 when $k=4n+1$, and it  is equal to
 $-1$ when $k=4n+3$, we deduce that
 \[\frac{1}{2}F\left(\frac{\pi}{2}\right)= 1-\frac{1}{3}+\frac{1}{5}-\frac{1}{7}+\cdots=\frac{1}{2}f\left(\frac{\pi}{2}\right)=\frac{\pi}{4},\]which is just the Newton-Gregory formula.
\end{example}
\begin{figure}[ht]
\centering
\includegraphics[scale=0.2]{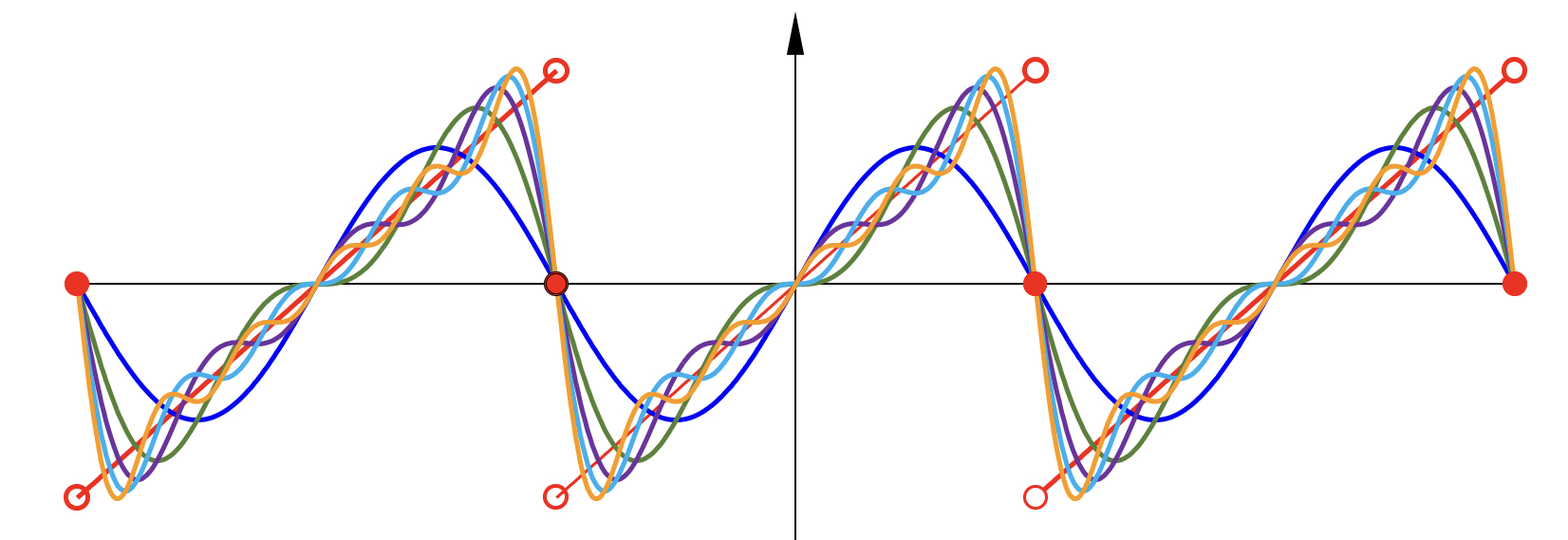}

\caption{Convergence of the Fourier series of the function $f(x)=x, -\pi < x<\pi$.}\label{figure69_2}
\end{figure}

\begin{example}[label=230815_10]{}
 The function $f:[0,2\pi]\to \mb{R}$, $\di f(x)=x(2\pi -x)$ considered in Example \ref{230815_6} is a strongly piecewise differentiable function. Hence, its Fourier series 
 \[F(x)= \frac{2}{3}\pi^2-4\sum_{k=1}^{\infty}\frac{\cos kx}{k^2}\]

 converges everywhere. Since $f:[-\pi,\pi]\to\mb{C}$ is continuous and $f(-\pi)=f(\pi)$, we find that 
   $F(x)=f(x)$ for all $x\in [-\pi,\pi]$.
 In particular, setting $x=0$ and $x=\pi$ respectively, we find that
 \[F(0)=\frac{2\pi^2}{3}-4\sum_{k=1}^{\infty}\frac{1}{k^2}=f(0)=0,\]
\[F(\pi)=\frac{2\pi^2}{3}-4\sum_{k=1}^{\infty}\frac{(-1)^k}{k^2}=f(\pi )=\pi^2.\]These give
\begin{gather*}
\sum_{k=1}^{\infty}\frac{1}{k^2}=1+\frac{1}{2^2}+\frac{1}{3^2}+\frac{1}{4^2}+\frac{1}{5^2}+\frac{1}{6^2}+\cdots =\frac{\pi^2}{6},\\ 
\sum_{k=1}^{\infty}\frac{(-1)^{k-1}}{k^2}=1-\frac{1}{2^2}+\frac{1}{3^2}-\frac{1}{4^2}+\frac{1}{5^2}-\frac{1}{6^2}+\cdots =\frac{\pi^2}{12}.
\end{gather*}
\end{example}
\begin{figure}[ht]
\centering
\includegraphics[scale=0.2]{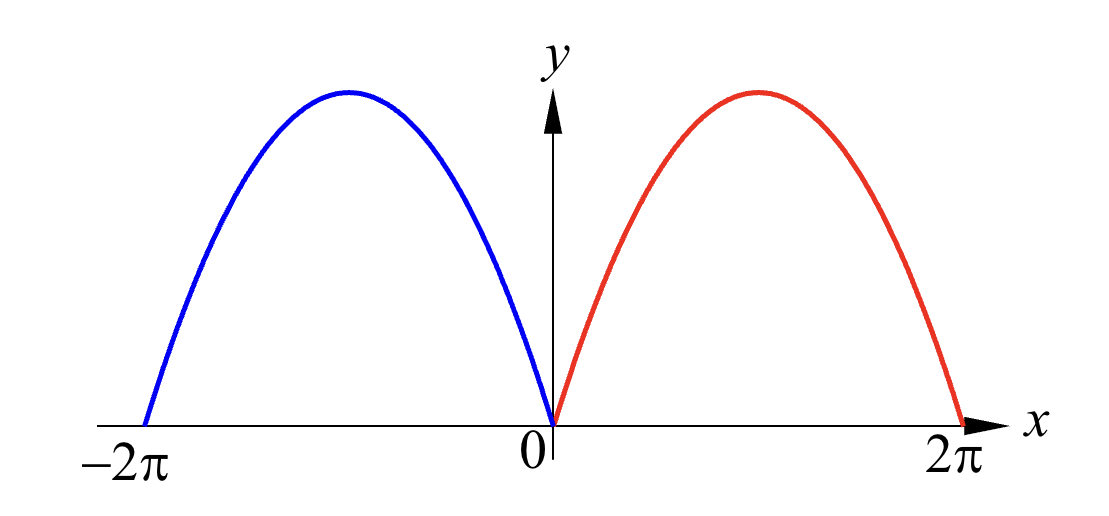}
\caption{The function $f:[0,2\pi]\to\mb{R}$, $f(x)=x(2\pi-x)$ and its entension.}\label{figure71_2}
\end{figure}
\begin{example}[label=230815_8]{}
In Example \ref{230814_1}, we consider the function $g:[-\pi, \pi]\to\mb{R}$   defined as
 \[g(x)=\begin{cases} 1,\quad &\text{if}\; a\leq x\leq b,\\
 0,\quad &\text{otherwise},\end{cases}\] where $a$ and $b$ are two numbers satisfying $-\pi\leq a<b\leq \pi$. Notice that $g$ is a strongly piecewise differentiable function. Thus, its Fourier series
 \[G(x)=\frac{b-a}{2\pi}+\frac{i}{2\pi }\sum_{k=1}^{\infty} \frac{(e^{-ikb}-e^{-ika})e^{ikx}-(e^{ikb}-e^{ika})e^{-ikx}}{  k}
 \]converges pointwise. 
  If $x\in (a, b)$, Dirichlet's theorem says that $G(x)=g(x)=1$. Hence, for any $x\in (a,b)\subset(-\pi, \pi)$,
 \[b-a=2\pi-i \sum_{k=1}^{\infty} \sum_{k=1}^{\infty} \frac{(e^{-ikb}-e^{-ika})e^{ikx}-(e^{ikb}-e^{ika})e^{-ikx}}{  k}. \]

\end{example}
 
 \begin{remark}{}
 If we scrutinize the proof of Theeorem \ref{230820_1}, we find that a necessary and sufficient condition for the Fourier series of a $2\pi$-periodic function $f:\mb{R}\to\mb{C}$ to converge at a point $x$ is that the limit 
\[\lim_{n\to\infty} \int_{0}^{\pi}\frac{f(x+t)+f(x-t)}{ t}\sin\left[ \left(n+\frac{1}{2}\right)t\right]dt\]should exists. This is known as the {\it Riemann's localization theorem}. Theeorem \ref{230820_1} says that if $f:[-\pi,\pi]\to\mb{R}$ is strongly piecewise differentiable, then this limit exists. This is sufficient for most of our applications.
 \end{remark}

\begin{remark}{Fourier Sine Series and Fourier Cosine Series}
Let $L>0$, and let $f:[0, L]\to\mb{C}$ be a Riemann integrable function defined on $[0,L]$. We can extend $f$ to be an odd function $f_o:[-L,L]\to\mb{C}$ by defining
\[
f_o(x)=\begin{cases}- f(-x),\quad &\text{if}\; -L\leq x<0,\\
0,\quad &\text{if}\quad x=0,\\
f(x),\quad &\text{if}\; \;\;0<x\leq L.\end{cases}
\] 
We can also extend $f$ to be an even function $f_e:[-L,L]\to \mb{C}$ by defining
\[f_e(x)=\begin{cases}  f(-x),\quad &\text{if}\; -L\leq x<0,\\
f(x),\quad &\text{if}\;\; \;0\leq x\leq L.\end{cases}
\] The Fourier series of $f_o:[-L, L]\to\mb{C}$ is called the Fourier sine series of $f:[0,L]\to\mb{C}$. The  Fourier series of $f_e:[-L, L]\to\mb{C}$ is called the Fourier cosine series of $f:[0,L]\to\mb{C}$. 
\end{remark}

\vp
\noindent
{\bf \large Exercises  \thesection}
\setcounter{myquestion}{1}
\begin{question}{\themyquestion}
Consider the function $f:[-\pi, \pi]\to\mb{R}$, $f(x)=x^3-\pi^2 x$.
\begin{enumerate}[(a)]
\item Find the Fourier series of $f$.
\item Use the Fourier series to find the sum
\[\sum_{k=1}^{\infty}\frac{(-1)^{k-1}}{(2k-1)^3}=1-\frac{1}{3^3}+\frac{1}{5^3}-\frac{1}{7^3}+\cdots.\]
\end{enumerate}
\end{question}
 
\atc

\begin{question}{\themyquestion}
Let $f:[-\pi, \pi]$ be the function defined as
\[f(x)=\begin{cases} x+\pi,\quad &\text{if}\; -\pi \leq x<0,\\x-\pi,\quad &\text{if}\quad 0\leq x\leq \pi.\end{cases}\]
\begin{enumerate}[(a)]
\item
Find the Fourier series of $f$.
\item Study the pointwise convergence of the Fourier series.

\end{enumerate}

\end{question}
\atc
\begin{question}{\themyquestion}
  Study the pointwise convergence of the Fourier series of the function $f:[-\pi, \pi]\to\mb{R}$, $f(x)=|x|$ obtained in Exercises \ref{sec7.1}.
\end{question}
\atc
\begin{question}{\themyquestion}
  Study the pointwise convergence of the Fourier series of the function $f:[-\pi, \pi]\to\mb{R}$, $f(x)=x^2$ obtained in Exercises \ref{sec7.1}.
\end{question}

 \atc
 
 \begin{question}{\themyquestion}
   Study the pointwise convergence of the  Fourier series of the function $f:[0, 2\pi]\to\mb{R}$, $f(x)=x^2$ obtained in Exercises \ref{sec7.1}.
\end{question}

\section[The $L^2$ Convergence of a Fourier Series ]{The $\pmb{L^2}$ Convergence of a Fourier Series}

In this section, we consider the $L^2$-convergence of a 
Fourier series. 
We first define $L^2$-converges for a sequence of Riemann integrable functions.
\begin{definition}{$\pmb{L^2}$-Convergence}
Let $I=[a,b]$ be an interval in $\mb{R}$, and let $\{f_n:I\to\mb{C}\}$ be a sequence of Riemann integrable functions. We say that $\{f_n:I\to\mb{C}\}$ converges in $L^2$ to a function $g:I\to\mb{C}$ in $\mathcal{R}(I,\mb{C})$ if 
\[\lim_{n\to\infty}\Vert f_n-g\Vert=0.\]
\end{definition}
In the vector space
$\mathcal{R}(I,\mb{C})$, we have nonzero functions $h:I\to\mb{C}$ which has zero norm. Hence, if $\{f_n:I\to\mb{C}\}$ converges in $L^2$ to a function $g:I\to\mb{C}$, the function $g$ is not unique. Nevertheless, we have the following.
\begin{theorem}[label=230814_2]{}
Let $I=[a,b]$ be an interval in $\mb{R}$, and let $\{f_n:I\to\mb{C}\}$ be a sequence of   functions in $\mathcal{R}(I,\mb{C})$ that converges in $L^2$ to the two functions $g_1:I\to\mb{C}$ and $g_2:I\to\mb{C}$ in $\mathcal{R}(I,\mb{C})$, then $g_1$ and $g_2$ are $L^2$-equivalent.
\end{theorem}

\begin{myproof}{Proof}
By triangle inequality,
\[\Vert g_1-g_2\Vert\leq\Vert f_n-g_1\Vert+\Vert f_n-g_2\Vert.\]
Since \[\lim_{n\to\infty}\Vert f_n-g_1\Vert=0\quad\text{and}\quad \lim_{n\to\infty}\Vert f_n-g_2\Vert=0,\]
we find that $\Vert g_1-g_2\Vert=0$. Thus, $g_1$ and $g_2$ are $L^2$-equivalent.
\end{myproof}

\begin{example}{}
Consider the sequence of functions $\{f_n:[0,1]\to\mb{R}\}$ defined as
\[f_n(x)=\begin{cases} 1,\quad & \text{if}\;x=\di \frac{m}{n}\;\text{for some integer $m$},\\0,\quad &\text{otherwise}.\end{cases}\]
Since $f_n:[0,1]\to\mb{R}$ is a function that is nonzero only for finitely many points, we find that 
$\Vert f_n\Vert =0$.
This implies that $\{f_n:I\to \mb{R}\}$ converges in $L^2$ to the function $f_0:I\to\mb{R}$ that is identically zero. However, $\{f_n:I\to\mb{R}\}$ does not converge pointwise. Take for example the point $x_0=1/2$. Then $f_n(x_0)=1$ if $n$ is even and $f_n(x_0)=0$ if $n$ is odd. Hence, the sequence $\{f_n(x_0)\}$ does not converge. In other words, for sequences of functions, pointwise convergence and $L^2$-convergence are different.
\end{example}
The Fourier series  of a Riemann integrable function $f:I\to\mb{C}$  converges in $L^2$ to  the function $f:I\to\mb{C}$ if 
\[\lim_{n\to\infty}\Vert s_n-f\Vert=0.\]Here $s_n(x)$ is the $n^{\text{th}}$-partial sum of the Fourier series. 
The main theorem we want to prove in this section is the Fourier series of any Riemann integrable function $f:I\to\mb{C}$ converges in $L^2$ to $f$ itself.
 We start with
 the following theorem which asserts that a Riemann integrable function  $f:[a,b]\to\mb{C}$   can be approximated in $L^2$ by   step functions.
 \begin{theorem}[label=230814_7_2]{}
 Let $f:[a,b]\to\mb{C}$ be a Riemann integrable function. 
  For every $\varepsilon>0$, there is a step function $g:[a,b]\to\mb{C}$ such that
 \[\Vert f-g\Vert<\varepsilon.\]
  
 \end{theorem}
 \begin{myproof}{Proof}
 As in the proof of Theorem \ref{230814_7}, it is sufficient to consider the case where the function $f$ is real-valued.
Since $f:[a,b]\to\mb{R}$ is Riemann integrable, it is bounded. Therefore, there exists $M>0$ such that
\[|f(x)|\leq M\hspace{1cm}\text{for all}\;x\in [a,b].\]

Given $\varepsilon>0$,  Theorem \ref{230814_7} says that there is a step function $g:[a,b]\to\mb{R}$ such that 
\[\int_a^b |f(x)-g(x)|dx<\frac{\varepsilon^2}{2M}.\]
 
By the construction of $g$ given in the proof of Theorem \ref{230814_7}, we find that 
$ |g(x)|\leq M$ for all $x\in [a,b]$.
Therefore,
\[(f(x)-g(x))^2\leq |f(x)-g(x)||f(x)+g(x)|\leq 2M|f(x)-g(x)|.\]
 
This implies that
\[\Vert f-g\Vert^2=\int_a^b (f(x)-g(x))^2dx\leq 2M\int_a^b |f(x)-g(x)|dx<\varepsilon^2.\]
Hence, $\Vert f-g\Vert<\varepsilon$.
 \end{myproof}

\begin{theorem}[label=230815_9]{}
Let $I=[-\pi, \pi]$. Given a Riemann integrable function $f:I\to\mb{C}$, let
$\di \sum_{k=-\infty}^{\infty}c_ke^{ikx}$
be its Fourier series, and let $\di s_n(x)=\sum_{k=-n}^nc_ke^{ikx}$ be the $n^{\text{th}}$-partial sum. We have the followings.
\begin{enumerate}[(a)]
\item For each $n\geq 0$, $\Vert s_n\Vert^2=\di 2\pi\sum_{k=-n}^n |c_k|^2$.
\item
For each $n\geq 0$, we have the Bessel's inequality $\Vert s_n\Vert\leq \Vert f\Vert$. 
\item The Fourier series converges in $L^2$ to  $f$ if and only if
\[\lim_{n\to\infty}\Vert s_n\Vert^2=\Vert f\Vert^2.\]
\end{enumerate}
\end{theorem}

\begin{myproof}{Proof}
For $n\in\mathbb{Z}$, let $\phi_n:\mb{R}\to\mb{C}$ be the function $\phi_n(x)=e^{inx}$. Then $\mathcal{S}=\left\{\phi_n\,|\,n\in\mb{Z}\right\}$ is an orthogonal system of functions in $\mathcal{R}(I,\mathbb{C})$, and 
$\Vert \phi_n\Vert=\sqrt{2\pi}$ for all $n\in\mathbb{Z}$. For $n\geq 0$, the set $\mathcal{S}_n=\left\{\phi_k\,|\, -n\leq k\leq n\right\}$ spans the subspace $W_n$, and 
\[\sum_{k=-n}^{n}c_k\phi_k(x)=s_n(x)=\left(\text{proj}_{W_n}f\right)(x).\]

Since $\mathcal{S}_n$ is an orthogonal system, the generalized Pythagoras theorem says that
\[\Vert s_n\Vert^2=\sum_{k=-n}^n|c_k|^2\Vert\phi_k\Vert^2=2\pi \sum_{k=-n}^n |c_k|^2.\]
This proves part (a). 
 
For part (b), recall that $f-s_n$ is orthogonal to $s_n$. By generalized Pythagoras theorem again,
\[\Vert f\Vert^2=\Vert s_n+(f-s_n)\Vert^2=\Vert s_n\Vert^2+\Vert f-s_n\Vert^2\geq \Vert s_n\Vert^2.\]
Hence, we find that $\Vert s_n\Vert\leq \Vert f\Vert$. This proves part (b).

Part (c) follows from  
\[\Vert f\Vert^2 =\Vert s_n\Vert^2+\Vert f-s_n\Vert^2.\]
 The Fourier series converges in $L^2$ to  $f$ if and only if $\di\lim_{n\to\infty}\Vert s_n-f\Vert=0$, if and only if
 \[\lim_{n\to\infty}\Vert s_n\Vert^2=\Vert f\Vert^2.\]
\end{myproof}

\begin{remark}[label=230816_9]{}
Part (b) of Theorem \ref{230815_9} says  for the trigonometric series $\di\sum_{k=-\infty}^{\infty}c_ke^{ikx}$ to be the Fourier series of a Riemann integrable function, it is necessary that the series $\di\sum_{k=-\infty}^{\infty}|c_k|^2$ is convergent.
\end{remark}

Now we will prove that the Fourier series of a special type of  step functions converges in $L^2$ to the function itself. 
\begin{theorem}[label=230815_11]{}
 Let $a$ and $b$ be two numbers satisfying $-\pi\leq a<b\leq \pi$, and let $g:[-\pi, \pi]\to\mb{R}$ be the function defined as
 \[g(x)=\begin{cases} 1,\quad &\text{if}\; a\leq x\leq b,\\
 0,\quad &\text{otherwise}.\end{cases}\] The Fourier series of $g$ converges in $L^2$ to the function $g$.
\end{theorem}
\begin{myproof}{Proof}
By Theorem \ref{230815_9}, 
it is sufficient to show that 
 
\[\lim_{n\to\infty} \Vert s_n\Vert^2 =\Vert g\Vert^2.\]

Now,
\[\Vert g\Vert^2=\int_{-\pi}^{\pi}|g(x)|^2 dx=\int_a^b dx=b-a.\]
In Example \ref{230814_1}, we have seen that the Fourier coefficients of $g$ is
\[c_0=\frac{b-a}{2\pi},\quad\text{and}\quad c_k=\frac{e^{-ikb}-e^{-ika}}{-2  \pi i k}\quad\text{when}\; k\neq 0.\]
By part (a) of Theorem \ref{230815_9},
\begin{align*}
\Vert s_n\Vert^2&=2\pi \left(|c_0|^2+\sum_{k=1}^n(|c_k|^2+|c_{-k}|^2)\right)\\
&=\frac{(b-a)^2}{2\pi} +4\pi \sum_{k=1}^n \frac{(e^{ikb}-e^{ika})(e^{-ikb}-e^{-ika})}{4\pi^2 k^2}\\
&=\frac{(b-a)^2}{2\pi} +\frac{2}{\pi}\sum_{k=1}^{n}\frac{1-\cos k(b-a)}{k^2}.
\end{align*}
\bp
By Example \ref{230815_10},
 \[  \frac{2}{3}\pi^2-4\sum_{k=1}^{\infty}\frac{\cos kx}{k^2}=x(2\pi-x)\hspace{1cm}\text{for all}\;x\in [0,2\pi],\] and 
 \[\sum_{k=1}^{\infty}\frac{1}{k^2}=\frac{\pi^2}{6}.\]
  
 Therefore,
 \begin{align*}
 \lim_{n\to\infty}\Vert s_n\Vert^2&= \frac{(b-a)^2}{2\pi}
 +\frac{2}{\pi}\sum_{k=1}^{\infty}\frac{1}{k^2}-\frac{2}{\pi}\sum_{k=1}^{\infty}\frac{\cos k(b-a)}{k^2}\\
 &=\frac{(b-a)^2}{2\pi}+\frac{\pi}{3}-\frac{1}{2\pi}\left( \frac{2}{3}\pi^2-2\pi(b-a)+(b-a)^2\right)\\
 &=b-a=\Vert g\Vert^2.
 \end{align*}
 This proves that the Fourier series of $g$ converges in $L^2$ to $g$.
\end{myproof}
 
Now we can prove our main theorem.
\begin{theorem}[label=230815_12]{$\pmb{L^2}$ Convergence of Fourier Series}
Let $I=[-\pi, \pi]$ and let $f:I\to\mb{C}$ be a Riemann integrable function. Then the Fourier series of $f$ converges in $L^2$ to   $f$ itself. 
\end{theorem}
 
\begin{myproof}{Proof}
Since we will be dealing with more than one functions here, we use $s_n(f)$ to denote the $n^{\text{th}}$-partial sum of the Fourier series of $f$.

We will show that given $\varepsilon>0$, there exists a positive integer $N$ such that for all $n\geq N$,
\[\Vert s_n(f)-f\Vert<\varepsilon.\]
Fixed $\varepsilon>0$. Theorem \ref{230814_7_2} says  that there is step function $g:[-\pi, \pi]\to\mb{R}$ such that
\[\Vert f-g\Vert<\frac{\varepsilon}{3}.\]
\bp
Let
 $P=\{x_0, x_1, \ldots, x_l\}$ be the partition of $[-\pi,\pi]$ such that for each $1\leq j\leq l$, $g$ is constant on $(x_{j-1}, x_j)$. Define $g_j:[-\pi, \pi]\to\mb{C}$ by
 
 \[g_j(x)=\begin{cases} g(x),\quad &\text{if}\; x_{j-1}<x<x_j,\\
 0,\quad &\text{otherwise}.\end{cases}\]
 
 Then 
 \[g(x)=g_1(x)+\cdots +g_l(x)\hspace{1cm}\text{for all}\;x\in [-\pi, \pi]\setminus P.\]
 Since Riemann integrals are not affected by function values at finitely many points, it follows that
 for each $n\geq 0$, 
 \[s_n(g)=\sum_{j=1}^ls_n(g_j).\]

 By Theorem \ref{230815_9}, 
 \[\Vert s_n(f)-s_n(g)\Vert=\Vert s_n(f-g)\Vert\leq \Vert f-g\Vert<\frac{\varepsilon}{3}.\]
 
 By triangle inequality,
 \[\Vert g-s_n(g)\Vert=\left\Vert \sum_{j=1}^l(g_j-s_n(g_j))\right\Vert\leq \sum_{j=1}^l\Vert g_j-s_n(g_j)\Vert.\]
 By Theorem \ref{230815_11}, $\di \lim_{n\to \infty} \Vert g_j-s_n(g_j)\Vert=0$ for $1\leq j\leq l$. Therefore,
 \[\lim_{n\to\infty}\Vert g-s_n(g)\Vert=0.\]
 This implies that there is a positive integer $N$ such that
 \[\Vert g-s_n(g)\Vert<\frac{\varepsilon}{3}\hspace{1cm}\text{for all}\;n\geq N.\]
 It follows that for all $n\geq N$,
 \[\Vert f-s_n(f)\Vert\leq \Vert f-g\Vert+\Vert g-s_n(g)\Vert+\Vert s_n(f)-s_n(g)\Vert<\varepsilon.\]
 This completes the proof.
\end{myproof}

A consequence of Theorem \ref{230815_12} is the following.
\begin{theorem}[label=230815_14]{Parseval's Identity I}
Let $f:[-\pi,\pi]\to\mb{C}$ be a Riemann integrable function, and let
\[c_k=\frac{1}{2\pi}\int_{-\pi}^{\pi}f(x)e^{-ikx}dx,\hspace{1cm}k\in\mb{Z}\]be its Fourier coefficients. Then
\[\int_{-\pi}^{\pi}|f(x)|^2 dx=\Vert f\Vert^2=2\pi \sum_{k=-\infty}^{\infty} |c_k|^2.\]

\end{theorem}
 
\begin{myproof}{Proof}
In Theorem   \ref{230815_12}, we have shown that the Fourier series of $f$ converges in $L^2$ to $f$. By Theorem \ref{230815_9}, this means that
\[\lim_{n\to\infty}\Vert s_n\Vert^2=\Vert f\Vert^2.\]
Since
\[\Vert s_n\Vert^2=2\pi\sum_{k=-n}^n |c_k|^2,\]
the Parseval's identity follows.
\end{myproof}

\begin{corollary}{Parseval's Identity II}
Let   $f:[-\pi,\pi]\to\mb{R}$ be a real-valued Riemann integrable function, and let
\[a_k=\frac{1}{ \pi}\int_{-\pi}^{\pi}f(x)\cos kx dx,\hspace{1cm}k\geq 0,\]
\[b_k=\frac{1}{ \pi}\int_{-\pi}^{\pi}f(x)\sin kx dx,\hspace{1cm}k\geq 1\]be its Fourier coefficients. Then
\[\int_{-\pi}^{\pi}f(x)^2 dx=\Vert f\Vert^2=\pi \left\{\frac{a_0^2}{2}+ \sum_{k=1}^{\infty} \left(a_k^2+b_k^2\right)\right\}.\]

\end{corollary}
\begin{myproof}{Proof}
This can be proved using $c_0=\di\frac{a_0}{2}\in\mb{R}$, and when $k\geq 1$,
\[c_k=\frac{a_k-ib_k}{2},\hspace{1cm} c_{-k}=\frac{a_k+ib_k}{2}, \]and $a_k$ and $b_k$ are real numbers.
 
\end{myproof}
Let us look at a few examples.
\begin{example}{}
In Example \ref{230815_5}, we have seen that the Fourier series of the function $f:[-\pi,\pi]\to\mb{R}$, $f(x)=x$ is 
\[F(x)=   \sum_{k=1}^{\infty}2(-1)^{k-1} \frac{\sin kx}{k}.\] 

Using Parseval's identity, we deduce that
\[\pi\sum_{k=1}^{\infty}\frac{4}{k^2}=\int_{-\pi}^{\pi}x^2dx=\frac{2\pi^3}{3}.\]
 
This gives
\[\sum_{k=1}^{\infty} \frac{1}{k^2}=\frac{\pi^2}{6},\]an identity we have obtained before.
\end{example}

\begin{example}{}
In Example \ref{230815_6}, we have seen that the Fourier series of the function $f:[0,2\pi]\to\mb{R}$, $f(x)=x(2\pi -x)$ is 
\[F(x)= \frac{2}{3}\pi^2-4\sum_{k=1}^{\infty}\frac{\cos kx}{k^2}.\]
 \be
Using Parseval's identity, we have
\begin{align*}
&2\pi \times \frac{4}{9}\pi^4+16\pi\sum_{k=1}^{\infty}\frac{1}{k^4}\\
&=\int_0^{2\pi}x^2(2\pi -x)^2dx=\int_0^{2\pi}\left(4\pi^2 x^2-4\pi x^3+x^4\right)dx\\
&=\left[\frac{4\pi^2 x^3}{3}- \pi x^4+\frac{x^5}{5}\right]_0^{2\pi}\\
&=\left(\frac{32}{3}-16+\frac{32}{5}\right)\pi^5= \frac{16}{15}\pi^5.
\end{align*} 
Therefore,
\[\sum_{k=1}^{\infty}\frac{1}{k^4}=\frac{\pi^4}{15}-\frac{ \pi^4}{18}=\frac{\pi^4}{90}.\]
\end{example2}

From the Parseval's identity, we can also obtain the following.
\begin{theorem}[label=230815_15]{}
Let $I=[-\pi, \pi]$. Given that $f:I\to \mb{C}$ and $g: I\to\mb{C}$ are Riemann integrable functions, let
\[c_k(f)=\frac{1}{2\pi}\int_{-\pi}^{\pi} f(x)e^{-ikx}dx,\hspace{1cm} c_k(g)=\frac{1}{2\pi}\int_{-\pi}^{\pi} g(x)e^{-ikx}dx\]
be their Fourier coefficients. Then
\[\int_{-\pi}^{\pi}f(x)\overline{g(x)}dx=\langle f,g\rangle=2\pi\sum_{k=-\infty}^{\infty}c_k(f)\overline{c_k(g)}.\]
\end{theorem}
\begin{myproof}{Proof}
This follows from Theorem \ref{230815_14} and the polarization formula
\begin{align*}
\langle f, g\rangle &=\frac{1}{4}\left( \langle f+g, f+g\rangle -\langle f-g, f-g\rangle\right.\\&\hspace{1cm}\left. +i\langle f+ig, f+ig\rangle - i\langle f-ig, f-ig\rangle\right).\end{align*}
\end{myproof}

We can use Theorem \ref{230815_15} to prove that Fourier series can be integrated term by term.
\begin{theorem}[label=230815_16]{Term-by-Term Integration of Fourier Series}
Let $I=[-\pi,\pi]$. Given that $f:I\to\mb{C}$ is a Riemann integrable function, let 
$\di \sum_{k=-\infty}^{\infty}c_k e^{ikx}$ be its Fourier series. On any compact interval $J=[a,b]$ that is contained in $I$, we can integrate term by term and obtain
\[\int_a^b f(x)dx=\sum_{k=-\infty}^{\infty} c_k\int_a^b e^{ikx}dx.\]
\end{theorem}

\begin{myproof}{Proof}
Let   $g:[-\pi, \pi]\to\mb{R}$ be the function defined as
 \[g(x)=\begin{cases} 1,\quad &\text{if}\; a\leq x\leq b,\\
 0,\quad &\text{otherwise}.\end{cases}\] Then
 \[c_k(g)= \frac{1}{2\pi}\int_a^b e^{-ikx}dx.\]
 Using Theorem \ref{230815_15}, we find that
 \begin{align*}
 \int_a^b f(x)dx&=\int_{-\pi}^{\pi}f(x)\overline{g(x)}dx\\
 &=2\pi \sum_{k=-\infty}^{\infty} c_k(f)\overline{c_k(g)}\\
 &=\sum_{k=\infty}^{\infty}c_k(f)\int_a^b e^{ikx}dx.
 \end{align*}This proves the assertion.
\end{myproof}

\begin{remark}{}
Theorem \ref{230815_16} is remarkable since we do not require the Fourier series of $f$ to converge uniformly.
\end{remark}

\begin{example}{}
In Example \ref{230815_5}, we have seen that the Fourier series of the function $f:[-\pi,\pi]\to\mb{R}$, $f(x)=x$ is 
\[F(x)=   \sum_{k=1}^{\infty}2(-1)^{k-1} \frac{\sin kx}{k}.\]  For $x\in [-\pi, \pi]$, term-by-term integration gives
\[\int_0^x  t dt =\sum_{k=1}^{\infty}(-1)^{k-1}\frac{2}{k}\int_0^x  \sin kt  dt.\]
 
This implies that
\[\frac{x^2}{2}=\sum_{k=1}^{\infty}(-1)^{k-1}\frac{2(1-\cos k x)}{k^2}.\]
Since \[ \sum_{k=1}^{\infty}(-1)^{k-1}\frac{1}{k^2}=\frac{\pi^2}{12},\]
we find that
\[x^2=\frac{\pi^2}{3}+\sum_{k=1}^{\infty}(-1)^{k}\frac{4}{k^2}\cos kx.\]
\end{example}

\begin{theorem}[label=230815_17]{General Term-by-Term Integration of Fourier Series}
Let $I=[-\pi,\pi]$. Given that $f:I\to\mb{C}$ is a Riemann integrable function, let 
$\di \sum_{k=-\infty}^{\infty}c_k e^{ikx}$ be its Fourier series. Let $g:I\to\mb{C}$ be any other Riemann integrable function. On any compact interval $J=[a,b]$ that is contained in $I$, we have
\[\int_a^b f(x)g(x)dx=\sum_{k=-\infty}^{\infty} c_k\int_a^b g(x) e^{ikx}dx.\]
\end{theorem}

\begin{myproof}{Proof}
Let   $  h:[-\pi, \pi]\to\mb{R}$ be the function defined as
 \[  h(x)=\begin{cases} \overline{g(x)},\quad &\text{if}\; a\leq x\leq b,\\
 0,\quad &\text{otherwise}.\end{cases}\] Then
 \[c_k(  h)= \frac{1}{2\pi}\int_a^b\overline{g(x)} e^{-ikx}dx.\]

 Hence,
 \[\overline{c_k(  h)}=\frac{1}{2\pi}\int_a^b g(x)e^{ikx}dx.\]
 
 Using Theorem \ref{230815_15}, we find that
 \begin{align*}
 \int_a^b f(x)g(x)dx&=\int_{-\pi}^{\pi}f(x)\overline{h(x)}dx\\
 &=2\pi \sum_{k=-\infty}^{\infty} c_k(f)\overline{c_k(  h)}\\
 &=\sum_{k=-\infty}^{\infty}c_k(f)\int_a^b g(x) e^{ikx}dx.
 \end{align*}This proves the assertion.
\end{myproof}

\vp
\noindent
{\bf \large Exercises  \thesection}
\setcounter{myquestion}{1}
\begin{question}{\themyquestion}
Consider the function $f:[-\pi, \pi]\to\mb{R}$, $f(x)=x^3-\pi^2 x$ whose Fourier series has been obtained in Excerses \ref{sec7.2}.
Use Parseval's identity to find the sum 
\[\sum_{k=1}^{\infty}\frac{1}{k^6}=1+\frac{1}{2^6}+\frac{1}{3^6}+\frac{1}{4^6}+\cdots.\]
 
\end{question}
 
\atc

\begin{question}{\themyquestion}
Use term by term integration to find the Fourier series of the function $f:[-\pi,\pi]\to\mb{R}$, $f(x)=x^3$.
\end{question}

\section{The  Uniform Convergence of a Trigonometric Series}

In this section, we consider  uniform convergence of a 
trigonometric series.  
In volume I, we have studied uniform convergence. If a series of continuous functions  converges uniformly, then it represents  a continuous function, and the series can be integrated term-by-term.

Applying to trigonometric series, we have the following.
\begin{theorem}[label=230816_2]{}
If the trigonometric series   
$\di \sum_{k=-\infty}^{\infty}c_ke^{ikx}$ converges uniformly, it  defines a continuous $2\pi$-periodic function $F:\mb{R}\to\mb{C}$. The Fourier series of $F:[-\pi,\pi]\to \mb{C}$ is the series $\di \sum_{k=-\infty}^{\infty}c_ke^{ikx}$ itself.

\end{theorem}

\begin{myproof}{Proof}
For $k\in\mb{Z}$, the function $\phi_k(x)=e^{ikx}$ is a  $2\pi$-periodic continuous function. The  assertion that the trigonometric series  $\di \sum_{k=-\infty}^{\infty}c_ke^{ikx}$ defines a continuous $2\pi$-periodic function $F:\mb{R}\to\mb{C}$ follows from the uniform convergence.
Since each $\phi_k(z)$, $k\in\mb{Z}$ is Riemann integrable, and the series $\di\sum_{k=-\infty}^{\infty}c_ke^{ikx}$ converges uniformly, we can integrate term by term to find that
\[c_k(F)=\frac{1}{2\pi}\int_{-\pi}^{\pi}F(x)e^{-ikx}dx=\frac{1}{2\pi}\sum_{l=-\infty}^{\infty}c_l\int_{-\pi}^{\pi}e^{i(l-k)x}dx=c_k.\]

This proves that the Fourier series of $F:[-\pi,\pi]\to \mb{C}$ is the series $\di \sum_{k=-\infty}^{\infty}c_ke^{ikx}$ itself.
\end{myproof}

We would like to have a sufficient condition for a trigonometric series \[ \sum_{k=-\infty}^{\infty}c_ke^{ikx}\] to converge uniformly.
In volume I, 
we have discussed Weierstrass $M$-test for real-valued functions.  This can be generalized  to complex-valued functions in a straightforward way. Let $\{f_n:D\to\mb{C}\}$ be a sequence of functions defined on a subset $D$ of $\mb{R}$. Assume that for each $n\in\mathbb{Z}^+$, there is a constant $M_n$ so that 
\[|f_n(x)|\leq M_n\hspace{1cm}\text{for all}\;x\in D.\]
The Weiertrass $M$-test states that if $\di\sum_{n=1}^{\infty}M_n$ is convergent, then the series of functions $\di\sum_{n=1}^{\infty}f_n(x)$ converges absolutely and uniformly.

Using Weierstrass $M$-test, we can deduce the following.

\begin{theorem}[label=230816_10]{}
If the series $\di\sum_{k=-\infty}^{\infty}|c_k|$ is convergent, then the  trigonometric series   
$\di \sum_{k=-\infty}^{\infty}c_ke^{ikx}$ converges absolutely and uniformly, and it defines a continuous function $F:\mb{R}\to\mb{C}$ whose Fourier series is itself.
\end{theorem}

\begin{myproof}{Proof}
For any $k\in\mb{Z}$,
\[|c_ke^{ikx}|\leq |c_k|\hspace{1cm}\text{for all}\;x\in\mb{R}.\]
By Weierstrass $M$-test, the series
$\di \sum_{k=-\infty}^{\infty}c_ke^{ikx}$  converges absolutely and uniformly. The rest follows from Theorem \ref{230816_2}.
\end{myproof}

\begin{example}{}
Consider the series 
\[\sum_{k=1}^{\infty}\frac{\cos kx}{k^{\frac{3}{2}}}.\]\be
Since $\di\sum_{k=1}^{\infty}\frac{1}{k^{\frac{3}{2}}}$ is a $p$-series with $p=\frac{3}{2}>1$, it is convergent. 
Hence, the function $F:\mb{R}\to\mb{R}$ defined as
\[F(x)=\sum_{k=1}^{\infty}\frac{\cos kx}{k^{\frac{3}{2}}}\]is a continuous function whose Fourier series is itself.
\end{example2}

\begin{remark}{}
By Remark \ref{230816_9}, in order for the series $\di\sum_{k=-\infty}^{\infty}c_ke^{ikx}$ to be the Fourier series of a Riemann integrable function, it is necessary that the series $\di\sum_{k=-\infty}^{\infty}|c_k|^2$ is convergent. However, the convergence of $\di\sum_{k=-\infty}^{\infty}|c_k|^2$  does not imply the convergence of  $\di\sum_{k=-\infty}^{\infty}|c_k|$.  
\end{remark}
Theorem \ref{230816_2} gives a criterion for a function defined as a trigonometric series to have Fourier series that is equal to itself.
However, we will usually start by a Riemann integrable function.

\begin{theorem}[label=230816_1]{}
Let $I=[-\pi,\pi]$, and let $f:I\to\mb{C}$ be a Riemann integrable function.
If the  Fourier series
$\di\sum_{k=-\infty}^{\infty}c_ke^{ikx}$ of $f:I\to\mb{C}$ converges uniformly, it defines a continuous $2\pi$-periodic function $F:\mb{R}\to\mb{C}$  whose restriction to $I$ is  $L^2$ equivalent to $f:I\to\mb{C}$. If the periodic extension $\widetilde{f}:\mb{R}\to\mb{C}$ of $f$ is continuous at $x_0$, then $F(x_0)=\widetilde{f}(x_0)$. 
\end{theorem}
\begin{myproof}{Proof}
The fact that $\di\sum_{k=-\infty}^{\infty}c_ke^{ikx}$ defines a continuous periodic function $F:\mb{R}\to\mb{C}$ has been asserted in Theorem \ref{230816_2}. 
Since $\di\sum_{k=-\infty}^{\infty}c_ke^{ikx}$ is the Fourier series for both $f:I\to\mb{C}$ and $F:I\to\mb{C}$, Theorem \ref{230815_12} says that it  converges in $L^2$ to $f:I\to\mb{C}$ and $F:I\to\mb{C}$.  Theorem \ref{230814_2} then asserts that $F:I\to\mb{C}$ and $f:I\to\mb{C}$ are   $L^2$ equivalent. The values of two $L^2$-equivalent functions agree at a point where both of them are continuous.
\end{myproof}

\begin{example}{}
In Example \ref{230815_6}, we have seen that the Fourier series of the function $f:[0,2\pi]\to\mb{R}$, $f(x)=x(2\pi -x)$ is 
\[F(x)= \frac{2}{3}\pi^2-4\sum_{k=1}^{\infty}\frac{\cos kx}{k^2}.\]

Since the series $\di\sum_{k=1}^{\infty}\frac{1}{k^2}$ is convergent, and $f:[0,2\pi]\to\mb{R}$ is continuous with $f(0)=f(2\pi)=0$, Theorem \ref{230816_10} and Theorem \ref{230816_1} imply that the Fourier series $F(x)$ converges uniformly to  the periodic extension of the function $f(x)$. 
\end{example}

\begin{example}{}
In Example \ref{230815_5}, we have seen that the Fourier series of the function
$f:[-\pi, \pi]\to\mb{R}$, $f(x)=x$ is
\[F(x)= \sum_{k=1}^{\infty}2(-1)^{k-1}\frac{\sin kx}{k}.\]  Since the harmonic series $\di \sum_{k=1}^{\infty}\frac{1}{k}$ is divergent, we cannot apply Theorem \ref{230816_10}.
However, we can argue that the Fourier series does not converge uniformly in the following way. 
 \be
In Example \ref{230816_3}, we have used   Dirichlet's theorem to conclude that $F(x)$ converges to $f(x)$ for any $x \in (-\pi, \pi)$. It is obvious that $F(\pi)=0$.

Now
\[\lim_{x\to \pi^-}F(x)=\lim_{x\to \pi^-}f(x)=\lim_{x\to\pi^-}x=\pi\neq F(\pi).\]
This shows that $F$ is not continuous at $x=\pi$. Hence, the convergence of the Fourier series is not uniform.
\end{example2}

\begin{remark}{}
 Theorem \ref{230816_2} and Theorem \ref{230816_1} can be regarded as uniquessness of Fourier series.
\end{remark}

Next we want to consider term-by-term differentiation.

\begin{example}[label=230816_8]{}
Let us consider again  the function
$f:[-\pi, \pi]\to\mb{R}$, $f(x)=x$, whose Fourier series  is given by
\[F(x)= \sum_{k=1}^{\infty}2(-1)^{k-1}\frac{\sin kx}{k}.\]  Since $f:(-\pi,\pi)\to\mb{R}$ is continuously differentiable with $f'(x)=1$, the extension function $\widetilde{f}:\mb{R}\to\mb{R}$ is continuously differentiable at all the points where $x\notin (2n+1)\mb{Z}$, and
\[\widetilde{f}'(x)=1,\quad \text{when}\;\; x\notin (2n+1)\mb{Z}.\]
 
Hence, for the function $f':[-\pi,\pi]\to\mb{R}$, regardless of how it is defined at $\pm \pi$,  its Fourier series   is the constant 
$G(x)=1$. However, if we differentiate the Fourier series of $f:[-\pi, \pi]\to\mb{R}$ term-by-term, we obtain the series
\[\sum_{k=1}^{\infty}2(-1)^{k-1}\cos kx.\] 
\be
As  $|a_k|=2$ for all $k\in \mathbb{Z}$, $\di\sum_{k=1}^{\infty}|a_k|^2$ is divergent. Thus, the series
\[\sum_{k=1}^{\infty}2(-1)^{k-1}\cos kx.\]
 cannot be the Fourier series of any Riemann-integrable function.
\end{example2}
Example \ref{230816_8} shows that   term-by-term differentiation can   fail even though the function $f:[-\pi,\pi]\to\mb{R}$ and its derivative are strongly piecewise differentiable. 

 \begin{theorem}[label=230816_11]{}
 Let $I=[-\pi,\pi]$ and let $f:I\to\mb{R}$ be a Riemann integrable function with Fourier series
$\di\sum_{k=-\infty}^{\infty} c_ke^{ikx}$.
 If the series $\di\sum_{k=-\infty}^{\infty}c_k$ is convergent, and the series $\di\sum_{k=-\infty}^{\infty}kc_ke^{ikx}$ converges uniformly, then $f:I\to\mb{R}$ is $L^2$-equivalent to the continuously differentiable function $F:I\to\mb{R}$,
$\di F(x)=\sum_{k=-\infty}^{\infty} c_ke^{ikx}$. Moreover,  the Fourier series of  $F':I\to\mb{R}$   is
\[\sum_{k=-\infty}^{\infty}ikc_ke^{ikx},\] which converges to $F'(x)$ for all $x\in I$. 
 \end{theorem}
\begin{myproof}{Proof}
By Theorem \ref{230816_2}, the uniformly convergent series  \[\sum_{k=-\infty}^{\infty}ikc_ke^{ikx} \] defines a continuous function $G:I\to\mb{C}$, whose Fourier series is itself. 

\bp
 Let 
\[H(x)=\int_0^xG(x)dx+\sum_{k=-\infty}^{\infty}c_k.\]
By fundamental theorem of calculus, $H$ is differentiable and $H'(x)=G(x)$. Thus, $H(x)$ is continuously differentiable.
 Since the series \[\sum_{k=-\infty}^{\infty}ikc_ke^{ikx} \]converges uniformly, we can do term-by-term integration to obtain
\[H(x)=\sum_{k=-\infty}^{\infty}ikc_k\int_0^x e^{ikt}dt +\sum_{k=-\infty}^{\infty}c_k=\sum_{k=-\infty}^{\infty}c_ke^{ikx}=F(x).\] 
Hence, the function $F:I\to\mb{C}$, 
\[F(x)=\sum_{k=-\infty}^{\infty}c_ke^{ikx}\] is continuously differentiable, with derivative
\[F'(x)=H'(x)=G(x)=\sum_{k=-\infty}^{\infty}ikc_ke^{ikx}.\]
This completes the proof.

\end{myproof}
Let us look at an example.
\begin{example}{}
Consider the trigonometric series 
\[\sum_{k=1}^{\infty}\frac{\sin kx}{k^3}\hspace{1cm}\text{and}\hspace{1cm}\sum_{k=1}^{\infty}\frac{\cos kx}{k^2}.\]
Since  the series $\di\sum_{k=1}^{\infty}\frac{1}{k^3}$ and  $\di\sum_{k=1}^{\infty}\frac{1}{k^2}$ are convergent, both the series
\be
\[\sum_{k=1}^{\infty}\frac{\sin kx}{k^3}\hspace{1cm}\text{and}\hspace{1cm}\sum_{k=1}^{\infty}\frac{\cos kx}{k^2}\]

converge uniformly, and so they define  continuous functions. Let 
\[F(x)=\sum_{k=1}^{\infty}\frac{\sin kx}{k^3}\hspace{1cm}\text{and}\hspace{1cm}G(x)=\sum_{k=1}^{\infty}\frac{\cos kx}{k^2}.\]
Since \[\frac{d}{dx}\frac{\sin kx}{k^3}=\frac{\cos kx}{k^2}\hspace{1cm}\text{for all}\;k\in\mb{Z}^+,\]
  $F$ is continuously differentiable and $F'(x)=G(x)$.
\end{example2}
At the end of this section, we give a brief discussion about the Ces$\grave{\text{a}}$ro mean of a Fourier series.
As an application, we give another proof of the Weiestrass approximation theorem.

\begin{definition}{Ces$\grave{\textbf{a}}$ro Mean of a Fourier Series}
  Given a Riemann integrable function $f:[-\pi,\pi]\to\mb{C}$, let its Fourier series be
\[\sum_{k=-\infty}^{\infty}c_ke^{ikx},\hspace{1cm}\text{where  }\;c_k=\frac{1}{2\pi}\int_{-\pi}^{\pi} f(x)e^{-ikx}dx.\]
For $n\geq 1$, the $n^{\text{th}}$ Ces$\grave{\text{a}}$ro mean of the Fourier series is
\[\sigma_n(x)=\frac{s_0(x)+s_1(x)+\cdots+s_{n-1}(x)}{n},\]where
\[s_n(x)=\sum_{k=-n}^n c_ke^{ikx}\] is the $n^{\text{th}}$-partial sum of  the Fourier series.
\end{definition}

\begin{proposition}{}
Let $I=[-\pi,\pi]$, and let  $f:I\to\mb{C}$ be a Riemann integrable function. The $n^{\text{th}}$ Ces$\grave{\text{a}}$ro mean  $\sigma_n(x)$ of the Fourier series of $f$ has an integral representation given by
\[\sigma_n(x)=\frac{1}{2\pi}\int_{-\pi}^{\pi}f(t)\mathcal{F}_n(x-t)dt,\]
where $\mathcal{F}_n:\mb{R}\to\mb{R}$ is the   kernel  
\[\mathcal{F}_n(x) =\begin{cases} n,\quad &\text{if}\; x\in 2\pi\mb{Z},
\\\di
\frac{ \sin^2 \frac{nx}{2}}{ n \sin^2\frac{x}{2}},\quad &\text{otherwise}.\end{cases}\]
\end{proposition}
\begin{myproof}{Proof}
By Proposition \ref{230815_2},
\[s_n(x)=\frac{1}{2\pi}\int_{-\pi}^{\pi}f(t)D_n(x-t)dt,\]
 where
$D_n(t)$ is the Dirichlet kernel. This gives
\[\sigma_{n}(x)=\frac{1}{2\pi}\int_{-\pi}^{\pi}f(t)\mathcal{F}_n(x-t)dt,\]
where  
\[\mathcal{F}_n(x)=\frac{D_0(x)+D_1(x)+\cdots+D_{n-1}(x)}{n}.\]
Using the fact that
\[D_n(x)=\sum_{k=-n}^ne^{ikx} =\begin{cases} 2n+1,\quad &\text{if}\; x\in 2\pi\mb{Z},
\\\di
\frac{\sin \left(n+\frac{1}{2}\right)x}{ \sin\frac{x}{2}},\quad &\text{otherwise},\end{cases}\]
we find that when $x\in 2\pi\mb{Z}$, 
\[\mathcal{F}_n(x)=\frac{1+3+\ldots+(2n-1)}{n}=n.\]
\bp
When $x\notin 2\pi\mb{Z}$,
\begin{align*}
\mathcal{F}_n(x)&=\frac{1}{n\sin \frac{x}{2}}\sum_{k=0}^{n-1}\sin\left(k+\frac{1}{2}\right)x
\\
&=\frac{1}{ n\sin \frac{x}{2}}\;\text{Im}\,\left\{\sum_{k=0}^{n-1} \exp\left(i \left[k+\frac{1}{2}\right]x\right)\right\}.\\&=\frac{1}{ n\sin \frac{x}{2}}\;\text{Im}\, \left\{e^{\frac{ix}{2}}\frac{e^{inx}-1}{e^{ix}-1}\right\}\\&=\frac{1}{n\sin \frac{x}{2}}\;\text{Im}\, \left\{ \frac{e^{inx}-1}{  2i\sin \frac{x}{2} }\right\}\\
&=\frac{1-\cos nx}{2n\sin^2 \frac{x}{2} }\\
&=\frac{ \sin^2\frac{nx}{2}}{n\sin^2 \frac{x}{2} }.
\end{align*}
This completes the proof.
\end{myproof}

\begin{definition}{}For $n\geq 1$, the Fej\'er kernel
$\mathcal{F}_n:\mb{R}\to\mb{R}$ is the kernel given by
\[\mathcal{F}_n(x) =\begin{cases} n,\quad &\text{if}\; x\in 2\pi\mb{Z},
\\\di \frac{ \sin^2 \frac{nx}{2}}{ n \sin^2\frac{x}{2}},\quad &\text{otherwise}.\end{cases}\]
\end{definition}
A good property about the Fej\'er kernel is $\mathcal{F}_n(t)\geq 0$ for all $t\in\mb{R}$.

  \begin{figure}[ht]
\centering
\includegraphics[scale=0.2]{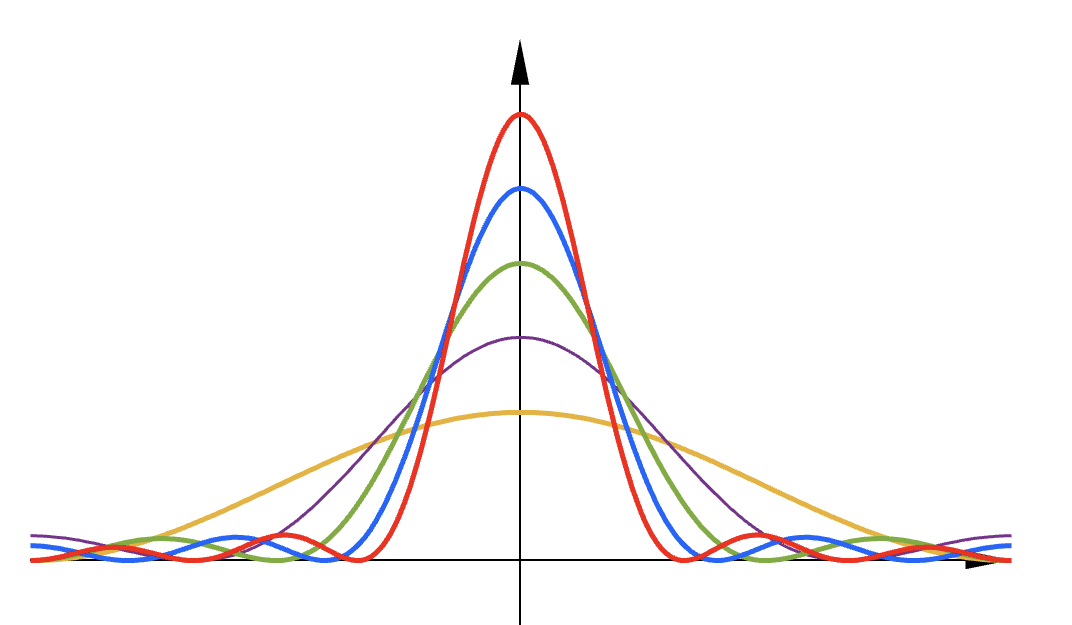}

\caption{The Fej\'er kernels $\mathcal{F}_n:[-\pi, \pi]\to\mb{R}$ for $2\leq n\leq 6$.}\label{figure77}
\end{figure}
Now we can prove the following theorem.

\begin{theorem}[label=230817_8]{}
Let $f:[-\pi,\pi]\to\mb{C}$ be a continuous function with $f(-\pi)=f(\pi)$, and let  
$\sigma_n(x)$
be the $n^{\text{th}}$  Ces$\grave{\text{a}}$ro mean of the Fourier series of $f$. Then the sequence of functions $\{\sigma_n:[-\pi,\pi]\to\mb{C}\}$ converges uniformly to the function $f:[-\pi,\pi]\to\mb{C}$.
\end{theorem}
\begin{myproof}{Proof}
As in the proof of the Dirichlet's theorem, we find that
\[\sigma_n(x)=\frac{1}{2\pi}\int_{-\pi}^{\pi}\widetilde{f}(x-t)\mathcal{F}_n(t)dt,\]
where $\widetilde{f}:\mb{R}\to\mb{C}$ be its $2\pi$-periodic extension of $f$.
Since 
\[\frac{1}{2\pi}\int_{-\pi}^{\pi} D_n(t)dt=1\hspace{1cm}\text{for all}\;n\geq 0,\]

we have
\[\frac{1}{2\pi}\int_{-\pi}^{\pi} \mathcal{F}_n(t)dt=1\hspace{1cm}\text{for all}\;n\geq 1.\]
It follows that for $x\in [-\pi,\pi]$,
\[\sigma_n(x)-f(x)=\frac{1}{2\pi}\int_{-\pi}^{\pi}\left(\widetilde{f}(x-t)-\widetilde{f}(x)\right)\mathcal{F}_n(t)dt.\]

For $x\in [-\pi, \pi]$ and $t\in [-\pi, \pi]$, $x-t\in [-2\pi, 2\pi]$. Since $ f:[-\pi,\pi]\to\mb{C}$ is a continuous function with $f(-\pi)=f(\pi)$, $\widetilde{f}:[-2\pi,2\pi]\to\mb{C}$ is continuous. Hence, it is uniformly continuous.

Given $\varepsilon>0$, there exists $\delta>0$ such that if $u$ and $v$ are in $[-2\pi ,2\pi]$ and $|u-v|<\delta$, then
\[|\widetilde{f}(u)-\widetilde{f}(v)|<\frac{\varepsilon}{2}.\]

By continuity, if $|u-v|\leq\delta$, then
\[|\widetilde{f}(u)-\widetilde{f}(v)|\leq \frac{\varepsilon}{2}.\]

 \bp
Being continuous on a compact interval, the function $\widetilde{f}:[-2\pi,2\pi]\to\mb{C}$ is also bounded. Therefore, there exists $M>0$ such that
\[|\widetilde{f}(x)|\leq M\hspace{1cm}\text{for all}\;x\in [-2\pi, 2\pi].\]

This implies that for any $x\in [-\pi, \pi]$ and $t\in [-\pi, \pi]$,
\[\left|\widetilde{f}(x-t)-\widetilde{f}(x)\right|\leq 2M.\]
 
On the other hand, if $\delta \leq |t|\leq \pi$,
\[ \sin^2 \frac{t}{2} \geq \sin^2\frac{\delta}{2}>0.\]

Therefore,
\[0\leq \mathcal{F}_n(t)\leq  \frac{1}{n  \sin^2\frac{\delta}{2}}\hspace{1cm}\text{when}\;\delta\leq|t|\leq \pi.\]
Let $N$ be  a positive integer such that
\[N>\frac{4M}{\varepsilon \sin^2\frac{\delta}{2}}.\]
For $n\geq N$ and $x\in [-\pi,\pi]$, we have
\begin{align*}
\left|\sigma_n(x)-f(x)\right|&\leq \frac{1}{2\pi}\int_{|t|\leq\delta }\left|\widetilde{f}(x-t)-\widetilde{f}(x)\right|\mathcal{F}_n(t)dt
\\
&\quad+\frac{1}{2\pi}\int_{ \delta \leq |t|\leq \pi}\left|\widetilde{f}(x-t)-\widetilde{f}(x)\right|\mathcal{F}_n(t)dt.
\end{align*}We estimate the two terms separately. Since $\mathcal{F}_n(t)\geq 0$ for all $t\in[-\pi,\pi]$,
\[0\leq \frac{1}{2\pi}\int_{|t|\leq\delta } \mathcal{F}_n(t)dt\leq \frac{1}{2\pi}\int_{-\pi }^{\pi} \mathcal{F}_n(t)dt=1.\]

Hence,
\[\frac{1}{2\pi}\int_{|t|\leq\delta }\left|\widetilde{f}(x-t)-\widetilde{f}(x)\right|\mathcal{F}_n(t)dt\leq \frac{\varepsilon}{2}\times
\frac{1}{2\pi}\int_{|t|\leq\delta } \mathcal{F}_n(t)dt\leq \frac{\varepsilon}{2}.\]
\bp
For the second term, we have
\begin{align*}
 \frac{1}{2\pi}\int_{ \delta \leq |t|\leq \pi}\left|\widetilde{f}(x-t)-\widetilde{f}(x)\right|\mathcal{F}_n(t)dt 
&\leq \frac{1}{2\pi}\int_{ \delta \leq |t|\leq \pi}\frac{2M}{n  \sin^2\frac{\delta}{2}}dt\\
&\leq \frac{2M}{N  \sin^2\frac{\delta}{2}}<\frac{\varepsilon}{2}.
\end{align*}
This shows that for all $n\geq N$, 
\[\left|\sigma_n(x)-f(x)\right|<\varepsilon\hspace{1cm}\text{for all}\;x\in [-\pi, \pi].\]
Thus, the sequence of functions $\{\sigma_n:[-\pi,\pi]\to\mb{C}\}$ converges uniformly to the function $f:[-\pi,\pi]\to\mb{C}$.
\end{myproof}

Notice that since $s_n(x)$ is in the span of $\mathcal{S}_n=\{e^{ikx}\,|\, -n\leq n\leq n\}$, $\sigma_{n+1}(x)$ is in the span of $\mathcal{S}_{n}(x)$. 
Now we apply Theorem \ref{230817_8} to give another proof of the Weierstrass approximation.
\begin{theorem}{Weierstrass Approximation Theorem}
Let $f:[a,b]\to\mathbb{R}$ be a continuous function defined on $[a, b]$. Given $\varepsilon>0$, there is a polynomial $p(x)$ such that
\[|f(x)-p(x)|<\varepsilon\hspace{1cm}\text{for all}\;x\in [a,b].\]
\end{theorem}
\begin{myproof}{Proof}
It is sufficient to prove the theorem for  a specific $[a,b]$. We take $[a,b]=[0,1]$. 
Given $f:[0,1]\to\mb{R}$ is a real-valued continuous function, we extend it to be an even function $f_e:[-1,1]\to\mb{R}$, and let $g:[-\pi,\pi]\to\mb{R}$ be the function defined as
\[g(x)=f_e(\cos x).\] 

This is well-defined since the range of $\cos x$ is $[-1,1]$. Since $\cos x$ and $f_e:[-1,1]\to\mb{R}$ are continuous even functions, $g:[-\pi,\pi]\to\mb{R}$ is a continuous even function.  Hence, we also have $g(\pi)=g(-\pi)$. 

\bp
Given $\varepsilon>0$,  Theorem \ref{230817_8} implies that there is a positive integer $n$ such that
\begin{equation}\label{230817_9}|g(x)-\sigma_{n+1}(x)|<\varepsilon.\end{equation}
Here
$\sigma_{n+1}(x)$ is the $(n+1)^{\text{th}}$  Ces$\grave{\text{a}}$ro mean of the Fourier series of $g$.  Since $g:[-\pi,\pi]\to\mb{R}$ is a real-valued even function, the Fourier series of $g$ has the form
\[\frac{a_0}{2}+\sum_{k=1}^{\infty}a_k\cos kx,\] where $a_k, k\geq 0$ are real.
This implies that 
\[\sigma_{n+1}(x)=\sum_{k=0}^{n}\alpha_k\cos kx\]
 for some real constants $\alpha_0, \alpha_1, \ldots, \alpha_{n}$.
For any $m\geq 1$, $\cos mx $ can be written as a linear combination of $1, \cos x, \cos^2x, \ldots, \cos^m x$. This shows that there are real constants $\beta_0, \beta_1, \ldots, \beta_{n}$ such that
\[\sigma_{n+1}(x)=\sum_{k=0}^{n}\beta_k\cos^kx.\]

Let
\[p(x)=\sum_{k=0}^{n} \beta_kx^k.\]
Then $\sigma_{n+1}(x)=p(\cos x)$.  Thus, \eqref{230817_9} says that
\[|f_e(\cos x)-p(\cos x)|<\varepsilon\hspace{1cm}\text{for all}\;x\in [-\pi,\pi].\]
This implies that
\[|f(x)-p(x)|<\varepsilon\hspace{1cm}\text{for all}\;x\in [0,1],\]which completes the proof of the theorem.
\end{myproof}
 \begin{remark}{}In the proof of the Weierstrass approximation theorem given above, we do not use Fourier series since
 the Fourier series of a $2\pi$-periodic continuous function does not necessary converge uniformly. 
 An example is given in \cite{Stein_Fourier}. 
 However,
  there are other approaches to prove the Weierstrass approximation theorem using Fourier series. For example, one can  approximate a continuous function uniformly by a continuous piecewise linear   function first. The Fourier series of a continuous piecewise linear function does converge uniformly to the function itself. 
  
  In the proof given above, we used the even extension $f_e$ of the given function $f$. The  Fourier series of $f_e(\cos x)$ is a cosine series, so that the Ces$\grave{\text{a}}$ro mean is a polynomial in $\cos x$. One can also bypass the even extension and the composition with the cosine function, using directly uniform approximation of trigonometric functions by Taylor polynomials, as asserted by the general theory of power series.
 \end{remark}

\vp
\noindent
{\bf \large Exercises  \thesection}
\setcounter{myquestion}{1}

\begin{question}{\themyquestion}
Consider the function $f:[-\pi, \pi]$   defined as
\[f(x)=\begin{cases} x+\pi,\quad &\text{if}\; -\pi \leq x<0,\\x-\pi,\quad &\text{if}\quad 0\leq x\leq \pi.\end{cases}\]
 The  Fourier series of this function has been obtained in Exercises \ref{sec7.2}. Does the Fourier series converge uniformly? Justify your answer.
 
\end{question}
\atc
\begin{question}{\themyquestion}
  Study the uniform convergence of the Fourier series of the function $f:[-\pi, \pi]\to\mb{R}$, $f(x)=x^2$ obtained in Exercises \ref{sec7.1}.
\end{question}

  \atc

\begin{question}{\themyquestion}
Show that the trigonometric series
\[\sum_{k=1}^{\infty}\left( \frac{2k\cos kx +3 \sin kx}{k^4}\right)\] defines a continuously differentiable function $F:\mb{R}\to\mb{R}$, and find the Fourier series of the function $F':[-\pi,\pi]\to\mb{R}$.
\end{question}

\section{Fourier Transforms} 
We have seen that the Fourier series of a function $f:[-L,L]\to\mb{C}$ defined on $[-L,L]$ is  
\[\sum_{k=-\infty}^{\infty} c_k\exp\left(\frac{i\pi kx}{L}\right),\]
where the Fourier coefficients $c_k$, $k\in\mb{Z}$ are given by
\[c_k=\frac{1}{2L}\int_{-L}^{L}f(t)\exp\left(-\frac{i\pi kt}{L}\right)dt.\]
This is also the Fourier series of the $2L$-periodic extension of the function $f$. Substitute the expression for $c_k$, we find that the Fourier series can be written as
\begin{equation}
\label{230818_1}\frac{1}{2\pi}\sum_{k=-\infty}^{\infty}\frac{\pi}{L}
\int_{-L}^{L}f(t)\exp\left( \frac{i\pi k(x-t)}{L}\right)dt.\end{equation}
Heuristically, 
\[\sum_{k=-\infty}^{\infty}\frac{\pi}{L}
 \exp\left( \frac{i\pi kt}{L}\right)\]
 can be regarded as a Riemann sum for the function $g:\mb{R}\to \mb{C}$,
 \[g(\omega)=e^{i\omega t}.\]In the limit $L\to\infty$, one obtain heuristically the integral
 \[\int_{-\infty}^{\infty}e^{i\omega t}d\omega,\]
 so that \eqref{230818_1} becomes
 \[\frac{1}{2\pi}\int_{-\infty}^{\infty}\int_{-\infty}^{\infty}f(t)e^{i\omega (x-t)}dtd\omega.\]
 This motivates us to define the Fourier transform of a function $f:\mb{R}\to\mb{C}$ as
 \[\widehat{f}(\omega)=\int_{-\infty}^{\infty}f(t)e^{-i\omega t}dt.\]
 We know that under certain conditions, the Fourier series of a function would converge to the function itself. Hence, we can also explore the conditions in which 
 \begin{equation}\label{230818_8}f(x)=\frac{1}{2\pi}\int_{-\infty}^{\infty}\int_{-\infty}^{\infty}f(t)e^{i\omega (x-t)}dtd\omega=\frac{1}{2\pi}\int_{-\infty}^{\infty}\widehat{f}(\omega)e^{i\omega x}d\omega.\end{equation}
 However, now the integrals we are working with are improper integrals. Therefore, there is another  convergence issue that we need to deal with. In this section, we only give a brief discussion about Fourier transforms. An in-depth  analysis would require advanced tools.
 
We say that a function $f:\mb{R}\to\mb{C}$ is Riemann integrable if it is Riemann integrable on any compact intervals.
 
 \begin{definition}{$\pmb{L^1}$ and $\pmb{L^2}$ Functions}
 Let $f:\mb{R}\to\mb{C}$ be a   Riemann integrable function.
 We say that 
 $f$ is $L^1$ if 
 the improper integral
 \[\int_{-\infty}^{\infty}|f(x)|dx\]is convergent. 
 In this case, we define the $L^1$-norm of $f$ as
 \[\Vert f\Vert_1=\int_{-\infty}^{\infty}|f(x)|dx.\]
 We say that $f$ is $L^2$ if 
 the improper integral
 \[\int_{-\infty}^{\infty}|f(x)|^2dx\]is convergent. 
 In this case, we define the $L^2$-norm of $f$ as
 \[\Vert f\Vert_2=\sqrt{\int_{-\infty}^{\infty}|f(x)|^2dx}.\]
 \end{definition}
 If $f:[a,b]\to\mb{C}$ is any Riemann integrable function, the zero-extension of $f$ to $\mb{R}$ is both a  $L^1$ and a $L^2$ function. As before, the $L^1$ and $L^2$ norms  are semi-norms which are positive semi-definite, where there are nonzero functions that have zero norms.

 \begin{example}{}
 Consider the function $f:\mb{R}\to\mb{R}$ defined as
 \[f(x)=\frac{1}{\sqrt{x^2+1}}.\]
 \be
 The integral 
 \[\int_{-\infty}^{\infty}\frac{1}{\sqrt{x^2+1}}dx\]
  is not convergent, but the integral
 \[\int_{-\infty}^{\infty}\frac{1}{ x^2+1 }dx\] is convergent.  Hence, $f:\mb{R}\to\mb{R}$ is $L^2$  but not $L^1$. 
 \end{example2}

   \begin{definition}{Fourier transform}
 Let $f:\mb{R}\to\mb{C}$ be a    Riemann integrable  function.    The  Fourier transform of $f$, denoted by $\mathcal{F}[f]$ or $\widehat{f}$, is defined as
  \[\mathcal{F}[f](\omega)=\widehat{f}(\omega)=\int_{-\infty}^{\infty}f(t)e^{-i\omega t}dt,  \] for all the $\omega\in \mb{R}$ which this improper integral is convergent.
 \end{definition}
 
 \begin{example}{}
 If $f:\mb{R}\to\mb{C}$ is a $L^1$-function, for any $\omega\in\mb{R}$, the integral
 \[\int_{-\infty}^{\infty}f(t)e^{-i\omega t}dt\] converges absolutely.
 Hence, a $L^1$ function has Fourier transform $\widehat{f}$ which is defined on $\mb{R}$. In particular, a function that vanishes outside a bounded interval has a  Fourier transform that is defined for all $\omega\in\mb{R}$.
 \end{example}
 \begin{proposition}{}
 Fourier transform is a linear operation. Namely, if $f:\mb{R}\to\mb{C}$ and $g:\mb{R}\to\mb{C}$ are functions that have Fourier transforms, then for any complex numbers $\alpha$ and $\beta$, the function $\alpha f+\beta g:\mb{R}\to\mb{C}$ also has Fourier transform, and
 \[\mathcal{F}[\alpha f+\beta g]=\alpha\mathcal{F}[f]+\beta\mathcal{F}[g].\]
 \end{proposition}
 \begin{remark}{}
 In engineering, it is customary to use $t$ as the independent variable for the function $f:\mb{R}\to\mb{C}$, and $\omega$ as the independent variable for its Fourier transform $\widehat{f}:\mb{R}\to\mb{C}$. The function $f$ is usually a function of time $t$, and its Fourier transform is a function of frequency $\omega$. Hence, the Fourier transform is a transform from the time domain to the frequency domain.
 \end{remark}
 
 \begin{example}[label=230818_2]{}
 Let $a$ and $b$ be two real numbers with $a<b$. Define the function $g:\mb{R}\to\mb{R}$ by
 \[g(t)=\begin{cases}1,\quad &\text{if}\;a\leq t\leq b,\\0,\quad &\text{otherwise}.\end{cases}\]
 Find the Fourier transform of $g$.
 \end{example}
 \begin{solution}{Solution}
 The Fourier transform of $g$ is
 \[\widehat{g}(\omega)=\int_a^b e^{-i\omega t} dt=\begin{cases}\di \frac{i(e^{-i\omega b}-e^{-i\omega a})}{\omega},\quad &\text{if}\;\omega\neq 0,\\
 b-a,\quad &\text{if}\;\omega=0.\end{cases}\]
 \end{solution}
 
 Of special interest is when 
 $g:\mb{R}\to\mb{R}$ is given by
 \begin{equation}\label{230818_3}g(t)=\begin{cases}1,\quad &\text{if}\;-a\leq t\leq a,\\0,\quad &\text{otherwise}, \end{cases}\end{equation}
 which is an even function. Example \ref{230818_2} shows that its Fourier transform is
 \[\widehat{g}(w)=\frac{2\sin a\omega}{\omega}.\]
 One can show that this function is not $L^1$ but is $L^2$.
 
\begin{figure}[ht]
\centering
\includegraphics[scale=0.2]{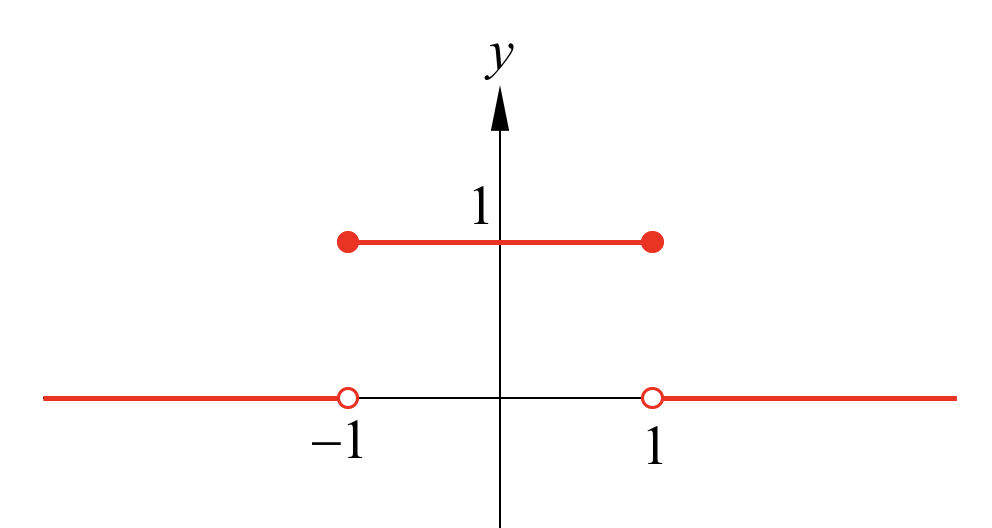}

\caption{The function $g:\mb{R}\to\mb{R}$ defined by \eqref{230818_3} with $a=1$.}\label{figure78}
\end{figure}
\begin{figure}[ht]
\centering
\includegraphics[scale=0.2]{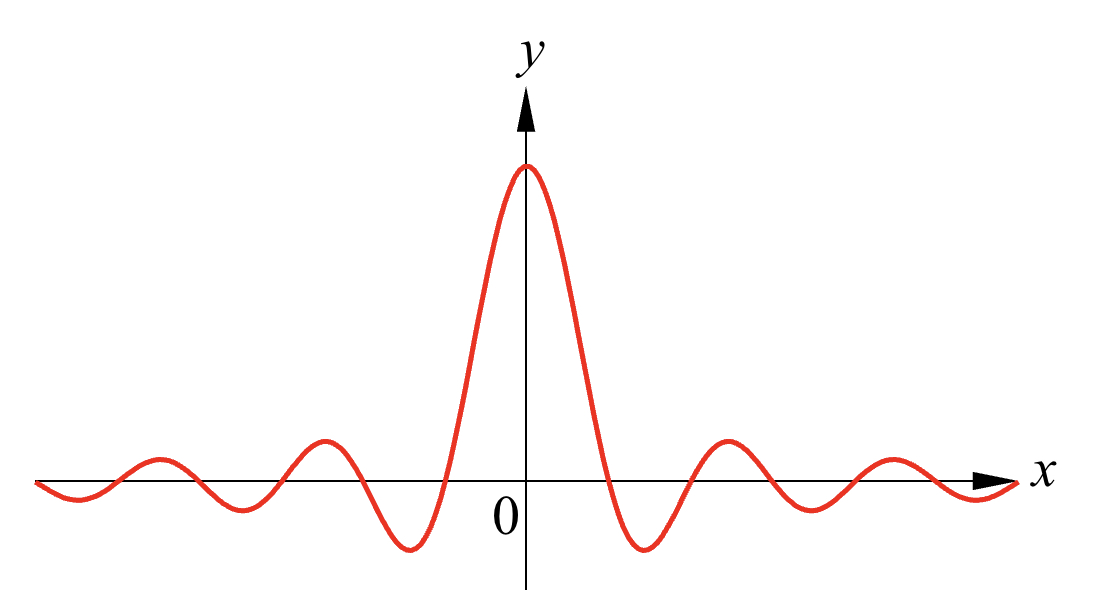}

\caption{The Fourier transform of the function $g:\mb{R}\to\mb{R}$ defined by \eqref{230818_3} with $a=1$.}\label{figure79}
\end{figure}
 \begin{remark}{}
 A function that vanishes outside a bounded interval is said to have compact support. In general, the support of a function $f:\mb{R}\to\mb{C}$ is defined to be the closure of the set of those points $x$ such that $f(x)\neq 0$. Namely,
 \[(\text{support}\; f)=\overline{\{x\in\mb{R}\,|\, f(x)\neq 0\}}.\]
 Since a set is bounded if and only if its closure is bounded, a function $f$ has compact support if and only if the set of points where $f$ does not vanish is bounded.
 \end{remark}
 
 Let us look at Fourier transforms of functions that does not have compact support.
 \begin{example}[label=230819_3]{}
 Let  $a$ be a positive number, and let $f:\mb{R}\to\mb{R}$ be the function defined as $f(t)=e^{-a|t|}$. Find the Fourier transform of $f$.
 \end{example}
 \begin{solution}
 {Solution}
 The Fourier transform of $f$ is given by
 \begin{align*}
 \widehat{f}(\omega) &=\int_{-\infty}^{\infty}e^{-a|t|}e^{-i\omega t}dt\\
 &=\lim_{L\to\infty}\int_0^{L}e^{-at}\left(e^{-i\omega t}+e^{i\omega t}\right)dt
\\&=\lim_{L\to\infty}\int_0^{L}\left(e^{-(a+i\omega)t}+e^{-(a-i\omega)t}\right)dt\\
 &=\lim_{L\to\infty}\left[-\frac{e^{-(a+i\omega)t}}{a+i\omega}-\frac{e^{-(a-i\omega)t}}{a-i\omega}\right]_0^L\\
 &= \frac{1}{a+i\omega}+\frac{1}{a-i\omega}\\
& = \frac{2a}{a^2+\omega^2}.
 \end{align*}
 \end{solution}
 Notice that the function $\widehat{f}:\mb{R}\to\mb{C}$, $\widehat{f}(\omega)= \di \frac{2a}{a^2+\omega^2}$ is a   function that is both $L^1$ and $L^2$.
 
 \begin{figure}[ht]
\centering
\includegraphics[scale=0.2]{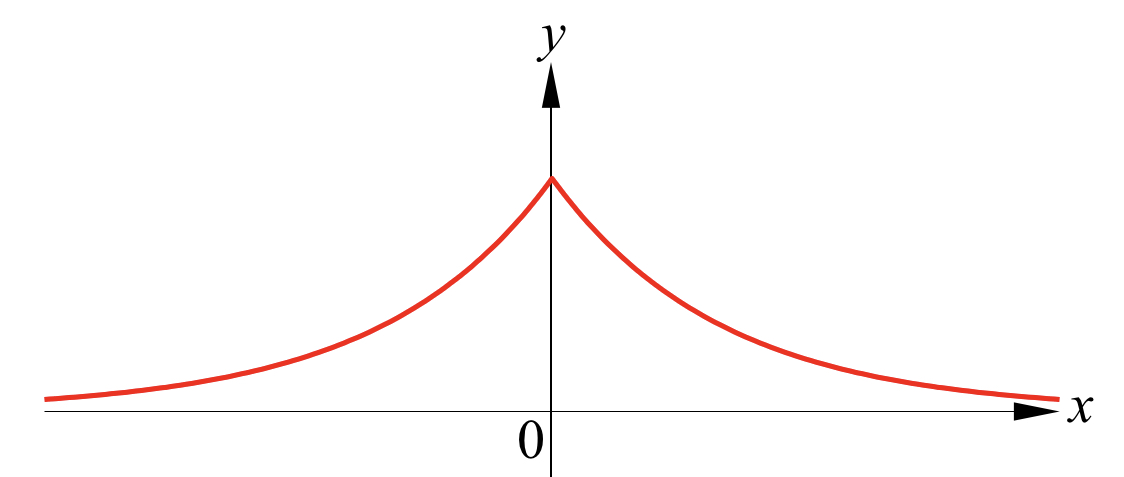}

\caption{The   function $f:\mb{R}\to\mb{R}$, $f(t)=e^{-|t|}$.}\label{figure80}
\end{figure}

\begin{figure}[ht]
\centering
\includegraphics[scale=0.2]{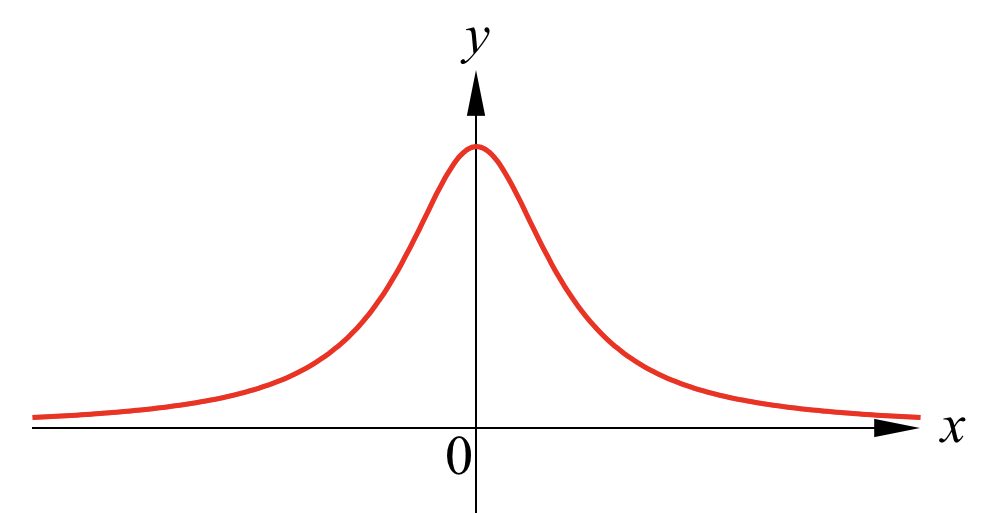}

\caption{The   function $\widehat{f}:\mb{R}\to\mb{R}$, $\widehat{f}(\omega)=\di \frac{2}{1+\omega^2}$, which is the Fourier transform of  $f:\mb{R}\to\mb{R}$, $f(t)=e^{-|t|}$.}\label{figure81}
\end{figure}

  A function $f:\mb{R}\to\mb{R}$ of the form \begin{equation}\label{230818_7}
  f(x)=c\exp\left(-\frac{(x-\mu)^2}{2\sigma^2}\right),\end{equation}
  is the probability density function of a normal distribution with mean $\mu$ and standard deviation $\sigma$ when
  \[c=\frac{1}{\sqrt{2\pi}\sigma}.\] 
   It  is also known as a {\it Gaussian function}. These functions are infinitely differentiable and they decay exponentially to 0 when $x$ gets large. 
   When $\mu=0$ and $\sigma=1$, 
   \[\Phi(x)=\frac{1}{\sqrt{2\pi}}\exp\left(-\frac{x^2}{2}\right)\]
   is the probability density of the standard normal distribution. 
   The Fourier transform of the Gaussian function $f(t)=\di \exp\left(-\frac{t^2}{2}\right)$   is
 \begin{align*}
 \widehat{f}(\omega) &= \int_{-\infty}^{\infty}e^{-\frac{t^2}{2}}e^{-i\omega t}dt\\
 &= e^{-\frac{\omega^2}{2}}\int_{-\infty}^{\infty}
 \exp\left(-\frac{1}{2}\left(t+ i\omega \right)^2\right)dt\\
  &=e^{-\frac{\omega^2}{2}}\int_{-\infty}^{\infty}e^{-\frac{t^2}{2}}dt\\
  &=\sqrt{2\pi }e^{-\frac{\omega^2}{2}}.
  \end{align*}
  In the computation, the equality
 \[ \int_{-\infty}^{\infty}
 \exp\left(-\frac{1}{2}\left(t+ i\omega \right)^2\right)dt=\int_{-\infty}^{\infty}e^{-\frac{t^2}{2}}dt\] can be understood in complex analysis as shifting contours of integrations. We leave the details to the students.
 
 Notice that for  the function  $f(t)=\di \exp\left(-\frac{t^2}{2}\right)$, its Fourier transform $\widehat{f}(\omega)$ is 
  equal to $f(\omega)$ multiplied by $\sqrt{2\pi}$.  Namely,
  \[\widehat{f}(\omega)=\sqrt{2\pi}f(\omega).\] 
  The factor $\sqrt{2\pi}$ here is due to our normalization.
  Different textbooks use different conventions for Fourier transforms. Among them are the followings:
  \begin{gather*}
  \int_{-\infty}^{\infty}f(t)e^{-i\omega t}dt,\hspace{1cm} \int_{-\infty}^{\infty}f(t)e^{i\omega t}dt,\\
 \frac{1}{\sqrt{2\pi}}  \int_{-\infty}^{\infty}f(t)e^{-i\omega t}dt,\hspace{1cm} \frac{1}{\sqrt{2\pi}}\int_{-\infty}^{\infty}f(t)e^{i\omega t}dt,\\
\frac{1}{2\pi}  \int_{-\infty}^{\infty}f(t)e^{-i\omega t}dt,\hspace{1cm} \frac{1}{2\pi}\int_{-\infty}^{\infty}f(t)e^{i\omega t}dt.
  \end{gather*}Some might also replace $i\omega t$ by $2\pi i \omega t$. When one is reading about Fourier transforms, it is important to check the definition of Fourier transform that is being used.
  
  One can show that in our definition,  the Fourier transform of the Gaussian function $f(t)=e^{-at^2}$ with $a>0$ is
  \[\widehat{f}(\omega)=\int_{-\infty}^{\infty}e^{-at^2}e^{-i\omega t}dt=\sqrt{\frac{\pi}{a}}e^{-\frac{\omega^2}{4a}},\]so that
  \[\frac{1}{\sqrt{2\pi}}\widehat{f}(\omega)=\frac{1}{\sqrt{2a}}e^{-\frac{\omega^2}{4a}}.\]
  
  \begin{figure}[ht]
\centering
\includegraphics[scale=0.2]{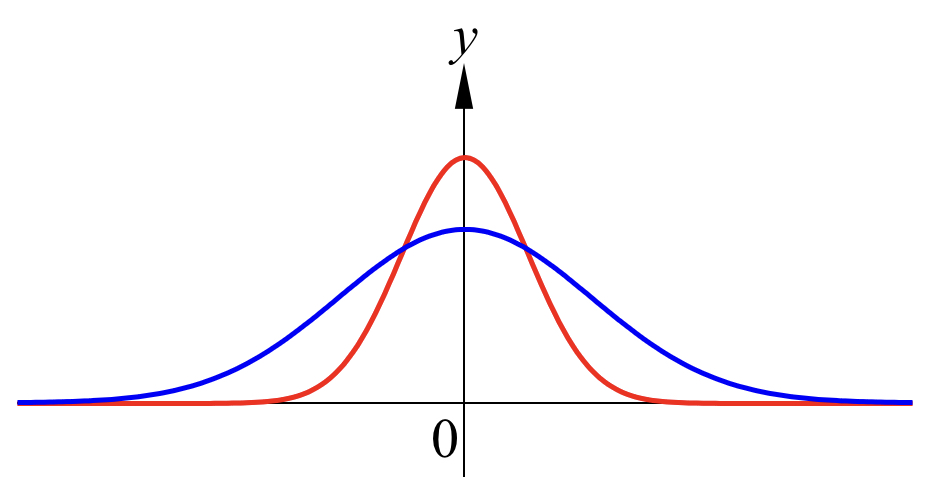}

\caption{The   function  $f(t)=e^{-t^2}$ and the function $g(\omega)=\di\frac{1}{\sqrt{2\pi}}\widehat{f}(\omega)=\frac{1}{\sqrt{2}}e^{-\frac{\omega^2}{4}}$.}\label{figure82}
\end{figure}
\begin{remark}{}
If   $f:\mb{R}\to\mb{C}$ is a $L^1$ function and 
\[\Vert f\Vert_1=\int_{-\infty}^{\infty}|f(t)|dt=0,\]
we say that $f$ is $L^1$-equivalent to the zero function. If  $f:\mb{R}\to\mb{C}$ is $L^1$-equivalent to the zero function, then for any $\omega\in\mb{R}$,
\[\left|\widehat{f}(\omega)\right|=\left|\int_{-\infty}^{\infty}f(t)e^{-i\omega t}dt\right|
\leq \int_{-\infty}^{\infty}|f(t)e^{-i\omega t}|dt=\int_{-\infty}^{\infty}|f(t)|dt=0.\]
Thus, the Fourier transform of $f$ is identically zero.
\end{remark}

Example \ref{230818_2} shows that the Fourier transform of a $L^1$ function is not necessary $L^1$. Nevertheless, we have the following, which is an extension of the Riemann-Lebesgue lemma to $L^1$ functions on $\mb{R}$.
\begin{theorem}{Extended Riemann-Lebesgue Lemma}
If the function $f:\mb{R}\to\mb{C}$ is $L^1$, then 
\[\lim_{\beta\to\infty}\int_{-\infty}^{\infty}f(t)e^{i\beta t}dt=0.\]
In other words, the Fourier transform $\widehat{f}:\mb{R}\to\mb{C}$ of $f$ is a function satisfying
\[\lim_{\omega\to \pm\infty} \widehat{f}(\omega)=0.\]
\end{theorem}
\begin{myproof}{Proof}
We are given that \[J=\int_{-\infty}^{\infty}|f(t)|dt <\infty.\]
Given   $\varepsilon>0$, there is a $L>0$ such that
\[\int_{|t|\geq L}|f(t)|dt<\frac{\varepsilon}{2}.\]
\bp
By triangle inequality, we have
\[\left|\int_{-\infty}^{\infty}f(t)e^{i\beta t}dt\right|\leq \left|\int_{-L}^{L}f(t)e^{i\beta t}dt\right|+\left|\int_{|t|\geq L}f(t)e^{i\beta t}dt\right|.\]

For the second term, we have
\[\left|\int_{|t|\geq L}f(t)e^{i\beta t}dt\right|\leq \int_{|t|\geq L}|f(t)e^{i\beta t}|dt=\int_{|t|\geq L}|f(t)|dt<\frac{\varepsilon}{2}.\]
By the Riemann-Lebesgue lemma,
\[\lim_{\beta\to 0}\int_{-L}^{L}f(t)e^{i\beta t}dt=0.\]Therefore, there exists $M>0$ such that if $\beta>M$, then
\[\left|\int_{-L}^{L}f(t)e^{i\beta t}dt\right|<\frac{\varepsilon}{2}.\]
It follows that for all $\beta>M$,
\begin{align*}
\left|\int_{-\infty}^{\infty}f(t)e^{i\beta t}dt\right|<\varepsilon.
\end{align*}This proves the assertion.
\end{myproof}

The following theorem imposes a strong condition on a function $g:\mb{R}\to\mb{C}$ to be the Fourier transform of a $L^1$ function $f:\mb{R}\to\mb{C}$.
\begin{theorem}[label=230819_2]{}
If $f:\mb{R}\to\mb{C}$ is a $L^1$ function, then its Fourier transform $\widehat{f}:\mb{R}\to\mb{C}$ is  uniformly continuous.
\end{theorem}
\begin{myproof}{Proof}

We are given that \[J=\int_{-\infty}^{\infty}|f(t)|dt <\infty.\]
\bp
Without loss of generality, we can assume that $J>0$.
Notice that for any $\omega_1$ and $\omega_2$ in $\mb{R}$,
\[\widehat{f}(\omega_1)-\widehat{f}(\omega_2)=\int_{-\infty}^{\infty}f(t)\left(e^{-i\omega_1 t}-e^{-i\omega_2t}\right)dt.\]
Given   $\varepsilon>0$, there is a $L>0$ such that
\[\int_{|t|\geq L}|f(t)|dt<\frac{\varepsilon}{3}.\]
By triangle inequality, we have
\begin{align*}
&\left|\widehat{f}(\omega_1)-\widehat{f}(\omega_2)\right|\leq \int_{-\infty}^{\infty}|f(t)|\left|e^{-i\omega_1 t}-e^{-i\omega_2t}\right|dt\\
&=\int_{|t|\leq L}|f(t)|\left|e^{-i\omega_1 t}-e^{-i\omega_2t}\right|dt+
\int_{|t|\geq L}|f(t)|\left|e^{-i\omega_1 t}-e^{-i\omega_2t}\right|dt.
\end{align*}
The second term is easy to estimate since $\left|e^{-i\omega_1 t}-e^{-i\omega_2t}\right|\leq 2$. We have
\[\int_{|t|\geq L}|f(t)|\left|e^{-i\omega_1 t}-e^{-i\omega_2t}\right|dt\leq 
2\int_{|t|\geq L}|f(t)|dt<\frac{2\varepsilon}{3}.\]

Since the function $g:\mb{R}\to \mb{C}$, $g(u)=e^{iu}$ is continuous at $u=0$, there exists a $\delta>0$ such that if $|u|<\delta$, then
\[|e^{iu}-1|<\frac{\varepsilon}{3J}.\]
Thus, given $\omega_1$ and $\omega_2$ in $\mb{R}$, if $|\omega_1-\omega_2|<\di\frac{\delta}{L}$, then for any $t\in [-L, L]$, 
\[|(\omega_1-\omega_2)t|\leq L|\omega_1-\omega_2|<\delta.\]
It follows that
\[|e^{-i\omega_1 t}-e^{-i\omega_2t}|=|e^{i(\omega_1-\omega_2)t}-1|<\frac{\varepsilon}{3J}.\]
Therefore,
\[ \int_{|t|\leq L}|f(t)|\left|e^{-i\omega_1 t}-e^{-i\omega_2t}\right|dt  \leq \frac{\varepsilon}{3J}\int_{|t|\leq L}|f(t)|dt\leq \frac{\varepsilon}{3}.\]
\bp
This proves that whenever $|\omega_1-\omega_2|< \di\frac{\delta}{L}$, then
\[\left|\widehat{f}(\omega_1)-\widehat{f}(\omega_2)\right|<\varepsilon.\]
Hence, $\widehat{f}:\mb{R}\to\mb{C}$ is uniformly continuous.
\end{myproof}

\begin{example}{}
By Example \ref{230818_2}, the Fourier transform of the $L^1$ function $f:\mb{R}\to\mb{C}$, 
\[f(t)=\begin{cases}1,\quad &\text{if}\; -1\leq t\leq 1,\\0,\quad &\quad  \text{otherwise},\end{cases}\]
is $\di \widehat{f}(\omega)=\frac{2\sin\omega}{\omega}$. 
Theorem \ref{230819_2} then implies that the function $g:\mb{R}\to\mb{R}$, 
\[g(x)=\frac{\sin x}{x}\]
is uniformly continuous.
\end{example}
  
  Motivated by the heuristics  \eqref{230818_8} from the theory of Fourier series, we make  the following  definition.
  \begin{definition}{Inverse Fourier Transform}
  Given a Riemann integrable function $\widehat{f}:\mb{R}\to\mb{C}$, we define its inverse Fourier transform by
  \[\mathcal{F}^{-1}[\widehat{f}](t)=\frac{1}{2\pi}\int_{-\infty}^{\infty} \widehat{f}(\omega)e^{i\omega t}d\omega\]for all the $t\in\mb{R}$ where this integral is convergent.
  \end{definition}
   Notice that if $\widehat{f}:\mb{R}\to\mb{C}$ is a $L^1$ function, then  $\mathcal{F}^{-1}[\widehat{f}](t) $ exists for all $t\in\mb{R}$, and 
 \[ \mathcal{F}^{-1}[\widehat{f}](t)=\frac{1}{2\pi} \mathcal{F}[\widehat{f}](-t).\]
 
 In other words, we have the following.
  
 \begin{proposition}{}
 Let $f:\mb{R}\to\mb{C}$  be a Riemann integrable function. If $f:\mb{R}\to\mb{C}$ has Fourier transform $\widehat{f}:\mb{R}\to\mb{C}$, and the function $\widehat{f}:\mb{R}\to\mb{C}$ has inverse Fourier transform given by $h:\mb{R}\to\mb{C}$, then the function $g:\mb{R}\to\mb{C}$,  $g(t)=\widehat{f}(t)$   has a Fourier transform given by
 \[\widehat{g}(\omega)= 2\pi h(-\omega).\]
 \end{proposition}
   \begin{example}[label=230819_1]{}
  For the function $f:\mb{R}\to\mb{C}$, $f(t)=e^{-at^2}$, its Fourier transform is $\widehat{f}:\mb{R}\to\mb{C}$,
  \[\hat{f}(\omega)=\sqrt{\frac{\pi}{a}}e^{-\frac{\omega^2}{4a}}.\]

  Therefore,
  \[\mathcal{F}^{-1}[\widehat{f}](t)=\frac{1}{2\sqrt{\pi a}}\int_{-\infty}^{\infty} e^{-\frac{\omega^2}{4a}}e^{i\omega t}d\omega =e^{-at^2}.\]
 For $g(t)=\di \sqrt{\frac{\pi}{a}}e^{-\frac{t^2}{4a}}$, 
 \[\widehat{g}(\omega)=\sqrt{\frac{\pi}{a}}\times \sqrt{4\pi a}\exp\left(-\frac{\omega^2}{\di 4\times\frac{1}{4a}}\right)=2\pi e^{-a\omega^2}=2\pi\mathcal{F}^{-1}[\widehat{f}](-\omega).\]
  \end{example}
  
Example \ref{230819_1} shows that for the function  $f:\mb{R}\to\mb{C}$, $f(t)=e^{-at^2}$, we have
   \[ \mathcal{F}^{-1}[\widehat{f}](t)=f(t).\]
  In general, we are interested in the following. If $f:\mb{R}\to\mb{C}$ is a Riemann integrable function with Fourier transform $\widehat{f}:\mb{R}\to\mb{C}$, under what conditions does $\mathcal{F}^{-1}[\widehat{f}](t)$ exists and 
 \begin{equation}\label{230818_10} \mathcal{F}^{-1}[\widehat{f}](t)=f(t)?\end{equation}

 \begin{example}[label=230820_12]{}
 For the function $g:\mb{R}\to\mb{C}$, 
 \[g(\omega)=\frac{2a}{a^2+\omega^2}\hspace{1cm}\text{with}\;a>0,\]
 one can use contour integration techniques in complex analysis to show that when $t\in\mb{R}$,
 \[\mathcal{F}^{-1}[g](t)=\frac{1}{2\pi}\int_{-\infty}^{\infty} \frac{2a e^{i\omega t}}{a^2+\omega^2}d\omega=e^{-a|t|}.\]Hence, for the function $f:\mb{R}\to\mb{C}$, $f(t)=e^{-a|t|}$, we also have
 \[ \mathcal{F}^{-1}[\widehat{f}](t)=f(t)\hspace{1cm}\text{for all}\;t\in\mb{R}.\]
 \end{example}
 
 \begin{definition}{Fourier Transform Pairs}
 If $f:\mb{R}\to\mb{C}$ and $g:\mb{R}\to\mb{C}$ are Riemann integrable functions, and 
 \[\mathcal{F}[f](\omega)=g(\omega),\hspace{1cm}\mathcal{F}^{-1}[g](t)=f(t),\]
 then we call the pair of functions $(f,g)$ a Fourier transform pair.
 \end{definition}
 
 \begin{example}{}
 For $a>0$, let $f:\mb{R}\to\mb{C}$ be the function $f(t)=e^{-a|t|}$, and let $g:\mb{R}\to\mb{C}$ be the function $\di g(t)=\frac{2a}{a^2+t^2}$. Then Example \ref{230820_12}  says that $(f,g)$ is a Fourier transform pair.
 \end{example}

 The following is important for the proofs later. 
 \begin{theorem}[label=230819_7]{}
The function $f:\mb{R}\to\mb{R}$,  $f(x)=\di\frac{\sin  x}{x}$ is an infinitely differentiable even function that satisfies 
 \[ \int_{0}^{\infty}f(x)dx= \int_{0}^{\infty}\frac{\sin  x}{x}dx =\frac{\pi}{2}.\]
 \end{theorem}The fact that  $f:\mb{R}\to\mb{R}$,  $f(x)=\di\frac{\sin  x}{x}$ is infinitely differentiable has been established earlier. 
 The formula for the improper integral can be proved using contour integration techniques and the fact that the function $g(z)=\di \frac{e^{iz}}{z}$ has a simple pole at $z=0$ with residue $1$. See for example \cite{Churchill}.
 
 \begin{corollary}[label=230819_8]{}
 For any $a>0$, 
 \[\lim_{L\to\infty}\int_0^a\frac{\sin Lx}{x}dx=\frac{\pi}{2}.\]
 
 \end{corollary}
 
 \begin{myproof}{Proof}
 Making a change of variables, we have
 \[
\int_0^a\frac{\sin Lx}{x}dx= \int_0^{aL}\frac{\sin x }{x}dx.
 \]Therefore,
 \[\lim_{L\to\infty}\int_0^a\frac{\sin Lx}{x}dx=\lim_{L\to\infty} \int_0^{aL}\frac{\sin x }{x}dx=\int_0^{\infty}\frac{\sin x}{x}dx=\frac{\pi}{2}.\]
 \end{myproof}

 Now we can prove our main theorem.  A function  $f:\mb{R}\to\mb{C}$  is said to be  strongly piecewise differentiable if it is strongly piecewise differentiable on any compact intervals.
 \begin{theorem}[label=230819_5]{Fourier Inversion Theorem}
 Let $f:\mb{R}\to\mb{C}$ be a $L^1$-function that is strongly piecewise differentiable, and let 
 \[\widehat{f}(\omega)=\int_{-\infty}^{\infty}f(t)e^{-i\omega t}dt\]
 be its Fourier transform. Then for any $x\in\mb{R}$,
 \[\lim_{L\to\infty}\frac{1}{2\pi}\int_{-L}^{L}\widehat{f}(\omega)e^{i\omega x}d\omega=\frac{f_+(x)+f_-(x)}{2}.\]
 \end{theorem}
 \begin{myproof}{Proof}
 Notice that

 \begin{align*}\int_{-L}^{L}\widehat{f}(\omega)e^{i\omega x}d\omega&= \int_{-L}^{L}\int_{-\infty}^{\infty}f(t)e^{i\omega (x-t)}dtd\omega\\
 &= \int_{-L}^{L}\int_{-\infty}^{\infty}f(x-t)e^{i \omega t}dtd\omega.\end{align*}
 \end{myproof}
 To continue, we need a technical lemma which guarantees we can interchange the order of integrations.
 \begin{lemma}[label=230819_6]{}
 Let $f:\mb{R}\to\mb{C}$ be a function that satisfies the conditions in Theorem \ref{230819_5}. Then for any $L>0$, we have
 \[\int_{-L}^{L}\int_{-\infty}^{\infty}f(x-t)e^{i\omega t}dtd\omega=\int_{-\infty}^{\infty}\int_{-L}^{L}f(x-t)e^{i\omega t}d\omega dt.\]
 \end{lemma}
  Assuming this lemma, we can continue with the proof of Theorem \ref{230819_5}.
  \begin{myproof}{\linkt Proof of Theorem \ref{230819_5} Continued}
By Lemma \ref{230819_6}, we have
\[\int_{-L}^{L}\widehat{f}(\omega)e^{i\omega x}d\omega= \int_{-\infty}^{\infty}\int_{-L}^{L}f(x-t)e^{i\omega t}d\omega dt.\]
Now we can integrate the integral with respect to $\omega$ and obtain
\[ \int_{-L}^{L}\widehat{f}(\omega)e^{i\omega x}d\omega= 2\int_{-\infty}^{\infty} f(x-t)\frac{\sin Lt }{t}dt.\]
Using the fact that $\di\frac{\sin Lt}{t}$ is an even function, we find that
\[\int_{-L}^{L}\widehat{f}(\omega)e^{i\omega x}d\omega=2\int_{0}^{\infty} \frac{f(x+t)+f(x-t)}{t}\sin Lt dt.\]
\bp
Split the integral into two parts, we have
\begin{align*}
\int_{-L}^{L}\widehat{f}(\omega)e^{i\omega x}d\omega 
&=2\int_{0}^{1} \frac{f(x+t)+f(x-t)}{t}\sin Lt dt\\&\quad +2\int_{1}^{\infty} \frac{f(x+t)+f(x-t)}{t}\sin Lt dt.
\end{align*}
Let
\[u=\frac{f_+(x)+f_-(x)}{2}.\]
As in the proof of Lemma \ref{230820_4},
 the function $h:[0,1]\to\mb{C}$ with
\[h(t)=\frac{f(x+t)+f(x-t)-2u}{t}\hspace{1cm}\text{when}\;t\in (0,1]\]
is a Riemann integrable function. Thererfore, the Riemann-Lebesgue lemma implies that
\[\lim_{L\to\infty}\int_0^1 h(t)\sin Lt dt=0.\]
It follows from Corollary \ref{230819_8} that
\[\lim_{L\to\infty} 2\int_{0}^{1} \frac{f(x+t)+f(x-t)}{t}\sin Lt dt=\lim_{L\to\infty} 4u\int_{0}^{1} \frac{\sin Lt}{t} dt=2\pi u.\]
On the other hand, 
\[\int_1^{\infty}\left| \frac{f(x+t)+f(x-t)}{t}\right|dt\leq 2 \int_{-\infty}^{\infty}|f(t)|dt<\infty.\] By the extended Riemann-Lebesgue lemma,
\begin{align*}
\lim_{L\to\infty}2\int_{1}^{\infty} \frac{f(x+t)+f(x-t)}{t}\sin Lt dt=0.
\end{align*}
This completes the proof that
\[\lim_{L\to\infty}\frac{1}{2\pi}\int_{-L}^{L}\widehat{f}(\omega)e^{i\omega x}d\omega=u=\frac{f_+(x)+f_-(x)}{2}.\]
\end{myproof}

Now we prove Lemma \ref{230819_6}.
\begin{myproof}{\linkt Proof of Lemma \ref{230819_6}}

Given $\varepsilon>0$, since
\[ \int_{-\infty}^{\infty}|f(t)|dt <\infty,\]
there is an $M>0$ such that
\[\int_{|t|\geq M}|f(t)|dt<\frac{\varepsilon}{ 4L }.\]
Since $e^{i\omega t} $ is an infinitely differentiable function, and $f(t)$ is a piecewice continuous function on any compact intervals, Fubini's theorem implies that
\[\int_{-L}^L\int_{x-M}^{x+M}f(x-t)e^{i \omega t} dtd\omega=\int_{x-M}^{x+M}\int_{-L}^Lf(x-t)e^{i \omega t} d\omega dt.\]
Now,
\begin{align*}
&\left|\int_{-L}^L\int_{-\infty}^{\infty}f(x-t)e^{i \omega t }dtd\omega-\int_{-L}^L\int_{x-M}^{x+M}f(x-t)e^{i \omega t} dtd\omega\right|
\\&\leq \int_{-L}^L\int_{|x-t|\geq M}|f(x-t)|dt d\omega \leq 2 L\int_{|t|\geq M}|f(t)|dt<\frac{  \varepsilon}{2}.\end{align*}
On the other hand, since $|\sin Lt|\leq L|t|$ for all $t\in\mb{R}$, we have
\begin{align*}
&\left|\int_{-\infty}^{\infty}\int_{-L}^Lf(x-t)e^{i \omega t } d\omega dt-\int_{x-M}^{x+M}\int_{-L}^Lf(x-t)e^{i\omega t } d\omega dt\right|
\\&\leq 2\int_{|x-t|\geq M}|f(x-t)| \left|\frac{\sin Lt}{t}\right| dt \leq2L\int_{|t|\geq M}|f(t)|dt<\frac{ \varepsilon}{2}.\end{align*}
This proves that
\[\left|\int_{-L}^{L}\int_{-\infty}^{\infty}f(x-t)e^{i \omega t}dtd\omega-\int_{-\infty}^{\infty}\int_{-L}^{L}f(x-t)e^{i \omega t}d\omega dt\right|<\varepsilon.\]Since $\varepsilon>0$ is arbitrary, the assertion follows.
\end{myproof}
\begin{corollary}{}
 Let $f:\mb{R}\to\mb{C}$ be a $L^1$ function that is continuous and strongly piecewise differentiable, and let 
 \[\widehat{f}(\omega)=\int_{-\infty}^{\infty}f(t)e^{-i\omega t}dt\]
 be its Fourier transform. If $\widehat{f}:\mb{R}\to\mb{C}$ is also a $L^1$ function, then for any $t\in \mb{R}$,
 \[\mathcal{F}^{-1}[\widehat{f}](t)=f(t).\]
\end{corollary}

\begin{example}{}
Since the function $g:\mb{R}\to\mb{R}$, 
\[g(t)=\begin{cases}1,\quad &\text{if}\;-a\leq t\leq a,\\0,\quad &\quad \text{otherwise}, \end{cases}\] is strongly piecewise differentiable $L^1$ function with Fourier transform
\[\widehat{g}(\omega)=\frac{2\sin a\omega}{\omega},\]
  the Fourier inversion theorem implies that for $|t|<a$,
  \[\lim_{L\to\infty}\frac{1}{ \pi}\int_{-L}^{L}  \frac{\sin a\omega}{\omega}e^{i\omega t}d\omega =1,\]
 
  while if $|t|>a$,
  \[\lim_{L\to\infty}\frac{1}{ \pi}\int_{-L}^{L}\frac{\sin a\omega}{\omega}e^{i\omega t}d\omega =0,\]and for $|t|=a$,
  \[\lim_{L\to\infty}\frac{1}{ \pi}\int_{-L}^{L} \frac{\sin a\omega}{\omega}e^{i\omega t}d\omega =\frac{1}{2}.\]
\end{example}

If $f:\mb{R}\to\mb{C}$ and $g:\mb{R}\to\mb{C}$ are $L^2$ functions, then the Cauchy Schwarz inequality implies that for any $L>0$,
\begin{align*}
\left(\int_{-L}^{L}f(t)\overline{g(t)}dt\right)^2&\leq \left(\int_{-L}^L|f(t)|^2dt \right)\left(\int_{-L}^L|g(t)|^2dt \right)\\&\leq \left(\int_{-\infty}^{\infty}|f(t)|^2dt \right)\left(\int_{-\infty}^{\infty}|g(t)|^2dt \right).\end{align*}
This implies that the improper integral 
\[\int_{-\infty}^{\infty}f(t)\overline{g(t)}dt\] converges absolutely. Thus, we can define a positive semi-definite inner product on the space of $L^2$ functions on $\mb{R}$ by
\[\langle f, g\rangle=\int_{-\infty}^{\infty}f(t)\overline{g(t)}dt.\]
The $L^2$ semi-norm is the norm induced by this inner product.

The following is a generalization of the Parseval's identity to Fourier transforms.
\begin{theorem}{Parseval- Plancherel Identity}
If $f:\mb{R}\to\mb{C}$ is a   Riemann integrable function that is both $L^1$ and $L^2$, then its Fourier transform $\widehat{f}:\mb{R}\to\mb{C}$ is  a $L^2$ function. Moreover,
\[\Vert f\Vert_2^2=\di \frac{1}{2\pi}\Vert \widehat{f}\Vert_2^2.\]
Namely,
\begin{equation}\label{230819_10}\int_{-\infty}^{\infty}|f(t)|^2dt=\frac{1}{2\pi}\int_{-\infty}^{\infty}|\widehat{f}(\omega)|^2d\omega.\end{equation}
\end{theorem}
\begin{myproof}{Sketch of Proof}
A rigorous proof of this theorem requires advanced tools in analysis. We give a heuristic argument for the validity of the  formula \eqref{230819_10}  under the additional assumption that $f:\mb{R}\to\mb{C}$ is continuous and  strongly piecewise differentiable, and $\widehat{f}:\mb{R}\to\mb{C}$ is also $L^1$.
\bp
Since $\widehat{f}:\mb{R}\to\mb{C}$ is a continuous and strongly piecewise differentiable  $L^1$ function, the Fourier inversion theorem implies that for all $t\in \mb{R}$,
\[f(t)=\mathcal{F}^{-1}[\widehat{f}](t)=\frac{1}{2\pi}\int_{-\infty}^{\infty}\widehat{f}(\omega)e^{i\omega t}dt.\]  
 Notice that
$\widehat{f}:\mb{R}\to\mb{C}$ is a $L^2$-function if  the limit
\[\lim_{L\to\infty}\int_{-L}^L\left|\widehat{f}(\omega)\right|^2d\omega \]exists. By the definition of Fourier transform,
\[\int_{-L}^L\left|\widehat{f}(\omega)\right|^2d\omega =\int_{-L}^L\overline{\widehat{f}(\omega)}\int_{-\infty}^{\infty}f(t)e^{-it\omega}dtd\omega.\]
By Theorem \ref{230819_2}, $\widehat{f}(\omega)$ is uniformly continuous. By the Riemann-Lebesgue lemma, $\di\lim_{\omega\to\pm\infty}\widehat{f}(\omega)=0$. These imply that the function $\widehat{f}:\mb{R}\to\mb{C}$ is bounded. Using the same reasoning as in the proof of Lemma \ref{230819_6}, we can interchange the order of integrations and obtain
\[\int_{-L}^L\left|\widehat{f}(\omega)\right|^2d\omega =\int_{-\infty}^{\infty}f(t)\int_{-L}^L\overline{\widehat{f}(\omega)}e^{-it\omega}d\omega dt.\]

Since $\widehat{f}:\mb{R}\to\mb{C}$ is   $L^1$, we can take the $L\to\infty$ limit under the integral sign.  Since
\[\lim_{L\to\infty}\int_{-L}^L\overline{\widehat{f}(\omega)}e^{-it\omega}d\omega=\overline{\int_{-\infty}^{\infty}\widehat{f}(\omega)e^{i\omega t}d\omega},\]
we conclude that
\[\lim_{L\to\infty}\int_{-L}^L\left|\widehat{f}(\omega)\right|^2d\omega =
2\pi\int_{-\infty}^{\infty}|f(t)|^2dt.\]
\end{myproof}

\begin{example}{}
For the function $f:\mb{R}\to\mb{C}$, $f(t)=e^{-a|t|}$ with $a>0$, its Fourier transform is $\widehat{f}:\mb{R}\to\mb{C}$, $\widehat{f}(\omega)=\di \frac{2a}{a^2+\omega^2}$.
Notice that
\[\int_{-\infty}^{\infty}|f(t)|^2 dt=2\int_0^{\infty}e^{-2at}dt=\frac{1}{a}.\]
The Parseval-Plancherel formula implies that
\[\int_{-\infty}^{\infty}\frac{1}{(a^2+\omega^2)^2}d\omega=\frac{1}{4a^2}\int_{-\infty}^{\infty} \left|\widehat{f}(\omega)\right|^2d\omega =\frac{\pi}{2a^2} \int_{-\infty}^{\infty}|f(t)|^2 dt=\frac{\pi}{2a^3}.\]
\end{example}

 One of the applications of Fourier transforms is to solve differential equations. For this we need the following.

\begin{theorem}[label=230820_2]{}
Let $f:\mb{R}\to\mb{C}$ be a continuously differentiable $L^1$ function such that 
\[\lim_{t\to\pm \infty}f(t)=0,\]
and
its derivative $f':\mb{R}\to\mb{C}$ is a Riemann integrable function that has Fourier transform. Then
\[\mathcal{F}[ f'](\omega)=i\omega \mathcal{F}[ f](\omega).\]
\end{theorem}
\begin{myproof}{Proof}
This follows from integration by parts.   For any $a$ and $b$ with $a<b$,
\begin{align*}
\int_{a}^{b}f'(t)e^{-it\omega} dt=\left[f(t)e^{-it\omega}\right]_{a}^b+i\omega\int_{a}^b f(t)e^{-it\omega}dt.
\end{align*}The assertion follows by taking the limit $a\to -\infty$ and $b\to\infty$.
\end{myproof}
\begin{example}{}
Find the Fourier transform of the function $f:\mb{R}\to\mb{C}$, $f(t)=t^2e^{-t^2}$.

\end{example}
\begin{solution}{Solution}
Let $g:\mb{R}\to\mb{C}$ be the function
\[g(t)=e^{-t^2}.\]
Then 
\[g'(t)=-2te^{-t^2},\hspace{1cm} g''(t)=-2e^{-t^2}+4t^2e^{-t^2}=-2g(t)+4f(t).\]
 
Since \[\lim_{t\to\pm\infty}g(t)=0,\hspace{1cm}\di\lim_{t\to\pm\infty}g'(t)=0,\] and
\[\mathcal{F}[g](\omega)=\sqrt{\pi}e^{-\frac{\omega^2}{4}},\]
  Theorem \ref{230820_2} implies that
  \[\mathcal{F}[g''](\omega)=-\sqrt{\pi}\omega^2e^{-\frac{\omega^2}{4}}.\]
  By linearity,
  \[\mathcal{F}[g''] =-2\mathcal{F}[g]+4\mathcal{F}[f].\]
 Therefore,
  \[\mathcal{F}[f](\omega)=  \sqrt{\pi}\left(\frac{1}{2}-\frac{\omega^2}{4}\right)e^{-\frac{\omega^2}{4}}.\]

\end{solution}

In the following, we consider an operation called convolution. 
\begin{definition}{Convolution}
Let $f:\mb{R}\to\mb{C}$ and $g:\mb{R}\to\mb{C}$ be Riemann integrable functions. The convolution of $f$ and $g$ is the function  defined as
\[(f*g)(x)=\int_{-\infty}^{\infty}f(x-t)g(t)dt=\int_{-\infty}^{\infty}f(t)g(x-t)dt,\] whenever this integral is convergent.
\end{definition}
Notice that the improper integral defining   $f*g$ is convergent for any $x$ in $\mb{R}$ when $f$ and $g$ are $L^2$ functions.
Convolutions can be defined for a wider class of functions. For example, if the supports of the functions $f$ and $g$ are both contained in $[0, \infty)$, then the integral is only nonzero when $0\leq t\leq x$. This gives
\[(f*g)(x)=\int_0^x f(t)g(x-t)dt,\] which is also well-defined for any $x\in\mb{R}$. In fact, this is the convolution one sees in the theory of Laplace transforms.

\begin{example}[label=230820_6]{}
Let $f:\mb{R}\to\mb{R}$ be the function defined as
\begin{equation}\label{230820_7}f(x)=\begin{cases}1,\quad &\text{if}\;0\leq x\leq 1,\\0,\quad &\;\;\text{otherwise},\end{cases}\end{equation} and let $g:\mb{R}\to\mb{R}$ be the function $g(x)=x$. Then
\[(f*g)(x)=\int_{-\infty}^{\infty}f(t)g(x-t)dt=\int_0^1 (x-t)dt=x-\frac{1}{2}.\]
For $f*f$, we have
\[(f*f)(x)=\int_0^1 f(x-t)dt=\int_{x-1}^xf(t) dt=\begin{cases} 0,\quad &\text{if}\;\quad x<0,\\
x, &\text{if}\;0\leq x<1,\\
2-x,\quad &\text{if} \;1\leq x\leq 2,\\
0, &\text{if}\;\quad x>2.\end{cases}\] 
\end{example}
\begin{figure}[ht]
\centering
\includegraphics[scale=0.2]{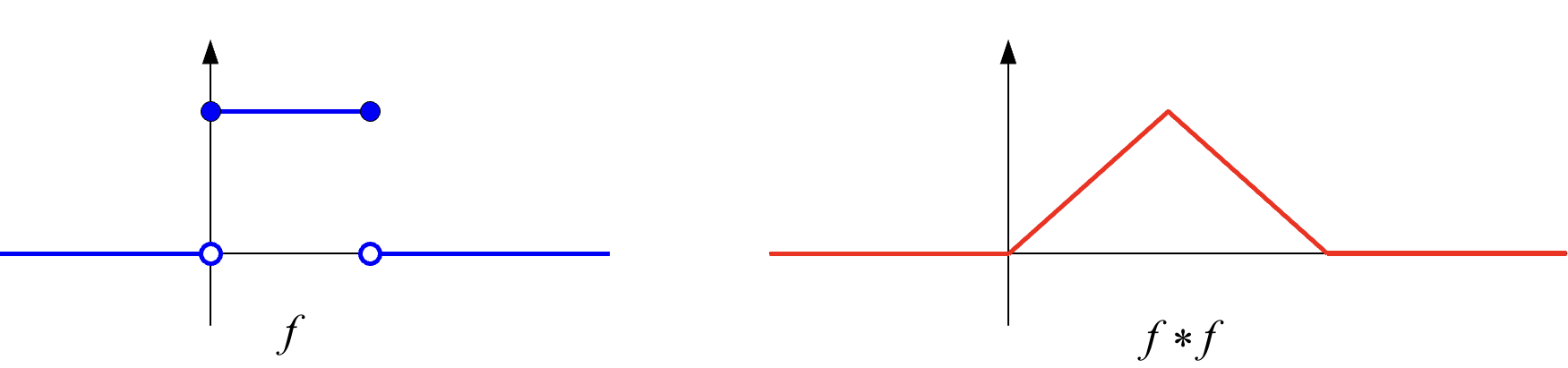}

\caption{The function $f:\mb{R}\to\mb{R}$ defined by \eqref{230820_7}  and the function $f*f$.}\label{figure83}
\end{figure}
 
 Convolution usually smooths up a function, as shown in Figure \ref{figure83}.
 
In the theory of Fourier transforms, convolution plays an important role because of the following.
\begin{theorem}{}
Let $f:\mb{R}\to\mb{C}$ and $g:\mb{R}\to\mb{C}$ be    functions that are both $L^1$ and $L^2$. Then $(f*g):\mb{R}\to\mb{C}$ is  an $L^1$ function and 
\[\mathcal{F}[f*g]=\mathcal{F}[f]\mathcal{F}[g].\]
\end{theorem}
\begin{myproof}{Sketch of Proof}
 Fubini's theorem implies that
\begin{align*}
\int_{-\infty}^{\infty}|(f*g)(x)|dx &\leq \int_{-\infty}^{\infty}\int_{-\infty}^{\infty} |f(x-t)||g(t)|dt dx\\&
\leq\int_{-\infty}^{\infty}|g(t)|\int_{-\infty}^{\infty}|f(x-t)|dx dt\\
&=\Vert f\Vert_1\int_{-\infty}^{\infty}|g(t)|dt=\Vert f\Vert_1\Vert g\Vert_1.
\end{align*}This shows that $f*g$ is an $L^1$ function.
By Fubini's theorem again, we have
\begin{align*}
\mathcal{F}[f*g](\omega)&=\int_{-\infty}^{\infty} \int_{-\infty}^{\infty}f(x-t)g(t)dt e^{-i\omega x}dx\\
&=\int_{-\infty}^{\infty}g(t)\int_{-\infty}^{\infty}f(x-t)e^{-i\omega(x-t)}dx e^{-i\omega t}dt\\
&=\mathcal{F}[f](\omega)\int_{-\infty}^{\infty}g(t)e^{-i\omega t}dt=\mathcal{F}[f](\omega)\mathcal{F}[g](\omega).
\end{align*}
\end{myproof}
\begin{example}{}
Find the Fourier transform of the function $g:\mb{R}\to\mb{R}$, 
\[g(t)=\begin{cases} 0,\quad &\text{if}\;\quad x<0,\\
t, &\text{if}\;0\leq t<1,\\
2-t,\quad &\text{if} \;1\leq t\leq 2,\\
0, &\text{if}\;\quad t>2.\end{cases}\] 
\end{example}
\begin{solution}{Solution}
By Example \ref{230820_6}, $g=f*f$, where $f$ is the function given by \ref{230820_7}. The Fourier transform of $f$ is
\[\widehat{f}(\omega)=\frac{1-e^{-i\omega}}{i\omega}=e^{-\frac{i\omega}{2}}\frac{2\di\sin\frac{\omega}{2}}{\omega}.\]
Therefore, the Fourier transform of $g$ is
\[\widehat{g}(\omega)=\widehat{f}(\omega)\times \widehat{f}(\omega)=e^{-i\omega}\frac{4\di\sin^2\frac{\omega}{2}}{\omega^2}.\]

\end{solution}

Now we list down some other useful properties of Fourier transforms. The proofs are left as exercises.
\begin{theorem}{}
Let $f:\mb{R}\to\mb{C}$ be a $L^1$ function and let $a$ be a real number. \begin{enumerate}[(a)]
\item If  $g:\mb{R}\to\mb{C}$ is the function $g(t)=f(t-a)$, then
 $\widehat{g}(\omega)=e^{-ia\omega}\widehat{f}(\omega)$.
 \item If  $h:\mb{R}\to\mb{C}$ is the function $h(t)=f(t)e^{iat}$, then 
 $\widehat{h}(\omega)= \widehat{f}(\omega-a)$.
 \end{enumerate}
\end{theorem}

\begin{example}{}
Let us consider solving a partial differential equation of the form 
\begin{equation}\label{230820_14}\frac{\pa^2u}{\pa t^2}-c^2\frac{\pa^2 u}{\pa x^2}=0,\end{equation} where $c$ is a positive constant.  This is the called the wave equation. The function $u$ is a function in $(t,x)\in \mb{R}^2$. 
For simplicity, we assume that $u$, $u_t$, $u_{tt}$  are infinitely differentiable bounded $L^1$ functions which decays to 0 when $t\to\pm\infty$.

Let $\widehat{u}(\omega, x)$ be the Fourier transform of $u$ with respect to the variable $t$. Then 
\[\mathcal{F}[u_{tt}](\omega, x)=-\omega^2\widehat{u}(\omega, x).\]
\be
 It can be justified that
\[\mathcal{F}[u_{xx}]=\frac{\pa^2}{\pa x^2}\mathcal{F}[u].\]

Thus, under Fourier transform with respect to $t$, the partial differential equation \eqref{230820_14} is transformed to
  a second order ordinary differential equation  
\begin{equation}\label{230820_11}\widehat{u}_{xx}(\omega, x)+\frac{\omega^2}{c^2} \widehat{u}(\omega, x)=0 \end{equation}
with respect to the variable $x$. The general solution is
\[ \widehat{u}(\omega, x)=A(\omega)e^{\frac{i\omega }{c}x}+B(\omega)e^{-\frac{i\omega}{c}x}\]for some infinitely differentiable functions $A(\omega)$ and $B(\omega)$.
Assume that
\[\mathcal{F}^{-1}[A](t)=\widetilde{A}(t),\hspace{1cm} \mathcal{F}^{-1}[B](t)=\widetilde{B}(t).\]
Then $A(\omega)e^{\frac{i\omega}{c}x}$ and $B(\omega)e^{-\frac{i\omega}{c}x}$ are the Fourier transforms of the functions
\[\widetilde{A}\left(t+\frac{x}{c}\right)\quad \text{and}\quad \widetilde{B}\left(t-\frac{ x}{c}\right)\]
respectively. 
These give
\[u(t,x)=\widetilde{A}\left(t+\frac{x}{c}\right)+\widetilde{B}\left(t-\frac{ x}{c}\right).\]Let
\[\phi(t)=\widetilde{A}\left(\frac{t}{c}\right)\quad\text{and}\quad \psi(t)=\widetilde{B}\left(-\frac{t}{c}\right).\]Then
\[u(t,x)=\phi(x+ct)+\psi(x-ct).\]
 This shows that the solution of the wave equation can be written as a sum of a left-travelling wave $\phi(x+ct)$ and a right-travelling wave $\psi(x-ct)$.
 
\end{example2}
 
\vp
\noindent
{\bf \large Exercises  \thesection}
\setcounter{myquestion}{1}
\begin{question}{\themyquestion}
 Let $f:\mb{R}\to\mb{C}$ be a $L^1$ function and let $a$ be a real number. Define the function $g:\mb{R}\to\mb{C}$ by 
 \[g(t)=f(t-a).\]
 Show that
 \[\widehat{g}(\omega)=e^{-ia\omega}\widehat{f}(\omega).\]
\end{question}
 \atc
 \begin{question}{\themyquestion}
 Let $f:\mb{R}\to\mb{C}$ be a $L^1$ function and let $a$ be a real number. Define the function $g:\mb{R}\to\mb{C}$ by 
 \[g(t)=f(t)e^{iat}.\]
 Show that
 \[\widehat{g}(\omega)= \widehat{f}(\omega-a).\]
\end{question}
 \atc
 \begin{question}{\themyquestion}
Find the Fourier transform of the function $f:\mb{R}\to\mb{C}$.
\begin{enumerate}[(a)]
\item
$\di f(t)=\frac{1}{t^2+4}$
\item $\di f(t)=\frac{1}{t^2+4t+13}$
\item $\di f(t)=\frac{\sin t}{t^2+4t+13}$
\end{enumerate}
\end{question}

\atc
 \begin{question}{\themyquestion}
Let  $f:\mb{R}\to\mb{R}$ be the function $f(t)= e^{-3|t|}$, and let $g:\mb{R}\to\mb{R}$ be the function $g(t)=(f*f)(t)$. Use convolution theorem to find the Fourier transform of the function $g:\mb{R}\to\mb{R}$.
\end{question}

 \atc
 \begin{question}{\themyquestion}
 Let $a$ and $b$ two distinct positive numbers, and let $f:\mb{R}\to\mb{R}$ and $g:\mb{R}\to\mb{R}$ be the functions $f(t)=e^{-at^2}$ and $g(t)=e^{-bt^2}$. Find the   function $h:\mb{R}\to\mb{R}$ defined as   $h(t)=(f*g)(t)$.
\end{question}

\atc
 \begin{question}{\themyquestion}
 Let $f:\mb{R}\to\mb{C}$ be a bounded $L^1$ function. Show that $f$ is $L^2$.
\end{question}

\appendix

\chapter{Sylvester's Criterion} \label{appA}
In this section, we give a proof of the Sylvester's criterion, which gives a necessary and sufficient condition for a symmetric matrix to be positive definite. The proof uses the $LDU$ factorization of a matrix.

Given an $n\times n$ matrix $A$ and an integer $1\leq k\leq n$, the $k^{\text{th}}$ principal submatrix of $A$, denoted by $M_k(A)$, is the $k\times k$ matrix consists of the first $k$ rows and  first $k$ columns of $A$. 
The Sylvester's criterion is the following. 
\begin{theorem}{Sylvester's Criterion for Positive Definiteness}
An $n\times n$ symmetric matrix $A$  is positive definite if and only if $\det M_k>0$ for all $1\leq k\leq n$, where $M_k$ is its $k^{\text{th}}$ principal submatrix.
\end{theorem}

For a positive integer $n$, let $\mathcal{M}_{n}$ be the vector space of $n\times n$ matrices, and let
$\mathcal{L}_{n}$, $\mathcal{U}_{n}$ and $\mathcal{D}_n$ be respectively the subspaces that consist of lower triangular, upper triangular, and diagonal matrices. Also, let
\begin{align*}
\widetilde{\mathcal{L}}_n&=\left\{L\in \mathcal{L}_n\,|\, \text{all the diagonal entries of $L$ are equal to 1} \right\},\\
\widetilde{\mathcal{U}}_n&=\left\{U\in \mathcal{U}_n\,|\, \text{all the diagonal entries of $U$ are equal to 1} \right\}.
\end{align*}
Notice that $L$ is in $\widetilde{\mathcal{L}}_n$ if and only if its transpose $L^T$ is in $\widetilde{\mathcal{U}}_n$.

The set of $n\times n$ invertible matrices is a group under matrix multiplication. This group is denoted by  $\text{GL}(n, \mb{R})$, and is called the general linear group. As a set, it is the subset of $\mathcal{M}_n$ that consists of all the matrices $A$ with $\det A\neq 0$.  The group $\text{GL}(n, \mb{R})$ has a subgroup that contains all the invertible matrices with determinant 1, deonoted by $\text{SL}(n, \mb{R})$, and is called the special linear group. 
The sets
$\widetilde{\mathcal{L}}_n$ and $\widetilde{\mathcal{U}}_n$ are subgroups of $\text{SL}(n, \mb{R})$.

If $A$ is an $n\times n$ matrix, an $LDU$ factorization of $A$ is a factorization of the form
\[A=LDU,\]
where $L\in  \widetilde{\mathcal{L}}_n$, $D\in\mathcal{D}_n$, and  $U\in \widetilde{\mathcal{U}}_n$. Notice that $\det A=\det D$. Hence, $A$ is invertible if and only if all the diagonal entries of $D$ are nonzero.

The following proposition says that the $LDU$ decomposition of an invertible matrix is unique.
\begin{proposition}[label=230806_11]{Uniqueness of $\pmb{LDU}$ Factorization}
If $A$ is an  $n\times n$ invertible matrix that has an $LDU$ factorization, then the factorization is unique.
\end{proposition}
\begin{myproof}{Proof}
We need to prove that if $L_1, L_2$ are in $\widetilde{\mathcal{L}}_n$, $U_1, U_2$ are in $\widetilde{\mathcal{U}}_n$, $D_1, D_2$ are in $\mathcal{D}_n$, and
\[L_1D_1U_1=L_2D_2U_2,\]
then $L_1=L_2$,    $U_1=U_2$ and $D_1=D_2$.

Let $L=L_2^{-1}L_1$ and $U=U_2U_1^{-1}$. Then
\[LD_1=D_2U.\]
Notice that 
$L$ is in $\widetilde{\mathcal{L}}_n$ and $LD_1$ is in $\mathcal{L}_n$. Similarly, $U$ is in $\widetilde{\mathcal{U}}_n$ and $D_2U$ is in $\mathcal{U}_n$.  The intersection of $\mathcal{L}_n$ and $\mathcal{U}_n$ is $\mathcal{D}_n$. Thus,   there exists $D\in\mathcal{D}_n$ such that
\[LD_1=D_2U=D.\]
Since $A$ is invertible, $D_1$ and $D_2$ are invertible. Hence,
\[L=DD_1^{-1}\quad\text{and}\quad U=D_2^{-1}D\] are diagonal matrices.
Since all the diagonal entries of $L$ and $U$ are 1, we find that $DD_1^{-1}=I_n$ and $D_2^{-1}D=I_n$, where $I_n$ is the $n\times n$ identity matrix. This proves that
\[D_1=D=D_2.\]
But then $L=I_n=U$, which imply that $L_1=L_2$ and $U_1=U_2$.
 
\end{myproof}

\begin{corollary}[label=230806_13]
{}
\begin{enumerate}[(i)]
\item
Given $L_0\in \mathcal{L}_n$, if $L_0$ is invertible, it has a unique $LDU$ decomposition with $U=I_n$ the $n\times n$ identity matrix.
\item Given $U_0\in \mathcal{U}_n$, if $U_0$ is invertible, it has a unique $LDU$ decomposition with $L=I_n$ the $n\times n$ identity matrix.
\end{enumerate}
\end{corollary}
\begin{myproof}{Proof}
It suffices to establish (i).
The uniqueness is asserted in Proposition \ref{230806_11}.  For the existence, let $L_0=[a_{ij}]$, where $a_{ij}=0$ if $i<j$. Since $L_0$ is invertible, $a_{ii}\neq 0$ for all $1\leq i\leq n$. Let $D=[d_{ij}]$ be the diagonal matrix with $d_{ii}=a_{ii}$ for $1\leq i\leq n$.  Then $D$ is invertible. Define $L=L_0D^{-1}$. Then $L$ is a lower triangular matrix and 
for $1\leq i\leq n$,
\[L_{ii}=a_{ii}d_{ii}^{-1}=1.\]
This shows that $L$ is in $\widetilde{\mathcal{L}}_n$. Thus, $L_0=LD$ is the $LDU$ decomposition of $L_0$ with $U=I_n$.
\end{myproof}

The following lemma says that multiplying by a matrix $L$ in $\widetilde{\mathcal{L}}_n$ does not affect the determinants of the principal submatrices.
\begin{lemma}[label=230806_12]{}
Let $A$ be an $n\times n$ matrix, and let $L$ be a matrix in $\widetilde{\mathcal{L}}_n$. If $B=LA$, then for $1\leq k\leq n$, 
\[\det M_k(B)=\det M_k(A).\]
\end{lemma}
\begin{myproof}{Proof}
For an $n\times n$ matrix $C$, we partition it into four blocks 
\[C=\begin{bmatrix}\;\begin{matrix} M_k(C) &\rvline & N_k(C)\\
\hline P_k(C) &\rvline & Q_k(C)\end{matrix}\;\end{bmatrix}.\]
\bp
For $L\in \widetilde{\mathcal{L}}_n$, $M_k(L)$ is in $\widetilde{\mathcal{L}}_k$, and $N_k(L)$ is the zero matrix.
Now $B=LA$ implies  that
\[\begin{bmatrix}\;\begin{matrix} M_k(B) &\rvline & N_k(B)\\
\hline P_k(B) &\rvline & Q_k(B)\end{matrix}\;\end{bmatrix}=\begin{bmatrix}\;\begin{matrix} M_k(L) &\rvline & N_k(L)\\
\hline P_k(L) &\rvline & Q_k(L)\end{matrix}\;\end{bmatrix}\begin{bmatrix}\;\begin{matrix} M_k(A) &\rvline & N_k(A)\\
\hline P_k(A) &\rvline & Q_k(A)\end{matrix}\;\end{bmatrix}.\]
This implies that
\[M_k(B)=M_k(L)M_k(A)+N_k(L)P_k(A)=M_k(L)M_k(A).\]
Since $M_k(L)\in \widetilde{\mathcal{L}}_k$, $\det M_k(L)=1$. Therefore,
\[\det M_k(B)=\det M_k(A).\]
\end{myproof}
Lemma \ref{230806_12} has an upper triangular counterpart. 
\begin{corollary}[label=230806_16]{}
Let $A$ be an $n\times n$ matrix, and let $U$ be a matrix in $\widetilde{\mathcal{U}}_n$. If $B=AU$, then for $1\leq k\leq n$, 
\[\det M_k(B)=\det M_k(A).\]
\end{corollary}
\begin{myproof}{Sketch of Proof}
Notice that  that $M_k(B^T)=M_k(B)^T$, and $B^T=U^TA^T$, where $U^T$ is in $\widetilde{\mathcal{L}}_n$. The result follows from the fact that $\det C^T=\det C$ for any $k\times k$ matrix $C$.
\end{myproof}

Now we prove the following theorem which asserts the existence of $LDU$ decomposition for a matrix $A$ with $\det M_k(A)\neq 0$ for all $1\leq k\leq n$.
\begin{theorem}[label=230806_15]{}
Let $A=[a_{ij}]$ be an $n\times n$ matrix such that $\det M_k(A)\neq 0$ for all $1\leq k\leq n$. Then $A$ has a unique $LDU$ decomposition.
\end{theorem}
\begin{myproof}{Proof}

Notice that $M_n(A)=A$. Since we assume that $\det M_n(A)\neq 0$, $A$ is invertible. The uniqueness of the $LDU$ decomposition of $A$ is asserted in Proposition \ref{230806_11}. 

We prove the statement by induction on $n$. When $n=1$, take $L=U=[1]$ and $D=A=[a]$ itself. Then $A=LDU$ is the $LDU$ decomposition of $A$.

Let $n\geq 2$. 
Suppose we have proved that any $(n-1)\times (n-1)$ matrix $B$ that satisfies $\det M_k(B)\neq 0$ for $1\leq k\leq n-1$ has a unique $LDU$ decomposition.
 
Now assume that $A$ is  an $n\times n$ matrix with $\det M_k(A)\neq 0$ for all $1\leq k\leq n$.
Since $\det M_1(A)=a_{11}$, $a=a_{11}\neq 0$. Let $L_1=[L_{ij}]$ be the matrix in $\widetilde{\mathcal{L}}_n$ such that for $2\leq i\leq n$,
\[L_{i1}=\frac{a_{i1}}{a},\]
and for $2\leq j<i\leq n$, $L_{ij}=0$. Namely,
\[L_1=\begin{bmatrix}\;\begin{matrix} 1 & \rvline & 0 \\\hline
  P_{1}(L_1) &\rvline & I_{n-1}\end{matrix}\;\end{bmatrix},\]
  where
  \[P_1(L_1)=\frac{1}{a}P_1(A).\]
  Notice that
  \[L_1^{-1}=\begin{bmatrix}\;\begin{matrix}1 & \rvline & 0 \\\hline
  -P_{1}(L_1) &\rvline & I_{n-1}\end{matrix}\;\end{bmatrix},\]
  and
  \[C=L_1^{-1}A=\begin{bmatrix}\;\begin{matrix}1 & \rvline & 0 \\\hline
  -P_{1}(L_1) &\rvline & I_{n-1}\end{matrix}\;\end{bmatrix}\begin{bmatrix}\;\begin{matrix} a &\rvline & N_1(A)\\
\hline P_1(A) &\rvline & Q_1(A)\end{matrix}\;\end{bmatrix}=\begin{bmatrix}\;\begin{matrix} a &\rvline & N_1(C)\\
\hline P_1(C) &\rvline & Q_1(C)\end{matrix}\;\end{bmatrix}\]is a matrix with \[P_1(C)=-aP_1(L_1)+P_1(A)=0.\]

\bp
By Lemma \ref{230806_12},
\[\det M_k(C)=\det M_k(A)\hspace{1cm}\text{for all}\;1\leq k\leq n.\]
Let $B=Q_1(C)$. Then $B$ is an $(n-1)\times (n-1)$ matrix. Since $P_1(C)=0$, we find that for $1\leq k\leq n-1$,
\[\det M_{k+1}(C)=a\det M_k(B).\]This shows that 
\[\det M_k(B)\neq 0\hspace{1cm}\text{for all}\;1\leq k\leq n-1.\]
By inductive hypothesis, $B$ has a unique $LDU$ decomposition given by
\[B=L_BD_BU_B.\]

Now let $L_2$ be the matrix in $\widetilde{\mathcal{L}}_{n}$ given by
\[L_2=\begin{bmatrix}\;\begin{matrix}1 & \rvline & 0 \\\hline
 0 &\rvline & L_B^{-1}\end{matrix}\;\end{bmatrix}.\]
 One can check that  
 \[(L_1L_2)^{-1}A=L_2^{-1}C=\begin{bmatrix}\;\begin{matrix} a &\rvline & N_1(C)\\\hline 0 &\rvline & D_BU_B\end{matrix}\;\end{bmatrix}.\]Let $L=L_1L_2$.  Then $L$ is in $ \widetilde{\mathcal{L}}_n$.  Since $D_BU_B$ is an upper triangular $(n-1)\times (n-1)$ matrix,
 $L^{-1}A$ is an upper triangular $n\times n$ matrix. By Corollary \ref{230806_13}, $L^{-1}A$ has a decomposition
 \[L^{-1}A=DU,\]
 where $D\in \mathcal{D}_n$ and $U\in\widetilde{\mathcal{U}}_n$. Thus, $A=LDU$ is the $LDU$ decomposition of $A$.

\end{myproof}

Now  we can complete the proof of the Sylvester's criterion for a  symmetric matrix to be positive definite.
\begin{myproof}{Proof of Sylvester's Criterion}
Let $A$ be an $n\times n$ symmetric matrix. First we prove that if $A$ is positive definite, then for $1\leq k\leq n$, $\det M_k(A)>0$. Notice that $M_k(A)$  is also a symmetric matrix. For $\mf{u}\in \mb{R}^k$, let $\mf{v}$ be the vector in $\mb{R}^n$ given by $\mf{v}=(\mf{u}, 0, \ldots, 0)$. Then
\[\mf{v}^TA\mf{v}=\mf{u}^TM_k(A)\mf{u}.\]
This shows that $M_k(A)$ is also positive definite. Hence, all the eigenvalues of $M_k(A)$ must be positive. This implies that $\det M_k(A)>0$.

Conversely, assume that $\det M_k(A)>0$ for all $1\leq k\leq n$.  By Theorem \ref{230806_15}, $A$ has a $LDU$ decomposition given by
\[A=LDU.\]Since $A$ is symmetric, $A^T=A$.
This gives
\[U^TD^TL^T=A^T=A=LDU.\]
Since $U^T$ is in $\widetilde{\mathcal{L}}_n$ and $L^T$ is in $\widetilde{\mathcal{U}}_n$, the uniqueness of $LDU$ decomposition implies that $U=L^T$. 
Hence, 
\[A=LDL^T.\]
By Lemma \ref{230806_12} and Corollary \ref{230806_16}, 
\[\det M_k(A)=\det M_k(D).\]
 If $D=[d_{ij}]$, let $\tau_i=d_{ii}$. Then $M_k(D)=\tau_1\tau_2\ldots \tau_k$. Since $\det M_k(A)>0$ for all $1\leq k\leq n$, $\tau_i>0$ for all $1\leq i\leq n$. By the invertible change of coordinates $\mf{y}=L^T\mf{x}$, we find that if $\mf{x}\in\mb{R}^n\setminus \{\mf{0}\}$,
 \[\mf{x}^TA\mf{x}=\mf{y}^TD\mf{y}=\tau_1y_1^2+\tau_2y_2^2+\cdots+\tau_ny_n^2>0.\]
 This proves that $A$ is positive definite.

\end{myproof}

\chapter{Volumes of Parallelepipeds} \label{appB}

In this appendix, we give a geometric proof   of the formula for the volume of a parallelepiped in $\mb{R}^n$.
\begin{theorem}[label=230903_11]{}
Let $\mathscr{P}$ be a parallelepiped in $\mb{R}^n$ spanned by the linearly independent vectors $\mf{v}_1, \ldots, \mf{v}_n$. Then the volume of $\mathscr{P}$ is equal to $|\det A|$, where $A$ is the matrix whose column vectors are $\mf{v}_1, \ldots, \mf{v}_n$.
\end{theorem}

Let us look at 
a special case of parallelepiped where this theorem is easy to prove by simple geometric consideration. 

\begin{definition}{Generalized Rectangles}
  A parallelepiped that is spanned by $n$ nonzero orthogonal vectors $\mf{w}_1$, $\ldots$, $\mf{w}_n$ is called a generalized rectangle.
  \end{definition}
 A generalized rectangle $R$ based at the origin and spanned by the $n$ nonzero orthogonal vectors $\mf{w}_1$, $\ldots$, $\mf{w}_n$ is equal to  $B(\mb{Q}_n)$, where $Q_n=[0,1]^n$ is the standard unit cube, and $B$ is the matrix
 \[B=\begin{bmatrix}\mf{w}_1 &\rvline & \cdots &\rvline & \mf{w}_n\end{bmatrix}.\] By geometric consideration, the volume of $R$ is given by the product of the lengths of its edges. Namely, 
\[\text{vol}\,(R)=\Vert\mf{w}_1\Vert\,\cdots\,\Vert\mf{w}_n\Vert.\] To see that this is equal to $\det B$, let $\mf{u}_1$, $\ldots$, $\mf{u}_n$  be the unit vectors in the directions of $\mf{w}_1$, $\ldots$, $\mf{w}_n$. Namely,
\[\mf{u}_i=\frac{\mf{w}_i}{\Vert\mf{w}_i\Vert}\hspace{1cm}1\leq i\leq n.\]
Then $B=PD$, where $P$ is an orthogonal matrix and $D$ is a diagonal matrix given respectively by
\begin{equation}\label{230905_1}P=\begin{bmatrix}\mf{u}_1 &\rvline & \cdots &\rvline & \mf{u}_n\end{bmatrix},\hspace{1cm}D=\begin{bmatrix} \Vert\mf{w}_1\Vert & 0 & \cdots & 0\\
0&\Vert\mf{w}_2\Vert & \cdots & 0\\\vdots & \vdots & \ddots & \vdots \\
0 & 0 & \cdots & \Vert\mathbf{w}_n\Vert\end{bmatrix}.\end{equation}

An $n\times n$ matrix $P$ is called an {\it orthogonal matrix} if \[P^TP=PP^T=I_n,\]
where $I_n$ is the $n\times n$ identity matrix. A matrix $P$ is orthogonal if and only if the column vectors of $P$ form an orthonormal basis of $\mb{R}^n$. If $P$ is orthogonal, $P^{-1}=P^T$, and  $P^{-1}$ is also orthogonal.
From $P^TP=I_n$, we find that 
\[\det(P)\det(P^T)=\det(I_n)=1.\]
Since $\det(P^T)=\det(P)$, we have $\det(P)^2=1$. Hence, the determinant of an orthogonal matrix can only be 1 or $-1$. 
Therefore, when $B=PD$, with $P$ and $D$ as given in \eqref{230905_1}, we have
\[|\det B|=|\det P\det D|=|\det D|=\Vert\mf{w}_1\Vert\,\cdots\,\Vert\mf{w}_n\Vert.\]
\begin{remark}{}
In the argument above, we do not show that the volume of a generalized rectangle spanned by the $n$ nonzero orthogonal vectors $\mf{w}_1$, $\ldots$, $\mf{w}_n$ is equal to $\Vert\mf{w}_1\Vert\,\cdots\,\Vert\mf{w}_n\Vert$ using the definition of $\text{vol}\,(R)$ in terms of a Riemann integral  $\di\int_R \,d\mf{x}$.  This is elementary but tedious.
\end{remark}

A linear transformation $\mf{T}:\mb{R}^n\to\mb{R}^n$, $\mf{T}(\mf{x})=P\mf{x}$ defined by an orthogonal matrix $P$ is called an {\it orthogonal transformation}. The significance of an orthognal transformation is as follows.
For any $\mf{u}$ and $\mf{v}$ in $\mb{R}^n$, 
\[\langle \mf{T}(\mf{u}), \mf{T}(\mf{v})\rangle =(P\mf{u})^T(P\mf{v})=\mf{u}^TP^TP\mf{v}=\mf{u}^T\mf{v}=\langle\mf{u},\mf{v}\rangle.\]Namely, $\mf{T}$ preserves inner products. Since lengths and angles are defined in terms of the inner product, this implies that an orthogonal transformation preserves lengths and angles.

Under an orthogonal transformation, the image of a rectangle $R$ is a rectangle that is congruent to  $R$. Since the volume of a Jordan measurable set $\mk{D}$  is obtained by taking the limit of a sequence of Darboux lower sums, and each Darboux lower sum is a sum of volumes of  rectangles with disjoint interiors that lie in $\mk{D}$, we find that orthogonal transformations also preserve the volumes of Jordan measurable sets.

\begin{theorem}{}
If $\mf{T}:\mb{R}^n\to\mb{R}^n$, $\mf{T}(\mf{x})=P\mf{x}$ is an orthogonal transformation, and $\mk{D}$ is a Jordan measurable set, then $\mf{T}(\mk{D})$ is also Jordan measurable and
\[\text{vol}\,(\mf{T}(\mk{D}))=\text{vol}\,(\mk{D}).\]
\end{theorem}  

To finish the proof of  Theorem \ref{230903_11}, we also need the following fact.
\begin{proposition}[label=230905_2]{}
Let $\mathscr{P}$ be a  parallelepiped based at the origin and spanned by the vectors $\mf{v}_1$, $\ldots$, $\mf{v}_n$. Assume that
\[\pi_n(\mf{v}_i)=0\hspace{1cm}\text{for}\;1\leq i\leq n-1,\]
or equivalently, $\mf{v}_1, \ldots, \mf{v}_{n-1}$ lies in the plane $x_n=0$. For $1\leq i\leq n-1$, let $\mf{z}_i\in\mb{R}^{n-1}$ be such that $\mf{v}_i=(\mf{z}_i, 0)$. If $\mathcal{Q}$ is the parallelepiped in $\mb{R}^{n-1}$ based at the origin and spanned by $\mf{z}_1$, $\ldots$, $\mf{z}_{n-1}$, then 
\[\text{vol}\,(\mathscr{P})=\text{vol}\,(\mathcal{Q})h,\]
where $h$ is the distance from $\mf{v}_n$ to the $x_n=0 $ plane, which is given explicitly by
\[h=\left|\,\text{proj}_{\mf{e}_n}\mf{v}_n\right|.\]
\end{proposition}
When $n=3$, this can be argued geometrically. For general $n$, let us give a proof using the definition of volume as a Riemann integral.
\begin{myproof}{Proof}
Recall that
\[\mathscr{P}=\left\{t_1\mf{v}_1+\cdots+t_{n-1}\mf{v}_{n-1}+t_n\mf{v}_n\,|\,\mf{t}\in [0,1]^n\right\}.\]
Notice that $\mf{v}_n$ can be written as $\mf{v}_n=(\mf{a}, h)$ for some $\mf{a}\in \mb{R}^{n-1}$. Hence, if a point $\mf{x}$ is in $\mathscr{P}$, then \[\mf{x}=\left(\di \frac{t}{h}\mathbf{a}+\mf{z}, t\right),\] where $0\leq t\leq h$, and $\mf{z}$ is a point in $\mathcal{Q}$. For $0\leq t\leq h$, let
\[\mathcal{Q}_t=\left\{\left.\left(\di \frac{t}{h}\mathbf{a}+\mf{z}, t\right)\,\right|\, 0\leq t\leq h\right\}.\]
Then it is a $(n-1)$-dimensional parallelepiped contained in the hyperplane $x_n=t$, which is a translate of the $(n-1)$-dimensional parallelepiped $\mathcal{Q}_0$. By Fubini's theorem,
\begin{align*}
\text{vol}\,\left(\mathscr{P}\right)&=\int_{\mathscr{P}}\,d\mf{x}=\int_0^h\left( \int_{\mathcal{Q}_t}\,dx_1\cdots dx_{n-1}\right) dx_n\\
&=\int_0^h \text{vol}\,(\mathcal{Q}_t)dt=\int_0^h \text{vol}\,(\mathcal{Q}_0)dt= \text{vol}\,(\mathcal{Q}) h.
\end{align*}
\end{myproof}

Now we can prove Theorem \ref{230903_11}.

\begin{myproof}{\linkt Proof of Theorem \ref{230903_11} } We prove by induction on $n$. The $n=1$ case is obvious. 
Assume that we have proved the $n-1$ case. Now given that $\mathscr{P}$ is a parallelepiped in $\mb{R}^n$ which is spanned by $\mf{v}_1, \ldots, \mf{v}_n$, we can assume that $\mathscr{P}$ is based at the origin $\mf{0}$ because translations preserve volumes. Let 
\[A=\begin{bmatrix}\mf{v}_1 &\rvline & \cdots & \rvline & \mf{v}_n\end{bmatrix}.\] We want to show that
\[\text{vol}\,(\mathscr{P})=|\det A|.\]
\bp
Let $W$ be the subspace of $\mb{R}^{n-1}$ that is spanned by $\mf{v}_1, \ldots, \mf{v}_{n-1}$. Applying the Gram-Schmidt process to the basis $\{\mf{v}_1, \ldots, \mf{v}_n\}$ of $\mb{R}^n$, we obtain an orthonormal basis $\{\mf{u}_1, \ldots, \mf{u}_n\}$. By the algorithm, the unit vector  $\mf{u}_n$ is orthogonal to the subspace $W$. 
Let 
\[P=\begin{bmatrix}\mf{u}_1 &\rvline & \cdots & \rvline & \mf{u}_n\end{bmatrix}\]
be the orthogonal matrix whose column vectors are $\mf{u}_1, \ldots, \mf{u}_n$, and consider the orthogonal transformation $\mf{T}:\mb{R}^n\to\mb{R}^n$, $\mf{T}(\mf{x})=P^{-1}\mf{x}=P^T\mf{x}$. 
For $1\leq i\leq n$, let \[\widetilde{\mf{v}}_i=\mf{T}(\mf{v}_i).\]
Then $\widetilde{\mathscr{P}}=\mf{T}(\mathscr{P})$ is a parallelepiped that has the same volume as $\mathscr{P}$, and it is spanned by $\widetilde{\mf{v}}_1, \ldots, \widetilde{\mf{v}}_n$. Notice that
\begin{align*}
\widetilde{A}&=\begin{bmatrix}\widetilde{\mf{v}}_1 &\rvline & \cdots &\rvline &\widetilde{\mf{v}}_n\end{bmatrix}=P^T \begin{bmatrix}\mf{v}_1 &\rvline & \cdots & \rvline & \mf{v}_n\end{bmatrix}\\&=\begin{bmatrix}\;\begin{matrix}\hspace{0.5cm}\mf{u}_1^T\hspace{0.5cm}\\\hline \vdots\\\hline \mf{u}_n^T\end{matrix}\;\end{bmatrix} \begin{bmatrix}\mf{v}_1 &\rvline & \cdots & \rvline & \mf{v}_n\end{bmatrix}
 =\begin{bmatrix} B & \rvline &\begin{matrix} \langle \mf{u}_1, \mf{v}_n\rangle \\\vdots \\  \langle \mf{u}_{n-1}, \mf{v}_n\rangle\end{matrix}\\\hline
\hspace{0.8cm}\mf{0} \hspace{0.8cm} & \rvline & \langle \mf{u}_{n}, \mf{v}_n\rangle\end{bmatrix}.\end{align*}
From this, we find that 
\[\det (\widetilde{A})=\det(B)\times \langle \mf{u}_{n}, \mf{v}_n\rangle.\]
Comparing the columns, we also have
\[\widetilde{\mf{v}}_i=(\mf{z}_i, 0)\hspace{1cm}\text{for}\; 1\leq i\leq n-1,\]
where $\mf{z}_1, \ldots, \mf{z}_{n-1}$ are the column vectors of $B$, which are vectors in $\mb{R}^{n-1}$; and 
\[\widetilde{\mf{v}}_n=\begin{bmatrix} \;\;\begin{matrix} \langle \mf{u}_1, \mf{v}_n\rangle \\\vdots \\  \langle \mf{u}_{n-1}, \mf{v}_n\rangle\\\hline\end{matrix}\;\; \\\langle \mf{u}_{n}, \mf{v}_n\rangle\end{bmatrix}. \]

\bp
The transformation $\mf{T}$ maps the subspace $W$ to the hyperplane $x_n=0$, which can be identified with $\mb{R}^{n-1}$. Let $\mathcal{Q}$ be the parallelepiped in $\mb{R}^{n-1}$ based at the origin and spanned by the vectors $\mf{z}_1, \ldots, \mf{z}_{n-1}$. The volume of the parallepiped $\widetilde{\mathscr{P}}$ is equal to the volume of $\mathcal{Q}$ times the distance $h$ from the tip of the vector $\widetilde{\mf{v}}_n$ to the plane $x_n=0$. 
By definition,
\[h=\left\Vert\,\text{proj}_{\mf{e}_n}\widetilde{\mf{v}}_n\right\Vert=\left|\langle \mf{u}_n, \mf{v}_n\rangle\right|.\]
 
Proposition \ref{230905_2} gives
\[
\text{vol}\,(\widetilde{\mathscr{P}})=\text{vol}\,(\mathcal{Q})\times \left|\langle \mf{u}_n, \mf{v}_n\rangle\right|.\]
By inductive hypothesis,
\[\text{vol}\,(\mathcal{Q})=|\det (B)|.\]
Therefore,
\[ \text{vol}\,(\widetilde{\mathscr{P}})=\left|\det(B)\times \langle \mf{u}_n, \mf{v}_n\rangle\right|=|\det (\widetilde{A})|.\]
 Since $\widetilde{A}=P^TA$, we find that
 \[\det (\widetilde{A})=\det(P^T)\det(A)=\pm \det (A).\]
 Hence, 
\[   \text{vol}\,( \mathscr{P})=\text{vol}\,(\widetilde{\mathscr{P}})=|\det (A)|.\]
 This completes the proof of Theorem \ref{230903_11}.
 \end{myproof}
 As a corollary , we have the following. 

\begin{theorem}[label=230903_12]{}
Let $\mf{I}$ be a closed rectangle in $\mb{R}^n$, and let $\mf{T}:\mb{R}^n\to\mb{R}^n$, $\mf{T}(\mf{x})=A\mf{x}$ be an invertible linear transformation. Then
\begin{equation*}\text{vol}\,(\mf{T}(\mf{I}))=|\det A|\,\text{vol}\,(\mf{I}).\end{equation*}
\end{theorem}
\begin{myproof}{Proof}
Let $\mf{I}=\di\prod_{i=1}^n[a_i, b_i]$. Then $\mf{I}=\mf{S}(Q_n)+\mf{a}$, where $\mf{a}=(a_1, \ldots, a_n)$, $Q_n$ is the standard unit cube $[0,1]^n$, and $\mf{S}:\mb{R}^n\to\mb{R}^n$ is the linear transformation defined by the diagonal matrix $B$ with diagonal entries $b_1-a_1, b_2-a_2, \ldots, b_n-a_n$. Therefore,
\[\mf{T}(\mf{I})=(\mf{T}\circ \mf{S})(Q_n)+\mf{T}(\mf{a}).\]
Since the matrix associated with the linear transformation $(\mf{T}\circ \mf{S}):\mb{R}^n\to\mb{R}^n$ is $AB$,  $\mf{T}(\mf{I})$ is a parallelepiped based at $\mf{T}(\mf{a})$ and spanned by the columnn vectors of $AB$. By Theorem \ref{230903_11}, 
 \begin{align*}
 \text{vol}\,\left(\mf{T}(\mf{I})\right)=|\det(AB)|=|\det(A)||\det(B)|.
 \end{align*}
Obviously, \[|\det B|=\di\prod_{i=1}^n (b_i-a_i)=\text{vol}\,(\mf{I}).\]
This proves that
\[\text{vol}\,(\mf{T}(\mf{I}))=|\det A|\,\text{vol}\,(\mf{I}).\]
\end{myproof}
\begin{remark}{}
The formula 
\[\text{vol}\,(\mf{T}(\mf{I}))=|\det A|\,\text{vol}\,(\mf{I})\]still holds even though the matrix $A$ is not invertible. In this case, $\det A=0$, and  the column vectors of $A$ are not linearly independent. Therefore,  $\mf{T}(\mf{I})$ lies in a plane in $\mb{R}^n$, and so  $\mf{T}(\mf{I})$ has zero volume.
\end{remark}

\chapter[Riemann Integrability]{Necessary and Sufficient Condition for Riemann Integrability} \label{appC}

In this appendix, we want to prove the Lebesgue-Vitali theorem which gives a necessary and sufficient condition for a   bounded function $f:\mk{D}\to\mb{R}$ to be Riemann integrable. We will introduce the concept of Lebesgue measure zero without introducing the concept of general Lebesgue measure. The latter is often covered in a standard course in real analysis.

Recall that the volume of a closed rectangle $\mf{I}=\di\prod_{i=1}^n[a_i,b_i]$ or its interior $\di \text{int}\,(\mf{I})=\prod_{i=1}^n(a_i,b_i)$ is
\[\text{vol}\,(\mf{I})=\text{vol}\,(\text{int}\,\mf{I}) =\prod_{i=1}^n(b_i-a_i).\] 

 If $A$ is a subset of $\mb{R}^n$, we say that $A$ has Jordan content zero if 
 \begin{enumerate}[(i)]
 \item for every $\varepsilon>0$, there are finitely many closed rectangles $\mf{I}_1$, $\ldots$, $\mf{I}_k$ such that 
\[A\subset\bigcup_{j=1}^k \mf{I}_j\quad\text{and}\quad\sum_{j=1}^k\text{vol}\,(\mf{I}_j)<\varepsilon.\]
\end{enumerate}
This is equivalent to any of the followings.
\begin{enumerate}[(i)]
\item[(ii)] For every $\varepsilon>0$, there are finitely many closed cubes $Q_1$, $\ldots$, $Q_k$ such that 
\[A\subset\bigcup_{j=1}^k Q_j\quad\text{and}\quad\sum_{j=1}^k\text{vol}\,(Q_j)<\varepsilon.\]

\item[(iii)] For every $\varepsilon>0$, there are finitely many open rectangles $U_1$, $\ldots$, $U_k$ such that 
\[A\subset\bigcup_{j=1}^k U_j\quad\text{and}\quad\sum_{j=1}^k\text{vol}\,(U_j)<\varepsilon.\]
\item[(iv)]
For every $\varepsilon>0$, there are finitely many open cubes $V_1$, $\ldots$, $V_k$ such that 
\[A\subset\bigcup_{j=1}^k V_j\quad\text{and}\quad\sum_{j=1}^k\text{vol}\,(V_j)<\varepsilon.\]
\end{enumerate}

A set     has Jordan content zero if and only if it is Jordan measurable and its volume is zero. Hence, we also call  a set   that has Jordan content zero as a set that has {\it Jordan measure zero}. The Jordan measure of a Jordan measurable set $A$   is the volume of $A$ defined as the  Riemann integral of the characteristic function $\chi_A:A\to\mb{R}$.

In Lebesgue measure, instead of a covering by finitely many rectangles, we allow a covering  by countably many rectangles.
A set $S$ is countable if it is finite or it is countably infinite. The latter means that there is a one-to-one correspondence between $S$ and the set $\mb{Z}^+$. In any case, a set $S$ is countable if and only if there is a surjection $h:\mb{Z}^+\to S$, which allows us to write 
\[S=\left\{s_k\,|\, k\in\mb{Z}^+\right\},\hspace{1cm}\text{where}\;s_k=h(k).\] 

\begin{definition}
{Lebesgue Measure Zero}
Let $A$ be a subset of $\mb{R}^n$. We say that $A$ has Lebesgue measure zero if for every $\varepsilon>0$, there is a countable collection of open rectangles $\{U_k\,|\,k\in\mb{Z}^+\}$ that covers $A$, the sum of whose volumes is less than $\varepsilon$. Namely,
\[A\subset\bigcup_{k=1}^{\infty}U_k\quad\text{and}\quad \sum_{k=1}^{\infty}\text{vol}\,(U_k)<\varepsilon.\]
\end{definition}

The following is obvious.
\begin{proposition}[label=230912_10]{}
Let $A$  be a subset of $\mb{R}^n$. If $A$ has Jordan content zero, then it has Lebesgue measure zero.
\end{proposition}
The converse is not true. There are sets with Lebesgue measure zero, but they do not have Jordan content zero. The following gives an example of such sets.
\begin{example}[label=230912_1]{}
Let $A=\mb{Q}\cap [0,1]$. The function $\chi_A:[0,1]\to\mb{R}$ is the Dirichlet's function, which is not Riemann integrable. Hence, $A$ is not Jordan measurable. Nevertheless, we claim that $A$ has Lebesgue measure zero.
 
Recall that $\mb{Q}$ is a countable set. As a subset of $\mb{Q}$, $A$ is also countable.   Hence, we can write $A$ as \[A=\left\{a_k\,|\,k\in\mb{Z}^+\right\}.\]

Given $\varepsilon>0$ and $k\in\mb{Z}^+$, let $U_k$ be the open rectangle
\[U_k=\left(a_k-\frac{\varepsilon}{2^{k+2}}, a_k+\frac{\varepsilon}{2^{k+2}}\right).\]
Then $a_k\in U_k$ for each $k\in\mb{Z}^+$. Thus, 
\[A\subset \bigcup_{k=1}^{\infty} U_k.\]
Now,
\[\sum_{k=1}^{\infty}\text{vol}\,(U_k)=\sum_{k=1}^{\infty}\frac{\varepsilon}{2^{k+1}}=\frac{\varepsilon}{2}<\varepsilon.\]
Therefore, $A$ has Lebesgue measure zero.
\end{example}

The converse to Proposition \ref{230912_10} is true if $A$ is compact.

\begin{proposition}[label=230912_9]{}
Let $A$ be a compact subset of $\mb{R}^n$. If $A$ has Lebesgue measure zero, then it has Jordan content zero.
\end{proposition}
\begin{myproof}{Proof}
Given $\varepsilon>0$, since $A$ has Lebesgue measure zero, there is a countable collection $\{U_{\alpha}\,|\,\alpha\in \mb{Z}^+\}$ of open rectangles that covers $A$, and
\[\sum_{\alpha\in\mb{Z}^+}\text{vol}\, (U_{\alpha})<\varepsilon.\]
\bp
Since $A$ is compact, there is a finite subcollection $\left\{U_{\alpha_l}\,|\,1\leq l\leq m\right\}$ that covers $A$. Obviously, we also have
\[\sum_{l=1}^m\text{vol}\,(U_{\alpha_l})<\varepsilon.\]
Hence, $A$ has Jordan content 0.
\end{myproof}

\begin{example}{}
Using the same reasoning as in Example \ref{230912_1}, one can show that any countable subset of $\mb{R}^n$ has Lebesgue measure zero.
\end{example}

 We have seen that if $A$ is a subset of $\mb{R}^n$ that has Jordan content zero, then its closure $\overline{A}$ also has Jordan content zero.  However, the same is not true for Lebesgue measure.
\begin{example}{}
Example \ref{230912_1} shows that the set $A=\mb{Q}\cap [0,1]$ has Lebesgue measure zero. Notice that $\overline{A}=[0,1]$. It  cannot have Lebesgue measure zero.
\end{example}

As in the case of Jordan content zero, we have the following equivalences for  a set $A$ in $\mb{R}^n$ to have Lebesgue measure zero.

\begin{enumerate}[(i)]
\item
For every $\varepsilon>0$, there is a countable collection of open rectangles $\{U_k\,|\,k\in\mb{Z}^+\}$ such that
\[A\subset\bigcup_{k=1}^{\infty}U_k\quad\text{and}\quad \sum_{k=1}^{\infty}\text{vol}\,(U_k)<\varepsilon.\]

\item[(ii)]
For every $\varepsilon>0$, there is a countable collection of closed rectangles $\{\mf{I}_k\,|\,k\in\mb{Z}^+\}$ such that
\[A\subset\bigcup_{k=1}^{\infty}\mf{I}_k\quad\text{and}\quad \sum_{k=1}^{\infty}\text{vol}\,(\mf{I}_k)<\varepsilon.\]

\item[(iii)]
For every $\varepsilon>0$, there is a countable collection of open cubes $\{V_k\,|\,k\in\mb{Z}^+\}$ such that
\[A\subset\bigcup_{k=1}^{\infty}V_k\quad\text{and}\quad \sum_{k=1}^{\infty}\text{vol}\,(V_k)<\varepsilon.\]
\item[(iv)]
For every $\varepsilon>0$, there is a countable collection of closed cubes $\{Q_k\,|\,k\in\mb{Z}^+\}$ such that
\[A\subset\bigcup_{k=1}^{\infty}Q_k\quad\text{and}\quad \sum_{k=1}^{\infty}\text{vol}\,(Q_k)<\varepsilon.\]
\end{enumerate}

The following is obvious.
\begin{proposition}{}
Let $A$ be a subset of $\mb{R}^n$. If $A$ has Lebesgue measure zero, and $B$ is a subset of $A$, then $B$ also has Lebesgue measure zero.
\end{proposition}

Using the fact that the set $\mb{Z}^+\times\mb{Z}^+$ is countable, we find that a countable union of coutable sets is countable. This gives the following.
\begin{proposition}[label=230912_6]{}
Let $\{A_m\,|\,m\in\mb{Z}^+\}$ be a countable collection of subsets of $\mb{R}^n$. If each of the $A_m, m\in\mb{Z}^+$ has Lebesgue measure zero, then the set
$\di A=\bigcup_{m=1}^{\infty}A_m$
also has Lebesgue measure zero. 
\end{proposition}
\begin{myproof}{Proof}
Fixed $\varepsilon>0$. For each $m\in \mb{Z}^+$, since $A_m$ has Lebesgue measure zero, there is a countable collection $\mathscr{B}_m=\di\left\{U_{m,k}\,|\,k\in\mb{Z}^+\right\}$ of open rectangles such that
\[A_m\subset \bigcup_{k=1}^{\infty}U_{m,k}\quad\text{and}\quad\sum_{k=1}^{\infty}\text{vol}\,(U_{m,k})<\frac{\varepsilon}{2^m}.\]
 
It follows that
\[A\subset \bigcup_{m=1}^{\infty} \bigcup_{k=1}^{\infty}U_{m,k}\quad\text{and}\quad \sum_{m=1}^{\infty}\sum_{k=1}^{\infty}\text{vol}\,(U_{m,k})<\sum_{m=1}^{\infty}\frac{\varepsilon}{2^m}=\varepsilon.\]
\bp
Notice that the collection 
\[\mathscr{B}=\bigcup_{m=1}^{\infty}\mathscr{B}_m=\left\{U_{m,k}\,|\,m\times k\in \mb{Z}^+\times \mb{Z}^+\right\}\] is countable. This proves that $A$ has Lebesgue measure zero.
\end{myproof}

Now we   proceed to  the main theorem.
\begin{theorem}[label=230912_3]{Lebesgue-Vitali Theorem}
Let $\di\mf{I}=\prod_{i=1}^{n}[a_i,b_i]$ be a closed rectangle in $\mb{R}^n$.  Given a bounded function $f:\mf{I}\to\mb{R}$, let $\mathcal{N}$ be its set of discontinuities. Then $f:\mf{I}\to\mb{R}$ is Riemann integrable if and only if $\mathcal{N} $ has Lebesgue measure zero.
\end{theorem}
In this theorem, we only consider functions defined on closed rectangles. This is because the Riemann integrability of a function $f:\mk{D}\to\mb{R}$ defined on a bounded set $\mk{D}$ is defined in terms of the Riemann integrability of its zero extension $\check{f}:\mf{I}\to\mb{R}$ to a closed rectangle $\mf{I}$ that contains $\mk{D}$.
 
 To prove the Lebesgue-Vitali theorem, we need a few  lemmas.
 Given a bounded function $f:\mk{D}\to\mb{R}$, we define the oscillation of 
 $f$ at a point $\mf{x}_0\in\mk{D}$ as
 \[\omega_f(\mf{x}_0)=\lim_{r\to 0^+}\sup\left\{f(\mf{u})-f(\mf{v})\,|\, \mf{u} , \mf{v}\in B(\mf{x}_0, r)\cap\mk{D}\right\}.\]\
 Notice that the set 
 \[\mathcal{F}_r=\left\{f(\mf{u})-f(\mf{v})\,|\, \mf{u} , \mf{v}\in B(\mf{x}_0, r)\cap\mk{D}\right\}\] is a bounded subset of real numbers   and $-\mathcal{F}_r=\mathcal{F}_r$. Thus, the supremum of $\mathcal{F}_r$  always exists and is  nonnegative. 
 It is easy to see that
 \[\mathcal{F}_{r_1}\subset\mathcal{F}_{r_2}\hspace{1cm}
 \text{if}\; r_1<r_2.\]
 Therefore, 
 $\sup\mathcal{F}_r$   decreases as $r\to 0^+$
 This implies that 
 \[\omega_f(\mf{x}_0)=\lim_{r\to 0^+}\sup \mathcal{F}_r=\inf_{r>0}\sup_{\mf{u},\mf{v}\in B(\mf{x}_0, r)\cap\mk{D}}(f(\mf{u})-f(\mf{v}))\]
 exists and is nonnegative.
 
 \begin{lemma}[label=230912_5]{}
 Let $\mk{D}$ be a subset of $\mb{R}^n$, and let $\mf{x}_0$ be a point in $\mk{D}$. Assume that $f:\mk{D}\to\mb{R}$ is a bounded function. Then $f$ is continuous at $\mf{x}_0$ if and only if $\omega_f(\mf{x}_0)=0$.
 \end{lemma}
 \begin{myproof}{Proof}
 First assume that $f$ is continuous at $\mf{x}_0$. Given $\varepsilon>0$, there is a $\delta>0$ such that for all $\mf{x}\in B(\mf{x}_0, \delta)\cap\mk{D}$,
 \[|f(\mf{x})-f(\mf{x}_0)|<\frac{\varepsilon}{3}.\]
 It follows that for all $\mf{u}, \mf{v} \in B(\mf{x}_0, \delta)\cap\mk{D}$,
 \[|f(\mf{u})-f(\mf{v})|<\frac{2\varepsilon}{3}.\]
 
 Thus, if $r<\delta$, 
 \[0\leq \sup\left\{f(\mf{u})-f(\mf{v})\,|\, \mf{u} , \mf{v}\in B(\mf{x}_0, r)\cap\mk{D}\right\}\leq \frac{2\varepsilon}{3}<\varepsilon.\]
This shows that  
\[  \omega_f(\mf{x}_0)=\lim_{r\to 0^+}  \sup\left\{f(\mf{u})-f(\mf{v})\,|\, \mf{u} , \mf{v}\in B(\mf{x}_0, r)\cap\mk{D}\right\}=0.\]
Conversely, assume that $\omega_f(\mf{x}_0)=0$. Given $\varepsilon>0$, there is a $\delta>0$ such that for all $0<r<\delta$,
\[\sup\left\{f(\mf{u})-f(\mf{v})\,|\, \mf{u} , \mf{v}\in B(\mf{x}_0, r)\cap\mk{D}\right\}<\varepsilon.\]
If $\mf{x}$ is in $ B(\mf{x}_0, \delta/2)\cap\mk{D}$, 
\[|f(\mf{x})-f(\mf{x}_0)|\leq\sup \left\{f(\mf{u})-f(\mf{v})\,|\, \mf{u} , \mf{v}\in B(\mf{x}_0, \delta/2)\cap\mk{D}\right\}<\varepsilon.\]
This proves that $f$ is continuous at $\mf{x}_0$.

\end{myproof}
 
\begin{corollary}[label=230912_7]{}
 Let $\mk{D}$ be a subset of $\mb{R}^n$, and let   $f:\mk{D}\to\mb{R}$ be a bounded function defined on $\mk{D}$. If $\mathcal{N}$ is the set of discontinuities of $f$, then
  
 \[\mathcal{N} =\left\{\mf{x}\in\mk{D}\,\left|\, \omega_f(\mf{x})>0\right.\right\}.\]
  
\end{corollary}
 
 We also  need the following proposition.

 \begin{proposition}[label=230912_8]{}
 Let $\mk{D}$ be a compact subset of $\mb{R}^n$, and let   $f:\mk{D}\to\mb{R}$ be a bounded function defined on $\mk{D}$. 
 \begin{enumerate}[(a)]
 \item For any $a>0$, the set
 \[A=\left\{\mf{x}\in\mk{D}\,|\, \omega_f(\mf{x})\geq a\right\}\] is a compact subset of $\mb{R}^n$.
 \item If $\mathcal{N}$ is the set of discontinuities of $f$, then
 \[\mathcal{N}=\bigcup_{k=1}^{\infty}\mathcal{N}_k,\]
 where 
 \[\mathcal{N}_k=\left\{\mf{x}\in\mk{D}\,\left|\, \omega_f(\mf{x})\geq \frac{1}{k}\right.\right\}.\]
 \item 
 The set $\mathcal{N}$ has Lebesgue measure zero if and only if $\mathcal{N}_k$ has Jordan content zero for each $k\in\mb{Z}^+$.
 
 \end{enumerate}
 \end{proposition}
 \begin{myproof}{Proof}
 Since $\mk{D}$ is compact, it is closed and bounded.
 For part (a),  $A\subset\mk{D}$ implies $A$ is bounded. 
To prove that $A$ is compact, we only need to show that $A$ is closed. This is equivalent to  $\mb{R}^n\setminus A$ is open.
Notice that \[ \mb{R}^n\setminus A=U_1\cup U_2,\]
where 
\[U_1=\mb{R}^n\setminus \mk{D}\quad\text{and}\quad U_2=\mk{D}\setminus A.\]
Since $\mk{D}$ is closed, $U_1$ is open. If $\mf{x}_0\in U_2$, then $\omega_f(\mf{x}_0)<a$. Let $\varepsilon=a-\omega_f(\mf{x}_0)$. Then $\varepsilon>0$. By definition of $\omega_f(\mf{x}_0)$, there is a $\delta>0$ such that for all $0<r<\delta$, 
\[\sup\{f(\mf{u})-f(\mf{v})\,|\,\mf{u}, \mf{v}\in B(\mf{x}_0, r)\cap\mk{D}\}<\omega_f(\mf{x}_0)+\varepsilon=a.\]

Take $c=\delta/3$. If $\mf{x}$ is in $ B(\mf{x}_0, c)$, and $\mf{u}, \mf{v}$ are in $B(\mf{x}, c)$, then $\mf{u}, \mf{v}$ are in $B(\mf{x}_0, 2c)$. Since  $2c<\delta$, we find that    
\begin{align*}
\omega_f(\mf{x})&\leq \sup\{f(\mf{u})-f(\mf{v})\,|\,\mf{u}, \mf{v}\in B(\mf{x}, c)\cap\mk{D}\}\\&\leq \sup\{f(\mf{u})-f(\mf{v})\,|\,\mf{u}, \mf{v}\in B(\mf{x}_0, 2c)\cap\mk{D}\}<a.\end{align*}
 
This shows that $B(\mf{x}_0, c)\subset U_2$. Hence, $U_2$ is open. Since $\mb{R}^n\setminus A$ is a union of two open sets, it is open. This completes the proof.

Part (b) follows   from Corollary \ref{230912_7} and the identity
 \[(0,\infty)=\bigcup_{k=1}^{\infty}\left[\frac{1}{k}, \infty\right).\]
 
 For part (c), 
 if $\mathcal{N}$ has Lebesgue measure zero, then for any $k\in\mb{Z}^+$, $\mathcal{N}_k$ also has Lebesgue measure zero. By part (a) and  Proposition \ref{230912_9}, $\mathcal{N}_k$ has Jordan content zero. Conversely, assume that  $\mathcal{N}_k$  has Jordan content zero for each $k\in\mb{Z}^+$. Then $\mathcal{N}_k$ has Lebesgue measure zero for each $k\in\mb{Z}^+$. Part (b) and Proposition \ref{230912_6} implies that $\mathcal{N}$ also has Lebesgue measure zero.
\end{myproof}

Now we can prove the Lebesgue-Vitali theorem.
\begin{myproof}{Proof of the Lebesgue-Vitali Theorem}
First we assume that $f:\mf{I}\to\mb{R}$ is Riemann integrable. Given $k\in \mb{Z}^+$,  we will show that  the set 
\[\mathcal{N}_k= \left\{\mf{x}\in\mk{D}\,\left|\, \omega_f(\mf{x})\geq \frac{1}{k}\right.\right\}\]
has Jordan content zero. By Proposition \ref{230912_8}, this implies that the set $\mathcal{N}$ of discontiuities of $f:\mf{I}\to\mb{R}$ has Lebesgue measure zero.

Fixed $k\in\mb{Z}^+$. Given $\varepsilon>0$, since $f:\mf{I}\to\mb{R}$ is Riemann integrable, there is a partition $\mf{P}$ of $\mf{I}$ such that
\[U(f,\mf{P})-L(f,\mf{P})<\frac{\varepsilon}{2k}.\]
 
Let
\[\mathscr{A}=\left\{\mf{J}\in\mathcal{J}_{\mf{P}}\,|\,(\text{Int}\, \mf{J})\cap\mathcal{N}_k\neq\emptyset\right\}.\]

Then 
\[\mathcal{N}_k=A_1\cup A_2,\] where
\[A_1= \left(\bigcup_{\mf{J}\in \mathscr{A}}\text{int}\,\mf{J}\right)\,\cap\,\mathcal{N}_k,\]

and
\[ A_2=\left(\bigcup_{\mf{J}\in\mathcal{J}_{\mf{P}}}\,\pa\mf{J}\right)\cap \mathcal{N}_k.\]

Notice that the set $A_2$ has Jordan content zero. Therefore, there are finitely many open rectangles $U_1, \ldots, U_m$ such that
\[A_2\subset \bigcup_{l=1}^m U_l\quad\text{and}\quad \sum_{l=1}^m\text{vol}\,(U)<\frac{\varepsilon}{2}.\]
The set $A_1$ itself is contained in a finite union of open rectangles $\text{int}\,\mf{J}$ with $\mf{J}\in \mathscr{A}$. Notice that
\[\sum_{\mf{J}\in\mathscr{A}}\left(M_{\mf{J}}(f)-m_{\mf{J}}(f)\right)\,\text{vol}\,(\mf{J})\leq U(f,\mf{P})-L(f,\mf{P})<\frac{\varepsilon}{2k}.\]
If  $\mf{J}$ is in $\mathscr{A}$, there is an $\mf{x}_0\in A_1$ such that $\mf{x}_0\in\text{int}\,\mf{J}$. Since $\text{int}\,\mf{J}$ is an open set, there is a $\delta>0$ such that $B(\mf{x}_0, \delta)\subset\mf{J}$. 
\bp
Now,
\[\sup\{f(\mf{u})-f(\mf{v})\,|\, \mathbf{u}, \mf{v}\in B(\mf{x}_0, \delta)\}\geq \omega_f(\mf{x}_0)\geq \frac{1}{k}.\]
Therefore, 
\[M_{\mf{J}}(f)-m_{\mf{J}}(f)\geq \frac{1}{k}\]
 
for each $\mf{J}$ in $\mathscr{A}$. This implies that
\[\frac{1}{k}\sum_{\mf{J}\in\mathscr{A}}\text{vol}\,(\text{int}\,\mf{J})\leq \sum_{\mf{J}\in\mathscr{A}}\left(M_{\mf{J}}(f)-m_{\mf{J}}(f)\right)\text{vol}\,(\text{int}\,\mf{J})<\frac{\varepsilon}{2k}.\]
Thus,
\[\sum_{\mf{J}\in\mathscr{A}}\text{vol}\,(\text{int}\,\mf{J})<\frac{\varepsilon}{2}.\]
 
Hence,
\[\mathscr{B}=\left\{\text{int}\,\mf{J}\,|\,\mf{J}\in \mathscr{A}\right\}\cup \left\{U_l\,|\, 1\leq l\leq m\right\}\] is a finite collection of open rectangles that covers $\mathcal{N}_k$, and the sum of the volumes of the rectangles in $\mathscr{B}$ is less than $\varepsilon$. This proves that $\mathcal{N}_k$ indeed has Jordan content zero.

 Conversely, assume that $\mathcal{N}$ has Lebesgue measure zero. Since $f:\mf{I}\to\mb{R}$ is bounded, there is a positive number $M$ such that
 \[|f(\mf{x})|\leq M\hspace{1cm}\text{for all}\;\mf{x}\in\mf{I}.\]
Given $\varepsilon>0$, let $k$ be a positive integer such that 
 \[k \geq \frac{2\text{vol}\,(\mf{I})}{\varepsilon}.\]
  Proposition \ref{230912_8} says that $\mathcal{N}_k$   has Jordan content zero. 
 Thus,  there is a finite collection  of open rectangles  $\mathscr{B}_k=\left\{U_{l}\,|\,1\leq l\leq m\right\}$ such that
 \[\mathcal{N}_k\subset \bigcup_{l=1}^{m} U_{l}\quad\hspace{1cm}\quad \sum_{l=1}^{m}\text{vol}\, (U_{l})<\frac{\varepsilon}{4M}.\]
  Let $\mf{P}_0$ be a partition of $\mf{I}$ such that each rectangle $\mf{J}$ in $\mathcal{J}_{\mf{P}_0}$ lies entirely in the closure of one of the rectangles $U_l$, $1\leq l\leq m$ or it is disjoint from all the $U_l$, $1\leq l\leq m$.
 Let 
 \[\mathscr{C}=\left\{\mf{J}\in\mathcal{J}_{\mf{P}_0}\,|\, \mf{J}\cap U_l= \emptyset\;\text{for all}\;1\leq l\leq m\right\}.\]
 \bp
 Then 
 \[\mathcal{N}_k \subset \bigcup_{\mf{J}\in\mathcal{J}_{\mf{P}_0}\setminus\mathscr{C}}\text{int}\,\mf{J}.\]
 For each point $\mf{x}$ that is in $\mathbf{I}\setminus \mathcal{N}_k$, there is an $ r_{\mf{x}}>0$ such that the open cube $\mf{x}+(-r_{\mf{x}},r_{\mf{x}})^n$  is contained in the open set $\mb{R}\setminus \mathcal{N}_k$. By taking a smaller $r_{\mf{x}}$, we can assume that
 \[\sup\left\{f(\mf{u})-f(\mf{v})\,|\, \mf{u}, \mf{v}\in B(\mf{x}, 2r_{\mf{x}})\right\}<\frac{1}{k}.\]
 The ball $B(\mf{x}, 2r_{\mf{x}})$ contains the cube  $Q=\mf{x}+[-r_{\mf{x}},r_{\mf{x}}]^n$. 
 Therefore,
 \[M_Q(f)-m_Q(f)\leq\frac{1}{k}.\]
  The collection
 \[\left\{\mf{x}+(-r_{\mf{x}},r_{\mf{x}})^n\,|\,\mf{x}\in \mathbf{I}\setminus \mathcal{N}_k\right\}\]
 is an open covering of the compact set
 \[K=\bigcup_{\mf{J}\in \mathscr{C}}\mf{J}.\]

 Thus, there is a finite subcover $\{V_1, \ldots, V_s\}$. For $1\leq j\leq s$, let $\mf{I}_j=\overline{V}_j\cap\mf{I}$. Then we still have
 \[K=\bigcup_{\mf{J}\in \mathscr{C}}\mf{J} \subset  \bigcup_{j=1}^s \mf{I}_j.\] After renaming the rectangles, let
 \[\left\{\overline{U}_l\,|\,1\leq l\leq m\right\}\cup \left\{\mf{I}_j\,|\, 1\leq j\leq s\right\}=\left\{W_1, W_2, \ldots, W_{q}\right\}.\]
 Now let $\mf{P}$ be a partition of $\mf{I}$ so that each rectangle in $\mathcal{J}_{\mf{P}}$ is either disjoint from the interior of all the $W_j$, $1\leq j\leq q$, or is contained in one of the $W_j$.
 Let 
 \[\mathscr{D}=\left\{\mf{J}\in\mathcal{J}_{\mathcal{P}}\,|\, \mf{J}\in \overline{U}_l \;\text{for some}\;1\leq l\leq m\right\}.\]
 
 Then
 \[\sum_{\mf{J}\in\mathscr{D}}(M_{\mf{J}}-m_{\mf{J}})\,\text{vol}\,(\mf{J})\leq 2M\sum_{l=1}^m\text{vol}\,\left(\overline{U}_l\right)= 2M\sum_{l=1}^m\text{vol}\,\left(U_l\right)<\frac{\varepsilon}{2}.\]
 \bp
 
 For those $\mf{J}$ that is not in $\mathscr{D}$, it is contained in one of the cubes $\mf{x}+[-r_{\mf{x}}, r_{\mf{x}}]^n$. Therefore,
 \[M_{\mf{J}}(f)-m_{\mf{J}}(f)\leq\frac{1}{k}.\]
 
 It follows that
 \[\sum_{\mf{J}\in\mathcal{J}_{\mf{P}}\setminus \mathscr{D}}(M_{\mf{J}}-m_{\mf{J}})\,\text{vol}\,(\mf{J})\leq \frac{1}{k}\sum_{\mf{J}\in\mathcal{J}_{\mf{P}}\setminus \mathscr{D}} \,\text{vol}\,(\mf{J})\leq \frac{\text{vol}\,(\mf{I})}{k}\leq \frac{\varepsilon}{2}.\]
 This proves that
 \[U(f,\mf{P})-L(f,\mf{P})<\varepsilon.\]
 Hence, $f:\mf{I}\to\mb{R}$ is Riemann integrable.
\end{myproof}
\backmatter

\chapter*{References}
\bibliographystyle{amsalpha}
\bibliography{ref}
 
\begin{coverpage}
~
 \end{coverpage}
\end{document}